%% file: exposes/introduction.tex
\ifx\danslelivre\undefined
\documentclass[10pt]{../fichiers-configuration-latex/smfart}
\usepackage{ifthen}
\newboolean{english}
\newboolean{francais}
\setboolean{english}{false}
\setboolean{francais}{true}
\input{../fichiers-configuration-latex/paquets}
\input{../fichiers-configuration-latex/lemme-proposition-theoreme}

\input{../fichiers-configuration-latex/numerotation}
\input{../fichiers-configuration-latex/macros}

\author{Luc Illusie}

\title{Introduction}

\begin{document}

\maketitle

\else
\setboolean{english}{true}
\setboolean{francais}{false}
\selectlanguage{francais}
\chapter{Introduction}
\addtocontents{toc}{par Luc Illusie}
\begin{center}
Luc Illusie
\end{center}
\fi

\renewcommand\thesubsection{\textbf{\arabic{subsection}}}

Le présent volume rassemble les exposés d'un groupe de travail qui s'est tenu à l'École
polytechnique du printemps 2006 au printemps 2008 sur les travaux de Gabber présentés dans
\cite{Conference-Illusie@Gabber} et \cite{Conference-Deligne@Gabber}. Ceux-ci portent sur la
cohomologie étale et l'uniformisation des schémas quasi-excellents. 

En ce qui concerne la cohomologie étale, un résultat central est le théorème de finitude suivant :

\begin{theoreme}\label{0-finitude1} Soient $Y$ un schéma noethérien quasi-excellent (\textit{cf.}
\refext{I}{1def-excellent}), $f : X \ra Y$ un morphisme de type fini, $n$ un entier inversible sur
$Y$, et $F$ un faisceau constructible de $\ZZ/n\ZZ$-modules sur X. Alors, pour tout $q$, $R^qf_*F$
est constructible, et il existe un entier $N$ tel que $R^qf_*F = 0$ pour $q \ge N$.  \end{theoreme} 

Rappelons que, sans hypothèse sur $Y$ ni sur $n$, mais lorsqu'on suppose $f$ propre et de
présentation finie, les faisceaux $R^qf_*F$ sont constructibles, et nuls pour $q > 2d$ si $d$ majore
la dimension des fibres de $f$ \sga{4}{xiv}{1.1}. Si $f$ n'est pas supposé propre, l'hypothèse que
$n$ soit inversible sur $Y$ est essentielle : si $k$ est un corps algébriquement clos de
caractéristique $p >0$ et $X$ la droite affine sur $Y = \SP(k)$, $H^1(X,\ZZ/p\ZZ)$ est un
$\FF_p$-espace vectoriel de dimension infinie (\textit{cf.} \sga{4}{xiv}{1.3}). Pour $Y$ excellent
de caractéristique nulle, \ref{0-finitude1} est démontré dans l'exposé d'Artin \sgacourt{4}{xix}{}. Si
$S$ est un schéma noethérien régulier de dimension $\le 1$ (non nécessairement quasi-excellent) et
$f$ un morphisme de $S$-schémas de type fini,  la conclusion de \ref{0-finitude1} est encore vraie,
d'après \sgaquatreetdemi{Th. Finitude}{1.1} (ce résultat peut d'ailleurs se déduire de
\ref{0-finitude1}, \textit{cf.} XIII).  

La démonstration d'Artin dans \sgacourt{4}{xix}{} utilise la résolution des singularités de Hironaka pour
se ramener au cas où $f$ est l'inclusion du complément d'un diviseur régulier dans un schéma
régulier et $F$ un faisceau constant, auquel cas la conclusion découle du théorème de pureté
cohomologique absolu établi également dans (\textit{loc. cit.}). 

La démonstration de Gabber de \ref{0-finitude1} suit la même méthode, mais :

(a) on doit faire appel au théorème de pureté cohomologique absolu établi dans le cas général par
Gabber \cite{AZUMINO@Fujiwara}, 

(b) on ne dispose plus de la résolution des singularités sous la forme de Hironaka ; celle-ci est
remplacée par un théorème d'uniformisation locale, dû à Gabber (\cite{Conference-Deligne@Gabber}),
qui s'énonce ainsi (\refext{II}{2-théorème uniformisation faible}, \refext{III}{3-unif},
\refext{IX}{9-unif}) : \begin{theoreme}\label{0-uniformisation l'} Soient $X$ un schéma noethérien
quasi-excellent, $Z$ un fermé rare de $X$ et $\ell$ un nombre premier inversible sur $X$. Il existe
une famille finie de morphismes $(p_i : X_i \ra X)_{i \in I}$, couvrante pour la topologie des
$\ell'$-altérations et telle que, pour tout $i \in I$ :

(i) $X_i$ soit régulier et connexe,

(ii) $p_i^{-1}(Z)$ soit le support d'un diviseur à croisements normaux stricts.  \end{theoreme} La
topologie des $\ell'$-altérations est une topologie du type de celles considérées par Voevodsky
(\textit{cf.} \cite{Cohomological@Goodwillie-Lichtenbaum}, pour laquelle les $\ell'$-altérations
(i.~e. les morphismes propres surjectifs génériquement finis de degré résiduel générique premier à
$\ell$) et les recouvrements étales complètement décomposés (i.~e.  \textit{de Nisnevich}) sont des
familles couvrantes, voir (\refext{II}{2-alt}) pour une définition précise.

La première partie de ce volume est consacrée, après des rappels, dans I, sur les notions de schéma
quasi-excellent ou excellent, à la démonstration de \ref{0-uniformisation l'} et de certains
compléments et variantes. 

Il y a trois grandes étapes.

(A) \textit{Fibration en courbes}. La question étant locale pour la topologie de Nisnevich, et en
particulier pour la topologie de Zariski, on peut supposer $X$ de dimension finie. On raisonne par
récurrence sur la dimension $d$ de $X$. On peut supposer $X$ local hensélien. L'hypothèse de
quasi-excellence sur $X$ permet d'appliquer le théorème d'algébrisation de Popescu au complété
$\widehat X$ de $X$, et de se ramener, via des techniques d'approximation dues à Gabber, expliquées
dans III, au cas où $X$ est local complet de dimension $d$, et même intègre, normal, avec $d \ge 2$.
On dispose alors des classiques théorèmes de structure de Cohen. Tels quels, ils sont toutefois
insuffisants. Mais un délicat raffinement, dû à Gabber, démontré dans IV, permet de se ramener,
quitte à remplacer $X$ par une extension finie de degré générique premier à $\ell$, au cas où $X$
est le complété en un point fermé d'un schéma $X'$ de type fini et de dimension relative 1 sur un
schéma local noethérien régulier complet de dimension $d-1$, et le fermé $Z$ le complété d'un fermé
rare $Z'$ de $X'$. Ce théorème de « fibration », ou « d'algébrisation partielle », est
établi dans V. Après quelques nouvelles réductions faciles, on est ramené au cas où $X$ est un
schéma normal intègre, propre sur un schéma affine normal intègre excellent $Y$ de dimension $d-1$,
à fibre générique géométriquement intègre, lisse, et de dimension 1, et le fermé $Z$ un diviseur de
fibre générique étale. 

(B) \textit{de Jong et log régularité}. On peut alors appliquer le théorème de la courbe nodale de
de Jong, sous sa forme équivariante, \cite[2.4]{Families@de_Jong} : il existe un groupe fini $G$ et
une altération projective $G$-équivariante du morphisme $f : X \ra Y$ en un morphisme $f' : X' \ra
Y'$, qui est une courbe projective nodale $G$-équivariante, et l'image inverse $Z'$ de $Z$ un
diviseur de composantes dominantes étales et transverses au lieu lisse de $f'$ (\refext{IX}{9-dJ}).
Appliquant l'hypothèse de récurrence au quotient de $Y'$ par un $\ell$-Sylow de $G$, et utilisant le
lemme d'Abhyankar, on se ramène finalement au cas où $G$ est un $\ell$-groupe, $X = X'/G$, $Y =
Y'/G$, $Y'$ contient un fermé ($G$-équivariant) $T'$ tel que le couple $(Y',T')$ soit \textit{log
régulier}, ainsi que le couple $(X',f'^{-1}(T') \cup D)$, où $D$ est un diviseur ($G$-équivariant)
étale sur $Y'$, avec $Z' \subset f'^{-1}(T') \cup D$ (\textit{cf.} \refext{VI}{6-nodal implique log
lisse}, \refext{VI}{6-Kummer}). 

(C) \textit{Modifications équivariantes}. Le quotient d'un log schéma log régulier par l'action d'un
groupe fini n'est pas en général log régulier. En particulier, $X = X'/G$, muni du fermé
$(f'^{-1}(T') \cup D)/G$, n'est pas nécessairement log régulier. Si ce couple était log régulier, la
désingularisation de Kato-Niziol des schémas log réguliers (\cite{Toric@Kato}, \cite{Toric@Niziol}),
généralisant la classique désingularisation des variétés toriques \cite{Toroidal@KKMS}, terminerait
la démonstration de \ref{0-uniformisation l'}. Gabber a dégagé des conditions suffisantes assurant
que le quotient d'un log schéma log régulier par un groupe fini est encore log régulier. Il s'agit
de la notion d'\textit{action très modérée}, étudiée dans VI.  Gabber montre qu'il existe une
modification projective $G$-équivariante $p : X'' \ra X'$ et un fermé $G$-équivariant $D''$ de $X''$
contenant l'image inverse de $f'^{-1}(T') \cup D$ tels que le couple $(X'',D'')$ soit log régulier,
et l'action de $G$ sur $(X'',D'')$ très modérée. On conclut alors en appliquant la désingularisation
de Kato-Niziol au quotient de $(X'',D'')$ par $G$. La démonstration de ce théorème de modification
est donnée en VIII. Elle s'appuie sur la théorie des désingularisations canoniques en
caractéristique nulle (Hironaka, Bierstone-Milman, Villamayor, Temkin). 

Gabber établit dans X une \textit{variante relative} de ce théorème de modification, où le log
schéma $G$-équivariant considéré est non seulement log régulier, mais log lisse sur une base log
régulière $S$, avec action triviale de $G$ : on peut construire une modification équivariante
respectant la log lissité sur $S$. Il en déduit notamment les raffinements suivants de théorèmes de
de Jong :
\begin{theoreme}\label{0-deJong l'} (1) Soit $X$ un schéma séparé et de type fini sur un
corps $k$, $Z \subset X$ un fermé rare, $\ell$ un nombre premier $\ne \mathop{\rm car}(k)$. Il
existe alors une extension finie $k'$ de $k$ de degré premier à $\ell$ et une $\ell'$-altération $h
: X' \rightarrow X$ au-dessus de $\mathop{\rm Spec}k' \rightarrow \mathop{\rm Spec}k$, avec $X'$
lisse et quasi-projectif sur $k'$, et $h^{-1}(Z)$ le support d'un diviseur à croisements normaux
stricts.

(2) Soient $S$ un schéma noethérien séparé, intègre, excellent, régulier,  de dimension 1, de point
générique $\eta$, $X$ un schéma séparé, plat et de type fini sur  $S$, $\ell$ un nombre premier
inversible sur $S$, $Z \subset X$ un fermé rare. Alors il existe une extension finie $\eta'$ de
$\eta$ de degré premier à $\ell$, et une $\ell'$-altération projective $h : \tilde X \rightarrow X$
au-dessus de  $S'\rightarrow S$, où $S'$ est le normalisé de $S$ dans $\eta'$, avec $\tilde X$
régulier et quasi-projectif sur $S'$, un diviseur à croisements normaux stricts $\tilde T$ sur
$\tilde X$,  et une partie fermée finie $\Sigma$ de $S'$ tels que :

(i) en dehors de $\Sigma$, $\tilde X \ra S'$ est lisse et $\tilde T \ra S'$ est un diviseur relatif
à croisements normaux ;

(ii) localement pour la topologie étale autour de chaque point géométrique $x$ de $\tilde X_{s'}$,
où $s' = \SP \,k'$ appartient à $\Sigma$, le couple $(\tilde X, \tilde T)$ est isomorphe au couple
formé de $$ X'=S'[t_1,\cdots,t_n,u_1^{\pm 1}, \cdots, u_s^{\pm
1}]/(t_1^{a_1}...t_r^{a_r}u_1^{b_1}\cdots u_s^{b_s}-\pi), $$ $$ T' = V(t_{r+1} \cdots t_m) \subset
X' $$ autour du point $(u_i = 1)$, $(t_j = 0)$, avec $1 \le r \le m \le n$, pour des entiers $> 0$
$a_1, \cdots, a_r, b_1, \cdots, b_s$ tels que $\mathrm{pgcd}(p,a_1,\cdots,a_r,b_1,\cdots, b_s) = 1$,
$p$ désignant l'exposant caractéristique de $k'$, et $\pi$ une uniformisante locale en $s'$ ;

(iii) $h^{-1}(Z)_{\mathrm{red}}$ est un sous-diviseur de $\cup_{s' \in \Sigma} (\tilde X_{s'}) \cup
\tilde T$.  \end{theoreme}

Le théorème d'uniformisation \ref{0-uniformisation l'} a le corollaire suivant, dit théorème
d'uniformisation « faible », où n'apparaît plus de nombre premier $\ell$ :
\begin{corollaire}\label{0-uniformisation faible} Soient $X$ un schéma noethérien quasi-excellent,
$Z$ un fermé rare de $X$. Il existe une famille finie de morphismes $(p_i : X_i \ra X)_{i \in I}$,
couvrante pour la topologie des altérations et telle que, pour tout $i \in I$ :

(i) $X_i$ soit régulier et connexe,

(ii) $p_i^{-1}(Z)$ soit le support d'un diviseur à croisements normaux stricts.  \end{corollaire} La
topologie des altérations est définie de manière analogue à celle des $\ell'$-altérations. Elle est
plus fine que la topologie étale, et aucune contrainte n'est imposée sur le degré des extensions
résiduelles génériques, \textit{cf.} \refext{II}{2-alt}. On peut démontrer \ref{0-uniformisation
faible} indépendamment de \ref{0-uniformisation l'}, en suivant seulement les étapes (A) et (B)
décrites plus haut, à l'aide, dans (A), d'une forme faible du théorème de Cohen-Gabber. Le théorème
de modification de (C) est inutile, il n'y a pas besoin de faire appel à la résolution canonique des
singularités en caractéristique nulle. La résolution des singularités toriques de Kato-Niziol
suffit. Cette démonstration est exposée dans VII.

\medskip La démonstration de \ref{0-finitude1} est donnée dans (\refext{XIII}{13-section
constructibilité forte}), en application de \ref{0-uniformisation l'} et du théorème de pureté
cohomologique absolu. L'énoncé de \ref{0-finitude1} comporte deux assertions :

(i) la constructibilité des $R^qf_*F$ pour tout $q$,

(ii) l'existence d'un entier $N$ (dépendant de $(f,F)$) tel que $R^qf_*F = 0$ pour $q \ge N$,
autrement dit, le fait que $Rf_*$ envoie $D^b_c(X,\ZZ/n\ZZ)$ dans $D^b_c(Y,\ZZ/n\ZZ)$, l'indice $c$
désignant la sous-catégorie pleine de $D^b$ formée des complexes à cohomologie constructible. 

Dans (\refext{XIII}{13-section constructibilité forte}), ces deux assertions sont démontrées
simultanément. On peut toutefois démontrer (i) en n'invoquant que le théorème d'uniformisation
faible \ref{0-uniformisation faible} (et le théorème de pureté). Ceci est fait dans
(\refext{XIII}{13-constructibilité sans annulation}). L'idée est la suivante. Si, dans
\ref{0-uniformisation faible}, les morphismes $p_i$ étaient propres, on pourrait, après réduction au
cas où $f$ est une immersion ouverte et $F$ un faisceau constant, se ramener, par descente
cohomologique propre, au cas de l'immersion du complément d'un diviseur à croisements normaux
stricts dans un schéma régulier, justiciable du théorème de pureté. Cependant, les $p_i$ ne sont pas
propres en général. Gabber tourne la difficulté à l'aide du théorème de constructibilité générique
de Deligne \sgaquatreetdemi{Th. finitude}{1.9 (i)} et d'un théorème de descente cohomologique
« orientée », établi dans XII. Ce théorème utilise la notion de \textit{produit orienté de topos},
due à Deligne \cite{AG_Tokyo-Kyoto@Laumon}. Les définitions et propriétés de base sont rappelées
dans XI. Un résultat clé est un théorème de changement de
base pour des « tubes » (\refext{XI}{11-ch. base}).

D'après des exemples classiques de Nagata, un schéma noethérien quasi-excellent n'est pas
nécessairement de dimension finie. Si $Y$ est de dimension finie, alors $Rf_*$ est de dimension
cohomologique finie, d'après un théorème de Gabber, exposé dans XVIII, et par suite (ii) découle de
(i). Gabber prouve, plus précisément, que si $X$ est un schéma noethérien, strictement local, de
dimension $d >0$, et $\ell$ un nombre premier inversible sur $X$, alors, pour tout ouvert $U$ de
$X$, on a $\mathrm{cd}_{\ell}(U) \le 2d-1$.

L'hypothèse de quasi-excellence dans \ref{0-finitude1} est essentielle, comme le montre l'exemple
donné dans XIX, d'une immersion ouverte $j : U \ra X$ de schémas noethériens, avec 2 inversible sur
$X$, telle que $R^1j_*(\ZZ/2\ZZ)$ ne soit pas constructible. 

Gabber a démontré des variantes de \ref{0-finitude1} pour les faisceaux d'ensembles ou de groupes
non commutatifs \cite{Conference-Illusie@Gabber}~:

\begin{theoreme}\label{0-finitude2}
Soit $f : X \ra Y$ un morphisme de type fini entre schémas noethériens. Alors :

(1) Pour tout faisceau d'ensembles constructible $F$ sur $X$, $f_*F$ est constructible.

(2) Si $Y$ est quasi-excellent, et si $L$ est un ensemble de nombres premiers inversibles sur $Y$,
pour tout faisceau de groupes constructible et de $L$-torsion $F$ sur $X$, $R^1f_*F$ est
constructible.
\end{theoreme}

La démonstration est donnée dans XXI. Elle ne fait pas appel aux
théorèmes d'uniformisation précédents. Elle utilise, pour le point clé, une technique
d'ultra-produits, et un théorème de rigidité pour les coefficients non abéliens, dû également à
Gabber, qui est établi dans XX. Ce théorème est une variante d'un théorème de Fujiwara-Gabber pour
les coefficients abéliens \cite[6.6.4, 7.1.1]{TUBULAR@Fujiwara} (signalée dans (\textit{loc. cit.},
6.6.5). Il s'énonce ainsi :

\begin{theoreme}\label{0-rigidité}
Soient $(X,Y)$ un couple hensélien,
où $X = \SP A$, $Y = V(I)$, l'idéal $I$ étant supposé de type fini, $\widehat X$ le complété
$I$-adique de $X$, $U$ un ouvert de $X$ contenant $X-Y$, et $\widehat U = \widehat X \times_X U$.
Alors, pour tout champ en groupoïdes ind-fini $C$ sur $U$, la flèche de restriction $\Gamma(U,C) \ra
\Gamma(\widehat U, C \vert \widehat U)$ est une équivalence.
\end{theoreme}
La démonstration donnée dans XXI est indépendante de \cite{TUBULAR@Fujiwara}.

\medskip Le reste de la partie II est consacré à trois autres applications des théorèmes
d'uniformisation de la partie I.

(a) \textit{Lefschetz affine}. Il s'agit de généralisations des théorèmes de \sga{4}{xiv}{3.1} (pour
les morphismes affines entre schémas de type fini sur un corps) et \sga{4}{xix}{6.1} (pour les
morphismes affines de type fini de schémas excellents de caractéristique nulle), ainsi que du
théorème de Gabber pour les morphismes affines de schémas de type fini sur un trait
\cite{PV@Illusie}. L'énoncé principal est le suivant (\refext{XV}{15-theoreme-principal}) :

\begin{theoreme}\label{0-Lefschetz affine}
Soit $f : X \ra Y$ un morphisme affine de type fini, où
$Y$ est un schéma noethérien quasi-excellent, muni d'une fonction de dimension $\delta_Y$, et soit
$n$ un entier inversible sur $Y$. Alors, pour tout faisceau constructible $F$ de $\ZZ/n\ZZ$-modules
sur $X$, et tout entier $q$, on a :
\begin{equation}\label{0-Lefschetz affine1}
\delta_Y(R^qf_*F) \le \delta_X(F) - q
\end{equation}
\end{theoreme}
Une \textit{fonction de dimension} $\delta$ sur un
schéma $T$ est une fonction $\delta : T  \ra \ZZ$ telle que $\delta(y) = \delta(x) - 1$ si $y$ est
une spécialisation étale immédiate de $x$. Cette notion est due à Gabber.  Elle généralise celle de
\textit{dimension rectifiée} introduite dans \sgacourt{4}{xiv}{}. Elle est définie et étudiée dans XIV.
Dans \ref{0-Lefschetz affine1}, la fonction $\delta_X$ est reliée à $\delta_Y$ par $\delta_X(x) =
\delta_Y(f(x)) + \mathrm{degtr}(k(x)/k(f(x))$, et  pour un faisceau $G$ sur $X$ (resp.  $Y$)
$\delta_X(G)$ (resp. $\delta_Y(G)$) désigne la borne supérieure des $\delta_X(x)$ (resp.
$\delta_Y(x)$) pour $G_x \ne 0$. La démonstration de \ref{0-Lefschetz affine} se fait par réduction
au théorème de Gabber cité plus haut, pour les schémas de type fini sur un trait, à l'aide du
théorème d'uniformisation « faible » \ref{0-uniformisation faible} et du théorème de descente
cohomologique orientée de XII. De \ref{0-Lefschetz affine} résulte :
\begin{corollaire}\label{0-Lefschetz affine2} Si $X$ est local noethérien, hensélien,
quasi-excellent, et de dimension $d$, de corps résiduel $k$, et si $U$ est un ouvert affine de $X$,
alors, pour tout nombre premier $\ell$ inversible sur $X$, on a \begin{equation}\label{0-Lefschetz
affine3} \mathrm{cd}_{\ell}(U) \le d + \mathrm{cd}_{\ell}(k).  \end{equation} \end{corollaire} En
particulier, si $X$ est intègre, de corps des fractions $K$, on en déduit $\mathrm{cd}_{\ell}(K) = d
+ \mathrm{cd}_{\ell}(k)$, formule conjecturée dans \sga{4}{x}{3.1}. Les valeurs possibles de
$\mathrm{cd}_{\ell}(U)$, pour $U$ ouvert, non nécessairement affine, de $X$ sont étudiées dans
XVIII. Gabber donne également des contre-exemples à \ref{0-Lefschetz affine3} lorsqu'on omet
l'hypothèse de quasi-excellence. 

\medskip (b) \textit{Une nouvelle démonstration de la conjecture de pureté absolue}. La
démonstration de cette conjecture qui est donnée dans \cite{AZUMINO@Fujiwara} utilise, dans sa
dernière partie, des techniques de $K$-théorie algébrique (résultats de Thomason). Gabber a annoncé
dans \cite{Conference-Deligne@Gabber} qu'on peut éviter tout recours à la $K$-théorie algébrique, en
utilisant, à la place, la forme raffinée du théorème de de Jong \ref{0-deJong l'}{(2)}. Cette
nouvelle démonstration est exposée en détail dans XVI. Ce chapitre contient en outre une théorie de
classes fondamentales généralisées (due à Gabber), utilisée pour construire une théorie de
morphismes de Gysin pour les morphismes d'intersection complète lissifiables, généralisant les
constructions de \sgaquatreetdemicourt{Cycle}{}.

\medskip (c) \textit{Complexes dualisants}. La notion de complexe dualisant est due à Grothendieck.
L'unicité, l'existence et les propriétés générales des complexes dualisants sont étudiées dans
\sgacourt{5}{I}{}. Toutefois, dans (\textit{loc.  cit.}) l'existence n'est établie qu'en caractéristique
nulle, ou sous des hypothèses d'existence de résolution des singularités, et de validité du théorème
de pureté absolue (conjecturale à l'époque). Dans le cas de schémas de type fini sur un schéma
régulier de dimension $\le 1$, l'existence est prouvée inconditionnellement par Deligne dans
\sgaquatreetdemicourt{Dualité}{}. Dans le cas général, l'existence, et la théorie de dualité locale qui
en résulte, ont été annoncées par Gabber dans \cite{Conference-Deligne@Gabber}. Le chapitre XVII
expose cette théorie en détail. Si $X$ est un schéma noethérien, et $\Lambda = \ZZ/n\ZZ$, où $n$ est
un entier inversible sur $X$, un \textit{complexe dualisant} sur $X$ est un objet de
$D^b_c(X,\Lambda)$ tel que le foncteur $D_K = R\mathscr{H}om(-,K)$ envoie $D^b_c(X,\Lambda)$ dans
$D^b_c(X,\Lambda)$ et que, pour tout $L \in D^b_c(X,\Lambda)$, la flèche de bidualité $L \ra
D_KD_K(L)$ soit un isomorphisme. Cette définition diffère légèrement de celle de \sga{5}{I}{1.7},
voir (\refext{XVII}{17-definition-complexe-dualisant}).  L'unicité, à décalage, et torsion près par
un $\Lambda$-module inversible, est prouvée dans \sga{5}{I}{2.1}. Le résultat principal de XVII est
que, si $X$ est excellent et admet une fonction de dimension, au sens précisé plus haut, alors $X$
admet un complexe dualisant. De plus, ces complexes dualisants ont les propriétés de fonctorialité
attendues, et, si $X$ est régulier, le faisceau constant $\Lambda_X$ est dualisant. Cette dernière
assertion était une conjecture dans \sgacourt{5}{I}{}, démontrée dans le cas de caractéristique nulle.
Nous renvoyons le lecteur à l'introduction de XVII pour des énoncés plus complets, et des
indications sur la méthode de démonstration, dont les ingrédients essentiels sont le théorème de
finitude \ref{0-finitude1} et le théorème d'algébrisation partielle de V (voir (A) \textit{supra}).

\renewcommand\thesubsection{\textbf{\thesection.\arabic{subsection}}}

\ifx\danslelivre\undefined
\bibliography{../fichiers-configuration-latex/bibliographie}
\bibliographystyle{../fichiers-configuration-latex/smfalpha}
\end{document}
\fi


\begin{center}
I. Anneaux excellents\\
Michel Raynaud rédigé par Yves Laszlo
\end{center}

\tableofcontents

\else
\setboolean{english}{false}
\setboolean{francais}{true}
\selectlanguage{francais}
\chapter{Anneaux excellents}
\addtocontents{toc}{par Michel Raynaud, rédigé par Yves Laszlo}
\begin{center}
Michel Raynaud, rédigé par Yves Laszlo
\end{center}
\fi

\section{Introduction} Ce texte est une version  un peu modifiée d'un exposé donné par Michel
Raynaud au printemps 2006 dans le cadre d'un séminaire sur les travaux de Gabber sur la cohomologie
des schémas excellents. Le but est de familiariser le lecteur avec la notion d'excellence et de lui
donner un fil d'Ariane pour se repérer dans \textbf{ÉGA IV} où l'on trouve les principales
propriétés des anneaux excellents. Son ambition n'est certainement pas de donner une exposition
complète de la théorie, mais une idée de la stratégie qui ramène pour l'essentiel les preuves à des
énoncés, souvent difficiles, dans le cas complet. Dans un second temps, on montre que toutes les
propriétés définissant les anneaux excellents peuvent être mises en défaut, même en petite
dimension. Notamment, il existe des anneaux de valuation discrète non excellents ainsi que des
anneaux noethériens intègres de dimension $1$ dont le lieu régulier n'est pas ouvert. Ce dernier
exemple est un sous-produit d'une construction proposée par Gabber (\refext{XIX}{19-rem-cexgabber}).
Elle montre que le théorème de constructibilité des images directes
(\refext{XIII}{13-constructibilité}) n'est plus vrai si on omet la condition de quasi-excellence.

\section{Définitions} Soit $A$ un anneau noethérien et $X=\Spec(A)$ son spectre. On va s'intéresser
à des conditions sur $X$ de deux sortes.

\textbf{$\bullet$ Conditions globales :}

\subsection{Condition 1 : conditions d'ouverture} \emph{Tout schéma intègre $Y$ fini sur $X$
contient un ouvert dense \begin{enumerate} \item[1.a)]  régulier.  \item[1.b)] normal.
\end{enumerate}}

\begin{remarque}\label{11->ouvert}La condition 1.a)  entraîne d'après le critère d'ouverture de
Nagata que le lieu régulier de tout schéma fini sur $X$ est ouvert (\egalong{IV}{2}{6.12.4}). De même,
la condition 1.b) entraîne que le lieu normal de tout schéma fini sur $X$ est ouvert
(\egalong{IV}{2}{6.13.7})\footnote{Et en fait, 1.a) (resp. 1.b)) entraîne que le lieu régulier (resp.
normal) de tout schéma intègre de type fini sur $X$ est ouvert}. Ces critères d'ouverture  assurent
en outre que pour tester 1.a) ou 1.b) on peut se limiter à des schémas intègres $Y$  qui sont de
plus finis \emph{et} radiciels sur $X$.  \end{remarque}

\textbf{$\bullet$  Conditions locales.}

Elles sont de deux types.

\subsection{Condition 2 : Conditions sur les fibres formelles} \emph{Pour tout point \emph{fermé}
$x$ de $X$, le morphisme de complétion \footnotemark $\Spec(\chap{\OO_x})\ra\Spec(\OO_x)$ est
\begin{enumerate} \item[2.a)] régulier.  \item[2.b)] normal.  \item[2.c)] réduit.  \end{enumerate} }
\footnotetext{Ses fibres sont appelées les fibres formelles (de $X$ ou $A$) en $x$.}

Un anneau vérifiant 2.a) est dit « G-ring » en anglais, ce en l'honneur de Grothendieck qui a
dégagé l'importance de la notion et étudié ses propriétés.

\begin{remarque}\label{1rem-fibre} Rappelons qu'un morphisme de schémas noethériens est dit régulier
(resp. normal, réduit) respectivement s'il est plat et si les fibres géométriques en tout point de
la base sont régulières (resp. normales, réduites). On dit que les fibres formelles de $X$ en $x$
sont géométriquement régulières, géométriquement normales ou géométriquement réduites si le
morphisme de complétion $\Spec(\chap{\OO_x})\ra\Spec(\OO_x)$ est régulier, normal ou réduit. Bien
entendu, il suffit de tester la régularité, normalité, ou réduction des fibres après changement de
base radiciel fini (\egalong{IV}{2}{6.7.7}).  Notons que la fibre fermée de
$\Spec(\chap{\OO_x})\ra\Spec(\OO_x)$ est le spectre du corps résiduel $k(x)$ : elle est toujours
géométriquement régulière. La fibre formelle en $y\in\Spec(\OO_x)$ s'identifie à la fibre formelle
\emph{générique} du sous-schéma fermé $\overline{\{y\}}$ (muni de sa structure réduite), adhérence
de $y$ dans $\Spec(\OO_{x})$ ; ceci explique qu'on s'intéresse dans la littérature aux fibres
formelles génériques des anneaux intègres. Dans le cas où $A$ est local mais pas un corps, elles
peuvent avoir des dimensions arbitraires entre $0$ et $\dim(A)-1$ et contenir des points fermés de
hauteurs différentes, même dans le cas excellent (\ref{1def-excellent}) régulier (\cite{Rot91}).
Dans le cas où $A$ est un localisé d'une algèbre intègre de type fini sur un corps, la dimension de
la fibre formelle générique est bien $\dim(A)-1$ (\cite{Mat88}).  \end{remarque}

\subsection{Condition 3 : condition de formelle caténarité} \emph{Pour tout point fermé $y$ d'un
sous-schéma fermé irréductible $Y$ de $X$, le complété\footnotemark $\chap{\OO_{Y,y}}$ est
équidimensionnel.  }\footnotetext{Bien entendu, même si $\OO_{Y,y}$ est intègre, son complété n'est
en général pas intègre : penser à une courbe nodale.}

On dit alors que $X$ est formellement caténaire. Par exemple, si $X$ est de dimension $1$,  $X$ est
formellement caténaire.

\begin{exemple}\label{1reg-form-cat}\label{1complet-formel-cat} Tout anneau local noethérien
\emph{complet} est formellement caténaire.\end{exemple}

\bigskip

Rappelons que $X$ est dit \emph{caténaire} si toutes les chaînes saturées de fermés irréductibles
de $X$ ayant mêmes extrémités ont même longueur, \emph{universellement
caténaire}\footnote{\label{1form-dim}Cette dernière notion est utile en théorie de la dimension :
si $A$ est intègre universellement caténaire contenue dans $B$ intègre de type fini sur $A$, on a
pour tout $\p\in\Spec(B)$ au dessus de $\q\in\Spec(A)$ la formule $$\dim
B_\p+\text{degtr}_{k(\q)}k(\p)=\dim A_\q+\text{degtr}_AB.$$ Mais, en pratique, on teste plutôt la
formelle caténarité qui, comme on le voit juste après, entraîne l'universelle caténarité, et même
lui est équivalente (voir~(\ref{1Ratliff}) plus bas) !} si tout schéma affine de type fini sur $X$
est caténaire. La caténarité est une notion locale. La terminologie de \emph{formelle caténarité}
est alors justifiée par la proposition élémentaire suivante (\egalong{IV}{2}{7.1.4}), proposition qui
résulte de la fidèle platitude du morphisme de complétion

\begin{lemme}\label{1form-versus-cat} \label{1note-cat-non uni cat}Soit $A$ local noethérien de
complété équidimensionnel. Alors \begin{itemize} \item[i)] $A$ est équidimensionnel et caténaire.
\item[ii)]Pour tout idéal $I$ de $A$, le quotient $A/I$ est équidimensionnel si et seulement si son
complété l'est ; en particulier, $A/I$ est formellement caténaire.  \item[iii)] En particulier, un
schéma affine $X$ noethérien formellement caténaire est caténaire et même universellement caténaire.
\end{itemize} \end{lemme}

Notons que \textit{iii)} découle immédiatement de \textit{i)} puisque $X$ est caténaire si et
seulement si ses composantes irréductibles le sont. On verra plus bas dans la section~\ref{1inv-exc}
que la propriété de formelle caténarité est notamment stable par extension finie d'où l'universelle
caténarité annoncée (cf. la preuve de la proposition~\ref{1equi-excellent} et, pour une réciproque,
voir (\ref{1Ratliff})).

\begin{exemple}\label{1CM-form-cat}Soit alors $B\ra A$ un morphisme local surjectif d'anneaux
noethériens et supposons $B$ de Cohen-Macaulay (par exemple régulier).  Comme $\chap{B}$ est de
Cohen-Macaulay, il est équidimensionnel de sorte que $A$ est formellement caténaire d'après
(\ref{1form-versus-cat}).  \end{exemple}

Regardons ce qui se passe dans le cas complet. Rappelons pour mémoire le théorème de structure de
Cohen des anneaux locaux complets noethériens (\egalong{IV}{1}{0.19.8.8}):

\begin{theoreme}[Cohen]\label{1Co} Soit $A$ un anneau local noethérien complet de corps résiduel
$k$.\begin{itemize} \item [(i)] $A$ est isomorphe à un quotient d'un anneau de séries formelles sur
un anneau de Cohen\footnote{Rappelons (\egalong{IV}{2}{19.8.5}) que les anneaux de Cohen $C$ sont les
corps de caractéristique nulle et les anneaux de valuation discrète complets d'inégale
caractéristique $p$ non ramifiés.  Lorsque le corps résiduel $\kappa$ de $C$ est parfait, $C$ n'est
autre que l'anneau des vecteurs de Witt de $\kappa$.}. Si $A$ contient  un corps, il est isomorphe à
un quotient d'un anneau de séries formelles sur $k$.  \item [(ii)] Si $A$ est de plus intègre, il
existe un sous-anneau $B$ isomorphe à un anneau de séries formelles sur un anneau de Cohen ou un
corps\footnote{Voir (\ref{1Coh-Gab}) pour une amélioration.} de sorte que l'inclusion $B\ra A$ soit
locale,  finie et  induise un isomorphisme  des corps résiduels.  \end{itemize} \end{theoreme}

Tout anneau local noethérien complet est donc quotient d'un anneau régulier.

\begin{definition}\label{1def-excellent} Soit $X$ un schéma (resp. $X=\Spec(A)$ un schéma affine)
noethérien. On dit que $X$ (resp. $A$) est \begin{itemize} \item \emph{excellent} si $X$ vérifie
$1.a)+2.a)+3)$.  \item \emph{quasi-excellent} si $X$ vérifie $1.a)+2.a)$.  \item
\emph{universellement japonais}\footnote{Ou \emph{Nagata} en anglais, voire
\emph{pseudo-géométrique} (chez Nagata notamment).} si $X$ vérifie 1.b)+2.c).  \end{itemize}
\end{definition}

\subsection{}\label{1desing} L'existence d'une classe de schémas stable par extension finie pour
laquelle le théorème de désingularisation est vérifié impose de se limiter aux schémas
quasi-excellents. Précisément, si tous les schémas intègres et finis $Y$ sur $X$ admettent une
désingularisation (au sens de l'existence de $Y'\ra Y$ propre et birationnel avec $Y'$ régulier),
alors $X$ est quasi-excellent \egalong{IV}{2}{7.9.5})\footnote{Si de plus $X$ peut localement se
plonger dans un schéma régulier, alors $X$ vérifie 3) et est donc excellent.}.  Inversement, le
théorème de désingularisation d'Hironaka se généralise à tout schéma réduit quasi-excellent de
caractéristique nulle (\cite[3.4.3]{desingularization@Temkin})\footnote{Ce résultat a été longtemps
considéré comme « bien connu des experts » alors que sa preuve, tout à fait non triviale, date
de 2008.}

On regroupe plus bas (\ref{1pathologie}) des exemples de « méchants anneaux ». Commençons par un
regard plus positif.

\section{Exemples immédiats.}

\begin{proposition} Un corps, un anneau de Dedekind de corps des fractions \emph{de caractéristique
nulle} est excellent.\end{proposition}

\begin{proof} Vérifions qu'un corps est excellent. En effet, une algèbre finie et intègre sur un
corps est un corps : les propriétés 1.a), 2.a) et 3) sont donc vérifiées ce qui prouve que tout
corps est excellent.

\medskip

Soit $A$ un anneau de Dedekind de corps des fractions $K$ \emph{de caractéristique nulle} est
excellent.

\begin{itemize} \item Vérifions 1.a). Soit donc $B$ intègre finie sur $A$. Soit $B$ est un corps,
auquel cas 1.a) est vérifié, soit $A$ se plonge dans $B$. Comme $K$ est de caractéristique nulle,
$B$ est génériquement étale sur $A$, prouvant que le lieu régulier de $B$ contient un ouvert non
vide (le lieu étale par exemple).  \item Pour 2.a), considérons $x$ fermé dans $\Spec(A)$. La fibre
formelle non fermée en $x$ est le complété $\chap{ K_x}$ de $K$ pour la valuation définie par $x$.
Comme $K$ est de caractéristique nulle, le corps $\chap{ K_x}$ est séparable sur $K$ d'où 2.a).
\item La propriété 3) est claire puisque le complété de $A$ en $x$ est intègre donc
équidimensionnel.  \end{itemize} \end{proof} On verra plus bas (\ref{1patho-DVR}) qu'il existe de
nombreux anneaux de valuation discrète qui ne sont pas quasi-excellents.

\section{L'exemple de base : les anneaux locaux noethériens complets.} Expliquons avec  Nagata
pourquoi les anneaux locaux noethériens complets sont excellents\footnote{Ceci permet de construire
de nombreux exemples d'anneaux de valuation discrète excellents de caractéristique positive (par
complétion de schémas réguliers aux points de hauteur $1$).}.

\bigskip

La propriété 2.a) est tautologique. La formelle caténarité a été vue (\ref{1complet-formel-cat}).
Reste 1.a). Une  extension finie d'un anneau complet étant complet, on doit prouver le résultat
suivant (cf.~\egalong{IV}{2}{22.7.6}).

\begin{theoreme}[Nagata]\label{1complet-exc} Si $X$ est local noethérien intègre et complet
\footnote{D'après~(\ref{11->ouvert}), ceci entraîne que le lieu régulier d'un schéma local
noethérien complet est ouvert, qu'il soit intègre ou non.}, alors le lieu régulier est
ouvert.\end{theoreme} \begin{proof} On va distinguer les cas d'égales et d'inégales
caractéristiques.

\begin{itemize} \item[Cas I] (Cf.~\egalong{IV}{2}{0.22.7.6}.) Supposons que $A$ contienne un corps et
notons ${k_0}$ son corps premier (qui est parfait !) de sorte que le corps résiduel $k$ de $A_\p$
est séparable sur ${k_0}$ pour tout $\p\in\Spec(A)$. L'anneau $A_\p$ est régulier si et seulement si
$A_\p$ est formellement lisse sur ${k_0}$ (voir dans ce cas \egalong{IV}{2}{0.19.6.4}). D'autre part,
le théorème de structure de Cohen (\ref{1Co}) assure que $A$ est isomorphe à $k[[T_1,\cdots,T_n]]/I$
de sorte que $\p$ s'identifie à un idéal de $B=k[[T_1,\cdots,T_n]]$ contenant $I$. Le \emph{critère
jacobien de formelle lissité de Nagata} (\egalong{IV}{2}{0.22.7.3}) assure que $A_\p$ est régulier si
et seulement si il existe des ${k_0}$-dérivations $D_i,i=1,\cdots,m$ de $B$ dans $B$ et
$f_i,i=1,\cdots,m$ des éléments engendrant $I_\p$ tels que $\det(D_if_j)\not\in\p$. Cette condition
étant visiblement ouverte, le théorème suit.

\item[Cas II] Supposons que $A$ est d'inégale caractéristique, et donc de  corps des fractions $K$
de caractéristique nulle. D'après le théorème de structure de Cohen (\ref{1Co}), $A$ contient un
sous-anneau régulier (et complet) $B$ faisant de $A$ une $B$-algèbre de finie. Le corps des
fractions $L$ de $B$ est de caractéristique nulle comme $K$.  Quitte à remplacer $A$ par un localisé
$A[1/a]$, on peut supposer que $B$ est libre de rang fini sur $A$ de base $y_1,\cdots,y_m$.  Mais
$\Spec(A)\ra\Spec(B)$ est étale en dehors du fermé $d=\det_{A/B}(\mathrm{Tr}(y_iy_j))\not=0$ de
$\Spec(B)$, qui est non trivial puisqu'il contient le point générique, l'extension
$\Frac(B)/\Frac(A)$ étant séparable -- de caractéristique nulle -- ! Comme $B$ est régulier, le
théorème suit.  \end{itemize} \end{proof}

\begin{remarque}\label{1Coh-Gab} Ainsi, un anneau de séries formelles sur un corps est excellent.

Notons que la preuve se simplifie si on connaît l'amélioration de Gabber du théorème de structure de
Cohen (\refext{IV}{4-Cohen-Gabber égale} et \refext{IV}{4-Cohen-Gabber-Epp}) : si $A$ noethérien est
local complet et intègre, il contient un anneau $B$ isomorphe à un anneau de séries formelles sur un
anneau de Cohen ou un corps tel que $\Spec(A)\ra\Spec(B)$ est fini et \emph{génériquement étale}. On
n'a alors pas besoin de distinguer les caractéristiques des corps de fractions dans la preuve. Mais
la preuve de cette amélioration est difficile.  \end{remarque}

\section{Permanence par localisation et extension de type fini}\label{1inv-exc} La notion de (quasi)
excellence est remarquablement stable.  Précisément, on a

\begin{theoreme}\label{1th-permanent} Toute algèbre de type fini ou plus généralement
essentiellement de type fini sur un anneau excellent (resp.  quasi-excellent) est excellente (resp.
quasi-excellente). En particulier, tout localisé d'algèbre de type fini sur un corps ou sur un
anneau de Dedekind ($\ZZ$ par exemple) de corps des fractions de caractéristique nulle est
excellent.\end{theoreme}

(Rappelons que, dans ce contexte, un morphisme $\Spec(B) \ra \Spec(A)$
est dit \emph{essentiellement de type fini} si $B$ est une
localisation d'une $A$-algèbre de type fini par un système
multiplicatif.)

Expliquons les grandes lignes de la preuve.

\subsection{Condition 1)} Le passage au localisé ne pose pas de problème. Soit $B$ de type fini sur
$A$. Si $A$ vérifie 1.a) ou 1.b), les critères d'ouverture de Nagata (\egalong{IV}{2}{6.12.4 et
6.13.7}) entraînent qu'il en est de même de $B$.

\subsection{Condition 2)} C'est la partie la plus difficile de la théorie (\egalong{IV}{2}{7.4.4}),
entièrement due à Grothendieck. Le point le plus délicat est la localisation :

\begin{theoreme2} Si $A$ vérifie 2.a) (resp. 2.b) ou 2.c)), alors pour tout $\p\in\Spec(A)$,
l'anneau $A_\p$ vérifie 2.a) (resp. 2.b) ou 2.c)), autrement dit les fibres formelles en tout point
de $\Spec(A)$ sont géométriquement régulières (resp. géométriquement normales ou géométriquement
réduites).\end{theoreme2}

\begin{proof} La preuve se fait par réduction au cas complet. On se limite à la propriété 2.a), le
cas de 2.b) ou 2.c) se traitant de même. Soit $\m$ maximal contenant $\p\in\Spec(A)$ et soit $B$ le
complété $\m$-adique de $A$. Comme $A_\m\ra B$ est fidèlement plat, il existe $\q\in\Spec(B)$ au
dessus de $\p$. Par hypothèse, $A_\m\ra B$ est régulier. Les morphismes réguliers étant stables par
localisation, $A_\p\ra B_\q$ est régulier. On regarde alors le diagramme commutatif

$$\xymatrix{ \chap{A_\p}\ar[r]^{\hat f}&\chap{B_\q}\\ A_\p\ar[r]^f\ar[u]^\alpha&B_\q\ar[u]^\beta}.$$

Supposons que $\beta$ soit régulier. Alors, $\hat f\circ\alpha$ est régulier comme composé de deux
morphismes réguliers. Comme $\hat f$ est fidèlement plat (comme morphisme local  complété du
morphisme plat $f$), on déduit que $\alpha$ est régulier (exercice ou~\egalong{IV}{2}{6.6.1}) ce qu'on
voulait. On est donc ramené à $\beta$, donc au cas complet. La régularité de $\beta$ résulte alors
de

\begin{theoreme2}\label{1perm-fib-form-loc}Soit $B$ un anneau local noethérien complet $B$. Alors,
les fibres formelles de $B$ en $\q\in\Spec(B)$ sont géométriquement régulières.\end{theoreme2}

Ce théorème est le noyau dur de la théorie. On se ramène (\ref{1rem-fibre}) à étudier les fibres
formelles génériques. On montre donc dans un premier temps (\egalong{IV}{1}{0.22.3.3}) que si $\p$ est
un idéal premier de $A$  local noethérien complet intègre, la fibre formelle générique
$\chap{A_\p}\otimes_{A_\p}\Frac(A_\p)$  de $A_\p$ est formellement lisse sur $\Frac(A_\p)=\Frac(A)$
en tout point. Dans un second temps, on montre (\egalong{IV}{1}{0.22.5.8}) qu'une algèbre locale
noethérienne sur un corps est formellement lisse\footnote{Rappelons qu'une $k$-algèbre locale $B$
(muni de la topologie adique) est formellement lisse  sur $k$ si tout $k$-morphisme continu
d'algèbre $B\ra C/I$ avec $I^2=0$ se relève continûment à la $C$-algèbre discrète $C$.} si et
seulement si elle est géométriquement régulière\footnote{En fait, on n'a visiblement besoin que du
sens formellement lisse entraîne géométriquement régulier, qui est le plus facile. Notons que la
preuve de l'équivalence a été considérablement simplifiée par Faltings (\cite{Fal78} ou pour le
lecteur non germaniste \cite{CRT@Matsumura},~28.7).} \end{proof}

Une fois prouvée la permanence par localisation, on peut montrer :

\begin{theoreme2} Soit $B$ une $A$-algèbre de type fini. Si $A$ vérifie 2.a) (resp. 2.b) ou 2.c)),
alors $B$ vérifie 2.a) (resp. 2.b) ou 2.c)).\end{theoreme2}

La  preuve se fait par récurrence sur le nombre de générateurs de $B$. Grâce à l'invariance par
localisation, on se ramène aisément à l'étude des fibres formelles de $B$ en un idéal maximal dans
le cas où  $B$ engendré par \emph{un} élément et $A$ est complet. La preuve n'est pas facile, mais
beaucoup plus simple que celles de (\egalong{IV}{1}{0.22.3.3 et 0.22.5.8}).

\subsection{Condition 3).} De même que pour les conditions de type 2), la stabilité par localisation
et extension finie résulte comme plus haut (\egalong{IV}{2}{7.1.8}) du cas complet, la platitude du
morphisme de localisation permettant de descendre du complété à l'anneau -- ce n'est pas immédiat
malgré tout --. Le cas complet est facile comme on a vu (\ref{1reg-form-cat}).

\subsection{Application au cas local} Dans le cas local, la condition d'ouverture du lieu régulier
découle de 2.a). Précisons.

\begin{proposition2}\label{1loc-exc} \begin{itemize} \item [i)] Le lieu régulier d'un anneau local
noethérien vérifiant 2.a) est ouvert.  \item [ii)] En particulier, un anneau local noethérien est
quasi-excellent (resp.  excellent) si et seulement s'il vérifie 2.a) (resp.  s'il vérifie 2.a) et
3)).  \end{itemize} \end{proposition2}

\begin{proof}

Soit $f:X\ra Y$ un morphisme fidèlement plat de schémas noethériens à fibres régulières (resp.
normales ou réduites). Alors, $\OO_x$ est  régulier (resp.  normal ou réduit) si et seulement si
$\OO_{f(x)}$ l'est (\egalong{IV}{2}{6.4.2, 6.5.1}).  Notons $U_R(X)$ l'ensemble des $x\in X$ tel que
$R(\OO_x)$ est régulier (resp.  normal ou réduit). Autrement dit, on a $f^{-1}(U_R(X))=U_R(Y)$. Or
le lieu régulier ou normal d'un anneau complet intègre est ouvert  (\ref{1complet-exc}). De plus, le
morphisme de complétion d'un anneau local noethérien $A$ est régulier si et seulement si $A$ vérifie
2.a) d'après (\ref{1perm-fib-form-loc}). Il suit que 2.a) entraîne 1.a) (resp.  2.b) entraîne 1.b))
dans le cas local.  \end{proof}

\section{Comparaison avec ÉGA IV : le cas des anneaux universellement japonais}

Rappelons la définition usuelle des anneaux universellement japonais (\egalong{IV}{1}{0.23.1.1}).

\begin{definition}\label{1def-japonais} $X$  est dit \begin{enumerate} \item \emph{japonais} s'il
est intègre et si la clôture intégrale de $A$ dans toute extension finie\footnote{On peut se
contenter des extensions finies radicielles si l'on veut : exercice ou \egalong{IV}{1}{23.1.2}.} de son
corps des fractions est finie sur $A$\footnote{Comme module ou comme algèbre : c'est la même chose
car la clôture intégrale est entière sur $A$ par construction.}; \item \emph{universellement
japonais} si tout anneau intègre qui est extension de type fini de $A$ est japonais\footnote{Ou, ce
qui est équivalent (\egalong{IV}{2}{7.7.2}), si tout quotient intègre est japonais.}.  \end{enumerate}
\end{definition}

La définition d'anneau japonais n'est que technique en ce qu'elle ne sert qu'à définir la seule
notion véritablement utile (et vérifiable à vrai dire) : celle d'anneau universellement japonais.
Cette définition est compatible avec~\ref{1def-excellent}. Expliquons pourquoi. D'après Nagata, $X$
est universellement japonais au sens de~\ref{1def-japonais} si et seulement si $X$ vérifie 1.b) et
si tous les quotients intègres des localisés $\OO_{X,x}$ en les points fermés $x\in X$ sont japonais
(\egalong{IV}{2}{7.7.2}).  Or, le théorème de Zariski-Nagata (\egalong{IV}{2}{7.6.4}) assure  que les
quotients intègres de $\OO_{X,x}$ sont japonais si et seulement si les fibres formelles de
$\OO_{X,x}$  sont géométriquement réduites \footnote{\label{1Zar-Nag}Ou, de façon équivalente, que
le complété de toute $\OO_{X,x}$-algèbre finie et réduite est réduit. Comme d'habitude, la preuve se
fait par réduction au cas complet, et même régulier complet grâce au théorème de structure de Cohen.
Le caractère japonais de tels anneaux est garanti par le théorème de Nagata
(\egalong{IV}{1}{0.23.1.5}).}. D'où l'équivalence entre les deux définitions des anneaux
universellement japonais.

Si on renforce la condition 2.c) en 2.b) (fibres formelles géométriquement normales), le passage à
la clôture intégrale commute à la complétion. Précisément, on a (\egalong{IV}{2}{7.6.1 et 7.6.3})

\begin{proposition}\label{1prop-normal-versus-complet} Supposons que $A$ local noethérien vérifie
2.b) et soit réduit. Alors, la clôture intégrale $A'$ de $A$ dans son anneau total des fractions est
finie sur $A$ et son complété est isomorphe à la clôture intégrale de $\chap{A}$\footnote{Qui est
réduit puisque $A$ est japonais (cf.  note~(\ref{1Zar-Nag})).} dans son anneau total des fractions.
\end{proposition}

On déduit l'important critère d'intégrité du complété.

\begin{corollaire}\label{1-cor-normal-versus-complet} Soit $A$ local noethérien.  \begin{itemize}
\item[i)] Supposons $A$  intègre et vérifiant 2.b). Alors, $\chap{A}$ est intègre si et seulement si
$A$ est unibranche (\emph{i.e.} $A'$ local).  \item [ii)] Supposons $A$ hensélien. Alors $A$ est excellent
si et seulement s'il vérifie 2.a). Si $A$ est de plus intègre, il en est de même de son complété.
\end{itemize} \end{corollaire}

\begin{proof} Prouvons \textit{i)}. Comme $A$ est unibranche, la clôture intégrale $A'$ de $A$ est
locale : il en est de même de son complété $\chap{A'}$. D'après (\ref{1prop-normal-versus-complet}),
on a $\chap{A'}=(\chap{A})'$ et donc est normal. Or, un anneau normal et local est intègre. Comme
$\chap{A'}$ contient $\chap{A}$, le résultat suit.

Prouvons \textit{ii)}. D'après (\ref{1loc-exc}), on doit seulement se convaincre qu'un anneau local
hensélien vérifiant 2.a) vérifie aussi 3), \emph{i.e.} est formellement caténaire. On peut supposer $A$
intègre et on doit prouver que $\chap{A}$ est équidimensionnel. Mais comme $A$ est hensélien
intègre, il est unibranche~\egalong{IV}{4}{18.8.16}, donc $\chap{A}$ est intègre d'après le premier
point, ce qui assure l'équidimensionalité.\end{proof}

\section{Comparaison avec ÉGA IV : le cas des anneaux excellents}

La définition des anneaux noethériens excellents de Grothendieck est \textit{a priori} différente de
celle donnée ici. Notamment, elle fait intervenir, un peu bizarrement, l'\emph{universelle
caténarité} en lieu et place de la \emph{formelle caténarité}. Précisément, elle fait intervenir
trois propriétés. Dans cette partie $A$ désigne un anneau noethérien et $X=\Spec(A)$ le schéma
affine correspondant.

\begin{itemize} \item [1EGA)] Pour tout quotient intègre $B$ de $A$ et toute extension finie
radicielle $K'$ du corps des fractions $K$ de $B$, il existe une sous-$B$-algèbre finie $B'$ de $K'$
contenant $B$, de corps des fractions $K'$ telle que le lieu régulier de $\Spec(B')$ soit un ouvert
dense.  \item [2EGA)] Les fibres formelles de $X$ en tout point (fermé ou non) sont géométriquement
régulières.  \item [3EGA)] $A$ est universellement caténaire.  \end{itemize}

Les anneaux excellents au sens des \textbf{ÉGA} sont les anneaux noethériens vérifiant les 3
propriétés précédentes (\egalong{IV}{2}{7.8.2}).

Notons tout de suite, ce qui est élémentaire, que l'universelle caténarité de $A$ équivaut à celle
des anneaux locaux de $\OO_{X,x}$ en tous ses points fermés -- ou tous ses points si on préfère --
(\egalong{IV}{2}{5.6.3}).

Pour que la définition des anneaux excellents de Grothendieck (\egalong{IV}{2}{7.8.2}) soit la même que
(\ref{1def-excellent}), on doit prouver la proposition suivante.

\begin{proposition}\label{1equi-excellent}Pour tout anneau noethérien et $i=1,2,3$, les propriétés
$i$\emph{)} et $i$\emph{ÉGA)} son équivalentes. En particulier, les notions de quasi-excellence et
d'excellence de la première partie coïncident avec celles des \textbf{ÉGA}.\end{proposition}

\begin{proof}

La condition 1EGA) équivaut à 1) d'après~\egalong{IV}{2}{6.12.4} (seule la partie 1EGA) entraîne 1) est
délicate même si elle n'utilise pas le critère de régularité de Nagata mais seulement de l'algèbre
commutative standard -- essentiellement le critère de régularité par fibres et la non dégénérescence
de la trace des extensions finies séparables de corps --).

Pour l'équivalence de 2) et 2EGA), il faut se convaincre que la géométrique régularité des fibres
formelles en tout point fermé entraîne la géométrique régularité des fibres formelles en tout point
: c'est un cas particulier des propriétés de permanence (\ref{1inv-exc}).

Ceci prouve la compatibilité des définitions de la quasi-excellence.

Si $X$ vérifie 3), tous ses anneaux locaux sont formellement caténaires (permanence par
localisation, cf. la section \ref{1inv-exc}) et donc sont caténaires (\ref{1form-versus-cat}).
Comme tout schéma (affine) de type fini sur $X$ vérifie 3) (permanence par extension de type finie,
cf. la section \ref{1inv-exc}), on déduit que $X$ est universellement caténaire et donc $X$ vérifie
3EGA).

La réciproque est due à Ratliff :

\begin{proposition2}[Ratliff]\label{1Ratliff}\label{1form-cat->uni-cat} Un anneau noethérien
universellement caténaire est formellement caténaire.\end{proposition2}

Précisément, Ratliff prouve (\cite{Rat71}, 3.12) que si $A$ est caténaire, $A_\p$ est formellement
caténaire dès que $\p$ n'est pas maximal\footnote{Dans l'étrange terminologie de l'auteur, c'est la
condition $\text{depth}(\p)>0$, ce qui signifie donc que la \emph{dimension} de $A/\p$ est $>0$.}.
Pour montrer la proposition, on peut donc supposer $\p$ maximal et $A$ local intègre. Alors, $\p[X]$
est premier non maximal dans $A[X]$ de sorte que le complété $\p[X]$-adique $\chap{A[X]_{\p[X]}}$
est formellement équidimensionnel. Comme $\chap{A}\ra\chap{A[X]_{\p[X]}}$ est local et plat,
l'argument de platitude (\egalong{IV}{2}{7.1.3}) utilisé plus haut assure que $\chap{A}$ est
équidimensionnel.

\end{proof}

\section{Hensélisation et anneaux excellents}
\label{1-henselisation-anneaux-excellents}

Rappelons qu'un morphisme d'anneaux noethériens $A\ra B$ est dit absolument plat s'il est réduit à
fibres discrètes et si les extension résiduelles sont algébriques et séparables. Ou, de façon
équivalente, s'il est plat ainsi que le morphisme diagonal $B\otimes_AB\ra B$ (cf. \cite{Fer72},
prop. 4.1 et~\cite{Oli71}, 3.1). Lorsque $B$ est (localement) de type fini sur $A$, ceci équivaut au
fait que $B$ soit étale sur $A$. En particulier, les extensions résiduelles sont séparables de sorte
qu'un tel morphisme est en fait régulier. Par exemple, tout morphisme ind-étale est absolument plat.
On a alors le résultat suivant (\cite{Gre76}).

\begin{theoreme}\label{1th-greco} Soit $f:A\ra B$ un morphisme absolument plat d'anneaux
noethériens.  Alors \begin{itemize} \item[i)] Si $A$ est vérifie 2.a) (resp.  2.b) ou 2.c)), $B$ un
vérifie 2.a) (resp.  2.b) ou 2.c))\footnote{Ou plus généralement, si $A$ est un $\mathbf{P}$-anneau
au sens de Grothendieck (\egalong{IV}{2}{7.3}), $B$ est un $\mathbf{P}$-anneau}.  \item[ii)] Si $A$ est
universellement-japonais, $B$ est universellement japonais.  \item[iii)] Si $A$ est quasi-excellent,
$B$ est quasi-excellent.  \item[iv)] Si $A$ est excellent, $B$ est excellent.  \item[v)] Si $f$ est
fidèlement plat, la réciproque de i), ii) et iii) est vraie.  \item[vi)] Si $f$ est  fidèlement plat
et $B$ est localement intègre, la réciproque de iv) est vraie.  \end{itemize} \end{theoreme}

Comme les morphismes d'hensélisation et de stricte hensélisation sont absolument plats
(\egalong{IV}{4}{18.6.9 et 18.8.12}) et fidèlement plats (ils sont locaux), on trouve, en particulier,
que la quasi-excellence et l'excellence sont stables par hensélisation et hensélisation stricte et
que, le caractère quasi-excellent ou universellement japonais d'un anneau local se teste sur
l'hensélisé ou l'hensélisé strict. Dans le cas de l'hensélisé, ces résultats étaient connus de
Grothendieck (\egalong{IV}{4}{18.7}), notamment 18.7.6).

En revanche, on ne peut espérer une propriété de descente de l'excellence comme en vi) sans
condition d'intégrité locale (cf.\ref{1pathologie}).

\section{Complétion formelle et anneaux excellents}

Soit $I$ un idéal d'un anneau noethérien $A$ contenu dans son radical de Jacobson et $\gtilde A$ sa
complétion $I$-adique. On peut se demander si les propriétés d'excellence passent au complété. La
réponse est oui en général. Précisément, on a :

\begin{proposition} Soit $I$ un idéal d'un anneau noethérien $A$ contenu dans son radical de
Jacobson et $\gtilde A$ sa complétion $I$-adique.  \begin{itemize} \item [i)] Si $A$ est
(semi)-local quasi-excellent (resp. excellent), il en est de même de $\gtilde A$ ; \item [ii)] Si
$A$ est une $\QQ$-algèbre excellente, il en est de même de $\gtilde A$.  \end{itemize}
\end{proposition}

La permanence de la quasi-excellence dans le cas (semi)-local, \emph{i.e.} de la géométrique régularité des
fibres formelles, est due à Rotthaus (\cite{Rot80}), tandis que celle de l'universelle caténarité
est due à Seydi\footnote{Comme me l'a expliqué Christel Rotthaus (communication privée), si $A, B$
sont locaux tels que $A\subset B \subset \chap{A}$ et $\chap{B}=\chap{A}$, alors l'universelle
caténarité de $A$ entraîne celle de $B$.} (le théorème 1.12 de \cite{Sey70} prouve qu'un anneau de
série-formelles $A[[t_1,\cdots,t_n]]$ est universellement caténaire dès que $A$ l'est ; il suffit
alors de considérer des générateurs $i_1,\cdots,i_n$ de $I$ définissant une surjection
$A[[t_1,\cdots,t_n]]\twoheadrightarrow \gtilde A$). Pour \textit{ii)}, reste à étudier l'ouverture
du lieu régulier dans le.  C'est ce qui est fait dans \cite{Bro-Rot80}, en utilisant le théorème de
désingularisation d'Hironaka de façon cruciale. Les techniques de \cite{Bro-Rot80} ont d'ailleurs
permis de montrer que si le théorème de désingularisation était vrai dans le cas local excellent,
toute complétion $I$-adique d'anneau excellent comme plus haut serait excellente (\cite{Nis-Nis87}).

En fait, le résultat est général. Plus précisément, Gabber (\cite{Gab07}) peut remplacer le théorème
d'Hironaka par son théorème d'uniformisation (\refext{VII}{7-uniformisation}) dans les arguments de
\cite{Nis-Nis87} pour prouver le résultat suivant

\begin{theoreme}[Gabber]Soit $A$ un anneau noethérien $I$-adiquement complet.  Alors, si $A/I$ est
quasi-excellent, $A$ est quasi-excellent.\end{theoreme}

On ne peut remplacer quasi-excellent par excellent dans le théorème précédent. En effet, Greco
(\cite{Gre82}) a construit un idéal $I$ d'un anneau $A$ intègre de dimension $3$, noethérien
semi-local $I$-adiquement complet et séparé qui est quasi-excellent non excellent alors que $A/I$
est excellent. On peut même supposer que $A$ est une $\QQ$-algèbre. La construction se fait par
pincements d'idéaux maximaux de hauteurs différentes (cf.~\ref{1ex-non-uni-cat}). Malgré tout, comme
on vient de le voir, la formelle caténarité passe aux complétions partielles (\cite{Sey70}) de sorte
que le complété $I$-adique d'un anneau excellent $A$ est excellent dès lors que $I$ est contenu dans
le radical de Jacobson de $A$.

\section{Approximation d'Artin et anneaux excellents}

Rappelons la définition suivante (cf. \cite{Art69}).

\begin{definition}[M. Artin] Un anneau local noethérien $(A,\m)$ a la propriété d'approximation
(\emph{AP})  si pour toute variété affine $X$ de type fini sur $A$, l'ensemble $X(A)$ est dense dans
$X(\chap{A})$.  \end{definition}

Bien entendu, il revient au même de dire que pour tout $X$ comme plus haut, on a
$$X(\chap{A})\not=\emptyset\Rightarrow X(A)\not=\emptyset.$$

Si $A$ vérifie AP, $A$ est certainement hensélien. Mais l'exemple~\ref{1DVR-bad} prouve qu'il ne
suffit pas que $A$ soit hensélien pour qu'il possède la propriété d'approximation. En fait, Rotthaus
a observé que l'excellence était une condition nécessaire à l'approximation d'Artin :

\begin{lemme}[\cite{Rot90}] Un anneau local noethérien vérifiant AP est hensélien et excellent.
\end{lemme}

Soit $k$ un corps de caractéristique nulle muni d'une valuation non triviale à valeurs réelles.
Artin a prouvé (\cite{Art68}) que les anneaux de séries convergentes à coefficients dans $k$ (pas
nécessairement complets) ont la propriété d'approximation. Ils sont donc henséliens et
excellents.\footnote{Dans le même ordre d'idées, l'anneau $k\{x_1,\cdots,x_n\}$ des séries formelles
restreintes (séries dont la suite des coefficients tendent vers $0$) à coefficients dans un corps
valué complet non archimédien $k$ est excellent dès que $k$ est de caractéristique nulle ou que $k$
est de degré fini sur $k^p, p=\text{car}(k)$ (\cite{Gre-Val74}). Le cas général a été obtenu par
Kiehl (\cite{Kie69} et aussi \cite{Con99} pour une preuve et des développements). Ceci répond,
partiellement, à une question de Grothendieck (\egalong{IV}{2}{7.4.8 B)}). En revanche, si $k$ est
valué non archimédien non complet de caractéristique positive tel que le morphisme de complétion
$k\ra\chap{k}$ n'est pas séparable, Gabber sait prouver que $k\{x_1\}$ n'est pas excellent.}

La situation est maintenant complètement clarifiée grâce aux travaux de Popescu culminant avec le
résultat suivant (\cite{NPD@Swan}) :

\begin{theoreme}[Popescu] \label{1-Artin-Popescu}\begin{itemize} \item [i)]Soit $A\ra B$ un
morphisme   régulier d'anneaux noethériens. Alors, $B$ est limite inductive filtrante de
$A$-algèbres lisses.  \item [ii)] Tout anneau local noethérien hensélien et excellent satisfait la
propriété d'approximation \emph{AP}\footnotemark.  \end{itemize} \end{theoreme} \footnotetext{Voir
\cite{Spi99}, th. 11.3 pour un énoncé un peu plus général.}

Le fait que i) entraîne ii) est un simple exercice. En effet, si $A$ est quasi-excellent, le
morphisme de complétion $A\ra \chap{A}$ est régulier et donc on peut écrire $\chap{A}=\limind\ L$ où
$L$ est lisse sur $A$. Soit $X=\Spec(B)$ avec $B$ de type fini et $\hat a\in X(\chap{A})$ d'image
$\hat a(0)\in X(k)$ où $k$ corps résiduel de $A$.  Il existe donc $L$ lisse sur $A$ tel que $\hat a$
provienne de $l\in X(L)$. Comme $A$ est hensélien, il existe $a\in X(A)$ (tel que $a(0)=\hat
a(0))$\footnote{L'argument n'utilise que la géométrique régularité des fibres -- et le caractère
hensélien -- mais pas la formelle caténarité. Ce n'est pas paradoxal, car un anneau local hensélien
est excellent si et seulement si ses fibres formelles sont géométriquement régulières
(\ref{1loc-exc}).}.

\section{Exemples de méchants anneaux noethériens}\label{1pathologie} Soit $A$ un anneau noethérien
et $X=\Spec(A)$ le schéma affine correspondant.  Il ressortira de cet inventaire que les propriétés
désagréables des anneaux du point de vue de l'excellence n'ont en général pas seulement à voir avec
la caractéristique $>0$ mais peuvent aussi se produire pour des $\QQ$-algèbres.

\subsection{Formelle caténarité : condition 3)} \label{1ex-non-uni-cat} Regardons d'abord de mauvais
anneaux du point de vue de la formelle caténarité.  \subsubsection{La caténarité n'entraîne pas la
formelle caténarité.} Dans \egalong{IV}{2}{5.6.11}, Grothendieck construit un exemple d'anneau local
noethérien de dimension $2$ intègre, caténaire mais non universellement caténaire, donc non
formellement caténaire d'après~\ref{1form-cat->uni-cat} (\emph{i.e.} ne vérifiant pas 3)).

Expliquons la construction qui consiste à pincer une surface lisse sur un corps $k$, de
caractéristique nulle si l'on veut, le long de deux points de hauteurs différentes ayant des corps
résiduels isomorphes à $k$.

On part d'un corps $k$ extension transcendante pure de degré  infini sur son corps premier par
exemple, que l'on peut même supposer être $\QQ$. Soit $S$ une surface lisse munie d'un morphisme
projectif sur $T=\Spec(k[\tau])$ et $t\in S(T)$. On suppose qu'il existe un point $s\in S(k)$ de la
fibre $S_0$ de $S\ra T$ au dessus de $0\in T$ qui n'est pas dans l'image de $t$. Par exemple, on
peut prendre $S=\Spec(k[\sigma,\tau])$ avec $t$ la section d'image $\sigma=0$ et $s=(1,0)$. Les
corps $k(s)=k$ et $k(t)=\Frac(k[\tau])$ sont des extensions transcendantes pures de $\QQ$ de même
degré (infini) de sorte qu'on peut choisir un isomorphisme de corps $k(s)\isononcan k(t)$.  Ceci
permet de définir le sous-anneau $\OO_\Sigma$ de $\OO_S$ des fonctions qui coïncident en $s$ et $t$.
On dispose donc d'un morphisme $\pi:S\ra \Sigma$ qui envoie $s,t$ sur $\sigma\in\Sigma$.  Posons
$A=\OO_{\Sigma,\sigma}$ et soit $B$ l'anneau de coordonnées de
$S\times_\Sigma\Spec(\OO_{\Sigma,\sigma})$. Par construction, $\dim(B_s)=2$ et $\dim(B_t)=1$. Alors,
$A$ est noethérien, et $B$ est la normalisation de $A$ et est fini sur $A$. Comme $A$ est de
dimension $2$ et intègre il est évidemment caténaire. Si $A$ était universellement caténaire, la
formule de dimension (voir note~(\ref{1form-dim})) entraînerait $\dim(A)=\dim B_s=\dim B_t$, une
contradiction.

On peut même trouver pour tout $n\geq 2$ des anneaux locaux noethériens intègres de dimension $n$
vérifiant 2.a), caténaires non universellement caténaires donc ne vérifiant pas 3)
(\cite{Hei-Rot-Wie00}).

\subsubsection{La formelle caténarité ne se teste pas sur l'hensélisé.} Par des techniques de
pincements analogues de surfaces sur $k$, donc de caractéristique nulle si on veut, comme plus haut,
Grothendieck construit en effet un exemple d'anneau local non universellement caténaire (donc non
formellement caténaire) dont l'hensélisé  est excellent (\egalong{IV}{4}{18.7.7}). Quitte à changer de
base par la clôture séparable, on s'aperçoit que la formelle caténarité ne se teste pas plus sur
l'hensélisé strict, contrairement à la quasi-excellence (\ref{1th-greco}).

\subsubsection{La formelle caténarité n'entraîne certainement pas 2.a) (ni même 2.c)).} Par exemple,
un anneau de valuation discrète $A$ non excellent (cf. la partie~\ref{1patho-DVR}) a une fibre
formelle générique non géométriquement régulière (en effet, il est formellement caténaire
(\ref{1CM-form-cat}) et 2.a) entraîne 1.a) dans le cas local (\ref{1loc-exc})). Or, cette fibre
générique formelle est artinienne dans ce cas (elle ne contient pas l'idéal maximal de $\chap{A}$)
et donc la géométrique régularité équivaut ici à la géométrique réduction.

\subsubsection{Il existe des anneaux intègres \emph{normaux} non formellement caténaires.} Ogoma
(\cite{Ogo81}) a construit une $\QQ$-algèbre locale $A$ intègre \emph{normale} de dimension $3$ dont
le complété à une composante de dimension $2$ et une composante de dimension $3$ et donc n'est pas
équidimensionnel. Pire, cet anneau n'est même pas caténaire : il possède une infinité de chaînes
saturées d'idéaux premiers de longueur 2 ou 3.

\subsection{Quasi-excellence : conditions d'ouverture  1.a) et 1.b)}\label{1ex-non-1} On s'intéresse
ici à des anneaux ayant un lieu régulier ou normal non ouvert.

 Comme on le verra plus loin (\refext{XIX}{19-rem-cexgabber}), Gabber a construit un exemple de
schéma, qu'on peut même supposer être un $\QQ$-schéma,  intègre de dimension $1$, dont le lieu
régulier (ou normal, c'est la même chose ici) contient une infinité de points et en particulier
n'est pas ouvert. La construction assure que les fibres formelles sont géométriquement régulières.
Comme on est en dimension $1$, normalité et régularité coïncident de sorte qu'on a un exemple
vérifiant 2.a) et 3) mais pas 1.b).

 Dans~\cite{Rot79}, Rotthaus construit une $\QQ$-algèbre noethérienne locale intègre de dimension
$3$ qui est formellement caténaire, universellement japonaise mais dont le lieu régulier n'est pas
ouvert.

\subsection{Quasi-excellence. Fibres formelles : conditions 2a), 2b) et 2c)}\label{1ex-non-uni-2} On
s'intéresse ici à des anneaux ayant des fibres formelles non géométriquement régulières voire pire.

\begin{itemize} \item Rotthaus construit une $\QQ$-algèbre locale $A$ noethérienne de dimension $3$
régulière (donc formellement caténaire), universellement japonaise mais pas excellente
(\cite{Rot79}).  Précisément, la fibre formelle au dessus d'un point de hauteur $1$ n'est pas
régulière. Ainsi, elle vérifie 2.c), 3) car $A$ régulier mais pas 2.a).

Dans l'exemple d'Ogoma précédent, la fibre formelle générique est connexe non intègre (elle a une
composante de dimension $1$ et une de dimension $2$ qui se coupent), donc non normale. On a donc un
exemple de $\QQ$-algèbre locale (de dimension $3$) noethérienne (intègre et normale) ne vérifiant
pas 2.b).

On peut descendre d'une dimension : Nagata construit (\cite{Nag62}, ex. 7 de l'appendice A1) une
$\QQ$-algèbre locale $B$ qui est intègre normale, formellement caténaire et de dimension $2$ mais
dont le complété n'est pas intègre\footnote{La construction est la suivante : soient $x,y$
algébriquement indépendants sur $\QQ$ et $w=\sum_{i>0}a_ix^i\in\QQ[[x]]$ transcendant sur $K(x)$. On
pose $z_1=(y+w)^2$ et $z_{i+1}=(z-(y+\sum_{j<i}a_jx^j)^2)/x^i$. Soit $A$ le localisé de
$\QQ[x,y,z_i,i\geq 1]$ en $(x,y,z_i,i\geq 1)$. Alors, $B=A[X]/(X^2-z_1)$ est l'exemple cherché. On
vérifie facilement que la complétion de $A$ est $\QQ[[x,y]]$ de sorte que $\chap{B}=
\QQ[[x,y]][X]/(X^2-(y+w)^2)$ n'est pas intègre. Comme d'habitude dans ces constructions, c'est le
caractère noethérien de $A$ qui pose problème. Une fois ceci acquis, $A$ est régulier de dimension
$2$ et $B$ normal puisque que Cohen-Macaulay de dimension $2$ singulier uniquement à l'origine.
Notons que $B$ est formellement caténaire comme quotient d'un régulier.}.
D'après~\ref{1prop-normal-versus-complet}, ceci prouve que $B$ ne vérifie pas 2.b) (mais vérifie
3)).

\item En caractéristique $>0$, Rotthaus  construit également (\cite{Rot79}) une algèbre locale
noethérienne de dimension $2$ régulière universellement japonaise mais pas excellente. Dans ce
denier cas, comme les fibres formelles sont de dimension $<2$, elles ne sont pas non plus
géométriquement normales. Ainsi, elle vérifie 2.c), 3)  mais pas 2.b).

\item En caractéristique nulle, Ferrand et Raynaud ont construit (\cite{Fer-Ray70}, prop. 3.3 et 3.5)
une $\CC$-algèbre locale noethérienne $A$ intègre de dimension $2$ telle que

$\bullet$ le normalisé $A'$ de $A$ est l'anneau des séries convergentes noté $\CC\{x,y\}$ (ne pas
confondre avec l'hensélisé de $\CC[x,y]$) et donc est excellent.

$\bullet\ A$ n'est pas japonais (en fait, $A'$ n'est pas fini sur $A$).

$\bullet\ \chap{A}$ a des composantes immergées (de sorte que -- platitude -- la fibre formelle
générique a des composantes immergées et donc ne vérifie pas 2.c)).

$\bullet\ $ Le lieu normal de $A[[T]]$ n'est pas ouvert.  D'après~\egalong{IV}{2}{6.13.5} son lieu
régulier n'est donc pas ouvert non plus.

$\bullet$ L'anneau $A$ est formellement caténaire (le spectre de son complété est irréductible). Il
en est donc de même de $A[[T]]$ (\cite{Sey70}).

\item Dans~\cite{Nag62}, ex. 5 de l'appendice A1, Nagata construit même un anneau local noethérien
intègre de dimension $3$ (de caractéristique $>0$) dont la clôture intégrale n'est même pas
noethérienne ; en particulier, cet anneau n'est pas japonais\footnote{La construction est du même
type que celle d'un anneau de valuation discrète décrite dans la note (\ref{1note-DVR-bad}) dont on
reprend les notations. On considère cette fois-ci l'anneau $B=k^p[[X,Y,Z]][k]$ et
$d=Y\sum_{i>0}X_iX^i+Z\sum_{i>0}X_{2i+1}X^i$. L'anneau $B[d]$ convient.}.

\item Pire, à partir d'anneaux construits par Nagata, Seydi construit (\cite{Sey72}) un anneau
noethérien  intègre normal $A$ de dimension $3$ dont le corps des fractions est de caractéristique
nulle et dont le complété n'est pas réduit. En particulier, il est japonais mais pas universellement
japonais. Ogoma construit (\cite{Ogo81}) une $\QQ$-algèbre noethérienne normale, donc japonaise, qui
n'est ni universellement japonaise ni caténaire.

\item Les fibres formelles peuvent être épouvantables, même en dimension $1$ : Ferrand et Raynaud
construisent un $\CC$-schéma local intègre de dimension $1$ dont la fibre formelle générique est un
schéma artinien qui n'est même pas Gorenstein -- donc certainement non réduit -- (\cite{Fer-Ray70},
prop. 3.1) : $X$ ne vérifie pas 2.c). En particulier, $X$ n'est pas universellement japonais (et
donc pas quasi-excellent). Bien entendu, $X$ vérifie 1.a) et 3) pour des raisons de dimension.

\item  Les exemples d'anneaux de valuation discrète non excellents (donc de caractéristique
positive) donnent des exemples d'anneaux ne vérifiant pas 2.c) (2.a) et 2.c) sont équivalents en
dimension $\leq 1$) mais vérifiant 1.a) et 3).

\end{itemize}

\begin{remarque}\label{1DVR-bad}Nagata a construit (\cite{Nag62}, (E3.3)) un anneau de valuation
discrète dont la fibre formelle générique est une extension radicielle non triviale de son  corps
des fractions, donc non excellent\footnote{\label{1note-DVR-bad}Voici la construction : soit $k$ le
corps des fractions de $\FF_p[X_n,n> 0]$ et $K$ celui de $\chap{A}=k[[Y]]$. Soit $L$ le sous-corps
de $K=k((Y))$ corps des fractions de $A=k^p[[Y]][k]$. Le complété de $A$ est $\chap{A}$. On montre,
et c'est le point délicat, que $A$ est noethérien. L'outil est le critère de Cohen : un anneau
semi-local est noethérien si et seulement si les idéaux maximaux sont de type fini et les idéaux de
type fini fermés (\cite{Nag62}, 31.8). Son complété étant régulier, il est lui même régulier donc de
valuation discrète (dimension).  Soit $L$ le corps des fractions de $A$. On vérifie facilement que
$c=\sum_{n> 0}X_nY^n,$ n'est pas dans $L$. Choisissons une $p$-base $\{c_i\}$ de $K$ sur $L$
contenant $c$ (ce qui est possible car $c\not\in L^p$, cf.~\egalong{IV}{1}{21.4.3}). Soit $K_0$ le
corps engendré sur $L$ par les $c_i$ distincts de $c$. L'extension $K/K_0$ est radicielle de degré
$p$ par construction. L'anneau $A\cap K_0$ est un anneau de valuation discrète de complété $k[[Y]]$
de sorte que la fibre formelle générique n'est pas géométriquement réduite.}.  \end{remarque}

 On va voir maintenant que de tels anneaux se rencontrent très facilement.

\subsection{Méthode systématique de construction d'anneaux non quasi-excellents}\label{1patho-DVR}
En fait, on peut construire (Orgogozo) de façon systématique de très nombreux anneaux de valuation
discrète non quasi-excellents. Précisons\footnote{Cette construction généralise en fait, de façon
indépendante, un exemple obtenu par Rotthaus dans~\cite{Rot97}}.

\begin{proposition}\label{1lot-bad-DVR}Soit $k((t))$ le corps des séries de Laurent à coefficients
dans un corps $k$ de caractéristique $p>0$ muni de sa valuation $t$-adique et $L/k$ une
sous-extension de type fini de $k((t))/k$ de degré de transcendance $>1$ sur $k$. Alors, le
sous-anneau $A$ de $L$ des éléments de valuation $\geq 0$ est un anneau de valuation discrète non
excellent.  \end{proposition}

\begin{proof}

Soit ${\bar L}$ un corps de caractéristique $>0$. Le $p$-rang est la dimension, finie ou non, de
$\Omega_{\bar L}$, le module des différentielles absolues. C'est aussi $\log_p([{\bar L}:{\bar
L}^p])$ où $[{\bar L}:{\bar L}^p]$ est la dimension de ${\bar L}$ sur ${\bar L}^p$
(\egalong{IV}{1}{21.3.5}).  La remarque clef est que le $p$-rang croit par extension de corps
\emph{séparable} ${\bar K}/{\bar L}$,  finie ou non, puisqu'on a dans ce cas un plongement ${\bar
K}\otimes_{\bar L}\Omega_{\bar L}\hookrightarrow \Omega_{{\bar K}}$ (\egalong{IV}{2}{20.6.3})
\begin{equation}\label{1ineq-sep} [{\bar L}:{\bar L}^p]\leq [{\bar K}:{\bar K}^p].  \end{equation}

Par ailleurs, si ${\bar L}$ est de \emph{type fini} sur un corps $k$, on a (\cite{Bou50}, V.6.3)

\begin{equation}\label{1prang} [{\bar L}:{\bar L}^p]=p^{\text{degtr}_k({\bar L})}[k:k^p].
\end{equation}

Plaçons nous dans la situation du lemme. L'anneau $A$ est de valuation discrète par construction et
son corps des fractions est $L$. L'hensélisé $A^h$ est local régulier de dimension $1$ donc intègre
et son corps des fractions $K=\Frac(A^h)$ contient $L=\Frac(A)$.  Le  complété $\chap{A^h}$ est un
anneau de séries formelles $\chap{K}=k[[\varpi]],\varpi $ uniformisante de $A^h$ (comme complété
d'une $k$-algèbre locale régulière de dimension $1$) et son corps des fractions  $\chap{K}$ est la
fibre générique formelle de $\Spec(\chap{A^h})\ra\Spec(A^h)$.

Supposons que $A^h$ soit quasi-excellent (précisément vérifie 2.a) de sorte que l'extension
$\chap{K}/K$ est \emph{séparable}.

On a donc dans ce cas $$[K:K^p]\leq[\chap{K}:\chap{K}^p].$$ Comme $\chap{K}=k((\varpi)),$ on a
$$[\chap{K}:\chap{K}^p]=p[k:k^p].$$

On a donc (\ref{1ineq-sep}), l'extension $K/L$ étant séparable, $$[K:K^p]\geq [L:L^p].$$ de sorte
que, grâce à (\ref{1prang}), on a $$p[k:k^p]=[\chap{K}:\chap{K}^p]\geq [L:L^p]>p[k:k^p],$$ une
contradiction. Ceci interdit à $A$ également d'être quasi-excellent (\ref{1th-greco}).\end{proof}

\ifx\danslelivre\undefined
\bibliography{../fichiers-configuration-latex/bibliographie}
\bibliographystyle{../fichiers-configuration-latex/smfalpha}


\begin{center}
II. Topologies adaptées à l'uniformisation locale\\
Fabrice Orgogozo
\end{center}

\version
\tableofcontents

\else
\setboolean{english}{false}
\setboolean{francais}{true}
\selectlanguage{francais}
\chapter{Topologies adaptées à l'uniformisation locale}
\label{2-topologies}
\addtocontents{toc}{par Fabrice Orgogozo}
\begin{center}
Fabrice Orgogozo
\end{center}
\fi

Dans cet exposé, $ℓ$ est un nombre premier, l'entier $1$
ou bien le symbole $∞$ et l'on note $ℓ′$ l'ensemble
des entiers naturels premiers à $ℓ$ où, par convention, $∞′=\{1\}$.

\section{Morphismes maximalement dominants et la catégorie $\alt ∕ S$}

\subsection{Morphismes maximalement dominants}

\subsubsection{}
Rappelons (\ega{IV}{1.1.4}) qu'un point d'un schéma est dit \emph{maximal} s'il est le
point générique d'une composante irréductible ou, de façon
équivalente, s'il est maximal pour l'ordre sur l'ensemble des
points du schéma défini par la relation : $x≥y$ si et seulement si $y$ est
une spécialisation de $x$ (c'est-à-dire si $y∈\sur{\{x\}}$).
Les points maximaux d'un schéma affine correspondent aux idéaux
premiers minimaux. Tout ouvert dense d'un schéma contient
la totalité des points maximaux.

\begin{définition2}\label{2-maximalement-dominant}
Un morphisme de schémas est dit \emph{maximalement dominant}\index{maximalement dominant}
s'il envoie tout point maximal de la source sur un point maximal du but.
\end{définition2}

Un morphisme entre schémas irréductibles est maximalement dominant si et seulement si il est
dominant. Il est clair que le composé de deux morphismes maximalement dominants est
maximalement dominant.

\begin{exemple2}
D'après \ega{IV₂}{2.3.4}, un morphisme plat, ou plus généralement quasi-plat (\emph{op. cit.},
2.3.3), est générisant (\ega{I′}{3.9.1}) donc maximalement dominant (\emph{op. cit.} 3.9.5).
\end{exemple2}

\begin{proposition2}\label{2-image maximalement dominant}
Soit $f:X→Y$ un morphisme maximalement dominant. Tout point maximal de $Y$ appartenant à l'image de $f$
est l'image d'un point maximal de $X$.
\end{proposition2}

\begin{démo}
Cela résulte du fait que $f$ est croissante pour le préordre ci-dessus
et du fait que tout point de $X$ a une générisation maximale.
\end{démo}

\begin{proposition2}\label{2-cb maximalement dominant par plat}
Soient $f:X→Y$ un morphisme maximalement dominant et $Y'→Y$ un morphisme
\emph{plat}. Alors, le morphisme $X'=X×_Y Y'→Y'$ est maximalement dominant.
\end{proposition2}

\begin{démo}
Cf. \ega{IV₂}{2.3.7} (ii).
\end{démo}

Rappelons la proposition suivante.

\begin{proposition2}[\ega{iv}{20.3.5}]
Soit $f:X→Y$ un morphisme maximalement dominant. Supposons que $X$
est \emph{réduit} et que $Y$ n'a qu'un nombre fini de composantes
irréductibles (par exemple, $Y$ nœthérien).
Alors, pour tout ouvert $U$ de $Y$ et tout ouvert $V$
dense dans $U$, l'ouvert $f^{-1}(V)$ est dense dans $f^{-1}(U)$.
\end{proposition2}

L'hypothèse sur $Y$ assure que tout point maximal de $U$ appartient à $V$.

\begin{proposition2}
\label{2-génériquement fini}
Soit $Y$ un schéma nœthérien et soit $f:X→Y$ un morphisme de type fini, \emph{maximalement
dominant}. Les conditions suivantes sont équivalentes :
\begin{enumerate}
\item Il existe un ouvert dense de $Y$ au-dessus duquel $f$ est fini.
\item Pour tout point maximal $η$ de $X$, l'extension $κ(η)/κ(f(η))$
est finie.
\end{enumerate}
\end{proposition2}

Un morphisme de schémas satisfaisant la condition (i)
ci-dessus est souvent dit \emph{génériquement fini} (en bas).

\begin{démo}
L'implication (i)⇒(ii) résulte du Nullstellensatz
et n'utilise pas les hypothèses de finitude faites ($Y$ nœthérien, $f$ de type fini).
Démontrons la réciproque. On peut supposer $Y$ irréductible et $X,Y$ réduits : on utilise le fait que $f$ est
fini si et seulement si $f_{\red}$ l'est. Par passage à la limite
(\ega{iv}{8.10.5.(x)}), on peut également supposer
que $Y$ est le spectre d'un corps $k$. On peut également
supposer $X$ irréductible donc intègre. Le résultat
est alors conséquence du fait qu'une algèbre intègre
entière sur un corps est un corps.
(Alternativement, on peut utiliser
l'inégalité $\dim(X)≤\dim(Y)=0$,
qui est un cas particulier de \ega{iv}{5.6.6}.)
\end{démo}

\begin{définition2}\label{2-L-maximalement fini}
On dit qu'un morphisme $f:X → Y$ est \emph{maximalement $ℓ′$-fini}
si pour tout point maximal $η$ de $X$, l'extension $κ(η)/κ(f(η))$
est finie et si pour tout point maximal $μ$ de $Y$ dans l'image de $f$,
il existe un point maximal $η$ de $X$ au-dessus de $μ$
tel que l'extension $κ(η)/κ(μ)$ soit de degré appartenant à $ℓ′$.
\end{définition2}

\subsubsection{}Lorsque $ℓ=1$, la seconde condition est vide.
Il est utile de faire les conventions de langage suivantes :
un morphisme maximalement $ℓ′$-fini et maximalement dominant
est dit \emph{maximalement $ℓ′$-fini dominant},
et un morphisme maximalement $1′$-fini « \emph{maximalement fini} ».

\begin{proposition2}\label{2-morphisme entre maximalement dominants gen finis}
Soient $S$ un schéma et $f:X→Y$ un $S$-morphisme entre $S$-schémas
maximalement dominants.
Si $Y/S$ est maximalement fini, $f$ est maximalement dominant.
Si l'on suppose de plus $X ∕ S$ maximalement fini, le morphisme $f$ est maximalement fini dominant.
\end{proposition2}

\begin{démo}
Soient $x$ un point maximal de $X$ et $s$ (resp. $y$) son image dans $S$
(resp. $Y$). Soit $y'≥y$ une générisation maximale de $y$. Les schémas
$X$ et $Y$ étant maximalement dominants sur $S$, $s$ est un point maximal et $y'$
est d'image $s$. Enfin, si $Y/S$ est maximalement fini, 
l'extension $κ(y')/κ(s)$ est finie de sorte que $y'$ est isolé dans la
fibre $Y_s$. Le point $y$, appartenant à l'adhérence de $y'$ dans
$Y_s$, coïncide donc avec $y'$ : le morphisme $f$ est maximalement dominant.
Si l'on suppose de plus $X/S$ maximalement fini,
la finitude de l'extension $κ(x)/κ(s)$ entraîne celle de l'extension
intermédiaire $κ(x)/κ(y)$.
\end{démo}

\subsection{La catégorie $\alt ∕ S$}

\subsubsection{}Soit $S$ un schéma nœthérien. Notons $η_S$ le schéma
coproduit (fini) de ses points maximaux.

\begin{définition2}
\label{2-alt}
On note $\alt/S$ \index{alt} la catégorie des $S$-schémas de type fini,
maximalement finis dominants, de source un schéma \emph{réduit}.
Les morphismes dans $\alt/S$ sont les $S$-morphismes.
\end{définition2}

\subsubsection{}Notons les faits suivants :
\begin{itemize}
\item le $S$-schéma $S_{\red}$ est final dans la catégorie $\alt /S$ ;
\item tout morphisme de $\alt/S$ est maximalement fini dominant
(\ref{2-morphisme entre maximalement dominants gen finis}) ;
\item les images inverses de diviseurs existent pour tout morphisme de $\alt/S$ (\ega{iv}{21.4.5.(iii)}) ;
\item si $X∈\Ob\,\alt/S$ et $S'→S$ est un morphisme \emph{réduit} (\ega{IV₂}{6.8.1}) avec $S'$
nœthérien, le produit fibré usuel $X'=X×_S S'$ est naturellement un objet de $\alt/S'$. Il
en est plus généralement ainsi de $X ′_\red$ dès lors que $S ′ → S$ est plat.
\end{itemize}

\begin{remarques2}
\begin{enumerate}
\item Le produit fibré usuel de deux $S$-schémas maximalement dominants
n'est pas nécessairement maximalement dominant, comme on peut
le constater lorsque $S=\Aff²$ et $X=Y$ sont l'éclatement en l'origine.
\item La définition originale de la catégorie $\alt ∕ S$, due à O. Gabber,
est moins restrictive sur le schéma $S$ :
il est supposé cohérent et ayant un nombre fini de composantes
irréductibles. Les objets de $\alt ∕ S$ sont alors les
$S$-schémas de type fini quasi-séparés réduits, maximalement dominants,
génériquement finis. Le cadre nœthérien semble suffisant pour
nos besoins. Signalons cependant que les « localisés »
d'un schéma nœthérien pour la topologie des altérations
(introduite ci-dessous) ne sont pas nécessairement nœthériens (cf. \ref{2-points de alt}).
\end{enumerate}
\end{remarques2}

\subsubsection{}
\label{2-dom}
Soit $X$ un $S$-schéma de type fini.
On note $X_\maxdom$ l'adhérence de l'image (ensembliste) de $X_{η_S}$
dans $X$, muni de la structure réduite. C'est la
réunion des composantes irréductibles de $X$ dominant une composante
irréductible de $S$, munie de la structure réduite.
Le foncteur $T↦T_\maxdom$ est adjoint à droite au foncteur
d'inclusion de la catégorie des schémas réduits, de type fini et
maximalement dominants sur $S$ dans la catégorie des $S$-schémas de type fini.

\begin{proposition2}\label{2-existence produit}
Les produits fibrés existent dans $\alt/S$.
\end{proposition2}

\begin{démo}
Soient $X→S'←Y$ deux flèches dans $\alt/S$ ; d'après
\ref{2-morphisme entre maximalement dominants gen finis}, les
schémas $X$ et $Y$ sont naturellement des objets de $\alt/S'$.
Le composé de deux morphismes maximalement finis dominants de type fini
étant de même nature, un produit de $X$ et $Y$, vus dans
$\alt/S'$, est — s'il existe — un produit fibré dans $\alt/S$.
On peut donc supposer $S=S ′$ et $S$ réduit.
Soient $X$ et $Y$ deux objets de $\alt/S$.
Il résulte formellement de l'existence du produit dans la catégorie des
$S$-schémas de type fini et de la propriété d'adjonction de $X↦ X_\maxdom$
que le schéma $(X×_S Y)_\maxdom$, muni des deux projections évidentes, est le produit
de $X$ et $Y$ dans la catégorie des schémas réduits, de type fini et
maximalement dominants sur $S$. Il appartient à $\Ob\,\alt/S$ car
$((X×_S Y)_\maxdom)_{η_S}=(X_{η_S}×_{η_S} Y_{η_S})_{\red}$ est fini sur $η_S$.
\end{démo}

\begin{proposition2}
\label{2-lprime-fini stable par cb}
La propriété d'être maximalement $ℓ′$-fini est stable par
changement de base.
\end{proposition2}

\begin{démo}
Il faut montrer que si $k ′ ∕ k$ est une extension finie de
degré premier à $ℓ$, où $k$ est de caractéristique $p$
et $K ∕ k$ est une extension quelconque
l'un des corps résiduels du produit tensoriel $K′=k'⊗_k K$ est de
degré premier à $ℓ$ sur $K$. On peut pour cela
supposer $k′ ∕ k$ étale ou radicielle.
Dans le premier cas, on écrit $K′= ∏_i K_i$, où $K_i$
est une extension étale de degré $d_i$ sur $K$, et on remarque
que si la somme $∑_i d_i=[k′:k]$ est première à $ℓ$,
il en est de même de l'un des $d_i$.
Le second cas ne se produit pas si $ℓ=p$ et est
trivial lorsque $ℓ≠p$, toute extension finie radicielle
étant de degré une puissance de $p$.
\end{démo}

\begin{proposition2}\label{2-adjoint preserve imm ouv, ps}
Soit $f:X→Y$ une immersion ouverte (resp. un morphisme propre et
surjectif, resp. quasi-fini) entre deux $S$-schémas de type fini. Le morphisme
$f_\maxdom:X_\maxdom→Y_\maxdom$ est une immersion ouverte (resp.
un morphisme propre et surjectif, resp. quasi-fini).
\end{proposition2}

\begin{démo}
Le cas d'une immersion ouverte est conséquence du fait général
suivant sur les topologies induites : la trace sur un \emph{ouvert} de l'adhérence d'une partie
coïncide avec l'adhérence de la trace de cette partie (cf. p. ex.
\bbk{TG}{I}{3}{1}{prop. 1} pour une variante).
Considérons maintenant le cas d'un morphisme propre et surjectif
$f:X→Y$. Par construction, les morphismes $X_\maxdom→X$ et
$Y_\maxdom→Y$ sont des immersions fermées ; les schémas
$X_\maxdom$ et $Y_\maxdom$ sont donc propres sur $Y$. Le $Y$-morphisme
$f_\maxdom$ est donc propre. Son image contient $f_{η_S}(X_{η_S})=Y_{η_S}$ donc son adhérence
$Y_\maxdom$. Le dernier cas est tout aussi trivial et laissé au lecteur.
\end{démo}

\subsubsection{}Notons que si $(U_i)_{i∈I}$ est un recouvrement
par des ouverts de Zariski d'un $S$-schéma de type fini $X$,
les ouverts $({U_i}_\maxdom)_{i∈I}$ recouvrent le schéma $X_\maxdom$.

\section{Topologies : définitions}

Dans ce paragraphe, on fixe un schéma nœthérien $S$.

\subsection{Topologie étale $ℓ′$-décomposée}

\subsubsection{}
\label{2-L-décomposé}
Nous dirons qu'un recouvrement étale $(X_i→ X)_{i∈I}$
d'un schéma $X$ est $ℓ′$-\emph{décomposé}
si tout point de $X$ peut être relevé en un point $x_i$ d'un $X_i$
tel que le degré $[κ(x_i): κ(x)]$ appartienne à $ℓ′$.
Lorsque $ℓ=1$ la condition imposée est vide et
l'on dit simplement que la famille constitue
un recouvrement étale. Lorsque $\mbox{ℓ=∞}$, on retrouve la définition de \cite[§1]{Nisnevic}
et l'on dit plutôt que le recouvrement est \emph{complètement} décomposé.
Il résulte de \ref{2-lprime-fini stable par cb} que
la propriété d'être un recouvrement étale $ℓ′$-décomposé est stable par
changement de base.

\subsubsection{}On appelle \emph{topologie étale $ℓ′$-décomposée}
la topologie de Grothendieck sur $\alt ∕ S$, notée $\et_{ℓ′}$,
définie par la prétopologie constituée des recouvrements
étales $ℓ′$-décomposés.

\subsection{Sorites sur le lieu $ℓ′$-décomposé}
\label{2-sorites lieu cd}

\subsubsection{}Pour chaque morphisme de schémas $f:Y → X$,
posons
\[
\mathrm{d\acute{e}c}_{ℓ′}(f) = \{x∈X:∃ y∈Y\textrm{ tel que } f(y)=x, [κ(x):κ(y)]\text{ fini appartenant à }ℓ′\}.
\]

Lorsque $ℓ=∞$, on retrouve l'ensemble $\mathrm{cd}(f)$ introduit
par Nisnevič. Nous utiliserons également cette notation.

\begin{proposition2}
Soit $f:Y → X$ un morphisme étale avec $X$ nœthérien.
L'ensemble $\mathrm{d\acute{e}c}_{ℓ′}(f)$ est
\emph{ind-constructible}, c'est-à-dire — $X$ étant nœthérien — réunion de parties localement fermées.
\end{proposition2}

\begin{démo}
On peut supposer $X$ et $Y$ intègres,
de points génériques notés $η$ et $μ$ respectivement.
Par récurrence nœthérienne, il suffit de montrer
que si $η$ appartient à l'ensemble $\mathrm{d\acute{e}c}_{ℓ′}(f)$,
celui-ci contient un ouvert de $X$.
On peut supposer de plus $X$ et $Y$ affines d'anneaux $A$ et $B$
respectivement, et le morphisme $A → B$ \emph{fini}.
La fonction $𝔭↦ \dim_{κ(𝔭)} B/𝔭B$, $X → 𝐍$,
est localement constante pour la topologie de Zariski
car $B$ est plat de présentation finie — donc localement
libre — sur $A$. Elle prend la valeur $[κ(μ): κ(η)]$ — première à $ℓ$ — en $η$.
Il en est donc ainsi au voisinage de $η$ ;
la conclusion en résulte aussitôt.
\end{démo}

Précisons un peu ce résultat lorsque $ℓ=∞$.

\begin{proposition2}\label{2-section cd}
Soit $f:Y→X$ un morphisme localement de type fini, avec $X$ nœthérien. Un point $x∈X$ appartient à
$\mathrm{cd}(f)$ si et seulement si il existe un sous-schéma $Z$ de $X$
contenant $x$ au-dessus duquel $f$ a une section.
\end{proposition2}

\begin{démo}
La condition est bien entendu suffisante.
Considérons $x∈\mathrm{cd}(f)$ ; c'est le point
générique du sous-schéma fermé réduit $\sur{\{x\}}$.
Par hypothèse, il existe une section au-dessus de ce point.
Le morphisme $f$ étant localement de présentation finie, cette section
s'étend par passage à la limite à un ouvert $Z=U⋂\sur{\{x\}}$
de $\sur{\{x\}}$, où $U$ est un ouvert de $X$.
\end{démo}

\begin{corollaire2}\label{2-cd noetherien}
Soient $X$ un schéma nœthérien et $(U_i\dessusdessous{f_i}{→}X)_{i∈I}$ un recouvrement de $X$ pour
la topologie étale $ℓ′$-décomposée. Il existe alors un sous-ensemble
fini $I₀⊂I$ tel que la famille $(U_i→X)_{i∈I₀}$ soit également couvrante.
\end{corollaire2}

Rappelons qu'un morphisme étale est, par définition, localement
de présentation finie.

\begin{démo}
D'après \ega{iv}{1.9.15}, l'espace topologique $X^{\mathrm{cons}}$,
dont l'espace sous-jacent est $X$ et dont les ouverts sont les parties
ind-constructibles de $X$, est \emph{compact}. Les ensembles $\mathrm{d\acute{e}c}_{ℓ′} (f_i)$
en constituent un recouvrement par des ouverts.
\end{démo}

\begin{proposition2}\label{2-cd limite}
Soit $(U_α\dessusdessous{f_α}{→}X_α)_{α∈A}$ un système projectif, filtrant, cartésien, de
morphismes étales entre schémas nœthériens, à morphismes de transition affines. Notons $f_∞:U_∞→X_∞$ le morphisme
induit sur la limite projective et, pour chaque $α$, $π_α$ la projection $X_∞ →  X_α$. 
On a 
$$
\mathrm{cd}(f_∞)=⋃_α π_α^{-1}\big(\mathrm{cd}(f_α)\big). 
$$
\end{proposition2}

\begin{démo}
L'inclusion du terme de droite dans le terme de gauche est évidente. 
Considérons réciproquement un point $x_∞$ dans $\mathrm{cd}(f_∞)$.
Le morphisme $f_∞$ a une section sur un sous-schéma de présentation
finie $Z_∞$ contenant $x_∞$. Le morphisme et la section se descendent par passage
à la limite à un niveau fini $α$ (cf. \ega{iv}{8.6.3, 8.8.2}).
\end{démo}

Les quatre énoncés précédents sont valables, \emph{mutatis mutandis},
lorsque l'on suppose simplement $X$ cohérent et le morphisme $f$
localement de présentation finie.
Pour référence ultérieure, signalons le lemme de descente
d'une section suivant.

\begin{proposition2}
\label{2-descente section}
Soient $k ′ ∕ k$ une extension finie de corps de degré premier à $ℓ$ et
$K ∕ k$ une extension finie de degré une puissance de $ℓ$. Posons $K
′= k′ ⊗_k K$. Si le morphisme $\Spec(K ′) → \Spec(K)$
possède une section, le morphisme $\Spec(k ′) → \Spec(k)$
possède également une section : c'est un isomorphisme.
\end{proposition2}

(Notons qu'un morphisme $X → \Spec(k)$ a une section
si et seulement si le morphisme $X_\red → \Spec(k)$
qui s'en déduit a une section.)

\begin{démo}
Cela résulte du fait que l'image de $k'$ dans $K$
par le morphisme composé $k ′ → K' → K$, où la seconde
flèche est la rétraction dont on suppose l'existence, est
à la fois de degré premier à $ℓ$ et de degré une puissance de $ℓ$ sur $k$.
\end{démo}

\subsection{Topologie des $ℓ′$-altérations}
\label{2-définition altL}

\subsubsection{}
On appelle \emph{topologie des $ℓ′$-altérations}
la topologie de Grothendieck sur $\alt ∕ S$, notée $\alt_{ℓ′}$,
définie par la prétopologie engendrée par
\begin{enumerate}
\item les recouvrements étales \emph{$ℓ′$-décomposés} ;
\item les morphismes propres \emph{maximalement $ℓ′$-finis} surjectifs.
\end{enumerate}

Prendre garde au fait que la seconde condition (« maximalement
$ℓ′$-fini ») porte sur les points maximaux tandis que la première
(« $ℓ′$-décomposé ») sur tous les points.

\begin{remarque2}
Les familles précédentes ne constituent
pas une prétopologie au sens de \sga{4}{II}{1.3} : la condition
de stabilité par composition n'est pas satisfaite. Les autres
conditions le sont, notamment la quarrabilité des morphismes (\ref{2-existence produit}).
\end{remarque2}

\subsubsection{}La topologie $\alt_{1′}$ est appelée \emph{topologie des
altérations}, notée simplement $\alt$.

\subsubsection{}
\label{2-stabilité par cb}
Notons pour référence ultérieure que si $(X_i → X)$ est une famille $\alt_{ℓ′}$-couvrante dans $\alt ∕ S$
et $S′ → S$ est un morphisme plat de source
un schéma localement nœthérien,
la famille $(X_i ×_S S′ → X×_S S′)$ de morphismes
dans $\alt ∕ S′$ est également $\alt_{ℓ′}$-couvrante
(cf. \ref{2-L-décomposé}, \ref{2-cb maximalement dominant par plat}).

\section{Formes standards}
Dans ce paragraphe, on fixe un schéma nœthérien $S$
et $X$ un objet de $\alt ∕ S$.
\subsection{Topologie étale}Le cas $ℓ=1$ de l'énoncé ci-dessous est un prototype
bien connu  des résultats que nous allons établir.

\begin{proposition2}
\label{2-forme standard étale}
Toute famille couvrante $(U_i → X)_{i ∈ I}$ pour la topologie étale
$ℓ′$-décomposée est dominée par une famille $\alt_{ℓ′}$-couvrante
du type
\[
(V_i → Y → X)_{i ∈ I}
\]
où $Y → X$ est fini, maximalement $ℓ′$-fini, surjectif et $(V_i → Y)_{i ∈ I}$ est un recouvrement
étale \emph{complètement} décomposé. Si $ℓ=1$, on peut supposer
que $(V_i → Y)_{i ∈ I}$ est un recouvrement par des ouverts de Zariski.
\end{proposition2}

\begin{démo}
On peut supposer l'ensemble $I$ fini (\ref{2-cd noetherien})
et $X$ intègre. Par passage à la limite, on peut
supposer de plus $X$ normal, de corps des fractions ayant
un groupe de Galois (absolu) pro-$ℓ$. (Le schéma $X$ n'est donc
pas nécessairement nœthérien.) Or, un morphisme étale $ℓ′$-décomposé
d'un tel schéma est nécessairement complètement décomposé
car le groupe de Galois des corps résiduels est également pro-$ℓ$.
Pour le complément lorsque $ℓ=1$, cf. p. ex. \cite[lemme 10.3]{Cycles@Orgogozo}.
\end{démo}

\subsubsection{}Il résulte immédiatement de la proposition
précédente que dans la définition du §\ref{2-définition altL},
on peut remplacer dans (i) la condition d'être
$ℓ′$-décomposé par celle d'être \emph{complètement}
décomposé (resp. de Zariski, si $ℓ=1$).

\subsection{Topologie des altérations}

Dans ce sous-paragraphe, on fixe un nombre premier $ℓ$.

\begin{théorème2}
\label{2-forme standard l'-altérations}
Supposons $X$ irréductible et quasi-excellent. Toute famille couvrante pour la topologie des $ℓ'$-altérations
de $X$ est dominée par une famille couvrante du type suivant :
\[
(V_i→ Y → X)_{i ∈ I},
\]
où $Y$ est un schéma intègre, $Y → X$ est propre et surjectif de degré
générique premier à $ℓ$ et $(V_i → Y)_{i ∈ I}$ est un recouvrement
pour la topologie complètement décomposée.
Si de plus $ℓ=1$, on peut supposer que $(V_i → Y)_{i ∈ I}$
est un recouvrement par des ouverts de Zariski.
\end{théorème2}

L'hypothèse de quasi-excellence est très certainement
superflue (procéder par passage à la limite).

Commençons par la démonstration du cas particulier $ℓ=1$,
bien qu'il résulte du cas général (joint
à \ref{2-forme standard étale} pour le complément).

\begin{démo}[Démonstration dans le cas où $ℓ=1$]
On peut supposer l'ensemble $I$ fini.
D'après \ref{2-forme standard étale}, il suffit de montrer
que si $(U_i→ X)_{i∈I}$ est un recouvrement de $X$ par des ouverts
de Zariski et $(X_i→ U_i)_{i∈I}$ une collection
de morphismes propres et surjectifs, il existe un morphisme
propre et surjectif $Y→X$ dans $\alt ∕ S$ et un
recouvrement par des ouverts
de Zariski $(V_i→Y)_{i ∈ I}$ tels que chaque morphisme
composé $V_i→ X$ se factorise à travers $X_i→ X$.
Choisissons pour chaque $i$ une compactification $\sur{X_i}→ X$
de $X_i → X$ ; on a $\sur{X_i}_{|U_i}=X_i$. Posons $Y=\sur{X_1}×_X \cdots
×_X \sur{X_r}$, où $I=\{1,…,r\}$, et $V_1=X_1×_X \sur{X_2}×_X\cdots
×_X \sur{X_r}$, $V_2=\sur{X_1}×_X X_2 ×_X \sur{X_3} ×_X \cdots ×_X
\sur{X_r}$, etc. Les ouverts $V_i$ recouvrent le schéma $\sur{X}$, qui
est propre et surjectif sur $X$. Par projection sur le $i$-ième facteur,
chaque $V_i$ s'envoie sur $X_i$. Quitte à appliquer le foncteur
$T↦ T_\maxdom$ (\ref{2-dom}), qui transforme un morphisme propre et surjectif (resp. un recouvrement
de Zariski) en un morphisme propre et surjectif (resp. en un recouvrement de
Zariski) (cf. \ref{2-adjoint preserve imm ouv, ps}), on obtient un recouvrement du type désiré
dans $\alt ∕ S$. Si $X$ est irréductible, on peut supposer $Y$ intègre.
\end{démo}

\begin{démo}[Démonstration dans le cas général]
Il suffit de vérifier que
si $(U_i→ X)_{i∈I}$ est un recouvrement étale complètement décomposé de $X$
et $(X_i→U_i)_{i∈I}$ une collection de morphismes propres,
surjectifs, maximalement $ℓ′$-finis, il existe une famille comme
dans l'énoncé la dominant. On peut supposer l'ensemble $I$ fini
et le schéma $X$ intègre et même normal compte-tenu
de l'hypothèse de quasi-excellence. Il est également loisible de supposer les $U_i$ connexes.
Notons qu'ils sont normaux. Observons d'autre part que l'on peut supposer les morphismes $X_i→U_i$ \emph{finis},
surjectifs, plats et de degré générique premier à $ℓ$ : chaque $X_i→X$ étant
génériquement plat ($X$ est réduit), quitte à remplacer $X$ par une
modification $X'→X$, on peut les platifier (\cite{Platification@Raynaud-Gruson},
 I 5.2.2), ce qui les rend \emph{finis}.
Quitte à ne considérer qu'une seule composante irréductible de chaque $X_i$,
de degré générique premier à $ℓ$ et à la normaliser,
on peut supposer les $X_i$ normaux connexes.
Résumons :

— les schémas $X$, $U_i$ et $X_i$ sont normaux intègres
de corps des fractions notés respectivement $K$, $K_i$
et $K_i′$ ;

— les morphismes $X_i → U_i$ sont finis de degrés
génériques premiers à $ℓ$.

Soit $L$ une clôture quasi-galoisienne sur $K$
d'une extension composée des $K_i′$,  et considérons
$X_L$ le normalisé de $X$ dans $L$. De même, considérons pour chaque $i$, le produit fibré
$U_{iL}=U_i×_X X_L$ (resp. le produit fibré
réduit normalisé $X_{iL}=(X_i×_X X_L)_{\red}^ν$). Compte tenu du
choix de $L$, pour chaque $i$ le morphisme $X_{iL}→U_{iL}$
a une section au-dessus de $\SP(L)$. Ce dernier étant fini
et $U_{iL}$ étant \emph{normal}, la section s'étend en une section $σ_i:U_{iL}→X_{iL}$.
Soit maintenant un $ℓ$-Sylow $S_ℓ$ de $\mathrm{Aut}(L/K)=G$
et notons $p$ la caractéristique du corps $K$. Si $ℓ=p$, les extensions
$K'_i/K$ sont donc étales de sorte que l'on peut supposer $L/K$ étale donc
\emph{galoisienne} ; l'extension $L^{S_ℓ}/K$ est alors de degré $(G:S_ℓ)$, premier à $ℓ$.
Si $ℓ≠p$, l'extension $L^{S_ℓ}/K$ est de degré $(G:S_ℓ)$ multiplié par une
puissance de $p$ ; c'est donc également un entier premier à $ℓ$.
Comme ci-dessus, notons $X_{L^{S_ℓ}}$ le normalisé de $X$ dans $L^{S_ℓ}$,
$U_{iL^{S_ℓ}}$ (resp. $X_{iL^{S_ℓ}}$) le produit fibré (resp. réduit normalisé)
de $U_i$ (resp. $X_i$) avec $X_{L^{S_ℓ}}$ au-dessus de $X$. Le morphisme
$X_{L^{S_ℓ}}→X$ est fini, de degré générique premier à $ℓ$.
D'après ce qui précède on a pour chaque $i$ un diagramme commutatif de 
schémas normaux :
$$
\xymatrix{
X_{iL^{S_\ell}} \ar[d]_{\text{premier à }\ell} & X_{iL} \ar[d] \ar[l] \\
U_{iL^{S_\ell}} & X_{iL} \ar[l]^{\text{puiss. de }\ell} \ar@/_1pc/[u]_{\sigma_i}
}
$$
En considérant isolément les composantes irréductibles du schéma
normal $U_{iL^{S_ℓ}}$ et, pour chacune d'entre elles, un point maximal de $X_{iL^{S_ℓ}}$ de degré
générique premier à $ℓ$ au-dessus, on montre immédiatement (cf. proposition \ref{2-descente section})
qu'il existe une section $U_{iL^{S_ℓ}}→X_{iL^{S_ℓ}}$.
Cela achève la démonstration de la proposition car
les $U_{iL^{S_ℓ}}$ forment un recouvrement pour la topologie complètement
décomposé du schéma $X_{L^{S_ℓ}}$, irréductible, de degré générique premier à $ℓ$ sur
$X$.
\end{démo}

\begin{remarques2}
\label{2-altération plus fine que étale cd}
\begin{enumerate}
\item Une famille couvrante pour la topologie des $ℓ ′$-altération le reste après changement de
base dans $\alt ∕ X$.
\item Il résulte du théorème que l'on obtient la même topologie si l'on remplace
la condition (i) de \ref{2-définition altL} par la condition
d'être un recouvrement étale \emph{complètement décomposé}.
Esquissons une preuve légèrement différente.
En passant à la limite sur les $Y$
comme dans le théorème précédent ($ℓ′$-altération de $X$)
on obtient un schéma normal dont le corps des fonctions a un
groupe de Galois absolu pro-$ℓ$. Cette propriété passe
aux corps résiduels (vérification triviale) si bien
qu'un revêtement étale $ℓ′$-décomposé est même
nécessairement \emph{complètement} décomposé.
\end{enumerate}
\end{remarques2}

Nous ferons également usage de la variante suivante
du théorème précédent.

\begin{théorème2}\label{2-forme standard pour annulation}
Supposons $X$ irréductible et quasi-excellent. Toute famille couvrante
pour la topologie des $ℓ'$-altérations de $X$
est dominée par une famille couvrante du type suivant :
\[
(W_i → V_i→ Y → X)_{i ∈ I},
\]
où tous les schémas sont irréductibles, $Y→X$ est propre birationnel,
$(V_i→Y)_{i ∈ I}$ est un recouvrement pour la topologie
complètement décomposée, et les morphismes
$W_i→V_i$ sont finis, plats, de degré premier à $ℓ$.
\end{théorème2}

Un énoncé semblable est également valable si $X$ n'est pas
irréductible : considérer le coproduit de ses composantes irréductibles.

\begin{démo}[Esquisse de démonstration]
Par platification, il suffit de montrer le résultat
d'échange suivant : si $Y → X$ est fini, plat, de degré
générique premier à $ℓ$ et $(V_i → Y)_{i ∈ I}$ est un
recouvrement étale complètement décomposé, il existe un
recouvrement étale complètement décomposé $(U_j → X)_{j ∈ J}$
et des morphismes $Z_j → W_j$, finis, plats, de degré
premier à $ℓ$ tels que la famille de morphismes
composés $(Z_j → X)$ domine celle des $(V_i → X)$.
Par passage à la limite, on peut supposer $X$ hensélien.
Par hypothèse, il existe une composante connexe $Y⁰$ de $Y$
qui est plate de degré premier à $ℓ$ sur $X$. Le schéma $Y⁰$ étant local hensélien,
il existe un indice $i$ tel que la restriction du
morphisme $V_i→Y$ à $Y⁰$ ait une section,
de sorte que le morphisme composé $Y⁰→X$ se factorise
à travers $V_i$. Cela permet de conclure.
\end{démo}

\section{Applications}

\subsection{Sorites}

\begin{proposition2}\label{2-voisinages fibre restent couvrants}
Soient $X$ un schéma nœthérien, $x$ un point de $X$,
$(X_i→X)_{i=1,\dots,n}$ un recouvrement pour la topologie des altérations
et, pour chaque indice $i$, un ouvert $X⁰_i↪X_i$ contenant la fibre
$(X_i)_x$. Il existe un voisinage ouvert de Zariski $U$ de $x$
tels que la famille $({X⁰_i}_{|U}→U)_i$ soit $\alt$-couvrante.
\end{proposition2}

Un cas particulièrement utile — et auquel on pourrait se
ramener d'après la proposition suivante  — est celui où
$X$ est local de point fermé $x$, de sorte que $U=X$.

\begin{démo}
D'après le théorème \ref{2-forme standard l'-altérations},
la famille $(X_i→X)_i$ est dominée
par une famille $(V_j↪Y→X)_j$ où $f:Y→X$ est notamment
propre et surjectif, et les $(V_j→Y)_j$ sont un recouvrement ouvert de $Y$. Soit $Y⁰_j$ l'image
inverse de $X⁰_i↪X_i$ par une factorisation $V_j→X_i$.
Par hypothèse, au-dessus du point $x$ de $X$, $Y⁰_j$ et $V_j$
coïncident, de sorte que $(Y⁰_j)_j$ est une famille d'ouverts
de $Y$ recouvrant la fibre $Y_x$. Leur réunion $Y⁰$ est un ouvert
de $Y$, contenant cette fibre, et on vérifie aussitôt que l'ouvert $U=X-f(Y-Y⁰)$
de $X$ convient.
\end{démo}

\begin{proposition2}\label{2-famille alt-couvrante et passage à la limite}
Soit $(X_α)_{α∈A}$ un système projectif filtrant de schémas
nœthériens affines à morphismes de transitions dominants. 
On suppose que la limite $X=\lim_α X_α$ est un schéma nœthérien
irréductible sur lequel un nombre premier $ℓ$ est inversible. 
Alors, toute famille $\altl$-couvrante de $X$ est dominée par l'image
inverse d'une famille $\altl$-couvrante de l'un des $X_α$.
\end{proposition2}

\begin{démo}
D'après le théorème \ref{2-forme standard l'-altérations}, joint au fait que les
diagrammes de morphismes se descendent, il suffit de démontrer le théorème
dans les cas particuliers suivants : la famille couvrante est
constituée d'un unique morphisme propre surjectif de degré générique
premier à $ℓ$ ou bien elle est finie, couvrante pour la
topologie étale complètement décomposée.
Traitons le premier cas. D'après \ega{iv}{8.10.5} (vi, xii),
tout morphisme propre et surjectif $Y→X$ se descend en un
morphisme propre et surjectif $Y_α→X_α$ ($α$ suffisamment grand).
Vérifions que l'on peut supposer $Y_α∈\Ob\,\alt/X_α$. D'après \ega{iv}{8.4.2}, on peut supposer
$X_α$ irréductible et $X→X_α$ maximalement dominant. Il en résulte,
puisque la flèche $Y_\maxdom → Y$ est un isomorphisme, que le morphisme $Y→Y_α$ se factorise à
travers $(Y_α)_\maxdom↪Y_α$. On peut donc supposer le
morphisme $Y_α → X_α$ maximalement dominant. Soit $η_α$ le
point générique de $X_α$. On doit vérifier que
la fibre générique $(Y_α)_{η_α}$ est intègre, finie de degré premier
à $ℓ$ sur $η_α$. Cela résulte du fait qu'elle l'est après
changement de base par $η→η_α$, où $η$ est le point
générique de $X$. Traitons maintenant le cas d'un recouvrement étale, complètement décomposé
$(U_i→X)_{i∈I}$ ($I$ fini). D'après \ega{iv}{17.7.8} (ii)
et \ega{iv}{8.10.5} (vi), cette famille provient d'une famille
$(U_{iα}→X_α)_{i∈I}$ couvrante pour la topologie étale d'un $X_α$,
pour $α$ suffisamment grand. D'après \ref{2-cd limite}, si pour 
tout $β≥α$ on note $f_β$ le morphisme $X_β×_{X_α} ∐_i U_{iα}→X_β$ et $π_β$
le morphisme $X→X_β$, on a $X=⋃_{β≥α} π_β^{-1}\big(\cd(f_β)\big)$.
C'est une réunion croissante d'ouverts du compact $X^\cons$ ;
pour $β$ assez grand on a donc $X=π_β^{-1}\big(\cd(f_β)\big)$.
Il résulte de \ega{iv}{8.3.11} que $\cd(f_β)=X_β$ pour $β$ assez
grand : sur $X_β$ le recouvrement étale est donc complètement
décomposé.
\end{démo}

\subsection{Caractérisation ponctuelle}

Terminons par une caractérisation de la topologie
des altérations, semblable à la caractérisation de 
la topologie étale à l'aide des anneaux
locaux strictement henséliens.

\begin{théorème2}[\cite{Cohomological@Goodwillie-Lichtenbaum}, 3.5]\label{2-points de alt}
Soient $X$ un schéma nœthérien et $Y$ un objet de $\alt/X$.
La flèche $Y→X$ est couvrante pour la topologie des altérations 
si et seulement si elle est valuativement surjective au sens suivant : 
tout morphisme maximalement dominant $\SP(V)→X$, où $V$ est un anneau de valuation à corps
résiduel algébriquement clos se relève en un morphisme
$\Spec(V)→Y$.
\end{théorème2}

\begin{remarque2}
On dispose d'un analogue du théorème \ref{2-points de alt}, où l'on remplace le point générique géométrique
$\sur{η}$ (resp. la condition que le corps des fractions est algébriquement clos) par
le quotient $\sur{η}/H$ où $H$ est un $ℓ$-Sylow de $\Aut(κ(\sur{η})/κ(η))$ (resp. la condition
que le corps des fractions est parfait de groupe de Galois absolu un
pro-$ℓ$-groupe).
\end{remarque2}

\begin{démo}
Montrons que la condition est nécessaire.
D'après le théorème précédent, on peut supposer $Y→X$ propre et
surjectif. (Le cas où $Y$ est associé à un recouvrement par des
ouverts de Zariski est trivial car $\Spec(V)$ est local.)
Notons $η_V$ le point générique de $\Spec(V)$. Le morphisme
$η_V→X$ se relève en un morphisme (non unique) $η_V→Y$ 
d'après le Nullstellensatz, car $Y→X$ est surjectif 
et $κ(η_V)$ algébriquement clos.
Il résulte du critère valuatif de propreté que le morphisme
$η_V→Y$ s'étend en un $X$-morphisme $\Spec(V)→Y$.
(Remarquons qu'il n'est pas nécessaire de supposer $\Spec(V)→X$
maximalement dominant.)

Montrons que la condition est suffisante. 
On peut supposer pour simplifier $X$ affine intègre, de point
générique noté $η$. On peut également supposer $Y$ affine, 
de sorte qu'il existe une immersion ouverte $Y↪\sur{Y}$ dans un 
schéma propre et surjectif sur $Y$. 
Choisissons enfin un point générique géométrique $\sur{η}→X$ et 
considérons l'espace de Zariski-Riemann 
$\ZR_{\sur{η}}(X)$, limite des espaces annelés $X'$, où $X'$ est un schéma intègre,
propre et surjectif, muni d'un $X$-morphisme $\sur{η}→X'$.
On peut montrer qu'il est quasi-compact (cf. \emph{op. cit.} ou 
\cite{CA@Zariski-Samuel}, chap. VI, th. 40 pour une variante)
et que si $X=\Spec(A)$, l'application qui à un anneau de valuation $A⊆V⊆κ(\sur{η})$
associe le point de $\ZR_{\sur{η}}(X)$ correspondant par le critère
valuatif de propreté est une bijection. 

Tout relèvement $r:\sur{η}→Y$ du point générique géométrique de 
$X$ induit un $X$-morphisme $\sur{η}→\sur{Y}$ donc un morphisme
continu $π_r:\ZR_{\sur{η}}(X)→\sur{Y}$, qui se factorise à travers la composante
irréductible de $\sur{Y}$ atteinte par $r$. Par hypothèse,
les ouverts $π_r^{-1}(Y)$, pour $r$ variable, recouvrent 
$\ZR_{\sur{η}}(X)$. Par quasi-compacité, il existe donc un nombre 
fini de relèvements $r₁,\dots,r_n:\sur{η}→Y$
tels que $\ZR_{\sur{η}}→X$ se factorise à travers
la réunion des $n$ \emph{ouverts} images inverses de $Y$
dans le $X$-schéma propre et surjectif $\sur{Z}=\sur{Y}×_X 
\sur{Y}×\cdots×_X \sur{Y}$ ($n$ fois).
Par définition de l'espace de Zariski-Riemann, cette factorisation 
entraîne que $\sur{Z}$ est la \emph{réunion} de ces ouverts.
Ainsi, $Y→X$ peut-être raffiné en un recouvrement par des ouverts
de Zariski d'un schéma propre et surjectif. CQFD. 
\end{démo}

\subsection{Réduction des théorèmes d'uniformisation locale au cas hensélien}

Rappelons qu'un des objectifs de ce livre est de démontrer
les théorèmes d'uniformisation locale \refext{VII}{7-uniformisation}
et \refext{0}{0-uniformisation l'} dont nous reproduisons
l'énoncé.

\begin{théorème2}[\cite{Conference-Deligne@Gabber}, 1.1]\label{2-théorème uniformisation faible}
Soient $X$ un schéma nœthérien quasi-excellent et $Z$ un
fermé rare de $X$.
Il existe une famille finie de morphismes $(X_i→X)_{i∈I}$, couvrante pour la topologie des
altérations et telle que pour tout $i∈I$ on ait :
\begin{enumerate}
\item le schéma $X_i$ est régulier et connexe ;
\item l'image inverse de $Z$ dans $X_i$ est le support
d'un diviseur à croisements normaux stricts.
\end{enumerate}
\end{théorème2}

Par convention, l'ensemble vide est considéré comme
un diviseur strictement à croisements normaux : c'est
une somme indexée par l'ensemble vide.

\begin{théorème2}[\emph{op. cit.}, 1.3]\label{2-théorème uniformisation première à l}
Soient $X$ un schéma nœthérien quasi-excellent, $Z$ un
fermé rare de $X$ et $ℓ$ un nombre premier inversible sur $X$.
Il existe une famille finie de morphismes $(X_i\ra X)_{i∈I}$ , 
couvrante pour la topologie des $ℓ'$-altérations et telle que pour
tout $i\in I$ on ait :
\begin{enumerate}
\item le schéma $X_i$ est régulier et connexe ;
\item l'image inverse de $Z$ dans $X_i$ est le support d'un diviseur à croisements normaux stricts.
\end{enumerate}
\end{théorème2}

Dans la fin de cet exposé, nous allons nous contenter de
démontrer le fait suivant.

\begin{proposition2}\label{2-réduction cas hensélien}
Si l'un des théorèmes d'uniformisation est vrai pour tout schéma
$X$ local nœthérien hensélien excellent normal (resp. pour tout
schéma $X$ local nœthérien hensélien excellent normal de dimension au plus
un entier $d$ fixé), il est vrai en général (resp. pour tout schéma nœthérien
excellent de dimension finie inférieure ou égale à $d$).
\end{proposition2}

La réduction au cas où $X$ est local nœthérien \emph{complet}
est bien plus délicate ; elle fait l'objet de l'exposé suivant.

\begin{démo}
Supposons le théorème \ref{2-théorème uniformisation faible}
(resp. \ref{2-théorème uniformisation première à l})
démontré dans le cas local hensélien excellent.
Soit $X$ un schéma nœthérien quasi-excellent et $Z$ un fermé rare.
On peut supposer $X$ normal intègre car le morphisme de normalisation est couvrant
pour la topologie des $ℓ'$-altérations et l'image inverse de $Z$ reste
rare. Fixons $x∈X$. D'après \ega{iv}{18.7.6} l'hensélisé
$X_{(x)}$ de $X$ en $x$ est excellent et $Z_{(x)}$ est un fermé
rare de $X_{(x)}$. Il existe donc une famille finie 
de diagrammes
$$
\xymatrix{
Y \ar[d] & Y_i \ar[d] \ar[l] \\
X_{(x)} & X_i \ar[l]
}
$$
dans $\alt/X_{(x)}$, où $Y$ est intègre, propre et surjectif 
(resp. et de degré générique premier à $ℓ$) sur $X$, 
$(Y_i→Y)$ est une famille couvrante pour la topologie
de Zariski (resp. complètement décomposée) et la famille $(X_i→X)$
satisfait les conditions (i-ii).
Il résulte de la démonstration de la proposition \ref{2-famille alt-couvrante et passage à la limite}
que cette famille de diagrammes s'étend en
famille du même type sur un voisinage étale complètement
décomposé $U$ de $x$ dans $X$.
Il reste à vérifier que les propriétés (i) et (ii) sont bien
conservées. Si un morphisme $T→X_{(x)}$ de type fini,
avec $T$ \emph{régulier}, est le changement de base d'un morphisme
$V→U$ de type fini où $U$ est un voisinage étale de $x$,
le schéma $V$ est régulier en les points de l'image du morphisme
$T→V$. (Un schéma local est régulier si et seulement si son hensélisé
l'est.) En particulier, $V$ est régulier en les points de la
fibre $V_x$. Le lieu régulier étant ouvert, on peut supposer
d'après \ref{2-voisinages fibre restent couvrants} — quitte à
rétrécir le voisinage $U$ de $x$ — que $V$ est régulier.
Enfin, il résulte de \ega{iv}{19.8.1 (ii)} que la propriété d'être
un diviseur à croisements normaux stricts se descend si elle est satisfaite à la limite. 
On conclut par compacité de $X$ pour la topologie étale
complètement décomposée, moins fine que les topologies
$\alt$ et $\altl$. Le cas respé de la proposition \ref{2-réduction cas hensélien}
est un corollaire immédiat de la démonstration du cas non respé.
\end{démo}

\ifx\danslelivre\undefined
\bibliography{../fichiers-configuration-latex/bibliographie}
\bibliographystyle{../fichiers-configuration-latex/smfalpha}


\author{Luc Illusie et Yves Laszlo}

\title{III. Approximation}

\setcounter{tocdepth}{1}
\maketitle
\tableofcontents

\else
\setboolean{english}{false}
\setboolean{francais}{true}
\selectlanguage{francais}

\chapter{Approximation}
\addtocontents{toc}{par Luc Illusie et Yves Laszlo}
\begin{center}
Luc Illusie et Yves Laszlo
\end{center}
\fi

\begin{abstract}On montre ici comment des méthodes
d'approximation de complexes et le théorème de Popescu
(\refext{I}{1-Artin-Popescu}) permettent de réduire la preuve
du théorème d'uniformisation des schémas quasi-excellents
(\ref{3-unif}) au cas complet.\end{abstract}

\section{Introduction} On montre ici comment ramener la preuve
du théorème d'uniformisation (\ref{3-unif}) au cas local,
noethérien complet (\ref{3-reduc}). On utilise pour cela
le théorème de Popescu (qui implique que les anneaux locaux
noethériens, henséliens et excellents vérifient la propriété
d'approximation d'Artin, cf.~\refext{I}{1-Artin-Popescu}) et des
méthodes d'approximations de complexes de longueur $2$ adaptées
de \cite{Con-deJ02} (cf. section~\ref{3-section-appro}).

L'exposé oral donné par Alban Moreau utilisait des résultats
(dus à Ofer Gabber) d'approximations de complexes plus forts
que ceux utilisés ici (\ref{3-CdJ}). Une version écrite de son
exposé a été très utile pour la rédaction de ce texte : nous
l'en remercions. Nous remercions également Fabrice Orgogozo de
nous avoir signalé que l'énoncé \cite[3.1]{Con-deJ02} suffisait
pour les applications en vue.

\section{Modèles et approximations à la Artin-Popescu}\label{3-mod}

Soit $A$  un anneau local noethérien, $\m$ son idéal maximal,
$\hat A$ son complété. On suppose $A$ excellent et hensélien. Soit
$\pi:{\s}=\Spec(\hat A)\ra S=\Spec(A)$ le morphisme canonique.  Pour
tout $n\geq 0$, on  note $$i_n:S_n\hookrightarrow {\s}$$ l'immersion
fermée définie par d'idéal $\hat \m^{n+1}=\m^{n+1}\hat A$ de
$\hat A$. Le composé $$\pi i_n:S_n\ra{\s}\ra S$$ est l'immersion
fermée $S_n\hookrightarrow S$ définie par l'idéal $\m^{n+1}$.

\begin{definition}  Soient $g:{\s}\ra T$ et $f:S\ra T$ des
morphismes de schémas et $n\in\NN$. On dira que  $f$ et $g$
sont $(n+1)$-proches si leurs restrictions $f\pi i_n$ et $gi_n$
à $S_n$ coïncident.  \end{definition}

Si $\X$ est un ${\s}$-schéma, on note $\X_n$ le $S_n$-schéma
$\X\times_{{\s}}S_n\ra S_n$.

Écrivons $\hat A$ comme limite inductive suivant
un ensemble ordonné filtrant $E$ de $A$-algèbres
de type fini $A_\alpha$. On a  des diagrammes commutatifs
\begin{equation}\label{3-salpha}\xymatrix{&S_\alpha=\Spec(A_\alpha)\ar[d]^{t_\alpha}\\{\s}\ar[ru]^{s_\alpha}\ar[r]^\pi&S}\end{equation}
avec  $t_\alpha$ de type fini et un isomorphisme  ${\s}=\varprojlim
S_\alpha$ \ega{IV$_3$}{8.2.3}.

\begin{definition}\label{3-def-mod} Soit $\X$ un
${\s}$-schéma de type fini et $h:\X\ra\Y$ un morphisme
de $\s$-schémas de type fini.  \begin{itemize}\item[]
Un modèle de $\X$  sur $S_\alpha$ est un diagramme cartésien
$$\xymatrix{\X\ar[r]\ar[d]_f\ar@{}[rd]|\square&X_\alpha\ar[d]^{f_\alpha}\\{\s}\ar[r]&S_\alpha\\}$$
où $X_\alpha$ est de type fini sur $S_\alpha$.  \item[] Un
modèle de $h$ sur $X_\alpha$ est un $S_\alpha$-morphisme
$h_\alpha:X_\alpha\ra Y_\alpha$ muni d'un isomorphisme
$h\iso(h_\alpha)_{{\s}}$.\end{itemize} \end{definition}

Des modèles de $\X$ sur $S_\alpha$  existent pourvu que
$\alpha$ soit assez grand \ega{IV$_3$}{8.8.3}. De plus, si
$X_\alpha,X_\beta$ sont des modèles de $\X$ sur $S_\alpha,S_\beta$,
il existe $\gamma\geq \alpha,\beta$ et un $S_{\gamma}$-isomorphisme
$$X_\alpha\times_{S_\alpha}S_\gamma\xrightarrow{\sim}
X_\beta\times_{S_\beta}S_\gamma$$ (\textit{loc. cit.}). De même,
des modèles $h_\alpha$ de $h:\X\ra\Y$ sur $S_\alpha$  existent
pourvu que $\alpha$ soit assez grand et les images inverses de tels
modèles $h_\alpha,h_\beta$ sur $S_\gamma$ sont $S_\gamma$-isomorphes
pour $\gamma\geq\alpha,\beta$ assez grand.

Si $T$ est un $S$-schéma et $B$ une $A$-algèbre, on note
$T(B)=\Hom_S(\Spec(B),T)$ l'ensemble des $S$-points de $T$ à valeurs
dans $B$. D'après le théorème de Popescu \cite[1.3]{Pop86},
comme $A$ est excellent et hensélien, il vérifie la propriété
d'approximation d'Artin, cf.  \refext{I}{1-Artin-Popescu}. Donc,
comme  $S_\alpha\ra S$  est de type fini,  $S_\alpha(A)$ est dense
dans $S_\alpha(\hat A)$ (pour la topologie $\m$-adique). Il existe
donc, pour tout $n\geq 0$ une section $u:S\ra S_\alpha$ de $t_\alpha$
qui est $n$-proche de $s_\alpha:\s\ra S_\alpha$. On définit alors
$X_u$ par le diagramme cartésien

\begin{equation}\label{3-Xu}
\xymatrix{X_u\ar[r]\ar[d]_{f_u}\ar@{}[rd]|\square&X_\alpha\ar[d]^{f_\alpha}\\S\ar[r]^u&S_\alpha\\}
\end{equation}

Comme $u$ est $n$-proche de $s_\alpha$, on a par définition
l'égalité $$u\pi i_n=s_\alpha i_n$$ de sorte la restriction de
$X_u\ra S$ à $S_n$ s'identifie à $\X_n\ra S_n$, autrement dit on
a un carré cartésien \begin{equation}\label{3-Xn}
    \xymatrix{\X_n\ar[r]\ar[d]_{f_n}\ar@{}[rd]|\square&X_u\ar[d]^{f_u}\\S_n\ar[r]^{i_n}&S\\}
\end{equation}

De même, si $\X,{{\Y}}$ sont de type fini sur ${\s}$ et $h_\alpha$
est un modèle de $h\in\Hom_{{\s}}(\X,{{\Y}})$ sur $S_\alpha$,
l'image inverse $h_u:X_u\ra Y_u$ est un $S$-morphisme induisant la
restriction $h_n:\X_n\ra {{\Y}}_n$ de $h$ au dessus de $S_n$.

\section{Approximations et topologie des altérations} Commençons
par un rappel (cf. \refext{II}{2-topologies}) sur la topologie
des altérations. Soit $T$ un schéma noethérien.  La catégorie
${\pspf}/T$ est la sous-catégorie pleine de la catégorie des
$T$-schémas dont les objets sont les $T$-schémas réduits de
type fini $X$, dont tout point maximal s'envoie sur un point
maximal de $T$ avec extension résiduelle finie. Notons que les
morphismes de $\pspf/T$ envoient point maximal sur point maximal.
  On définit deux topologies sur $\pspf/T$.

\begin{itemize}
  \item[(i)] La \textit{topologie des altérations} est la
  moins fine pour laquelle les familles suivantes sont couvrantes
  \begin{itemize} \item[(a)]les recouvrements ouverts de Zariski ;
  \item[(b)] les morphismes propres et surjectifs .  \end{itemize}
Une famille couvrante pour la topologie des altérations sera dite
\pspf-couvrante.
  \item[(ii)] Soit $\ell$ un nombre premier. La \textit{topologie
  des $\ell'$-altérations} sur ${\pspf}/T$ est la moins fine pour
  laquelle les familles suivantes sont couvrantes

  \begin{itemize} \item[(a)]les recouvrements étales de Nisnevich ;
  \item[(b)] les morphismes  propres surjectifs $\X'\ra \X$ tels que
  pour tout point maximal $\eta$ de $X$, il existe un point maximal
  $\eta'$ de $X'$ au dessus de $\eta$ avec  $\ell$ ne divisant pas
  $\deg(k(\eta')/k(\eta))$.  \end{itemize} Une famille couvrante
  pour la topologie des $\ell'$-altérations
sera dite \altl-couvrante.  \end{itemize}

Pour tout $T$-schéma $X$ dominant, on note $X_\r$ le sous-schéma
fermé réduit de $X$  réunion des composantes irréductibles qui
dominent une composante irréductible de $T$.

\begin{proposition}\label{3-alt-mod}On reprend les notations de
\ref{3-mod} : soit $A$  un anneau local noethérien, $\m$ son
idéal maximal, $\hat A$ son complété. On suppose $A$ excellent
et hensélien. Soit $\pi:{\s}=\Spec(\hat A)\ra S=\Spec(A)$
le morphisme canonique.Soit $\X\ra {\s}$  un objet non vide de
${\pspf/{\s}}$. Supposons de plus $S$ intègre.
 \begin{itemize} \item[(i)]Alors, il existe $\alpha_0\in E,
 n_0\in\mathbf{N}$ tel que  pour tout $\alpha\geq\alpha_0$, tout
 entier $n\geq n_0$, toute section $u:S\ra S_\alpha$ de $t_\alpha$
 qui est $n$-proche de $s_\alpha:{\s}\ra S_\alpha$, tout modèle
 (cf. \ref{3-def-mod}) $X_\alpha$ de $\X$ sur $S_\alpha$, $X_u$ est
 à fibre générique finie et le morphisme composé $(X_u)_\r\ra
 X_u\ra S$ est un objet non vide de ${\pspf/S}$.
   \item[(ii)] Supposons $X \ra{\s}$ \pspf-couvrant. Alors,
   il existe $\alpha_0\in E, n_0\in\mathbf{N}$ tel que  pour
   tout $\alpha\geq\alpha_0$, tout entier $n\geq n_0$, toute
   section $u:S\ra S_\alpha$ de $t_\alpha$ qui est $n$-proche de
   $s_\alpha:{\s}\ra S_\alpha$, tout modèle $X_\alpha$ de $\X$
   sur $S_\alpha$, le morphisme composé $(X_u)_\r\ra X_u\ra
   S$ est $\pspf$-couvrant.  \item[(iii)] Supposons $X \ra{\s}$
   \altl-couvrant. Alors, il existe $\alpha_0\in E, n_0\in\mathbf{N}$
   tel que  pour tout $\alpha\geq\alpha_0$, tout entier $n\geq n_0$,
   toute section $u:S\ra S_\alpha$ de $t_\alpha$ qui est $n$-proche
   de $s_\alpha:{\s}\ra S_\alpha$, tout modèle $X_\alpha$ de $\X$
   sur $S_\alpha$, le morphisme composé $(X_u)_\r\ra X_u\ra S$
   est $\altl$-couvrant.
 \end{itemize}
\end{proposition}

\begin{proof} Observons d'abord que, $S$ étant
hensélien et excellent, $\s$ est intègre,
cf. \refext{I}{1-cor-normal-versus-complet}.

Prouvons (i). Comme $\X\ra\s$ est génériquement fini, il existe
$a\in\hat A-\{0\}$ tel que $\X$ soit fini, surjectif et libre
de rang $d>0$ au dessus de l'ouvert non vide  $\X-V(a)$. On peut
choisir $\alpha_0$ assez grand de sorte que \begin{itemize}
  \item $a$ provienne de $a_\alpha\in A_\alpha-\{0\}$ pour
  $\alpha\geq\alpha_0$ ; \item  $X_\alpha\ra S_\alpha$ soit fini,
  surjectif (\ega{IV$_3$}{8.10.5}) et libre de rang $d$ sur
  $X_\alpha-V(a_\alpha)$ (utiliser \ega{IV$_3$}{8.5.2}).
\end{itemize}

Choisissons alors un entier $n$ tel que $a\not\in\hat\m^{n+1}$. Pour
tout $\alpha\geq\alpha_0, m\geq n$, toute  section $u$ qui est
$m$-proche de $t_\alpha,$, on a $$u^*(a_\alpha)\not\in\m^{n+1}$$
et donc $u^*(a_\alpha)$ est non nul. Ceci assure que  $X_u$ est
fini, surjectif et libre de rang $d$ au dessus de l'ouvert non vide
$S-V(u^*(a_\alpha))$ image réciproque de de $S_\alpha-V(a_\alpha)$
par $u$. Le premier point en découle.

Prouvons (iii) [La preuve de (ii) est en tout point similaire]. On
suppose donc que  $\X\ra {\s}$ est $\altl$-couvrant. On sait
(\refext{II}{2-forme standard l'-altérations}) que $\X\ra {\s}$
est dominé dans $\pspf/\s$ par un recouvrement standard $${{\Y}}\ra
\X'\ra {\s}$$ avec

 \begin{itemize}
   \item ${{\Y}}\ra \X'$  Nisnevich couvrant \item $\X'\ra {\s}$
   propre et surjectif dont la restriction à chaque composante
   irréductible est dominante et génériquement, le degré
   générique de l'une d'elles étant premier à $\ell$.
 \end{itemize}

 Quitte à remplacer le schéma réduit $\X'$ par une  composante
 convenable et $\Y$ par le \altl-recouvrement Nisnevich induit,
 on peut supposer $\X'$ \emph{intègre  } de degré générique
 $\deg(\X'/\s)=\delta$ premier à $\ell$.

Soit $\eta$ le point générique de $S$. La construction
$X\mapsto X_\r$ est fonctorielle pour la sous-catégorie pleine des
$S$-schémas $X$ à fibre générique  finie. Or, d'après (i), pour
des choix de modèles et de section $u$ de $t_\alpha$ convenables,
on sait que $Y_u,X'_u$ et $X_u$ sont à fibre générique finie. On
a donc une factorisation
 $$\xymatrix{(Y_u)_\r\ar[r]\ar[d]&(X_u)_\r\ar[d]\\(X'_u)_\r\ar[r]&S}$$
 Or, toujours d'après (i), on peut en outre supposer que
 $(Y_u)_\r,(X'_u)_\r$ et $(X_u)_r$ sont des objets de $\pspf/S$.
 Pour conclure que $(X_u)_\r\ra S$ est $\altl$-couvrant, il
 suffit de prouver que pour $u$ convenable $(Y_u)_\r\ra S$ est
 $\altl$-couvrant.

 Tenant compte des propriétés de permanence usuelles des modèles
 \ega{IV$_3$}{8.8.3 et 8.10.5}, la preuve de (i) assure que pour
 des modèles et $u$ convenables le $X'_u\rightarrow S$ est propre
 et surjectif et que sa fibre générique est de degré premier
 à $\ell$. Ceci assure que la restriction de $X'_u\rightarrow S$
 à au moins une des composantes réduites de $X'_u$ dominant $S$
 est de degré premier à $\ell$. 
 Ainsi,  $(X'_u)_\r\ra S$ est bien \altl-couvrant.

 La propriété d'être un recouvrement Nisnevich (resp. propre et
 surjectif) étant stable par changement de base, reste à prouver
 le lemme suivant.

\begin{lemme2}\label{3-stabmod}Il existe $\alpha_0$ tel que
pour tout $\alpha\geq\alpha_0$, tout modèle ${ Y}_\alpha\ra {
X}'_\alpha\ra  S_\alpha$ de ${{\Y}}\ra \X'\ra {\s}$ vérifie
$\Y_\alpha\ra \X'_\alpha$ est Nisnevich couvrant.  \end{lemme2}

\begin{proof}Dire que ${{\Y}}\ra \X'$ est Nisnevich couvrant, c'est
dire qu'il est lisse, quasi-fini et qu'on a une stratification
$$\varnothing=\X'_0\subset \X'_1\cdots \subset \X'_n=\X'$$ avec
$\X'_i$ fermé de $\X'$ et ${{\Y}}/\X'$ a une section au dessus
de $\X'_{i+1}-\X'_i$. La conclusion découle immédiatement de
cette remarque et des propriétés de permanence usuelles des
modèles \ega{IV$_3$}{8.8.3 et 8.10.5} et \ega{IV$_4$}{17.7.8}.
\end{proof}\end{proof}

Le but de ce qui suit est d'améliorer les résultats topologiques de
la proposition~\ref{3-alt-mod} en montrant que des épaississements
convenables des cônes normaux des fibres spéciales de $\X$
(resp. $X_u$) dans $\X$ (resp. $X_u$) sont isomorphes. Ceci
permettra de prouver des énoncés de stabilité de propriétés
dans le passage de $\X$ à $X_u$, en l'occurrence la dimension et
la régularité (corollaire~\ref{3-perm}).

\section{Gradués supérieurs et approximations de
complexes}\label{3-section-appro}  Soient $I$ un idéal d'un
topos annelé $(\mathscr{X},\OO)$, $\F$ un ${\OO}$-module
de $\mathscr{X}$ et $a$ un entier $\geq 1$. On pose
$I^n=\OO$ si $n\leq 0$. On définit le module $\ZZ$-gradué
$$\gr_a(\F)=\mathop{\oplus}_{n\in\ZZ}I^n\F/I^{n+a}\F$$ qui est
donc la somme $$\gr_a(\F)=\F/{I}\F\oplus\cdots\oplus\F/{I}^aF\oplus
{I}\F/{I}^{a+1}\F\oplus {I}^2\F/{I}^{a+2}\F\oplus\cdots$$ concentrée
en degrés $\geq -(a-1)$.  C'est un ${\OO}/I^a$-module ; de plus,
le produit $$I^n\otimes I^m\ra I^{n+m}$$ induit une structure de
${\OO}/I^a$-algèbre $\ZZ$-graduée sur $\gr_a({\OO})$ et $\gr_a(\F)$
est un $\gr_a(\OO)$-module $\ZZ$-gradué.

On s'intéresse ici au cas où $\mathscr{X}$ est le topos de Zariski
d'un $S$-schéma $X$ annelé par son faisceau structural $\OO$
et $I=\hat m{\OO}$.

\begin{remarque}Le morphisme surjectif tautologique
$\gr_a(\OO)\ra\gr_1(\OO)$ a pour noyau $J=I\gr_a(\OO)$. On a donc
$J^a=0$ (puisque $J$ est un $\OO_{X_{a-1}}$-module) de sorte que
$C_a(\X)=\Spec(\gr_a(\OO))$ est un épaississement d'ordre $a-1$
du cône normal $\Spec(\gr_1(\OO))$.\end{remarque}

\begin{definition} Soient $\X,\Y$ des $S$-schémas (resp. des
$\hat S$-schémas). Un $a$-isomorphisme $\X\iso_a {{\Y}}$
est la donnée d'un $S$-isomorphisme $\phi:X_{a-1}\iso
Y_{a-1}$ et d'un isomorphisme de $\gr_a(A)$-algèbres graduées
$\phi^{-1}\gr_a(\OO_{{\Y}})\iso\gr_a(\OO_\X)$. On dit dans ce cas
que $\X,{{\Y}}$ sont $a$-proches. \end{definition}

On identifiera alors leurs fibres spéciales $X_0,Y_0$ grâce à
l'isomorphisme $X_{a-1}\iso Y_{a-1}$.

\subsection{} On adapte ici le théorème 3.2 de \cite{Con-deJ02}
(et le lemme clef 3.1 de \textit{loc. cit.}). Commençons par
une définition. Soient $B$ un anneau noethérien et $I$ un idéal
de $B$.

\begin{definition}\label{3-AR} Soit $f:M\ra N$ un morphisme de
$B$-modules de type fini. Un entier $c\geq 0$ est une constante
d'Artin-Rees de $f$ si pour tout $n\geq c$ on a $${{I}}^nN\cap
\Im(f)\subset {{I}}^{n-c}\Im(f).$$\end{definition}

Le lemme d'Artin-Rees assure l'existence d'une constante
d'Artin-Rees.

\begin{proposition}\label{3-CdJ} Soient
$(L^\bullet,\d_L^\bullet),(M^\bullet,\d_M^\bullet)$ des
complexes de $B$-modules libres de type fini concentrés
en degré $[-2,0]$ avec $L^i=M^i$ pour tout $i$. Soit $c$
une constante d'Artin-Rees pour $\d_L^{-2}$ et $\d_L^{-1}$
et $n$ un entier $\geq c$. Supposons $H^{-1}(L^\bullet)=0$ et
$$\d_L^\bullet=\d_M^\bullet\mod I^{n+1}.$$ Alors : \begin{itemize}
\item[(i)] $c$ est une constante d'Artin-Rees pour $\d_M^{-1}$
; \item[(ii)] $H^{-1}(M^\bullet)=0$ ; \item[(iii)] L'identité
de $L^0=M^0$induit un isomorphisme de $\gr_{n+1-c}(B)$-modules
$$\gr_{n+1-c}(H^{0}(L^\bullet))\iso\gr_{n+1-c}(H^{0}(M^\bullet)) ;$$
\item[(iv)] De plus, si $L^0=M^0=B$, l'isomorphisme précédent
est un isomorphisme de $\gr_{n+1-c}(B)$-algèbres, autrement
dit les algèbres $H^{0}(L^\bullet)$ et $H^{0}(M^\bullet)$ sont
$(n+1-c)$-isomorphes.  \end{itemize} \end{proposition}

\begin{proof} Les deux premiers points sont prouvés dans le lemme
3.1 de \textit{loc. cit.}. Le dernier est trivial. Reste le point
(iii).

Pour $n=c$, c'est le théorème 3.2  de \textit{loc. cit.}
dont on ne fait qu'adapter la preuve dans le
cas $n>c$. Soit $m\in\ZZ$. On écrit $\d_L,\d_M$
pour $\d_L^{-1},\d^{-1}_M$. Pour $\delta=\d_L,\d_M$, on a
$$\gr_{n+1-c}^m(\coker(\delta))={I}^mL^0/({I}^{m+n+1-c}L^0+{I}^mL^0\cap\Im(\delta))$$
de sorte qu'il s'agit de montrer l'égalité
$${I}^{m+n+1-c}L^0+{I}^mL^0\cap\Im(\d_L)={I}^{m+n+1-c}L^0+{I}^mL^0\cap\Im(\d_M)$$
pour tout $m\in\ZZ$. Soit $x\in L^{-1}$ tel que $\d_L(x)\in{I}^mL^0$.

Supposons $m\leq c$.  Comme $$\d_L(x)-\d_M(x)\in{I}^{n+1}L^0\text{
et }m\leq c\leq n,$$ on a $\d_L(x)-\d_M(x)\in{I}^{m}L^0$ de sorte
que $$\d_M(x)=\d_L(x)+\d_M(x)-\d_L(x)\in{I}^mL^0\cap\Im(d_M).$$
Comme $n+1\geq m+n+1-c$, on a également
$$\d_L(x)-\d_M(x)\in{I}^{n+1}L^0\subset{I}^{m+n+1-c}L^0$$
de sorte que $$\d_L(x)=\d_L(x)-\d_M(x)+\d_M(x)\in
{I}^{m+n+1-c}L^0+{I}^mL^0\cap\Im(\d_M)$$ et donc
$${I}^{m+n+1-c}L^0+{I}^mL^0\cap\Im(\d_L)\subset{I}^{m+n+1-c}L^0+{I}^mL^0\cap\Im(\d_M).$$
Par symétrie des rôles de $\d_L$ et $\d_M$, on a l'égalité
cherchée dans ce cas.

Si $m>c$, le calcul est analogue. On a (\ref{3-AR})
$${I}^mL^0\cap\Im(\d_L)\subset{I}^{m-c}\d_L(L^{-1})$$ de sorte
que $$\d_L(x)=\d_L(x')\text{ avec }x'\in{I}^{m-c}L^{-1}.$$
Comme $\d_L-\d_M=0\mod{I}^{n+1}$, la matrice de $\d_L-\d_M$
est à coefficients dans ${I}^{n+1}$ de sorte que
$$\d_L-\d_M\in{I}^{n+1}\Hom_B(L^{-2},L^{-1}).$$ On a donc
$$\d_L(x')-\d_M(x')\in{I}^{n+1}{I}^{m-c}L^{0}={I}^{n+1+m-c}L^{0}.$$
Comme
$$\d_M(x')=\d_L(x')+\d_M(x')-\d_L(x')=\d_L(x)+\d_M(x')-\d_L(x'),$$
on a d'une part $$\d_M(x')\in({I}^m
L^0+{I}^{n+1+m-c}L_0)\cap\Im(d_M)\stackrel{n\geq
c}{\subset}{I}^m L^0\cap\Im(d_M),$$ et, d'autre part,
$$\d_L(x)=\d_M(x')-\d_L(x')+\d_L(x')+\in{I}^{m+n+1-c}L^0+{I}^mL^0\cap\Im(\d_M).$$
On conclut comme plus haut  par symétrie.\end{proof}

\section{Modèles et $a$-isomorphismes}
\begin{theoreme}[Approximation]\label{3-theo} Soit $A$  un
anneau local noethérien, $\m$ son idéal maximal, $\hat A$
son complété. On suppose $A$ excellent et hensélien. Soit
$\pi:{\s}=\Spec(\hat A)\ra S=\Spec(A)$ le morphisme canonique. Soit
$\X$ de type fini sur $S$. On se donne de plus $\alpha_0\in E$
et  un modèle (cf. \ref{3-def-mod}) $X_{\alpha_0}$ de $\X$
sur $S_{\alpha_0}$.  Pour tout $\alpha\geq \alpha_0$ on note
$X_\alpha=X_\alpha\times_{S_{\alpha_0}}S_\alpha$ le modèle de
$\X$ sur $S_\alpha$ déduit par changement de base.  Il existe
$\alpha_1\geq\alpha_0$ et des entiers $n_0\geq c>0$ tels que
pour tout $n\geq n_0,\alpha\geq\alpha_1$ et toute section $u$ de
$t_\alpha$ qui est $(n+1)$-proche de $s_\alpha$, il existe un unique
$(n+1-c)$-isomorphisme $\X\iso_{n-c}X_u$ au dessus de l'isomorphisme
$X_{n-c}\ra (X_u)_{n-c}$ déduit de (\ref{3-Xn}).\end{theoreme}

\begin{definition}\label{3-rem} Dans les conditions précédentes,
on dit que $(X_\alpha,\alpha,u)$ (ou, si aucune conclusion n'est
à craindre, $X_u$) est une \emph{approximation} de $\X$ sur $S$
(à l'ordre $n-c$).  \end{definition}

L'assertion « Il existe $\alpha_0$, un entier  $n_0$ tels
que pour tout $n\geq n_0,\alpha\geq\alpha_0$ et toute section
$u$ de $t_\alpha$ qui est $(n+1)$-proche de $s_\alpha$, $X_u$
vérifie la propriété $P$ » pourra parfois être condensée
en « Toute approximation $X_u$ assez fine de $\X$ vérifie la
propriété $P$ ». On emploiera une terminologie analogue pour
les approximations de ${\s}$-morphismes.

\begin{proof}Deux $(n+1-c)$-isomorphismes diffèrent par un
automorphisme $$\iota:\gr_{n+1-c}(\OO_\X)\iso\gr_{n+1-c}(\OO_\X)$$
de $\OO_{X_{n-c}}$-algèbres graduées. Il est en particulier
$\OO_S$-linéaire. Comme $\gr_{n+1-c}(\OO_\X)$ est engendré sur
$\gr_a(\OO_S)$ par $\OO_{X_{n-c}}$, l'automorphisme $\iota$ est
l'identité. D'où l'unicité.

On peut donc supposer $\X$ affine. Comme $\X$ est de type fini sur
$\s$, $X$ se plonge dans l'espace affine $$\Aff^m_{\s}=\Spec(\hat
A[t])$$ de coordonnées $t=(t_1,\cdots,t_m)$ comme le sous-schéma
fermé d'idéal $$J=\langle\tilde P_1,\cdots,\tilde P_N\rangle$$ où
$\tilde P_i\in B=\hat A[t]$. Choisissons une résolution partielle
du $B$-module $C=B/J$ par des $B$-modules libres de type fini

 \begin{equation}\label{3-1}
    B^a\xrightarrow{\tilde R}B^b\xrightarrow{\tilde P=(\tilde
    P_i)}B\ra C\ra 0
\end{equation} où $\tilde R$ est une matrice à coefficients
dans $B$.

Pour $\alpha_0$ assez grand,  $\tilde P$  et $\tilde R$ proviennent
de matrices $P_{\alpha_0},R_{\alpha_0}$ à coefficients dans
$$ B_{\alpha_0}=A_{\alpha_0}[t]\text{ telles que }PR=0$$ de
sorte que le fermé $F$ de $\Aff^m_{A_{\alpha_0}}$ d'équations
$P_{\alpha_0,1}=\cdots=P_{\alpha_0,N}=0$ est un modèle de $\X$
sur $S_{\alpha_0}$. Comme rappelé dans la section~\ref{3-mod},
quitte à changer $\alpha_0$ en un indice plus grand, on peut
supposer qu'on a $F=X_{\alpha_0}$. Pour $\alpha\geq\alpha_0$, note
$P_\alpha,R_\alpha$ les matrices à coefficients dans $B_\alpha$
déduites de $P_{\alpha_0},\R_{\alpha_0}$ par le morphisme
$$B_{\alpha_0}=A_{\alpha_0}[t]\ra B_\alpha=A_\alpha[t].$$
Pour tout $\alpha\geq \alpha_0$, les matrices à coefficients
dans $B$ déduites de $P_{\alpha},\R_{\alpha}$ par le morphisme
$$B_{\alpha_0}=A_{\alpha_0}[t]\ra B=\hat A[t]$$ sont les mêmes :
on les note $P,R$.

On s'est ramené, pour $\alpha\geq \alpha_0$,
au cas où $$X_\alpha=\Spec(C_\alpha)\text{ avec
}C_\alpha=B_\alpha/(P_\alpha).$$ On dispose donc d'une part
d'un complexe (en degrés [-2,0]) de $B_\alpha$-modules libres
\begin{equation*}
  L_\alpha=(
  B_\alpha^a\xrightarrow{R_\alpha}B_\alpha^b\xrightarrow{P_\alpha=(P_{i,\alpha})}B_\alpha)
\end{equation*} avec $H^0(L_\alpha)=C_\alpha$. Le complexe de
$B$-modules libres de rang fini \begin{equation*}
  L=B\otimes_{B_\alpha}L_\alpha=(
  B^a\xrightarrow{R}B^b\xrightarrow{P=1\otimes P_\alpha}B)
\end{equation*}  est acyclique en degré $-1$ par construction.

\begin{remarque}\textit{A priori}, $L_\alpha$ n'a pas de
raison d'être acyclique en degré $-1$, même pour $\alpha$
grand. \end{remarque}

D'autre part, la section $u$ de $t_\alpha$ est définie par un
morphisme de $A$-algèbres $$u^*:A_\alpha\ra A$$ de sorte que
$$u^*\mod\m^{n+1}=s^*_\alpha\mod\hat\m^{n+1},$$
 où $s_\alpha^*: A_\alpha\ra\hat A$ est défini par $s_\alpha:\hat
 S\ra S_\alpha$ (\ref{3-salpha}). Par action sur les coefficients
 des polynômes, on obtient un morphisme d'anneau
$$\bar u:B_\alpha=A_\alpha[t]\ra A[t]\ra\hat A[t]=B$$ d'où un
complexe \begin{equation*}
  M=(  B^a\xrightarrow{\bar u(R)}B^b\xrightarrow{\bar u(P)}B)
\end{equation*}  Par construction, on a $$L/\m^{n+1}L=M/\m^{n+1}M.$$
On choisit alors une constante d'Artin-Rees $c$ pour
$B^{b}\xrightarrow{P}B$ et on invoque la proposition \ref{3-CdJ}
pour conclure.  \end{proof}

\begin{corollaire}\label{3-perm}Soient $\X,{{\Y}}$ des
$S$-schémas noethériens qui sont $a$-proches. Soit $x\in X_0=Y_0$.
\begin{itemize} \item[(i)] Si $a\geq 1$, les dimensions de $\X$ et
${{\Y}}$ en $x$ sont les mêmes.  \item[(ii)] Si $a\geq 2$ et $\X$
régulier en $x$, alors ${{\Y}}$ régulier en $x$.  \item[(iii)]
Supposons $\X\ra{\s}$ de type fini et $\X$ régulier. Alors,
il existe $\alpha_0\in E, n_0\in\mathbf{N}$ tels que  pour tout
$\alpha\geq\alpha_0$, tout entier $n\geq n_0$, toute section $u$
de $t_\alpha$ qui est $n$-proche de $s_\alpha:{\s}\ra S_\alpha$,
tout modèle $X_\alpha$ de $\X$ sur $S_\alpha$, le schéma $X_u$
soit régulier dans un voisinage ouvert de la fibre spéciale.
\end{itemize} \end{corollaire}

\begin{proof}

Par hypothèses, les cônes normaux de $X_0,Y_0$ dans $\X,{{\Y}}$
sont $S$-isomorphes. Comme la dimension de $\X$ en $x$ est égale
à celle de son cône normal \cite[15.9]{CRT@Matsumura}, le premier
point en découle.

Supposons maintenant que $\X,{{\Y}}$ soient $2$-proches. D'après
(i), on sait que $\X$ et ${{\Y}}$ ont même dimension en $x$. Comme
$\X,{{\Y}}$ sont $2$-proches, $X_1$ et $Y_1$ sont isomorphes. Puisque
l'espace tangent de Zariski à $\X$ en un point de $X_0$ ne dépend
que de son second voisinage infinitésimal $X_1$,  les $k(x)$-espaces
vectoriels cotangents de Zariski en $x$ à $\X$ et ${{\Y}}$ sont
isomorphes, d'où ii).

Pour le dernier point, il suffit d'invoquer les deux premiers et
le théorème~\ref{3-theo} pour conclure qu'une approximation assez
fine est régulière au voisinage de la fibre spéciale. Comme $X_u$
est excellent (puisque de type fini sur $S$ excellent), son lieu
régulier $R$ est ouvert de sorte que $R$ est un voisinage ouvert
régulier de la fibre spéciale.\end{proof}

 \begin{remarque}\label{3-remgabber} O. Gabber sait généraliser
 la proposition~\ref{3-CdJ} au cas où les complexes envisagés
 sont seulement de type fini sur un anneau noethérien pour
 obtenir  les proximités de la cohomologie également en
 degré $-2$ (et pas seulement en degré $0,-1$). Il peut plus
 précisément montrer des énoncés de proximité pour les images,
 noyaux des différentielles\footnotemark. Gabber en déduit de
 nombreux énoncés de permanence par approximation analogues au
 corollaire~\ref{3-perm} . Notamment, si $X,Y$ sont $a$-proches pour
 $a$ assez grand, alors $X$ réduit (resp. normal) le long de $X_0$ entraîne $Y$
 réduit (resp. normal) le long de $Y_0$. Cependant, plusieurs questions naturelles
 restent suspens comme par exemple la permanence des propriétés
 $S_n,R_n$. \end{remarque}
\footnotetext{La preuve de cette généralisation a été exposée
par A. Moreau lors du séminaire sur les travaux de Gabber.}
\section{Réduction au cas local noethérien complet}Rappelons
l'énoncé du théorème d'uniformisation

\begin{theoreme}[Uniformisation]\label{3-unif} Soient $T$ un schéma
noethérien quasi-excellent et $Z$ un fermé rare de $T$. Soit $\ell$
un nombre premier inversible sur $T$.  Il existe une famille finie de
morphismes $(X_i\ra T)_{i\in I}$  et telle que pour tout $i\in I$ on
ait \begin{itemize} \item[(i)]La famille finie de morphismes $(X_i\ra
T)_{i\in I}$   est {\pspf}-couvrante (resp. \altl-couvrante) ;
  \item[(ii)] $X_i$ est régulier et connexe ; \item[(iii)]l'image
  inverse de $Z$ dans $X_i$ est vide ou le support d'un diviseur
  à croisements normaux stricts;
\end{itemize} \end{theoreme}

Nous allons montrer l'énoncé de réduction suivant.

 \begin{proposition}\label{3-reduc} Si (\ref{3-unif}) est vrai pour
 tout $T$ noethérien, local, complet, alors (\ref{3-unif}) est vrai.
  \end{proposition}
\begin{proof} On peut d'abord supposer  $T$ local, excellent et
hensélien (rappelons (\ref{1-cor-normal-versus-complet}) qu'un
schéma local, hensélien et quasi-excellent est excellent).

En effet, supposons le théorème prouvé dans ce cas. D'après
\ega{IV$_4$}{18.7.6}, l'hensélisé $T_{(t)}$ de $T$ en $t\in T$
est excellent et $Z_{(t)}$ est un fermé rare de $T_{(t)}$. On peut
alors trouver une famille finie $(f_i:X_i\ra T_{(t)})_{i\in I}$
vérifiant les 3 propriétés (i)-(iii) du théorème. La famille
$f_i$ est dominée par un recouvrement standard $f:Y\ra T'\ra
T_{(t)}$ avec $Y\ra T'$ Zariski (resp. Nisnevich) couvrant et
$T'\ra T$ propre et surjectif (resp. propre et surjectif tels que
pour tout point maximal $\eta$ de $T$, il existe un point maximal
$\eta'$ de $T'$ au dessus de $\eta$ avec  $\ell$ ne divisant pas
$\deg(k(\eta')/k(\eta))$). D'après \ega{IV$_3$}{8.2.2}, les $f_i$
et $f$ proviennent de morphismes $\tilde f_i:\tilde X_i\ra \tilde
T$ où  $\tilde T\ra T$ est un voisinage Nisnevich étale de $T$ en
$t$. La topologie de Nisnevich étant moins fine que les topologies
$\pspf$ et $\altl$ (\refext{II}{2-altération plus fine que étale
cd} (ii)), les propriétés usuelles de passage à la limite
(\ref{3-stabmod} et\ega{IV$_3$}{8.4.3}, \ega{IV$_4$}{19.8.1})
permettent de conclure comme plus haut.

Supposons donc $T$, schéma local, hensélien, excellent.

Quitte à remplacer $S$ par la somme disjointe de ses composantes
réduites, on se ramène au cas où $S$ est de plus intègre.

On peut supposer de plus $S=\Spec(A)$ normal intègre. En
effet, comme $A$ est excellent, le morphisme de normalisation
est fini de degré générique $1$, donc est $\pspf$-couvrant
(resp. $\altl$-couvrant). Comme $A$ est local intègre et hensélien,
$A$ est unibranche de sorte que le normalisé de $A$ est local,
donc intègre puisque normal et est noethérien hensélien puisque
fini sur $A$.

Comme $A$ est excellent, la normalisation commute à la complétion
(\ref{1prop-normal-versus-complet}) de sorte que $\hat A$ est dès
lors normal comme $A$, donc également intègre puisque normal
et local.

On peut donc supposer $T=S$ avec $S$ schéma local intègre, normal,
hensélien et excellent.

Comme $\s$ est plat sur $S$, l'image inverse ${\Z}$ de $Z$ est
encore un fermé rare de ${\s}$.  Choisissons une uniformisation
$$(\tilde\X_i\ra {\s})_{i\in I'}$$ de $({\s},{\Z})$ comme
dans ~\ref{3-unif}.  D'après ~\ref{3-theo},~\ref{3-alt-mod}
et \ref{3-perm}, on peut trouver un entier $n\geq 0$, des
$n$-isomorphismes $\tilde\X_i\ra_{n}(\tilde X_i)_u$  de sorte  que

 a)  chaque $S$-schéma $(X_i)_u$ est régulier le long de sa fibre
 spéciale $(\tilde X_i)_0$, donc au voisinage (le lieu régulier
 étant ouvert puisque les schémas considérés sont excellents).

 b) la famille $((\tilde X_i)_u)_\r$ est \pspf-couvrante
 (resp. \altl-couvrante).

 D'après a), $(\tilde X_i)_u$ est régulier au voisinage de la
 fibre spéciale et y est la réunion disjointe de ses composantes
 connexes qui sont intègres. Ainsi, au voisinage de la fibre
 spéciale, $((\tilde X_i)_u)_\r$ est schématiquement la réunion
 disjointe des composantes de $(\tilde X_i)_u$ dominant $S$. Comme
 tout voisinage ouvert de la fibre spéciale $(\tilde X_i)_0$ dans
 $((\tilde X_i)_u)_\r$ est {\pspf}-couvrant (resp. \altl-couvrant)
 (\refext{II}{2-voisinages fibre restent couvrants}), la famille
 $(X_i\ra S)_{i\in I}$ des composantes connexes de voisinages
 convenables des $(\tilde X_i)_0$ dans  $(\tilde X_i)_u, i\in I'$
 vérifie les conditions (i) et (ii).

 Soit $D'$ l'image inverse de ${\Z}$ dans $\X=\sqcup_{i\in I}
 X_i$ qu'on peut supposer non vide. Par hypothèse, $D=D'_\red$
 est un diviseur à croisements normaux strict, c'est-à-dire
 $D=\sum_{j\in J} D_j$ avec  $$D_K=\cap_{j\in K} D_j$$ régulier de
 codimension $\text{card}(J)$ pour toute partie $K\subset J$. Quitte
 à augmenter $\alpha$, on peut supposer que les $D_j$ ont des
 modèles sur $S_\alpha$, ces modèles induisant des modèles des
 $D_J$. Comme $u$ est une section de $t_\alpha$, le schéma $D_u$
 réunion schématique des $(D_i)_u$ est, topologiquement, l'image
 inverse de $Z$ dans $X_u$. D'après ~\ref{3-perm}, on peut supposer
 que chaque $(D_u)_K$ est régulier de codimension $\text{card}(K)$
 (par caténarité de $X_i$ qui est excellent puisque de type fini
 sur $S$) de sorte que $D_u$ est un diviseur à croisements normaux
 stricts le long de la fibre spéciale. Les lieux réguliers de
 $(D_u)_K$ et $X_u$ étant ouverts, on peut supposer que $D_u$
 est un diviseur à croisements normaux stricts au voisinage de la
 fibre spéciale (caténarité de $X_u$).
\end{proof}

\ifx\danslelivre\undefined
\bibliography{../fichiers-configuration-latex/bibliographie}
\bibliographystyle{../fichiers-configuration-latex/smfalpha}


\begin{center}
IV. Le théorème de Cohen-Gabber\\
Fabrice Orgogozo
\end{center}

\version
\tableofcontents

\else
\setboolean{english}{false}
\setboolean{francais}{true}
\selectlanguage{francais}
\chapter{Le théorème de Cohen-Gabber}
\addtocontents{toc}{par Fabrice Orgogozo}
\begin{center}
Fabrice Orgogozo
\end{center}
\fi

\section{$p$-bases et différentielles (rappels)}

\subsection{Définition et caractérisation différentielle}

\subsubsection{}
Pour la commodité du lecteur, et pour fixer les
notations, nous rappelons ici quelques résultats bien
connus dont nous ferons usage ci-après.
Nous conseillons au lecteur de ne s'y reporter qu'en
cas de besoin.

\begin{définition2}
Soient $k$ un corps de caractéristique $p>0$,
$K$ une extension de $k$, et $(b_i)_{i∈I}$ une famille d'éléments de $K$.
On dit que les $(b_i)$ constituent une
\emph{$p$-base} de $K$ sur $k$ (resp. sont \emph{$p$-libres} sur $k$) 
si les monômes $\prod_i b_i^{n(i)}$ ($0≤n(i)<p$,
$(n(i))_{i∈I}$ de support fini) forment une base
du $k(K^p)$-espace vectoriel $K$ (resp. sont linéairement indépendants
sur $k(K^p)$).
\end{définition2}

Si $k=\FF_p$, on parle alors de $p$-base \emph{absolue},
ou de $p$-base s'il n'y a pas d'ambiguïté.
Enfin, on appelle parfois \emph{$p$-monôme} un produit comme
ci-dessus. Un lien entre cette notion et la structure des anneaux
locaux complets ressort du théorème suivant.

\begin{theoreme2}[\ac{v}{3}{3}{, th. 1 b)}]
Soient $A$ un anneau local séparé complet de caractéristique $p>0$
et $(β_i)_{i∈I}$ une famille d'éléments de $A$ dont les classes modulo
l'idéal maximal $𝔪 _A$ forment une $p$-base du corps résiduel $A ∕ 𝔪_A$.  Il existe
alors un unique corps de représentants de $A$ contenant les éléments $β_i$.
\end{theoreme2}

\begin{remarque2}
On peut étendre de façon évidente la notion de $p$-base
au cas d'un anneau quelconque de caractéristique $p>0$,
cf. \ega{$0_{\textsc{iv}}$}{21.1-4}. Nous n'en aurons pas besoin.
\end{remarque2}

\subsubsection{}On vérifie immédiatement que les $(b_i)_{i∈I}$
forment une $p$-base de $K$ sur $k$ si et seulement si, pour tout $i∈I$,
l'élément $b_i$ n'appartient pas au sous-corps $k(K^p,(b_j)_{j≠i})$
de $K$. 
(Voir p. ex. \ega{$0_{\textsc{iv}}$}{21.4.3}.)

\subsubsection{}Pour toute extension de corps $K/k$, nous noterons $d_{K/k}$ la différentielle $K→ Ω¹_{K/k}$.

\begin{proposition2}
\label{4-caractérisation différentielle}
Soient $k$ un corps de caractéristique $p>0$, et $K$ une extension de $k$.
Une famille $(b_i)_{i∈I}$ d'éléments de $K$
est une $p$-base de $K$ sur $k$ si et seulement si les différentielles
$d_{K/k}(b_i)$ forment une base du $K$-espace vectoriel $Ω¹_{K/k}$.
\end{proposition2}

\begin{démo}
Soit $B=(b_i)_{i∈I}$ une $p$-base de $K$ sur $k$. Tout morphisme
(ensembliste) $Δ:B→K$ s'étend de manière unique en une $k$-dérivation
$D$ de $K$ : il suffit de poser $D(b₁^{n₁}\cdots b_r^{n_r})= ∑_i n_i
b₁^{n₁}\cdots b_i^{n_i-1}\cdots b_r^{n_r}$ et de 
l'étendre par $k(K^p)$-linéarité. Cela est équivalent au fait que les $d_{K/k}(b_i)$
forment une base de $Ω¹_{K/k}$. (Le fait que $K$ soit un corps
n'est pas utilisé dans cette implication.) Réciproquement, si les $d_{K/k}(b_i)$ forment une base, on 
observe que les $p$-monômes sont $k(K^p)$-linéairement indépendants : 
dans le cas contraire on aurait, pour un indice $i$ convenable,
$b_i∈K^p((b_j)_{j≠i})$, ce qui se traduirait par une relation
linéaire entre les différentielles. Soit $B'$ une $p$-base de $K$ sur $k$ contenant les
$b_i$ (\emph{loc. cit.}, 21.4.2) ; d'après l'implication précédente, on a nécessairement
$B'=(b_i)_{i∈I}$.
\end{démo}

\begin{corollaire2}
\label{4-dérivée nulle}
Soient $k$ un corps de caractéristique $p>0$, $K$ une extension de $k$.
Un élément $x$ de $K$ appartient à $k(K^p)$ si et seulement si
$d_{K/k}(x)=0$.
\end{corollaire2}

\subsubsection{}Rappelons que le \emph{$p$-rang} d'un corps est le cardinal
d'une $p$-base absolue (bien défini en vertu de ce qui précède).
On vérifie immédiatement que ce \emph{cardinal} (fini ou non) est invariant par
extension \emph{finie} de corps.

\subsection{Stabilisation}

\begin{lemme2}[cf. p. ex. \ega{$0_{\textsc{iv}}$}{21.8.1}]
Soient $K$ un corps, $k$ un sous-corps, $(k_α)_{α∈I}$ une famille
de sous-corps de $K$ telle que $⋂k_α=k$ et
\emph{filtrante décroissante}, 
c'est-à-dire telle que pour toute paire d'indices $α,β$, il existe un indice
$γ$ tel que $k_γ ⊂ k_α ∩ k_β$.
Soient $V$ un $K$-espace
vectoriel,  et $(v_i)$ ($1≤i≤n$) une famille finie de vecteurs de $V$.
Si la famille $(v_i)$ est libre sur $k$, il existe un indice $γ$ telle qu'elle
soit aussi libre sur $k_γ$.
\end{lemme2}

\begin{lemme2}
Soient $K$ un corps de caractéristique $p>0$, $k$ un sous-corps et 
$(K_α)_{α∈I}$ une famille filtrante décroissante de sous-corps contenants
$k$.
Les conditions suivantes sont équivalentes :
\begin{enumerate}
\item $⋂_α K_α(K^p)=k(K^p)$ ; 
\item pour tout ensemble fini $\{b₁,\dots,b_n\}⊂K$, $p$-libre sur $k$, il existe
un indice $α$ tel qu'il soit $p$-libre sur $K_α$ ;  
\item il existe une $p$-base de $K$ sur $k$ telle que tout sous-ensemble fini 
soit $p$-libre sur un $k_α$ pour $α$ convenable ; 
\item le morphisme canonique $Ω¹_{K/k}→ \lim_α   Ω¹_{K/K_α}$ est \emph{injectif}.
\end{enumerate}
\end{lemme2}

\begin{démo}
(i)⇒(ii) est une conséquence immédiate du lemme précédent. (ii)⇒(iii) est
trivial (toute $p$-base convient). (iii)⇒(iv) trivial (utiliser
\ref{4-caractérisation différentielle}). Vérifions (iv)⇒(i). Soit $x\notin
k(K^p)$. D'après \ref{4-dérivée nulle}, $d_{K/k}(x)≠0$ de sorte qu'il existe
$α$ tel que $d_{K/k_α}(x)$ soit également non nul. D'après \emph{loc. cit.},
cela entraîne que $x\notin k_α(K^p)$.
\end{démo}

On en déduit le lemme suivant, qui est un cas particulier
de \egazero{iv}{21.8.5}.

\begin{lemme2}
Soient $K$ un corps de caractéristique $p$, $k$ un sous-corps
 et $(K_α)_{α∈I}$ une famille filtrante décroissante de
sous-corps de $K$ contenant $k$ telle que $⋂_α K_α(K^p)=k(K^p)$. Pour toute extension finie
$L$ de $K$, on a également \mbox{$⋂_α K_α(L^p)=k(L^p)$}.
\end{lemme2}

\begin{démo}
On se ramène immédiatement au cas où $L/K$ est monogène.
Si $L/K$ est (algébrique) séparable, la conclusion résulte immédiatement
de l'existence des isomorphismes canoniques $Ω¹_{L/k}\iso Ω¹_{K/k}⊗_K L$,
$Ω¹_{L/K_α}\iso Ω¹_{K/K_α}⊗_K L$ et du critère (iv) ci-dessus.
Dans le cas contraire, $L=K(a)$, où $b=a^p∈K-K^p$.
On distingue naturellement deux cas.
Premier cas : $d_{K/k}(b)=0$, c'est-à-dire $b ∈ k(K^p)$.
Il en résulte que pour toute sous-extension $M$ de $L ∕ k$, on a l'égalité
$M(L^p)=M(K^p,b)=M(K^p)$. Ainsi,
\[⋂_α K_α(L^p)=⋂_α K_α(K^p)=k(K^p)=k(L^p),\]
et la conclusion résulte du critère (i) ci-dessus.
Second cas : $d_{K/k}(b)≠0$. On peut alors compléter $\{b\}$ en une $p$-base
de $K$ sur $k$, que l'on note $(b,(b_i)_{i∈I})$. La famille
$(a,(b_i)_{i∈I})$ est alors une $p$-base de $L$ sur $k$
et on vérifie immédiatement le critère (iii) ci-dessus :
si $(b,b₁,\dots,b_n)$ est $p$-libre sur $k_α$, il en est de même
de $(a,b₁,\dots,b_n)$.
\end{démo}

\begin{proposition2}[\cite{KKK@Matsumura}, §30, lemme 6]\label{4-lemme Nagata}
Soient $K$ un corps de caractéristique $p>0$ et
$(K_α)$ une famille filtrante décroissante de sous-corps
\emph{cofinis} — c'est-à-dire tels que les degrés $[K:K_α]$
soient finis — telle que $⋂_α K_α=K^p$.
Alors, pour toute extension finie $L/K$, il existe un
indice $β$ tel que pour tout sous-corps cofini $K′⊂K_{β}$ on ait :
\[
\rang_L Ω¹_{L/K′}=\rang_{K} Ω¹_{K/K′}.
\]

\end{proposition2}

\begin{démo}
On souhaite se ramener au cas où $L$ est monogène sur $K$.
Pour cela considérons une sous-$K$-extension $M$ de $L$ et posons $M_α=K_α(M^p)⊂M$.
Ce sont des sous-corps cofinis de $M$ et, pour $α,β$ et $γ$ comme
dans l'énoncé, on a $M_γ⊂M_α∩M_β$. En vertu du lemme
précédent, appliqué dans le cas particulier où $k=K^p$,
on a l'égalité $⋂_α M_α=M^p$. On se ramène au cas où $L/M$ est monogène
(par récurrence) en remarquant également que les extensions $M_α/K_α$ sont
\emph{finies}, et que pour tout sous-corps $K′$ de $K$ on a
$Ω¹_{M/K′}=Ω¹_{M/K′(M^p)}$.

Supposons dorénavant l'extension $L ∕ K$ monogène.
Si elle est (algébrique) séparable, le théorème est trivial :
on a $Ω¹_{L/K′} ⭉  Ω¹_{K/K′}⊗_K L$ pour tout $K′⊆ K$.
Sinon, $L=K(a)$, où $a^p=b∈K-K^p$,
et, pour chaque $K′ ⊆ K$, $Ω¹_{L ∕ K ′}$ est naturellement isomorphe à
\[
\big(Ω¹_{K ∕ K′}/Kd_{K ∕ K′}(b)\big)\! ⊗_K \! L   ⊕ L d_{L ∕ K′}(a).
\]
De plus, $d_{K ∕ K^p}(b)$ et $d_{L ∕ L^p}(a)$ sont non nuls
car $b$ (resp. $a$) n'appartient pas à $K^p$ (resp. $L^p$).
Puisque $K^p=⋂K_α$ (resp. $L^p=⋂ K_α(L^p)$) , il existe
un $β$ tel que $d_{K ∕ K_β}(b) ≠ 0$ (resp. $d_{L ∕ K_β}(a)≠
0$). Il résulte de l'isomorphisme ci-dessus que pour chaque $K′ ⊆ K_β$, on a
l'égalité $\rang_L Ω¹_{L ∕ K ′}=\rang_K Ω¹_{K ∕ K ′}$. CQFD.
\end{démo}

\subsubsection{}\label{4-dérivations continues}
Rappelons enfin que si $A$ et $B$ sont deux anneaux linéairement topologisés
(\egazero{i}{7.1.1}), et $A→B$ un morphisme \emph{continu}, 
le $B$-module $Ω¹_{B/A}$ est un $B$-module
\emph{topologique}, la topologie étant déduite 
de celle de $B⊗_A B$ par restriction et passage au quotient 
(\egazero{iv}{20.4.3}). Le $B$-module sous-jacent ne dépend pas des topologies de $A$ et $B$.
On note $\chap{Ω}¹_{B/A}$ son \emph{séparé complété} ; 
il est isomorphe à une limite de $Ω¹$ de morphismes
entre anneaux topologiques 
discrets (\emph{loc. cit.}, 20.7.4).
Pour tout $B$-module topologique $L$, le morphisme canonique
induit par la dérivation universelle est un \emph{isomorphisme} :
$$
\mathrm{Hom.cont}_{B}(Ω¹_{B/A},L)\iso \mathrm{D\acute{e}r.cont}_{A}(B,L).
$$
Si $B$ et $L$ sont séparés et complets, le terme de gauche s'identifie
canoniquement à $\mathrm{Hom.cont}_{B}(\chap{Ω}¹_{B/A},L)$.
Comme on le constate dans le cas particulier très simple où $A$ est un corps
et $B$ un anneau de séries formelles, le $B$-module $\chap{Ω}¹_{B/A}$ a des propriétés de finitude
bien plus remarquables que $Ω¹_{B/A}$ (\emph{loc. cit.}, exemple 20.7.6 et prop.
20.7.5).

\section{Les théorèmes de Cohen-Gabber en caractéristique $>0$}

\subsection{Le théorème de Cohen-Gabber non équivariant en caractéristique $>0$}
Le but de ce paragraphe est de démontrer la variante suivante du théorème de structure
des anneaux locaux nœthériens complets 
\egazero{iv}{19.8.8} (ii), dû à Irving S. Cohen.

\begin{theoreme2}[théorème de Cohen-Gabber ;  \cite{Conference-Illusie@Gabber}, lemme 8.1]
\label{4-Cohen-Gabber égale}
Soit $A$ un anneau local complet nœthérien \emph{réduit}, d'égale
caractéristique $p>0$, équidimensionnel
de dimension $d$ et de corps résiduel $k$.
Il existe un sous-anneau $A₀$ de $A$, isomorphe à $k[[t₁,\dots,t_d]]$,
tel que $A$ soit fini sur $A₀$, sans torsion et génériquement étale.
De plus, le morphisme $A₀→A$ induit un isomorphisme sur les corps résiduels.
\end{theoreme2}

\begin{remarques2}Ce résultat apparaît explicitement comme hypothèse,
pour $A$ intègre, dans \ega{$0_{\textsc{iv}}$}{21.9.5}.
L'expression « génériquement étale » signifie ici qu'il existe un
ouvert dense de $\Spec(A₀)$ au-dessus duquel le morphisme
$\Spec(A)→\Spec(A₀)$ est étale.
\end{remarques2}

\subsubsection{}La démonstration du théorème, qui est une adaptation au cas non irréductible
de \cite{Conference-Illusie@Gabber}, occupe le reste de cette section.
Nous supposerons par la suite $d>0$, sans quoi l'énoncé est évident.
Dans les alinéas \ref{4-début Cohen-Gabber} à \ref{4-fin Cohen-Gabber},
nous allons montrer qu'il existe un corps de représentants $κ$ de $A$ tel
que le $A$-module des formes différentielles
complété $\chap{Ω}¹_{A/κ}$ soit de rang générique
égal à $d$ sur chaque composante irréductible.
En (\ref{4-fin démo}) nous verrons comment en déduire
rapidement le théorème.

\subsubsection{}\label{4-début Cohen-Gabber}Soit $(b_i)_{i∈E}$ une $p$-base de $k=A/\MM_A$.
Choisissons des relèvements arbitraires $β_i$
des $b_i$ dans $A$. Rappelons qu'il existe un unique corps de
représentants $κ\subset A$ contenant les $β_i$ et se surjectant sur $k$
(cf. \ac{ix}{2}{2}{, th. 1 a)}). Changer de corps de
représentants revient donc à changer les $β_i$.
Fixons également un système de paramètres $τ₁,\dots,τ_d$
de $A$ ;  nous ne le changerons qu'à la fin de la
démonstration (\ref{4-fin démo}).

\subsubsection{}Pour toute partie \emph{finie} $e\subset E$, posons $κ_e:=κ^p(β_i,\,i\notin e)\subset κ$.
Les trois propriétés suivantes sont évidentes :
$$\text{pour toute partie finie $e⊂E$}, [κ:κ_e]<+∞,$$
$$\text{pour toutes parties finies $e,e'⊂E$}, κ_{e\cup
e'} ⊂ κ_{e}\cap κ_{e'},$$
$$⋂_{e\subset E} κ_e = κ^p.$$

\subsubsection{}Soient $\Spec(\sur{A})$ une composante
irréductible de $\Spec(A)$, munie de la structure
réduite, et $\sur{τ₁},\dots,\sur{τ_d}$ les images des $τ_i$
dans $\sur{A}$ par la surjection canonique $A↠\sur{A}$.
Considérons le diagramme d'anneaux :
\[
\xymatrix{
\kappa_e[[\sur{\tau_1}^p,\dots,\sur{\tau_d}^p]] \ar[r] \ar[d] &
\kappa[[\sur{\tau_1},\dots ,\sur{\tau_d}]]  \ar[r]
\ar[d] & \sur{A} \ar[d] \\
L_{\kappa,e} \ar[r] & L_\kappa \ar[r] & L 
}
\]
où les flèches horizontales sont les homomorphismes
canoniques, et les flèches verticales les inclusions
dans les corps de fractions respectifs. Les flèches
horizontales sont injectives et correspondent à
des morphismes finis. Pour la seconde, cela résulte
du fait que le module $\sur{A}$ est \emph{quasi-fini}
(\ega{$0_{\textrm{I}}$}{7.4.1}) sur $κ[[\sur{τ₁},\dots
,\sur{τ_d}]]$ donc de type fini car l'idéal $(\sur{τ₁},\dots
,\sur{τ_d})\sur{A}$ est un idéal de définition
(\ega{$0_{\textrm{I}}$}{7.4.4}). Enfin, 
les $\sur{τ_i}$ sont analytiquement indépendants sur
$κ$ : le sous-anneau $κ[[\sur{τ₁},\dots
,\sur{τ_d}]]$ de $\sur{A}$ est bien un anneau de séries
formelles (\ega{$0_{\textrm{IV}}$}{16.3.10}).

On a observé ci-dessus que la famille des $κ_e⊂κ$, $e⊂E$,
satisfait aux hypothèses de la proposition \ref{4-lemme
Nagata}. On vérifie
immédiatement qu'il en est de même de la famille
des sous-corps $L_{κ,e}$ de $L_κ$ ;  on a donc l'égalité
\begin{equation}\label{4-equation1}
\rang_L Ω¹_{L/L_{κ,e}}=\rang_{L_κ} Ω¹_{L_κ/L_{κ,e}},
\end{equation}
dès que l'ensemble fini $e$ est suffisamment grand.

Posons $R_κ=\kappa[[\sur{τ_1},\dots
,\sur{τ_d}]]$ et
$R_{κ,e}=\kappa_e[[\sur{\tau_1}^p,\dots,\sur{\tau_d}^p]]$.
Le terme de gauche de (\ref{4-equation1})
est le rang générique du $\sur{A}$-module
$Ω¹_{\sur{A}/R_{κ,e}}$, c'est-à-dire le rang de son tensorisé
avec $L$.
Remarquons que d'après \ega{$0_{\textrm{IV}}$}{21.9.4}, 
$Ω¹_{\sur{A}/R_{κ,e}}$ s'identifie au module $\sur{A}$-module
$\chap{Ω}¹_{\sur{A}/κ_e}$ de formes différentielles complété.
Le terme de droite est quant à lui le 
rang du $R_{κ,e}$-module libre $Ω¹_{R_{κ,e}/R_κ}$.
Ce dernier est égal à $d+\rang_κ Ω¹_{κ/κ_e}=d+|e|$
(où $|-|$ désigne le cardinal d'un ensemble), de sorte que 
la formule \ref{4-equation1} se réécrit :
\begin{equation}\label{4-equation2}
\rang_{\sur{A}} \chap{Ω}¹_{\sur{A}/κ_e}=d+|e|.
\end{equation}

\subsubsection{}La proposition suivante va nous permettre de modifier
le corps des représentants de façon à pouvoir supposer
$e$ vide (de façon équivalente : $κ_e=κ$).

\begin{proposition2}\label{4-changement relèvements}
Il existe une partie finie $e$ de $E$ et des éléments
$β'_i$, pour $i∈e$, relevant les $b_i$ tels que,
pour chaque composante irréductible intègre $\Spec(\sur{A})$ de
$\Spec(A)$,
les conditions suivantes soient vérifiées :
\begin{enumerate}
\item $\rang_{\sur{A}} \chap{Ω}¹_{\sur{A}/κ_e}=d+|e|$,
\item les images des $dβ'_i$ dans
$\chap{Ω}¹_{\sur{A}/κ_e}⊗_{\sur{A}}
L$, où $L=\Frac(\sur{A})$, sont $L$-linéairement indépendantes.
\end{enumerate}
\end{proposition2}

L'égalité \ref{4-equation1} (et donc \ref{4-equation2}) étant valable, pour chaque
composante irréductible, dès que $e$ est suffisamment 
grand, on peut choisir un tel ensemble qui convient
pour chacune d'entre elles.
La propriété (i) en découle.

Pour démontrer la propriété (ii), nous utiliserons le lemme élémentaire suivant.

\begin{lemme2}
Soient $\sur{A}$ et $L$ comme ci-dessus. 
Pour tout idéal non nul $I$ de $\sur{A}$, l'ensemble
des $df\otimes_{\sur{A}} L$, pour $f∈I$, est une famille 
génératrice du $L$-espace vectoriel
$\chap{Ω}¹_{\sur{A}/κ_e}\otimes_{\sur{A}} L$.
\end{lemme2}

\begin{proof}
Soient $f₀∈I$ non nul, et $ω₀=df₀$. Pour tout 
$b∈\sur{A}$, $d(bf₀)=bω₀+f₀db$.
La famille des $d(bf₀)⊗1$ 
contient $ω₀⊗1$ ;  d'après la formule précédente,
le $L$-espace vectoriel qu'elle engendre
contient donc les $db⊗1$ pour chaque $b∈\sur{A}$.
\end{proof}

Soit $\{𝔭₁,\dots,𝔭_c\}$ l'ensemble des idéaux premiers
minimaux de $A$. 
Pour chaque $j∈\{1,\cdots,c\}$, posons $A_j=A/𝔭_j$ et
$X_j=\Spec(A_j)$
la composante irréductible intègre de $X=\Spec(A)$ correspondante.
Notons pour tous $i∈e$ et $j∈\{1,\dots,c\}$,
$β_{i,j}$ l'image dans $A_j$ de $β_i∈A$.
(Rappelons que les $β_i$ font partie d'une $p$-base de $κ\subset A$.)
Nous allons démontrer par récurrence sur $j$ ($0≤j≤c$) qu'il existe des éléments 
$\{m_{i,j}\}$ dans $\MM_A$, pour $i∈e$, tels que
les images des éléments $β_i+m_{i,j}$ dans chacun des
anneaux $A_1,\dots,A_j$ aient des différentielles
linéairement indépendantes dans chacun des
espaces vectoriels
$Ω¹_{A₁/R_{κ,e}}\otimes_{A₁} \Frac{\,A₁}$,
\dots, $Ω¹_{A_j/R_{κ,e}}\otimes_{A_j} \Frac{\,A_j}$. 
Pour $j=0$, cette condition est vide.
Supposons l'assertion démontrée pour un $j≤c-1$ et
montrons la pour $j+1$.
Quitte à remplacer $β_i$ par $β_i+m_{i,j}$, on peut 
supposer que $m_{i,j}=0$ pour tout $i∈e$.
L'anneau $A$ étant réduit, les $𝔭_α$ forment une
décomposition primaire \emph{réduite}
de $(0)$, de
sorte que l'idéal $\got{q}_j:=𝔭_1\cap \cdots \cap 𝔭_j$
($=\mathrm{Ker}(A→A_1×\cdots×A_j)$) 
n'est pas contenu dans $𝔭_{j+1}$.
Si $j>0$, notons $I_{j+1}$ son image dans
$\sur{A}=A_{j+1}(=A/𝔭_{j+1})$ ;  c'est un idéal non nul. Si $j=0$,
on considère $\MM_{\sur{A}}$. 
D'après (i), $\rang_{\sur{A}} \chap{Ω}¹_{\sur{A}/κ_e}=d+|e|≥|e|$ ; 
d'autre part, la famille $d(I_{j+1})$ est génératrice
dans $\chap{Ω}¹_{\sur{A}/κ_e}\otimes_{\sur{A}} L$
(où $L=\Frac{\,\sur{A}}$). 

\begin{lemme2}
Soient $V$ un espace vectoriel de dimension au moins
$n$, $b₁,\dots,b_n$ des vecteurs de $V$ et $W$ une
famille génératrice. Il existe une famille
$w₁,\dots,w_n$ d'éléments de $W ∪ \{0\}$ tels que 
les $b_i+w_i$ soient linéairement indépendants.
\end{lemme2}

\begin{démo}Par récurrence immédiate sur $n$.\end{démo}

Il existe donc des éléments $m'_{i,j+1}∈I_{j+1}$, $i∈e$,
tels que les différentielles des éléments
\mbox{$d((β_i \mod 𝔭_{j+1})+m'_{i,j+1})$}, $i∈e$, soient linéairement indépendantes
dans $\chap{Ω}¹_{\sur{A}/κ_e}\otimes_{\sur{A}} L$.

Relevons les $m'_{i,j+1}$ en des éléments $m_{i,j+1}$
de $\got{q}_j$ si $j>0$, ou de $\MM_A$ si $j=0$. Par construction,
ils satisfont la propriété escomptée au cran $j+1$.

\subsubsection{}\label{4-fin Cohen-Gabber}
Considérons le sous-corps $κ':=κ^p(β_i, i∉e;β'_i, i∈e)=κ_e(β'_i, i∈e)\subset A$, où les $β'_i$ ($i∈e$) sont
comme en \ref{4-changement relèvements}. Il s'envoie isomorphiquement sur
$k=A/\MM_A$ par réduction : son image contient $k^p$ et
les images des $β_i$ ($i∉e$), $β'_i$ ($i∈e$), qui constituent une $p$-base
de $k$.
Des égalités \ref{4-equation2} et de la propriété (ii) de 
\ref{4-changement relèvements},  on tire : 
$$\rang_{\sur{A}} \chap{Ω}¹_{\sur{A}/κ'}=d,$$ pour toute composante irréductible
intègre
$\Spec(\sur{A})$ de $X$. 
Par la suite, nous noterons encore $κ$ ce nouveau
corps de représentants.

\subsubsection{}\label{4-fin démo}Le $A$-module $\chap{Ω}¹_{A/κ}$ étant de rang générique $d$ sur chaque
composante irréductible, on montre en procédant comme précédemment, 
qu'il existe des éléments $f_1,\dots,f_d$ de $A$ tels que les $d(f_i \ \mathrm{mod}
\, 𝔭_α)\otimes_{A_j} \Frac{\,A_j}$ forment une base de 
$\chap{Ω}¹_{A_j/κ} \otimes_{A_j} \Frac{\,A_j}$ pour chaque composante
irréductible $\Spec(A_j)$ de $X$. Quitte à les
multiplier individuellement par une puissance $p$-ième
d'un élément appartenant à $\MM_A-⋃_j 𝔭_j$, on peut les supposer dans $\MM_A$.
Rappelons que l'on a choisi un système de paramètres
$τ_1,\dots,τ_d$ dans $A$, de sorte que le morphisme
$\Spec(A)→\Spec(k[[τ_1,\dots,τ_d]])$ soit fini.

Posons, pour $i∈\{1,\dots,d\}$,
$$
t_i:=τ_i^p(1+f_i).
$$
Soient $A₀$ le sous-anneau $κ[[t_1,\dots,t_d]]$ de
$A$, $X₀=\Spec(A₀)$. Le morphisme $X→X₀$ est fini : cela résulte du 
fait que les éléments $1+f_i$ sont des unités de $A$. Vérifions 
qu'il est génériquement étale.
L'anneau $A$ étant nœthérien complet, le $A$-module de type
fini $Ω¹_{A/A₀}$ est également complet et coïncide donc
avec le module des formes différentielles complété $\chap{Ω}¹_{A/A₀}$.
Les anneaux $A₀$ et $A$ étant métrisables, et tout sous-$A$-module
de $\chap{Ω}¹_{A/κ}$ étant fermé, la suite $$
\chap{Ω}¹_{A₀/κ}\chap{⊗}_{A₀} A →  \chap{Ω}¹_{A/κ} →
\chap{Ω}¹_{A/A₀}=Ω¹_{A/A₀}→0
$$
est \emph{exacte} (\ega{$0_{\textrm{IV}}$}{20.7.17}).
Il résulte de l'hypothèse sur les éléments $f_i$ et de la formule
$$
d(t_i)=τ_i^p df_i
$$
qu'au-dessus de chaque point maximal de $X=\Spec(A)$, la première flèche
est surjective. On en déduit que le $A$-module $Ω¹_{A/A₀}$ est génériquement nul,
CQFD.

\subsection{Le théorème de Cohen-Gabber équivariant en caractéristique $>0$}

\subsubsection{}Nous allons démontrer ici une généralisation
du théorème \ref{4-Cohen-Gabber égale} dans le cas d'un anneau non nécessairement
équidimensionnel, muni d'une action d'un groupe fini.

\begin{theoreme2}\label{4-Cohen-Gabber équivariant}
Soient $A$ un anneau local nœthérien complet réduit, de dimension $d$, de corps résiduel $κ$ et 
$G$ un groupe fini agissant sur $A$ avec $|G|$ inversible dans $κ$. Alors, il
existe un morphisme fini génériquement étale, $G$-équivariant, $κ[[t₁,\dots,t_d]]→ A$, où $κ→A$ relève l'identité
de $κ$ et $G$ agit trivialement sur les $t_i$.
\end{theoreme2}

Commençons par une proposition.

\begin{proposition2}
Soit $A$ un anneau muni d'une action d'un groupe fini $G$
d'ordre inversible sur $A$ et soit $B=\Fix_G(A)$ le sous-anneau
des invariants.
\begin{enumerate}
\item L'anneau $B$ est
\begin{enumerate}
\item nœthérien si $A$ l'est ;
\item réduit si $A$ l'est ;
\item local d'idéal maximal $𝔪 ∩ B$ si $A$ est local d'idéal
maximal $𝔪$, de corps résiduel isomorphe au sous-corps $\Fix_G(A ∕ 𝔪)$ de $κ=A ∕ 𝔪$.
\end{enumerate}
\item Le morphisme $\Spec(A) → \Spec(B)=\Spec(A) ∕ G$ est
\begin{enumerate}
\item fini si $A$ est nœthérien ;
\item génériquement étale si $A$ est de plus réduit.
\end{enumerate}
\end{enumerate}
\end{proposition2}

\begin{démo}
(i) Notons $\Tr$ le morphisme $B$-linéaire $\tr:A→  B$, $x\mapsto \frac{1}{|G|}∑_{g∈G}
g(x)$, parfois appelé « opérateur de Reynolds ».
Pour tout idéal $I$ de $B$, on a $IA ∩ B = I$. En effet,
l'inclusion $I ⊆ IA ∩ B$ est triviale et l'inclusion opposée
résulte du fait que si $x ∈ IA ∩ B$,
sa « trace » $x=\Tr(x)$ appartient, par $I$-linéarité, à $IB=I$.
On en déduit immédiatement l'énoncé a). L'énoncé b) est trivial.
Si $A$ est local, on a $A-𝔪=A^×$. Il résulte d'une part que $G$
stabilise globalement $𝔪$ et d'autre part que $\Fix_G(A) - \Fix_G(𝔪) = \Fix_G(A)^×$.
Ainsi, $B$ est maximal d'idéal $𝔫=\Fix_G(𝔪)$. Enfin, Le morphisme
canonique $B ∕ 𝔫 → \Fix_G(κ)$ déduit de l'inclusion canonique
$B ∕ 𝔫 → κ$ est un isomorphisme. En effet, si $a ∈A$ est un
relèvement arbitraire de $λ ∈ \Fix_G(κ)$, l'élément
$b=\Tr(a)$ en est un relèvement $G$-équivariant. Ceci achève
la démonstration du c).
(ii.a). Nous allons montrer que le morphisme entier $\Spec(A) →
\Spec(B)$ est fini par réduction au cas bien connu où $A$ est un corps.

❧ Réduction au cas réduit. Soient $N$ le nilradical de $A$ et $M=N∩B$ celui de $B$.
Pour chaque entier $i∈𝐍$, le $A/N$-module $N^i/N^{i+1}$ est de type
fini, car $A$ est supposé nœthérien, et nul pour $i ≫  0$.
Le module $\mathrm{gr}_{N}(A)=⊕_{n≥0}N^i/N^{i+1}$
est donc de type fini sur $\mathrm{gr}_{N}^0(A)=A/N$.
Si ce dernier est de type fini sur $B/M=\mathrm{gr}_{M}^0(B)$,
il en est de même de $\mathrm{gr}_{N}(A)$
sur $\mathrm{gr}_{M}(B)$ et finalement (\ac{iii}{2}{9}{, cor. 1})
de $A$ sur $B$, par complétude de l'anneau nœthérien $B$
pour la topologie $M$-adique.

❧ Réduction au cas d'un produit de corps.
Supposons $A$ réduit et considérons l'ensemble fini $\{𝔭_i\}_{i∈I}$
des idéaux premiers minimaux de $A$.
Pour chaque $i$, $𝔮_i=𝔭_i ∩ B$ est un idéal premier \emph{minimal} de $B$.
Cela résulte du théorème de Cohen-Seidenberg
(\ac{v}{2}{1}{, th. 1 et cor. 2}) et de
la transitivité de l'action de $G$ sur les fibres
de $\Spec(A) → \Spec(B)$ (\emph{op. cit.}, n°2, th. 2).
Soit $\Frac{ A}$ (resp. $\Frac{ B}$) l'anneau total des fractions de $A$ (resp. $B$) ;
c'est un produit de corps dans lequel $A$ (resp. $B$) s'injecte,
isomorphe au semi-localisé de $A$ en les $\{𝔭_i\}_{i∈I}$ (resp.
$\{𝔮_i\}_{i∈I}$). Soit $S=A-⋃_i 𝔭_i$ ; on a donc $\Frac{ A}=S^{-1}A$.
D'après (\emph{op. cit.}, §1, n°1, prop. 23), on a
$\Fix_G(S^{-1}A)=\Fix_G(S)^{-1}B$, de sorte que
$\Fix_G(\Frac{ A})=\Frac{ B}$ et $A⊗_{B} \Frac{ B}≅\Frac{ A}$.
Supposons $\Frac{\,A}$ fini sur $\Frac{\,B}$, de sorte qu'il existe d'après
l'isomorphisme précédent un nombre fini $n$ d'éléments de $A$ qui engendrent $\Frac{\,A}$
sur $\Frac{\,B}$. Observons que l'opérateur $\tr:A→B$ définit,
par composition avec le produit, un accouplement $A⊗_{B} A→ B$ 
qui est parfait sur les anneaux de fractions :
on se ramène à montrer que si $e_i$ est un idempotent
correspondant au facteur $K_i=\Frac{\,A/𝔭_i}$ de $\Frac{\,A}$, l'élément
$\tr(e_i)$ est non nul ;  il est en effet égal à $\frac{|G_i|}{|G|}$,
où $G_i$ est le stabilisateur de $e_i$.
Les $n$ éléments ci-dessus définissent donc un \emph{plongement} $A^G$-linéaire
de $A$ dans $B^n$. On peut conclure par nœthérianité.

❧ Réduction au cas d'un corps.
Soit donc $A=∏_i K_i$ un produit fini de corps et posons $X=\Spec(A)=∐_i η_i$.
Si $X=X₁∐X₂$, où $X₁$ et $X₂$ sont $G$-stables, $X/G=(X₁/G)∐(X₂/G)$ de sorte
que l'on se ramène immédiatement au cas où $X/G$ est connexe, c'est-à-dire où l'action
de $G$ est \emph{transitive}. Pour tout $i$, notons $G_i$ le groupe de
décomposition correspondant. D'après le cas classique (cas d'un corps),
$η_i → η_i/G_i$ est fini étale. Il en résulte que le morphisme
$X→ ∐ η_i / G_i$ est fini. Enfin, puisque pour tout $i$,
$η_i/G_i\iso X/G$ (\emph{loc. cit.}, §2, n°2, prop. 4), le résultat
(ii.a) en découle.

L'énoncé (ii.b) est désormais évident.
\end{démo}

\subsubsection{}Soient $A$ et $G$ comme dans l'énoncé
du théorème \ref{4-Cohen-Gabber équivariant}. Il résulte de la
proposition précédente l'on a l'égalité $\dim(B)=\dim(A)<+∞$, où
l'on note $B=\Fix_G(A)$. Nous noterons $d$ leur dimension commune.
Soit $B/I$ le quotient maximal $d$-équidimensionnel de $B$. D'après le théorème de
Cohen-Gabber \ref{4-Cohen-Gabber égale}, il existe un corps de représentants
$λ↪B/I$ et un système de paramètres $t₁,\dots,t_d$
de $B/I$ tel que $λ[[t₁,\dots,t_d]]→ B/I$ soit fini, génériquement étale.
On peut relever l'inclusion $λ↪B/I$ en une inclusion $λ↪B$ :
cela résulte par exemple, en caractéristique résiduelle positive
(seul cas non trivial), de la correspondance entre sous-corps de
représentants et relèvements d'une $p$-base donnée du corps résiduel.
Enfin, on peut relever le système de paramètres de $B/I$ en un système
de paramètres de $B$ : cela résulte, par dévissage, du lemme suivant.

\begin{lemme2}
Soient $A↠B$ une surjection d'anneaux locaux nœthériens et
$b∈\MM_B$ un élément \emph{sécant} pour $B$, c'est-à-dire
tel que $\dim(B/b)=\dim(B)-1$.
Il existe un relèvement de $b$ dans $A$ sécant pour $A$.
\end{lemme2}

Pour des généralités sur les suites sécantes,
voir par exemple \emph{op. cit.}, chap. VIII, §3. n°2.

\begin{démo}
On se ramène immédiatement au cas où $B=A/(f)$, $f∈A$.
Soit $a ∈ A$ un relèvement arbitraire de $b$ ;  par hypothèse, on a
$\dim(A/(f,a))=\dim(B)-1$.
Si $\dim(B)=\dim(A)-1$,
on a nécessairement $\dim(A/a)=\dim(A)-1$ car la dimension
chute d'au plus un par équation. Dans le cas contraire,
$f$ appartient à la réunion $⋃_{i=1}^{n}𝔭_i$, où les $𝔭_i$ sont les
idéaux premiers de $A$ de cohauteur
$\dim(A)$. Supposons que $f∈𝔭₁,\dots,𝔭_r$, et seulement
ces idéaux-ci. La conclusion ne peut être mise en défaut
que si $a+(f)⊆ ⋃_{i=1}^n 𝔭_i$, c'est-à-dire si tous les relèvements de $b$
sont non sécants. Pour chaque $i≤r$,
on a $a∉𝔭_i$ car $f$ appartenant à $𝔭_i$,
on aurait $\dim(A/(f,a))=\dim(A)$. Il en résulte notamment
que $r ≠ n$.
Il suffit donc de montrer que l'hypothèse
$a+(f) ⊆ ⋃_{i=r+1}^n 𝔭_i$ est absurde.
On aurait en effet $a+f^m=a+f ⋅ f^{m-1}∈  ⋃_{i=r+1}^n 𝔭_i$ pour tout $m$
et finalement $f^m(1-f^{m-m'})∈ 𝔭_i$ pour deux entiers $m>m'$ et un indice $r+1 ≤ i ≤ n$.
On en tire immédiatement $f∈𝔭_i$, ce qui est contraire à l'hypothèse.
\end{démo}

\subsubsection{}L'extension $κ/λ$ étant étale, car $λ=\Fix_G(κ)$,
le morphisme $κ→A/𝔪$ se relève \emph{uniquement} en un $λ$-homomorphisme $k→A$ ;
ce morphisme est $G$-équivariant. Le morphisme $A/B$ étant fini, génériquement étale,
ceci achève la démonstration du théorème \ref{4-Cohen-Gabber équivariant}.

\section{Autour du théorème de Epp}

\subsection{Énoncé (rappel)}

\subsubsection{}Si $X$ est un schéma \emph{réduit} n'ayant qu'un nombre fini
de composantes irréductibles, nous noterons $X_{\norm}$ son normalisé
(\ega{2}{6.3.6--8}).

\begin{theoreme2}[Helmut Epp, \cite{Eliminating@Epp}, théorème 1.9]\label{4-Epp}
Soit $T→S$ un morphisme local dominant de traits complets, de caractéristique
résiduelle $p>0$. Notons $κ_S$ et $κ_T$ leurs corps résiduels respectifs.
Supposons $κ_S$ \emph{parfait} et le sous-corps parfait maximal
de $κ_T$ \emph{algébrique} sur $κ_S$.
Il existe une extension finie de traits $S'→S$ telle que
le produit fibré réduit normalisé
$$
T':={(T\times_S S')_{\red}}_{\norm}
$$
ait une \emph{fibre spéciale réduite} au-dessus de $S'$.
\end{theoreme2}

\begin{remarque2}\label{4-remarque Epp}
En caractéristique mixte, le produit fibré $T\times_S S'$ 
est réduit. En effet, le morphisme $T'→S'$ (obtenu par
changement de base d'un plat) est plat, et $S'$ est
intègre si bien que l'anneau des fonctions de $T'$
s'injecte dans l'anneau des fonctions de sa fibre
générique. Il suffit donc de prouver que 
cette dernière est réduite. Or, en caractéristique nulle,
toute extension de corps est séparable. 
On vérifie
également sans difficulté que la conclusion du théorème est encore valable si 
l'on suppose seulement $S$ complet, mais pas
nécessairement $T$ (cf. \emph{loc. cit.}, §2).
\end{remarque2}

\subsection{Sorites}\label{4-stabilité E}

\subsubsection{}Nous dirons qu'une extension
de corps $K/k$ de caractéristique $p>0$ a la \emph{propriété de Epp} si tout
élément du sous-corps parfait maximal de $K$, $K^{p^∞}:=\cap_{i≥0} K^{p^{n}}$,
est algébrique séparable sur $k$. Pour $k$ parfait, c'est l'hypothèse
faite sur $κ_T/κ_S$ dans \ref{4-Epp}. Dans ce court
paragraphe, on rappelle quelques résultats élémentaires
de stabilité pour cette notion.
Commençons par un lemme.

\begin{lemme2}
Pour tout corps $K$ d'exposant caractéristique $p>1$, on a, dans une clôture séparable $K^{\sep}$ de $K$,
\[(K^{p^∞})^{\sep}=(K^{\sep})^{p^∞}.\]
\end{lemme2}

\begin{démo}
L'inclusion $(K^{p^{∞}})^{\sep}\subset (K^{\sep})^{p^∞}$ est
évidente : $K^{p^∞}$ est parfait donc toute extension algébrique,
en particulier sa clôture séparable
$(K^{p^{∞}})^{\sep}$, l'est également. Comme cette dernière est contenue
dans $K^{\sep}$, elle est également contenue dans son plus grand
sous-corps parfait $(K^{\sep})^{p^∞}$.

Réciproquement, considérons $x∈(K^{\sep})^{p^∞}$, et notons, pour chaque
entier $n≥0$, $x_n$ sa racine $p^n$-ième dans $K^{\sep}$ et $f_n$ son polynôme
minimal (\emph{unitaire}). Compte tenu d'une part de l'expression de $f_n$ en fonction
des polynômes symétriques en les conjugués galoisiens de $x_n$ et d'autre part
de l'injectivité et de l'additivité de l'élévation à la puissance $p^n$-ième,
on a l'égalité $f₀=f_n^{(p^n)}$, où $f_n^{(p^n)}$ est le polynôme
obtenu à partir de $f_n$ en élevant les coefficients à la puissance $p^n$-ième.
Il en résulte que les coefficients du polynôme minimal $f₀$ de $x$
appartiennent à $K^{p^∞}$.
\end{démo}

\begin{proposition2}[Cf. \cite{Eliminating@Epp}, §0.4]\label{4-prop Epp}
 Soit $k$ un corps d'exposant caractéristique $p$.
\begin{enumerate}
\item Soient $L/K$ et $K/k$ ayant la propriété
de Epp. Alors, $L/k$ a la propriété de Epp.
\item Toute extension finie de $k$ a la propriété de Epp.
\item Si $p>1$, pour tout entier naturel $d$, l'extension $\big(\Frac{\,k[[x₁,\dots,x_d]]}\big)/k$
a la propriété de Epp.
\item Si $p>1$, pour toute inclusion $k⊂A$,
où $A$ est un anneau local complet noethérien intègre,
induisant un isomorphisme sur les corps résiduels,
l'extension $(\Frac{\,A})/k$ a la propriété de Epp.
\end{enumerate}
\end{proposition2}

\begin{démo} Supposons immédiatement $p>1$ sans quoi
(i) et (ii) sont triviaux.\\
(i) Par hypothèse on a dans une clôture séparable de $L$ l'inclusion $L^{p^∞}\subset K^{\sep}$.
Comme le corps $L^{p^∞}$ est \emph{parfait}, on en déduit
que $L^{p^∞}\subset (K^{\sep})^{p^∞}=(K^{p^∞})^{\sep}\subset
k^{\sep}$, où l'égalité résulte du lemme précédent.

(ii) Toute extension étale a tautologiquement
la propriété de Epp. D'après (i), il reste à considérer
le cas d'une extension radicielle $K/k$.
Si elle est de hauteur $≤r$, on a $K^{p^r}\subset k$ et
en particulier $K^{p^∞}\subset k\subset k^{\sep}$.

(iii) Soit $A=k[[x₁,\dots,x_d]]$ et $K$ son corps
des fractions. Montrons que $K^{p^∞}=k^{p^∞}$.
Comme $K$ est contenu dans $k((x₁,\dots,x_{d-1}))((x_d))$,
on se ramène par récurrence au cas où $d=1$.
Tout élément non nul de $k((t))^{p^∞}$ a une valuation infiniment
$p$-divisible donc nulle, de sorte que
$k((t))^{p^∞}-\{0\}$ est contenu dans
$k[[t]]^×$ et finalement dans $k^{p^∞}$ par un calcul immédiat.
(iv) Cela résulte des observations précédentes et du théorème
de structure de·Cohen.
\end{démo}

\section{Le théorème de Cohen-Gabber en caractéristique mixte}

\subsection{Anneaux de Cohen et lissité formelle (rappels)}

\subsubsection{}Pour la commodité du lecteur, nous énonçons quelques résultats,
principalement dus à Cohen. Pour les démonstrations, nous renvoyons à \accourt{ix}{2} et \egazero{iv}{§19}.

\begin{définition2}[\egazero{iv}{19.3.1}]\label{4-19.3.1}
Soit $A$ un anneau topologique. Une $A$-algèbre topologique $B$
est dite \emph{formellement lisse} si pour toute $A$-algèbre topologique
\emph{discrète} $C$, et tout idéal \emph{nilpotent} $I$ de $C$,
tout $A$-morphisme continu $u:B→C/I$ se factorise en
$B\dessusdessous{v}{→}C\dessusdessous{φ}{↠}C/I$,
où $v$ est un $A$-morphisme continu et $φ$ l'homomorphisme canonique.
\end{définition2}

On dit aussi que $A→B$ est un \emph{morphisme formellement lisse}.
La proposition suivante énonce une propriété de relèvement
un peu plus générale que celle de la définition.

\begin{proposition2}[\emph{loc. cit.}, 19.6.1]\label{4-19.3.10}
Soient $A$ un anneau topologique, et $B$ un $A$-algèbre formellement lisse.
Soient $C$ une $A$-algèbre topologique, $I$ un idéal de $C$ vérifiant
les conditions suivantes :
\begin{enumerate}
\item $C$ est métrisable et complet ; 
\item $I$ est fermé et la suite $(I^n)_{n∈𝐍}$ tend vers zéro.
\end{enumerate}
Alors, tout $A$-morphisme continu $u:B→C/I$ se factorise en
$B\dessusdessous{v}{→}C↠C/I$,
où $v$ est un $A$-morphisme continu.
\end{proposition2}

Les deux théorèmes suivants donnent deux critères
importants de lissité formelle.

\begin{théorème2}[\emph{loc. cit.}, 19.6.1]
Une extension de corps munis de la topologie discrète
est formellement lisse si et seulement si l'extension
est \emph{séparable}.
\end{théorème2}

\begin{théorème2}[\emph{loc. cit.}, 19.7.1]\label{4-19.7.1}
Soient $A,B$ deux anneaux locaux nœthériens, $\got{m},\got{n}$
leurs idéaux maximaux respectifs et $k=A/\got{m}$ le corps
résiduel de $A$. Munissons $A$ et $B$ respectivement
des topologies $\got{m}$-adique et $\got{n}$-adique.
Soit $φ:A→B$ un morphisme \emph{local}, et posons $B₀=B⊗_A k$.
Les propriétés suivantes sont équivalentes :
\begin{enumerate}
\item $B$ est une $A$-algèbre formellement lisse ; 
\item $B$ est un $A$-module plat, et $B₀$ munie
de la topologie quotient est une $k$-algèbre formellement lisse.
\end{enumerate}
\end{théorème2}

Le théorème suivant, joint au précédent, est à la base
de la démonstration de l'existence des anneaux de Cohen
définis ci-après.
\begin{quote}
\emph{Dans les énoncés qui suivent, les anneaux locaux sont munis de la topologie 
de l'idéal maximal.}
\end{quote}

\begin{théorème2}[\emph{loc. cit.}, 19.7.2]\label{4-19.7.2}
Soient $A$ un anneau local nœthérien, $I$ un idéal strict, $A₀=A/I$, $B₀$
un anneau local nœthérien complet, $A₀→B₀$ un morphisme local formellement
lisse. Il existe alors un anneau local nœthérien complet $B$, un morphisme local
$A→B$ faisant de $B$ un $A$-module \emph{plat}, et un $A₀$-isomorphisme
$u:B⊗_A A₀\iso B₀$.
\end{théorème2}

\begin{définition2}[\emph{loc. cit.}, 19.8.4 et 5]\label{4-19.8.4}
\label{4-définition anneau Cohen}
On appelle \emph{anneau de Cohen} un anneau qui est soit
un corps de caractéristique nulle, soit un anneau de valuation discrète
complet, de corps résiduel de caractéristique $p>0$ et d'idéal
maximal engendré par $p$.
\end{définition2}

\begin{théorème2}[Cohen, \emph{loc. cit.}, 19.8.6 et 21.5.3]\label{4-19.8.6}
\begin{enumerate}
\item \label{1} Soient $W$ un anneau de Cohen de corps résiduel $K$, $C$ un anneau local nœthérien
complet, et $I$ un idéal strict de $C$. Alors, tout morphisme local
$u:W→C/I$ se factorise en $W\dessusdessous{v}{→}C→C/I$, où $v$ est local.
De plus, la factorisation est \emph{unique} si et seulement si 
$Ω¹_K=0$.
\item \label{2} Soit $K$ un \emph{corps}. Il existe un anneau de Cohen $W$ de corps 
résiduel isomorphe à $K$. Si $W′$ est un second anneau de Cohen, de corps 
résiduel $K′$, tout isomorphisme $u:K\iso K′$ provient par passage au
quotient d'un isomorphisme $v:W\iso W′$.
\end{enumerate}
\end{théorème2}

\subsubsection{}Rappelons que l'hypothèse $Ω¹_K=0$ est équivalente au fait
que $K$ est \emph{parfait} s'il est de caractéristique \mbox{$>0$} ou bien
est une extension algébrique de $𝐐$ s'il est de caractéristique nulle.

Signalons que si $K$ est \emph{parfait} de caractéristique $p>0$, le morphisme
$v$ de (\ref{2}) est \emph{unique}. Dans ce cas, $W$ est d'ailleurs 
isomorphe à l'anneau des vecteurs de Witt sur $K$.

\subsection{Le théorème de Cohen-Gabber en caractéristique mixte}

\subsubsection{}Soit $A$ un anneau local nœthérien complet de caractéristique résiduelle $p>0$.
Le schéma $X=\Spec(A)$ est de manière unique un $\Spec(𝐙_p)$-schéma. Notons
$X_p$ le sous-schéma fermé de $X$, fibre au-dessus du point fermé de $\Spec(𝐙_p)$.
Nous dirons qu'un ouvert $U⊂X$ est \emph{$p$-dense} si $U∩X_p$ est 
dense dans $X_p$.

\begin{théorème2}\label{4-Cohen-Gabber-Epp}
Soit $X=\Spec(A)$ un schéma local nœthérien complet normal de corps résiduel
$k$, de dimension $d≥2$ et de point générique de caractéristique nulle.
Il existe un morphisme fini surjectif $X'→X$, où $X'$
est normal intègre de corps résiduel $k'$, et
un morphisme fini surjectif $X'→\Spec(V[[t₁,\dots,t_{d-1}]])$,
où $V$ est un anneau de valuation discrète de corps résiduel $k'$,
\emph{étale} au-dessus d'un ouvert $p$-dense du but.
\end{théorème2}

La suite de ce paragraphe est consacrée à la démonstration
du théorème précédent.

\subsubsection{}Soit $X$ comme dans l'énoncé. Considérons 
le sous-corps parfait maximal $k₀=k^{p^∞}$ du corps résiduel $k$ de
$A=Γ(X,𝒪_X)$ et notons $W₀=W(k₀)$
l'anneau des vecteurs de Witt correspondant.
Il résulte du théorème de Cohen qu'il existe un unique morphisme
$X→S₀=\Spec(W₀)$ qui étende le morphisme 
$\Spec(k)→\Spec(k₀)$ entre les points fermés (\ref{4-19.8.6}, (i)).

Pour tout point maximal $\got{p}$ de la fibre spéciale $X_p$
de ce morphisme, l'anneau de valuation discrète $A_{\got{p}}$ a pour corps
résiduel $\Frac{\,A/\got{p}}$, où l'anneau $A/\got{p}$ est local
nœthérien complet intègre de corps résiduel $k$. D'après
\ref{4-prop Epp} (i) \& (iv), l'extension $\Frac(A/\got{p})/k₀$ a la propriété
de Epp. De tels idéaux $\got{p}$ étant
en nombre fini et la conclusion du théorème de Epp
(\ref{4-Epp}, \ref{4-remarque Epp}) étant
stable par changement de base fini car c'est un résultat
de lissité formelle, il existe donc un changement de base fini $S₀'=\Spec(W₀')→S₀$ 
tel que la fibre spéciale du produit fibré normalisé $X₀':=(X\times_{S₀}
S₀')_{\norm}=\Spec(A₀')$ soit réduite en ses points maximaux. (On utilise le fait que les points maximaux
de la fibre spéciale de $X₀'→ S₀'$ se trouvent au-dessus des points maximaux de la fibre
spéciale de $X→S₀$ ;  cf. p. ex. \egazero{iv}{16.1.6}.)

D'après le lemme suivant, la fibre spéciale du morphisme
$X₀'→ S₀'$ est alors réduite. 

\begin{lemme2}
Soit $X$ un schéma nœthérien normal. Tout diviseur de Cartier effectif 
génériquement réduit est réduit.
\end{lemme2}

\begin{démo}
On peut supposer $X$ affine égal à $\Spec(A)$ et le
diviseur de Cartier effectif défini par une
fonction $f∈A$. Soient $a∈A$ et $n≥1$ tels que 
$a^n∈fA$ ;  on souhaite montrer que $a∈fA$. 
L'anneau $A/f$ étant génériquement réduit,
l'élément $a/f$ de $\Frac{\,A}$ appartient à $A_𝔭$ 
pour tout idéal premier $𝔭$ contenant $f$ (un tel
idéal est de hauteur $1$ s'il est minimal
(\ac{viii}{3}{1}{\,prop. 1})). Il en est
évidemment de même pour $f∉𝔭$. L'anneau
$A$ étant normal, $⋂_𝔭 A_𝔭=A$ où $𝔭$ parcourt les idéaux de hauteur
$1$ (\emph{cf. loc. cit.} \textsc{vii}, §1, n°7, th. 4) 
de sorte que $a/f∈A$. CQFD.
\end{démo}

\subsubsection{}\label{4-fibration mixte} 
Notons $k₀'$ le corps résiduel de $W₀'$ 
, $\varpi'$ une uniformisante de $W₀'$, et considérons
une composante connexe $X'=\Spec(A')$ de $X₀'$ ; c'est un schéma fini
surjectif au-dessus de $X$. 
Soit $k'$ son corps résiduel.
L'inclusion $k₀'\inj k'$ déduite du morphisme $X'→S₀'$
est formellement lisse, car $k₀'$ est parfait, donc
se relève d'après \ref{4-19.7.1} et \ref{4-19.7.2} en
un morphisme \emph{formellement lisse} $W′₀→V$ où 
$V$ est un anneau local complet nœthérien. Cet anneau
est un anneau de valuation discrète.
L'anneau $A'/ϖ'$ étant \emph{réduit}, équidimensionnel de dimension $d-1$, de corps
résiduel $k'$, il existe d'après le théorème de
Cohen-Gabber (\ref{4-Cohen-Gabber égale}),
un relèvement $k₀'$-linéaire $k'↪ A'/ϖ'$ et 
des éléments $x_1,\dots,x_{d-1}$ dans l'idéal maximal de $A'/ϖ'$ 
tels que le morphisme induit
$k''[[t_1,\dots,t_{d-1}]]\ra A''/ϖ'$, envoyant
l'indéterminée $t_i$ sur $x_i$, soit
fini, \emph{génériquement étale} en haut et en bas.

Par lissité formelle de $W₀'→V$, le morphisme composé
$V→k'→A'/ϖ'$ se relève en un $W₀'$-morphisme $V→A'$.
En relevant les $x_i$ dans $A'$,
nous obtenons un morphisme $V[[t_1,\dots,t_{d-1}]]→A'$,
fini injectif (cf. p. ex. \egazero{iv}{19.8.8 (démonstration)}), 
\emph{étale} au-dessus du point générique de la fibre spéciale. CQFD.

\subsection{Le théorème de Cohen-Gabber premier à $ℓ$ en caractéristique mixte}

\begin{théorème2}\label{4-Cohen-Gabber-Epp équivariant}
Soit $X=\Spec(A)$ un schéma local nœthérien complet normal de
dimension $d≥2$, de corps résiduel $k$
de caractéristique $p>0$ et de point générique de caractéristique
nulle. Soit $ℓ$ un nombre premier différent de $p$.
Il existe alors :
\begin{enumerate}
\item un schéma local nœthérien intègre normal $Y$ muni d'une action
d'un $ℓ$-groupe fini $H$ et un morphisme fini surjectif
$H$-équivariant $Y → X$ tel que le quotient $Y ∕ H$ soit
de degré (générique) premier à $ℓ$ sur $X$ ;
\item un anneau de valuation discrète complet $V$ de même corps
résiduel $k ′$ que $Y$, de caractéristique mixte, muni d'une action de $H$
compatible avec son action sur $k ′$ ;
\item un morphisme local $Y→Y'=\Spec(V[[t_1,\dots,t_{d-1}]])$
qui soit fini, étale au-dessus d'un ouvert $p$-dense de $Y'$,
et $H$-équivariant avec action triviale de $H$ sur les $t_i$.
\end{enumerate}
\end{théorème2}

Ces morphismes sont représentés dans le diagramme ci-dessous,
où toutes les flèches sont des morphismes finis surjectifs.

$$
\xymatrix@C=3cm{
\Spec(V[[t_1,\dots,t_{d-1}]])=Y' & Y \ar[l]_(0.3){\text{$p$-génériquement étale}} \ar[d]
\ar[r] & Y/H \ar[dl]^{\textrm{ordre premier à } \ell} \\
& X & 
}
$$

\begin{remarque2}
Observons que les conditions (i) — (iii)
sur les morphismes $Y → X$ et $Y → Y ′$
n'entraînent pas que le schéma $Y/H$
soit étale au-dessus d'un ouvert $p$-dense de
$\Spec(\Fix_H(V)[[t_1,\dots,t_{d-1}]])=Y'/H$.
Voici un exemple, dû à Takeshi Saitô.
Soient $k$ un corps algébriquement clos de caractéristique $p>0$,
$W$ l'anneau des vecteurs de Witt sur $k$, $ℓ$ un nombre premier différent
de $p$, $A=W[[x, y]]/(x^ℓ y − p)$. Soient $W′=W[π]/(π^ℓ-p)$ et $B$ le
normalisé de $A⊗_W W′$, $W ′$-isomorphe à $W′[[x,z]]/(xz-π)$. Le groupe $H = μ_ℓ(k)$ agit sur $B$,
via son action sur $W′$ : $ζ\cdot x=x$ et $ζ\cdot z=ζz$.
Le morphisme $Y=\Spec(B)→X=\Spec(A)$ défini par $x↦x$ et $y↦z^ℓ$ satisfait les propriétés du théorème car
$Y/H⥲X$ et $Y→Y'=\Spec(W′[[x]])$, $x↦x$, est $p$-génériquement étale.
Cependant, $Y/H$ a une fibre spéciale isomorphe au schéma
non réduit $\Spec\big( k[[x,y]]/(x^ℓ\cdot y)\big)$.
Elle n'est donc pas étale au-dessus d'un
ouvert dense de la fibre spéciale de $Y ′$.
\end{remarque2}

La suite de ce paragraphe est consacrée à la démonstration du théorème
précédent. Notons que si la fibre spéciale $X_p$ de $X$ sur $\Spec(𝐙_p)$ est réduite,
ce théorème — comme le précédent — résulte simplement du
théorème \ref{4-Cohen-Gabber équivariant}, dans le cas particulier
où le groupe $G$ est trivial : on peut prendre $Y=X$ et $H$ trivial.

\subsubsection{}Considérons à nouveau le sous-corps parfait maximal $k₀$
du corps résiduel $k$ de $A$ et $W₀=W(k₀)↪A$ l'unique morphisme
relevant l'inclusion $k₀↪k$.
Soit $W₀^ν$ la clôture intégrale de $W₀$ dans $A$.

\begin{lemme2}
L'extension $W₀^ν/W₀$ est \emph{finie}, totalement ramifiée.
\end{lemme2}

\begin{démo}
Soit $W′/W₀$ une extension finie de traits, où $W′$ est
contenu dans $A$. Le corps résiduel de $W′$ est une
extension finie de $k₀$ ; c'est donc un corps parfait,
contenu dans $k$ et contenant $k₀$. Il est
donc égal à $k₀$ : l'extension est totalement ramifiée.
Le degré de l'extension $W′/W₀$ est par conséquent égal à son indice de
ramification, qui est \emph{majoré} par tout entier $N$ tel
que $p$ appartienne à $𝔪_A^N-𝔪_A^{N+1}$.
Si $W″$ est tel que le degré de l'extension $\Frac(W″) ∕ \Frac(W₀)$
soit maximal, on a nécessairement $W′ ⊆ W″$,
comme on le voit immédiatement en considérant
la sous-extension composée, dans $\Frac(A)$, des corps des fractions.
Ainsi, $W₀^ν=W″$ est fini sur $W₀$.
\end{démo}

\subsubsection{} D'après le théorème de Epp (\ref{4-Epp}),
il existe une extension finie d'anneaux de valuation
discrète $W₀^ν→W₀'$, que l'on peut supposer génériquement galoisienne
de groupe un groupe fini $G$, telle que la fibre spéciale sur $W₀'$ de la
normalisation $A'$ de $A⊗_{W₀^ν} W₀'$ soit \emph{réduite}. 
Observons que l'anneau $W₀^ν$ étant intégralement clos dans $A$,
l'anneau $A'$ est local.
Notons $k₀^ν$ (resp. $k₀'$) le corps résiduel de $W₀^ν$ 
(resp. $W₀'$) et $k'$ le corps résiduel de $A'$.
Choisissons des anneaux de Cohen $I(k')$ et $I({k'}^G)$ relatifs aux corps $k'$ et $k'^G$.
Il existe un morphisme $I({k'}^G)→I(k')$ relevant l'inclusion.
Ce morphisme étant fini étale entre anneaux locaux complets donc
henséliens, l'action du quotient $\Gal(k'/{k'}^G)$ de $G$ sur $k'$ 
se relève en une action $I({k'}^G)$-linéaire sur $I(k')$
(cf. p. ex. \cite{CL@Serre}, III, §5, th. 3).
Le corps $k₀'$ étant parfait, il existe d'après le théorème \ref{4-19.8.6}
un morphisme \emph{$G$-équivariant} $W(k₀')→I(k')$. 
Soient enfin $V=W₀'⊗_{W(k₀')} I(k')$, $ϖ$ une uniformisante de $W₀'$,
$\sur{A'}=A'/ϖA'$ et $H$ un $ℓ$-Sylow de $G$. D'après le théorème 
\ref{4-Cohen-Gabber équivariant}, il existe un morphisme \emph{fini, étale,
$H$-équivariant}, $φ:k'[[t₁,\dots,t_{d}]]→ \sur{A'}$, où les $t_i$
s'envoient dans $\sur{A'}^H$. Le morphisme ${A'}^H→ \sur{A'}^H$ étant
surjectif — comme cela se voit en utilisant la trace — on peut relever les images des $t_i$ 
en des $x'_i$ dans ${A'}^H$. De plus, par lissité formelle de $V/W₀'$, on peut relever
$k'→ \sur{A'}$ en un $W₀'$-morphisme $ψ:V→A'$ : cela résulte par exemple
de \egazero{iv}{19.3.10}. En procédant cran par cran,
et en considérant des isobarycentres dans les espaces affines
définis par le lemme bien connu suivant, on
constate qu'existe même un tel relèvement qui est \emph{$H$-invariant}.
\begin{lemme2}
Soit $(A,I)→ (B,J)$ un morphisme formellement lisse, où $A$ (resp. $B$) est
muni de la topologie définie par son idéal $I$ (resp. $J$). 
Soit $C↠C'$ une surjection de $A$-algèbres, de noyau $\mc{N}$ tel que
$\mc{N}²=I\mc{N}=0$. Alors, l'ensemble des $A$-linéaires relèvements d'un morphisme
$B→ C'$ à $C$ est soit vide soit un torseur sous $\mathrm{D\acute{e}r}_A (B,\mc{N})$.
\end{lemme2}

\subsubsection{}Le $X$-schéma $Y=\Spec(A')$ est bien fini $p$-génériquement étale 
sur $Y'=\Spec(V[[t₁,\dots,t_{d}]])$ si l'on envoie $V$ dans $A'$ par
$ψ$ comme ci-dessus et les variables $t_i$ sur les $x'_i$.
 Par construction $Y$ est, génériquement sur $X$, galoisien de groupe un sous-groupe de $G$ ;
son quotient $Y/H$ est donc génériquement d'ordre premier à $ℓ$.
Ceci achève la démonstration du théorème.

\ifx\danslelivre\undefined
\bibliography{../fichiers-configuration-latex/bibliographie}
\bibliographystyle{../fichiers-configuration-latex/nomcomplet-date}


\begin{center}
V. Algébrisation partielle\\
Fabrice Orgogozo
\end{center}
\version

\else
\setboolean{english}{false}
\setboolean{francais}{true}
\selectlanguage{francais}
\chapter{Algébrisation partielle}
\addtocontents{toc}{par Fabrice Orgogozo}
\begin{center}
Fabrice Orgogozo
\end{center}
\fi

\section{Préparatifs (rappels)}

\subsection{Le théorème de préparation de Weierstraß}
On trouvera dans \ac{vii}{3}{7-8}{} une démonstration du théorème
suivant.

\begin{théorème2}\label{5-Weierstrass}
Soient $A$ un anneau local séparé complet d'idéal maximal $𝔪$,
$d>0$ un entier et $f∈A[[\sous{X},T]]$ une série formelle,
où l'on pose $\sous{X}=(X₁,\dots,X_d)$.
\begin{enumerate}
\item Soit $ρ$ un entier naturel tel que
$f$ soit \emph{$ρ$-régulière} relativement à $T$, c'est-à-dire congrue à $(u∈A[[T]]^×)\cdot
T^ρ$ modulo $(𝔪,\sous{X})$. Alors, pour tout $g∈A[[\sous{X},T]]$, il existe un unique
couple $(q,r)∈A[[\sous{X},T]]×A[[\sous{X}]][T]$ tel que
$g=qf+r$ et $\deg_T(r)<ρ$. De plus, il existe
un unique polynôme $P=T^ρ+∑_{i<ρ} p_i T^i$, où
les coefficients $p_i$ appartiennent à $(𝔪,\sous{X})A[[\sous{X}]]$,
et une unité $u∈A[[\sous{X},T]]^×$ tels que $f=uP$.
\item Si $f$ est \emph{non nulle modulo $𝔪$},
il existe un entier naturel $ρ$ et un automorphisme
$A[[T]]$-linéaire $c$ de $A[[\sous{X},T]]$, tel que
$c(X_i)=X_i+T^{N_i}$ ($N_i>0$) et la série
entière $c(f)$ soit $ρ$-régulière.
\end{enumerate}
\end{théorème2}

\subsubsection{}Signalons que l'on peut satisfaire la condition (ii)
simultanément pour un nombre \emph{fini} d'éléments :
cf. \emph{loc. cit.}, n°7, lemme 2 où l'on considérera
un produit (fini) de séries formelles. 

Nous ferons usage de la propriété suivante des
polynômes comme en (ii) ci-dessus. 

\begin{lemme2}\label{5-lemme complétion}
Soient $B$ un anneau local complet nœthérien et $P∈B[X]$ un polynôme
de la forme $X^ρ+∑_{i<ρ}b_i X^i$, où $b_i∈\MM_{B}$ et
$ρ>0$. Alors, le complété $(P)$-adique de $B\{X\}$ s'identifie à
$B[[X]]$.
\end{lemme2}

Rappelons que $B\{X\}$ désigne l'hensélisé en l'origine de l'anneau $B[X]$.
Un polynôme $P$ comme ci-dessus est parfois dit \emph{de Weierstraß}.

\begin{démo}
Soient $N$ un entier naturel et $Q=P^N$.
Il résulte de \ref{5-Weierstrass} (i), que l'anneau quotient
$B[[X]]/(Q)$ est isomorphe comme $B$-module à $B[X]/(X^{\deg(Q)})$
et en particulier fini sur $B$. 
Par fidèle platitude du morphisme $B\{X\}→B[[X]]$, on a
$QB[[X]]∩B\{X\}=QB\{X\}$ de sorte que le $B$-morphisme
$B\{X\}/Q→B[[X]]/Q$ est \emph{injectif}. L'anneau $B\{X\}/Q$
est donc également fini sur l'anneau complet $B$ ;
il est donc isomorphe à son complété $B[[X]]/Q$. En faisant tendre
$N$ vers l'infini, on en déduit que le séparé-complété $(P)$-adique de $B\{X\}$ est 
donc isomorphe à celui de $B[[X]]$ ; ce dernier est isomorphe à $B[[X]]$ 
puisque $\deg(P)>0$.
\end{démo}

\subsection{Le théorème d'algébrisation d'Elkik}

\begin{définition2}
Une paire $(\Spec(C),V(J))$ est dite \emph{hensélienne}
si pour tout polynôme $f∈C[T]$, toute racine $β$ de $f$
dans $C/J$ telle que $f'(β)$ soit une unité de $C/J$
se relève en une racine dans $C$.
\end{définition2}

\subsubsection{}Il résulte de \cite{Caracterisation@Crepeaux}, prop. 2, (iii)⇒(i)
et \cite{Anneaux@Raynaud}, chap. XI, §2, prop. 1, 4)⇒1)
(voir aussi \cite{K-theory@Gabber}, déf. p. 59) que la
définition précédente est équivalente à la définition
\ega{iv}{18.5.5} (relèvement d'ouverts-fermés).
Remarquons que la notion de \emph{racine simple} introduite
dans la définition est plus forte que celle de \a{iv}{2}{1}{, déf. 1} et que le relèvement
ci-dessus est nécessairement unique. Notons également
que l'idéal $I$ est nécessairement contenu
dans le radical de Jacobson de l'anneau $C$.

\begin{lemme2}\label{5-paire hensélienne}
Soient $C$ un anneau \emph{local} hensélien d'idéal maximal
$𝔪$, et $J⊆𝔪$ un idéal. La paire $(\Spec(C),V(J))$ est hensélienne.
\end{lemme2}

En particulier, pour $B$ et $P$ comme dans le lemme \ref{5-lemme
complétion}, la paire $(\Spec(B\{X\}),V(P))$  est hensélienne.

\begin{démo}
Soient $f$ et $β$ comme ci-dessus. L'anneau $C$ étant local hensélien,
l'image $γ$ de $β$ dans le corps résiduel $C/𝔪$ se relève en
une racine $α$ de $P$. Notons $β'$ son image dans $C/J$ et vérifions que $β=β'$.
Remarquons tout d'abord que puisque $P'(α)$ est une unité de $C$,
$P'(β')$ est une unité de $C/J$. De plus, l'égalité
$P(β)=P(β')+(β-β')P'(β')+(β-β')²b$ où $b∈B/J$ se réduit
à $β-β'=(β-β')²\frac{-b}{P'(β')}$ ; si l'on pose $x=β-β'$, on a donc $x(1-ax)=0$
pour un $a∈C/J$. Comme $x$ appartient à $𝔪$ (car $β$ et $β'$ ont
pour image $γ$ dans $C/𝔪$), on a $x=0$.
\end{démo}

Terminons ces rappels par l'énoncé du théorème d'algébrisation de Renée Elkik (\cite{Solutions@Elkik}, théorème 5).

\begin{théorème2}\label{5-Elkik}
Soient $(X=\Spec(A),Y=V(I))$ une paire hensélienne avec $A$ nœthérien,
et $U$ le sous-schéma ouvert complémentaire de $Y$ dans $X$.
Notons $X_{\chap{Y}}$ le complété de $X$ le long de $Y$, $\chap{Y}$
le fermé correspondant à $Y$ et $\chap{U}$ son complémentaire dans 
$X_{\chap{Y}}$.
Le foncteur $X'↦X'×_X X_{\chap{Y}}$ induit une
équivalence de catégories entre la catégorie des
$X$-schémas finis, étales sur $U$, et la catégorie
des $X_{\chap{Y}}$-schémas finis, étales sur $\chap{U}$. 
\end{théorème2}

\section{Algébrisation partielle en égale caractéristique}\label{5-égale caractéristique}

\subsection{Énoncé}

\subsubsection{}\label{5-déf algébrisation partielle}
Soient $A$ un anneau local nœthérien complet et $\{I_e\}_{e∈E}$ une
collection d'idéaux de $A$. On dit que la paire $(A,\{I_e\}_{e∈E})$ est \emph{partiellement algébrisable} 
s'il existe un anneau local nœthérien complet $B$ de
dimension strictement inférieure à celle de $A$, une
$B$-algèbre de type fini $C$, un idéal maximal $𝔫$ au-dessus
de l'idéal maximal de $B$, et un isomorphisme
$A\isononcan \chap{C_{𝔫}}$ tel que les idéaux $I_e$ ($e∈E$) proviennent
d'idéaux de $C_{\got{n}}$.

\begin{théorème2}\label{5-algébrisation partielle égale caractéristique}
Soit $A$ un anneau local nœthérien complet réduit d'égale caractéristique
qui ne soit pas un corps. Alors, $A$ muni d'un ensemble \emph{fini} quelconque
d'idéaux est partiellement algébrisable.
\end{théorème2}

\subsection{Démonstration}
\subsubsection{}Soient $X=\Spec(A)$ et $I_1,\dots,I_n⊆A$ comme dans
l'énoncé. Il résulte de la définition \ref{5-déf algébrisation partielle}
que si un idéal $I$ de $A$ est de la forme $J₁∩\cdots∩J_r$ 
et que les $J_i$ sont partiellement algébrisables, 
l'idéal $I$ l'est également. D'après le théorème de
décomposition primaires des idéaux, on peut supposer les $I_i$
\emph{primaires}.

\subsubsection{}Notons $k$ le corps résiduel de $A$ et $d>0$ sa dimension. 
D'après \refext{IV}{4-Cohen-Gabber égale} si $X$ est équidimensionnel
ou bien \refext{IV}{4-Cohen-Gabber équivariant} (avec $G=\{1\}$) dans
le cas général, il existe un morphisme fini génériquement
étale $π:X→X₀$, où $X₀=\Spec(k[[t₁,\dots,t_d]])$.

\subsubsection{}Soit $I$ l'un des $I_i$. Deux cas se présentent.
\begin{enumerate}
\item $\dim(A/I)=d$. L'idéal $I$ est donc un idéal premier
minimal de $A$. 
\item $\dim(A/I)<d$. L'image de $V(I)$ dans $X₀$
est donc de dimension au plus $d-1$ donc contenue dans un fermé
$V(g_I)$ où $g_I∈A₀=Γ(X₀,𝒪_{X₀})-\{0\}$.
\end{enumerate}

Soient $g=∏_{I_i} g_{I_i}$ où $I_i∈\{I_1,\dots,I_n\}$ 
parcourt le sous-ensemble des idéaux du second type,
et $f∈A₀-\{0\}$ telle que le lieu de ramification de
$π$ soit contenu dans $V(f)$. Posons $h=gf$.
D'après \ref{5-Weierstrass} (ii) et (i),
quitte à changer de base par un automorphisme (c'est-à-dire changer les
coordonnées), on peut supposer que $h$ est un polynôme
unitaire en $t_d$. Considérons le sous-anneau $\gtilde{A₀}=k[[t₁,\dots,t_{d-1}]]\{t_d\}$
de $A₀$. Il est hensélien et contient $h$. 
D'après les lemmes \ref{5-lemme complétion} et \ref{5-paire hensélienne}, 
la paire $(\gtilde{X₀}=\Spec(\gtilde{A₀}),V(h))$ est hensélienne. On est donc en
mesure d'appliquer le théorème \ref{5-Elkik} et d'en déduire qu'il
existe un diagramme cartésien :
$$
\xymatrix{X \ar[d] \ar[r] & \gtilde{X} \ar[d] \\
X₀ \ar[r] & \gtilde{X₀}
}
$$
où la flèche verticale de gauche est, par hypothèse, étale hors de 
$V(h)$ et les flèches horizontales sont des morphismes de complétion
(à la fois pour la topologie $h$-adique et celle définie par leurs
idéaux maximaux respectifs).

Les idéaux $I$ du premier type (c'est-à-dire premier minimaux)
se descendent à $\gtilde{X}$ d'après le lemme suivant.

\begin{lemme2}\label{5-premier minimal}
Soit $B$ un anneau local hensélien quasi-excellent de complété
noté $\chap{B}$. Tout idéal premier minimal de $\chap{B}$ provient par
image inverse d'un idéal premier minimal de $B$.  
\end{lemme2}

\begin{démo}
Par restriction à l'adhérence du point générique du fermé,
il suffit de démontrer que le complété d'un intègre hensélien
quasi-excellent est intègre. Ce fait est bien connu et résulte
d'ailleurs immédiatement du théorème d'approximation de
Popescu, appliqué à l'équation \mbox{$xy=0$}.
\end{démo}

\subsubsection{}Quant aux idéaux $I$ du second type, il suffit d'observer que 
chaque $V(I)$ est fini sur $\Spec(k[[t₁,\dots,t_{d-1}]])$,
donc sur $\gtilde{X}$, et d'appliquer 
le

\begin{lemme2}\label{5-fini par complétion}
Soient $B$ un anneau local nœthérien, $J⊆𝔪_B$ un idéal, et 
$\chap{B}$ le complété $J$-adique de $B$.
Tout quotient de $\chap{B}$ fini sur $B$ se descend à $B$.
\end{lemme2}

\subsubsection{}Admettons momentanément ce lemme et achevons la démonstration 
de \ref{5-algébrisation partielle égale caractéristique}. Comme on l'a vu,
les idéaux $I_1,\dots,I_n$ proviennent d'idéaux de $\gtilde{A}$.
Cet anneau est fini — donc \emph{a fortiori} de présentation finie
par nœthérianité — sur l'anneau $k[[t₁,\dots,t_{d-1}]]\{t_{d}\}$. Ce dernier est
l'hensélisé en l'origine de $k[[t₁,\dots,t_{d-1}]][t_d]$ ;
il est isomorphe à la colimite filtrante d'anneaux de type fini sur
$k[[t₁,\dots,t_{d-1}]][t_d]$ donc sur l'anneau $B=k[[t₁,\dots,t_{d-1}]]$.
La conclusion résulte de \ega{iv}{8.8.2} qui assure l'existence d'une
$B$-algèbre $C$ comme en \ref{5-déf algébrisation partielle}
dont proviennent les idéaux $I_i$.

\subsubsection{}Revenons à la démonstration du lemme \ref{5-fini par
complétion} ci-dessus.
Soit $I⊆\chap{B}$ tel que $\chap{B}/I$ soit fini sur $B$.
Quitte à remplacer $B$ par $B/\Ker(B→\chap{B}/I)$,
c'est-à-dire $\Spec(B)$ par l'image schématique de $V(I)$, on peut supposer
$B→\chap{B}/I$ injectif, c'est-à-dire $V(I)→\Spec(B)$ schématiquement
dominant. Le $B$-module $\chap{B}/I$ étant fini, la topologie
$J$-adique sur $\chap{B}/I$ induit la topologie $J$-adique sur $B$.
Puisque l'application $B→\chap{B}/I$ est injective, d'image dense,
et continue, il en résulte que $\chap{B}/I$ est le
\emph{séparé-complété} de $B$ pour la topologie $J$-adique. On a donc
$I=(0)$ ; il se descend tautologiquement à $B$.

\section{Algébrisation partielle première à $ℓ$ en caractéristique mixte}

\subsection{Énoncé}

\begin{théorème2}
\label{5-algébrisation partielle mixte raffinée}
Soient $A$ un anneau local nœthérien complet normal de caractéristique mixte
$(0,p)$ de dimension $d≥2$ et $ℓ≠p$ un nombre premier. Il existe
un morphisme injectif fini $A→A'$ de degré générique premier à $ℓ$,
où $A'$ est un anneau normal intègre dont toute famille \emph{finie}
d'idéaux est partiellement algébrisable.
\end{théorème2}

\begin{remarque2}Signalons qu'il suffit pour démontrer \refext{XIII}{13-constructibilité}
d'établir la variante affaiblie de l'énoncé précédent selon laquelle
— en reprenant les notations de \ref{5-déf algébrisation partielle} —
tout fermé rare de $\Spec(A')≃\Spec(\chap{C_{𝔫}})$
est ensemblistement contenu dans l'image inverse d'un \emph{diviseur} de
$\Spec(C_𝔫)$.
\end{remarque2}

\subsubsection{Reformulation géométrique}
\label{5-algébrisation partielle mixte raffinée géométrique}
Soit $X$ un schéma local nœthérien complet de caractéristique mixte,
de dimension $d ≥ 2$ et soit $\{Z_i\}_{i ∈ I}$ une famille finie de
fermés de $X$. Pour tout nombre premier $ℓ$ inversible sur $X$,
il existe un diagramme
\[
\xymatrix{X & X ′ \ar[r]^{a} \ar[l]_{π} & Y \ar[d]^f  \\
& & S }
\]
où :
\begin{itemize}
\item $S$ est un schéma nœthérien régulier complet de caractéristique
mixte et de dimension $d-1$ ;
\item $X ′$ est un schéma local normal ;
\item $π$ est un morphisme fini de degré générique premier à $ℓ$ ;
\item $f$ est un morphisme de type fini ;
\item $a$ induit un isomorphisme entre $X ′$ et le complété de $Y$ en
un point fermé de la fibre spéciale de $f$,
\end{itemize}
et des fermés $F_i$ de $Y$ tels que $π^{-1}(Z_i)=a^{-1}(F_i)$ pour
tout $i ∈ I$.

\begin{remarque2}
Il découle de \ref{5-algébrisation partielle égale caractéristique}
que le résultat précédent est également vrai en égale caractéristique,
et que l'on peut alors supposer $X=X′$.
\end{remarque2}

\subsection{Démonstration}

\subsubsection{}Soient $X=\Spec(A)$ de dimension $d≥2$ et $ℓ$ comme
dans l'énoncé. D'après le théorème \refext{IV}{4-Cohen-Gabber-Epp équivariant},
il existe un diagramme commutatif
$$
\xymatrix{
V[[t_1,\dots,t_{d-1}]]=B₀ \ar[r] & B  & A'=\Fix_H(B) \ar[l]\\
& A \ar[u] \ar[ur] &
}
$$
où $V$ est le spectre d'un anneau de valuation discrète complet
d'idéal maximal $𝔪_V$, et $H$ est un $ℓ$-groupe agissant sur l'anneau
normal $B$, son sous-anneau $V$, et trivialement sur les variables $t_i$
($1≤i≤d-1$). De plus, $A→A'$ est une injection finie de rang générique premier à $ℓ$
et $π:\Spec(B)→\Spec(B₀)$ est fini, $p$-génériquement étale.

\subsubsection{}Nous allons montrer que toute famille finie d'idéaux de $A'$ est partiellement
algébrisable. Soit $I'₁,\cdots,I'_n$ une telle famille, que l'on peut 
supposer constituée d'idéaux primaires (cf. §\ref{5-égale caractéristique}).

\subsubsection{}Soit $I'$ l'un des $I_i$. Deux cas se présentent. 
\begin{enumerate}
\item $\dim(A'/I+p)=d-1$. Supposons $I'≠(0)$
et notons $𝔭$ l'idéal premier, nécessairement de hauteur un,
pour lequel $I'$ est primaire. Par hypothèse, l'idéal premier $𝔭$ contient $p$ ;
c'est un idéal premier minimal de $A'/p$. D'autre part, 
$A'$ est normal car $B$ l'est. Il résulte par exemple
de \cite{Algebre@Serre}, chap. \textsc{iii}, C, §1 que $I'$
est une \emph{puissance symbolique} de $𝔭$, c'est-à-dire
l'image inverse dans $A'$ d'une puissance de l'idéal (principal) $𝔭A'_𝔭$.
\item $\dim(A'/(I+p))<d-1$. L'image de $V(I)$ dans
$\Spec(V^H[[t₁,\dots,t_{d-1}]])$ est donc contenu dans un fermé
$V(g_{I'})$ où $g_{I'}∈V^H[[t₁,\dots,t_{d-1}]]$ est \emph{non nulle
modulo $𝔪_{V^H}$}.
\end{enumerate}

Soient $g=∏_{I'_i} g_{I'_i}$ où $I'_i∈\{I'_1,\dots,I'_n\}$ 
parcourt le sous-ensemble des idéaux du second type,
et $f∈V[[t₁,\dots,t_{d-1}]]-𝔪_{V}$ telle que le lieu de ramification de
$π$ soit contenu dans $V(f)$. Posons $h=gf$.
D'après le théorème de préparation (\ref{5-Weierstrass}) et l'observation lui faisant
immédiatement suite, on peut supposer que $h$ est un polynôme
de Weierstraß en $t_{d-1}$, \emph{$H$-équivariant}. (Rappelons que $H$ agit
trivialement sur les variables).
Comme en §\ref{5-égale caractéristique}, le morphisme $B₀→B$ se descend donc
d'après \ref{5-Elkik} en un morphisme
$\gtilde{B₀}→\gtilde{B}$, où $\gtilde{B₀}=V[[t_1,\dots,t_{d-2}]]\{t_{d-1}\}$. 
Le groupe $H$ préservant l'ouvert $D(f)$ de $\Spec(B₀)$, son action se
descend. Le diagramme ci-dessus se complète donc en un diagramme
$$
\xymatrix{
V[[t_1,\dots,t_{d-2}]]\{t_{d-1}\}=\gtilde{B₀} \ar[d] \ar[r] &
\gtilde{B} \ar[d] & \gtilde{A'}=\gtilde{B}^H \ar[l] \ar[d] \\
B₀ \ar[r] & B  & A' \ar[l]
}
$$
où les flèches verticales sont les morphismes de complétion
et les flèches horizontales sont finies.

\subsubsection{}Les idéaux $I'$ du second type se descendent de $A'$ à $\gtilde{A'}$
car $A'/I'$ est fini sur $V^H[[t_1,\dots,t_{d-2}]]$ donc \emph{a
fortiori} sur $\gtilde{A'}$ (cf. \ref{5-fini par complétion}).
Quant aux idéaux du premier type (puissances symboliques), il suffit d'appliquer
le lemme \ref{5-premier minimal} à la paire constituée de $\gtilde{A'}/p$ et de
son complété $A'/p$.
Comme en §\ref{5-égale caractéristique}, on utilise le fait que
$\gtilde{A'}$ soit fini — de type fini suffirait — sur $V^H[[t_1,\dots,t_{d-2}]]\{t_{d-2}\}$ pour
descendre, par passage à la limite, les idéaux à un anneau
de type fini sur $V^H[[t_1,\dots,t_{d-2}]]$.

\ifx\danslelivre\undefined
\bibliography{../fichiers-configuration-latex/bibliographie}
\bibliographystyle{../fichiers-configuration-latex/nomcomplet-date}


\begin{center}
VI Log régularité, actions très modérées\\
Luc Illusie
\end{center}

\else
\setboolean{english}{false}
\setboolean{francais}{true}
\selectlanguage{francais}
\chapter{Log régularité, actions très modérées}
\addtocontents{toc}{par Luc Illusie}
\begin{center}
Luc Illusie
\end{center}
\fi

\section{Log régularité} \subsection{}\label{6-1.1} Pour le langage des log schémas nous renvoyons
le lecteur à \cite{LSFI@Kato}, \cite{Toric@Niziol}, \cite{HT@GR}. Sauf mention du contraire, les log
structures considérées le seront au sens de la topologie étale. Un log schéma \textit{fin} (resp.
\textit{fs}, i. e. \textit{fin et saturé}) \cite{LSFI@Kato} est un schéma muni d'une log structure
admettant localement (pour la topologie étale) une carte sur un monoïde fin (resp. fin et saturé).
On note en général $M_X$ le faisceau de monoïdes d'un log schéma $X$, $\alpha : M_X \ra
\OO_X$ le morphisme structural, et $\overline{M}_X = M_X/\OO^*_X$. Sauf mention du
contraire, les log schémas considérés sont supposés localement noethériens.

Dans ce qui suit, $X$ désigne un fs log schéma.  \subsection{}\label{6-1.2} Soient $\overline{x}$ un
point géométrique de $X$, d'image $x \in X$, $\OO_{X,\overline{x}}$ le localisé strict de
$X$ en $\overline{x}$. Notons $I_{\overline{x}}$ l'idéal de $\OO_{X,\overline{x}}$ engendré
par $\alpha(M_{X,\overline{x}}-\OO_{X,\overline{x}}^*)$, $C_{X,\overline{x}}$ le sous-schéma
fermé de $X_{(\overline{x})} = \Spec \OO_{X,\overline{x}}$ défini par $I_{\overline{x}}$.
$C_{X,\overline{x}}$ est la trace sur $X_{(\overline{x})}$ de la strate de $X$ où le rang de
$\overline{M}^{\mathrm{gp}}$ est égal à $r(x) =
\mathrm{rg}(\overline{M}_{X,\overline{x}}^{\mathrm{gp}})$. 

On dit que $X$ est \textit{log régulier} en $x$ (ou $\overline{x}$) si $C_{X,\overline{x}}$ est
régulier et de codimension égale à $r(x)$ (cette condition ne dépend que de $x$). On dit que $X$ est
\textit{log régulier} si $X$ est log régulier en tout point. La définition analogue pour les log
schémas zariskiens est due à Kato \cite{Toric@Kato}. La variante dans le cadre étale a été traitée
par Niziol \cite{Toric@Niziol}. Nous rappelons ci-après quelques propriétés de cette notion.
\subsection{}\label{6-1.3} $X$ est log régulier en $x$ si et seulement si $X$, muni de la log
structure Zariski $M_X^{\mathrm{Zar}} := \varepsilon_*M_X$, où $\varepsilon : X_{\mathrm{\et}} \ra
X_{\mathrm{Zar}}$, est log régulier en $x$ au sens de Kato \cite{Toric@Kato},
(\cite[4.6]{Saturated@Tsuji}, \cite[2.4]{Toric@Niziol}. En particulier, si  $X$ est log régulier en
$x$, $X$ est log régulier en toute générisation $y$ de $X$ (\cite[7.1]{Toric@Kato}).

Si la log structure de $X$ est triviale, $X$ est log régulier si et seulement si $X$ est régulier au
sens usuel.  \subsection{}\label{6-1.4} Supposons $X$ log régulier. Soit $j : U \inj X$ l'inclusion
de l'ouvert de trivialité de sa log structure. Alors $U$ est un ouvert dense de $X$ et on a $$ M_X =
\OO_X \cap j_*\OO_U^* $$ (\cite[2.6]{Toric@Niziol}).

Nous dirons qu'un couple $(X,Z)$ formé d'un schéma $X$ et d'un fermé $Z$ est un couple \textit{log
régulier} si, pour la log structure sur $X$ définie par $M_X = \OO_X \cap
j_*\OO_U^*$, où $j : U \inj X$ est l'ouvert complémentaire de $Z$, $X$ est log régulier et
$Z$ est le complément de l'ouvert de trivialité de sa log structure. La log structure précédente sur
$X$ sera dite \textit{associée} au couple $(X,Z)$.  \subsection{}\label{6-1.5} Supposons $X$ log
régulier. Pour $i \in \NN$, soit $X^{(i)}$ l'ensemble des points $x$ de $X$ tels que $r(x) =
i$, avec la notation de \ref{6-1.2}. C'est une partie localement fermée, sous-jacente à un
sous-schéma régulier de $X$, de codimension $i$, dont la trace sur $X_{(\overline{x})}$, en chaque
point géométrique $\overline{x}$ localisé en  $x \in X^{(i)}$, est $C_{X,\overline{x}}$. On dit que
$X^{(i)}$ est la \textit{strate de codimension} $i$ définie par le rang de
$\overline{M}^{\mathrm{gp}}$. La stratification par les $X^{(i)}$ est appelée \textit{stratification
par le rang de $\overline{M}^{gp}$}, ou \textit{stratification canonique}. Voici deux exemples. 

(i) Si $X$ est un schéma noethérien régulier, muni de la log structure définie par un diviseur à
croisements normaux $D$, $X^{(i)}$ est l'ensemble des points où passent exactement $i$ branches de
$D$, i. e. tels que le normalisé de $D$ ait $i$ points au-dessus de $x$.

(ii) Si $X$ est une variété torique sur un corps $k$, de tore $T$, munie de sa log structure
canonique, $X$ est un log schéma log régulier, l'ouvert de trivialité de la log structure est $T$,
et $X^{(i)}$ est la réunion des orbites de $T$ de codimension $i$.  \subsection{}\label{6-str loc
log rég}\label{6-str loc log reg} Supposons $X$ log régulier en $x$, image du point géométrique
$\overline{x}$, soient $P = \overline{M}_{X,\overline{x}}$, $k = k(\overline{x)}$. Notons que $P$
est un monoïde fs \textit{saillant} (i. e. tel que $P^* = 0$). Soit $\widehat  X_{\overline{x}}$ le
complété de $X_{(\overline{x})}$ au point fermé. Alors, d'après \cite[3.2]{Toric@Kato}, $X$ admet
une carte modelée sur $P$ en $\overline{x}$, qui donne lieu à des isomorphismes

(i) $$ \widehat  X_{\overline{x}} \iso \Spec k[[P]][[t_1,\cdots,t_n]] $$ si $\OO_{X,x}$ est
d'égale caractéristique,

(ii) $$ \widehat  X_{\overline{x}} \iso \Spec C(k)[[P]][[t_1,\cdots,t_n]]/(f) $$ si
$\OO_{X,x}$ est d'inégale caractéristique $(0,p)$, où $C(k)$ est un anneau de Cohen pour k,
et $f$ est congru à $p$ modulo l'idéal engendré par $P- \{0\}$ et les $t_i$.
\subsection{}\label{6-1.7} Supposons $X$ log régulier. Alors $X$ est régulier en $x$, image de
$\overline{x}$, si et seulement si $\overline{M}_{\overline{x}} \isononcan \NN^r$
(\cite[5.2]{Toric@Niziol}, voir aussi \cite [1.8]{LR@Vidal}). Il en résulte que l'ensemble des
points de régularité de $X$ coïncide avec l'ensemble des points de régularité du monoïde
$\overline{M}_{\overline{x}}$, et en particulier est ouvert dans $X$ (\cite[5.3]{Toric@Niziol}). Si
$X$ est log régulier et régulier, l'ouvert de trivialité de la log structure est alors le complément
d'un diviseur à croisements normaux.  \subsection{} Soit $f : X \ra Y$ un morphisme de log schémas
fs. Si $Y$ est log régulier et $f$ log lisse, $X$ est log régulier.

L'analogue pour les log structures Zariski est \cite[8.2]{Toric@Kato}. La démonstration de
(\textit{loc. cit.}) s'applique, \textit{mutatis mutandis}, dans le cadre étale.

Le corollaire suivant jouera un rôle clé dans l'application des résultats de de Jong à la
démonstration du théorème d'uniformisation de Gabber. Rappelons que, si $S$ est un schéma, une
\textit{courbe nodale} $f : C \ra S$ est un morphisme propre et plat, purement de dimension relative
1, dont les fibres géométriques ont pour seules singularités des croisements normaux.
\begin{proposition} \label{6-nodal implique log lisse} Soit $(Y,T)$ un couple log régulier
(\ref{6-1.4}). Soit $f : X \ra Y$ une courbe nodale, lisse au-dessus de $Y-T$. Soit $D$ un diviseur
effectif sur $X$, de support contenu dans le lieu de lissité de $f$, et étale sur $Y$. Alors le
couple $(X,f^{-1}(T) \cup D)$ est log régulier, et pour les log structures associées, $f$ est un
morphisme de log schémas et est log lisse.  \end{proposition} La question est locale sur $X$ et
l'assertion est triviale sur l'ouvert de lissité de $f$ et sur $f^{-1}(Y-T)$. Soit $\overline{x}$ un
point géométrique de non lissité de $f$, d'image $\overline{y}$ dans $Y$. D'après \sga{7}{XV}{1.3.2}
(voir aussi \cite[2.23]{Alterations@de_Jong}, on a \begin{equation} \OO_{X,\overline{x}}
\simeq \OO_{Y,\overline{y}}\{u,v\}/(uv-h), \label{6-(1.9.1)} \end{equation} où
$\OO_{Y,\overline{y}}\{u,v\}$ désigne le localisé strict de $\Spec
\OO_{Y,\overline{y}}[u,v]$ au point $(u = v = 0)$ au-dessus du point fermé de
$Y_{(\overline{y})}$, et $h$ est un élément de $\OO_{Y,\overline{y}}$ inversible sur $Y-T$,
donc appartenant à $M_{Y,\overline{y}}$ (\ref{6-1.4}). Soit $c$ une carte locale de $Y$ en
$\overline{y}$ définie par une section de  $M_{Y,\overline{y}} \ra P :=
\overline{M}_{Y,\overline{y}}$. Écrivons $h = \varepsilon a$, avec $\varepsilon \in
\OO^*_{Y,\overline{y}}$ et $a \in P$. Soit $Q$ le monoïde fs défini par le carré cocartésien
$$ \xymatrix{\NN^{2} \ar[r] & Q \\ \NN \ar[u] \ar[r] & \ZZ \times P
\ar[u]^{g}}, $$ où la flèche $\NN \ra \NN^2$ (resp. $\NN \ra \ZZ \times
P $) est donnée par $1 \mapsto (1,1)$ (resp. $1 \mapsto (1,a)$). Le carré $$ \xymatrix{\Spec
\OO_{Y,\overline{y}}[u,v]/(uv-h) \ar[r]^-{d} \ar[d] & \Spec \ZZ[Q] \ar[d]^{\Spec
\ZZ[g]} \\ \Spec \OO_{Y,\overline{y}} \ar[r]^{c'} & \Spec \ZZ[\ZZ\times
P]}, $$ où la flèche $c'$ est donnée par $c$ et $\ZZ \ra \OO_{Y,\overline{y}}$, $1
\mapsto \varepsilon$, et $d$ par $(1,0) \mapsto u$, $(0,1) \mapsto v$ sur $\NN^2$ et $c'$ sur
$\ZZ \times P$, est cartésien. Comme $c'$ est une carte de $Y$ en $\overline{y}$, $d$ est une
carte en $\overline{x}$ de $X$ muni de la log structure image inverse par $f$, i. e. définie par
$\OO_X \cap \OO^*_{f^{-1}(Y-T)}$, et $(c',d,g)$ est une carte de $f$. Comme $g$ est
injectif et $\coker g^{gp} \simeq \ZZ$, $f$ est log lisse en $\overline{x}$, et le couple
$(X,f^{-1}(T) \cup D)$ est log régulier en $\overline{x}$.  \section{Revêtements Kummer étales} Dans
cette section et la suivante, nous aurons à considérer des action de groupes sur des schémas ou des
log schémas : sauf mention du contraire, il s'agira d'actions \textit{à droite}.
\subsection{}\label{6-Kummer} Rappelons quelques définitions (cf. \cite[3]{Overview@Illusie},
\cite[3.1]{Ket@Stix}, \cite{GFL@Vidal}). Un homomorphisme $h : P \ra Q$ de monoïdes intègres  est
dit \textit{de Kummer} si $h$ est injectif et, pour tout $q \in Q$, il existe $n \ge 1$ et $p \in P$
tel que $nq = h(p)$. On dit qu'un morphisme $f : X \ra Y$ de fs log schémas est \textit{de Kummer}
si, pour tout point géométrique $\overline{x}$ de $X$ d'image $\overline{y} = f(\overline{x})$ dans
$Y$, l'homomorphisme induit $\overline{M}_{\overline{y}} \ra \overline{M}_{\overline{x}}$ est de
Kummer. On dit qu'un morphisme $f : X \ra Y$ est \textit{Kummer étale} si $f$ est de Kummer et log
étale. On dit que $f$ est un \textit{revêtement Kummer étale} si $f$ est de Kummer étale et le
morphisme de schémas sous-jacent est fini.

Le \textit{site Kummer étale} d'un fs log schéma $X$ est la catégorie des fs log schémas de Kummer
étales au-dessus de $X$ munie de la topologie définie par les familles surjectives de $
X$-morphismes (lesquels sont automatiquement de Kummer étale). Pour $X$ connexe,  les revêtements
Kummer étales de $X$ forment une catégorie galoisienne, équivalente à la catégorie des
représentations d'un groupe profini, le \textit{groupe fondamental logarithmique} de $X$,
$\pi_1^{\mathrm log}(X,x)$, où $x$ est un \textit{point géométrique logarithmique} de $X$, cf.
(\textit{loc. cit.}).  \subsection{}\label{6-ket standard} Un morphisme $f : X \ra Y$ de log schémas
fs qui est déduit par changement de base par un morphisme strict $g : Y \ra \Spec \ZZ[P]$ d'un
morphisme $\Spec \ZZ[h] :\Spec \ZZ[Q] \ra \Spec \ZZ[P]$, où $h : P \ra Q$ est un homomorphisme de Kummer
entre monoïdes fs tels que le conoyau de $h^{\mathrm{gp}}$ soit annulé par un entier $n$ inversible
sur $Y$, est un revêtement Kummer étale. On dit qu'un tel revêtement est un revêtement Kummer étale
\textit{standard}. 

Un revêtement Kummer étale standard $f : X \ra Y$ comme ci-dessus est \textit{galoisien} dans le
sens suivant : le groupe diagonalisable (étale) $G =
Hom(Q^{\mathrm{gp}}/P^{\mathrm{gp}},\mathbf{G}_{m}{}_{Y})$ opère sur $X$ par automorphismes de
$Y$-log schémas, et le morphisme canonique $$ X \times_Y G \ra X \times_Y X, (x,g) \mapsto (x,xg),
$$ où le produit au second membre est pris dans la catégorie des fs log schémas, est un isomorphisme
(cf. \cite[3.2]{Overview@Illusie}). 

Tout morphisme de Kummer étale est, localement pour la topologie étale, isomorphe à un revêtement
Kummer étale standard (cf. \cite[3.1.4]{Ket@Stix}). Plus précisément, si $f : X \ra Y$ est un
morphisme de Kummer étale surjectif, avec $X$ (resp. $Y$) strictement local de point fermé $x$
(resp. $y$), on a un carré cartésien $$ \xymatrix{X \ar[r]^-a \ar[d]^f & \Spec \ZZ[Q] \ar[d]^{\Spec
\ZZ[h]} \\ Y \ar[r]^-b & \Spec \ZZ[P]}, $$ où $a$ (resp. $b$) est une carte de $X$ (resp. $Y$), et $h
: P \ra Q$ un morphisme de Kummer tel que $\coker h^{\mathrm{gp}}$ soit annulé par un entier $n$
inversible sur $Y$. Alors $f$ est un revêtement Kummer étale galoisien de groupe $G =
\mathrm{Hom}(\coker h^{gp},\mu_n(k(y)))$.

\subsection{}\label{6-Kummer standard local} Dans ce cas, l'action de $G$ sur $X$ peut se décrire de
la manière suivante. Soit $C = \coker h^{gp}$. La suite exacte $$ 0 \ra P^{\mathrm{gp}} \ra
Q^{\mathrm{gp}} \ra C \ra 0 $$ donne, par application du dual de Cartier
$\mathrm{Hom}(-,\GG_m)$, une suite exacte de groupes diagonalisables $$ 0 \ra \Gamma \ra T_Q
\ra T_P \ra 0, $$ où $T_P$ (resp. $T_Q$) est le tore sur $\ZZ$ dual de $P$ (resp. $Q$), et $G$ est
la valeur du groupe fini constant $\Gamma_Y$. Posons $Z_P = \Spec \ZZ[P]$, $Z_Q = \Spec \ZZ[Q]$.
L'accouplement canonique $$ T_P \otimes P^{\mathrm{gp}} \ra \GG_m $$ définit une famille de
caractères $$ (\chi_p : T_P \ra \GG_m)_{p \in P^{\mathrm{gp}}}, $$ qui détermine l'action de
$T_P$ sur $Z_P$ par $$ g.p = \chi_p(g)p, $$ pour un point $g$ de $T_P$ (à valeurs dans $S$), et $p
\in P$, vu comme un point de $Z_P$ à valeurs dans $S$. L'action de $T_Q$ sur $Z_Q$ est décrite de
façon similaire, et la carte $a : X \ra Z_Q$ est équivariante relativement aux actions de $G$ et
$T_Q$ sur $X$ et $Z_Q$ respectivement : pour $g \in G$ et $q \in Q$, \begin{equation} g.a^*(q) =
\chi_q(g)a^*(q), \label{6-g.q} \end{equation} où $a^*(q) \in M_{X,x}$ est l'image de $q$ par $a$. En
particulier, $G$ opère \textit{trivialement} sur $\overline{M}_{X,x}= M_{X,x}/\OO^*_{X,x}$. 

Par ailleurs, $P = Q \cap P^{\mathrm{gp}}$, donc $$ P = \{q \in Q ; \chi_q(g) = 1 \  \forall g \in G
\}.  $$ Ainsi, $Y$ est un quotient de $X$ par $G$, en tant que log schéma, et en tant que schéma :
$$ \OO_Y = (f_*\OO_X)^G \ \ ; \ \ M_Y = (f_*M_X)^G.  $$

Supposons de plus $Y$ log régulier, de sorte que $X$ l'est également. Alors $G$ opère librement
au-dessus de l'ouvert (dense) de $Y$ où la log structure est triviale, et fait de $X$ un revêtement
étale galoisien de $Y$ de groupe $G$. De plus, $G$ opère \textit{trivialement} sur la
\textit{strate} $C_{X,x}$ (\ref{6-1.5}). Pour le voir, on peut remplacer $X$ (resp. $Y$) par son
complété en $x$ (resp. $y$). On peut supposer de plus $P$ saillant, i. e.  $P^* = \{0\}$ (et même $P
= \overline{M}_{Y,y}$). Alors $Q$ est également saillant : on a $Q = Q^* \oplus Q_1$ avec $Q_1$
saillant, et $Q^* \cap P = \{0\}$, donc comme $Y \times_{\Spec \ZZ [P]} \Spec \ZZ [Q_1]$ est local, $Q^*
= \{0\}$.  D'après \ref{6-str loc log rég}, la carte $(h,a,b)$ de $f$ permet d'identifier
$\OO_{Y,y} \ra \OO_{X,x}$ à un homomorphisme de la forme $$ k[[t_1,\cdots,t_r]][[P]]
\ra k[[t_1,\cdots,t_r]][[Q]], $$ dans le cas d'égale caractéristique (resp.  $$
C[[t_1,\cdots,t_r]][[P]]/(u) \ra C[[t_1,\cdots,t_r]][[Q]]/(u) $$ dans le cas d'inégale
caractéristique $(p,0)$), où $k = k(x)  = k(y)$ (resp. $C$ est un anneau de Cohen de $k$, et $u$ est
congru à $p$ modulo l'idéal engendré par $P - \{0\}$ et les $t_i$), les $t_i$ et $u$ étant fixes par
$G$. Les strates $C_{X,x}$  et $C_{Y,y}$ sont respectivement $\Spec B$ et $\Spec A$, où $A$ et $B$ sont
les réductions de $\OO_{Y,y}$ et $\OO_{X,x}$ modulo les idéaux engendrés par
$P-\{0\}$ et $Q - \{0\}$. Comme $G$ opère par \ref{6-g.q} sur $Q$, $G$ opère donc trivialement sur
$C_{X,x}$, et $f : X \ra Y$ induit un isomorphisme $C_{X,x} \iso C_{Y,y}$.  \section{Actions très
modérées} \subsection{}\label{6-3.1} Soit $X$ un log schéma fs, muni d'une action d'un groupe fini
$G$. On se propose de dégager des conditions suffisantes sur l'action de $G$ pour que, lorsque $X$
est log régulier, le quotient de $X$ par $G$ existe comme log schéma et soit log régulier. 

On dit que $G$ opère \textit{modérément} sur $X$ en un point géométrique $\overline{x}$ de $X$
localisé en $x$, si le stabilisateur $G_{\overline{x}}$ de $\overline{x}$ (\textit{groupe d'inertie
en $x$}) est d'ordre premier à la caractéristique de $k(x)$. On dit que $G$ opère modérément sur $X$
si $G$ opère modérément en $\overline{x}$ pour tout $\overline{x}$. Ces définitions ne font pas
intervenir la log structure de $X$.

La définition et les résultats qui suivent sont dues à Gabber. On dit que $G$ opère \textit{très
modérément} sur $X$ \textit{en} $\overline{x}$ si les trois conditions suivantes sont satisfaites :

(i) $G$ opère modérément en $\overline{x}$ ;

(ii) $G$ opère trivialement sur $\overline{M}_{X,\overline{x}}$ ;

(iii) $G_{\overline{x}}$ opère trivialement sur la strate $C_{X,\overline{x}}$ (\ref{6-1.2}).

On dit que $G$ opère \textit{très modérément} sur $X$ si $G$ opère très modérément sur $X$ en tout
point géométrique.

Le résultat principal est le suivant : \begin{theoreme}\label{6-très modéré}\label{6-tres modere}\label{6-3.2}
Soit $X$ un log schéma fs log régulier, muni d'une action d'un groupe fini $G$. On suppose que $G$
opère de façon admissible sur le schéma sous-jacent à $X$, librement sur un ouvert dense, et très
modérément. Alors le quotient $Y = X/G$ existe comme log schéma, est un log schéma fs log régulier,
et la projection $f : X \ra Y$ est un revêtement Kummer étale de groupe $G$.  \end{theoreme}
Rappelons (\sga{1}{V}{1}) que dire que $G$ opère de façon admissible signifie que $X$ est réunion
d'ouverts affines stables par $G$, de sorte que le quotient $X/G$ est un schéma. 

Compte tenu de \ref{6-ket standard}, on en déduit : \begin{corollaire}\label{6-inertie abélienne}
Sous les hypothèses de \ref{6-très modéré}, pour tout point géométrique $\overline{x}$ de $X$
localisé en $x$, d'image $\overline{y}$ dans $Y$, le groupe d'inertie $H = G_{\overline{x}}$ en $x$
est abélien, d'ordre premier à la caractéristique de $k(x)$. Le localisé strict $X_{(\overline{x)}}$
est un revêtement Kummer étale du localisé strict $Y_{(\overline{y})}$ de groupe $H$, et le
revêtement Kummer étale $Y_{(\overline{y})} \times_Y X$ induit sur $Y_{(\overline{y})}$ s'en déduit
par extension du groupe $H$ à $G$.  \end{corollaire} \subsection{}\label{6-carte équivariante} Nous
démontrerons un résultat plus précis que \ref{6-très modéré}, pour lequel nous aurons besoin de la
notion suivante. Soient $n$ un entier $\ge 1$, $G$ un groupe fini, et $Q$ un monoïde fs. Supposons
donné un homomorphisme $$ \chi : G^{ab} \otimes Q^{gp} \ra \mu_n := \mu_n(\CC), \ \ g \otimes q
\mapsto \chi_q(g).  $$ Soit $\Lambda = \ZZ[\mu_n][1/n]$. On déduit de $\chi$  une action de $G$ sur
le log schéma $\Spec \Lambda[Q]$, caractérisée par $$ g.q = \chi_q(g)q $$ pour $g \in G$, $q \in Q$.

Soit $X$ un log schéma fs muni d'une action de $G$. Par une \textit{carte $G$-équivariante de $X$
modelée sur $\chi$ et Q}, on entend un morphisme strict, $G$-équivariant $$ c : X \ra \Spec
\Lambda[Q].  $$ \begin{proposition}\label{6-très modéré ponctuel}\label{6-tres modere ponctuel} Soit
$X$ un fs-log schéma, muni d'une action d'un groupe fini $G$, et soit $\overline{x}$ un point
géométrique de $X$ localisé en $x$. On note $H = G_{\overline{x}}$ le groupe d'inertie en $x$. 

(a) On suppose que les conditions (i) et (ii) de \ref{6-3.1} sont vérifiées en $\overline{x}$. Soit
$n$ l'exposant de $|H|$. Il existe, localement pour la topologie étale au voisinage de
$\overline{x}$, une carte $H$-équivariante de $X$ modelée sur $Q =\overline{M}_{\overline{x}}$ et un
homomorphisme $\chi : H^{ab} \otimes Q^{gp} \ra \mu_n$. 

(b) On suppose de plus que $G$ agit très modérément en $\overline{x}$, i. e. que $H$ agit
trivialement sur la strate $C_{\overline{x}}(X)$. Alors, localement pour la topologie étale au
voisinage de $\overline{x}$ et de son image $\overline{y}$ dans $Y$, le quotient $Y = X/H$ existe
comme log schéma, et $X$ est un $Y$-sous log schéma fermé strict d'un revêtement Kummer étale
standard de $Y$.

(c) Si, sous les hypothèses de (ii), $X$ est supposé en outre log régulier, alors
$Y_{(\overline{y})}$ est log régulier, et $X_{(\overline{x})}$ est un revêtement Kummer étale
standard de $Y_{(\overline{y})}$ de groupe $H$.  \end{proposition} Sous les hypothèses de
\ref{6-très modéré}, le log schéma de schéma sous-jacent $Y = X/G$, muni de la log structure $M_Y =
(f_*M_{X})^{G}$ est log régulier, et s'identifie, par localisation stricte en $\overline{y}$ au
quotient $Y_{(\overline{y})} = X_{(\overline{x})}/H$. Le fait que $G$ opère génériquement librement
implique la dernière assertion de \ref{6-inertie abélienne}, et \ref{6-très modéré} s'en déduit.

Par ailleurs, compte tenu de la description locale \ref{6-Kummer standard local} des revêtements
Kummer étales standard, \ref{6-très modéré ponctuel} implique : \begin{corollaire} \label{6-tres
modere dans voisinage} Sous les hypothèses de \ref{6-très modéré ponctuel}, si $X$ est log régulier
et si $G$ agit très modérément en $\overline{x}$, alors $G$ agit très modérément sur $X$ dans un
voisinage étale de $\overline{x}$.  \end{corollaire} \subsection{} \label{6-preuve très modéré
ponctuel} \label{6-preuve tres modere ponctuel} \textit{Preuve de  \ref{6-très modéré ponctuel}} 

(a) Considérons la suite exacte de groupes abéliens \begin{equation} 0 \ra
\OO_{X,\overline{x}}^* \ra M_{X,\overline{x}}^{gp} \ra \overline{M}_{X,\overline{x}}^{gp}
\ra 0.  \label{6-exacte} \end{equation} Elle est $H$-équivariante, et si $Q =
\overline{M}_{\overline{x}}$, $Q^{gp}$ est de type fini sans torsion. Choisissons un scindage $s :
Q^{gp} \ra M_{X,\overline{x}}^{gp}$ de \ref{6-exacte} (comme suite exacte de groupes abéliens).
Comme $H$ agit trivialement sur $Q^{gp}$, on a, pour $a \in Q^{gp}$ et $g \in H$ $$ g.s(a) =
z(g,a)s(a) $$ avec $z(g,a) \in \OO_{X,\overline{x}}^*$. Pour $g_1$, $g_2$ dans $H$, on a $$
z(g_1g_2,a) = (g_1z(g_2,a)).z(g_2,a), $$ en d'autres termes, $g \mapsto (a \mapsto z(g,a))$ est un
1-cocycle $$ z \in Z^1(H,\mathrm{Hom}(Q^{gp},\OO_{X,\overline{x}}^*)).  $$ L'image $[z]$ de
$z$ dans $H^1(H,\mathrm{Hom}(Q^{gp},\OO_{X,\overline{x}}^*))$ est la classe de cohomologie
de \ref{6-exacte}. Dans le carré commutatif de flèches canoniques $$
\xymatrix{Z^1(H,\mathrm{Hom}(Q^{gp},\OO_{X,\overline{x}}^*)) \ar[r] \ar[d]
&Z^1(H,\mathrm{Hom}(Q^{gp},k(\overline{x})^*)) \ar[d] \\
H^1(H,\mathrm{Hom}(Q^{gp},\OO_{X,\overline{x}}^*)) \ar[r] &
H^1(H,\mathrm{Hom}(Q^{gp},k(\overline{x})^*))}, $$ la flèche verticale de droite est (trivialement)
un isomorphisme. Comme $n$ est premier à la caractéristique de $k(x)$,  $z$ est l'unique relèvement
de son image $\chi$ dans $Z^1(H,\mathrm{Hom}(Q^{gp},k(\overline{x})^*))$, et la flèche horizontale
inférieure est un isomorphisme, puisque $(1 + \mathbf{m}_{X,\overline{x}})^*$ est $n$-divisible.
Identifiant $Z^1(H,\mathrm{Hom}(Q^{gp},k(\overline{x})^*))$ à
$\mathrm{Hom}(H,\mathrm{Hom}(Q^{gp},k(\overline{x})^*))$, on trouve ainsi un homomorphisme $$ \chi :
H^{ab} \otimes Q^{gp} \ra \mu_n, $$ et la restriction de $s$ à $Q$ est, au voisinage de
$\overline{x}$, une carte $H$-équivariante $a$ de $X$ modelée sur $Q$ et $\chi$.

(b) On peut remplacer $X$ par son localisé strict $X_{(\overline{x})}$ et $G$ par $H$. Le quotient
$Y = X/G$ existe alors comme schéma, est strictement local, de point fermé $\overline{y}$, avec
$k(\overline{y}) = k(\overline{x})$, et $X \ra Y$ est fini. Soit $P'$ le sous-groupe de $Q^{gp}$
défini par $$ P' = \{q \in Q^{gp} | \chi_q(g) = 1 \ \ \forall g \in G\}.  $$ Il est d'indice fini,
premier à la caractéristique de $k(\overline{x})$. Soit $$ P = P' \cap Q.  $$ La carte équivariante
$a : X \ra \Spec \Lambda[Q]$ de (a) (où $\Lambda = \ZZ [\mu_n][1/n]$) définit, par passage au
quotient, un morphisme $b : Y \ra \Spec \Lambda[P]$. Ce morphisme $b$ définit une log structure sur
$Y$ pour laquelle $(f_*M_X)^G = M_Y$ ($f : X \ra Y$ désignant la projection), et une carte de $f$
\begin{equation} \xymatrix{X \ar[r] \ar[d]^f & \Spec \Lambda[Q] \ar[d] \\ Y \ar[r] & \Spec \Lambda[P]},
\label{6-carré} \end{equation} modelée sur l'homomorphisme de Kummer $P \subset Q$. Comme $G$ agit
trivialement sur la strate $C_{\overline{x}}(X) = \Spec
\OO_{X,\overline{x}}/I_{\overline{x}}$, et que $G$ est d'ordre premier à la caractéristique
de $k(\overline{x})$, l'homomorphisme $$ \OO_{Y,\overline{y}} (=
\OO_{X,\overline{x}}^G) \ra \OO_{X,\overline{x}}/I_{\overline{x}} (=
(\OO_{X,\overline{x}}/I_{\overline{x}})^G) $$ est surjectif, donc il en est de même de
l'homomorphisme $$ k(\overline{y})[Q] \  = \ k(\overline{x})[Q] \ra
\OO_{X,\overline{x}}/\mathbf{m}_{Y,\overline{y}}\OO_{X,\overline{x}}.  $$ Par
Nakayama, il en résulte que le morphisme \begin{equation} i : X \ra Y \times_{\Spec \Lambda[P]} \Spec
\Lambda[Q] \label{6-morphisme} \end{equation} déduit de \ref{6-carré} est une immersion fermée
stricte.

(c) Dans ce cas, la strate $C_{\overline{x}}(X)$ est régulière, et se projette isomorphiquement sur
$C_{\overline{y}}(Y)$. Comme $\mathrm{rg}(\overline{M}^{\mathrm{gp}}_{\overline{x}}) =
\mathrm{rg}(\overline{M}^{\mathrm{gp}}_{\overline{y}}) =
\mathrm{codim}(C_{\overline{y}}(Y),Y_{(\overline{y})})$, $Y$ est log régulier en $y$, et l'immersion
$i$ est un isomorphisme, de sorte que $X_{(\overline{x})}$ est un revêtement Kummer étale de
$Y_{(\overline{y})}$ de groupe $H$.  \section{Points fixes} \subsection{}\label{6-stratification}
Soit $X$ un schéma noethérien séparé sur lequel opère un groupe fini $G$. Pour chaque sous-groupe
$H$ de $G$, on note $X^H$ le schéma des points fixes de $H$. C'est un sous-schéma fermé de $X$,
représentant le foncteur $S \mapsto X(S)^H$, intersection des graphes des translations  $h : X \ra
X$ pour $h \in H$. En chaque point géométrique $\overline{x}$ de $X^H$ d'image $x$, le sous-groupe
d'inertie $G_{\overline{x}}$ de $G$ contient $H$, et est égal à $H$ si et seulement si $x$
appartient au sous-schéma (localement fermé) $$ X_H = X^H - \cup_{H' \supset H, H' \ne H}X^{H'}.  $$
Pour $g \in G$, on a $$ gX^H = X^{gHg^{-1}}, $$ et de même $$ gX_h = X_{gHg^{-1}} $$ de sorte que la
réunion $X^C$ (resp. $X_C$) des $X^H$ (resp. $X_H$) pour $H$ dans une classe de conjugaison $C$ de
sous-groupes de $G$ est $G$-stable. Les $X_C$, pour $C$ parcourant les classes de conjugaison de
sous-groupes de $G$, forment une stratification de $X$ par des sous-$G$-schémas, avec la propriété
que pour tout point géométrique $\overline{x}$ localisé en $X_C$, le groupe d'inertie
$G_{\overline{x}}$ appartient à $C$. Nous appellerons cette stratification la \textit{stratification
par l'inertie}.

Le but de ce numéro est de donner des exemples d'actions très modérées de groupes finis $G$ sur des
log schéma log réguliers et réguliers $X$, où un raffinement de la stratification par l'inertie est
déduite de la stratification canonique d'un diviseur à croisements normaux $G$-stable.

Le résultat suivant est classique, nous en donnons une démonstration, faute de référence.
\begin{proposition}\label{6-pts fixes reguliers} Soit $X$ un schéma noethérien régulier sur lequel
opère un groupe fini $G$ de façon modérée \ref{6-3.1}. Alors le schéma des points fixes $X^G$ est
régulier.  \end{proposition} On peut supposer $X$ local, $X = \Spec A$, de sorte que $X^G = \Spec A_G$,
où $A_G$ est l'algèbre des co-invariants $A/I$, $I$ l'idéal engendré par les $ga - a$, pour $g$ dans
$G$ et $a$ dans l'idéal maximal $\mathbf{m}$ de $A$. On peut supposer de plus que $G$ est contenu
dans le groupe d'inertie au point fermé, sinon $X^G$ est vide. Comme $\widehat X^G = \Spec \hat A/I
\hat A = (\hat X)^G$, on peut supposer $X $local et complet. On démontre le lemme en
\textit{linéarisant} l'action de $G$.

Supposons d'abord que $G$ soit d'égale caractéristique $p$. Choisissons une base $t = (t_i)_{1 \le i
\le r}$ de l'espace cotangent $T = \mathbf{m}/\mathbf{m}^2$ et des éléments $x_i$ dans $\mathbf{m}$
relevant les $t_i$. Choisissons un corps de représentants, noté encore $k$, de $k$ dans $A$, et
notons $$ \varphi : k[[T]] \ra A $$ l'homomorphisme envoyant $t_i$ sur $x_i$. L'homomorphisme $$ f :
k[[T]] \ra A $$ envoyant $t$ sur le système $$ y = (1/\mathrm{card}(G))\sum_{g \in G} g\varphi
g^{-1}t $$ est $G$-équivariant. En outre, $y$ est congru à $x$ modulo $\mathbf{m}^2$, donc est un
système régulier de paramètres de $A$, et donc $f$ est un isomorphisme. Par suite $f$ induit un
isomorphisme $$ k[[T]]_G = k[[T_G]] \ra A_G, $$ où $T_G$ est l'espace des co-invariants de $G$ sur
$T$, ce qui démontre le lemme dans ce cas.

Supposons maintenant que $A$ soit de caractéristique mixte $(0,p)$. Soit $C$ un anneau de Cohen pour
$k$. Choisissons un homomorphisme $C \ra A$ relevant $C \ra k$. Comme $G$ est d'ordre premier à $p$,
il existe un unique $C[G]$-module $V$, libre de type fini sur $C$, relevant $T$. Soit $t = (t_i)_{1
\le i \le r}$ une base de $T$, $v = (v_i)$ une base de $V$ relevant $t$. Comme plus haut,
choisissons des relèvements $x_i$ dans $\mathbf{m}$ des $t_i$, qui forment donc une suite régulière
de paramètres dans $A$. Prolongeons l'homomorphisme $C \ra A$ en $$ \varphi : C[[V]] \ra A $$ en
envoyant $v_i$ sur $x_i$. L'homomorphisme $$ f : C[[V]] \ra A $$ envoyant $v = (v_i)$ sur $$ y =
(1/\mathrm{card}(G))\sum_{g \in G} g\varphi g^{-1}v $$ est $G$-équivariant, et $y$ est congru à $x$
modulo $\mathbf{m^2}$, donc est un système régulier de paramètres de $A$. L'homomorphisme $f$ est
donc surjectif, et son noyau est défini par un élément $F$ congru à $p$ modulo l'idéal engendré par
les $v_i$. Ainsi $X$ est un diviseur régulier $G$-équivariant dans $X' = \Spec C[[V]]$, d'équation
$F=0$. L'image de $F$ dans $C[[V_G]]$ est un paramètre régulier, donc, comme $X'^G = \Spec C[[V_G]]$,
$X^G = X \times_{X'} X'^G$ est régulier.  \begin{corollaire}\label{6-stratification régulière} Soit
$X$ un schéma noethérien séparé, régulier, muni d'une action modérée d'un groupe fini $G$. Alors la
stratification de $X $par l'inertie est formée de schémas réguliers.  \end{corollaire}
\subsection{}\label{6-dcn stables} Sous les hypothèses de \ref{6-stratification régulière}, soit $Y$
un diviseur à croisements normaux stricts $G$-stable, réunion de composantes irréductibles $Y_i$,
$1\le i \le m$. On munit $X$ de la log structure définie par $Y$. Rappelons (cf.
\cite[7.1]{Alterations@de_Jong}) qu'on dit que $Y$ est \textit{$G$-strict} si la condition suivante
est réalisée :  pour tout $i$ et pour tout $g \in G$, si $gY_i \cap Y_i \ne \emptyset$, alors $gY_i
= Y_i$. Si $Y$ est $G$-strict, alors la condition (ii) de \ref{6-3.1} est vérifiée en chaque point
géométrique $\overline{x}$ de $X$. En effet, si $(D_i)_{1 \le i \le r}$ est l'ensemble des branches
de $Y$ passant par $\overline{x}$, alors $gD_i = D_i$ pour tout $i$. Comme
$\overline{M}^{\mathrm{gp}}_{\overline{x}} = \oplus_{1 \le i \le r}\ZZ e_i$, $e_i$ correspondant à
$D_i$, le groupe d'inertie $G_{\overline{x}}$ opère trivialement sur
$\overline{M}^{\mathrm{gp}}_{\overline{x}}$. Rappelons également (\cite[7.2]{Alterations@de_Jong})
qu'il existe une modification $G$-équivariante canonique $f : \tilde X \ra X$ telle que
$f^{-1}(Y)_{\mathrm{red}}$ soit $G$-strict.

\begin{corollaire}\label{6-dcn stricts} Supposons que $Y$ soit $G$-strict, que $G$ opère de façon
modérée, admissible et génériquement libre, et que la stratification canonique de $X$ soit plus fine
que la stratification par l'inertie. Alors $G$ opère très modérément sur $X$ (et donc la conclusion
de \ref{6-très modéré} s'applique). Le groupe d'inertie le long de la strate $X^{(i)}$ (\ref{6-1.5}
(i)) est constant, de valeur $G_i$, avec $\mathrm{rg}(G_i) = \mathrm{codim}(X^{(i)},X)$. En
particulier, $G$ opère librement sur la strate $X^{(0)} = X-Y$.  \end{corollaire}
\subsection{}\label{6-exemples}\textit{Exemples}. 

(a) Soient $k$ un corps algébriquement clos, $n$ un entier $\ge 2$ premier à la caractéristique
résiduelle de $k$, $G$ le groupe $\mu_n = \mu_n(k)$. On fait opérer $G$ sur $X = \mathbf{A}^2_k$ par
homothéties ($(\lambda,x) \mapsto \lambda x$ pour $\lambda \in G$, $x \in X(k)$). La stratification
par l'inertie comporte deux strates, $X-\{0\}$, où $G$ opère librement, et $X^G = \{0\}$. La donnée
de deux droites $Y_1$, $Y_2$ telles que $Y_1 \cap Y_2 = \{0\}$ définit un diviseur à croisements
normaux $G$-strict $Y = Y_1 \cup Y_2$, et le couple $(X,Y)$ vérifie les conditions de \ref{6-dcn
stricts}. Le choix de paramètres $t_1$, $t_2$ tels que $Y_i = V(t_i)$ permet de définir une carte
équivariante \ref{6-carte équivariante} $c : \Spec \ZZ[\NN^2] \la X$, $e_i \mapsto t_i$, associée à
l'homomorphisme $\chi : G \otimes \ZZ^2 \ra \mu_n$ tel que $\chi(\lambda \otimes e_i) = \lambda$. Le
quotient $X/G$ est le schéma torique $\Spec k[x^n,x^{n-1}y, \cdots, xy^{n-1}, y^n]$.

Plus généralement, soient $n$ un entier $\ge 1$, $G$ un groupe abélien d'ordre $n$, $S$ un schéma
noethérien régulier au-dessus de $\Spec \ZZ[1/n,\mu_n]$, $E$ un $\OO_{S}$-module localement
libre de rang fini, muni d'une action linéaire de $G$, $X$ le fibré vectoriel $V(E)$. Pour chaque
caractère $\chi : G \ra \mu_n$, notons $L_{\chi}$ le $G$-$\OO_S$-module correspondant.
L'homomorphisme canonique $G$-équivariant $$ \bigoplus_{\chi} L_{\chi} \otimes E_{\chi} \ra E, $$ où
$E_{\chi} = \mathscr{H}om_G(L_{\chi},E)$ et $\chi$ parcourt le groupe des caractères de $G$, est un
isomorphisme. Il définit une décomposition $G$-équivariante $$ X \iso \oplus_{\chi} X_{\chi}, $$ où
$X_{\chi} = V(E_{\chi})$, muni de l'action de $G$ à travers $\chi$. En particulier, $X^G = X_1$, où
$1 : G \ra \mu_n$ est le caractère trivial. Supposons $S$ local, $S = \Spec A$. Pour chaque $\chi \in
\mathrm{Hom}(G,\mu_n)$, choisissons une base $(t_i)_{i \in I_{\chi}}$ de $\mathrm{Hom}(E_{\chi},A)$,
de sorte que $X_{\chi} = \Spec A[(t_i)_{i \in I_{\chi}}]$, avec $gt_i = \chi(g)t_i$ pour $g \in G$, $i
\in I_{\chi}$. Le couple formé de $X$ et du diviseur à croisements normaux (relatifs à $S$) $Y =
\sum_{\chi, i \in I_{\chi}}Y_i$, où $Y_i = (t_i = 0)$ pour $i \in I_{\chi}$, vérifie les conditions
de \ref{6-dcn stricts} (et les vérifie d'ailleurs fibre à fibre).  

(b) Soient $k$ un corps algébriquement clos d'exposant caractéristique $p$, $n$ un entier $\ge 2$
tel que $(2n,p) = 1$, $G$ le groupe diédral $D_n = ⟨s,r : s^2 = 1, r^n = 1, srs = r^{-1}⟩$. Soit
$\zeta \in k$ une racine primitive $n$-ième de 1. Soit $\rho : G \ra \mathrm{GL}(E)$  la
représentation de degré 2 induite du caractère $\chi$ de $\mu_n \subset G$ tel que $\chi(r) = \zeta$
: $\rho(s) = \begin{pmatrix} 0 & 1 \\ 1 & 0 \\  \end{pmatrix} $, $\rho(r) = \begin{pmatrix} \zeta &
0 \\ 0 & \zeta^{-1} \\ \end{pmatrix}$. Soit $X$ le $G$-schéma $V(E) = \Spec k[u,v]$, $s(u) = v, r(u) =
\zeta u$. Pour $0 \le i \le n-1$, notons $Z_i \subset X$ la droite $v = \zeta^iu$, et $Z = \bigcup_{0
\le i \le n-1}Z_i$. La stratification par l'inertie comporte $n+2$ strates : $X-Z$, où $G$ opère
librement, $Z_i - \{0\}$ ($0 \le i \le n-1$), où le groupe d'inertie est d'ordre 2 (de générateur
$sr^i$), et $\{0\} = X^G$. 

Pour $n = 2$, $G = (\ZZ/2\ZZ)^2$, $Z$ est un diviseur à croisements normaux $G$-strict, et le couple
$(X,Y)$ vérifie les conditions de \ref{6-dcn stricts}. Pour $n > 2$, $Z$ n'est plus un diviseur à
croisements normaux, et l'inertie en $\{0\}$ n'est plus abélienne. Soient $f : X' \ra X$ l'éclaté de
$\{0\}$ dans $X$, $E = f^{-1}(0)$ le diviseur exceptionnel, $Z' = \bigcup_{0 \le i \le n-1}Z'_i$ le
transformé strict de $Z$. Alors $G$ opère de façon naturelle sur $X'$, la projection $f$ est
$G$-équivariante, et $Y' = f^{-1}(Z) = E \cup Z'$ est un diviseur à croisements normaux $G$-strict.
Le couple $(X',Y')$ vérifie les conditions de \ref{6-dcn stricts}. La stratification de $X'$ par
l'inertie se compose des strates $Z'_i$, où l'inertie est un groupe à deux éléments, et de $X'-Z'$,
où $G$ opère librement. La stratification canonique associée à $Y'$ la raffine : $X'^{(0)} = X' -
Y'$, $X'^{(1)} = Y' - \bigcup_{0 \le i \le n-1}(E \cap Z'_i)$, $X'^{(2)} = \bigcup_{0 \le i \le
n-1}(E \cap Z'_i)$.  \subsection{}\label{6-éclatements}\label{6-eclatements} La construction
précédente, qui rend les inerties abéliennes, se généralise. Soit $X$ un schéma noethérien régulier,
séparé, muni d'une action modérée d'un groupe fini $G$, et soit $Y$ un diviseur à croisements
normaux $G$-strict. Si $H$ est un sous-groupe de $G$, $X^H$ est régulier (et séparé), donc il en est
de même de l'éclaté $X' = \Ecl_{X^H}(X)$ de $X$ le long de $X^H$. Le normalisateur $N = N_G(H)$ de $H$
dans $G$ stabilise $X^H$, donc agit sur $X'$, et le morphisme $f : X' \ra X$ est équivariant
relativement à $N \ra G$. De plus, $f^{-1}(X^H)$ est un diviseur régulier dans $X'$. Si $D$ est une
composante de $Y$, comme $D$ est $H$-stable, $D \times_X X^H = D^H$ est régulier, et le transformé
strict $\tilde D = \Ecl_{D^H}(D)$ est un diviseur régulier croisant $f^{-1}(X^H)$ transversalement. Il
s'ensuit que le transformé total réduit $Y' = f^{-1}(Y)_{\mathrm{red}}$ est un diviseur à
croisements normaux $N$-strict dans $X'$.  \begin{proposition}\label{6-inerties abéliennes} Sous les
hypothèses de \ref{6-éclatements}, soit $\overline{x}$ un point géométrique de $X$ en lequel le
groupe d'inertie $G_{\overline{x}}$ n'est pas abélien, et soit $H$ le sous-groupe des commutateurs
$(G_{\overline{x}},G_{\overline{x}})$. Alors $G_{\overline{x}} = N_{G_{\overline{x}}}(H)$ agit sur
$X' = \Ecl_{X^H}(X)$, et en chaque point géométrique $\overline{y}$ de $X'$ au-dessus de
$\overline{x}$, le groupe d'inertie $(G_{\overline{x}})_{\overline{y}}$ est strictement plus petit
que $G_{\overline{x}}$.  \end{proposition} En effet, le point $\overline{y}$ correspond à une droite
$L$ dans la fibre en $\overline{x}$ du fibré normal $T_{\overline{x}}/T_{\overline{x}}^H$ de $X^H$
dans $X$, où $T_{\overline{x}} = T_{\overline{x}}(X)$. Supposons que $\overline{y}$ soit fixe sous
$G_{\overline{x}}$. Alors $G_{\overline{x}}$ agit sur $L$ par un caractère, donc $H$ agit
trivialement sur $L$. Or $(T_{\overline{x}}/T_{\overline{x}}^H)^H = 0$, contradiction. (Noter que
cet argument montre en particulier que, si $H \ne \{1\}$, $X^H$ est de codimension $\ge 2$ dans $X$
en $\overline{x}$.)

\ifx\danslelivre\undefined
\bibliography{../fichiers-configuration-latex/bibliographie}
\bibliographystyle{../fichiers-configuration-latex/smfalpha}


\begin{center}
VII. Démonstration du théorème d'uniformisation locale (faible)\\
Fabrice Orgogozo
\end{center}
\version

\else
\setboolean{english}{false}
\setboolean{francais}{true}
\selectlanguage{francais}
\chapter{Démonstration du théorème d'uniformisation
locale (faible)}
\addtocontents{toc}{par Fabrice Orgogozo}
\begin{center}
Fabrice Orgogozo
\end{center}
\fi

\section{Énoncé}

L'objet de cet exposé est de démontrer le théorème
\refext{II}{2-théorème uniformisation faible}
(voir aussi \refext{0}{0-uniformisation faible}),
dont nous rappelons l'énoncé ci-dessous :

\begin{théorème}\label{7-uniformisation}
Soient $X$ un schéma nœthérien quasi-excellent
et $Z$ un fermé rare de $X$.
Il existe une famille finie de morphismes 
$(X_i→X)_{i∈I}$, couvrante pour la topologie des
altérations et telle que pour tout $i∈I$ on ait :
\begin{enumerate}
\item le schéma $X_i$ est régulier et connexe ;
\item l'image inverse de $Z$ dans $X_i$ est le support
d'un diviseur à croisements normaux stricts.
\end{enumerate}
\end{théorème}

\section{Réductions : rappel des résultats antérieurs}\label{7-réductions}

\subsection{Réduction au cas local, normal de dimension finie}
Nous avons vu en \refext{II}{2-réduction cas hensélien}
qu'il suffit de démontrer le théorème lorsque
le schéma $X$ est local nœthérien normal hensélien excellent.
Faisons cette hypothèse supplémentaire. Un tel schéma est nécessairement de dimension finie,
que nous noterons ici $d$. De plus,
on a vu en \emph{loc. cit.} que si le théorème
est établi pour chaque schéma local nœthérien hensélien
excellent de dimension au plus $d$, il en est
de même pour les schémas nœthériens quasi-excellents
de dimension au plus $d$.

\subsection{Réduction au cas complet}
Il résulte de la proposition \refext{III}{3-reduc}
qu'il suffit de démontrer le théorème pour le schéma
local nœthérien complet $\chap{X}$, ce dernier
étant de même dimension que $X$ et également normal.

\subsection{Récurrence}
\label{7-récurrence}
Il résulte de ce qui précède que l'on peut supposer
le schéma $X$ local nœthérien complet normal
de dimension $d$ et le théorème connu pour chaque
schéma nœthérien quasi-excellent de dimension au plus $d-1$.
Lorsque $d=1$, le théorème est bien connu ; nous supposerons
dorénavant $d ≥ 2$.

\section{Fibration en courbes et application d'un théorème de A. J.
de Jong}

\subsection{}Soit $X=\Spec(A)$ un schéma local nœthérien complet normal
comme en \ref{7-récurrence} et $Z$ un fermé rare.
Quitte à remplacer $X$ (resp. $Z$) par un $X$-schéma fini
également local nœthérien normal excellent de dimension $d$
(resp. par son image inverse), on peut supposer
d'après \refext{V}{5-algébrisation partielle mixte raffinée géométrique},
qu'il existe un schéma local nœthérien régulier $S$ de dimension $d-1$,
un $S$-schéma de type fini dominant $X′$ intègre et affine, un point fermé $x′$ de la
fibre spéciale de $f:X′ → S$, et un fermé rare $Z′$ de $X′$
satisfaisant les conditions suivantes :

— le morphisme $c$ induit un isomorphisme $X ⥲ \Spec(\chap{𝒪_{X′,x′}})$ ;

— l'image inverse $c^{-1}(Z′)$ de $Z′$ coïncide avec $Z$.

\begin{center}
\begin{tikzpicture}[auto]
\matrix(diag)[matrix of math nodes]
{|(X)| X & |(Xprime)| X′ \\
   & |(S)| S  \\
};
 \draw[->] (X) -- node{$c$}  (Xprime);
 \draw[->] (Xprime) -- node{$f$} (S);
\end{tikzpicture}
\end{center}

\subsection{}Supposons l'existence d'une famille $(X′_i → X′)$
couvrante pour la topologie des altérations (\refext{II}{2-définition altL})
telle que chaque $X_i′$ soit régulier et chaque image inverse $Z_i ′$
de $Z′$ dans $X_i′$ soit le support d'un diviseur à croisements
normaux. Il résulte de \refext{II}{2-stabilité par cb}
que la famille $(X_i → X)$ obtenue par changement de base
(plat) $X → X′$ est également $\alt$-couvrante.
D'autre part, l'hypothèse d'excellence faite sur les schémas
garantit que le morphisme de complétion $c$ est \emph{régulier}
(\refext{I}{1def-excellent}).
La régularité d'un morphisme étant stable par changement de base
localement de type fini (\ega{IV₂}{6.8.3}),
et préservant la régularité des schémas (\ega{IV₂}{6.5.2 (ii)})
il en résulte que chaque $X_i$ est régulier. De même, l'image
inverse $Z_i$ de $Z_i′$ dans $X_i$ — qui coïncide
avec l'image inverse de $Z$ dans $X_i$ par le morphisme
évident — est le support d'un diviseur à croisements normaux
pour chaque indice $i$. 

\subsection{}Quitte à remplacer $X$ par $X′$, ce
qui est licite d'après ce qui précède, nous pouvons
supposer le schéma $X$ intègre de dimension $d$,
équipé d'un morphisme dominant de type fini $f:X → S$ où $S$
est local nœthérien régulier de dimension $d-1$.
Quitte à compactifier $f$, on peut le supposer \emph{propre} ;
quitte à éclater, on peut supposer que le fermé $Z$ est un
\emph{diviseur}.

\subsection{}
Nous sommes dans les conditions d'application
du théorème \cite[2.4]{Families@de_Jong},
d'après lequel, quitte à altérer $S$ et $X$, on peut supposer les faits suivants :

— le morphisme $f$ est une courbe nodale ;

— le diviseur $Z$ est contenu dans la réunion
d'un diviseur $D$ étale sur $S$, contenu dans le lieu
lisse de $f$, et de l'image inverse $f^{-1}(T)$ d'un
fermé rare $T$ de $S$.

\section{Résolution des singularités}
\label{7-utilisation log-résolution}

\subsection{Résolution des singularités de la base}
Les altérations précédentes conduisent à une situation
où le schéma $X$ n'est plus local (ni même affine) et
le schéma $S$ n'est plus nécessairement régulier.
Il est cependant excellent de dimension $d-1$
donc justiciable de l'hypothèse de récurrence \ref{7-récurrence}.
Ainsi, on peut supposer que la paire $(S,T)$ est régulière,
c'est-à-dire que le schéma $S$ est régulier et que $T$
est un diviseur à croisements normaux. En effet,
il en est ainsi localement pour la topologie des altérations.

\subsection{}
D'après \refext{VI}{6-nodal implique log lisse},
la paire $(X,D ∪ f^{-1}(T))$ est log-régulière
au sens de \refext{VI}{6-1.2}.
Observons qu'un diviseur contenu dans un
diviseur à croisements normaux est lui-même à croisements normaux.
Ceci nous permet de supposer d'agrandir $Z$ si
nécessaire et auquel on applique alors le théorème suivant
de Katô K. (\cite[10.3, 10.4]{Toric@Kato}), complété par W. Nizioł
(\cite[5.7]{Toric@Niziol}). (Voir aussi \cite[§5.6]{HT@GR} et \refext{VIII}{8-resolution of log schemes1}.)

\begin{théorème}
Soit $(X,Z)$ une paire log-régulière, où $X$ est un schéma nœthérien.
Il existe un schéma noethérien régulier $Y$, un diviseur à
croisements normaux $D⊂Y$ et un morphisme projectif birationnel $π:Y→ X$
tel que l'image inverse $π^{−1}(Z)$ soit contenue dans $D$.
\end{théorème}

\ifx\danslelivre\undefined
\bibliography{../fichiers-configuration-latex/bibliographie}
\bibliographystyle{../fichiers-configuration-latex/nomcomplet-date}


\begin{center} XVIII. Gabber's modification theorem (absolute case)\\
Luc Illusie and Michael Temkin\footnote{The research of M.T. was partially supported by the
European Union Seventh Framework Programme (FP7/2007-2013) under grant agreement 268182.}
\end{center}

\tableofcontents

\else
\setboolean{english}{true}
\setboolean{francais}{false}
\selectlanguage{english}
\chapter{Gabber's modification theorem
(absolute case)}
\addtocontents{toc}{par Luc Illusie et Michael Temkin}
\begin{center}
Luc Illusie and Michael Temkin
\footnote{The research of M.T. was partially supported by the European Union Seventh
Framework Programme (FP7/2007-2013) under grant agreement 268182.}

\end{center}
\fi

In this expos\'e we state and prove Gabber's modification theorem
mentioned in the introduction (see step (C)). Its main application is
to Gabber's refined -- i.e. prime to $\ell$ -- local uniformization
theorem. This is treated in expos\'e IX. A relative variant of
the modification theorem, also due to Gabber, has applications to
prime to $\ell$ refinements of theorems of de Jong on alterations of
schemes of finite type over a field or a trait. This is discussed
in expos\'e X.  In \S1, we state Gabber's modification theorem in
its absolute form (Theorem \ref{8-TH. MODIFICATION}). The proof
of this theorem occupies \S\S\ref{8-prelimsec}--\ref{8-endsec}. A
key ingredient is the existence of functorial (with respect to
regular morphisms) resolutions in characteristic zero; the relevant
material is collected in \S\ref{8-functorressec}. We apply it in
\S\ref{8-logregsec} to get resolutions of log regular log schemes,
using the language of Kato's fans and Ogus's monoschemes. The main
results, on which the proof of \ref{8-TH. MODIFICATION} is based,
are Theorems \ref{8-monoidaldesth} and \ref{8-Morel 2.4.2}. \S\S
\ref{8-functorressec} and \ref{8-logregsec} can be read independently
of \S\S \ref{8-mainsec}, \ref{8-prelimsec}, \ref{8-endsec}.

We wish to thank Sophie Morel for sharing with us her notes
on resolution of log regular log schemes and Gabber's magic box.
They were quite useful.

\section{Statement of the main theorem}\label{8-mainsec}
\begin{theorem}[Gabber]\label{8-TH. MODIFICATION} Let $X$ be a
noetherian, qe, separated, log regular fs log scheme (\refext{VI}{6-1.2}),
endowed with an action of a finite group $G$. We assume that $G$
acts tamely (\refext{VI}{6-3.1}) and generically freely on $X$
(i.e. there exists a $G$-stable, dense open subset of $X$ where
the inertia groups $G_{\overline{x}}$ are trivial). Let $Z$ be the
complement of the open subset of triviality of the log structure of
$X$, and let $T$ be the complement of the largest $G$-stable open
subset of $X$ over which $G$ acts freely. Then there exists an fs log
scheme $X'$ and a  $G$-equivariant morphism $f=f_{(G,X,Z)}{\colon}
X' \to X$ of log schemes having the following properties:

(i) As a morphism of schemes, $f$ is a projective modification,
i.e. $f$ is projective and induces an isomorphism of dense open
subsets.

(ii) $X'$ is log regular and $Z' = f^{-1}(Z \cup T)$ is the
complement of the open subset of triviality of the log structure
of $X'$.

(iii) The action of $G$ on $X'$ is very tame (\refext{VI}{6-3.1}).
\end{theorem}

When proving the theorem we will construct $f_{(G,X,Z)}$ that
satisfies a few more nice properties that will be listed in Theorem
\ref{8-addendum}. We remark that Gabber also proves the theorem, 
more generally, when $X$ is not assumed to be qe. However, 
the quasi-excellence assumption simplifies the proof so we impose it here.
Most of the proof works for a general noetherian $X$, so we will assume 
that $X$ is qe only when this will be needed in \S\ref{8-endsec}.

\subsection{}\label{8-remarks main th.} (a) Note that we do
not demand that $f$ is log smooth. In general, it is not. Here
is an example. Let $k$ be an algebraically closed field of
characteristic $\ne 2$. Let $G = \{\pm 1\}$ act on the affine plane
$X = \mathbf{A}^2_k$, endowed with the trivial log structure, by
$x \mapsto \pm x$. Then $X$ is regular and log regular, and $T=
\{0\}$. The action of $G$ on $X$ is tame, but not very tame, as
$G_{\{0\}}$ (= $G$) does not act trivially on the (only) stratum $X$
of the stratification by the rank of $M^{gp}/\OO^*$. Let
$f\colon X' \to X$ be the blow up of $T$, with its natural action
of $G$. Then the pair $(X',Z' = f^{-1}(T))$ is log regular, $f$
is a $G$-equivariant morphism of log schemes, $X'-Z'$ is at the
same time the open subset of triviality of the log structure and
the largest $G$-stable open subset of $X'$ where $G$ acts freely,
and $G$ acts very tamely on $X'$. However, $f$ is not log smooth
(the fiber of $f$ at $\{0\}$ is the line $Z'$ with the log structure
associated to $\mathbf{N} \to \OO, 1 \mapsto 0$, which is
not log smooth over $\Spec\ \!  k$ with the trivial log structure).

(b) In the above example, let $D_1$, $D_2$ be distinct lines
in $X$ crossing at $\{0\}$, and put the log structure $M(D)$ on
$X$ defined by the divisor with normal crossings $D = D_1 \cup
D_2$. Then $\tilde X = (X,M(D))$ is log smooth over $\Spec\ \! k$
endowed with the trivial log structure, and $G$ acts very tamely
on $(X,M(D))$ (\refext{VI}{6-dcn stricts}).  Moreover, the modification $f$
considered above underlies the log blow up $\tilde f\colon \tilde X'
\to \tilde X$ of $X$ at (the ideal in $M(D)$) of $\{0\}$. While
$f$ depends only on $X$, the log \'etale morphism $\tilde f$
is not canonical, as it depends on the choice of $D$. However,
one can recover $f$ from the \textit{canonical} resolutions of
toric singularities (discussed in the next section). Namely,
as $G$ acts very tamely on $\tilde X$, the quotient $\tilde Y =
\tilde X/G$ is log regular (\refext{VI}{6-3.2}): $\tilde Y = \Spec\
\! k[P]$, where $P$ is the submonoid of $\mathbf{Z}^2$ generated
by $(2,0)$, $(1,1)$ and $(0,2)$, and the projection $p\colon
\tilde X \to \tilde Y$ is a Kummer \'etale cover of group $G$,
in particular is a $G$-\'etale cover of $V = p(U)$, where $U = X -
D$. Let $g\colon \tilde Y' \to \tilde Y$ be the log blow up of $\{0\}
= p(\{0\})$ in $\tilde Y$. We then have a cartesian diagram of log
schemes $$ \xymatrix{\tilde X' \ar[d] \ar[r]& \tilde Y' \ar[d]^g \\
\tilde X \ar[r]^{p} & \tilde Y}, $$ where the horizontal maps are
Kummer \'etale covers of group $G$. Now, as a morphism of schemes,
$\tilde Y' \to \tilde Y$ is the \textit{canonical resolution}
of $\tilde Y$, and the underlying scheme $X'$ of $\tilde X'$
is the \textit{normalization} of $\tilde Y'$ in the $G$-\'etale
cover $p\colon U \to V$. This observation, suitably generalized,
plays a key role in the proof of \ref{8-TH. MODIFICATION}.

\section{Functorial resolutions}\label{8-functorressec} All schemes
considered from now on will be assumed to be noetherian.

\subsection{Towers of blow ups}\label{8-blowtowersec} In this
section we review various known results on the following related
topics: blow ups and their towers, various operations on towers,
such as strict transforms and pushforwards, associated points of
schemes and schematic closure.

\subsubsection{Blow ups}\label{8-blowsec} We start with
recalling basic properties of blow ups; a good reference is
\cite[\S1]{Nagata@Conrad}. Let $X$ be a scheme. By a \textit{blow
up} of $X$ we mean a triple consisting of a morphism $f\colon Y
\to X$, a closed subscheme $V$ of $X$ (the \textit{center}), and
an $X$-isomorphism $\alpha\colon Y\iso \mathrm{Proj}(\oplus_{n
\in \mathbf{N}}\calI^n)$, where $\calI=\calI(V)$ is the ideal
of $V$.  We will write $Y=\Bl_V(X)$. When there is no risk of
confusion we will omit $V$ and $\alpha$ from the notation. A blow
up $(f,V,\alpha)$ is said to be \textit{empty} if $V=\emptyset$. In
this case, $Y=X$ and $f$ is the identity. The only $X$-automorphism
of a blow up is the identity. Also, it is well known that $Y$ is
the universal $X$-scheme such that $V\times_XY$ is a Cartier divisor
(i.e. the ideal $\calI\calO_Y$ is invertible).

\subsubsection{Total and strict transforms}\label{8-transformsec}
Given a blow up $f{\colon}Y=\Bl_V(X)\to X$, there are two natural
ways to pullback closed subschemes $i{\colon}Z\into X$. The {\em
total transform} of $Z$ under $f$ is the scheme-theoretic preimage
$f^\tot(Z)=Z\times_XY$. The {\em strict transform} $f^\st(Z)$
is defined as the schematic closure of $f^{-1}(Z-V)\iso Z-V$.

\begin{remark2}\label{8-strictrem} (i) The strict transform depends
on the centers and not only on $Z$ and the morphism $Y\to X$. For
example, if $D\into X$ is a Cartier divisor then the morphism
$\Bl_D(X)\to X$ is an isomorphism but the strict transform of $D$
is empty.

(ii) While $f^\tot(Z)\to Z$ is just a proper morphism, the morphism
$f^\st(Z)\to Z$ can be provided with the blow up structure because
$f^\st(Z)\iso\Bl_{V\times_XZ}(Z)$ (e.g., if $Z\into V$ then
$f^\st(Z)=\emptyset=\Bl_Z(Z)$). Thus, the strict transform can be
viewed as a genuine blow up pullback of $f$ with respect to $i$.
\end{remark2}

\subsubsection{Towers of blow ups}\label{8-towerssec} Next we
introduce blow up towers and study various operations with them
(see also \cite[\S2.2]{Temkin@qe}). By a \textit{tower of blow ups}
of $X$ we mean a finite sequence of length $n\ge 0$ $$ X_{\bullet}
= (\xymatrix{ X_n \ar[r]^{f_{n-1}} & X_{n-1} \ar[r] &\cdots
\ar[r]^{f_0} & X_0 = X}) $$ of blow ups. In particular, this data
includes the centers $V_i\into X_i$ for $0\le i\le n-1$. Usually,
we will denote the tower as $X_\bullet$ or $(X_\bullet,V_\bullet)$. Also,
we will often use notation $X_n\dashto X_0$ to denote a sequence
of morphisms.

If $n=0$ then we say that the tower is {\em trivial}. Note that the
morphism $X_n\to X$ is projective, and it is a modification if and
only if the centers $V_i$ are nowhere dense. If $X_{\bullet}$ is a
tower of blow ups, we denote by $(X_\bullet)_c$ the {\em contracted}
tower deduced from $X_{\bullet}$ by omitting the empty blow ups.

\subsubsection{Strict transform of a tower}\label{8-strictsec}
Assume that $\calX=(X_\bullet,V_\bullet)$ is a blow up tower
of $X$ and $h{\colon}Y\to X$ is a morphism. We claim that
there exists a unique blow up tower $\calY=(Y_\bullet,W_\bullet)$
of $Y$ such that $Y_i\to Y\to X$ factors through $X_i$ and
$W_i=V_i\times_{X_i}Y_i$. Indeed, this defines $Y_0,W_0$ and $Y_1$
uniquely. Since $V_0\times_XY_1=W_0\times_{Y_0}Y_1$ is a Cartier
divisor, $Y_1\to X$ factors uniquely through $X_1$. The morphism
$Y_1\to X_1$ uniquely defines $W_1$ and $Y_2$, etc. We call $\calY$
the {\em strict transform} of $\calX$ with respect to $h$ and denote
it $h^\st(\calX)$.

\begin{remark2}\label{8-strtowerrem} The following observation
motivates our terminology: if $h{\colon}Y\into X$ is a closed
immersion then $Y_i\into X_i$ is a closed immersion and $Y_{i+1}$
is the strict transform of $Y_i$ under the blow up $X_{i+1}\to X_i$.
\end{remark2}

\subsubsection{Pullbacks}\label{8-pullsec} One can also define
a naive base change of $\calX$ with respect to $h$ simply as
$(Y_\bullet,W_\bullet)=\calX\times_XY$. This produces a sequence
of proper morphisms $Y_n\dashto Y_0$ and closed subschemes
$W_i\into Y_i$ for $0\le i\le n-1$. If this datum is a blow up
sequence, i.e. $Y_{i+1}\iso\Bl_{W_i}(Y_i)$, then we say that
$\calX\times_XY$ is the {\em pullback} of $\calX$ and use the
notation $h^*(\calX)=\calX\times_XY$.

\begin{remark2}\label{8-pullrem} The pullback exists if and only if
$\calX\times_XY\iso h^\st(\calX)$. Indeed, this is obvious for
towers of length one, and the general case follows by induction on
the length.  \end{remark2}

\subsubsection{Flat pullbacks}\label{8-flatpullsec} Blow ups are
compatible with flat base changes $h\colon Y\to X$ in the sense that
$\Bl_{V\times_XY}(Y)\iso\Bl_V(X)\times_XY$ (e.g. just compute
these blow ups in the terms of $\Proj$). By induction on length
of blow up towers it follows that pullbacks of blow up towers with
respect to flat morphisms always exist. One can slightly strengthen
this fact as follows.

\begin{remark2}\label{8-nonflatrem} Assume that $X_\bullet$ is a
blow up tower of $X$ and $h{\colon}Y\to X$ is a morphism. If there
exists a flat morphism $f{\colon}X\to S$ such that the composition
$g{\colon}Y\to S$ is flat and the blow up tower $X_\bullet$ is
the pullback of a blow up tower $S_\bullet$ then the pullback
$\h^*(X_\bullet)$ exists and equals to $g^*(S_\bullet)$.  \end{remark2}

\subsubsection{Equivariant blow ups}\label{8-equivsec} Assume
that $X$ is an $S$-scheme acted on by a flat $S$-group scheme
$G$. We will denote by $p,m\colon X_0=G\times_SX\to X$ the
projection and the action morphisms. Assume that $V\into X$ is
a $G$-equivariant closed subscheme (i.e., $V\times_{X}(X_0;m)$
coincides with $V_0=V\times_{X}(X_0;p)$) then the action
of $G$ lifts to the blow up $Y=\Bl_V(X)$. Indeed, the blow up
$Y_0=\Bl_{V_0}(X_0)\to X_0$ is the pullback of $Y\to X$ with respect
to both $m$ and $p$, i.e. there is a pair of cartesian squares $$
\xymatrix{Y_0 \ar@<0.5ex>[d]^{m'}\ar@<-0.5ex>[d]_{p'} \ar[r]& X_0
\ar@<0.5ex>[d]^m\ar@<-0.5ex>[d]_p \\ Y \ar[r] & X} $$ So, we obtain
an isomorphism $Y_0\iso G\times_SY$ (giving rise to the projection
$p'\colon Y_0\to Y$) and a group action morphism $m'\colon Y_0\to Y$
compatible with $m$. Furthermore, the unit map $e{\colon}X\to X_0$
satisfies the condition of Remark \ref{8-nonflatrem} (with $X=S$),
hence we obtain the base change $e'{\colon}Y\to Y_0$ of $e$.

It is now straightforward to check that $m'$ and $e'$ satisfy
the group action axioms, but let us briefly spell this out
using simplicial nerves. The action of $G$ on $V$ defines a
cartesian sub-simplicial scheme $\Ner(G,V)\into\Ner(G,X)$. By
the flatness of $G$ over $S$ and Remark \ref{8-nonflatrem},
$\Bl_{\Ner(G,V)}(\Ner(G,X))$ is cartesian over $\Ner(G)$, hence
corresponds to an action of $G$ on $\Bl_V(X)$.

\subsubsection{Flat monomorphism}\label{8-flatmonsec} Flat
monomorphisms are studied in \cite{Ray@Raynaud}. In particular, it is
proved in \cite[Prop. 1.1]{Ray@Raynaud} that $i:Y\into X$ is a flat
monomorphism if and only if $i$ is injective and for any $y\in Y$ the
homomorphism $\calO_{X,y}\to\calO_{X,x}$ is an isomorphism. Moreover,
it is proved in \cite[Prop. 1.2]{Ray@Raynaud} that in this case $i$
is a topological embedding and $\calO_Y=\calO_X|_{Y}$. In addition
to open immersions, the main source of flat monomorphisms for us
will be morphisms of the form $\Spec(\calO_{X,x})\into X$ and their
base changes.

\subsubsection{Pushouts of ideals}\label{8-pushidealsec}
Let $i{\colon}Y\to X$ be a flat monomorphism,
e.g. $\Spec(\calO_{X,x})\into X$. By the {\em pushout} $U=i_*(V)$
of a closed subscheme $V\into Y$ we mean its schematic image in
$X$, i.e. $U$ is the minimal closed subscheme such that $V\into X$
factors through $U$. It exists by \ega{I}{9.5.1}.

\begin{lemma2}\label{8-pushideallem} The pushout $U=i_*(V)$ extends
$V$ in the sense that $U\times_XY=V$.  \end{lemma2} \begin{proof} By
minimality of $U$, it suffices to show that $V$ admits some extension
to a closed subscheme of $X$. Furthermore, it suffices to extend $V$
to a neighborhood $X'$ of $i(Y)$ because closed subschemes of $X'$
can be extended to $X$ by the schematic image. Thus we can replace
$X$ with $X'$, in particular, by \cite[1.5(i)]{Ray@Raynaud} we
can assume that $i{\colon}Y\to X$ is affine. Then $i$ is a limit
of affine finite type morphisms $X_i\to X$ by \ega{I}{9.6.6},
and by \egalong{IV}{3}{8.8.2 and 8.10.5(iv)} $V$ is the preimage of a
closed subscheme $V_i\into X_i$ for large enough $i$. We claim that
the morphism $Y\to X_i$ factors through a small enough neighborhood
$X'$ of $i(Y)$. Indeed, by quasi-compactness of $Y$ it suffices to
check that for any $y\in Y$ the morphism $\Spec(\calO_{X,y})\to X_i$
factors through a neighborhood of $y$. Since $\Spec(\calO_{X,y})$ is
the limit of all affine neighborhoods of $y$ the latter follows from
\egalong{IV}{3}{8.13.2}. Having such a factoring $Y\to X'\to X_i$ it
remains to note that $U$ is the preimage of $V'=V_i\times_{X_i}X'$,
hence $U$ extends to $V'$ in $X'$, and as we observed earlier this
implies the lemma.  \end{proof}

\subsubsection{Pushouts of blow up towers}\label{8-pushblowsec}
Given a blow up $f{\colon}\Bl_V(Y)\to Y$ and a flat monomorphism
$i{\colon}Y\into X$ we define the pushforward $i_*(f)$ as the blow
up along $U=i_*(V)$. Using Lemma \ref{8-pushideallem} and flat
pullbacks we see that $i^*i_*(f)=f$ and $\Bl_V(Y)\into\Bl_U(X)$ is
a flat monomorphism. So, we can iterate this procedure to construct
pushforward with respect to $i$ of any blow up tower $Y_\bullet$
of $Y$. It will be denoted  $(X_\bullet,U_\bullet)=i_*(Y_\bullet,V_\bullet)$.

\begin{remark2}\label{8-pushblowrem} (i) Clearly,
$i^*i_*(Y_\bullet)=Y_\bullet$.

(ii) In the opposite direction, a blow up tower $(X_\bullet,U_\bullet)$
of $X$ satisfies  $i_*i^*(X_\bullet)=X_\bullet$ if and only if the
preimage of $Y$ in each center $U_i$ of the tower is schematically
dense.  \end{remark2}

\subsubsection{Associated points of a scheme}\label{8-asschemesec}
Recall that a point $x\in X$ is called associated if $m_x$ is an
associated prime of $\calO_{X,x}$, i.e. $\calO_{X,x}$ contains an
element whose annihilator is $m_x$, see \egalong{IV}{2}{3.1.1}. The
set of all such points will be denoted $\Ass(X)$. The following
result is well known but difficult to find in the literature.

\begin{lemma2}\label{8-asschemelem} Let $i{\colon}Y\into X$ be
a flat monomorphism. Then the schematic image of $i$ coincides
with $X$ if and only if $\Ass(X)\subseteq i(Y)$.  \end{lemma2}
\begin{proof} Note that the schematic image of $i$ can be described
as $\Spec(\calF)$, where $\calF$ is the image of the homomorphism
$\phi{\colon}\calO_X\to i_*(\calO_Y)$. Thus, the schematic image
coincides with $X$ if and only if $\Ker(\phi)=0$.

If $i(Y)$ omits a point $x$ then any $m_x$-torsion element
$s\in\calO_{X,x}$ is in the kernel of $\phi_x{\colon}\calO_{X,x}\to
i_*(\calO_Y)_x$ (we use that $\calO_Y=\calO_X|_Y$ by
\S\ref{8-flatmonsec} and since $\ox$ is the support of $s$, the
restriction $s|_Y$ vanishes). So, if there exists $x\in\Ass(X)$
with $x\notin i(Y)$ then $\Ker(\phi)\neq 0$.

Conversely, if the kernel is non-zero then we take $x$ to
be any maximal point of its support and choose any non-zero
$s\in\Ker(\phi_x)$. In particular, $s|_{Y\cap X_x}=0$ and hence
$x\notin i(Y)$. For any non-trivial generalization $y$ of $x$ the
image of $s$ in $\calO_{X,y}$ vanishes because $\Ker(\calO_{X,y}\to
i_*(\calO_Y)_y)=0$ by maximality of $x$. Thus, $x$ is the support
of $s$, and hence $s$ is annihilated by a power of $m_x$. Since $X$
is noetherian, we can find a multiple of $s$ whose annihilator is
$m_x$, thereby obtaining that $x\in\Ass(X)$.  \end{proof}

\subsubsection{Associated points of blow up
towers}\label{8-astowersec} If $(X_\bullet,V_\bullet)$ is a blow up tower
of $X$ then by the set $\Ass(X_\bullet)$ of its associated points
we mean the union of the images of $\Ass(V_i)$ in $X$. Combining
Remark \ref{8-pushblowrem}(i) and Lemma \ref{8-asschemelem} we
obtain the following:

\begin{lemma2}\label{8-astowerlem} Let $i{\colon}Y\into X$ be a
flat monomorphism and let $X_\bullet$ be a blow up tower of $X$. Then
$i_*i^*(X_\bullet)=X_\bullet$ if and only if $\Ass(X_\bullet)\subset i(Y)$.
\end{lemma2}

\subsection{Normalized blow up towers}\label{8-normaltowersec} For
reduced schemes most of the notions, constructions and results of
\S\ref{8-blowtowersec} have normalized analogs. We develop such a
"normalized" theory in this section.

\subsubsection{Normalization}\label{8-normsec} The normalization of
a reduced noetherian scheme $X$, as defined in \ega{III}{6.3.8},
will be denoted $X^\nor$. Recall that normalization is compatible
with open immersions and for an affine $X=\Spec(A)$ its normalization
is $X^\nor=\Spec(B)$ where $B$ is the integral closure of $A$
in its total ring of fractions (which is a finite product of
fields). The normalization morphism $X^\nor\to X$ is integral but
not necessarily finite (since we do not impose the quasi-excellence
assumption on schemes).

\subsubsection{Functoriality}\label{8-functornormsec} Recall
(see expos\'e II) that a morphism $f{\colon}Y\to X$ is called
{\em maximally dominating} if it takes generic points of $Y$ to
generic points of $X$. Normalization is a functor on the category of
reduced schemes with maximally dominating morphisms. Furthermore,
it possesses the following universal property: any maximally
dominating morphism $Y\to X$ with normal $Y$ factors uniquely
through $X^\nor$. (Both claims are local on $Y$ and $X$ and are
obvious for affine schemes.)

\subsubsection{Normalized blow ups}\label{8-normblowsec} By the {\em
normalized blow up} of a reduced scheme $X$ along a closed subscheme
$V$ we mean the morphism $f{\colon}\Bl_V(X)^\nor\to X$. Note that $f$
is universally closed but does not have to be of finite type. As
in the case of usual blow ups, $V$ is a part of the structure. In
particular, $\Bl_V(X)$ has no $X$-automorphisms and we can talk
about equality of normalized blow ups (as opposed to an isomorphism).

\begin{proposition2}\label{8-normblowprop} (i) Keep the above
notation. Then $\Bl_V(X)^\nor\to X$ is the universal maximally
dominating morphism $Y\to X$ such that $Y$ is normal and $V\times_XY$
is a Cartier divisor.

(ii) For any blow up $f{\colon}Y=\Bl_V(X)\to X$ its normalization
$f^\nor{\colon}Y^\nor\to X^\nor$ is the normalized blow up along
$V\times_XX^\nor$.  \end{proposition2} \begin{proof} Combining the
universal properties of blow ups and normalizations we obtain (i),
and (ii) is its immediate corollary.  \end{proof}

Towers of normalized blow ups and their transforms can now be
defined similarly to their non-normalized analogs.

\subsubsection{Towers of normalized blow ups}\label{8-normtowersec}
A {\em tower of normalized blow ups} is a finite sequence $X_n\dashto
X_0$ of normalized blow ups with centers $V_i\into X_i$ for $0\le
i\le n-1$. The centers are part of the datum.

\subsubsection{Normalization of a blow up tower}\label{8-norblowsec}
Using induction on length and Proposition \ref{8-normblowprop}(ii),
we can associate to a blow up tower $\calX=(X_\bullet,V_\bullet)$
of a reduced scheme $X$ a normalized blow up tower
$\calX^\nor=(Y_\bullet,W_\bullet)$, where $Y_0=X$, $Y_i=X_i^\nor$ for
$i\ge 1$, and $W_i=V_i\times_{X_i}X^\nor_i$. We call $\calX^\nor$
the {\em normalization} of $\calX$.

\subsubsection{Strict transforms}\label{8-strnorsec}
If $\calX=(X_\bullet,V_\bullet)$ is a normalized blow up
tower of $X=X_0$ and $f{\colon}Y\to X$ is a morphism
then we define the strict transform $f^\st(\calX)$ as
the normalized blow up tower $(Y_\bullet,W_\bullet)$ such that
$Y_0=Y$ and $W_i=V_i\times_{X_i}Y_i$. Using induction on the
length and the universal property of normalized blow ups, see
\ref{8-normblowprop}(i), one shows that such a tower exists and is
the universal normalized blow up tower of $Y$ such that $f=f_0$
extends to a compatible sequence of morphisms $f_i{\colon}Y_i\to
X_i$.

\subsubsection{Pullbacks}\label{8-normpullsec} The strict transform
$f^\st(\calX)$ as above will be called the {\em pullback} and denoted
$f^*(\calX)$ if $Y_i\iso X_i\times_XY$ for any $0\le i\le n$.

\begin{lemma2}\label{8-smoothpulllem} If $f{\colon}Y\to X$ is
regular then any normalized blow up tower $\calX$ of $X$ admits
a pullback $f^*(\calX)$.  \end{lemma2} \begin{proof} Blow ups
are compatible with flat morphisms hence we should only show that
normalization is compatible with regular morphisms: if $f{\colon}Y\to
X$ is a regular morphism of reduced schemes then the morphism
$h{\colon}Y^\nor\to X^\nor\times_XY$ is an isomorphism. Note that
$h$ is an integral morphism which is generically an isomorphism and
the target is normal because it is regular over a normal scheme,
see \cite[21.E(iii)]{CA@Matsumura}. Hence $h$ is an isomorphism.
\end{proof}

\begin{remark2} Using \egalong{IV}{2}{6.14.1} instead of the
reference in the proof, one obtains that the claim of the lemma
holds, more generally, whenever $f$ is a normal morphism (i.e. flat
with geometrically normal fibres). We will use in the sequel only
the case when $f$ is regular.  \end{remark2}

\subsubsection{Fpqc descent of blow up towers}\label{8-descsubsec}
The classical fpqc descent of ideals (and modules) implies that
there is also an fpqc descent for blow up towers. Namely, if $Y\to
X$ is an fpqc covering and $Y_\bullet$ is a blow up tower of $Y$
whose both pullbacks to $Y\times_XY$ are equal then $Y_\bullet$
canonically descends to a blow up tower of $X$ because the centers
descend. In the same way, normalized blow up towers descend with
respect to quasi-compact surjective regular morphisms.

\subsubsection{Associated points}\label{8-assnorsec} The
material of \S\S\ref{8-pushblowsec}--\ref{8-astowersec} extends
to normalized blow up towers almost verbatim. In particular,
if $\calX=(X_\bullet,V_\bullet)$ is such a tower then $\Ass(\calX)$
is the union of the images of $\Ass(V_j)$ and for any flat
monomorphism $i{\colon}Y\into X$ (which is a regular morphism by
\S\ref{8-flatmonsec}) with a blow up tower $\calY$ of $Y$ we always
have that $i^*i_*\calY=\calY$, and we have that $i_*i^*\calX=\calX$
if and only if $\Ass(\calX)\subset i(Y)$.

\subsection{Functorial desingularization}\label{8-desingsec} In
this section we will formulate the desingularization result about
toric varieties that will be used later in the proof of Theorem
\ref{8-TH. MODIFICATION}. Then we will show how it is obtained from
known desingularization results.

\subsubsection{Desingularization of a scheme}\label{8-desschsec}
By a {\em resolution} (or {\em desingularization}) {\em tower} of
a scheme $X$ we mean a tower of blow ups with nowhere dense centers
$X_{\bullet}$ such that $X=X_0$, $X_n$ is regular and no $f_i$ is an
empty blow up. For example, the trivial tower is a desingularization
if and only if $X$ itself is regular.

\subsubsection{Normalized desingularization}\label{8-normdessec} We
will also consider normalized blow up towers such that each center
is non-empty and nowhere dense, $X=X_0$ and $X_n$ is regular. Such
a tower will be called a {\em normalized desingularization tower}
of $X$.

\begin{remark2}\label{8-normdesrem} (i) For any desingularization
tower $\calX$ of $X$ its normalization $\calX^\nor$ is a normalized
desingularization tower of $X$.

(ii) Usually one works with non-normalized towers; they are subtler
objects that possess more good properties. All known constructions of
functorial desingularization (see below) produce blow up towers by
an inductive procedure, and one cannot work with normalized towers
instead. However, it will be easier for us to deal with normalized
towers in log geometry because in this case one may work only with
fs log schemes.  \end{remark2}

\subsubsection{Functoriality of
desingularization}\label{8-functorsec} For concreteness, we will
consider desingularizations in \S\ref{8-functorsec}, but all what
we say holds for normalized desingularizations too. Assume that a
family $\calS^0$ of schemes is provided with desingularizations
$\calF(X)=X_\bullet$ for any $X\in\calS^0$. We say that the
desingularization (family) $\calF$ is {\em functorial} with respect
to a family $\calS^1$ of morphisms between the elements of $\calS^0$
if for any $f{\colon}Y\to X$ from $\calS^1$ the desingularization of
$X$ induces that of $Y$ in the sense that $f^*\calF(X)$ is defined
and its contraction coincides with $\calF(Y)$ (so, $\mathscr{F}(Y)
= (Y \times_X \mathscr{F}(X))_c$). Note that we put the $=$ sign
instead of an isomorphism sign, which causes no ambiguity by the
fact that any automorphism of a blow up is the identity as we
observed above.

\begin{remark2}\label{8-functorrem} (i) Contractions in the pulled
back tower appear when some centers of $\calF(X)$ are mapped to
the complement of $f(Y)$ in $X$. In particular, if $f\in\calS^1$
is surjective then the precise equality $\mathscr{F}(Y) = Y \times_X
\mathscr{F}(X)$ holds.

(ii) Assume that $X=\cup_{i=1}^nX_i$ is a Zariski covering and
the morphisms $X_i\into X$ and $\coprod_{i=1}^nX_i\to X$ are
in $\calS^1$. In general, one cannot reconstruct $\calF(X)$
from $\calF(X_i)$'s because the latter are contracted
pullbacks and it is not clear how to glue them with correct
synchronization. However, all information about $\calF(X)$ is
kept in $\calF(\coprod_{i=1}^nX_i)$. The latter is the pullback
of $\calF(X)$ hence we can reconstruct $\calF(X)$ by gluing the
restricted blow up towers $\calF(\coprod_{i=1}^nX_i)|_{X_i}$. Note
that $\calF(\coprod_{i=1}^nX_i)|_{X_i}$ can be obtained from
$\calF(X_i)$ by inserting empty blow ups, and these empty blow ups
make the gluing possible. This trick with synchronization of the
towers $\calF(X_i)$ by desingularizing disjoint unions is often
used in the modern desingularization theory, and one can formally
show (see \cite[Rem. 2.3.4(iv)]{Temkin@qe}) that such approach is
equivalent to the classical synchronization of the algorithm with
an invariant.

(iii) Assume that $\calS^1$ contains all identities $\Id_X$
with $X\in\calS^0$ and a morphism $Y\coprod Z\to X$ is in
$\calS^1$ whenever its restrictions $Y\to X$ and $Z\to X$ are in
$\calS^1$. As an illustration of the above trick, let us show that
even if $f,g{\colon}Y\to X$ are in $\calS^1$ but not surjective,
we have an equality $\calF(X)\times_X(Y,f)=\calF(X)\times_X(Y,g)$
of non-contracted towers. Indeed, set $Y'=Y\coprod X$ and consider
the morphisms $f',g'{\colon}Y'\to X$ that agree with $f$ and
$g$ and map $X$ by identity. Then $\calF(X)\times_X(Y',f')$
and $\calF(X)\times_X(Y',g')$ are equal because $f'$ and $g'$
are surjective, hence their restrictions onto $Y$ are also
equal, but these are precisely $\calF(X)\times_X(Y,f)$ and
$\calF(X)\times_X(Y,g)$.  \end{remark2}

\subsubsection{Gabber's magic box}\label{8-magicboxsec} Now we have
tools to formulate the aforementioned desingularization result.

\begin{theorem2}\label{8-magicboxth} Let $\calS^0$
denote the family of disjoint unions of affine toric
varieties over $\QQ$, i.e. $X\in\calS^0$ if it is of the form
$\coprod_{i=1}^n\Spec(\QQ[P_i])$, where $P_1\. P_n$ are fs torsion
free monoids. Let $\calS^1$ denote the family of smooth morphisms
$$f\colon\coprod_{j=1}^m\Spec(\QQ[Q_j])\to\coprod_{i=1}^n\Spec(\QQ[P_i])$$
such that for each $1\le j\le m$ there exists $1\le i=i(j)\le n$ and a
homomorphism of monoids $\phi_j\colon P_i\to Q_j$ so that the restriction
of $f$ onto $\Spec(\QQ[Q_j])$ factors through the toric morphism
$\Spec(\QQ[\phi_j])$. Then there exists a desingularization $\calF$
on $\calS^0$ which is functorial with respect to $\calS^1$ and,
in addition, satisfies the following compatibility condition:
if $\calO_1\.\calO_l$ are complete noetherian rings containing
$\QQ$, $Z=\coprod_{i=1}^l\Spec(\calO_i)$, and $g,h{\colon}Z\to X$
are two regular morphisms with $X\in\calS^0$ then \begin{equation}
\label{8-equalizer2a} (Z,g) \times_X \mathscr{F}(X) = (Z,h) \times_X
\mathscr{F}(X) \end{equation} \end{theorem2}

Before showing how this theorem follows from known desingularization
results, let us make few comments.

\begin{remark2} (i) Gabber's original magic box also requires that
the centers are smooth schemes. This (and much more) can also be
achieved as will be explained later, but we prefer to emphasize
the minimal list of properties that will be used in the proof of
Theorem \ref{8-TH. MODIFICATION}.

(ii) It is very important to allow disjoint unions in the theorem
in order to deal with synchronization issues, as explained in
Remark \ref{8-functorrem}(ii). This theme will show up repeatedly
throughout the paper.  \end{remark2}

\subsubsection{Desingularization of qe schemes
over $\mathbf{Q}$}\label{8-desingQsec} Gabber's magic box
\ref{8-magicboxth} is a particular case of the following
theorem, see \cite[Th. 1.2.1]{Temkin@qe}. Indeed, due to Remark
\ref{8-functorrem}(iii), functoriality with respect to regular
morphisms implies (\ref{8-equalizer2a}).

\begin{theorem2}\label{8-qeth} There exists a desingularization
algorithm $\calF$ defined for all reduced quasi-excellent schemes
over $\QQ$ and functorial with respect to all {\em regular}
morphisms. In addition, $\calF$ blows up only regular centers.
\end{theorem2}

\begin{remark2}
Although this is not stated in \cite{Temkin@qe}, one can strengthen
Theorem \ref{8-qeth} by requiring that $\calF$ blows up only regular centers
contained in the singular locus. An algorithm $\calF$ is constructed in
\cite{Temkin@qe} from an algorithm $\calF_\Var$ that desingularizes varieties
of characteristic zero, and one can check that if the centers of $\calF_\Var$ lie
in the singular loci (of the intermediate varieties) then the same is true
for $\calF$. Let us explain how one can choose an appropriate $\calF_\Var$.
In \cite{Temkin@qe}, one uses the algorithm of Bierstone-Milman to construct $\calF$,
see Theorem 6.1 and its Addendum in \cite{BMT@Bierstone-Milman-Temkin}
for a description of this algorithm and its properties.
It follows from the Addendum that the algorithm blows up centers lying in the
singular loci until $X$ becomes smooth, and then it performs some additional blow ups
to make the exceptional divisor snc. Eliminating the latter blow ups we obtain
a desingularization algorithm $\calF_\Var$ which only blows up regular
centers lying in the singular locus.
\end{remark2}

It will be convenient for us to use the algorithm $\calF$
from Theorem \ref{8-qeth} in the sequel. Also, to simplify the
exposition we will freely use all properties of $\calF$ but the
careful reader will notice that only the properties of Gabber's
magic box will be crucial in the end. Also, instead of working with
$\calF$ itself we will work with its normalization $\calF^\nor$
which assigns to a reduced qe scheme over $\QQ$ the normalized
blow up tower $\calF(X)^\nor$. It will be convenient to use the
notation $\tilcalF=\calF^\nor$ in the sequel.

\begin{remark2}\label{8-qenorth} (i) Since normalized blow ups
are compatible with regular morphisms, it follows from Theorem
\ref{8-qeth} that the normalized desingularization $\tilcalF$
is functorial with respect to all regular morphisms.

(ii) The feature which is lost under normalization (and which is
not needed for our purposes) is some control on the centers. The
centers $\tilV_i$ of $\tilcalF(X)$ are preimages of the centers
$V_i\into X_i$ of $\calF(X)$ under the normalization morphisms
$X_i^\nor\to X_i$, so they do not have to be even reduced. It will
only be important that $\tilV_i$'s are equivariant when a smooth
group acts on $X$. In the original Gabber's argument it was important
to blow up only regular centers because they were not part of the
blow up data, and one used that a regular center without codimension one
components intersecting the regular locus is determined already by the
underlying morphism of the blow up.  \end{remark2}

\subsubsection{Alternative desingularization
inputs}\label{8-altersec} For the sake of completeness, we discuss
how other algorithms could be used instead of $\calF$. Some
desingularization algorithms for reduced varieties over $\QQ$
are constructed in \cite{BM@resolution}, \cite{Wlodarczyk} and
\cite{Kollar@resolution}. They all are functorial with respect to
equidimensional smooth morphisms (though usually one "forgets"
to mention the equidimensionality restriction). It is shown in
\cite[\S6.3]{BMT@Bierstone-Milman-Temkin} how to make the algorithm
of \cite{BM@resolution} fully functorial by a slight adjusting of
the synchronization of its blow ups. There are also a few other canonical
constructions (e.g. due to Villamayor), which are probably functorial
with respect to equidimensional smooth morphisms too. All these
algorithms can be used to produce a desingularization of log regular
schemes (see \S\ref{8-logregsec}), so the only difficulty is in
establishing the compatibility (\ref{8-equalizer2a}).

For the algorithm of \cite{BMT@Bierstone-Milman-Temkin}
it was shown by Bierstone-Milman (unpublished, see
\cite[Rem. 7.1(2)]{BMT@Bierstone-Milman-Temkin}) that the induced
desingularization of a formal completion at a point depends only on
the formal completion as a scheme. This is precisely what we need in
(\ref{8-equalizer2a}).

Finally, there is a much more general result by Gabber, see
Theorem \ref{8-invarth}, whose proof uses Popescu's theorem and the
cotangent complex. It implies that, actually, any desingularization
of reduced varieties over $\QQ$ which is functorial with respect to
smooth morphisms automatically satisfies (\ref{8-equalizer2a}). So,
in principle, any functorial desingularization of varieties over
$\QQ$ could be used for our purposes. Since Gabber's result and
its proof are powerful and novel for the desingularization theory
(and were missed in \cite{BMT@Bierstone-Milman-Temkin}, mainly
due to a not so trivial involvement of the cotangent complex),
we include them in \S\ref{8-invarsec}.

\subsubsection{Invariance of the regular locus}\label{8-consfunsec}
Until the end of \S\ref{8-desingsec} we consider only qe schemes
of characteristic zero, and our aim is to establish a few useful
properties of $\calF$ (and $\tilcalF$) that are consequences of
the functoriality property $\calF$ satisfies. First, we claim that
$\calF$ does not modify the regular locus of $X$, and even slightly
more than that:

\begin{corollary2}\label{8-reglocuscor} All centers of $\calF(X)$
and $\tilcalF(X)$ sit over the singular locus of $X$. In particular,
$X$ is regular if and only if $\calF(X)$ is the trivial tower.
\end{corollary2} \begin{proof} It suffices to study $\calF$. The
claim is obvious for $S=\Spec(\QQ)$ because $S$ does not contain
non-dense non-empty subschemes. By functoriality, $\calF(T)$ is
trivial for any regular $T$ of characteristic zero, because it is
regular over $S$. Finally, if $T$ is the regular locus of $X$ then
$\calF(T)=(\calF(X)\times_XT)_c$ and hence any center $V_i\into X_i$
of $\calF(X)$ does not intersect the preimage of $T$.  \end{proof}

\subsubsection{Equivariance of the
desingularization}\label{8-equivdessec} It is well known that
functorial desingularization is equivariant with respect to
any smooth group action (and, moreover, extends to functorial
desingularization of stacks). For the reader's convenience we
provide an elementary argument.

\begin{corollary2}\label{8-cor BMT group action} Let $S$ be
a qe scheme over $\QQ$, $G$ be a smooth $S$-group and $X$
be a reduced $S$-scheme of finite type acted on by $G$. Then
the action of $G$ on $X$ extends naturally to an action of
$G$ on $\mathscr{F}(X)$ and $\tilcalF(X)$.  \end{corollary2}
\begin{proof} Again, it suffices to study $\calF$. Let $\calF(X)$
be given by $X_n\dashto X_0=X$ and $V_i\into X_i$ for $0\le i\le
n-1$. By $p,m\colon Y=G\times_S X\to X$ we denote the projection
and the action morphisms. Note that $m$ is smooth (e.g. $m$ is the
composition of the automorphism $(g,x)\to(g,gx)$ of $G\times_S X$ and
$p$). Therefore, $\calF(X)\times_{X}(Y;m)=\calF(Y)=\calF(X)\times_SG$
by Theorem \ref{8-qeth} and Remark \ref{8-functorrem}(i). In
particular, $V_0\times_X (Y; p) = V_0\times_X(Y; m)$, i.e. $V_0$
is $G$-equivariant. By \S\ref{8-equivsec}, $X_1$ inherits
a $G$-action. Then the same argument implies that $V_1$ is
$G$-equivariant and $X_2$ inherits a $G$-action, etc.  \end{proof}

\subsection{Complements on functorial
desingularizations}\label{8-invarsec} This section is devoted to
Gabber's result on a certain non-trivial compatibility property
that any functorial desingularization satisfies. It will not be
used in the sequel, so an uninterested reader may safely skip it.

\begin{theorem2}\label{8-invarth} Assume that $S$ is a scheme,
$\calS^0$ is a family of reduced $S$-schemes of finite type and
$\calS^1$ is a family of morphisms between elements of $\calS^0$
such that if $f{\colon}Y\to X$ is smooth and $X\in\calS^0$
then $Y\in\calS^0$ and $f\in\calS^1$. Let $\mathscr{F}$
be a desingularization on $\calS^0$ which is functorial with
respect to all morphisms of $\calS^1$. Then any pair of regular
morphisms $g{\colon}Z\to X$ and $h{\colon}Z\to Y$ with targets
in $\calS^0$ induces the same desingularization of $Z$; namely,
$\calF(X)\times_XZ=\calF(Y)\times_YZ$.  \end{theorem2}

Note that the theorem has no restrictions on the characteristic
(because no such restriction appears in Popescu's theorem). Before
proving the theorem let us formulate its important corollary,
whose main case is when $S=\Spec(k)$ for a field $k$ and $\calS^0$
is the family of all reduced $k$-schemes of finite type.

\begin{corollary2}\label{8-invarcor} Keep the notation of Theorem
\ref{8-invarth}. Then $\calF$ canonically extends to the family
$\widehat{\calS}^0$ of all schemes that admit a regular morphism to
a scheme from $\calS^0$ and the extension is functorial with respect
to all regular morphisms between schemes of $\widehat{\calS}^0$.
\end{corollary2}

The main ingredient of the proof will be the following result that
we are going to establish first.

\begin{proposition2}\label{8-invarprop} Consider a commutative
diagram of schemes $$ \xymatrix{{} & Z \ar[dl]_g \ar[dr]^h  \ar[d]^f
& {} \\ X \ar[dr]_a & Z' \ar[l]_{g'} \ar[r]^{h' }& Y \ar[dl]^b \\
{} & S & {}} $$ such that $a$ and $b$ are of finite type, $g$ and
$h$ are regular and $g'$ is smooth. Then $h'$ is smooth around the
image of $f$.  \end{proposition2}

For the proof we will need the following three lemmas. In the first
one we recall the Jacobian criterion of smoothness, rephrased in
terms of the cotangent complex.

\begin{lemma2}\label{8-invarlem1} Let $f{\colon} X \ra S$ be a
morphism which is locally of finite presentation, and let $x \in
X$. Then the following conditions are equivalent:

(i) $f$ is smooth at $x$;

(ii) $H_1(L_{X/S} \otimes^L k(x)) = 0$.  \end{lemma2} In the lemma
we use the convention $H_i = H^{-i}$, and $L_{X/S}$ denotes the
cotangent complex of $X/S$.  \begin{proof} (i) $\Rightarrow$
(ii) is trivial: as $f$ is smooth at $x$, up to shrinking $X$
we may assume $f$ smooth, then $L_{X/S}$ is cohomologically
concentrated in degree zero and locally free \cite[III
3.1.2]{Illusie@cotangent}. Let us prove (ii) $\Rightarrow$ (i).
We may assume that we have a factorization $$ \xymatrix{X \ar[r]^i
\ar[d]^f & Z \ar[dl]^g \\ S}, $$ where $i$ is a closed immersion of
ideal $I$ and $g$ is smooth. Consider the standard exact sequence
$$ I/I^2 \ra i^*\Omega^1_{Z/S} \ra \Omega^1_{X/S} \ra 0. \leqno
(*) $$ By the Jacobian criterion \egalong{IV}{4}{17.12.1} and
\egazero{IV}{19.1.12}, the smoothness of $f$ at $x$
is equivalent to the fact that the morphism $$ (I/I^2) \otimes k(x)
\ra \Omega^1_{Z/S} \otimes k(x) \leqno (**) $$ deduced from the left
one in (*) is injective.  Now, $(I/I^2) \otimes k(x) = H_1(L_{X/Z}
\otimes^L k(x))$ \cite[III 3.1.3]{Illusie@cotangent}, and (**)
is a morphism in the exact sequence associated with the triangle
deduced from the transitivity triangle $Li^*L_{Z/S} \ra L_{X/S}
\ra L_{X/Z} \ra Li^*L_{Z/S}[1]$ by applying $\otimes^L k(x)$ :
$$ H_1(L_{X/S} \otimes^L k(x)) \ra H_1(L_{X/Z} \otimes^L k(x))
(= (I/I^2) \otimes k(x)) \ra \Omega^1_{Z/S} \otimes k(x).  $$ By
(ii), $H_1(L_{X/S} \otimes^L k(x)) = 0$, hence (**) is injective,
which completes the proof.  \end{proof}

\begin{lemma2}\label{8-invarlem2} Consider morphisms $f{\colon}
X \ra Y$, $g{\colon} Y \ra S$, $h = gf{\colon} X \ra S$, and let
$x \in X$, $y=f(x)\in Y$. Assume that

(i) $H_1(L_{X/S} \otimes ^L k(x)) = 0$

(ii) $H_2(L_{X/Y} \otimes^L k(x)) = 0$.

Then $H_1(L_{Y/S} \otimes^L k(y)) = 0$. In particular, if $g$
is locally of finite presentation then $g$ is smooth at $y$.
\end{lemma2} \begin{proof}It is equivalent to show that $H_1(L_{Y/S}
\otimes^L k(x)) = 0$, and this follows trivially from the exact
sequence $$ H_2(L_{X/Y} \otimes^L k(x)) \ra H_1(L_{Y/S} \otimes^L
k(x)) \ra H_1(L_{X/S} \otimes^L k(x)).  $$ \end{proof}

\begin{lemma2}\label{8-invarlem3} Let $f {\colon} X \ra S$ be a
regular morphism between locally noetherian schemes. Then $L_{X/S}$
is cohomologically concentrated in degree zero and $H_0(L_{X/S})
= \Omega^1_{X/S}$ is flat.  \end{lemma2} \begin{proof} We may
assume $X = \SP B$ and $S = \SP A$ affine. Then, by Popescu's
theorem \cite[1.1]{NPD@Swan}, $X$ is a filtering projective limit
of smooth affine $S$-schemes $X_{\alpha} = \SP B_{\alpha}$. By
\cite[II (1.2.3.4)]{Illusie@cotangent}, we have $$ L_{B/A} =
\colim_{\alpha} L_{B_{\alpha}/A}.  $$ By \cite[III 3.1.2 and II
2.3.6.3]{Illusie@cotangent}, $L_{B_{\alpha}/A}$ is cohomologically
concentrated in degree zero and $H_0(L_{B_{\alpha}/A}) =
\Omega^1_{B_{\alpha}/A}$ is projective of finite type over
$B_{\alpha}$, so the conclusion follows.  \end{proof}

\begin{proof}[Proof of Proposition \ref{8-invarprop}] The question is
local around a point $y = f(x) \in Z'$, $x \in Z$. In view of Lemma
\ref{8-invarlem1}, by Lemma \ref{8-invarlem2} applied to $Z \ra Z'
\ra Y$ it suffices to show that $H_1(L_{Z/Y} \otimes^L k(x)) = 0$
and $H_2(L_{Z/Z'} \otimes^L k(x)) = 0$. As $Z$ is regular over $Y$,
the first vanishing follows from Lemma \ref{8-invarlem3}. For the
second one, consider the exact sequence
 $$
H_2(L_{Z/X} \otimes^L k(x)) \ra H_2(L_{Z/Z'} \otimes^L k(x)) \ra
H_1(L_{Z'/X} \otimes^L k(x)).  $$ By the regularity of $Z/X$ and
Lemma \ref{8-invarlem3}, $H_2(L_{Z/X} \otimes^L k(x)) = 0$. As $Z'$
is smooth over $X$, $H_1(L_{Z'/X} \otimes^L k(x)) = 0$ by Lemma
\ref{8-invarlem1}, which proves the desired vanishing and finishes
the proof.  \end{proof}

\begin{proof}[Proof of Theorem \ref{8-invarth}] Find finite
affine coverings $X=\cup_i X_i$, $Y=\cup_i Y_i$ and $Z=\cup_i
Z_i$ such that $g(Z_i)\subseteq X_i$ and $h(Z_i)\subseteq
Y_i$. Set $X'=\coprod_i X_i$, $Y'=\coprod_i Y_i$ and
$Z'=\coprod_i Z_i$ and let $Z'\to X'$ and $Z'\to Y'$ be the
induced morphisms. It suffices to check that $\calF(X)\times_XZ$
and $\calF(Y)\times_YZ$ become equal after pulling them back to
$Z'$. So, we should check that $(\calF(X)\times_XX')\times_{X'}Z'$
coincides with $(\calF(Y)\times_YY')\times_{Y'}Z'$. The
morphisms $X'\to X$ and $Y'\to Y$ are smooth and hence
contained in $\calS^1$. So, $\calF(X)\times_XX'=\calF(X')$
and similarly for $Y$. In particular, it suffices to prove that
$\calF(X')\times_{X'}Z'=\calF(Y')\times_{Y'}Z'$. This reduces the
problem to the case when all schemes are affine, so in the sequel
we assume that $X,Y$ and $Z$ are affine.

Next, note that it suffices to find factorizations $g = g_0f$ and $h
= h_0f$, where $f{\colon} Z \ra Z_0$ is a morphism with target in
$\calS^0$ and $g_0{\colon} Z_0 \ra X$, $h_0{\colon} Z_0 \ra X$ are
smooth. By Popescu's theorem, one can write $ g {\colon} Z \ra X$ as
a filtering projective limit of affine smooth morphisms $g_{\alpha}
{\colon} Z_{\alpha} \ra X$, $\alpha \in A$. As $Y$ is of finite
type over $S$, $h$ will factor through one of the $Z_{\alpha}$'s
(\egalong{IV}{3}{8.8.2.3}): there exists $\alpha \in A$, $f_{\alpha}
{\colon}Z \ra Z_{\alpha} $, $h_{\alpha} {\colon} Z_{\alpha} \ra Y$
such that $g = g_{\alpha}f_{\alpha}$, $h = h_{\alpha}f_{\alpha}$. By
Proposition \ref{8-invarprop}, $h_{\alpha}$ is smooth around the
image of $f_{\alpha}$, so we can take $Z_0$ to be a sufficiently
small neighborhood of the image of $f_\alpha$.  \end{proof}

\section{Resolution of log regular log schemes}\label{8-logregsec}
Unless said to the contrary, by log structure we mean
a log structure with respect to the \'etale topology. We
will say that a log structure $M_X$ on a scheme $X$ is {\em
Zariski} if $\varepsilon^*\varepsilon_*M_X\iso M_X$, where
$\varepsilon{\colon}X_\et\to X$ is the morphism between the \'etale
and Zariski sites. In this case, we can safely view the log structure
as Zariski log structure $\varepsilon_*M_X$. A similar convention
will hold also for log schemes.

\subsection{Fans}\label{8-fansec} Many definitions/constructions on
log schemes are of "combinatorial nature". Roughly speaking, these
constructions use only multiplication and ignore addition. Naturally,
there exists a category of geometric spaces whose structure
sheaves are monoids, and most of combinatorial constructions
can be described as "pullbacks" of analogous "monoidal"
operations. The first definition of such a category was done by
Kato in \cite{Toric@Kato}. Kato called his spaces {\em fans} to
stress their relation to the classical combinatorial fans obtained
by gluing polyhedral cones. For example, to any combinatorial fan
$C$ one can naturally associate a fan $F(C)$ whose set of points
is the set of faces of $C$. The main motivation for the definition
is that fans can be naturally associated to various log schemes.

It took some time to discover that fans are sort of "piecewise
schemes" rather than a monoidal version of schemes. A more geometric
version of combinatorial schemes was introduced by Deitmar in
\cite{Deitmar@F1}. He called them $\FF_1$-schemes, but we prefer
the terminology of monoschemes introduced by Ogus in his book
in preparation \cite{Ogus@book}. Note that when working with a
log scheme $X$, we use the sheaf $M_X$ in some constructions and
we use its sharpening $\oM_X$ (see \S\ref{8-notmonsec}) in other
constructions. Roughly speaking, monoschemes naturally arise when we
work with $M_X$ while fans naturally arise when we work with $\oM_X$.

In \S\ref{8-logregsec}, we will show that: (a) a functorial
desingularization of toric varieties over $\QQ$ descends to
a desingularization of monoschemes, (b) to give the latter is
more or less equivalent to give a desingularization of fans,
(c) desingularization of fans can be used to induce a monoidal
desingularization of log schemes, (d) the latter induces a
desingularization of log regular schemes, which (at least in some
cases) depends only on the underlying scheme.

In principle, we could work locally, using desingularization of
disjoint unions of all charts for synchronization. In this case, we
could almost ignore the intermediate categories by working only with
fine monoids and blow up towers of their spectra. However, we decided
to emphasize the actual geometric objects beyond the constructions,
and, especially, stress the difference between fans and monoschemes.

\subsubsection{Sharpening}\label{8-notmonsec} For a monoid $M$,
by $M^\times$ (or $M^*$) we denote the group of its invertible
elements, and its {\em sharpening} $\oM$ is defined as $M/M^\times$.

\subsubsection{Localization}\label{8-localmonsec} By {\em
localization} of a monoid $M$ along a subset $S$ we mean the
universal $M$-monoid $M_S$ such that the image of $S$ in $M_S$
is contained in $M_S^\times$. If $M$ is integral then $M_S$ is
simply the submonoid $M[S^{-1}]\subseteq M^\gp$ generated by $M$
and $S^{-1}$. If $M$ is a fine then any localization is isomorphic
to a localization at a single element $f$, and will be denoted $M_f$.

\subsubsection{Spectra of fine monoids}\label{8-specmonoidsec} All
our combinatorial objects will be glued from finitely many spectra of
fine monoids. Recall that with any fine monoid $P$ one can associate
the set $\Spec(P)$ of prime ideals (with the convention that
$\emptyset$ is also a prime ideal) equipped with the Zariski topology
whose basis is formed by the sets $D(f)=\{p\in\Spec(P)|\ f\notin
p\}$ for $f\in P$, see, for example, \cite[\S9]{Toric@Kato}. The
structure sheaf ${M}_P$ is defined by ${M}_P(D(f))=P_f$, and the
sharp structure sheaf $\oM_P={M}_P/{M}_P^\times$ is the sharpening
of ${M}_P$ (we will see in Remark \ref{8-specrem}(iii) that actually
$\oM_P(D(f))=\overline{P_f}=P_f/(P_f^\times)$, i.e. no sheafification
is needed).

\begin{remark2}\label{8-specrem} (i) Since $P\setminus P^\times$
and $\emptyset$ are the maximal and the minimal prime ideals of
$P$, $\Spec(P)$ possesses unique closed and generic points $s$ and
$\eta$. The latter is the only point whose stalk $M_{P,\eta}=P^\gp$
is a group.

(ii) The set $\Spec(P)$ is finite and its topology is the
specialization topology, i.e. $U$ is open if and only if it is
closed under generalizations. (More generally, this is true for any
finite sober topological space, such as a scheme that has finitely
many points.)

(iii) A subset $U\subseteq\Spec(P)$ is affine (and even of the form
$D(f)$) if and only if it is the localization of $\Spec(P)$ at a
point $x$ (i.e. the set of all generalizations of $x$). Any open
covering $U=\cup_i U_i$ of an affine set is trivial (i.e. $U$
is equal to some $U_i$), therefore any functor $\calF(U)$
on affine sets uniquely extends to a sheaf on $\Spec(P)$. In
particular, this explains why no sheafification is needed when
defining $\oM_P$. Furthermore, we see that, roughly speaking, any
notion/construction that is "defined in terms of" localizations
$X_x$ and stalks $M_x$ or $\oM_x$ is Zariski local. This is very
different from the situation with schemes.  \end{remark2}

\subsubsection{Local homomorphisms of monoids}\label{8-localhomsec}
Any monoid $M$ is local because $M\setminus M^\times$ is its unique
maximal ideal. A homomorphism $f{\colon}M\to N$ of monoids is {\em
local} if it takes the maximal ideal of $M$ to the maximal ideal
of $N$. This happens if and only if $f^{-1}(N^\times)=M^\times$.

\subsubsection{Monoidal spaces}\label{8-monoidspsec} A {\em
monoidal space} is a topological space $X$ provided with
a sheaf of monoids $M_X$. A morphism of monoidal spaces
$(f,f^\#){\colon}(Y,M_Y)\to(X,M_X)$ is a continuous map
$f{\colon}Y\to X$ and a homomorphism $f^\#{\colon}f^{-1}(M_X)\to
M_Y$ such that for any $y\in Y$ the homomorphism of monoids
$f_y^\#{\colon}{M}_{X,f(y)}\to{M}_{Y,y}$ is local.

\begin{remark2} Strictly speaking one should have called the above
category the category of locally monoidal spaces and allow non-local
homomorphisms in the general category of monoidal spaces. However,
we will not use the larger category, so we prefer to abuse the
notation slightly.  \end{remark2}

Spectra of monoids possess the usual universal property, namely:

\begin{lemma2}\label{8-monoidsplem} Let $(X,M_X)$ be a monoidal
space and $P$ be a monoid.

(i) The global sections functor $\Gamma$ induces
a bijection between morphisms of monoidal spaces
$(f,f^\#){\colon}(X,M_X)\to(\Spec(P),M_P)$ and homomorphisms
$\phi{\colon}P\to\Gamma(M_X)$.

(ii) If $M_X$ has sharp stalks then $\Gamma$ induces
a bijection between morphisms of monoidal spaces
$(f,f^\#){\colon}(X,M_X)\to(\Spec(P),\oM_P)$ and homomorphisms
$\phi{\colon}P\to\Gamma(M_X)$.  \end{lemma2} \begin{proof} (i)
Let us construct the opposite map. Given a homomorphism $\phi$,
for any $x\in X$ we obtain a homomorphism $\phi_x{\colon}P\to
M_{X,x}$. Clearly, $m=P\setminus\phi_x^{-1}(M^\times_{X,x})$ is a
prime ideal and hence $\phi_x$ factors through a uniquely defined
local homomorphism $P_m\to M_{X,x}$. Setting $f(x)=m$ we obtain
a map $f{\colon}X\to\Spec(P)$, and the rest of the proof of (i)
is straightforward.

If the stalks of $M_X$ are sharp then any morphism
$(X,M_X)\to(\Spec(P),M_P)$ factors uniquely through
$(\Spec(P),\oM_P)$. Also, $\Gamma(M_X)$ is sharp, hence
any homomorphism to it from $P$ factors uniquely through
$\oP$. Therefore, (ii) follows from (i).  \end{proof}

\subsubsection{Fine fans and monoschemes}\label{8-finfansec} A
fine {\em monoscheme} (resp. a fine {\em fan}) is a monoidal space
$(X,{M}_X)$ that is locally isomorphic to $A_P=(\Spec(P),{M}_P)$
(resp. $\oA_P=(\Spec(P),\oM_P)$), where $P$ is a fine monoid. We
say that $(X,M_X)$ is {\em affine} if it is isomorphic to $A_P$
(resp.f $\oA_P$). A morphism of monoschemes (resp. fans) is a
morphism of monoidal spaces. A monoscheme (resp. a fan) is called
{\em torsion free} if it is covered by spectra of $P$'s with torsion
free $P^\gp$'s. It follows from Remark \ref{8-specrem}(iii) that this
happens if and only if all groups $M^\gp_{X,x}$ are torsion free.

\begin{remark2}\label{8-fanrem} (i) Any fs fan is torsion
free because if an fs monoid is torsion free then so is any its
localization. This is not true for general fine fans. For example,
if $\mu_2=\{\pm 1\}$ then $P=\NN\oplus\mu_2\setminus\{(0,-1)\}$
is a sharp monoid with $P^\gp=\ZZ\oplus\mu_2$.

(ii) For any point $x$ of a fine monoscheme (resp. fan) $X$ the
localization $X_x$ that consists of all generalizations of $x$ is
affine. In particular, by Remark \ref{8-specrem}(i) there exists
a unique maximal point generalizing $x$, and hence $X$ is a finite
disjoint union of irreducible components.  \end{remark2}

\subsubsection{Comparison of monoschemes and
fans}\label{8-comparesec} There is an obvious sharpening functor
$(X,{M}_X)\mapsto(X,\oM_X)$ from monoschemes to fans, and there is
a natural morphism of monoidal spaces $(X,\oM_X)\to(X,M_X)$. The
sharpening functor loses information, and one needs to know
${M}_X^\gp$ to reconstruct ${M}_X$ from $\oM_X$ as a fibred product
(see \cite[2.1.8.3]{Ogus@book}). Actually, there are much more fans
than monoschemes. For example, the generic point $\eta\in\Spec(P)$
is open and $\oM_{P,\eta}$ is trivial hence for any pair of fine
monoids $P$ and $Q$ we can glue their sharpened spectra along the
generic points. What one gets is sort of "piecewise scheme" and,
in general, it does not correspond to standard geometric objects,
such as schemes or monoschemes. We conclude that, in general,
fans can be lifted to monoschemes only locally.

\begin{remark2} As a side remark we note that sharpened monoids
naturally appear as structure sheaves of piecewise linear spaces
(a work in progress of the second author on skeletons of Berkovich
spaces). In particular, PL functions can be naturally interpreted
as sections of the sharpened sheaf of linear functions on polytopes.
\end{remark2}

\subsubsection{Local smoothness}\label{8-localsmoothfansec}
A local homomorphism of fine monoids $\phi{\colon}P\to Q$
is called {\em smooth} if it can be extended to an isomorphism
$P\oplus\NN^r\oplus\ZZ^s\iso Q$. The following lemma checks
that this property is stable under localizations.

\begin{lemma2}\label{8-smoothmonoidlem} Assume that $\phi{\colon}P\to
Q$ is smooth and $P',Q'$ are localizations of $P,Q$ such that $\phi$
extends to a local homomorphism $\phi'{\colon}P'\to Q'$. Then $\phi'$
is smooth.  \end{lemma2} \begin{proof} Recall that $P'=P_a$ for $a\in
P$ (notation of \S\ref{8-localmonsec}), and $\phi'$ factors through
the homomorphism $\phi_a{\colon}P'\to Q_a$, which is obviously
smooth. Therefore, replacing $\phi$ with $\phi_a$ we can assume
that $P=P'$. Let $b=(p,n,z)\in Q$ be such that $Q'=Q_b$. Then $p\in
P^\times$ because $P\to Q'$ is local, and hence $Q'$ is isomorphic
to $P\oplus(\NN^r)_n\oplus\ZZ^s$. It remains to note that any
localization of $\NN^r$ is of the form $\NN^{r-t}\oplus\ZZ^t$.
\end{proof}

\subsubsection{Smoothness}\label{8-smoothfansec} The lemma
allows to globalize smoothness: a morphism $f{\colon}Y\to X$
of monoschemes is called {\em smooth} if the homomorphisms of
stalks ${M}_{X,f(y)}\to{M}_{Y,y}$ are smooth. In particular, $X$
is {\em smooth} if its morphism to $\Spec(1)$ is smooth, that is,
the stalks ${M}_{X,x}$ are of the form $\NN^r\oplus\ZZ^s$. In
particular, a smooth monoscheme is torsion free.

Analogous smoothness definitions are given for fans. Moreover,
in this case we can consider only sharp monoids, and then the
group component $\ZZ^s$ is automatically trivial. It follows
that we can rewrite the above paragraph almost verbatim but with
$s=0$. Obviously, a morphism of monoschemes is smooth if and only
if its sharpening is a smooth morphism of fans.

\begin{remark2}\label{8-chartsmoothrem} (i) Recall that
any fine monoscheme $X$ admits an open affine covering
$X=\cup_{x\in X}\Spec(M_{X,x})$. It follows that a morphism
of fine monoschemes $f{\colon}Y\to X$ is smooth if and
only if it is covered by open affine charts of the form
$\Spec(P\oplus\NN^r\oplus\ZZ^s)\to\Spec(P)$.

(ii) Smooth morphisms of fine fans admit a similar local description,
and we leave the details to the reader.  \end{remark2}

\subsubsection{Saturation}\label{8-satursec} As usually, for a
fine monoid $P$ we denote its saturation by $P^\sat$ (it consists
of all $x\in P^\gp$ with $x^n\in P$ for some $n>0$). Saturation is
compatible with localizations and sharpening and hence extends to
a saturation functor $F\mapsto F^\sat$ on the categories of fine
monoschemes (resp. fine fans). We also have a natural morphism
$F^\sat\to F$, which is easily seen to be bijective. So, actually,
$(F,M_F)^\sat=(F,M_F^\sat)$.

\subsubsection{Ideals}\label{8-idealsec} A subsheaf of ideals
$\calI\subseteq M_X$ on a monoscheme $(X,{M}_X)$ is called a
{\em coherent ideal} if for any point $x\in X$ the restriction of
$\calI$ on $X_x$ coincides with $\calI_xM_{X_x}$. (Due to Remark
\ref{8-specrem}(iii) this means that $\calI$ is coherent in the usual
sense, i.e. its restriction on an open affine submonoscheme $U$ is
generated by the global sections over $U$.) We will consider only
coherent ideals, so we will omit the word "coherent" as a rule. An
ideal $\calI\subseteq{M}_X$ is {\em invertible} if it is locally
generated by a single element.

\subsubsection{Blow ups}\label{8-blowmonosec} Similarly to schemes,
for any non-empty ideal $\calI\subseteq\calO_X$ there exists a
universal morphism of monoschemes $h{\colon}X'\to X$ such that the
{\em pullback ideal} $h^{-1}\calI=\calI{M}_{X'}$ is invertible. We
call $\calI$ the {\em center} of the blow up. (One does not have
an adequate notion of closed submonoscheme, so unlike blow ups
of the scheme it would not make sense that "$V(\calI)$" is the
center.) An explicit construction of $X'$ copies its scheme analog:
it is local on the base and for an affine $X=\Spec(P)$ with an ideal
$I\subseteq P$ corresponding to $\calI$ one glues $X'$ from the
charts $\Spec(P[a^{-1}I])$, where $a\in I$ and $P[a^{-1}I]$ is the
submonoid of $P^\gp$ generated by the fractions $b/a$ for $b\in I$
(see \cite[Ch. 2, \S1.6]{Ogus@book} for details).

\begin{remark2}\label{8-blowmonorem} Blow ups induce isomorphisms on
the stalks of $M^\gp$; this is an analog of the fact that blow ups
of schemes along nowhere dense subschemes are birational morphisms.
\end{remark2}

\subsubsection{Saturated blow ups}\label{8-satblowsec} Analogously
to normalized blow ups, one defines {\em saturated blow up} of a
monoscheme $X$ along an ideal $\calI\subset M_X$ as the saturation
of $\Bl_\calI(X)$. The same argument as for schemes shows that
$\Bl_\calI(X)^\sat$ is the universal saturated $X$-monoscheme such
that the pullback of $\calI$ is invertible.

\subsubsection{Towers and pullbacks}\label{8-towpullsec} Towers
of (saturated) blow ups of a monoscheme $X$ are defined in the
obvious way. Given such a tower $X_\bullet$ with $X=X_0$ and a
morphism $f{\colon}Y\to X$ we define the {\em pullback tower}
$Y_\bullet=f^*(X_\bullet)$ as follows: $Y_0=Y$ and $Y_{i+1}$ is the
(saturated) blow up of $Y_i$ along the pullback of the center
$\calI_i$ of $X_{i+1}\to X_i$. Due to the universal property of
(saturated) blow ups this definition makes sense and $Y_\bullet$
is the universal (saturated) blow up tower of $Y$ that admits a
morphism to $X_\bullet$ extending $f$.

\begin{remark2}\label{8-towpullrem} Unlike pullbacks of
(normalized) blow up towers of schemes, see \S\S\ref{8-pullsec}
and \ref{8-normpullsec}, we do not distinguish strict transforms
and pullbacks. The above definition of pullback covers our needs,
and we do not have to study the base change of monoschemes. For the
sake of completeness, we note that fibred products of monoschemes
exist and in the affine case are defined by amalgamated sums of
monoids, see \cite{Deitmar@F1}. Also, it is easy to check that for
a smooth $f$ (which is the only case we will use) one indeed has
that $f^*(X_\bullet)=X_\bullet\times_XY$ for any (saturated) blow up
tower $X_\bullet$. For blow ups one checks this with charts and in
the saturated case one also uses that saturation is compatible
with a smooth morphism $f{\colon}Y\to X$, i.e. $Y^\sat\iso
X^\sat\times_XY$.  \end{remark2}

\subsubsection{Compatibility with sharpening}\label{8-compsharpsec}
Ideals and blow ups of fans are defined in the same way, but with
$\oM_X$ used instead of ${M}_X$ (Kato defines their saturated version
in \cite[9.7]{Toric@Kato}). Towers of blow ups of fans are defined
in the obvious way.

\begin{lemma2}\label{8-idealfanlem} Let $X=(X,{M}_X)$ be a
monoscheme, let $(F,M_F)=(X,\oM_X)$ be the corresponding fan and
let $\lambda{\colon}M_X\to M_F$ denote the sharpening homomorphism.

(i) $\calI\mapsto\lambda(\calI)$ induces a natural bijection between
the ideals on $X$ and on $F$.

(ii) Blow ups are compatible with sharpening, that is,
the sharpening of $\Bl_\calI(X))$ is naturally isomorphic to
$\Bl_{\lambda(\calI)}(F)$. The same statement holds for saturated
blow ups.

(iii) Sharpening induces a natural bijection between the (saturated)
blow up towers of $X$ and $F$.  \end{lemma2} \begin{proof} (i) is
obvious. (ii) is shown by comparing the blow up charts. Combining
(i) and (ii), we obtain (iii).  \end{proof}

\subsubsection{Desingularization}\label{8-fandesingsec} Using
the above notions of smoothness and blow ups, one can copy
other definitions of the desingularization theory. By a {\em
desingularization} (resp. {\em saturated desingularization}) of
a fine monoscheme $X$ we mean a blow up tower (resp. saturated
blow up tower) $X_n\dashto X_0=X$ with smooth $X_n$. By Remark
\ref{8-blowmonorem}, if $X$ admits a desingularization then it is
torsion free, and we will later see that the converse is also true.

For concreteness, we consider below non-saturated desingularizations,
but everything extends to the saturated case verbatim. A family
$\calF^\mono(X)$ of desingularizations of torsion free monoschemes
is called {\em functorial} (with respect to smooth morphisms) if for
any smooth $f{\colon}Y\to X$ the desingularization $\calF^\mono(Y)$
is the contracted pullback of $\calF^\mono(X)$. The same argument
as for schemes (see Remark \ref{8-functorrem}(ii)) shows that
$\calF^\mono$ is already determined by its restriction onto the
family of finite disjoint unions of affine monoschemes.

The definition of a functorial desingularization $\calF^\fan$ of fine
torsion free fans is similar. Since blow up towers and smoothness
are compatible with the sharpening functor, it follows that a
desingularization of a monoscheme $X$ induces a desingularization of
its sharpening. Moreover, any affine fan can be lifted to an affine
monoscheme and $\calF^\fan$ is determined by its restriction onto
disjoint unions of affine fans, hence we obtain the following result.

\begin{theorem2}\label{8-fandesingth} The sharpening functor
induces a natural bijection between functorial desingularizations
of fine torsion free monoschemes and functorial desingularizations
of fine torsion free fans. A similar statement holds for saturated
desingularizations.  \end{theorem2}

\begin{remark2}\label{8-satdesrem} Similarly to the normalization
of a desingularization tower, to any desingularization $\calF$ of
monoschemes or fans one can associate a saturated desingularization
$\calF^\sat$: one replaces all levels of the towers, except the zero
level, with their saturations. In this case blow ups are replaced
with saturated blow ups along the pulled back ideals. If $\calF$
is functorial with respect to all smooth morphisms then the same is
true for $\calF^\sat$. Indeed, for any smooth $Y\to X$ the centers
of $\calF(Y)$, $\calF^\sat(X)$ and $\calF^\sat(Y)$ are the pullbacks
of those of $\calF(X)$. In addition, the saturation construction
is compatible with the bijections from Theorem \ref{8-fandesingth}
in the obvious way.  \end{remark2}

\begin{remark2}\label{8-fandesrem} In principle, (saturated)
desingularization of fans or monoschemes can be described in
purely combinatorial terms of fans and their subdivisions (e.g. see
\cite[9.6]{Toric@Kato} or \cite[\S4]{Toric@Niziol}). However, it is
not easy to construct a functorial one directly. We will instead
use a relation between monoschemes and toric varieties to descend
desingularization of toric varieties to monoschemes and fans.
\end{remark2}

\subsection{Monoschemes and toric varieties}\label{8-monotoricsec}

\subsubsection{Base change from monoschemes to
schemes}\label{8-monoliftsec} Let $S$ be a scheme. The following
proposition introduces a functor from monoschemes to $S$-schemes
that can be intuitively viewed as a base change with respect to a
"morphism" $S\to\Spec(1)$.

\begin{proposition2}\label{8-monoliftprop} Let $S$ be a scheme
and $F$ be a monoscheme. Then there exists an $S$-scheme $X=S[F]$
with a morphism of monoidal spaces $f{\colon}(X,\calO_X)\to(F,M_F)$
such that any morphism $(Y,\calO_Y)\to(F,M_F)$, where $Y$ is an
$S$-scheme, factors uniquely through $f$.  \end{proposition2}
\begin{proof} Assume, first, that $F=\Spec(P)$ is affine. By Lemma
\ref{8-monoidsplem}(i), to give a morphism $(Y,\calO_Y)\to(F,M_F)$
is equivalent to give a homomorphism of monoids
$\phi{\colon}P\to\Gamma(\calO_Y)$, and the latter factors uniquely
through a homomorphism of sheaves of rings $\calO_S[P]\to\calO_Y$. It
follows that $S[F]=\mathbf{\Spec}(\calO_S[P])$ in this case. Since
the above formula is compatible with localizations by elements $a\in
P$, i.e. $S[F_a]\iso S[F]_{\phi(a)}$, it globalizes to the case of
an arbitrary monoscheme. Thus, for a general monoscheme $F$ covered
by $F_i=\Spec(P_i)$, the scheme $S[F]$ is glued from $S[F_i]$.
\end{proof}

\begin{remark2}\label{8-monoliftremrem} Note that if $S=\Spec(R)$
and $F=\Spec(P)$ then $S[F]=\Spec(R[P])$. However, we will often
consider an "intermediate" situation where $S=\Spec(R)$ is affine
and $F$ is a general monoscheme. To simplify notation, we will
abuse them by writing $R[F]$ instead of $\Spec(R[F])$. Such "mixed"
notation will always refer to a scheme.  \end{remark2}

\subsubsection{Toric schemes}\label{8-toricfansec} If $F$ is torsion
free and connected then we call $S[F]$ a {\em toric scheme} over
$S$. Recall that by Remark \ref{8-fanrem}(ii), $F$ possesses a unique
maximal point $\eta=\Spec(P^\gp)$, where $\Spec(P)$ is any affine
open submonoscheme. Hence $X=S[F]$ possesses a dense open subscheme
$T=S[\eta]$, which is a split torus over $S$, and the action of $T$
on itself naturally extends to the action of $T$ on $X$.

\begin{remark2} Assume that $k$ is a field. Classically, a toric
$k$-variety is defined as a normal finite type separated $k$-scheme
$X$ that contains a split torus $T$ as a dense open subscheme such
that the action of $T$ on itself extends to the whole $X$. If
$F$ is saturated then our definition above is equivalent to the
classical one. However, we also consider non-normal toric varieties
corresponding to non-saturated monoids.  \end{remark2}

\subsubsection{Canonical log structures}\label{8-canonsec}
For any monoscheme $F$, the $S$-scheme $X=S[F]$ possesses
a natural log structure induced by the universal morphism
$f{\colon}(X,\calO_X)\to(F,M_F)$. Namely, $M_X$ is the log structure
associated with the pre-log structure $g^{-1}M_F\to\calO_{X_\et}$,
where $g{\colon}(X_\et,\calO_{X_\et})\to(X,\calO_X)\to(F,M_F)$ is
the composition. We call $M_X$ the {\em canonical} log structure
of $X=S[F]$.

\begin{remark2}\label{8-canonrem} (i) The canonical log
structure is Zariski, as $\varepsilon_* M_X$ coincides with
the Zariski log structure associated with the pre-log structure
$f^{-1}M_F\to\calO_X$.

(ii) The log scheme $(X=S[F],M_X)$ is log smooth over the scheme
$S$ provided with the trivial log structure. In particular,
if $S$ is regular and $F$ is saturated then $(X,M_X)$ is fs and
log regular. Without the saturation assumption we still have that
$X^\sat\iso S[F^\sat]$ is log regular, hence $X$ is log regular
in the sense of Gabber (see \S\ref{8-nonsatursec}).

(iii) If $\eta$ is the set of generic points of $F$ then $T=S[\eta]$
is the open subset of $X$ which is the triviality set of its
log structure. However, the map $M_{X} \to \OO_X
\cap j_* \OO_{T}^\times$ is not an isomorphism in
general, as the case where $T = \Spec \mathbf{C}[t,t^{-1}]
\subset X = \Spec \mathbf{C}[t^2,t^3]$ already shows: the
image of $t^2 + t^3$ in $\OO_{X,\{0\}}$ belongs to $(
j_*\OO_{T}^\times)_{\{0\}}$, but does not belong to
$M_{\{0\}} = t^{P}\OO_{X,\{0\}}^\times$, where $P$ is the
(fine, but not saturated) submonoid of $\mathbf{N}$ generated by 2
and 3. (We use that $\OO_{X,\{0\}}$ is strictly smaller than
its normalization $\mathbf{C}[t]_{(t)}=\OO_{X,\{0\}}[t]$
and hence $\frac{t^2+t^3}{t^2}=1+t$ is not contained in
$\OO_{X,\{0\}}$.)  \end{remark2}

\subsubsection{Toric saturation}\label{8-toricsatsec} Saturation of
monoschemes corresponds to normalization of schemes. This will play
an essential role later, since we get a combinatorial description
of the normalization.

\begin{lemma2}\label{8-toricsatlem} If $S$ is a normal scheme
and $F$ is a fine monoscheme then there is a natural isomorphism
$S[F]^\nor\iso S[F^\sat]$.  \end{lemma2} \begin{proof} Note
that $f{\colon}\ZZ[F^\sat]\to\Spec(\ZZ)$ is a flat morphism and
its fibers are normal because  they are classical toric varieties
$\FF_p[F^\sat]$. So, $f\times S{\colon}S[F^\sat]\to S$ is a flat
morphism with normal fibers and normal target, and we obtain that
its source is normal by \cite[21.E(iii)]{CA@Matsumura}. It remains
to note that $S[F^\sat]\to S[F]$ is a finite morphism inducing
isomorphism of dense open subschemes $S[F^\gp]$, hence $S[F^\sat]$
is the normalization of $S[F]$.  \end{proof}

\begin{remark2} The same argument shows that if $S$ is Cohen-Macaulay
then so is $S[F^\sat]$.  \end{remark2}

\subsubsection{Toric smoothness}\label{8-toricsmoothsec} Next,
let us compare smoothness of morphisms of monoschemes as defined
in \S\ref{8-smoothfansec} and classical smoothness of toric
morphisms. The following lemma slightly extends the classical result
(e.g. see \cite[\S2.1]{Toric@Fulton}) that if $P$ is fs and $\CC[P]$
is regular then $P\iso\NN^r\oplus\ZZ^s$.

\begin{lemma2}\label{8-smoothfanlem} Let $f{\colon}F\to F'$ be a
morphism of fine monoschemes and let $S$ be a scheme.

(i) If $f$ is smooth then $S[f]$ is smooth (as a morphism of
schemes).

(ii) If $F$ is torsion free and the morphism $S[F]\to S$ is
smooth then $F$ is smooth.  \end{lemma2} \begin{proof} Part
(i) is obvious, so let us check (ii). We can also assume that
$F=\Spec(P)$ is affine. Also, we can replace $S$ with any of its
points achieving that $S=\Spec(k)$. Then $P$ is a fine torsion
free monoid and $k[P]\subseteq k[P^\gp]\iso k[\ZZ^n]$. It
follows that $X=\Spec(k[P])$ is an integral smooth $k$-variety
of dimension $n$. Note that $\Spec(k[P^\sat])$ is a finite
modification of $X$ which is generically an isomorphism. Since $X$
is normal we have that $\Spec(k[P^\sat])\iso X$, and it follows
that $P$ is saturated. Now, $P\iso\oP\oplus\ZZ^l$ and hence
$X\iso\Spec(k[\oP])\times_k\GG_m^l$. Obviously, $\Spec(k[\oP])$
is smooth of dimension $r=n-l$ and our task reduces to showing
that $\oP\iso\NN^{r}$.

Let $m=\oP\setminus\{1\}$ be the maximal ideal of $\oP$. Then
$I=k[m]$ is a maximal ideal of $k[\oP]$ with residue field
$k$. In particular, by $k$-smoothness of $k[\oP]$ we have that
$\dim_k(I/I^2)=r$. On the other hand, $I=\oplus_{x\in m}xk$
and $I^2=\oplus_{x\in m^2}xk$, hence $I/I^2\iso\oplus_{x\in
m\setminus m^2}xk$ and we obtain that $m\setminus m^2$ consists of
$r$ elements $t_1\.t_r$. Note that these elements generate $\oP$
as a monoid and hence they generate $\oP^\gp$ as a group. Since
$\oP^\gp\iso\ZZ^r$, the elements $t_1\.t_r$ are linearly
independent in $\oP^\gp$, and we obtain that the surjection
$\oplus_{i=1}^rt_i^\NN\to\oP$ is an isomorphism.  \end{proof}

\begin{remark2} It seems very probable that, much more
generally, $f$ is smooth whenever $S[f]$ is smooth as a
morphism of schemes and one of the following conditions holds:
(a) $S$ has points in all characteristics, (b) the homomorphisms
${M}^\gp_{F',x'}\to{M}^\gp_{F,x}$ induced by $f$ have torsion free
kernels and cokernels. We could prove this either in the saturated
case or under some milder but unnatural restrictions. The main ideas
are similar but the proof becomes more technical. We do not develop
this direction here since the lemma covers our needs.  \end{remark2}

\subsubsection{Toric ideals}\label{8-toricidealsec} Let $F$ be
a torsion free monoscheme, $k$ be a field and $X=k[F]$. For any
ideal $\calI$ on $F$ one naturally defines an ideal $k[\calI]$
on $X$: in local charts, an ideal $I\subseteq P$ goes to the
ideal $Ik[P]=k[I]=\oplus_{a\in I}ak$ in $k[P]$. We say that
$\calJ=k[\calI]$ is a {\em monoidal ideal} on $k[F]$. Note that
$\calJ$ determines $\calI$ uniquely and, actually, $\calI=\calJ\cap
M_X$ (this is obvious in local charts).

\begin{lemma2}\label{8-toricideallem} Assume that $F$ is connected,
$\eta$ is its maximal point, $X=k[F]$, and $T=k[\eta]$ is the torus
of the toric scheme $X$. A coherent ideal $\calJ\subseteq{M}_X$
is $T$-equivariant if and only if it is monoidal.  \end{lemma2}
\begin{proof} Any monoidal ideal is obviously $T$-equivariant,
so let us prove the inverse implication. The claim is local, so
we should prove that any $T$-equivariant ideal $J\subset A=k[P]$
is of the form $k[I]$ for an ideal $I$ of $P$. Consider the
coaction homomorphism $\mu{\colon}A\to A\otimes_k k[P^\gp]=B$. The
equivariance of $J$ means that $JB$ (with respect to the embedding
$A\into A\otimes_k k[P^\gp]$) is equal to $\mu(J)B$. In particular,
$\mu|_J{\colon}J\to JB=J\otimes_k k[P^\gp]$ induces a $P^\gp$-grading
on $J$ compatible with the $P^\gp$-grading $A=\oplus_{\gamma\in P}
A_\gamma$. Thus, $J$ is a homogeneous ideal in $A$ and, since $A_1=k$
is a field and each $k$-module $A_\gamma$ is of rank one, we obtain
that $J=\oplus_{\gamma\in I} A_\gamma$ for a subset $I\subseteq
P$. Thus, $J=k[I]$, and clearly $I$ is an ideal.  \end{proof}

\subsubsection{Toric blow ups}\label{8-toricblowsec} We will also
need the well known fact that toric blow ups are of combinatorial
origin, i.e. they are induced from blow ups of monoschemes.

\begin{lemma2}\label{8-blowfanlem} Assume that $F$ is a torsion
free monoscheme, $\calI\subseteq{M}_F$ is an ideal, $X=k[F]$,
and $\calJ=k[\calI]$. Then there is a canonical isomorphism
$\Bl_\calJ(X)\iso k[\Bl_\calI(F)]$.  \end{lemma2} \begin{proof}
Assume first that $F=\Spec(P)$. Then $\calI$ corresponds to an ideal
$I\subseteq P$ and we can simply compare charts: $\Bl_\calI(F)$
is covered by the charts $\Spec(P[a^{-1}I])$ for $a\in I$, and,
since $I$ generates $J$, the charts $k[a^{-1}J]=k[P[a^{-1}I]]$
cover $\Bl_J(X)$. This construction is compatible with localizations
$(P,I)\mapsto (P_b,b^{-1}I)$ hence it globalizes to the case of a
general fine monoscheme with an ideal.  \end{proof}

Using Lemma \ref{8-toricsatlem} we obtain a similar relation between
saturated blow ups and normalized toric blow ups.

\begin{corollary2}\label{8-blowfancor} Keep notation of
Lemma \ref{8-blowfanlem}. Then $\Bl_\calJ(X)^\nor\iso
k[\Bl_\calI(F)^\sat]$.  \end{corollary2}

\subsubsection{Desingularization of
monoschemes}\label{8-desingmonosec} Let $F$ be a torsion free
monoscheme and $X=k[F]$ for a field $k$ of characteristic zero
(e.g. $k=\QQ$). Recall that the normalized desingularization functor
$\tilcalF$ from \S\ref{8-desingQsec} is compatible with the action
of any smooth $k$-group, hence the centers of $\tilcalF(X)\colon
X_n\dashto X_0=X$ are $T$-equivariant ideals. By Lemma
\ref{8-toricideallem}, the blown up ideal of $X_0$ is of the form
$k[\calI]$ for an ideal $\calI\subseteq{M}_F$, hence $X_1=k[F_1]$
for $F_1=\Bl_\calI(F)^\sat$ by Lemma \ref{8-blowfancor}. Applying
this argument inductively we obtain that the entire normalized
blow up tower $\tilcalF(X)$ descends to a saturated blow up tower
of $F$, which we denote as $\tilcalF^\mono(F)$ (in other words,
$\tilcalF(X)=k[\tilcalF^\mono(F)]$). Since $F_n$ is smooth by
Lemma \ref{8-smoothfanlem}(ii), the tower $\tilcalF^\mono(F)$
is a desingularization of $F$. Moreover, part (i) of the same
lemma implies that $\tilcalF^\mono$ is functorial with respect to
smooth morphisms of monoschemes. Namely, for any smooth morphism
$F'\to F$, $\tilcalF^\mono(F')$ is the contracted pullback of
$\tilcalF^\mono(F)$ (see \S\ref{8-towpullsec}). Summarizing, we
have obtained:

\begin{theorem2}\label{8-monodesinth} Let $k$ be a field and let
$\tilcalF$ be a normalized desingularization of disjoint unions
of toric $k$-varieties which is functorial with respect to smooth
morphisms. Then each normalized blow up tower $\tilcalF(k[F])$
is the pullback of a uniquely defined saturated blow up tower
of monoschemes $\tilcalF^\mono(F)$. This construction produces a
saturated desingularization of monoschemes which is functorial with
respect to smooth morphisms.  \end{theorem2}

Combining theorems \ref{8-fandesingth} and \ref{8-monodesinth}
we obtain a functorial saturated desingularization of fans that
will be denoted $\tilcalF^\fan$.

\begin{remark2} We will work with normalized and saturated
desingularizations, so we formulated the theorem for $\tilcalF$. The
same argument shows that $\calF$ induces desingularizations
$\calF^\mono$ and $\calF^\fan$ that are functorial with respect
to all smooth morphisms. Moreover, the descent from toric
desingularization is compatible with normalization/saturation,
i.e. $\tilcalF^\mono=(\calF^\mono)^\sat$ and similarly for fans.
\end{remark2}

\subsection{Monoidal desingularization}\label{8-logfansec} In this
section we will establish, what we call, monoidal desingularization
of fine log schemes $(X,M_X)$. This operation "resolves" the sheaf
$\oM_X$ but does not "improve" the log strata of $X$.

\subsubsection{Log stratification}\label{8-logstratsec}
Using charts one immediately checks that for any fine log scheme $(X,M_X)$
the function $x\mapsto {\rm rank}(\oM^\gp_\ox)$ is constructible.
The corresponding stratification of $X$, whose strata are the maximal locally
closed subsets on which this function is constant, will be called
the {\em log stratification} (the analogous stratification in \refext{VI}{6-1.5} was called
canonical or stratification by rank).

\subsubsection{Monoscheme charts of log
schemes}\label{8-monochartsec} A (global) {\em monoscheme
chart} of a Zariski log scheme $(X,M_X)$ is a morphism of
monoidal spaces $c{\colon}(X,\varepsilon_*M_X)\to(F,{M}_F)$
such that the target is a monoscheme and $\varepsilon_*M_X$
is isomorphic to the Zariski log structure associated with the
pre-log structure $c^{-1}M_F\to\calO_X$ (obtained as $c^{-1}M_F\to
M_X\to\calO_X$). In particular, $M_X$ is the log structure associated
with $(c\circ\varepsilon)^{-1}M_F\to\calO_{X_\et}$. We say that
the chart is {\em fine} if $(F,M_F)$ is so. For example, any toric
scheme $R[F]$, where $F$ is a monoscheme, possesses a canonical
chart $R[F]\to F$.

\begin{lemma2}\label{8-monochartlem} Let $(X,M_X)$ be a log scheme
and let $(F,M_F)$ be a monoscheme. Then any morphism of monoidal
spaces $f{\colon}(X,M_X)\to(F,M_F)$ factors uniquely into the
composition of a morphism of log schemes $(X,M_X)\to\ZZ[F]$ and
the canonical chart $\ZZ[F]\to F$.  \end{lemma2} \begin{proof}
Note that $(Id_X,\alpha){\colon}(X,\calO_X)\to(X,M_X)$ is
a morphism of monoidal spaces, hence so is the composition
$h{\colon}(X,\calO_X)\to(F,M_F)$. If $F=\Spec(P)$ then $h$ is
determined by the homomorphism $P\to\Gamma(\calO_X)$ by Lemma
\ref{8-monoidsplem}. Since the latter factors uniquely into the
composition of the homomorphism of monoids $P\to\ZZ[P]$ and the
homomorphism of rings $\ZZ[P]\to\Gamma(\calO_X)$, we obtain a
canonical factoring $X\to\Spec(\ZZ[P])\to\Spec(P)$. Furthermore,
this affine construction is compatible with localizations of $P$,
hence it globalizes to the case when the monoscheme $F$ is arbitrary.
\end{proof}

\begin{remark2}\label{8-usualchartrem} (i) Usually, one works with
log schemes using local charts $(X,M_X)\to\Spec(\ZZ[P])$. By Lemma
\ref{8-monochartlem} this is equivalent to working with affine
monoscheme charts.

(ii) In particular, any fine log scheme $(X,M_X)$ admits a fine
monoscheme chart \'etale-locally, i.e., there exists a strict (in
the log sense) \'etale covering  $(Y,M_Y)\to(X,M_X)$ whose source
possesses a fine monoscheme chart. Similarly, any Zariski fine log scheme
admits a fine monoscheme chart Zariski locally. \end{remark2}

\subsubsection{Chart base change}\label{8-chartchangesec}
Given a fine monoscheme chart $(X,M_X)\to F$ and a morphism of
monoschemes $F'\to F$ we will write $(X,M_X)\times_FF'$ instead
of $(X,M_X)\times_{\ZZ[F]}\ZZ[F']$, where the second product
is taken in the category of fine log schemes. This notation is
partially justified by the following result.

\begin{lemma2}\label{8-chartbaselem} Keep the above notation and
let $(X',M_{X'})=(X,M_X)\times_FF'$.

(i) The morphism $c'{\colon}(X',M_{X'})\to F'$ is a monoscheme chart.

(ii) If $(Y,M_Y)$ is a log scheme over $(X,M_X)$ and
$d{\colon}(Y,M_Y)\to F$ is the induced morphism of monoidal spaces,
then any lifting of $d$ to a morphism $(Y,M_Y)\to F'$ factors
uniquely through $c'$.  \end{lemma2} \begin{proof} Strictness is
stable under base changes, hence $(X',M_{X'})\to\ZZ[F']$ is strict
and we obtain (i). The assertion of (ii) is a consequence of Lemma
\ref{8-monochartlem}.  \end{proof}

\subsubsection{Log ideals}\label{8-logidealsec} By {\em log
ideal} on a fine log scheme $(X,M_X)$ we mean any ideal
$\calI\subseteq M_X$ that \'etale-locally on $X$ admits
a coherent chart as follows: there exists a strict \'etale
covering $f{\colon}(Y,M_Y)\to(X,M_X)$, a fine monoscheme chart
$c{\colon}(Y,M_Y)\to F$ and a coherent ideal $\calI_F\subseteq M_F$
such that $f^{-1}(\calI)M_Y=c^{-1}(\calI_F)M_Y$.

\subsubsection{Log blow ups}\label{8-logblowsec} It is proved
in \cite[4.2]{Toric@Niziol} that there exists a universal
morphism $f{\colon}(X',M_{X'})\to(X,M_X)$ such that the ideal
$f^{-1}(\calI)M_{X'}$ is invertible, i.e. locally (in the \'etale
topology) generated by one element. (We use here that, unlike rings,
any principal ideal $aM$ of an integral monoid $M$ is invertible
in the usual sense, i.e. $M\iso aM$ as $M$-sets.) Actually,
the formulation in \cite{Toric@Niziol} refers only to saturated
blow ups, but the proof deals also with the non-saturated ones.

The construction of log blow us is standard and it also shows
that they are compatible with arbitrary strict morphisms. If
$(X,M_X)$ and $\calI$ admit a chart $F,\calJ\subseteq M_F$ then
$(X',M_{X'})=(X,M_X)\times_F\Bl_\calJ(F')$ is as required. If
$(Y,M_Y)\to(X,M_X)$ is strict then $F$ is also a chart of $(Y,M_Y)$,
hence the local construction is compatible with strict morphisms. The
general case now follows by descent because any fine log scheme
admits a chart  \'etale locally. We call $f$ the {\em log blow up}
of $(X,M_X)$ along $\calI$ and denote it $\LogBl_{\calI}(X,M_X)$
(it is called unsaturated log blow up in \cite{Toric@Niziol}). Log
blow up towers are defined obviously. As usually, contraction
of such a tower is obtained by removing all empty log blow ups
(i.e. blow ups along $\calI=M_X$).

\subsubsection{Saturated log blow ups}\label{8-satlogblowsec}
Saturated log blow up along a log ideal $\calI$ is defined
as $(\LogBl_{\calI}(X,M_X))^\sat$. It satisfies an obvious
universal property too. (It is called log blow up in
\cite{Toric@Niziol}). Towers of saturated log blow ups, their
pullbacks, and saturation of a tower of log blow ups are defined
in the obvious way.

\subsubsection{Pullbacks}\label{8-logpullsec} Let
$f{\colon}(Y,M_Y\to(X,M_X)$ be a morphism of log schemes. By
{\em pullback} of the log blow up $\LogBl_\calI(X,M_X)$
along a log ideal $\calI\subseteq M_X$ we mean the log blow up
$\LogBl_{\calJ}(Y,M_Y)$, where $\calJ=f^{-1}(\calI)M_Y$. This is
the universal log scheme over $(Y,M_Y)$ whose morphism to $(X,M_X)$
factors through $\LogBl_\calI(X,M_X)$. The pullback of saturated blow
ups is defined similarly, and these definitions extend inductively
to pullbacks of towers of (saturated) log blow ups.

\subsubsection{Basic properties}\label{8-basiclogblowsec}
Despite the similarity with usual blow ups of schemes,
log blow ups (resp. saturated log blow ups) have nice
properties that are not satisfied by usual blow ups. First,
it is proved in \cite[4.8]{Toric@Niziol} that log blow ups
are compatible with any log base change $f{\colon}Y\to X$,
i.e. $\LogBl_{f^{-1}\calI}(Y)\iso\LogBl_\calI(X)\times_XY$
for a monoidal ideal $\calI$ on $X$. In particular, saturated
blow ups are compatible with saturated base changes. Second, log
blow ups (resp. saturated log blow ups) are log \'etale morphisms
because so are both saturation morphisms and charts of the form
$\ZZ[\Bl_\calI(F)]\to\ZZ[F]$.

\subsubsection{Fan charts}\label{8-fanchart} A {\em
fan chart} of a Zariski log scheme $(X,M_X)$ is a
morphism $d{\colon}(X,\varepsilon_*\oM_X)\to(F,M_F)$
of monoidal spaces such that the target is a fan and
$d^{-1}(M_F)\iso\varepsilon_*\oM_X$. For example, for any
monoscheme chart $c{\colon}(X,\varepsilon_*M_X)\to(F,M_F)$,
its sharpening $\oc{\colon}(X,\varepsilon_*\oM_X)\to(F,\oM_F)$
is a fan chart. Fan charts were considered by Kato (e.g.,
in \cite[9.9]{Toric@Kato}). They contain less information than
monoscheme charts, but "remember everything about ideals and blow
ups" because there is a one-to-one correspondence between ideals and
blow up towers of $M_F$ and $M_\oF$. Let us make this observation
rigorous. For concreteness, we discuss only non-saturated (log)
blow ups, but everything easily extends to the saturated case.

\begin{remark2}\label{8-fanchartrem} (i) Assume that
$\oc{\colon}(X,\oM_X)\to(F,M_F)$ is a fan chart. Any ideal
$\calI_F\subseteq M_F$ induces a log ideal $\calI\subset
M_X$, which is the preimage of $\oc^{-1}(\calI_F)\oM_X$
under $M_X\to\oM_X$. We say that the blow up $F'=\Bl(\calI_F)$
induces the log blow up $(X',\oM_{X'})=\LogBl_\calI(X,\oM_X)$
or that the latter log blow up is the {\em pullback} of
$\Bl(\calI_F)$. Furthermore, $(X',\oM_{X'})\to F'$ is also a fan
chart (see \cite[9.9]{Toric@Kato}, where the fs case is treated),
hence this definition iterates to a tower $F_\bullet$ of blow ups
of $F$. We will denote the pullback tower of log blow ups as
$\oc^*(F_\bullet)$.

(ii) By a slight abuse of language, Kato and Niziol denote
$\oc^*(F_\bullet)$ as $(X,M_X)\times_FF_\bullet$. One should be very
careful with this notation because, in general, there is no morphism
$(X,M_X)\to F$ that lifts $\oc$. Also, one cannot define analogous
"base change" for an arbitrary morphism of fans $F'\to F$. The reason
is that there are many "unnatural" gluings in the category of fans
(e.g. along generic points), and such gluings cannot be lifted to
log schemes (and even to monoschemes).

(iii) For blow up towers, however, the base change
notation is safe and agrees with the base change from
the monoschemes. Namely, if $\oc$ is the sharpening of a
monoscheme chart $c{\colon}(X,M_X)\to(F,M_F)$ then there exists
a one-to-one correspondence between blow up towers of the
monoscheme $F=(F,M_F)$ and the fan $\oF=(F,\oM_F)$, see Lemma
\ref{8-idealfanlem}. Clearly, the matching towers induce the same
log blow up tower of $(X,M_X)$. In particular, $\calF^\mono(F,M_F)$
and $\calF^\fan(F,\oM_F)$ (see Theorem \S\ref{8-fandesingth})
induce the same log blow up tower of $(X,M_X)$.  \end{remark2}

\subsubsection{Monoidal desingularization of log
schemes}\label{8-monoidaldessec} Let $(X,M)$ be a fine log
scheme and assume that $(X,M)$ is {\em monoidally torsion free}
in the sense that the groups $\oM^\gp_\ox$ are torsion free. By
{\em monoidal desingularization} (resp. {\em saturated monoidal
desingularization}) of a fine log scheme $(X,M)$ we mean a
tower of log blow ups (resp. tower of saturated log blow ups)
$(X_n,M_n)\dashto (X_0,M_0)=(X,M)$ such that for any geometric
point $\ox\to X_n$ the stalk of $\oM_n$ at $\ox$ is a free monoid. A
morphism $(Y,N)\to(X,M)$ is called {\em monoidally smooth} if each
induced homomorphism of stalks of monoids $\oM_\ox\to\oN_\oy$ can
be extended to an isomorphism $\oM_\ox\oplus\NN^r\iso\oN_\oy$.

\begin{theorem2}\label{8-monoidaldesth} Let $\tilcalF^\fan$ be a
saturated desingularization of fine torsion free fans which is
functorial with respect to smooth morphisms. Then there exists
unique saturated monoidal desingularization $\tilcalF^\rmlog(X,M)$
of monoidally torsion free fine log schemes $(X,M)$, such that
$\tilcalF^\rmlog$ is functorial with respect to all monoidally
smooth morphisms and $\tilcalF^\rmlog(X,M)$ is the contraction of
$c^*(\tilcalF^\fan(F))$ for any log scheme $(X,M)$ that admits a
fan chart $c{\colon}(X,\oM_X)\to F$. In the same way, a functorial
desingularization $\calF^\fan$ induces a monoidal desingularization
$\calF^\rmlog$.  \end{theorem2}

\begin{remark2}\label{8-monoidaldesrem} Since any monoscheme
chart induces a fan chart, it then follows from Remark
\ref{8-fanchartrem}(iii) that $\tilcalF^\rmlog(X,M)$ is the
contraction of $d^*(\tilcalF^\mono(F))$ for any log scheme
$(X,M)$ that admits a monoscheme chart $d{\colon}(X,\oM_X)\to F$.
\end{remark2}

\begin{proof}[Proof of Theorem \ref{8-monoidaldesth}] Both
cases are established similarly, so we prefer to deal with
$\calF^\fan$ (to avoid mentioning saturations at any step of
the proof). By descent, it suffices to show that the pullback
from fans induces a functorial monoidal desingularization of
those fine log schemes that admit a global fan chart. Thus,
if $\calF^\rmlog$ exists then it is unique, and our aim is to
establish existence and functoriality. Both are consequences of
the following claim: assume that $f{\colon}(Y,M_Y)\to(X,M_X)$
is a monoidally smooth morphism whose source and target admit fan
charts $d{\colon}(Y,\oM_Y)\to G$ and $d'{\colon}(X,\oM_X)\to F$,
then contractions of $d^*(\calF^\fan(G))$ and $c^*(\calF^\fan(F))$
are equal, where $c=d'\circ\of{\colon}(Y,\oM_Y)\to(X,\oM_X)\to
F$. Note that $d^{-1}(\oM_G)\iso\oM_Y$ and the homomorphism
$c^{-1}(\oM_F)\iso\oM_Y$ is smooth.

Choose a point $y\in Y$ and consider the localizations
$Y'=\Spec(\calO_{Y,y})$, $F'=\Spec(M_{F,c(y)})$ and
$G'=\Spec(M_{G,d(y)})$ at $y$ and its images in the fans. Since
$\phi{\colon}M_{G,d(y)}\to\oM_{Y,y}$ is an isomorphism and the
homomorphism $\psi{\colon}M_{F,c(y)}\to\oM_{Y,y}$ is smooth,
we obtain a factorization of $\psi$ into a composition of
a homomorphism $\lambda{\colon}M_{F,c(y)}\to M_{G,d(y)}$
and $\phi$, where $\lambda$ is smooth. Set $U=c^{-1}(F')\cap
d^{-1}(G')$. Then $U$ is a neighborhood of $y$, and $c$ and $d$
induce homomorphisms $\phi_U{\colon}M_{G,d(y)}\to\oM_Y(U)$ and
$\psi_U{\colon}M_{F,c(y)}\to\oM_Y(U)$. Sine the monoids are fine,
we can shrink $U$ so that the equality $\psi_U=\phi_U\circ\lambda$
holds. It then follows from Lemma \ref{8-monoidsplem}(ii) that
$c|_U$ factors into a composition of $d|_U$ and the smooth morphism
$\Spec(\lambda){\colon}G'\to F'$.

By quasi-compactness of $Y$ we can now find finite coverings
$Y=\cup_{i=1}^nY_i$, $F=\cup_{i=1}^nF_i$ and $G=\cup_{i=1}^nG_i$,
and smooth morphisms $\lambda_i{\colon}G_i\to F_i$ such that
$Y_i$ is mapped to $F_i$ and $G_i$ by $c$ and $d$, respectively,
and the induced maps of monoidal spaces $c_i{\colon}Y_i\to
F_i$ and $d_i{\colon}Y_i\to G_i$ satisfy $c_i=\lambda_i\circ
d_i$. Set $Y'=\coprod_{i=1}^nY_i$, $F'=\coprod_{i=1}^n F_i$,
$G'=\coprod_{i=1}^nG_i$, $c'{\colon}Y'\to F$ and $d'{\colon}Y'\to
G$. By descent, it suffices to check that contractions of
$c'^*(\calF^\fan(F))$ and $d'^*(\calF^\fan(G))$ are equal. Since,
the morphism $Y'\to F$ factors through the surjective smooth morphism
$F'\to F$, and similarly for $G$, these two pullbacks are equal to
the contracted pullbacks of $\calF^\fan(F')$ and $\calF^\fan(G')$,
respectively. It remains to note that $Y'\to F'$ factors through the
smooth morphism $\coprod_{i=1}^n\lambda_i{\colon}G'\to F'$. Hence
$\calF^\fan(G')$ is the contracted pullback of $\calF^\fan(F')$,
and their contracted pullbacks to $Y'$ coincide.  \end{proof}

\subsection{Desingularization of log regular log
schemes}\label{8-resolution of log schemes1} In this section we will
see how saturated monoidal desingularization leads to normalized
desingularization of log regular log schemes. Up to now we freely
considered saturated and unsaturated cases simultaneously, and did
not feel any essential difference. This will not be the case in the
present section because the notion of log regularity was developed
by Kato and Niziol in the saturated case. Gabber generalized the
definition for non-saturated case and extended to that case all
main results about log regular log schemes. This was necessary for
his original approach, but can be by passed by use of saturated
monoidal desingularization. So, we prefer to stick to the saturated
case and simply refer to all foundational results about log regular
fs log schemes to \cite{Toric@Kato} and \cite{Toric@Niziol}. For
the sake of completeness, we will outline Gabber's results about
the general case in \S\ref{8-nonsatursec}.

\subsubsection{Conventions}\label{8-logregconvsec} Recall
that Kato's notion of log regular fs log schemes was already
used in expos\'e VI, \S1.2. Throughout \S\ref{8-resolution of
log schemes1} we assume that $(X,M_X)$ is a log regular fs log
scheme. Note that the homomorphism $\alpha_X{\colon}M_X\to\calO_X$
of $X_\et$-sheaves is injective by \cite[2.6]{Toric@Niziol}, and
actually $M_X=\calO^\times_U\cap\calO_X$, where $U\subseteq X$ is
the triviality locus of $M_X$. So, we will freely identify $M_X$
with a multiplicative submonoid of $\calO_X$.

\subsubsection{Monoidal ideals}\label{8-logregidealsec}
For any log ideal $\calI\subseteq M_X$ consider the ideal
$\calJ=\alpha(\calI)\calO_X$ it generates. We call $\calJ$ a {\em
monoidal ideal} and by a slight abuse of language, we will write
$\calJ=\calI\calO_X$.

\begin{lemma2}\label{8-logregideallem} Let $X=(X,M_X)$ be as in
\S\ref{8-logregconvsec}. The rules $\calI\mapsto\calI\calO_X$
and $\calK\mapsto\calK\cap M_X$ give rise to a one-to-one
correspondence between log ideals $\calI\subseteq M_X$
and monoidal ideals $\calK\subseteq\calO_X$.  \end{lemma2}
\begin{proof} It suffices to show that any log ideal $\calI$
coincides with $\calJ=\calI\calO_X\cap M_X$. The question easily
reduces to the case when $X$ admits a chart $X\to\Spec(\ZZ[P])$
and $\calI=IM_X$, $\calJ=JM_X$ for ideals $I\subseteq J$
of $P$. Assume that the lemma fails, i.e. $I\subsetneq
J$, and consider the exact sequence \begin{equation}
\label{8-logregeq} I\ZZ[P]\otimes_{\ZZ[P]}\calO_X\to
J\ZZ[P]\otimes_{\ZZ[P]}\calO_X\to
J\ZZ[P]/I\ZZ[P]\otimes_{\ZZ[P]}\calO_X\to 0 \end{equation}
Since $\Tor_1^{\ZZ[P]}(\ZZ[P]/I\ZZ[P],\calO_X)=0$
by \cite[6.1(ii)]{Toric@Kato},
$I\ZZ[P]\otimes_{\ZZ[P]}\calO_X=I\calO_X$ and similarly
for $J$, and we obtain that the first morphism in the
sequence (\ref{8-logregeq}) is $I\calO_X\to J\calO_X$. To
obtain a contradiction, it suffices to show that
$J\ZZ[P]/I\ZZ[P]\otimes_{\ZZ[P]}\calO_X\neq 0$. Note that
$\ZZ[P]/m_P\ZZ[P]$ is a quotient of $J\ZZ[P]/\ZZ[I]$,
so it remains to note that $m_P\calO_X\neq \calO_X$ and hence
$\ZZ[P]/m_p\ZZ[P]\otimes_{\ZZ[P]}\calO_X\neq 0$.  \end{proof}

\subsubsection{Interpretation of monoidal smoothness}\label{8-regsec}
Note that by \refext{VI}{6-1.7} $(X,M_X)$ is monoidally smooth if and only if
$X$ is regular and in this case the non-triviality locus of $M_X$
is a normal crossings divisor $D$.

\subsubsection{Saturated log blow ups of log regular
log schemes}\label{8-logregsatblowsec} Using Kato's
$\Tor$-independence result \cite[6.1(ii)]{Toric@Kato} Niziol
proved in \cite[4.3]{Toric@Niziol} that saturated log blow
ups of $(X,M_X)$ are compatible with normalized blow ups along
monoidal ideals. Namely, if $(Y,M_Y)=\LogBl_\calI(X,M_X)^\sat$
then $Y\iso\Bl_{\calI\calO_X}(X)^\nor$. We will also need more
specific results that showed up in the proof of loc.cit., so we
collect them altogether in the following lemma.

\begin{lemma2}\label{8-logregblowlem} Let
$f{\colon}(X,M_X)\to(Y,M_Y)$ be a strict morphism of fs log
regular log schemes, let $\calI\subseteq M_Y$ be a log ideal and
$\calJ=f^{-1}\calI M_X$. Set $(X',M_{X'})=\LogBl_\calJ(X,M_X)$,
$(X'',M_{X''})=(X',M_{X'})^\sat$, $(Y',M_{Y'})=\LogBl_\calI(Y,M_Y)$
and $(Y'',M_{Y''})=(Y',M_{Y'})^\sat$. Then

(i) The (saturated) log blow up of $(X,M_X)$ is compatible with
(normalized) blow up of $X$: $X'\iso\Bl_{\calI\calO_X}(X)$
and $X''\iso\Bl_{\calI\calO_X}(X)^\nor$.

(ii) The (normalized) blow up of $X$ along $\calJ$ is the pullback
of the (normalized) blow up of $Y$ along $\calI$. In particular,
$X'\iso X\times_YY'$ and $X''\iso X\times_YY''$.  \end{lemma2}
\begin{proof} All claims can be checked \'etale locally, hence we can
assume that there exists a chart $g{\colon}(Y,M_Y)\to(\ZZ[P],P)$
and $\calI=g^{-1}(I_0)M_X$ for an ideal $I_0\subseteq P$. Then it
suffices to prove (ii) for $g$ and the induced chart $g\circ f$ of
$(X,M_X)$. In particular, this reduces the lemma to the particular
case when $X=\Spec(A)$ and $f$ is a chart $(X,M_X)\to(\ZZ[P],P)$. It
is shown in the first part of the proof of \cite[4.3]{Toric@Niziol}
that $$X'\iso\Proj(A\otimes_{\ZZ[P]}(\oplus_{n=0}^\infty
I_0^n))\iso\Proj(\oplus_{n=0}^\infty\calJ^n)$$ The first
isomorphism implies that $X'\to X$ is the base change of
$\Proj(\oplus_{n=0}^\infty I_0^n)=Y'\to Y$, and the second
isomorphism means that $X'\iso\Bl_{\calI\calO_X}(X)$. This
establishes the unsaturated and unnormalized part of the Lemma,
and the second part follows in the same way from the second part
of the proof of \cite[4.3]{Toric@Niziol}.  \end{proof}

\begin{remark2} It follows from the lemma that the unsaturated log
blow up $(X',M_{X'})$ is log regular in the sense of Gabber, see
\S\ref{8-nonsatursec}. Thus, once log regularity is correctly defined
in full generality, it becomes a property preserved by log blow ups
(as it should be, since log blow ups are log smooth).  \end{remark2}

\subsubsection{Desingularization of log regular log
schemes}\label{8-desingloregsec} By Lemmas \ref{8-logregideallem}
and \ref{8-logregblowlem}, any saturated log blow up tower
$(X_n,M_n)\dashto (X,M_X)$ induces a normalized blow up tower
$X_n\dashto X$, which completely determines it. Furthermore,
by \S\ref{8-regsec}, the first tower is a saturated monoidal
desingularization if and only if the second tower is a normalized
desingularization. In particular, the saturated monoidal
desingularization $\tilcalF^\rmlog(X,M_X)$ induces a normalized
desingularization of the scheme $X$ that depends on $(X,M_X)$
and will be denoted $\tilcalF(X,M_X)$.

\begin{theorem2}\label{8-deslogregth} The saturated
monoidal desingularization $\tilcalF^\rmlog$ (see Theorem
\ref{8-monoidaldesth}) gives rise to desingularization
$\tilcalF$ of log regular log schemes that possesses the same
functoriality properties: if $\phi{\colon}(Y,M_Y)\to(X,M_X)$
is a monoidally smooth morphism of log regular log
schemes then $\tilcalF(Y,M_Y)$ is the contraction of
$\phi^\st(\tilcalF(X,M_X))$. Furthermore, if $\phi$ is strict then
$\phi^\st(\tilcalF(X,M_X))=\tilcalF(X,M_X)\times_XY$.  \end{theorem2}
Strictness of $\phi$ in the last claim is not needed. To remove it
one should work out the assertion of Remark \ref{8-towpullrem}.
\begin{proof} We should only establish functoriality. Let
$(X',M_{X'})=\LogBl_\calI(X,M_X)^\sat$ be the first saturated blow
up of $\tilcalF^\rmlog(X,M_X)$. Set $\calJ=\phi^{-1}\calI M_Y$, then
$(Y',M_{Y'})=\LogBl_\calJ(Y,M_Y)^\sat$ is the first saturated blow up
of $\tilcalF^\rmlog(Y,M_Y)$ by functoriality of $\tilcalF^\rmlog$. By
Lemma \ref{8-logregblowlem}, $X'=\Bl_{\calI\calO_X}(X)^\nor$
and $Y'=\Bl_{\calJ\calO_Y}(Y)^\nor$, and using that
$\phi^{-1}(\calI\calO_X)\calO_Y=\calJ\calO_Y$ we obtain that $Y'$
is the strict transform of $X'$. It remains to inductively apply
the same argument to the other levels of the towers. The last claim
follows from Lemma \ref{8-logregblowlem}(ii).  \end{proof}

\begin{remark2} The same results hold for (non-saturated) monoidal
desingularization, which induces a (usual) desingularization of
log regular log schemes. For non-saturated log regular log schemes
(see \S\ref{8-nonsatursec}) one should first establish analogs of
Lemmas \ref{8-logregideallem} and \ref{8-logregblowlem} (where the
input in the second one does not have to be saturated). Then the
same proof as above applies.  \end{remark2}

\subsubsection{Canonical fans and associated
points}\label{8-canonfansec} By the {\em canonical fan} $\Fan(X)$ of $(X,M_X)$
we mean the set of maximal points of the log stratification (i.e. the maximal
points of the log strata). Alternatively, $\Fan(X)$ can be described as in
\cite[\S2.2]{Toric@Niziol} as the set of points $x\in X$ such that
$m_\ox$ coincides with the ideal $I_\ox\subseteq\calO_\ox$
generated by $\alpha(M_\ox-M_\ox^\times)$.

We provide $F=\Fan(X)$ with the induced topology and define $M_F$ to be the
restriction of $\oM_X$ onto $F$. For example, for a toric $k$-variety
$X=\Spec(k[Z])$, where $Z$ is a monoscheme, $(F,M_F)$ is isomorphic
to the sharpening of $Z$. More generally, if log scheme $(X,M_X)$
is Zariski then $(F,M_F)$ is a fan and the map $c{\colon}X\to F$ sending any
point to the maximal point of its log stratum is a fan chart of $X$.
This follows easily from the fact that such $(X,M_X)$ admits monoscheme
charts Zariski locally.

\begin{remark2}
In general, $(F,M_F)$ does not have to be a fan, but it seems probable that
it can be extended to a meaningful object playing the role of algebraic
spaces in the category of fans. We will not investigate this direction here.
\end{remark2}

\begin{lemma2}\label{8-assoclem} Let $X=(X,M_X))$ be an fs log
regular log scheme with a monoidal ideal $\calI\subseteq\calO_X$. Then:

(i) The set of associated points of $\calO_X/\calI$ is contained
in $\Fan(X)$.

(ii) The fans of $\Bl_\calI(X)$ and $\Bl_\calI(X)^\nor$ are contained
in the preimage of $\Fan(X)$.

(iii) For any tower of monoidal blow ups (resp. normalized
monoidal blow ups) $X_n\dashto X$, its set of associated points
is contained in $\Fan(X)$.  \end{lemma2} \begin{proof} (i) It
suffices to check that any $x\in X-\Fan(X)$ is not an associated
point of $\calO_X/\calI$. Since associated points are compatible
with flat morphisms, we can pass to the formal completion
$\hatX_x=\Spec(\hatcalO_{X,x})$. Let us consider first the more
difficult case when $A=\hatcalO_{X,x}$ is of mixed characteristic
$(0,p)$. By \refext{VI}{6-str loc log reg}, $A\iso B/(f)$ where $B=C(k)[[Q]][[\ut]]$,
$Q$ is a sharp monoid defining the log structure, $\ut=(t_1\. t_n)$,
and $f\in B$ equals to $p$ modulo $(Q\setminus\{1\},\ut)$. Note that
$n\ge 1$ as otherwise $Q\setminus\{1\}$ would generate the maximal
ideal of $A$. In this case, $x$ is the log stratum of $\hatX_x$
and hence $x$ is the maximal point of its log stratum in $X_x$,
which contradicts that $x\notin\Fan(X)$. The completion $\hatcalI$
of $\calI$ is of the form $IA$ for an ideal $I\subset Q$. We should
prove that $x\in\hatX_x$ is not an associated point of $A/\hatcalI$,
or, equivalently, $\depth(A/\hatcalI)\ge 1$. Since $B/(fB+\hatcalI
B)=A/\hatcalI$, it suffices to show that $\depth(B/\hatcalI B)\ge 2$
and $f$ is a regular element of $B/\hatcalI B$.

Regularity of $f$ follows from the following easy claim by
taking $C=C(k)$, $J=I$ and $R=Qt_1^\NN\dots t_n^\NN$: if $C$
is a domain, $R$ is a sharp monoid with an ideal $J$, $f$ is an
element of $C[[R]]$ with a non-zero free term, and $g\in C[[R]]$
satisfies $fg\in JC[[R]]$ then $g\in JC[[R]]$. Next, let us bound
the depth of $B/\hatcalI B$. In view of \cite[21.C]{CA@Matsumura},
depth is preserved by completions of local rings hence it suffices
to show that $\depth(D_r/ID_r)\ge 2$, where $D=C(k)[Q][\ut]$ and
$r=(m_Q,p,\ut)$ is the ideal generated by $m_Q=Q\setminus\{1\}$,
$p$ and $\ut$. Note that $D/ID\iso C(k)[\ut][Q\setminus I]$
as a $C(k)[\ut]$-module, hence it is a flat $C(k)[\ut]$-module
and the local homomorphism $C(k)[\ut]_{(p,\ut)}\to
D_r/ID_r$ is flat. By \cite[21.C]{CA@Matsumura},
$\depth(D_r/ID_r)\ge\depth(C(k)[\ut]_{(p,\ut)})=n+1\ge 2$, so
we have established the mixed characteristic case. In the equal
characteristic case we have that $A=k[[Q]][[\ut]]$, and the same
argument shows that $\depth(A/\hatcalI)\ge\depth(k[\ut]_{(\ut)})=n\ge
1$.

To prove (ii) we should check that if $x\in X-\Fan(X)$ then no
point of $\Fan(\Bl_\calI(X))$ sits over $x$. We will only check the
mixed characteristic case since it is more difficult. As earlier,
$\hatX_x=\Spec(A)$ where $A=B/(f)$ and $B=C(k)[[Q]][[\ut]]$. Note
that $\psi{\colon}\hatX_x\to X$ is a flat strict morphism of
log schemes. Hence $\psi$ is compatible with blow ups and it
maps the fans of $\hatX_x$ and $\Bl_\hatcalI(\hatX_x)$ to the
fans of $X_x$ and $\Bl_\calI(X)$, respectively. Therefore,
we should only check that $\Fan(\Bl_\hatcalI(\hatX_x))$ is
disjoint from the preimage of $x$. The latter blow up is covered
by the charts $V_a=\Spec(A[a^{-1}\hatcalI])$ with $a\in I$. Set
$P'=Q[a^{-1}I]$, $B'=C(k)[[P']][[t_1\. t_n]]/(f)$ and $A'=B'/(f)$
(where $f$ is as above). Then the $m_x$-adic completion of $V_a$
is $\hatV_a=\Spec(A')$. The ideal defining the closed point of
$\Fan(\hatV_a)$ is generated by the maximal ideal $m'$ of $P'$. This
ideal does not contain any $t_i$. Indeed, $t_i\notin m'B'+fB'$
because $t_i\notin fC(k)[[t_1\.t_n]]$ as $f-p\in(t_1\.t_n)$. Thus,
$\Fan(\hatV_a)$ is disjoint from the preimage of $x$, and hence the
same is true for $\Fan(V_a)$. This proves (ii) in the non-saturated
case, and the saturated case is dealt with similarly but with $P'$
replaced by its saturation. Finally, (iii) follows immediately from
(i) and (ii).  \end{proof}

\subsubsection{Independence of the log structure}\label{8-indepsec}
Dependence of $\tilcalF(X,M_X)$ on $M_X$ is a subtle question. In
this section we will use functoriality of $\tilcalF$ to prove that
$\tilcalF(X,M_X)$ is independent of $M_X$ in characteristic zero
(i.e., it is compatible with any automorphism of $X$). This will
cover our needs, but the restriction on the characteristic will
complicate our arguments later. Conjecturally, $\tilcalF(X,M_X)$
does not depend on $M_X$ at all and the following result of
Gabber supports this conjecture: if $P$ and $Q$ are two fine
sharp monoids and $k[[P]]\iso k[[Q]]$ for a field $k$ (of any
characteristic!) then $P\iso Q$. For completeness, we will give
a proof of this in \S\ref{8-reconstrsec}.

\begin{theorem2}\label{8-Morel 2.4.2} Let $\tilcalF$ be
a functorial normalized desingularization of regular qe
schemes of characteristic zero, and let $\tilcalF(X,M_X)$ be
the normalized desingularization of log regular log schemes
it induces (see Theorem \ref{8-deslogregth}). Assume that
$(X,M_X)$ and $(Y,M_Y)$ are saturated log regular log schemes
such that there exists an isomorphism $\phi{\colon}Y\iso X$
of the underlying schemes. Assume also that the maximal points
of the strata of the stratifications of $X$ and $Y$ by the
rank of $\overline{M}^{gp}$ are of characteristic zero. Then
$\tilcalF(X,M_X)$ and $\tilcalF(Y,M_Y)$ are compatible with $\phi$,
that is, $\tilcalF(Y,M_Y)=\phi^*(\tilcalF(X,M_X))$.  \end{theorem2}
\begin{proof} We can check the assertion of the theorem \'etale
locally. Namely, we can replace $X$ with a strict \'etale covering
$X'$ and replace $Y$ with $Y'=Y\times_XX'$ with the log structure
induced from $Y$. In particular, we can assume that the log
structures are Zariski, and so the canonical fans $\Fan(X)$ and
$\Fan(Y)$ are defined. Our assumption  on the maximal points actually
means that $\Fan(X)$ is contained in $X_\QQ=X\otimes_\ZZ\QQ$. By
Lemma \ref{8-assoclem}(iii) and \S\ref{8-assnorsec}, $\tilcalF(X,M)$
is the pushforward of its restriction onto $X_\QQ$, and
similarly for $Y$. So, it suffices to prove that the normalized
desingularizations of $X_\QQ$ and $Y_\QQ$ are compatible with
respect to $\phi\otimes_\ZZ\QQ{\colon}X_\QQ\iso Y_\QQ$. Thus,
it suffices to prove the theorem for $X$ and $Y$ of characteristic
zero, and, in the sequel, we assume that this is the case.

To simplify notation we identify $X$ and $Y$, and set $N_X=M_Y$. It
suffices to check that the blow up towers $\tilcalF(X,M_X)$ and
$\tilcalF(X,N_X)$ coincide after the base change to each completion
$\hatX_\ox=\Spec(\hatcalO_{X,\ox})$ at a geometric point $\ox$. By
\refext{VI}{6-str loc log reg}, we have that $\hatX_\ox\iso\Spec(k[[P]][[t_1\. t_n]])$,
where $P\iso\oM_{X,\ox}$, and the morphism of fs log
schemes $(\hatX_\ox,P)\to(X,M_X)$ is strict. By Theorem
\ref{8-monoidaldesth} the contracted pullback of $\tilcalF(X,M_X)$
to $\hatX_\ox$ coincides with $\tilcalF^\rmlog(\hatX_\ox,P)$. In
the same way, the contracted pullback of $\tilcalF(X,N_X)$
coincides with $\tilcalF^\rmlog(\hatX_\ox,Q)$, where
$\hatcalO_{X,\ox}\iso\Spec(k[[Q]][[t_1\. t_m]])$ (we use that $k$
is isomorphic to the residue field of this ring and hence depends
only on the ring $\hatcalO_{X,\ox}$).

Let us now recall how $\tilcalF^\rmlog(\hatX_\ox,P)$ is
constructed (Theorems \ref{8-fandesingth}, \ref{8-monodesinth}
and \ref{8-monoidaldesth}). We have the obvious strict
morphism $\hatX_\ox\to Z:=\Spec(\QQ[P][t_1\. t_n])$,
hence $\tilcalF^\rmlog(\hatX_\ox,P)$ is the pullback
of $\tilcalF^\rmlog(Z,P)$. The latter is the pullback of
$\tilcalF^\mono(P)=\tilcalF^\fan(P)$, which, in its turn, is
induced from $\tilcalF(Z)$. Therefore, $\tilcalF(Z,P)=\tilcalF(Z)$
and, by the functoriality of $\tilcalF$, its pullback to
$\hatX_\ox$ is $\tilcalF(\hatX_\ox)$. The same argument shows
that $\tilcalF(\hatX_\ox)$ is the contracted pullback of
$\tilcalF(X,N_X)$.

In order to prove that $\tilcalF(X,M_X)=\tilcalF(X,N_X)$ it only
remains to resolve the synchronization issues, i.e. to prove
equality without contractions. For this one should take the
union $S$ of the fans of $(X,M_X)$ and $(X,N_X)$, and consider
the morphism $\hatX_S:=\coprod_{x\in S}\hatX_\ox$ rather than
the individual completions. The pullbacks of $\tilcalF(X,M_X)$
and $\tilcalF(X,N_X)$ to $\hatX_S$ have no empty blow ups because
the fans contain all associated points of the towers by Lemma
\ref{8-assoclem}(iii). Hence the same argument as above shows that
they both coincide with $\tilcalF(\hatX_S)$.  \end{proof}

\begin{remark2} (i) Without taking the completion, $\tilcalF(X)$
does not even have to be defined as $X$ may be non-qe. In order
to pass to the completion we used functoriality of the monoidal
desingularization with respect to strict morphisms, which may be very
bad (e.g. non-flat) on the level of usual morphisms of schemes. Even
when $(X,M)$ is log regular, the formal completion morphism
$\hatX_x\to X$ can be non-regular in the non-qe case. However,
one can show that it is regular over the fan, and this is enough
to relate the (log) desingularization of $X$ and $\hatX_x$.

(ii) We used the very strong desingularization $\tilcalF$ from
Theorem \ref{8-qeth}. However, it is easy to see that only the
properties listed in Theorem \ref{8-magicboxth} were essential.
\end{remark2}

\subsection{Complements on non-saturated log regular log
schemes}\label{8-nonsatursec} For the sake of completeness, we
mention Gabber's results on non-saturated log regular log schemes
that will not be used. We only formulate results but do not give
proofs. Gabber defines a fine log scheme $(X,M_X)$ to be {\em log
regular} if its saturation is log regular in the usual sense. Assume
now that $(X,M_X)$ is fine and log regular. The key result that
lifts the theory off the ground is that Kato's Tor independence
extends to non-saturated log regular log schemes. Namely, if
$(X,M_X)$ admits a chart $\ZZ[P]$ and $I\subseteq P$ is an
ideal then ${\rm Tor}_1^{\ZZ[P]}(\calO_X,\ZZ[P]/I)=0$. For fs
log schemes this is due to Kato, and Gabber deduces the general
case using a non-flat descent. It then follows similarly to
the saturated case that if $(Y,M_Y)=\LogBl_\calI(X,M_X)$ then
$Y\iso\Bl_{\calI\calO_X}(X)$ and $(Y,M_Y)$ is log regular. In
addition, one shows that $(X,M_X)$ is saturated if and only if $X$ is
normal, and if $(Y,M_Y)=(X,M_X)^\sat$ then $Y\iso X^\nor$. Using
these foundational results on log regular fine log schemes one can
imitate the method of \S\ref{8-resolution of log schemes1} to extend
Lemma \ref{8-assoclem} to the non-saturated case. As a corollary,
one obtains an analog of Theorem \ref{8-Morel 2.4.2} for $\calF$
and $\calF^\rmlog$.

\subsection{Reconstruction of the monoid}\label{8-reconstrsec}
This section will not be used in the sequel. Its aim is to prove
that a fine torsion free monoid $P$ can be reconstructed from a ring
$A=k[[P]]$ where $k$ is a field. The main idea of the proof is that
an isomorphism $k[[P]]\iso A$ defines an action of the torus
$\Spec(k[P^\gp])$ on $A$, and any two maximal tori in $\Aut_k(A)$
are conjugate.

\subsubsection{Automorphism groups of complete
rings}\label{8-automsec} Let $k$ be a field and $A$ be a complete
local noetherian $k$-algebra with residue field $k$. Let $m$ denote
the maximal ideal and set $A_n=A/m^{n+1}$. Consider the algebraic
$k$-groups $G_n=\Aut_k(A_n)$, i.e. $G_n(B)=\Aut_B(A_n\otimes_kB)$
for any $k$-algebra $B$. They form a filtering projective family
$\dots G_2\to G_1\to G_0$, which we call algebraic pro-group
$G_\bullet$. Note that $G_\bullet$ induces a functor $G(B)=\lim_n G_n(B)$
and $G(k)=\Aut_k(A)$.

\begin{remark2} Gabber also considered more complicated filtering
families, but we stick to the simplest case we need.  \end{remark2}

\subsubsection{Stabilization}\label{8-stabsec} We say that an
algebraic pro-group $G_\bullet$ is {\em stable} if the homomorphisms
$G_{n+1}\to G_n$ are surjective for large enough $n$. Any algebraic
group can be stabilized as follows. For each $G_n$ and $i\ge
0$ let $G_{n,i}$ denote the image of $G_{n+i}$ in $G_n$. Then
$G_{n,0}\supseteq G_{n,1}\supseteq\dots$ is a decreasing sequence
of algebraic subgroups of $G$, hence it stabilizes on a subgroup
$G^\st_n\subseteq G_n$. The family $G^\st_\bullet$ with obvious
transition morphisms is an algebraic pro-subgroup of $G_\bullet$, and
it is clear from the definition that $G^\st_\bullet$ is stable. Note
also that $G_\bullet^\st$ is isomorphic to $G_\bullet$ at least in the
sense that $G(B)\iso\lim_n G^\st_n(B)$ for any $k$-algebra $B$.

\subsubsection{Maximal pro-tori}\label{8-maxtorsec} By a {\em
pro-torus} $T_\bullet$ in $G_\bullet$ we mean a {\em compatible}
family of tori $T_n\into G_n$  for $n\gg 0$, in the sense that
$\pi_n(T_{n+1})=T_n$. It is  called a {\em torus} if $\pi_n$'s
are isomorphisms for $n\gg 0$. A pro-torus is {\em split} or {\em
maximal split} if so are $T_n$ for $n\gg 0$. Pro-tori $T_\bullet$
and $T'_\bullet$ are conjugate if for $n\gg 0$ there exists a {\em
compatible} family of conjugations $G_n\to G_n$ taking $T_n$
to $T'_n$.

\begin{proposition2}\label{8-maxtorprop} Assume that $k$
is a field, $G_\bullet$ is a stable pro-algebraic $k$-group,
and $T_\bullet,T'_\bullet\into G_\bullet$ are split pro-tori. If
$T_\bullet$ is maximal then $T'_\bullet$ is conjugate to a
sub-pro-torus of $T_\bullet$. In particular, any two maximal
split pro-tori are conjugate.  \end{proposition2} \begin{proof}
It is a classical result that maximal split tori in algebraic
$k$-groups are conjugate. In particular, for each $n$ we can
move $T'_n$ into $T_n$ by a conjugation, and the only issue
is compatibility of the conjugations. Naturally, to achieve
compatibility we should lift conjugations inductively from $G_n$
to $G_{n+1}$. It suffices to prove that if $\pi_n$ is surjective and
$c_n{\colon}G_n\to G_n$ conjugates $T'_n$ into $T_n$ then it lifts
to $c_{n+1}{\colon}G_{n+1}\to G_{n+1}$ that conjugates $T'_{n+1}$
into $T_{n+1}$. By the stability assumption, we can lift $c_n$
to a conjugation $c'$ of $G_{n+1}$. It takes $T'_{n+1}$ to the
subgroup $H=KT_{n+1}$, where $K$ is the kernel of $\pi_n$. Since
maximal split tori in $H$ are conjugate and conjugation by elements
of $T_{n+1}$ preserves $T_{n+1}$, we can find a conjugation $c''$
by an element of $K$ that takes $c'(T'_{n+1})$ to $T_{n+1}$. Then
$c_{n+1}=c''c'$ is a lifting of $c_n$ as required \end{proof}

\begin{corollary2}\label{8-maxtorcor} Let $k$ be a
field and let $G_\bullet$ be an pro-algebraic $k$-group. If
$T_\bullet,T'_\bullet\subseteq G_\bullet$ are stable pro-tori and $T_\bullet$
is maximal as a split pro-torus of $G_\bullet$, then $T'_\bullet$ is
conjugate to a subtorus of $T_\bullet$.  \end{corollary2} \begin{proof}
Obviously, $T_n,T'_n\subseteq G^\st_n$ for $n\gg 0$. So, $T'_\bullet$
is conjugate to a subtorus of $T_\bullet$ already inside of $G^\st$
by Proposition \ref{8-maxtorprop}.  \end{proof}

\subsubsection{Certain tori in $\Aut_k(A)$}\label{8-torisec}
Any isomorphism $C[[P]]\iso A$, where $C$ is a complete
local $k$-algebra and $P$ is a sharp fine monoid, induces an
algebraic action of the split torus $T_P=\Spec(k[P^\gp])$ on $A$:
a character $\chi{\colon}P\to k$ acts on $C$ trivially and acts
on $p\in P$ by $p\mapsto \chi(p)p$ (and the action of $B$-points
$\chi{\colon}P\to B$ is analogous). Thus we obtain homomorphisms
$\psi{\colon}T_P\to G'_n$ which are injective for $n>0$. In
particular, the image is a split torus of $G$. Furthermore, we claim
that if $C=k$ then the torus is maximal (as a pro-torus). Indeed, if
$\phi\in\Aut_k(A)$ commutes with $T_P$ then its action on $k[[P]]$
preserves the $T_P$-eigenspaces $pk$ for $p\in P$ and on each
$pk$ it acts by multiplication by a number $\lambda(p)$. Clearly,
$\lambda{\colon}P\to k$ is a homomorphism and we obtain that $\phi$
belongs to $\psi(T_P(k))$ and corresponds to $\lambda\in T_P(k)$.

\begin{theorem2}\label{8-toritheor} Assume that $k$ is a field,
$P$ is a sharp fine monoid and $A=k[[P]]$. If $C$ is a complete
local $k$-algebra, $Q$ is a sharp fine monoid and $C[[Q]]\iso A$
then $C\iso k[[R]]$ and $P\iso Q\times R$ for a sharp fine
monoid $R$. In particular, $P$ is uniquely determined by $A$.
\end{theorem2} \begin{proof} Consider $G=\Aut_k(A)$ with split tori
$T_P,T_Q\into G$ corresponding to these isomorphisms. By maximality
of $T_P$ and Corollary \ref{8-maxtorcor} there exists a conjugation
of $G$ that maps $T_Q$ into $T_P$. This produces a new isomorphism
$C[[Q]]\iso A=k[[P]]$ that respects the grading, i.e. each $pk$
lies in some $qC$, and we obtain a surjective map $f{\colon}P\to Q$,
which is clearly a homomorphism. If $C=k$ then $f$ is an isomorphism
and we obtain that $P$ is determined by $A$.

Set $R_q=\prod_{p\in f^{-1}(q)}pk$. We have natural embeddings
$\prod_{p\in f^{-1}(q)}pk\into qC$ which are all isomorphisms
because $A=\prod_{q\in Q} R_q$ is isomorphic to $C[[Q]]=\prod_{q\in
Q}qC$. In particular, for $R=R_1$ we have that $C=\prod_{q\in R}
qk=k[[R]]$. Therefore,\linebreak[4] $A\iso k[[R]][[Q]]\iso k[[R\times Q]]$,
and since the monoid is determined by $A$ we obtain that $P\iso
Q\times R$.  \end{proof}

\section{Proof of Theorem \ref{8-TH. MODIFICATION} -- preliminary
steps}\label{8-prelimsec} The goal of \S\ref{8-prelimsec} is to
reduce the proof of Theorem \ref{8-TH. MODIFICATION} to the case
when the following conditions are satisfied: (1) $X$ is regular,
(2) the log structure is given by an snc divisor $Z$ which is {\em
$G$-strict} in the sense that for any $g\in G$ and a component
$Z_i$ either $gZ_i=Z_i$ or $gZ_i\cap Z_i=\emptyset$, (3) $G$ acts
freely on $X\setminus Z$ and for any geometric point $\ox\to X$
the inertia group $G_\ox$ is abelian.

\subsection{Plan} A general method for constructing a $G$-equivariant
morphism $f$ as in Theorem \ref{8-TH. MODIFICATION} is to construct
a tower of $G$-equivariant morphisms of log schemes $X'=X_n\dashto
X_0=X$, where the underlying morphisms of schemes are normalized
blow ups along $G$-stable centers sitting over $Z\cup T$, such that
various properties of the log scheme $X_i$ with the action of $G$
gradually improve to match all assertions of the Theorem. To simplify
notation, we will, as a rule, replace $X$ with the new log scheme
after each step. The three conditions above will be achieved in
three steps as follows.

\subsubsection{Step 1. Making $X$ regular.}\label{8-step1}
This is achieved by the saturated log blow up
tower $\tilcalF(X,Z){\colon}X'\dashto X$ from Theorem
\ref{8-deslogregth}. In particular, the morphism $X'\to X$ is
even log smooth. In the sequel, we assume that $X$ is regular, in
particular, $Z$ is a normal crossings divisor by \refext{VI}{6-1.7}. We will
call $Z$ the {\em boundary} of $X$. In the sequel, these conditions
will be preserved, so let us describe an appropriate restriction
on further modifications.

\subsubsection{Permissible blow ups}\label{8-permissec} After
Step 1 any modification in the remaining tower will be of the form
$f{\colon}(X',Z')\to(X,Z)$ where $X'=\Bl_V(X)$, $Z'=f^{-1}(Z\cup
V)$ and $V$ has {\em normal crossings} with $Z$, i.e. \'etale
locally on $X$ there exist regular parameters $t_1\. t_d$ such that
$Z=V(\prod_{i=1}^lt_i)$ and $V=V(t_{i_1}\. t_{i_m})$. We call such
modification {\em permissible} and a blow up $f{\colon}X'\to X$
is called {\em permissible} (with respect to the boundary $Z$)
if it underlies a permissible modification. Since there is an
obvious bijective correspondence between permissible modifications
and blow ups we will freely pass from one to another. Note that
$Z'=f^\st(Z)\cup E'$, where $E'=f^{-1}(V)$ is the exceptional
divisor.

\subsubsection{Permissible towers}\label{8-permistowersec} A {\em
permissible modification tower} $(X_d,Z_d)\dashto(X_0,Z_0)=(X,Z)$
is defined in the obvious way and we say that a blow up tower
$X_d\dashto X$ is permissible if it underlies such a modification
tower. Again, we will freely pass between permissible towers of
these types. Note that $Z_i=Z^\st_i\cup E_i$, where $Z_i^\st$
is the strict transform of $Z$ under $h_i{\colon}X_i\dashto X$
and $E_i$ is the exceptional divisor of $h_i$ (i.e. the union of
the preimages of the centers of $h_i$).

\begin{remark2}\label{8-permisrem} (i) Consider a permissible tower
as above. It is well known that for any $i$ one has that $X_i$
is regular, $Z_i$ is normal crossings and $E_i$ is even snc. For
completeness, let us outline the proof. Both claims follow from
the following: if $Z$ is snc then $Z_i$ is snc. Indeed, the claim
about $E_i$ follows by taking $Z=\emptyset$ and the claim about
$Z_i$ can be checked \'etale locally, so we can assume that
$Z$ is snc. Finally, if $Z$ is snc then $Z_i$ is snc by Lemma
\ref{8-regblowlem} below.

(ii) Permissible towers are very common in embedded desingularization
because they do not destroy regularity of the ambient scheme and the
normal crossings (or snc) property of the boundary (or accumulated
exceptional divisor). Even when one starts with an empty boundary,
a non-trivial boundary appears after the first step, and this
restricts the choice of further centers. Actually, any known
self-contained proof of embedded desingularization constructs a
permissible resolution tower.  \end{remark2}

\subsubsection{$G$-permissible towers}\label{8-Gpermissec} In
addition, we will only blow up $G$-equivariant centers $V$. So,
$f{\colon}X'=\Bl_V(X)\to X$ is $G$-equivariant and the exceptional
divisor $E=f^{-1}(V)$ is regular and $G$-equivariant and hence
$G$-strict. Such a blow up (or their tower) will be called {\em
$G$-permissible}. It follows by induction that the exceptional
divisor of such a tower is $G$-strict.

\subsubsection{Step 2. Making $Z$ snc and $G$-strict}\label{8-step2}
Consider the stratification of $Z$ by multiplicity: a point $z\in Z$
is in $Z^n$ if it has exactly $n$ preimages in the normalization
of $Z$. Note that $\{Z^n\}$ is precisely the log stratification
as defined in \S\ref{8-logstratsec}. By
{\em depth} of the stratification we mean the maximal $d$ such
that $Z^d\neq\emptyset$. In particular, $Z^d$ is the only closed
stratum. Step 2 proceeds as follows: $X_{i+1}\to X_i$ is the blow
up along the closed stratum of $Z_i^\st$.

\begin{remark2}\label{8-step2rem} What we use above is
the standard algorithm that achieves the following two
things: $Z'$ is snc and $Z^\st=\emptyset$ (see, for example,
\cite[7.2]{Alterations@de_Jong}). Even when $Z$ is snc, the second
property is often used in the embedded desingularization algorithms
to get rid of the old components of the boundary.  \end{remark2}

\subsubsection{Justification of Step 2}\label{8-step2proof}
Since the construction is well known, we just sketch the
argument. First, observe that $Z^d$ has normal crossings with
$Z$, that is, $X'=\Bl_{Z^d}(X)\to X$ is permissible. Thus, $Z'$
is normal crossings and hence $Z^\st$ is also normal crossings. A
simple computation with blow up charts shows that the depth of
$Z^\st$ is $d-1$ (for example, one can work \'etale-locally, and
then this follows from Lemma \ref{8-weaksnclem} below). It follows
by induction that the tower produced by Step 2 is permissible,
of length $d$ and with $Z_d^\st=\emptyset$. So, $Z_d=E_d$ is snc
by Remark \ref{8-permisrem} and $G$-strict by \S\ref{8-Gpermissec}.

\subsubsection{Step 3. Making the inertia groups abelian
and the action of $G$ on $X\setminus Z$ free}\label{8-step3}
Recall (\refext{VI}{6-stratification}) that the inertia strata are of the
form $X_H=X^H\setminus\cup_{H\subsetneq H'}X^{H'}$. Step 3 runs
analogously to Step 2, but this time we will work with the inertia
stratification of $X$ instead of the log stratification, and will
have to apply the same operation a few times. Let us first describe
the blow up algorithm used in this step; its justification will be
given in \S\ref{8-step3proof}.

Let $f_{\{X^H\}}{\colon}X'\dashto X$ denote the following blow
up tower. First we blow up the union of all closed strata
$X_H$. In other words, $V_0$ is the union of all non-empty
minimal $X^H$, i.e. non-empty $X^H$ that do not contain $X^K$
with $\emptyset\subsetneq X^K\subsetneq X^H$. Next, we consider
the family of all strict transforms of $X^H$ and blow up the union
of the non-empty minimal ones, etc. Obviously, the construction is
$G$-equivariant. We will prove in Proposition \ref{8-towerprop} that
$f_{\{X^H\}}$ is permissible of length bounded by the length of the
maximal chains of subgroups. Also, we will show in \S\ref{8-step3sec}
that $f_{\{X^H\}}$ decreases all non-abelian inertia groups,
so the algorithm of Step 3 goes as follows: until all inertia
groups become abelian, apply $f_{\{X^H\}}{\colon}X'\dashto X$
(i.e. replace $(X,Z)$ with $(X',Z')$).

\subsection{Justification of Step 3}\label{8-step3proof}

\subsubsection{Weakly snc families}\label{8-weaksncsec} Assume that
$X$ is regular, $Z\into X$ is an snc divisor, and $\{X_i\}_{i\in
I}$ is a finite collection of closed subschemes of $X$. For any
$J\subseteq I$ we denote by $X_J$ the scheme-theoretic intersection
$\cap_{j\in J}X_j$. The family $\{X_i\}$ is called {\em weakly
snc} if each $X_i$ is nowhere dense and $X_J$ is regular. The
family $\{X_i\}$ is called {\em weakly $Z$-snc} if it is weakly
snc and each $X_J$ has normal crossings with $Z$. In particular,
$\{X_i\}_{i\in I}$ is weakly snc (resp. weakly $Z$-snc) if and only
if the family $\{X_J\}_{\emptyset\neq J\subseteq I}$ is weakly snc
(resp. weakly $Z$-snc).

\begin{remark2}\label{8-weaksncrem} (i) Here is a standard criterion
of being an snc divisor, which is often taken as a definition. Let
$D\into X$ be a divisor with irreducible components $\{D_i\}_{i\in
I}$. Then $D$ is snc if and only if it is weakly snc and each
irreducible component of $D_J$ is of codimension $|J|$ in $X$.

(ii) The condition on the codimension is essential. For example,
$xy(x+y)=0$ defines a weakly snc but not snc divisor in
$\AA^2_k=\Spec(k[x,y])$.

(iii) The criterion from (i) implies that if $X$ is qe then the snc
locus of $D$ is open -- it is the complement of the union of singular
loci of $D_J$'s. (Note that this makes sense for all points of $X$
because $D$ is snc at a point $x\in X\setminus D$ if and only if
$X=D_\emptyset$ is regular at $x$.)
\end{remark2}

\begin{lemma2}\label{8-Zsnclem} Let $(X,Z)$ and $G$ be as achieved
in Step 2. Then the family $\{X^H\}_{1\neq H\subseteq G}$ is
weakly $Z$-snc.
\end{lemma2}

\begin{proof} Recall that for any
subgroup $H\subseteq G$ the fixed point subscheme $X^H$ is regular
by Proposition \refext{VI}{6-pts fixes reguliers}, and $X^H$ is nowhere dense for $H\neq 1$
by generic freeness of the action of $G$. Since for any pair of
subgroups $K,H\subseteq G$ we have that $X^H\times_XX^K=X^{KH}$,
the family is weakly snc. It remains to show that $Y=X^H$ has normal
crossings with $Z$. Note that it is enough to consider the case when
$X$ is local with closed point $x$ and $G=H$. The cotangent spaces
at $x$ will be denoted $T^*X=m_{X,x}/m^2_{X,x}$, $T^*Y$, etc. Their
dual spaces will be called the tangent spaces, and denoted $TX$,
$TY$, etc. Let $\phi^*:T^*X\onto T^*Y$ denote the natural map
and let $\phi:TY\into TX$ denote its dual. We will systematically
use without mention that $|G|$ is coprime to ${\rm char} k(x)$,
in particular, the action of $G$ on $T^*X$ is semi-simple.

The proof of \refext{VI}{6-pts fixes reguliers} also shows that for any point $x\in X^H$, the
tangent space $T_x(X^H)$ is isomorphic to $(T_xX)^H$. In particular,
$TY\iso(TX)^G$ and hence $\phi^*$ maps $(T^*X)^G\subset
T^*X$ isomorphically onto $T^*Y$. Therefore, $U=\Ker(\phi^*)$
is the $G$-orthogonal complement to $(T^*X)^G$, i.e. the
only $G$-invariant subspace such that $(T^*X)^G\oplus U\iso
T^*X$. Let $Z_i=V(t_i)$, $1\le i\le n$ be the components of $Z$ and
let $dt_i\in T^*X$ denote the image of $t_i$. By $Z$-strictness
of $G$, each line $L_i=\Span(dt_i)$ is $G$-invariant, so $G$
acts on $dt_i$ by a character $\chi_i$. Without restriction of
generality, $\chi_1\.\chi_l$ for some $0\le l\le n$ are the only
trivial characters. In particular, $L=\Span(dt_1\. dt_n)$ is the
direct sum of $L^G=\Span(dt_1\. dt_l)$ and its $G$-orthogonal
complement $L\cap U$, which (by uniqueness of the complement)
coincides with $\Span(dt_{l+1}\.dt_n)$. Complete the basis of
$L$ to a basis $\{dt_1\.dt_n,e_1\. e_m\}$ of $T^*X$ such that
$\{dt_{l+1}\.dt_n,e_1\. e_r\}$ for some $r\le m$ is a basis of
$U$ and choose functions $s_1\.s_m$ on $X$ so that $ds_j=e_j$
and $s_1\. s_r$ vanish on $Y$. Clearly, $t_1\. t_n,s_1\.s_m$ is a
regular family of parameters of $\calO_{X,x}$, so the lemma would
follow if we prove that $Y=V(t_{l+1}\.t_n,s_1\.s_r)$.

Since $Y$ is regular and $\Ker(\phi^*)$ is spanned by the images of
$t_{l+1}\. t_n,s_1\.s_r$, we should only check that these functions
vanish on $Y$. The $s_j$'s vanish on $Y$ by the construction, so we
should check that $t_i$ vanishes on $Y$ whenever $l<i\le n$. Using
the functorial definition from \refext{VI}{6-stratification} of the subscheme of fixed
points, we obtain that $Z_i^G=Z_i\times_XX^G$, hence $Z_i\times_XY$
is regular by \refext{VI}{6-pts fixes reguliers}. However, $TY$ is contained $TZ_i$, which
is the vanishing space of $dt_i$, hence we necessarily have that
$Y\into Z_i$.  \end{proof}

\subsubsection{Snc families}\label{8-sncsec} A family of nowhere
dense closed subschemes $\{X_i\}_{i\in I}$ is called {\em snc}
(resp. {\em $Z$-snc}) at a point $x$ if $X$ is regular at $x$ and
there exists a regular family of parameters $t_j\in\calO_{X,x}$ such
that in a neighborhood of $x$ each $X_i$ (resp. and each irreducible
component $Z_k$ of $Z$) passing through $x$ is given by the vanishing of a subfamily
$t_{j_1},\dots,t_{j_l}$. Note that the family $\{X_i\}_{i\in I}$
is $Z$-snc if and only if the union $\{X_i\}\cup\{Z_k\}$ is snc. A
family is {\em snc} if it is so at any point of $X$ (in particular,
$X$ is regular).

\begin{remark2}\label{8-sncrem} (i) It is easy to see that the family
$\{X^H\}_{H\subseteq G}$ is snc whenever $G$ is abelian. Indeed,
it suffices to show that for any point $x$ there exists a basis of
$T_xX$ such that each $T_x(X^H)$ is given by vanishing of some of the
coordinates. But this is so because the action of $G$ on $T_xX$ is
(geometrically) diagonalizable and $T_x(X^H)=(T_xX)^H$. In general,
the family $\{X^H\}_{H\subseteq G}$ does not have to be snc, as
the example of a dihedral group $D_n$ with $n\ge 3$ acting on the
plane shows.

(ii) If $Z\into X$ is an snc divisor with components $Z_i$ and
$V\into X$ is a closed subscheme then the family $\{Z_i,V\}$ is
snc if and only if $V$ has normal crossings with $Z$.  \end{remark2}

The transversal case of the following lemma can be deduced from
\egalong{IV}{4}{\S19.1}, but we could not find the general case in
the literature (although it seems very probable that it should have
appeared somewhere).  \begin{lemma2}\label{8-snclem} Any weakly snc
family with $|I|=2$ is snc.  \end{lemma2} \begin{proof} We should
prove that if $X=\Spec(A)$ is a regular local scheme and $Y$, $Z$
are regular closed subschemes such that $T=Y\times_XZ$ is regular
then there exists a regular family of parameters $t_1\.t_n\in A$
such that $Y$ and $Z$ are given by vanishing of some set of these
parameters. Let $m$ be the maximal ideal of $A$, and let $I$, $J$
and $K=I+J$ be the ideals defining $Y$, $Z$ and $T$, respectively. By
$T^*X=m/m^2$, $T^*Y$, etc., we denote the cotangent spaces at the
closed point of $X$. Note that $I/mI\iso\Ker(T^*X\to T^*Y)$,
and similar formulas hold for $J/mJ$ and $K/mK$. Indeed, we can
choose the parameters so that $Y=V(t_1\.t_l)$ and then the images
of $t_1\.t_l$ form a basis both of $I/mI$ and $\Ker(T^*X\to T^*Y)$.

Now, let us prove the lemma. Assume first that $Y$ and $Z$ are
transversal, i.e. $T^*X\into T^*Y\oplus T^*Z$. Choose elements
$t_1\. t_{l+k}$ such that $Y=V(t_1\.t_l)$, $Z=V(t_{l+1}\.t_{l+k})$,
$l=\codim(Y)$ and $k=\codim(Z)$. Then the images $dt_i\in T^*X$
of $t_i$ are linearly independent because $dt_1\.dt_l$ span
$\Ker(T^*X\to T^*Y)$ and $dt_{l+1}\.dt_{l+k}$ span $\Ker(T^*X\to
T^*Z)$. Hence we can complete $t_i$'s to a regular family of
parameters by choosing $t_{l+k+1}\.t_n$ such that $dt_1\.dt_n$
is a basis of $T^*X$. This proves the transversal case, and to
establish the general case it now suffices to show that if $Y$
and $Z$ are not transversal then there exists an element $t_1\in
m\setminus m^2$ which vanishes both on $Y$ and $Z$. (The we can
replace $X$ with $X_1=V(t_1)$ and repeat this process until $Y$
and $Z$ are transversal in $X_a=V(t_1\.t_a)$.) Tensoring the
exact sequence $0\to I\cap J\to I\oplus J\to K\to 0$ with $A/m$
we obtain an exact sequence $$(I\cap J)/m(I\cap J)\to I/mI\oplus
J/mJ\stackrel\phi\to K/mK\to 0$$ The failure of transversality
is equivalent to non-injectivity of $\phi$, hence there exists
an element $f\in I\cap J$ with a non-zero image in $I/mI\oplus
J/mJ$. Thus, $f\in m\setminus m^2$ and we are done.  \end{proof}

\subsubsection{Blowing up the minimal strata of a weakly
snc family}\label{8-weakblowsec} Given a weakly snc family
$\{X_i\}_{i\in I}$, we say that a scheme $X_J$ with $J\subseteq I$
is {\em minimal} if it is non-empty and any $X_{J'}\subsetneq X_J$
is empty. Also, we will need the following notation: if $Z\into X$
is a closed subscheme and $D\into X$ is a Cartier divisor with the
corresponding ideals $\calI_Z,\calI_D\subset\calO_X$, then $Z+D$
is the closed subscheme defined by the ideal $\calI_Z\calI_D$.

\begin{proposition2}\label{8-weaksnclem} Assume that $X$ is regular,
$Z\into X$ is an snc divisor with irreducible components $Z_1\. Z_l$,
and $\{X_i\}_{i\in I}$ is a $Z$-snc (resp. weakly $Z$-snc) family
of subschemes. Let $V$ be the union of all non-empty minimal
subschemes $X_J$, $f{\colon}X'=\Bl_V(X)\to X$, $X'_i=f^\st(X_i)$
and $Z'=f^{-1}(Z\cup V)$. Then

(i) $X'$ is regular and $Z'$ is snc.

(ii) The family $\{X'_i\}_{i\in I}$ is $Z'$-snc (resp. weakly
$Z'$-snc).

(iii) For any $J\subseteq I$, the scheme-theoretical intersection
$X'_J=\cap_{j\in J}X'_j$ coincides with $f^\st(X_J)$.

(iv) For any $J\subseteq I$ the total transform $X_J\times_XX'$ is of
the form $X'_J+D'_J$ where $D'_J$ is the divisor consisting of all
connected components of $E'=f^{-1}(V)$ contained in $f^{-1}(X_J)$.
\end{proposition2} \begin{proof} We start with the following lemma.

\begin{lemma2}\label{8-regblowlem} Assume that $X$ is regular,
$Z$ is an snc divisor and $V\into Y$ are closed subschemes
having normal crossings with $Z$. Let $f{\colon}X'\to X$ be
the blow up along $V$, $Y'=f^\st(Y)$, $Z'=f^{-1}(Z\cup V)$, and
$E'=f^{-1}(V)$. Then $Z'$ is snc, $Y'$ has normal crossings with
$Z'$ and $Y\times_XX'=Y'+E'$.  \end{lemma2} \begin{proof} The proof
is a usual local computation with charts. Take any point $u\in V$
and choose a regular family of parameters $t_1\.t_n$ at $u$ such
that $Y$ (resp. $V$, resp. $Z$) are given by the vanishing of
$t_1\.t_m$ (resp. $t_1\.t_l$, resp. $\prod_{i\in I}t_i$), where
$0\le m\le l\le n$ and $I\subseteq\{1\. n\}$. Locally over $u$
the blow up is covered by $l$ charts, and the local coordinates on
the $i$-th chart are $t'_j$ such that $t'_j=t_j$ for $j>l$ or $j=i$
and $t'_j=\frac{t_j}{t_i}$ otherwise. On this chart, $Y\times_XX'$
(resp. $Y'$, resp. $E'$, resp. $Z'$) is given by the vanishing of
$t_1\.t_m$ (resp. $t'_1\.t'_m$, resp. $t'_i$, resp. $\prod_{j\in
I\cup\{i\}}t'_j$), hence the lemma follows.  \end{proof}

The lemma implies (i). In addition, it follows from the
lemma that $f^\st(X_J)$ has normal crossings with $Z'$ and
$X_J\times_XX'=f^\st(X_J)+D'_J$. Thus, (iii) implies (iv), and
(iii) implies (ii) in the case when the family $\{X_i\}_{i\in I}$
is weakly $Z$-snc. Note also that if this family is $Z$-snc then
locally at any point $x\in X$ there exists a family of regular
parameters $\underline{t}=\{t_1\. t_n\}$ such that each $X_J$ and
each component of $Z$ is given by the vanishing of a subfamily of
$\underline{t}$ locally in a neighborhood of $x$. Then the same local
computation as was used in the proof of Lemma \ref{8-regblowlem}
proves also claims (ii) and (iii) of the proposition. So, it remains
to prove (iii) when the family is weakly snc. It suffices to prove
that if (iii) holds for $X_J$ and $X_K$ then it holds for $X_{J\cup
K}$. Moreover, (iii) does not involve the boundary so we can assume
that $Z=\emptyset$. It remains to note that $\{X_J,X_S\}$ is an
snc family by Lemma \ref{8-snclem}, hence our claim follows from
the snc case.  \end{proof}

\subsubsection{Blow up tower of a weakly $Z$-snc
family}\label{8-weaktowersec} Let $X$ be a regular scheme, $Z$
be an snc divisor and $\{X_i\}_{i\in I}$ be a weakly $Z$-snc
family. By the blow up tower $f_{\{X_i\}}$ of $\{X_i\}$ we mean
the following tower: the first blow up $h_1:X_1\to X_0=X$ is
along the union of all non-empty minimal schemes of the form $X_J$
for $\emptyset\neq J\subseteq I$, the second blow up is along the
union of all non-empty minimal schemes of the form $h_1^\st(X_J)$
for $\emptyset\neq J\subseteq I$, etc.

\begin{proposition2}\label{8-towerprop} Keep the above notation. Then
the tower $f_{\{X_i\}}$ is permissible with respect to $Z$ and
its length equals to the maximal length of chains $\emptyset\neq
X_{J_1}\subsetneq\dots\subsetneq X_{J_d}$ with $\emptyset\neq
J_d\subsetneq\dots\subsetneq J_1\subseteq I$. Furthermore, the strict
transform of any scheme $X_J$ is empty and the total transform of
$X_J$ is a Cartier divisor.  \end{proposition2} \begin{proof} The
claim about the length is obvious. Let $f_{\{X_i\}}:X_d\dashto X_0=X$
and let $h_n:X_n\dashto X_0=X$ be its $n$-th truncation. For each
$i\in I$ set $X_{n,i}=h_n^\st(X_i)$ and for each $J\subseteq I$ set
$X_{n,J}=\cap_{i\in J}X_{n,i}$. Using Proposition \ref{8-weaksnclem}
and straightforward induction on length, we obtain that the family
$\{X_{n,i}\}_{i\in I}$ is weakly $Z_n$-snc, $X_{n,J}=h_n^\st(X_J)$,
the blow up $X_{n+1}\to X_n$ is along the union of non-empty
minimal $X_{n,J}$'s, and $X_J\times_XX_n=X_{n,J}+D_{J,n}$, where
$D_{J,n}$ is a divisor. So, the tower is permissible, and since
$X_{d,J}=\emptyset$ we also have that $X_J\times_XX_d$ is a divisor.
\end{proof}

\begin{remark2} We will not need this, but it is easy to deduce from
the proposition that on the level of morphisms the modifications
$X_d\to X$ is isomorphic to the blow up along $\prod_{\emptyset\neq
J\subseteq I}\calI_{X_J}$.  \end{remark2}

\subsubsection{Justification of Step 3}\label{8-step3sec} The
blow up tower $f_{\{X^H\}}{\colon}X'\dashto X$ from Step 3 is
$G$-equivariant in an obvious way, and it is permissible by
Proposition \ref{8-towerprop}. In addition, $Z'=f^{-1}(Z)\cup
f^{-1}(\cup_{1\neq H\subseteq G} X^H)$ and, since $G$ acts freely
on $X'\setminus f^{-1}(\cup_{1\neq H\subseteq G} X^H)\iso
X\setminus\cup_{1\neq H\subseteq G}X^H$, it also acts freely on
$X'\setminus Z'$. It remains to show that applying $f_{\{X^H\}}$ we
decrease all non-abelian inertia groups. Namely, for any $x'\in X'$
mapped to $x\in X$ we want to show that either $G_\ox$ is abelian
or the inclusion $G_{\ox'}\subsetneq G_\ox$ is strict.

Let $H\subseteq G$ be any non-abelian subgroup with commutator
$K=[H,H]$. Since $X^K\times_XX'$ is a divisor by Proposition
\ref{8-towerprop}, the universal property of blow ups implies that
$X'\to X$ factors through $Y=\Bl_{X^K}(X)$. On the other hand, it
is proved in \refext{VI}{6-eclatements} that $Y\to X$ is an $H$-equivariant
blow up, and if a geometric point $\ox\to X$ with $G_\ox=H$ lifts to
a geometric point $\oy\to Y$ then $G_\oy\subsetneq H$. Therefore,
the same is true for the $G$-equivariant modification $X'\to X$,
and we are done.

\section{Proof of Theorem \ref{8-TH. MODIFICATION} -- abelian
inertia}\label{8-endsec}

\subsection{Conventions}\label{8-notationendsec} Throughout
\S\ref{8-endsec} we assume that $(X,Z)$ and $G$ satisfy
all conditions achieved at Steps 1, 2, 3, and our aim is to
construct a modification $f_{(G,X,Z)}{\colon}X'\to X$ as in Theorem
\ref{8-TH. MODIFICATION}. Unless specially mentioned, we do not assume 
that $X$ is qe. This is done in order to isolate the only place where
this assumption is needed (existence of rigidifications).

\subsection{Outline of our method and other
approaches}\label{8-outlinesec}

\subsubsection{Combinatorial nature of the problem}\label{8-combsec}
On the intuitive level it is natural to expect that "everything
relevant to our problem" should be determined by the following
"combinatorial" data: the log structure of $X$, the inertia
stratification of $X$ by $X_H:=X^H\setminus\cup_{H'\subsetneq
H}X^{H'}$ and the representations of the inertia groups on the
tangent spaces (which are essentially constant along $X_H$). This
combinatorial nature is manifested in both approaches to the problem
that we describe below.

\subsubsection{Combinatorial algorithm}\label{8-algorsec}
The most natural approach is to seek for a "combinatorial
algorithm" that iteratively blows up the closures of {\em
log-inertia strata} (i.e. intersections of a log stratum with
an inertia stratum) or some generalization thereof. Actually, our (very simple) algorithms in
\S\ref{8-prelimsec} were of this type. The algorithm should be
governed by a combinatorial data, such as the number of components
of $Z$ through a point $x$ and the history of their appearance
(similarly to the desingularization algorithms), the representation
of $G_\ox$ on $T_x$ plus some history (e.g. representation of $G_\ox$
on the tangent spaces to the components of $Z$), etc.

It is natural to expect that building such an algorithm would lead to
a relatively simple proof of Theorem \ref{8-TH. MODIFICATION}. In
particular, it would be non-sensitive to quasi-excellence
issues. Unfortunately, despite partial positive results, we could
not construct such an algorithm.
Thus, the question whether such an algorithm exists remains
open.\footnote{F. Pop told to the second author that he has a plan
of constructing such a combinatorial algorithm.}

\subsubsection{Our method}\label{8-ourmethodsec} A general plan
of our method is as follows. In \S\ref{8-localsec} we will show
that such a modification $f$ exists \'etale-locally on the base
if $X$ is qe. A priori, our construction will be canonical up to
an auxiliary choice, but then we will prove in \S\ref{8-globalsec}
that actually it is independent of the choice and hence descends to a
modification $f$ as required. To prove independence we will show in
\S\ref{8-functsec} that the construction is functorial with respect
to strict {\em inert} morphisms (i.e. morphisms that preserve both
the log and the inertia structures, see \S\ref{8-inertsec}). The
latter is a manifestation of the "combinatorial nature" of our
algorithm.

\subsection{Local construction}\label{8-localsec}

\subsubsection{Very tame action and Zariski
topology}\label{8-verytamesec} Note that the log structure on
$(X,Z)$ is Zariski since $Z$ is snc. Thus, it will be convenient
to describe very tame action in terms of Zariski topology. We say
that the action of $G$ is {\em very tame at} a point $x\in X$ if
for any geometric point $\ox$ over $x$ the action is very tame at
$\ox$. Since the log structure is Zariski, the $G_\ox$-equivariant
log scheme $\Spec(\calO_{X,\ox})$ is independent of the choice of
$\ox$ up to an isomorphism. In particular, the action is very tame
at $x$ if and only if it is very tame at a single geometric point
$\ox$ above $x$.

\begin{lemma2}\label{8-tamelocallem} Assume that $(G,X,Z)$ is as
in \S\ref{8-notationendsec}, $x\in X$ is a point, and $T$ is the
set of points of $X$ at which the action is very tame. Then

(i) $T$ is open,

(ii) $x\in T$ if and only if locally at $x$ the log stratification of
$X$ is finer than the inertia stratification (i.e. any log stratum
is contained in an inertia stratum).  \end{lemma2} \begin{proof}
Choose a geometric point $\ox$ above $x$. If $x\in T$ then by
Corollary \refext{VI}{6-tres modere dans voisinage} the action is very tame on a suitable \'etale
neighborhood $X'$ of $\ox$. The image of $X'$ in $X$ is an open
neighborhood of $x$ which is contained in $T$. Thus, $T$ is open. To
prove the second claim we note that the conditions (i) and (ii) from
\refext{VI}{6-3.1} are automatically satisfied at $\ox$. Indeed, (i) is satisfied
because the action is tame by assumption of \ref{8-TH. MODIFICATION}
and (ii) is satisfied because $Z$ is snc and $G$-strict. Since
log and inertia stratifications are compatible with the strict
henselization morphism, condition (iii) from \refext{VI}{6-3.1} is satisfied
at $\ox$ if and only if the log stratification is finer than the
inertia stratification at $x$. This proves (ii).  \end{proof}

\subsubsection{Admissibility}\label{8-admissec} In the sequel,
by saying that $X$ is {\em admissible} we mean that the action
of $G$ on $X$ is admissible in the sense of \sga{1}{V}{1.7}
(e.g. $X$ is affine). This is needed to ensure that $X/G$ exists as
a scheme. An alternative would be to allow $X/G$ to be an algebraic
space (and $(X/G,Z/G)$ to be a log algebraic space).

\subsubsection{Rigidification}\label{8-rigidsec} By a {\em
rigidification} of $X$ we mean a $G$-equivariant normal crossings
divisor $\oZ$ that contains $Z$ and such that the action of $G$ on
the log regular log scheme $\oX=(X,\oZ)$ is very tame. If $\oZ$ is
snc then we say that the rigidification is {\em strict}. Sometimes,
by a rigidification of $X=(X,Z)$ we will mean the log scheme $\oX$
itself. Our construction of a modification $f$ uses a rigidification,
so let us first establish local results on existence of the latter.

\begin{lemma2}\label{8-rigidlem} Let $X,Z,G$ be as
in \S\ref{8-notationendsec}. Assume that $X=\Spec(A)$ is a
local scheme with closed point $x$, $G=G_\ox$, and $\mu_N\subset
A$, where $N$ is the order of $G$. Then $X$ possesses a strict
rigidification.  \end{lemma2} \begin{proof} Choose $t_1\. t_n\in A$
such that $Z_i=(t_i)$ are the components of $Z$. By $G$-strictness
of $Z$, for any $g\in G$ we have that $Z_i=(gt_i)$ locally at
$x$, in particular, the tangent space to each $Z_i$ at $x$ is
$G$-invariant. Now, we can use averaging by the $G$-action to make
the parameters $G$-equivariant. Namely, $G$ acts on $dt_i$ by a
character $\chi_i$ and replacing $t_i$ with $\frac{1}{|G|}\sum_{g\in
G}\frac{gt_i}{\chi_i(g)}$ we do not change $dt_i$ and achieve that
$gt_i=\chi_i(g)t_i$.

The action of $G$ on the cotangent space at $x$ is diagonalizable
because $G$ is abelian and $\mu_N\subset k(x)$. In particular,
we can complete the family $dt_1\. dt_n$ to a basis $dt_1\. dt_l$
such that $t_i\in\calO_{X,x}$ and $G$ acts on each $dt_i$ by
a character $\chi_i$. Then $t_1\. t_l$ is a regular family of
parameters of $\calO_{X,x}$ and using the same averaging procedure
as above we can make them $G$-equivariant. Take now $\oZ$ to be the
union of all divisors $(t_i)$ with $1\le i\le l$. The action of $G$
on $(X,\oZ)$ is very tame at $x$ because it is the only point of
its log stratum. So, it remains to use Lemma \ref{8-tamelocallem}.
\end{proof}

\subsubsection{Inert morphisms}\label{8-inertsec} Let $(X,Z)$
and $G$ be as in \S\ref{8-notationendsec}. Our next aim is to
find an \'etale cover $f:(Y,T)\to(X,Z)$ which "preserves" the
log-inertia structure of $(X,Z)$ and such that $(Y,T)$ admits a
rigidification. The condition on the log structure is obvious:
we want $f$ to be strict, i.e. $f^{-1}(Z)=T$. Let us introduce a
restriction related to the inertia groups.

Assume that $(Y,T)$ with an action of $H$ is another such triple,
and let $\lambda:H\to G$ be a homomorphism. A $\lambda$-equivariant
morphism $f{\colon}(Y,T)\to(X,Z)$ will be called {\em inert} if
for any point $y\in Y$ with $x\in X$ the induced homomorphism of
inertia groups $G_\oy\to G_\ox$ is an isomorphism. In particular,
the inertia stratification of $Y$ is the preimage of the inertia
stratification of $X$.

\begin{lemma2}\label{8-rigidlem2} Let $X$, $Z$ and $G$ be as above,
and assume that $X$ is qe. Then there exists a $G$-equivariant
surjective \'etale inert strict morphism $h{\colon}(Y,T)\to(X,Z)$
such that $Y$ is affine and $(Y,T)$ possesses a strict
rigidification.  \end{lemma2} \begin{proof} First, we note that the
problem is local on $X$. Namely, it suffices for any point $x\in
X$ to find a $G_\ox$-equivariant \'etale inert strict morphism
$h{\colon}(Y,T)\to(X,Z)$ such that $Y$ is affine, $x\in h(Y)$, any
point $x'\in h(Y)$ satisfies $G_{\ox'}\subseteq G_\ox$, and $(Y,T)$
admits a strict rigidification. Indeed, $h$ can be extended to a
$G$-equivariant morphism $Y\times_X(X\times G/G_\ox)=\coprod_{g\in
G/G_\ox}Y_g\to X$, where each $Y_g$ is isomorphic to $Y$ and
the morphism $Y_g\to X$ is obtained by composing $Y\to X$ with
$g{\colon}X\to X$. Clearly, the latter morphism is \'etale, inert,
and strict, and by quasi-compactness of $X$ we can combine finitely
many such morphism to obtain a required cover of $(X,Z)$.

Now, fix $x\in X$ and consider the $G_\ox$-invariant neighborhood
$X'=X\setminus\cup_{H\nsubseteq G_\ox}X^H$ of $x$. We
will work over $X'$, so the condition  $G_{\ox'}\subseteq
G_{\ox}$ will be automatic. Let $N$ be the order of
$G_\ox$, and consider the $G_\ox$-equivariant morphism
$f{\colon}Y=X'\times\Spec(\ZZ[\frac{1}{N},\mu_N])\to X$ (with
$G_\ox$ acting trivially on the second factor). Let $T=f^{-1}(Z)$
and let $y$ be any lift of $x$. It suffices to show that
$(Y,T)$ admits a rigidification in an affine $G_\ox$-invariant
neighborhood of $y$ (such neighborhoods form a fundamental
family of neighborhoods of $y$). By Lemma \ref{8-rigidlem}, the
localization $Y_y=\Spec(\calO_{Y,y})$ with the restriction $T_y$
of $T$ possesses a strict rigidification $\oT_y$. Clearly, $\oT_y$
extends to a divisor $\oT$ with $T\into\oT\into Y$ and we claim
that it is a rigidification in a neighborhood of $y$. Indeed, $\oT$
is snc at $y$, hence it is snc in a neighborhood of $y$ by Remark
\ref{8-weaksncrem}, and it remains to use Lemma \ref{8-tamelocallem}.
\end{proof}

\subsubsection{Main construction}\label{8-constrsec} Assume,
now, that $X=(X,Z)$ is admissible and admits a rigidification
$\oZ$. We are going to construct a $G$-equivariant modification
$$f_{(G,X,Z,\oZ)}{\colon}(X',Z')\to(X,Z)$$ such that $G$ acts very
tamely on the target and $f_{(G,X,Z,\oZ)}$ is independent of the
rigidification. The latter is a subtle property (missing in the
obvious modification $(X,\oZ)\to(X,Z)$), and it will take us a
couple of pages to establish it.

The quotient log scheme $\oY=(Y,\oT)=(\oX/G,\oZ/G)$ is log regular
by Theorem \refext{VI}{6-tres modere}, hence by Theorem \ref{8-monoidaldesth}
there exists a functorial saturated log blow up tower
$\oh=\tilcalF^\rmlog(\oY){\colon}\oY'=(Y',\oT')\to\oY$ with a regular
and log regular source. Let $\of{\colon}\oX'=(X',\oZ')\to\oX$
be the pullback of $\oh$ (as a saturated log blow up tower, see
\S\ref{8-logpullsec}), then $\alpha'{\colon}\oX'\to\oY'$ is a Kummer
\'etale $G$-cover because $\alpha{\colon}\oX\to\oY$ is so by \refext{VI}{6-tres modere}
as the square $$ \xymatrix{(X',\oZ') \ar[d]^{\alpha'} \ar[r]^\of&
(X,\oZ) \ar[d]^\alpha \\ (Y',\oT') \ar[r]^\oh & (Y,\oT)} $$ is
cartesian in the category of fs log scheme.

Since $G$ acts freely on $U=X-Z$, $V=U/G$ is regular and $\oT|_V$
is snc. In particular, $\oh$ is an isomorphism over $V$ and hence
$\of$ is an isomorphism over $U$.  We claim that the Weil divisor
$T=Z/G$ of $Y$ is $\QQ$-Cartier (it does not have to be Cartier,
as the orbifold case with $X=\AA^2$, $Z=\AA^1$ and $G=\{\pm 1\}$
shows). Indeed, it suffices to check this \'etale-locally at a
point $y\in Y$. In particular, we can assume that $\mu_N\subset
\OO_{Y,y}$, where $N=|G|$. Then, as we showed in the proof of Lemma
\ref{8-rigidlem}, $Z$ can be locally defined by equivariant
parameters, in particular, $Z=V(f)$ where $G$ acts on $f$ by
characters. Therefore, $f^N$ is $G$-fixed, and we obtain that
$T$ is the reduction of the Cartier divisor $C$ of $Y$ given by
$f^N=0$. (The same argument applied to $T\times_YX$ shows that $NT$
is Cartier, so $T$ is $\QQ$-Cartier.) So, $C'=C\times_YY'$ is a
Cartier divisor whose reduction is $T'=\oh^{-1}(T)$. Since $Y'$
is regular and $T'$ lies in the snc divisor $\oT'$, we obtain that
$T'$ is itself an snc divisor.

Let $Z'$ denote the divisor $\alpha'^{-1}(T')=\of^{-1}(Z)$. Since
$X'\to Y'$ is \'etale over $V=Y'-T'$, the morphism of log schemes
$(X',Z')\to(Y',T')$ is a Kummer \'etale $G$-cover. This follows
from a variant of the classical Abhyankar's lemma (IX 2.1), which
is independent of the results of the present expos\'e.

In particular, $X'=(X',Z')$ is log regular and it follows that the
action of $G$ on $X'$ is very tame. We define $f_{(G,X,Z,\oZ)}$
to be the modification $X'\to X$.

\begin{remark2}\label{8-constrrem} (i) Note that $X'\to X$ satisfies
all conditions of Theorem \ref{8-TH. MODIFICATION} because the action
is very tame and $G$ acts freely on $X'\setminus f^{-1}(Z)$. So, we
completed the proof in the case when $(X,Z)$ admits a rigidification
$\oZ$. Our last task will be to get rid of the rigidification.

(ii) The only dependence of our construction on the rigidification
is when we construct the resolution of $(Y,\oT)$. Conjecturally,
it depends only on the scheme $Y$, and then $(Y',T')$, and hence
also $(X',Z')$, would depend only on $(X,Z)$. Recall that we
established in Theorem \ref{8-Morel 2.4.2} the particular case
of this conjecture when all maximal points of the log strata of
$(Y,\oT)$ are of characteristic zero. Hence independence of the
rigidification is unconditional in this case, and, fortunately,
this will suffice.  \end{remark2}

\subsubsection{Finer structure of
$f_{(G,X,Z,\overline{Z})}$}\label{8-finersec} Obviously, the saturated
log blow up tower $\of{\colon}(X',\oZ')\to(X,\oZ)$ depends on the
rigidification, and this is the reason why we prefer to consider
the modification $f_{(G,X,Z,\oZ)}{\colon}(X',Z')\to(X,Z)$
instead. However, there is an additional structure on
$f_{(G,X,Z,\oZ)}$ that has a chance to be independent of $\oZ$, and
which should be taken into account. By \S\ref{8-desingloregsec}, the
modification of schemes $f{\colon}X'\to X$ has a natural structure of
a normalized blow up tower $X_\bullet$ with $X=X_0$ and $X'=X_n$. Note
also that the tower contains no empty blow ups because this is true
for $\tilcalF^\rmlog(X/G,\oZ/G)$ and $f_{(G,X,Z,\oZ)}$ is its strict
transform with respect to the surjective morphism $X\to X/G$.

Note also that the log structure on $(X',Z')$ is reconstructed
uniquely from $f$ because $Z'=f^{-1}(Z)$ and $(X',Z')$ is
saturated and log regular. Therefore, it is safe from now on to
view $f_{(G,X,Z,\oZ)}$ as a normalized blow up tower of $X$, but
the modification of log schemes $(X',Z')\to(X,Z)$ will also be
denoted as $f_{(G,X,Z,\oZ)}$.

\begin{remark2}\label{8-finnormrem} Although we do not assume that
$X$ is qe, all normalizations in the tower $f_{(G,X,Z,\oZ)}$ are
finite. This happens because they underly saturations of fine log
schemes, which are always finite morphisms.  \end{remark2}

\subsection{Functoriality}\label{8-functsec} Clearly, the
construction of $f$ depends canonically on $(G,X,Z,\oZ)$,
i.e. is compatible with any automorphism of such
quadruple. Our next aim is to establish functoriality with
respect to strict inert $\lambda$-equivariant morphisms
$\phi{\colon}(H,Y,T,\oT)\to(G,X,Z,\oZ)$ (i.e. morphisms
that "preserve the combinatorial structure"). For this one
has first to study the quotient morphism of log schemes
$\tilphi{\colon}(Y/H,\oT/H)\to(X/G,\oZ/G)$.

\subsubsection{Log structure of the quotients}\label{8-logquotsec}
Recall the following facts from Proposition \refext{VI}{6-tres modere ponctuel}(b)
and its proof. Assume that $X=(X,M_X)$ is an fs log scheme provided with a
very tame action of a group $G$. After replacing $X$ with its strict
localization at a geometric point $\ox$, it admits an equivariant
chart $X\to\Spec(\Lambda[Q])$, where $\Lambda=\ZZ[1/N,\mu_N]$
for the order $N$ of $G_\ox$, $Q$ is an fs monoid and the
action of $\oG_x$ is via a pairing $\chi{\colon}G_\ox\otimes
Q\to\mu_N$. Moreover, if $P\subseteq Q$ is the maximal submonoid with
$\chi(G_\ox\otimes P)=1$ then $\Spec(\Lambda[Q])\to\Spec(\Lambda[P])$
is a chart of $X\to X/G_\ox$. Now, let us apply this description
to the study of $\tilphi$.

\begin{proposition2}\label{8-logquotprop} Assume that fs log schemes
$X$, $Y$ are provided with admissible very tame actions of groups $G$
and $H$, respectively,  $\lambda{\colon}H\to G$ is a homomorphism,
and $\phi{\colon}Y\to X$ is a strict inert $\lambda$-equivariant
morphism. Then the quotient morphism $\tilphi{\colon}Y/H\to X/G$
is strict.  \end{proposition2} \begin{proof} Fix a geometric point
$\oy$ of $Y$ and let $\ox$ be its image in $X$. It suffices to
show that $\tilphi$ is strict at the image of $\oy$ in $Y/H$. The
morphism $Y/H_\oy\to Y/H$ is strict (and \'etale) over the image of
$\oy$, and the same is true for $X$. Therefore we can replace $H$
and $G$ with $H_\oy\iso G_\ox$, and then we can also replace $X$
and $Y$ with their strict localizations at $\ox$ and $\oy$. Now,
the morphism $X\to\tilX=X/G_\ox$ admits an equivariant chart
$h{\colon}\Spec(\Lambda[Q])\to\Spec(\Lambda[P])$ as explained before
the proposition. Since $\phi$ is strict, the induced morphism
$Y\to\Spec(\Lambda[Q])$ is also a chart and hence $h$ is also a
chart of $Y\to\Y/G_\oy$. Thus, $\tilphi$ is strict.  \end{proof}

\subsubsection{An application to functoriality of $f_\bullet$} Assume
that $(G,X,Z,\oZ)$ is as earlier, and let $(H,Y,T,\oT)$ be another
such quadruple (i.e. $(Y,T)$ with the action of $H$ satisfies
conditions of Steps 1, 2, 3 and $(Y,\oT)$ is its rigidification).

\begin{corollary2}\label{8-logquotcor} Assume that
$\lambda{\colon}H\to G$ is a homomorphism and $\phi{\colon}Y\to
X$ is a $\lambda$-equivariant inert morphism such that
$T=Z\times_XY$ and $\oT=\oZ\times_XY$. Then $f_\bullet$ is
compatible with $\phi$ in the sense that $f_{(H,Y,T,\oT)}$ is
the contraction of $\phi^\st(f_{(G,X,Z,\oZ)})$. In addition,
$\phi^\st(f_{(G,X,Z,\oZ)})=f_{(G,X,Z,\oZ)}\times_XY$.
\end{corollary2} \begin{proof} The morphism
$\psi{\colon}(Y,\oT)\to(X,\oZ)$ is strict, hence the morphism of
quotients is strict by Proposition \ref{8-logquotprop}, and by
functoriality of saturated monoidal desingularization we obtain
that $\tilcalF^\rmlog(Y/H,\oT/H)$ is the contracted pullback
of $\tilcalF^\rmlog(X/G,\oZ/G)$. So, both $f_{(H,Y,T,\oT)}$ and
$f_{(G,X,Z,\oZ)}$ are obtained as the contraction of the strict
transform of $\tilcalF^\rmlog(X/G,\oZ/G)$. The first claim of the
Corollary follows.

Furthermore, $f_{(G,X,Z,\oZ)}$ underlies a log blow
up tower of $(X,\oZ)$ which is the strict transform of
$\tilcalF^\rmlog(X/G,\oZ/G)$, and the same is true for
$f_{(H,Y,T,\oT)}$. Since $\psi$ is strict it follows from Lemma
\ref{8-logregblowlem}(ii) that the strict transform is a pullback,
i.e. $\phi^\st(f_{(G,X,Z,\oZ)})=f_{(G,X,Z,\oZ)}\times_XY$.
\end{proof}

\subsubsection{Localizations and completions}\label{8-generalizsec}
In particular, it follows that the construction of $f_\bullet$
is compatible with localizations and completions. Namely, if
$x\in X$ is a point, $X_x=\Spec(\calO_{X,x})$, $Z_x=Z\times_XX_x$
and $\oZ_z=\oZ\times_XX_x$, then $f_{(G_\ox,X_x,Z_z,\oZ_x)}$ is
the contraction of $f_{(G,X,Z,\oZ)}\times_XX_x$. Similarly, if
$\hatX_x=\Spec(\hatcalO_{X,x})$, $\hatZ_x=Z\times_X\hatX_x$
and $\widehat{\oZ}_x=\oZ\times_X\hatX_x$, then
$f_{(G_\ox,\hatX_x,\widehat{Z}_x,\widehat{\oZ}_x)}$ is the
contraction of $f_{(G,X,Z,\oZ)}\times_X\hatX_x$.

\subsection{Globalization}\label{8-globalsec} To complete the
proof of Theorem \ref{8-TH. MODIFICATION} it suffices to show that
$f_{(G,X,Z,\oZ)}$ is independent of $\oZ$, and hence the local
constructions glue to a global normalized blow up tower. The main
idea is to simultaneously lift two rigidifications to characteristic
zero and apply Theorem \ref{8-Morel 2.4.2}.

\subsubsection{Independence of rigidification}\label{8-indeprigsec}
We start with the case of complete local rings. Then the problem
is solved by lifting to characteristic zero and referencing to
\ref{8-Morel 2.4.2}. The general case will follow rather easily.

\begin{lemma2}\label{8-formallem} Keep assumptions on $(X,Z)$
and $G$ as in \S\ref{8-notationendsec} and assume, in addition,
that $X=\coprod_{i=1}^m\Spec(A_i)$ where each $A_i$ is a complete
noetherian regular local ring with a separably closed residue
field. Then for any pair of rigidifications $\oZ$ and $\oZ'$ the
equality $f_{(G,X,Z,\oZ)}=f_{(G,X,Z,\oZ')}$ holds.  \end{lemma2}
\begin{proof} Almost the whole argument runs independently on
each irreducible component, so assume first that $X=\Spec(A)$
is irreducible. By Remark \ref{8-constrrem}(ii), it suffices
to consider the case when ${\rm char}(k)=p>0$, so let $C(k)$
be a Cohen ring of $k$. Note that we can work with $H=G_\ox$
instead of $G$ because $f_{(H,X,Z,\oZ)}=f_{(G,X,Z,\oZ)}$ by
Corollary \ref{8-logquotcor}. Since $H$ acts trivially on $k$,
for any element $t\in A$ its $H$-averaging is an element of $A^H$
with the same image in the residue field. Hence $k$ is the residue
field of $A^H$ and the usual theory of Cohen rings provides a
homomorphism $C(k)\to A^H$ that lifts $C(k)\to A^H/m_{A^H}$. Note
that $\oZ$ and $\oZ'$ are snc because each $A_i$ is strictly
henselian. Using averaging on the action of $H$ again, we can find
regular families of $H$-equivariant parameters $\uz=(z_1\. z_d)$
and $\uz'=(z'_1\. z'_d)$ such that $Z=V(\prod_{i=1}^lz_i)$,
$z'_i=z_i$ for $1\le i\le l$, $\oZ=V(\prod_{i=1}^nz_i)$ and
$\oZ'=V(\prod_{i=1}^{n'}z'_i)$. Explicitly, the action on $z_i$
(resp. $z'_i$) is by a character $\chi_i{\colon}H\to k(x)^\times$
(resp. $\chi'_i{\colon}H\to k(x)^\times$).

Since the image of $\uz$ is a basis of the cotangent space at $x$,
we obtain a surjective homomorphism $f{\colon}B=C(k)[[t_1\. t_d]]\to
A$ taking $t_i$ to $z_i$. Provide $B$ with the action of $H$
which is trivial on $C(k)$ and acts on $t_i$ via $\chi_i$, in
particular, $f$ is $H$-equivariant. Let us also lift each $z'_i$
to an $H$-equivariant parameter $t'_i\in B$. For $i\le l$ we
take $t'_i=t_i$, and for $i>l$ we first choose any lift and then
replace it with its $\chi_i$-weighted $H$-averaging. Consider
the regular scheme $Y=\Spec(B)$ with $H$-equivariant snc
divisors $T=V(\prod_{i=1}^lz_i)$, $\oT=V(\prod_{i=1}^nz_i)$ and
$\oT'=V(\prod_{i=1}^{n'}z'_i)$.

Since $H$ acts vary tamely on $(X,\oZ)$, it acts trivially on
$V(z_1\. z_n)=\Spec(k[[z_{n+1}\. z_d]])$ and we obtain that
$\chi_i=1$ for $i>n$. Therefore, $H$ also acts trivially on
$\Spec(B/(t_1\. t_n))=\Spec(C(k)[[t_{n+1}\. t_d]])$ and we
obtain that the action on $(Y,\oT)$ is very tame. Since the
closed immersion $j:X\to Y$ is $H$-equivariant and strict, and
$\oZ=\oT\times_YX$, Corollary \ref{8-logquotcor} implies that
$f_\bullet$ is compatible with $j$, i.e., $f_{(H,X,Z,\oZ)}$ is the
contracted strict transform of $f_{(H,Y,T,\oT)}$. The same argument
applies to the rigidifications $\oZ'$ and $\oT'$, so it now suffices
to show that $f_{(H,Y,T,\oT)}=f_{(H,Y,T,\oT')}$. For this we observe
that maximal points of log strata of the log schemes $(Y/H,\oT/H)$
and $(Y/H,\oT'/H)$ are of characteristic zero, hence the latter
equality holds by Theorem \ref{8-Morel 2.4.2}.

Finally, let us explain how one deals with the case of $m>1$. First
one finds an $H$-equivariant strict closed immersion $i:X\to Y$
such that $\oZ$ and $\oZ'$ extend to rigidifications $\oT$ and
$\oT'$ of $(Y,T)$, and the maximal points of the log strata of
$(Y/H,\oT/H)$ and $(Y/H,\oT'/H)$ are of characteristic zero. For
this we apply independently the above construction to the connected
components of $X$. Once $i$ is constructed, the same reference to
\ref{8-Morel 2.4.2} shows that $f_{(H,Y,T,\oT)}=f_{(H,Y,T,\oT')}$
and hence $f_{(H,X,Z,\oZ)}=f_{(H,X,Z,\oZ')}$ \end{proof}

\begin{corollary2}\label{8-indepprop} Let $(X,Z)$ and $G$
be as in \S\ref{8-notationendsec} and assume that the action
is admissible. Then for any choice of rigidifications $\oZ$
and $\oZ'$ we have that $f_{(G,X,Z,\oZ)}=f_{(G,X,Z,\oZ')}$.
\end{corollary2} \begin{proof} For a point $x\in X$ let
$\widehat{\calO^\sh_{X,x}}$ denote the completion of the strict
henselization of $\calO_{X,x}$. It suffices to check that for
any point $x$ the normalized blow up towers $f_{(G,X,Z,\oZ)}$
and $f_{(G,X,Z,\oZ')}$ pull back to the same normalized blow up
towers of $\hatX^\sh_x=\Spec(\widehat{\calO^\sh_{X,x}})$ (with
respect to the morphism $\hatX^\sh_x\to X$). Indeed, any normalized
blow up tower $\calX=(X_\bullet,V_\bullet)$ is uniquely determined by
its centers $V_i$. For each $i$ the morphism $Y_i=\coprod_{x\in
X}X_i\times_X\hatX^\sh_x\to X_i$ is faithfully flat, hence $V_i$
is uniquely determined by $V_i\times_{X_i}Y_i$, which is the center
of $\calX\times_X\coprod_{x\in X}\hatX^\sh_x$.

By \S\ref{8-generalizsec},
$f_{(G_\ox,\hatX^\sh_x,Z\times_X\hatX^\sh_x,\oZ\times_X\hatX^\sh_x)}$
is the contracted pullback of $f_{(G,X,Z,\oZ)}$, and an analogous
result is true for $f_{(G,X,Z,\oZ')}$. Thus the contracted
pullbacks are equal by Lemma \ref{8-formallem}. We have, however,
to worry also for the synchronization, i.e. to establish equality
of non-contracted pullbacks. For this we will use the following
trick. The towers $f_{(G,X,Z,\oZ)}$ and $f_{(G,X,Z,\oZ')}$ are
of finite length, hence there exists a finite subset $S\subset X$
such that the image of any center of either of these towers has a
non-empty intersection with $S$. Set $\hatX^\sh_S=\coprod_{s\in
S}\hatX^\sh_s$, then the pullbacks of $f_{(G,X,Z,\oZ)}$ and
$f_{(G,X,Z,\oZ')}$ to $\hatX^\sh_S$ are already contracted. Now,
in order to compare the pullbacks to $\hatX^\sh_x$, consider the
pullbacks to $\hatX^\sh_x\coprod\hatX^\sh_S$. They are contracted,
so Lemma \ref{8-formallem} (which covers disjoint unions) implies
that these pullbacks are equal. Restricting them onto $\hatX^\sh_x$
we obtain equality of non-contracted pullbacks to $\hatX^\sh_x$.
\end{proof}

\begin{remark2}\label{8-indeprem} (i) The above corollary implies
that the modification $f_{(G,X,Z,\oZ)}$ depends only on $(G,X,Z)$,
so it will be denoted $f_{(G,X,Z)}$ in the sequel. At this stage,
$f_{(G,X,Z)}$ is defined only when $X$ is admissible and $(X,Z)$
admits a rigidification.

(ii) Corollaries \ref{8-logquotcor} and \ref{8-indepprop} imply
that $f_{(G,X,Z)}$ is functorial with respect to equivariant strict
inert morphisms.  \end{remark2}

\subsubsection{Theorem \ref{8-TH. MODIFICATION} -- end of proof
for}\label{8-endproofsec} Let $X=(X,Z)$ be as assumed in
\S\ref{8-notationendsec}, and suppose that $X$ is qe. By Lemma
\ref{8-rigidlem2} there exists a surjective \'etale inert strict
morphism $h{\colon}X_0\to X$ such that $X_0$ is affine and possesses
a rigidification. Then $X_1=X_0\times_XX_0$ is affine and also
admits a rigidification (e.g. the preimage of that of $X_0$ by
one of the canonical projections). By Remark \ref{8-indeprem}(i),
$X_0$ and $X_1$ possess normalized blow up towers $f_{(G,X_0,Z_0)}$
and $f_{(G,X_1,Z_1)}$, which are  compatible with both projections
$X_1\to X_0$ by Remark \ref{8-indeprem}(ii). It follows that
$f_{(G,X_0,Z_0)}$ is induced from a unique normalized blow up
tower of $X$ that we denote as $f_{(G,X,Z)}$. This modification
satisfies all assertions of Theorem \ref{8-TH. MODIFICATION}
because $f_{(G,X_0,Z_0)}$ does so by Remark \ref{8-constrrem}(i).

\subsection{Additional properties of $f_{(G,X,Z)}$}\label{8-addsec}
Finally, let us formulate an addendum to Theorem
\ref{8-TH. MODIFICATION} where we summarize additional properties
of the constructed modification of $(X,Z)$. At this stage
we drop any assumptions on $(X,Z)$ beyond the assumptions of
\ref{8-TH. MODIFICATION}. By $f_{(G,X,Z)}$ we denote below the
entire modification from Theorem \ref{8-TH. MODIFICATION} that also
involves the modifications of Steps 1, 2, 3.

\begin{theorem2}\label{8-addendum} Keep assumptions of Theorem
\ref{8-TH. MODIFICATION}. In addition to assertions of the theorem,
the modifications $f_{(G,X,Z)}$ can be constructed uniformly for all
triples $(G,X,Z)$ such that the following properties are satisfied:

(i) Each $f_{(G,X,Z)}$ is provided with a structure of a normalized
blow up tower and its centers are contained in the preimages of
$Z\cup T$.

(ii) For any homomorphism $\lambda:H\to G$ the construction is
functorial with respect to $\lambda$-equivariant inert strict
morphisms $(Y,T)\to(X,Z)$.
\end{theorem2}
\begin{proof} The total modification $f_{(G,X,Z)}$ is obtained
by composing four modifications $f_1$, $f_2$, $f_3$ and $f_4$:
the modifications from Steps 1, 2, 3 and the modification we have
constructed in \S\ref{8-endsec}. Recall that $f_1$ and $f_4$ are constructed
as normalized blow up towers. Modifications $f_2$ and $f_3$ are
permissible blow up towers, hence they are also normalized blow up
towers with the same centers. This establishes the first part of (i).

Claim (ii) follows from the following functorialities: $f_1$
and $f_2$ are functorial with respect to all strict morphisms,
$f_3$ is functorial with respect to all inert morphisms, and $f_4$
is functorial with respect to strict inert morphisms. To prove
the second part of (i) we use (ii) to restrict $f_{(G,X,Z)}$ onto
$U=X\setminus{Z\cup T}$. Then $U$ is a regular scheme with a trivial
log structure which is acted freely by $G$. It follows from the
definitions of $f_1$, $f_2$, $f_3$ and $f_4$ that they are trivial
for such $U$. So, $f_{(G,U,\emptyset)}$ is the trivial tower, and
hence all centers of $f_{(G,X,Z)}$ are disjoint from the preimage
of $U$.
\end{proof}

\ifx\danslelivre\undefined
\bibliography{../fichiers-configuration-latex/bibliographie}
\bibliographystyle{../fichiers-configuration-latex/smfalpha}


\begin{center}
IX. Uniformisation locale première à $\ell$\\
Luc Illusie.
\end{center}

\else
\setboolean{english}{false}
\setboolean{francais}{true}
\selectlanguage{francais}
\chapter{Uniformisation locale première à $\ell$}
\addtocontents{toc}{par Luc Illusie}
\begin{center}
Luc Illusie
\end{center}
\fi

\section{Rappel de l'énoncé et premières réductions}
Rappelons l'énoncé du théorème d'uniformisation locale première à $\ell$ (\refext{II}{2-théorème
uniformisation faible}, \refext{III}{3-unif}) : \begin{theoreme}\label{9-unif} Soient $X$ un schéma
noethérien quasi-excellent, $Z$ un fermé rare de $X$ et $\ell$ un nombre premier inversible sur $X$.
Il existe une famille finie de morphismes $(p_i : X_i \ra X)_{i \in I}$, couvrante pour la topologie
des $\ell'$-altérations et telle que, pour tout $i \in I$ :

(i) $X_i$ soit régulier et connexe,

(ii) $p_i^{-1}(Z)$ soit le support d'un diviseur à croisements normaux stricts.  \end{theoreme} Le
premier ingrédient essentiel de la démonstration de \ref{9-unif} est le résultat suivant, forme
faible d'un résultat de de Jong \cite[2.4]{Families@de_Jong} : \begin{theoreme}\label{9-dJ} Soient
$f : X \ra Y$ un morphisme propre de schémas noethériens excellents intègres, et $Z$ un sous-schéma
fermé rare de $X$. Soit $\eta$ le point générique de $Y$. On suppose que $X_{\eta}$ est lisse,
irréductible et de dimension 1, et que $Z_{\eta}$ est étale. Il existe alors un groupe fini $G$, un
diagramme commutatif de $G$-schémas $$ \xymatrix{X' \ar[r]^a \ar[d]^{f'} & X \ar[d]^f \\ Y' \ar[r]^b
& Y}, $$ un diviseur effectif $G$-équivariant $D$ dans $X'$, et un fermé rare $G$-équivariant $T'$
de $Y'$ possédant les propriétés suivantes :

(i) $f'$ est projectif ;

(ii) $G$ agit trivialement sur $X$ et $Y$, librement sur $Y'-T'$ ;

(iii) $a$ et $b$ sont des altérations projectives génériquement étales, et $Y'/G \ra Y$ (resp. $X'/G
\ra X$) induit un isomorphisme en $\eta$ (resp. au point générique de $X_{\eta}$) ;

(iv) $f'$ est une courbe nodale, lisse hors de $T'$ ;

(v) $D$ est étale sur $Y'$, et contenu dans le lieu lisse de $f'$ ;

(vi) $Z' := a^{-1}(Z)$ est contenu dans $D \cup f'^{-1}(T')$.  \end{theoreme} Rappelons que dire que
$f'$ est une courbe nodale signifie que $f'$ est plat, à fibres géométriques connexes de dimension
1, ayant pour seules singularités des points quadratiques ordinaires. 

Il suffit en effet d'appliquer (\textit{loc. cit.}) au couple $(f,Z)$, avec le groupe $G$ de
(\textit{loc. cit.}) égal à $\{1\}$. Les hypothèses faites sur $X_{\eta}$ et $Z_{\eta}$ assurent que
$(f,Z)$ vérifie la condition (2.1.1) de \cite{Families@de_Jong}, et donc que le couple $(a,b)$
vérifie (2.2.1) et (2.2.5) de (\textit{loc. cit.});ce qui implique (iii). Le fermé $T'$ est donné
par le fermé noté $D$ dans (\textit{loc. cit.}, 2.5), éventuellement agrandi pour que $G$ opère
librement sur $Y'-T'$. Noter que, si $Y$ est séparé, il en est de même de $Y'$, et (ii) entraîne que
$G$ opère fidèlement sur $X'$ et $Y'$. 

\subsection{}\label{9-réductions} Les premières réductions de la démonstration de \ref{9-unif} sont
analogues à celles de la démonstration du théorème d'uniformisation locale faible. Il suffit de
prouver \ref{9-unif} pour $X$ de dimension finie. On raisonne par récurrence sur la dimension de
$X$. Le théorème est connu en dimension $\le 1$ (normalisation). Soit $d$ un entier $\ge 2$.
Supposons le théorème établi en dimension $< d$. D'après (\refext{III}{3-reduc}), on peut supposer
$X$ local noethérien complet, et même normal. D'après \refext{V}{5-algébrisation partielle mixte
raffinée géométrique}, quitte à faire une extension finie de $X$ de degré générique premier à
$\ell$, on peut supposer qu'il existe un diagramme $$ \xymatrix{X \ar[r]^g & X' \ar[d]^f \\ {} & Y},
$$ avec $Y$ local noethérien régulier complet de dimension $d-1$ et $f$ de type fini, à fibres de
dimension 1, et un point fermé $x'$ de $X'$ et un fermé rare $Z'$ de $X'$ tels que $f(x')$ soit le
point fermé de $Y$, $g$ induise un isomorphisme de $X$ sur le complété de $X'$ en $x'$, et qu'enfin
$Z = g^{-1}(Z')$. Comme $X'$ est excellent, $g$ est régulier. Le changement de base par $g$ préserve
régularité, diviseurs à croisements normaux, et familles couvrantes pour la topologie des
$\ell'$-altérations (\refext{II}{2-définition altL}).  On peut donc remplacer $X$ par $X'$, donc,
quitte à changer les notations, supposer $X$ de dimension $d$, muni d'un morphisme de type fini $f :
X \ra Y$, à fibres de dimension 1. Le problème étant local pour la topologie des
$\ell'$-altérations, donc \textit{a fortiori} pour la topologie de Zariski, on peut supposer $X$
affine. Compactifiant $f$, on se ramène à supposer $f$ propre. Quitte à éclater dans $X$ un
sous-schéma fermé ayant $Z$ pour espace sous-jacent, on peut supposer que $Z$ est un diviseur dans
$X$. Le morphisme $f$ n'est plus nécessairement une courbe relative, mais sa fibre générique reste
de dimension 1. Soit $\eta$ le point générique de $Y$. D'après (\ega{IV}{4.6.6}) il existe une
extension radicielle finie $\eta'$ de $\eta$ telle que $(X_{\eta'})_{\mathrm{red}}$ soit
géométriquement réduit, et $(Z_{\eta'})_{\mathrm{red}}$ étale sur $\eta'$. Quitte à remplacer $Y$
par son normalisé dans $\eta'$, $X$ par son normalisé dans le corps des fractions de
$(X_{\eta'})_{\mathrm{red}}$, et $Z$ par son image inverse réduite, on peut donc supposer que
$X_{\eta}$ est lisse et que $Z_{\eta}$ est étale. Le schéma $Y$ n'est plus nécessairement régulier,
mais reste affine, normal, intègre et excellent. Considérons la factorisation de Stein $\xymatrix{X
\ar[r]^{f_1} & Y_1 \ar[r]^q & Y}$ de $f$ : $q$ est fini surjectif, génériquement étale, $f_1$ est
propre et surjectif, et ses fibres géométriques sont connexes. Comme $f_1{}_*\OO_X =
\OO_{Y_1}$ et que $X$ est intègre, $Y_1$ est intègre également. Quitte à remplacer $Y$ par
$Y_1$ (et $f$ par $f_1$), on peut donc supposer que la fibre générique de $f$ est lisse,
géométriquement connexe, et que $Z_{\eta}$ est étale. On se trouve alors dans la situation de
\ref{9-dJ}, avec $Z$ un diviseur, et $Y$ affine.  \subsection{}\label{9-récurrence} Appliquons
\ref{9-dJ} à la situation que nous venons d'obtenir. D'après (iii), le morphisme $X'/G \ra X$ est
$\mathrm{alt}_{\ell'}$ couvrant. Remplaçant $X$ par $X'/G$, $Z$ par son image inverse, et $Y$ par
$Y'/G$, et changeant les notations, on peut donc supposer que $X = X'/G$, $Y = Y'/G$. Noter que,
comme $Y$ est affine et $b$, $f'$ projectifs, les actions de $G$ sur $X'$ et $Y'$ sont admissibles,
et $X'/G$ et $Y'/G$ sont encore intègres et excellents. Soit $H$ un sous-groupe de $\ell$-Sylow de
$G$. Considérons la factorisation $$ \xymatrix{X' \ar[r]^{a_1} \ar[d]^{f'} & X'/H \ar[r]^{a_2}
\ar[d] & X \ar[d]^f \\ Y' \ar[r]^{b_1} & Y'/H \ar[r]^{b_2}& Y }.  $$ Comme $a_2$ est
$\mathrm{alt}_{\ell'}$-couvrant, on peut (utilisant l'admissibilité de l'action de $H$) remplacer
$X$ par $X'/H$, $Z$ par son image inverse dans $X'/H$, $Y$ par $Y'/H$, et enfin $G$ par $H$, de
sorte qu'on peut supposer que $G$ est un $\ell$-groupe.

Appliquons l'hypothèse de récurrence au couple $(Y = Y'/G,T := T'/G)$. Il existe une famille finie
$\mathrm{alt}_{\ell'}$-couvrante$(Y_i \ra Y)_{i \in I}$, avec $Y_i$ régulier connexe, et (l'espace
sous-jacent à) $T_i = Y_i \times_Y T$ le support d'un diviseur à croisements normaux stricts. Pour
chaque $i \in I$, soient $Y'_i$ le normalisé d'une composante de $Y' \times_Y Y_i$ et $G_i \subset
G$ le groupe de décomposition de $Y_i$. Remplaçant $Y$ par $Y_i$, $Y'$ par $Y'_i$, $G$ par $G_i$, et
les autres données par leurs images inverses par $Y_i \ra Y$, $Y'_i \ra Y'$, et travaillant
séparément sur chaque $Y_i$, on se ramène à supposer que, dans le diagramme de \ref{9-dJ}, on a les
propriétés additionnelles suivantes :

(*) $Y = Y'/G$ est affine, régulier, connexe, $T = T'/G$ est un diviseur à croisements normaux
stricts dans $Y$, $Y'-T' = Y' \times _Y (Y-T)$ est un revêtement étale galoisien de $Y-T$ de groupe
$G$, $Y'$ est le normalisé de $Y$ dans $Y'-T'$, $X = X'/G$.  \section{Log régularité, fin de la
démonstration} Nous aurons besoin du résultat suivant, cas particulier d'un théorème de
Fujiwara-Kato \cite[3.1]{FK@Fujiwara-Kato}: \begin{proposition}\label{9-FK} Soient $Y$ un schéma
noethérien régulier, $T \subset  Y$ un diviseur à croisements normaux stricts, $V = Y-T$. Munissons
$Y$ de la log structure telle que le couple $(Y,T)$ soit log régulier (\refext{VI}{6-1.4}). Alors :

(i) Le foncteur de restriction de la catégorie des revêtements Kummer étales de $Y$ dans celle des
revêtements étales de $V$ modérément ramifiés le long de $T$ est une équivalence de catégories. 

(ii) Si $Y'$ est un revêtement Kummer étale de $Y$, $Y'$ est le normalisé de $Y$ dans $Y' \times_Y
V$. 

(iii) Si $V'$ est un revêtement étale de $V$, modérément ramifié le long de $T$, il existe une
unique log structure fs sur le normalisé $Y'$ de $Y$ dans $V'$ faisant de $Y'$ un revêtement Kummer
étale de $Y$.  \end{proposition} Il suffit de prouver (i) et (ii). La question est locale pour la
topologie étale sur $Y$. On peut donc supposer $Y$ strictement local, $Y = \SP A$, de point fermé $y
= \SP k$, et $T = \sum_{1 \le i \le r}\mathrm{div}(t_i)$, où les $t_i$ font partie d'un système
régulier de paramètres de $A$.  Le log schéma $Y$ admet la carte $\overline{M}_{y} = \NN^r \ra A$,
$e_i \mapsto t_i$. D'après le théorème de structure locale des revêtements Kummer étales (VI 2.2),
tout revêtement Kummer étale $Y'$ de $Y$ est somme de revêtements Kummer étale standard, de la forme
$Z = Y \times_{\SP \ZZ[\NN^r]} \SP \ZZ[Q]$, où $\NN^r \ra Q$ est un morphisme de Kummer tel que $nQ
\subset \NN^r$ pour un entier $n$ premier à la caractéristique de $k$. On a donc $Q = \NN^r \cap L$,
pour un sous-groupe $L$ de $\ZZ^r$ tel que $n[\ZZ^r : L] = 0$, et $Z = (\SP
A[x_1,\cdots,x_r]/(x_1^n-t_1,\cdots,x_r^n-t_r))^G$, où $G = \mathrm{Hom}(\ZZ^r/ L,\QQ/\ZZ)$ est un
sous-groupe de $\mu_n^r$ opérant sur $\SP A[x_1,\cdots,x_r]/(x_1^n-t_1,\cdots,x_r^n-t_r)$ de la
manière naturelle. Il en résulte que $Z$ est le normalisé de $Y$ dans le revêtement étale $Z
\times_Y V$ de $V$. On en déduit (ii), et la pleine fidélité en (i), par la considération de graphes
de morphismes entre revêtements Kummer étales de $Y$. L'essentielle surjectivité découle du lemme
d'Abhyankar (\sga{1}{xiii}{5.3}) donnant la structure du groupe fondamental modéré de $V$, qui
implique que le normalisé de $Y$ dans un revêtement étale modéré connexe de $V$ est de la forme $Z$
décrite précédemment.  \subsection{}\label{9-Kummer} Partons de la situation obtenue à la fin de
\ref{9-récurrence}. D'après \ref{9-FK}, compte tenu de (*), il existe sur $Y'$ une unique log
structure faisant de $Y'$ un revêtement Kummer étale de $Y$ (muni de la log structure définie par
$T$), galoisien de groupe $G$. En particulier, le couple $(Y',T')$ est log régulier. D'après
(\refext{VI}{6-nodal implique log lisse}), le couple $(X',f'^{-1}(T') \cup D)$ est log régulier, et
pour la log structure correspondante sur $X'$, $f'$ est log lisse. De plus, l'image inverse $Z'$ de
$Z$ dans $X'$ est un diviseur contenu dans $D' = f'^{-1}(T') \cup D$. 

L'action de $G$ sur $X'$ est modérée ($G$ est un $\ell$-groupe), mais pas nécessairement très
modérée (\refext{VI}{6-3.1}). Si elle l'était, le couple $(X = X'/G, D'/G)$ serait alors log
régulier (\refext{VI}{6-très modéré}), et l'on pourrait terminer la démonstration de \ref{9-unif}
comme dans (\refext{VII}{7-utilisation log-résolution}), à l'aide de la résolution des singularités
des couples log réguliers.  On se ramène à ce cas grâce au théorème de modification
(\refext{VIII}{8-TH. MODIFICATION}), dont nous rappelons l'énoncé :
\begin{theoreme}\label{9-modification} Soit $(X,Z)$ un couple log régulier (\refext{VI}{6-1.4}),
muni d'une action d'un groupe fini $G$. On suppose que $X$ est noethérien, séparé, et que l'action
de $G$ sur $X$ est modérée et génériquement libre. Soit $T$ le complément du plus  grand ouvert
$G$-stable de $X$ où $G$ opère librement. Il existe alors une modification projective
$G$-équivariante $f : X' \ra X$ telle que, si $Z' = f^{-1}(Z \cup T)$, le couple $(X',Z')$ soit log
régulier, et l'action de $G$ sur $X'$ très modérée.  \end{theoreme} \subsection{}\label{9-fin} Le
couple $(X',D')$ de \ref{9-Kummer} vérifie les hypothèses sur $(X,Z)$ de \ref{9-modification} :
$(X',D')$ est log régulier, $X'$ est séparé (car projectif sur $Y$), l'action de $G$ sur $X'$ est
modérée, et libre sur $X'-f'^{-1}(T')$, en particulier génériquement libre. De plus, l'action de $G$
sur $X'$ est admissible. Il existe donc un diagramme commutatif $G$-équivariant $$
\xymatrix{(X'',D'') \ar[r] \ar[d]^p & (X''/G,D''/G) \ar[d]^q \\ (X',D') \ar[r] & (X,D) =
(X'/G,D'/G)}, $$ où les flèches horizontales sont les projections canoniques, $p$ est une
modification projective ($G$-équivariante), $(X'',D'')$ est log régulier, avec $X'' - D'' \subset
p^{-1}(X'-D')$, et l'action de $G$ sur $(X'',D'')$ est très modérée. D'après (VI 3.2), le couple
$(X''/G,D''/G)$ est donc log régulier. Appliquons à ce couple le théorème de désingularisation de
Kato-Niziol (\cite[10.3, 10.4]{Toric@Kato}, \cite[5.7]{Toric@Niziol}, \cite[4.5]{HT@GR}) : il existe
un log éclatement $e : \tilde X \ra X''/G$, avec $\tilde X$ régulier, et un diviseur à croisements
normaux stricts $\tilde D$ tels que $\tilde X - \tilde D \subset e^{-1}(X''/G - D''/G)$. Alors
$\tilde Z := (qe)^{-1}(Z)$ a pour support un diviseur contenu dans $\tilde D$, donc strictement à
croisement normaux. Comme $q$ est une modification, $qe$ en est une également, et \ref{9-unif} est
démontré.

\ifx\danslelivre\undefined
\bibliography{../fichiers-configuration-latex/bibliographie}
\bibliographystyle{../fichiers-configuration-latex/smfalpha}


\begin{center}
X. Gabber's modification theorem (log smooth case)\\
Luc Illusie and Michael Temkin\footnote{The research of M.T. was partially supported by the
European Union Seventh Framework Programme
(FP7/2007-2013) under grant agreement 268182.}
\end{center}

\tableofcontents

\else
\setboolean{english}{true}
\setboolean{francais}{false}
\selectlanguage{english}
\chapter{Gabber's modification theorem (log smooth case)}
\addtocontents{toc}{par Luc Illusie et Michael Temkin}
\begin{center}
Luc Illusie and Michael Temkin\footnote{The research of M.T. was partially supported by the
European Union Seventh Framework Programme
(FP7/2007-2013) under grant agreement 268182.}
\end{center}
\fi

In this chapter we state and prove a variant of the main theorem of VIII
(see \refext{VIII}{8-TH. MODIFICATION}) for schemes $X$ which are log smooth
over a base $S$ with trivial $G$-action. See \ref{10-modif-logsmooth} for
a precise statement. The proof is given in §1 and in the remaining part
of the exposé we deduce refinements of classical theorems of de Jong,
for schemes of finite type over a field or a trait, where the degree of
the alteration is made prime to a prime $\ell$ invertible on the base.
Sections 2 and 3 are independent and contain two different proofs of
such a refinement, so let us outline the methods briefly.

For concreteness, assume that $k$ is a field, $S=\Spec(k)$, and $X$ is a
separated $S$-scheme of finite type. Two methods to construct regular
$l'$-alterations of $X$ are: (1) use a pluri-nodal fibration to construct
a regular $G$-alteration $X'\to X$ and then factor $X'$ by an $l$-Sylow
subgroup of $G$, and (2) construct a regular $l'$-alteration by induction
on $\dim(S)$ so that one factors by an $l$-Sylow subgroup at each step of
the induction. The first approach is presented in \S\ref{10-primesec}.
It is close in spirit to the approach of \cite{Families@de_Jong} and its
strengthening by Gabber-Vidal, see \cite[\S4]{Brauer@Vidal}. The weak point
of this method is that one uses inseparable Galois alterations. In particular,
even when $k$ is perfect, one cannot obtain a separable alteration of $S$.

In order to deal with the latter case, Gabber suggested to try the second
approach, which is close in spirit to the original de Jong's method of
\cite{Alterations@de_Jong}. Note that the results of \cite{Alterations@de_Jong}
do not provide equivariant alterations, which was the main reason to switch to
the method of \cite{Families@de_Jong}. However, the recent work \cite{Temkin@stable}
resolves this issue, and Gabber's idea was to use \cite{Temkin@stable} instead of
\cite{Alterations@de_Jong}. This approach is realized in \S\ref{10-complsec}
(which is due to the second author) and it outperforms the method of \S\ref{10-primesec}
when $k$ is perfect. Moreover, developing this method the second author discovered
Theorem \ref{10-l'theor} that generalizes Gabber's theorems \ref{10-de Jong1} and
\ref{10-de Jong2} to the case of a general base $S$ satisfying a certain resolvability
assumption (see \S\ref{10-lresolvesec}). In addition, if $S$ is of characteristic zero
then the same method allows to use modifications instead of $l'$-alterations, see
Theorem \ref{10-char0theor}. As an application, in Theorem \ref{10-weaksemistable}
we generalize Abramovich-Karu's weak semistable reduction theorem. Finally,
we minimize separatedness assumptions in \S\ref{10-complsec}, and for this we show
in \S\ref{10-sepelimsec} how to weaken the separatedness assumptions in Theorems
\refext{VIII}{8-TH. MODIFICATION} and \ref{10-modif-logsmooth}.

\section{The main theorem}

\begin{theorem}[Gabber]\label{10-modif-logsmooth}
Let $f : X \ra S$ be an equivariant log smooth map between fs log schemes
endowed with an action of a finite group $G$. Assume that :

(i) $G$ acts trivially on $S$ ;

(ii) $X$ and $S$ are noetherian, qe, separated, log regular, and $f$ defines a
map of log regular pairs $(X,Z) \ra (S,W)$ ;

(iii) $G$ acts tamely and generically freely on $X$.

Let $T$ be the complement of the largest open subset of $X$ over which $G$
acts freely. Then there exists an equivariant projective modification $h :
X' \ra X$ such that, if $Z' = h^{-1}(Z \cup T)$, the pair $(X',Z')$ is log
regular, the action of $G$ on $X'$ is very tame, and $(X',Z')$ is log smooth
over $(S,W)$ as well as the quotient $(X'/G,Z'/G)$ when $G$ acts admissibly
on $X$.
\end{theorem}

\begin{remark2}\label{10-remark main th.} (a) In the absence
of the hypothesis (i) it may not be possible to find a modification $h$
satisfying the properties of \ref{10-modif-logsmooth}, as the example at
the end of \refext{VIII}{8-remarks main th.} shows.

(b) By \cite[8.2]{Toric@Kato} the log smoothness of $f$ and the log regularity
of $S$ imply the log regularity of $X$. Conversely, according to Gabber
(private communication), if $X$ is log regular and $f$ is log smooth and
surjective, then $S$ is log regular.

(c) We will deduce Theorem \ref{10-modif-logsmooth} from Theorem 
\refext{VIII}{8-TH. MODIFICATION}. Recall that in the latter theorem we assumed
that $X$ is qe, though Gabber has a subtler argument that works for a general $X$. 
This forces us to assume that $S$ (and hence $X$) is qe in Theorem 
\ref{10-modif-logsmooth}. However, our argument also shows that
once one removes the quasi-excellence assumption from \refext{VIII}{8-TH. MODIFICATION},
one also obtains the analogous strengthening of Theorem \ref{10-modif-logsmooth}.
\end{remark2}

\medskip For the proof of \ref{10-modif-logsmooth} we will use the following
result on the local structure of equivariant log smooth maps.

\begin{proposition}\label{10-local-structure} (Gabber's preparation
lemma). Let $f : X \ra Y$ be an equivariant log smooth map between fine
log schemes endowed with an action of a finite group $G$. Let $x$ be a
geometric point of $X$, with image $y$ in $Y$. Assume that $G$ is the inertia
group at $x$ and is of order invertible on $Y$. Assume furthermore that $G$
acts trivially on $\overline{M}_{x}$ and $\overline{M}_{y}$\footnote{If $M$
is the sheaf of monoids of a log scheme, $\overline{M}$ denotes, as usual,
the quotient $M/\OO^*$.} and we are given an equivariant chart $a
: Y \ra \Spec \,\Lambda[Q]$ at $y$, modeled on some pairing $\chi : G^{ab}
\otimes Q^{gp} \ra \mu = \mu_N(\mathbf{C})$ (in the sense of (\refext{VI}{6-carte équivariante})),
where $Q$ is fine, $\Lambda = \mathbf{Z}[1/N, \mu]$, with $N$ the exponent of
$G$. Then, up to replacing $X$ by an étale equivariant neighborhood of $x$,
there is an equivariant chart $b : X \ra \Spec \,\Lambda[P]$ extending $a$,
such that $Q^{gp} \ra P^{gp}$ is injective, the torsion of its cokernel is
annihilated by an integer invertible on $X$, and the resulting map $b' :
X \ra X' = Y \times_{\Spec \,\Lambda[Q]}\Spec \,\Lambda[P]$ is smooth. Moreover,
up to further shrinking $X$ around $x$, $b'$ lifts to an equivariant étale
map $c : X \ra X' \times_{\Spec \,\Lambda} \Spec \, \mathrm{Sym}_{\Lambda}(V)$,
where $V$ is a finitely generated projective $\Lambda$-module equipped with a
$G$-action. If $X$, $Y$, and $Q$ are fs, with $Q$ sharp, then $P$ can be chosen
to be fs with its subgroup of units $P^*$ torsionfree.
\end{proposition}

\begin{proof}[Proof of \ref{10-local-structure}.] This is an adaptation of the proof
of \cite[3.5]{LSFI@Kato} to the equivariant case. Consider the canonical
homomorphism of \textit{loc. cit.} $$ k(x) \otimes_{\OO_{X,x}}
\Omega^1_{X/Y,x} \ra k(x) \otimes_{\mathbf{Z}} \overline{M}^{gp}_{X/Y,x}
\leqno(1) $$ sending $1 \otimes d\mathrm{log}t$ to the class of $1 \otimes
t$, where $$ \overline{M}^{gp}_{X/Y,x} = M^{gp}_{X,x}/(\OO_{X,x}^*
+ \mathrm{Im} f^{-1}(M^{gp}_{Y,y})).  $$ It is surjective, and as $G$ fixes
$x$, it is $G$-equivariant. As $G$ is of exponent invertible in $k(x)$
and acts trivially on the right hand side, (1) admits a $G$-equivariant
decomposition $$ k(x) \otimes_{\OO_{X,x}} \Omega^1_{X/Y,x} = V_0
\oplus( k(x) \otimes_{\mathbf{Z}} \overline{M}^{gp}_{X/Y,x}), \leqno(2)
$$ where $V_0$ is a finite dimensional $k(x)$-vector space, endowed with an
action of $G$. Let $(t_i)_{1 \le i \le r}$ be elements of $M^{gp}_x$ such that
the classes of $1 \otimes t_i$ form a basis of $k(x) \otimes_{\mathbf{Z}}
\overline{M}^{gp}_{X/Y,x}$. Let $Z$ be the free abelian group with basis
$(e_i)_{1 \le i \le r}$, and $h : Z \ra M^{gp}_x$ the homomorphism sending
$e_i$ to $t_i$. As $G$ acts trivially on $\overline{M}^{gp}_x$, by the
method of (\refext{VI}{6-preuve très modéré ponctuel}) we get characters
$$ \psi_i : G \ra \mu $$ such that $gh(e_i) = \psi_i(g)h(e_i)$. As in
the proof of \cite[3.5]{LSFI@Kato}, consider the homomorphism $$ u : Z
\oplus Q^{gp} \ra M^{gp}_x $$ defined by $h$ on $Z$ and the composition
$Q^{gp} \ra M^{gp}_y \ra M^{gp}_x$ on the second factor. We have $$ gu(a)
= \psi(g \otimes a))u(a) $$ for some homomorphism $$ \psi : G^{ab} \otimes
(Z \oplus Q^{gp}) \ra \mu $$ extending $\chi$ and such that $\psi(g \otimes
e_i))u(e_i) = \psi_i(g)h(e_i)$. As in \textit{loc. cit.}, if $\overline{u}$
denotes the composition $$ \overline{u} : Z \oplus Q^{gp} \ra M^{gp}_x \ra
\overline{M}^{gp}_x (= M^{gp}_x/\OO^*_x) $$ we see that $k(x) \otimes
\overline{u}$ is surjective, hence the cokernel $C$ of $\overline{u}$ is killed
by an integer $m$ invertible in $k(x)$. Using that $\OO_{X,x}^*$
is $m$-divisible, one can therefore choose elements $a_i  \in M^{gp}_x$ and
$b_i \in Z \oplus Q^{gp}$ ($1 \le i \le n$) such that the images of $a_i$
generate $\overline{M}^{gp}_x$ and $a_i^m = u(b_i)$. Let $E$ be the free
abelian group with basis $e_i$ ($1 \le i \le n$), and let $F$ be the abelian
group defined by the push-out diagram $$ \xymatrix{E \ar[r]^m \ar[d] & E \ar[d]
\\ Z \oplus Q^{gp} \ar[r]^w &F}, \leqno (3) $$ where the left vertical arrow
sends $e_i$ to $b_i$. The lower horizontal map is injective and its cokernel
is killed by $m$. The relation $a_i^m = u(b_i)$ implies that $u$ extends
to a homorphism $$ v : F \ra M^{gp}_x $$ whose composition $\overline{v} :
F \ra M^{gp}_x \ra \overline{M}^{gp}_x$ is surjective. Associated with $v$
is a morphism $$ \varphi : G^{ab} \otimes F \ra \mu $$ extending $\psi$, such
that $gv(a) = \varphi(g \otimes a)v(a)$ for $a \in F$. Let $P := v^{-1}(M_x)
\subset F$. Then $P$ is a fine monoid containing $Q$, $P^{gp} = F$, and $v$
sends $P$ to $M_x$. By \cite[2.10]{LSFI@Kato}, using the ($G$-equivariant)
isomorphism $\mathscr{H}om(\Lambda[P],M_X)_{x} \iso Hom(\Lambda[P],M_{X,x})$,
up to shrinking $X$ at $x$, the equivariant map $$ \Lambda[P] \ra M_x $$
defined by $(v,\varphi)$ extends to an equivariant chart $$ b : X \ra
\Spec \,\Lambda[P] $$ extending the chart $a : Y \ra \Spec \,\Lambda[Q]$. The
homomorphism $Q^{gp} \ra P^{gp}$ is injective, and the torsion part of its
cokernel, which is the same as the torsion part of the cokernel of $w :
Z \oplus Q^{gp} \ra F$ in (3), is killed by $m$. Consider the resulting map
$$ b' : X \ra X' = Y \times_{\Spec \,\Lambda[Q]} \Spec \, \Lambda[P].  $$ This
map is strict. Showing that the underlying schematic map is smooth at $x$ is
equivalent to showing that $b'$ is log smooth at $x$. To do this, as $X$ and
$X'$ are log smooth over $Y$, by the jacobian criterion \cite[3.12]{LSFI@Kato}
it suffices to show that the map $$ k(x) \otimes \Omega^1_{X'/Y} \ra k(x)
\otimes \Omega^1_{X/Y} $$ induced by $b'$ is injective. We have $$ k(x)
\otimes \Omega^1_{X'/Y} = k(x) \otimes P^{gp}/Q^{gp} = k(x) \otimes Z $$
(the last equality by the fact that $F/(Z \oplus Q^{gp})$ is killed by $m$),
and by construction (cf. (2)), we have $$ k(x) \otimes Z = k(x) \otimes
\overline{M}^{gp}_{X/Y,x}, $$ which by the map induced by $b'$ injects into
$k(x) \otimes \Omega^1_{X/Y}$.

Let us now prove the second assertion. For this, as $b'$ is strict, we may
forget the log structures of $X$ and $X'$, and by changing notations, we may
assume that $X' = Y$ and the log structures of $X$ and $Y$ are trivial. In
particular, we have $$ k(x) \otimes \Omega^1_{X/Y} = V_0, $$ with the notation
of (2). As the question is étale local on $X$, and closed points are very
dense in the fiber $X_y$, in particular, any point has a specialization at
a closed point of $X_y$, we may assume that $x$ is localized at a closed
point of $X_y$, and even, up to base changing $Y$ by a finite radicial
extension, that $x$ is a rational point of $X_y$. We then have $$ k(x) \otimes
\Omega^1_{X/Y} = \mathbf{m}_x/(\mathbf{m}_x^2 + \mathbf{m}_y\OO_x),
\leqno (4) $$ where $\mathbf{m}$ denotes a maximal ideal. By \cite[12 2.3,
13 14.4]{Representations@Serre} there is a finitely generated projective
$\Lambda[G]$-module $V$ such that $V_0 = k(x) \otimes V$. The homomorphism $V
\ra \mathbf{m}_x/(\mathbf{m}_x^2 + \mathbf{m}_y\OO_x)$ therefore lifts
to a homomorphism of $\Lambda[G]$-modules $$ V \ra \mathbf{m}_x, \leqno (5) $$
inducing an isomorphism $k(x) \otimes V \ra k(x) \otimes \Omega^1_{X/Y}$. By
the jacobian criterion, it follows that the ($G$-equivariant) map $$ X \ra
Y \times_{\Spec \, \Lambda} \Spec \, Sym_{\Lambda}(V) $$ is étale at $x$.

Let us prove the last assertion. First of all, as $M_x$ is fs, $v : P \ra
M_x$ factors through the saturation $P^{\mathrm{sat}}$ of $P$ in $P^{gp}$,
so we may assume that $P$ is fs. Then (cf.\cite[5.3.42]{HT@GR}) we have
a split exact sequence $$ 0 \ra H \ra P \ra P_0 \ra 0 $$ with $P_0^*$
torsionfree and $H$ a finite group. As $Q$ is fs and sharp, $Q^{gp}$ is
torsionfree, so the composition $Q^{gp} \ra P^{gp} \ra P_0^{gp}$ is still
injective, as well as the composition $H \ra P^{gp} \ra (P^{gp}/Q^{gp})$,
hence $H$ is contained in the torsion part of $(P^{gp}/Q^{gp})$, and
we have an exact sequence $$ 0 \ra H \ra (P^{gp}/Q^{gp})_{tors} \ra
(P_0^{gp}/Q^{gp})_{tors} \ra 0, $$ where the subscript \textit{tors} denotes
the torsion part. Thus $(P_0^{gp}/Q^{gp})_{tors}$ is killed by an integer
invertible on $X$. As $\overline{M}_x$ is torsionfree, the composition
$P \ra M_x \ra \overline{M}_x$ factors through $P_0$, into a map $v_0 :
P \ra \overline{M}_x$. Consider the diagram $$ \xymatrix{{} & M_x \ar[d]
\ar[r] & M^{gp}_x \ar[d] \\ P_0 \ar[r]^{v_0} & \overline{M}_x \ar[r] &
\overline{M}^{gp}_x}, $$ where the square is cartesian. As $P_0^{gp}$
is torsionfree, the map $P^{gp}_0 \ra \overline{M}^{gp}_x$ defined by the
lower row admits a lifting $s : P_0^{gp} \ra M^{gp}_x$, sending $P_0$ to
$M_x$. As $v$ is a chart, $P/v^{-1}(\OO_x^*)  \ra \overline{M}_x$
is an isomorphism, and since $H$ is contained in $v^{-1}(\OO_x^*)$,
$P_0/s^{-1}(\OO_x^*) \ra \overline{M}_x$ is an isomorphism as well,
hence $s$ is a chart at $x$. Associated with $s$ is a homomorphism $$
\varphi_0 : G^{ab} \otimes P_0^{gp} \ra \mu $$ and the map $$ b_0 : X \ra
\Spec \, \Lambda[P_0] $$ defined by the pair $(s,\varphi_0)$ is an equivariant
chart of $X$ at $x$ (extending $a$).
\end{proof}

\medskip \textit{Proof of 1.1 (beginning)}.

\medskip The strategy is to check that, at each step of the proof of the
absolute modification theorem (\refext{VIII}{8-TH. MODIFICATION}), the log
smoothness of $X/S$ is preserved, and, at the end, that of the quotient
$(X/G)/S$ as well. For some of them, this is trivial, as the modifications
performed are log blow ups. Others require a closer inspection.

\subsection{}\label{10-preliminary} \textit{Preliminary reductions.} We
may assume that conditions (1) and (2) at the beginning of  (\refext{VIII}{8-prelimsec})
are satisfied, namely :

(1)\textit{$X$ is regular,}

(2) \textit{$Z$ is a $G$-strict snc divisor in $X$.}

Indeed, these conditions are achieved by $G$-equivariant saturated log blow
up towers (\refext{VIII}{8-step1}, \refext{VIII}{8-step2}).

We will now exploit Gabber's preparation lemma 1.3 to give a local picture
of $f$ displaying both the log stratification and the inertia stratification
of $X$. We work étale locally at a geometric point $x$ in $X$ with image $s$
in $S$. Up to replacing $G$ by the inertia group $G_x$ at $x$, we may assume
that $G = G_x$.

We now apply 1.3. Let $N$ be the exponent of $G$. Assume $S$ strictly
local at $s$. We may replace $\Lambda = \mathbf{Z}[1/N,\mu]$ by its
localization at the (Zariski) image of $s$, so that $\Lambda$ is either the
cyclotomic field $\mathbf{Q}(\mu)$ or its localization at a finite place
of its ring of integers, of characteristic $ p = \mathrm{char}(k(s))$
not dividing $n$. Choose a chart $$ a : S \ra \Spec \,\Lambda[Q] $$ with
$Q$ fs and the inverse image of $\OO_{S,s}^*$ in $Q$ equal
to $\{1\}$, so that $Q$ is sharp and $Q \iso \overline{M}_s$. Let $C$
denote $k(s)$ if $\OO_{S,s}$ contains a field, and a Cohen ring of
$k(s)$ otherwise. Let $(y_i)_{1 \le i \le m}$ be a family of elements of
$\mathbf{m}_s$ such that the images of the $y_i$'s in $\OO_{S,s}/I_s$
form a regular system of parameters, where $I_s = I(s,M_s)$ is the ideal
generated by the image of $M_s - \OO_{S,s}^*$ by the canonical
map $\alpha : M_s \ra \OO_{S,s}$. By \cite[3.2]{Toric@Kato}, the
chart $a$ extends to an isomorphism
\begin{equation}\label{10-Kato-Cohen}
C[[y_1,\cdots,y_m]][[Q]])/(g) \iso \widehat{\OO}_{S,s},
\end{equation}
where $g \in C[[y_1,\cdots,y_m]][[Q]]$ is $0$ if $k(s)$ contains a field,
and congruent to $p = \mathrm{char}(k(s)) > 0$ modulo the ideal generated
by $Q - \{1\}$ and $(y_1,\cdots,y_m)$ otherwise. By 1.3, up to shrinking $X$
around $x$, we can find a $G$-equivariant commutative diagram (with trivial
action of $G$ on the bottom row)
\begin{equation}\label{10-chart-diagram}
\xymatrix{X \ar[dr] \ar[r]^{c} & X' \ar[d] \ar[r]^<<<{b} & \Spec \,(\Lambda[P]
\otimes_{\Lambda}\mathrm{Sym}_{\Lambda}(V))\ar[d] \\ {} & S \ar[r]^a &
\Spec \,\Lambda[Q]},
\end{equation} where :

(i) the square is cartesian ;

(ii) $a$, $b$, and $c$ are strict, where the log
structure on $\Spec \,\Lambda[Q]$ (resp. $\Spec \,(\Lambda[P]
\otimes_{\Lambda}\mathrm{Sym}_{\Lambda}(V))$) is the canonical one, given by
$Q$ (resp. $P$) ; $P$ is an fs monoid, with $P^*$ torsionfree ; $G$ acts on
$\Lambda[P]$ by $g(\lambda p) = \lambda\chi(g,p)p$, for some homomorphism $$
\chi : G^{ab} \otimes P^{gp} \ra \mu $$

(iii) $V$ is a free, finitely generated $\Lambda$-module, equipped with a
$G$-action ;

(iv) the right vertical arrow is the composition of the projection onto the
factor $\Spec \, \Lambda[P]$ and $\Spec \, \Lambda[h]$, for a homomorphisme
$h : Q \ra P$ such that $h^{gp}$ is injective and the torsion part of
$\mathrm{Coker}h^{gp}$ is annihilated by an integer invertible on $X$ ;

(v) $c$ is étale.

\medskip (vi) Consider the map $$ v : P \ra M_x $$ defined by the chart
$X \ra \Spec \, \Lambda[P]$ induced by $bc$. Up to localizing on $X'$
around $x$, we may assume that $v$ factors through the localization
$P_{(\frak p)}$ of $P$ at the prime ideal $\frak p$ complementary of the
face $v^{-1}(\OO_{X,x}^*)$. Replacing $P$ by $P_{(\frak p)}$,
$P$ decomposes into
\begin{equation}\label{10-P1} P = P^* \oplus P_1,
\end{equation} with $P^* = v^{-1}(\OO_{X,x}^*)$ free finitely generated
over $\mathbf{Z}$, and $P_1$ sharp, and the image of $x$ by $bc$ into the
factor $\Spec \, \Lambda[P]$ is the rational point at the origin. Then $v$
induces an isomorphism $P_1 \iso \overline{M}_x$. By the assumptions (1),
(2), we have $\overline{M}_x \iso \mathbf{N}^r$. One can therefore choose
$(e_i \in P_1)$ $(1 \le i \le r)$ forming a basis of $P_1$.  Then $v(e_i) =
t_i \in M_x \subset \OO_{X,x}$ is a local equation for a branch $Z_i$
of $Z$ at $x$, $(Z_i)_{1 \le i \le r}$ is the set of branches of $Z$ at $x$,
and $G$ acts on $t_i$ through the character $\chi_i = \chi | \mathbf{Z}e_i :
G \ra \mu$.

\medskip Furthermore :

\medskip (vii) The square in (\ref{10-chart-diagram}) is tor-independent.

\medskip Indeed, by the log regularity of $S$ and the
choice of the chart $a$, we have, by \cite[6.1]{Toric@Kato},
$Tor_i^{\mathbf{Z[Q]}}(\OO_{S,s},\mathbf{Z}[P]) = 0$ for $i > 0$.

\medskip Though this will not be needed, one can describe the
local structure of (\ref{10-chart-diagram}) more precisely as
follows. Let
\begin{equation}\label{10-Y} Y := \Spec \,(\Lambda[P]
\otimes_{\Lambda}\mathrm{Sym}_{\Lambda}(V)) = \Spec \,(\Lambda[P^*]
\otimes_{\Lambda} \Lambda[P_1]  \otimes_{\Lambda}\mathrm{Sym}_{\Lambda}(V))
\end{equation} and let $Y' := \Spec \,C[[y_1,\cdots, y_m]][[Q]] \times_{\Spec
\Lambda} Y$, with the notation of \ref{10-Kato-Cohen}. Replacing $X'$
by its strict localization at $x$ we may assume that $X = X'$. Then the
completion of $X$ at $x$ is either isomorphic to the completion of $Y'$ at
$x$, or a regular divisor in it, defined by the equation $g' = 0$, where $g'$
is the image of  $g$ in $\widehat{\OO}_{Y',x}$, with the notation
of \ref{10-Kato-Cohen}.

\subsection{}\label{10-Step 3} \textit{Step 3 and log smoothness
(beginning)}. We will now analyze the modifications performed in the
proof of Step 3 in \refext{VIII}{8-step3}, \refext{VIII}{8-step3sec}. The permissible towers used in
\textit{loc. cit.} are iterations of operations of the form : for a subgroup
$H$ of $G$, blow up the fixed point (regular) subscheme $X^H$, and replace $Z$
by the union of its strict transform $Z^{\mathrm{st}}$ and the exceptional
divisor $E$. Though such a blow up is not a log blow up in general, we will
see that it still preserves the log smoothness of $X$ over $S$.

We work étale locally around $x$, so we can assume $X =
X'$ in \ref{10-chart-diagram}. We then have a cartesian square
\begin{equation}\label{10-XH} \xymatrix{X^H \ar[r]^{b^H} \ar[d]^{f} &
Y^H \ar[d] \\ X \ar[r]^b & Y},
\end{equation}
with $Y$ as in (\ref{10-Y}). We
also have cartesian squares
\begin{equation}\label{10-Z} \xymatrix{Z \ar[r]^b
\ar[d]^f & T \ar[d] \\ X \ar[r]^{b} & Y},
\end{equation}
where $T \subset Y$
is the snc divisor $\sum T_i$, $T_i $ defined by the equation $e_i \in P_1$
(\ref{10-P1}), and
\begin{equation}\label{10-Z inter XH} \xymatrix{Z
\times_X X^H \ar[r] \ar[d] & T \times_Y Y^H \ar[d] \\ X \ar[r] & Y}.
\end{equation}
\begin{lemma}\label{10-tor-ind} The squares (\ref{10-XH}),
(\ref{10-Z}), and (\ref{10-Z inter XH}) are tor-independent.
\end{lemma}
\textit{Proof}. For (\ref{10-Z}), this is because $Z$ (resp. $T$) is a
divisor with normal crossings in $X$ (resp. $Y$) (cf. \sga{6}{vii}{1.2}). For
(\ref{10-XH}), as the square (\ref{10-chart-diagram}) is tor-independent
(by \ref{10-preliminary} (vii)), it is enough to show that the composite
(cartesian) square
\begin{equation}\label{10-composite} \xymatrix{X^H
\ar[d] \ar[r] & Y^H \ar[d] \\ S \ar[r] &\Spec \, \Lambda[Q]}
\end{equation}
is tor-independent. We have a decomposition
\begin{equation}\label{10-YH}
Y^H = (\Spec \,\Lambda[P^*]) ^H \times (\Spec \, \Lambda[P_1])^H \times (\Spec \,
\mathrm{Sym}_{\Lambda}(V))^H,
\end{equation}
(products taken over $\Spec \,
\Lambda$), and the map to $\Spec \, \Lambda[Q]$ is the composition of the
projection onto $(\Spec \,\Lambda[P^*]) ^H \times (\Spec \, \Lambda[P_1])^H$
and the canonical map induced by $\Spec \, \Lambda[Q] \ra \Spec \, \Lambda[P]$,
which factors through the fixed points of $H$, $G$ acting trivially on the
base. Let us examine the three factors.

(a) We have $$ (\Spec \,\mathrm{Sym}_{\Lambda}(V))^H = \Spec \,
\mathrm{Sym}_{\Lambda}(V_H), $$ where $V_H$ is the module of coinvariants,
a free module of finite type over $\Lambda$, as $H$ is of order invertible
in $\Lambda$. Therefore $\Spec \, \Lambda[Q] \times_{\Spec \,\Lambda} (\Spec
\,\mathrm{Sym}_{\Lambda}(V))^H$ is flat over $\Spec \, \Lambda[Q]$, and its
enough to check that $(\Spec \,\Lambda[P^*]) ^H \times (\Spec \, \Lambda[P_1])^H$
is tor-independent of $S$ over $\Spec \, \Lambda[Q]$.

(b) The restriction to $P^* = v^{-1}(\OO_{X,x}^*)$ of the
1-cocycle $z(v) \in Z^1(H,Hom(P,k(x)^*)$ associated with $v : P \ra M_x$
($hv(a) = z(v)(h,a)v(a)$ for $h \in H$, $a \in P$, see the proof of
\ref{10-local-structure} and (\refext{VI}{6-preuve très modéré ponctuel})) is a 1-coboundary, hence trivial,
as $B^1(H, Hom(P,k(x)^*) = 0$. Therefore $$
 (\Spec \,\Lambda[P^*]) ^H =  \Spec \,\Lambda[P^*].
$$

(c) Recall that $$ P_1 = \oplus_{1 \le i \le r} \mathbf{N}e_i, $$ with $e_i$
sent by $v$ to a local equation of the branch $Z_i$ of $Z$, and that $G$ acts
on $\Lambda[\mathbf{N}e_i]$ through the character $\chi_i : G \ra \mu$. Let
$A \subset \{1,\cdots,r\}$ be the set of indices $i$ such that $\chi_i | H$
is trivial. Then $$ (\Spec \, \Lambda[P_1])^H = \Spec \, \Lambda[\oplus_{i \in
A} \mathbf{N}e_i].  $$

Let $I$ be the ideal of $P$ generated by $\{e_i \}_{i \notin A}$. It follows
from (b) and (c) that $$ (\Spec \, \Lambda[P])^H = \Spec \, \Lambda[P]/(I),
$$ where $(I)$ is the ideal of $\Lambda[P]$ generated by $I$. By
\cite[6.1]{Toric@Kato}, $Tor_i^{\Lambda[Q]}(\OO_S,\Lambda[P]/(I)) =
0$ for $i > 0$, and therefore (\ref{10-composite}), hence (\ref{10-XH})
is tor-independent. It remains to show the tor-independence of (\ref{10-Z inter XH}).
For this, again it is enough to show the tor-independence
of
\begin{equation}\label{10-new composite} \xymatrix{Z \times_X X^H
\ar[r] \ar[d] & T \times_Y Y^H \ar[d] \\ S \ar[r] & \Spec \, \Lambda[Q]}.
\end{equation} By (a), (b), (c), we have $$ T \times_Y Y^H = \sum_{i \in A}
\Spec \, \Lambda[P]/(J_i) \times \Spec \, \mathrm{Sym}_{\Lambda}(V_H), $$ where
$J_i \subset P$ is the ideal generated by $e_i \in P_1$, and $(J_i)$ the ideal
generated by $J_i$ in $\Lambda[P]$.  The desired tor-independence follows
from the vanishing of $Tor_i^{\Lambda[Q]}(\OO_S,\Lambda[P]/(J_B))$,
where for a subset $B$ of $A$, $J_B$ denotes the ideal generated by the
$e_i$'s for $i \in B$.

\begin{lemma}\label{10-pull-back blow-up} Consider a cartesian square
\begin{equation}\label{10-pull-back 1} \xymatrix{V' \ar[d] \ar[r] & V \ar[d]
\\ X' \ar[r]^g & X},
\end{equation}
where the right vertical arrow is a
regular immersion. If (\ref{10-pull-back 1}) is tor-independent, then the
left vertical arrow is a regular immersion, and $$ \mathrm{Bl}_{V'}(X') = X'
\times_X \mathrm{Bl}_V(X).  $$ Let $W \ra X$ be a second regular immersion,
such that $V \times_X W \ra W$ is a regular immersion, and let $W' = X'
\times_X W$. If moreover the squares
\begin{equation}\label{10-pull-back 2}
\xymatrix{V' \ar[d] \ar[r] & V \ar[d] \\ X' \ar[r]^g & X},
\end{equation}
and
\begin{equation}\label{10-pull-back 3} \xymatrix{V' \times_{X'} W'
\ar[d] \ar[r] & V \times_X W \ar[d] \\ X' \ar[r]^g & X},
\end{equation}
are
tor-independent, then the left vertical arrows are regular immersions, and $$
W'^{\mathrm{st}} = X' \times_X W^{\mathrm{st}}, $$ where $W^{\mathrm{st}}$
(resp. $W'^{\mathrm{st}}$) is the strict transform of $W$ (resp. $W'$)
in $\mathrm{Bl}_V(X)$ (resp. $\mathrm{Bl}_{V'}(X')$).
\end{lemma}
\textit{Proof.} Let $I$ (resp. $I'$) be the ideal of $V$ (resp. $V'$) in $X$
(resp. $X'$). By the tor-independence of (\ref{10-pull-back 1}), if $u :
E \ra I$ is a local surjective regular homomorphism \sga{6}{vii}{1.4}, the
Koszul complex $g^*K(u)$ is a resolution of $\OO_{V'}$, hence $V'
\ra X'$ is a regular immersion. Moreover, by \sga{6}{vii}{1.2}, for any $n
\ge 0 $, the natural map $g^*I^n \ra I'^n$ is an isomorphism,  and therefore
$\mathrm{Bl}_{V'}(X') = X' \times_X \mathrm{Bl}_V(X)$. The tor-independence of
(\ref{10-pull-back 2}) and (\ref{10-pull-back 3}) imply that of $$ \xymatrix{
V' \times_{X'} W' \ar[r] \ar[d] & V \times_X W \ar[d] \\ W' \ar[r] & W}.
$$ The second assertion then follows from the first one and the formulas
(\refext{VIII}{8-strictrem} (ii)) $$ W^{\mathrm{st}} = \mathrm{Bl}_{V \times_X
W} W, $$ $$ W'^{\mathrm{st}} = \mathrm{Bl}_{V' \times_{X'} W'} W'.  $$

\subsection{}\label{10-Step 3 (end)} \textit{Step 3 and log smoothness
(end)}. As recalled at the beginning of \ref{10-Step 3}, we have to show
that, if $H$ is a subgroup of $G$, then the log regular pair $(X_1,Z_1)$
is log smooth over $S$, where $X_1 := \mathrm{Bl}_{X^H}(X)$ and $Z_1$ is
the snc divisor $Z^{\mathrm{st}} \cup E$, $Z^{\mathrm{st}}$ (resp. $E$)
denoting the strict transform of $Z$ (resp. the exceptional divisor) in the
blow-up $h : X_1 \ra X$.

The question is again étale local above $X$ around $x$, so we may assume
that $X = X'$ and we look at the cartesian square (\ref{10-chart-diagram}) $$
\xymatrix{X \ar[r] \ar[d] & Y \ar[d] \\ S \ar[r] & \Spec \, \Lambda[Q]}, $$ with
$Y$ as in (\ref{10-Y}), and the associated cartesian squares (\ref{10-XH}),
(\ref{10-Z}), and (\ref{10-Z inter XH}).

\textit{Claim.} We have
\begin{equation}\label{10-blow-up}
\mathrm{Bl}_{X^H}(X) = X \times_Y \mathrm{Bl}_{Y^H}(Y),
\end{equation}
\begin{equation}\label{10-strict transform} Z^{\mathrm{st}} = X \times_Y
T^{\mathrm{st}}.
\end{equation}

\textit{Proof.} In view of \ref{10-tor-ind} and \ref{10-pull-back blow-up},
(\ref{10-blow-up}) follows from the fact that the immersion $Y^H \ra
Y$ is regular. For (\ref{10-strict transform}) recall that $$ T = T_0
\times_{\Spec \, \Lambda} \Spec \, \mathrm{Sym}_{\Lambda}(V), $$ where $T_0
\subset \Spec \, \Lambda[P]$ is the snc divisor $$ T_0 = \sum_{1 \le i \le r}
\mathrm{div}(z_i) $$ with $z_i \in \Lambda[P]$ the image of $e_i \in P_1$
as in \ref{10-P1}. Hence
\begin{equation}\label{10-T} T = \sum_{1 \le
i \le r} T_i,
\end{equation}
where $T_i = \mathrm{div}(z_i) \times_{\Spec \,
\Lambda}\mathrm{Sym}_{\Lambda}(V)$, and $T^{\mathrm{st}} = \sum_{1 \le i \le
r}T_i^{\mathrm{st}}$. We have (\ref{10-YH}) $$ Y^H =  (\Spec \, \Lambda[P])^H
\times \Spec \, \mathrm{Sym}_{\Lambda}(V_H), $$ with $(\Spec \, \Lambda[P])^H$
defined by the equations $(z_i = 0)_{i \notin A}$, with the notations of
\ref{10-tor-ind} (c). In particular, the immersion $Y^H \times_Y Z_i \ra Z_i$
is regular, hence, by \ref{10-pull-back blow-up}, we have $Z_i^{\mathrm{st}}
= X \times_Y T_i^{\mathrm{st}}$, hence (\ref{10-strict transform}), which
finishes the proof of the claim.

Since the map $S \ra \Spec \, \Lambda[Q]$ is strict, in order to prove the
desired log smoothness, we may, by this claim, replace the triple $(X, X^H,
Z)$ over $S$ by $(Y,Y^H,T)$ over $\Spec \, \Lambda[Q]$. We choose coordinates on
$P^*$, $P_1 = \mathbf{N}^r$, $V$ : $$ P^* = \oplus_{1 \le i \le t}\mathbf{Z}f_i
, \, \, P_1 = \oplus_{1 \le i \le r}\mathbf{N}e_i, \, \, V = \oplus_{1 \le i
\le s} \Lambda y_i $$ $$ \Lambda[P] = \Lambda[u_1^{\pm 1}, \cdots, u_t^{\pm 1},
z_1, \cdots, z_r], \, \, \mathrm{Sym}_{\Lambda}(V) = \Lambda[y_1,\cdots, y_s],
$$ with $u_i$ (resp. $z_i$) the image of $f_i$ (resp. $e_i$) in $\Lambda[P]$,
in such a way that $$ \Lambda[P]^H = \Lambda[u_1^{\pm 1}, \cdots, u_t^{\pm 1},
z_{m+1}, \cdots, z_r], $$ i. e. is defined in $\Lambda[P]$ by the equations
$(z_1 = \cdots = z_m = 0$, for some $m$,  $1 \le m \le r$, and $$ \Lambda[V_H]
= \Lambda[y_{n+1}, \cdots, y_s], $$ i. e. is defined in $\Lambda[V]$ by the
equations $y_1 = \cdots = y_n = 0$ for some $n$, $1 \le n \le s$. Then $$
Y^H \subset Y = \Spec \, \Lambda[u_1^{\pm 1}, \cdots, u_t^{\pm 1},z_1, \cdots,
z_r, y_1, \cdots, y_s] $$ is defined by the equations $$ z_1 = \cdots =
z_m = y_1 = \cdots = y_n = 0.  $$ Then$$ Y' := \mathrm{Bl}_{Y^H}(Y) $$ is
covered by affine open pieces : $$ U_i = \Spec \, \Lambda[(u_j^{\pm 1})_{1
\le j  \le t}, z'_1,\cdots, z'_{i-1}, z_i, z'_{i+1}, \cdots, z'_m, z_{m+1},
\cdots, z_r, y'_1, \cdots, y'_n, y_{n+1}, \cdots, y_s] $$ ($1 \le i \le m$),
with $U_i \ra Y$ given by $z_j \ra z_iz'_j$ for $1 \le j \le m$, $j \ne i$,
$y_j \ra z_iy'_j$, $1 \le j \le n$, and the other coordinates unchanged,
and $$ V_i = \Spec \, \Lambda[(u_j^{\pm 1})_{1 \le j  \le t}, z'_1, \cdots,
z'_m,  z_{m+1}, \cdots, z_r, y'_1, \cdots, y'_{i-1},y_i,y'_{i+1}, \cdots,
y'_n, y_{n+1}, \cdots, y_s] $$ ($1 \le i \le n$), with $V_i \ra Y$ given by
$z_j \mapsto y_iz'_j$ for $1 \le j \le m$, $y_j \mapsto y_iy'_j$, $1 \le j \le
n$, $j \ne i$, and the other coordinates unchanged.  Recall that $Y$ has the
log structure defined by the log regular pair $(Y,T)$, where $T$ is the snc
divisor $$ T = (z_1 \cdots z_r = 0), $$ and $Y'$ is given the log structure
defined by the log regular pair $(Y',T')$, where $T'$ is the snc divisor $$
T' = F \cup T^{\mathrm{st}}, $$ where $F$ is the exceptional divisor of the
blow up of $Y^H$ and $T^{\mathrm{st}}$ the strict transform of $T$. Consider
the canonical morphisms $$ \xymatrix{Y' \ar[r]^b & Y \ar[r]^<<<{g} & \Sigma :=
\Spec \, \Lambda[Q]}.  $$ They are both morphisms of log schemes. The morphism
$g$ is given by the homomorphism of monoids $\gamma : Q \ra P$, i. e.  $$
q \in Q \mapsto (\gamma_1(q), \cdots, \gamma_t(q), \gamma_{t+1}(q), \cdots,
\gamma_{t+r}(q), 0, \cdots, 0) \in \Lambda[u_1^{\pm 1}, \cdots, u_t^{\pm 1},
z_1, \cdots, z_r, y_1, \cdots, y_s].  $$ The blow up $b$ has been described
above in the various charts. Note that $b$ is not log étale, or even log
smooth, in general. However, the composition $gb : Y' \ra \Sigma$ is log
smooth. We'll check this on the charts $(U_i)$, $(V_i)$.

(a) \textit{Chart of type $U_i$}. We have $F = (z_i = 0)$, $T^{\mathrm{st}}
= (\prod_{1 \le j \le r, j \ne i}z'_i = 0)$. Hence the log strucure of
$U_i$ is given by the canonical log structure of $\Lambda[\mathbf{N}^r]$
in the decomposition $$ U_i = \Spec \, \Lambda[\mathbf{Z}^t] \times \Spec \,
\Lambda[\mathbf{N}^r] \times \Spec \, \Lambda[y'_1, \cdots, y'_n, y_{n+1},
\cdots, y_s] $$ with the basis element $e_k$ of $\mathbf{N}^r$ sent to the
$k$-th place in $(z'_1,\cdots, z'_{i-1}, z_i, z'_{i+1}, \cdots, z'_m, z_{m+1},
\cdots, z_r)$ (and the basis element $f_k$ of $\mathbf{Z}^t$ sent to $u_k$),
the third factor having the trivial log structure. Checking the log smoothness
of $gb : U_i \ra \Sigma$ amounts to checking the log smoothness of its factor
$\Spec \, \Lambda[P] \ra \Sigma = \Spec \, \Lambda(Q)$, which is defined by the
composition of homomorphisms of monoids $$ \xymatrix{Q \ar[r]^{\gamma} &
\mathbf{Z}^t \oplus  \mathbf{N}^r \ar[r]^{Id \oplus \beta} & \mathbf{Z}^t
\oplus \mathbf{N}^r}, $$ where $\beta$ is the homomorphism $\mathbf{N}^r
\ra \mathbf{N}^r$ sending $e_j$ to $e_j + e_i$ for $1 \le j \le m, j
\ne i$, $e_i$ to $e_i$, and $e_j$ to $e_j$ for $m+1 \le j \le r$. Recall
((\ref{10-chart-diagram}), (iv)) that $\gamma^{gp}$ is injective and the
torsion part of its cokernel is invertible in $\Lambda$. As $\beta^{gp}$ is
an isomorphism, the same holds for the composition $(Id \oplus \beta)\gamma$,
hence $gb : U_i \ra \Sigma$ is log smooth.

\medskip (b) \textit{Chart of type $V_i$}. We have $F = (y_i = 0)$,
$T^{\mathrm{st}} = \prod_{1 \le j \le m}z'_j\prod_{j \ge m+1}z_i$. Hence
the log structure of $V_i$ is given by the canonical log structure
of $\Lambda[\mathbf{N}^{r+1}]$ in the decomposition $$ V_i = \Spec \,
\Lambda[\mathbf{Z}^t] \times \Spec \, \Lambda[\mathbf{N}^{r+1}] \times \Spec
\, \Lambda[(y'_j)_{1 \le j \le n, j \ne i}, y_{n+1}, \cdots, y_s] $$
with the basis element $e_k$ of $\mathbf{N}^{r+1}$ sent to the $k$-th
place in $(z'_1,\cdots, z'_m, z_{m+1}, \cdots, z_r)$ if $k \le r$, and
$e_{r+1}$ sent to $y_i$ (and the basis element $f_k$ of $\mathbf{Z}^t$
sent to $u_i$), the third factor having the trivial log structure. Again,
Checking the log smoothness of $gb : V_i \ra \Sigma$ amounts to checking
the log smoothness of its factor $\Spec \, \Lambda[\mathbf{Z}^t] \times \Spec \,
\Lambda[\mathbf{N}^{r+1}] \ra \Spec \, \Lambda(Q)$. This factor is defined by
the composition of homomorphisms of monoids $$ \xymatrix{Q \ar[r]^{\gamma}
& \mathbf{Z}^t \oplus \mathbf{N}^r \ar[r]^{Id \oplus \beta} & \mathbf{Z}^t
\oplus \mathbf{N}^{r+1}} $$ where $\beta : \mathbf{N}^r \ra \mathbf{N}^{r+1}$
sends $e_j$ to $e_j + e_{r+1}$ for $1 \le j \le m$, and to $e_j$ for  $m+1
\le j \le r$. Then $\beta^{gp}$ is injective, and its cokernel is isomorphic
to $\mathbf{Z}$, hence $(\beta \gamma)^{gp}$ is injective, and we have an
exact sequence $$ 0 \ra \mathrm{Coker}\,  \gamma^{gp} \ra \mathrm{Coker} \,
(\beta \gamma)^{gp} \ra \mathbf{Z} \ra 0.  $$ In particular, the torsion
part of $\mathrm{Coker} (\beta \gamma)^{gp}$ is isomorphic to that of
$\mathrm{Coker} \, \gamma^{gp}$, hence of order invertible in $\Lambda$,
which implies that $gb : V_i \ra \Sigma$ is log smooth.

This finishes the proof that Step 3 preserves log smoothness.

\subsection{}\label{10-end of proof} \textit{End of proof of
\ref{10-modif-logsmooth}.} We may now assume that in addition to conditions
(1) and (2) of \ref{10-preliminary}, condition (3) is satisfied as well, namely

(3) \textit{$G$ acts freely on $X - Z$ (i. e. $Z = Z \cup T$ in the notation
of \ref{10-modif-logsmooth} or (\refext{VIII}{8-TH. MODIFICATION})),  and,
for any geometric point $x \ra X$, the inertia group $G_x$ is abelian.}

We have to check :

\medskip \textit{Claim.} \textit{If $f_{(G,X,Z)} : (X',Z') \ra (X,Z)$ is the
modification of (\refext{VIII}{8-logquotcor}),
then $(X',Z')$ and $(X'/G,Z'/G)$ are log smooth over $S$.}

\medskip Working étale locally around a geometric point $x$ of $X$, we will
first choose a strict rigidification $(X,\overline{Z})$ of $(X,Z)$ such that
$(X,\overline{Z})$ is log smooth over $S$.  We will define $(X,\overline{Z})$
as the pull-back by $S \ra\Sigma = \Spec \, \Lambda[Q]$ of a rigidification
$(Y,\overline{T})$ of $(Y,T)$ which is log smooth over $\Sigma$, with the
notation of $(\ref{10-Y})$. Using that $G$  $(= G_x)$ is abelian, one can
decompose $V$ into a sum of $G$-stable lines, according to the characters
of $G$ : $$ V = \oplus_{1 \le i \le s} \Lambda y_i $$ with $G$ acting on
$\Lambda y_i$ through a character $\chi_i : G \ra \mu_N$, i. e. $gy_i  =
\chi_i(g)y_i$. We define $\overline{T}$ to be the divisor $z_1 \cdots z_r y_1
\cdots y_s = 0$ in $Y = \Spec \, \Lambda[u_1^{\pm 1}, \cdots, u_t^{\pm 1},z_1,
\cdots, z_r, y_1, \cdots, y_s]$. The action of $G$ on $(Y,\overline{T})$ is
very tame at $x$ because the log stratum at $x$ is $\Spec \, \Lambda[u_1^{\pm
1}, \cdots, u_t^{\pm 1}]$, hence very tame in a neighborhood of $x$ by
(\refext{VIII}{8-tamelocallem}) (actually on the whole of $Y$,
cf. (\refext{VIII}{6-dcn stricts}, \refext{VIII}{6-exemples}(a)).
On the other hand, $(\Spec \, \Lambda[y_1,\cdots, y_s], y_1 \cdots y_s
= 0)$ is log smooth over $\Spec \,\Lambda$, and as $\Spec \, \Lambda[P]$ is log
smooth over $\Sigma$, $(Y,\overline{T}) $ is log smooth over $\Sigma$. Since
$f_{(G,X,Z)}$ is compatible with base change by strict inert morphisms,
it is enough to check that if $f_{(G,Y,T)} = f_{(G,Y,T,\overline{T})}
: (Y',T') \ra (Y,T)$ is the modification of (\refext{VIII}{8-logquotcor})
then $(Y',T')$ is log smooth over $\Sigma$. Recall
(\refext{VIII}{8-constrsec}) that we have a cartesian $G$-equivariant diagram
\begin{equation}\label{10-cartesian diagram} \xymatrix{(Y',\overline{T'})
\ar[r]^{\overline{h'}} \ar[d]^{\alpha'} & (Y,\overline{T}) \ar[d]^{\alpha}
\\ (Y'_1,\overline{T'_1}) \ar[r]^{\overline{h_1}} & (Y_1,\overline{T_1})},
\end{equation} where the horizontal maps are the compositions of saturated
log blow up towers, and the vertical ones Kummer étale $G$-covers. From
(\ref{10-cartesian diagram}) is extracted the relevant diagram involving $h
:= f_{(G,Y,T,\overline{T})}$, $$ \xymatrix{(Y',T') \ar[d]^{\beta} \ar[r]^h &
(Y,T) \\ (Y'_1,T'_1)}, $$ where $T'_1 = \overline{h}^{-1}_1(T_1)$, with $T_1 =
T/G$, $T' = \alpha'^{-1}(T'_1)$, and $h$ (resp. $\beta$) is the restriction
of $\overline{h'}$ (resp. $\alpha'$) over $(Y,T)$ (resp. $(Y_1',T_1')$). In
particular, $\beta$ is a Kummer étale $G$-cover (as Kummer étale
$G$-covers are stable under any fs base change). As $G$ acts trivially
on $S$, this diagram can be uniquely completed into a commutative diagram
\begin{equation}\label{10-relevant diagram} \xymatrix{(Y',T') \ar[d]^{\beta}
\ar[r]^h & (Y,T) \ar[d]^f \\ (Y'_1,T'_1) \ar[r]^g & \Sigma}.
\end{equation}
Here $f$ is log smooth and $\beta$ is a Kummer étale $G$-cover. Though
$\overline{h'}$ and $\overline{h_1}$ are log smooth, $h$ and $h_1$ are not,
in general. However, it turns out that :

\medskip (*) \textit{$g : (Y'_1,T'_1) \ra \Sigma$, hence $g\beta = fh :
(Y',T') \ra \Sigma$, are log smooth},

\medskip \noindent which will finish the proof of the claim, hence of
\ref{10-modif-logsmooth}. We first prove

\medskip (**) \textit{With the notation of (\ref{10-relevant diagram}),
$(Y_1,\overline{T_1})$ is log smooth over $\Sigma$}.

\medskip Let us write $Y = \Spec \, \Lambda[\overline{P}]$, with
\begin{equation}\label{10-Peq} \overline{P} = P \times \mathbf{N}^s =
\mathbf{Z}^t \times \mathbf{N}^r \times \mathbf{N}^s.
\end{equation} As
$G$ acts very tamely on $(Y,\overline{T})$, the quotient pair $(Y_1 = Y/G,
\overline{T_1} = \overline{T}/G)$ is log regular. More precisely, by the
calculation in (\refext{VI}{6-très modéré ponctuel}(b)), this pair consists of the log scheme $Y_1 =
\Spec \, \Lambda[\overline{R}]$ with its canonical log structure, where $$
\overline{R} = \mathrm{Ker}(\overline{P}^{gp} \ra \mathrm{Hom}(G,\mu_N))
\cap \overline{P}, $$ $\overline{P}^{gp} \ra \mathrm{Hom}(G,\mu_N)$ being
the homomorphism defined by the pairing $\chi : G \otimes \overline{P}^{gp}
\ra \mu_N$. The inclusion $\overline{R} \subset \overline{P}$ is a Kummer
morphism, and $\overline{P}^{gp}/\overline{R}^{gp}$ is annihilated by
an integer invertible in $\Lambda$. As $Q^{gp} \ra P^{gp}$ is injective,
with the torsion part of its cokernel annihilated by an integer invertible
in $\Lambda$, the same is true of $Q^{gp} \ra \overline{P}^{gp}$, hence
of $Q^{gp} \ra \overline{R}^{gp}$. Thus $(Y_1,\overline{T_1}) = \Spec \,
\Lambda[\overline{R}]$ is log smooth over $\Sigma$.

Finally, let us prove (*). It is enough to work locally on $Y'_1$ so we
can replace the log blow up sequence $(Y'_1,\oT'_1)\to(Y_1,\oT_1)$ with
an affine chart (i.e. we replace the first log blow up with a chart, then
do the same for the second one, etc.). Then $Y'_1=\Spec\Lambda[\oR']$,
and $\oR^\gp\iso\oR'^\gp$ by \refext{VIII}{8-blowmonosec}. Note that
$\oR'\iso\ZZ^a\times\NN^b$ where $D_1,\dots,D_b$ are the components of $\oT'_1$.
We can assume that $D_1,\dots,D_c\subseteq T'_1$ and $D_{c+1},\dots,D_b$
are not contained in $T'_1$. Let $R'\iso\ZZ^a\times\NN^c$ denote the
submonoid $\oR'$ that defines the log structure of $(Y'_1,T'_1)$. Note that
$R'$ consists of all elements $g'\in\oR'$ such that $(g'=0)\subseteq T'_1$
(as a set). Also, by $\nu:\oR\to\oR'$ we will denote the homomorphism that
defines $(Y'_1,\oT'_1)\to(Y_1,\oT_1)$.

We showed in \refext{VIII}{8-constrsec} that $T_1=T/G$ is a $\QQ$-Cartier divisor in $Y_1$
and observed that therefore $T'_1$ is a Cartier divisor in $Y'_1$. Note that
the inclusion $R \subset \overline{R}$, where $$ R = \mathrm{Ker}(P^{gp}
\ra \mathrm{Hom}(G,\mu_N)) \cap \overline{P} $$ defines a log structure
on $Y_1$. Denote the corresponding log scheme $(Y_1,T_1)$. We obtain the
following diagram of log schemes (on the left). The corresponding diagram of
groups is placed on the right; we will use it to establish log smoothness
of $g$. Existence of dashed arrows requires an argument; we will construct
them later.

\begin{equation}\label{10-doublediagram} \xymatrix{(Y'_1,\oT'_1)
\ar[d] \ar[r]^{\oh_1} & (Y_1,\oT_1) \ar[d]\ar[dr] &
& & \oR'^\gp&\oR^\gp\ar[l]_\sim^{\nu^\gp}\\ (Y'_1,T'_1)
\ar@{-->}[r] & (Y_1,T_1)\ar[r]& \Sigma& & R'^\gp\ar@{^{(}->}[u]&
R^\gp\ar@{_{(}-->}[l]\ar@{^{(}->}[u]&Q^\gp\ar[l]\ar[lu]}
\end{equation}

Part (ii) of the following remark clarifies the notation $(Y_1,T_1)$. It
will not be used so we only sketch the argument.

\begin{remark} (i) Note that $(Y_1,T_1)$ may be not log smooth over
$\Sigma$. For example, even when $\Sigma$ is log regular, e.g. $\Spec k$ with
trivial log structure, $(Y_1,T_1)$ does not have to be log regular, as $T_1$
may even be non-Cartier. Nevertheless, as $\overline{h_1}$ is log smooth (even
log étale), $(Y'_1,\overline{T'_1})$ is log smooth over $\Sigma$. Moreover,
$Y'_1$ is regular, and $\overline{T'_1}$ an snc divisor in it.

(ii) Although $T_1$ may be bad, one does have that
$R\calO_{Y_1}^\times=\calO_{Y_1}\cap i_*\calO^\times_{Y_1\setminus
T_1}$ for the embedding $i:Y_1\setminus T_1\into Y_1$. This can
be deduced from the formulas for $R$ and $\oR$ and the fact that
$\oR\calO_{Y_1}^\times=\calO_{Y_1}\cap j_*\calO^\times_{Y_1\setminus\oT_1}$
by log regularity of $(Y_1,\oT_1)$.
\end{remark}

Note that $Q\to\oP$ factors through $P$, hence $Q\to\oR$ factors through
$R=\oR\cap P$. It follows from (\ref{10-Peq}) that $P$ consists of all elements
$f\in\oP$ whose divisor $(f=0)$ is contained in $T$ (as a set). Therefore
$g\in\oR$ lies in $R$ if and only if $(g=0)\subseteq T_1$ (as a set). This
fact and the analogous description of $R'$ observed earlier imply that
$\nu:\oR\to\oR'$ takes $R$ to $R'$. Thus, we have established the dashed
arrows in (\ref{10-doublediagram}).

Let $\varphi:Q\to\oR'$ be the homomorphism defining the composition
$(Y'_1,\overline{T'_1}) \ra (Y_1,\overline{T_1})\to\Sigma$. Since the
latter is log smooth, $\varphi$ is injective, and the torsion part
of $\Coker(\varphi^\gp)$ is annihilated by an integer $m$ invertible
in $\Lambda$.  Note that $R^\gp\into R'^\gp\into\oR'^\gp$, and therefore
we also have that $Q^\gp\into R'^\gp$ and the torsion of its cokernel is
annihilated by $m$. Therefore, $(Y'_1,T'_1)$ is log smooth over $\Sigma$,
which finishes the proof of (*), hence of $\ref{10-modif-logsmooth}$.

\begin{remark}\label{10-remark proof of (*)} In the proof of (*) above, we
first proved that $g$ is log smooth, and deduced that $g\beta$ is, too. In
fact, as $\beta$ is a Kummer étale $G$-cover, the log smoothness of $g\beta$
implies that of $g$. More generally, we have the following descent result,
due to Kato-Nakayama (\cite[3.4]{quasisections@Nakayama}) :
\end{remark}

\begin{theorem} Let $\xymatrix{X' \ar[r]^g & X \ar[r]^f & Y}$ be morphisms
of fs log schemes. If $g$ is surjective, log étale and exact,  and $fg$
is log smooth, then $f$ is log smooth.
\end{theorem} The assumption
on $g$ is equivalent to saying that $g$ is a Kummer étale cover
(cf. \cite[1.6]{Overview@Illusie}).

\section{Prime to $\ell$ variants of de Jong's alteration theorems}\label{10-primesec}

Let $X$ be a noetherian scheme, and $\ell$  be a prime number. A morphism
$h : X' \ra X$ is called an \textit{$\ell'$-alteration} if $h$ is proper,
surjective, generically finite, maximally dominating (i.e., (\refext{II}{2-maximalement-dominant})
sends each maximal point to a maximal point) and the degrees of
the residual extensions $k(x')/k(x)$ over each maximal point $x$ of $X$
are prime to $\ell$. The next theorem was stated in (Introduction 0.3 (1)) :

\begin{theorem}[Gabber]\label{10-de Jong1} Let $k$ be a field, $\ell$ a
prime number different from the characteristic of $k$, $X$ a separated and
finite type $k$-scheme, $Z \subset X$ a nowhere dense closed subset. Then
there exists a finite extension $k'$ of $k$, of degree prime to $\ell$,
and a projective $\ell'$-alteration $h : \tilde X \ra X$ above $\Spec \,
k' \ra \Spec \, k$, with $\tilde X$ smooth and quasi-projective over $k'$,
and $h^{-1}(Z)$ is the support of a relative, strict normal crossings divisor.
\end{theorem} Recall that a relative, strict normal crossings divisor in
a smooth scheme $T/S$ is a divisor $D = \sum_{i \in I}D_i$, where $I$ is
finite, $D_i \subset T$ is an $S$-smooth closed subscheme of codimension 1,
and for every subset $J$ of $I$ the scheme-theoretic intersection $\cap_{i
\in J}D_i$ is smooth over $S$.

\medskip We will need the following variant, due to Gabber-Vidal (proof
of \cite[4.4.1]{BS@Vidal}), of de Jong's alteration theorems \cite[5.7,
5.9, 5.11]{Families@de_Jong}, cf. \cite[3.8]{Independance@Zheng} :
\begin{lemma}\label{10-de Jong-Gabber-Vidal1} Let $X$ be a proper scheme over
$S = \Spec \, k$, normal and geometrically reduced and irreducible, $Z \subset X$
a nowhere dense closed subset. We assume that a finite group $H$ acts on $X \ra
S$, faithfully on $X$, and that $Z$ is $H$-stable. Then there exists a finite
extension $k_1$ of $k$, a finite group $H_1$, a surjective homomorphism $H_1
\ra H$, and an $H_1$-equivariant diagram with a cartesian square (where $S =
\Spec \, k$, $S_1 = \Spec \, k_1$)
\begin{equation}\label{10-Vidal-diagram1}
\xymatrix{X \ar[d] & X_1 \ar[d] \ar[l]_b & X_2 \ar[l]_a \ar[dl] \\ S &
S_1 \ar[l] & {}}
\end{equation} satisfying the following properties :

(i) $S_1/\mathrm{Ker}(H_1 \ra H) \ra S$ is a radicial extension ;

(ii) $X_2$ is projective and smooth over $S_1$ ;

(iii) $a : X_2 \ra X_1$ is projective and surjective, maximally dominating
and generically finite and flat,  and there exists an $H_1$-admissible dense
open subset $W \subset X_2$ over a dense open subset $U$ of $X$, such that
if $U_1 = S_1 \times_S U$ and $K = \mathrm{Ker}(H_1 \ra \mathrm{Aut}(U_1))$,
$W \ra W/K$ is a Galois étale cover of group $K$ and the morphism $W/K \ra
U_1$ induced by $a$ is a universal homeomorphism ;

(iv) $(ba)^{-1}(Z)$ is the support of a strict normal crossings divisor
in $X_2$.
\end{lemma}

\textit{Proof.} We may assume $X$ of dimension $d \ge 1$. We apply
\cite[4.4.3]{BS@Vidal} to $X/S$, $Z$, and $G = H$. We get the data of
(\textit{loc. cit.}), namely an equivariant finite extension of fields
$(S_1, H_1) \ra (S,H)$ such that $S_1/\mathrm{Ker}(H_1 \ra H) \ra S$ is
radicial, an $H_1$-equivariant pluri-nodal fibration $(Y_d \ra \cdots \ra
Y_1 \ra S_1, \{ \sigma_{ij} \}, Z_0 = \emptyset)$, and an $H_1$-equivariant
alteration $a_1 : Y_d \ra X$ over $S$, satisfying the conditions (i),
(ii), (iii) of (\textit{loc. cit.}) (in particular $a_1^{-1}(Z) \subset
Z_d$). Then, as in the proof of \cite[4.4.1]{BS@Vidal}, successively applying
\cite[4.4.4]{BS@Vidal} to each nodal curve $f_i : Y_i \ra Y_{i-1}$, one can
replace $Y_i$ by an $H_1$-equivariant projective modification $Y'_i$ of it
such that $Y'_i$ is regular, and the inverse image $Z'_i$ of $Z_i := \cup_j
\sigma_{ij}(Y_{i-1}) \cup f_i^{-1}(Z_{i-1})$ in $Y'_i$ is an $H_1$-equivariant
strict snc divisor. Then,  $X_2 := Y'_d$ is smooth over $S_1$ and $Z'_d$
is a relative snc divisor over $S_1$. This follows from the analog of
the remark following \cite[4.4.4]{BS@Vidal} with ``semistable pair over
a trait" replaced by ``pair consisting of a smooth scheme and a relative
snc divisor over a field". In particular, $(ba)^{-1}(Z)_{\mathrm{red}}$
is a subdivisor of $Z'_d$, hence an snc divisor. After replacing $H_1$ by
$H_1/\mathrm{Ker}(H_1 \ra \mathrm{Aut}(X_2))$ the open subsets $U$ and $V$
of (iii) are obtained as at the end of the proof of \cite[4.4.1]{BS@Vidal}.

\subsection{}\label{10-proof of de Jong1} \textit{Proof of \ref{10-de
Jong1}}. There are three steps.

\medskip \textit{Step 1}. \textit{Preliminary reductions.} By Nagata's
compactification theorem \cite{Nagata@Conrad}, there exists a dense open
immersion $X \subset \overline{X}$ with $\overline{X}$ proper over $S$. Up
to replacing $X$ by $\overline{X}$ and $Z$ by its closure $\overline{Z}$,
we may assume $X$ proper over $S$. By replacing $X$ by the disjoint sum
of its irreducible components, we may further assume $X$ irreducible, and
geometrically reduced (up to base changing by a finite radicial extension
of $k$). Up to blowing up $Z$ in $X$ me may further assume that $Z$ is a
(Cartier) divisor in $X$. Finally, replacing $X$ by its normalization $X'$,
which is finite over $X$, and $Z$ by its inverse image in $X'$, we may assume
$X$ normal.

\medskip \textit{Step 2}. \textit{Use of \ref{10-de Jong-Gabber-Vidal1}}.
Choose a finite Galois extension $k_0$ of $k$ such that the irreducible
components of $X_0 = X \times_S S_0$ ($S_0 = \Spec \, k_0$) are geometrically
irreducible. Let $G = \mathrm{Gal}(k_0/k)$ and $H \subset G$ the
decomposition subgroup of a component $Y_0$ of $X_0$. We apply \ref{10-de
Jong-Gabber-Vidal1} to ($Y_0/S_0$, $Z_0 \cap Y_0$), where $Z_0 = S_0 \times_S
Z$.  We find a surjection $H_1 \ra H$ and an $H_1$-equivariant diagram of
type \ref{10-Vidal-diagram1} :
\begin{equation}\label{10-Vidal-diagram2}
\xymatrix{Y_0 \ar[d] & Y_1 \ar[d] \ar[l]_{b'} & Y_2 \ar[l]_{a'} \ar[dl]
\\ S_0 & S_1 \ar[l] & {}},
\end{equation}
satisfying conditions (i), (ii),
(iii), (iv) with $S$ replaced by $S_0$, and $X_2 \ra X_1 \ra X$ by $Y_2 \ra
Y_1 \ra Y_0$. As $G$ transitively permutes the components of $X_0$, $X_0$ is,
as a $G$-scheme over $S_0$, the contracted product $$ X_0 = Y_0 \times^H G,
$$ i. e. the quotient of $Y_0 \times G$ by $H$ acting on $Y_0$ on the right
and on $G$ on the left (cf. proof of \refext{VIII}{8-rigidlem2}), and similarly
$Z = Z_0\times^H G$. Choose an extension of the diagram $\xymatrix{H_1 \ar[r]^{u}
& H \ar[r]^{i} & G}$ into a commutative diagram of finite groups $$
\xymatrix{H_1 \ar[d]^{u} \ar[r]^{i_1} & G_1 \ar[d]^v \\ H \ar[r]^i & G} $$
with $i_1$ injective and $v$ surjective (for example, take $i_1$ to be the
graph of $iu$ and $v$ the projection). Define $$ X_1 := Y_1 \times^{H_1} G_1,
\, \, \, X_2 := Y_2 \times^{H_1} G_1.  $$ Then (\ref{10-Vidal-diagram2})
extends to a $G_1$-equivariant diagram of type \ref{10-Vidal-diagram1}
\begin{equation}\label{10-Vidal-diagram3} \xymatrix{X_0 \ar[d] & X_1 \ar[d]
\ar[l]_b & X_2 \ar[l]_a \ar[dl] \\ S_0 & S_1 \ar[l] & {}},
\end{equation}
satisfying again (i), (ii), (iii), (iv). In particular, the composition $h
: X_2 \ra X_1 \ra X_0 \ra X$ is an alteration above the composition $S_1
\ra S_0 \ra S$, $X_2$ is projective and smooth over $S_1$, and $h^{-1}(Z)$
is the support of an snc divisor. However, as regard to \ref{10-de Jong1},
the diagram $$ \xymatrix{X \ar[d] & X_2 \ar[l]_h \ar[d] \\ S & S_1 \ar[l]}
$$ deduced from (\ref{10-Vidal-diagram3}) has two defects :

(a) the extension $S_1 \ra S$ is not necessarily of degree prime to $\ell$,

(b) the alteration $h$ is not necessarily an $\ell'$-alteration.

We will first repair (a) and (b) at the cost of temporarily losing the
smoothness of $X_2/S_1$ and the snc property of $h^{-1}(Z)$. By (i),
$S_1/\mathrm{Ker}(G_1 \ra G) \ra S_0$ is a radicial extension, hence
$S_1/G_1 \ra S = S_0/G$ is a radicial extension, too. Similarly, by (iii),
$X_2/G_1 \ra X$ is an alteration over $S_1/G_1 \ra S$,  which is a universal
homeomorphism over a dense open subset. Now let $L$ be an $\ell$-Sylow subgroup
of $G_1$. Then $S_1/L \ra S_1/G_1$ is of degree prime to $\ell$, and $X_2/L \ra
X_2/G_1$ is a finite surjective morphism of generic degree prime to $\ell$. Let
$S' := \Spec \, k' = S_1/L$, $X' = S' \times_S X$. We get a commutative diagram
with cartesian square
\begin{equation}\label{10-Vidal-diagram4} \xymatrix{X
\ar[d] & X' \ar[d] \ar[l] & X_2/L \ar[l] \ar[dl] & X_2 \ar[l] \ar[dl]\\ S &
S' \ar[l] & S_1\ar[l]},
\end{equation}
where $S' \ra S$ is an extension of
degree prime to $\ell$, $S_1 \ra S'$ a Galois extension of group $L$, $$ h_2
: X_2/L \ra X $$ an $\ell'$-alteration, $X_2/S_1$ is projective and smooth,
and if now $h$ denotes the composition $X_2 \ra X$, $Z_1 := h^{-1}(Z)$ is
an snc divisor in $X_2$. If $X_2/L$ was smooth over $S'$ and $Z_1/L$ an snc
divisor in $X_2/L$, we would be finished. However, this is not the case in
general. We will use Gabber's theorem \ref{10-modif-logsmooth}  to fix this.

\medskip \textit{Step 3.} \textit{Use of \ref{10-modif-logsmooth}.} Let $Y$
be a connected component of $X_2$, $(Z_1)_Y = h^{-1}(Z) \cap Y$, $D$ the
stabilizer of $Y$ in $L$, $I \subset D$ the inertia group at the generic
point of $Y$. Then $D$ acts on $Y$ through $K := D/I$, and this action
is generically free. As $Y$ is smooth over $S_1$ and $(Z_1)_Y$ is an snc
divisor in $Y$, $(Y,(Z_1)_Y)$ makes a log regular pair, log smooth over $S_1$,
hence over $S' = S_1/L$ (equipped with the trivial log structure). We have a
$K$-equivariant commutative diagram
\begin{equation}\label{10-Gabber-diagram1}
\xymatrix{(Y/K,(Z_1)_Y/K) \ar[d] & (Y,(Z_1)_Y) \ar[l] \ar[dl]_f \\ S'},
\end{equation} where $K$ acts trivially on $S'$, and $f$ is projective,
smooth, and log smooth ($S'$ having the trivial log structure). We now apply
\ref{10-modif-logsmooth} to ($f : (Y,(Z_1)_Y) \ra S' = (S',\emptyset),
K$), which satisfies conditions (i) - (iii) of (\textit{loc. cit.}). We
get a $D$-equivariant projective modification $g : Y_1 \ra Y$ ($D$ acting
through $K$) and a $D$-strict snc divisor $E_1$ on $Y_1$, containing
$Z_1 := g^{-1}((Z_1)_Y)$ as a subdivisor, such that the action of $D$
on $(Y_1,E_1)$ is very tame, and $(Y_1,E_1)$ and $(Y_1/D,E_1/D)$ are log
smooth over $S'$. Pulling back $g$ to the orbit $Y \times^D L$ of $Y$
under $L$, i. e. replacing $g$ by $g \times^D L$, and working separately
over each orbit, we eventually get an $L$-equivariant commutative square
\begin{equation}\label{10-Gabber-diagram2} \xymatrix{(Y_2/L,E_2/L) \ar[d]_v &
(Y_2,E_2) \ar[l] \ar[d]_u\\ (X_2/L,Z_1/L) & (X_2,Z_1) \ar[l]},
\end{equation}
where $u$, $v$ are projective modifications (and $Z_1 = h^{-1}(Z)$, $Z_1/L =
h_2^{-1}(Z)$ as above), with the property that the pair $(Y_2/L,E_2/L)$ is
an fs log scheme log smooth over $S'$ ($= S_1/L$), and $v^{-1}(h_2^{-1}(Z))
\subset E_2/L)$. Let $w : (\tilde X, \tilde E) \ra (Y_2/L,E_2/L)$ be a
projective, log étale modification such that $\tilde X$ is regular, and
$\tilde E = w^{-1}(E_2/L)$ is an snc divisor in $\tilde X$. For example,
one can take for $w$ the saturated monoidal desingularization $\tilde
{\mathscr{F}}^{\mathrm{log}}$ of (\refext{VIII}{8-deslogregth}). We then apply
\ref{10-local-structure} to the log smooth morphism $\tilde X \ra S'$. By a
special case of the (\ref{10-chart-diagram}), with $Q = \{1\}$, $G = \{1\}$,
as $P^*$ is torsion free, $\tilde X$ is not only regular, but smooth over $S'$,
and $\tilde E$ a relative snc on $\tilde X$. Let $$ \tilde h : \tilde X \ra X
$$ be the composition $$ \xymatrix{\tilde X \ar[r]^w & Y_2/L \ar[r]^v & X_2/L
\ar[r]^{h_2} & X}.  $$ This is a projective $\ell'$-alteration, and it fits
in the commutative diagram $$ \xymatrix{X \ar[d] & \tilde X \ar[l]_{\tilde
h} \ar[d] \\ S & S' \ar[l]}, $$ where $S'$ is an extension of $S$ of degree
prime to $\ell$, $\tilde X$ is projective and smooth over $S'$, and $\tilde
h^{-1}(Z)_{\mathrm{red}}$ is a sub-divisor of the $S'$-relative snc divisor
$\tilde E$, hence a relative snc divisor as well. This finishes the proof
of \ref{10-de Jong1}.

Recall now the theorem stated in (Introduction 0.3 (2)) :
\begin{theorem}[Gabber]\label{10-de Jong2} Let $S$ be a separated, integral,
noetherian, excellent, regular scheme of dimension 1, with generic point
$\eta$, $X$ a scheme separated, flat and of finite type over $S$, $\ell$
a prime number invertible on $S$, $Z \subset X$ a nowhere dense closed
subset. Then there exists a finite extension $\eta'$ of $\eta$ of degree prime
to $\ell$ and a projective $\ell$'-alteration $h : \tilde X \ra X$ above $S'
\ra S$, where $S'$ is the normalization of $S$ in $\eta'$, with $\tilde X$
regular and quasi-projective over $S'$, a strict normal crossings divisor
$\tilde T$ on $\tilde X$, and a finite closed subset $\Sigma$ of $S'$ such
that :

(i) outside $\Sigma$, $\tilde X \ra S'$ is smooth and $\tilde T \ra S'$
is a relative divisor with normal crossings ;

(ii) étale locally around each geometric point $x$ of $\tilde X_{s'}$,
where $s' = \Spec \, k'$ belongs to $\Sigma$, the pair $(\tilde X, \tilde T)$
is isomorphic to the pair consisting of
$$ X'=S'[u_1^{\pm 1}, \cdots, u_s^{\pm 1}, t_1, \cdots, t_n]/
(u_1^{b_1} \cdots u_s^{b_s}t_1^{a_1} \cdots t_r^{a_r}- \pi), $$
$$ T'=V(t_{r+1}\cdots t_m)\subset X'$$
around the point $(u_i = 1), (t_j = 0)$, with $1
\le r \le m \le n$, for positive integers $a_1, \cdots, a_r, b_1, \cdots, b_s$
satisfying $\mathrm{gcd}(p,a_1, \cdots, a_r, b_1, \cdots, b_s) = 1$, $p$ the
characteristic exponent of $k'$, $\pi$ a local uniformizing parameter at $s'$ ;

(iii) $\tilde h^{-1}(Z)_{\mathrm{red}}$ is a sub-divisor of $\cup_{s' \in
\Sigma}(\tilde X_{s'})_{\mathrm{red}} \cup \tilde T$.
\end{theorem}

The proof follows the same lines as that of \ref{10-de Jong1}. We
need again a Gabber-Vidal variant of de Jong's alteration
theorems (cf. \cite[3.8]{Independance@Zheng}). This is essentially
\cite[4.4.1]{BS@Vidal}), except for the additional data of $Z \subset X$
and the removal of the hypothesis that $S$ is a strictly local trait :

\begin{lemma}\label{10-de Jong-Gabber-Vidal2} Let $X$ be a normal,
proper scheme over $S$, whose generic fiber is geometrically reduced and
irreducible, $Z \subset X$ a nowhere dense closed subset. We assume that
a finite group $H$ acts on $X \ra S$, faithfully on $X$, and that $Z$ is
$H$-stable. Then there exists a finite group $H_1$, a surjective homomorphism
$H_1 \ra H$, and an $H_1$-equivariant diagram with a cartesian square
\begin{equation}\label{10-Vidal-diagram5} \xymatrix{X \ar[d] & X_1 \ar[d]
\ar[l]_b & X_2 \ar[l]_a \ar[dl] \\ S & S_1 \ar[l] & {}},
\end{equation}
satisfying the following properties :

(i) $S_1 \ra S$ is the normalization of $S$ in a finite extension $\eta_1$
of $\eta$ such that $\eta_1/(\mathrm{Ker}(H_1 \ra H) \ra \eta$ is a radicial
extension (where $\eta$ is the generic point of $S$) ;

(ii) $X_2$ is projective and strictly semistable over $S_1$ (i. e. is
strictly semistable over the localizations of $S_1$ at closed points
\cite[2.16]{Alterations@de_Jong}) ;

(iii) $a : X_2 \ra X_1$ is  projective and surjective, maximally dominating and
generically finite and flat, and  there exists an $H_1$-admissible dense open
subset $W \subset X_2{}_{\eta_1}$ over a dense open subset $U$ of $X_{\eta}$,
such that if $U_1 = \eta_1 \times_{\eta} U$ and $K = \mathrm{Ker}(H_1 \ra
\mathrm{Aut}(U_1))$, $W \ra W/K$ is a Galois étale cover of group $K$
and the morphism $W/K \ra U_1$ induced by $a$ is a universal homeomorphism ;

(iv) $(ba)^{-1}(Z)$ is the support of a strict normal crossings
divisor in $X_2$, and $(X_2,(ba)^{-1}(Z))$ is a strict semistable pair
over $S_1$ (i. e. over the localizations of $S_1$ at closed points
\cite[6.3]{Alterations@de_Jong}).
\end{lemma}

\medskip Note that (ii) and (iv) imply that there exists a finite closed subset
$\Sigma$ of $S_1$ such that, outside $\Sigma$, the pair $(X_2,(ba)^{-1}(Z))$
consists of a smooth morphism and a relative strict normal crossings divisor.

\medskip \textit{Proof.} Up to minute modifications the proof is the
same as that of \ref{10-de Jong-Gabber-Vidal1}. We may assume the generic
fiber $X_{\eta}$ of dimension $d \ge 1$. We apply \cite[4.4.3]{BS@Vidal}
to $X/S$, $Z$, and $G = H$. We get the data of (\textit{loc. cit.}), namely
an equivariant finite surjective morphism $(S_1, H_1) \ra (S,H)$, with $S_1$
regular (hence equal to the normalization of $S$ in the extension $\eta_1$
of the generic point $\eta$ of $S$) such that $\eta_1/\mathrm{Ker}(H_1 \ra H)
\ra \eta$ is radicial, an $H_1$-equivariant pluri-nodal fibration $(Y_d \ra
\cdots \ra Y_1 \ra S_1, \{ \sigma_{ij} \}, Z_0)$, and an $H_1$-equivariant
alteration $a_1 : Y_d \ra X$ over $S$, satisfying the conditions (i),
(ii), (iii) of (\textit{loc. cit.}) (in particular $a_1^{-1}(Z) \subset
Z_d$). Then, as in the proof of \cite[4.4.1]{BS@Vidal}, successively applying
\cite[4.4.4]{BS@Vidal} to each nodal curve $f_i : Y_i \ra Y_{i-1}$, one
can replace $Y_i$ by an $H_1$-equivariant projective modification $Y'_i$
of it such that $Y'_i$ is regular, and the inverse image of $Z_i := \cup_j
\sigma_{ij}(Y_{i-1}) \cup f_i^{-1}(Z_{i-1})$ in $Y'_i$ is an $H_1$-equivariant
strict snc divisor. Then, by the remark following \cite[4.4.4]{BS@Vidal} $X_2
:= Y'_d$ is strict semistable over $S_1$ and $(X_2,Z_d)$ is a strict semistable
pair over $S_1$. In particular, $(ba)^{-1}(Z)_{\mathrm{red}}$ is a subdivisor
of $Z_d$, hence an snc divisor, and $(X_2,(ba)^{-1}(Z))_{\mathrm{red}}$ is
a strict semistable pair over $S_1$. The open subsets $U$ and $V$ of (iii)
are constructed as at the end of the proof of \cite[4.4.1]{BS@Vidal}.

\subsection{}\label{10-proof of de Jong2} \textit{Proof of \ref{10-de
Jong2}}. It is similar to that of \ref{10-de Jong1}. There are again three
steps. We will indicate which modifications should be made.

\medskip \textit{Step 1.} \textit{Preliminary reductions}. Up to replacing $X$
by the disjoint union of the schematic closures of the reduced components of
its generic fiber, and working separately with each of them, we may assume $X$
integral (and $X_{\eta} \ne \emptyset$). Applying Nagata's compactification
theorem, we may further assume $X$ proper and integral.  Base changing by
the normalization of $S$ in a finite radicial extension of $\eta$ and taking
the reduced scheme, we reduce to the case where, in addition, $X_{\eta}$
is irreducible and geometrically reduced. Then we blow up $Z$ in $X$ and
normalize as in the previous step 1. Here we used the excellency of $S$
to guarantee the finiteness of the normalizations.

\medskip \textit{Step 2.} \textit{Use of \ref{10-de Jong-Gabber-Vidal2}}. Let
$S_0$  be the normalization of $S$ in a finite Galois extension $\eta_0$ of
$\eta$ such that the irreducible components of the generic fiber of $X_0 =
X \times_S S_0$ ($S_0 = \Spec \, k_0$) are geometrically irreducible. Let $G
= \mathrm{Gal}(\eta_0/\eta)$ and $H \subset G$ the decomposition subgroup
of a component $Y_0$ of $X_0$. We apply \ref{10-de Jong-Gabber-Vidal2}
to ($Y_0/S_0$, $Z_0 \cap Y_0$), where $Z_0 = S_0 \times_S Z$.  We find
a surjection $H_1 \ra H$ and an $H_1$-equivariant diagram of type
(\ref{10-Vidal-diagram5}) satisfying conditions (i), (ii), (iii), (iv)
with $S$ replaced by $S_0$, and $X_2 \ra X_1 \ra X$ by $Y_2 \ra Y_1 \ra
Y_0$. We then, as above,  extend $H_1 \ra H$ to a surjection $G_1 \ra G$
and obtain a $G_1$-equivariant diagram of type (\ref{10-Vidal-diagram5})
\begin{equation}\label{10-Vidal-diagram6} \xymatrix{X_0 \ar[d] & X_1 \ar[d]
\ar[l]_b & X_2 \ar[l]_a \ar[dl] \\ S_0 & S_1 \ar[l] & {}},
\end{equation}
satisfying again (i), (ii), (iii), (iv). In particular, the composition $h :
X_2 \ra X_1 \ra X_0 \ra X$ is an alteration above the composition $S_1 \ra S_0
\ra S$, $X_2$ is projective over $S_1$, with strict semistable reduction, and
$h^{-1}(Z)$ is the support of an snc divisor forming a strict semistable pair
with $X_2/S_1$. It follows from (i) that $S_1/G_1 \ra S = S_0/G$ is generically
radicial, and by (iii) that $X_2/G_1 \ra X$ is an alteration over $S_1/G_1 \ra
S$,  which is a universal homeomorphism over a dense open subset. As above,
choose an $\ell$-Sylow subroup $L$ of $G_1$. Then $S_1/L$ is regular, $S_1/L
\ra S_1/G_1$ is finite surjective of generic degree prime to $\ell$, and
$X_2/L \ra X_2/G_1$ is a finite surjective morphism of generic degree prime
to $\ell$. Putting $S' = S_1/L$, $X' = S' \times_S X$, we get a commutative
diagram with cartesian square
\begin{equation}\label{10-Vidal-diagram7}
\xymatrix{X \ar[d] & X' \ar[d] \ar[l] & X_2/L \ar[l] \ar[dl] & X_2 \ar[l]
\ar[dl]\\ S & S' \ar[l] & S_1\ar[l]},
\end{equation}
where $S'$ is regular,
$S' \ra S$ is finite surjective of generic degree prime to $\ell$, $S_1 \ra
S'$ generically étale of degree the order of $L$, $$ h_2 : X_2/L \ra X $$
an $\ell'$-alteration, $X_2/S_1$ is projective and strictly semistable, and
if $h$ denotes the composition $X_2 \ra X$, $Z_1 := h^{-1}(Z)_{\mathrm{red}}$
is an snc divisor in $X_2$, forming a strictly semistable pair with $X_2/S_1$.

\medskip \textit{Step 3}. \textit{Use of \ref{10-modif-logsmooth}}. Defining
$Y$, $(Z_1)_Y$, $I \subset D$, $K = D/I$ as in the former step 3, $K$ acts
generically freely on $Y$. As the pair $(Y,(Z_1)_Y)$ is strictly semistable
over $S_1$, there exists a finite closed subset $\Sigma_1$ of $S_1$ such that
$(Y, Y_{\Sigma_1} \cup (Z_1)_Y)$ forms a log regular pair, log smooth over
$S_1$ equipped with the log structure defined by $\Sigma_1$. As $S_1 \ra S'
= S_1/L$ is Kummer étale, $(Y,(Y_{\Sigma'})_{\mathrm{red}} \cup (Z_1)_Y)$
(where $\Sigma'$ is the inverse image of $\Sigma_1$) is also log smooth
over $S'$ (equipped with the log structure given by $\Sigma'$), and we get a
$K$-equivariant commutative diagram (\ref{10-Gabber-diagram1}), with trivial
action of $K$ on $S'$ and $f$ projective and log smooth over $S'$. We then
apply \ref{10-modif-logsmooth} to $f : (Y,(Y_{\Sigma'})_{\mathrm{red}} \cup
(Z_1)_Y) \ra S'$, and the proof runs as above. We get a $D$-equivariant
projective modification $g : Y_1 \ra Y$ ($D$ acting through $K$) and a
$D$-strict snc divisor $E_1$ on $Y_1$, containing $(g^{-1}((Z_1)_Y) \cup
(Y_1)_{\Sigma'})_{\mathrm{red}}$ as a subdivisor, such that the action
of $D$ on $(Y_1,E_1)$ is very tame, and $(Y_1,E_1)$ and $(Y_1/D,E_1/D)$
are log smooth over $S'$. After extending from $D$ to $L$ we get an
$L$-equivariant commutative square of type (\ref{10-Gabber-diagram2}), with
$(Y_2/L,E_2/L)$ log smooth over $S'$ ($= S_1/L$), and $(v^{-1}(h_2^{-1}(Z))
\cup (Y_2/L)_{\Sigma'})_{\mathrm{red}} \subset E_2/L$. As above, we take
a projective, log étale modification such that $\tilde X$ is regular,
and $\tilde E = w^{-1}(E_2/L)$ is an snc divisor in $\tilde X$.

We now apply \ref{10-local-structure} to the log smooth morphism
$(\tilde X, \tilde E) \ra S'$. It's enough to work étale locally on
$\tilde X$ around some geometric point $x$ of $\tilde X_{s'}$, with $s'
\in \Sigma'$. We replace $S'$ by its strict localization at the image
of $x$, and consider (\ref{10-chart-diagram}), with $Q = \mathbf{N}$,
$G = \{1\}$, $\Lambda = \mathbf{Z}_{(p)}$ if $p > 1$ and $\mathbf{Q}$
otherwise, and the chart $a : \mathbf{N} \ra M_{S'}, 1 \mapsto \pi$,
$\pi$ a uniformizing parameter of $S'$. In (\ref{10-P1}) we have $P^*=
\mathbf{Z}^{\mu}$, $P_1 = \mathbf{N}^{\nu}$, for nonnegative integers
$\mu$, $\nu$, hence $$ P = \mathbf{Z}^{\nu} \oplus \mathbf{N}^{\mu}.
$$ Let $((b_i),(a_i)$ be the image of $1 \in \mathbf{N}$ in $P$ in the
above decomposition, and let $(a_1,\cdots,a_r)$, $(b_1,\cdots,b_s)$ be the
sets of those $a_i$'s and $b_i$'s which are $\ne 0$. We may assume $b_i
>0$ if $b_i \ne 0$. As the torsion part of $\mathrm{Coker}(\mathbf{Z}
\ra P^{gp})$ is annihilated by an integer invertible on $\tilde X$, we
have $\mathrm{gcd}(p,a_1,\cdots,a_r,b_1,\cdots,b_s) = 1$, where $p$ is the
characteristic exponent of $k$. We have $P = \mathbf{Z}^s \oplus \mathbf{N}^r
\oplus \mathbf{Z}^{\nu-s} \oplus \mathbf{N}^{\mu - r}$. Choosing a basis of
$V$, we get that étale locally around $x$, $\tilde X$ is given by a cartesian
square $$ \xymatrix{\tilde X \ar[d] \ar[r] & \Spec \, \Lambda[u_1^{\pm 1},
\cdots, u_s^{\pm 1}, t_1, \cdots, t_n] \ar[d]\\ S' \ar[r] & \Spec \, \Lambda[z]}
$$ with $x$ going to the point $(u_i = 1), (t_j = 0)$, and $z \mapsto \pi$,
$z \mapsto u_1^{b_1} \cdots u_s^{b_s}t_1^{a_1}\cdots t_r^{a_r}$, in other
words, $$ \tilde X = S'[u_1^{\pm 1}, \cdots, u_s^{\pm 1}, t_1, \cdots,
t_n]/(u_1^{b_1} \cdots u_s^{b_s}t_1^{a_1} \cdots t_r^{a_r} - \pi), $$ Finally,
$\tilde E$ is the union of the special fiber $\tilde X_{s'}$ and a horizontal
divisor $\tilde T$, étale locally given by the equation $t_{r+1}\cdots
t_m = 0$, for an integer $m$ such that $1 \le r \le m \le n$. As $\tilde
h^{-1}(Z)_{\mathrm{red}}$ is a sub-divisor of $(\tilde X_{s'})_{\mathrm{red}}
\cup \tilde T$, this finishes the proof of \ref{10-de Jong2}.

\section{Resolvability, log smoothness, and weak semistable reduction}\label{10-complsec}

\subsection{Elimination of separatedness assumptions}\label{10-sepelimsec}
The main aim of \S\ref{10-sepelimsec} is to weaken the separatedness assumptions
in Theorems \refext{VIII}{8-TH. MODIFICATION} and \ref{10-modif-logsmooth}.
We start with an example showing how things can go wrong without such an assumption.
Recall, see \refext{VI}{6-stratification}, that if a finite group $G$ acts on
a scheme $X$ then the fixed point subscheme $X^G$ is the intersection of graphs
of the translations $g{\colon}X\to X$.  In particular, $X^G$ is closed whenever $X$
is separated.

\begin{example2}\label{10-nonsepex}
Let $G=\{1,\sigma\}$ be a cyclic group of order two. Consider an affine plane
$\AA^2_k=\Spec(k[x,y])$ over a field $k$ of characteristic different from $2$
and let $X$ be the affine plane with doubled origin, i.e., $X$ is glued from two
copies of $\AA^2_k$ along $U=\AA^2_k\setminus\{0\}$. We provide $U$ with the
action of $G$ by the rules $\sigma x=x$ and $\sigma y=-y$, and extend this action
to $X$ so that $\sigma$ intertwines the origins. Then $X^G$ is the closed subscheme
of $U$ given by $y=0$. In particular, it is not closed in $X$, and the inertia
stratification (used in the proof of Theorem \refext{VIII}{8-TH. MODIFICATION})
does not make sense. Note that the closure $Z$ of $X^G$ is an affine line with
doubled origin on which $G$ acts by permuting the origins. For any $G$-equivariant
modification $Y\to X$ the proper transform of $Z$ is a modification of $Z$ and hence
is isomorphic to $Z$. Thus, $Y^G$ contains the generic point of $Z$ but not its origins,
and we obtain that $Y^G$ is not closed. Therefore, the action on $Y$ is not very tame
for any choice of a log structure, and we obtain that the assertion of Theorem
\refext{VIII}{8-TH. MODIFICATION} fails for such $(G,X,Z=\emptyset)$.
\end{example2}

\subsubsection{Inertia specializing actions}
Example \ref{10-nonsepex} motivates the following definition: an action of
a finite group $G$ on a scheme $X$ is {\em inertia specializing} if for any
point $x\in X$ with a specialization $y\in X$ one has that $G_x\subseteq G_y$.
Since a subscheme of $X$ is closed if and only if it is closed under specializations,
we obtain the following result.

\begin{lemma2}
An action of $G$ on $X$ is inertia specializing if and only if for each subgroup
$H\subseteq G$ the subscheme $X^H$ is closed.
\end{lemma2}

\begin{remark2}\label{10-inspecrem}
(i) A large class of examples of inertia specializing actions can be described as follows.
The following conditions are equivalent and imply that the action is inertia specializing:
(a) any $G$ orbit is contained in an open separated subscheme of $X$, (b) $X$ admits
a covering by $G$-equivariant separated open subschemes $X_i$. In particular, any
admissible action is inertia specializing.

(ii) If $(G,X,Z)$ is as in Theorem \refext{VIII}{8-TH. MODIFICATION}, but instead of
separatedness of $X$ one only assumes that it possesses a covering by $G$-equivariant
separated open subschemes $X_i$, then the assertion of the theorem still holds true.
Indeed, the theorem applies to the $G$-equivariant log schemes $(X_i,Z_i=Z|_{X_i})$,
and by Theorem \refext{VIII}{8-addendum} the modifications $f_{(G,X_i,Z_i)}$ agree
on the intersections and hence glue to a required modification $f_{(G,X,Z)}$ of $X$.
\end{remark2}

A quick analysis of the proof of \refext{VIII}{8-TH. MODIFICATION} is required
to obtain the following stronger result.

\begin{theorem2}\label{10-elimth}
(i) Theorem \refext{VIII}{8-TH. MODIFICATION} and its complement \refext{VIII}{8-addendum}
hold true if the assumption that $X$ is separated is weakened to the assumption that
the action of $G$ on $X$ is inertia specializing.

(ii) Theorem \ref{10-modif-logsmooth} holds true if the assumptions that $X$ and $S$ are
separated are replaced with the single assumption that the action of $G$ on $X$ is inertia specializing.
\end{theorem2}
\begin{proof}
The proof of \refext{VIII}{8-TH. MODIFICATION} in the general (i.e. not necessarily qe)
case runs as follows. First, for any point $x\in X$ and localizations $X_x$ and $Z_x$
one constructs $f_x=f_{(G_x,X_x,Z_x)}$. This stage does not use any separatedness assumption.
Second, for any specialization $y$ of $x$ one shows that after eliminating empty normalized
blow ups, the restriction  $f_y|_{X_x}$ coincides with $f_x$. For this stage to work we only
need that $G_x\subseteq G_y$, i.e. the assumption that the action is inertia specializing is
precisely what one uses. The third stage is to show that the family $f_y$ forms an atomic
normalized blow up tower and hence descends to a normalized blow up tower $f_{(G,X,Z)}$.
It uses the log-inertia stratification in order to control the associated points of the
normalized blow up towers $f_x$. Since the inertia stratification is well defined for
inertia specializing actions, the log-inertia stratification is well defined too, and
this stage works fine. Finally, the last stage is to uniformly bound the normalization
threshold of $f_{(G,X,Z)}$ in terms of the combinatorial datum. The argument is local
and the combinatorial datum involves only the log-inertia stratification and the size
of $G$. Hence this stage also extends to our case without any changes.

The proof of Theorem \ref{10-modif-logsmooth} runs as follows. One considers the modification
$f_{(G,X,Z)}$ from \refext{VIII}{8-TH. MODIFICATION} and checks that it satisfies the additional
properties asserted by \ref{10-modif-logsmooth}. This check is local on $X$ and hence applies to
non-separated schemes as well. Since  by part (i) of \ref{10-elimth}, $f_{(G,X,Z)}$ is well defined
whenever $G$ acts inertia specializing on $X$, we obtain (ii).
\end{proof}

\subsubsection{Pseudo-projective morphisms and non-separated Chow's lemma}\label{10-chowsec}
We conclude \S\ref{10-sepelimsec} with recalling a non-separated version of Chow's lemma
due to Artin-Raynaud-Gruson. It will be needed to avoid unnecessary separatedness
assumptions in the future. We prefer to use the following non-standard terminology:
a finite type morphism $f{\colon}X\to S$ is {\em pseudo-projective} if it can be factored
into a composition of a {\em local isomorphism} $X\to\oX$ (i.e. $X$ admits an open covering
$X=\cup_i X_i$ such that the morphisms $X_i\to\oX$ are open immersions) and a projective morphism $\oX\to S$.

\begin{remark2}\label{10-chowrem}
(i) We introduce pseudo-projective morphisms mainly for terminological convenience.
Although pseudo-projectivity is preserved by base changes, it can be lost under
compositions. Moreover, even if $X$ is pseudo-projective over a field $k$, its
blow up $X'$ does not have to be pseudo-projective over $k$ (thus giving an example
of a projective morphism $f{\colon}X'\to X$ and a pseudo-projective one $X\to\Spec(k)$
so that the composition is not pseudo-projective). Indeed, let $X$ be an affine plane
with a doubled origin $\{o_1,o_2\}$, and let $X'$ be obtained by blowing up $o_1$. By
$\eta$ we denote the generic point of $C_1=f^{-1}(o_1)$. The ring $\Spec(\calO_{X',\eta})$
is a DVR and its spectrum has two different $k$-morphisms to $X'$: one takes the closed point
to $\eta$ and another one takes it to $o_2$. It then follows from the valuative criterion of
separatedness that any $k$-morphism $g{\colon}X'\to Y$ with a separated target takes $o_2$
and $\eta$ to the same point of $Y$. In particular, such $g$ cannot be a local isomorphism,
and hence $X'$ is not pseudo-projective over $k$.

(ii) Note that a morphism $f$ is separated (resp. proper) and pseudo-projective if and only
if it is quasi-projective (resp. projective). So, the following result extends the classical
Chow's lemma to non-separated morphisms.
\end{remark2}

\begin{proposition2}\label{10-chowlem}
Let $f{\colon}X\to S$ be a finite type morphism of quasi-compact and quasi-separated schemes,
and assume that $X$ has finitely many maximal points. Then there exists a projective
modification $g{\colon}X'\to X$ (even a blow up along a finitely generated ideal with
a nowhere dense support) such that the morphism $X'\to S$ is pseudo-projective.
\end{proposition2}
\begin{proof}
As a simple corollary of the flattening theorem, it is proved in \cite[5.7.13]{Platification@Raynaud-Gruson}
that there exists a modification $X'\to X$ such that $X'\to S$ factors as a composition of an
étale morphism $X'\to\oX$ that induces an isomorphism of dense open subschemes and a projective
morphism $\oX\to S$. (In loc.cit. one works with algebraic spaces and assumes that $f$ is locally
separated, but the latter is automatic for any morphism of schemes.) Our claim now follows from
the following lemma (which fails for locally separated morphisms between algebraic spaces).
\end{proof}

\begin{lemma2}
Assume that $\phi{\colon}Y\to Z$ is an étale morphism of schemes that restricts to an open
embedding on a dense open subscheme $Y_0\into Y$. Then $\phi$ is a local isomorphism.
\end{lemma2}
\begin{proof}
Let us prove that if, in addition, $\phi$ is separated then $\phi$ is an open immersion.
Since $Y$ possesses an open covering by separated subschemes, this will imply the lemma.
The diagonal $\Delta_\phi{\colon}Y\to Y\times_ZY$ is an open immersion, and by the
separatedness of $\phi$, it is also a closed immersion. Thus, $Y$ is a connected
component of $Y\times_ZY$, and since both $Y$ and $Y\times_ZY$ have dense open
subschemes that map isomorphically onto $Y_0$, $\Delta_\phi$ is an isomorphism.
This implies that $\phi$ is a monomorphism, but any étale monomorphism is an
open immersion by \cite[17.9.1]{EGAIV-4}.
\end{proof}

\subsection{Semistable curves and log smoothness}

\subsubsection{Log structure associated to a closed subset}
Let $S$ be a reduced scheme. Any closed nowhere dense subset $W\subset S$ induces
a log structure $j_*\calO^*_U\cap\calO_S\into\calO_S$, where $j{\colon}U\into S$
is the embedding of the complement of $W$. The associated log scheme will be denoted
$(S,W)$. By \refext{VI}{6-1.4}, any log regular log scheme is of the form $(S,W)$,
where $W$ is the non-triviality locus of the log structure.

\subsubsection{Semistable relative curves}\label{10-relcurvesec}
Following the terminology of \cite{Temkin@stable}, by a semistable multipointed
relative curve over a scheme $S$ we mean a pair $(C,D)$, where $C$ is a flat
$S$-scheme of relative dimension one and with fibers having only ordinary nodes
as singularities, and $D\into C$ is a closed subscheme which is étale over $S$
and disjoint from the singular locus of $C\to S$. We do not assume $C$ to be neither
proper nor even separated over $S$.

\begin{proposition2}\label{10-logsmprop}
Assume that $(S,W)$ is a log regular log scheme and $(C,D)$ is a semistable multipointed
relative $S$-curve such that the morphism $f{\colon}C\to S$ is smooth over $S\setminus W$.
Then the morphism of log schemes $(C,D\cup f^{-1}(W))\to(S,W)$ is log smooth.
\end{proposition2}
\begin{proof}
If $C$ is proper over $S$ then this claim is nothing else than \refext{VI}{6-nodal implique
log lisse}. Moreover, the proof of \refext{VI}{6-nodal implique log lisse}. is local on $C$,
hence it applies to our case as well.
\end{proof}

\subsection{$l'$-resolvability}\label{10-lresolvesec}

\subsubsection{Alterations}
Assume that $S'$ and $S$ are reduced schemes with finitely many maximal points and let
$\eta'\subseteq S'$ and $\eta\subseteq S$ denote the subschemes of maximal points.
Let $f{\colon}S'\to S$ be an {\em alteration}, i.e. a proper, surjective, generically finite,
and maximally dominating morphism. Recall that $f$ is an {\em $l'$-alteration} if
$([k(x):k(f(x))],l)=1$ for any $x\in\eta'$. We say that $f$ is {\em separable} if
$k(\eta')$ is a separable $k(\eta)$-algebra (i.e. $\eta'\to\eta$ has geometrically
reduced fibers). If, in addition, $S'$ and $S$ are provided with an action of
a finite group $G$ such that $f$ is $G$-equivariant, the action on $S$ is trivial,
the action on $\eta'$ is free, and $\eta'/G\iso\eta$, then we say that $f$ is a
{\em separable Galois alteration} of group $G$, or just {\em separable $G$-alteration}.

\begin{remark2}
We add the word "separable" to distinguish our definition from Galois alterations
in the sense of de Jong (see \cite{Alterations@de_Jong}) or Gabber-Vidal
(see \cite[p. 370]{Brauer@Vidal}). In the latter cases, one allows alterations
that factor as $S'\to S''\to S$, where $S'\to S''$ is a separable Galois alteration
and $S''\to S$ is purely inseparable.
\end{remark2}

\subsubsection{Universal $l'$-resolvability}\label{10-resolvsec}
Let $X$ be a noetherian scheme and let $l$ be a prime invertible on $X$. Assume that
for any alteration $Y\to X$ and nowhere dense closed subset $Z\subset Y$ there exists
a surjective projective morphism $f{\colon}Y'\to Y$ such that $Y'$ is regular and
$Z'=f^{-1}(Z)$ is an snc divisor. (By a slight abuse of language, by saying that
a closed subset is an snc divisor we mean that it is the support of an snc divisor.)
If, furthermore, one can always choose such $f$ to be a modification, separable
$l'$-alteration, $l'$-alteration, or alteration, then we say that $X$ is
{\em universally resolvable}, {\em universally separably $l'$-resolvable},
{\em universally $l'$-resolvable}, or {\em universally $\QQ$-resolvable}, respectively.

\begin{remark2}
(i) Due to resolution of singularities in characteristic zero, any qe scheme over
$\Spec(\QQ)$ is universally resolvable for any $l$. This is essentially due to
Hironaka, \cite{Hironaka}, though an additional work was required to treat qe
schemes that are not algebraic in Hironaka's sense, see \cite{Temkin@des1} for
the noetherian case and \cite{Temkin@des2} for the general case.

(ii) It is hoped that all qe schemes admit resolution of singularities
(in particular, are universally resolvable). However, we are, probably,
very far from proving this. Currently, it is known that any qe scheme of
dimension at most two admits resolution of singularities (see \cite{CJS@surfaces}
for a modern treatment). In particular, any qe scheme of dimension at most two is
universally resolvable.

(iii) One can show that any universally $\QQ$-resolvable scheme is qe, but we
prefer not to include this proof here, and will simply add quasi-excellence
assumption to the theorems below.

(iv) On the negative side, we note that there exist regular (hence resolvable)
but not universally $\QQ$-resolvable schemes $X$. They can be constructed analogously
to examples from \refext{I}{1patho-DVR}. For instance, there exists a discretely valued
field $K$ whose completion $\hatK$ contains a finite purely inseparable extension $K'/K$
(e.g. take an element $y\in k((x))$ which is transcendental over $k(x)$ and set
$K=k(x,y^p)\subset  K' = k(x,y) \subset k((x))$ with the induced valuation).
The valued extension $K'/K$ has a defect in the sense that $e_{K'/K}=f_{K'/K}=1$.
In other words, the DVR's $A'$ and $A$ of $K'$ and $K$, have the same residue field
and satisfy $m_{A'}=m_AA'$. Since $A'$ is $A$-flat, it cannot be $A$-finite. On the
other hand, $A'$ is the integral closure of $A$ in $K'$, and we obtain that $A$ is
not qe. In addition, although $X=\Spec(A)$ is regular, any $X$-finite integral scheme
$X'$ with $K'\subseteq k(X')$ possesses a non-finite normalization and hence does not
admit a desingularization. Thus, $X$ is not universally $\QQ$-resolvable.
\end{remark2}

Our main goal will be to show that universal $l'$-resolvability of a qe base scheme
$S$ is inherited by finite type $S$-schemes whose structure morphism $X\to S$ is
maximally dominating (see Theorem \ref{10-l'theor} below, where a more precise
result is formulated). The proof will be by induction on the relative dimension,
and the main work is done when dealing with the case of generically smooth relative curves.

\begin{theorem}\label{10-relcurveth}
Let $S$ be an integral, noetherian, qe scheme with generic point $\eta=\Spec(K)$,
let $f{\colon}X\to S$ be a maximally dominating (\refext{II}{2-maximalement-dominant}) morphism
of finite type, and let $Z\subset X$ be a nowhere dense closed subset. Assume that $S$ is universally
$l'$-resolvable (resp. universally separably $l'$-resolvable), $X_\eta=X\times_S\eta$
is a smooth curve over $K$, and $Z_\eta=Z\times_S\eta$ is étale over $K$. Then there
exist a projective $l'$-alteration (resp. a separable projective $l'$-alteration)
$a{\colon}S'\to S$, a projective modification $b{\colon}X'\to (X\times_SS')^{\rm pr}$,
where $(X\times_SS')^{\rm pr}$ is the proper transform of $X$, i.e. the schematic closure
of $X_\eta\times_SS'$ in $X\times_SS'$,
\begin{equation*}
\xymatrix{X'\ar[r]^b\ar[rrd]_{f'} & (X\times_SS')^{\rm pr}\ar@{^{(}->}[r] & X\times_SS' \ar[d]\ar[r]& X\ar[d]^f\\
& & S' \ar[r]^a & S}
\end{equation*}
and divisors $W'\subset S'$ and $Z'\subset X'$ such that $S'$ and $X'$ are regular,
$W'$ and $Z'$ are snc, the morphism $f'{\colon}X'\to S'$ is pseudo-projective
(\S\ref{10-chowsec}), $(X',Z')\to(S',W')$ is log smooth, and $Z'=c^{-1}(Z)\cup f'^{-1}(W')$,
where $c$ denotes the alteration $X'\to X$.
\end{theorem}
We also note if $f$ is separated (resp. proper) then $f'$ is even quasi-projective
(resp. projective) by Remark \ref{10-chowrem}(ii).
\begin{proof}
It will be convenient to represent $Z$ as $Z_h\cup Z_v$, where the horizontal component
$Z_h$ is the closure of $Z_\eta$ and the vertical component $Z_v$ is the closure of
$Z\setminus Z_h$. The following observation will be used freely: if $a_1{\colon}S_1\to S$
is a (resp. separable) projective $l'$-alteration with an integral $S_1$ and
$b_1{\colon}X_1\to(X\times_SS_1)^{\rm pr}$ is a projective modification, then
it suffices to prove the theorem for $f_1{\colon}X_1\to S_1$ and the preimage
$Z_1\subset X_1$ of $Z$ (note that the generic fiber of $f_1$ is smooth because
it is a base change of that of $f$). So, in such a situation we can freely replace
$f$ by $f_1$, and $Z$ will be updated automatically without mentioning, as a rule.
We will change $S$ and $X$ a few times during the proof. We start with some preliminary steps.

Step 1. {\it We can assume that $f$ is quasi-projective.} By Proposition \ref{10-chowlem},
replacing $X$ with its projective modification we can achieve that $f$ factors through
a local isomorphism $X\to\oX$, where $\oX$ is $S$-projective. Let $X_1\subseteq\oX$ be
the image of $X$. Then the induced morphism $X_1\to S$ is quasi-projective and with
smooth generic fiber. If the theorem holds for $f_1$ and the image $Z_1\subset X_1$
of $Z$, i.e., there exist $a{\colon}S'\to S$ and $b_1{\colon}X'_1\to(X_1\times_SS')^{\rm pr}$
that satisfy all assertions of the theorem, then the theorem also holds for $f$ and $Z$ because
we can keep the same $a$ and take $b=b_1\times_{X_1}X$. This completes the step, and in the
sequel we assume that $f$ is quasi-projective. As we will only use projective modifications
$b_1{\colon}X_1\to(X\times_SS_1)^{\rm pr}$, the quasi-projectivity of $f$ will be preserved
automatically.

Step 2. {\it We can assume that $f$ and $Z_h\to S$ are flat.} Indeed, due to the flattening
theorem of Raynaud-Gruson, see \cite[5.7.9]{Platification@Raynaud-Gruson}, this can be achieved
by replacing $S$ with an appropriate projective modification $S'$, replacing $X$ with the proper
transform, and replacing $Z$ with its preimage. From now on, the proper transforms of $X$ will
coincide with the base changes.

Step 3. {\it Use of the stable modification theorem.} By the stable modification theorem
\cite[1.5 and 1.1]{Temkin@stable} there exist a separable alteration $\oa{\colon}\oS\to S$
with an integral $\oS$ and a projective modification $\oX\to X\times_S\oS$ such that
$(\oX,\oZ_h)$ is a {\em semistable multipointed} $\oS$-curve (see \S\ref{10-relcurvesec}),
where $\oZ_h\subset\oX$ is the horizontal part of the preimage $\oZ$ of $Z$. Enlarging
$\oS$ we can assume that it is integral and normal.
\begin{equation*}
\xymatrix{\oX \ar[rd]_\of\ar[r] & X\times_S\oS \ar[d]\ar[r]& X\ar[d]^f\\
& \oS \ar[r]^\oa & S}
\end{equation*}

Step 4. {\it We can assume that $\oa$ is a separable projective $G$-alteration, where $G$
is an $l$-group.} Since semistable multipointed relative curves are preserved by base changes,
we can just enlarge $\oS$ by replacing it with any separable projective Galois alteration that
factors through $\oS$. Once $\oS\to S$ is Galois, let $\oG$ denote its Galois group and let
$G\subseteq\oG$ be any Sylow $l$-subgroup. Since $\oS\to\oS/G$ is a separable $G$-alteration
and $\oS/G\to S$ is a separable projective $l'$-alteration, we can replace $S$ with $\oS/G$,
replace $X$ with $X\times_S(\oS/G)$, and update $Z$ accordingly, accomplishing the step.

Step 5. {\it The action of $G$ on $X\times_S\oS$ via $\oS$ lifts equivariantly to $\oX$.
In particular, $\of$ becomes $G$-equivariant and $\oX\to X$ becomes a separable projective
$G$-alteration.} This follows from \cite[1.6]{Temkin@stable}.

Step 6. {\it The action of $G$ on $\oX$ is inertia specializing.} Indeed, any covering of
$X$ by separated open subschemes induces a covering of $\oX$ by $G$-equivariant separated
open subschemes. So, it remains to use Remark \ref{10-inspecrem}(i).

Step 7. {\it We can assume that $\oS\to S$ is finite.} By Raynaud-Gruson there exists
a projective modification $S'\to S$ such that the proper transform $\oS'$ of $\oS$ is
flat over $S'$. Let $\eta$ denote the generic point of $S$ and $S'$ and let $\eta'$
denote the generic point of $\oS$ and $\oS'$. Since the morphisms $\oS\times_SS'\to S'$
and $\eta'\to\eta$ are $G$-equivariant and $\oS'$ is the schematic closure of $\eta'$ in
$\oS\times_SS'$, we obtain that the morphism $\oS'\to S'$ is a separable projective $G$-alteration.
Replacing $\oS\to S$ with $\oS'\to S'$, and updating $X$ and $\oX$ accordingly, we achieve that
$\oS\to S$ becomes flat, and hence finite. All conditions of steps 1--6 are preserved with the
only exception that $\oS$ may be non-normal. So, we replace $\oS$ with its normalization and
update $\oX$. This operation preserves the finiteness of $\oS\to S$, so we complete the step.

Step 8. {\it Choice of $W$.} Fix a closed subset $W\subsetneq S$ such that $\oS\to S$ is étale
over $S\setminus W$, $\of(\oZ_v)\subseteq\oW$, where $\oZ_v$ is the vertical part of $\oZ$ and
$\oW=\oa^{-1}(W)$, and $\of$ is smooth over $\oS\setminus\oW$.

Step 9. {\it We can assume that $S$ is regular and $W$ is snc.} Indeed, by our assumptions on
$S$ there exists a projective $l'$-alteration (resp. a separable projective $l'$-alteration)
$a{\colon}S'\to S$ such that $S'$ is regular and $a^{-1}(W)$ is snc. Choose any preimage of
$\eta$ in $S'\times_S\oS$ and let $\oS'$ be the normalization of its closure. Then $\oS'\to S'$
is a separable projective Galois alteration with Galois group $G'\subseteq G$, so we can replace
$S,\oS, G$ and $X$ with $S',\oS', G'$ and $X\times_SS'$, respectively, and update $W$, $\oW$ and
$Z$ accordingly (i.e. replace them with their preimages). Note that step 9 is the only step where
a non-separable alteration of $S$ may occur.

Step 10. {\it The morphism $(\oS,\oW)\to(S,W)$ is Kummer étale.} Indeed, $\oS\to S$ is
an étale $G$-covering outside of $W$, and $\oS$ is the normalization of $S$ in this
covering, so the assertion follows from \refext{IX}{9-FK}.

Consider the $G$-equivariant subscheme $\oT=\oZ\cup\of^{-1}(\oW)$ of $\oX$. The morphism
$(\oX,\oT)\to(\oS,\oW)$ is log smooth by Proposition \ref{10-logsmprop}, hence so is the
composition $(\oX,\oT)\to(S,W)$ and we obtain that $(\oX,\oT)$ is log regular. The group
$G$ acts freely on $\oS\setminus\oW$ and hence also on $\oX\setminus\oT$. Also, its action
on $\oX$ is tame and inertia specializing (step 6), hence we can apply Theorem \ref{10-modif-logsmooth}
to $(\oX,\oT)\to(S,W)$. As a result, we obtain a projective $G$-equivariant modification
$(\oX',\oT')\to(\oX,\oT)$ such that $\oT'$ is the preimage of $\oT$, $G$ acts very tamely
on $(\oX',\oT')$, and $(X',Z')=(\oX'/G,\oT'/G)$ is log smooth over $(S,W)$ (the quotient
exists as a scheme as $f$ is quasi-projective by Step 1 and the morphisms $\oS\to S$ and
$\oX'\to\oX\to X\times_S\oS$ are projective). Clearly, $X'$ is a projective modification
of $X$ and $Z'$ is the union of the preimages of $W$ and $Z$, hence it only remains to
achieve that $X'$ is regular and $Z'$ is snc. For this it suffices to replace $(X',Z')$
with its projective modification $\tilcalF^\rmlog(X',Z')$ introduced in  \refext{VIII}{8-deslogregth}.
\end{proof}

\begin{remark2}
It is natural to compare Theorem \ref{10-relcurveth} and the classical de Jong's result
recalled in \refext{IX}{9-dJ}. The main differences are as follows.

(i) One considers non-proper relative curves in \ref{10-relcurveth}, and this is the only
point that requires to use the stable modification theorem instead of de Jong's result.
The reason is that although the problem easily reduces to the case of a quasi-projective
$f$ (see Step 2), one cannot make $f$ projective, as the compactified generic fiber
$\oX_\eta$ does not have to be smooth (i.e. {\em geometrically} regular) at the added points.

(ii) One uses $l'$-alterations in \ref{10-relcurveth}. This is more restrictive than
in \refext{IX}{9-dJ}, but the price one has to pay is that the obtained log smooth
morphism $(X',Z')\to(S',W')$ does not have to be a nodal curve (e.g. $X'\to S'$ may
have non-reduced fibers). The construction of such $b{\colon}X'\to X$ involves a
quotient by a Sylow subgroup, and is based on Theorem \ref{10-modif-logsmooth}.
(Note also that it seems probable that instead of \ref{10-modif-logsmooth} one
could use a torification argument of Abramovich-de Jong, see \cite[\S1.4.2]{AdJ@Resolution}.)
\end{remark2}

Now, we are going to use Theorem \ref{10-relcurveth} to prove the main result of \S\ref{10-complsec}.

\begin{theorem}\label{10-l'theor}
Let $f{\colon}X\to S$ be a maximally dominating (\refext{II}{2-maximalement-dominant})
morphism of finite type between noetherian qe schemes, let $Z\subset X$ be a nowhere
dense closed subset, and assume that $S$ is universally $l'$-resolvable, then:

(i) $X$ is universally $l'$-resolvable.

(ii) There exist projective $l'$-alterations $a{\colon}S'\to S$ and $b{\colon}X'\to X$
with regular sources, a pseudo-projective (\S\ref{10-chowsec}) morphism $f'{\colon}X'\to S'$
compatible with $f$
\begin{equation*}
\xymatrix{X' \ar[d]^{f'}\ar[r]^b & X \ar[d]^f\\
S' \ar[r]^a & S&}
\end{equation*}
and snc divisors $W'\subset S'$ and $Z'\subset X'$ such that $Z'=b^{-1}(Z)\cup f'^{-1}(W')$
and the morphism $(X',Z')\to(S',W')$ is log smooth.

(iii) If $S=\Spec(k)$, where $k$ is a perfect field, then one can achieve in addition to (ii)
that $a$ is an isomorphism and the alteration $b$ is separable. In particular, $X$ is universally
separably $l'$-resolvable in this case.
\end{theorem}
\begin{proof}
Note that (i) follows from (ii) because any alteration $X_1$ of $X$ is also of finite type
over $S$, so we can apply (ii) to $X_1$ as well. Thus, our aim is to prove (ii) and its
complement (iii). We will view $Z$ both as a closed subset and a reduced closed subscheme.
We start with a few preliminary steps, that reduce the theorem to a special case. We will
tacitly use that if $S_1\to S$ and $X_1\to X$ are projective $l'$-alterations, separable in
case (iii), and $f_1{\colon}X_1\to S_1$ is compatible with $f$, then it suffices to prove the
theorem for $f_1$ and the preimage $Z_1\subset X_1$ of $Z$. So, in such situation we can freely
replace $f$ with $f_1$, and $Z$ will be updated automatically without mentioning, as a rule.

Step 1. {\it We can assume that $X$ and $S$ are integral and normal.} For a noetherian scheme
$Y$ let $\Nor(Y)$ denote the normalization of its reduction. Since $f$ is maximally dominating,
it induces a morphism $\Nor(f){\colon}\Nor(X)\to\Nor(S)$, and replacing $f$ with $\Nor(f)$ we can
assume that $S$ and $X$ are normal. Since we can work separately with the connected components,
we can now assume that $S$ and $X$ are integral.

Step 2. {\it We can assume that $f$ is projective.} By Proposition \ref{10-chowlem} there exists
a projective modification $X_1\to X$ such that the morphism $X_1\to S$ factors into a composition
of a local isomorphism $X_1\to\oX$ and a projective morphism $\of{\colon}\oX\to S$. Replacing $X$
with $X_1$ we can assume that $X$ itself admits a local isomorphism $g{\colon}X\to\oX$ with an
$S$-projective target. Let $\oZ$ be the closure of $g(Z)$. Then it suffices to solve our problem
for $\of$ and $\oZ$, as the corresponding alteration of $\oX$ will induce an alteration of $X$
as required. Thus, replacing $X$ and $Z$ with $\oX$ and $\oZ$, we can assume that $f$ is projective.

Step 3. {\it It suffices to find $f'$ which satisfies all assertions of the theorem except the
formula for $Z'$, while the latter is weakened as $b^{-1}(Z)\cup f'^{-1}(W')\subseteq Z'$.} Given
such an $f'$ note that $Z''=b^{-1}(Z)\cup f'^{-1}(W')$ is a subdivisor of $Z'$, hence it is an snc
divisor too. We claim that $X',Z'',S',W'$ satisfy all assertions of the theorem, and the only thing
one has to check is that the morphism $(X',Z'')\to(S',W')$ is log smooth. The latter follows from
Lemma \ref{10-logsmlem} whose proof will be given below.

Step 4. {\it In the situation of (iii) we can assume that the field $k$ is infinite.} Assume that
$S=\Spec(k)$ where $k$ is a finite field and fix an infinite algebraic $l'$-prime extension $\ok/k$
(i.e. it does not contain the extension of $k$ of degree $l$). We claim that it suffices to prove
(ii) and (iii) for $\oS=\Spec(\ok)$ and the base changes $\oX=X\times_S\oS$ and $\oZ=Z\times_S\oS$.
Indeed, assume that $a{\colon}\oX'\to\oX$ is a separable $l'$-alteration with a regular source and
such that $\oZ'=a^{-1}(\oZ)$ is an snc divisor (obviously, we can take $\oS'=\oS$ and $\oW'=\emptyset$).
Since $\oS=\lim_i S_i$ where $S_i=\Spec(k_i)$ and $k_i/k$ run over finite subextensions of $\ok/k$,
\cite[8.8.2(ii)]{EGAIV-3} implies that there exists $i$ and a finite type morphism
$X'_i\to X_i=X\times_SS_i$ such that $\oX'\iso X'_i\times_{S_i}\oS$. For any finite
subextension $k_i\subseteq k_j\subset\ok$ set $X'_j=X'_i\times_{S_i}S_j$ and $X_j=X_i\times_{S_i}S_j$,
and let $Z'_j\subset X'_j$ be the preimage of $Z$. Then it follows easily from \cite[8.10.5]{EGAIV-3}
and \cite[17.7.8]{EGAIV-4} that $X'_j\to X_j$ is an $l'$-alteration and $X'_j\to S_j$ is smooth for
large enough $k_j$. In the same manner one achieves that $Z'_j$ is an snc divisor. Now, it is obvious
that $(X'_j,Z'_j)\to(S,\emptyset)$ is log smooth and $X'_j\to X_j\to X$ is an $l'$-alteration.

Now we are in a position to prove the theorem. We will use induction on $d=\trdeg(k(X)/k(S))$, with
the case of $d=0$ being obvious. Assume that $d\ge 1$ and the theorem holds for smaller values of $d$.

Step 5. {\it Factorizing $f$ through a relative curve.} After replacing $X$ by a projective modification
we can factor $f$ through an integral scheme $Y$ such that $g{\colon}X\to Y$ is maximally dominating,
$h{\colon}Y\to S$ is projective, and $\trdeg(k(X)/k(Y))=1$. Indeed, one can obviously construct such
a rational map $h'{\colon}X\dasharrow Y$ even without modifying $X$ (i.e. $h'$ is well defined on a
non-empty open subscheme $U\subseteq X$). Let $X'$ be the schematic image of the morphism $U\into X\times Y$.
Then $X'\to X$ is a projective modification (an isomorphism over $U$), and the morphism $X'\to S$
factors through $Y$.

Let $\eta=\Spec(k(Y))$ denote the generic point of $Y$, $X_\eta=X\times_{Y}\eta$ and
$Z_\eta=Z\times_Y\eta$. If $k$ is an infinite perfect field and $S=\Spec(k)$ then
a Bertini type argument from the end of proof of \cite[4.11]{Alterations@de_Jong}
implies that we can choose $Y$ so that $X_\eta$ and $Z_\eta$ are smooth over $\eta$.
(One assumes that $k$ is algebraically closed in \cite{Alterations@de_Jong}, but to
use Bertini's theorem one only needs that $k$ is perfect and infinite.) In the general
case, we can achieve that $X_\eta$ and $Z_\eta$ are smooth at cost of a purely inseparable
alteration of $Y$. Indeed, pick up any finite purely inseparable extension $K/k(Y)$ such
that $Z_K=\Nor(Z\times_\eta\Spec(K))$ (i.e. just the reduction) and $X_K=\Nor(X\times_\eta\Spec(K))$
are smooth, extend $K$ to a projective alteration $Y'\to Y$, and replace $Y$ and $X$ with $Y'$
and the schematic closure of $X_K$ in $X\times_YY'$, respectively.

Step 6. {\it Use of Theorem \ref{10-relcurveth}.} So far, we have constructed the right column
of the following diagram
$$
\xymatrix{(X',Z') \ar[rrdd]_{f'}\ar[rr]^{\tilcalF^\rmlog(L,M_L)} & &
(L,M_L) \ar[r]\ar[d]^{g'}\ar@{}[rd]|\square & (X'',Z'') \ar[d]^{g''}\ar[r]^{b''} & X \ar[d]^g\ar@/^1pc/[dd]^f\\
 & & (Y',V')\ar[r]^{c'}\ar[d]^{h'} & (Y'',V'')\ar[r]^{c''}\ar[d]^{h''} & Y\ar[d]^h\\
 & & (S',W')\ar[r]^a & S\ar@{=}[r] & S
}
$$
By Theorem \ref{10-relcurveth} there exists a projective $l'$-alteration $c''{\colon}Y''\to Y$
with regular source, a projective modification $X''\to(X\times_YY'')^{\rm pr}$ with regular
source, a projective morphism $g''{\colon}X''\to Y''$ compatible with $g$, and snc divisors
$V''\subset Y''$ and $Z''\subset X''$ such that $(X'',Z'')\to(Y'',V'')$ is log smooth and
$b''^{-1}(Z)\subseteq Z''$. In case (iii), $Y$ is universally separably $l'$-resolvable by
the induction assumption, hence we can take $c''$ to be separable, and then $b''{\colon}X''\to X$
is also separable. In addition, by the induction assumption applied to $h''{\colon}Y''\to S'$ and
$V''\subset Y''$ there exist projective $l'$-alterations $a{\colon}S'\to S$ and
$c'{\colon}Y'\to Y''$ with regular sources and snc divisors $W'\subset S'$ and
$V'\subset Y'$ and a projective morphism $h'{\colon}Y'\to S'$ compatible with $h''$
such that $(Y',V')\to(S',W')$ is log smooth, $c'^{-1}(V'')\subseteq V'$, and $c'$ is
separable if the assumption of (iii) is satisfied.

Set $(L,M_L)=(X'',Z'')\times^{\rm fs}_{(Y'',V'')}(Y',V')$, where the product is
taking place in the category of fs log schemes. To simplify notation we will write
$P^\sat$ instead of $(P^{\rm int})^\sat$ for monoids and log schemes. Recall that
$(L,M_L)=(F,M_F)^\sat$, where $(F,M_F)$ is the usual log fibered product, and
$F=X''\times_{Y''}Y'$ by \cite[1.6]{LSFI@Kato}. Furthermore, we have local Zariski
charts for $c'$ and $g''$ modeled, say, on $P_i\to P'_i$ and $P_i\to Q_i$. Hence
$(F,M_F)$ is a Zariski log scheme with charts modeled on $R_i=P'_i\oplus_{P_i}Q_i$,
and $(L,M_L)$ is a Zariski log scheme with charts modeled on $R_i^\sat$. Furthermore,
the saturation morphism $L\to F$ is finite hence the composition $L\to F\to X''$ is
projective. The morphism $g'{\colon}(L,M_L)\to (Y',V')$ is a saturated base change of
the log smooth morphism $g''{\colon}(X'',Z'')\to(Y'',V'')$, hence it is log smooth.
As $(Y',V')$ is log regular, $(L,M_L)$ is also log regular. Applying to $(L,M_L)$ the
saturated monoidal desingularization functor $\tilcalF^\rmlog$ from \refext{VIII}{8-deslogregth}
we obtain a log regular Zariski log scheme $(X',Z')$ with a regular $X'$. Then $Z'$ is a normal
crossings divisor, which is even an snc divisor since the log structure is Zariski.

We claim that $(X',Z')$ and $(S',W')$ are as asserted by the theorem except of the
weakening dealt with in Step 3. Indeed, the morphism $(X',Z')\to(S',W')$ is log
smooth because it is the composition $(X',Z')\to (L,M_L)\to(Y',V')\to(S',W')$ of
log smooth morphisms. The preimage of $Z$ in $X''$ is contained in $Z''$, which is
the non-triviality locus of the log structure of $(X'',Z'')$, hence its preimage in
$X'$ is also contained in the non-triviality locus of the log structure of $(X',Z')$,
which is $Z'$. Clearly, $Z'$ also contains the preimage of $W'$. By the construction,
$S'\to S$ is a projective $l'$-alteration, and it remains to check that $X'\to X$ is
also a projective $l'$-alteration. Since $\tilcalF^\rmlog(L,M_L)$ is a saturated log
blow up tower and $(L,M_L)$ is log regular, the underlying morphism of schemes $X'\to L$
is a projective modification by \refext{VIII}{8-logregblowlem}(i). The projective morphism
$L\to X''$ is an $l'$-alteration because generically (where the log structures are trivial)
it is a base change of the projective $l'$-alteration $Y'\to Y''$. And $X''\to X$ is
a projective $l'$-alteration by the construction. It remains to recall that in the
situation of (iii) the alterations $c'{\colon}Y'\to Y''$ and $b''{\colon}X''\to X$
are separable, hence so are $(L,M_L)\to X''$ and the total composition $X'\to X$.
\end{proof}

\begin{remark2}
The only place where inseparable alterations are used is the argument at step 5,
where we had to choose an inseparable extension $K/k(Y)$ when $X_\eta$ or $Z_\eta$
is not geometrically regular.
\end{remark2}

\begin{lemma2}\label{10-logsmlem}
Assume that $S$ and $X$ are regular schemes, $W\subset S$ and $Z\subset X$ are snc
divisors, and $f:X\to S$ is a morphism such that $f^{-1}(W)\subseteq Z$ and the
induced morphism of log schemes $h{\colon}(X,Z)\to(S,W)$ is log smooth. Then for
any intermediate divisor $f^{-1}(W)\subseteq Z'\subseteq Z$ the morphism $h'{\colon}(X,Z')\to(S,W)$
is log smooth.
\end{lemma2}
\begin{proof}
We can work locally at a geometric point $\ox\to X$. Let $x\in X$ and $s\in S$
be the images of $\ox$, and let $q_1\. q_r\in\calO_{S,s}$ define the irreducible
components of $W$ through $s$. Set $Q=\oplus_{i=1}^r q_i^\NN$. Shrinking $S$ we
obtain a chart $c{\colon}S\to\Spec(\ZZ[Q])$ of $(S,W)$. By Proposition \ref{10-local-structure}
applied to $c$, $h$, and $G=1$, after localizing $X$ along $\ox$ one can find an fs
chart of $h$ consisting of $c$, $X\to\Spec(\ZZ[P])$, and $\phi{\colon}Q\to P$ such
that the morphism $X\to S\times_{\Spec(\ZZ[Q])}\Spec(\ZZ[P])$ is smooth, $P^*$ is
torsion free, $\phi$ is injective, and the torsion of $\Coker(\phi^\gp)$ is annihilated
by an integer $n$ invertible on $S$.

Let $p_1\. p_t\in\calO_{X,x}$ define the irreducible components of $Z$ through $x$.
Our next aim is to adjust the chart similarly to \ref{10-chart-diagram}(vi) to achieve
that $\oP\iso N=\oplus_{j=1}^t p_j^\NN$. Note that $\oM_{X,x}\iso N$, where
$M_X\into\calO_X$ is the log structure of $(X,Z)$. The homomorphism $\psi{\colon}P\to M_{X,x}$
factors through the fs monoid $R=P[T^{-1}]$ where $T=\psi^{-1}(M^*_{X,x})$. Clearly,
$\R^*$ is torsion free, $\oR\iso N$, and shrinking $X$ around $x$ we obtain a chart
$X\to\Spec(\ZZ[R])$. Since $\Spec(\ZZ[R])$ is open in $\Spec(\ZZ[P])$ the morphism
$X\to S\times_{\Spec(\ZZ[Q])}\Spec(\ZZ[R])$ is smooth. So, we can replace $P$ with
$R$ achieving that $\oP\iso N$, and hence $P\iso\NN^t\oplus\ZZ^u$.

Without restriction of generality $Z'$ is defined by the vanishing of $\prod_{j=1}^{t'}p_j$
for $t'\le t$. Since $f^{-1}(W)\subseteq Z'$, the image of $Q$ in $\oP$ is contained in
$\oP'=\oplus_{j=1}^{t'}p_j^\NN$. Hence $\phi$ factors through a homomorphism
$\phi'{\colon}Q\to P'=\oP'\oplus P^*$, and we obtain a chart of $h'$ consisting
of $c$, $X\to\Spec(\ZZ[P'])$, and $\phi'$. By \cite[3.5,3.6]{Toric@Kato},
to prove that $h'$ is log smooth it remains to observe that $\phi'$ is injective,
the torsion of $\Coker(\phi'^\gp)$ is annihilated by $n$, and the morphism
$$X\to S\times_{\Spec(\ZZ[Q])}\Spec(\ZZ[P])\to S\times_{\Spec(\ZZ[Q])}\Spec(\ZZ[P'])$$
is smooth because $\Spec(\ZZ[P])\to\Spec(\ZZ[P'])$ is so.
\end{proof}

\subsubsection{Comparison with Theorems \ref{10-de Jong1} and \ref{10-de Jong2}}
Theorems \ref{10-de Jong1} and \ref{10-de Jong2} follow by applying Theorem \ref{10-l'theor}(ii)
to $X\to S$ and $Z$ (where one takes $S=\Spec(k)$ in \ref{10-de Jong1}). Indeed,
the main part of the proofs of \ref{10-de Jong1} and \ref{10-de Jong2} was to
construct $l'$-alterations $X'\to X$ and $S'\to S$ with regular sources, snc
divisors $Z'\subset X'$ and $W'\subset S'$, and a log smooth morphism
$f'{\colon}(X',Z')\to(S',W')$ compatible with $f$. Then, in the last paragraphs
of both proofs, Proposition \ref{10-local-structure} was used to obtain a more
detailed description of $X'$, $Z'$, and $f'$. In particular, for a zero-dimensional
base this amounted to saying that $X'$ is $S$-smooth and $Z'$ is relatively snc over
$S$, and for a one-dimensional base this amounted to the conditions (i) and (ii)
of \ref{10-de Jong2}.

Conversely, Theorem \ref{10-de Jong1} (resp. \ref{10-de Jong2}) implies assertion (ii)
of Theorem \ref{10-l'theor} under the assumptions of \ref{10-de Jong1} (resp. \ref{10-de Jong2})
on $X$ and $S$. Moreover, the non-separated Chow's lemma could be used in their proofs as well,
so the separatedness assumption there could be easily removed. In such case,
Theorems \ref{10-de Jong1} and \ref{10-de Jong2} would simply become the low
dimensional (with respect to $S$) cases of Theorem  \ref{10-l'theor}(ii) plus
an explicit local description of the log smooth morphism $f'$. The strengthening
\ref{10-l'theor}(iii), however, was not achieved in \ref{10-de Jong1}, and required
a different proof of the whole theorem.

\subsection{Saturation}\label{10-satursec} In Theorems \ref{10-relcurveth} and \ref{10-l'theor}
we resolve certain morphisms $f{\colon}X\to S$ with divisors $Z\subset X$ by log smooth morphisms
$f'{\colon}(X',Z')\to(S',W')$. However, as we insisted to use only $l'$-alterations and to obtain
regular $X'$ and snc $Z'$, we had to compromise a little on the "quality" of $f'$. For example,
our $f'$ may have non-reduced fibers. Due to de Jong's theorem, if the relative dimension is one,
then one can make $f'$ a nodal curve. We will see that a similar improvement of $f'$ is possible
in general if one uses arbitrary alterations and allows non-regular $X'$. The procedure reduces
to saturating $f'$ and is essentially due to Tsuji and Illusie-Kato-Nakayama.

\subsubsection{Saturated morphisms} Recall that a homomorphism $P\to Q$ of fs (resp. fine)
monoids is {\em saturated} (resp. {\em integral}) if for any homomorphism $P\to P'$ with
fs (resp. fine) target the pushout $Q\oplus_{P}P'$ is fs (resp. fine). A morphism of fs (resp. fine)
log schemes $f{\colon}(Y,M_Y)\to(X,M_X)$ is {\em saturated} (resp. {\em integral}) if so are
the homomorphisms $\oM_{X,f(y)}\to\oM_{Y,y}$.

\begin{remark2}\label{10-satrem}
(i) Integral morphism were introduced already by Kato in \cite[\S4]{LSFI@Kato}. Kato also
introduced the notion of saturated morphisms, which was first seriously studied by Tsuji
in \cite{Saturated@Tsuji}. Actually, one can define saturated morphisms for arbitrary fine
log schemes, but the definition is more involved than we use. For fs log schemes our definition
coincides with the usual one due to \cite[II 2.13(2)]{Saturated@Tsuji}.

(ii) The following two basic properties of saturated morphisms follow from the definition:
(a) a composition of saturated morphisms between fs log schemes is saturated,
(b) if $f{\colon}Y\to X$ is a saturated morphism between fs log schemes, then,
for any morphism of fs log schemes $Y'\to Y$, the base change $f'{\colon} Y'\to X'$ of $f$
in the category of log schemes is a saturated morphism of fs log schemes. (Also, it is proved
in \cite[II 2.11]{Saturated@Tsuji} that analogous properties hold for saturated morphisms between
arbitrary integral log schemes.)

(iii) Let $f{\colon}Y\to X$ be a morphism of fs log schemes. It is shown in \cite[II 3.5]{Saturated@Tsuji}
that if $f$ can be modeled on charts corresponding to saturated homomorphisms of fs monoids $P_i\to Q_i$
then $f$ is saturated. Let us remark that the converse is also true: if $f$ is saturated then it can be
modeled on charts corresponding to $P_i\to Q_i$ as above.
\end{remark2}

\subsubsection{Integrality and saturatedness for log smooth morphisms}
 We recall the following result that relates the notions of integral and saturated
 morphisms to certain properties of the underlying morphisms of schemes.

\begin{proposition2}\label{10-saturprop1}
Let $f{\colon}(Y,M_Y)\to(X,M_X)$ be a log smooth morphism between fs log schemes.

(i) If $f$ is integral then $Y\to X$ is flat.

(ii) Assume that $f$ is integral. Then $f$ saturated if and only if $Y\to X$ has reduced fibers.
\end{proposition2}
\begin{proof}
The first claim is proved in \cite[4.5]{LSFI@Kato} and the second one is proved in \cite[II 4.2]{Saturated@Tsuji}.
\end{proof}

Sometimes one can also go in the opposite direction: from flatness to integrality.

\begin{proposition2}\label{10-saturprop2}
Let $f{\colon}(Y,M_Y)\to(X,M_X)$ be a log smooth morphism between fs log schemes and
 assume that the morphism $Y\to X$ is flat, then $f$ is integral.
\end{proposition2}
\begin{proof}
It suffices to show that if $\oy\to Y$ is a geometric point and $\ox=f(\oy)$ then
the homomorphism $\ophi{\colon}\oM_{X,\ox}\to\oM_{Y,\oy}$ is integral. By Proposition
\ref{10-local-structure} and the argument in \ref{10-preliminary}(vi), localizing $X$
and $Y$ along these points we can assume that $f$ is modeled on a chart
$Y_0=\Spec(\ZZ[P])\to X_0=\Spec(\ZZ[Q])$ corresponding to a homomorphism
$\phi{\colon}Q\to P$ so that the morphism $g{\colon}Y\to Z=X\times_{X_0}Y_0$
is smooth (in particular, flat), and $\phi$ has the following properties:
$Q =\oM_{X,\ox}$, $P$ is fs, $P^*$ is torsion free, the composition
$Q\to P \to\oP = P/P^*$ coincides with $\ophi$, the kernel of $\phi^{\gp}$ is finite,
killed by an integer invertible at $x$, as well as the torsion part of its cokernel
(but we will not need these last two properties). Since $Q$ is sharp and saturated,
$Q^{\gp}$ is torsion free, so $\phi$ is injective. We claim that $\ophi$ is is integral
if and only if $\phi$ is integral. To see this note that if $Q\to R$ is a homomorphism
of monoids, then $R\oplus_Q\oP$ is isomorphic to the quotient of $R\oplus_Q P$ by the
image of $P^*$, and hence either both pushouts are integral or neither of them is integral.
Thus, we need only prove that $\phi$ is integral.

Note that the morphism $h{\colon}Z\to X$ is flat at the (Zariski) image $z\in Z$ of $\oy$
because $f$ and $g$ are flat. Set $x=h(z)$ and $k=k(x)$, then the fiber
$h_x{\colon}\Spec(k[P])\to\Spec(k[Q])$ of $h$ over $x$ is flat at $z$.
In other words, if $I\subset k[P]$ is the ideal corresponding to $z$ then
the homomorphism $k[Q]\to k[P]_I$ is flat. The preimage of $m_y$ under
$\ZZ[P]\to\calO_{Y,y}$ contains the set $m_P=P\setminus P^*$, hence $m_P\subset I$
and we obtain that $I$ contains $J=k[m_P]$. Note that the ideal $J$ is prime as
$k[P]/J\iso k[P^*]$ is a domain due to $P^*$ being torsion free. Thus, the
localization $k[P]_J$ makes sense, and we obtain a flat homomorphism $\psi{\colon}k[Q]\to k[P]_J$.

It is proved in \cite[4.1]{LSFI@Kato}, that if the homomorphisms $K[\phi]{\colon}K[Q]\to K[P]$
are flat for any field $K$ then $\phi$ is integral. The proof consists of two parts. First one
checks that $\phi$ is injective, which is automatic in our case. This is the only argument in
loc.cit. where a play with different fields is needed. We claim that the second part of the
proof of the implication (iii)$\implies$(v) in \cite[4.1]{LSFI@Kato} works fine with a single
field $k$, and, moreover, it suffices to only use that $k[Q]\to k[P]_J$ is flat. Let us indicate
how the argument in loc.cit. should be adjusted.

Assume that, as in the proof of \cite[4.1]{LSFI@Kato}, we are given $a_1,a_2\in Q$ and
$b_1,b_2\in P$ such that $\phi(a_1)b_1=\phi(a_2)b_2$. Let $S$ be the kernel of the
homomorphism of $k[Q]$-modules $k[Q]\oplus k[Q]\to k[Q]$ given by $(x,y)\mapsto a_1x-a_2y$.
By the flatness, the kernel of $k[P]_J\oplus k[P]_J\to k[P]_J$, $(x,y)\mapsto\phi(a_1)x-\phi(a_2)y$
is generated by the image of $S$. Hence there exist representations
$b_1=\sum_{i=1}^r\phi(c_i)\frac{f_i}{s}$ and $b_2=\sum_{i=1}^r\phi(d_i)\frac{f_i}{s}$
with $c_i,d_i\in k[Q]$, $f_i\in k[P]$, $s\in k[P]\setminus J$, and $a_1c_i=a_2d_i$.
Moreover, multiplying $s$ and $f_i$'s by an appropriate unit $u\in P^*$ we can assume
that $s=1+s'$ for $s'\in{\rm Span}_k(P\setminus\{1\})$. Then
$b_1+\sum_{1\le\alpha\le m}\lambda_\alpha t_\alpha=\sum_{1\le i \le r}\phi(c_i)f_i$,
with $\lambda_\alpha\in k^*$, and the $t_\alpha\in P$ pairwise distinct and distinct
from $b_1$, so we see that there exist $a_3\in Q$, $b\in P$, and $1\le i\le r$, such
that $a_3$ appears in $c_i$, $b$ appears in $f_i$, and $b_1=\phi(a_3)b$. The remaining
argument copies that of the loc.cit. verbatim, and one obtains in the end that $\phi$
satisfies the condition (v) from \cite[4.1]{LSFI@Kato}. Thus, $\phi$ is integral and
we are done.
\end{proof}

Before going further, let us discuss an incarnation of saturated morphisms in (more classical)
toroidal geometry.

\begin{remark2}
In toroidal geometry an analog of saturated morphisms was introduced by Abramovich and Karu
in \cite{AK@Toroidal}. In the language of log schemes toroidal morphisms can be interpreted
as log smooth morphisms $f{\colon}(X,Z)\to(S,W)$ between log regular schemes (with the toroidal
structure given by the triviality loci of the log structures). If  $f$ is a toroidal morphism
as above then Abramovich-Karu called it weakly semistable when the following conditions hold:
$S$ is regular, $f$ is locally equidimensional, and the fibers of $f$ are reduced. Furthermore,
they remarked that the equidimensionality condition is equivalent to flatness of $f$ whenever
$S$ is regular, see \cite[4.6]{AK@Toroidal}. Thus, the weak semistability condition is nothing
else but saturatedness of $f$ and regularity of the target. In particular, saturated log smooth
morphisms between log regular log schemes may be viewed as the generalization of weakly semistable
morphisms to the case of an arbitrary log regular (or toroidal) base.
\end{remark2}

Now, we are going to prove the main result of \S\ref{10-satursec}.

\begin{theorem}\label{10-saturth}
Assume that $f{\colon}(X,Z)\to(S,W)$ is a log smooth morphism such that $(S,W)$ is log regular
and $S$ is universally $\QQ$-resolvable (\S\ref{10-resolvsec}). Then there exists an alteration
$h{\colon}S'\to S$ such that $S'$ is regular, $W'=g^{-1}(W)$ is an snc divisor, and the fs base
change $f'{\colon}(X',Z')\to(S',W')$ is a saturated morphism.
\end{theorem}
Recall that $(X',Z')=(X,Z)\times_{(S,W)}^{\rm fs}(S',W')$ and $f'$ is log smooth because the
saturation morphism is log smooth.
\begin{proof}
By \refext{VIII}{8-deslogregth}, applying to $(S,W)$ an appropriate saturated log blow up
tower and replacing $(X,Z)$ with the fs base change we can achieve that $S$ is regular and
$W$ is normal crossings. By an additional sequence of log blow ups we can also make $W$ snc
(see \refext{VIII}{8-step2}), so $(S,W)$ becomes a Zariski log scheme. Now, we can étale-locally
cover $f$ by charts $f_i{\colon}(X_i,Z_i)\to(S_i,W_i)$ modeled on $P_i\to Q_i$ such that $S_i$ are
open subschemes in $S$. By \cite[A.4.4, A.4.3]{IKN@Quasi-unipotent}, for each $i$ there exists
a morphism $h_i{\colon}(S'_i,W'_i)\to(S_i,W_i)$ such that $h_i$ is a composition of a Kummer
morphism and a log blow up, and the fs base change of $f_i$ is saturated. (Although the proof
in loc.cit. is written in the context of log analytic spaces, it translates to our situation
almost verbatim. The only changes are that we have to distinguish étale and Zariski topology
on the base (in order to construct log blow ups), and $h_i$ does not have to be log étale as
there might be inseparable Kummer morphisms.)

Note that $W'_i=h_i^{-1}(W_i)$. In addition, $S'_i\to S_i$ is a projective alteration by
\refext{VIII}{8-logregblowlem}. Extend each $h_i$ to a projective alteration $g_i{\colon}T_i\to S$,
and let $h{\colon}S'\to S$ be a projective alteration that factors through each $T_i$.
By the universal $\QQ$-resolvability assumption we can enlarge $S'$ so that it becomes
regular and $Z'=h^{-1}(Z)$ becomes snc. We claim that $h$ is as claimed. It suffices to
check that the fs base change of each morphism $f'_i{\colon}(X_i,Z_i)\to(S,W)$ is saturated.
However, already the fs base change of $f'_i$ to $(T_i,g_i^{-1}(W))$ is saturated by the construction,
hence so is its further base change to $(S',W')$.
\end{proof}

\begin{remark2}
Our proof is an easy consequence of \cite[A.4.4 and A.4.3]{IKN@Quasi-unipotent}. The first cited
result shows that (locally) any log smooth morphism can be made exact by an appropriate log blow
up of the base. This result is somewhat analogous to the flattening theorem of Raynaud-Gruson.
The second cited result shows that by a Kummer extension of the base one can (locally) saturate
an exact log smooth morphism. It is somewhat analogous to the reduced fiber theorem of Bosch-L\"utkebohmert-Raynaud
(\cite{BLR@Reduced_Fiber}) which implies that if $f{\colon}Y\to X$ is a finite type morphism between
reduced noetherian schemes then there exists an alteration $X'\to X$ such that the normalized base
change $f'{\colon}\Nor(Y\times_XX')\to X'$ has reduced fibers. Although the proof of the latter is
far more difficult.
\end{remark2}

\subsection{Characteristic zero case}\label{10-char0sec} Theorem \ref{10-l'theor} can be substantially
strengthened when $S$ is of characteristic zero, i.e., the morphism $S\to\Spec(\ZZ)$ factors through $\Spec(\QQ)$.

\begin{theorem}\label{10-char0theor}
Assume that $S$ is a reduced, noetherian, qe scheme of characteristic zero, $f{\colon}X\to S$
is a maximally dominating morphism of finite type with reduced source, and $Z\subset X$ is a
nowhere dense closed subset. Then there exist projective modifications $a{\colon}S'\to S$ and
$b{\colon}X'\to X$ with regular sources, a pseudo-projective morphism $f'{\colon}X'\to S'$
compatible with $f$, and snc divisors $W'\subset S'$ and $Z'\subset X'$ such that
$Z'=b^{-1}(Z)\cup f'^{-1}(W')$ and the morphism $(X',Z')\to(S',W')$ is log smooth.
\end{theorem}
\begin{proof}
The proof is very close to the proof of Theorem \ref{10-l'theor}, so we will just say which changes
in that proof should be made. First, we note that any $S$-scheme $Y$ of finite type is noetherian and qe.
Thus, if $Y$ is reduced and $T\subset Y$ is a nowhere dense closed subset then the pair $(Y,T)$
can be desingularized by \cite{Temkin@des1} in the following sense: there exists a projective
modification $h{\colon}Y'\to Y$ with regular source and such that $h^{-1}(T)$ is an snc divisor.
This result replaces the $l'$-resolvability assumption in Theorem \ref{10-l'theor}, and it allows
to apply the proof of that theorem to our situation with the only changes that one always uses
projective modifications instead of projective $l'$-alterations, and Theorem \ref{10-relcurveth}
is replaced with Lemma \ref{10-relcurvelem} below. (Note that Lemma \ref{10-relcurvelem}
is weaker than Theorem \ref{10-char0theor}, while Theorem \ref{10-relcurveth} does not follow
from Theorem \ref{10-l'theor}.)
\end{proof}

\begin{lemma2}\label{10-relcurvelem}
Let $S$ be an integral, noetherian, qe scheme with generic point $\eta=\Spec(K)$,
let $f{\colon}X\to S$ be a maximally dominating morphism of finite type, and let
$Z\subset X$ be a nowhere dense closed subset. Assume that $X_\eta=X\times_S\eta$
is a smooth curve over $K$, and $Z_\eta=Z\times_S\eta$ is étale over $K$. Then
there exist projective modifications $a{\colon}S'\to S$ and $b{\colon}X'\to X$
with regular sources, a pseudo-projective morphism $f'{\colon}X'\to S'$ compatible
with $f$ and snc divisors $W'\subset S'$ and $Z'\subset X'$ such that
$Z'=b^{-1}(Z)\cup f'^{-1}(W')$ and the morphism $(X',Z')\to(S',W')$ is log smooth.
\end{lemma2}
\begin{proof}
The proof copies the proof of Theorem \ref{10-relcurveth} with the only difference
that instead of an $l$-Sylow subgroup $G\subseteq\oG$ one simply takes $G=\oG$.
The latter is possible because the schemes are of characteristic zero and hence
any action of $\oG$ is tame.
\end{proof}

Combining Theorem \ref{10-char0theor} and \ref{10-saturth} we obtain the following
weak semistable reduction theorem.

\begin{theorem}\label{10-weaksemistable}
Assume that $S$ is a reduced, noetherian, qe scheme of characteristic zero,
$f{\colon}X\to S$ is a maximally dominating morphism of finite type with
reduced source, and $Z\subset X$ is a nowhere dense closed subset. Then
there exists an alteration $S'\to S$, a modification $X'\to X\times_SS'$
of the proper transform of $X$, a pseudo-projective morphism $f'{\colon}X'\to S'$
compatible with $f$, and divisors $W'\subset S'$ and $Z'\subset X'$ such that $S'$
is regular, $W'$ is snc, $Z'=b^{-1}(Z)\cup f'^{-1}(W')$, and the morphism $(X',Z')\to(S',W')$
is log smooth and saturated (i.e. $X'\to S'$ is weakly semistable).
\end{theorem}

\begin{remark2}
(i) In the case when $X$ and $S$ are integral proper varieties over an algebraically
closed field $k$ of characteristic zero, this theorem becomes the weak semistable
reduction theorem of Abramovich-Karu. Our proof has many common lines with their
arguments. In particular, the first step of their proof was to make $f$ toroidal,
and it was based on de Jong's theorem. (Note also that in a recent work \cite{ADK@Toroidal}
of Abramovich-Denef-Karu, the toroidalization theorem was extended to separated schemes
of finite type over an arbitrary ground field of characteristic zero.)
Our Theorem \ref{10-char0theor} can be viewed as a generalization
of the toroidalization theorem of Abramovich-Karu.

(ii) The second stage in the proof of the weak semistable reduction theorem of Abramovich-Karu
(the combinatorial stage) is analogous to Theorem \ref{10-saturth}. It obtains as an input
a toroidal morphism $f{\colon}(X,Z)\to(S,W)$ between proper varieties of characteristic zero
and outputs an alteration $h{\colon}S'\to S$ such that $S'$ is regular, $W'=h^{-1}(W)$ is snc,
and the saturated base change of $f$ is weakly semistable. The proof is similar to the arguments
used in the proofs of \cite[A.4.4 and A.4.3]{IKN@Quasi-unipotent}. First, one constructs a toroidal
blow up of the base that makes the fibers equidimensional (i.e. makes the log morphism integral),
and then an appropriate normalized finite base change is used to make the fibers reduced.
\end{remark2}

\ifx\danslelivre\undefined
\bibliography{../fichiers-configuration-latex/bibliographie}
\bibliographystyle{../fichiers-configuration-latex/smfalpha}


\begin{center}
XI. Produits orientés\\
Luc Illusie
\end{center}

\else
\setboolean{english}{false}
\setboolean{francais}{true}
\selectlanguage{francais}
\chapter{Produits orientés}
\addtocontents{toc}{par Luc Illusie}
\begin{center}
Luc Illusie
\end{center}
\fi

On fixe un univers $\mathscr{U}$. Sauf mention du contraire, les sites (resp. topos) considérés
seront des $\mathscr{U}$-sites (resp. topos), et « petit  » signifiera $\mathscr{U}$-petit.

\section{Construction des produits orientés} La construction suivante est due à Deligne
\cite{AG_Tokyo-Kyoto@Laumon} :

\subsection{}\label{11-site} Soient $f : X \ra S$, $g : Y \ra S$ des
morphismes de topos. On suppose que $X$, $Y$, $S$ ont des sites de définition  $C_1$, $C_2$, $D$,
admettant des limites projectives finies, et associés à des prétopologies, et que $f^{*}$, $g^{*}$
prolongent des foncteurs continus entre sites, et commutant aux limites projectives finies. Soit $C$
le site suivant :

(i) $C$ est la catégorie des couples de morphismes $U \ra V \la W$ au-dessus de $X \ra S \la Y$, où
$U \ra V$ (resp. $V \la W$) désigne un morphisme $U \ra f^{*}V$ (resp. $g^{*}V \la W$) de $C_1$
(resp. $C_2$) et $V$ est un objet de $D$.

(ii) $C$ est muni de la topologie définie par la prétopologie engendrée par les familles couvrantes
$(U_i \ra V_i \la W_i) \ra (U \ra V \la W)$ $(i \in I)$ du type suivant :

(a) $V_i = V$, $W_i = W$ pour tout $i$, et la famille $(U_i \ra U)$ est couvrante ;

(b) $U_i = U$, $V_i = V$ pour tout $i$, et la famille $(W_i \ra W)$ est couvrante ;

(c) $(U' \ra V' \la W') \ra (U \ra V \la W)$, où $U' = U$ et $W' \ra W$ est déduit par changement de
base d'un morphisme $V' \ra V$ de $D$. 

On note $\gtilde C$ le topos des faisceaux sur $C$.

\begin{lemme}\label{11-faisceau}
Soit $F$ un préfaisceau sur $C$. Pour que $F$ soit un
faisceau il faut et il suffit que les deux conditions
suivantes soient vérifiées :

(i) pour toute famille couvrante $(Z_i \ra Z)$ de $C$ du type (a) ou (b), la suite $F(Z) \ra
\prod_{i \in I} F(Z_i) \rightrightarrows \prod_{(i,j) \in I \times I} F(Z_i \times_Z Z_j)$ est
exacte ;

(ii) pour toute famille couvrante $(U' \ra V' \la W') \ra (U \ra V \la W)$ du type (c),
l'application $$ F(U \ra V \la W) \ra F(U' \ra V' \la W') $$ est bijective. En particulier, si l'on
note $(-)^a$ le foncteur faisceau associé, pour toute famille couvrante $(Z' \ra Z)$ du type (c), le
morphisme de faisceaux associés $Z'^a \ra Z^a$ est un isomorphisme.
\end{lemme}

La nécessité est triviale pour (i), et pour (ii), il suffit d'observer que le morphisme diagonal $$
(U \ra V' \la W') \ra (U \ra V' \times_V V' \la W' \times_W W') $$ est un morphisme couvrant (du
type (c)), qui égalise la double flèche $$ (U \ra V' \times_V V' \la W' \times_W W')
\rightrightarrows (U \ra V' \la W').  $$ Pour la suffisance, on note que les familles couvrantes
$(Z_i \ra Z)$ de type (a), (b), ou (c) sont stables par changement de base $Z' \ra Z$, et on
applique \sga{4}{ii}{2.3}.

\subsection{}\label{11-tau} Notons $e_X$ (resp. $e_Y$, $e_S$) l'objet final de $C_1$ (resp. $C_2$,
$D$). On a des projections naturelles $$ p_1 : \gtilde C \ra X, p_2 : \gtilde C \ra Y $$ données par
$$ p_1^*(U) = (U \ra e_S \la e_Y), p_2^*(W) = (e_X \ra e_S \la W).  $$ On a par ailleurs un
morphisme canonique $$ \tau : gp_2 \ra fp_1 $$ donné par le morphisme de foncteurs $\tau : (gp_2)_*
\ra (fp_1)_*$ défini de la façon suivante : pour un faisceau $F$ sur $C$, et un objet $V$ de $S$, $$
\tau : ((gp_2)_*F)(V) \ra ((fp_1)_*F)(V) $$ est le composé $$ F(e_X \ra e_S \la g^*V) \ra F(f^*V \ra
V \la g^*V) \ra F(f^*V \ra e_S \la e_Y), $$ où la première flèche est induite par la localisation
$(f^*V \ra V \la g^*V) \ra (e_X \ra e_S \la g^*V)$, et la seconde est l'inverse de l'isomorphisme
donné par \ref{11-faisceau}, relativement au morphisme de type (c) $(f^*V \ra V \la g^*V \ra (f^*V
\ra e_S \la e_Y)$.

\begin{theoreme}\label{11-existence} Soit $T$ un topos muni de morphismes $a : T
\ra X$, $b : T \ra Y$ et d'un morphisme $t : gb \ra fa$. Il existe alors un triplet $(h : T \ra
\gtilde C, \alpha : p_1h \iso a, \beta : p_2h \iso b$, unique à isomorphisme unique près, tel que le
composé $$ \xymatrix{gb \ar[r]^{\beta^{-1}} & gp_2h \ar[r]^{\tau} & fp_1h \ar[r]^{\alpha} & fa} $$
soit égal à $t$.
\end{theoreme}

Nous aurons besoin, pour la démonstration, du lemme suivant :

\begin{lemme}\label{11-lemme existence}
Soit $Z = (U \ra V \la W)$ un objet de $C$. Avec la notation
de \ref{11-faisceau}, le carré suivant est cartésien :
\begin{equation}
\xymatrix{Z^a \ar[d] \ar[r]
& (p_2^*W)^a \ar[d]^v \\ (p_1^*U)^a \ar[r]^u & ((gp_2^*)V)^a}.  \label{11-carré}
\end{equation}
Dans ce carré, $v$ et les flèches issues de $Z^a$ sont les flèches évidentes, et $u$ est la flèche
composée $$ \xymatrix{(p_1^*U)^a \ar[r] & ((fp_1)^*V)^a \ar[r]^{\tau} & ((gp_2^*)V)^a}, $$ où $\tau$
est le composé $$ \xymatrix{(f^*V \ra e_S \la e_Y)^a \ar[r]^{r^{-1}} & (f^*V \ra V \la g^*V)^a
\ar[r]  & (e_X \ra e_S \la g^*V)^a}, $$ $r$ désignant l'isomorphisme $(f^*V \ra V \la g^*V)^a \iso
(f^*V \ra e_S \la e_Y)^a$ de \ref{11-faisceau}.
\end{lemme}

Soit $z : Z \ra Z' = (f^*V \ra V \la
g^*V)$ la projection canonique. Le composé $Z \ra p_1^*U \ra (fp_1)^*V$ se factorise à travers $z$.
Par suite, et par définition de $\tau$, le diagramme $$ \xymatrix{Z^a \ar[r]^z \ar[d] & Z'^a
\ar[d]^{r} \ar[dr] & {}\\ (p_1^*U)^a \ar[r] & ((fp_1^*)V)^a \ar[r]^{\tau} & ((gp_2)^*V)^a } $$ est
commutatif. Comme le composé  $Z \ra p_2^*W \ra (gp_2)^*W$ se factorise aussi à travers $z$, le
carré  \ref{11-carré} est donc commutatif. Celui-ci est le pourtour du diagramme suivant, où les
flèches autres que $\tau$ sont les flèches évidentes : $$ \xymatrix {(f^*U \rightarrow V \rightarrow
g^*W)^a \ar[r]  \ar[d]&(f^*V \rightarrow V \leftarrow g^*W)^a \ar[r] \ar[d]& (e_X \rightarrow e_S
\leftarrow g^*W)^a \ar[d] \\ (f^*U \rightarrow V \rightarrow  g^*V) ^a\ar[r]  \ar[d]&(f^*V
\rightarrow V \leftarrow g^*V)^a \ar[r] \ar[d]^r & (e_X \rightarrow e_S \leftarrow g^*V)^a
\ar[d]^{\Id} \\ (f^*U \rightarrow e_S \leftarrow e_Y)^a \ar[r] &(f^*V \rightarrow e_S \leftarrow
e_Y)^a \ar[r]^{\tau} & (e_X \rightarrow e_S \leftarrow g^*V)^a}.  $$ Chacun des carrés qui le
composent est cartésien. Il en est donc de même de \ref{11-carré}.

\subsection{}\label{11-preuve du th.}
Prouvons \ref{11-existence}. On peut supposer que $a$ et $b$ sont donnés par des morphismes de
sites, la topologie du site de définition de $T$ étant moins fine que la topologie canonique
(d'après \sga{4}{iv}{4.9.4}, on pourrait prendre $T$ lui-même comme site de définition, avec sa
topologie canonique). Par \ref{11-lemme existence} l'unicité est claire : pour $Z = (U \ra V \la W)$
dans $C$, on doit avoir
\begin{equation}
h^*Z = a^*U \times_{(gb)^*V}b^*W, \label{11-déf de h}
\end{equation}
où $a^*U \ra (gb)^*V$ est le composé $\xymatrix{a^*U \ar[r] & (fa)^*V \ar[r]^t &
(gb)^*V}$. Les isomorphismes $\alpha$ et $\beta$ sont alors tautologiques, nous les négligerons dans
le reste de la démonstration. Vérifions que le foncteur $h^*$ donné par \ref{11-déf de h} définit un
morphisme de topos $h$ vérifiant la propriété énoncée en \ref{11-existence}. Comme $h^*$ commute aux
limites projectives finies, pour vérifier que $h^*$ induit un morphisme de topos, il suffit de
vérifier que $h^*$ est continu (\sga{4}{iv}{4.9.1, 4.9.2}). Il est trivial que $h^*$ transforme
familles couvrantes du type (a) ou (b) en familles couvrantes. Par ailleurs, si $(U' \ra V' \la W')
\ra (U \ra V \la W)$ est une famille couvrante du type (c), le carré $$ \xymatrix {b^*W' \ar[r]
\ar[d] &b^*W \ar[d]  \\ (gb)^*V' \ar[r] &(gb)^*V} $$ est cartésien, et par suite $$ a^*U
\times_{(gb)^*V'} b^*W' \rightarrow a^*U \times_{(gb)^*V} b^*W $$ est un isomorphisme. Il reste à
vérifier que $\tau$ induit $t$. Mais par définition, le morphisme de faisceaux défini par
$h^*(\tau)$ appliqué à $V$ est le composé $$ \xymatrix{((fa)^*V)^a \ar[r]^-{r^{-1}} & ((fa)^*V
\times_{(gb)^*V} (gb)^*V)^a \ar[r] & ((gb)^*V)^a}, $$ donc est égal à celui défini par $t$ appliqué
à $(fa)^*V$, ce qui achève la démonstration. 

\begin{definition}
\label{11-produit orienté} Le topos
$\gtilde C$ construit en \ref{11-site} s'appelle le \textit{produit orienté (gauche) de $X$ et $Y$
au-dessus de $S$}, et se note $X \ori_S Y$.  Les morphismes du diagramme $$ \xymatrix{{} & X \ori_S
Y \ar[dl]_{p_1} \ar[dr]^{p_2} & {} \\ X \ar[dr]^f  & {} & Y \ar[dl]_g \\ {} & S & {}} $$ sont reliés
par la 2-flèche $\tau : gp_2 \ra fp_1$. Il découle de \ref{11-existence} que le quadruplet $(X
\ori_S Y, p_1,p_2,\tau)$ est indépendant (à isomorphisme unique près) du choix des sites de
définition $C_1$, $C_2$, $D$.
\end{definition}

On définit de même le produit orienté droit $X \stackrel{\ra}{\times}_S Y$, avec ses projections
canoniques $p_1 : X \stackrel{\ra}{\times}_S Y \ra X$, $p_2 : X \stackrel{\ra}{\times}_S Y \ra Y$ et
la 2-flèche $\tau' : fp_1 \ra gp_2$, qui possède la propriété universelle de \ref{11-existence},
avec $X$ et $Y$ échangés.

\subsection{}\label{11-points} Désignons par $pt$ un topos ponctuel
(catégorie des faisceaux d'ensembles sur un espace réduit à un point). Soient $x : pt \ra X$, $y :
pt \ra Y$ des points de $X$ et $Y$ respectivement, et $u : gy \ra fx$ une 2-flèche. Par
\ref{11-existence}, le triplet $(x,y,u)$ définit un point $z : \mathrm{pt} \ra X \ori_S Y$ tel que
$p_1z \isononcan x$, $p_2z \isononcan y$. Ce point sera noté $(x,y,u)$ (ou parfois $(x,y)$ s'il n'y
a pas de confusion à craindre). Tout point de $X \ori_S Y$ est de cette forme.

\subsection{}\label{11-fonctorialité}
Considérons un diagramme de 1-morphismes de topos
\begin{equation}
\xymatrix{X' \ar[r]^{f'} \ar[d]^u & S' \ar[d]^h & Y' \ar[l]_{g'} \ar[d]^v \\ X
\ar[r]^f & S & Y \ar[l]_{g}} \label{11-fonct1}
\end{equation}
et des 2-flèches $a : hf' \ra fu$, $b
: gv \ra hg'$. Notons $T = X \ori_S Y$, $T' = X' \ori_{S'}Y'$, $p_1 : T \ra X$, $p_2 : T \ra Y$,
$p'_1 : T' \ra X'$, $p'_2 : T' \ra Y'$ les projections canoniques, $\tau : gp_2 \ra fp_1$, $\tau' :
g'p'_2 \ra f'p'_1$ les 2-flèches canoniques. Considérons la 2-flèche composée $$ c :
\xymatrix{gvp'_2 \ar[r]^b & hg'p'_2 \ar[r]^{\tau'} &hf'p'_1 \ar[r]^a &fup'_1}.  $$ D'après
\ref{11-existence}, $c$ définit un diagramme de 1-morphismes
\begin{equation}
\xymatrix{X' \ar[d]^u
&T' \ar[d]^t \ar[l]_{p'_1} \ar[r]^{p'_2} & Y' \ar[d]^v \\ X & T \ar[l]_{p_1}  \ar[r]^{p_2} & Y}
\label{11-fonct2}
\end{equation}
et des 2-isomorphismes $\alpha : p_1t \iso up'_1$, $\beta : p_2t
\iso vp'_2$ rendant commutatif le carré
\begin{equation}
\xymatrix{gp_2t \ar[r]^{\tau}
\ar[d]^{\beta} &  fp_1t \ar[d]^{\alpha} \\ gvp'_2 \ar[r]^{a\tau'b} &fup'_1}.  \label{11-fonct3}
\end{equation}
On dit que le triplet $(t,\alpha,\beta)$ (ou simplement $t : T' \ra T$) de
\ref{11-fonct2} est déduit de \ref{11-fonct1} par \textit{fonctorialité}. On notera $$ t = u \ori_h
v.  $$ On a une compatibilité évidente pour un composé de deux données \ref{11-fonct1}.

\subsection{}\label{11-exemples} Voici quelques exemples.

(a) Dans la situation de \ref{11-existence}, le triplet $(a,b,t)$ définit un diagramme de type
\ref{11-fonct1} $$ \xymatrix{T \ar[d]^a \ar[r]^{\Id} & T \ar[d]^{fa} & T \ar[l]_{\Id} \ar[d]^b \\ X
\ar[r]^f & S & Y \ar[l]_{g}}, $$ avec $t : gb \ra fa$, d'où un morphisme $$ a \ori_{fa} b : T \ori_T
T \ra X \ori_S T.  $$ Par ailleurs, d'après \ref{11-existence}, les flèches identiques de $T$
définissent un morphisme canonique, dit \textit{diagonal} $$ \Delta : T \ra T \ori_T T.  $$ La
1-flèche $h$ de \ref{11-existence} est la composée $$ h = (a \ori_{fa} b)\Delta : T \ra X \ori_S Y.
$$ En particulier, prenant pour $T$ un topos ponctuel, de sorte que $\Delta$ est un isomorphisme, on
a, avec les notations de \ref{11-points} $$ (x,y,u) = x \ori_{fx} y : pt \ra X \ori_S Y.  $$

(b) Dans la situation de \ref{11-produit orienté}, soient $X'$, $S'$, $Y'$ des objets de $X$, $S$,
$Y$ respectivement, et $f' : X' \ra S'$, (resp. $g' : Y' \ra S'$) une flèche au-dessus de $f$ (resp.
$g$), i. e. une flèche $f' : X' \ra f^*(S')$, (resp. $g' : Y' \ra g^*(Y')$). Notons $X' \ori_{S'}
Y'$ l'objet $(X' \ra S' \la Y') = p_1^*(X') \times_{(gp_2)^*(S')} p_2^*Y'$ de $X \ori_S Y$,
\textit{cf.} \ref{11-lemme existence}, On en déduit un diagramme 2-commutatif de1- flèches
naturelles

\begin{equation}
\xymatrix{{} & (X \ori_S Y)_{/(X' \ori_{S'}Y')} \ar[dl] \ar[dr] & {} \\
X_{/X'} \ar[r] \ar[d] & S_{/S'} \ar[d] & Y_{/Y'} \ar[l] \ar[d] \\ X \ar[r] & S & Y \ar[l]},
\label{11-localisation1}
\end{equation}
où la notation $(-)_{/-}$ désigne un topos localisé.
D'après \ref{11-existence}, la partie supérieure de \ref{11-localisation1} définit un 1-morphisme
\begin{equation} m :  (X \ori_S Y)_{/(X' \ori_{S'}Y')} \ra X_{/X'} \ori_{S_{/S'}} Y_{/Y'}.
\label{11-localisation2}
\end{equation}
Il résulte de \ref{11-déf de h} que $m$ est une équivalence,
par laquelle, dans la suite, nous identifierons les deux membres. D'autre part, les carrés
2-commutatifs de \ref{11-localisation1} définissent, d'après \ref{11-fonctorialité}, une flèche de
fonctorialité $$ X_{/X'} \ori_{S_{/S'}} Y_{/Y'} \ra X \ori_S Y, $$ Celle-ci, ou son composé avec
$m$,
\begin{equation}
(X \ori_S Y)_{/(X' \ori_{S'}Y')} \ra X \ori_S Y \label{11-localisation3}
\end{equation}
s'appelle \textit{flèche de localisation}.

\begin{proposition}\label{11-localisation4} Supposons que $X'$ soit l'objet final de $X$ et que $g'$
soit cartésien au-dessus de $g$, i. e. $g' : Y' \iso g^*(S')$. Alors la flèche
\ref{11-localisation3} est une équivalence. 
\end{proposition}

En effet, avec les notations de
\ref{11-localisation1}, il résulte de \ref{11-faisceau} que la flèche $e_X \ori_{S'} Y' \ra e_X
\ori_{e_S} e_Y$ de $X \ori_{S} Y$ est un isomorphisme.

\subsection{}\label{11-tubes} Considérons en
particulier le cas où $S = Y$ est un schéma muni de la topologie étale, $g = \Id$ et $f : X \ra S =
Y$ est l'inclusion d'un sous-schéma fermé de $Y$. Le topos $T = X \ori_Y Y$ joue le rôle d'un
\textit{voisinage tubulaire étale} de $Y$ dans $X$. Les points de $T$ sont les triplets $(x,y,t)$,
où $x$ (resp. $y$) est un point géométrique de $X$ (resp. $Y$) et $t : y \ra x$ une flèche de
spécialisation (\textit{cf.  (\sga{4}{viii}{7.9})}. En d'autres termes, $(x,y,t)$ est la donnée d'un
point géométrique $x$ de $X$, d'une générisation $y_0$ du point fermé (noté encore par abus $x$) du
localisé strict  $X_{(x)}$ de $X$ en $x$ et d'un point géométrique de  $X_{(x)}$ localisé en $y_0$,
ou encore d'une extension séparablement close $y \ra y_0$ du point générique de
$\overline{\{y_0\}}$. Par ailleurs, si $v : Y' \ra Y$ est un voisinage étale de $X$ dans $Y$, i. e.
un diagramme commutatif $$ \xymatrix{{} & Y' \ar[d]^v \\ X \ar[ur] \ar[r] & Y}, $$ où $v$ est étale,
alors, d'après \ref{11-localisation4} le morphisme canonique $$ X \ori_{Y'} Y' \ra X \ori_Y Y $$ est
une équivalence. Ainsi, $T$ ne dépend que du hensélisé de $X$ le long de $Y$ (lorsque celui-ci est
défini, en particulier, pour $Y$ affine, \textit{cf.} \cite{Anneaux@Raynaud}). Nous verrons au
numéro suivant et dans l'exposé XII d'autres propriétés de $T$ précisant cette analogie avec un
voisinage tubulaire.  \section{Tubes et changement de base}

\subsection{}\label{11-topos locaux}
Soit $(S,s)$ un topos ponctué, i. e. un couple formé d'un topos $S$ et d'un point $s : pt \ra S$ de
$S$. Si $(S,s)$ et $(T,t)$ sont des topos ponctués, un morphisme (ponctué) de $(S,s)$ dans $(T,t)$
est un couple $(f,a)$ d'un morphisme $f : S \ra T$ et d'une 2-flèche $a : fs \ra t$. Une 2-flèche $c
: (f,a) \ra (g,b)$ est une 2-flèche $c : f \ra g$ telle que $b(cs) = a$. Si $(S,s)$ est un topos
ponctué, on note $F \mapsto F_s = s^*F$ le foncteur fibre en $s$.

Rappelons qu'un topos ponctué $(S,s)$ est dit \textit{local de centre} $s$ (\sga{4}{vi}{8.4.6}) si,
pour tout objet $F$ de $S$, la flèche naturelle $\Gamma(S,F) \ra F_s$ est bijective. Un morphisme
$(f,a) : (S,s) \ra (T,t)$ de topos ponctués est dit \textit{local} si la 2-flèche $a : fs \ra t$ est
un isomorphisme.

\subsection{}\label{11-section canonique} La construction qui suit est due à
Gabber. Soit $(S,s)$ un topos local de centre $s$. Notons $\varepsilon : S \ra pt$ la projection.
Par définition la flèche canonique $\varepsilon_* \ra s^*$ est un isomorphisme. On en déduit un
isomorphisme
\begin{equation}
\varepsilon^*\varepsilon_* \iso (s\varepsilon)^*.  \label{11-can1}
\end{equation} La flèche d'adjonction $\varepsilon^*\varepsilon_* \ra \Id$ s'identifie donc, par
\ref{11-can1}, à un morphisme $(s\varepsilon)^* \ra \Id$, i. e. à une 2-flèche
\begin{equation}
c_s :
\Id  \ra s\varepsilon \label{11-can2}
\end{equation}
entre les 1-morphismes $\Id : S \ra S$ et
$s\varepsilon : S \ra S$. Si $F$ est un objet de $S$, $(s\varepsilon)^*F$ est le faisceau constant
sur $S$ de valeur $F_s = \varepsilon_*F = \Gamma(S,F)$. Si $U$ est un objet connexe de $S$, le
morphisme $(s\varepsilon)^*F \ra F$ induit sur $\Gamma(U,-)$ le morphisme de restriction
$\Gamma(S,F) \ra \Gamma(U,F)$. Le composé $c_ss : s \ra s$ est l'identité : $(s\varepsilon)^*F \ra
F$ induit l'identité sur les fibres en $s$.

Soient $f : X \ra S$, $g : Y \ra S$ des morphismes de topos, $x : pt \ra X$ un point de $X$, $s = fx
: pt \ra S$ son image dans $S$. Le diagramme
\begin{equation}
\xymatrix{pt \ar[r]^{\Id} \ar[d]^x & pt
\ar[d]^{s} & Y \ar[l]_{\varepsilon g} \ar[d]^{\Id} \\ X \ar[r]^f & S & Y \ar[l]_g}, \label{11-can3}
\end{equation} (où le carré de gauche est 2-commutatif) et la 2-flèche
\begin{equation}
c_sg : g \ra
s\varepsilon g \label{11-can4}
\end{equation}
sont une donnée de type \ref{11-fonct1}. Pour un objet
$F$ de $S$, $(s\varepsilon g)^*F$ est le faisceau constant sur $Y$ de valeur $F_s = \Gamma(S,F)$, et
la flèche $(s\varepsilon g)^*F \ra g^*F$ est la composée $\Gamma(S,F)_Y \ra \Gamma(Y,g^*F)_y \ra
g^*F$. Notons que, par \ref{11-existence}, le produit orienté $pt \ori_{pt} Y$ s'identifie
canoniquement à $Y$, avec $p_1 = \Id : Y \ra Y$. De \ref{11-can3} et \ref{11-can4} on déduit donc un
diagramme de type \ref{11-fonct2} :
\begin{equation}
\xymatrix{pt \ar[d]^x & Y \ar[d]^{\sigma}
\ar[l]_{\varepsilon g} \ar[r]^{\Id} & Y \ar[d]^{\Id} \\ X & X \ori_S Y \ar[l]_{p_1} \ar[r]^{p_2} & Y},
\label{11-can5}
\end{equation}
en d'autres termes, une section $\sigma : Y \ra X \ori_S Y$ de $p_2$
telle que $p_1\sigma = x\varepsilon g$. On peut voir cette section comme étant définie, via
\ref{11-existence} par le couple de morphismes $x\varepsilon g : Y \ra X$, $\Id : Y \ra Y$ et la
2-flèche $c_sg : g \ra fx\varepsilon g = s\varepsilon g$. On dit que $\sigma$ est la \textit{section
canonique} définie par le point $x$. Par composition avec $p_2{}_*$, la flèche d'adjonction $\Id \ra
\sigma_*\sigma^*$ donne une flèche canonique
\begin{equation}
\gamma : p_2{}_* \ra \sigma^*.
\label{11-can6}
\end{equation}
Le résultat suivant est dû à Gabber :

\begin{proposition}\label{11-can7} Soient $f : (X,x) \ra (S,s)$ un morphisme local de topos locaux
($fx = s$), et $g : Y \ra S$ un morphisme de topos. Soit $y : pt \ra Y$ un point de $Y$. Pour tout
objet $F$ de $X \ori_S Y$, $\gamma$ (\ref{11-can6}) induit un isomorphisme $$ \gamma_y :
(p_2{}_*F)_y \iso (\sigma^*F)_y.  $$
\end{proposition}

Soit $t = \sigma y : pt \ra T $ le point de
$T = X \ori_S Y$ image de $y$ par $\sigma$. Ce point est défini (\textit{cf.} 1.4) par le triplet
$(x,y,u)$, où $u : gy \ra fx = fx\varepsilon gy = s\varepsilon g y = s$ est déduit de \ref{11-can4}.
On a $(\sigma^*F)_y = F_t$. Par définition, $$ F_t = \mathrm{colim}_z F(U \ra V \la W), $$ où $z :
pt \ra (U \ra V \la W)$ parcourt les voisinages de $t$ dans $T$. Comme $X$ et $S$ sont locaux, les
voisinages de $t$ de la forme $\sigma w : pt \ra (e_X \ra e_S \la W)$, où $w : pt \ra W$ est un
voisinage de $y$ dans $Y$ forment un système cofinal. Donc $$ F_t = \mathrm{colim}_w F(e_X \ra e_S
\la W), $$ où $z = \sigma w : pt \ra (e_X \ra e_S \la W)$ parcourt les voisinages précédents, avec
$U = e_X$, $V = e_S$.  Par ailleurs, $$ (p_2{}_*F)_y = \mathrm{colim}_w F(e_X \ra e_S \la W), $$ où
$w : pt \ra W$ parcourt les voisinages de $y$ dans $Y$. La flèche $\gamma_y$ est la restriction
naturelle. C'est donc un isomorphisme. 

\begin{corollaire2}
\label{11-can8} Sous les hypothèses de
\ref{11-can7}, si $S$ a assez de points, en particulier, si $S$ est localement cohérent
(\sga{4}{vi}{9.0}), $\gamma$ (\ref{11-can6}) est un isomorphisme.
\end{corollaire2}

Il est plausible que l'hypothèse d'existence d'assez de points soit superflue. Celle-ci sera
cependant satisfaite dans les applications que nous avons en vue.

\begin{corollaire2}\label{11-can9} Sous les hypothèses de \ref{11-can7}, supposons $Y$ local de
centre $y$. Alors $X \ori_S Y$ est local de centre $\sigma(y)$. 
\end{corollaire2}

En effet, on a
alors $(p_2{}_*F)_y = \Gamma(Y,p_2{}_*F) = \Gamma(X \ori_S Y,F)$, et $\gamma_y$ s'identifie à la
restriction $ \Gamma(X \ori_S Y,F) \ra F_{\sigma(y)}$.

Notons que, si $f : X \ra S$ est un morphisme local de schémas strictement locaux, $g : Y \ra S$ un
morphisme de schémas strictement locaux, le produit fibré schématique $X \times_S Y$ n'est pas en
général strictement local, ni même local, même si $g$ est local.

Le résultat ci-après est dû également à Gabber :

\begin{theoreme}
\label{11-ch. base} Soient $f : X
\ra S$, $g : Y \ra S$ des morphismes de topos, $T = X \ori_S Y$, $p_1 : T \ra X$, $p_2 : T \ra Y$
les projections canoniques, $\tau = gp_2 \ra fp_1$ la 2-flèche canonique. On suppose $X$, $Y$, $S$
cohérents et $f$, $g$ cohérents (\sga{4}{vi}{2.3,2.4.5,3.1}). Soit $\Lambda$ un anneau. Alors, pour
tout $F \in D^+(Y,\Lambda)$, la flèche de changement de base, déduite de $\tau$,
\begin{equation}
f^*Rg_*F \ra Rp_1{}_*p_2^*F \label{11-ch. base1}
\end{equation}
est un isomorphisme (de
$D^+(X,\Lambda)$).
\end{theoreme}

Nous aurons besoin du lemme suivant, qui généralise \cite[9.1]{Cycles@Orgogozo} :

\begin{lemme}\label{11-cohérent} Sous les hypothèses de \ref{11-ch. base}, $T$ est cohérent, et les
projections $p_1$, $p_2$ sont des morphismes cohérents.
\end{lemme}

Le topos $X$ (resp. $Y$, resp. $S$) admet une (petite) famille génératrice $C_1$ (resp. $C_2$, resp.
$D$) formée d'objets cohérents, stable par limites projectives finies (\sga4{vi}{2.4.5}). Comme $f$
et $g$ sont cohérents, $f^*V$ (resp. $g^*V$) est cohérent si $V$ est dans $D$ (\sga4{vi}{3.2}). La
sous-catégorie pleine de $X$ (resp. $Y$, resp.  $S$ correspondante, munie de la topologie induite,
est un site de définition de $X$ (resp. $Y$, resp. $S$). Soit $C$ la catégorie définie comme en
\ref{11-site}, munie de la topologie définie par la prétopologie engendrée par les familles finies
de type (a) et (b) et les familles de type (c). Elle est stable par limites projectives finies, et
est un site de définition de $T$. Il suffit donc de montrer que tout objet de $C$ est quasi-compact
(\sga4{vi}{2.4.5}). Pour cela, notons $\mathscr{P}$ la prétopologie définie en \ref{11-site}.
Décrivons $\mathscr{P}$. Pour chaque objet $Z = (U \ra V \la W)$ de $C$, notons $\mathrm{Cov}(Z)$
l'ensemble des familles $(Z_i \ra Z)_{i \in I}$ obtenues par composition d'un nombre fini de
familles de type (c) et de familles (finies) de type (a) et (b). En particulier, l'ensemble $I$ est
fini. Par définition, la donnée des $\mathrm{Cov}(Z)$ vérifie les axiomes PT0, PT2 et PT3 de
\sga4{ii}{1.3}.  L'axiome PT1 (stabilité par changement de base) est également vérifié, les familles
de type (c), ainsi que les familles finies de type (a) (resp. (b)) étant stables par changement de
base, et le changement de base commutant à la composition des familles. La donnée des
$\mathrm{Cov}(Z)$ est donc une prétopologie, et par définition, c'est la prétopologie $\mathscr{P}$.
Comme les familles appartenant à $\mathrm{Cov}(Z)$ sont finies, tout objet de $C$ est
automatiquement quasi-compact, comme annoncé. La cohérence des projections $p_1$ et $p_2$ en
découle.

\begin{remarque}\label{11-g suffit} Gabber sait montrer que la conclusion de \ref{11-cohérent} vaut
sous les seules hypothèses que $X$, $Y$, $S$ et $g$ sont cohérents. Nous n'aurons pas besoin de
cette généralisation.
\end{remarque}

\subsection{}\label{11-preuve de ch. base} Prouvons
\ref{11-ch. base}. Comme $X$ est cohérent, donc possède assez de points, il suffit de vérifier que,
pour tout point $x : pt \ra X$ de $X$, la fibre en $x$ de \ref{11-ch. base1}
\begin{equation}
(f^*Rg_*F)_x \ra (Rp_1{}_*p_2^*F)_x \label{11-(*)}
\end{equation}
est un isomorphisme. Soit $s : pt
\ra S$ l'image de $x$ par $f$. Soit $X_{(x)}$ (resp. $S_{(s)}$) le localisé de $X$ (resp. $S$) en
$x$ (resp. $s$). Rappelons (\sga{4}{vi}{8.4.2}) que $X_{(x)}$ (resp. $Y_{(y)}$) (noté
$\mathrm{Loc}_x(x)$ (resp. $\mathrm{Loc}_s(S)$) dans \textit{loc. cit.}) est la limite projective $$
X_{(x)} = \mathrm{Limtop}_{U \in Vois(x)} X_{/U}, $$ (resp.  $$ S_{(s)} = \mathrm{Limtop}_{V \in
Vois(s)} S_{/V}), $$ où $U$ (resp. $V$) parcourt la catégorie cofiltrante $Vois(x)$ (resp. $Vois(s))
$des voisinages de $x$ (resp. $s$) dans $X$ (resp. $S$). Comme $X$ (resp. $S$) est cohérent, on peut
d'ailleurs se borner aux $U$ (resp. $V$) qui sont cohérents, les morphismes de transition étant
alors automatiquement cohérents. C'est ce que nous ferons dans la suite, notant encore $Vois(x)$
(resp. $Vois(s)$) la sous-catégorie pleine formée des $U$ (resp. $V$) cohérents. Le topos $X_{(x)}$
(resp. $S_{(s)}$) est un topos local, au-dessus de $X$ (resp. $S)$, dont l'image du centre est $x$
(resp. $s$) (\sga{4}{vi}{8.4.6}). Le morphisme $f$ induit un morphisme local $f_{(x)} : X_{(x)} \ra
S_{(s)}$. Définissons de même $$ T_{(x)} = \mathrm{Limtop}_{U \in Vois(x)} T_{/p_1^*U}, $$ $$
Y_{(s)} = \mathrm{Limtop}_{V \in Vois(s)} Y_{/g^*V}), $$ de sorte qu'on obtient un carré
\begin{equation} \xymatrix{T_{(x)} \ar[d]^{p_1} \ar[r]^{p_2} & Y_{(s)} \ar[d]^g \\ X_{(x)} \ar[r]^f
& S_{(s)}}, \label{11-carré localisé}
\end{equation}
avec une 2-flèche $\tau : gp_2 \ra fp_1$. Par
la compatibilité de la formation des produits orientés à la localisation (\ref{11-localisation2}),
la flèche
\begin{equation}
T_{(x)} \ra X_{(x)} \ori_{S_{(s)}} Y_{(s)}, \label{11-(0)}
\end{equation}
déduite de ce carré par \ref{11-existence} est une équivalence. D'après \ref{11-cohérent}, les topos
$X_{/U}$, $S_{/V}$, $T_{/p_1^*U}$, $Y_{/g^*V}$ sont cohérents, les flèches de transition des
systèmes projectifs $X_{/U}$, $S_{/V}$ sont cohérents, et les morphismes $g : Y_{/g^*V} \ra S_{/V}$,
$p_1 : T_{/p_1^*U} \ra X_{/U}$ sont cohérents. On est donc dans les conditions d'application de
\sga4{vi}{8.7.3}, qui, compte tenu de ce que $S_{(s)}$ et $X_{(x)}$ sont locaux, implique que les
flèches canoniques
\begin{equation}
(Rg_*F)_s \ra  R\Gamma(S_{(s)},Rg_*F) \ra R\Gamma(Y_{(s)},F),
\label{11-(1)} \end{equation}
\begin{equation}
(Rp_1{}_*p_2^*F)_x \ra
R\Gamma(X_{(x)},Rp_1{}_*p_2^*F) \ra R\Gamma(T_{(x)},p_2^*F) \label{11-(2)}
\end{equation}
sont des
isomorphismes. Avec les identifications \ref{11-(0)}, \ref{11-(1)} et \ref{11-(2)}, la flèche
\ref{11-(*)} s'identifie à la fibre en $x$ de la flèche de changement de base (déduite de $\tau$) du
carré \ref{11-carré localisé}. Cette flèche s'écrit
\begin{equation}
R\Gamma(Y_{(s)},F)  \ra
R\Gamma(T_{(x)},p_2^*F).  \label{11-(3)}
\end{equation}
On a : 

\medskip (*) : \textit{La flèche \ref{11-(3)} s'identifie canoniquement à la flèche de fonctorialité
définie par $p_2$}.  \medskip

Pour le vérifier, on peut supposer $S$ et $X$ locaux, de centres respectifs $s$ et $x$, et $f$
local. Par définition, \ref{11-(3)} est la flèche composée $$ \xymatrix{R\Gamma(S,Rg_*F)
\ar[r]^{\ref{11-(3)}} \ar[d]^{\alpha} & R\Gamma(T,p_2^*F) \\ R\Gamma(S,Rg_*Rp_2{}_*p_2^*F)
\ar[r]^{R\Gamma(S,\tau)} & R\Gamma(S,Rf_*Rp_1{}_*p_2^*F) \ar[u]_{\beta}}, $$ où la flèche $\alpha$
est définie par la flèche d'adjonction $adj : F \ra Rp_2{}_*p_2^*F$ (et donc est la flèche de
fonctorialité $R\Gamma(Y,F)  \ra R\Gamma(T,p_2^*F)$ définie par $p_2$), et $\beta$ est
l'isomorphisme canonique de transitivité relatif à $fp_1 : T \ra S$. Or, par définition de $\tau$
(\ref{11-tau}), pour tout faisceau $G$ sur $T$, la flèche $$ \Gamma(S,\tau) : \Gamma(S,(gp_2)_*G)
\ra \Gamma(S,(fp_1)_*G) $$ est l'identité. Il en est donc de même de la flèche horizontale
inférieure du diagramme ci-dessus, ce qui prouve (*). Il reste à prouver que \ref{11-(3)} est un
isomorphisme. En fait, la flèche
\begin{equation}
adj : F \ra Rp_2{}_*p_2^*F \label{11-adj}
\end{equation} est un isomorphisme. Pour le voir, il suffit d'observer que, compte tenu de la
description de $\gamma$ donnée en \ref{11-can9}, le composé $$ \xymatrix{F \ar[r]^{adj}  &
Rp_2{}_*p_2^*F \ar[r]^{\gamma} & \sigma^*p_2^*F = F}, $$ où $\sigma$ est la section de $p_2$ définie
en \ref{11-can5} et $\gamma$ l'isomorphisme de \ref{11-can8}, est l'identité. Ceci achève la
démonstration de \ref{11-ch. base}.

\begin{remarques}\label{11-remarques}
\renewcommand{\theenumi}{(\arabic{enumi})}
\begin{enumerate}
\item Supposons que les données de \ref{11-ch. base} proviennent de morphismes de schémas, munis de la
topologie étale, avec $X$, $S$, $Y$ cohérents et $f$ et $g$ cohérents. Les points $x$, $s$ sont des
points géométriques, et les localisés $X_{(x)}$, $S_{(s)}$ des localisés stricts. Si $f$ est une
immersion fermée, la flèche d'adjonction \ref{11-adj} est un isomorphisme (on peut en effet supposer
$S$ strictement local, et il en et alors de même de $X$). Supposons de plus que $Y = S$, $g = \Id_S$
comme en \ref{11-tubes}. On a vu en \textit{loc. cit.} que $T = X \ori_S S$ joue le rôle d'un
voisinage tubulaire de $X$ dans $S$. Soient $j : S^*  = S - X \ra S$ l'ouvert complémentaire de $X$,
et $T^* = X \ori_S S^* = T_{/(e_X \ra e_S \la S^*)}$ le topos induit. Alors  $T^* $ joue le rôle
d'un \textit{voisinage tubulaire épointé} de $X$ dans $S$ : pour $F \in D^+(S^*,\Lambda)$, on a, par
\ref{11-ch. base}, $$ f^*Rj_*F \iso Rp_1{}_*p_2^*F.  $$

\item Sans l'hypothèse de cohérence sur $g$, la conclusion de \ref{11-ch. base} peut être en
défaut, comme le montre l'exemple suivant, dû à Gabber. Soient $X$ un espace topologique connexe,
non vide, $i : Y \ra X$ l'inclusion d'un fermé non vide distinct de $X$, $j : U= X - Y \ra X$
l'inclusion de l'ouvert complémentaire. Alors $i^*j_*\ZZ$ est non nul. Mais, si tout point de $U$ a
un voisinage dont l'adhérence dans $X$ ne rencontre pas $Y$, alors le produit orienté $Y \ori_X U$
est vide. C'est le cas par exemple, si $X$ est le segment $[0,1]$ et $Y$ le point $\{0\}$. 

\item Sous les hypothèses de \ref{11-ch. base}, on montre de manière analogue que :
\begin{itemize}

\item[(a)] Pour tout faisceau d'ensembles $F$ sur $Y$, la flèche de changement de base $$ f^*g_*F \ra
p_1{}_*p_2^*F $$ est un isomorphisme.

On peut espérer des variantes non abéliennes supérieures :

\item[(b)] Pour tout faisceau en groupes $F$ sur $Y$, la flèche de changement de base $$ f^*R^1g_*F \ra
R^1p_1{}_*p_2^*F, $$ est un isomorphisme (de faisceaux d'ensembles pointés).

\item[(c)] Plus généralement, pour tout champ $F$ sur $Y$, la flèche de changement de base $$ f^*g_*F \ra
p_1{}_*p_2^*F $$ est une équivalence.

\end{itemize}

La vérification de (b) et (c) semble requérir, outre les techniques de réduction des champs aux
gerbes de \cite[III 2.1.5]{CNA@Giraud}, des résultats de passage à la limite pour la cohomologie non
abélienne analogues à ceux de \sga4{vi}{8.7}, pour lesquels nous ne connaissons pas de référence.

\item Gabber sait démontrer la généralisation suivante de \ref{11-ch. base}. Soient $f : X
\ra S$, $g : Y \ra S$ des morphismes de topos, $T = X \ori_S Y$. On suppose que $Y$ et $S$ sont
localement cohérents, et que, pour tout objet cohérent algébrique $V$ de $S$, $g^*V$ est cohérent
algébrique (\sga4{vi}{2.1,2.3}). Alors, pour tout $F \in D^+(Y,\Lambda)$, la flèche de changement de
base \ref{11-ch. base1} est un isomorphisme, et on devrait avoir des résultats analogues dans le cas
non abélien, comme en (3) (a), (b), (c) ci-dessus. Gabber déduit ces résultats d'un théorème général
de changement de base pour certains topos fibrés.
\end{enumerate}
\end{remarques}

\section{Produits fibrés} Les compléments donnés
dans ce numéro et le suivant ne seront pas utilisés dans le reste du volume.

\subsection{}\label{11-site produit fibré} Les produits fibrés de topos ont été construits par
Giraud \cite[3.4]{classtop@Giraud}. La construction suivante est due à Gabber. Soient $f : X \ra S$,
$g : Y \ra S$ des morphismes de topos comme en \ref{11-site}. Soit $D$ le site suivant :

\renewcommand{\theenumi}{(\roman{enumi})}
\begin{enumerate}
\item La catégorie sous-jacente à $D$ est la catégorie $C$ considérée en \ref{11-site} (i). 

\item $D$ est munie de la topologie définie par la prétopologie engendrée par les familles couvrantes
$(U_i \ra V_i \la W_i) \ra (U \ra V \la W)$ $(i \in I)$ de la forme (a), (b), (c) de \ref{11-site}
(ii) et de la forme

(d) $(U' \ra V' \la W') \ra (U \ra V \la W)$, où $W' = W$ et $U' \ra U$ est déduit par changement de
base d'un morphisme $V' \ra V$ du site de définition de $S$.
\end{enumerate}

En d'autres termes, la topologie sur $D$ est la borne supérieure des topologies sur $C$ définissant
les produits orientés $X \ori_S Y$ et $X \stackrel{\ra}{\times}_S Y$. 

D'après \ref{11-faisceau}, pour qu'un préfaisceau $F$ sur $D$ soit un faisceau, il faut et il suffit
que $F$ vérifie les conditions d'exactitude habituelles relatives aux familles couvrantes de type
(a) et (b), et que, pour toute famille couvrante $Z' \ra Z$ de type (c) ou (d), $F(Z) \ra F(Z')$
soit un isomorphisme. 

Soit $\gtilde D$ le topos des faisceaux sur $D$. On a des projections naturelles $$ p_1 : \gtilde D
\ra X, p_2 : \gtilde D \ra Y $$ données par les mêmes formules qu'en \ref{11-tau}, et la construction
de $\tau$ en (\textit{loc. cit.}) donne un \textit{isomorphisme}
\begin{equation}
\varepsilon : gp_2
\iso fp_1.  \label{11-epsilon} \end{equation}

\begin{theoreme}\label{11-existence produit fibré}
Soit $T$ un topos muni de morphismes $a : T \ra X$, $b : T \ra Y$ et d'un isomorphisme $t : gb \iso
fa$. Il existe alors un triplet $(h : T \ra \gtilde D, \alpha : p_1h \iso a, \beta : p_2h \iso b$,
unique à isomorphisme unique près, tel que le composé $$ \xymatrix{gb \ar[r]^{\beta^{-1}} & gp_2h
\ar[r]^{\varepsilon} & fp_1h \ar[r]^{\alpha} & fa} $$ soit égal à $t$.
\end{theoreme}

La démonstration est analogue à celle de \ref{11-existence}. Le foncteur $h^*$ est encore donné par la
formule \ref{11-déf de h}. Comme $t$ est un isomorphisme, on a $$ a^*U \times_{(gb)^*V}b^*W = a^*U
\times_{(fa)^*V} b^*W, $$ où $b^*V \ra (fa)^*V$ est le composé $\xymatrix{b^*V \ar[r] & (gb)^*V
\ar[r]^{t^{-1}} & (fa)^*V}$. Il s'ensuit que $h^*$ transforme famille couvrante de type (d) en
famille couvrante, et l'on conclut comme dans \ref{11-preuve du th.}.

\begin{definition}\label{11-produit fibré} Le topos $\gtilde D$ s'appelle \textit{produit fibré de
$X$ et $Y$ au-dessus de $S$}, et se note $X \times_S Y$. Les morphismes du diagramme $$ \xymatrix{{}
& X \times_S Y \ar[dl]_{p_1} \ar[dr]^{p_2} & {} \\ X \ar[dr]^f  & {} & Y \ar[dl]_g \\ {} & S & {}}
$$ sont reliés par le 2-isomorphisme $\varepsilon : gp_2 \ra fp_1$.
\end{definition}

\begin{exemples}\label{11-exemples produits fibrés}
\renewcommand{\theenumi}{(\arabic{enumi})}
\begin{enumerate}
\item \textit{Espaces topologiques}. Soient $f : X \ra S$, $g : Y \ra S$ des applications continues
entre espaces topologiques. Alors le topos $\widetilde{X \times_S Y} $ des faisceaux sur le produit
fibré usuel $X \times_S Y$ représente le produit fibré $\gtilde X \times_{\gtilde S} \gtilde Y$ : le
morphisme naturel $\widetilde{X \times_S Y} \ra \gtilde X \times_{\gtilde S} \gtilde Y$ est une
équivalence. On définit en effet un quasi-inverse en observant que les ouverts de $X \times_S Y$ de
la forme $U \times_V W$ forment une base. 

\item \textit{Schémas}. Soient $f : X \ra S$, $g : Y \ra S$ des morphismes de schémas. Désignons par
l'indice $\zar$ (resp. $\et$) le topos zariskien resp. étale) associé. Du fait de (1),
le morphisme naturel $(X \times_S Y)_{\zar} \ra X_{\zar} \times_{S_{\zar}} Y_{\zar}$ n'est pas une
équivalence en général. De même, le morphisme naturel $(X \times_S Y)_{\et} \ra
X_{\et} \times_{S_{\et}} Y_{\et}$ n'est pas une équivalence en général, même
si $X$, $Y$, $S$ sont les spectres de corps : si $S = \SP k$, $S_{\et}$ est équivalent au
topos classifiant $BG$ du groupe profini $G = \mathrm{Gal}(\overline{k}/k)$, où $\overline{k}$ est
une clôture séparable de $k$, et la formation de $BG$ commute aux produits fibrés.  \section{Topos
évanescents et co-évanescents}
\end{enumerate}
\end{exemples}

\subsection{}\label{11-déf. topos évanescents} Soient $f : X \ra S$,
$g : Y \ra S$ des morphismes de topos comme en \ref{11-site}. Le produit orienté
\begin{equation}
X \ori_S S, \label{11-évanescent}
\end{equation}
où $S \ra S$ est le morphisme identique, s'appelle le
\textit{topos évanescent} de $f$. Il est étudié dans \cite{AG_Tokyo-Kyoto@Laumon} et
\cite{Cycles@Orgogozo}. Le produit orienté
\begin{equation}
S \ori_S Y, \label{11-co-évanescent}
\end{equation}
où $S \ra S$ est le morphisme identique, joue un rôle important dans les travaux de
Faltings sur les théorèmes de comparaison $p$-adiques et la correspondance de Simpson $p$-adique
(\textit{topos de Faltings}) (cf. \cite{almost-etale@Faltings}, \cite{Simpson@Faltings},
\cite{Simpson@Abbes-Gros}, \cite{toposdefaltings@Deligne}, \cite{co-evanescent@AG}). On propose ici de
l'appeler \textit{topos co-évanescent} de $g$. Du topos évanescent de $\Id_S$,
\begin{equation}
\fl S = S \ori_S S, \label{11-S fléché}
\end{equation}
qui est aussi le topos co-évanescent de $\Id_S$, les
produits orientés $X \ori_S Y$ se déduisent par changement de base. Considérons en effet le produit
fibré itéré
\begin{equation}
Z = X \times_S \fl S \times_S Y,
\end{equation} où la flèche de gauche
(resp. droite) de $\fl S$ vers $S$ est $p_1$ (resp. $p_2$). On a des projections naturelles $q_1 : Z
\ra X$, $q_2 : Z \ra Y$ et $m : Z \ra \fl S$, avec des isomorphismes $p_1m \iso fq_1$, $gq_2 \iso
p_2m$. Par composition avec la flèche structurale $\tau : p_2 \ra p_1$ de $\fl S$, on en déduit une
flèche $z : gq_2 \ra fp_1$. D'après\ref{11-existence}, le triplet $(q_1,q_2,z)$ définit donc une
flèche
\begin{equation}
h : Z \ra X \ori_S Y \label{11-h1}
\end{equation}
et des isomorphismes $p_1h
\iso q_1$, $p_2h \iso q_2$, par lesquels $z$ s'identifie à $\tau h$, où $p_1$, $p_2$ désignent les
projections canoniques de $X \ori_S Y$ sur $X$ et $Y$. 

\begin{proposition}
\label{11-produit orienté
et produit fibré} Le morphisme $h$ (\ref{11-h1}) est un isomorphisme. En particulier, il définit des
isomorphismes canoniques
\begin{equation}
X \times_S \fl S \iso X \ori_S S,
\end{equation}
\begin{equation} \fl S \times_S Y \iso S \ori_S Y. 
\end{equation}
\end{proposition}

Il suffit de montrer que $Z$, muni de $(q_1,q_2,z)$ vérifie la
propriété universelle du produit orienté. Soit $T$
un topos muni de morphismes $a : T \ra X$, $b : T \ra Y$, et d'une 2-flèche $t : gb \ra fa$. Par la
propriété universelle de $\fl S$, on en déduit d'abord un unique triplet, formé d'un morphisme $k :
T \ra \fl S$ et d'isomorphismes $p_1k \ra fa$, $p_2k \ra gb$ tels que $t = \tau k$ modulo ces
identifications. Puis, par la propriété universelle des produits fibrés, on en déduit un unique
quadruplet formé d'un morphisme $s : T \ra Z$ et d'isomorphismes $ms \iso k$, $q_1s \iso a$, $q_2s
\iso b$ tel que $zs = t$ modulo ces identifications.  \subsection{}\label{11-psi} Soit $f : X \ra S$
un morphisme de topos. D'après \ref{11-existence}, les morphismes $\Id_X$ et $f$ définissent un
morphisme
\begin{equation}
\Psi : X \ra X \ori_S S \label{11-psi1}
\end{equation}
tel que
\begin{equation}
p_1\Psi = \Id_X, p_2\Psi = f, \tau \Psi = \Id_f, \label{11-psi2}
\end{equation}
où
$\tau : p_2 \ra fp_1$ est la 2-flèche structurale de $X \ori_S S$ :
\begin{equation}
\xymatrix{{} &
X \ar[d]^{\Psi} \ar[dl]_{\Id_X} \ar[dr]^f & {} \\ X \ar[dr]_f & X \ori_S S \ar[l]_{p_1} \ar[r]^{p_2}
& S \ar[dl]^{\Id_S} \\ {} & S & {}} \label{11-psi1bis}
\end{equation}

 Le foncteur $\Psi_*$ s'appelle \textit{foncteur cycles proches}. Pour un objet $(U \ra V \la W)$ du
site $C$ définissant $X \ori_S S$, on a $\Psi^*(U \ra V \la W) = U \times_V W$, où $U \times_V W :=
U \times_{f^*(V)} f^*W$. Si $\Lambda$ est un anneau et $F \in D^+(X,\Lambda)$, le complexe $R\Psi_*F
\in D^+(X \ori_S S,\Lambda)$ (noté aussi $R\Psi F$) s'appelle \textit{complexe des cycles proches}
(de $f$ relatif à $F$). 

L'identité $p_1{}_*\Psi_* = \Id$ définit, par adjonction, un morphisme canonique
\begin{equation}
p_1^* \ra \Psi_*.  \label{11-psi3}
\end{equation}
Pour $F \in D^+(X,\Lambda)$, le cône du morphisme
$p_1^*F \ra R\Psi F$ qui s'en déduit s'appelle le \textit{complexe des cycles évanescents} (de $f$
relatif à $F$) et se note $R\Phi F$. Dans le cas des schémas (munis de la topologie étale), ces
foncteurs, qui généralisent les foncteurs $R\Psi$ et $R\Phi$ de Grothendieck (\sga{7}{i, xiii}),
sont étudiés dans \cite{AG_Tokyo-Kyoto@Laumon} et \cite{Cycles@Orgogozo}. 

Considérons le morphisme de changement de base
\begin{equation}
p_1{}_* \ra \Psi^* \label{11-psi4}
\end{equation} déduit de l'identité $p_1\Psi = \Id_X$, en d'autres termes, le morphisme déduit, par
application de $p_1{}_*$, de la flèche d'adjonction $\Id \ra \Psi_* \Psi^*$, compte tenu de ce que
$p_1{}_*\Psi_* = \Id$. Le résultat suivant est donné sans démonstration dans
\cite{AG_Tokyo-Kyoto@Laumon} :

\begin{proposition}
\label{11-isocan1} Le morphisme \ref{11-psi4} est
un isomorphisme. 
\end{proposition}

On va définir un morphisme
\begin{equation}
\Psi^* \ra p_1{}_*,
\label{11-psi5}
\end{equation}
dont on montrera qu'il est inverse de \ref{11-psi4}. Pour cela, on
définit un morphisme
\begin{equation}
\Id \ra \Psi_*p_1{}_* \label{11-psi6}
\end{equation}
de la
façon suivante. Pour un faisceau $F$ sur $X \ori_S S$ et un objet $Z = (U \ra V \la W)$ du site $C$
de \ref{11-site}, la flèche $F(Z) \ra (\Psi_*p_1{}_*F)(Z)$ est la composée
\begin{equation}
F(Z) \ra
F(U \times_{f^*V} f^*W \ra W \la W) \ra F(U \times_{f^*V} f^*W \ra e_S \la e_S), \label{11-psi7}
\end{equation}
où $W \ra W$ est l'identité, la première flèche est la restriction et la seconde,
l'inverse de l'isomorphisme relatif au recouvrement de type (c) $$ (U \times_{f^*V} f^*W \ra W \la
W) \ra (U \times_{f^*V} f^*W \ra e_S \la e_S).  $$ Le morphisme \ref{11-psi5} est adjoint de
\ref{11-psi6}. Notons $u$ (resp. $v$) le morphisme \ref{11-psi4} (resp. \ref{11-psi5}). On va
montrer que $u$ et $v$ sont inverses l'un de l'autre. L'argument qui suit est dû à Orgogozo. Il
s'agit de montrer que, pour tout faisceau $F$ sur $X \ori_S S$ et tout faisceau $G$ sur $X$, les
applications $\alpha(F,G) = \mathrm{Hom}(u(F),G) : \mathrm{Hom}(\Psi^*F,G) \ra
\mathrm{Hom}(p_1{}_*F,G)$ et $\beta(F,G) = \mathrm{Hom}(v(F),G) :\mathrm{Hom}(p_1{}_*F,G) \ra
\mathrm{Hom}(\Psi^*F,G)$ sont inverses l'une de l'autre.

 L'application $$ \alpha(F,G) : \mathrm{Hom}(\Psi^*F,G) = \mathrm{Hom}(F,\Psi_*G) \ra
\mathrm{Hom}(p_1{}_*F,G) $$ envoie $a : F \ra \Psi_*G$ sur $p_1{}_*a : p_1{}_*F \ra p_1{}_*\Psi_*G =
G$. L'application $a$ est la donnée d'une famille compatible d'applications $a_{(U \ra V \la W)} :
F(U \ra V \la W) \ra G(U \times_V W)$, pour $(U \ra V \la W)$ parcourant les objets de $C$,
``compatible" voulant dire compatible aux flèches de restriction. L'application $\alpha(a)$ est la
famille $a_{(U \ra e_S \la e_S)} : F(U \ra e_S \la e_S) \ra G(U)$, $U$ parcourant les objets de $X$.

L'application $$ \beta(F,G) : \mathrm{Hom}(p_1{}_*F,G) \ra \mathrm{Hom}(\Psi_*F,G) =
\mathrm{Hom}(F,\Psi_*G) $$ envoie $b : p_1{}_*F \ra G$ sur le composé $F \ra \Psi_*p_1{}_*F \ra
\Psi_*G$, où la première flèche est \ref{11-psi6} et et la seconde $\Psi_*b$. L'application
$\beta(b)$ est la famille $\beta(b)_{(U \ra V \la W)} : F(U \ra V \la W) \ra G(U \times_V W)$ est la
composée de \ref{11-psi7} et de $$ (\Psi_*b)(U \times_V W \ra e_S \la e_S) : F(U \times_V W \ra e_S
\la e_S) = (p_1{}_*F)(U \times_V W) \ra G(U \times_V W).  $$

Pour chaque $a : F \ra \Psi_*G$, on a un diagramme commutatif $$ \xymatrix{F(U \ra V \la W) \ar[d]
\ar[r]^{a_{(U \ra V \la W)}} & G(U \times_V W) \\ F(U \times_V W \ra W \la W) \ar[r] & F(U \times_V
W \ra e_S \la e_S) \ar[u]_{a_{(U\times_V W \ra e_S \la e_S)}}}, $$ où la flèche horizontale
inférieure est l'isomorphisme figurant dans \ref{11-psi7}. Dans ce diagramme, le composé des flèches
autres que la flèche horizontale supérieure est $\beta \alpha(a)_{(U \ra V \la W)}$, donc $\beta
\alpha = \Id$. La vérification de $\alpha \beta = \Id$ est triviale également. 

\begin{remarque}\label{11-isocan2} Lorsque les topos $X$ et $S$ sont localement
cohérents, donc en particulier dans le cas des schémas, on peut prouver \ref{11-isocan1} plus
simplement, en se ramenant au cas local. Comme $X$ a assez de points, il suffit de montrer que, pour
tout point $x : pt \ra X$, et tout faisceau $F$ sur $X \ori_S S$, la fibre en $x$ de \ref{11-psi4},
est un isomorphisme. Quitte à remplacer $X$ par son localisé en $x$ (\textit{cf.}
\ref{11-localisation1}), on peut supposer $X$ local de centre $x$. Soit $s : pt \ra S$ l'image de
$x$ par $f$. Soit $F$ un faisceau sur $T = X \ori_S S$. On doit montrer que $$ (p_1{}_*F)_x =
\Gamma(T,F) \ra (\Psi_*F)_x = F_{(x,s)} $$ est un isomorphisme. Soient $S_{(s)}$ le localisé de $S$
en $s$ et $T_{(s)} = T \ori_{S_{(s)}} S_{(s)}$. D'après les résultats de passage à limite invoqués
dans \ref{11-preuve de ch. base}, on a $$ \Gamma(T_{(s)},F) = \mathrm{colim} \Gamma(T_U,F), $$ où
$U$ parcourt les voisinages cohérents de $s$ et $T_U = X \ori_{S_{/U}} S_{/U}$. D'après
\ref{11-localisation4}, les flèches de restriction $\Gamma(T,F) \ra \Gamma(T_U,F)$ sont des
isomorphismes. On peut donc supposer $S$ local de centre $s$. D'après \ref{11-can9}, $T$ est alors
local, de centre $(x,s)$, et $\Gamma(T,F) = F_{(x,s)} = (\Psi^*F)_x$.
\end{remarque}

\subsection{}\label{11-topos de Faltings}
Soit $g : Y \ra S$ un morphisme de topos. Les faisceaux sur le topos co-évanescent $T =
S \ori_S Y$ (\ref{11-co-évanescent} ont une description simple, due à Deligne
\cite{toposdefaltings@Deligne}. Pour un faisceau $F$ sur $T$, la flèche de restriction $F(U \ra V
\la W) \ra F(U \ra U \ra U \times_V W)$ est un isomorphisme. Cela suggère de considérer le site
$C_0$ suivant. La catégorie $C_0$ est celle des flèches $(V \la W)$ au-dessus de $f : Y \la S$, i.
e. des flèches $W \ra g^*V$ de $Y$. On munit $C_0$ de la topologie définie par la prétopologie
engendrée par les familles couvrantes des types (a) et (b) ci-après :

(a) $(V \la W_i)_{i\in I} \ra (V \la W)$, où la famille $(W_i \ra W)_{i \in I}$ est couvrante,

(b) $(V_i \la W_i)_{i \in I} \ra (V \la W)$, où la famille $(V_i \ra V)_{i \in I}$ est couvrante, et
$W_i = V_i \times_V W$.

On montre (\cite[4.10]{co-evanescent@AG}) que $C_0$  est un site de définition de $T$. Un objet $F$
de $S \ori_S Y$, décrit comme un faisceau sur $C_0$, s'interprète comme la donnée d'une famille de
faisceaux $F_V : W \mapsto F(V \la W)$ sur $g^*V$ et de flèches de restriction $F_V \ra
j_{V'V}{}_*F_{V'}$ pour $V' \ra V$, définissant $j_{V'V} : g^*V' \ra g^*V$, satisfaisant la
condition de descente que, pour une famille couvrante $(V_i \ra V)_{i \in I}$, la suite $$ F_V \ra
\prod_{i} j_{V_iV}{}_*F_{V_i} \rightrightarrows \prod_{ii'} j_{V_{ii'}V}{}_*F_{V_{ii'}} $$ soit
exacte, où $V_{ii'} = V_i \times_V V_{i'}$. D'après \ref{11-existence}, les morphismes $g : Y \ra S$
et $\Id_Y$ définissent un morphisme
\begin{equation}
\Psi : Y \ra S \ori_S Y \label{11-copsi1}
\end{equation} tel que

\begin{equation}
p_1\Psi = g , p_2\Psi = \Id_Y , \tau\Psi = \Id_g,
\label{11-copsi2}
\end{equation}
où $\tau : gp_2 \ra p_1$ est la 1-flèche structurale de $S \ori_S Y$ :
\begin{equation}
\xymatrix{{} & Y \ar[d]^{\Psi} \ar[dl]_{g} \ar[dr]^{\Id_Y} & {} \\ S
\ar[dr]_{\Id_S} & S \ori_S Y \ar[l]_{p_1} \ar[r]^{p_2} & Y \ar[dl]^{g} \\ {} & S & {}}.
\label{11-copsi1bis}
\end{equation}
Le foncteur $\Psi_*$, qu'on pourrait appeler \textit{foncteur
cycles co-proches}, se comporte de manière très différente du foncteur cycles proches de
\ref{11-psi1}. En effet, de l'identité $p_2\Psi = \Id_Y$ on déduit, par adjonction, un morphisme
canonique
\begin{equation}
p_2^* \ra \Psi_*, \label{11-copsi3}
\end{equation}
analogue de
\ref{11-psi3}, et l'on a :

\begin{proposition}
\label{11-copsi4} Le morphisme \ref{11-copsi3} est un
isomorphisme. 
\end{proposition}
En particulier, \textit{le foncteur $\Psi_*$ est exact.} Ici, c'est
la flèche de changement de base, déduite de l'identité $p_2\Psi = \Id_Y$,
\begin{equation}
p_2{}_*
\ra \Psi^*, \label{11-copsi5}
\end{equation}
analogue de \ref{11-psi4}, qui n'est pas, en général,
un isomorphisme. On peut donner de \ref{11-copsi4} une démonstration analogue à celle de
\ref{11-isocan1}. Il est plus simple de déduire ce résultat des descriptions explicites suivantes
des foncteurs $p_1^*$, $p_2^*$, $\Psi$ et du morphisme $\tau$. Ces descriptions sont dues à Deligne
\cite{toposdefaltings@Deligne}.  \subsection{}\label{11-descriptions} (a) \textit{Description de
$p_1^*$}. On a $p_1^*V = (V \la g^*V)$ (l'objet $(V \ra e_S \la e_Y)$ de $C$ correspondant à l'objet
$(V \la g^*V)$ de $C_0$). Si $F$ est un faisceau sur $S$, $p_1^*F$ est le faisceau associé au
préfaisceau dont la valeur en $(V \la W)$ est la limite inductive des $F(V')$ suivant la catégorie
des flèches $(V \la W) \ra p_1^*V'$. Cette catégorie ayant $(V \la W) \ra (V \la g^*V)$ pour objet
initial, cette limite est égale à $F(V)$. En d'autres termes, $p_1^*F$ est le faisceau associé au
préfaisceau dont la valeur en $(V \la W)$ est $F(V)$. Dans la description donnée plus haut d'un
faisceau sur $S \ori_S Y$ en termes d'une famille de faisceaux sur les $g^*V$, $p_1^*F$ est la
famille des faisceaux constants $G_V$ sur $g^*V$ de valeur $F(V)$.

\medskip (b) \textit{Description de $p_2^*$.} On a $p_2^*W = (e_S \la W)$. Si $F$ est un faisceau
sur $Y$, $p_2^*F$ est le faisceau associé au préfaisceau dont la valeur en $(V \la W)$ est la limite
inductive des $F(W')$ suivant la catégorie des flèches $(V \la W) \ra p_2^*W'$. Cette catégorie
ayant $(V \la W) \ra (e_S \la W)$ pour objet initial, cette limite est égale à $F(W)$. En d'autres
termes, $p_2^*F$ est la famille des faisceaux $H_V$ sur $g^*V$, où $H_V$ est le faisceau $W \mapsto
F(W)$. 

\medskip (c) \textit{Description de $\tau$.} Si $F$ est un faisceau sur $S$, le morphisme $\tau :
p_1^*F \ra (gp_2)^*F$ est déduit du morphisme de préfaisceaux qui, pour $(V \la W)$ dans $C_0$,
envoie $F(V)$ dans $(gp_2)^*F(V \la W) = (g^*F)(W)$ par le morphisme composé $F(V) \ra (g^*F)(g^*V)
\ra (g^*F)(W)$. On le voit à l'aide de (a) et (b), en explicitant le morphisme $(fp_1)^* \ra
(gp_2)^*$ adjoint du morphisme $\tau$ décrit en \ref{11-tau}.

\medskip (d) \textit{Description de $\Psi_*$.} Si $F$ est un faisceau sur $Y$, et $(U \ra V \la W)$
un objet de $C$, on a $(\Psi_*F)(U \ra V \la W) = F(U \times_V W$. Dans la description de $S \ori_S
Y$ à l'aide du site $C_0$, on a donc $$ (\Psi_*F)(V \la W) = F(W).  $$ Compte tenu de (b), on a donc
$$ \Psi_*F = p_2^*F $$ Cette identification est celle donnée par \ref{11-copsi3}, ce qui prouve
\ref{11-copsi4}.

Notons encore que les points de $T = S \ori_S Y$ sont les flèches $s \la y$, où $s$ (resp. $y$) est
un point de $S$ (resp. $Y$), et que, pour un faisceau $F$ sur $Y$, si $(s \la y)$ est un point de
$T$, on a $$ (p_2^*F)_{(s \la y)} = F_y = (\Psi_*F)_{(s \la y)}.  $$

\begin{remarque}
Comme l'observe Gabber, l'isomorphisme \ref{11-copsi3} implique le théorème de
changement de base \ref{11-ch. base} pour $X = S, f = \Id_S$ sans hypothèse de cohérence sur $S$,
$Y$, et $g$.
\end{remarque}

\ifx\danslelivre\undefined
\bibliography{../fichiers-configuration-latex/bibliographie}
\bibliographystyle{../fichiers-configuration-latex/smfalpha}


\begin{center}
XII. Descente cohomologique orientée\\
Fabrice Orgogozo
\end{center}

\version

\else
\setboolean{english}{false}
\setboolean{francais}{true}
\selectlanguage{francais}
\chapter{Descente cohomologique orientée}
\addtocontents{toc}{par Fabrice Orgogozo}
\begin{center}
Fabrice Orgogozo
\end{center}
\fi

\section{Acyclicité orientée des morphismes propres}

\subsection{}
\label{12-notations invcoho}
L'objet de cette section est de démontrer le théorème \ref{12-invariance cohomologique}
ci-dessous, qui généralise l'invariance par éclatement admissible
du voisinage tubulaire défini à l'aide du produit fibré orienté
(cf. \refext{XI}{11-produit orienté} et \refext{XI}{11-remarques}).
Avant d'énoncer le théorème, fixons quelques notations.
On considère un schéma $S$ cohérent, $g:Y → S$, $π: S′ → S$,
et $a:X → S′$ des morphismes de schémas avec $π$
\emph{propre} et $a$ \emph{cohérent}. 
À ces données sont associés :
\begin{itemize}
\item le morphisme composé $b= π ∘ a$ ;
\item le morphisme $g′: Y′ → S′$ déduit de $g$ par
le changement de base $π$, où $Y′=Y ×_S S′$, et le morphisme de projection $ρ:Y′ → Y$ ;
\item les topos $T=Y\po_S X$ et $T'=Y'\po_{S'} X$
tels que définis en \emph{loc. cit.} par produits
orientés des topos étales associés aux schémas considérés,
ainsi enfin que le morphisme $\fl{ρ}:T'→T$ déduit par fonctorialité
des morphismes $ρ$ et $π$ (\refext{XI}{11-fonctorialité}).
\end{itemize}

\[
\xymatrix{
Y \cartesien \ar[d]_g & Y' \ar[d]^{g'} \ar[l]_{ρ} &  & T'=Y'\po_{S'} X
\ar[d]_{\fl{ρ}} \\
S & S' \ar[l]^{π} & X \ar[l]^{a}
\ar@/^1pc/[ll]^b & T=Y\po_S X
}
\]

\subsection{Énoncés}

\begin{théorème2}
\label{12-invariance cohomologique}
Sous les hypothèses du \ref{12-notations invcoho},
le morphisme $\fl{ρ}$ est acyclique pour les faisceaux de torsion :
pour tout entier $n≥1$ et tout objet $𝒦$
de $\dpluscoefs{T}{𝐙/n}$, le morphisme
d'adjonction $𝒦→\R \fl{ρ}_* \fl{ρ}^*𝒦$ est un isomorphisme.
\end{théorème2}

\begin{théorème2}\label{12-invariance cohomologique non abélienne}
Sous les hypothèses du \ref{12-notations invcoho},
le morphisme $\fl{ρ}$ est $1$-acyclique
pour tout champ ind-fini $𝒞$ sur $T$,
le $2$-morphisme d'adjonction $𝒞→\fl{ρ}_* \fl{ρ}^*𝒞$ est une équivalence de catégories.
\end{théorème2}

\subsection{Démonstrations}

\subsubsection{Réductions}
\label{12-réductions démonstration invcoho}
Soit $𝒦$ comme dans l'énoncé \ref{12-invariance cohomologique}.
Notons $𝒦'$ son image inverse $\fl{ρ}^*𝒦$ sur $T'$.
Le topos $T$ est \emph{cohérent}, cf. \refext{XI}{11-cohérent}.
Il a donc assez de points (\sga{4}{VI}{9.0}) et
ceux-ci sont comme décrits en \refext{XI}{11-points},
c'est-à-dire associés à une paire de points (géométriques) $(y,x)$
de $Y$ et $X$ et à une spécialisation $b(x) ⤳ g(y)$.
En calculant la fibre du morphisme $𝒦→\R \fl{ρ}_* \fl{ρ}^*𝒦$ en un tel
point, on vérifie comme en \refext{XI}{11-ch. base}
que l'on peut supposer $Y,X,S$ et $g$ locaux (pour la topologie étale),
de sorte que le topos $T$ l'est également (cf. \refext{XI}{11-topos locaux}, \refext{XI}{11-can9}).
Sous ces hypothèses supplémentaires, il nous faut montrer que le
morphisme d'adjonction (image inverse)
\[
α : \RΓ(T,𝒦) → \RΓ(T',𝒦')
\]
est isomorphisme.
Considérons les projections $p₁',p₁,p₂$ et $p₂'$
telles que ci-dessous :
$$
\xymatrix{
Y' \ar[d]_{ρ} &  T'=Y'\po_{S'} X \ar[l]_{p₁'} \ar[d]_{\fl{ρ}}
\ar[rd]^{p₂'} & \\
Y & T=Y\po_S X \ar[l]_{p₁} \ar[r]^{p₂} & X
}
$$
de sorte qu'on a en particulier l'égalité tautologique
$\RΓ(T',𝒦')=\RΓ(Y',\R{p₁'}_*𝒦')$.

\subsubsection{Cas où le complexe $𝒦$ provient du schéma $X$}
Supposons qu'il existe un complexe $𝒢 ∈ \dpluscoefs{X}{𝐙/n}$
tel que
\[
𝒦 = p₂^* 𝒢.
\]
Par commutativité du triangle de droite ci-dessus,
on a également $𝒦′={p₂'}^*𝒢$.
Le but $\RΓ(Y',R{p₁'}_*𝒦')$ du morphisme d'adjonction $α$
est donc le complexe $RΓ(Y',\R{p₁'}_*{p₂'}^*𝒢)$,
lui-même isomorphe, d'après \refext{XI}{11-ch. base}
à $\RΓ(Y',{g'}^*\R a_*𝒢)$.
Par propreté de $a$ et le théorème de changement de base
propre, le morphisme d'adjonction
\[
\RΓ(Y,g^*\R π_*\R a_* 𝒢 ) → \RΓ(Y,\R ρ_*{g'}^*\R a_*𝒢)
= \RΓ(Y',{g'}^*\R a_*𝒢)
\]
est un isomorphisme.
Les schémas $Y$ et $S$ étant locaux pour la topologie
étale et $g$ étant un morphisme local, on
a successivement $\RΓ(Y,g^*\R π_*\R a_* 𝒢 )=\RΓ(S,\R b_* 𝒢)$
et $\RΓ(S,\R b_* 𝒢)=𝒢_x$, où $x$ est le point
de $X$ évident. Enfin, le topos $T$ et le
morphisme $p₂$ étant locaux également,
on a $𝒢_x = \RΓ(T,𝒦)=𝒦_t$.
Nous avons exprimé la flèche $α$ comme composée
de morphismes d'adjonction qui sont des isomorphismes. CQFD.

\subsubsection{Cas général}
\label{12-cas général démonstration invcoho}
On ne suppose dorénavant plus que $𝒦$ se descend à $X$.
Considérons cependant le complexe $𝒢=σ^*𝒦$, où $σ$
est la section du morphisme $p₂$ définie en \refext{XI}{11-section canonique}.
(Si $𝒦$ se descend, on retrouve le complexe $𝒢$ ci-dessus.)
D'après \refext{XI}{11-can7}, le morphisme d'adjonction $\R {p₂}_*𝒦→σ^*𝒦$
est un isomorphisme, de sorte que la source $𝒦_t=\RΓ(T,𝒦)=\RΓ(X, \R {p₂}_* 𝒦)$
du morphisme d'adjonction $α$ est isomorphe à $\RΓ(X,𝒢)=𝒢_x=\fl{𝒢}_t$
où l'on pose $\fl{𝒢}={p₂}^*𝒢$. En d'autre termes,
la coünité $\fl{𝒢}→𝒦$ de l'adjonction induit un isomorphisme sur les sections globales
du topos \emph{local} $T=Y\po_S X$. D'après ce qui précède (cas où $𝒦$ se descend),
le morphisme d'adjonction $\RΓ(T,\fl{𝒢}) → \RΓ(T′,\fl{𝒢}′)$ est un
isomorphisme, où $\fl{𝒢}'$ désigne le tiré en arrière
par $p₁'$ de $\fl{𝒢}$ sur le topos $T'$.
Il serait donc suffisant de montrer que le morphisme
$c:\R{p₁'}_*\fl{𝒢}' → \R{p₁'}_*𝒦'$ est un isomorphisme
car en appliquant le foncteur $\RΓ(Y', —)$ on obtient
le morphisme $\RΓ(T′,\fl{𝒢}') → \RΓ(T′,\fl{𝒦}')$ déduit
de la coünité. Par propreté de $π$ donc de $Y′$ sur $Y$,
il suffit même de montrer que la restriction de $c$
à la fibre spéciale $Y′_y$ est un isomorphisme. C'est
ce que nous allons vérifier. (On utilise ici le théorème
de changement de base propre pour les faisceaux de torsion.)
Soient donc $y'$ un point (« géométrique ») de $Y'_y$, $s'$ son image dans $S'$
et calculons la fibre en $y'$ du morphisme $c$.
Notons $Y″$ et $S″$ les localisés $Y'_{(y')}$ et $S'_{(s')}$ ainsi
que $X″$ le produit fibré $X ×_{S′} S″$.
La fibre en $y′$ d'une image directe par $p₁′$
s'identifie à la cohomologie du topos $T″= Y″ \po_{S″} X″$
(à valeurs dans le tiré en arrière).
Si $σ″$ désigne la section canonique (\refext{XI}{11-section canonique})
de la seconde projection $p₂″:T″ → X″$, on a un isomorphisme
canonique (\refext{XI}{11-can7})
\[
\RΓ(T″,-) = \RΓ(X″,{σ″}^*-).
\]
Il en résulte qu'il suffit de montrer que la coünité $\fl{𝒢}″ → 𝒦″$
devient un isomorphisme après application du foncteur $σ″$,
où l'on note $𝒦″$ et $\fl{𝒢}″$ les tirés en arrière sur $T″$.
Ceci résulte immédiatement de la transitivité des images inverses
et de l'égalité $σp₂pσ″=pσ″$ où $p$ désigne
le morphisme $T″ → T$. Cette égalité résulte
à son tour du fait que le morphisme $Y″ → Y$ est
\emph{local}.

\subsubsection{Cas non abélien}
\label{12-remarque-cas-non-abélien}
La démonstration est identique. Pour le théorème de changement de base
propre non abélien, on fait appel à \cite[VII.2.2.2]{CNA@Giraud}.
Signalons que, comme signalé en \refext{XI}{11-remarques} (3)
nous n'avons pas connaissance d'une référence publiée
permettant de justifier le passage à la limite nécessaire au
calcul des fibres.

\section{Descente cohomologique orientée}

\subsection{Topologie orientée des altérations}
\subsubsection{}
Soient $S$ un schéma nœthérien et $ℬ$ la catégorie des diagrammes de schémas
nœthériens $X\dessusdessous{f}{→} S \dessusdessous{g}{←} Y$ et des
morphismes rendant les carrés commutatifs.
Considérons le pseudo-foncteur de $ℬ$ vers la $2$-catégorie
des topos, qui envoie l'objet précédent sur le produit
fibré orienté $Y\po_S X$, où l'on note abusivement $Y$ pour $Y_\et$, etc.
Remarquons que les limites finies (resp. les coproduits) existent
(resp. existent et sont disjoints, universels) dans la
catégorie $ℬ$ ; ils se calculent « terme à terme ». Par exemple,
le produit fibré de $X₁→S₁←Y₁$ et $X₂→S₂←Y₂$ au-dessus de $X→S←Y$
est $(X₁×_X X₂)→ (S₁×_S S₂) ← (Y₁×_Y Y₂)$.

\subsubsection{}
\label{12-topologie alt-orientée}
On considère la topologie orientée des altérations sur $ℬ$ engendrée
par les familles ci-dessous :

\begin{enumerate}
\item $\big((X→S←Y_i) → (X→S←Y)\big)_{i∈I}$, où $(Y_i→Y)_{i∈I}$ est une famille
$\alt$-couvrante ;
\item $\big((X_i→S←Y) → (X→S←Y)\big)_{i∈I}$, où $(X_i→X)_{i∈I}$ est
une famille $\alt$-couvrante ;
\item $\big((X×_S S_i → S_i ← Y×_S S_i) → (X→S←Y)\big)_{i∈I}$, où
$(S_i→S)_{i∈I}$ est une famille $\alt$-couvrante ;
\item $\big((X → S' ← Y×_S S') → (X→S←Y)\big)$, où $S'→S$ est un morphisme propre ;
\item $\big((X×_S S' → S' ← Y)→ (X→S←Y)\big)$, où $S'→S$ est un morphisme étale.
\end{enumerate}

Remarquons que les propriétés des familles
de morphismes (i-v) sont stables par changement de base dans
la catégorie $ℬ$.

\subsection{Énoncés}

\begin{théorème2}\label{12-descente cohomologique}
Soit $(X_{∙}→S_{∙}←Y_{∙})→ (X→S←Y)$ un hyperrecouvrement
pour la topologie orientée des altérations, c'est-à-dire un objet simplicial de $ℬ$
au-dessus de $(X→S←Y)$ tel que pour tout entier $i≥0$, le morphisme
canonique $(X_{i+1}→S_{i+1}←Y_{i+1})→\big(\cosq_i(X_{∙}→S_{∙}←Y_{∙})\big)_{i+1}$
de la catégorie $ℬ$ soit couvrant pour la topologie orientée des
altérations. Sous cette hypothèse, et si l'on note
$\fl{ε}:T_∙=Y_{∙}\po_{S_{∙}} X_{∙} → T=Y\po_S X$ le morphisme d'augmentation, on a les résultats
de descente suivants :
\begin{enumerate}
\item pour tout complexe de torsion borné inférieurement $𝒦$ sur $Y\po_S X$, le morphisme
d'adjonction $𝒦→\R \fl{ε}_*\fl{ε}^*𝒦$ est un isomorphisme.
\item pour tout champ $𝒞$ ind-fini sur $Y\po_S X$, le morphisme
d'adjonction $𝒞→\fl{ε}_*\fl{ε}^*𝒞$ est une équivalence.
\end{enumerate}
\end{théorème2}

\subsubsection{}En d'autres termes, le morphisme $\fl{ε}$ est de descente
cohomologique pour les faisceaux de torsion et les champs ind-finis.
Pour le sens à donner à l'énoncé (ii), on renvoie à \cite{CNA@Giraud} (spécialement
chap. VII, §2.2) et \cite{Alterations@Orgogozo}, §2.

\subsubsection{}
La principale application que nous ferons du théorème précédent
est la formule de changement de base suivante.

\begin{théorème2}\label{12-cb orienté}
Soient $X\dessusdessous{f}{→ }S\dessusdessous{g}{←}Y$ un diagramme de
schémas nœthériens et $ε :S_{∙}→ S$ un hyperrecouvrement pour la topologie des altérations.
Notons $f_{∙}:X_{∙}=X×_{S} S_{∙} →  S_{∙}$
et $g_{∙}:Y_{∙}=X×_{S} S_{∙}→S_{∙}$
(resp. $ε_Y:Y_{∙}=Y×_{S} S_{∙} → Y$)
les hyperrecouvrements pour la topologie des altérations
(resp. pour la topologie orientée des altérations) qui s'en déduisent.
\begin{enumerate}
\item Pour tout complexe de torsion borné inférieurement
$𝒦$ sur $X$, le morphisme
$$
g^*\R f_*𝒦 →  \R{ε_Y}_*\big( {g_{∙}}^*\R{f_{∙}}_{*}
𝒦_{|X_{∙}}\big)
$$
est un isomorphisme.
\item Pour tout champ ind-fini $𝒞$ sur $X$, le morphisme
$$
g^*f_*𝒞 →  {ε_Y}_*\big( {g_{∙}}^*{f_{∙}}_{*}
𝒞_{|X_{∙}}\big)
$$
est un isomorphisme.
\end{enumerate}
\end{théorème2}

\subsection{Démonstration du théorème \ref{12-descente cohomologique}}

Nous commençons par démontrer le (i). Pour la variante
non abélienne (ii), cf. \ref{12-variante non abélienne descente}.

\subsubsection{Réductions}
D'après la théorie générale de la descente — cf. p. ex. \cite[5.3.5]{Hodge3@Deligne} —,
il suffit de démontrer le théorème pour les hyperrecouvrements correspondant à un $0$-cosquelette
associé à une des familles de morphismes du type (i) à (v) (\ref{12-topologie alt-orientée}).
De plus, et pour la même raison, on peut
se contenter du cas particulier où les familles $\alt$-couvrantes
apparaissant dans la description de recouvrement sont,
soit un recouvrement ouvert de Zariski, soit un morphisme
propre et surjectif. Il suffit donc de démontrer les cinq propositions suivantes,
où l'on fixe un objet $(X\dessusdessous{f}{→}S\dessusdessous{g}{←}Y)$
de $ℬ$. Rappelons que $f$ est automatiquement cohérent.
Dans les paragraphes ci-dessous, on reprend les notations
de l'énoncé. En particulier, on note $T$ le topos
$Y\po_S X$, on fixe un entier $n ≥ 1$, un complexe $𝒦 ∈ \Ob  \dpluscoefs{T}{𝐙/n}$
et on note $\fl{ε}$ le morphisme du topos simplicial $T_∙$
vers $T$. On veut montrer que la flèche d'adjonction
\[
α: 𝒦 → \R \fl{ε}_* \fl{ε}^* 𝒦 = \R \fl{ε}_* 𝒦_∙
\]
est un isomorphisme.

\subsubsection{Famille de type (i) : $(X→S←Y′) → (X→S←Y)$}
Soit $Y'→Y$ un morphisme couvrant pour la topologie de Zariski
ou bien propre et surjectif, et $Y_{∙}=\cosq_0^Y (Y')$ son cosquelette.
Par passage aux fibres, on peut supposer les schémas $X$, $Y$ et $S$
ainsi que morphisme $g : Y → S$ locaux. Comme on l'a déjà
dit, le topos $T$ est alors un topos local.
Si $Y'→Y$ est Zariski couvrant, il possède alors une section.
Il en est donc de même du morphisme $T′=Y′ \po_S X → T$, auquel cas
le résultat est connu (cf. \sga{4}{V}{bis 3.3.1}).
Supposons donc le morphisme $Y'→Y$ propre et surjectif.
La flèche $α$ s'identifie au morphisme
\[
𝒦_t=\RΓ(T,𝒦)→\R Γ(T_{∙},𝒦_{∙})=\R Γ(Y_{∙},\R {p_{1 ∙}}_* 𝒦_∙).
\]
Reprenons les notations de \ref{12-réductions démonstration invcoho}.
Par descente cohomologique classique, l'adjonction
$𝒦_t = \RΓ(Y,\R{p₁}_* 𝒦) → \RΓ(Y_∙,(\R{p₁}_* 𝒦)_∙))$ est un isomorphisme.
Par propreté de $ρ:Y′ → Y$, il nous suffit de montrer
que la \emph{restriction au-dessus de $y$} de l'adjonction
\[β:{(\R{p₁}_* 𝒦)}_∙ → {\R p_{1 ∙}}_* 𝒦_∙\]
est un isomorphisme. En effet, le théorème de changement de base
propre nous dit que la source et le but de $α$ se calculent
par restriction aux fibres spéciales.
Que $β_y$ soit un isomorphisme résulte immédiatement du lemme ci-dessous (i),
appliqué aux schémas $Z=Y_i$.

\begin{lemme2}
\begin{enumerate}
Soient $X → S$, $Y → S$ et $Z → S$  des morphismes cohérents.
\item Soit $γ:Z → Y$ un morphisme cohérent.
Notons $\fl{γ}:Z\po_S X→ Y\po_S X$ le morphisme induit.
Le morphisme d'adjonction
\[
γ^*\big(\R{p₁^Y}_*𝒦\big)→ \R{p₁^Z}_*\fl{γ}^*𝒦
\]
est un isomorphisme en tout point de $Z_y$.
\item Soit $δ:Z→X$ un morphisme cohérent.
Notons $\fl{δ}:Y\po_S Z → Y\po_S X$ le morphisme induit.
Le morphisme d'adjonction
\[
δ^*\big(\R{p₂^X}_*𝒦\big)→ \R{p₂^Z}_*\fl{δ}^*𝒦
\]
est un isomorphisme en tout point de $Z$.
\end{enumerate}
\end{lemme2}

\begin{démo}
(i) Si $γ$ est un morphisme local de schéma locaux,
le morphisme $\fl{γ}$ est un morphisme local de topos locaux.
La conclusion en résulte par passage aux fibres en observant
que $Z_{(z)} → Y$ est local si $z$ est un point (géométrique)
localisé sur $Z_y$.
(ii) Même argument : pour tout point (géométrique) $z$
de $Z$, d'image $x$ par $δ$, le morphisme $Y \po_S Z_{(z)} →
Y \po_S X_{(x)}$ est un morphisme local de topos locaux.
\end{démo}

\subsubsection{Famille de type (ii) : $(X′→S←Y) → (X→S←Y)$}
\label{12-X propre}
Quitte à remplacer $p₁$ par $p₂$ et utiliser
le (ii) du lemme ci-dessus, la même démonstration s'applique.
Notons le théorème de changement de base propre
n'apparaît pas directement ici mais est néanmoins
utilisé dans la démonstration du théorème de descente classique
pour $X_{∙}→X$.

\subsubsection{Famille de type (iii) : $(X′=X×_S S′ → S′ ←Y′=Y×_S S′) → (X→S←Y)$}
Notons $S_{∙}=\cosq_0^S (S')$ le cosquelette du morphisme $S′ → S$
et $X_{∙}$, $Y_{∙}$ les schémas simpliciaux qui s'en déduisent.
Par localisation, on se ramène au cas où $S$ est local de sorte que
le cas de la topologie de Zariski est un corollaire de l'existence
d'une section au morphisme $S'→S$. Supposons donc $S'→S$ propre et
surjectif et factorisons $\fl{ε}$ en
\[
Y_{∙}\po_{S_{∙}} X_{∙} \dessusdessous{φ}{→ } Y\po_S X_{∙}
\dessusdessous{π}{→ } Y\po_S X.
\]
On a vu en \ref{12-X propre} que le morphisme $𝒦→\R π_*π^*𝒦$
est un isomorphisme. D'autre part, il résulte du théorème \ref{12-invariance cohomologique} et du fait
que les images directes se calculent cran par cran que 
pour tout complexe $𝒢_{∙} ∈ \dpluscoefs{Y\po_S X_{∙}}{𝐙/n}$ sur $Y\po_S X_{∙}$ — par exemple $π^*𝒦$ —
le morphisme d'adjonction $𝒢_{∙}→\R φ_*φ^*𝒢_{∙}$ est un isomorphisme.
Le résultat en découle.

\subsubsection{Famille de type (iv) : $(X → S' ← Y′= Y×_S S') → (X→S←Y)$}
On suppose ici que le morphisme $X → S$ se factorise
à travers un morphisme \emph{propre} $S′ → S$, et on définit $S_∙$ et $Y_∙$ comme ci-dessus.
Remarquons que pour chaque entier $i ≥ 0$, le topos $Y_i\po_{S_i} X$ est
le topos associé au produit fibré itéré dans la catégorie $ℬ$ du
morphisme $(X→S'←Y′)→ (X→S←Y)$.
Factorisons le morphisme $\fl{ε} : Y_{∙}\po_{S_{∙}} X → Y\po_S X$
en
\[
Y_{∙}\po_{S_{∙}} X \dessusdessous{φ}{→ } (Y\po_S X)_{\mathrm{const}} \dessusdessous{π}{→ } Y\po_S X,
\]
où $(Y\po_S X)_{\mathrm{const}}$ est le topos simplicial constant.
Il résulte du théorème \ref{12-invariance cohomologique} appliqué pour
chaque $i$ aux morphismes $a:X\dessusdessous{\text{diag}}{→}S_i$
et $b:S_i→S$ que l'adjonction $\Id→\R φ_*φ^*$ est un isomorphisme. On s'est donc ramené
au cas trivial où $S'=S$.

\subsubsection{Famille de type (v) : $(X′= X×_S S' → S' ← Y)→ (X→S←Y)$}
On suppose ici que le morphisme $Y → S$ se factorise
à travers un morphisme \emph{étale} $S′ → S$.
Notons à nouveau $S_{∙}$ le $0$-cosquelette de ce morphisme
et $X_{∙}=X ×_S S_{∙}$ le schéma simplicial qui s'en déduit.
Comme ci-dessus, notons que pour chaque entier $i ≥ 0$ le topos $Y\po_{S_i} X_i$
est le topos associé au produit fibré itéré dans la catégorie $ℬ$ du
morphisme $(X′ →S'←Y)→ (X→S←Y)$. Pour montrer que la flèche
d'adjonction $α: 𝒦 → \R \fl{ε}_* \fl{ε}^* 𝒦$ est un isomorphisme, on peut supposer — par passage aux
fibres — que le schéma $S$ est strictement hensélien de
sorte que le morphisme $S'→S$, et par conséquent $(X×_S S' → S' ← Y)→
(X→S←Y)$, aient une section. Le résultat est alors évident.

\subsubsection{Variante non abélienne}
\label{12-variante non abélienne descente}
Modulo la difficulté signalée ci-dessus (\ref{12-remarque-cas-non-abélien}), il suffit de remplacer les références
\cite[5.3.5]{Hodge3@Deligne} et \sga{4}{V}{bis 3.3.1}
par \cite[2.5,2.8]{Alterations@Orgogozo}. 

\subsection{Démonstration du théorème \ref{12-cb orienté}}

\subsubsection{Cas des complexes}
Soit $𝒦$ un complexe sur $X$ comme dans l'énoncé,
et notons $\fl{𝒦}$ son image inverse par la seconde
projection $p₂:Y\po_S X → X$.
Il résulte du théorème \ref{12-descente cohomologique} (i)
que l'adjonction $α: 𝒦 → \R \fl{ε}_* \fl{ε}^*𝒦$
est un isomorphisme. (On rappelle que $\fl{ε}$
désigne le morphisme d'augmentation $Y_{∙}\po_{S_{∙}} X_{∙} → Y\po_S
X$.)

D'autre part, on a la chaîne d'isomorphismes :
\begin{align*}
g^*\R f_*𝒦					& ⥲ \R{p₁}_*\fl{𝒦} \\
\R{p₁}_*α : \R{p₁}_*\fl{𝒦}			& ⥲ \R{p₁}_*{\R \fl{ε}}_* {\fl{ε}}^*\fl{𝒦} \\
\R{p₁}_*\R{\fl{ε}}_* {\fl{ε}}^*\fl{𝒦}	& = \R {ε_Y}_* \R {p_{1∙}}_* \fl{𝒦}_{∙}  \text{, où } \fl{𝒦}_{∙}={\fl{ε}}^*\fl{𝒦} \\
\R {ε_Y}_* \R {p_{1∙}}_* \fl{𝒦}_{∙}		& ⭉ \R{ε_Y}_* g_{∙}^* \R {f_{∙}}_* \fl{𝒦}_{∙}
\end{align*}
Les premier et dernier isomorphismes résultent de \refext{XI}{11-ch. base}
et le troisième de la fonctorialité des images directes.

\subsubsection{Cas des champs}
Même démonstration. 

\ifx\danslelivre\undefined
\bibliography{../fichiers-configuration-latex/bibliographie}
\bibliographystyle{../fichiers-configuration-latex/nomcomplet-date}


\begin{center}
XIII. Le théorème de finitude\\
Fabrice Orgogozo
\end{center}

\else
\setboolean{english}{false}
\setboolean{francais}{true}
\selectlanguage{francais}
\chapter{Le théorème de finitude}
\addtocontents{toc}{par Fabrice Orgogozo}
\begin{center}
Fabrice Orgogozo
\end{center}
\fi

\section{Introduction}

\subsection{}L'objet de cet exposé est de démontrer le théorème
suivant (\refext{0}{0-finitude1}).

\begin{théorème2}
\label{13-constructibilité}
Soient $X$ un schéma nœthérien \emph{quasi-excellent} (\refext{I}{1def-excellent}),
$f:Y→X$ un morphisme de type fini, $n≥1$ un entier inversible sur $X$
et $ℱ$ un faisceau constructible de $𝐙/n$-modules sur $Y$.
Alors :
\begin{enumerate}
\item Pour tout entier $q ≥ 0$ le faisceau $\R^q f_* ℱ$ est
constructible.
\item Il existe un entier $N$ tel que $\R^q f_* ℱ=0$ pour $q ≥ N$.
\end{enumerate}
\end{théorème2}

\subsubsection{}De façon équivalente, le morphisme
$\R f_*:\detalepluscoefs{Y}{𝐙/n}→\detalepluscoefs{X}{𝐙/n}$
induit un morphisme $\detalebccoefs{Y}{𝐙/n}→\detalebccoefs{X}{𝐙/n}$
entre les sous-catégories de complexes à cohomologie bornée et
constructible.

\subsection{Remarques}
\subsubsection{Organisation de l'exposé}L'énoncé ci-dessus est la conjonction d'un résultat de
\emph{constructibilité} (i) et d'un résultat d'\emph{annulation} (ii).
Dans le §\ref{13-constructibilité sans annulation}, nous présentons une démonstration
de la constructibilité qui ne requiert pas
la forme forte du théorème d'uniformisation mais seulement
la forme faible (\refext{VII}{7-uniformisation}).
Les ingrédients clefs supplémentaires sont le théorème de pureté
absolu, le théorème de constructibilité générique (dû à P. Deligne)
et la descente cohomologie orientée.
Au paragraphe \ref{13-compléments annulation}, nous donnons une
démonstration de résultat d'annulation pour les schémas \emph{de
dimension finie}, qui complète la démonstration du théorème \ref{13-constructibilité} pour ces schémas.
Le cas général est traité en §\ref{13-section constructibilité forte}, en s'appuyant sur le
théorème d'uniformisation premier à $ℓ$ (\refext{IX}{9-unif}), où $ℓ$ est un nombre premier divisant $n$.
Enfin, nous étendons ce résultat d'abord au cas des coefficients
$ℓ$-adiques, où $ℓ$ est un nombre premier inversible sur
les schémas considérés, puis au cas des champs (comme coefficients).

\subsubsection{Terminologie et notations}
\label{13-terminologie et notations}
Nous dirons d'un complexe $𝒦∈\ob\,\detalepluscoefs{Y}{Λ}$,
où $Λ$ est un anneau fini, est \emph{constructible} s'il appartient à
$\ob\,\detalebccoefs{Y}{Λ}$, c'est-à-dire si ses faisceaux de cohomologie sont
constructibles, nuls en grands degrés. Lorsque $n ≥ 1$ est fixé
et que cela ne semble pas créer de confusion nous notons $Λ$
l'anneau $𝐙/n$. De même, un complexe $𝒦$ sur un schéma $X$
étant donné, nous noterons souvent encore $𝒦$ ses images
inverses sur différents $X$-schémas.

\section{Constructibilité via l'uniformisation locale faible}
\label{13-constructibilité sans annulation}

Dans cette section, on démontre \ref{13-constructibilité} (i),
dont on reprend les notations.

\subsection{Réductions}
Les réductions suivantes sont classiques : cf. p. ex. \sga{4}{XVI}{4.5}.

\subsubsection{Réduction au cas où le faisceau $ℱ$ est constant}
\label{13-réduction cas constant}
D'après \sga{4}{IX}{2.14 (ii)}, le faisceau $ℱ$ s'injecte
dans une somme finie $𝒢= ⨁_{i ∈ I} {g_i}_* C_i$ d'images directes par
des morphismes finis $g_i$ de faisceaux en $𝐙/n$-modules constants constructibles $C_i$.
On peut supposer $ℱ=𝒢$. Cela résulte d'une part du fait qu'un
sous-quotient d'un faisceau constructible est
constructible\footnote{Pour le voir, on peut utiliser le
fait qu'un faisceau est constructible si et seulement si il est nœthérien, cf.
\sga{4}{IX}{2.9 (i)}.} et d'autre part de la suite exacte longue
de cohomologie associée au triangle
\[
\R f_*ℱ→\R f_* 𝒢→\R f_* (𝒢/ℱ)\dessusdessous{+1}{→}.
\]
Enfin, on peut supposer $ℱ$ constant constructible car
on peut supposer l'ensemble $I$ être un singleton
et l'égalité $\R f_*(g_* C)=\R (f ∘ g)_* g$,
où $g$ est un morphisme fini, nous permet de supposer $g=\Id$.
Décomposant $n$ en produit, on se ramène au cas où $ℱ$ est
un faisceau constant $𝐅_ℓ$, le nombre premier $ℓ$ étant inversible
sur $X$ (cf. p. ex. \sgaquatreetdemi{Th. finitude}{2.2 b)}).

\subsubsection{Réduction au cas où le morphisme $f$ est une immersion ouverte}
\label{13-réduction cas immersion ouverte}
Un faisceau sur le schéma nœthérien donc cohérent $X$
étant constructible si et seulement si il l'est localement
pour la topologie de Zariski (\sga{4}{IX}{2.4 (i)}),
on peut supposer $X$ affine. On utilise ici le fait
trivial que la formation des images directes commute
au changement de base par un ouvert de Zariski.
On peut également supposer $Y$ affine ; cela résulte
par exemple de l'analogue faisceautique
\[
E₁^{p,q}=\R^q {f_p}* (ℱ_{|Y_p}) ⇒ \R^{p+q} f_* ℱ
\]
de la suite spectrale de Leray (\cite{Hodge3@Deligne}, 5.2.7.1),
où les $f_p:Y_p → X$ sont déduits de $f$ et d'un hyperrecouvrement
Zariski $Y_{∙} → Y$.
Le morphisme $f:Y→X$ est alors affine donc quasi-projectif,
et le théorème de constructibilité étant connu pour les morphismes
propres (\sga{4}{XIV}{1.1}), on peut supposer que $f$ est une
\emph{immersion ouverte} dominante. (On pourrait également utiliser
le théorème de compactification de Nagata.)
Conformément à l'usage, nous noterons dorénavant $j:U→X$ le morphisme $f$.

\subsection{Fin de la démonstration du théorème \ref{13-constructibilité} (i)}
\label{13-fin première démo}

\subsubsection{}Soit $q≥0$ un indice pour lequel on souhaite montrer
que le faisceau $\R^q j_* 𝐅_ℓ$ est constructible. On rappelle que $j$
est une immersion ouverte $U ↪ X$ et $ℓ$ est un nombre premier
inversible sur $X$.
Il résulte du critère de constructibilité \sga{4}{ix}{2.4.(v)} qu'il suffit de démontrer que
pour toute immersion fermée $g:Z↪ X$, le faisceau $g^*\R^q j_*𝐅_ℓ$
est constructible \emph{sur un ouvert dense} de $Z$.
Le théorème d'uniformisation locale (\refext{VII}{7-uniformisation}),
joint à la méthode classique de construction d'hyperrecouvrements
(\cite{Hodge3@Deligne}, §6.2), a pour corollaire immédiat
le fait suivant.

\begin{théorème2}
\label{13-cb}
Il existe un hyperrecouvrement pour la topologie des altérations
$ε_{∙} : X_{∙}→ X$ satisfaisant les conditions suivantes :
\begin{enumerate}
\item pour chaque $i≤q+1$, le schéma $X_i$ est régulier connexe,
\item pour chaque $i≤q+1$, l'image inverse $U_i$ de $U$ dans $X_i$
est le complémentaire du support d'un diviseur strictement à croisements normaux,
\end{enumerate}
\end{théorème2}

\subsubsection{}
Notons $Z_{∙}$, $g_{∙}$ et $j_{∙}$ le schéma et les morphismes
simpliciaux qui se déduisent de $Z$, $g$ et $j$ respectivement
par le changement de base $X_{∙}→X$.
Il résulte du théorème de pureté absolue
(\refext{XVI}{16-theoreme-purete-absolue}) que le complexe
$g_{∙}^*\R {j_{∙}}_* Λ$ sur $Y_{∙}$
est à cohomologie constructible en degré $≤q+1$.
Par ailleurs, il résulte du théorème de constructibilité
générique \sgaquatreetdemi{Th. finitude}{1.9 (i)}
— appliqué aux morphismes $ε_p:Z_p → Z$ et aux
complexes $g_p^*\R {j_p}_* 𝐅_ℓ$ — et de la suite spectrale
rappelée ci-dessus qu'il existe un ouvert
dense de $Z$ au-dessus duquel le tronqué en degrés inférieurs
ou égaux à $q$ de l'image directe $\R {ε_{Z ∙}}_*
\big(g_{∙}^*\R {j_{∙}}_* 𝐅_ℓ\big)$ est \emph{constructible}.
D'après \refext{XII}{12-cb orienté} (i)
cette image directe est isomorphe à $g^*\R j_*𝐅_ℓ$.
Le faisceau $g^* \R j_* 𝐅_ℓ$ est donc constructible
sur un ouvert dense de $Z$. CQFD.

\subsection{Compléments}
\label{13-compléments annulation}

\begin{théorème2}
\label{13-finitude sur base dimension 1}
Soient $S$ un schéma nœthérien, $f: X → Y$ un morphisme de type fini
entre $S$-schémas de type fini et $n$ un entier inversible sur $S$.
Supposons l'une des deux conditions suivantes satisfaite :
\begin{enumerate}
\item le schéma $S$ est de dimension $1$ ;
\item le schéma $S$ est \emph{local} de dimension $2$.
\end{enumerate}
Alors, pour faisceau de $𝐙/n$-modules constructible $ℱ$ sur $X$
les faisceaux image directe $\R^q f_* ℱ$ sont constructibles
et nuls pour $q ≫ 0$.
\end{théorème2}

\begin{remarque2}
On verra en \refext{XIX}{19-prop-cexgabber} qu'il existe un contre-exemple
à la constructibilité lorsque $S$ est nœthérien
de dimension $2$ (non local). Ceci résulte de l'existence
d'une surface régulière et d'un diviseur possédant une
infinité de points doubles.
Il serait intéressant de construire un contre-exemple
à l'énoncé de constructibilité précédent lorsque $S$
est \emph{local} (nœthérien) de dimension $3$,
ou bien de montrer qu'il n'en existe pas.
\end{remarque2}

\begin{démo}[Esquisse de démonstration]
(i) D'après le théorème de constructibilité générique
(\sgaquatreetdemi{Th. finitude}{1.9 (i)}), il existe un
ouvert dense de $S$ au-dessus duquel le résultat est acquis.
On peut donc supposer le schéma $S$ \emph{local}.
Il est également loisible de le supposer strictement
hensélien. Par restriction à ses composantes irréductibles,
on peut finalement supposer $S$ local \emph{intègre} (de
dimension $1$). Soit $S′ → S$ le morphisme de normalisation.
C'est un homéomorphisme universel de sorte que l'on peut
remplacer $S$ par $S′$. Or, ce dernier schéma est nœthérien régulier, de dimension $1$.
La conclusion résulte alors du théorème de constructibilité
\sgaquatreetdemi{Th. finitude}{1.1}
et du théorème de finitude de la dimension
cohomologique \sga{4}{X}{3.2, 4.4}. (Voir aussi \cite[2.4]{PV@Illusie}
et \refext{XVIII}{18A-GABBER'S BOUND}.)

(ii) Soit $s$ le point fermé de $S$.
Le schéma $S-\{s\}$ étant de dimension $1$ le résultat
est acquis au-dessus ce cet ouvert.
Soit $ℱ$ un $𝐙/n$-Module constructible sur $X$ et considérons
un triangle distingué
\[
𝒦 → ℱ → \R j_* j^* ℱ,
\]
où $j$ est l'immersion ouverte $X -X_s ↪ X$.
Notons que $𝒦$ est à support dans $X_s$ et
constructible si $\R j_* j^* ℱ$. Appliquant
le foncteur $\R f_*$ au triangle précédent
et utilisant la finitude sur $X-X_s$ (resp. $X_s$),
on est ramené à montrer la constructibilité
des images directes par les immersions ouvertes
$X-X_s ↪ X$ et $Y-Y_s ↪ Y$.
On utilise alors le morphisme de complétion $\chap{S} → S$
et le théorème de comparaison de Gabber-Fujiwara
(\cite[6.6.4]{TUBULAR@Fujiwara}) pour se ramener au cas
où le schéma local $S$ est complet, donc excellent.
\end{démo}

\subsubsection{}Reprenons les notations du théorème \ref{13-constructibilité}
et supposons le schéma $X$ de dimension finie.
Le (i) de \emph{loc. cit.} joint au théorème de
Lefschetz affine \refext{XV}{15-theoreme-principal}\footnote{Le lecteur
constatera que cette référence à un exposé ultérieur ne génère pas de cercle vicieux.}
entraînent le complément suivant, qui sera amélioré dans la section suivante.

\begin{proposition2}
\label{13-annulation si dimension finie}
Soit $X$ un schéma nœthérien quasi-excellent
de \emph{dimension finie}. Soit $f:Y→X$ un morphisme de type fini.
Pour tout entier $n≥1$ inversible sur $X$, le foncteur $\R f_*:
\detalepluscoefs{Y}{𝐙/n}→\detalepluscoefs{X}{𝐙/n}$ est de dimension
cohomologie finie. En particulier, il induit un foncteur
de $\detalepluscoefs{Y}{𝐙/n}$ dans $\detalepluscoefs{X}{𝐙/n}$.
\end{proposition2}

\begin{démo}
Supposons $X$ de dimension finie $d$. On peut supposer $X$ affine.
Quitte à considérer un hyperrecouvrement Zariski de $Y$ par
des schémas affines, on peut également supposer ce schéma
affine. Soit maintenant $N$ un majorant de la dimension des fibres
de $f$. La dimension cohomologique du foncteur
image directe par $f$ est au plus $d+N$. En effet, si $\sur{x}$ est
un point géométrique de $X$, et $ℱ$ un $𝐙/n$-faisceau
constructible sur $Y$, on a $(\R f_* ℱ)_{\sur{x}}=\R Γ(Y×_X
X_{(\sur{x})},ℱ)$. Les schémas $X′= X_{(\sur{x})}$
et $Y′= Y×_X X_{(\sur{x})}$ admettent respectivement les fonctions
de dimension $δ_{X′}: x′↦ \dim(\sur{\{x′\}})$ et la fonction
induite $δ_{Y′}$ définie en \refext{XIV}{14-coro-fonction-dim-induite-tf}.
Notons que $δ_{X′}$ est bornée par $d$ et $δ_{Y′}$
par $d+N$. Il résulte donc du théorème de Lefschetz affine
(sous la forme \refext{XV}{15-prop-devissage-classique})
que $H^q(Y′, ℱ)=0$ pour $q>d+N$.
\end{démo}

\begin{remarque2}
On verra en \refext{XVIII}{18A-GABBER'S BOUND} que l'on a un résultat d'annulation
sous la seule hypothèse que $X$ est nœthérien de dimension finie :
si $X$ est un schéma nœthérien strictement local hensélien
de dimension $d>0$ et $n$ est inversible sur $X$, alors
tout ouvert de $X$ est de $n$-dimension cohomologique au plus
$2d-1$.
\end{remarque2}

\subsubsection{Constructibilité des images directes dans le cas non abélien}
Quitte à remplacer la réduction \ref{13-réduction cas constant}
par \sga{1}{XIII}{§3, (4)}, l'usage de
\sga{4}{XIV}{1.1} en \ref{13-réduction cas immersion ouverte}
par \sga{1}{XIII}{6.2}, le théorème de pureté absolu
par  \sga{1}{XIII}{2.4}, le théorème de finitude \sgaquatreetdemi{Th. finitude}{1.9 (i)}
par \cite[2.2]{Alterations@Orgogozo} et enfin
\refext{XII}{12-cb orienté} (i)
par \refext{XII}{12-cb orienté} (ii), on obtient
essentiellement par la même méthode une démonstration du
théorème suivant.

\begin{théorème}[\refext{XXI}{21-groupes}]\label{13-constructibilité non abélienne}
Soient $X$ un schéma nœthérien \emph{quasi-excellent}, $f:Y→X$ 
un morphisme de type fini, et $\mathsf{L}$ l'ensemble des nombres
premiers inversibles sur $X$. 
Pour tout champ en groupoïdes constructible ind-$\mathsf{L}$-fini sur $Y_{\et}$
le champ $f_*𝒞$ est constructible.
\end{théorème}

\section{Constructibilité et annulation via l'uniformisation locale première à $ℓ$}
\label{13-section constructibilité forte}

Dans cette section, on démontre le théorème \ref{13-constructibilité}.
Constructibilité (i) et annulation (ii) sont établis simultanément.

\subsection{Réduction au cas d'une immersion ouverte et de la finitude hors d'un lieu de codimension donnée}
\label{13-hr}

\subsubsection{}Comme en \ref{13-terminologie et notations},
posons $Λ=𝐙/n$, où $n$ est l'entier inversible sur $X$ de l'énoncé.
Pour chaque entier $c≥0$, considérons la propriété $(P_c)$ suivante :
\begin{quote}
\emph{Pour tout schéma quasi-excellent nœthérien $X$, toute immersion ouverte
dominante $j:U↪X$ et tout complexe $𝒦∈\ob\,\detalebccoefs{U}{Λ}$,
il existe un fermé $T↪X$ de codimension strictement
supérieure à $c$ tel que $(R j_*𝒦)_{|X-T}$ appartienne à $\ob\,\detalebccoefs{(X-T)}{Λ}$.}
\end{quote}

\subsubsection{}
\label{13-Pc impliquent théorème}
La conjonction des énoncés $(P_c)$ pour
chaque $c ≥ 0$, entraîne le théorème.
Soit en effet une paire $(f,ℱ)$ comme dans l'énoncé du théorème. D'après
le théorème de compactification de Nagata, il existe
une immersion ouverte $j:X↪\sur{X}$ et un morphisme propre
$\sur{f}:\sur{X}→Y$ tels que $f=\sur{f}j$. La formule
de composition $\R f_*=\R \sur{f}_* \R j_{*}$ et le théorème de finitude pour les
morphismes propres, nous ramènent à démontrer la
constructibilité de complexe $𝒦=\R j_*ℱ$.
(Voir aussi \ref{13-réduction cas immersion ouverte}.)
La conclusion résulte alors du lemme suivant.

\begin{lemme2}
Soient $X$ un schéma nœthérien et $𝒦∈\ob\,\detalepluscoefs{X}{Λ}$.
Supposons que pour tout entier $c≥0$, il existe un fermé $T_c$ de
codimension strictement supérieure à $c$ tel que
$\mc{K}_{|X-T_c}∈\ob\,\detalepluscoefs{(X-T_c)}{Λ}$.
Alors, $\mc{K}∈\ob\,\detalebccoefs{X}{Λ}$. 
\end{lemme2}

\begin{démo}
Le schéma $X$ étant
nœthérien, ses localisés sont de dimension finie et pour toute suite
de fermés $(T_c)_{c∈𝐍}$ comme dans l'énoncé,
on a $X=⋃_c(X-T_c)$. D'autre part, le schéma $X$ étant quasi-compact,
il est recouvert par un nombre \emph{fini} des ouverts $X-T_c$.
La conclusion résulte alors du fait que si $U,U'$ sont deux ouverts de $X$ tels que
$\mc{K}_{|U}∈\ob\,\detalebccoefs{U}{Λ}$,
$\mc{K}_{|U'}∈\ob\,\detalebccoefs{U'}{Λ}$, on a également 
$\mc{K}_{|U∪U'}∈\ob\,\detalebccoefs{(U∪U')}{Λ}$. 
\end{démo}

\subsubsection{}
Nous allons démontrer la propriété
$(P_c)$ ci-dessus par récurrence sur $c$. Insistons sur le fait que
\emph{le schéma $X$ et le complexe $𝒦$ sont
variables}. Pour $c=0$, cette propriété est triviale : prendre $T=X-U$.
Soit $c≥1$ et supposons la propriété établie au cran $c-1$.
On souhaite la démontrer au cran $c$.

\subsection{Récurrence : l'ingrédient clef et une première réduction}

\subsubsection{}
\label{13-uniformisation}
D'après le théorème d'uniformisation première à $ℓ$
(\refext{IX}{9-unif}) et le théorème de la forme standard (\refext{II}{2-forme standard pour annulation}),
il existe une famille finie indexée, par un ensemble $I$ d'éléments $i$, de diagrammes commutatifs
\begin{center}
\begin{tikzpicture}[auto]
 \matrix(diag)[matrix of math nodes,column sep=4em,row sep=3ex]
{ |(X3)| X‴_i & & |(Y)| Y=∐_{j ∈ J} Y_j & \\
  |(X2)| X″_i & |(X1)| X′ & |(X)| X & |(U)| U\\};
 \draw[->] (X3) -- node[swap]{fini, plat, surjectif} (X2);
 \draw[->] (X3) -- node{degré premier à $ℓ$} (X2);
 \draw[->] (X3) -- (Y);
 \draw[->] (Y) -- (X);
 \draw[->] (X2) -- node[swap]{étale} (X1);
 \draw[->] (X1) -- node{$p$} (X);
 \draw[->] (X1) -- node[swap]{propre, birationnel} (X);
 \draw[left hook->] (U) -- node[swap]{$j$} (X);
 \end{tikzpicture}
\end{center}
où, en plus des propriétés indiquées ci-dessus,

— la famille $(X″_i→X')$ est couvrante pour la topologie
étale complètement décomposée ;

— les schémas $Y_j$, $j ∈ J$, sont \emph{réguliers} ;

— l'image inverse de $U$ dans $Y_j$ est le complémentaire
d'un diviseur strictement à croisements normaux.

\subsubsection{}
Soit $(j,𝒦)$ une paire comme dans l'énoncé
de la propriété $P_c$ (\ref{13-hr}).
Nous verrons seulement plus tard que l'on peut supposer $𝒦=Λ$.
D'après l'hypothèse de récurrence appliquée aux paires
$(j,𝒦)$ et $(j',𝒦)$, où $j'$ est l'immersion ouverte
de $U'=U×_X X'$ dans $X'$, il existe deux fermés $T↪X$ et $T'↪X'$ 
de codimension $≥c$ tels que les complexes
$\R j_* 𝒦$ et $\R j'_*𝒦$ soient constructibles sur les ouverts complémentaires
correspondants. Le fermé $T$ n'ayant qu'un nombre fini de points maximaux
et l'énoncé à démontrer — la constructibilité hors d'un
fermé de codimension $>c$ — étant un problème local au
voisinage de ces points, on peut supposer $T$ irréductible,
de codimension $c$, de point générique noté $η_T$. Soit $η'$
un point maximal de $T'$.
Si l'image par $p$ de $η'$ n'est pas égale à $η_T$, la composante irréductible
correspondante de $T'$ disparaît après localisation
(Zariski) au voisinage de $η_T$.
Compte tenu du fait que $T$ est de codimension $c$ et $T'$ de codimension au
moins égale, toute composante irréductible $T'_α$ de $T'$ dominant $T$ est nécessairement de
dimension égale à celle de $T$, et le morphisme induit $T'_α→T$ est
génériquement fini. Quitte à se restreindre à un voisinage
ouvert de $η_T$ dans $X$, on peut finalement supposer que $T'$
est une somme $\coprod_α T'_α$,
où les $T'_α$ sont irréductibles et les morphismes $T'_α→T$ sont \emph{finis},
surjectifs.

\subsection{Notation : le complexe $ψ_f(g,𝒦)$}

\subsubsection{}
Pour tout $X$-schéma $f:X₁→X$ et tout $X₁$-schéma $g:X₂ → X₁$, notons
$h$ le morphisme composé $X₂→X₁→X$ et $j₁$ l'immersion
ouverte $U₁=X₁×_X U↪X₁$ déduite de $j$ par changement de base.
Considérons le complexe de faisceaux sur $X$,

\begin{minipage}[c]{.65\linewidth}
\[
ψ_f(g,𝒦) ≔ \R g_*\big(h^*(\R {j_Y}_*𝒦)\big).
\]
\end{minipage} \hfill
\begin{minipage}[c]{.46\linewidth}
\begin{tikzpicture}[auto]
 \matrix(diag)[matrix of math nodes]
{ |(X2)| X₂ & \\
  |(X1)| X₁ & |(U1)| U₁ \\
  |(X)| X & |(U)| U \\};
 \draw[->] (X2) -- node{$g$} (X1);
 \draw[->] (X1) -- node{$f$} (X);
 \draw[left hook->] (U1) -- node[swap]{$j₁$} (X1);
 \draw[left hook->] (U) -- node[swap]{$j$} (X);
 \draw[->] (U1) -- (U);
 \draw[->] (X2) to [in=130,out=210] node[swap]{$h$} (X);
 \end{tikzpicture}
\end{minipage}

Ci-dessous, le morphisme $g$ sera le plus souvent une immersion
fermée, qui sera parfois supprimée de la notation, ainsi que $f$,
si cela ne semble pas induire de confusion. Par exemple, $ψ(X,𝒦)=\R j_*𝒦$.

\subsubsection{}
La formation du complexe $ψ$ est fonctorielle en le sens suivant :
pour tout diagramme commutatif
\begin{center}
\begin{tikzpicture}[auto]
 \matrix(diag)[matrix of math nodes]
{ & |(X2)| X₂ & & |(X2p)| X₂′ & \\
  |(U1)| U₁ & |(X1)| X₁ & & |(X1p)| X₁′ & |(U1p)| U₁′ \\
        & & |(X)| X & & \\};
 \draw[->] (X2) -- node[swap]{$g$} (X1);
 \draw[->] (X1) -- node[swap]{$f$} (X);
 \draw[->] (X2p) -- node{$g′$} (X1p);
 \draw[->] (X1p) -- node{$f′$} (X);
 \draw[<-] (X2) -- node{$m$} (X2p);
 \draw[<-] (X1) -- node{$n$} (X1p);
 \draw[right hook->] (U1) -- node{$j₁$} (X1);
 \draw[left hook->] (U1p) -- node[swap]{$j₁′$} (X1p);
 \end{tikzpicture}
\end{center}
le morphisme de changement de base (adjonction)
$n^* \R {j_Y}_*𝒦 → \R {j_{Y'}}_*𝒦$ induit un morphisme
\[
ψ_f(g,𝒦)→ ψ_{f'}(g',𝒦).
\]

\subsection{Seconde localisation}

\subsubsection{}Nous dirons qu'un morphisme d'une catégorie dérivée $\mathsf{D}^+(𝒯,Λ)$,
où $𝒯$ est le topos étale d'un schéma, est un
\emph{$\dbc$-isomorphisme} ou isomorphisme \emph{modulo} $\dbc$,
s'il a un cône dans $\dbc(𝒯,Λ)$. Cela revient
d'après \cite[2.1.35]{TC@Neeman} à supposer que
la flèche induite dans la catégorie triangulée quotient
$\mathsf{D}^+(𝒯,Λ) ∕ \dbc(𝒯,Λ)$ est un \emph{isomorphisme}.
Notons que dans la terminologie d'\emph{op. cit.},
la sous-catégorie $\dbc(𝒯,Λ)$ est \emph{épaisse}.
La localisation considérée ici (due à J.-L. Verdier)
est l'analogue triangulé de celle considérée par J.-P. Serre
dans le cas des catégories abéliennes.

\begin{proposition2}
Quitte à se restreindre au voisinage de $η_T$, on peut supposer que le morphisme
d'adjonction
\[
ψ_{\Id}(T↪X,ℱ)→ψ_p(T'↪X',ℱ)
\]
est un $\dbc$-isomorphisme.
\end{proposition2}

Notons que le terme de droite, $ψ_p(T'↪X',ℱ)$, est isomorphe
à la somme directe $⊕_α ψ_p(T'_α↪X',ℱ)$.

\begin{démo}
Soit $p_U$ le morphisme induit par $p$ au-dessus de l'ouvert $U$ de $X$ ;
c'est un isomorphisme au-dessus d'un ouvert $W$ de $U$.
Notons $i$ l'immersion fermé du complémentaire $Z=U-W$ dans $U$.
On a sur $U$ un triangle distingué
\[
𝒦→ \R {p_{U}}_* {p_{U}}^* 𝒦 → i_*ℋ \dessusdessous{+1}{→},
\]
où $ℋ$ est constructible sur $Z$, d'après le théorème de finitude pour
le morphisme propre $p_U$.
Il résulte du théorème de changement de base propre pour $p$ que le triangle
distingué précédent devient, après application du foncteur $ψ_{\Id}(T↪X,-)$,
le triangle distingué de complexes supportés sur $T$ suivant :
\[
ψ_{\Id}(T↪X,𝒦)→ ψ_p(p^{-1}(T)↪X',𝒦) → ψ_{\Id}(T↪X,i_*ℋ)  \dessusdessous{+1}{→}.
\]

Première étape. Nous allons commencer par montrer que la première flèche
est génériquement sur $T$ un $\dbc$-isomorphisme.
(« Génériquement sur $T$ » : quitte à se restreindre à un
voisinage Zariski convenable de $η_T$.)
Soient en effet $\sur{Z}$ l'adhérence de $Z$ dans $X$, $\sur{j}:Z↪\sur{Z}$
l'immersion ouverte et $\sur{i}:\sur{Z}↪X$ l'immersion
fermée, représentés dans le diagramme ci-dessous.

\begin{center}
\begin{tikzpicture}[auto]
 \matrix(diag)[matrix of math nodes,column sep=5em,row sep=7ex]
{ |(U)| U & |(Z)| Z=U-W    \\
  |(X)| X & |(Zb)| \sur{Z} \\};
 \draw[left hook->] (Z) -- node[swap]{$i$} (U);
 \draw[left hook->] (Zb) -- node[swap]{$\sur{i}$} (X);
 \draw[left hook->] (U) -- node[swap]{$j$} (X);
 \draw[left hook->] (Z) -- node{$\sur{j}$} (Zb);
 \end{tikzpicture}
\end{center}
La restriction à $T$ du complexe $ψ(T,i_*ℋ)$ — dont on veut montrer
qu'elle est génériquement $\dbc$-nulle — est isomorphe à la
restriction du complexe ${\sur{i}}_*\R \sur{j}_* ℋ$.
Le fermé $\sur{Z}$ étant de codimension $≥1$ dans $X$, car $W$ est partout
dense dans $X$, l'hypothèse de récurrence
pour la paire $(\sur{j},ℋ)$ entraîne
immédiatement le résultat.

Deuxième étape. Pour conclure, il nous faut maintenant montrer que le morphisme d'adjonction
$ψ_p(p^{-1}(T),𝒦)→ψ_p(T',𝒦)$, à travers lequel le morphisme
$ψ_{\Id}(T,𝒦)→ψ_p(T',𝒦)$ de l'énoncé se factorise est,
génériquement sur $T$, un $\dbc$-isomorphisme.
Sur le fermé $T_p′=p^{-1}(T)$ de $X′$, considérons
la restriction $ℒ=(\R j'_{*}𝒦)_{|p^{-1}(T)}$ de l'image directe
par $j′$ de $𝒦$, et le triangle distingué
\[
(T_p′ - T′ ↪ T_p′)_! ℒ_{|T_p′-T′} → ℒ → (T′ ↪ T_p′)_* ℒ_{|T′} \dessusdessous{+1}{→}
\]
constitué de ses prolongements par zéro.
Rappelons que $j'$ désigne l'immersion ouverte de $U'$
dans $X'$. Par définition de $T'$, le premier complexe est
constructible ; il en est donc de même de son image directe
(dérivée) par le morphisme \emph{propre} $p_T$.
Or, l'image directe de la seconde flèche par $p_T$ n'est autre
que le morphisme d'adjonction $ψ_{p}(T′_p,𝒦)→ ψ_p(T',𝒦)$. CQFD.
\end{démo}

\subsection{Construction d'une rétraction}

\subsubsection{}
Quitte à rétrécir $X$ un peu plus encore, on peut supposer que pour tout $α$
(on rappelle que $T′= ∐_α T′_α$), il existe un indice $i_α$ tel que le morphisme étale $X″_{i_α}→X'$
ait une section $σ_α$ au-dessus de $T'_α$. Cela résulte
du fait que la famille $(X″_i→X')_i$ est complètement
décomposée, de sorte qu'une section existe au voisinage du point
générique de $T'_α$ (\refext{II}{2-section cd}).
La propreté du morphisme dominant $X'→X$ permet de
déduire l'existence d'un ouvert convenable de $X$
de celle d'un ouvert de $X′$.

\subsubsection{}
Pour simplifier les notations, on pose pour chaque indice $α$,
$X″_α=X″_{i_α}$, $X‴_α=X‴_{i_α}$ et on note $T″_α⊂X″_α$
l'image de $T'_α$ par une section $σ_α$ comme ci-dessus,
et enfin $T‴_α ⊆ X‴_α$ l'image inverse de $T″_α$
par le morphisme fini $X‴_α → X″_α$.

\begin{proposition2}
Le morphisme d'adjonction $ψ(T'_α ↪ X′_α,𝒦)→ ψ(T″_α ↪ X″_α,𝒦)$ est un \emph{isomorphisme}.
\end{proposition2}

Bien entendu, les complexes ci-dessus sont calculés en
munissant les schémas $X′_α$ et $X″_α$ de la structure
de $X$-schéma évidente. Nous nous autoriserons dorénavant
cet abus de notation.

\begin{démo}
Résulte du fait que le morphisme $X″_α→X'$ est étale.
\end{démo}

\begin{proposition2}
Le morphisme d'adjonction $ψ(T″_α↪ X′_α,𝒦)→ ψ(T‴_α↪ X″_α,𝒦)$ a un inverse à gauche.
\end{proposition2}

\begin{démo}
Considérons le diagramme à carrés cartésiens suivant :
\begin{center}
\begin{tikzpicture}[auto]
\matrix(diag)[matrix of math nodes,column sep=4em,row sep=5ex]
{|(T3)| T‴_α & |(X3)| X‴_α & |(U3)| U‴_α \\
 |(T2)| T″_α & |(X2)| X″_α & |(U2)| U″_α \\};
\draw[->] (T3) -- node[swap]{$π_T$} (T2) ;
\draw[->] (T3) -- (X3) ;
\draw[->] (T2) -- (X2) ;
\draw[->] (X3) -- node{$π$} (X2) ;
\draw[->] (U3) -- node[swap]{$j‴_α$} (X3) ;
\draw[->] (U2) -- node{$j″_α$} (X2) ;
\draw[->] (U3) -- node{$π_U$} (U2) ;
\end{tikzpicture}
\end{center}
où $U″_α=U ×_X X″_α$, de même pour $U‴_α$, et $π: X‴_α → X″_α$ est
comme en \ref{13-uniformisation}.
En particulier, le morphisme $π_U$ est fini, plat, et de
degré générique premier à $ℓ$, de sorte que le morphisme
composé
\[𝒦→{π_U}_*{π_U}^*𝒦\dessusdessous{\tr}{→}𝒦\]
est la multiplication par le degré, donc inversible.
Appliquons le foncteur $\R j″_α$.
Par composition des images directes, le terme du milieu
est $π_* \R {j‴_α}_*𝒦$, où l'omet le foncteur image
inverse de la notation (\ref{13-terminologie et notations}).
D'après le théorème de changement de base pour les morphismes finis, sa
restriction au fermé $T″_α$ est isomorphe à
${π_T}_* \big((\R {j‴_α}_*𝒦)_{|T‴_α}\big)$.
En poussant les faisceaux sur $X$ par le morphisme
$T″_α→X$, la suite précédente devient donc
\[
ψ(T″_α↪ X″_α,𝒦)→ ψ(T‴_α↪ X‴_α,𝒦)→ψ(T″_α↪ X″_α,𝒦)
\]
et la composition de ces flèches est un isomorphisme.
\end{démo}

\subsection{Cas des coefficients constants : utilisation du théorème de pureté}

\subsubsection{}
Posons $T‴=\coprod T‴_α$, $X‴=∐X‴_{α}$ et considérons le diagramme
commutatif de morphismes d'adjonction, complété du morphisme trace :

\begin{center}
{\small
\begin{tikzpicture}[auto]
\matrix(diag)[matrix of math nodes,column sep=1.5em,row sep=3ex]
{|(X)| ψ(T ↪ X,𝒦) & |(X1)| ψ(T' ↪ X', 𝒦) &
|(X2)| ψ(T ″ ↪ X″,𝒦) & |(X3)| ψ(T‴ ↪ X‴,𝒦) & |(X2bis)| ψ_p(T″ ↪ X ″, 𝒦).\\
 & & & |(Y)| ψ(T‴ → Y,𝒦) & & \\};
\draw[->,densely dashed] (X) -- (X1) ;
\draw[->,densely dashed] (X1) -- (X2) ;
\draw[->] (X2) -- (X3) ;
\draw[->,densely dashed] (X2) to [out=20,in=160] (X2bis) ;
\draw[->] (X3) -- node[swap]{$\tr$} (X2bis) ;
\draw[->] (X) to [out=-20,in=160] (Y) ;
\draw[->] (Y) -- (X3) ;
\end{tikzpicture}
}
\end{center}
D'après les trois propositions précédentes,
les flèches en tirets deviennent des isomorphismes
modulo $\dbc$. Si le complexe $ψ(T‴ → Y,𝒦)$
est \emph{constructible}, c'est-à-dire nul
modulo $\dbc$, il en résulte que $ψ(T ↪ X,𝒦)$
— ou, de façon équivalente, $(\R j_* 𝒦)_{|T}$ —
est également constructible.

\begin{proposition2}
Le complexe $ψ(T‴ → Y, Λ)$ est constructible.
\end{proposition2}

\begin{démo}
Le morphisme composé $T‴ → X$ étant fini, il suffit
de démontrer que le complexe $\R {j_Y}* Λ$ est constructible.
Cela résulte des hypothèses faites en \ref{13-uniformisation}
et du théorème de pureté \refext{XVI}{16-theoreme-purete-absolue}.
\end{démo}

\subsection{Réduction au cas des coefficients constants}
\label{13-ops constant}

\subsubsection{}
Pour achever la démonstration du théorème \ref{13-constructibilité}, il nous faut maintenant
montrer que la propriété $(P_c)$ de §\ref{13-hr},
où $c$ est fixé, résulte du cas particulier
où $𝒦=Λ$ et de l'énoncé $P_{c-1}$.

\subsubsection{}Commençons par observer que l'on peut supposer $𝒦$
concentré en degré $0$, c'est-à-dire être un \emph{faisceau} constructible,
que nous noterons dorénavant $ℱ$. L'ensemble des faisceaux constructibles
satisfaisant à la propriété à établir au rang $c$ est, à $X$ fixé, stable
par extension et facteur direct. D'après
\sga{5}{i}{3.1.2}, on peut supposer $ℱ=π_*k_!Λ$ où $π:U'→U$ est un morphisme
\emph{fini} et $k:W↪U'$ une immersion ouverte, avec $U'$ \emph{intègre}.
D'après le théorème principal de Zariski (\egalong{iv}{3}{8.12.6}),
le morphisme composé $U'→X$, \emph{quasi-fini}, se factorise en une immersion
ouverte $j':U'↪X'$ suivie d'un morphisme fini $\sur{π}:X'→X$.
\begin{center}
\begin{tikzpicture}[auto]
 \matrix(diag)[matrix of math nodes,column sep=4em,row sep=5ex]
{ |(W)| W & |(Up)| U′ & |(Xp)| X′ \\
          & |(U)|  U  & |(X)|  X \\};
 \draw[right hook->] (W)  -- node{$k$} (Up);
 \draw[right hook->] (Up) -- node{$j′$} (Xp);
 \draw[->] (Up) -- node[swap]{$π$} (U);
 \draw[->] (Xp) -- node{$\sur{π}$} (X);
 \draw[right hook->] (U)  -- node{$j$} (X);
 \end{tikzpicture}
\end{center}
Le complexe $\R j_* π_* k_! Λ$, dont on s'interroge
sur la constructibilité, est isomorphe au complexe $\sur{π}_*\R j'_* k_!Λ$.
En vertu du lemme suivant, on peut supposer $X'=X$.

\begin{lemme2}
Soient $f:Y→X$ un morphisme fini de schémas, $T_Y$ un fermé de $Y$
et $T_X=f(T_Y)$ son image.
\begin{enumerate}
\item On a l'inégalité : $\mathrm{codim}(T_X,X)≥\mathrm{codim}(T_Y,Y)$.
\item Soit $K∈\ob\,\detalepluscoefs{Y}{Λ}$ tel que
$K_{|Y-T_Y}∈\ob\,\detalepluscoefs{(Y-T_Y)}{Λ}$.
Alors, $(f_*K)_{|X-T_X}∈\ob\,\detalepluscoefs{(X-T_X)}{Λ}$.
\end{enumerate}
\end{lemme2}

\begin{démo}
Le premier énoncé est bien connu. 
Le second est un corollaire immédiat de la préservation de la constructibilité 
par le morphisme composé, fini, $Y-f^{-1}(T_X)↪Y-T_Y→X-T_X$.
\end{démo}

\subsubsection{}Soient $j:U→X$ et $k:W→U$ sont deux immersions
ouvertes, avec $U$ intègre. Nous souhaitons maintenant
déduire la constructibilité du complexe $\R j_* k_!Λ$
hors d'un fermé de codimension au moins $c$
de la propriété analogue pour les complexes $\R j_* Λ$.
Admettant ce résultat pour ces derniers, il résulte du
triangle distingué $k_!Λ→Λ→i_*Λ\dessusdessous{+1}{→}$, où $i$ est l'immersion 
fermée du complémentaire $F$ de $W$ dans $U$,
qu'il suffit de démontrer la constructibilité
de $R j_* i_*Λ$ hors d'un fermé de codimension au moins $c$.
Le schéma $U$ étant intègre, l'adhérence $\sur{F}$ de $F$ dans $X$ est de codimension
strictement positive $X$. Soit $m:F↪\sur{F}$ l'immersion ouverte
correspondante et $n:\sur{F}↪X$ l'immersion fermée. On a tautologiquement :
\[
R j_* i_*Λ=n_* \R m_*Λ,
\]
par commutativité du diagramme
\begin{center}
\begin{tikzpicture}[auto]
 \matrix(diag)[matrix of math nodes,column sep=4em,row sep=5ex]
{ |(U)| U & |(X)|  X \\
  |(F)| F & |(Fb)| \sur{F} \\};
 \draw[right hook->] (U)  -- node{$j$} (X);
 \draw[right hook->] (F) -- node[description]{$m$} (Fb);
 \draw[->] (F) -- node[swap]{$i$} (U);
 \draw[->] (Fb) -- node{$\sur{π}$} (X);
 \end{tikzpicture}
\end{center}
Par hypothèse de récurrence $(P_{c-1})$, il existe un fermé
$T_{\sur{F}}$ de $\sur{F}$, de codimension
au moins $c$ dans $\sur{F}$, tel que la restriction
de $\R m_*Λ$ à l'ouvert $\sur{F}-T_{\sur{F}}$ de $\sur{F}$ soit constructible.
L'image directe par l'immersion fermée $n$ du complexe $\R m_*Λ$
est donc constructible sur l'ouvert $X-T_{\sur{F}}$ de $X$.
La conclusion résulte maintenant
du fait que la codimension de $T_{\sur{F}}$ dans $X$ est
\emph{strictement supérieure} à $c$.

\begin{remarque2}
O. Gabber sait également démontrer un résultat
de \emph{constructibilité uniforme}, dans l'esprit
de ceux de \cite[§3]{Katz-Laumon} mais sans hypothèse
sur la caractéristique.
Cf. courriel à Luc Illusie, du 3 avril 2007 ;
voir aussi \cite{uniformite@Orgogozo}.
\end{remarque2}

\section{Coefficients $ℓ$-adiques}

\subsection{Définitions}
On rappelle ici la construction, due à Torsten Ekedahl
(\cite{adic@Ekedahl}), de la catégorie triangulée des complexes
bornés constructibles $ℓ$-adiques. Voir aussi
\cite{monodromie@Fargues}, §5 pour un résumé et quelques améliorations.
Pour un formalisme également valable pour les champs
et sous des hypothèses moins restrictives,
cf. « Enhanced six operations and base change theorem for sheaves on
Artin stacks » de Zheng Weizhe et Liu Y.

On fixe ici un schéma nœthérien $X$, sur lequel un nombre premier $ℓ$
est inversible.

\subsubsection{Systèmes projectifs}
Notons $X^𝐍$ le topos des systèmes projectifs
indicés par $𝐍$ de faisceaux étales sur $X$ ; on
en fait un topos annelé via $𝐙/ℓ^𝐍:=(𝐙/ℓ^n)_n$.
Un \emph{système projectif $ℓ$-adique} de faisceaux
sur $X$ est un $𝐙/ℓ^𝐍$-module sur $X^𝐍$,
c'est-à-dire un système projectif
de faisceaux abéliens $ℱ=(… → ℱ_{n+1} → ℱ_n → …)$ sur $X$, où $ℱ_n$
un $𝐙/ℓ^n 𝐙$-modules sur $X$. Ils constituent
une catégorie abélienne, dont on note
$\dcoefs{X^𝐍}{𝐙/ℓ^𝐍}$ la catégorie dérivée.
Un tel système projectif de faisceaux est dit \emph{essentiellement nul} si pour tout $n$
il existe un entier $m ≥ n$ tel que le morphisme de
transition $ℱ_m→ℱ_n$ correspondant soit nul.
De même, un complexe $𝒦 ∈ \Ob\,\dcoefs{X^𝐍}{𝐙/ℓ^𝐍}$
est essentiellement nul si chaque système projectif $H^i(𝒦)$ de faisceau l'est.

\renewcommand{\mod}{\,\mathrm{mod}\,}

\subsubsection{$𝐙_ℓ$-complexes}
Soit $\dbccoefs{X^𝐍}{𝐙/ℓ^𝐍}$
la sous-catégorie pleine de $\dcoefs{X^𝐍}{𝐙/ℓ^𝐍}$
des complexes
\[
𝒦=\big(𝒦_n ∈ \Ob\,\dbcoefs{X^𝐍}{𝐙/ℓ^n}\big)_n
\]
dont la « réduction modulo $ℓ$ »
\[
𝐅_ℓ \dessusdessous{𝕃}{⊗}_{𝐙/ℓ^𝐍} 𝒦 = (𝐅_ℓ \dessusdessous{𝕃}{⊗}_{𝐙/ℓ^n} 𝒦_n)_n
\]
soit un \emph{pro-objet essentiellement
constant} de $\dbccoefs{X}{𝐙/ℓ}$.
(« Essentiellement constant » : isomorphe, modulo les complexes essentiellement nuls, à un système
projectif provenant de $\dbccoefs{X}{𝐙/ℓ}$.)
Un tel objet est appelé un \emph{$𝐙_ℓ$-complexe borné
constructible} ; ils forment une catégorie triangulée.

\subsubsection{Catégorie triangulée des $𝐙_ℓ$-faisceaux}
On note $𝒟^{𝖻}_{𝖼}(X,𝐙_ℓ)$ la catégorie triangulée
obtenue à partir de la catégorie $\dbccoefs{X^𝐍}{𝐙/ℓ^𝐍}$
en inversant les \emph{isomorphismes essentiels modulo $ℓ$},
c'est-à-dire les flèches $u$ tel que
$𝐅_ℓ \dessusdessous{𝕃}{⊗}_{𝐙/ℓ^𝐍} u$ ait un cône essentiellement nul.
De même, on peut définir variantes non bornées et non
constructibles.

\subsubsection{}Comme expliqué en \cite[§5.9]{monodromie@Fargues},
lorsque $X$ est de type fini sur un corps séparablement clos
ou fini, la catégorie obtenue est équivalente
à la catégorie $2-\lim_n \dbctfcoefs{X}{𝐙/ℓ^n}$ considérée
par P. Deligne dans Weil II. Notons qu'en toute généralité,
les constituants d'un $𝐙_ℓ$-complexe borné constructible
sont de tor-dimensions finies.
On rappelle (\sga{4}{XVII}{4.1.9}) qu'un complexe
$𝒦∈\ob\,\db(T,A)$, où $T$ est un topos et $A$ un Anneau
commutatif, est dit de \emph{tor-dimension} inférieure ou égale
à $n$ si pour tout complexe $ℒ$ de $A$-modules concentré en degrés
positifs ou nuls, $\HH^i(𝒦⊗^L ℒ)=0$ pour tout $i<-n$.

\subsubsection{}Un des points clefs de la théorie
est le \emph{lemme de \japonais{中山}-Ekedahl} d'après
lequel le foncteur triangulé noté $𝐅_ℓ \dessusdessous{𝕃}{⊗}_{𝐙_ℓ} -$,
déduit du foncteur de réduction modulo $ℓ$ ci-dessus,
est \emph{conservatif} : $𝒦 ∈ \Ob\,𝒟^{𝖻}_{𝖼}(X,𝐙_ℓ)$
est nul si et seulement si le complexe $ 𝐅_ℓ
\dessusdessous{𝕃}{⊗}_{𝐙_ℓ} 𝒦 ∈ \Ob\,\dbccoefs{X}{𝐙/ℓ}$ l'est.
Nous renvoyons le lecteur à \cite[prop. 1.1]{multiplicativeII@Ekedahl}
ou \cite[2.3.7,2.4.5]{Finiteness@Illusie} pour
une première apparition de ce lemme, et \cite[3.6.(ii)]{adic@Ekedahl}
pour le résultat précédent.

\subsubsection{}D'après ce lemme, et ces corollaires
(\emph{op. cit.}, th. 5.1.(ii) et th. 6.3), le théorème de
finitude $ℓ$-adique ci-dessous résulte d'un théorème de finitude
pour les coefficients finis.

\subsection{Théorèmes : énoncés}

\begin{théorème2}
Soient $X$ un schéma nœthérien quasi-excellent
de dimension finie, $ℓ$ un nombre premier inversible
sur $X$ et $f:Y→X$ un morphisme de type fini.
Pour tout entier $n≥1$, le foncteur $\R
f_*:\detalepluscoefs{Y}{𝐙/ℓ^n}→\detalepluscoefs{X}{𝐙/ℓ^n}$
envoie $\detalebctfcoefs{Y}{𝐙/ℓ^n}$ dans $\detalebctfcoefs{X}{𝐙/ℓ^n}$.
\end{théorème2}

\begin{démo}
D'après le résultat du paragraphe §\ref{13-compléments annulation},
le foncteur $\R f_*$ est de dimension cohomologique finie.
La conclusion résulte de \sga{4}{XVII}{5.2.11} (tor-dimension finie)
et du théorème \ref{13-constructibilité} (constructibilité).
\end{démo}

\begin{remarque2}
On devrait pouvoir montrer, par réduction au cas des schémas de type
fini sur $𝐙$, que pour tout morphisme cohérent $f$, le
foncteur $\R f_*$ envoie un faisceau plat constructible sur un complexe de de tor-dimension 
au plus zéro.
\end{remarque2}

D'après les résultats esquissés dans
\cite[§5-6]{adic@Ekedahl}, on peut déduire du théorème
précédent le théorème $ℓ$-adique suivant.

\begin{théorème2}
Soient $X$ un schéma nœthérien quasi-excellent de dimension finie,
$ℓ$ un nombre premier inversible
sur $X$ et $f:Y→X$ un morphisme de type fini.
Le foncteur $\R f_*:\dpluscoefs{Y^𝐍}{𝐙/ℓ^𝐍}→
\dpluscoefs{X^𝐍}{𝐙/ℓ^𝐍}$, $𝒦=(𝒦_n)_n↦ \R f_* 𝒦 = (\R f_* 𝒦_n)_n$
  induit un foncteur
$\R f_*:\dbccoefs{Y}{𝐙_ℓ}→\dbccoefs{X}{𝐙_ℓ}$.
\end{théorème2}

\ifx\danslelivre\undefined
\bibliography{../fichiers-configuration-latex/bibliographie}
\bibliographystyle{../fichiers-configuration-latex/style-GTG}


\begin{center}
XIV. Fonctions de dimension\\
Vincent Pilloni et Benoît Stroh
\end{center}

\else
\setboolean{english}{false}
\setboolean{francais}{true}
\selectlanguage{francais}
\chapter{Fonctions de dimension}
\addtocontents{toc}{par Vincent Pilloni et Benoît Stroh}
\begin{center}
Vincent Pilloni et Benoît Stroh
\end{center}
\fi

Nous définissons la notion de fonction de dimension sur un schéma $X$ et nous montrons l'existence
de telles fonctions localement pour la topologie étale si $X$ est quasi-excellent.

\section{Universelle caténarité des schémas henséliens}

Dans cette partie, nous rappelons les notions de \emph{caténarité} et d'\emph{universelle
caténarité}.
Le lecteur pourra consulter l'exposé~I pour plus de détails.

\subsection{Schémas universellement caténaires}

Soient $S$ un espace topologique et $X\subset Y$ des fermés irréductibles de $S$. Notons
$\codim(X,Y)$ la borne supérieure de l'ensemble des longueurs des chaînes strictement croissantes
de fermés irréductibles $X\varsubsetneq Z \subset Y$. Si $S$ est un schéma,  $X$ et $Y$ des
sous-schémas fermés intègres et $x$ le point générique de $X$, on a
$$\codim(X,Y)\:=\:\dim(\OO_{Y,x})\: .$$

\begin{definition2}[\egazero{iv}{14.3.2}] \label{14-definition-catenaire} Un schéma $S$ est
\emph{caténaire} s'il est nœthérien et si pour toute chaîne $X\subset Y\subset Z$ de fermés
irréductibles de $S$, on a $$\codim(X,Z)=\codim(Y,Z)+\codim(X,Y).$$ Un schéma $S$ est
\emph{universellement caténaire} si tout schéma de type fini sur $S$ est caténaire.
\end{definition2}

La notion de caténarité est stable par localisation et par restriction à des sous-schémas fermés.
Ainsi, $S$ est universellement caténaire si et seulement si pour tout entier $n\geq 0$, le schéma
$\Aff^n_S$ est caténaire.

\begin{lemme2} \label{14-lemme-CohenMac-implique-univ-cat}  Un schéma de Cohen-Macaulay est
universellement caténaire.  \end{lemme2}

\begin{démo} Si $S$ est Cohen-Macaulay, il est caténaire d'après \cite{CA@Matsumura}~16.B. Comme
pour tout $n\geq 0$, le schéma $\Aff^n_S$ reste Cohen-Macaulay, le schéma $S$ est bien
universellement caténaire.  \end{démo}

\begin{exemple2}\label{14-exemple-regulier-implique-univ-cat}  Tout schéma régulier est
universellement caténaire car Cohen-Macaulay. En particulier, le spectre d'un corps, un trait et le
spectre d'une algèbre de séries formelles sur un corps ou sur un anneau de valuation discrète sont
universellement caténaires. Tout schéma de type fini sur un schéma régulier est universellement
caténaire.  \end{exemple2}

\begin{proposition2}[\cite{CA@Matsumura}~28.P] \label{14-prop-complet-univ-cat} Un schéma local
complet nœthérien est universellement caténaire.  \end{proposition2}

\begin{démo} Le théorème de structure de Cohen~\egazero{iv}{19.8.8} permet d'écrire tout schéma
local complet nœthérien comme fermé dans le spectre d'une algèbre de séries formelles sur un anneau
de Cohen. L'universelle caténarité résulte de l'exemple précédent et de la stabilité de cette notion
par passage à un fermé.  \end{démo}

\subsection{Un théorème de Ratliff}

On dit qu'un schéma nœthérien est \emph{équidimensionnel} si toutes ses composantes irréductibles
ont même dimension (finie). Soit $S$ un schéma local nœthérien. On note $\widehat{S}$ le spectre du
complété de l'anneau de $S$ en son idéal maximal.

\begin{definition2} \label{14-defi-formellement-catenaire} Le schéma local $S$ est
\emph{formellement équidimensionnel} si  $\widehat{S}$ est équidimensionnel. Il est
\emph{formellement caténaire} si pour tout $s\in S$, l'adhérence $\overline{\{s\}}$ est formellement
équidimensionnelle.  \end{definition2}

Soit $S$ un schéma local nœthérien. Ratliff a démontré le théorème fondamental suivant, qui a déjà
été mentionné dans la proposition~\refext{I}{1Ratliff}.

\begin{theoreme2}[\cite{CRT@Matsumura} 31.7] \label{14-theoreme-Ratliff} Pour un schéma local
nœthérien $S$, les conditions suivantes sont équivalentes~: \begin{itemize} \item $S$ est
formellement caténaire, \item $S$ est universellement caténaire, \item $\Aff^1_S$ est caténaire,
\item $S$ est caténaire et pour tout $s \in S$, tout schéma intègre $S'$ muni d'une flèche finie et
dominante $S'\rightarrow \overline{\{s\}}$ et tout point fermé $s'$ de $S'$, on a $\dim(\OO_{S',s'})
= \dim(\overline{\{s\}})$.  \end{itemize} \end{theoreme2}

On a ajouté une quatrième condition équivalente à l'énoncé \cite[Theorem~31.7]{CRT@Matsumura}. Il
résulte de~\egalong{iv}{2}{5.6.10} que les trois premières conditions équivalentes impliquent la
quatrième. La réciproque est démontrée au cours de la démonstration de
\cite[Theorem~31.7]{CRT@Matsumura} (au second paragraphe de la page~255).

\begin{corollaire2}[\cite{CRT@Matsumura} 31.2] \label{14-corollaire-dim2-cat} Tout schéma nœthérien
de dimension~$\leq 2$ est caténaire. Tout schéma nœthérien de dimension~$\leq 1$ est universellement
caténaire.  \end{corollaire2}

\subsection{Schémas henséliens et caténarité}

Nous avons vu  que tout schéma local complet nœthérien est universellement caténaire dans la
proposition~\ref{14-prop-complet-univ-cat}. Les schémas locaux henséliens  jouissent également de
bonnes propriétés de caténarité.

\begin{proposition2} \label{14-prop-henselien-cat-univ-cat} Tout schéma local hensélien caténaire
est universellement caténaire.  \end{proposition2}

\begin{démo} Soit $S=\SP(A)$ un schéma local hensélien caténaire, soit $P$ un idéal premier de~$A$,
soit $L$ une extension finie de $\mathrm{Frac}(A/P)$ et soit $B$ une extension finie de $A/P$
contenue dans~$L$. D'après le théorème \ref{14-theoreme-Ratliff}, il suffit de prouver que la
dimension du localisé de $B$ en chacun de ses idéaux maximaux est égale à la dimension de $A/P$.
Toute algèbre finie sur un anneau hensélien est semi-locale d'après \egalong{iv}{4}{18.5} et
\egalong{iv}{4}{18.6}. Comme le schéma $B$ est intègre, il est local. Le théorème du « going-up »
(\cite{CRT@Matsumura}~9.3 et~9.4) montre qu'on a bien $\dim(B)=\dim(A/P)$.  \end{démo}

\medskip Rappelons également le résultat suivant, conséquence du
corollaire~\refext{I}{1-cor-normal-versus-complet}~ii).

\begin{proposition2} \label{14-prop-henselien-qexc-univ-cat} Tout schéma local hensélien
quasi-excellent est universellement caténaire.  \end{proposition2} Ainsi, tout schéma local
hensélien quasi-excellent est excellent.

\section{Spécialisations immédiates et fonctions de dimension}

\subsection{Définitions}

Soit $X$ un schéma. Pour tout point $x$ de $X$ et tout point géométrique $\bar{x}$ au-dessus de $x$,
on note $X_{(x)}$, $X_{(x)}^\h$ et $\widehat{X}_{(x)}$ le localisé, l'hensélisé et le complété de
$X$ en $x$. De même, l'on note $X_{(\bar{x})}$ l'hensélisé strict de $X$ en $\bar{x}$.

Soient $x$ et $y$ deux points de $X$, et $\bar{x}$ et $\bar{y}$ deux points géométriques au dessus
de $x$ et $y$. 

\begin{definition2}[\sga{4}{vii}{7.2}]\label{14-definition-specialisation} Un morphisme de
spécialisation $\bar{x} \spim\bar{y}$ est la donnée d'un $X$-morphisme $X_{(\bar x )} \rightarrow
X_{(\bar{y})}$ entre hensélisés stricts.  \end{definition2} D'après \sga{4}{vii}{7.4}, la donnée
d'une spécialisation  $\bar{x } \spim \bar{y}$ est équivalente à la donnée d'un $X$-morphisme
$\bar{x} \rightarrow X_{(\bar{y})}$.

\begin{definition2} \label{14-definition-specialisation-immediate} Soit $r \in \NN$. On dit qu'une
spécialisation $\bar{x} \spim \bar{y}$ est une  \emph{spécialisation de codimension $r$} si
l'adhérence de l'image de $\bar{x}$ dans $X_{(\bar{y})}$ est un schéma de dimension $r$.  \\ On dit
que $y$ est une \emph{spécialisation étale immédiate} de $x$ s'il existe une spécialisation $\bar{x}
\spim \bar{y}$ qui soit de codimension~$1$.  \\ On dit que $y$ est une \emph{spécialisation Zariski
immédiate} de $x$  si $y \in \overline{\{ x \} }$ et si le localisé en $y$ de l'adhérence de $x$
est de dimension $1$.  \end{definition2}

Si $y$ est une {spécialisation étale immédiate} de $x$, on dit également que $x$ est une
\emph{générisation étale immédiate} de $y$. Désignons par $f:X_{(\bar y)}\ra X_{(y)}$ le morphisme
d'hensélisation stricte. Les générisations  étales immédiates de $y$ sont alors les images par $f$
des points $x'\in X_{(\bar y)}$ tels que $\dim(\overline{\{x'\}})=1$.

Avant d'examiner plus en détail ces notions, on rappelle le fait facile suivant que nous utiliserons
implicitement plus bas : si $f\colon X\to S$ est un morphisme plat, $f$ envoie les points maximaux
de $X$ sur des points maximaux de $S$, autrement dit toute composante irréductible de $X$ domine une
composante irréductible de $S$.

Si $x$ et $y$ sont deux points d'un schéma nœthérien $X$ tels que $y\in \overline{\{x\}}$ (au sens
habituel, c'est dire que $y$ est une spécialisation de $x$ ou encore que $x$ est une générisation de
$y$), alors $y$ est une spécialisation Zariski (resp. étale) immédiate de $x$ si et seulement si
c'est le cas dans $\overline{\{x\}}_{(y)}$. Pour certaines considérations, ceci permet de supposer
que $X$ est local intègre de point générique $x$ et de point fermé $y$. Dans ce cas, $y$ est
spécialisation Zariski immédiate de $x$ si et seulement si $\dim(X)=1$. Dans le cas étale, cela se
lit sur l'hensélisé strict :

\begin{proposition2} \label{14-proposition-Laszlo-1} Si $x$ et $y$ sont deux points d'un schéma
nœthérien $X$, le point $y$ est une {spécialisation étale immédiate} de $x$ si et seulement si $y
\in \overline{\{ x \} }$ et que le hensélisé strict en un point géométrique au-dessus de $y$ de
l'adhérence de $x$ possède une composante irréductible de dimension $1$.  \end{proposition2}

\begin{démo} On se ramène au cas particulier $X=\overline{\{x\}}_{(y)}$ envisagé plus haut. Le point
$y$ est une spécialisation étale immédiate si et seulement s'il existe un point $\tilde{x}$ de
$X_{(\bar{y})}$ au-dessus de $x$ tel que l'adhérence de $\tilde{x}$ dans $X_{(\bar{y})}$ soit de
dimension $1$. Cela équivaut comme il est énoncé ici que $X_{(\bar{y})}$ possède une composante
irréductible de dimension $1$. En effet, si on note $\tilde{x}$ le point générique de $C$, par
l'argument de platitude énoncé plus haut, $C$ domine $X$, c'est-à-dire que $\tilde{x}$ est au-dessus
de $x$. Inversement, si $\tilde{x}$ est un point au-dessus de $x$ dont l'adhérence dans $X_{(\bar
y)}$ soit de dimension $1$, on peut noter $C$ une composante irréductible de $X_{(\bar y)}$
contenant $\tilde{x}$. Le point générique de $C$ et $\tilde{x}$ étant tous les deux au-dessus de
$x$, ils sont égaux puisque l'un est une générisation de l'autre et que les fibres de $X_{(\bar
y)}\to X$ sont discrètes.  \end{démo}

\medskip

\begin{proposition2} \label{14-prop-specialisation-zar-et} Soit $X$ un schéma nœthérien. Une
spécialisation Zariski immédiate entre points de $X$ est une spécialisation étale immédiate, et la
réciproque est vraie si $X$ est universellement caténaire.  \end{proposition2}

On peut supposer que $X=\overline{\{x\}}_{(y)}$ comme précédemment. Pour l'implication, on suppose
que $\dim(X)=1$ et on veut montrer que $X_{(\bar y)}$ possède une composante irréductible de
dimension $1$. En fait, toutes les composantes irréductibles de $X_{(\bar y)}$ sont de dimension
$1$. En effet, le schéma $X_{(\bar y)}$ est local et tout point de $X_{(\bar y)}$ autre que $\bar y$
appartient à la fibre de $x$ pour le morphisme \guil{pro-étale} $X_{(\bar y)}\to X$. L'espace
topologique sous-jacent à $X_{(\bar y)}-\{\bar y\}$ est donc fini, discret (et non vide). Ceci
implique que toute composante irréductible de $X_{(\bar y)}$ est de dimension $1$. Pour la
réciproque, nous utiliserons deux lemmes :

\begin{lemme2}\label{14-lemme-riou-1} Soit $X$ un schéma local nœthérien hensélien de point fermé
$y$. Soit $\bar y$ un point géométrique au-dessus de $y$. Alors, $X$ possède une composante
irréductible de dimension $1$ si et seulement si le hensélisé strict $X_{(\bar y)}$ en possède une.
\end{lemme2}

\begin{démo} Si $C$ est une composante irréductible de dimension $1$ de $X_{(\bar y)}$, son image
ensembliste dans $X$ est fermée car $p\colon X_{(\bar y)}\to X$ est entier. Comme $p$ est plat,
$p(C)$ est une composante irréductible de $X$ contenant exactement deux points donc $\dim(p(C))=1$.
Inversement, la surjectivité et la platitude de $p$ impliquent que si $D\subset X$ est une
composante irréductible de dimension $1$, il existe une composante irréductible $C$ de $X_{(\bar
y)}$ telle que $p(C)=D$. On a bien sûr $\dim(C)\geq 1$. Soit $z\in C$ un point qui ne soit pas le
point générique de $C$. Le point $p(z)$ ne peut pas être le point générique de $D$ car sinon la
fibre générique de $p$ ne serait pas discrète. C'est donc que $p(z)$ est le point fermé de $C$. Le
fait que $p^{-1}(y)$ soit discret implique alors que $z$ ne peut être que le point fermé de $C$.  Le
schéma local intègre $C$ possède donc exactement deux points : $\dim(C)=1$.  \end{démo}

\begin{lemme2}\label{14-lemme-riou-2} Soit $X$ un schéma local nœthérien. Si $X$ possède une
composante irréductible de dimension $1$, alors son complété $\widehat{X}$ aussi et la réciproque
est vraie si $X$ est universellement caténaire.  \end{lemme2}

\begin{démo} Commençons par le cas où $X$ est intègre. Comme $X$ et $\widehat{X}$ ont la même
dimension (cf. \cite[chap.~VIII, \S{}9]{CA@Zariski-Samuel}), si $\dim(X)=1$, alors
$\dim(\widehat{X})=1$. Un schéma local de dimension $1$ n'a bien sûr que des composantes
irréductibles de dimension $1$. Inversement, si $X$ est universellement caténaire, d'après le
théorème~\ref{14-theoreme-Ratliff}, les composantes irréductibles de $\widehat{X}$ ont toutes la
même dimension. Si l'une d'entre elles est de dimension $1$, le schéma $\widehat{X}$ est lui aussi
de dimension $1$, et alors $\dim(X)=\dim(\widehat{X})=1$.

Dans le cas général, notons $X_i$ les composantes irréductibles de $X$.  Pour tout $i$, le produit
fibré $X_i\times_X \widehat{X}$ s'identifie à $\widehat{X_i}$ (voir \sga{1}{iv}{3}). C'est un fait
que les composantes irréductibles des différents $\widehat{X_i}$ sont exactement les composantes
irréductibles de $\widehat{X}$ : ce sont des parties fermées irréductibles recouvrant $\widehat{X}$
et aucune d'entre elles n'est contenue dans une autre (ceci se déduit du fait que chaque composante
irréductible de $\widehat{X_i}$ domine $X_i$). Il est dès lors évident que l'énoncé pour $X$ résulte
de l'énoncé pour les schémas locaux intègres $X_i$.  \end{démo}

Montrons la réciproque énoncée dans la proposition~\ref{14-prop-specialisation-zar-et}. Soit $X$ un
schéma local nœthérien intègre universellement caténaire de point générique $x$ et de point fermé
$y$. On suppose que $y$ est une spécialisation étale immédiate de $x$. D'après la
proposition~\ref{14-proposition-Laszlo-1}, cela signifie que $X_{(\bar y)}$ possède une composante
irréductible de dimension $1$, ce qui équivaut d'après le lemme~\ref{14-lemme-riou-1} à dire que le
hensélisé $X_{(y)}^\h$ possède une composante irréductible de dimension $1$. Le complété
$\widehat{X}$ de $X$ est aussi celui de $X_{(y)}^\h$, le sens facile du lemme~\ref{14-lemme-riou-2}
appliqué à $X_{(y)}^\h$ montre que $\widehat{X}$ possède une composante irréductible de dimension
$1$. La réciproque de ce lemme appliquée au schéma universellement caténaire $X$ montre que $X$
possède une composante irréductible de dimension $1$. Ceci montre que $y$ est une spécialisation
Zariski immédiate de $x$.

\bigskip

Si le schéma $X$ est quasi-excellent, on peut alors lire les spécialisations étales d'un point $x$
de $X$ dans le complété de $X$ en $x$ :

\begin{proposition2}\label{14-prop-Laszlo-2} Soit $X$ un schéma quasi-excellent. Soient $x$ et $y$
deux points de $X$. On suppose que $y\in\overline{\{x\}}$. Notons $c:\widehat X_{(y)}\ra X_{(y)}$ le
morphisme de complétion. Le point $y$ est une spécialisation étale immédiate de $x$ si et seulement
si $c^{-1}(\overline{\{x\}})$ possède une composante irréductible de dimension $1$.
\end{proposition2}

\begin{démo} On peut supposer que $X=\overline{\{x\}}_{(y)}$. Le point $y$ est une spécialisation
étale immédiate de $X$ si et seulement si $X_{(\bar y)}$ possède une composante irréductible de
dimension $1$, c'est-à-dire, d'après le lemme~\ref{14-lemme-riou-1} que $X_{(y)}^\h$ en possède une.
Comme $X_{(y)}^\h$ est universellement caténaire (cf.
proposition~\ref{14-prop-henselien-qexc-univ-cat}), cela équivaut encore d'après le
lemme~\ref{14-lemme-riou-2} à dire que $\widehat{X}$ possède une composante irréductible de
dimension $1$.  \end{démo}

\medskip

\begin{definition2} \label{ddefi4} On appelle \emph{fonction de dimension} sur $X$ toute fonction
$\delta  : X \ra \ZZ $ telle que pour toute spécialisation étale immédiate $x \spim y$ entre points
de $X$, on ait $$\delta(y) = \delta(x) - 1.$$ \end{definition2}

La notion de fonction de dimension ne voit pas les éléments nilpotents~: $\delta$ est une fonction
de dimension sur~$X$ si et seulement si elle induit une fonction de dimension sur le sous-schéma
réduit~$X^{red}$. De plus si $U \stackrel{i}{\hookrightarrow } X$ est un morphisme étale et $\delta$
est une fonction de dimension sur $X$, $\delta\circ i$ définit une fonction de dimension sur $U$.
Plus précisément, l'ensemble des fonctions de dimension sur les $X$-schémas étales définit un
faisceau étale sur $X$. La différence entre deux fonctions de dimension sur $X$ est une fonction
invariante par spécialisations quelconques, donc une fonction localement constante. Nous montrerons
plus loin que si $X$ est quasi-excellent, des fonctions de dimension existent localement pour la
topologie étale sur $X$ et que les fonctions de dimension forment un $\ZZ$-torseur étale.

\subsection{Fonctions de dimension et universelle caténarité}

Le but de ce paragraphe est de démontrer le résultat suivant.

\begin{theoreme2}[Gabber]\label{14-theoreme-univ-cat-fonction-dim} Un schéma  nœthérien est
universellement caténaire si et seulement s'il possède une fonction de dimension localement pour la
topologie de Zariski.  \end{theoreme2}

Le théorème résulte de la conjonction du
corollaire~\ref{14-corollaire-univ-cat-fonction-dim-loc-zar} et de la
proposition~\ref{14-proposition-fonc-dim-univ-cat} ci-dessous. 

\begin{proposition2} \label{14-proposition--codim-fonc-dim} Soit $X$ un schéma intègre
universellement caténaire. La fonction $\delta:X\ra\ZZ$ définie  par $\delta(x)=-\dim(\OO_{X,x})$
est une fonction de dimension sur $X$.  \end{proposition2}

\begin{démo} En vertu de la proposition~\ref{14-prop-specialisation-zar-et}, comme $X$ est
universellement caténaire, il suffit de montrer que $\delta(y) = \delta(x) - 1$ pour toute
spécialisation Zariski immédiate~$x \spim y$. Comme $X$ est caténaire intègre, on a
$$\dim(\OO_{X,y}) = \dim(\OO_{X,x}) + \dim(\OO_{\overline{\{y\} },x}) = \dim(\OO_{X,x}) +1\: .$$
\end{démo}

\begin{remarque2} Si $X$ n'est pas supposé intègre, la fonction $\delta(x)=-\dim(\OO_{X,x})$ n'est
pas forcément une fonction de dimension comme le montre l'exemple où $X$ est obtenu par recollement
en un point d'une droite et d'un plan.  \end{remarque2}

\begin{corollaire2} \label{14-corollaire-univ-cat-fonction-dim-loc-zar} Tout schéma universellement
caténaire admet des fonctions de dimension localement pour la topologie de Zariski.
\end{corollaire2}

\begin{démo} Soit $X$ un schéma universellement caténaire. Soit $x\in X$. Il s'agit de montrer qu'il
existe un voisinage ouvert de $x$ pouvant être muni d'une fonction de dimension. L'espace
topologique $X$ est réunion de ses composantes irréductibles  $X_1,...,X_n$. Quitte à remplacer $X$
par l'ouvert complémentaire des  composantes $X_i$ ne contenant pas $x$, on peut supposer que $x$
appartient à toutes les composantes $X_i$.   Pour tout $0\leq i\leq n$,  notons $\mathscr{F}_i$
l'ensemble des fonctions de dimension sur $X_i$. D'après la
proposition~\ref{14-proposition--codim-fonc-dim}, cet ensemble est non vide et est un torseur sous
$\ZZ$. On choisit  un élément $\delta_i\in\mathscr{F}_i$ qui vaut $0$ sur le point $x$. Pour tous
$1\leq i , \; j \leq n$, la fonction $\delta_i-\delta_j$ est localement constante sur $X_i\cap X_j$
et vaut $0$ au point $x$. Soit $F_{i,j}$ le fermé de $X$, réunion des composantes connexes de
$X_i\cap X_j$ ne contenant pas $x$.  Soit $U$ le complémentaire dans $X$ de la réunion des
$F_{i,j}$.  Les fonctions $\delta_i$ se recollent en une fonction de dimension sur $U$.  \end{démo}

\medskip

Démontrons une réciproque  partielle du
corollaire~\ref{14-corollaire-univ-cat-fonction-dim-loc-zar}.

\begin{lemme2}\label{14-lemme-fonc-dim-cat} Un schéma nœthérien qui possède  une fonction de
dimension localement pour la topologie de Zariski est caténaire.  \end{lemme2}

\begin{démo} Pour montrer la caténarité, on peut supposer que le schéma $S$ possède une fonction de
dimension $\delta$. Supposons que $X\subset Y$ sont des fermés irréductibles de points génériques
respectifs $x$ et $y$. Choisissons une chaîne de spécialisations Zariski immédiates $y=x_0\spim
x_1\spim \dots x_d=x$ de longueur maximale. Par définition de la codimension, on a $\codim(X,Y)=d$
et par définition des fonctions de dimension, compte tenu du sens facile de la
proposition~\ref{14-prop-specialisation-zar-et}, on obtient $\delta(x)=\delta(y)-d$, d'où
$\codim(X,Y)=\delta(y)-\delta(x)$.

Maintenant, si $X\subset Y\subset Z$ sont des fermés irréductibles, on a :

\begin{eqnarray*} \delta(y)-\delta(x) &=& \codim(X,Y) , \\ \delta(z)-\delta(y) &=& \codim(Y,Z) , \\
\delta(z)-\delta(x) &=& \codim(X,Z) .  \end{eqnarray*}On en déduit
$\codim(X,Z)=\codim(X,Y)+\codim(Y,Z)$ d'où la caténarité.  \end{démo}

\medskip

Grâce au théorème~\ref{14-theoreme-Ratliff}, on peut remplacer \guil{caténaire} par
\guil{universellement caténaire} dans  le lemme~\ref{14-lemme-fonc-dim-cat} :

\begin{proposition2} \label{14-proposition-fonc-dim-univ-cat} Un schéma nœthérien qui possède une
fonction de dimension localement pour la topologie de Zariski est universellement caténaire.
\end{proposition2}

\begin{démo} On peut supposer qu'un tel schéma $S$ est local et possède une fonction de dimension
$\delta$. Notons $S^h$ son hensélisé. La fonction $\delta$ induit une fonction de
dimension~$\delta^h$ sur $S^h$. D'après le lemme \ref{14-lemme-fonc-dim-cat} appliqué à $S^h$ et à
$\delta^h$, le schéma $S^h$ est caténaire. D'après la proposition
\ref{14-prop-henselien-cat-univ-cat}, il est universellement caténaire. Le théorème
\ref{14-theoreme-Ratliff} montre alors que $S^h$ est formellement caténaire. Par conséquent, $S$ est
formellement caténaire donc également universellement caténaire.  \end{démo}

\subsection{Existence locale pour la topologie étale}

Dans ce paragraphe nous allons démontrer le théorème suivant.

\begin{theoreme2}[Gabber]\label{14-theoreme-exist-local-et-fn-dim} Tout schéma quasi-excellent
possède des fonctions de dimension localement pour la topologie étale.  \end{theoreme2}

Une application répétée du lemme suivant (variante de l'argument du
corollaire~\ref{14-corollaire-univ-cat-fonction-dim-loc-zar}) permet de montrer que si l'énoncé du
théorème est vrai pour les composantes irréductibles d'un schéma nœthérien $X$, alors le théorème
vaut aussi pour $X$. Plus loin, on pourra ainsi supposer que $X$ est intègre.

\begin{lemme2} Soit $X$ un schéma nœthérien dont l'espace topologique sous-jacent soit réunion de
deux sous-schémas fermés $X_1$ et $X_2$. Soit $\overline x$ un point géométrique de $X_1\cap X_2$.
On suppose que pour tout $i\in\{1,2\}$, il existe un voisinage étale $U_i$ de $\overline x$ dans
$X_i$ tel que $U_i$ admette une fonction de dimension.  Alors, il existe un voisinage étale $U$ de
$\overline x$ dans $X$ tel que $U$ admette une fonction de dimension.  \end{lemme2}

\begin{démo} Pour tout $i\in\{1,2\}$, on choisit un voisinage étale $U_i$ de $\overline x$ dans
$X_i$ tel que $U_i$ admette une fonction de dimension $\delta_i$. On se donne un point géométrique
distingué $\overline u_i$ au-dessus de $\overline x$ et on peut supposer que $\delta_i(u_i)=0$.
D'après \egalong{iv}{4}{18.1.1}, quitte à remplacer $U_i$ par un voisinage ouvert de $u_i$, on peut
supposer qu'il existe un morphisme étale $\widetilde{U_i}\to X$ et un isomorphisme
$\widetilde{U_i}\times_X{X_i}\simeq U_i$. On peut former le produit fibré $V=\widetilde{U_1}\times_X
\widetilde{U_2}$. Notons $\pi\colon V\to X$ la projection et $\overline{v}$ un point géométrique de
$V$ au-dessus de $\overline{u}_1$ et $\overline{u}_2$. Pour tout $i\in\{1,2\}$, la projection de $V$
sur le facteur $\widetilde{U_i}$ induit un morphisme étale $\pi^{-1}(X_i)\to U_i$. Par composition
avec ce morphisme étale, la fonction de dimension $\delta_i$ sur $U_i$ induit une fonction de
dimension $\tilde{\delta}_i$ sur le sous-schéma fermé $\pi^{-1}(X_i)$ de $V$ et elle vérifie
$\tilde{\delta}_i(v)=0$. Ces fonctions de dimensions $\tilde{\delta}_i$ pour $i\in\{1,2\}$ se
recollent sur l'ouvert $U$ complémentaire dans $V$ de la réunion des composantes connexes de
$\pi^{-1}(X_1\cap X_2)$ ne contenant pas $v$.  \end{démo}

Avant de traiter le cas des schémas intègres, commençons par celui des schémas normaux :

\begin{proposition2}\label{14-proposition-codim-fonction-de-dimension-sur-normal-quasi-excellent}
Soit $X$ un schéma normal quasi-excellent. La fonction $\delta : X \rightarrow \ZZ$ définie par
$\delta(x) = -\dim(\OO_{X,x})$ est une fonction de dimension.  \end{proposition2}

\begin{démo} On peut supposer de plus que $X$ est local. Notons $Y$ son hensélisé et $h\colon Y\to
X$ le morphisme de hensélisation.  D'après le théorème~\refext{I}{1th-greco} et les commentaires
subséquents, $Y$ est lui aussi quasi-excellent. D'après la
proposition~\ref{14-prop-henselien-qexc-univ-cat}, $Y$ est universellement caténaire. Par ailleurs,
comme le morphisme $Y\to X$ est régulier, la normalité de $X$ implique celle de $Y$ (cf.
\egalong{iv}{2}{6.5.4}). Le schéma local $Y$ est donc intègre et universellement caténaire, l'opposé
de la codimension définit une fonction de dimension $\delta'\colon Y\rightarrow \ZZ$. Comme une
spécialisation étale immédiate entre points de $X$ se relève pour ainsi dire par définition en une
spécialisation étale immédiate de points de $Y$, pour montrer que $\delta$ est une fonction de
dimension sur $X$, il suffit de montrer que pour tout $y\in Y$, si on note $x=h(y)$, on a
$\delta(x)=\delta'(y)$. Il s'agit donc de montrer que $\dim(Y_{(y)})=\dim(X_{(x)})$. Comme la
dimension d'un anneau local nœthérien ne change pas par complétion (et donc par hensélisation), il
suffit de montrer que $\dim(Y^\h_{(y)})=\dim(X^\h_{(x)})$. Les deux schémas locaux excellents
$Y^\h_{(y)}$ et $X^\h_{(x)}$ ayant le même hensélisé strict, ils ont bien la même dimension d'après
le lemme suivant : \end{démo}

\begin{lemme2} Soit $X$ un schéma local excellent. Notons $\overline x$ un point géométrique
au-dessus du point fermé $x$ de $X$. Alors, $\dim(X)=\dim X_{(\bar x)}$.  \end{lemme2}

\begin{démo} D'après le théorème~\ref{14-theoreme-univ-cat-fonction-dim}, il existe une fonction de
dimension $\delta$ sur $X$. On peut supposer qu'elle vérifie $\delta(x)=0$. Par composition avec le
morphisme $X_{(\bar x)}\to X$, on obtient une fonction de dimension $\delta'$ sur $X_{(\bar x)}$. En
utilisant le lien entre fonction de dimension et codimension énoncé dans la démonstration du
lemme~\ref{14-lemme-fonc-dim-cat}, il est évident que $\dim(X)$ (resp. $\dim(X_{(\bar x)})$) est le
maximum atteint par la fonction $\delta$ sur $X$ (resp. par $\delta'$ sur $X_{(\bar x)}$). Comme
$X_{(\bar x)}\to X$ est surjectif, le maximum de $\delta$ est le même que celui de $\delta'$, donc
$\dim(X)=\dim(X_{(\bar x)})$.  \end{démo}

\medskip

Revenons au cas du théorème~\ref{14-theoreme-exist-local-et-fn-dim} où $X$ est supposé intègre et
quasi-excellent, et notons $Y$ son normalisé. Le morphisme $p : Y \rightarrow  X$ est fini et
surjectif, donc de descente cohomologique universelle. Notons~$\delta$ une fonction de dimension sur
$Y$ ; son existence est assurée par la
proposition~\ref{14-proposition-codim-fonction-de-dimension-sur-normal-quasi-excellent}. Soit~$p_1$
et~$p_2$ les deux projections $Y\times_X Y \ra Y$. La fonction $$p_1^\star \delta -p_2^\star\delta :
Y\times_X Y \longrightarrow \ZZ$$ définit un $1$-cocycle de Čech, donc une classe~$[p_1^\star
\delta -p_2^\star\delta]$ dans $\HH^1_{\v{C}ech}(Y \ra X, \ZZ)$.  D'après la théorie de la descente
cohomologique, il existe une injection naturelle $$\HH^1_{\v{C}ech}(Y \ra X, \ZZ) \hookrightarrow
\HH^1(X,\ZZ)\:.$$ La classe~$[p_1^\star \delta -p_2^\star\delta]$ définit donc une classe
d'isomorphisme de $\ZZ$-torseurs étales sur~$X$. Il résulte alors immédiatement de la proposition
suivante que $X$ admet une fonction de dimension localement pour la topologie étale :

\begin{proposition2} Soit $U$ un schéma étale sur $X$. L'annulation de la classe $[p_1^\star \delta
-p_2^\star\delta]\vert_U$ dans $\HH^1(U,\ZZ)$ entraîne l'existence d'une fonction de dimension sur
$U$.  \end{proposition2}

\begin{démo} En utilisant la compatibilité des constructions au changement de base étale $U\to X$,
on peut supposer que $U=X$.  L'annulation de $[p_1^\star \delta -p_2^\star\delta]$ dans
$\HH^1_{\v{C}ech}(Y \ra X,\ZZ) \hookrightarrow \HH^1(X,\ZZ)$ signifie qu'il existe une fonction
\emph{localement constante} $\gamma: Y \ra \ZZ$ telle que $p_1^\star \delta-p_2^\star
\delta=p_1^\star \gamma-p_2^\star \gamma$. Autrement dit, quitte à remplacer $\delta$ par
$\delta-\gamma$, on peut supposer que $p_1^\star \delta=p_2^\star \delta$. Ainsi, $\delta:Y\ra \ZZ$
se descend en une fonction $\delta':X\ra \ZZ $.

Pour conclure, il s'agit de montrer que si $p: Y\ra X$ est un morphisme fini surjectif entre schémas
quasi-excellents, que $\delta'\colon X\to \ZZ$ est une fonction et $\delta=\delta'\circ p$, alors
$\delta'$ est une fonction de dimension sur $X$ si $\delta$ est une fonction de dimension sur $Y$.
Pour cet énoncé, on peut supposer que $X$ est local hensélien quasi-excellent, donc universellement
caténaire (cf.  proposition~\ref{14-prop-henselien-qexc-univ-cat}), ce qui permet d'appliquer la
proposition~\ref{14-prop-specialisation-zar-et}. Pour montrer que $\delta'$ est une fonction de
dimension sur $X$, il suffit alors de savoir qu'au-dessus d'une spécialisation Zariski immédiate
$s\spim t$ de points de $X$, il existe une spécialisation Zariski immédiate $s'\spim t'$ de points
de $Y$.  \end{démo}

\subsection{Existence globale de fonctions de dimension}

Suivant \egazero{iv}{14.2.1}, on dit qu'un schéma nœthérien $X$ est équicodimensionnel si ses points
fermés ont tous la même codimension (qui est alors égale à $\dim(X)$).

\begin{exemple2} Les schémas de type fini équidimensionnels sur un corps $k$ ou sur $\ZZ$ sont
équicodimensionnels~: il est classique que dans cette situation, on a $\dim(X) = \dim(\OO_{X,x})$
pour tout point fermé $x$.  Les schémas locaux sont équicodimensionnels car ils possèdent un unique
point fermé.  Si $S=\mathrm{\SP(R)}$ est un trait d'uniformisante $\pi$, le schéma $\Aff^1_S$ n'est
pas équicodimensionnel. En effet, il existe un point fermé de $\Aff^1_S$ au-dessus du point
générique de $S$~:~il suffit d'écrire $\Aff^1_S=\SP(R[t])$ et de considérer $\MM=(\pi t -1 )$, qui
est un idéal maximal de corps résiduel $\mathrm{Frac}(R)$.  \end{exemple2}

Le lemme suivant est inspiré de~\egazero{iv}{14.3.3}\;\footnote{Gabber remarque que la proposition
\egazero{iv}{14.3.3} est fausse. Les assertions $a$, $c$ et $d$ de \emph{loc. cit.} sont
équivalentes entre elles et impliquent $b$ mais ne lui sont pas équivalentes. Il faut remplacer $b$
par la condition « $X$ est caténaire équidimensionnel et ses composantes irréductibles sont
équicodimensionnelles ». Gabber donne comme contre-exemple le spectre du localisé de
$k[x,y,z,w]/(xz,xw)$ en le complémentaire de l'union des idéaux premiers $(x-1,y)$ et $(x,z,w)$ avec
$k$ un corps.}.

\begin{lemme2} \label{14-lemme-EGA-faux} Soit $X$ un schéma équidimensionnel caténaire dont les
composantes irréductibles sont équicodimensionnelles. Pour tout $x\in X$, on a $$\dim(X) =
\dim\left(\overline{\{x\}}\right)+\dim(\OO_{X,x})\: .$$ \end{lemme2}

\begin{remarque2} En particulier, cette égalité est vérifiée pour tout schéma intègre local
caténaire. D'après~\cite[th.~31.4]{CRT@Matsumura}, si $X$ est intègre local nœthérien et si pour
tout $x\in X$, on a  $\dim(X)=\dim(\overline{\{x\}} )+\dim(\OO_{X,x})$, alors $X$ est caténaire.
\end{remarque2}

Le lemme~\ref{14-lemme-EGA-faux} et la proposition~\ref{14-proposition--codim-fonc-dim} impliquent
le résultat suivant.

\begin{corollaire2} \label{14-coro-equicodim-fonction-dim} Soit $X$ un schéma intègre,
équicodimensionnel et universellement caténaire. La fonction $\delta:X\ra\ZZ$ définie par $\delta(x)
= \dim(\overline{\{x\}})$ est une fonction de dimension sur $X$.  \end{corollaire2}

Les conclusions du corollaire sont prises en défaut si $X$ n'est pas équicodimensionnel. Soient par
exemple $S=\SP(R)$ un trait d'uniformisante $\pi$ et $X=\Aff^1_S=\SP(R[t])$. Si l'on note $x$ le
point fermé de $X$ correspondant à l'idéal maximal $(\pi t -1)$ et  $\eta$  le point générique de
$\Aff^1_S$, alors la spécialisation $\eta\spim x$ est immédiate et pourtant
$\dim(\overline{\{x\}})=0$ et $\dim(\overline{\{\eta\}})=2$.

\begin{corollaire2}\label{14-coro-schema-tf-sur-corps-fonc-dim} Soit $X$ un schéma qui est soit de
type fini sur un corps, soit de type fini sur $\ZZ$, ou soit local universellement caténaire. La
fonction définie par $\delta(x) = \dim(\overline{\{x\}})$ est une fonction de dimension sur $X$.
\end{corollaire2}

\begin{démo} Le schéma $X$ est universellement caténaire. D'après le
corollaire~\ref{14-coro-equicodim-fonction-dim}, la fonction $\delta$ est une fonction de dimension
sur chaque composante irréductible de $X$. Cette fonction est définie globalement donc est une
fonction de dimension sur $X$.  \end{démo}

\subsection{Fonction de dimension induite}

Soient $Y\ra X$ un morphisme de schémas et $\delta_X$ une fonction de dimension sur $X$. Dans
certains cas nous  pouvons construire une fonction de dimension~$\delta_Y$ induite sur $Y$. On admet
la proposition suivante.

\begin{proposition2}[\cite{CA@Matsumura},~14.C]\label{14-Matsumura-14-C} Soient $X$ un schéma
nœthérien intègre universellement caténaire, $Y$ un schéma intègre et $Y\rightarrow X$ un morphisme
de type fini dominant. Soient $k(X)$ et $k(Y)$ les corps de fractions respectifs de $X$ et $Y$,
soient $y$ un point de $Y$ et $x$ son image dans $X$, et soient $k(y)$ et $k(x)$ leurs corps
résiduels. On a $$\dim(\OO_{Y, y}) - \mathrm{degtr}\left(k(Y)/k(X)\right)  =   \dim(\OO_{X, x})-
\mathrm{degtr}\left(k(y)/k(x)\right) \: .$$ \end{proposition2}

\begin{corollaire2} \label{14-coro-fonction-dim-induite-tf} Soient $X$ un schéma nœthérien qui
possède une fonction de dimension $\delta_X$ et $f  :  Y \ra X$ un morphisme de type fini. La
fonction $\delta_Y  :  Y \ra \ZZ$ définie par $$\delta_Y(y) =  \delta_X(f(y))  +
\mathrm{degtr}\left(k(f(y))/k(y)\right)$$ est une fonction de dimension sur $Y$.  \end{corollaire2}

\begin{démo} On peut supposer que $X$ et $Y$ sont intègres et que $f$ est dominant. D'après la
proposition \ref{14-proposition-fonc-dim-univ-cat}, $X$ est universellement caténaire et d'après la
proposition  \ref{14-proposition--codim-fonc-dim}, $x\mapsto -\dim(\OO_{X,x})$ est une fonction de
dimension sur $X$. Comme les fonctions de dimension forment un $\ZZ$-torseur, on peut supposer que
$\delta_X(x)  =  -\dim(\OO_{X,x})$ pour tout $x\in X$.

Le corollaire~\ref{14-Matsumura-14-C} montre que $\delta_Y(y) =-\dim(\OO_{Y,y}) +
\mathrm{degtr}(k(Y)/k(X))$ et la proposition~\ref{14-proposition--codim-fonc-dim} montre que $y
\mapsto -\dim(\OO_{Y,y})$ est une fonction de dimension sur $Y$. Ainsi, $\delta_Y$ est une fonction
de dimension sur $Y$.  \end{démo}

\medskip Avant d'établir la fonctorialité des fonctions de dimension vis-à-vis des morphismes
réguliers entre schémas excellents, démontrons un énoncé de changement de base par un morphisme
régulier en cohomologie étale. Ce lemme est une simple conséquence du théorème de Popescu
\refext{I}{1-Artin-Popescu} et du théorème de changement de base par un morphisme
lisse~\sga{4}{xvi}{1.2}.

\begin{lemme2} \label{14-lemme-cb-regulier} Soient $$\xymatrix{T' \ar[r]^{g'} \ar[d]_{f'} & T
\ar[d]^f \\ S' \ar[r]^g & S}$$ un diagramme cartésien de schémas, $n$ un entier inversible sur~$S$
et~$\mc{F}$ un faisceau étale en~$\ZZ/n\ZZ$-modules sur~$T$. Supposons que $f$ est quasi-compact et
quasi-séparé et que $g$ est un morphisme régulier entre schémas nœthériens. La flèche naturelle de
changement de base $$g^* \R f_*  (\mc{F})  \isolong   \R f'_* g'^*  (\mc{F})$$ est un isomorphisme.
\end{lemme2}

\begin{démo} D'après le théorème de Popescu, il existe un ensemble ordonné filtrant~$I$ et une
famille de schémas~$S_i$ indexée par~$I$, tels que~$S_i$ soit lisse sur~$S$ pour tout~$i\in I$ et
que $S' = \lim_i \: S_i$. Il existe donc pour tout~$i\in I$ un diagramme commutatif à carrés
cartésiens $$\xymatrix{T' \ar[r]^{h'_i} \ar[d]^{f'} & T_i \ar[r]^{g'_i} \ar[d]^{f_i} & T \ar[d]^f \\
S' \ar[r]^{h_i} & S_i \ar[r]^{g_i} & S}$$ On conclut grâce à la suite d'isomorphismes suivante pour
tout $q\geq 0$ : \begin{eqnarray*} \R^q f'_*  g'^*  (\mc{F})  & \isolongg & {\colim}_{i} \: h_i^*
\R^q {f_i}_*  g_i'^* (\mc{F}) \\ & \isolongg &  {\colim}_{i} \: g^*  \R^q f_* (\mc{F}) \\ &
\isolongg & g^*  \R^q f_* (\mc{F}) \end{eqnarray*} Le premier de ces isomorphismes résulte du
théorème de passage à la limite~\sga{4}{vii}{5.11}, et le second du théorème de changement de base
par le morphisme lisse~$g_i$~\sga{4}{xvi}{1.2}.  \end{démo}

Nous prouvons à présent qu'un morphisme régulier entre schémas excellents permet d'induire des
fonctions de dimension.

\begin{proposition2} \label{14-prop-fonction-dim-induite-regulier} Soient $f : Y \ra X$ un morphisme
régulier entre schémas excellents et~$\delta_X$ une fonction de dimension sur~$X$. La fonction
$\delta_Y : Y \ra \ZZ$ définie par $$\delta_Y(y) =  \delta_X(f(y))  - \dim \left( \OO_{Y_{f(y)},y}
\right) $$ est une fonction de dimension sur $Y$.  \end{proposition2}

\begin{démo} Comme la vérification est locale, il n'y a pas de mal à supposer $X$ et $Y$ strictement
locaux et $f$ local. Soit $\delta$ une fonction de dimension sur $Y$ ; son existence est assurée par
le théorème \ref{14-theoreme-univ-cat-fonction-dim}. Il suffit de montrer que $\delta_Y-\delta$ est
une fonction constante sur $Y$. Les fibres de~$f$ sont  régulières donc universellement caténaires
d'après~\ref{14-exemple-regulier-implique-univ-cat}. La
proposition~\ref{14-proposition--codim-fonc-dim} montre que la fonction qui à~$y\in Y$ associe  $$-
\dim \left( \OO_{Y_{f(y)},\: y} \right) $$ induit une fonction de dimension sur chacune des fibres
de~$f$. La fonction $\delta_Y-\delta$ est donc localement constante sur chaque fibre de~$f$. Il
résulte du lemme~\ref{14-lemme-cb-regulier} que ces fibres sont connexes : en effet, on a
$\HH^0(f^{-1}(x),\ZZ/n\ZZ)=\HH^0(x,\ZZ/n\ZZ)=\ZZ/n\ZZ$ pour tout $x\in X$ et tout entier~$n$
inversible sur~$X$. La fonction $\delta_Y-\delta$ est donc constante sur les fibres de~$f$ et
descend à $X$. Il suffit de montrer que~$\gamma=\delta_Y-\delta$ est localement constante sur $X$.
Une façon de calculer la valeur de $\gamma$ en un point $s$ de $X$ consiste à considérer le point
générique $\eta_s$ du schéma régulier connexe $f^{-1}(s)$, de sorte que
$\gamma(s)=\delta_X(s)-\delta(\eta_s)$. Soit $s' \spim s$ une spécialisation Zariski immédiate entre
deux points de $X$. Il s'agit de montrer que $\gamma(s)=\gamma(s')$. Vu que $\delta_X$ et $\delta$
sont des fonctions de dimension sur $X$ et $Y$ respectivement, pour montrer cela, il suffit de
savoir que $\eta_{s'}$ est une spécialisation immédiate de $\eta_s$. Pour montrer cela, quitte à
remplacer $X$ par le localisé en $s$ de l'adhérence de $s'$, on peut supposer que $X$ est local
intègre de dimension $1$, de point générique $s'$ et de point fermé $s$. Il s'agit alors de montrer
que la fibre $f^{-1}(s)$ est de codimension $1$ dans $Y$, ce qui résulte facilement du
\emph{Hauptidealsatz}.  \end{démo}

\subsection{Contre-exemple}

Rappelons l'exemple de~\ega{ii}{5.6.11} d'un schéma caténaire non universellement caténaire.  Soient
$k_0$ un corps et $k$ une extension purement transcendante de $k_0$  de degré de transcendance
infini. Notons $S=k[X]_{(X)}$ le localisé de l'anneau de la droite affine sur~$k$ en l'origine et
$V=S[T]$. Les idéaux maximaux $\MM=(X,T)$ et $\MM'=(XT-1)$ de $V$ sont respectivement de hauteur $2$
et $1$, et il existe un isomorphisme $\phi \: : \: V/\mathfrak{m} \iso V/\mathfrak{m}'$. Posons
$C=\{\:v\in V \: | \: \phi(\: v \: \mathrm{mod} \: \mathfrak{m}'\:) = v \: \mathrm{mod} \:
\mathfrak{m} \:\}$. C'est un sous-anneau de $V$ qui n'est pas de type fini sur $k$ et dont le
spectre est le quotient de $\SP(V)$ par la relation d'équivalence identifiant $\MM$ et $\MM'$. Il
est intègre et contient l'idéal maximal $\mathfrak{n}$ image réciproque de $\MM$.

\begin{proposition2}\label{14-propo-contre-ex} Le schéma $\SP(C)$ est nœthérien, quasi-excellent,
caténaire mais  non universellement caténaire. Le point fermé correspondant à l'idéal maximal
$\mathfrak{n}$ de $C$ est une spécialisation étale immédiate mais non Zariski immédiate du point
générique de $\SP(C)$.  \end{proposition2}

\begin{démo} Le caractère nœthérien est montré dans \ega{ii}{5.6.11} et le caractère quasi-excellent
dans \ega{ii}{7.8.4}. Le schéma $\SP(C)$ est caténaire d'après le
corollaire~\ref{14-corollaire-dim2-cat} car il est de dimension~$2$. On montre qu'il n'est pas
universellement caténaire en vérifiant que la dernière condition du
théorème~\ref{14-theoreme-Ratliff} est prise en défaut. En effet, le morphisme
$$\SP(V_{(\mathfrak{n})}) \ra \SP(C_{(\mathfrak{n})})$$ est fini  et la fibre du point fermé
correspondant à $\mathfrak{n}$ est constituée des points fermés correspondants à $\mathfrak{m}$ et
$\mathfrak{m}'$, qui sont de hauteur respective $1$ et $2$. Le schéma $\SP(C)$ est intègre de
dimension $2$ donc $\mathfrak{n}$ n'est pas une spécialisation Zariski immédiate du point générique.
\end{démo}

\ifx\danslelivre\undefined
\bibliography{../fichiers-configuration-latex/bibliographie}
\bibliographystyle{../fichiers-configuration-latex/smfalpha}


\begin{center}
XV. Théorème de Lefschetz affine\\
Vincent Pilloni et Benoît Stroh
\end{center}

\else
\setboolean{english}{false}
\setboolean{francais}{true}
\selectlanguage{francais}
\chapter{Théorème de Lefschetz affine}
\addtocontents{toc}{par Vincent Pilloni et Benoît Stroh}
\begin{center}
Vincent Pilloni et Benoît Stroh
\end{center}
\fi

\section{Énoncé du théorème et premières réductions}

\subsection{Énoncé}
\label{15-énoncés}

\subsubsection{}
Soient $X$ un schéma muni d'une fonction de dimension $\delta_X$
(\refext{XIV}{ddefi4}) et $n$ un entier inversible sur $X$.
Pour tout faisceau étale $ℱ$ de $𝐙/n$-modules sur $X$,
$$\delta_X(\mc F) = \mathrm{sup}\: \big\{ \delta_X(x) \ , \ x\in X \ | \
{\mc F}_{\xb} \neq 0 \big\}.$$
Rappelons (\refext{XIV}{14-coro-fonction-dim-induite-tf})
qu'un morphisme de type fini $f: Y → X$ induit une
fonction de dimension sur $Y$ ; nous la noterons ici $f^{⋆}δ_X$.
Le théorème principal de cet exposé est le suivant (voir aussi
\refext{0}{0-Lefschetz affine}).

\begin{theoreme2}
\label{15-theoreme-principal}
Supposons le schéma $X$ quasi-excellent et le morphisme $f:Y \ra X$ 
\emph{affine} de type fini. Alors, pour tout faisceau constructible $ℱ$
de $𝐙/n$-modules sur $Y$, on a :
\[
\delta_X(\R^{q}f_⋆(\mc{F})) ≤  f^⋆\delta_X(\mc{F}) - q .
\]
\end{theoreme2}

\begin{remarque2}
\label{15-LA sur trait}
Ce théorème a été déjà démontré en 1994 par O. Gabber
lorsque $X$ est de type fini sur un trait,
cf. \cite{PV@Illusie}. La démonstration du théorème
précédent procède notamment par réduction à ce cas.
\end{remarque2}

\subsection{Reformulation et réductions}

\subsubsection{}
Soient $f$ et $ℱ$ comme ci-dessus.
La conclusion du théorème signifie que pour
tout point géométrique $\xb$ de $X$ localisé en un
point $x$ et tout entier $q>f^⋆δ_X(ℱ)-δ_X(x)$, on a
\[
(\R f_⋆ ℱ)_{\xb}=H^q(Y_{(\xb)},ℱ)=0,
\]
où l'on note $Y_{(\xb)}$ le produit fibré $Y×_X X_{(\xb)}$.
Rappelons (\refext{XIV}{14-coro-schema-tf-sur-corps-fonc-dim})
que le schéma strictement local $X_{(\xb)}$
peut être muni la fonction de dimension $δ_{X_{(\xb)}}:t↦ \dim(\sur{\{t\}})$
(\refext{XIV}{14-coro-schema-tf-sur-corps-fonc-dim}) ;
c'est l'unique fonction de dimension nulle en $x$.
Notons l'inégalité $f^⋆δ_X(ℱ)-δ_X(x) ≥ f_{(\xb)}^⋆δ_{X_{(\xb)}}(ℱ)$,
triviale dans le cas où $δ_X(x)=0$, auquel on peut se ramener.
Ainsi, le théorème \ref{15-theoreme-principal} est équivalent
à l'énoncé suivant.

\begin{corollaire2}
\label{15-prop-principale-locale}
Supposons de plus le schéma $X$ quasi-excellent, \emph{strictement local},
muni de la fonction de dimension $δ_X:t↦ \dim(\sur{\{t\}})$.
Alors, pour tout faisceau constructible $ℱ$ de $𝐙/n$-modules
sur $Y$, on a :
\[
\HH^q(Y,\mc{F})=0 \text{ si } q > f^⋆\delta_X(\mc{F}).
\]
\end{corollaire2}

\subsubsection{}Procédant comme en \sga{4}{xiv}{4.4} pour se ramener au cas d'une
immersion ouverte affine puis utilisant la méthode de la trace
(\sga{4}{ix}{5.5} ou \sga{5}{I}{3.1.2}) pour se ramener au cas
des coefficients constants (voir aussi \refext{XIII}{13-ops constant}),
on montre que le théorème est également équivalent au corollaire suivant.
(On pourrait d'ailleurs supposer l'entier $n$ premier.)

\begin{corollaire2}
\label{15-prop-devissage-classique}
Soient $X$ un schéma strictement local quasi-excellent de dimension $d$,
$U$ un ouvert affine de $X$, et $n$ un entier inversible sur $X$.
Alors,
\[
\HH^q(U,𝐙/n)=0 \text{ si } q>d.
\]
\end{corollaire2}

\subsubsection{Réduction au cas complet}
Supposons dorénavant $X$ strictement local quasi-excellent,
de complété $\chap{X}$ en le point fermé, et
fixons un ouvert affine $U$ de $X$, dont on
note $\chap{U}$ l'image inverse sur $\chap{X}$.
Le morphisme naturel de~$\widehat{X}$ dans~$X$ est régulier car~$X$ est quasi-excellent.
En appliquant le lemme de changement de base par
un morphisme régulier (\refext{XIV}{14-lemme-cb-regulier}) au diagramme cartésien
$$\xymatrix{\widehat{U} \ar[d] \ar[r] & U\ar[d] \\
\widehat{X}\ar[r] & X}$$
on obtient $\HH^q(\widehat{U},\ZZ/n\ZZ)=\HH^q({U},\ZZ/n\ZZ)$
pour tout~$q\geq 0$. Voir également \cite[7.1.1]{TUBULAR@Fujiwara}
pour une autre démonstration.

\subsubsection{}Dans les deux sections qui vont suivre, nous
allons démontrer l'énoncé \ref{15-prop-devissage-classique}
dans le cas particulier où $X$ est local nœthérien complet à
corps résiduel séparablement clos.

\section{Pureté, combinatoire des branches et descente}
\label{15-section pureté et combinatoire}

\subsection{Pureté : rappel et une application}

\subsubsection{}Nous rappelons le théorème de pureté absolue démontré par O. Gabber
(\cite{AZUMINO@Fujiwara}). Par convention, on considère le schéma vide
comme un diviseur strictement à croisements normaux dont
l'ensemble des branches est indexé par l'ensemble vide.

\begin{theoreme2}[\refext{XVI}{16-theoreme-purete-absolue}]\label{15-theoreme-purete-rappel}
Soient $X$ un schéma régulier, $Z$ un diviseur strictement à
croisements normaux de complémentaire $j:U = X\backslash Z
\hookrightarrow X$ et de branches~$\{Z_i\}_{i\in I}$, et~$n$ un entier inversible sur~$X$. Il existe des isomorphismes canoniques
\begin{eqnarray*}
\R^1j_{\star}\: (\ZZ/n\ZZ) & \isolong & \bigoplus_{i\in I} \: (\ZZ/n\ZZ)_{Z_i}(-1) \\
\R^qj_{\star}\: (\ZZ/n\ZZ) & \isolong & \bigwedge^q\: \R^1 j_{\star}\: (\ZZ/n\ZZ)
\end{eqnarray*}
\end{theoreme2}

\subsubsection{Combinatoire des branches : notations}
Soient $g:X′ → X$ un morphisme entre schémas, et $U$ un ouvert de $X$.
Notons $j:U ↪ X$ l'immersion ouverte, $j′:U′ ↪ X′$
l'immersion ouverte qui s'en déduit par changement de base
et $Z$ et $Z′$ les fermés complémentaires respectifs.
Enfin on se donne un fermé $F ⊆ Z$ dont on note $F′$
l'image inverse.

\begin{définition2}
\label{15-définition même combinatoire}
On dit que $(Z ↪ X)$ et $(Z′ ↪ X′)$ ont
\emph{même combinatoire le long de $F$}
si pour tout point géométrique $\sur{z}′$ de $F′$
d'image le point géométrique $\sur{z}$ de $F$, les propriétés suivantes sont satisfaites :
\begin{enumerate}
\item les schémas $X_{(\sur{z})}$ et $X′_{(\sur{z}′)}$ sont
réguliers ;
\item le fermé $Z_{(\sur{z})}$ est un diviseur à croisements normaux
définis par des équations $f₁,…,f_r$ ;
\item les fonctions $g^⋆f₁,…,g^⋆f_r$ forment une famille
libre de $𝔪/𝔪²$, où $𝔪$ est l'idéal maximal de $Z′_{(\sur{z}′)}$.
\end{enumerate}
\end{définition2}

\subsubsection{}Il résulte de (i) et (iii) ci-dessus
que le fermé $Z′_{(\sur{z}′)}$ est un diviseur à croisements
normaux dans $X′_{(\sur{z}′)}$, de même codimension que
$Z_{(\sur{z})}$ dans $X_{(\sur{z})}$.

\subsubsection{}
\label{15-abus terminologie}
Lorsque $X$ est un schéma sur une
base $S$, et $F$ est un fermé de ce dernier,
on s'autorise à dire « … le long de $F$ »
pour « … le long de l'image inverse $F×_S X$ ».

\begin{proposition2}
\label{15-pureté et combinatoire}
Supposons que $(Z ↪ X)$ et $(Z′ ↪ X′)$ aient
même combinatoire le long d'un fermé $F$ de $X$.
Alors, pour tout entier $n$ inversible sur $X$, le morphisme d'adjonction
\[
(\R j_⋆ 𝐙/n)_{|F′} → (\R j′_⋆ 𝐙/n)_{|F′}
\]
est un \emph{isomorphisme}.
\end{proposition2}

\begin{démo}
Quitte à localiser en des points géométriques $\sur{z}′$ et $\sur{z}$,
on peut supposer les schémas strictement locaux et le morphisme
$X′ → X$ local. Il faut alors montrer que
$\R Γ(U, 𝐙/n) → \R Γ(U′, 𝐙/n)$ est un isomorphisme. D'après
le théorème de pureté \ref{15-theoreme-purete-rappel},
il suffit de montrer que le morphisme induit sur le $\HH¹$
est un isomorphisme, ce qui résulte aussitôt de la structure de ces
groupes et de ce que la classe associée à une branche $Z_i=V(f_i)$ de $Z$
est envoyée par restriction sur la classe de la
branche $g^{-1}(Z_i)=V(g^⋆f_i)$ (cf. \refext{XVI}{16-section-gysin}).
\end{démo}

\subsection{Application du théorème de descente fléchée}
\label{15-paragraphe-descente}

\subsubsection{}
Soient $X$ un schéma strictement local nœthérien, $j:U ↪ X$
un ouvert, $i_x:x ↪ X$ l'immersion fermée, et $ε:X_∙ → X$
un hyperrecouvrement pour la topologie des altérations (\refext{II}{2-définition altL}).
La proposition suivante — où les morphismes obtenus
par changement de base sont notés de façon évidente —
est un corollaire immédiat du théorème \refext{XII}{12-cb orienté}
et du fait que la cohomologie de $U$ est la fibre en $x$ de
l'image directe par $j$.

\begin{proposition2}
\label{15-cohomologie via hyperrecouvrement}
Sous les hypothèses précédentes, le morphisme d'adjonction
\[
\R Γ(U,𝐙/n) → \R {ε_x}_⋆ \big( i_{x ∙}^⋆ \R {j_∙}_⋆ 𝐙/n \big)
\]
est un \emph{isomorphisme}.
\end{proposition2}

\subsubsection{}
Supposons maintenant donné un morphisme \emph{local} $X′ → X$ de schémas strictement
locaux nœthériens. Comme précédemment, on note $U$ un
ouvert de $X$, $Z$ son complémentaire et $x$ le point fermé de $X$.
À tout hyperrecouvrement pour la topologie des altérations
$X_∙ → X$ de $X$ est associé par changement de base un
hyperrecouvrement de $X′$.

\begin{proposition2}
\label{15-isomorphisme via HR et combinatoire}
Supposons que pour chaque entier $q ≥ 0$
les fermés $(Z_q ↪ X_q)$ et $(Z′_q ↪ X′_q)$ aient
même combinatoire le long de la fibre spéciale $(X_{q})_x$.
Alors le morphisme d'adjonction
\[
\RΓ(U,𝐙/n) → \R Γ(U′,𝐙/n)
\]
est un isomorphisme.
De plus, si l'on fait l'hypothèse précédente
pour les seuls entiers $q ≤ N+1$, où $N$ est un
entier quelconque, les morphismes
$\HH^q(U,𝐙/n) → \HH^q(U′,𝐙/n)$ sont des isomorphismes pour
chaque $q ≤ N$.
\end{proposition2}

\begin{démo}
Le premier point est conséquence immédiate
des deux propositions précédentes.
En effet, d'après \ref{15-cohomologie via hyperrecouvrement}
et l'invariance de la cohomologie par changement de corps
séparablement clos, soit
\[
\Big( \R {ε_x}_⋆ \big( i_{x ∙}^⋆ \R {j_∙}_⋆ 𝐙/n \big) \Big)_{|x′} = \R
{ε′_{x′}}_⋆ \Big( \big({i_{x∙}}^⋆ \R {j_∙}_⋆ 𝐙/n\big)_{|x′_∙}\Big),
\]
il suffit de démontrer que pour chaque $q$,
l'adjonction $(\R {j_q}_⋆ 𝐙/n)_{|{X′_q}_x} → (\R {j′_q}_⋆ 𝐙/n)_{|{X′_q}_x}$
est un isomorphisme. Cela résulte de l'hypothèse
combinatoire et de \ref{15-pureté et combinatoire}.
La variante tronquée est un corollaire de la démonstration
précédente et du fait que la cohomologie en degré $q$
ne dépend que des étages $≤q+1$.
\end{démo}

\begin{remarque2}
Dans ce critère, on ne fait d'hypothèses qu'en les points des fibres spéciales des hyper-recouvrements ; c'est ce qui en fait toute sa force.
\end{remarque2}

\section{Uniformisation et approximation des données}\label{15-section-5}

\subsection{Notations}

\subsubsection{}Soit $X,U,Z$ et $n$ comme à la fin
du ❡\ref{15-énoncés} : le schéma $X$ est strictement local
complet, $U$ en est un ouvert \emph{affine}, $Z=X-U$ (muni de la
structure réduite), et $n$ est un entier inversible sur $X$.
On veut démontrer \ref{15-prop-devissage-classique} dans ce cas.
Fixons un entier $N$.

\subsubsection{}Il résulte du théorème d'uniformisation
(\refext{VII}{7-uniformisation}), complété par
\refext{XIII}{13-cb}, qu'il existe un hyperrecouvrement
pour la topologie des altérations $ε:X_∙ → X$ tel que
chaque $X_q$ soit régulier et chaque $Z_q=Z×_X X_q$ soit
un diviseur strictement à croisements normaux. Il est
loisible de supposer les schémas $X_q$ affines.

\subsubsection{}Soient $k$ le corps résiduel de $X$,
soit $C$ un anneau de Cohen de corps résiduel $k$
(\refext{IV}{4-19.8.4}) et $S=\Spec(C)$.
Il résulte du théorème de structure des anneaux locaux
nœthériens (\egazero{iv}{19.8.8}) qu'il existe un entier $m$
et une immersion fermée de $X$ dans le complété en l'origine
(de la fibre spéciale sur $S$) de l'espace affine $𝐀^m_S$.
Notons $\chap{E}$ ce complété, $E$ son analogue hensélien
(hensélisé de l'espace affine en l'origine) et enfin $e$
le point fermé de ce dernier.

\subsubsection{}Le schéma $E$ étant quasi-excellent,
le morphisme de complétion $\chap{E} → E$ est (local) régulier
de sorte que, d'après le théorème de Popescu,
on peut écrire :
\[
\chap{E} = \lim_α E_α,
\]
où les $E_α → E$ sont des morphismes essentiellement lisses
entre schémas strictement locaux, la limite étant filtrante.
Notons que les schémas $E_α$ sont essentiellement de type
fini sur le \emph{trait} $S$.

\subsection{Passage à la limite}

\subsubsection{}
Quitte à restreindre l'ensemble d'indices, c'est-à-dire à supposer
$α ≫ 0$, les principes généraux de \ega{IV₃}{§8}, joints au
fait que les schémas $X,Z,U$ et les $X_q$ pour $q ≤ N$
sont de présentation finie sur $\chap{E}$,
entraînent l'existence de diagrammes à carrés cartésiens
\begin{eqnarray*}
\xymatrix{ U_{≤ N, α} \ar[r] \ar[d] & X_{≤ N, α}   \ar[d] & Z_{≤ N, α}\ar[l] \ar[d] \\
U_{\alpha} \ar[r] \ar[dr] & X_{\alpha} \ar[d] & Z_{\alpha} \ar[l] \ar[dl] \\
& E_α & }
\end{eqnarray*}
dont l'analogue sur $\chap{E}$ se déduit par changement de
base $\chap{E} → E_α$. De plus, on peut supposer que pour
chaque $α$, $X_α → E_α$ est une immersion fermée — de sorte que $X_α$ est strictement local —,
et $U_α → X_α$ une immersion ouverte affine de
complémentaire $Z_α$.

\begin{remarque2}
Les schémas $X_q$ et $X_{q α}$ ont même fibre spéciale pour tout~$q ≤  N$.
\end{remarque2}

\subsubsection{}Il résulte du théorème de la forme standard
(\refext{II}{2-forme standard l'-altérations})\footnote{La
démonstration dans le cas « $ℓ=1$ », qui est celui
utilisé ici, est valable sans supposer $X$ irréductible.}
et de \ega{IV₃}{8.10.5} que l'on peut supposer — quitte à
considérer de grands indices — que les $X_{≤ N, α} → X$
sont des hyperrecouvrements (tronqués) pour la topologie des
altérations.

\subsubsection{}Vérifions que l'on peut supposer
que pour chaque $α$ et chaque $q ≤ N$, le « modèle » $X_{q α}$ de $X_q$
est régulier le long de sa fibre spéciale.
Fixons $q$ puis oublions le. Le schéma $X$
est donc maintenant supposé régulier, affine, de type fini
sur $\chap{E}$. Le problème est local pour la topologie de Zariski ;
on peut donc supposer que $X$ est un fermé $V(f₁,…,f_c)$
de codimension $c$ dans un espace affine $Y=𝐀^m_{\chap{E}}$.
Quitte à supposer les $α$ suffisamment grands,
il existe des fonctions $f_{i α}$
telles que $f_i=g^⋆f_{i α}$ pour chaque $i$, où $g$
est le morphisme évident $Y → Y_α=𝐀^m_{E_α}$.
Soit $x$ un point de $Y$ appartenant à la fibre
spéciale de $X$ et soit $x_α=g(x)$ son image.
Notons $𝔪_x$ (resp. $𝔪_{x_α}$) l'idéal maximal de $𝒪_{Y,x}$
(resp. $𝒪_{Y_α,x_α}$). Par hypothèse de régularité de $X$ en $x$,
les images $\sur{f₁},…,\sur{f_c}$ des $f_i$
dans $𝔪_x/ 𝔪²_x$ sont linéairement indépendantes
sur $κ(x)=𝒪_{Y,x}/𝔪_x$. Le diagramme
\begin{center}
 \begin{tikzpicture}[auto]
  \matrix(diag)[matrix of math nodes,column sep=3em,row sep=5ex]
  {|(f)|  f_i     & |(m)|  𝔪_x     & |(t)|  𝔪_x/𝔪²_x    \\
   |(fa)| f_{i α} & |(ma)| 𝔪_{x_α} & |(ta)| 𝔪_{x_α}/𝔪²_{x_α} \\};
  \draw[->>] (m) -- (t);
  \draw[->] (ma) -- node{$g^⋆$} (m);
  \draw[->>] (ma) -- (ta);
  \draw[->] (ta) --  (t);
  \draw[|->] (fa) --  (f);
 \end{tikzpicture}
\end{center}
étant commutatif, il résulte de l'égalité $κ(x)=κ(x_α)$
que les images $\sur{f_{i α}}$ des $f_{i α}$ dans $ 𝔪_{x_α}/𝔪²_{x_α}$
sont linéairement indépendantes sur $κ(x_α)$.
Il en résulte que le sous-schéma fermé $X_α=V(f_{1 α},…,f_{c α})$
de $Y_α$ est régulier en $x_α$.

\subsubsection{}On montre de même que
l'on peut supposer que pour chaque $α$ et chaque $q ≤ N$,
les immersions $(Z_q ↪ X_q)$ et $(Z_{q α} ↪ X_{q α})$ ont même combinatoire le long
du point fermé $e ∈ E$, c'est-à-dire le long des fibres spéciales.

\subsubsection{}
Il résulte de la proposition \ref{15-isomorphisme via HR et combinatoire}
que les morphismes d'adjonction
\[
\HH^q(U_α,𝐙/n) → \HH^q(U,𝐙/n)
\]
sont des isomorphismes pour $q<N$. Nous allons montrer
que si $q>d=\dim(X)$ et $α$ est suffisamment grand,
on a $\HH^q(U_α,𝐙/n)=0$. Ceci achèvera la démonstration
du théorème de Lefschetz affine. Notons que l'ouvert
affine $U_α$ du schéma strictement local $X_α$
n'est en général pas de dimension $d$.

\subsection{Utilisation d'une section}

\subsubsection{}
Soit $σ : E → E_α$ une section de $E_α → E$
et $X^σ_α$ (resp. $U^σ_α,Z^σ_α$) le $E$-schéma déduit de $X_α$
(resp. $U_α$, $Z_α$)
par changement de base.
Notons également $X^σ_{≤ N,α}$ l'hyperrecouvrement
de $X^σ_α$ pour la topologie des altérations obtenu
à partir de $X_{≤ N,α} → X_α$ par changement de base.
On rappelle que l'on remplace les schémas obtenus
par produit fibré usuel par la réunion des adhérences
de leurs composantes irréductibles dominant une composante
irréductible de $X^σ_α$. (En d'autres termes, on applique le
foncteur $T↦ T_{\maxdom}$ de \refext{II}{2-dom}.)
Enfin $U^σ_{≤ N,α}$ (resp. $Z^σ_{≤ N,α}$) est l'ouvert (resp. le fermé)
simplicial évident.

\subsubsection{}Il résulte de \refext{III}{3-perm} et
\refext{III}{3-reduc} (démonstration)
que si $α$ est suffisamment grand et $σ$ est
suffisamment proche de l'identité,
alors, les immersions fermées
$(Z^σ_{qα} ↪ X^σ_{qα})$ et
$(Z_{qα} ↪ X_{qα})$ ont \emph{même combinatoire}
le long de la fibre spéciale au-dessus de $E$ pour
chaque $q ≤ N$.
Il en résulte comme ci-dessus que le morphisme
\[
\HH^q(U_α,𝐙/n) → \HH^q(U^σ_α,𝐙/n)
\]
est un isomorphisme. Comme l'ouvert $U^σ_α$
est affine dans $X^σ_α$ de dimension $d$
et essentiellement de type fini sur le trait $S$,
il résulte par passage à la limite
du théorème de Lefschetz (dû à O. Gabber)
\cite[2.4]{PV@Illusie}, on a
\[
\HH^q(U^σ_α,𝐙/n)=0 \text{ si } d<q<N.
\]
Finalement, $\HH^q(U,𝐙/n)=0$ si $q>d=\dim(X)$. CQFD.

\ifx\danslelivre\undefined
\bibliography{../fichiers-configuration-latex/bibliographie}
\bibliographystyle{../fichiers-configuration-latex/smfalpha}


\begin{center}
XVI. Classes de Chern, morphismes de Gysin, pureté absolue\\
Joël Riou\\
exposés oraux des 16 juin 2006 et 13 avril 2007\\
\end{center}

\begin{center}
\verb$Id: gysin.tex,v 1.98 2009-08-20 10:06:45 cvs Exp $
\end{center}

\tableofcontents

\else
\setboolean{english}{false}
\setboolean{francais}{true}
\selectlanguage{francais}
\chapter{Classes de Chern, morphismes de Gysin, pureté absolue}
\addtocontents{toc}{par Joël Riou}
\begin{center}
Joël Riou
\end{center}
\fi

Dans ces notes, on présente la nouvelle démonstration par Ofer Gabber du
théorème de pureté cohomologique absolue, annoncée dans
\cite{ABELIEN@Gabber}. La
section~\ref{16-section-classes-chern} rappelle la construction des classes de
Chern en cohomologie étale. Celles-ci servent dans la
section~\ref{16-section-gysin} qui consiste en la construction et l'étude des
propriétés des morphismes de Gysin associés aux morphismes d'intersection
complète lissifiables. Dans la section~\ref{16-section-purete}, ces morphismes
de Gysin sont utilisés pour donner une formulation précise du théorème de
pureté absolue (théorème~\ref{16-theoreme-purete-absolue}). La démonstration
du théorème de pureté (différente de celle rédigée dans
\cite{AZUMINO@Fujiwara})
s'appuie notamment sur les résultats de géométrie logarithmique établis
dans l'exposé VI.

\bigskip

Dans tout cet exposé, on fixe un entier naturel $n\geq 1$. Tous les schémas
seront supposés être des schémas sur $\SP\left(\ZZ\cro{\frac 1 n}\right)$. On note
$\Lambda$ le faisceau d'anneaux constant de valeur $\ZZ/n\ZZ$, $\Lambda(1)$
le faisceau des racines $n$-ièmes de l'unité (pour la topologie étale)
et $\Lambda(r)$ ses
puissances tensorielles, auxquelles on peut donner un sens pour tout
$r\in\ZZ$.

\section{Classes de Chern}
\label{16-section-classes-chern}

Dans cette section, on rappelle la construction des classes de Chern
des fibrés vectoriels sur des schémas généraux à valeurs dans la 
cohomologie étale.
On s'appuie sur le calcul de la cohomologie étale des fibrés projectifs de
\sga{5}{vii}{2} et sur la méthode de \cite{CHERN@Grothendieck}.
Les démonstrations sont parfois différentes de celles
de \sga{5}{vii}{3} : on s'est efforcé de donner une présentation aussi
\guil{économique} que possible.

À la différence de l'exposé oral qui utilisait un langage géométrique,
dans ces notes, un fibré vectoriel est un Module $\mathscr E$
localement libre de rang fini et le fibré projectif de $\mathscr E$ est
le fibré des hyperplans défini dans \ega{ii}{4.1.1} :
$\PP(\mathscr E)=\Proj(\SS^\star\mathscr E)$ où $\SS^\star\mathscr E$
est l'Algèbre symétrique de $\mathscr E$.

\begin{definition}\label{16-definition-c-1}
Soit $X$ un $\ZZn$-schéma. Soit $\mathscr L$ un fibré en droites sur $X$.
Le faisceau des sections inversibles de $\mathscr L$ est naturellement
muni d'une structure de torseur sous le schéma en groupes $\Gm$. La classe
d'isomorphisme de $\mathscr L$ définit donc un élément dans $H^1_\et(X,\Gm)$.
On note $c_1(L)\in \HH^2_\et(X,\Lambda(1))$ l'image de cet élément par le
morphisme de bord déduit de la suite exacte courte de Kummer :
\[0\to \Lambda(1)\to \Gm\vers {\cro n}\Gm \to 0\;\footnote{%
Une grande liberté d'interprétation est laissée à l'imagination du lecteur
au sujet des conventions de signes à utiliser dans cette définition de
$c_1$. La même liberté me semblant exister dans \sga{4}{xviii}{1.1}
pour le morphisme
trace, il me paraît illusoire de faire davantage de zèle. Cependant, à un
endroit de ces notes, ce choix ne sera pas neutre : il faut que les choix
de signes soient mutuellement cohérents de sorte que l'on ait une
compatibilité entre le degré des fibrés en droites sur la droite
projective, le morphisme trace et la première classe de Chern (cf.
démonstration du lemme~\ref{16-lemme-compatibilite-trace-c-un}).}\text{.}\]
\end{definition}

Si $\mathscr L$ et $\mathscr L'$ sont deux 
fibrés en droites sur $X$, on a la relation
d'additivité :
\[c_1(\mathscr L\otimes \mathscr L')=
c_1(\mathscr L)+c_1(\mathscr L')\in \HH^2(X,\Lambda(1))\;%
\footnote{Il existe des théories cohomologiques
\guil{orientées} pour lesquelles cette propriété de la première classe de
Chern n'est pas satisfaite, cf.~\cite{CA1@Morel-Levine}.}\text{.}\]

Notons que les classes de Chern de fibrés en droites résident dans les degrés
pairs de la cohomologie étale, elles commutent donc avec toutes les classes
de cohomologie. Notons
aussi que si $f\colon Y\to X$ est un morphisme et 
$\mathscr L$ un fibré en droites
sur $X$, alors $f^\star(c_1(\mathscr L))=c_1(f^\star \mathscr L)$.

\begin{theoreme}[Formule du fibré projectif]%
\label{16-theoreme-formule-du-fibre-projectif}
Soit $X$ un $\ZZn$-schéma. Soit $\mathscr E$ un fibré 
vectoriel de rang constant $r$
sur $X$.
On note $\pi\colon\PP(\mathscr E)\to X$
le fibré projectif de $\mathscr E$.
On pose $\xi=c_1(\OO(1))\in \HH^2(X,\Lambda(1))$\;\footnote{Le faisceau 
$\OO(1)$ est le faisceau fondamental sur $\PP(\mathscr E)$ : c'est le 
quotient inversible de $\pi^\star \mathscr E$ par 
l'hyperplan universel.}.
Alors, les puissances $\xi^i\in \HH^{2i}(X,\Lambda(i))$ de $\xi$
définissent un isomorphisme dans $\detaleplus X$ :
\[(1,\xi,\dots,\xi^{r-1})\colon
\bigoplus_{i=0}^{r-1} \Lambda(-i)[-2i]
\iso \R\pi_\star \Lambda\]
\end{theoreme}

D'après le théorème de changement de base pour un morphisme propre, on peut
supposer que $X$ est le spectre d'un corps algébriquement clos $k$. On se
ramène ainsi au calcul de l'algèbre de cohomologie étale des espaces
projectifs sur $k$, cf.~\sga{5}{vii}{2}.

\begin{theoreme}\label{16-theoreme-classes-de-chern}
Il existe une unique manière de définir, pour tout $\ZZn$-schéma $X$ et
tout fibré vectoriel $\mathscr E$, 
des éléments $c_i(\mathscr E)\in \HH^{2i}_\et(X,\Lambda(i))$
pour tout $i\in\NN$ de sorte que si l'on définit la série formelle
$c_t(\mathscr E)=\sum_{i\geq 0} c_i(\mathscr E)t^i$,
on ait les propriétés suivantes :
\begin{itemize}
\item la série formelle $c_t(\mathscr E)$ ne dépend que 
de la classe d'isomorphisme
du fibré vectoriel $\mathscr E$ sur le $\ZZn$-schéma $X$ ;
\item si $f\colon Y\to X$ est un morphisme de $\ZZn$-schémas
et $\mathscr E$ un fibré vectoriel sur
$X$, alors $f^\star(c_t(\mathscr E))=c_t(f^\star \mathscr E)$ ;
\item si $0\to \mathscr E'\to \mathscr E\to \mathscr E''\to 0$ est une
suite exacte courte de fibrés vectoriels sur un $\ZZn$-schéma $X$, on a la
relation de Cartan-Whitney :
\[c_t(\mathscr E)=c_t(\mathscr E')c_t(\mathscr E'')\;\text{;}\]
\item si $\mathscr L$ est un fibré en droites 
sur un $\ZZn$-schéma $X$, la classe
$c_1(\mathscr L)$ est celle de la définition~\ref{16-definition-c-1} et 
\[c_t(\mathscr L)=1+c_1(\mathscr L)t\;\text{.}\]
\end{itemize}

On a alors les relations $c_0(\mathscr E)=1$ et 
$c_i(\mathscr E)=0$ pour $i>\rang \mathscr E$ pour tout
fibré vectoriel $\mathscr E$ sur un $\ZZn$-schéma $X$.
\end{theoreme}

La démonstration utilise plusieurs constructions géométriques :

\begin{proposition}[Principe de scindage I]\label{16-proposition-scindage-1}
Soit $X$ un $\ZZn$-schéma. Soit $\mathscr E$ un 
fibré vectoriel de rang $r$. On note
$\pi\colon \Drap(\mathscr E)\to X$ le fibré des drapeaux complets de
$\mathscr E$. Les
propriétés suivantes sont satisfaites :
\begin{itemize}
\item le fibré vectoriel $\pi^\star \mathscr E$ admet une filtration (canonique)
$0=\mathscr M_0\subset \mathscr M_1\subset \dots\subset \mathscr M_r=
\pi^\star \mathscr E$ par des fibrés
vectoriels de sorte que pour tout entier $1\leq i\leq r$, le quotient
$\mathscr L_i=\mathscr M_i/\mathscr M_{i-1}$ soit un fibré en droites ;
\item le morphisme canonique $\Lambda\to \R\pi_\star\Lambda$ est un
monomorphisme scindé dans $\detaleplus X$.
\end{itemize}
\end{proposition}

La seule propriété non triviale réside dans le fait que $\Lambda
\to\R\pi_\star\Lambda$ soit un monomorphisme scindé. En remarquant
que la projection $\Drap(\mathscr E)\to X$ peut s'écrire comme un composé
de $r$ projections de fibrés projectifs, ceci se déduit de la formule du
fibré projectif\;\footnote{Plus précisément, Grothendieck a montré
(cf.~\cite{PFTI@Grothendieck}, ou \sga{6}{vi}{4.6} pour le même argument
dans le cas de la $K$-théorie algébrique) que la théorie des classes de
Chern permettait de calculer l'algèbre de cohomologie des fibrés de
drapeaux, fussent-ils incomplets.}.

\begin{proposition}[Principe de scindage II]\label{16-proposition-scindage-2}
Soit $X$ un $\ZZn$-schéma. Soit $(E):0\to \mathscr E' \to \mathscr E \vers p 
\mathscr E''\to 0$
une suite exacte courte de fibrés vectoriels sur $X$. On note $\Sect(E)$
le $X$-schéma défini par le fait que pour tout $X$-schéma $f\colon Y\to
X$, l'ensemble des $X$-morphismes $Y\to \Sect(E)$ s'identifie
naturellement à l'ensemble des sections de la surjection de fibrés
vectoriels 
$f^\star(p)\colon f^\star \mathscr E\to f^\star \mathscr E''$ sur
$Y$\;\footnote{Je remercie Dennis Eriksson de m'avoir signalé cette
construction.}.
Le $Y$-schéma $\Sect(E)$ est naturellement muni d'une structure de 
torseur 
sous le $Y$-schéma en groupes vectoriel d'homomorphismes 
$\SheafHom(\mathscr E'',\mathscr E')$.
Notons $\pi\colon \Sect(E)\to X$ la projection. Les propriétés
suivantes sont satisfaites :
\begin{itemize}
\item l'image inverse par $\pi\colon \Sect(E)\to X$
de la suite exacte de fibrés vectoriels
$E$ est (canoniquement) scindée ;
\item le morphisme canonique $\Lambda\to \R\pi_\star\Lambda$
est un isomorphisme dans $\detaleplus X$.
\end{itemize}
\end{proposition}

L'existence de $\Sect(E)$ est évidente, la question étant de
nature locale sur $X$. Localement pour la topologie de Zariski sur $X$, la
projection $\pi$ est la projection depuis un espace affine, l'isomorphisme
$\Lambda\iso \R\pi_\star\Lambda$ résulte alors de l'invariance par
homotopie de la cohomologie étale pour les $\ZZn$-schémas (\sga{4}{xv}{2.2}).

\medskip

Démontrons le théorème~\ref{16-theoreme-classes-de-chern}. Grâce aux
propositions~\ref{16-proposition-scindage-1} et
\ref{16-proposition-scindage-2}, l'unicité est évidente. Il s'agit donc de
construire une théorie des classes de Chern satisfaisant les propriétés
demandées. Soit $\mathscr E$ un fibré vectoriel (que l'on peut
supposer de rang constant $r$) sur un $\ZZn$-schéma $X$. On considère le
fibré projectif $\PP(\mathscr E)$ sur $X$.
On note $\xi=c_1(\OO(1))$. D'après la
formule du fibré projectif
(théorème~\ref{16-theoreme-formule-du-fibre-projectif}), il existe d'uniques éléments, notés
$c_i(\mathscr E)\in \HH^{2i}(X,\Lambda(i))$ pour $1\leq i\leq r$ tels 
que l'on ait la
relation
\[\xi^r - c_1(\mathscr E)\xi^{r-1}+c_2(\mathscr E)\xi^{r-2}
+\dots+(-1)^r c_r(\mathscr E)=0\in
\HH^{2r}(\PP(\mathscr E),\Lambda(r))\;\text{.}\]
On pose $c_0(\mathscr E)=1$ et $c_i(\mathscr E)=0$ pour $i>r$.
Dans le cas où $\mathscr E$ est un
fibré en droites, $\PP(\mathscr E)\simeq X$
et $\OO(1)\simeq \mathscr E$, ce qui 
montre que cette définition étend la précédente pour les fibrés en droites.
La seule propriété non évidente est la formule de
Cartan-Whitney. Par principe de scindage
(propositions~\ref{16-proposition-scindage-1} et \ref{16-proposition-scindage-2}),
il suffit d'établir la formule
suivante :

\begin{lemme}
Soit $X$ un $\ZZn$-schéma. Soit $(\mathscr L_i)_{1\leq i\leq r}$ 
une famille finie
de fibrés en droites sur $X$, soit $\mathscr E=\bigoplus_{1\leq i\leq r}
\mathscr L_i$ leur
somme directe. Dans $\HH^{2r}(\PP(\mathscr E),\Lambda(r))$, on a la relation :
\[\prod_{i=1}^r(\xi-c_1(\mathscr L_i))=0\]
où $\xi=c_1(\OO(1))$. Autrement dit, 
\[c_t(\mathscr E)=\prod_{i=1}^r c_t(\mathscr L_i)\;\text{.}\]
\end{lemme}

L'argument qui suit est inspiré de \cite{PF@Panin-Smyrnov}. 
Pour $1\leq i\leq r$, on note $H_i$ l'hyperplan projectif de $\PP(\mathscr E)$
défini par l'inclusion $\mathscr L_i\to \mathscr E$. Notons $\pi\colon
\PP(\mathscr E)\to X$ la projection. Le morphisme canonique
$\pi^\star\mathscr L_i\to \OO(1)$ induit un isomorphisme sur l'ouvert
complémentaire de $H_i$ dans $\PP(\mathscr E)$.
On en déduit que l'élément $\xi-c_1(\mathscr L_i)$
de $\HH^2(X,\Lambda(1))$ peut être relevé en un élément $x_i$ du groupe de 
cohomologie à supports $\HH^2_{H_i}(X,\Lambda(1))$\;\footnote{Pour 
le moment,
peu importe de fixer un relèvement canonique.}.
Le produit des éléments $x_i$ vit naturellement dans le
groupe de cohomologie à support
$\HH^{2r}_{\cap_{1\leq i \leq r} H_i}(\PP(E),\Lambda(i))$ qui est nul puisque
l'intersection de ces $r$ hyperplans est vide ; on en
déduit la formule voulue par oubli du support.

\begin{proposition}
Soit $\mathscr E$ un fibré vectoriel sur un $\ZZn$-schéma $X$.
Pour tout entier
naturel $i$, on a l'égalité :
\[c_i(\mathscr E^\dual)=(-1)^i c_i(\mathscr E)\;\text{;}\]
autrement dit, on a une formule de changement de variables :
\[c_t(\mathscr E^\dual)=c_{-t}(\mathscr E)\;\text{.}\]
\end{proposition}

Grâce à la relation de Cartan-Whitney et au principe de scindage, on peut 
se ramener au cas où $\mathscr E$ est un fibré en droites. Cela résulte alors 
du fait que $c_1\colon \Pic(X)\to \HH^2(X,\Lambda(1))$ soit un homomorphisme
de groupes.

\section{Morphismes de Gysin}
\label{16-section-gysin}

Étant donné un morphisme d'intersection complète $X\vers f S$ entre
$\ZZn$-schémas vérifiant certaines hypothèses techniques, on va construire
un morphisme de Gysin $\Cl f\colon \Lambda\to f^?\Lambda$ où
$f^?=f^!(-d)[-2d]$ ($d$ est la dimension relative virtuelle de $f$).
Ces morphismes de Gysin seront
compatibles à la composition des morphismes d'intersection complète.

L'essentiel de cette section est consacrée à la construction de ces
morphismes de Gysin dans le cas des immersions régulières. Le morphisme trace
permettra de faire la construction dans le cas des morphismes lisses.
Ces deux définitions se recolleront pour donner la
définition~\ref{16-definition-cl-f-lissifiable} dans le cas général et le
théorème~\ref{16-theoreme-cl-f-transitivite} établira la compatibilité à la
composition de ces morphismes de Gysin.

\subsection{Première classe de Chern d'un pseudo-diviseur}

Soit $\mathscr L$ un fibré en droites sur un $\ZZn$-schéma $X$, $Z$ un
fermé de $X$ et $U$ l'ouvert complémentaire. On suppose donnée une section
inversible $s\colon \OO_U\iso {\mathscr L}_{|U}$. Au couple $(\mathscr
L,s)$ est canoniquement associée une classe $c_1(\mathscr L,s)\in
\HH^2_Z(X,\Lambda(1))$ induisant $c_1(\mathscr L)\in \HH^2(X,\Lambda(1))$ par
oubli du support (construire un élément de $\HH^1_Z(X,\Gm)$ et
utiliser la suite exacte de Kummer).

La classe $c_1(\mathscr L,s)$ définit un morphisme
$\Lambda_Z=\Lambda_X/\Lambda_U\to \Lambda_X(1)\cro 2$ dans $\detaleplus X$. En
\guil{composant} un tel morphisme avec une classe de cohomologie de $Z$
représentée par un morphisme $\Lambda_Z\to\Lambda_Z(q)\cro p$, il vient que
$c_1(\mathscr L,s)$ induit des morphismes \[\Gys_{(\mathscr L,s)}\colon
\HH^p(Z,\Lambda(q))\to \HH^{p+2}_Z(X,\Lambda(q+1))\;\text{.}\]

\begin{definition2}
Si $Z\to X$ est une immersion régulière de codimension $1$ définie par un
Idéal (inversible) $\II$, on pose $\Gys_{Z\subset X}=\Gys_{(\II
,1_{X-Z})}$.
\end{definition2}

On prendra garde au fait que \emph{via} les identifications usuelles, le
faisceau inversible $\II$ et le diviseur effectif $Z$ ont des
classes opposées dans le groupe de Picard de $X$.

\subsection{Classes fondamentales généralisées}

Pour étudier la compatibilité à la composition des classes
fondamentales définies dans \cite[§1]{AZUMINO@Fujiwara} dans le cas des
immersions régulières (cf.~\sga{6}{vii}{1.4}),
Ofer Gabber définit une classe fondamentale
généralisée pour une immersion fermée $Y\to X$ définie par un Idéal de
type fini $\II$. Cette construction n'est plus limitée aux
immersions régulières et est compatible aux 
changements de bases arbitraires, mais elle 
dépend d'une donnée supplémentaire, à savoir celle d'un fibré vectoriel
sur $Y$ se surjectant sur le faisceau conormal $\Normal Y X=\II/\II^2$.

\subsubsection{Éclatement modifié}

Soit $Y\to X$ une immersion fermée entre $\ZZn$-schémas définie 
par un Idéal de type fini $\II$. On note $U$ l'ouvert complémentaire.
Soit $\mathscr E
\to \Normal Y X$ un épimorphisme de Modules sur $Y$ où 
$\mathscr E$ est un Module localement libre de rang fini.
On définit une $\mathscr O_X$-Algèbre graduée quasi-cohérente $\mathscr
A_\star$ par produit fibré de façon à 
avoir un carré cartésien de $\mathscr O_X$-Modules,
pour tout entier naturel $n$ :
\[
\xymatrix{
\ar[d]{\mathscr A}_n\ar[r] & \II^n \ar[d]\\
\SS^n\mathscr E\ar[r] & \II^n/\II^{n+1}}
\]
où l'algèbre symétrique $\SS^\star\mathscr E$ est prise sur le faisceau
d'anneaux $\mathscr O_Y=\mathscr O_X/\II$.

\begin{definition3}
On pose $\Ecl_{Y,\mathscr E}(X)=\Proj({\mathscr A}_\star)$ et on note
$\pi\colon \Ecl_{Y,\mathscr E}(X)\to X$ la projection.
\end{definition3}

\begin{remarque3}
Si $Y\to X$ est une immersion fermée
régulière et que $\mathscr E\to \Normal Y X$
est un isomorphisme, $\Ecl_{Y,\mathscr E}(X)$ s'identifie à l'éclaté de $Y$
dans $X$. C'est ce cas particulier que l'on généralise ici en vue d'obtenir
une construction compatible aux changements de base.
\end{remarque3}

\begin{proposition3}\label{16-proposition-proprietes-eclatement-modifie}
L'Algèbre ${\mathscr A}_0$ est isomorphe à ${\mathscr O}_X$, les Modules
${\mathscr A}_n$ sont de type fini, l'Algèbre graduée
${\mathscr A}_\star$ est engendrée par ${\mathscr A}_1$ et on a un 
isomorphisme canonique de $\mathscr O_Y$-Algèbres graduées 
$\mathscr A_\star \otimes_{\mathscr O_X}{\mathscr O_X/\II}\iso
\SS^\star \mathscr E$.
\end{proposition3}

L'assertion concernant $\mathscr A_0$ est tautologique.
Soit $n$ un entier naturel. Comme $\II^n \to \II^n/\II^{n+1}$ est un
épimorphisme, la projection $\mathscr A_n\to \SS^n\mathscr E$ est aussi un
épimorphisme et si on note $\mathscr K_n$ son noyau, on a un isomorphisme
$\mathscr K_n\iso \II^{n+1}$. Par dévissage, il en résulte que $\mathscr
A_n$ est un Module de type fini.

Puisque $\SS^n\mathscr E$ est un $\mathscr O_X/\II$-Module, $\mathscr K_n$
contient $\II\cdot\mathscr A_n$. Comme $\mathscr E\to \II/\II^2$ est un
épimorphisme, le morphisme $\SS^n\mathscr E\to \II^n/\II^{n+1}$ est aussi
un épimorphisme, ce qui implique que la projection $\mathscr A_n\to \II^n$
est un épimorphisme. Par conséquent, l'inclusion
$\II\cdot\mathscr A_n\to\II\cdot
\II^n=\II^{n+1}\simeq \mathscr K_n$ est un isomorphisme : $\mathscr
K_n=\II\cdot \mathscr A_n$. Ceci permet d'obtenir l'isomorphisme
$\mathscr A_\star\otimes_{\mathscr O_X} {\mathscr O_X/\II}\iso
\SS^\star \mathscr E$.

Pour montrer que le morphisme évident $\mathscr A_1^{\otimes n}\to
\mathscr A_n$ de Modules
est un épimorphisme, il suffit, d'après le lemme de
Nakayama, de le tester après passage aux corps résiduels de $X$.
Au-dessus de l'ouvert $U$, c'est évident ; au-dessus de $Y$, cela résulte
de l'isomorphisme 
$\mathscr A_\star\otimes_{\mathscr O_X} {\mathscr O_X/\II}\iso 
\SS^\star \mathscr E$.

\begin{corollaire3}\label{16-corollaire-structure-eclatement-modifie}
Le morphisme $\pi\colon \Ecl_{Y,\mathscr E}(X)\to X$ est projectif et on
dispose d'isomorphismes canoniques $\pi^{-1}(U)\iso U$ et
$\pi^{-1}(Y)\simeq \PP(\mathscr E)$.
\end{corollaire3}

L'isomorphisme au-dessus de $U$ est évident. Compte tenu de \ega{ii}{3.5.3},
celui décrivant $\pi^{-1}(Y)$
se déduit de l'isomorphisme de $\mathscr O_Y$-Algèbres graduées $\mathscr
A_\star \otimes_{\mathscr O_X}\mathscr O_Y\iso \SS^\star\mathscr E$.

\begin{proposition3}\label{16-proposition-cb-eclatements-modifies}
Soit $p\colon X'\to X$ un morphisme.
On pose $Y'=Y\times_X X'$ et $\mathscr E'=p^\star \mathscr E$. On dispose
d'un épimorphisme évident $\mathscr E'\to\Normal {Y'} {X'}$. Le
morphisme canonique
\[\Ecl_{Y',\mathscr E'}(X')\to \Ecl_{Y,\mathscr E}(X)\times_X X'\]
est une nil-immersion.
\end{proposition3}

Notons $\mathscr A'_\star$ la $\mathscr O_{X'}$-Algèbre 
graduée quasi-cohérente donnant
naissance à $\Ecl_{Y',\mathscr E'}(X')$. On dispose d'un morphisme évident 
$p^\star \mathscr A_\star\to\mathscr A'_\star$ de $\mathscr
O_{X'}$-Algèbres graduées quasi-cohérentes. Pour tout entier, le morphisme
$p^\star \mathscr A_n\to \mathscr A'_n$ est un morphisme entre $\mathscr
O_{X'}$-Modules de type fini ; pour montrer qu'il s'agit d'un épimorphisme,
d'après le lemme de Nakayama, il suffit de vérifier que ce morphisme
induit un isomorphisme d'une part au-dessus de $U'=X'-Y'$ (c'est
évident) et d'autre part modulo l'idéal $\II'$ définissant $Y'$ dans $X'$
(cela résulte de la description donnée dans la 
proposition~\ref{16-proposition-proprietes-eclatement-modifie}). Le morphisme 
\[\Ecl_{Y',\mathscr E'}(X')\to \Ecl_{Y,\mathscr E}(X)\times_X X'\]
s'identifie au $X'$-morphisme évident
$\Proj(\mathscr A'_\star)\to \Proj (p^\star \mathscr A_\star)$
(\ega{ii}{3.5.3}) ; d'après ce
qui précède, il s'agit d'une immersion fermée. Le fait que ce morphisme
induise un isomorphisme au-dessus de $p^{-1}(U)$ et de $p^{-1}(Y)$
permet d'en déduire aussitôt que le morphisme induit au niveau des schémas
réduits associés
\[\Ecl_{Y',\mathscr E'}(X')_{\red}\to (\Ecl_{Y,\mathscr E}(X)\times_X
X')_\red\]
est un isomorphisme.

\subsubsection{Définition des classes}

On se donne toujours une immersion fermée $i\colon Y\to X$ définie par un
Idéal $\II$ de type fini. On note $j\colon U\to X$ l'inclusion de l'ouvert
complémentaire ($\II$ étant de type fini, $j$ est un morphisme de 
type fini). On suppose donné un épimorphisme de $\mathscr O_Y$-Modules
$\mathscr E \to \Normal Y X$
avec $\mathscr E$ est localement libre de rang
fini. On note $\pi\colon \Ecl_{Y,\mathscr E}(X)\to X$ la projection de
l'éclatement modifié, $j'\colon U\to \Ecl_{Y,\mathscr E}(X)$ l'immersion
ouverte évidente et $i'\colon \PP(\mathscr E)\to \Ecl_{Y,\mathscr E}(X)$
l'immersion fermée donnée par le 
corollaire~\ref{16-corollaire-structure-eclatement-modifie}. On note $r$ le
rang du fibré vectoriel $\mathscr E$ que l'on suppose de rang constant pour
simplifier et on suppose $r>0$.

On a ainsi le diagramme suivant de schémas, où les carrés sont cartésiens :

\[
\xymatrix{
\ar[d]^{\pi'}\PP(\mathscr E)\ar[r]^{i'} & 
\Ecl_{Y,\mathscr E}(X)\ar[d]^\pi &\ar[l]_{j'} U\ar@{=}[d]\\
Y \ar[r]^i & X & \ar[l]_{i'} U}
\]

\begin{proposition3}
\label{16-proposition-cohomologie-eclatement}
Le morphisme évident $\Lambda\to \R\pi'_\star\Lambda$ dans $\detaleplus Y$
est un monomorphisme scindé : la formule du fibré projectif identifie son
conoyau à 
\[\bigoplus_{k=1}^{r-1} \Lambda(-k)\cro{-2k}\;\text{.}\]
Les morphismes évidents définissent un triangle distingué :
\[\Lambda\to \R\pi_\star\Lambda \to
i_\star\coker(\Lambda\to \R\pi'_\star \Lambda)\vers 0 \Lambda\cro 1\]
dans $\detaleplus X$. On peut le récrire sous la forme
\[\Lambda\to \R\pi_\star\Lambda \vers \rho
\bigoplus_{k=1}^{r-1}i_\star\Lambda(-k)\cro{-2k}\vers 0 \Lambda\cro
1\;\text{,}\]
le morphisme $\rho$ admettant une section canonique
donnée par les éléments
$c_1(\OO(1),1_U)^k$ de $\HH^{2k}_{\PP(\mathscr E)}(\Ecl_{Y,\mathscr
E}(X),\Lambda(k))$, identifiés à des morphismes 
$i_\star \Lambda(-k)\cro {-2k}\to \R\pi_\star \Lambda$ dans $\detaleplus X$.
\end{proposition3}

On note $L$ une résolution injective du faisceau constant $\Lambda$ vu
comme faisceau de $\Lambda$-modules sur le grand site étale
des schémas de type fini sur $X$. Pour tout morphisme de type fini $W\vers
p X$, on note $L_{|W}$ le complexe de faisceau de $\Lambda$-modules sur
$W_\et$ induit par $L$ ; on peut le voir comme un objet de $\detaleplus W$
isomorphe à $\Lambda$.

\begin{lemme3}
Le carré commutatif évident de complexes de faisceaux sur $X_{\et}$
est homotopiquement bicartésien :
\[
\xymatrix{
\ar[d]L_{|X} \ar[r] & i_\star L_{|Y} \ar[d]\\
\pi_\star L_{|\Ecl_{Y,\mathscr E}(X)}\ar[r] & \pi_\star i'_\star
L_{|\PP(\mathscr E)}}
\]
(ceci signifie par exemple
que le complexe simple associé à ce diagramme, identifié à
un complexe $3$-uple, est acyclique).
\end{lemme3}

Les complexes simples associés aux complexes doubles 
\[j_!L_{|U}\to L_{|X} \to i_\star L_{|Y}\]
et 
\[j'_!L_{|U}\to L_{|\Ecl_{Y,\mathscr E}(X)} \to i'_\star L_{|\PP(\mathscr E)}\]
de faisceaux sur $X$ et $\Ecl_{Y,\mathscr E}(X)$ respectivement sont
acycliques.
Choisissons un foncteur de résolution \guil{flasque} additif
$r$ sur la catégorie des faisceaux de $\Lambda$-modules sur 
${\Ecl_{Y,\mathscr E}(X)}_\et$ et notons abusivement
$\R\pi_\star$ le foncteur (additif)
de la catégorie des
complexes (bornés inférieurement) de faisceaux de $\Lambda$-modules
sur $\Ecl_{Y,\mathscr E}(X)$ vers la catégorie des
complexes de faisceaux de $\Lambda$-modules
sur $X$ défini par la formule $\R\pi_\star K = \Tot(\pi_\star r K)$,
ce foncteur préserve les
quasi-isomorphismes et induit
le foncteur $\R\pi_\star\colon \detaleplus {\Ecl_{Y,\mathscr E}(X)}\to
\detaleplus X$ usuel.

On obtient ainsi un diagramme commutatif de complexes de faisceaux de
$\Lambda$-modules sur $X$ :
\[
\xymatrix{
j_!L_{|U} \ar[r]\ar[d] & L_{|X} \ar[r] \ar[d]& i_\star L_{|Y} \ar[d]\\
\R\pi_\star j'_! L_{|U} \ar[r] & \R\pi_\star
L_{|\Ecl_{Y,\mathscr E}(X)}\ar[r] & \R\pi_\star i'_\star
L_{|\PP(\mathscr E)}}
\]
Les lignes de ce diagramme constituent des complexes doubles dont les
complexes simples associés sont acycliques. D'après le théorème de
changement de base pour un morphisme propre, le morphisme 
$j_!L_{|U}\to \R\pi_\star j'_! L_{|U}$ est un quasi-isomorphisme.
On en déduit que le carré de droite est homotopiquement bicartésien, ce qui
permet de conclure.

\medskip

Revenons à la démonstration de la
proposition~\ref{16-proposition-cohomologie-eclatement},
la formule du fibré projectif pour $\PP(\mathscr E)$ implique que l'on a un
triangle distingué dans $\detaleplus X$ :
\[i_\star \Lambda\to \R\pi_\star i'_\star \Lambda
\to \bigoplus_{i=1}^{r-1}i_\star \Lambda(-i)\cro{-2i} \vers 0 i_\star
\Lambda\cro 1\;\text{.}\]
En considérant les colonnes du carré homotopiquement bicartésien donné par
le lemme, on peut conclure à l'existence d'un triangle distingué
\[
\xymatrix{
\Lambda \to \R\pi_\star \Lambda \to 
\displaystyle\bigoplus_{i=1}^{r-1} i_\star
\Lambda(-i)\cro{-2i}\to\Lambda\cro 1\;\text{.}}
\]
Ce triangle est scindé par les puissances de l'élément $c_1(\OO(1),1_U)$ ; le
morphisme de droite est donc nul, ce qui achève la démonstration de la
proposition.

\begin{corollaire3}\label{16-corollaire-caracterisation-image-restriction}
La suite suivante, dont les morphismes sont évidents, est exacte :
\[0\to \HH^{2r}_Y(X,\Lambda(r)) \to 
\HH^{2r}_{\PP(\mathscr E)}(\Ecl_{Y,\mathscr E}(X),\Lambda(r))\to
\coker(\HH^{2r}(Y,\Lambda(r))\to \HH^{2r}(\PP(\mathscr E),\Lambda(r)))\to 0\]
\end{corollaire3}

L'énoncé de ce corollaire serait bien évidemment juste en tout bidegré
$(p,q)$ au lieu de $(2r,r)$, mais nous n'utiliserons que ce cas
particulier.

\bigskip

On note $\Gys\colon \HH^p(\PP(\mathscr E),\Lambda(q))\to
\HH^{p+2}_{\PP(\mathscr E)}(\Ecl_{Y,\mathscr E}(X),\Lambda(q+1))$ le
morphisme de Gysin associé au pseudo-diviseur $(\OO(1),1_U)$ sur 
$\Ecl_{Y,\mathscr E}(X)$ et $\xi=c_1(\OO(1))\in \HH^2(\PP(\mathscr
E),\Lambda(1))$. Le lemme suivant est évident :

\begin{lemme3}\label{16-lemme-res-gys-egale-xi}
Le morphisme composé
\[\HH^p(\PP(\mathscr E),\Lambda(q))\vers {\Gys} \HH^{p+2}_{\PP(\mathscr
E)}(\Ecl_{Y,\mathscr E}(X),\Lambda(q+1))\to  \HH^{p+2}(\PP(\mathscr
E),\Lambda(q+1))\;\text{,}\]
où la flèche de droite est le morphisme de restriction, est la
multiplication par $\xi$.
\end{lemme3}

\begin{definition3}
On définit un élément $\Clf {i,\mathscr E}$ de $\HH^{2r-2}(\PP(\mathscr
E),\Lambda(r-1))$ par la formule :
\[\Clf {i,\mathscr E}=-\cro{\xi^{r-1}-c_1(\mathscr
E)\xi^{r-2}+\dots+ (-1)^{r-1}c_{r-1}(\mathscr E)}\;\text{.}\]
\end{definition3}

Le lemme suivant résulte aussitôt de la construction des classes de Chern :

\begin{lemme3}\label{16-lemme-xi-clf}
Dans $\HH^{2r}(\PP(\mathscr E),\Lambda(r))$, on a l'égalité
\[\xi\Clf{i,\mathscr E}=(-1)^rc_r({\mathscr E})\;\text{.}\]
\end{lemme3}

\begin{definition3}
Compte tenu du corollaire~\ref{16-corollaire-caracterisation-image-restriction}, 
les lemmes~\ref{16-lemme-res-gys-egale-xi}
et \ref{16-lemme-xi-clf} 
montrent que l'élément
$\Gys(\Clf{i,\mathscr E})\in \HH^{2r}_{\PP(\mathscr E)}(\Ecl_{Y,\mathscr
E}(X),\Lambda(r))$ provient par restriction d'un unique élément de
$\HH^{2r}_Y(X,\Lambda(r))$, noté $\Cl{i,\mathscr E}$.
\end{definition3}

\subsubsection{Propriétés des classes généralisées}

\begin{proposition3}
\label{16-proposition-classes-generalisees-cb}
La formation des classes généralisées $\Cl{i,\mathscr E}$ 
et $\Clf {i,\mathscr E}$ est compatible à tout changement de base $X'\to X$.
\end{proposition3}

Compte tenu de la 
proposition~\ref{16-proposition-cb-eclatements-modifies}, ceci résulte
aussitôt des définitions.

\begin{proposition3}
\label{16-proposition-classes-generalisees-changement-de-fibre}
Soit ${\mathscr E}'\to {\mathscr E}$ un épimorphisme
de Modules localement libres sur $Y$. Soit $\mathscr K$ le noyau de cet
épimorphisme. On suppose que ${\mathscr E}'$ est de rang constant $r'$.
On a alors la relation
\[\Cl {i,\mathscr E'}= (-1)^{r'-r}c_{r'-r}({\mathscr K})\cdot\Cl
{i,\mathscr E}\]
dans $\HH^{2r'}_Y(X,\Lambda(r'))$ où on a utilisé les accouplements
canoniques \[\HH^a(Y,\Lambda(b))\otimes 
\HH^{a'}_Y(X,\Lambda(b'))\to  \HH^{a+a'}_Y(X,\Lambda(b+b'))\;\text{.}\]
\end{proposition3}

On dispose d'une immersion fermée de 
$\Ecl_{Y,\mathscr E}(X)$ dans 
$\Ecl_{Y,{\mathscr E}'}(X)$, ce qui permet de considérer la composition
suivante de flèches de restriction :
\[\HH^{2r'}_Y(X,\Lambda(r'))\to \HH^{2r'}_{\PP(\mathscr E')}(\Ecl_{Y,\mathscr
E'}(X),\Lambda(r'))\to \HH^{2r'}_{\PP(\mathscr E)}(\Ecl_{Y,\mathscr
E}(X),\Lambda(r'))\;\text{.}\]
Cette composée étant injective, il s'agit de montrer que les images des
deux éléments considérés dans $\HH^{2r'}_{\PP(\mathscr E)}(\Ecl_{Y,\mathscr
E}(X),\Lambda(r'))$
sont égales, mais comme ces deux éléments sont naturellement définis
comme étant des images d'éléments de $\HH^{2r'-2}(\PP(\mathscr
E),\Lambda(r'-1))$ par le morphisme $\Gys$ associé au fibré en droites
$\OO(1)$ sur $\Ecl_{Y,\mathscr E}(X)$ canoniquement trivialisé
sur $X-Y$, on se ramène à montrer l'égalité
\[{\Clf{i,{\mathscr E}'}}_{|\PP(\mathscr E)}=
(-1)^{r'-r}c_{r'-r}({\mathscr K})\cdot\Clf{i,\mathscr
E}\]
dans $\HH^{2r'}(\PP(\mathscr E),\Lambda(r'))$. Ceci ne fait plus intervenir
que le schéma $Y$ et l'épimorphisme ${\mathscr E}'\to \mathscr E$ de fibrés
vectoriels sur $Y$. Il s'agit d'une identité \guil{universelle} dans cette
situation ; on peut appliquer le principe de scindage et faire une
récurrence sur la différence $r'-r$ pour se ramener au cas où $r'=r+1$,
c'est-à-dire que $\mathscr K$ est un fibré en droites.
Les relations de Cartan-Whitney entre les classes de Chern de
$\mathscr E$, ${\mathscr E}'$ et $\mathscr K$ donnent aussitôt l'égalité
\[{\Clf{i,\mathscr E'}}_{|\PP(\mathscr E)} = \cro{\xi \Clf{i,\mathscr
E}-(-1)^rc_r({\mathscr E})}-c_1(\mathscr K)\cdot\Clf{i,\mathscr E}\;\text{,}\]
ce qui donne bien la relation voulue puisque $\xi\Clf{i,\mathscr
E}=(-1)^rc_r({\mathscr E})$ (cf.~lemme~\ref{16-lemme-xi-clf}).

\subsection{Immersions régulières}

On rappelle que la notion d'immersion régulière est définie dans
\sga{6}{vii}{1.4}.

\begin{definition2}\label{16-definition-cl-i}
Soit $i\colon Y\to X$ une immersion régulière entre $\ZZn$-schémas.
On pose $i^?=i^\star(c)\cro{2c}\colon \detaleplus X\to \detaleplus Y$ où $c$ est la
codimension de $i$.
On définit un morphisme 
$\Cl i \colon \Lambda\to i^?\Lambda$ dans $\detaleplus Y$ de la façon suivante.
Quitte à décomposer $Y$ en réunion disjointe d'ouverts-fermés, on peut
supposer que la codimension $c$ de $i$ est constante.
Si $c=0$, $i$ est l'inclusion d'un ouvert, $\Cl i$ est
l'isomorphisme évident. Dans le cas où $c>0$,
choisissons un ouvert $U$ de $X$ dans lequel $Y$ est un sous-schéma fermé,
notons $i'\colon Y\to U$ cette immersion fermée. Le faisceau conormal
$\Normal Y X$ de $Y$ dans $X$ est un fibré vectoriel de rang $c$ sur $Y$
muni de l'épimorphisme tautologique $\Normal Y X\to\Normal Y X$ ; on peut
donc considérer la classe $\Cl {i'}=\Cl {i',\Normal Y X}\in
\HH^{2c}_Y(U,\Lambda(c))$, que l'on identifie à un morphisme $\Cl i\colon
\Lambda\to
i'^?\Lambda\simeq i^?\Lambda$ dans $\detaleplus Y$ ; il est évident que la
construction ne dépend pas de l'ouvert intermédiaire $U$.
\end{definition2}

Le théorème suivant généralise l'énoncé établi dans
\cite[proposition~1.2.1]{AZUMINO@Fujiwara} :

\begin{theoreme2}\label{16-theoreme-composition-cl-immersions-regulieres}
Si $Z\vers i Y$ et $Y \vers j X$ sont deux 
immersions régulières composables, le
diagramme suivant est commutatif dans $\detaleplus Z$ :
\[
\xymatrix{
\ar[rd]_{\Cl {j\circ i}}\Lambda\ar[r]^{\Cl i}  &
   i^?\Lambda \ar[d]^{i^?(\Cl j)} \\
&   i^?j^?\Lambda}
\]
\end{theoreme2}

On peut évidemment supposer que les immersions $i$ et $j$ sont des
immersions fermées et que les codimensions de $i$ et de $j$ sont constantes,
de valeurs respectives $m$ et $n$. Si $m=0$ ou $n=0$, c'est trivial ; on
suppose donc que $m>0$ et $n>0$.

\begin{lemme2}
On peut supposer que $n=1$ (\emph{i.e.} $j$ est de codimension $1$).
\end{lemme2}

On éclate $Y$ dans $X$ pour obtenir le diagramme suivant où les carrés sont
cartésiens :
\[
\xymatrix{
P'\ar[r]^{i'} \ar[d]^{p'}& P\ar[r]^{j'}\ar[d]^p & \Ecl_Y(X)\ar[d]^\pi \\
Z\ar[r]^i& Y\ar[r]^j & X}
\]

L'idée de la démonstration est d'utiliser une sorte de formule d'excès
d'intersection (cf.~\cite[theorem~6.3]{IT@Fulton} pour une formulation dans la
théorie de Chow) pour les immersions $j\circ i\colon Z\to X$ et $j\colon Y
\to X$ relativement au changement de base $\pi\colon \Ecl_Y(X)
\to X$ qui va faire
chuter la codimension de ces immersions fermées régulières.

On a des isomorphismes canoniques $P=\PP(\Normal Y X)$ et $P'=\PP({\Normal
Y X}_{|Z})$.
On vérifie facilement que $P\to \Ecl_Y(X)$ est une 
immersion fermée régulière de
codimension $1$. Par changement de base lisse, $P'\to P$ est une immersion
fermée régulière de codimension $m$. On suppose que $i'^?(\Cl {j'})\circ
\Cl {i'}=\Cl{j'\circ i'}$ et on veut montrer que $i^?(\Cl {j})\circ
\Cl {i}=\Cl{j\circ i}$. Les morphismes à comparer s'identifient à des
éléments de $\HH^{2(m+n)}_Z(X,\Lambda(m+n))$ (on fera ce type
d'identifications jusqu'à la fin de la démonstration). La
proposition~\ref{16-proposition-cohomologie-eclatement} implique que
l'application
\[\pi^\star\colon \HH^{2(m+n)}_Z(X,\Lambda(m+n))\to
\HH^{2(m+n)}_{P'}(\Ecl_Y(X),\Lambda(m+n))\]
est injective, il suffit donc de comparer les classes après application de
$\pi^\star$.

Considérons $\pi^\star(\Cl {j\circ i})\in
\HH^{2(m+n)}_{P'}(\Ecl_Y(X),\Lambda(m+n))$. La classe $\Cl {j\circ i}$ est la
classe généralisée $\Cl {j\circ i,\Normal Z X}$, la 
proposition~\ref{16-proposition-classes-generalisees-cb} implique l'égalité
\[\pi^\star \Cl {j\circ i}=\Cl {j'\circ i',\pi^\star\Normal Z X}\]
où est sous-entendu l'épimorphisme de fibrés vectoriels $p'^\star \Normal
Z X\to \Normal {P'} {\Ecl_Y(X)}$ dont on note $\mathscr E'$ le noyau, qui est un
$\mathscr O_{P'}$-Module localement libre de rang $n-1$. La
proposition~\ref{16-proposition-classes-generalisees-changement-de-fibre}
donne alors l'égalité
\[\pi^\star\Cl{j\circ i}=c_{n-1}({\mathscr E'}^\dual)\cdot \Cl {j'\circ i'}\]
où l'on a utilisé l'accouplement
\[\HH^{2(n-1)}(P',\Lambda(n-1))\times \HH^{2(m+1)}_{P'}(\Ecl_Y(X),\Lambda(m+1))
\to \HH^{2(m+n)}_{P'}(\Ecl_Y(X),\Lambda(m+n))\;\text{.}\]
La composition des classes admise provisoirement pour les immersions $j'$
et $i'$ donne l'égalité
\[\Cl {j'\circ i'}=\Cl {i'}\cdot \Cl{j'}\]
\emph{via} l'accouplement
\[\HH^{2m}_{P'}(P,\Lambda(m))\times \HH^2_P(\Ecl_Y(X),\Lambda(1))\to
\HH^{2(m+1)}_{P'}(\Ecl_Y(X),\Lambda(n+1))\;\text{.}\]
On a ainsi obtenu :
\[\pi^\star \Cl{j\circ i}=c_{n-1}({\mathscr E'}^\dual)\cdot \Cl {i'}\cdot \Cl
{j'}\;\text{.}\]
Notons $\mathscr E$ le noyau de l'épimorphisme $p^\star \Normal Y
X\to\Normal {P'} {\Ecl_Y(X)}$. Il vient aussitôt que dans le diagramme évident
de Modules sur $P'$ qui suit, les lignes et les colonnes sont exactes :
\[
\xymatrix{
& 0 \ar[d] & 0 \ar[d] &   & \\
0 \ar[r] & \ar[d]i'^\star\mathscr E\ar[r] & \ar[d]\mathscr E'\ar[r] & 0
\ar[d] & \\
0 \ar[r] & \ar[d]i'^\star p^\star \Normal Y X \ar[r] & \ar[d]
p'^\star \Normal Z X\ar[r]
&\ar[d] p'^\star \Normal Z Y \ar[r] & 0 \\
0 \ar[r] & \ar[d]i'^\star \Normal {P} {\Ecl_Y(X)} \ar[r] & \ar[d]
\Normal {P'} {\Ecl_Y(X)}\ar[r] & \ar[d] \Normal {P'} P \ar[r] & 0 \\
&0 &0 &0 &}
\]
En particulier, on obtient un isomorphisme canonique $i'^\star \mathscr
E\iso \mathscr E'$, d'où $i'^\star c_{n-1}({\mathscr E}^\dual)
=c_{n-1}({\mathscr E'}^\dual)\in \HH^{2(n-1)}(P',\Lambda(n-1))$. On en déduit :
\[\pi^\star \Cl{j\circ i}=c_{n-1}({\mathscr E'}^\dual)
\cdot \Cl {i'}\cdot \Cl{j'}=
\Cl {i'}\cdot c_{n-1}({\mathscr E}^\dual)\cdot \Cl {j'}\;\text{.}\]
On utilise implicitement dans ces notations l'associativité des structures
multiplicatives permettant par exemple de définir une application
\[\HH^{2m}_{P'}(P,\Lambda(m))\times \HH^{2(n-1)}(P,\Lambda(n-1))\times
\HH^2_P(\Ecl_Y(X),\Lambda(1))\to \HH^{2(m+n)}(\Ecl_Y(X),\Lambda(m+n))\]
sans qu'il y ait à s'inquiéter de l'ordre dans lequel les multiplications
sont faites.
Les propositions~\ref{16-proposition-classes-generalisees-changement-de-fibre}
et \ref{16-proposition-classes-generalisees-cb} impliquent les égalités
suivantes :
\[c_{n-1}({\mathscr E}^\dual)\cdot \Cl{j'}=\Cl{j',p^\star \Normal Y X}
=\pi^\star \Cl{j}\in \HH^{2n}_P(\Ecl_Y(X),\Lambda(n))\;\text{.}\]
Le morphisme $p$ étant lisse, on a aussitôt $\Cl {i'}=\pi^\star \Cl i$. On
a ainsi obtenu l'égalité voulue :
\[\pi^\star \Cl{j\circ i}=\pi^\star \Cl i\cdot \pi^\star \Cl j\;\text{,}\]
ce qui achève la démonstration du lemme.

\bigskip

On est ramené à établir le
théorème~\ref{16-theoreme-composition-cl-immersions-regulieres} dans le cas où
$j$ est de codimension $1$. On pose maintenant $P=\PP(\Normal Z X)$ et
$P'=\PP(\Normal Z Y)$. Le diagramme suivant récapitule la situation :

\[
\xymatrix{
P\ar[r] & \pi^{-1}(Y)\ar[r] & \Ecl_Z(X) \ar[dd]^\pi\\
P'\ar[u]\ar[r]\ar[d] & \ar[d]\Ecl_Z(Y)\ar[u] & \\
Z \ar[r]^i & Y \ar[r]^j & X}
\]

On veut établir l'égalité suivante dans $\HH^{2m+2}_Z(X,\Lambda(m+1))$ :
\[\Cl {j\circ i}=\Cl i\cdot \Cl j\;\text{.}\]
D'après la proposition~\ref{16-proposition-cohomologie-eclatement}, il suffit
de vérifier cette égalité dans $\HH^{2m+2}_{P}(\Ecl_Z(X),\Lambda(m+1))$ après
application de $\pi^\star$.

Par définition, la classe $\Cl{j\circ i}\in \HH^{2m+2}_Z(X,\Lambda(m+1))$ se
\guil{restreint} en un élément
\[\gamma=\Gys_{P\subset \Ecl_Z(X)}(\Clf{j\circ i})\]
dans $\HH^{2m+2}_P(\Ecl_Z(X),\Lambda(m+1))$ où $\Clf{j\circ i}\in
\HH^{2m}(P,\Lambda(m))$.

Notons $\II$ l'Idéal de $Y$ dans $X$, ${\II}_P$ celui de $P$
dans $\Ecl_Z(X)$ et $\tilde{\II}$ celui de $\Ecl_Z(Y)$ dans $\Ecl_Z(X)$.
On a un isomorphisme canonique de faisceaux inversibles sur $\Ecl_Z(X)$ :
\[\pi^\star\II\simeq {\II}_P\otimes \tilde{\mathscr
I}\;\text{.}\]
Cet isomorphisme est compatible aux trivialisations données sur
$\pi^{-1}(V)$ où $V=X-Y$. On obtient ainsi une égalité dans le 
groupe des classes d'équivalences de tels pseudo-diviseurs, ce qui 
permet de décomposer $c_1(\pi^\star\II, 1_{\pi^{-1}(V)})\in
\HH^2_{\pi^{-1}(V)}(\Ecl_Z(X),\Lambda(1))$
en une somme de deux composantes :
\[\pi^\star(c_1(\II,1_V))=c_1(\pi^\star\II, 1_{\pi^{-1}(V)})
= c_1({\II}_P,1_{\Ecl_Z(X)-P})+c_1(\tilde{\mathscr
I},1_{\Ecl_Z(X)-\Ecl_Z(Y)})\;\text{.}\]

On en déduit une décomposition
\[\pi^\star\Cl i\cdot\pi^\star \Cl j=\alpha+\beta\] dans
$\HH^{2m+2}_P(\Ecl_Z(X),\Lambda(m+1))$ où
\begin{eqnarray*}
\alpha &=&-\Gys_{P\subset \Ecl_Z(X)}({\Cl i}_{|P})\;\text{,} \\
\beta &=& -\Gys_{\Ecl_Z(Y)\subset \Ecl_Z(X)}(\Gys_{P'\subset \Ecl_Z(Y)}(\Clf i))\;\text{.}
\end{eqnarray*}
Le calcul de $\Cl k$ où $k$ est l'inclusion de l'intersection de diviseurs
de Cartier s'intersectant transversalement dans le schéma ambiant réalisé dans 
\cite[proposition~1.1.4]{AZUMINO@Fujiwara}
permet d'obtenir l'égalité d'opérateurs suivante :
\[\Gys_{\Ecl_Z(Y)\subset \Ecl_Z(X)}\circ \Gys_{P'\subset \Ecl_Z(Y)} = 
\Gys_{P\subset \Ecl_Z(X)}\circ \Gys_{P'\subset P}\;\text{.}\]
Notre but est d'établir l'égalité $\gamma=\alpha+\beta$. Les calculs
précédents permettent d'écrire chacun des éléments $\alpha$, $\beta$ et
$\gamma$ comme des images par le morphisme $\Gys_{P\subset \Ecl_Z(X)}$
de classes
$\tilde{\alpha}$, $\tilde{\beta}$ et $\tilde{\gamma}$ dans
$\HH^{2m}(P,\Lambda(m))$ :
\begin{eqnarray*}
\tilde{\alpha} &=& -{\Cl i}_{|P}\;\text{,}\\
\tilde{\beta} &=& -\Gys_{P'\subset P}(\Clf i)\;\text{,}\\
\tilde{\gamma} &=& \Clf {j\circ i}\;\text{.}
\end{eqnarray*}
On est ainsi ramené à établir l'égalité
$\tilde{\gamma}=\tilde{\alpha}+\tilde{\beta}$
dans $\HH^{2m}(P,\Lambda(m))$.

D'après les lemmes~\ref{16-lemme-res-gys-egale-xi} et \ref{16-lemme-xi-clf}, on a 
${\Cl i}_{|Z}=(-1)^{m}c_m(\Normal Z Y)$. On en déduit l'égalité
\[\tilde{\alpha}=(-1)^{m+1}c_m(\Normal Z Y)\;\text{.}\]
Pour calculer $\tilde{\beta}$, on observe que l'Idéal de $P'$ dans $P$
s'identifie au faisceau inversible $\mathscr K
\otimes_{{\mathscr O}_Z} \OO(-1)$
où $\mathscr K={\Normal Y X}_{|Z}$ est le noyau de l'épimorphisme
$\Normal Z X\to \Normal Z Y$. On en déduit 
\[\tilde{\beta}= (c_1(\mathscr K)-\xi)\cdot\cro{\xi^{m-1}-c_1(\Normal Z
Y)\xi^{m-2}+\dots + (-1)^{m-1}c_{m-1}(\Normal Z Y)}\;\text{.}\]
Par ailleurs, la définition de $\tilde{\gamma}$ donne l'égalité :
\[\tilde{\gamma} = -\cro{\xi^m-c_1(\Normal Z X)\xi^{m-1}+(-1)^mc_m(\Normal
Z X)}\;\text{.}\]
La formule de Cartan-Whitney appliquée à la suite exacte courte
\[0\to \mathscr K\to \Normal Z X\to\Normal Z Y\to 0\]
de fibrés vectoriels sur $Z$ permet d'obtenir aussitôt la relation voulue
$\tilde{\gamma}=\tilde{\alpha}+\tilde{\beta}$, ce qui achève la
démonstration du théorème.

\bigskip

La classe que l'on a définie est évidemment compatible avec celle de
\sgaquatreetdemi{Cycle}{2.2} :

\begin{proposition2}\label{16-proposition-compatibilite-sga-quatre-et-demi}
Soit $i\colon Y\to X$ une immersion régulière de codimension $c$
entre $\ZZn$-schémas.
Le morphisme de faisceaux $\Lambda\to {\mathscr H}^{2c}(i^!\Lambda(c))$
induit par le morphisme $\Cl i\colon \Lambda\to i^?\Lambda$ 
est donné par la classe $\clsga Y$ de \sgaquatreetdemi{Cycle}{2.2}.
\end{proposition2}

\subsection{Morphismes lisses}

Soit $p\colon X\to S$ un morphisme lisse compactifiable
de $\ZZn$-schémas de dimension
relative $d$. D'après
\sga{4}{xviii}{2.9}, on dispose d'un morphisme trace
\[\Tr_p\colon \R^{2d}p_!\Lambda(d)\to \Lambda\;\text{,}\]
que l'on peut réinterpréter sous la forme d'un morphisme
\[\R p_!\Lambda(d)\cro{2d}\to \Lambda\]
dans $\detaleplus S$ (en effet, d'après le théorème de changement de base pour un
morphisme propre et \sga{4}{x}{4.3}, les faisceaux $\R^ip_!\Lambda$
sont nuls pour $i>2d$). 

\begin{definition2}\label{16-definition-cl-p}
Soit $p\colon X\to S$ un morphisme lisse compactifiable
de $\ZZn$-schémas. Le morphisme 
$\Cl p\colon \Lambda\to p^?\Lambda$\;\footnote{On rappelle que l'on a posé
$p^?=p^!(-d)\cro{-2d}$ où $d$ est la dimension relative de $p$.}
dans $\detaleplus X$ est le morphisme déduit
par adjonction du morphisme $\R p_!\Lambda(d)\cro{2d}\to\Lambda$
défini ci-dessus.
\end{definition2}

D'après \sga{4}{xviii}{3.2.4}, ce morphisme $\Cl p$ est un isomorphisme :
c'est la dualité de Poincaré.

\begin{proposition2}\label{16-proposition-composition-cl-lisses}
Si $f\colon Z\to Y$ et $g\colon Y\to X$ sont des morphismes lisses
compactifiables composables, le diagramme suivant
est commutatif dans $\detaleplus Z$ :
\[
\xymatrix{
\ar[rd]_{\Cl {j\circ i}}\Lambda\ar[r]^{\Cl i}  &
   i^?\Lambda \ar[d]^{i^?(\Cl j)} \\
&   i^?j^?\Lambda}
\]
\end{proposition2}

Ceci est énoncé en \sga{4}{xviii}{3.2.4} et résulte de la
compatibilité des morphismes traces à la composition,
cf.~propriété (Var~3) dans \sga{4}{xviii}{2.9}.

\begin{remarque2}
Si cette théorie avait été à notre disposition, il eût peut-être été plus
commode d'utiliser ici la construction des foncteurs $f^!$ pour $f$
lissifiable mentionnée dans l'introduction de \sga{4}{xviii}{0.4}.
Dans le cadre axiomatique des \guil{foncteurs homotopiques stables}, ceci
est réalisé dans \cite{THESE@Ayoub}.
\end{remarque2}

\subsection{Morphismes d'intersection complète lissifiables}

\begin{definition2}\label{16-definition-intersection-complete}
Un morphisme d'intersection complète
est un morphisme 
$X\vers f S$ admettant localement une factorisation sous la forme
$X\vers i T\vers p S$ où $p$
est lisse et $i$ une immersion régulière (cf.~\sga{6}{vii}{1.4}).
On pose
$\dimrelvirt f = \dim p-\codim i$ : c'est la dimension relative virtuelle
de $f$ (cf.~\sga{6}{viii}{1.9}).
\end{definition2}

\begin{definition2}
On note $\Sch$ la catégorie dont les objets sont les $\ZZn$-schémas
quasi-compacts
admettant un faisceau inversible ample et dont les morphismes sont
les morphismes de type fini entre de tels schémas.
On note $\Schic$ la sous-catégorie de $\Sch$ ayant les mêmes objets mais
dont les morphismes sont les morphismes d'intersection complète.
\end{definition2}

Dans $\Sch$, tout morphisme $X\to Y$ peut se factoriser sous la forme
$X\vers i \PP^n_Y\vers \pi Y$ où $i$ est une immersion et $\pi$ la
projection canonique.
Tous les morphismes de $\Sch$ sont donc compactifiables, on
peut leur appliquer le formalisme des foncteurs $\R f_!$ et $f^!$.

Les morphismes de $\Schic$ admettent des factorisations globales dans
$\Schic$ sous la forme d'une immersion fermée régulière suivie d'un
morphisme lisse. 

\begin{definition2}
Pour tout morphisme $f\colon X\to Y$ dans $\Schic$, on peut
définir un foncteur
\[f^?\colon \detaleplus Y\to \detaleplus X\]
par la formule $f^?=f^!(-d)\cro{-2d}$
où $d=\dimrelvirt f$. 
\end{definition2}

Les foncteurs $f^?$ sont les foncteurs image inverse pour une structure de
catégorie fibrée convenable au-dessus de la catégorie $\Schic$ : on
utilisera implicitement les isomorphismes de transitivité $f^?g^?\simeq
(gf)^?$ associés à la composition de deux morphismes composables dans
$\Schic$.

\begin{definition2}\label{16-definition-cl-f}
Soit $f\colon X\to S$ un morphisme dans $\Schic$.
On suppose donnée une factorisation de $f$ dans $\Schic$ sous la forme
$X\vers i Y\vers p S$ où $i$ est une immersion régulière et $p$ un
morphisme lisse. On définit un morphisme
\[\Cl {p,i}\colon \Lambda\to f^?\Lambda\]
dans $\detaleplus X$ comme étant le morphisme composé
\[
\xymatrix{
\Lambda\ar[rd]_{\Cl {p,i}}\ar[r]^{\Cl i}&i^?\Lambda\ar[d]^{i^?(\Cl p)} \\
& f^?\Lambda}
\]
où $\Cl i$ est le morphisme de la définition~\ref{16-definition-cl-i}
et $\Cl p$ celui de la définition~\ref{16-definition-cl-p}.
\end{definition2}

\begin{theoreme2}\label{16-theoreme-independance-factorisation}
Soit $f\colon X\to S$ un morphisme dans $\Schic$.
Si $X\vers i Y\vers p S$ et $X\vers {i'}Y'\vers {p'}S$ sont
deux factorisations du type envisagé dans la
définition~\ref{16-definition-cl-f}, alors les deux morphismes suivants
dans la catégorie $\detaleplus X$ sont égaux :
\[\Cl {p,i}=\Cl{p',i'}\colon \Lambda\to f^?\Lambda\;\text{.}\]
\end{theoreme2}

La notation suivante s'avère assez commode pour cette démonstration :

\begin{definition2}
Si $f\colon Z\to Y$ et $g\colon Y\to X$ sont des morphismes composables
dans $\Schic$, $a\colon \Lambda\to g^?\Lambda$ et 
$b\colon \Lambda\to f^?\Lambda$
des morphismes dans $\detaleplus Y$ et $\detaleplus Z$ respectivement, on pose $a\bigstar
b=f^?(a)\circ b\colon \Lambda\to (g\circ f)^?\Lambda$.
\end{definition2}

Cette loi $\bigstar$ vérifiant une propriété d'associativité évidente, on
omettra les parenthèses.

Par définition, on a ainsi : $\Cl {p,i}=\Cl p\bigstar \Cl i$.
On veut vérifier l'égalité $\Cl p\bigstar\Cl i=\Cl {p'}\bigstar \Cl {i'}$.
Quitte à introduire le produit fibré de $Y$ et de $Y'$ au-dessus de $S$, on
peut supposer que \guil{$Y'$ coiffe $Y$}, à savoir qu'il existe un
morphisme lisse $q\colon Y'\to Y$ tel que $i=q\circ i'$ et $p'=p\circ q$ :
\[
\xymatrix{
& Y'\ar[rd]^{p'}\ar[d]^q & \\
X\ar[ru]^{i'}\ar[r]^i & Y\ar[r]^p &S}
\]
On a ainsi
\[\Cl{p',i'}=\Cl {p'}\bigstar \Cl{i'}=\Cl p\bigstar\Cl q\bigstar \Cl{i'}\;
\text{,}\]
la dernière égalité résultant de la
proposition~\ref{16-proposition-composition-cl-lisses}. On est ramené à montrer
l'égalité $\Cl i=\Cl q\bigstar\Cl {i'}$. Pour cela, on introduit le produit
fibré $X'$ de $X$ et $Y'$ au-dessus de $Y$ :
\[
\xymatrix{
\ar[d]^{q'}X'\ar[r]^j& Y'\ar[d]^q \\
X\ar@/^1pc/@{.>}[u]^s\ar[ru]_{i'}\ar[r]_i & Y}
\]
Le morphisme $i'$ donne naissance à la section $s$ de la projection
$q'\colon X'\to X$. Le morphisme $q$ étant lisse, l'immersion $j\colon
X'\to Y'$ est régulière. Admettons provisoirement les
égalités suivantes :
\[\Cl {q'}\bigstar \Cl s=\Id_\Lambda\;\text{,}\qquad 
\Cl q\bigstar \Cl j=\Cl i\bigstar \Cl {q'}\;\text{.}\]
Il vient :
\begin{eqnarray*}
\Cl i&=&\Cl i\bigstar \Cl {q'}\bigstar \Cl s\\
&=& \Cl q\bigstar \Cl j\bigstar \Cl s\;\text{.}
\end{eqnarray*}
On utilise alors la composition des morphismes de Gysin associés aux
immersions régulières 
(cf.~théorème~\ref{16-theoreme-composition-cl-immersions-regulieres}).
Celle-ci donne l'égalité $\Cl j\bigstar \Cl s=\Cl {i'}$ qui permet de
conclure que $\Cl i=\Cl q\bigstar \Cl {i'}$. Les deux lemmes qui suivent
permettent d'obtenir les deux égalités admises ci-dessus :

\begin{lemme2}\label{16-lemme-cb-lisse}
Soit un diagramme cartésien dans $\Sch$ :
\[
\xymatrix{
X'\ar[d]^{q'}\ar[r]^j & Y'\ar[d]^q \\
X\ar[r]^i & Y}
\]
On suppose que $q$ est lisse et que $i$ est une immersion régulière (donc
$j$ aussi). Alors on a l'égalité
\[\Cl q\bigstar \Cl j=\Cl i\bigstar\Cl {q'}\;\text{.}\]
\end{lemme2}

On peut supposer que $i$ est une immersion fermée.
On identifie $\Cl i$ (resp. $\Cl {j}$) à une classe dans 
$\HH^{2d}_X(Y,\Lambda(d))$ (resp. $\HH^{2c}_{X'}(Y',\Lambda(c))$) où $c$ est
la codimension de l'immersion régulière $i$. D'après la
proposition~\ref{16-proposition-classes-generalisees-cb}, on a
$q^\star(\Cl i)=\Cl j$. Bien que la vérification soit abracadabrante, 
la compatibilité du morphisme trace aux changements
de base (propriété (Var~2) de \sga{4}{xviii}{2.9}) permet de conclure.

\medskip

\begin{lemme2}
Soit $p\colon X\to S$ un morphisme lisse dans $\Sch$
admettant une section $s\colon
S\to X$ (qui est une immersion régulière). Alors, $\Cl p\bigstar \Cl
s=\Id_\Lambda$ dans $\detaleplus S$.
\end{lemme2}

Les endomorphismes de $\Lambda$ dans $\detaleplus S$ sont donnés par des sections
du faisceau $\Lambda$ dans $S$, il suffit de vérifier que les nombres
obtenus en passant aux points génériques de $S$ sont égaux à $1$. Comme on
peut supposer que $S$ est réduit et que la construction est compatible avec
le passage aux points génériques, on peut supposer que $S$ est le
spectre d'un corps $k$. Notons $x$ l'image de $\SP(k)$ dans $X$. Quitte à
remplacer $X$ par un voisinage ouvert, on peut supposer qu'il existe un
morphisme étale $\pi\colon X\to \Aff^d_k$ identifiant $x$ à l'image inverse
de l'origine dans $\Aff^d_k$. En utilisant l'isomorphisme évident
$\HH^{2d}_{(0,\dots,0)}(\Aff^d_k,\Lambda(d))\iso
\HH^{2d}_{x}(X,\Lambda(d))$, on se ramène au lemme suivant :

\begin{lemme2}
Pour tout entier naturel $d$ et tout schéma $S\in \Sch$, si on note $p\colon
\Aff^d_S\to S$ la projection et $s\colon S\to \Aff^d_S$ l'inclusion de
l'origine, on a l'égalité
\[\Cl p\bigstar \Cl s=\Id_\Lambda \]
dans $\detaleplus S$.
\end{lemme2}

L'énoncé est évident pour $d=0$. Une récurrence évidente s'appuyant sur le
théorème~\ref{16-theoreme-composition-cl-immersions-regulieres} et la
proposition~\ref{16-proposition-composition-cl-lisses} permet de se ramener au
cas où $d=1$, et comme précédemment, on peut supposer que $S=\SP(k)$ où
$k$ est un corps que l'on peut supposer séparablement clos. On se ramène
finalement au lemme suivant :

\begin{lemme2}\label{16-lemme-compatibilite-trace-c-un}
Pour tout corps séparablement clos $k$, si on note $p\colon \PP^1_k\to
\SP(k)$ la projection et $s\colon \SP(k)\to \PP^1_k$ l'inclusion de
$0$, on a l'égalité
\[\Cl p\bigstar \Cl s=\Id_\Lambda \]
dans $\detaleplus {\SP(k)}$.
\end{lemme2}

L'idéal de l'immersion fermée $s$ s'identifie au faisceau inversible
$\OO(-1)$. Par définition, l'image ${\Cl s}_{|\PP^1_k}$
de $\Cl s$ dans $\HH^2(\PP^1_k,\Lambda(1))$
est $c_1(\OO(-1))\cdot \Clf s$. Mais $\Clf s=-1$. Ainsi,
${\Cl s}_{|\PP^1_k}=c_1(\OO(1))$. Le degré du fibré en droites
$\OO(1)$ étant $1$, on peut conclure en utilisant la commutativité
du diagramme suivant (cf.~\sga{4}{xviii}{1.1.6}) :
\[
\xymatrix{
\ar[rd]_{\deg} \Pic(\PP^1_k) \ar[r]^-{c_1} &
\HH^2(\PP^1_k,\Lambda(1)) \ar[d]^{\Tr_p}_{\sim} \\
& \Lambda}
\]

\bigskip

\begin{definition2}\label{16-definition-cl-f-lissifiable}
Soit $f\colon X\to S$ un morphisme dans $\Schic$.
On note $\Cl f\colon
\Lambda\to f^?\Lambda$ le morphisme $\Cl {p,i}$ dans $\detaleplus X$ défini à
partir d'une factorisation de $f$ dans $\Schic$
sous la forme $f=p\circ i$ avec $i$ une
immersion régulière et $p$ un morphisme lisse. D'après le 
théorème~\ref{16-theoreme-independance-factorisation}, cette définition est
indépendante de la factorisation.
\end{definition2}

\begin{theoreme2}\label{16-theoreme-cl-f-transitivite}
Si $X\vers f Y$ et $Y\vers g Z$ sont des morphismes composables dans
$\Schic$, le diagramme suivant est commutatif dans $\detaleplus X$.
\[
\xymatrix{
\Lambda\ar[rd]_{\Cl {g\circ f}}\ar[r]^{\Cl f}&f^?\Lambda\ar[d]^{f^?(\Cl g)} \\
& (g\circ f)^?\Lambda }
\]
\end{theoreme2}

Paraphrasant \sga{6}{viii}{2.6}, on choisit une factorisation $Y\vers j
V'\vers p' Z$ dans $\Schic$ avec $j$ une immersion régulière et $p'$ lisse,
et une immersion régulière $X\vers i \PP^n_Y$, de façon à obtenir le
diagramme suivant :
\[
\xymatrix{
\ar[rd]_fX\ar[r]^i &\ar[d]^p \PP^n_Y\ar[r]^{j'} & \PP^n_{V'}\ar[d]^{p''}\\
& \ar[rd]_g Y\ar[r]^j & V' \ar[d]^{p'}\\
& & Z}
\]
En utilisant le
théorème~\ref{16-theoreme-composition-cl-immersions-regulieres} et la
proposition~\ref{16-proposition-composition-cl-lisses}, on obtient
\[\Cl{g\circ f} = (\Cl{p'}\bigstar\Cl{p''})\bigstar (\Cl {j'}\bigstar \Cl
i)
\;\text{.}\]
Le lemme~\ref{16-lemme-cb-lisse} donne l'égalité :
\[\Cl{p''}\bigstar \Cl {j'}=\Cl j\bigstar\Cl p\;\text{,}\]
ce qui permet d'obtenir :
\[\Cl {g\circ f}=(\Cl{p'}\bigstar \Cl j)\bigstar(\Cl p\bigstar \Cl i)
\;\text{,}\]
où l'on reconnaît l'égalité $\Cl {g\circ f}=\Cl g\bigstar \Cl f$.

\bigskip

\begin{proposition2}
Soit $f\colon X\to S$ un morphisme dans $\Schic$. On suppose que $f$ est
plat de dimension relative $d$. Alors le morphisme $\Cl f\colon \Lambda
\to f^?\Lambda$ correspond par adjonction au morphisme
$\R f_!\Lambda(d)[2d]\to \Lambda$ donné par le morphisme trace $\Tr_f\colon
\R^{2d}f_!\Lambda(d)\to \Lambda$.
\end{proposition2}

Compte tenu de la
proposition~\ref{16-proposition-compatibilite-sga-quatre-et-demi}, cela
résulte de \sgaquatreetdemi{Cycle}{2.3.8~(i)}.

\begin{remarque2}
Si $f\colon X\to Y$ est un morphisme propre
dans $\Schic$ de dimension relative
virtuelle $d$, le morphisme $\Cl f$ permet de définir, pour tout $K\in
\detaleplus Y$, un morphisme $f_\star\colon \HH^{p}(X,f^\star K)\to
\HH^{p-2d}(Y,K(-d))$, compatible à la composition. On peut aussi en définir une
version à supports $f_\star\colon \HH^{p}_Z(X,f^\star K)\to
\HH^{p-2d}_{Z'}(Y,K(-d))$ dès que $Z$ et $Z'$ sont des fermés de
$X$ et $Y$ respectivement tels que $f(Z)\subset Z'$.
\end{remarque2}

\section{Théorème de pureté}
\label{16-section-purete}
\subsection{Énoncés}
L'objectif de cette section est de donner une démonstration
du théorème suivant :

\begin{theoreme2}\label{16-theoreme-purete-absolue}
Soit $X$ un $\ZZn$-schéma régulier. Soit $Y$ un sous-schéma (fermé)
de $X$ qui est aussi
régulier. On note $i\colon  Y\to X$ l'immersion, et $c$ sa codimension.
Alors, le morphisme de
Gysin $\Cl i\colon
\Lambda\to i^?\Lambda=i^!\Lambda(c)\cro{2c}$ est un isomorphisme dans
$\detaleplus Y$.
\end{theoreme2}

\begin{corollaire2}\label{16-corollaire-purete-absolue}
Soit $f\colon X\to S$ un morphisme de type fini entre $\ZZn$-schémas
réguliers. On suppose que $X$ et $S$ admettent un faisceau ample. Alors, le
morphisme de Gysin $\Cl f\colon \Lambda(d)[2d]\to f^!\Lambda$ est un
isomorphisme dans $\detaleplus X$, où $d$ désigne la dimension relative virtuelle
de $f$.
\end{corollaire2}

\begin{corollaire2}\label{16-corollaire-purete-absolue-dcn}
Soit $X$ un $\ZZn$-schéma régulier. Soit $D$ un diviseur à croisements
normaux dans $X$. On note $j\colon X-D\to X$ l'inclusion de son 
complémentaire. Alors, $\R j_\star \Lambda$ appartient à $\detalebctf X$.
Plus précisément, si $D=D_1+\dots+D_n$ est un diviseur à croisements
normaux strict, alors $\R^1 j_\star \Lambda$ s'identifie à $\bigoplus_{1\leq
i\leq n} \Lambda_{D_i}(-1)$ et $\R^\star j_\star \Lambda$ est l'algèbre
extérieure sur $R^1 j_\star \Lambda$.
\end{corollaire2}

Ce corollaire mérite une démonstration. Pour la première assertion, on peut
travailler localement pour la topologie étale sur $X$ ; il suffit donc
d'établir la deuxième assertion. On suppose que $D=D_1+\dots+D_n$ est un
diviseur à croisements normaux strict. On note $j_i\colon X-D_i\to X$
l'inclusion du complémentaire de $D_i$ pour tout $i$. Nous allons montrer
que le morphisme de Künneth
\[\R j_{1\star}\Lambda \otimes^{\L}\dots \otimes^{\L} \R
j_{n\star}\Lambda\to \R j_\star\Lambda\]
est un isomorphisme dans $\detale X$, ce qui impliquera le résultat vu que
les faisceaux $\R^q j_{i\star}\Lambda$ sont connus par pureté ($\Lambda$
pour $q=0$, $\Lambda_{D_i}(-1)$ pour $q=1$ et $0$ sinon) et qu'ils sont
plats.

On procède par récurrence sur $n$. Les cas $n=0$ et $n=1$ sont évidents. On
suppose $n\geq 2$, on pose $D'=D_2+\dots+D_n$ et on fait l'hypothèse que le
résultat est connu pour $D'$. Il s'agit donc de monter que si on note
$j'\colon X-D'\to X$ l'inclusion du complémentaire de $D'$, alors le
morphisme de Künneth
\[\R j_{1\star}\Lambda\otimes^{\L} \R j'_\star \Lambda\to \R j_\star\Lambda\]
est un isomorphisme. Autrement dit, le morphisme canonique
\[\R\SheafHom(\Lambda_{X-D_1},\Lambda)\otimes^{\L}
\R\SheafHom(\Lambda_{X-D'},\Lambda)\to \R\SheafHom(\Lambda_{X-D_1}\otimes
\Lambda_{X-D'},\Lambda)\]
est un isomorphisme dans $\detale X$. À $K$ (resp. $L$) fixé dans $\detale
X$, la famille des $L$ (resp. $K$) tels que le morphisme
\[\R\SheafHom(K,\Lambda)\otimes^{\L}
\R\SheafHom(L,\Lambda)\to \R\SheafHom(K\otimes^{\L} L,\Lambda)\]
soit un isomorphisme, propriété que nous appellerons (Kü), est une
sous-catégorie triangulée de $\detale X$.

Pour $K=\Lambda$ ou $L=\Lambda$, la condition (Kü) est évidemment vérifiée,
la montrer pour $(\Lambda_{X-D_1},\Lambda_{X-D'})$ revient donc, par
dévissage, à la montrer pour $(\Lambda_{D_1},\Lambda_{X-D'})$ ou encore
pour $(\Lambda_{D_1},\Lambda_{D'})$. Il résulte aussitôt du théorème de
pureté et des compatibilités obtenues que si $Y$ et $Z$ sont deux
sous-schémas fermés réguliers de $X$ s'intersectant transversalement
(\emph{i.e.} $Y\cap Z$ est régulier de codimension la somme des
codimensions de $Y$ et de $Z$), alors $(\Lambda_Y,\Lambda_Z)$ vérifie (Kü).
En particulier, $(\Lambda_{D_1},\Lambda_{D_i})$ vérifie (Kü) pour $i\geq 2$
et plus généralement, pour tout sous-ensemble non vide $I$ de
$\{2,\dots,n\}$, $(\Lambda_{D_1},\Lambda_{D_I})$ vérifie (Kü) où $D_I$ est
l'intersection des $D_i$ pour $i\in I$. En utilisant la suite exacte
standard
\[0\to \Lambda_{D'}\to \bigoplus_{2\leq i\leq n}\Lambda_{D_i}\to
\bigoplus_{2\leq i<j\leq n} \Lambda_{D_{ij}}\to \dots\;\text{,}\]
on en déduit par dévissage la condition (Kü) pour
$(\Lambda_{D_1},\Lambda_{D'})$, ce qu'il fallait démontrer.

\begin{definition2}
Un couple régulier est un couple $(X,Y)$ où $X$ est un $\ZZn$-schéma régulier
et $Y$ un sous-schéma fermé de $X$ qui est régulier. On dit que $(X,Y)$ est
pur si la conclusion du théorème~\ref{16-theoreme-purete-absolue} est vraie
pour l'inclusion de $Y$ dans $X$. Si $\overline y\to Y$ est un point
géométrique de $Y$, on dira que $(X,Y)$ est pur en $\overline y$ si
$(\Cl i)_{\overline y}$ est un isomorphisme dans $\detaleplus {\overline y}$.
\end{definition2}

Le théorème~\ref{16-theoreme-purete-absolue} peut ainsi se reformuler en
disant que tout couple régulier est pur. Dans la
sous-section~\ref{16-subsection-purete-ponctuelle} sera introduite la notion de
pureté ponctuelle qui consiste à étudier les couples réguliers de la forme
$(X,x)$ où $X$ est un schéma local régulier de
point fermé $x$. Pour
démontrer le théorème de pureté, il suffira de savoir que les couples
réguliers de cette forme sont purs. Dans la
sous-section~\ref{16-subsection-changement-de-coefficients}, on se ramènera
au cas où l'anneau de coefficients $\Lambda$ est $\ZZ/\ell\ZZ$ avec $\ell$
un nombre premier inversible sur les schémas réguliers considérés. Dans
la sous-section~\ref{16-subsection-diviseurs-reguliers}, on établira quelques
propriétés utiles concernant la pureté des couples réguliers donnés par des
diviseurs. 
Comme dans la démonstration de \cite{AZUMINO@Fujiwara},
la démonstration de la pureté ponctuelle pour des schémas réguliers
arbitraires se ramènera à celle des
schémas réguliers qui sont de type fini sur
un trait $S$ (d'inégale caractéristique). Dans la
sous-section~\ref{16-subsection-schemas-sur-un-trait}, on obtiendra des
conditions suffisantes pour montrer que des schémas réguliers de type fini
sur $S$ sont ponctuellement purs. La
sous-section~\ref{16-subsection-geometrie-logarithmique} donnera les énoncés
de géométrie logarithmique permettant d'établir que si $(X,M)$ est un
log-schéma log-lisse sur un trait (muni de sa log-structure canonique)
et que le schéma $X$ est régulier, alors $X$ est ponctuellement pur.
La démonstration du théorème~\ref{16-theoreme-purete-absolue}
sera donnée dans la
sous-section~\ref{16-subsection-demonstration-theoreme-purete}. Elle utilisera
les résultats des sous-sections précédentes ainsi que trois théorèmes de
résolution des singularités que l'on peut résumer ainsi :
\begin{itemize}
\item utilisation d'altérations pour obtenir un schéma à réduction
semi-stable à partir d'un schéma (normal) sur $S$
(cf.~\cite[proposition~4.4.1]{BS@Vidal}) ;
\item résolution des singularités d'une action modérée
d'un groupe fini sur un log-schéma log-régulier de façon à
obtenir une action très modérée (cf. exposé VI) ;
\item résolution des log-singularités des log-schémas log-réguliers
(cf.~exposé VI).
\end{itemize}

\subsection{Pureté ponctuelle}
\label{16-subsection-purete-ponctuelle}

\begin{definition2}
Soit $X$ un $\ZZn$-schéma local régulier.
On dit que $X$ est ponctuellement pur en
son point fermé $x$ si le morphisme $\Cl i\colon \Lambda\to i^?\Lambda$
est un isomorphisme dans $\detaleplus x$ où $i\colon x\to X$ est
l'inclusion du point fermé de $X$.
\end{definition2}

Un schéma local régulier est ponctuellement pur en son point fermé si et
seulement si son hensélisé (resp. son hensélisé strict) l'est.

\begin{definition2}
Soit $X$ un $\ZZn$-schéma. Si $x\in X$, on dit que $X$ est ponctuellement
pur au point $x$ si le localisé de $X$ en $x$ est ponctuellement pur en son
point fermé. On dit que $X$ est ponctuellement pur s'il l'est en tous ses
points.
\end{definition2}

La proposition suivante est \cite[proposition~2.2.4]{AZUMINO@Fujiwara}. La
démonstration de cet article semble compliquée puisqu'elle passe par des
résultats plus fins que ceux dont nous avons besoin. On en redonne donc
une démonstration plus courte.

\begin{proposition2}\label{16-proposition-deux-sur-trois}
Soit $i\colon Y\to X$ une immersion fermée entre schémas réguliers. Le
nombre de conditions satisfaites parmi les trois suivantes ne peut pas être
deux :
\begin{enumerate}
\item[(a)] Le couple régulier $(X,Y)$ est pur ;
\item[(b)] Le schéma $Y$ est ponctuellement pur ;
\item[(c)] Le schéma $X$
est ponctuellement pur aux points situés dans l'image de $Y$.
\end{enumerate}
\end{proposition2}

Soit $y\in Y$, notons $V(y)$ le localisé de $Y$ en $y$ et $V(x)$ celui de
l'image $x$ de $y$ dans $X$. On a un diagramme de schémas :
\[
\xymatrix{
y\ar[rd]_{i_x}\ar[r]^{i_y} & V(y)\ar[d]^{i'} \\
& V(x)}
\]
La composition des morphismes de Gysin donne le diagramme commutatif
suivant dans $\detaleplus y$ :
\[
\xymatrix{
\ar[rrd]_{\Cl {i_x}}\Lambda\ar[rr]^{\Cl {i_y}}
&   & i_y^?\Lambda\ar[d]^{i^?_y\Cl {i'}} \\
&  & i_x^?\Lambda}
\]

Sur ce diagramme, on voit aussitôt que (a) et (b) impliquent (c) et que
(a) et (c) impliquent (b). Montrons que (b) et (c) impliquent (a). Il
s'agit de montrer que pour tout point $y$ de $Y$, le morphisme 
$i_y^\star \Cl {i'}$ est un isomorphisme.
On peut procéder par récurrence sur la
dimension de $V(y)$. On peut ainsi supposer que le support d'un
cône $C$ du morphisme
$\Cl {i'}$ dans $\detaleplus {V(y)}$ est contenu dans $\left\{y\right\}$.
Mézalor, le morphisme
canonique $i_y^!C\to i_y^\star C$ est un isomorphisme ; le diagramme
ci-dessus montre que $i_y^!C=0$, ce qui permet de conclure que $C=0$ et
finalement d'obtenir (a).

\bigskip

Rappelons quelques propriétés importantes concernant la pureté ponctuelle :

\begin{proposition2}[{\cite[proposition~2.2.2]{AZUMINO@Fujiwara}}]%
\label{16-proposition-purete-ponctuelle-completion}
Soit $X$ un schéma local strictement hensélien régulier.
Le complété $\hat{X}$ est
ponctuellement pur en son point fermé si et seulement si $X$ l'est.
\end{proposition2}

\begin{proposition2}[{\cite[corollary~2.2.3]{AZUMINO@Fujiwara}}]%
\label{16-proposition-purete-egale-caracteristique}
Soit $k$ un corps premier. Si $X$ est schéma régulier qui est
un $k$-schéma, alors $X$ est ponctuellement pur.
\end{proposition2}

\subsection{Changement de coefficients}
\label{16-subsection-changement-de-coefficients}

\begin{proposition2}\label{16-proposition-lemme-chinois}
Soit $n$ un entier naturel non nul.
Soit $n=\prod_{j=1}^k \ell_j^{\nu_j}$ la factorisation de $n$ en produit de
puissances de nombres premiers distincts. Un couple régulier $(X,Y)$ est
pur relativement à l'anneau de coefficients $\ZZ/n\ZZ$ si et seulement
s'il l'est relativement à l'anneau de coefficients $\ZZ/\ell_j^{\nu_j}\ZZ$
pour tout $j\in\{1,\dots,k\}$.
\end{proposition2}

Cela résulte aussitôt du lemme chinois et du fait que si $m$ est un
entier naturel divisant $n$, alors pour toute immersion fermée régulière
$i\colon Y \to X$, le diagramme évident commute dans $\detalepluscoefs{Y}{\ZZ/n\ZZ}$ :
\[
\xymatrix{
\ZZ/n\ZZ\ar[r]^{\Cl i} \ar[d]& i^?\ZZ/n\ZZ\ar[d]\\
\ZZ/m\ZZ\ar[r]^{\Cl i}       & i^?\ZZ/m\ZZ}
\]

\bigskip

\begin{proposition2}\label{16-proposition-ops-n-egale-ell}
Soit $\ell$ un nombre premier. Pour tout entier $\nu\geq 1$, un couple
régulier $(X,Y)$ est pur relativement à l'anneau de coefficients
$\ZZ/\ell\ZZ$ si et seulement s'il l'est relativement à l'anneau de
coefficients $\ZZ/\ell^\nu\ZZ$.
\end{proposition2}

En utilisant la résolution de Godement des faisceaux $\ZZ/\ell^\nu\ZZ(c)$
(où $c$ est la codimension de l'immersion $i\colon Y\to X$) pour tout
$\nu$, on peut représenter les morphismes de Gysin
$\Cl i\colon \ZZ/\ell^\nu\ZZ\to i^?\ZZ/\ell^\nu\ZZ$ dans
$\detalepluscoefs{Y}{\ZZ/\ell^\nu\ZZ}$ par des cocycles. Un tel cocycle pour
$\nu_0$ fixé induit pour tout entier $\nu\leq \nu_0$ un cocycle
représentant le morphisme de Gysin à coefficients dans $\ZZ/\ell^\nu\ZZ$.
Les propriétés élémentaires de la résolution de Godement font que, si on le
souhaite, on peut en fait trouver une famille compatible de cocycles pour
tout $\nu\in\NN$.

Compte tenu de ces observations, une fois ces cocycles convenablement
choisis, on dispose d'un cône privilégié $C(\nu)$ du morphisme $\Cl
i\colon\ZZ/\ell^\nu\ZZ\to i^?\ZZ/\ell^\nu\ZZ$ dans
$\detalepluscoefs{Y}{\ZZ/\ell^\nu\ZZ}$ pour tout $\nu\in\NN$ et de triangles
\[
\xymatrix{
C(\mu)\ar[r] & C(\mu+\nu)\ar[r] & C(\nu) \ar[r] & C(\mu)[1] }
\]
dans
$\detalepluscoefs{Y}{\ZZ/\ell^{\mu+\nu}\ZZ}$ pour tous $(\mu,\nu)\in\NN^2$.

Par conséquent, si $C(1)=0$, il vient que pour tout $\nu\geq 1$,
$C(\nu)=0$. Inversement, si $C(1)$ est non nul, son premier objet de
cohomologie non nul s'injecte dans celui de $C(\nu)$ pour tout $\nu\geq 1$.

\subsection{Diviseurs réguliers}
\label{16-subsection-diviseurs-reguliers}

\begin{definition2}
Si $X$ est un schéma et $\overline x \to X$ un point géométrique, on note
$V(\overline x)$ l'hensélisé strict de $X$ en $\overline x$ et 
$i_{\overline x}\colon V(\overline x)\to X$ le morphisme canonique.
\end{definition2}

\begin{proposition2}\label{16-proposition-caracterisation-purete-diviseur}
Soit $X$ un schéma régulier. Soit $D$ un diviseur régulier de $X$. Le
couple régulier $(X,D)$ est pur si et seulement si pour tout point
géométrique $\overline x\to D$, on a
$\HH^q_\et(V(\overline x)-i_{\overline x}^{-1}(D),\Lambda)=0$ pour tout
$q\geq 2$.
\end{proposition2}

Cela résulte du calcul de
$\HH^q_{\et}(V(\overline x)-i_{\overline x}^{-1}(D),\Lambda)$ pour
$q\in\{0,1\}$ (cf.~\sgaquatreetdemi{Cycle}{2.1.4}).

\bigskip

\begin{proposition2}\label{16-proposition-purete-argument-de-transfert-facile}
On suppose que l'anneau de coefficients est $\ZZ/\ell\ZZ$ où $\ell$ est un
nombre premier.
Soit $f\colon Y\to X$ un morphisme fini et plat de degré constant premier à
$\ell$ entre $\ZZn[\ell]$-schémas réguliers. Soit $D$ un diviseur régulier
de $X$. On suppose que $D'=f^{-1}(D)_{\red}$ est un diviseur régulier de
$Y$. Si le couple régulier $(Y,D')$ est pur, alors $(X,D)$ aussi.
\end{proposition2}

Grâce à la proposition~\ref{16-proposition-caracterisation-purete-diviseur},
on peut choisir un point géométrique de $D$ et remplacer $X$ par son
hensélisé strict en ce point. On suppose donc que $X$ et $D$ sont locaux
strictement henséliens et on se concentre sur la pureté du couple $(X,D)$
en le point fermé de $D$.
Le schéma $Y$ est
alors réunion disjointe finie de schémas locaux strictement henséliens ; au
moins un de ceux-ci est de degré premier à $\ell$ sur $X$. On peut donc
supposer que $Y$ aussi est local strictement hensélien.
Il suffit alors de montrer que $\HH^q(X-D,\ZZ/\ell\ZZ)$ s'injecte dans
$\HH^q(Y-D',\ZZ/\ell\ZZ)$, ce qui résulte du lemme suivant :

\begin{lemme2}
On suppose que l'anneau de coefficients $\Lambda$
est $\ZZ/\ell\ZZ$ où $\ell$ est un
nombre premier.
Soit $f\colon Y\to X$ un morphisme de présentation finie, fini et plat de
rang $d$ premier à $\ell$ entre $\ZZn[\ell]$-schémas.
Alors,
le morphisme canonique $\Lambda\to f_\star \Lambda$ est un monomorphisme
scindé dans $\detaleplus X$.
\end{lemme2}

D'après \sga{4}{xvii}{6.2.3}, on a un morphisme $\Tr_f\colon
f_\star\Lambda\to\Lambda$ tel que la composée
\[\Lambda\to f_\star\Lambda\to\Lambda\]
soit la multiplication par $d$, ce qui donne le scindage voulu puisque $d$
est inversible dans $\Lambda$.

\bigskip

\begin{proposition2}\label{16-proposition-up-and-down}
On suppose que l'anneau de coefficients $\Lambda$
est $\ZZ/\ell\ZZ$ où $\ell$ est un
nombre premier. Soit $X$ un $\ZZn[\ell]$-schéma régulier.
Soit $f$ une fonction sur $X$
dont le lieu des zéros $D=V(f)$ soit un diviseur régulier de $X$. On pose
$X'=\SP\left(\OO_X[T]/(T^\ell-f)\right)$. On note $\pi\colon X'\to X$ la projection,
$D'=\pi^{-1}(D)_{\red}$ (noter que $D'\to D$ est un isomorphisme). Alors,
$X'$ est un schéma régulier, et le couple régulier $(X',D')$ est pur si et
seulement si le couple régulier $(X,D)$ l'est.
\end{proposition2}

Soit $\overline x$ un point géométrique de $D$ (on identifiera aussi
$\overline x$ à un point géométrique de $D'$).
On va en fait montrer que
$(X,D)$ est pur en $\overline x$ si et seulement si $(X',D')$
l'est. On peut supposer que $X$ est le spectre premier d'un anneau local
strictement hensélien $A$ d'idéal maximal $\got{m}$ et que $\overline
x$ est au-dessus du point fermé de $X$. On a évidemment $f\in\got{m}$ ;
le fait que $D=V(f)$ soit régulier revient à dire que $f\not\in 
{\got{m}}^2$.

Notons $A'=A[T]/(T^\ell-f)$. En considérant le déterminant de
l'endomorphisme de $A'$ comme $A$-module donné par la multiplication par un
élément $b\in A'$, on observe que $b$ est inversible dans $A'$ si et
seulement si son image dans l'algèbre locale $(A/\got{m})[T]/(T^\ell)$
est inversible. Il en résulte que $A'$ est local d'idéal maximal
$\got{m}'=(T)+\got{m}A'$. Par ailleurs, on a un isomorphisme $A/(f)\iso
A'/(T)$ (\emph{i.e.} $D'\to D$ est un isomorphisme). On construit facilement un
isomorphisme $A_{(f)}[T]/(T^d-f)\iso A'_{(T)}$, ce qui montre que le
localisé de $A'$ en $(T)$ est un anneau de valuation discrète. La
codimension de l'idéal premier $(T)$ dans $A'$ est donc $1$. Compte tenu du
fait que $A'/(T)$ soit régulier, il en résulte que $X'$ est régulier.

On peut considérer, pour tout entier $\nu\geq 0$, le
$X$-schéma affine $X^\nu=\SP\left(\OO_X[T]/(T^{\ell^\nu}-f)\right)$. En élevant $T$ à
la puissance $\ell$, on obtient une tour de morphismes
\[\dots \to X^{\nu+1}\to X^{\nu} \to \dots \to X^1\to X^0\;\text{,}\]
le dernier morphisme $X^1\to X^0$ s'identifiant à $\pi\colon X'\to X$.
Au-dessus de $X-D$, cette tour de morphismes définit une tour de revêtements
étales galoisiens de $X-D$ de groupe de Galois
$\ZZ_\ell(1)=\underset{\nu}\lim\,\mu_{\ell^\nu}$ où pour tout entier $n$
inversible dans $A$, on note simplement $\mu_n=\mu_n(A)$.

Cette tour de revêtements définit un morphisme (surjectif) de groupes
profinis $\pi_1^\et(X-D)^\ab\to \ZZ_\ell(1)$. Par conséquent, on a un morphisme
de topos $\rho_f\colon (X-D)_\et\to \BB\ZZ_\ell(1)$ où $\BB\ZZ_\ell(1)$
désigne le topos des $\ZZ_\ell(1)$-ensembles discrets.

\begin{lemme2}\label{16-lemme-interpretation-purete-diviseur}
Le couple régulier $(X,D)$ est pur en $\overline x$
si et seulement si le morphisme
\[\RG(\BB\ZZ_\ell(1),\mu_\ell)\to\RG((X-D)_\et,
\mu_\ell)\]
induit par le morphisme de topos $\rho_f$ est un isomorphisme dans la
catégorie dérivée des groupes abéliens.
\end{lemme2}

Ce lemme découle des deux lemmes suivants :

\begin{lemme2}
Pour tout entier $q\geq 2$, $\HH^q(\BB\ZZ_\ell(1),\mu_\ell)=0$ et on a
des isomorphismes canoniques
\[\HH^0(\BB\ZZ_\ell(1),\mu_\ell)\simeq \mu_\ell\;\text{,}\quad
\HH^1(\BB\ZZ_\ell(1),\mu_\ell)\simeq  \Hom(\ZZ_\ell(1),\mu_\ell)
\simeq \ZZ/\ell\ZZ\;\text{.}\]
\end{lemme2}

Il s'agit de montrer que $\ZZ_\ell(1)$ est de $\ell$-dimension
cohomologique $1$. Pour cela, voir par exemple
\cite[§3.4, chapitre~I]{CG@Serre}.

\begin{lemme2}
Le morphisme composé \[\ZZ/\ell\ZZ\simeq
\HH^1(\BB\ZZ_\ell(1),\mu_\ell)\vers{\rho_f^\star}\HH^1((X-D)_\et,\mu_\ell)\iso
\HH^2_D(X,\mu_\ell)\] est donné (au signe près) par la classe de Gysin
$\Cl{D\subset X}$.
\end{lemme2}

Un relèvement dans $\HH^1(X-D,\mu_\ell)$ de $\Cl{D\subset X}$ est donné par
l'image de $f$ par le morphisme de bord $\HH^0((X-D)_\et,\Gm)\to
\HH^1((X-D)_\et,\mu_\ell)$ associé à la suite exacte de Kummer
\[0\to \mu_\ell\to \Gm\vers{p} \Gm\to 0\;\text{,}\]
autrement dit par le $\mu_\ell$-torseur $p^{-1}(f)\subset \Gm$.
Géométriquement, ce
torseur s'identifie tautologiquement au revêtement galoisien $X'-D'\to X-D$
de groupe de Galois $\mu_\ell$. Bien sûr,
la classe de ce $\mu_\ell$-torseur est 
donnée par le morphisme évident $\pi_1^\et(X-D)^\ab\to\ZZ_\ell(1)\to
\mu_\ell$ donnant l'image par $\rho_f^\star$ du générateur canonique de
$\HH^1(\BB\ZZ_\ell(1),\mu_\ell)$.

\bigskip

On peut appliquer le lemme~\ref{16-lemme-interpretation-purete-diviseur} à
$X'$ : il vient que le couple régulier $(X',D')$ est pur en $\overline x$
si et seulement si
le morphisme
\[\RG(\BB\ZZ_\ell(1),\mu_\ell)\to\RG((X'-D')_\et,
\mu_\ell)\]
induit par le morphisme de topos $\rho_T\colon (X'-D')_\et\to
\BB\ZZ_\ell(1)$ est un isomorphisme.

On dispose d'un carré commutatif de topos :
\[
\xymatrix{
\ar[d]^g (X'-D')_\et\ar[r]^{\rho_T} & \BB\ZZ_\ell(1)
 \ar[d]_{g'}\\
(X-D)_\et\ar[r]^{\rho_f} & \BB\ZZ_\ell(1)}
\]
où $g$ est induit par $\pi\colon X'\to X$ et $g'$ par la multiplication par
$\ell$ sur $\ZZ_\ell(1)$. Le faisceau $g_\star \ZZ/\ell\ZZ$
s'identifie canoniquement à
$\rho_f^\star g'_\star(\ZZ/\ell\ZZ)$. Il en découle aisément que le couple
régulier $(X',D')$ est pur en $\overline x$
si et seulement si le morphisme canonique
\[\RG(\BB\ZZ_\ell(1),g'_\star(\ZZ/\ell\ZZ))\to
\RG((X-D)_\et,\rho_f^\star
g'_\star(\ZZ/\ell\ZZ)) \]
est un isomorphisme.

La cohomologie relative du morphisme de topos $\rho_f\colon (X-D)_\et\to
\BB\ZZ_\ell(1)$ définit un foncteur triangulé 
\[F\colon \dplus(\BB\ZZ_\ell(1),\ZZ/\ell\ZZ)\to
\dplus(\ZZ/\ell\ZZ)\]
tel que pour tout $K\in\dplus(\BB\ZZ_\ell(1),\ZZ/\ell\ZZ)$, $F(K)$ soit
isomorphe à un cône du morphisme canonique $\RG(\BB\ZZ_\ell(1),K)\to
\RG((X-D)_\et,\rho_f^\star K)$.

Le lemme suivant découle de ce qui précède :

\begin{lemme2}
Le couple régulier $(X,D)$ est pur en $\overline x$
si et seulement si $F(\ZZ/\ell\ZZ)=0$
tandis que $(X',D')$ est pur en $\overline x$
si et seulement si $F(g'_\star(\ZZ/\ell\ZZ))=0$.
\end{lemme2}

Comme $\ZZ_\ell(1)$ est un pro-$\ell$-groupe, le faisceau
$g'_\star(\ZZ/\ell\ZZ)$ est une extension successive de $l$ copies de
$\ZZ/\ell\ZZ$. Le foncteur $F$ étant triangulé, on en déduit aussitôt que
si $F(\ZZ/\ell\ZZ)$ est nul, alors $F(g'_\star\ZZ/\ell\ZZ)$ aussi et que si
$F(\ZZ/\ell\ZZ)$ est non nul, son premier objet de cohomologie non nul
s'injecte dans celui de $F(g'_\star(\ZZ/\ell\ZZ))$. Ceci achève la
démonstration de la proposition~\ref{16-proposition-up-and-down}.

\subsection{Schémas sur un trait}
\label{16-subsection-schemas-sur-un-trait}

Soit $S$ un trait. On note $s$ son point fermé, $\eta$ son point
générique et $\pi$ une uniformisante.

\begin{proposition2}\label{16-proposition-un-trait-est-ponctuellement-pur}
Le trait $S$ est ponctuellement pur.
\end{proposition2}

Il s'agit de montrer que $S$ est ponctuellement pur en son point fermé.
On peut supposer que $S$ est strictement hensélien ; cela résulte alors
facilement du fait que le corps des fractions de $S$ soit de $\ell$-dimension
cohomologique $1$ pour tout nombre premier inversible sur $S$
(cf.~\sga{4}{x}{2.2}).

\begin{proposition2}
Pour tout entier naturel $n$, l'espace affine $\Aff^n_S$ est ponctuellement
pur.
\end{proposition2}

D'après la proposition~\ref{16-proposition-purete-egale-caracteristique}, les
schémas $\Aff^n_s$ et $\Aff^n_\eta$ sont ponctuellement purs. Ainsi,
$\Aff^n_S$ est ponctuellement pur en les points de la fibre générique.
Pour établir la pureté ponctuelle de $\Aff^n_S$ en les
points de la fibre spéciale, on utilise la
proposition~\ref{16-proposition-deux-sur-trois} : il suffit de montrer que
le couple régulier $(\Aff^n_S,\Aff^n_s)$ est pur.
Le cas $n=0$ résulte de la
proposition~\ref{16-proposition-un-trait-est-ponctuellement-pur} et le cas
général en découle en vertu du théorème de changement de base lisse.

\begin{corollaire2}\label{16-corollaire-lisse-implique-ponctuellement-pur}
Un $S$-schéma lisse est ponctuellement pur.
\end{corollaire2}

\begin{definition2}
Soit $p\colon X\to S$ un morphisme de type fini, avec $X$
régulier et admettant un faisceau ample.
On pose $K_X=p^?\Lambda_S$ et on dispose d'un morphisme de
Gysin $\Cl{X/S} \colon \Lambda_X\to K_X$ dans $\detaleplus X$
(cf. définition~\ref{16-definition-cl-f-lissifiable}).
\end{definition2}

\begin{proposition2}
Soit $p\colon X\to S$ un morphisme de type fini, avec $X$
régulier et admettant un faisceau ample.
Le schéma $X$ est ponctuellement pur si et seulement si le morphisme
$\Cl{X/S}\colon \Lambda_X\to K_X$ est un isomorphisme dans $\detaleplus X$.
\end{proposition2}

On choisit une factorisation $X\vers i Y \vers q S$ de $p$ (dans la
catégorie $\Schic$) avec $Y$ lisse
sur $S$ et $i$ une immersion fermée (régulière). D'après le
théorème~\ref{16-theoreme-cl-f-transitivite} (ou plutôt par définition de
$\Cl{X/S}$), le diagramme suivant est
commutatif :
\[
\xymatrix{
\ar[rd]_{\Cl{X/S}}\Lambda\ar[r]^{\Cl i} & i^?\Lambda
\ar[d]^{i^?\Cl q}_\sim\\
& i^?q^?\Lambda}
\]
Le morphisme $q$ étant lisse, le morphisme de Gysin $\Cl q$ est un
isomorphisme. Par conséquent, $\Cl {X/S}$ est un isomorphisme si et
seulement si $\Cl i\colon \Lambda\to i^?\Lambda$ en est un. D'après la
proposition~\ref{16-proposition-deux-sur-trois} et compte tenu du fait que $Y$
soit ponctuellement pur (cf. corollaire~%
\ref{16-corollaire-lisse-implique-ponctuellement-pur}), ceci équivaut encore
à dire que $X$ est ponctuellement pur.

\begin{corollaire2}\label{16-corollaire-lisse-sur-ponctuellement-pur}
Soit $X$ un $S$-schéma de type fini qui est régulier. Soit $Y$ un
$X$-schéma lisse. Si $X$ est ponctuellement pur, alors $Y$ aussi.
\end{corollaire2}

\begin{proposition2}\label{16-proposition-purete-argument-de-transfert}
Soit $f\colon X\to Y$ un morphisme propre et dominant de $S$-schémas où $X$
et $Y$ sont supposés de type fini sur $S$, intègres, réguliers et
admettant des faisceaux amples.
On suppose de plus que $f$ est génériquement étale de degré $d$ inversible
dans $\Lambda$. Alors, la pureté ponctuelle de $X$ implique celle de $Y$.
\end{proposition2}

Le morphisme $f$ est localement d'intersection complète lissifiable de
dimension relative virtuelle zéro, d'où $f^?=f^!$. Le morphisme de Gysin
relatif à $f$ est un morphisme $\Cl f\colon \Lambda\to f^!\Lambda$.

\begin{lemme2}
On peut généraliser le morphisme $\Cl f\colon \Lambda\to f^!\Lambda$ en des
morphismes $f^\star M\to f^!M$, fonctoriellement en $M\in\detaleplus Y$.
\end{lemme2}

Le morphisme $\Cl f\colon \Lambda\to f^!\Lambda$ correspond par adjonction
à un morphisme $\R f_!\Lambda\to \Lambda$,
que l'on peut tensoriser avec $M$
pour obtenir (\emph{via} la formule de projection) un morphisme
$\R f_!f^\star M\to M$
qui correspond lui-même par adjonction au morphisme
$f^\star M\to f^!M$ du type recherché.

\medskip

En appliquant la fonctorialité de la construction du lemme
au morphisme $\Cl{Y/S}\colon \Lambda_Y
\to K_Y$, on obtient un diagramme commutatif dans $\detaleplus X$.
\[
\xymatrix{
\ar[d]_{f^\star(\Cl {Y/S})}f^\star \Lambda_Y\ar[r] & f^!\Lambda_Y
   \ar[d]^{f^!(\Cl {Y/S})} \\
f^\star K_Y\ar[r] & f^!K_Y }
\]

\emph{Via} l'isomorphisme canonique $f^\star\Lambda_Y\simeq \Lambda_X$, le
morphisme du haut s'identifie au morphisme
$\Cl f\colon\Lambda_X\to f^!\Lambda_Y$ ;
celui de droite est $f^!(\Cl {Y/S})$. D'après le
théorème~\ref{16-theoreme-cl-f-transitivite}, il vient que le morphisme
composé
$\Lambda_X\simeq f^\star \Lambda_Y\to f^!K_Y\simeq K_X$ est le morphisme de
Gysin $\Cl{X/S}$. On déduit de ceci un diagramme commutatif de la forme
suivante dans $\detaleplus X$ :

\[
\xymatrix{
\ar[d]_{f^\star(\Cl {Y/S})}f^\star \Lambda_Y\ar[r]^\sim &
\Lambda_X\ar[d]^{\Cl{X/S}} \ar[r] & f^!\Lambda_Y\ar[d]^{f^!\Cl{Y/S}} \\
f^\star K_Y\ar[r] & K_X\ar[r]^\sim & f^!K_Y}
\]

Comme $f$ est propre, on obtient par adjonction un nouveau diagramme
commutatif dans $\detaleplus Y$ :
\[
\xymatrix{
\Lambda_Y\ar[d]_{\Cl{Y/S}} \ar[r] & \R f_\star
\Lambda_X\ar[d]^{\R f_\star (\Cl{X/S})}\ar[r] &\Lambda_Y 
\ar[d]^{\Cl{Y/S}}\\
K_Y\ar[r] & \R f_\star K_X\ar[r]& K_Y}
\]

Le diagramme ci-dessus met en évidence une relation entre les morphismes
$\Cl{Y/S}$ et $\R f_\star(\Cl{X/S})$. Comme va le montrer le lemme
suivant, le premier morphisme est un facteur direct du second, ce qui
montre que la pureté ponctuelle de $X$ implique celle de $Y$, achevant la
démonstration de la
proposition~\ref{16-proposition-purete-argument-de-transfert}.

\begin{lemme2}
Sur le diagramme précédent,
les morphismes composés $\Lambda_Y\to \Lambda_Y$ et
$K_Y\to K_Y$ sont les multiplications par le degré $d$ (en particulier, ce
sont des isomorphismes).
\end{lemme2}

Comme $Y$ est connexe (non vide), on a un isomorphisme évident 
$\Lambda\iso \End_{\detaleplus Y}(\Lambda_Y)$. D'après le théorème de bidualité
locale (cf. \sgaquatreetdemi{Th.~finitude}{4.3}),
on a aussi un isomorphisme $\Lambda\iso
\End_{\detaleplus Y}(K_Y)$. Il suffit donc d'obtenir la conclusion au-dessus d'un
ouvert non vide de $Y$.
Quitte à remplacer $Y$ par un ouvert non vide convenable, on peut supposer
que $f$ est un revêtement étale.
On est ainsi ramené au lemme suivant :

\begin{lemme2}
Soit $f\colon X\to Y$ un morphisme de schémas fini étale de degré constant
$d$. Pour tout objet $M\in \detaleplus Y$, le morphisme composé
\[M\to f_\star f^\star M\to M\]
déduit des adjonctions canoniques $(f^\star,f_\star)$ et
$(f_\star,f^\star)$ est la multiplication par $d$.
\end{lemme2}

Grâce aux formules de projection, on peut supposer que $M=\Lambda_Y$. Il
suffit alors d'établir le résultat après un changement de base étale (non
vide) trivialisant le revêtement $X\to Y$ (par exemple une clôture
galoisienne de ce revêtement). Bref, on peut supposer que $X$ est une
réunion disjointe de $d$ copies de $Y$, auquel cas le résultat est trivial.

\bigskip

\begin{definition2}\label{16-definition-schemas-v}
Soit $(e_1,\dots,e_n)\in\NN^n$.
On définit un $S$-schéma :
\[V(S,\pi,e_1,\dots,e_n)=
\SP\left(\OO_S[T_1,\dots,T_n]/
\left(\textstyle\prod_{i=1}^n T_i^{e_i}-\pi\right)\right)\;\text{.}\]
Pour tout $i$, on note $H_i$ le sous-schéma fermé de
$V(S,\pi,e_1,\dots,e_n)$ défini par l'équation $T_i=0$.
\end{definition2}

\begin{proposition2}\label{16-proposition-standards-reguliers}
Soit $(e_1,\dots,e_n)$ un $n$-uplet d'entiers naturels non tous nuls.
Alors,
le schéma $V(S,\pi,e_1,\dots,e_n)$ est régulier et ponctuellement pur.
\end{proposition2}

On peut supposer que l'anneau des coefficients est $\ZZ/\ell\ZZ$ où $\ell$
est un nombre premier inversible sur $S$.

\begin{lemme2}~

\begin{enumerate}
\item[(i)] Si $(e_1,\dots,e_n)$ est un $n$-uplet d'entiers non tous nuls
dont au moins un est inversible dans $\eta$, le $S$-schéma
$V(S,\pi,e_1,\dots,e_n)$ est intègre, régulier et de fibre générique
lisse ;
\item[(ii)] Soit $d\geq 1$, si $S'$ est le trait obtenu en
extrayant une racine $d$-ième $\pi'$
de l'uniformisante $\pi$\;\footnote{Rappelons brièvement pourquoi $S'$ est
bien un trait. Si $A$ est un anneau de valuation discrète, $\pi$ une
uniformisante et $d\geq 1$, notons $A'=A[X]/(X^d-\pi)$. Il est évident que
$A'$ est local noethérien et que son idéal maximal est
engendré par $X$. On vérifie aussitôt que la
suite évidente $0\to A'\vers X A' \to A'/(X) \to 0$ de $A'$-modules est
exacte, ce qui donne une résolution projective du corps résiduel de $A'$.
D'après \cite[§D.1, chapitre~IV]{Algebre@Serre},
$A'$ est régulier de dimension
$0$ ou $1$ et comme $A'$ n'est pas un corps, c'est un anneau de valuation
discrète.},
pour tout $n$-uplet $(e_1,\dots,e_n)$, on a un
isomorphisme de schémas
\[V(S',\pi',e_1,\dots,e_n)=V(S,\pi,de_1,\dots,de_n)\;\text{;}\]
\item[(iii)] Si $(e_1,\dots,e_n)$ est un $n$-uplet d'entiers non tous nuls,
le schéma $V(S,\pi,e_1,\dots,e_n)$ est régulier et intègre ;
\item[(iv)] Soit $(e_1,\dots,e_n)$ un $n$-uplet d'entiers non tous nuls,
le schéma $V(S,\pi,e_1,\dots,e_n)$ est ponctuellement pur si et seulement si
pour tout $i$ tel que $e_i>0$, le couple régulier
$(V(S,\pi,e_1,\dots,e_n),H_i)$ est pur\;\footnote{Si
$e_i=0$, c'est vrai aussi : c'est un cas particulier du théorème de pureté
relatif, cf.~\sga{4}{xvi}{3.7}.} ;
\item[(v)] Soit $(e_1,\dots,e_n)$ un $n$-uplet d'entiers non nuls,
soit $e$ le p.p.c.m. des $e_i$ ;
on suppose que $\ell$ ne divise pas $e$ ; si
$V(S,\pi,e,\dots,e)$ est
ponctuellement pur, alors $V(S,\pi,e_1,\dots,e_n)$ aussi ;
\item[(vi)] Soit $(e_1,\dots,e_n)$ un $n$-uplet d'entiers non tous nuls,
$V(S,\pi,e_1,\dots,e_n)$ est ponctuellement pur si et seulement si
$V(S,\pi,\ell e_1,e_2,\dots,e_n)$ est ponctuellement pur.
\end{enumerate}
\end{lemme2}

Les assertions (i) et (ii) sont laissées en exercice au lecteur.
L'assertion (iii) résulte aussitôt de (i) et de (ii).

Pour montrer l'assertion~(iv),
il suffit d'observer que les diviseurs $H_i$ pour
$e_i>0$ sont
ponctuellement purs (ce sont des espaces affines sur le corps résiduel de
$S$) et forment un recouvrement de la fibre spéciale de
$V(S,\pi,e_1,\dots,e_n)$. La fibre générique du schéma
$V(S_,\pi,e_1,\dots,e_n)$
étant ponctuellement pure (puisque lisse sur une extension
de $\eta$), on peut conclure en
utilisant la proposition~\ref{16-proposition-deux-sur-trois}.

Concernant l'assertion~(v),
l'élévation des $T_i$ à la puissance $\frac e
{e_i}$ définit un morphisme fini et plat $V(S,\pi,e,\dots,e)\to
V(S,\pi,e_1,\dots,e_n)$ de degré $\frac{e^n}{e_1\dots e_n}$ (premier à
$\ell$) ; compte tenu du critère~(iv), la
proposition~\ref{16-proposition-purete-argument-de-transfert-facile} permet de
conclure.

Pour établir~(vi), remarquons que l'élévation de $T_1$ à la puissance
$\ell$ définit un morphisme fini et plat 
$V(S,\pi,\ell e_1,e_2,\dots,e_n)\to V(S,\pi,e_1,\dots,e_n)$
de degré $\ell$ et étale en dehors du lieu d'annulation de $T_1$.
Il suffit
donc de montrer que $(V(S,\pi,\ell e_1,e_2,\dots,e_n),H_1)$
est pur si et seulement si
$(V(S,\pi,e_1,\dots,e_n),H_1)$ l'est, ce qui résulte
de la proposition~\ref{16-proposition-up-and-down}.

\medskip

Établissons la proposition~\ref{16-proposition-standards-reguliers}. D'après
l'assertion (iii), les schémas considérés sont réguliers. Pour établir leur
pureté ponctuelle, d'après
le corollaire~\ref{16-corollaire-lisse-sur-ponctuellement-pur}, on peut
supposer qu'aucun des exposants $e_i$ n'est nul. Dans le cas où les tous
les entiers $e_i$ valent $1$, le résultat est établi dans 
\cite[theorem~1.4]{DWORK@Illusie} (voir aussi
\cite[Satz~2.21]{UEBER@Rapoport-Zink}). Grâce à l'utilisation d'un trait auxiliaire,
l'assertion~(ii) permet d'en déduire
que pour tout entier $d\geq 1$, $V(S,\pi,d,\dots,d)$
est ponctuellement pur. En utilisant l'assertion~(v), on obtient que
$V(S,\pi,e_1,\dots,e_n)$ est ponctuellement pur si $\ell$ ne divise aucun
des entiers $e_i$. L'assertion~(vi) permet de passer au cas général.

\subsection{Géométrie logarithmique}
\label{16-subsection-geometrie-logarithmique}

\begin{definition2}
Soit $S$ un trait, de point générique $\eta$. La log-structure canonique
sur $S$ est la log-structure image directe de la log-structure triviale sur
$\eta$. Toute uniformisante de $S$ définit un morphisme de monoïdes
$\NN\to\Gamma(S,\OO_S)$ donnant naissance à une carte 
$S\to \SP(\ZZ[\NN])$ du log-schéma $S$.
\end{definition2}

L'objectif de cette sous-section d'établir le théorème suivant :

\begin{theoreme2}\label{16-theoreme-log-lisse-regulier-implique-p-pur}
Soit $S$ un trait muni de sa log-structure canonique. Soit $(X,M)\to S$  un
morphisme log-lisse de log-schémas fs. Si le schéma $X$ est régulier, alors
il est ponctuellement pur.
\end{theoreme2}

La proposition suivante précise \cite[theorem~3.5]{LSFI@Kato} dans le cas des
log-schémas fs :

\begin{proposition2}\label{16-proposition-critere-log-lissite-fs}
Soit $(X,M)\to (Y,N)$ un morphisme log-lisse entre log-schémas fs. On
suppose donnée une carte $Y\to\SP(\ZZ[Q])$ de $(Y,N)$ où $Q$ est un monoïde
fs sans torsion\;\footnote{Si $\overline y$ est un point géométrique de
$Y$, il existe un voisinage étale de $\overline y$ admettant une telle
carte avec $Q=M_{\overline y}/\OO_{Y,\overline y}^\times$ qui est fs 
saillant.}.
Pour tout point géométrique $\overline x$ de $X$, il
existe un voisinage étale $U$ de $\overline x$, un morphisme injectif de
monoïde $Q\to P$ avec $P$ fs sans torsion et une carte $U\to \ZZ[P]$ tels
que la partie de torsion de $\coker(Q^{\gp}\to P^{\gp})$ soit 
d'ordre inversible sur $U$ et que le morphisme de schémas
$U\to Y\times_{\SP(\ZZ[Q])}\SP(\ZZ[P])$ soit étale.
\end{proposition2}

Dans la démonstration du critère de log-lissité de
\cite[theorem~3.5]{LSFI@Kato}, des éléments $t_1,\dots,t_r$ de $M_{\overline x}$
sont choisis de sorte que la famille 
$(d\log t_1,\dots,d\log t_r)$ forme une base du faisceau des
log-différentielles
$\omega^1_{X/Y,\overline x}$. On considère ensuite le morphisme de monoïdes
évident $\NN^r\oplus Q\to M_{\overline x}$ donné sur la composante $\NN^r$
par les $t_1,\dots,t_r$. Il est tel que le conoyau de $\ZZ^r\oplus
Q^{\gp}\to M_{\overline x}^{\gp}$ soit fini d'exposant $n$ inversible dans
l'anneau $\OO_{X,\overline x}$ (en particulier, $\OO_{X,\overline
x}^\times$ est $n$-divisible). Il existe un
morphisme injectif $\ZZ^r\oplus Q^{\gp}\to G$ de conoyau tué par une
puissance de $n$ et un
prolongement $h\colon G\to M_{\overline x}^{\gp}$
de $\ZZ^r\oplus Q^{\gp}\to M_{\overline x}^{\gp}$
tel que $G\to M_{\overline
x}^{\gp}/\OO^\times_{X,\overline x}$ soit surjectif. Comme $M_{\overline
x}^{\gp}/\OO^\times_{X,\overline x}$ est un groupe abélien de type fini et
sans torsion, le lemme suivant montre que l'on peut s'arranger pour que $G$
soit un groupe abélien libre.
Dans la démonstration de \cite[theorem~3.5]{LSFI@Kato}, on pose ensuite
$P=h^{-1}(M_{\overline x})$ et il est montré que sur un voisinage étale $U$
de $\overline x$, $P$ engendre la log-structure de $(X,M)$ et que
le morphisme
de schémas
$U\to S\times_{\SP(\ZZ[Q])}\SP(\ZZ[P])$ est étale en $\overline x$. Le
monoïde $P$ ainsi construit est fs et sans torsion.

\begin{lemme2}
Soit $n$ un entier naturel non nul.
Soit $A$ un groupe abélien libre de type fini.
Soit $\varphi\colon A\to B$ un morphisme de groupes abéliens.
Soit $U\subset B$ un sous-groupe $n$-divisible. On suppose que $B/U$ est
sans torsion et que $\coker(A\to
B/U)$ est fini et tué par $n$.
Alors, il existe un groupe abélien $A'$ libre de type fini,
un morphisme injectif $A\to A'$ de groupes abéliens tel que $A'/A$ soit tué
par une puissance de 
$n$ et une extension $A'\to B$ du morphisme $A\to B$ telle que le
morphisme composé $A'\to B/U$ soit surjectif.
\end{lemme2}

Grâce à une récurrence sur l'ordre de $\coker(A\to B/U)$, on peut supposer
que $\coker(A\to B/U)$ est cyclique d'ordre $d\geq 2$, engendré par la
classe d'un élément $b\in B$. Il existe donc $a\in A$ et $u\in U$ tels que
$db=\varphi(a)+u$. Comme $u$ est $n$-divisible, il existe $\tilde{u}\in U$
tel que $u=d\tilde{u}$. Quitte à remplacer $b$ par $b-\tilde{u}$, on peut
supposer que $u=0$. On forme le carré cocartésien suivant dans la catégorie
des groupes abéliens :

\[
\xymatrix{ \mathbf{Z}\ar[r]^a\ar[d] & A \ar[d]\\
\frac 1 d \mathbf{Z}\ar[r] & A'}
\]

En raison de la relation $db=\varphi(a)$, on peut définir un unique
morphisme de groupes abéliens $\varphi'\colon A'\to B$ induisant
$\varphi\colon A\to B$ et envoyant $\frac 1 d$ sur $b$. On obtient ainsi
une surjection $A'\to B/U$ induisant un isomorphisme $A'/A\iso
\coker(A\to B/U)$. Il reste à vérifier que $A'$ est sans torsion. Soit $a'$
un élément de torsion de $A'$. L'image de $a'$ dans $B/U$ \emph{via}
$\varphi'$ est de torsion, mais $B/U$ étant sans torsion, on a
$\varphi(a')\in U$. Comme $\varphi'$ induit un isomorphisme $A'/A\iso
\coker(A\to B/U)$, on en déduit que $a'\in A$, mais $A$ est sans torsion,
donc $a'=0$.

\begin{proposition2}\label{16-proposition-structure-locale-log-lisse-regulier}
Soit $(X,M)\to S$ un log-schéma fs log-lisse sur un trait $S$ (muni de sa
log-structure canonique). On suppose que le schéma $X$ est régulier. Alors,
localement pour la topologie étale, $X$ admet un morphisme étale vers un
schéma $V(S,\pi,e_1,\dots,e_n)$ où $(e_1,\dots,e_n)$ est un $n$-uplet
d'entiers non tous nuls (cf.~définition~\ref{16-definition-schemas-v}).
\end{proposition2}

Soit $\pi$ une uniformisante de $S$ ; elle donne naissance
à une carte $S\to \SP(\ZZ[\NN])$.
D'après la
proposition~\ref{16-proposition-critere-log-lissite-fs}, on peut supposer qu'il
existe un monoïde $P$ fs sans torsion, un morphisme injectif $\NN\to P$, une
carte $X\to\SP(\ZZ[P])$ telle que le morphisme de schémas 
$X\to S\times_{\SP(\ZZ[\NN])}\SP(\ZZ[P])$ soit étale. Soit $\overline x$ un
point géométrique de $X$. On note $P'$ le sous-monoïde de $P$ formé des
éléments dont l'image dans $\Gamma(X,\OO_X)$ soit inversible au point
$\overline x$.

On peut supposer que $P'$ est un groupe. En effet, si $A$
est un sous-ensemble fini de $P'$ qui engendre le
groupe abélien (libre de type fini) $P'^{\gp}$\;\footnote{En fait,
on peut montrer que $P'$ est un monoïde de type fini
(c'est une face de $P$).},
on peut remplacer $X$ par
le voisinage ouvert de $\overline x$ sur lequel les images des éléments de
$A$ (et donc de $P'$) sont inversibles dans le faisceau structural et par
suite, remplacer $P$ par $P[-P']$ qui est encore fs et sans torsion.

Le fait que $X\to \SP(\ZZ[P])$ soit une carte implique alors que $P'$
est le noyau de
$P^{\gp}\to M_{\overline x}^{\gp}/\OO^\times_{X,\overline x}$. En
particulier, on obtient un isomorphisme
\[P/P'\iso M_{\overline
x}/\OO^\times_{X,\overline x}\;\text{.}\]
Comme $X$ est log-régulier, on reconnaît que $X$ est régulier au fait que
$M_{\overline x}/\OO^\times_{X,\overline x}$ soit un monoïde libre. Par
conséquent, il existe un entier $r$ et un isomorphisme de monoïdes
$\NN^r\iso P/P'$. On peut relever ce morphisme en un morphisme $\NN^r\to
P$, ce qui permet de construire un isomorphisme $\NN^r\oplus P'\iso
P$.

Il en résulte que le morphisme de carte $X\to \SP(\ZZ[P])$ a pour but un
schéma isomorphe à
$\SP(\ZZ[\NN^r\oplus P'])$ qui est le produit d'un espace affine
et d'un tore déployé (dont $P'$ est le groupe des caractères). Dans la
carte du morphisme $(X,M)\to S$ qui est donnée, l'image de $1$ par le
morphisme de monoïdes $\NN\to P$ peut s'écrire $(e_1,\dots,e_r,p')$ dans
$\NN^r\oplus P'$ \emph{via} les identifications ci-dessus. On peut choisir
une base $a_1,\dots,a_s$ de $P'$ comme groupe abélien telle que
$p'=\sum_{i=1}^s f_i a_i$ avec $f_i\in\NN$. On a ainsi construit un
morphisme étale $X\to V(S,\pi,e_1,\dots,e_r,f_1,\dots,f_s)$ (avec les
$e_1,\dots,e_r,f_1,\dots,f_s$ non tous nuls).

\medskip

Compte tenu de la proposition~\ref{16-proposition-standards-reguliers},
le théorème~\ref{16-theoreme-log-lisse-regulier-implique-p-pur} résulte
aussitôt de la
proposition~\ref{16-proposition-structure-locale-log-lisse-regulier}.

\subsection{Démonstration du théorème de pureté}
\label{16-subsection-demonstration-theoreme-purete}

Démontrons le théorème~\ref{16-theoreme-purete-absolue}.
D'après les propositions~\ref{16-proposition-lemme-chinois} et
\ref{16-proposition-ops-n-egale-ell}, on peut supposer que l'anneau des
coefficients $\Lambda$ est $\ZZ/\ell\ZZ$ où $\ell$ est un nombre premier.
D'après la proposition~\ref{16-proposition-deux-sur-trois}, il s'agit de montrer
que tout $\ZZn[\ell]$-schéma régulier est ponctuellement pur. D'après
\cite[corollary~6.1.5]{AZUMINO@Fujiwara}, on peut supposer que $X$ est un schéma
régulier intègre, quasi-projectif et plat 
sur un trait (strictement hensélien) $S$, que
l'on peut supposer d'inégale caractéristique d'après
la proposition~\ref{16-proposition-purete-egale-caracteristique}. On peut
utiliser les notations de la
sous-section~\ref{16-subsection-schemas-sur-un-trait}.
Quitte à étendre le
trait $S$, on peut supposer que l'anneau sous-jacent à $S$ est
intégralement fermé dans le corps des fonctions rationnelles sur $X$. La
fibre générique $X_\eta$ de $X$ est donc géométriquement intègre.

En appliquant \cite[proposition~4.4.1]{BS@Vidal} à la normalisation de
l'adhérence de $X$ dans un plongement projectif, on obtient qu'il existe un
groupe fini $G$ et un diagramme $G$-équivariant :
\[
\xymatrix{X'\ar[r]\ar[d] & \ar[d] X \\ S'\ar[r] & S}
\]
tels que :
\begin{itemize}
\item $G$ agisse trivialement sur $X$ et $S$ ;
\item $S'\to S$ soit une extension finie de traits ;
\item $X'\to X$ soit projectif, $X'$ soit régulier, connexe
et à réduction semi-stable sur $S'$ ;
\item $G$ agisse fidèlement sur $X'$ et $X'\to X$ soit génériquement un
revêtement étale galoisien de groupe $G$.
\end{itemize}

On munit $X'$ de la log-structure dont l'ouvert de trivialité est la fibre
générique de $X'\to S'$. Soit $H$ un $\ell$-Sylow de $G$. On note $T=S'/H$.
L'extension de traits (strictement henséliens) $S'\to T$ est d'ordre une
puissance de $\ell$, donc modérément ramifiée. Par conséquent, pour les
log-structures canoniques, $S'\to T$ est log-étale. Comme on sait que $X'$
est log-lisse sur $S'$, il l'est donc aussi sur $T$. Comme $H$ agit
trivialement sur $T$ et que son action sur $X'$ est modérée, on peut
appliquer le théorème de résolution équivariante qui
donne un morphisme projectif et birationnel $H$-équivariant
$X''\to X'$ de log-schémas tel que $X''$ soit log-lisse sur $T$ et que $H$
agisse très modérément sur $X''$. Le log-schéma quotient $X''/H$ est donc
log-lisse sur $T$ (en particulier, $X''/H$ est 
log-régulier). D'après le théorème de résolution
des log-singularités, il existe un log-éclatement
(en particulier, log-étale,
projectif et birationnel) $X'''\to X''/H$ tel que $X'''$
soit régulier. La situation est résumée sur le diagramme suivant :

\[
\xymatrix{& \ar@(lu,ru)^{H} X''\ar[r]^{\text{birat.}}\ar[d] 
& \ar@(lu,ru)^{G}X'\ar[r]\ar[d] & X\ar[dd] \\
X'''\ar[r]^{\text{birat.}}&\ar@{.>}[rru] X''/H\ar[rd] & \ar@(dl,ul)^{G}S'\ar[d] & \\
& & T\ar[r] &S}
\]

Le log-schéma $X'''$ est régulier et log-lisse sur $T$ ;
d'après le théorème~\ref{16-theoreme-log-lisse-regulier-implique-p-pur},
$X'''$ est ponctuellement pur. Le morphisme évident $X'''\to X$ est
projectif et génériquement un revêtement
étale de degré premier à $\ell$, d'après la 
proposition~\ref{16-proposition-purete-argument-de-transfert}, on peut
conclure que $X$ est ponctuellement pur, ce qui achève la démonstration du
théorème de pureté.


\ifx\danslelivre\undefined
\bibliography{../fichiers-configuration-latex/bibliographie}
\bibliographystyle{../fichiers-configuration-latex/smfalpha}


\begin{center}
XVII. Dualité\\
Joël Riou
\end{center}

\begin{center}
\verb$Id: dualite.tex,v 1.124 2008-05-09 13:52:48 cvs Exp $
\end{center}

\tableofcontents

\else
\setboolean{english}{false}
\setboolean{francais}{true}
\selectlanguage{francais}
\chapter{Dualité}
\addtocontents{toc}{par Joël Riou}
\begin{center}
Joël Riou
\end{center}
\fi


Ce texte vise à fournir une rédaction des résultats annoncés par Ofer
Gabber dans \cite{ABELIEN@Gabber} concernant les complexes dualisants dans
le contexte étale sur les schémas noethériens excellents. 

On fixe un entier naturel $n\geq 2$, on note
$\Lambda=\mathbf{Z}/n\mathbf{Z}$ : ce sera notre anneau de coefficients.
Si $X$ est un $\ZZn$-schéma noethérien, on note $\detalebc X$ la sous-catégorie
de $\detaleb X$ formée des complexes ayant des faisceaux de cohomologie
constructibles et $\detalebctf X$ la sous-catégorie pleine de $\detalebc X$ formée
des complexes de tor-dimension finie.

\begin{definition}\label{17-definition-complexe-dualisant}
Soit $X$ un $\ZZn$-schéma noethérien. Un complexe dualisant sur
$X$ est un objet $K\in\detalebc X$ tel que le foncteur de dualité 
$\Dual_K=\R\SheafHom(-,K)$ préserve $\detalebc X$ et que pour tout
$L\in\detalebc X$, le morphisme de bidualité $L\to\Dual_K\Dual_K L$ soit un
isomorphisme.
\end{definition}

Cette définition diffère de celle de \sga{5}{i}{1.7} dans la mesure où on ne
demande pas à un complexe dualisant d'être de dimension quasi-injective
finie.

Le théorème suivant récapitule l'essentiel des résultats de Gabber que nous
allons établir dans ces notes :

\begin{theoreme}[Gabber]\label{17-theoreme-gabber}
Soit $X$ un $\ZZn$-schéma noethérien excellent muni d'une fonction de
dimension (cf.~définition~\ref{17-definition-fonction-de-dimension}).
Alors, $X$ admet un complexe dualisant $K$, unique à produit
tensoriel avec un objet inversible près (cf.
proposition~\ref{17-proposition-objets-inversibles-d-t-a}).
Ce complexe dualisant $K$ appartient à $\detalebctf X$ ; il est de dimension
quasi-injective finie (autrement dit est un complexe dualisant au sens de
\sga{5}{i}{1.7}) si et seulement si $X$ est de dimension de Krull finie.

Par ailleurs, on a les résultats suivants :
\begin{itemize}
\item si $X$ est régulier, le faisceau constant $\Lambda$ est un complexe 
dualisant ;
\item si $f\colon Y\to X$ est un morphisme plat et à fibres géométriquement
régulières avec $Y$ noethérien excellent, alors $f^!K$ est un complexe
dualisant ;
\item si $f\colon Y\to X$ est un morphisme de type fini compactifiable,
alors $f^!K$ est un complexe dualisant.
\end{itemize}
\end{theoreme}

La démonstration s'appuie sur la notion de complexe dualisant potentiel
(définition~\ref{17-definition-complexe-dualisant-potentiel}) sur un
$\ZZn$-schéma noethérien excellent $X$ muni d'une fonction de dimension
$\delta$ : il s'agit de la donnée d'un complexe $K\in\detaleplus X$ muni
d'isomorphismes (appelés épinglages) $\SheafRG_x(K)\iso
\Lambda(\delta(x))[2\delta(x)]$ pour tout point $x\in X$, qui soient
compatibles aux morphismes de transition associés aux spécialisations
immédiates de points géométriques de $X$ (cf.
section~\ref{17-section-transition-codimension-un}). Gabber montre 
dans \cite[lemma~8.1]{T@Gabber}
que si un $\ZZn$-schéma noethérien excellent admet un
complexe dualisant, alors il admet aussi une fonction de dimension globale ;
cette hypothèse du théorème~\ref{17-theoreme-gabber}
est donc bien nécessaire.

Dans la section~\ref{17-section-premieres-proprietes}, nous verrons notamment
que sur un $\ZZn$-schéma noethérien régulier, le faisceau constant
$\Lambda$ est un complexe dualisant potentiel pour la fonction de dimension
$-\codim$ (cf. proposition~\ref{17-proposition-cdpotentiel-regulier}).
Pour ce faire, nous utiliserons de manière essentielle le
théorème de pureté cohomologique absolue démontré par Gabber, ainsi que les
propriétés des morphismes de Gysin établies
dans~\refext{XVI}{16-section-gysin}.

Dans la section~\ref{17-section-classe-de-cohomologie}, nous construirons un
isomorphisme $\Lambda\iso \HH^{2d}_x(X,\Lambda(d))$ où $X$ est un 
$\ZZn$-schéma noethérien
excellent local strictement hensélien normal de dimension $d$ et de
point fermé $x$, vérifierons que cet isomorphisme est compatible aux
spécialisations immédiates et nous servirons de ce résultat pour
construire des morphismes de transition $\HH^{i}_{\overline y}(X,K)\to
\HH^{i+2c}_{\overline x}(X,K(c))$ pour une spécialisation $\overline
y\to\overline x$ de codimension $c$ arbitraire entre points géométriques d'un
$\ZZn$-schéma noethérien excellent $X$, pour tout
$K\in\detaleplus X$. On utilisera
notamment la résolution des singularités pour les schémas noethériens
excellents de dimension $2$.

Dans la section~\ref{17-section-existence-unicite-cdp}, nous montrerons
l'existence et l'unicité à isomorphisme unique près d'un complexe dualisant
potentiel sur un $\ZZn$-schéma noethérien excellent muni d'une fonction de
dimension. Dans le cas d'un schéma normal, nous nous appuierons sur les
résultats de la section~\ref{17-section-classe-de-cohomologie} et sur les
résultats généraux de Gabber sur l'existence de $t$-structures définies par
des fonctions de perversité. Nous montrerons aussi que les complexes
dualisants potentiels vérifient de bonnes propriétés de stabilité par
rapport aux morphismes plats et à fibres géométriquement régulières et aux
morphismes de type fini.

Dans la section~\ref{17-section-dualite}, nous montrerons qu'un
complexe dualisant potentiel est un complexe dualisant. Une fois les
propriétés de finitude établies, nous procéderons par récurrence sur la
dimension, en utilisant d'une part une généralisation d'un argument de
\sgaquatreetdemi{Th.~finitude}{4.3} et d'autre part le théorème
d'algébrisation partielle~\refext{V}{5-algébrisation partielle mixte raffinée géométrique}.

Dans la section~\ref{17-section-coefficients-generaux}, nous montrons qu'à
partir d'un complexe dualisant à coefficients $\Lambda$, on peut construire
des complexes dualisants pour des anneaux de coefficients plus généraux.
Ces résultats sont essentiellement indépendants des sections précédentes.
Cependant, ce n'est qu'en combinant les résultats des sections précédentes
sur les complexes dualisants potentiels avec le résultat d'unicité des
complexes dualisants de la
proposition~\ref{17-proposition-unicite-complexe-dualisant} que l'on peut
déduire le théorème~\ref{17-theoreme-gabber}. En vertu de ce
théorème~\ref{17-theoreme-gabber}, les constructions de cette section donnent
en particulier des complexes dualisants pour des anneaux de coefficients
généraux sur les schémas noethériens excellents munis de fonctions de
dimension. Il semble très probable qu'il soit possible d'étendre les
résultats d'Ofer Gabber à des énoncés de dualité avec des coefficients
$\ell$-adiques. Toutefois, le rédacteur a renoncé à les rédiger.

Enfin, quelques résultats nécessaires à ce qui précède ont été rejetés en
appendice. On y décrit notamment une construction du produit tensoriel
dérivé sur la catégorie dérivée toute entière des faisceaux de Modules sur
un topos annelé (commutatif).

\section{Le morphisme de transition en \texorpdfstring{codimension~$1$}{codimension 1}}
\label{17-section-transition-codimension-un}
\subsection{Notations}

\begin{definition2}
Soit $X$ un $\ZZn$-schéma, soit $x\in X$, soit $X\in \detaleplus X$. On pose
$\SheafRG_x(K)=i_x^! K_{|X_{(x)}}\in \detaleplus x$ où $i_x$ est
l'inclusion du point fermé du localisé $X_{(x)}$\;\footnote{On obtiendrait
une définition équivalente en remplaçant le schéma local $X_{(x)}$ par son
hensélisé.}. Si $\overline x$ est un point géométrique de $X$ au-dessus de
$x$\;\footnote{Dans la suite, pour ne pas alourdir inutilement le texte, si
on fixe un point $x$ d'un schéma $X$, $\overline x$ désignera un point
géométrique au-dessus de $x$ et inversement, si on fixe un point
géométrique $\overline x\to X$, on notera $x$ le point de l'espace
topologique sous-jacent à $X$ au-dessus duquel $\overline x$ se trouve.},
on note $\RG_{\overline
x}(K)=\SheafRG_x(K)_{\overline x}\in \dplus(\Lambda)$ ; cet objet
s'identifie à la cohomologie à supports dans le point fermé de la
restriction de $K$ à l'hensélisé strict $X_{(\overline x)}$. Les objets de
cohomologie de $\RG_{\overline x}(K)$ seront notés $\HH^i_{\overline
x}(K)$ pour tout $i\in\mathbf{Z}$.
\end{definition2}

\begin{definition2}[\sga{4}{viii}{7.2}]
Si $\overline y$ et
$\overline x$ sont deux points géométriques d'un schéma $X$, une
spécialisation $\overline y\to \overline x$ est un $X$-morphisme
$X_{(\overline y)}\to X_{(\overline x)}$ entre les localisés stricts
correspondants, ce qui revient à la donnée d'un $X$-morphisme $\overline
y\to X_{(\overline x)}$. On définit la codimension d'une spécialisation
comme étant la dimension de l'adhérence du point de $X_{(\overline x)}$
en-dessous de $\overline y$. On dit qu'une spécialisation est immédiate si
elle est de codimension $1$.
\end{definition2}

On se propose ici de définir un morphisme de transition
$\transition^X_{\overline y\to\overline x}\colon
\RG_{\overline y}(K)\to \RG_{\overline
x}(K)(1)[2]$ dans $\dplus(\Lambda)$ pour toute spécialisation immédiate
$\overline y\to \overline x$ de points géométriques sur un $\ZZn$-schéma
excellent $X$, quel que soit $K\in\detaleplus X$.

\subsection{Cas d'un trait strictement hensélien}

Soit $X$ un trait strictement hensélien de point générique $\eta$ et de
point fermé $s$. 
Soit $\overline \eta$ un point géométrique au-dessus de
$\eta$. On va définir le morphisme de transition
\[\RG_{\overline \eta}(K)\to \RG_s(K)(1)[2]\]
pour tout $K\in\detaleplus X$. 

On note $p$ l'exposant caractéristique du corps résiduel de $X$, que l'on
suppose inversible dans $\Lambda$. On dispose
d'une suite exacte canonique de groupes profinis :
\[1\to S\to \Gal(\overline \eta/\eta)\to G\to 1\;\text{,}\]
où $G$ est le groupe d'inertie modérée, canoniquement isomorphe à
$\chap{\ZZ}'(1)$ et où $S$, le groupe de ramification sauvage, est un
pro-$p$-groupe (cf.~\cite[proposition~6.2.12]{ART@Gabber-Ramero}).
Remarquons que l'ordre de $G$ est multiple de
$(\#\Lambda)^\infty$ : ce fait sera utile au lecteur scrupuleux qui
voudrait vérifier en exercice les détails passés sous silence dans cette
sous-section.

\subsubsection{Algèbre du groupe \texorpdfstring{$\chap{\ZZ}'(1)$}{Z'(1)}}

\begin{definition3}
L'algèbre de groupe $\Lambda[[G]]$ est l'anneau des endomorphismes du
foncteur d'oubli de la catégorie des $\Lambda$-modules discrets munis
d'une action continue de $G$ vers celle des $\Lambda$-modules. Cette
algèbre est naturellement topologisée : elle est munie de la topologie la
moins fine qui soit telle que pour tout $\Lambda$-module discret $M$ muni
d'une action continue de $G$ et tout élément $m\in M$, l'application
$\Lambda[[G]]\to M$ qui à $a\in\Lambda[[G]]$ associe le résultat $a.m$ de
son action sur $m$ soit continue.
\end{definition3}

On a un isomorphisme canonique d'anneaux topologiques
\[\Lambda[[G]]\iso \lim \Lambda[G/H]\;\text{,}\]
où $H$ parcourt l'ensemble ordonné des sous-groupes ouverts distingués de
$G$ et où $\Lambda[G/H]$ est l'algèbre de groupe usuelle (discrète)
du groupe fini
$G/H$. On peut identifier les $\Lambda$-modules discrets munis d'une action
continue de $G$ aux $\Lambda[[G]]$-modules discrets.

L'action naturelle de $\Lambda[[G]]$ sur $\Lambda$ muni de l'action
triviale de $G$ définit un morphisme continu d'augmentation
$\varepsilon\colon\Lambda[[G]]\to \Lambda$. On note $I_G$ le noyau de
$\varepsilon$ : c'est l'idéal d'augmentation.

\begin{proposition3}\label{17-proposition-structure-i-g}
Le $\Lambda[[G]]$-module $I_G$ est libre de rang $1$, engendré par
$1-\sigma$ si $\sigma$ est un générateur topologique de $G$.
\end{proposition3}

La démonstration de cette proposition est laissée en exercice au lecteur.

\subsubsection{Description de la cohomologie à supports}

Il est clair que la catégorie des faisceaux (d'ensembles) sur $X_\et$ est
naturellement équivalente à la catégorie des flèches de $\Gal(\overline
\eta/\eta)$-ensembles (discrets) ${\mathscr M}_s\to 
{\mathscr M}_{\overline \eta}$ telles que 
l'action de $\Gal(\overline \eta/\eta)$ sur ${\mathscr M}_s$ soit triviale.

Pour tout faisceau de $\Lambda$-modules $\mathscr M$ sur $X_\et$, on note
$F(\mathscr M)$ le complexe évident de $\Lambda$-modules :
\[{\mathscr M}_s\to {\mathscr M}_{\overline \eta}^S\to
\Hom_{\Lambda[[G]]}(I_G,{\mathscr M}_{\overline \eta}^S)\]
où ${\mathscr M}_s$ est placé en degré $0$. Le foncteur qui à $\mathscr M$
associe $F(\mathscr M)$ étant additif, en utilisant la construction du
complexe simple associé à un complexe double, on peut naturellement définir
un complexe de $\Lambda$-modules $F(K)$  pour tout complexe $K$ de
faisceaux de $\Lambda$-modules sur $X_\et$. En observant en outre que le
foncteur qui à $\mathscr M$ associe $F(\mathscr M)$ est exact, on voit que 
le foncteur $F$ défini au niveau des complexes préserve les
quasi-isomorphismes. On a ainsi obtenu un foncteur $F\colon \detaleplus X\to
\dplus(\Lambda)$.

\begin{proposition3}\label{17-proposition-description-h-2-support-trait}
Pour tout $K\in\detaleplus X$, le morphisme évident
\[\RG_s(X,K)\to F(K)\]
est un isomorphisme dans $\dplus(\Lambda)$.
\end{proposition3}

On se ramène au cas où $K$ se réduit à un faisceau injectif $\mathscr M$
placé en degré $0$. Il est évident que le noyau de 
${\mathscr M}_{\overline \eta}^S\to 
\Hom_{\Lambda[[G]]}(I_G,{\mathscr M}_{\overline \eta}^S)$ s'identifie à
$\Gamma(\eta,\mathscr M)$.
On a donc bien un isomorphisme $\HH^0_s(X,\mathscr M)\iso
\HH^0(F(\mathscr M))$. Il reste à montrer que $\HH^1(F(\mathscr M))$ est nul.
Cela revient à dire qu'après application du foncteur
$\Hom_{\Lambda[[G]]}(-,{\mathscr M}_{\overline \eta}^S)$, l'inclusion de
l'idéal $I_G$ dans $\Lambda[[G]]$ induit une surjection.
Pour cela, il suffit de montrer que le $\Lambda[[G]]$-module discret
${\mathscr M}_{\overline \eta}^S$ est injectif :
par restriction à un ouvert, ${\mathscr M}_{|\eta}$
est un faisceau injectif de $\Lambda$-modules sur $\eta_\et$ et ${\mathscr
M}_{\overline \eta}^S$ est obtenu à partir de ${\mathscr M}_{|\eta}$ par
application du foncteur image directe associé au morphisme
de topos évident $\eta_\et \to \BB G$ (cf. \sga{4}{iv}{4.5.2}).

\subsubsection{Définition du morphisme de transition}

Soit $\sigma$ un générateur topologique de $G$. Soit $M$ un
$\Lambda[[G]]$-module discret. On observe que l'on a un isomorphisme
canonique de groupes abéliens $M(-1)\simeq \Hom(G,M)$. On définit un
morphisme de groupes abéliens
$M(-1)\to \Hom_{\Lambda[[G]]}(I_G,M)$ \emph{via} les
isomorphismes suivants :
\[\xymatrix{M(-1)\ar[r]^-{\sim} &
\Hom(G,M)\ar[r]_-\sim^-{\ev_\sigma} & M &
\ar[l]^-\sim_-{\ev_{1-\sigma}} \Hom_{\Lambda[[G]]}(I_G,M)\;\text{.}}\]

Si $\mathscr M$ est un faisceau de $\Lambda$-modules sur $X_\et$, en
composant l'application ci-dessus dans le cas de $M={\mathscr M}_{\overline
\eta}^S$ et le projecteur canonique 
${\mathscr M}_{\overline
\eta}\to {\mathscr M}_{\overline
\eta}^S$ ($S$ étant un pro-$p$-groupe et $p$ inversible dans $\Lambda$), on
définit un morphisme ${\mathscr M}_{\overline \eta}(-1)\to
\Hom_{\Lambda[[G]]}(I_G,{\mathscr M}_{\overline \eta}^S)$ et donc un
morphisme de complexes $s_\sigma\colon {\mathscr M}_{\overline \eta}(-1)\to
F(\mathscr M)[2]$, fonctoriel en $\mathscr M$.

\begin{definition3}\label{17-definition-transition-cas-trait}
On note $\transition^X_{\overline\eta\to s}\colon \RG_{\overline \eta}(K)\to
\RG_s(K)(1)[2]$ le morphisme dans $\dplus(\Lambda)$ défini par 
$s_\sigma$ fonctoriellement pour tout objet $K\in\detaleplus X$. D'après le
lemme suivant, ce morphisme $\transition^X_{\overline\eta\to s}$ 
est indépendant du générateur $\sigma$ de $G$ : c'est le morphisme de
transition associé à la spécialisation $\overline\eta\to s$ de points
géométriques de $X$.
\end{definition3}

\begin{lemme3}
Si $\sigma$ et $\sigma'$ sont deux générateurs topologiques de $G$, il
existe une unique homotopie fonctorielle en $\mathscr M$ entre $s_\sigma$
et $s_{\sigma'}$.
\end{lemme3}

On voit aussitôt qu'on peut se limiter aux $\mathscr M$ tels que 
${\mathscr M}_s=0$ et que $S$ agisse trivialement sur ${\mathscr
M}_{\overline \eta}$. On peut identifier cette catégorie de faisceaux à
celle des $\Lambda[[G]]$-modules discrets.

Soit $M$ un $\Lambda[[G]]$-module. Notons $\mathscr M$ le faisceau sur
$X_\et$ correspondant. On note $F_\sigma(M)$ le complexe
\[\xymatrix{\dots\ar[r] & 0\ar[r] & M\ar[r]^{1-\sigma} & M\ar[r] & 0 \ar[r] & \dots}\]
concentré en les degrés $1$ et $2$. On note $\Psi_\sigma\colon F(\mathscr
M)\iso F_\sigma(M)$ l'isomorphisme de complexes défini de façon évidente
à partir de $\sigma$ :
\[
\xymatrix{
\dots \ar[r] & \ar[d]0 \ar[r] &\ar@{=}[d] {\mathscr M}_{\overline \eta}\ar[r] &
\Hom_{\Lambda[[G]]}(I_G,{\mathscr M}_{\overline \eta})
\ar[d]^{\ev_{1-\sigma}}_\sim  \ar[r] & 0 \ar[r]\ar[d] & \dots\\ 
\dots \ar[r] & 0 \ar[r] & M\ar[r]^{1-\sigma} & M \ar[r] & 0 \ar[r] &\dots}
\]

Notons $\varphi_\sigma\colon {\mathscr M}_{\overline \eta}(-1)\iso M$
l'isomorphisme défini par l'évaluation en $\sigma$ \emph{via}
l'isomorphisme canonique ${\mathscr M}_{\overline\eta}(-1)\simeq
\Hom(G,{\mathscr M}_{\overline \eta})$. Notons $t_\sigma\colon M\to
F_\sigma(M)[2]$ le morphisme de complexe représenté par les flèches
verticales ci-dessous :
\[
\xymatrix{\dots \ar[r] & \ar[d]0 \ar[r] & \ar[d]0 \ar[r] &\ar@{=}[d] M\ar[r]
&\ar[d] 0 \ar[r] & \dots\\
\dots \ar[r] & 0 \ar[r] & M\ar[r]^{1-\sigma} & M\ar[r] & 0 \ar[r] & \dots}
\]

On dispose ainsi d'un carré commutatif de complexes, fonctoriel en $M$ :
\[
\xymatrix{{\mathscr M}_{\overline\eta}(-1)\ar[r]^-{s_\sigma}
\ar[d]^{\varphi_\sigma}_\sim& F(\mathscr
M)[2] \ar[d]^{\Psi_\sigma}_\sim\\
M\ar[r]^-{t_\sigma} & F_\sigma(M)[2]}
\]
Posons $f_{\sigma,\sigma'}=\Psi_{\sigma}\circ
s_{\sigma'}\circ\varphi_{\sigma}^{-1}$.
Les flèches verticales étant des isomorphismes de complexes, montrer que
les morphismes $s_\sigma,s_{\sigma'}\colon {\mathscr
M}_{\overline\eta}(-1)\to F(\mathscr M)[2]$ sont (fonctoriellement)
homotopes revient à vérifier que les deux morphismes
$t_\sigma,f_{\sigma,\sigma'}\colon M \to F_\sigma(M)[2]$ le sont.

On peut ainsi représenter la situation de façon plus concrète :
\[
\xymatrix{
\ar[d] 0\ar[r] & M\ar@/_/[d]_{\Id}\ar@/^/[d]^{g}\\
\ar[r]^{1-\sigma} M & M
}\]
où $g$ est la transformation naturelle induite par $f_{\sigma,\sigma'}$.

Comme $\Lambda[[G]]$ désigne précisément l'anneau des transformations
naturelles $M\to M$, on peut identifier $\Id$ et
$g$ à des éléments $1$ et $g$ de $\Lambda[[G]]$
respectivement.
Montrer l'existence et l'unicité de l'homotopie fonctorielle entre
$s_\sigma$ et $s_{\sigma'}$ se ramène donc à montrer
l'existence et l'unicité de $h\in\Lambda[[G]]$ tel que $(1-\sigma)\cdot
h=1-g$. D'après la proposition~\ref{17-proposition-structure-i-g}, cela
revient à montrer que $\varepsilon(g)=1$.
Si on note $u$ l'unité de $\Lambda[[G]]$ telle que 
$(1-\sigma')=u\cdot (1-\sigma)$ et qu'on note $\alpha$ l'élément de 
${\chap{\ZZ}{'}}^\times$
tel que $\sigma'=\sigma^\alpha$, alors on a la relation
$u\cdot g=\alpha$.
On est donc ramené à montrer que
$\varepsilon(u)=\varepsilon(\alpha)$. Pour cela, on utilise la formule
suivante :
\[
\frac{1-\sigma^\beta}{1-\sigma}=\sum_{i=0}^{\beta-1}\sigma^i\;\text{.}
\]
Cette formule est évidemment juste pour $\beta\in\mathbf{N}$ ; on peut lui
donner un sens pour tout $\beta\in\chap{\ZZ}$ en prolongeant chacun des
membres par continuité. En appliquant cette formule avec
$\beta=\alpha$, on obtient le
résultat voulu :
\[\varepsilon(u)=\varepsilon\left(\sum_{i=0}^{\alpha-1}\sigma^i\right)
=\sum_{i=0}^{\alpha-1}1=\varepsilon(\alpha)\;\text{,}\]
ce qui achève la démonstration du lemme.

\subsection{Cas d'un schéma local strictement hensélien intègre excellent de dimension~1}

Soit $X$ un $\ZZn$-schéma local strictement hensélien intègre excellent de
dimension $1$. Soit $\eta$ le point générique de $X$. Soit $\overline \eta$
un point géométrique au-dessus de $\eta$. Soit $s$ le point fermé de $X$.
Soit $\tilde X\vers f X$ la normalisation de $X$. Le schéma $\tilde X$ est
un trait strictement hensélien (de point fermé $\tilde s$)
et $f$ est un homéomorphisme
universel (en particulier, $\tilde s/s$ est une extension purement inséparable),
donc le foncteur image inverse $f^\star$ induit une équivalence
entre la catégorie des faisceaux sur $X_\et$ et sur ${\tilde X}_\et$.

\begin{definition2}\label{17-definition-transition-cas-local-dimension-un}
\emph{Via} les identifications ci-dessus, 
le morphisme de transition
$\transition^X_{\overline\eta \to s}$ est induit par $\frac 1 {[\tilde s:s]}
\transition^{\tilde X}_{\overline
\eta\to\tilde{s}}$ (cf. définition~\ref{17-definition-transition-cas-trait}).
\end{definition2}

\subsection{Cas général}

\begin{definition2}\label{17-defi-transition-codimension-un}
Soit $X$ un $\ZZn$-schéma excellent. Soit $\overline y\to\overline x$
une spécialisation immédiate de points géométriques de $X$. Pour définir le
morphisme de transition $\transition^X_{\overline y\to\overline x}\colon
\RG_{\overline y}(K)\to\RG_{\overline x}(K)(1)[2]$ pour tout
$K\in\detaleplus X$, quitte à remplacer $K$ par son image inverse \emph{via} le
morphisme canonique $X_{(\overline x)}\to X$, on peut supposer que $X$ est
local strictement hensélien (excellent) de point fermé $\overline x$. On
note alors $Z$ l'adhérence du point de $X$ en-dessous de $\overline y$ et
$i\colon Z\to X$ son immersion dans $X$. Ce $\ZZn$-schéma $Z$ est local
strictement hensélien intègre excellent de dimension $1$, le morphisme de
transition $\transition^{Z}_{\overline y\to\overline z}$ a été introduit dans la
définition~\ref{17-definition-transition-cas-local-dimension-un}.
Pour tout $K\in\detaleplus X$, on
définit le morphisme de transition $\transition^X_{\overline y\to\overline x}$ de
façon à faire commuter le diagramme suivant où les flèches verticales sont
les isomorphismes évidents :
\[
\xymatrix{\ar[d]^\sim \RG_{\overline y}(K)\ar[r]^-{\transition^X_{\overline
y\to\overline x}} & \RG_{\overline
x}(K)(1)[2]\ar[d]^\sim \\
\RG_{\overline y}(i^!K)\ar[r]^-{\transition^Z_{\overline
y\to\overline x}} & \RG_{\overline x}(i^!K)(1)[2]}
\]
\end{definition2}

\begin{remarque2}\label{17-remarque-composition-transition-codimension-0-ou-1}
Pour toute spécialisation $\overline x'\to\overline x$ de codimension $0$
entre points géométriques de $X$ (essentiellement, un élément du groupoïde
de Galois absolu du corps résiduel d'un des points de $X$), on a un
isomorphisme évident $\RG_{\overline x'}(K)\iso \RG_{\overline
x}(K)$, que l'on note $\transition^X_{\overline x'\to \overline x}$. Il est
évident que si $\overline z\to\overline y$ et $\overline y\to \overline x$
sont des spécialisations composables telles que la codimension $c$ de
$\overline z\to\overline x$ soit $0$ ou $1$, on a une égalité de morphismes
\[\transition^X_{\overline z\to\overline x}=\transition^X_{\overline y\to\overline
x}\circ \transition^X_{\overline z\to\overline y}\colon \RG_{\overline
z}(K)\to \RG_{\overline x}(K)(c)[2c]\;\text{.}\]
Ainsi, les morphismes de transition associés aux
spécialisations se composent bien dans l'étendue où ces constructions ont
été faites jusqu'à présent. Nous définirons plus tard des morphismes de
transition en codimension arbitraire et ce de façon compatible à la
composition, mais seulement au niveau des groupes de cohomologie
(cf. théorème~\ref{17-theoreme-transition-toute-codimension}).
\end{remarque2}

\section{Complexes dualisants putatifs et potentiels}
\label{17-section-premieres-proprietes}
\subsection{Définition des complexes dualisants putatifs et potentiels}

\begin{definition2}
\label{17-definition-fonction-de-dimension}
Une fonction de dimension sur un schéma localement
noethérien $X$ est une fonction $\delta\colon X\to\mathbf{Z}$ telle que
pour toute spécialisation immédiate $\overline y\to\overline x$ de points
géométriques de $X$, on ait $\delta(y)=\delta(x)+1$.
\end{definition2}

Localement pour la topologie de Zariski, un schéma excellent admet une
fonction de dimension, cf~\refext{XIV}{14-theoreme-univ-cat-fonction-dim}.

\begin{definition2}\label{17-definition-complexe-dualisant-potentiel}
Soit $X$ un $\ZZn$-schéma noethérien
excellent muni d'une fonction de dimension $\delta$. Un
complexe dualisant putatif consiste en la donnée de $K\in\detaleplus X$ et pour
tout $x\in X$ d'un isomorphisme (appelé épinglage en $x$)
$\SheafRG_x(K)\iso \Lambda(\delta(x))[2\delta(x)]$ dans $\detaleplus
x$. Un complexe dualisant putatif est un complexe dualisant potentiel si
pour toute spécialisation immédiate $\overline y\to\overline x$, le
diagramme suivant est commutatif :
\[\xymatrix{
\ar[rd]^\sim\RG_{\overline y}(K)
\ar[r]^-{\transition^X_{\overline y\to\overline x}}& 
\RG_{\overline x}(K)(1)[2]\ar[d]^\sim \\
& \Lambda(\delta(y))[2\delta(y)] 
}\]
Autrement dit, les épinglages sont compatibles aux morphismes de transition
associés aux spécialisations immédiates.
\end{definition2}

Les notions de complexes dualisants putatifs et potentiels (à coefficients
$\Lambda$) ne sont définis que pour les $\ZZn$-schémas noethériens
excellents munis d'une fonction de dimension. Certains des énoncés à venir
contiendront donc implicitement ces hypothèses sur les schémas. La fonction
de dimension ne sera pas non plus systématiquement mentionnée dans les
énoncés.

\begin{remarque2}
Si $X$ est connexe et que $K\in\detaleplus X$ est muni de deux 
structures de complexe dualisant
potentiel, pour vérifier que les épinglages sont les mêmes en tous les
points de $X$, il suffit de le faire en un seul point.
\end{remarque2}

L'objectif de cette section est de montrer que sur un schéma régulier
excellent muni de la fonction de dimension $-\codim$, le faisceau constant
$\Lambda$ est naturellement muni d'une structure de complexe dualisant
potentiel.

\subsection{Fonctorialité par rapport aux morphismes étales}

Des propriétés de stabilité importantes des complexes dualisants potentiels
par rapport à certaines classes de morphismes seront obtenues dans la
section~\ref{17-section-complements-cdp}. Pour le moment, mentionnons
simplement la compatibilité suivante pour les morphismes étales :

\begin{proposition2}\label{17-proposition-stabilite-cdp-morphisme-etale}
Soit $f\colon Y\to X$ un morphisme étale entre $\ZZn$-schémas noethériens
excellents. On suppose que $X$ est muni d'une fonction de dimension
$\delta_X$. On définit une fonction de dimension $\delta_Y$ sur $Y$ en
posant $\delta_Y(y)=\delta_X(f(y))$ pour tout $y\in Y$. Soit $K$ un
complexe dualisant putatif sur $X$. Alors, $f^\star K$ est naturellement
muni d'une structure de complexe dualisant putatif et c'est un complexe
dualisant potentiel si $K$ en est un.
\end{proposition2}

Si $y$ est un point de $Y$, que $x=f(y)$ et que l'on note $g\colon y\to x$
le morphisme induit par $f$, on a un isomorphisme canonique $g^\star
\SheafRG_x(K)\simeq \SheafRG_y(f^\star K)$. Ceci permet de
définir les épinglages sur $f^\star K$. On vérifie aussitôt que si $K$ est
un complexe dualisant potentiel, alors $f^\star K$ aussi.

\bigskip

La construction de cette proposition passe évidemment à la limite : le
résultat vaut aussi pour des localisations $X_{(x)}\to X$,
$X_{(x)}^\hensel\to X$ ou $X_{(\overline x)}\to X$. Nous utiliserons
librement ces observations simples dans la suite.

\subsection{Compléments sur les spécialisations}

\begin{definition2}\label{17-definition-specialisation-au-dessus-d-une-autre}
Soit $X'\to X$ un morphisme de schémas. Soit $\overline{y'}\to \overline{x'}$
et $\overline y\to \overline x$ des spécialisations de points géométriques 
de $X'$ et $X$ respectivement. Si l'on se donne des $X$-morphismes 
$\overline {y'}\to\overline y$ et $\overline{x'}\to \overline x$ tels que
le diagramme évident ci-dessous commute, alors on dit que 
$\overline{y'}\to\overline{x'}$
est au-dessus de $\overline y\to\overline x$.
\[\xymatrix{\overline {y'}\ar[d]\ar[r] & X'_{(\overline x')} \ar[d]\\
\overline y\ar[r] & X_{(\overline x)}}
\]
\end{definition2}

\begin{proposition2}\label{17-proposition-specialisation-going-down}
Soit $X'\to X$ un morphisme de schémas. Soit
$\overline {y'}\to\overline {x'}$ une
spécialisation de points géométriques de $X'$. Alors, à des isomorphismes
uniques près, il existe une unique spécialisation $\overline y\to\overline
x$ de points géométriques de $X$ en-dessous de $\overline {y'}\to\overline
{x'}$.
\end{proposition2}

C'est évident.

\begin{proposition2}\label{17-proposition-specialisation-going-up}
Soit $X'\to X$ un morphisme fini et $\overline y\to\overline x$ une
spécialisation de points géométriques de $X$. Soit $\overline
{y'}\to X'$ un point géométrique de $X'$ au-dessus de $\overline y$
(\emph{i.e.}
on s'est donné un $X$-morphisme $\overline {y'}\to\overline y$). Alors, à
des isomorphismes uniques près, il existe une unique spécialisation
de points géométriques de $X'$ au-dessus de $\overline y\to\overline x$ de
la forme $\overline {y'}\to\overline {x'}$.
\end{proposition2}

On peut supposer que $X$ est local strictement hensélien de point fermé
$\overline x$. Le schéma $X'$ étant fini sur $X$, c'est une réunion disjointe
finie de schémas locaux strictement henséliens. Quitte à remplacer $X'$ par
la composante connexe contenant $\overline {y'}$, on peut supposer que $X'$
est lui-aussi local strictement hensélien. Il n'y a alors manifestement
plus d'alternative : $\overline {x'}$ est le point fermé de $X'$.

\subsection{Construction d'un complexe dualisant potentiel dans le cas régulier}
\subsubsection{Un complexe dualisant putatif}
\begin{proposition3}\label{17-proposition-cdputatif-regulier}
Soit $X$ un $\ZZn$-schéma noethérien régulier excellent muni de la fonction de
dimension $-\codim$. Alors $\Lambda$ est naturellement muni d'une structure
de complexe dualisant putatif.
\end{proposition3}

Soit $x\in X$. Pour définir l'épinglage en $x$, on peut supposer que $X$
est local de point fermé $x$. On note $i\colon x\to X$ l'inclusion de ce
point fermé. Le morphisme de Gysin $\Cl i\colon
\Lambda(\delta(x))[2\delta(x)]\to i^!\Lambda$ est un isomorphisme d'après
le théorème de pureté cohomologique absolue
(cf.~\refext{XVI}{16-theoreme-purete-absolue}). L'épinglage en $x$ est
l'isomorphisme inverse.

\subsubsection{Cas d'un trait}

\begin{proposition3}\label{17-proposition-cdp-trait}
Soit $X$ un trait excellent, muni de la fonction de dimension
$\delta=-\codim$. On suppose que $X$ est un $\ZZn$-schéma.
\begin{enumerate}
\item[(a)] Le complexe dualisant putatif $\Lambda$ de la
proposition~\ref{17-proposition-cdputatif-regulier} est un complexe dualisant
potentiel.
\item[(b)] Si $K$ est un complexe dualisant putatif, il existe un unique
morphisme $\Lambda\to K$ compatible aux épinglages au point générique.
\item[(c)] Si $K$ est un complexe dualisant potentiel, le morphisme
$\Lambda\to K$ défini ci-dessus est un isomorphisme (compatible aux
épinglages).
\end{enumerate}
\end{proposition3}

La compatiblité (a) est facile (au moins au signe près). Établissons (b).
Soit $K$ un complexe dualisant putatif. Notons $i\colon s\to X$ l'inclusion
du point fermé $s$ de $X$ et $j\colon\eta\to X$ l'inclusion de son point
générique $\eta$. On a un triangle distingué canonique, que l'on peut
récrire en présence d'épinglages :
\[\xymatrix{\ar[d]^{\sim}i_\star i^! K\ar[r] & \ar@{=}[d] K \ar[r] &
\ar[d]^\sim \R j_\star j^\star
K\ar[r] & i_\star i^! K[1]\ar[d]^\sim \\
i_\star \Lambda(-1)[-2]\ar[r] & K \ar[r] & \R j_\star \Lambda\ar[r] & 
i_\star \Lambda(-1)[-1]}\]
En appliquant le foncteur $i$-ème faisceau de cohomologie ${\mathscr H}^i$,
on obtient l'annulation de ${\mathscr H}^i K$ pour $i<0$ et un 
isomorphisme canonique ${\mathscr H}^0K\simeq \Lambda$. En vertu des propriétés
élémentaires de la $t$-structure canonique sur $\detaleplus X$, on obtient un
unique morphisme $\Lambda\simeq {\mathscr H}^0 K\to K$ compatible à
l'épinglage au point générique.

Pour obtenir (c), supposons que $K$ soit un complexe dualisant potentiel. On
considère le diagramme commutatif induit par le morphisme canonique
$\Lambda\to K$ de (b) et la fonctorialité du morphisme de spécialisation
associé à un choix de  spécialisation $\overline \eta\to \overline s$ 
au-dessus des points $\eta$ et $s$ :
\[
\xymatrix{
\RG_{\overline\eta}(\Lambda)\ar[r]\ar[d]^{\transition^X_{\overline \eta\to\overline s}} & \RG_{\overline
\eta}(K)\ar[d]^{\transition^X_{\overline \eta\to\overline s}} \\
\RG_{\overline s}(\Lambda)(1)[2]\ar[r] & \RG_{\overline s}(K)(1)[2]}
\]
Comme $\Lambda$ et $K$ sont des complexes dualisants potentiels, les
flèches verticales sont des isomorphismes. Par construction, le morphisme
du haut est un isomorphisme. Il en résulte que le morphisme du bas aussi.
Par conséquent le morphisme $\Lambda\to K$ induit un isomorphisme après
application de $j^\star$ et de $i^!$ : c'est un isomorphisme.

\subsubsection{Fonctorialité par rapport aux morphismes quasi-finis}

\begin{proposition3}\label{17-proposition-cdputatif-quasi-fini}
Soit $f\colon Y\to X$ un morphisme quasi-fini. Soit $K$ un complexe dualisant
putatif sur $X$ pour une certaine fonction de dimension $\delta_X$ sur $X$.
Alors $f^!K$ est naturellement muni d'une structure de complexe dualisant
putatif pour la fonction de dimension $\delta_Y$ sur $Y$ définie par
$\delta_Y(y)=\delta_X(f(y))$ pour tout $y\in Y$.
\end{proposition3}

Soit $y\in Y$. Notons $x=f(y)$ et $\pi\colon y\to x$ le morphisme fini
induit par $f$. On a un isomorphisme canonique dans $\detaleplus y$ :
\[\SheafRG_y(f^!K)\simeq \pi^! \SheafRG_x(K)\;\text{.}\]
L'épinglage en $x$ donne un isomorphisme $\SheafRG_x(K)\simeq
\Lambda(\delta_X(x))[2\delta_X(x)]$. Pour obtenir l'isomorphisme voulu
$\SheafRG_y(f^!K)\simeq \Lambda(\delta_Y(y))[2\delta_Y(y)]$, il
suffit de définir un isomorphisme $\Lambda\iso \pi^!\Lambda$ : on
utilise le morphisme de Gysin $\Cl \pi$.

\begin{remarque3}\label{17-remarque-transitivite-construction-f-exclam-haut-quasi-fini}
Soit $g\colon Z\to Y$ un autre morphisme quasi-fini. \emph{Via}
l'isomorphisme de transitivité $g^!f^!\simeq (f\circ g)^!$, la structure de
complexe dualisant putatif sur $g^!(f^!K)$ obtenue en appliquant cette
construction à $f$ puis à $g$ est la même que celle obtenue en appliquant
directement la construction à $f\circ g$ : cela résulte aussitôt des
propriétés de composition des morphismes de Gysin.
\end{remarque3}

\begin{proposition3}\label{17-proposition-cdpotentiel-morphisme-fini}
Soit $f\colon Y\to X$ un morphisme fini et surjectif et $K$ un complexe
dualisant putatif sur $X$. Alors $K$ est un complexe dualisant potentiel si
et seulement si le complexe dualisant putatif $f^!K$ de la 
proposition~\ref{17-proposition-cdputatif-quasi-fini} est
un complexe dualisant potentiel.
\end{proposition3}

En vertu des propositions~\ref{17-proposition-specialisation-going-down} et
\ref{17-proposition-specialisation-going-up}, on peut supposer que $X$ et $Y$
sont des schémas locaux strictement henséliens intègres de dimension $1$
et que les fonctions de dimension prennent les valeurs $0$ et $-1$.
Compte tenu de la
remarque~\ref{17-remarque-transitivite-construction-f-exclam-haut-quasi-fini}, il
vient alors qu'il suffit de traiter deux cas :
\begin{enumerate}
\item[(1)] $Y$ est le normalisé de $X$ ;
\item[(2)] $X$ et $Y$ sont des traits.
\end{enumerate}

On obtient la conclusion dans le cas (1) en utilisant le fait que le
morphisme de transition pour $X$ et $K$ est défini à partir de celui pour
$Y$ et $f^\star K$ (cf.
définition~\ref{17-definition-transition-cas-local-dimension-un}) et que si
$M/L$ est une extension finie purement inséparable de corps,
le morphisme de Gysin
$\Cl \pi\colon \Lambda\to\pi^!\Lambda$ associé au morphisme $\pi\colon
\SP(M)\to\SP(L)$ s'identifie à la multiplication par le degré de $M/L$
\emph{via} les isomorphismes tautologiques $\pi^!\Lambda\simeq
\pi^\star\Lambda\simeq \Lambda$.

La démonstration dans le cas (2) va utiliser le lemme général suivant :

\begin{lemme3}\label{17-lemme-compatibilite-cdp-f-exclamation-schemas-reguliers}
Soit $f\colon Y\to X$ un morphisme quasi-fini entre $\ZZn$-schémas
noethériens excellents réguliers de
dimension relative virtuelle $-c$. On munit $X$ de la fonction de dimension
$\delta_X=-\codim$ et $Y$ de la fonction de dimension $\delta_Y$
définie par composition avec $f$ comme dans la 
proposition~\ref{17-proposition-cdputatif-quasi-fini}.
Compte tenu de la relation
$\delta_Y(y)=-\codim_Y(y)-c$ pour tout $y\in Y$,
la proposition~\ref{17-proposition-cdputatif-regulier} munit 
$\Lambda(-c)[-2c]$ d'une structure de complexe dualisant putatif pour la
fonction de dimension $\delta_Y$ sur $Y$. La
proposition~\ref{17-proposition-cdputatif-regulier} donne une structure de
complexe dualisant putatif sur $\Lambda$ sur $X$, dont on déduit, par la
proposition~\ref{17-proposition-cdputatif-quasi-fini}, une structure de
complexe dualisant putatif sur $f^!\Lambda$ sur $Y$ pour la fonction de
dimension $\delta_Y$. Alors, l'isomorphisme de pureté $\Cl
f\colon \Lambda(-c)[-2c]\iso f^!\Lambda$ est compatible aux épinglages
de ces deux complexes dualisants putatifs.
\end{lemme3}

Étudions les épinglages en un point $y\in Y$. Notons $x=f(y)$. On peut
supposer que $X$ et $Y$ sont locaux (strictement henséliens) de points
fermés respectifs $x$ et $y$. On a un diagramme commutatif :
\[
\xymatrix{y\ar[r]^j\ar[d]_g \ar[rd]^h& Y\ar[d]^f \\ x\ar[r]^i & X}
\]
Les schémas apparaissant sur ce diagramme sont affines et réguliers et 
les morphismes entre eux sont de type fini. Ces morphismes sont donc
localement d'intersection complète (lissifiables), on peut leur appliquer
la théorie des morphismes de Gysin. Le résultat du lemme découle alors
aussitôt de leur compatibilité à la composition, puisqu'elle
donne un diagramme commutatif dans $\detaleplus y$, où l'on a
noté $c'$ la codimension de $y$ dans $Y$ :
\[
\xymatrix{\Lambda\ar[r]^{\Cl g} \ar[d]^{\Cl j} & g^!\Lambda
\ar[d]^{g^!(\Cl i)} \\
j^!\Lambda(c')[2c']\ar[r]^-{j^!(\Cl f)}&h^!\Lambda(c+c')[2c+2c']
}
\]

\bigskip

Revenons à la démonstration de la
proposition~\ref{17-proposition-cdpotentiel-morphisme-fini}. Supposons que $K$
soit un complexe dualisant potentiel. D'après la
proposition~\ref{17-proposition-cdp-trait} appliquée à $X$, $K$ est
canoniquement isomorphe à $\Lambda$ (avec les épinglages de la
proposition~\ref{17-proposition-cdputatif-regulier}). On suppose donc
$K=\Lambda$. On a un isomorphisme de pureté
$\Cl f\colon \Lambda\iso f^!\Lambda$.
D'après le lemme~\ref{17-lemme-compatibilite-cdp-f-exclamation-schemas-reguliers}, cet isomorphisme est compatible avec les épinglages de la
proposition~\ref{17-proposition-cdputatif-regulier} pour $\Lambda$ et ceux de
la proposition~\ref{17-proposition-cdputatif-quasi-fini} pour $f^!\Lambda$.
D'après la proposition~\ref{17-proposition-cdp-trait} appliquée à $Y$,
il vient que
$f^!\Lambda$ est bien un complexe dualisant potentiel.

Inversement, supposons que $f^!K$ soit un complexe dualisant potentiel. La
proposition~\ref{17-proposition-cdp-trait} pour $X$
donne un morphisme canonique
$\Lambda\to K$ compatible aux épinglages au point générique de $X$.
Le morphisme $f^!\Lambda\to f^!K$ qui s'en déduit est un morphisme entre
deux complexes dualisants potentiels compatible aux épinglages au point
générique de $Y$.
D'après la proposition~\ref{17-proposition-cdp-trait},
$f^!\Lambda\to f^!K$ est un isomorphisme. Le foncteur $f^!$ étant
conservatif, il en résulte que le morphisme canonique $\Lambda\to K$ est un
isomorphisme. On peut donc supposer que $K=\Lambda$ (de façon compatible
aux épinglages au 
point générique). 
L'épinglage de $K$ au point fermé $x$ ne peut qu'être de la forme
$\lambda\cdot \Cl {x\subset X}^{-1}\colon \SheafRG_x(\Lambda)\iso
\Lambda(-1)[-2]$ pour $\lambda\in\Lambda^\times$. Compte tenu de ce qui
précède, il vient aussitôt
que le complexe dualisant putatif $f^!K$ sur $Y$ s'identifie à $\Lambda$,
épinglé trivialement au point générique et par $\lambda\cdot \Cl {y\subset
Y}^{-1}$ au point fermé $y$ de $Y$. Compte tenu de la
proposition~\ref{17-proposition-cdp-trait}, ce complexe dualisant putatif sur
$Y$ ne peut évidemment être un complexe dualisant potentiel que si
$\lambda=1$. Ainsi, $K$ est bien un complexe dualisant potentiel sur $X$,
puisqu'il s'identifie à $\Lambda$ de façon compatible aux épinglages.

\subsubsection{Un complexe dualisant potentiel}

\begin{proposition3}\label{17-proposition-cdpotentiel-regulier}
Soit $X$ un schéma régulier excellent, muni de la fonction de dimension
$-\codim$. Le complexe dualisant putatif $\Lambda$ de la
proposition~\ref{17-proposition-cdputatif-regulier} est un complexe dualisant
potentiel.
\end{proposition3}

Soit $\overline y\to\overline x$ une spécialisation immédiate de points
géométriques de $X$. On veut montrer que les épinglages sur $\Lambda$ sont
compatibles au morphisme de transition $\transition^X_{\overline y\to\overline
x}$. Pour cela, on peut supposer que $X$ est local strictement hensélien de
point fermé $\overline x$. Soit $C$ l'adhérence de l'image de $\overline
y$ dans $X$. Notons $i\colon C\to X$ l'immersion fermée de $C$ dans $X$.
Soit $n\colon \tilde C\to C$ le normalisé de $C$. Montrer la compatibilité
des épinglages avec le morphisme de transition associé à $\overline y\to
\overline x$ revient à montrer que le complexe dualisant putatif
$i^!\Lambda$ de la proposition~\ref{17-proposition-cdputatif-quasi-fini} est
un complexe dualisant potentiel, ce qui, d'après la
proposition~\ref{17-proposition-cdpotentiel-morphisme-fini},
revient encore à dire que
$n^!i^!\Lambda$ en est un. D'après le 
lemme~\ref{17-lemme-compatibilite-cdp-f-exclamation-schemas-reguliers},
le complexe dualisant
putatif $n^!i^!\Lambda$ s'identifie au complexe dualisant putatif
$\Lambda(-c)[-2c]$ obtenu par torsion et décalage à partir de $\Lambda$ sur
le trait $\tilde C$. En vertu de la
proposition~\ref{17-proposition-cdp-trait}, il s'agit bien d'un complexe
dualisant potentiel, ce qui achève la démonstration.

\section{Morphismes de transition généraux et classe de cohomologie en degré maximal}
\label{17-section-classe-de-cohomologie}
\subsection{Énoncés des théorèmes principaux}

\begin{theoreme2}\label{17-theoreme-classe-degre-maximal}
Soit $X$ un $\ZZn$-schéma local normal strictement 
hensélien excellent de dimension $d$ et de point fermé $x$.
Alors, $\HH^q_x(X,\Lambda(d))=0$ pour $q>2d$ et on a un isomorphisme 
$[x]\colon \Lambda\iso \HH^{2d}_x(X,\Lambda(d))$
compatible aux morphismes de transition associés aux spécialisations
immédiates.
\end{theoreme2}

\begin{theoreme2}\label{17-theoreme-transition-toute-codimension}
Soit $X$ un $\ZZn$-schéma noethérien
excellent. Pour toute spécialisation $\overline
y\to\overline x$ de points géométriques de $X$ de codimension $c$, on peut
définir, pour tout $K\in\detaleplus X$ et $i\in\mathbf{Z}$, un morphisme de
transition $\transition^X_{\overline y\to\overline x}\colon \HH^i_{\overline
y}(K)\to \HH^{i+2c}_{\overline x}(K(c))$, compatible à la composition des
spécialisations et induit par les définitions de la 
section~\ref{17-section-transition-codimension-un} pour $c\leq 1$. Par
ailleurs, ces morphismes de transition généralisés vérifient une
compatibilité avec les morphismes finis, énoncée dans la
proposition~\ref{17-proposition-compatibilite-transition-specialisation}.
\end{theoreme2}

\begin{remarque2}
Localement pour la topologie étale, un schéma quasi-excellent est
excellent. Il est donc évident que dans l'énoncé du
théorème~\ref{17-theoreme-transition-toute-codimension}, on peut remplacer
l'hypothèse \guil{excellent} par \guil{quasi-excellent}.
\end{remarque2}

\subsection{Dimension cohomologique}

On énonce ici deux résultats de dimension cohomologique qui se déduisent du
théorème de Lefschetz affine de Gabber,
cf.~\refext{XV}{15-prop-principale-locale}.

\begin{proposition2}\label{17-proposition-dimension-cohomologique}
Soit $\ell$ un nombre premier.
Soit $X$ un $\ZZn[\ell]$-schéma local strictement hensélien excellent 
de point fermé
$x$. On note $U=X-x$. Soit $\mathscr M$ un faisceau de
$\mathbf{Z}/\ell\mathbf{Z}$-modules
sur $U$. Soit $d\in\mathbf{N}$. On suppose que si $u\in U$ est tel que
${\mathscr M}_{\overline u}\neq 0$, alors la dimension de l'adhérence de
$u$ dans $X$ est $\leq d$. Alors $\HH^q(U,\mathscr M)=0$ pour $q\geq 2d$.

En particulier, si $d=\dim X$, on a $\cd_\ell U\leq 2d-1$ et $\cd_\ell
\Gamma_x\leq 2d$.
\end{proposition2}

Par passage à la limite inductive filtrante sur les sous-faisceaux
constructibles de $\mathscr M$, on peut supposer que $\mathscr M$ est
constructible. Quitte à remplacer $X$ par 
l'adhérence dans $X$ du support de $\mathscr M$, on peut supposer que
$d=\dim X$. On est alors ramené à montrer que $\cd_\ell U\leq 2d-1$.
Comme le schéma (séparé)
$U$ peut être recouvert par $d$ ouverts affines (cf.
\cite[§B.3, Chapitre~III]{Algebre@Serre}), et qu'un
ouvert affine de $X$ est de dimension cohomologique au plus 
$d$, en utilisant convenablement les suites exactes de Mayer-Vietoris,
on obtient bien que $\cd_\ell U\leq 2d-1$.

\begin{proposition2}\label{17-proposition-dimension-cohomologique-galois}
Soit $\ell$ un nombre premier.
Soit $X$ un $\ZZn[\ell]$-schéma local strictement hensélien intègre excellent
de dimension $d$, de point générique $\eta$. Alors $\cd_\ell\eta\leq d$.
\end{proposition2}

Le point générique $\eta$ s'identifie à limite projective du système
projectif formé par les ouverts affines non vides de $X$. Comme chacun de
ces ouverts affines est de $\ell$-dimension cohomologique au plus $d$,
il en va de même pour $\eta$.

\subsection{Morphismes finis}

\begin{proposition2}\label{17-proposition-compatibilite-transition-spec-immediate}
Soit $f\colon X'\to X$ un morphisme fini entre $\ZZn$-schémas noethériens
excellents. Soit $\overline{y'}\to\overline{x'}$ une spécialisation
immédiate de points géométriques de $X'$ au-dessus d'une spécialisation
immédiate $\overline y\to\overline x$ de points géométriques de $X$
(cf. définition~\ref{17-definition-specialisation-au-dessus-d-une-autre}).
Pour tout $K\in\detaleplus X$ et $i\in\mathbf{Z}$, le diagramme suivant 
commute :
\[
\xymatrix{\HH^i_{\overline{y'}}(f^!K)
     \ar[r]^-{\transition^{X'}_{\overline{y'}\to\overline{x'}}}
  & \HH^{i+2}_{\overline{x'}}(f^!K(1)) \\
\HH^i_{\overline y}(K)
     \ar[r]^-{\transition^{X}_{\overline{y}\to\overline{x}}}
     \ar[u]_{\sim}^{\Cl {y'\to y}}
  & \HH^{i+2}_{\overline{x}}(K(1))\ar[u]^{\sim}_{\Cl {x'\to x}}}
\]
\end{proposition2}

On peut supposer que $f$ est un morphisme dominant entre schémas locaux
strictement henséliens intègres de dimension $1$. Comme dans la
démonstration de la
proposition~\ref{17-proposition-cdpotentiel-morphisme-fini}, il y a deux cas à
traiter :
\begin{enumerate}
\item[(1)] $X'$ est le normalisé de $X$ ;
\item[(2)] $X$ et $X'$ sont des traits.
\end{enumerate}

Le cas (1) étant trivial, on se concentre sur le cas où $X$ et $X'$ sont
des traits. Pour vérifier la compatibilité, on peut évidemment supposer que
$i=0$. Pour tout $K\in\detaleplus X$, comme les flèches verticales sont des
isomorphismes, on peut noter 
$\delta_K\colon \HH^0_{\overline y}(K)\to \HH^2_{\overline x}(K(1))$ la
différence des flèches obtenues en suivant les deux chemins possibles. Il
s'agit de monter que pour tout $K\in\detaleplus X$, on a $\delta_K=0$.

Les résultats de la section~\ref{17-section-premieres-proprietes} montrent que
$\delta_\Lambda=0$. La suite de la démonstration va consister à se ramener
à ce cas-là.

Comme $\overline y$ est au-dessus du
point générique de $X$, $\HH^0_{\overline y}$ s'identifie à la fibre en
$\overline y$ du faisceau de cohomologie de $K$ en degré $0$. Si on note
$\tau_{\leq 0} K\to K$ le morphisme canonique déduit de la $t$-structure
canonique sur $\detaleplus X$\;\footnote{L'auteur pense qu'il eût été plus 
cohérent de noter cette troncature $\tau^{\leq 0}$ pour respecter 
la convention qui veut que les degrés cohomologiques soient indiqués
en exposant, mais la tradition ayant consacré l'usage inverse,
il s'y plie avec répugnance.}, on obtient un carré commutatif :
\[
\xymatrix{
\ar[d]^{\sim}\HH^0_{\overline y}(\tau_{\leq 0}K)\ar[r]^{\delta_{\tau_{\leq
0}K}} &\HH^2_{\overline
x}(\tau_{\leq 0}K(1)) \ar[d]\\
\HH^0_{\overline y}(K)\ar[r]^{\delta_K}& \HH^2_{\overline x}(K(1))
}
\]
Il résulte de ce diagramme que si $\delta_{\tau_{\leq 0}K}=0$, alors
$\delta_K=0$.

Notons ${\mathscr H}^0K$ le faisceau de cohomologie de $K$ en degré zéro.
Pour des raisons de dimension cohomologique, le morphisme canonique
${\mathscr H}^0K\to \tau_{\leq 0}K$ induit un isomorphisme après
application du foncteur $\HH^2_{\overline x}(-(1))$ (et aussi du foncteur
$\HH^0_{\overline y}$). Par conséquent, en considérant un carré commutatif du
type précédent, on obtient cette fois-ci que $\delta_{{\mathscr H}^0K}=0$
équivaut à $\delta_{\tau_{\leq 0}K}=0$.

Il résulte de ces remarques que pour montrer que
$\delta_K=0$  pour tout $K\in\detaleplus X$, on peut supposer que $K$ est
concentré en degré $0$.

On suppose maintenant que $K=\mathscr M$ où $\mathscr M$ est un faisceau de
$\Lambda$-modules sur $X$. Notons $j\colon y\to X$ l'inclusion du point
générique de $X$. En utilisant la description du foncteur $\HH^2_{\overline x}$
de la proposition~\ref{17-proposition-description-h-2-support-trait}, on
obtient que le morphisme canonique $j_!j^\star \mathscr M\to\mathscr M$
induit un isomorphisme après application de $\HH^0_{\overline y}$ et
$\HH^2_{\overline x}(-(1))$ : $\delta_{K}=0$ équivaut à $\delta_{j_!j^\star
K}=0$. La propriété, pour un faisceau $K$ sur $X$, d'être tel que
$\delta_K=0$ ne dépend donc que de sa restriction au point générique
$y$.

Pour tout faisceau de $\Lambda$-modules $\mathscr L$ sur $y$, on
note $\delta_{\mathscr L}=\delta_{j_!\mathscr L}$. Il s'agit de montrer que
$\delta_{\mathscr L}=0$ pour tout faisceau de $\Lambda$-modules sur $y$.
Notons $\mathscr L\to {\mathscr L}_{\constant}$ le plus grand
quotient constant de $\mathscr L$ : formellement, le foncteur $\mathscr
L\longmapsto {\mathscr
L}_{\constant}$ est le foncteur adjoint à gauche du foncteur d'inclusion
de la catégorie des faisceaux constants de $\Lambda$-modules dans la
catégorie des faisceaux de $\Lambda$-modules sur $y$.
En utilisant la
proposition~\ref{17-proposition-description-h-2-support-trait}, on montre
facilement que le morphisme canonique ${\mathscr L}\to{\mathscr
L}_{\constant}$ induit un isomorphisme après application de
$\HH^2_{\overline x}(j_!(1))$ (et une surjection après application
de $\HH^0_{\overline y}$). D'où $\delta_{\mathscr
L}=0$ si et seulement si $\delta_{{\mathscr L}_{\constant}}=0$.
Il en résulte que l'on peut supposer que $\mathscr L$ est un faisceau
constant de $\Lambda$-modules. Par ailleurs, on peut évidemment supposer
que $\mathscr L$ est constructible. Si ${\mathscr L}\to{\mathscr L}'$ est
un épimorphisme de faisceaux constants constructibles de $\Lambda$-modules,
il vient aussitôt que $\delta_{\mathscr L}=0$ implique $\delta_{{\mathscr
L}'}=0$. Il en résulte que l'on peut supposer que $\mathscr L=\Lambda^r$
pour un certain $r\in\mathbf{N}$, puis, par additivité, que $r=1$. Bref, on
s'est bien ramené au cas du faisceau constant $\Lambda$, ce qui achève la
démonstration de cette proposition.

\subsection{Le cas de la \texorpdfstring{dimension~$2$}{dimension 2}}
\label{17-subsection-le-cas-de-la-dimension-2}

Compte tenu des résultats établis jusqu'à présent, les
théorèmes~\ref{17-theoreme-classe-degre-maximal} et
\ref{17-theoreme-transition-toute-codimension} peuvent être considérés comme
ayant été établis en dimension $0$ et $1$. Cette sous-section se concentre
sur le cas de la dimension~$2$ qui est l'étape cruciale pour passer au cas
général.

On fixe un $\ZZn$-schéma local strictement
hensélien excellent normal $X$ de dimension~$2$. On note $x$ le point fermé
de $X$ et on pose $U=X-x$. En s'aidant d'une résolution des singularités
$X'\to X$, nous allons définir au
paragraphe~\ref{17-subsubsection-classe-construction} une classe dans
$\HH^4_x(X,\Lambda(2))$, puis nous démontrerons au
paragraphe~\ref{17-subsubsection-classe-independance} qu'elle est indépendante
de la résolution. Enfin, nous établirons au
paragraphe~\ref{17-subsubsection-classe-transition} une compatibilité entre
cette classe, les classes définies pour les localisés stricts en $X$ en
les points de codimension~$1$ et les morphismes de transition associés aux
spécialisations immédiates correspondantes.

\subsubsection{Construction d'une classe}
\label{17-subsubsection-classe-construction}

D'après le résultat principal de 
\cite{D2DS@Lipman}, il existe un morphisme birationnel $p\colon X'\to X$ avec
$X'$ régulier. Le morphisme $p$ induit alors automatiquement un
isomorphisme $p^{-1}(U)\iso U$. On peut
supposer de plus que $p^{-1}(x)_{\red}$ est un diviseur à croisements
normaux stricts (dont on note $D_1,\dots,D_n$ les composantes irréductibles).

En effet, on peut dans un premier temps supposer que les composantes
irréductibles $D_1,\dots,D_n$ sont de dimension~$1$ : si ce n'est pas le
cas, on peut éclater le point fermé correspondant (ceci ne se peut produire
que si $X$ est déjà régulier et que $p$ est un isomorphisme). Dans un
deuxième temps, on peut s'arranger pour que les composantes irréductibles
de $p^{-1}(x)_{\red}$ soient régulières en itérant le processus consistant
à éclater les points singuliers (cf. \cite[page~38]{TIFR@Shafarevich}). Dans un
dernier temps, on peut contraindre les croisements à devenir normaux en
éclatant les points fermés récalcitrants ; comme des invariants numériques
décroissent strictement dans cette opération (cf.
\cite[page~21]{TIFR@Shafarevich}), ce processus termine.

\begin{proposition3}\label{17-proposition-classe-sorites}
\begin{enumerate}
\item[(a)] Le morphisme de bord $\HH^{q-1}(U,\Lambda(2))\to
\HH^q_x(X,\Lambda(2))$ est un isomorphisme pour $q\geq 2$ et ces groupes sont
nuls pour $q\geq 5$ ;
\item[(b)] Le morphisme évident $\HH^q(X',\Lambda(2))\to
\HH^q(p^{-1}(x),\Lambda(2))$ est un isomorphisme pour tout $q\in\mathbf{Z}$
et ces groupes sont nuls pour $q\geq 3$ ;
\item[(c)] Le morphisme de bord $\HH^{q-1}(U,\Lambda(2))\to
\HH^q_{p^{-1}(x)}(X',\Lambda(2))$ est un isomorphisme pour $q\geq 4$ ;
\item[(d)] Le morphisme évident $\HH^q_x(X,\Lambda(2))\to
\HH^q_{p^{-1}(x)}(X',\Lambda(2))$ est un isomorphisme pour $q\geq 4$ et ces
groupes sont nuls pour $q\geq 5$.
\end{enumerate}
\end{proposition3}

(a) s'obtient en utilisant la suite exacte canonique :
\[\HH^{q-1}(X,\Lambda(2))\to \HH^{q-1}(U,\Lambda(2))\to \HH^q_x(X,\Lambda(2))\to
\HH^q(X,\Lambda(2))\;\text{.}\]
En effet, $X$ étant local strictement hensélien, on a $\HH^i(X,\Lambda(2))$
pour $i>0$. On conclut en utilisant le fait que $\cd_\ell U\leq 3$ (cf.
proposition~\ref{17-proposition-dimension-cohomologique}).

(b) résulte du théorème de changement de base pour un morphisme propre et
du fait que $p^{-1}(x)$ soit une courbe propre et donc de dimension
cohomologique $2$.

(c) se déduit aussitôt de (b).

(d) résulte de (a), (c) et de la commutativité du diagramme évident
suivant :
\[
\xymatrix{\HH^{q-1}(U,\Lambda(2))\ar[r]\ar@{=}[d] & \HH^q_x(X,\Lambda(2))\ar[d] \\
\HH^{q-1}(U,\Lambda(2))\ar[r] & \HH^q_{p^{-1}(x)}(X',\Lambda(2))}
\]

\begin{proposition3}\label{17-proposition-classe-suite-exacte-1}
On a une suite exacte 
\[
\bigoplus_{i<j} \HH^4_{D_i\cap D_j}(X',\Lambda(2))\to 
\bigoplus_i \HH^4_{D_i}(X',\Lambda(2))\to
\HH^4_{p^{-1}(x)}(X',\Lambda(2))\to 0\;\text{,}
\]
où les flèches sont induites par des morphismes d'agrandissement du
support, et leurs différences.
\end{proposition3}

Pour tout fermé $F$ de $X'$, on note $\Lambda_F$ le faisceau
$i_{F\star}\Lambda$
où $i_F$ est l'immersion de $F$ dans $X'$. On a une suite exacte courte
évidente de faisceaux de $\Lambda$-modules sur $X'$ :
\[0\to \bigoplus_{i<j} \Lambda_{D_i\cap D_j}\to 
\bigoplus_i \Lambda_{D_i}\to \Lambda_{p^{-1}(x)}\to 0\;\text{.}\]

En appliquant $\R\SheafHom(-,\Lambda(2))$ au triangle distingué de
$\detaleplus {X'}$ associé à cette suite exacte courte, on peut obtenir la
suite exacte voulue, pourvu que l'on sache montrer que $\HH^5_{D_i\cap
D_j}(X',\Lambda(2))$ est nul pour $i<j$.
D'après le théorème de pureté absolue, ce
groupe s'identifie à $\HH^1(D_i\cap D_j,\Lambda)$ qui est bien nul puisque
$D_i\cap D_j$ est une union disjointe finie de spectres de corps
séparablement clos.

\begin{proposition3}\label{17-proposition-classe-suite-exacte-2}
La suite exacte de la proposition~\ref{17-proposition-classe-suite-exacte-1}
se récrit sous la forme :
\[\bigoplus_{i<j}\bigoplus_{y\in D_i\cap D_j}\Lambda\vers {\delta}
\Lambda^n \vers\varepsilon 
\HH^4_{p^{-1}(x)}(X',\Lambda(2))\to 0\;\text{.}\]
Si on note $(\chi_y)_{i<j,y\in D_i\cap D_j}$ la base canonique du groupe de
gauche et $\chi_{D_1},\dots,\chi_{D_n}$ celle de $\Lambda^n$, la
différentielle $\delta$ vérifie la formule
\[\delta(\chi_y)=[y:x]\cdot \left(\chi_{D_i}-\chi_{D_j}\right)\;\text{.}\]
En outre, $\varepsilon(\chi_{D_i})\in \HH^4_{p^{-1}(x)}(X',\Lambda(2))$ est
obtenu par agrandissement du support à partir de la classe
$\frac{1}{[y:x]}\Cl{y\to X}\in \HH^4_y(X',\Lambda(2))$
pour tout point fermé $y$ d'un $D_i$.
\end{proposition3}

Le théorème de pureté cohomologique absolue donne des isomorphismes
\[\Cl{D_i\cap D_j\to X}\colon \HH^0(D_i\cap D_j,\Lambda)\iso \HH^4_{D_i\cap
D_j}(X',\Lambda(2))\;\text{,}\]
ce qui permet de décrire le groupe de gauche dans la suite exacte de la
proposition~\ref{17-proposition-classe-suite-exacte-1}.
Le même théorème donne aussi des isomorphismes \[\Cl{D_i\to X}\colon
\HH^2(D_i,\Lambda(1))\iso \HH^4_{D_i}(X',\Lambda(2))\;\text{.}\]
En outre, on dispose
du morphisme trace (relativement à $x$) $\tr\colon \HH^2(D_i,\Lambda(1))
\iso \Lambda$. Celui-ci est caractérisé par le fait que pour tout point fermé
$y$ de $D_i$, l'image de $\Cl{y\to D_i}\in \HH^2_y(D_i,\Lambda(1))$ par le
morphisme composé \[\HH^2_y(D_i,\Lambda(1))\to \HH^2(D_i,\Lambda(1))\vers{\tr}
\Lambda\]
soit $[y:x]$.
La description des morphismes $\delta$ et $\varepsilon$ donnée dans
l'énoncé s'obtient alors aussitôt à partir des propriétés de composition des
morphismes de Gysin.

\begin{proposition3}\label{17-proposition-classe-suite-exacte-3}
La suite suivante est exacte :
\[\bigoplus_{i<j}\bigoplus_{y\in D_i\cap D_j}\Lambda\vers {\delta}
\bigoplus_i \Lambda \vers\Sigma \Lambda\to 0\;\text{,}\]
où $\Sigma$ est défini par $\Sigma(\chi_{D_i})=1$ pour tout $i$.
\end{proposition3}

Tout d'abord, à partir de la formule donnée pour $\delta$ dans la 
proposition~\ref{17-proposition-classe-suite-exacte-2}, il est clair que
$\Sigma\circ \delta=0$. $\Sigma$ induit donc un morphisme
$\coker \delta\to \Lambda$.
Comme les degrés $[y:x]$ qui
interviennent sont inversibles dans $\Lambda$, on voit que si $i$ et $j$
sont tels que $D_i\cap D_j$ soit non vide, alors $\chi_{D_i}$ et
$\chi_{D_j}$ ont la même classe dans $\coker \delta$. La fibre $p^{-1}(x)$
étant connexe (Main Theorem), on en déduit que tous les éléments
$\chi_{D_i}$ ont la même classe dans $\coker \delta$. L'inclusion entre
l'image de $\delta$ et le noyau de $\Sigma$ est donc une égalité.

\begin{corollaire3}
Les suites exactes des propositions~\ref{17-proposition-classe-suite-exacte-2}
et \ref{17-proposition-classe-suite-exacte-3} donnent naissance à un
isomorphisme $\Lambda\iso \HH^4_{p^{-1}(x)}(X';\Lambda(2))$.
\end{corollaire3}

\begin{corollaire3}
\emph{Via} l'isomorphisme canonique $\HH^4_{x}(X,\Lambda(2))\iso
\HH^4_{p^{-1}(x)}(X',\Lambda(2))$ de la
proposition~\ref{17-proposition-classe-sorites}~(d), l'isomorphisme du
corollaire précédent donne un isomorphisme $\Lambda\iso
\HH^4_x(X,\Lambda(2))$.
\end{corollaire3}

\begin{definition3}
On note $[x]_{X'}\in \HH^4_x(X,\Lambda(2))$ le générateur défini par
l'isomorphisme du corollaire précédent.
\end{definition3}

\subsubsection{Indépendance en la résolution}
\label{17-subsubsection-classe-independance}

Dans ce paragraphe, nous allons montrer que si $q\colon X''\to X$ est une
autre résolution du type envisagé dans le
paragraphe~\ref{17-subsubsection-classe-construction}, alors
$[x]_{X'}=[x]_{X''}$. Quitte à introduire une désingularisation de la
composante irréductible dominant $X$ du produit
fibré de $X''$ de $X'$ au-dessus de $X$, on peut supposer qu'une des deux
désingularisations $X'$ et $X''$ considérées coiffe l'autre. On suppose
donc par exemple qu'il existe un (unique) $X$-morphisme $\pi\colon X''\to
X'$. Le morphisme $\pi$ est projectif et birationnel entre deux schémas
réguliers de dimension $2$, il induit donc un isomorphisme au-dessus d'un
ouvert $U'$ de $X'$ tel que le fermé $X'-U'$ soit de dimension $0$. Il
existe donc certainement un point fermé $y'$ de $p^{-1}(x)$ tel que le
morphisme induit soit un isomorphisme $\pi^{-1}(y')\to y'$. On note $y''$
l'unique point fermé de $\pi^{-1}(y')$. La compatibilité des classes de
Gysin au changement de base implique que la classe $\Cl{y'\to X'}$ est
envoyée sur $\Cl{y''\to X''}$ par le morphisme de restriction
$\pi^\star\colon \HH^4_{y'}(X',\Lambda(2))\to \HH^4_{y''}(X'',\Lambda(2))$.
Cette compatibilité vaut encore après élargissement du support à
$p^{-1}(x)$ et à $q^{-1}(x)$. En considérant la composition suivante :
\[\HH^4_{x}(X,\Lambda(2))\iso \HH^4_{p^{-1}(x)}(X',\Lambda(2))\iso
\HH^4_{q^{-1}(x)}(X'',\Lambda(2))\;\text{,} \]
on obtient aussitôt que $[y':x]\cdot [x]_{X'}=[y'':x]\cdot [x]_{X''}$, ce
qui permet de conclure que $[x]_{X'}=[x]_{X''}$.

\begin{remarque3}
On peut montrer que si $p\colon X'\to X$ est un morphisme projectif
birationnel avec $X'$ régulier, alors, même sans supposer que
$p^{-1}(x)_\red$ soit un diviseur à croisements normaux strict, 
l'application $\HH^4_x(X,\Lambda(2))\to
\HH^4_{p^{-1}(x)}(X',\Lambda(2))$ est un isomorphisme et la classe $[x]$
s'envoie bien sur l'élément induit par 
$\frac 1 {[x':x]}\Cl{x'\to X'}$ pour tout
point fermé $x'$ de $p^{-1}(x)$.
\end{remarque3}

\subsubsection{Compatibilité aux morphismes de transition}
\label{17-subsubsection-classe-transition}

Soit $\overline y\to X$ un point géométrique au-dessus d'un point $y$ de
codimension $1$. Autrement dit, on a une spécialisation immédiate
$\overline y\to x$ de points géométriques de $X$. Nous allons montrer que
si $\overline\eta\to\overline y$ est une spécialisation immédiate, alors on
a l'égalité
\[ [x]=\transition^X_{\overline y\to x}(\transition^X_{\overline\eta\to\overline x}(1))\]
dans $\HH^4_x(X,\Lambda(2))$, ce qui achèvera la démonstration du
théorème~\ref{17-theoreme-classe-degre-maximal} jusqu'en dimension~$2$.

On note $C$ l'adhérence de $y$ dans $X$. Il existe une désingularisation
$X'\to X$ du type envisagé dans le 
paragraphe~\ref{17-subsubsection-classe-construction} telle que l'adhérence
$C'$ de $y$ dans $X'$ soit un trait. Plus précisément, le morphisme évident
$C'\to C$ identifie $C'$ au normalisé de $C$. Le morphisme $C'\to C$ étant
un homéomorphisme universel, on peut noter $x'$
le point fermé de $C'$ (celui de
$C$ est bien entendu $x$) et $U'=X'-x'$.
On dispose d'immersions fermées évidentes
$i\colon C\to X$, $i'\colon C'\to X'$,
$k\colon y\to U$ et $k'\colon y\to U'$.

\begin{lemme3}\label{17-lemme-hexagone-commutatif}
Avec les notations ci-dessus, le diagramme évident qui suit
est commutatif (les
flèches marquées comme étant des isomorphismes devant être considérées
comme bidirectionnelles) :
\[
\xymatrix{
&\ar[ld]^{\sim} \HH^4_x(X,\Lambda(2))\ar@{=}[r] & \HH^4_x(C,i^!\Lambda(2)) & \\ 
\HH^4_{p^{-1}(x)}(X',\Lambda(2)) & & & \HH^3(y,k^!\Lambda(2))
\ar[ul]^{\sim}\ar[ld]\\
& \ar[lu]^{\sim}
\HH^4_{x'}(X',\Lambda(2))\ar@{=}[r] & \HH^4_{x'}(C',{i'}^!\Lambda(2)) & \\
}
\]
\end{lemme3}

Ceci équivaut à la commutativité du diagramme évident :
\[
\xymatrix{
\HH^3(U,\Lambda(2))\ar[r]^\sim & \HH^4_x(X,\Lambda(2))\ar@{=}[r] &
   \HH^4_x(C,i^!\Lambda(2)) & \ar[l] \HH^3(y,k^!\Lambda(2))\ar@{=}[d] \\
\ar[u]_{\sim}\HH^3(U',\Lambda(2))\ar[r]^\sim & \HH^4_{x'}(X',\Lambda(2))\ar@{=}[r] &
   \HH^4_{x'}(C',{i'}^!\Lambda(2)) & \ar[l] \HH^3(y,{k'}^!\Lambda(2)) 
}
\]

On peut identifier le carré externe de ce diagramme-ci au suivant, où les
flèches horizontales sont les flèches d'oubli du support :
\[
\xymatrix{
\HH^3(U,\Lambda(2)) & \ar[l] \HH^3_y(U,\Lambda(2)) \ar@{=}[d] \\
\ar[u]_{\sim} \HH^3(U',\Lambda(2)) & \ar[l] \HH^3_y(U',\Lambda(2))
}
\]

Ce diagramme-là est bien évidemment commutatif, ce qui achève la
démonstration du lemme~\ref{17-lemme-hexagone-commutatif}.

\bigskip

On peut prolonger sur la droite le diagramme commutatif du lemme :
\[
\xymatrix{
\HH^4_x(C,i^!\Lambda(2)) &
 \ar[l]_-{\transition^C_{\overline y\to x}} \HH^2_{\overline y}(C,i^!\Lambda(1)) & &\\
\HH^3(y,k^!\Lambda(2))\ar[u]^\sim\ar[d]^\sim  & &
({\mathscr H}^2 k^!\Lambda(1))_{\overline
y}\ar[dl]_{\sim}^1\ar[ul]_{\sim}^{[x':x]} &
\ar[l]_-{\transition^X_{\overline\eta\to\overline y}}^-{\sim}\Lambda\\
\HH^4_{x'}(C',{i'}^!\Lambda(2)) 
& \ar[l]_-{\transition^{C'}_{\overline y\to x}} \HH^2_{\overline
 y}(C',{i'}^!\Lambda(1)) & & }
\]
Du groupe $({\mathscr H}^2 k^!\Lambda(1))_{\overline y}$ partent deux
isomorphismes canoniques vers des groupes de cohomologie à supports de $C$
et de $C'$ ; ici, on les a multipliés respectivement par $[x':x]$ et $1$ de
façon à rendre le diagramme commutatif (cf.
définition~\ref{17-definition-transition-cas-local-dimension-un}).

On peut partir de l'élément $1$ dans le groupe $\Lambda$ tout à droite et
considérer son image dans le groupe $\HH^4_x(X,\Lambda(2))$ figurant sur le
diagramme du lemme. En suivant le chemin du bas, on obtient
$[x':x]\cdot [x]_{X'}$.
En suivant le chemin du haut, on obtient
$[x':x]\cdot \transition^X_{\overline y\to x}
(\transition^X_{\overline\eta\to\overline y}(1))$. On peut ainsi conclure que
l'on a bien l'egalité
\[[x]_{X'}=\transition^X_{\overline y\to x}
(\transition^X_{\overline\eta\to\overline y}(1))\]
dans $\HH^4_x(X,\Lambda(2))$.

\subsection{Morphismes de transition en codimension arbitraire}

Nous allons maintenant démontrer le
théorème~\ref{17-theoreme-transition-toute-codimension} en nous appuyant sur
le théorème~\ref{17-theoreme-classe-degre-maximal} établi pour le moment
jusqu'en dimension $2$. Ce résultat nous permettra ensuite d'établir le
théorème~\ref{17-theoreme-classe-degre-maximal} en toute généralité.

\begin{definition2}
Soit $X$ un $\ZZn$-schéma noethérien excellent.
Soit $\overline y=\overline {x_0}\to \dots \to \overline {x_n}=\overline x$ une
suite de spécialisations de points géométriques de $X$ telle que pour tout
$0\leq i<n$, la spécialisation $\overline x_i\to \overline x_{i+1}$ soit de
codimension $0$ ou $1$. On note $c$ la codimension de la spécialisation
$\overline y\to \overline x$.
Pour tout $K\in\detaleplus X$, on
note $\transition^X_{\overline x_0\to\dots \to
\overline x_n}\colon \HH^p_{\overline y}(K)\to 
\HH^{p+2c}_{\overline x}(K(c))$ le morphisme de transition obtenu par la
composition
$\transition^X_{\overline x_{n-1}\to\overline x_n}\circ\dots\circ \transition^X_{\overline x_0\to\overline x_1}$.
\end{definition2}

\begin{definition2}\label{17-definition-propriete-c}
Soit $X$ un $\ZZn$-schéma noethérien excellent. Soit $\overline y\to
\overline x$ une spécialisation de points géométriques de $X$ de
codimension $c$. Soit $K\in\detaleplus X$. On dit que la propriété
$\propC^c_{\overline y\to \overline x,K}$ est satisfaite si le morphisme
$\transition^X_{\overline {x_0}\to\dots\to\overline {x_n}}\colon
\HH^0_{\overline y}(K)\to \HH^{2c}_{\overline x}(K(c))$
ne dépend pas du choix de la factorisation de $\overline y\to\overline x$
en $\overline {x_0}\to\dots\to\overline {x_n}$.
\end{definition2}

\begin{definition2}
On dira que la propriété $\propC^{\leq c}$ est satisfaite si toutes les
propriétés $\propC^{c'}_{\overline y\to\overline x,K}$ envisagés dans
la définition~\ref{17-definition-propriete-c} le sont pour $c'\leq c$.
On dira que la
propriété $\propC^{\leq c}_{\loc}$ est vérifiée si les propriétés
$\propC^{c'}_{\overline \eta\to x,K}$ sont vérifiées 
dans la situation, dite locale, où le schéma $X$ est supposé 
local strictement hensélien intègre de dimension $\leq c$
de point fermé $x$ et où $\overline \eta$ est au-dessus du point générique
de $X$. Enfin, on dira que la propriété $\propC^{\leq
c}_{\loc,\normal,\Lambda}$ est satisfaite si on suppose de plus que $X$ est
normal et que $K=\Lambda$.
\end{definition2}

Il s'agit donc d'établir la propriété $\propC^{\leq c}$ pour tout $c\geq
0$. D'après la
remarque~\ref{17-remarque-composition-transition-codimension-0-ou-1}, on peut
considérer que la propriété $\propC^{\leq 1}$ est connue. En outre, notons
que le théorème~\ref{17-theoreme-transition-toute-codimension} affirme en
particulier $\propC^{\leq c}_{\loc,\normal,\Lambda}$ pour tout $c\geq 0$. Les
résultats de la sous-section~\ref{17-subsection-le-cas-de-la-dimension-2}
montrent que $\propC^{\leq 2}_{\loc,\normal,\Lambda}$ est satisfaite.

\begin{proposition2}\label{17-proposition-compatibilite-transition-specialisation}
Soit $f\colon X'\to X$ un morphisme fini entre $\ZZn$-schémas noethériens
excellents. Soit $\overline{y'}\to\overline{x'}$ une spécialisation
de points géométriques de $X'$ de codimension $c$
au-dessus d'une spécialisation
$\overline y\to\overline x$ de points géométriques de $X$
(cf. définition~\ref{17-definition-specialisation-au-dessus-d-une-autre}).
Pour toute factorisation $\overline{x'_0}\to\dots\overline{x'_n}$ de
$\overline {y'}\to\overline{x'}$ en une suite de spécialisations de
codimension $0$ ou $1$, si on note $\overline{x_0}\to\dots\overline{x_n}$
la factorisation de $\overline y\to\overline x$ en-dessous de la
précédente, pour tout $K\in\detaleplus X$ et $i\in\mathbf{Z}$,
le diagramme suivant commute :
\[
\xymatrix{\HH^i_{\overline{y'}}(f^!K)
     \ar[rr]^-{\transition^{X'}_{\overline{x'_0}\to\dots\to\overline{x'_n}}}
&  & \HH^{i+2c}_{\overline{x'}}(f^!K(c)) \\
\HH^i_{\overline y}(K)
     \ar[rr]^-{\transition^{X}_{\overline{x_0}\to\dots\to\overline{x_n}}}
     \ar[u]_{\sim}^{\Cl {y'\to y}}
&  & \HH^{i+2c}_{\overline{x}}(K(c))\ar[u]^{\sim}_{\Cl {x'\to x}}}
\]
En outre, on a l'équivalence
$\propC^c_{\overline{y}\to\overline{x},K}\Longleftrightarrow
\propC^c_{\overline{y'}\to\overline{x'},f^!K}$.
\end{proposition2}

La commutativité du diagramme se déduit aussitôt de la 
proposition~\ref{17-proposition-compatibilite-transition-spec-immediate}.
L'équivalence annoncée résulte de la commutativité du diagramme et du fait
qu'il est essentiellement équivalent de 
se donner une factorisation de la spécialisation
$\overline{y'}\to\overline{x'}$ ou de s'en donner une de
$\overline y\to\overline x$ (cf.
propositions~\ref{17-proposition-specialisation-going-down} et
\ref{17-proposition-specialisation-going-up}).

\begin{lemme2}\label{17-lemme-reduction-propriete-c-situation-locale}
Pour tout entier $c\geq 0$, on a l'équivalence $\propC^{\leq
c}\Longleftrightarrow \propC^{\leq c}_{\loc}$.
\end{lemme2}

Supposons $\propC^{\leq c}_{\loc}$ et montrons $\propC^{\leq c}$.
Il suffit évidemment de montrer
$\propC^{c'}_{\overline y\to\overline x,K}$ pour tout $K\in\detaleplus X$, avec $X$
un schéma local strictement hensélien de point fermé $\overline x$ et une
spécialisation $\overline y\to\overline x$ de
codimension $c'\leq c$. Il s'agit de vérifier que
l'on peut supposer que $X$ est intègre et que
$\overline y$ est au-dessus du point générique. Pour cela, on introduit
l'immersion fermée $i\colon Z\to X$ où $Z=\overline{\left\{y\right\}}$. La
proposition~\ref{17-proposition-compatibilite-transition-specialisation}
montre que $\propC^{c'}_{\overline y\to\overline x,i^!K}$ implique
$\propC^{c'}_{\overline y\to\overline x,K}$. Pour conclure, il suffit
d'observer que $\propC^{c'}_{\overline y\to\overline x,i^!K}$ est un cas
particuler de $\propC^{\leq c}_{\loc}$.

\begin{lemme2}\label{17-lemme-propriete-c-reductions}
Soit $X$ un $\ZZn$-schéma noethérien excellent local strictement hensélien
intègre de dimension $c$,
de point fermé $x$ et de point générique $\eta$. Soit $\overline
\eta$ un point géométrique au-dessus de $\eta$. Soit $K\in\detaleplus X$.
\begin{enumerate}
\item[(a)] $\propC^c_{\overline\eta\to x,\tau_{\leq 0}K}\Longrightarrow
\propC^c_{\overline\eta\to x,K}$ ;
\item[(b)] $\propC^c_{\overline\eta\to x,\tau_{\leq 0}K}\Longleftrightarrow
\propC^c_{\overline\eta\to x,{\mathscr H}^0K}$ ;
\item[(c)] Soit $j\colon U\to X$ l'inclusion d'un ouvert dense
et $\mathscr M$ et $\mathscr N$ deux faisceaux de
$\Lambda$-modules sur $X$ tels que $j^\star \mathscr M\simeq j^\star
\mathscr N$, alors 
$\propC^c_{\overline\eta\to x,\mathscr M}\Longleftrightarrow
\propC^c_{\overline\eta\to x,\mathscr N}$.
\end{enumerate}
\end{lemme2}

L'implication (a) résulte du fait que le morphisme canonique $\tau_{\leq
0}K\to K$ induise un isomorphisme $\HH^0_{\overline \eta}(\tau_{\leq
0}K)\iso \HH^0_{\overline \eta}(K)$.

On considère ensuite le morphisme canonique
$\tau_{\leq 0}K\to {\mathscr H}^0 K$. Il induit évidemment un isomorphisme
après application de $\HH^0_{\overline \eta}$, mais aussi après celle de
$\HH^{2c}_{x}$, pour des raisons de dimension cohomologique (cf.
proposition~\ref{17-proposition-dimension-cohomologique}). L'équivalence (b)
en résulte aussitôt.

Pour montrer (c), on peut supposer que $\mathscr N=j_!j^\star \mathscr M$.
On a alors un monomorphisme canonique $j_!j^\star \mathscr M\to\mathscr M$,
dont le conoyau est $i_\star i^\star\mathscr M$, où $i\colon Z\to X$
désigne une immersion fermée complémentaire. La
proposition~\ref{17-proposition-dimension-cohomologique} montre alors que
$\HH^q(X-x,(i_\star i^\star\mathscr M)_{|X-x})=0$ pour $q\geq 2c-2$. Ainsi,
le morphisme $j_!j^\star \mathscr M\to\mathscr M$ induit un isomorphisme
non seulement après application de $\HH^0_{\overline \eta}$, mais aussi celle
de $\HH^{2c}_x$, ce qui permet de conclure.

\begin{lemme2}\label{17-lemme-reduction-propriete-c-faisceau-constant}
Pour tout entier $c\geq 0$, on a l'implication $\propC^{\leq
c}_{\loc,\normal,\Lambda}\Longrightarrow \propC^{\leq c}_{\loc}$.
\end{lemme2}

Soit $X$ un $\ZZn$-schéma noethérien excellent
local strictement hensélien intègre de dimension $c'\leq c$, de point
fermé $x$ et de point générique $\eta$. Soit $\overline\eta$ un point
géométrique au-dessus de $\eta$. D'après le 
lemme~\ref{17-lemme-propriete-c-reductions}, pour montrer 
$\propC^{c'}_{\overline\eta\to x,K}$ pour tout $K\in\detaleplus X$, il suffit
de montrer $\propC^{c'}_{\overline\eta\to x,\mathscr M}$ pour tout faisceau
de $\Lambda$-modules $\mathscr M$ sur $X$. On peut évidemment supposer que
$\mathscr M$ est constructible. Mais alors, il existe un ouvert dense de
$X$ au-dessus duquel $\mathscr M$ soit localement constant. Il existe
évidemment un morphisme fini surjectif $f\colon \tilde X\to X$, avec
$\tilde X$ normal intègre, et un ouvert dense $U$ de $X$ tel que $f$
induise un revêtement fini étale $\tilde U=f^{-1}(U)\to U$ et que
$(f^\star \mathscr M)_{|\tilde U}$ soit isomorphe à un faisceau constant de
valeur un certain $\Lambda$-module $N$. Choisissons une spécialisation
$\overline {\tilde \eta}\to \tilde x$ au-dessus de
$\overline\eta\to x$. La propriété $\propC^{\leq c}_{\loc,\normal,\Lambda}$
contient $\propC^{c'}_{\overline {\tilde \eta}\to \tilde x,\Lambda}$
comme cas particulier. Cette dernière propriété
implique à son tour la propriété 
$\propC^{c'}_{\overline {\tilde \eta}\to \tilde x,N}$. On observe que l'on
a un isomorphisme canonique $(f^!\mathscr M)_{|\tilde U}\simeq N$.
L'assertion (c) du lemme~\ref{17-lemme-propriete-c-reductions} permet
d'obtenir la propriété $\propC^{c'}_{\overline {\tilde \eta}\to \tilde
x,{\mathscr H}^0(f^!\mathscr M)}$, et les assertions (a) et (b) d'en
déduire $\propC^{c'}_{\overline {\tilde \eta}\to \tilde
x,f^!\mathscr M}$. Enfin, la 
proposition~\ref{17-proposition-compatibilite-transition-specialisation}
permet de vérifier $\propC^{c'}_{\overline\eta\to x,\mathscr M}$, ce
qui achève la démonstration du lemme.

\begin{lemme2}\label{17-lemme-propriete-c-etape-de-recurrence}
Pour tout entier $c\geq 3$, on a l'implication $\propC^{\leq
c-1}\Longrightarrow \propC^{\leq c}_{\loc,\normal,\Lambda}$.
\end{lemme2}

Soit $X=\SP(A)$
un schéma noethérien excellent local strictement hensélien normal
de dimension $c$, de point fermé $x$ et de point générique $\eta$. On
choisit un point géométrique $\overline \eta$ au-dessus de $\eta$. On note
$s\colon \overline\eta\to x$ la spécialiation canonique.
Nous
allons montrer la propriété $\propC^c_{s,\Lambda}$.

Soit $\overline z\to X$ un point géométrique tel que $z\neq x$ et $z\neq
\eta$. On a une spécialisation canonique $t\colon \overline z\to x$. On peut
choisir une spécialisation $u\colon \overline \eta\to\overline z$. On obtient
peut-être pas ainsi une factorisation de la spécialisation
$\overline\eta\to x$ fixée plus haut, mais il existe certainement
$\sigma\in\Gal(\overline\eta/\eta)$ faisant commuter le diagramme suivant
de spécialisations de points géométriques de $X$ :
\[
\xymatrix{\overline\eta\ar[r]^s\ar[d]^{\sigma}_\sim & x\\ 
\overline\eta \ar[r]^u & \overline z\ar[u]^t}
\]
Les spécialisations $\sigma$, $u$ et $t$
sont de codimension $\leq c-1$, les morphismes de transitions
associés sont donc bien définis.
On peut considérer l'image de $1\in \HH^0_{\overline\eta}(X,\Lambda)$ par le
composé de ces morphismes de transition :
\[\gamma_{\sigma,u,t}=(\transition^X_{t}\circ\transition^X_u\circ\transition^X_\sigma)(1)\in
\HH^{2c}_x(X,\Lambda(c))\;\text{.}\]
Bien entendu, on a $\transition^X_\sigma(1)=1$. On note donc
simplement $\gamma_{u,t}$ l'élément
$\gamma_{\sigma,u,t}=(\transition^X_{t}\circ\transition^X_u)(1)$. 
J'affirme que cette classe $\gamma_{u,t}$
ne dépend que du point $z$ de $X$ en-dessous duquel $\overline z$ se
trouve. En effet, soit $\overline {z'}$ un autre point géométrique
au-dessus de $z$, et $t'\colon\overline{z'}\to x$ la spécialisation
canoniquement associée. On choisit une spécialisation $u'\colon
\overline\eta\to \overline{z'}$. On peut choisir un $z$-isomorphisme 
$\sigma'\colon \overline z\iso \overline {z'}$.
Le schéma $X$ étant géométriquement unibranche, on peut montrer facilement
qu'il existe un $\eta$-isomorphisme $\sigma''\colon
\overline\eta\to\overline\eta$ induisant 
un diagramme commutatif de spécialisations de points
géométriques de $X$ :
\[
\xymatrix{
\overline\eta\ar[r]^u\ar[d]^{\sigma''}_\sim&
\overline{z}\ar[d]^{\sigma'}_\sim\ar[r]^t & x\ar@{=}[d] \\
\overline\eta\ar[r]^{u'}& \overline{z'}\ar[r]^{t'} & x
}
\]
En utilisant cette fois-ci que $\sigma''$ agit trivialement sur $1\in
\HH^0_{\overline \eta}(X,\Lambda)$, on montre que 
$\gamma_{u,t}=\gamma_{u',t'}$. On peut donc noter
simplement $\gamma_z$ cette classe.

Il s'agit de montrer que la classe $\gamma_z$ est indépendante de $z\in
X-\left\{\eta,x\right\}$. Par construction, si 
$z'\in X-\left\{\eta,x\right\}$ est tel qu'il existe une spécialisation
(Zariski) $z\to z'$, on a $\gamma_z=\gamma_{z'}$.

Comme $A$ est local normal excellent de dimension $\geq 3$, d'après
\sga{2}{xiii}{2.1}, si $f\in A-\left\{0\right\}$, le schéma $\SP\left(
A/(f)-\left\{x\right\}\right)$ est connexe. Il en résulte aussitôt que les
classes $\gamma_z$ pour $z\in \SP\left(A/(f)-\left\{x\right\}\right)$
sont égales.
Maintenant, si $z$ et $z'$ sont deux éléments de $X-\left\{\eta,x\right\}$,
il existe bien évidemment $f\in A-\left\{0\right\}$ tel que $z$ et $z'$
appartiennent à l'hypersurface définie par $f$. On obtient donc bien
$\gamma_z=\gamma_{z'}$, ce qui achève la démonstration du lemme.

\bigskip

Nous sommes maintenant en mesure de démontrer le
théorème~\ref{17-theoreme-transition-toute-codimension}. Il s'agit d'établir
la propriété $\propC^{\leq c}$ pour tout $c\geq 0$. Elle équivaut à la
propriété $\propC^{\leq c}_{\loc}$ d'après le
lemme~\ref{17-lemme-reduction-propriete-c-situation-locale}. D'après le
lemme~\ref{17-lemme-reduction-propriete-c-faisceau-constant}, elle équivaut
encore à la propriété 
$\propC^{\leq c}_{\loc,\normal,\Lambda}$. Compte tenu de ces
équivalences, le
lemme~\ref{17-lemme-propriete-c-etape-de-recurrence} montre que pour $c\geq
3$, $\propC^{\leq c-1}_{\loc,\normal,\Lambda}$ implique 
$\propC^{\leq c}_{\loc,\normal,\Lambda}$. Un raisonnement par récurrence
permet donc de conclure puisque la propriété $\propC^{\leq
2}_{\loc,\normal,\Lambda}$ a été établie plus haut.

\subsection{Fin de la démonstration}

\begin{proposition2}
Soit $X$ un $\ZZn$-schéma noethérien excellent 
local hensélien normal de dimension $d$, de point fermé $x$ et
de point générique $\eta$. Soit $\overline x$ un point géométrique
au-dessus de $x$. Soit $\overline \eta$ un point géométrique au-dessus de
$\eta$. Soit $\overline\eta\to\overline x$ une spécialisation.
Alors, l'image de $1$ par le morphisme de transition
$\transition^X_{\overline\eta \to \overline x}\colon
\HH^0_{\overline \eta}(X,\Lambda)\to \HH^{2d}_{\overline x}(X,\Lambda(d))$ est
invariante par $\Gal(\overline x/x)$, et indépendante de la spécialisation
choisie $\overline \eta\to \overline x$, on la note $[x]$.
On obtient ainsi un morphisme
$\Lambda[-2d]\to \tau_{\geq 2d}\SheafRG_x(\Lambda(d))$
dans $\detale x$.
\end{proposition2}

Ceci résulte aussitôt du
théorème~\ref{17-theoreme-transition-toute-codimension} et du fait que $1\in
\HH^0_{\overline\eta}(X,\Lambda)$ est fixé par $\Gal(\overline\eta/\eta)$.

\bigskip

Par construction, les classes $[x]$ sont compatibles aux morphismes de
transition. Pour finir de démontrer le
théorème~\ref{17-theoreme-classe-degre-maximal}, il reste à montrer que dans
la situation de la proposition précédente, le morphisme $\Lambda[-2d]\to
\tau_{\geq 2d}\SheafRG_x(\Lambda(d))$ est un isomorphisme, ou
encore, dans la situation où $X$ est strictement hensélien, que le
morphisme $\Lambda\to \HH^{2d}_x(X,\Lambda(d))$ induit par $[x]$ est un
isomorphisme. En effet, des considérations de dimension cohomologique (cf.
proposition~\ref{17-proposition-dimension-cohomologique}) expliquent
l'annulation de $\HH^q_x(X,\Lambda(d))$ pour $q>2d$.

\begin{lemme2}
Soit $d\geq 3$.
On suppose le résultat du théorème~\ref{17-theoreme-classe-degre-maximal}
connu jusqu'en dimension $d-1$. Alors, pour tout $\ZZn$-schéma noethérien
excellent $X$ local strictement hensélien normal de dimension $d$, de point
fermé $x$, le morphisme $\Lambda\to \HH^{2d}_x(X,\Lambda(d))$
induit par $[x]$ est surjectif.
\end{lemme2}

On utilise la suite spectrale de coniveau calculant la cohomologie de
l'ouvert $U=X-x$ à coefficients dans $\Lambda(d)$ :
\[E_1^{pq}=\bigoplus_{y\in U^p} \HH^{p+q}_y(X_{(y)},\Lambda(d))\Longrightarrow
\HH^{p+q}(U,\Lambda(d))\;\text{,}\]
où $U^p$ désigne l'ensemble des points de codimension $p$ de $U$. La
topologie étale étant plus fine que la topologie de Nisnevich, dans
l'expression du terme $E_1$, on peut remplacer les localisés $X_{(y)}$ par
leurs hensélisés $X_{(y)}^\hensel$ pour tout $y\in U^p$ avec $0\leq p\leq
d-1$. On obtient ainsi un isomorphisme
canonique :
\[\RG_y(X_{(y)},\Lambda(d))\simeq
\RG(y,\SheafRG_y(\Lambda(p))(d-p))\;\text{.}\]
La structure de $\tau_{\geq 2p}\SheafRG_y(\Lambda(p))$ est
connue par notre connaissance limitée du
théorème~\ref{17-theoreme-classe-degre-maximal} puisque $p\leq d-1$. D'après
la proposition~\ref{17-proposition-dimension-cohomologique-galois}, on a une
majoration de la $\ell$-dimension cohomologique de $y$ pour tout nombre
premier $\ell$ divisant $n$ : $\cd_\ell y\leq d-p$.
En appliquant la suite spectrale de composition des
foncteurs dérivés au calcul de $\RG_y(X_{(y)},\Lambda(d))$, on
obtient
d'une part que $\HH^i_y(X_{(y)},\Lambda(d))=0$ pour $i>d+p$, c'est-à-dire que
$\HH^{p+q}(X_{(y)},\Lambda(d))=0$ pour $q>d$, et d'autre part que 
$\HH^{p+d}(X_{(y)},\Lambda(d))\simeq \HH^{d-p}(y,\Lambda(d-p))$.

On a donc montré que $E^{p,q}_1=0$ pour $q\geq d+1$, et par ailleurs, il est
évident que $E^{p,q}_1=0$ si $p\geq d$.
On en déduit aussitôt que l'on a
des isomorphismes canoniques $\coker(E^{d-2,d}_1\to E^{d-1,d}_1)\iso
\HH^{2d-1}(U,\Lambda(d))\iso \HH^{2d}_x(X,\Lambda(d))$. En particulier, le
morphisme canonique $E^{d-1,d}_1\to \HH^{2d}_x(X,\Lambda(d))$ est surjectif.
D'après le calcul ci-dessus, on connaît la structure de $E^{d-1,d}_1$ :
\[E^{d-1,d}_1 \simeq  \bigoplus_{y\in U^{d-1}} \HH^1(y,\Lambda(1))\;\text{.}\]
Pour tout $y\in U^{d-1}$, on a 
un isomorphisme évident $\HH^1(y,\Lambda(1))\simeq \Lambda$ (induit par une
spécialisation immédiate $\overline y\to x$ de points géométriques de $X$).
Il n'est pas difficile de vérifier que l'image du
morphisme canonique $\Lambda\simeq \HH^1(y,\Lambda(1))\subset
E_1^{d-1,d}\to \HH^{2d}_x(X,\Lambda(d))$ est 
le sous-groupe engendré par $[x]$. Par conséquent, l'image du morphisme
surjectif $E^{d-1,d}_1\to \HH^{2d}_x(X,\Lambda(d))$ est le sous-groupe
engendré par $[x]$, ce qui achève la démonstration du lemme.

\begin{lemme2}\label{17-lemme-recurrence-theoreme-classe-degre-maximal}
Soit $d\geq 3$.
On suppose le résultat du théorème~\ref{17-theoreme-classe-degre-maximal}
connu jusqu'en dimension $d-1$. Alors, pour tout $\ZZn$-schéma noethérien
excellent $X$ local strictement hensélien normal de dimension $d$, de point
fermé $x$, le morphisme $\Lambda\to \HH^{2d}_x(X,\Lambda(d))$
induit par $[x]$ est un isomorphisme, autrement dit, l'énoncé du
théorème~\ref{17-theoreme-classe-degre-maximal} vaut jusqu'en dimension $d$.
\end{lemme2}

En vertu du lemme précédent, il ne s'agit plus que de déterminer la
structure du $\Lambda$-module $\HH^{2d}_x(X,\Lambda(d))$. Compte tenu du
théorème de changement de base formel (cf.
\cite[corollary~6.6.4]{TUBULAR@Fujiwara}),
on peut supposer que $X$ est complet.
Les théorèmes de structure des anneaux locaux noethériens
complets (cf. \egazero{iv}{19.8.8})
montrent que $X$ est alors isomorphe à un
sous-schéma fermé d'un schéma régulier. Pour la fonction de dimension
$\delta$ sur $X$ telle que $\delta(\eta)=0$ où $\eta$ est le point
générique de $X$, compte tenu de ce qu'on sait déjà sur les complexes
dualisants potentiels (cf.
proposition~\ref{17-proposition-cdpotentiel-regulier}),
il existe un complexe dualisant potentiel
$K$ pour $(X,\delta)$. Nous allons utiliser le lemme général suivant :

\begin{lemme2}\label{17-lemme-structure-cohomologie-negative-cdp-normal}
Soit $X$ un $\ZZn$-schéma noethérien excellent normal intègre de point
générique $\eta$. On suppose $X$ muni de la fonction de dimension $\delta$
telle que $\delta(\eta)=0$. Soit $K$ un complexe dualisant potentiel pour
$(X,\delta)$. Alors, les faisceaux de cohomologie ${\mathscr H}^qK$  sont
nuls pour $q<0$ et l'épinglage en $\eta$ s'étend en un isomorphisme
${\mathscr H}^0K\simeq \Lambda$. Autrement dit, on a un isomorphisme
canonique $\Lambda\iso \tau_{\leq 0}K$.
\end{lemme2}

Pour obtenir le résultat pour $X$, il suffit de l'avoir pour ses hensélisés
stricts. On peut donc supposer que $X$ est strictement hensélien de
dimension $d$ et de point fermé $x$. On procède par récurrence sur $d$.
Si $d=0$, le résultat est évident. On suppose donc que $d\geq 1$ et que le
résultat est connu pour l'ouvert $U=X-x$. Notons $j\colon U\to X$ et
$i\colon x\to X$ les immersions évidentes. On a un triangle distingué
dans $\detaleplus X$ :
\[i_\star i^!K\to K\to \R j_\star j^\star K\to i_\star
i^! K[1]\;\text{.}\]
Grâce à l'épinglage en $x$, il vient que ${\mathscr H}^q i^!K=0$ pour
$q\leq 1$. L'hypothèse sur $U$ montre que $\Lambda\iso \tau_{\leq
0}j^\star K$. On en déduit que ${\mathscr H}^q K=0$ pour $q<0$ et que
${\mathscr H}^0K\simeq j_\star \Lambda$. Comme $X$ est normal, le morphisme
canonique $\Lambda\to j_\star \Lambda$ est un isomorphisme, ce qui achève
la démonstration du lemme.

\bigskip

Revenons à la démonstration du 
lemme~\ref{17-lemme-recurrence-theoreme-classe-degre-maximal}, on considère le
morphisme canonique $\Lambda\to K$ déduit du lemme ci-dessus. Soit
$\overline \eta$ un point géométrique au-dessus de $\eta$.
Choisissons une spécialisation $\overline\eta\to x$. On considère le
diagramme commutatif suivant, où les flèches verticales sont induites par
$\Lambda\to K$ et les flèches horizontales par les morphismes de transition
associés à la spécialisation $\overline\eta\to x$ :
\[
\xymatrix{
\ar[d]\HH^0_{\overline\eta}(\Lambda)\ar[r] & \HH^{2d}_x(X,\Lambda(d)) \ar[d]\\
\HH^0_{\overline\eta}(K)\ar[r] & \HH^{2d}_x(X,K(d))
}
\]
Le morphisme de gauche est évidemment un isomorphisme. Celui du bas aussi
puisque $K$ est un complexe dualisant potentiel. Le morphisme du haut est
donc injectif, mais on sait déjà qu'il est surjectif. Le $\Lambda$-module
$\HH^{2d}_x(X,\Lambda(d))$ est donc isomorphe à $\Lambda$, ce qui permet de
conclure.

\section{Compléments sur les complexes dualisants potentiels}
\label{17-section-complements-cdp}
\subsection{Énoncés}

Grâce aux résultats de la section~\ref{17-section-classe-de-cohomologie}, nous
allons pouvoir poursuivre l'étude de certaines propriétés des complexes
dualisants potentiels.
Le résultat principal de la
section~\ref{17-section-premieres-proprietes} était la construction d'un
complexe dualisant potentiel sur les schémas réguliers (cf.
proposition~\ref{17-proposition-cdputatif-regulier}). Les deux propositions
suivantes établissent des propriétés de stabilité des complexes dualisants
potentiels par rapport aux morphismes de type fini et aux morphismes
réguliers \ega{iv}{6.8.1}.

\begin{proposition2}\label{17-proposition-stabilite-cdp-morphisme-regulier}
Soit $f\colon Y\to X$ un morphisme régulier
entre $\ZZn$-schémas noethériens excellents. On suppose $X$ muni
d'une fonction de dimension $\delta_X$. On munit $Y$ de la fonction de
dimension $\delta_Y$ définie par l'égalité
$\delta_Y(y)=\delta_X(f(y))-\codim_{f^{-1}(f(y))}(y)$
pour tout $y\in Y$.

Si $K$ un complexe dualisant putatif sur $X$, alors $f^\star K$ est
naturellement muni d'une structure de complexe dualisant putatif sur $Y$,
et c'est un complexe dualisant potentiel si $K$ en est un.
\end{proposition2}

\begin{proposition2}\label{17-proposition-stabilite-cdp-morphisme-de-type-fini}
Soit $f\colon Y\to X$ un morphisme de type fini compactifiable
entre $\ZZn$-schémas noethériens excellents. On suppose $X$ muni
d'une fonction de dimension $\delta_X$. On munit $Y$ de la fonction de
dimension $\delta_Y$ définie par l'égalité
$\delta_Y(y)=\delta_X(f(y))+\degtr(y/f(y))$
pour tout $y\in Y$.

Si $K$ un complexe dualisant putatif sur $X$, alors $f^! K$ est
naturellement muni d'une structure de complexe dualisant putatif sur $Y$,
et c'est un complexe dualisant potentiel si $K$ en est un.
\end{proposition2}

\begin{remarque2}\label{17-remarque-cdp-compatibilite-composition}
Au cours des démonstrations, on observera que les constructions des
propositions~\ref{17-proposition-stabilite-cdp-morphisme-regulier} et
\ref{17-proposition-stabilite-cdp-morphisme-de-type-fini} sont compatibles à
la composition des morphismes : si $g\colon Z\to Y$ est un autre morphisme
du type envisagé et $K$ un complexe dualisant putatif sur $X$,
l'isomorphisme de transitivité $g^\star f^\star K\simeq (f\circ g)^\star K$
(resp. $g^! f^! K\simeq (f\circ g)^! K$) est compatible aux épinglages.
La proposition~\ref{17-proposition-stabilite-cdp-morphisme-regulier}
généralise la construction de la
proposition~\ref{17-proposition-stabilite-cdp-morphisme-etale} (cas où $f$ est
étale).
\end{remarque2}

Nous verrons plus bas que la
proposition~\ref{17-proposition-stabilite-cdp-morphisme-de-type-fini} résulte
facilement de la
proposition~\ref{17-proposition-stabilite-cdp-morphisme-regulier}. Nous allons
donc nous intéresser plus particulièrement aux morphismes réguliers. La
sous-section suivante sur le théorème de changement de base par un
morphisme régulier nous permettra de définir une structure de complexe
dualisant putatif sur $f^\star K$ dans la situation de la
proposition~\ref{17-proposition-stabilite-cdp-morphisme-regulier}. Dans le cas
d'un complexe dualisant potentiel, ce sont les résultats de la
section~\ref{17-section-classe-de-cohomologie} qui vont nous permettre de
vérifier que les épinglages sur $f^\star K$ sont compatibles aux
spécialisations.

\subsection{Changement de base par un morphisme régulier}

Dans ce numéro, on étudie une généralisation du théorème de changement de
base par un morphisme lisse au cas du changement de base par un
morphisme régulier.
La démonstration utilise le théorème de
Popescu (cf.~\cite{NPD@Swan}).
On examine ensuite quelques conséquences de ce résultat.

\begin{proposition2}\label{17-proposition-changement-de-base-regulier}
Soit $f\colon Y\to X$ un morphisme régulier entre $\ZZn$-schémas
localement noethériens. Soit $g\colon X'\to X$ un
morphisme quasi-compact et quasi-séparé. On constitue le carré cartésien de
schémas :
\[
\xymatrix{
Y'\ar[r]^{f'}\ar[d]^{g'} & X'\ar[d]^{g} \\
Y\ar[r]^f & X
}
\]
Alors, pour tout $K\in \detaleplus {X'}$, le morphisme de changement de base
$f^\star \R g_\star K\to \R {g'}_\star {f'}^\star K$ est un
isomorphisme dans $\detaleplus Y$.
\end{proposition2}

Les arguments habituels montrent que l'on peut supposer que $X$, $Y$, puis
$X'$, et donc $Y'$, sont affines. D'après le théorème de Popescu, on peut
écrire $Y$ comme la limite projective d'un système projectif filtrant
$(Y_{\alpha}\vers {f_\alpha} X)_\alpha$ de $X$-schémas 
affines lisses (voir \ega{iv}{8}). Le théorème de changement de
base par un morphisme lisse \sga{4}{xvi}{1.2} montre que les morphismes de
changement de base $f_\alpha^\star \R g_\star K\to 
\R {g'_\alpha}_\star {f'}_\alpha^\star K$ sont des isomorphismes.
Ce résultat pour les morphismes
$f_\alpha$ s'étend au morphisme $f$ grâce à la technique de passage à la
limite de \sga{4}{vii}{5.8}, ce qui permet de conclure.

\begin{proposition2}\label{17-proposition-rhom-morphismes-reguliers}
Soit $f\colon Y\to X$ un morphisme régulier
entre $\ZZn$-schémas noethériens.
Soit $K\in\detalebc X$. Soit $L\in \detaleplus X$. Alors, le morphisme évident est
un isomorphisme dans $\detaleplus Y$ :
\[f^\star\R\SheafHom(K,L)\iso \R\SheafHom(f^\star K,f^\star
L)\;\text{.}\]
\end{proposition2}

Avant de la démontrer, signalons un corollaire de cette proposition :

\begin{corollaire2}\label{17-corollaire-cb-regulier-pour-cohomologie-a-supports}
Soit $f\colon Y\to X$ un morphisme régulier
entre $\ZZn$-schémas noethériens. Soit $i\colon Z\to X$ une 
immersion
fermée. Notons $i'\colon f^{-1}(Z)\to Y$ l'immersion fermée déduite par
changement de base, et $f'\colon f^{-1}(Z)\to Z$ la projection.
Alors, pour tout $K\in \detaleplus X$, le morphisme évident
${f'}^\star i^! K\to {i'}^!f^\star K$ est un isomorphisme dans $\detaleplus
{f^{-1}(Z)}$.
\end{corollaire2}

Pour démontrer la proposition~\ref{17-proposition-rhom-morphismes-reguliers},
commençons par constater que ce corollaire est un cas
particulier de la proposition~\ref{17-proposition-changement-de-base-regulier}
(l'appliquer avec pour $g$ l'immersion ouverte complémentaire de $i$).

\medskip

\begin{lemme2}\label{17-lemme-rhom-morphismes-reguliers-reduction-ferme}
Soit $f\colon Y\to X$ un morphisme régulier
entre $\ZZn$-schémas noethériens. Soit $i\colon Z\to X$ une
immersion fermée. On constitue le carré cartésien suivant :
\[\xymatrix{Z'\ar[r]^{i'}\ar[d]^{f'} & Y\ar[d]^f \\Z\ar[r]^i & X}\]
Soit $K\in \detalebc Z$. On suppose que pour tout $L\in\detaleplus Z$, le
morphisme canonique
${f'}^\star\R\SheafHom(K,L)\to\R\SheafHom({f'}^\star
K,{f'}^\star L)$ est un isomorphisme. Alors, pour tout $L\in\detaleplus X$, le
morphisme canonique $f^\star\R\SheafHom(i_\star K,L)\to
\R\SheafHom(f^\star i_\star K,f^\star L)$ est un isomorphisme.
\end{lemme2}

Ce lemme-ci résulte aussitôt du 
corollaire~\ref{17-corollaire-cb-regulier-pour-cohomologie-a-supports}.

\begin{lemme2}\label{17-lemme-rhom-morphismes-reguliers-reduction-ouvert}
Soit $f\colon Y\to X$ un morphisme régulier
entre $\ZZn$-schémas noethériens. Soit $j\colon U\to X$ une
immersion ouverte. On constitue le carré cartésien suivant :
\[\xymatrix{U'\ar[r]^{j'}\ar[d]^{f'} & Y\ar[d]^f \\U\ar[r]^j & X}\]
Soit $K\in \detalebc U$. On suppose que pour tout $L\in\detaleplus U$, le
morphisme canonique
${f'}^\star\R\SheafHom(K,L)\to\R\SheafHom({f'}^\star
K,{f'}^\star L)$ est un isomorphisme. Alors, pour tout $L\in\detaleplus X$, le
morphisme canonique $f^\star\R\SheafHom(j_! K,L)\to
\R\SheafHom(f^\star j_! K,f^\star L)$ est un isomorphisme.
\end{lemme2}

Ce lemme-là résulte de l'isomorphisme $f^\star \R j_\star \simeq
\R j'_\star f'^\star$ qui est un cas particulier de la 
proposition~\ref{17-proposition-changement-de-base-regulier}.

\bigskip

Démontrons la proposition~\ref{17-proposition-rhom-morphismes-reguliers}.
Les lemmes~\ref{17-lemme-rhom-morphismes-reguliers-reduction-ferme} et 
\ref{17-lemme-rhom-morphismes-reguliers-reduction-ouvert} permettent de
supposer que $K$ est un faisceau localement constant constructible de
$\Lambda$-modules sur $X$. En utilisant que le résultat est
trivial si $f$ est étale, on se ramène finalement au cas facile
où $K$ est un faisceau constant.

\subsection{Démonstration de la proposition~\ref{17-proposition-stabilite-cdp-morphisme-regulier}}

\subsubsection{Structure de complexe dualisant putatif sur \texorpdfstring{$f^\star K$}{f*K}}

\begin{lemme3}\label{17-lemme-construction-cdputatif-f-star}
Dans la situation de la
proposition~\ref{17-proposition-stabilite-cdp-morphisme-regulier}, si $K$ est
un complexe dualisant putatif sur $X$,
on peut munir $f^\star K$ d'une structure de
complexe dualisant putatif sur $Y$.
\end{lemme3}

Soit $K$ un complexe dualisant putatif sur $X$.
Soit $x\in X$. Introduisons le localisé $X'=X_{(x)}$ de $X$ en $x$ et
$Y'=Y\times_X X'$. Notons
$i\colon f^{-1}(x)\to Y'$ l'immersion de la fibre au-dessus de $x$,
$g\colon f^{-1}(x)\to x$ la projection et $j\colon Y'\to Y$ le morphisme
canonique.
D'après le
corollaire~\ref{17-corollaire-cb-regulier-pour-cohomologie-a-supports},
on a un isomorphisme canonique $g^\star
\SheafRG_x(K)\simeq i^!j^\star f^\star K$. Définir des épinglages
sur $f^\star K$ en les points de $f^{-1}(x)$ revient à définir une
structure de complexe dualisant putatif sur $i^!j^\star f^\star K$ pour
la fonction de dimension ${\delta_Y}_{|f^{-1}(x)}$. On vient de voir que
cet objet s'identifie à $g^\star \SheafRG_x(K)$ qui s'identifie
lui-même à $\Lambda(\delta_X(x))[2\delta_X(x)]$ en vertu de l'épinglage
donné de $K$ en $x$. La fibre $f^{-1}(x)$ étant régulière, le faisceau
constant $\Lambda$ est naturellement muni d'une structure de complexe
dualisant potentiel pour la fonction de dimension $-\codim$ (cf.
proposition~\ref{17-proposition-cdputatif-regulier}). Vu la définition de
$\delta_Y$, on en déduit aussitôt une structure de complexe dualisant
potentiel sur $\Lambda(\delta_X(x))[2\delta_X(x)]\in\detaleplus {f^{-1}(x)}$
pour la fonction de dimension ${\delta_Y}_{|f^{-1}(x)}$.
On a ainsi obtenu des
épinglages pour $f^\star K$ en tous les points de
$f^{-1}(x)$. Grâce à cette construction fibre à fibre, on a défini
une structure de complexe dualisant putatif sur $f^\star K$.

\begin{lemme3}\label{17-lemme-compatibilite-cdp-f-star-schemas-reguliers}
Soit $f\colon Y\to X$ un morphisme régulier entre $\ZZn$-schémas
noethériens excellents réguliers. On munit $X$ et $Y$ des fonctions de
dimension $\delta_X=-\codim$ et $\delta_Y=-\codim$. On munit les faisceaux
constants $\Lambda$ sur $X$ et $Y$ des structures de complexes dualisants
putatifs définies dans la proposition~\ref{17-proposition-cdputatif-regulier} ;
le lemme~\ref{17-lemme-construction-cdputatif-f-star} munit $f^\star\Lambda$
d'une structure de complexe dualisant putatif sur $Y$ pour la fonction de
dimension $\delta_Y$. Alors, l'isomorphisme canonique $\Lambda\simeq
f^\star \Lambda$ sur $Y$ est compatible aux épinglages.
\end{lemme3}

Vérifions la compatibilité de l'isomorphisme canonique $\Lambda\simeq
f^\star\Lambda$ aux épinglages en un point $y$ de $Y$. Quitte à
localiser, on peut supposer que $X$ est local de point fermé $x=f(y)$ et
que $y$ est un point fermé dans la fibre $F=f^{-1}(y)$. Notons $c=\codim_X
x$ et $d=\codim_F y$. L'épinglage de $f^\star\Lambda$ en $y$ est donné par
le produit des classes $f^\star (\Cl{x\to X})\in \HH^{2c}_F(Y,\Lambda(c))$ et
$\Cl{y\to F}\in \HH^{2d}_y(F,\Lambda(d))$, ce produit trouvant demeure dans
$\HH^{2c+2d}_y(Y,\Lambda(c+d))$. La compatibilité des classes de Gysin au
changement de base implique que $f^\star(\Cl{x\to X})$ est la classe de
Gysin $\Cl{F\to Y}\in \HH^{2c}_F(Y,\Lambda(c))$. En utilisant la
compatibilité des classes de Gysin à la composition, le produit considéré
plus haut est $\Cl{y\to Y}$, qui est précisément la classe qui définit
l'épinglage de $\Lambda$ en $y$.

\begin{remarque3}
Supposons que l'on dispose de deux morphismes réguliers composables $Z\vers
g Y\vers f X$ entre  $\ZZn$-schémas noethériens excellents. Si $K$ est un
complexe dualisant putatif sur $X$, alors appliquer la construction du
lemme~\ref{17-lemme-construction-cdputatif-f-star}
à $f$, puis à $g$, revient à l'appliquer directement à $f\circ
g$, autrement dit l'isomorphisme évident $g^\star
f^\star K\simeq (f\circ g)^\star K$ dans $\detaleplus Z$ est
compatible aux épinglages. En effet, les constructions s'effectuant fibre à
fibre, on peut supposer que $X$ est le spectre d'un corps et, quitte à
modifier la fonction de dimension, que $K=\Lambda$ ;
en particulier, $X$, $Y$ et $Z$ sont réguliers, les faisceaux
constants $\Lambda$ sur $X$, $Y$ et $Z$ sont naturellement munis de
structures de complexes dualisants putatifs (et même potentiels,
cf.~proposition~\ref{17-proposition-cdpotentiel-regulier}) pour les
fonctions de dimension considérées, la compatibilité requise est
obtenue en appliquant le
lemme~\ref{17-lemme-compatibilite-cdp-f-star-schemas-reguliers} aux
morphismes $f$, $g$ et $f\circ g$.
\end{remarque3}

\subsubsection{Compatibilité aux spécialisations}

Supposons maintenant que $K$ soit un complexe dualisant potentiel.
Par définition des complexes dualisants potentiels, 
il s'agit maintenant de montrer que les épinglages sur $f^\star K$ sont
compatibles aux morphismes de transition associés aux spécialisations
immédiates. D'après le
théorème~\ref{17-theoreme-transition-toute-codimension}, on peut donner un
sens à cette compatiblité pour toute spécialisation
$\overline{y'}\to\overline y$ de points géométriques de $Y$, quelle que
soit sa codimension.

\begin{lemme3}\label{17-lemme-stabilite-cdp-morphisme-regulier}
Plaçons nous dans la situation de la
proposition~\ref{17-proposition-stabilite-cdp-morphisme-regulier}. Soit $K$ 
un complexe dualisant potentiel sur $X$. Soit ${\overline {y'}\to\overline y}$
une spécialisation de points géométriques de $Y$ au-dessus d'une
spécialisation $\overline {x'}\to\overline x$ de points géométriques de
$X$.
Si $\overline{x'}\to\overline x$ est de codimension $0$ ou $1$, alors les
épinglages sur $f^\star K$ sont compatibles à la spécialisation $\overline
{y'}\to\overline y$.
\end{lemme3}

Étant entendu que la construction de la structure de complexe dualisant
putatif sur $f^\star K$ a été réalisée fibre à fibre, on peut 
supposer que $X$ est local strictement hensélien intègre, de point
générique $x'$ et de point fermé $x$. On peut supposer que
$\delta_X(x')=0$.
Le schéma $X$ étant de dimension
$\leq 1$, si on note $n\colon \tilde X\to X$ la normalisation de $X$, le schéma
$\tilde X$ est régulier. On constitue le carré cartésien suivant :
\[
\xymatrix{\tilde Y\ar[r]^{n'}\ar[d]^{\tilde f} & Y\ar[d]^f \\
\tilde{X}\ar[r]^n & X
}
\]
Les morphismes $n$ et $n'$ sont finis surjectifs radiciels. Pour eux, on
dispose des constructions $n^!$ et ${n'}^!$ sur les complexes dualisants
putatif (cf. proposition~\ref{17-proposition-cdputatif-quasi-fini}).
Pour les morphismes
réguliers $f$ et $\tilde f$, on a les
constructions $f^\star$ et ${\tilde f}^\star$. On peut donc considérer les
complexes dualisants putatifs ${\tilde f}^\star n^!K$ et ${n'}^!f^\star K$.
D'après le lemme~\ref{17-lemme-compatibilite-f-star-homeo-universel} à suivre,
on a un isomorphisme canonique de complexes dualisant putatifs 
${\tilde f}^\star n^!K\simeq {n'}^!f^\star K$. En outre, la
proposition~\ref{17-proposition-cdpotentiel-morphisme-fini} montre que pour
montrer que $f^\star K$ est un complexe dualisant potentiel, il suffit de
montrer que ${n'}^!f^\star K$ en est un. Comme il s'identifie
${\tilde f}^\star n^!K$ et que la 
proposition~\ref{17-proposition-cdpotentiel-morphisme-fini} nous dit aussi que 
$n^!K$ est un complexe dualisant potentiel, on peut finalement remplacer
$f$ par $\tilde f$ et supposer que $X$ est
régulier de dimension $\leq 1$.
On peut alors supposer que $K=\Lambda$,
épinglé comme il convient de le faire. Le complexe $f^\star \Lambda$
s'identifie à $\Lambda$ dont on sait qu'il peut-être muni d'une structure
de complexe dualisant potentiel
(cf. proposition~\ref{17-proposition-cdpotentiel-regulier}). Il s'agit de
montrer que les épinglages sur $f^\star \Lambda$ et sur $\Lambda$ sont
compatibles : c'est le sens du
lemme~\ref{17-lemme-compatibilite-cdp-f-star-schemas-reguliers}.

\begin{lemme3}\label{17-lemme-compatibilite-f-star-homeo-universel}
Plaçons nous dans la situation de la
proposition~\ref{17-proposition-stabilite-cdp-morphisme-regulier}. Soit
$g\colon X'\to X$ un morphisme fini surjectif radiciel. On constitue le
carré cartésien :
\[
\xymatrix{
Y'\ar[r]^{g'}\ar[d]^{f'} & Y\ar[d]^f \\
X'\ar[r]^g & X
}
\]
Soit $K$ un complexe dualisant putatif sur $X$. Alors, l'isomorphisme
évident ${f'}^\star g^! K\simeq {g'}^!f^\star K$ dans $\detaleplus {Y'}$ est
compatible aux épinglages.
\end{lemme3}

Les constructions envisagées se réalisant fibre à fibre, on peut supposer
que $X$ et $X'$ sont les spectres de corps notés
respectivement $E$ et $E'$. Les morphismes $f$ et $f'$ étant plats et à
fibres géométriquement régulières, les schémas $Y$ et $Y'$ sont réguliers.
Dans cette situation, compte tenu de la 
proposition~\ref{17-proposition-cdpotentiel-morphisme-fini} et de 
la construction du lemme~\ref{17-lemme-construction-cdputatif-f-star},
il est manifeste que les complexes dualisants putatifs 
${f'}^\star g^! K$ et ${g'}^!f^\star K$ sont des complexes dualisant
potentiels. Pour montrer que l'isomorphisme évident 
${f'}^\star g^! K\simeq {g'}^!f^\star K$ est compatible aux épinglages, il
suffit donc de le faire aux points maximaux de $Y'$.
Bref, on peut supposer que
$Y$ est lui aussi le spectre d'un corps $F$. Comme $Y'$ est régulier et
homéomorphe à $Y$, le schéma $Y'$ est à son tour
le spectre d'un corps $F'$. Compte tenu
de la construction de la proposition~\ref{17-proposition-cdputatif-quasi-fini}
faisant intervenir les morphismes de Gysin associés aux morphismes
d'intersection complète $Y'\to Y$ et $X'\to X$,
la comparaison des deux structures de complexes dualisants putatifs
envisagées sur $Y'$ se ramène à l'égalité des degrés $[F':F]=[E':E]$,
qui résulte aussitôt de la définition de $Y'$ : $F'=E'\otimes_E F$.

\bigskip

Pour finir la démonstration de la
proposition~\ref{17-proposition-stabilite-cdp-morphisme-regulier}, il nous
reste à établir des compatibilités entre les épinglages sur $f^\star K$ et
les spécialisations de points géométriques de $Y$. Les lemmes suivants sur
les spécialisations nous seront utiles.

\begin{lemme3}\label{17-lemme-morphismes-reguliers-connexite}
Soit $f\colon Y\to X$ un morphisme régulier entre
schémas noethériens. On suppose que $f$ est un
morphisme local entre schémas locaux strictements henséliens. Soit
$\overline x\to X$ un point géométrique de $X$. Alors, $f^{-1}(\overline
x)=Y_{\overline x}=Y\times_X\overline x$ est un schéma intègre.
\end{lemme3}

La fibre $Y_x$ est géométriquement régulière et géométriquement connexe
d'après la proposition~\ref{17-proposition-changement-de-base-regulier}
appliquée à un faisceau constant. Le schéma $Y_{\overline x}$ est donc
intègre.

\begin{lemme3}\label{17-lemme-morphismes-reguliers-specialisations-generiques}
Soit $f\colon Y\to X$ un morphisme régulier
entre schémas noethériens excellents. Soit $\overline
{x'}\to\overline x$ une spécialisation de points géométriques de $X$. Soit
$\overline{\eta}$ un point géométrique de $X$ au-dessus de $\overline x$
tel que $\eta$ soit de codimension $0$ dans sa fibre pour $f$. Alors, à un
isomorphisme non nécessairement unique près, il existe une unique
spécialisation $\overline{\eta'}\to\overline\eta$ au-dessus de
$\overline{x'}\to\overline{x}$ telle que $\eta'$ soit de codimension $0$
dans sa fibre. La codimension de $\overline{\eta'}\to\overline\eta$
est la même que celle de $\overline {x'}\to\overline x$.
\end{lemme3}

On peut supposer que $f$ est un morphisme local entre schémas locaux
strictement henséliens $Y$ et $X$ de points fermés respectifs
$\overline\eta$ et $\overline x$. Établir le lemme dans ce cas précis
revient à montrer qu'à isomorphisme près, la fibre géométrique
$Y_{\overline x}$ n'a qu'un seul point géométrique au-dessus d'un point
maximal de $Y_{\overline x}$, ce qui résulte du 
lemme~\ref{17-lemme-morphismes-reguliers-connexite}.
Les formules liant les fonctions de dimension sur $X$ et $Y$ montrent que
la codimension de la spécialisation de 
$\overline{\eta'}\to\overline\eta$ est la même que celle de 
$\overline{x'}\to\overline{x}$.

\begin{lemme3}\label{17-lemme-morphismes-reguliers-relevement-specialisations}
Soit $f\colon Y\to X$ un morphisme régulier
entre schémas noethériens excellents. Soit
$\overline{\eta'}\to\overline\eta$ une spécialisation de points
géométriques de $Y$ au-dessus d'une spécialisation
$\overline{x'}\to\overline{x}$. On suppose que $\eta$ et $\eta'$ sont de
codimension $0$ dans leur fibre pour $f$. Soit
$\overline{x_0}\to\dots\to\overline{x_n}$ une factorisation de la
spécialisation $\overline{x'}\to\overline{x}$ en une suite de
spécialisations immédiates. Alors, on peut décomposer
$\overline{\eta'}\to\overline\eta$ en une suite de spécialisations
immédiates $\overline{\eta_0}\to\dots\to\overline{\eta_n}$ au-dessus de 
$\overline{x_0}\to\dots\to\overline{x_n}$, chacun des points $\eta_i$ étant
de codimension $0$ dans sa fibre pour $f$.
\end{lemme3}

Ceci résulte aussitôt du
lemme~\ref{17-lemme-morphismes-reguliers-specialisations-generiques},
propriété d'unicité comprise.

\begin{lemme3}\label{17-lemme-morphismes-reguliers-generisation-specialisations}
Soit $f\colon Y\to X$ un morphisme régulier
entre schémas noethériens excellents. Soit
$\overline{y'}\to\overline{y}$ une spécialisation de points géométriques de
$Y$ au-dessus de $\overline{x'}\to\overline x$.
Alors, il existe un diagramme commutatif de spécialisations de points
géométriques de $Y$ :
\[
\xymatrix{
\overline{\eta'}\ar[r]\ar[d] & \overline\eta\ar[d] \\
\overline{y'}\ar[r] & \overline y
}
\]
où les points géométriques $\overline{\eta'}$ et $\overline{\eta}$ sont
respectivement au-dessus de $\overline{x'}$ et de $\overline x$, et où
$\eta'$ et $\eta$ sont de codimension $0$ dans leur fibre pour $f$.
\end{lemme3}

On peut supposer que $X$ et $Y$ sont locaux strictement henséliens de
points fermés respectifs $\overline y$ et $\overline x$. On choisit une
spécialisation $\overline\eta\to \overline y$ de points géométriques de la
fibre $Y_{\overline x}$ avec $\eta$ de codimension $0$ dans cette fibre.
D'après le
lemme~\ref{17-lemme-morphismes-reguliers-specialisations-generiques}, il
existe une spécialisation $\overline {\eta'}\to\overline \eta$ de points
géométriques de $Y$ au-dessus de $\overline{x'}\to\overline x$, avec
$\eta'$ de codimension $0$ dans sa fibre pour $f$. Par ailleurs, la fibre
géométrique
$Y_{\overline{x'}}$ étant intègre, on peut en choisir un point géométrique
$\overline{\eta''}$ au-dessus du point générique.
Le point géométrique $\overline{y'}$ de $Y$ étant au-dessus de
$\overline{x'}$, il définit un point géométrique de $Y_{\overline{x'}}$ ;
on dispose donc d'une spécialisation $\overline{\eta''}\to\overline{y'}$ de
points géométriques au-dessus de $\overline{x'}$.
La fibre géométrique $Y_{\overline{x'}}$ étant intègre, il existe un
isomorphisme $\overline{\eta'}\simeq\overline{\eta''}$, ce qui donne le
diagramme commutatif souhaité.

\subsubsection{Fin de la démonstration}

Finissons la démonstration de la
proposition~\ref{17-proposition-stabilite-cdp-morphisme-regulier}, il s'agit
de montrer que pour tout complexe dualisant potentiel $K$ sur $X$, les
épinglages sur le complexe dualisant putatif $f^\star K$ sont compatibles
aux spécialisations de points géométriques de $Y$. Soit
$\overline{y'}\to\overline y$ une telle spécialisation. Le
lemme~\ref{17-lemme-morphismes-reguliers-generisation-specialisations}
s'applique et on
obtient un diagramme commutatif de spécialisations comme ci-dessus. Pour
montrer la compatibilité pour la spécialisation $\overline{y'}\to\overline
y$, il suffit de l'obtenir pour les trois autres spécialisations qui
interviennent dans le carré commutatif. Pour les spécialisations
$\overline{\eta'}\to\overline{y'}$ et $\overline{\eta}\to\overline{y}$,
cela résulte aussitôt des faits observés dans 
la construction même du lemme~\ref{17-lemme-construction-cdputatif-f-star}, à
savoir que chaque fibre de $f$ est munie d'un complexe dualisant potentiel.
On peut appliquer le
lemme~\ref{17-lemme-morphismes-reguliers-relevement-specialisations}
pour obtenir une
factorisation de $\overline{\eta'}\to\overline\eta$ en une suite de
spécialisations immédiates $\overline{\eta_0}\to\dots\to\overline{\eta_n}$
au-dessus d'une composition de spécialisations immédiates
$\overline{x_0}\to\dots\to\overline{\eta_n}$ de points géométriques de $X$.
On peut appliquer le lemme~\ref{17-lemme-stabilite-cdp-morphisme-regulier} aux
spécialisations $\overline{\eta_i}\to\overline{\eta_{i+1}}$ pour obtenir la
compatibilité souhaitée pour $\overline{\eta'}\to\overline\eta$, ce qui
achève la démonstration de la
proposition~\ref{17-proposition-stabilite-cdp-morphisme-regulier}.

\subsection{Démonstration de la proposition~\ref{17-proposition-stabilite-cdp-morphisme-de-type-fini}}
\label{17-subsection-morphismes-de-type-fini}

\begin{lemme2}\label{17-lemme-construction-cdputatif-f-bang}
Dans la situation de la
proposition~\ref{17-proposition-stabilite-cdp-morphisme-de-type-fini},
si $K$ est un complexe dualisant putatif sur $X$,
on peut munir $f^!K$ d'une structure de
complexe dualisant putatif sur $Y$.
\end{lemme2}

Définir la structure de complexe dualisant putatif sur $f^!K$ peut se faire
fibre à fibre. On suppose donc que $X=x$ est le spectre d'un corps.
Soit $y\in Y$. Soit $U$ un ouvert non vide de régularité de l'adhérence
(réduite) $\overline{\{y\}}$. On note $j\colon U\to Y$ l'immersion de $U$
dans $Y$ et $\pi\colon U\to x$ la projection. On a des isomorphisme
canoniques induits par le morphisme de Gysin $\Cl \pi$ et l'épinglage de
$K$ en $x$ :
\[\Lambda(\delta_Y(y))[2\delta_Y(y)]\iso
\pi^!\Lambda(\delta_X(x))[2\delta_X(x)]\giso \pi^!\SheafRG_x(K)
\simeq j^!f^!K\;\text{.}\]
En passant au point générique de $U$, on obtient l'isomorphisme voulu :
$\SheafRG_y(f^!K)\simeq \Lambda(\delta_Y(y))[2\delta_Y(y)]$. Il ne
dépend évidemment pas de l'ouvert $U$, ce qui permet de définir l'épinglage
souhaité de $f^!K$ en $y$.

\medskip

On peut vérifier que la construction de cette structure de complexe
dualisant putatif sur $f^!K$ est compatible à la composition des morphismes
de type fini (pour le sens de cette affirmation, cf.
remarque~\ref{17-remarque-cdp-compatibilite-composition}). Par ailleurs, il
est évident que dans le cas où $f$ est quasi-fini, les épinglages définis
ici sont les mêmes que ceux de la
proposition~\ref{17-proposition-cdputatif-quasi-fini}. Enfin, dans le cas où
$f$ est lisse de dimension relative $d$,
\emph{via} l'isomorphisme canonique $f^\star K\simeq
f^!K(d)[2d]$, les épinglages sont compatibles avec ceux de la
proposition~\ref{17-proposition-stabilite-cdp-morphisme-regulier}, compte tenu
du décalage entre les deux fonctions de dimension envisagées sur $Y$.

\bigskip

Montrons la
proposition~\ref{17-proposition-stabilite-cdp-morphisme-de-type-fini}. La
question étant de nature locale, on peut supposer que $f\colon X\to Y$ est
un morphisme de type fini entre schémas affines. Il existe donc une
factorisation $f=p\circ i$ avec $i$ une immersion fermée et $p$ un
morphisme lisse. Soit $K$ un complexe dualisant potentiel sur $X$.
Les remarques précédentes montrent que $p^!K$,
puis $i^!p^!K$ sont des complexes dualisants
potentiels et que ce dernier s'identifie à $f^!K$.
Par conséquent, $f^!K$ est un
complexe dualisant potentiel, ce qui achève la démonstration.

\section{Existence et unicité des complexes dualisants potentiels}
\label{17-section-existence-unicite-cdp}
\subsection{Énoncé du théorème}

L'objectif de cette section est d'établir le théorème suivant :

\begin{theoreme2}\label{17-theoreme-existence-unicite-cdp}
Soit $X$ un $\ZZn$-schéma excellent muni d'une fonction de dimension
$\delta$. Alors $(X,\delta)$ admet un complexe dualisant potentiel $K_X$,
unique à isomorphisme unique près, et le morphisme évident $
\Lambda\to \tau_{\leq 0}\R\SheafHom(K_X,K_X)$
est un isomorphisme dans $\detale X$. De plus,
$K_X\in\Perv^{-2\delta}(X,\Lambda)$ (cf.
sous-section~\ref{17-subsection-faisceaux-pervers}).
\end{theoreme2}

\subsection{Préliminaires sur les faisceaux pervers}
\label{17-subsection-faisceaux-pervers}

Si $X$ est un schéma noethérien et $p\colon
X\to\mathbf{Z}\cup \{+\infty\}$ une fonction de perversité
(c'est-à-dire que pour toute spécialisation $\overline y\to \overline x$ de
points géométriques de $X$, on a $p(x)\geq p(y)$),
Gabber a défini dans \cite{T@Gabber} une $t$-structure
$(\detaletronque X {\leq p},\detaletronque X {\geq p})$ sur $\detaleplus X$ de sorte que
pour tout $K\in\detaleplus X$, on ait :
\[K\in\detaletronque X {\leq p}\Longleftrightarrow \forall x\in X, K_{|x}\in
\detaletronque x {\leq p(x)}\;\text{,}\]
\[K\in\detaletronque X {\geq p}\Longleftrightarrow \forall x\in X,
\SheafRG_x(K)\in
\detaletronque x {\geq p(x)}\;\text{,}\]
où l'on a muni chacune des catégories $\detaleplus x$ de sa $t$-structure
canonique.

On note $\Perv^p(X,\Lambda)\subset \detaleplus X$ le cœur de cette
$t$-structure, les foncteurs de troncatures étant notés $\tau_{\leq p}$ et
$\tau_{\geq p}$.

\begin{definition2}
Soit $f\colon Y\to X$ un morphisme de schémas. On dit que $f$ est une
pseudo-immersion ouverte si $f$ induit un homéomorphisme sur son image et
que le morphisme induit $f^{-1}\OO_X\to \OO_Y$ est un isomorphisme.
\end{definition2}

Cette classe de morphismes est stable par composition ; elle contient les
immersions ouvertes et les localisations.

\begin{proposition2}\label{17-proposition-faisceaux-pervers-immersions-ouvertes}
Soit $p\colon X\to\mathbf{Z}\cup \{+\infty\}$ une fonction de perversité.
Soit $f\colon Y\to X$ une pseudo-immersion ouverte.
\begin{enumerate}
\item[(a)] La fonction $p\circ f\colon Y\to \mathbf{Z}\cup\{+\infty\}$
est une fonction de perversité (encore notée $p$)
et pour les $t$-structures définies par $p$, $f^\star\colon\detaleplus X\to
\detaleplus Y$ est $t$-exact et $\R f_\star\colon \detaleplus Y\to
\detaleplus X$
$t$-exact à gauche ;
\item[(b)] Le foncteur $f^\star$ induit un foncteur exact $f^\star\colon
\Perv^p(X,\Lambda)\to\Perv^p(Y,\Lambda)$ qui admet un adjoint à droite
$f^p_\star\colon \Perv^p(Y,\Lambda)\to \Perv^p(X,\Lambda)$ défini par la
formule $f_\star^p K=\tau_{\leq p}\R f_\star K$ ;
\item[(c)] Le morphisme d'adjonction $f^\star f^p_\star\to
\Id_{\Perv^p(Y,\Lambda)}$ est un
isomorphisme et le foncteur $f_\star^p\colon \Perv^p(Y,\Lambda)\to
\Perv^p(X,\Lambda)$ est pleinement fidèle.
\item[(d)] Si $g\colon Z\to Y$ est une pseudo-immersion ouverte composable
avec $f$, on dispose d'un isomorphisme de transitivité $f_\star^p\circ
g_\star^p\simeq (f\circ g)_\star^p$.
\end{enumerate}
\end{proposition2}

La $t$-exactitude de $f^\star$ est triviale. En particulier, $f^\star$ est
$t$-exact à droite ; par adjonction, $\R f_\star$ est $t$-exact à
gauche. (b) résulte aussitôt de (a). (c) en résulte aussi compte
tenu du fait que le morphisme d'adjonction $f^\star \R f_\star\to
\Id_{\Perv^p(Y,\Lambda)}$
est un isomorphisme. L'isomorphisme du (d) se déduit par adjonction de
l'isomorphisme de transitivité des foncteurs images inverses.

\subsection{Cas d'un schéma normal}

\begin{proposition2}\label{17-proposition-theoreme-existence-unicite-cdp-cas-normal}
L'énoncé du théorème~\ref{17-theoreme-existence-unicite-cdp} est vrai si on
suppose de plus que le schéma $X$ est normal.
Plus précisément, soit $X$ un $\ZZn$-schéma
noethérien excellent irréductible normal
de point générique $\eta$,
muni d'une fonction de dimension $\delta$ telle que $\delta(\eta)=0$. On
note $j\colon \eta\to X$ l'inclusion du point générique et on pose
$T=j^\varphi_\star \Lambda$ où $\varphi\colon X\to\mathbf{N}$ est la
fonction de perversité définie par l'égalité :
\[\varphi(x)=\max(0,2\dim \OO_{X,x}-2)\]
pour tout $x\in X$.
\begin{enumerate}
\item[(a)] Le morphisme d'adjonction $\Lambda\to j_\star^\varphi
j^\star\Lambda$ définit un morphisme évident $\Lambda\to T $ tel que pour
tout point géométrique $\overline x\to X$, l'application induite
\[\HH^{-2\delta(x)}_{\overline x}(\Lambda)\to \HH^{-2\delta(x)}_{\overline x}(T)\]
soit un isomorphisme et que $\HH^q_{\overline x}(T)=0$ si $q\neq -2\delta(x)$.
Le théorème~\ref{17-theoreme-classe-degre-maximal} donne alors un épinglage
de $T$ en $x$. Avec ces épinglages, $T$ est un complexe dualisant
potentiel pour $(X,\delta)$.
\item[(b)] Si $K$ est un complexe dualisant potentiel pour $(X,\delta)$,
alors $K$ appartient à $\Perv^{\varphi}(X,\Lambda)$
(et à $\Perv^{-2\delta}(X,\Lambda)$)
et le morphisme $K\to T$ qui s'en déduit aussitôt est un isomorphisme
compatible aux épinglages.
\item[(c)] Si $K$ est un complexe dualisant potentiel sur $X$, alors le
morphisme évident $\Lambda\to\tau_{\leq 0}\R \SheafHom(K,K)$ est un
isomorphisme.
\end{enumerate}
\end{proposition2}

Établissons (a). Le point essentiel est de montrer que pour tout $x\in X$,
le morphisme évident
$\HH^{-2\delta(x)}_{\overline x}(\Lambda)\to \HH^{-2\delta(x)}_{\overline x}(T)$
est un isomorphisme et que $\HH^q_{\overline x}(T)$ est nul si $q\neq
-2\delta(x)$. Le fait que $T$ soit un complexe dualisant
potentiel résultera alors aussitôt du
théorème~\ref{17-theoreme-classe-degre-maximal}. Pour établir le résultat
voulu, on peut supposer que $X$ est local strictement hensélien de
point fermé $x$. Par récurrence sur $d=\dim X$, on peut supposer (a) connu
sur l'ouvert $U=X-x$. Si $d=0$, le résultat est évident. On suppose donc
$d\geq 1$, de sorte que $U$ contienne le point générique de $X$.
Posons $L=T_{|U}$. D'après la
proposition~\ref{17-proposition-faisceaux-pervers-immersions-ouvertes}, si on
note $g\colon U\to X$ l'immersion ouverte évidente, on a un isomorphisme
canonique $T=g_\star^\varphi L$. On en déduit aussitôt un isomorphisme
canonique $T=\tau_{\leq \varphi(x)}\R g_\star L$.

Si $d=1$, $\varphi(\eta)=0$, donc $T=g_\star\Lambda=\Lambda$. La
proposition~\ref{17-proposition-cdp-trait} permet de conclure que $T$ est bien
un complexe dualisant potentiel avec les épinglages envisagés ici. On
suppose donc que $d\geq 2$. Dans ce cas, on a $\varphi(x)=2d-2$ et
$T=\tau_{\leq 2d-2}\R g_\star L$.

D'après le lemme~\ref{17-lemme-structure-cohomologie-negative-cdp-normal}, la
structure des objets de cohomologie ${\mathscr H}^q L$ pour $q\leq 0$ est
connue. Notons $i\colon x\to X$ l'immersion du point fermé de $X$.
Utilisons le triangle distingué canonique :
\[i_\star i^! T\to T\to \R g_\star L\to i_\star i^!T[1]\;\text{.}\]
Il en résulte une suite exacte longue :
\[\dots \to \HH^{q-1}(U,L)\to \HH^q_x(X,T)\to
({\mathscr H}^qT)_x \to \HH^q(U,L)\to \dots\]
Par construction, $({\mathscr H}^q T)_x\iso \HH^q(U,L)$ si $q\leq 2d-2$ et 
$({\mathscr H}^q T)_x=0$ si $q\geq 2d-1$. Il en résulte que $\HH^q_x(X,T)=0$
si $q\leq 2d-1$ et que l'on a des isomorphismes $\HH^{q-1}(U,L)\iso
\HH^q_x(X,T)$ pour $q\geq 2d$.
Il vient aussi que ${\mathscr H}^q T=0$ si $q<0$ et que l'on a un
isomorphisme canonique $\Lambda\iso {\mathscr H}^0T$. Ainsi, on a un
morphisme canonique
$\Lambda \to T$ induisant un isomorphisme $\Lambda\iso
\tau_{\leq 0}T$.

\begin{lemme2}\label{17-lemme-technique-dimcohom-existence-unicite-cdp}
Soit $U$ le complémentaire du point fermé dans un schéma local strictement
hensélien excellent normal $X$ de dimension $d\geq 2$.
Soit $M\in\detaletronque U {\leq \varphi}$. On suppose qu'il existe un isomorphisme 
$\Lambda\iso \tau_{\leq 0}M$. Alors, pour tout $q\geq 2d$, le morphisme
évident est un isomorphisme :
\[\HH^{q-1}(U,\Lambda)\iso \HH^{q-1}(U,M)\;\text{.}\]
\end{lemme2}

On note $M^+$ un cône du morphisme évident $\Lambda\to M$. Il suffit de
montrer que $\HH^q(U,M^+)=0$ pour $q\geq 2d-2$. Les hypothèses
impliquent que les objets de cohomologie de $M^+$ sont nuls en dehors
de l'intervalle $[1,2d-4]$. On peut aussi observer que pour tout 
$1\leq i\leq d-2$, si $y\in U$ est tel que $({\mathscr
H}^{2i-1} M^+)_{\overline y}$ ou $({\mathscr H}^{2i} M^+)_{\overline y}$
soit non nul,
alors l'adhérence (dans $X$) de $y$ est de dimension $\leq d-i-1$.
D'après la proposition~\ref{17-proposition-dimension-cohomologique} et la
suite spectrale d'hypercohomologie, on obtient bien que $\HH^q(U,M^+)=0$ si
$q\geq 2d-2$.

\bigskip

On peut appliquer le
lemme~\ref{17-lemme-technique-dimcohom-existence-unicite-cdp}
avec $M=L$. Il résulte alors de ce qui précède et du
théorème~\ref{17-theoreme-classe-degre-maximal} que $\HH^q_x(X,T)=0$ si $q\neq
2d$, que l'on a un isomorphisme canonique 
$\HH^{2d}_x(X,\Lambda(d))\iso \HH^{2d}_x(X,T(d))$, et que l'on peut ainsi 
définir des épinglages sur $T$ qui en font un complexe dualisant potentiel
sur $X$.

\bigskip

Montrons maintenant (b). Soit $K$ un complexe dualisant potentiel sur $X$.
Il est évident que $K\in\detaletronque X {\geq -2\delta}\subset
\detaletronque X {\geq
\varphi}$. Pour montrer que $K\in \detaletronque X {\leq \varphi}\subset
\detaletronque X
{\leq -2\delta}$, on peut supposer que $X$ est local strictement hensélien
de dimension $d$.
On note $i\colon x\to X$ l'inclusion du point fermé
et $g\colon U\to X$ l'inclusion de l'ouvert complémentaire $X-x$. Le cas où
$d=0$ étant trivial et celui où $d=1$ ayant été traité dans la 
proposition~\ref{17-proposition-cdp-trait}, on peut supposer que $d\geq 2$.
Par récurrence sur $d$, on peut supposer (b) connu pour l'ouvert $U$. On
note $L=K_{|U}\in\detaletronque U {\leq\varphi}$.

Comme $K$ est un complexe dualisant potentiel, 
la structure des faisceaux de cohomologie ${\mathscr H}^qK$ pour $q\leq 0$
est connue (cf.
lemme~\ref{17-lemme-structure-cohomologie-negative-cdp-normal}). On peut 
donc appliquer le
lemme~\ref{17-lemme-technique-dimcohom-existence-unicite-cdp} avec $M=L$.
Ainsi, le morphisme évident est un isomorphisme $\HH^{2d-1}(U,\Lambda)\iso
\HH^{2d-1}(U,L)$, et $\HH^q(U,L)=0$ pour $q\geq 2d$.

De même que pour établir (a), on utilise le triangle distingué
canonique :
\[i_\star i^! T\to T\to \R g_\star L\to i_\star i^!T[1]\;\text{,}\]
et la suite exacte longue qui s'en déduit :
\[\dots \to \HH^{q-1}(U,L)\to \HH^q_x(X,K)\to
({\mathscr H}^qK)_x \to \HH^q(U,L)\to \dots\]
La structure de $\HH^q_x(X,K)$ étant connue pour tout $q\in\mathbf{Z}$,
il vient aussitôt que $({\mathscr H}^qK)_x=0$ pour $q\geq 2d+1$ et que
l'on a une suite exacte
\[0\to ({\mathscr H}^{2d-1}K)_x\to \HH^{2d-1}(U,L)\to \HH^{2d}_x(X,K)
\to ({\mathscr H}^{2d}K)_x \to 0\;\text{.}\]
Pour montrer que $K\in\detaletronque X {\leq \varphi}$, il reste donc à montrer que
le morphisme canonique $\HH^{2d-1}(U,L)\to \HH^{2d}_x(X,K)$ est un
isomorphisme. Choisissons une spécialisation $\overline{\eta}\to x$. On
peut considérer le diagramme commutatif suivant, où les flèches verticales
sont induites par le morphisme canonique $\Lambda\to K$ :
\[
\xymatrix{
\ar[d]^{\sim}
\HH^0_{\overline {\eta}}(X,\Lambda)\ar[r]^-{\transition^X_{\overline{\eta}\to x}} &
   \HH^{2d}_x(X,\Lambda(d))\ar[d] & \ar[l]_{\sim}
   \HH^{2d-1}(U,\Lambda(d))\ar[d]^{\sim} \\
\HH^0_{\overline {\eta}}(X,K)\ar[r]^-{\transition^X_{\overline{\eta}\to x}} &
   \HH^{2d}_x(X,K(d)) & \ar[l] \HH^{2d-1}(U,L(d))
}\]
Les flèches de la colonne de gauche vers celle du milieu sont les 
morphismes de transition introduits au
théorème~\ref{17-theoreme-transition-toute-codimension}. Ici, ce sont des
isomorphismes : pour la flèche du haut, cela résulte du
théorème~\ref{17-theoreme-classe-degre-maximal} et pour la flèche du bas, du
fait que $K$ soit un complexe dualisant potentiel. On en déduit que le
morphisme du milieu $\HH^{2d}_x(X,\Lambda(d))\to \HH^{2d}_x(X,K(d))$ est un
isomorphisme. Compte tenu des autres isomorphismes connus, il vient que le
morphisme évident $\HH^{2d-1}(U,L(d))\to \HH^{2d}_x(X,K(d))$ est un
isomorphisme, ce qui achève de montrer que $K\in \Perv^\varphi(X,\Lambda)$.

On a alors un morphisme d'adjonction $K\to j_\star^\varphi j^\star K=T$
dans $\Perv^\varphi(X,\Lambda)$. Pour montrer que c'est un isomorphisme, on
peut se placer dans la situation locale précédente, et faire une récurrence
sur la dimension pour pouvoir supposer que le morphisme induit $K_{|U}\to
T_{|U}$ est un isomorphisme. Choisissons une spécialisation
$\overline{\eta}\to x$. On en déduit un diagramme commutatif :
\[
\xymatrix{\ar[d]\HH^0_{\overline
\eta}(X,K)\ar[r]^-{\transition^X_{\overline{\eta}\to x}} & \HH^{2d}_x(X,K(d)) \ar[d]\\
\HH^0_{\overline
\eta}(X,T)\ar[r]^-{\transition^X_{\overline{\eta}\to x}} & \HH^{2d}_x(X,T(d))
}
\]
La flèche de gauche est évidemment un isomorphisme. Les flèches
horizontales aussi puisque $K$ et $T$ sont des complexes dualisants
potentiels. Il vient donc que le morphisme induit $i^!K\to i^!T$ est un
isomorphisme. Il en découle que $K\to T$ est un isomorphisme (compatible
aux épinglages au point générique, donc à tous les épinglages),
ce qui achève la démonstration de (b).

Montrons (c). Soit $K$ un complexe dualisant potentiel sur $X$. Pour
montrer que le morphisme canonique $\Lambda\to\tau_{\leq
0}\R\SheafHom(K,K)$ est un isomorphisme, quitte à remplacer $X$ par
des schémas connexes étales sur lui,
il suffit de montrer que $\Lambda\iso
\Hom_{\detaleplus X}(K,K)$ et que $\Hom_{\detaleplus X}(K,K[q])=0$ pour $q<0$.
L'annulation de $\Hom_{\detaleplus X}(K,K[q])$ pour $q<0$ résulte aussitôt du
fait que $K$ appartienne au cœur de la $t$-structure définie par $\varphi$.
D'après (b), on a un isomorphisme canonique $K=j_\star^\varphi \Lambda$.
L'isomorphisme $\Lambda\iso\Hom_{\detaleplus X}(K,K)$ résulte alors de ce
que $j_\star^\varphi\colon \Perv^\varphi(\eta,\Lambda)\to
\Perv^\varphi(X,\Lambda)$ soit pleinement fidèle
(cf. proposition~\ref{17-proposition-faisceaux-pervers-immersions-ouvertes}).

\subsection{Un résultat de recollement}
\label{17-subsection-recollement-cdp}

La proposition suivante est un résultat de recollement qui nous permettra
de passer du cas normal au cas général :

\begin{proposition2}\label{17-proposition-existence-unicite-cdp-recollement}
Soit $X$ un $\ZZn$-schéma noethérien excellent muni d'une fonction de
dimension $\delta$. On suppose donné un carré cartésien
\[
\xymatrix{
Y'\ar[d]^{p'}\ar[r]^{i'} & X'\ar[d]^p \\
Y\ar[r]^i & X }
\]
où $i$ est une immersion fermée d'ouvert complémentaire $U$
et où $p$ est un morphisme fini induisant
un isomorphisme $p^{-1}(U)\iso U$. On note $q=p\circ i'$. On suppose que
l'énoncé du théorème~\ref{17-theoreme-existence-unicite-cdp} est connu pour
$X'$, $Y$ et $Y'$ (relativement aux fonctions de dimensions déduites de
celle sur $X$ par le procédé de la
proposition~\ref{17-proposition-cdputatif-quasi-fini}). Alors, cet énoncé est
également vrai pour $X$, et si on note $K_X$, $K_{X'}$,
$K_Y$ et $K_{Y'}$ les complexes dualisants potentiels de $X$, $X'$, $Y$ et
$Y'$ respectivement, on a un triangle distingué canonique dans $\detaleplus X$
:
\[q_\star K_{Y'}\to i_\star K_Y\oplus p_\star K_{X'}\to K_X\to q_\star
K_{Y'}[1]\;\text{.}\]
\end{proposition2}

Dans un premier temps, supposons que $X$ admette un complexe dualisant
potentiel $K_X$ et montrons que l'on peut définir un triangle distingué
canonique de la forme souhaitée. Pour cela, considérons la suite exacte
courte évidente de faisceaux sur $X$ :
\[0\to \Lambda\vers {(+,+)} i_\star\Lambda\oplus p_\star\Lambda\vers{(+,-)}
q_\star\Lambda\to 0\;\text{.}\]
En appliquant $\R\SheafHom(-,K_X)$ au triangle distingué
correspondant, on obtient un triangle distingué :
\[q_\star q^!K_X\vers {(+,-)}
i_\star i^! K_X\oplus p_\star p^!K_X\vers{(+,+)}
K_X\to q_\star q^!K_X[1]\;\text{,}\]
qui, compte tenu de la
proposition~\ref{17-proposition-stabilite-cdp-morphisme-de-type-fini} et de la
vertu d'unicité des complexes dualisants potentiels sur $X'$, $Y$ et $Y'$,
se récrit sous la forme :
\[q_\star K_{Y'}\vers {(+,-)} i_\star K_Y\oplus p_\star K_{X'}\vers{(+,+)}
K_X\to q_\star K_{Y'}[1]\;\text{.}\]

Revenant aux hypothèses de la proposition, nous allons montrer
qu'inversement, si on définit $K_X$ de façon à avoir un tel triangle
distingué (mais \emph{a priori} pas de façon canonique),
on obtient bien un complexe dualisant potentiel sur $X$. On
suppose donc le théorème~\ref{17-theoreme-existence-unicite-cdp} connu
seulement pour $X'$, $Y$ et $Y'$ et on note $K_{X'}$, $K_Y$ et $K_{Y'}$ les
complexes dualisants potentiels correspondants. La propriété d'unicité pour
les complexes dualisants potentiels sur $Y'$ donne des isomorphismes
canoniques $K_{Y'}\simeq {p'}^!K_Y$ et $K_{Y'}\simeq {i'}^!K_{X'}$. Par
adjonction, on en déduit des morphismes canoniques 
$p'_\star K_{Y'}\to K_Y$ et $i'_\star K_{Y'}\to K_{X'}$, puis en appliquant
respectivement $i_\star$ et $p_\star$, on obtient des morphismes canoniques
$q_\star K_{Y'}\to i_\star K_Y$ et $q_\star K_{Y'}\to p_\star
K_{X'}$. On peut considérer leur différence et constituer
un triangle distingué dans $\detaleplus X$ :
\[q_\star K_{Y'}\vers {(+,-)} i_\star K_Y\oplus p_\star K_{X'}\to
K_X\to q_\star K_{Y'}[1]\;\text{.}\]
On obtient ainsi un objet $K_X\in\detaleplus X$ et deux morphismes privilégiés
$i_\star K_Y\to K_X$ et $p_\star K_{X'}\to K_X$.

Notons $j\colon U\to X$ et $j'\colon U\to X'$ les immersions ouvertes
évidentes. En appliquant $j^\star$ au triangle distingué ci-dessus, on peut
commencer par observer que le morphisme évident ${j'}^\star K_{X'}\simeq
j^\star p_\star K_{X'}\to j^\star K_X$ est un isomorphisme. Par conséquent,
on peut munir $K_X$ d'épinglages en les points de l'ouvert $U$ de façon
compatible avec la structure de complexe dualisant potentiel obtenue sur
${j'}^\star K_{X'}$.

Considérons le morphisme canonique $q_\star K_{Y'}\to p_\star K_{X'}$.
J'affirme qu'il induit un isomorphisme après application du foncteur $i^!$.
En effet, on a des isomorphismes évidents de foncteurs $i^!q_\star
{i'}^!\simeq i^!i_\star p'_\star {i'}^\star\simeq p'_\star {i'}^\star\simeq
i^!p_\star$ et, compte tenu de l'isomorphisme canonique $K_{Y'}\simeq
{i'}^!K_{X'}$, on obtient bien que le morphisme évident $i^!q_\star
K_{Y'}\to i^!p_\star K_{X'}$ est un isomorphisme. En appliquant $i^!$ au
triangle distingué de définition de $K_X$, il vient alors que le morphisme
canonique $i_\star K_Y\to K_X$ induit un isomorphisme $K_Y\iso i^! K_X$
après application de $i^!$. Ceci permet de définir des épinglages pour
$K_X$ en tous les points de $Y$.

Finalement, on a obtenu une structure de complexe dualisant putatif sur
$K_X$. D'après la proposition~\ref{17-proposition-cdpotentiel-morphisme-fini},
pour montrer que $K_X$ est un complexe dualisant potentiel, il suffit de
montrer que $p^!K_X$ en est un. Nous allons bien évidemment le comparer à
$K_{X'}$. Par construction de $K_X$, on a un diagramme commutatif :
\[
\xymatrix{
& i_\star K_Y \ar[rd]& \\
q_\star K_{Y'}\ar[ur]\ar[dr] & & K_X \\
& p_\star K_{X'}\ar[ur] & 
}
\]
Par adjonction, on obtient que deux définitions concurrentes d'un morphisme
$K_{Y'}\to q^!K_X$ coïncident :
\[
\xymatrix{
& {p'}^! K_Y\ar[r] & {p'}^!i^!K_X & \\
K_{Y'}\ar@{=}[ru]\ar@{=}[rd] & & & q^!K_X\ar@{=}[lu]\ar@{=}[ld] \\
& {i'}^! K_{X'}\ar[r] & {i'}^!p^!K_X & 
}
\]
On a déjà montré que le morphisme canonique $K_Y\to i^! K_X$ était un
isomorphisme. Par conséquent, sur le diagramme ci-dessus, toutes les
flèches sont des isomorphismes. Ainsi, le morphisme évident 
$ K_{X'}\to p^!K_X$ induit un isomorphisme non seulement après
application de $j^\star$, mais aussi après celle de $i^!$. Il en résulte
que ce morphisme $K_{X'}\to p^! K_X$ est un isomorphisme. En outre, sur le
diagramme ci-dessus, tous les objets sont naturellement munis d'une
structure de complexe dualisant putatif et tous les isomorphismes, sauf
peut-être celui du bas, sont compatibles aux épinglages. Cet isomorphisme 
${i'}^! K_{X'}\to {i'}^!p^!K_X$ est donc lui aussi compatible aux
épinglages. Par conséquent, l'isomorphisme $K_{X'}\iso p^!K_X$ est
compatible aux épinglages non seulement sur $U$ mais aussi sur $Y'$.
Il en résulte que
$p^!K_X$ est un complexe dualisant potentiel, et on a vu que l'on pouvait
en conclure que $K_X$ en était un aussi.

En outre, il résulte aussitôt de la construction que l'hypothèse selon
laquelle les complexes dualisants potentiels sur $X'$, $Y$ et $Y'$ sont
pervers pour la fonction de perversité $-2\delta$ que $K_X$ est aussi
pervers pour $-2\delta$.

Pour conclure, il s'agit de montrer que si
$K$ et $L$ deux complexes dualisants potentiels sur $X$,
on a un isomorphisme 
privilégié $\Lambda\iso \tau_{\leq 0}\R\SheafHom(K,L)$ qui
donne naissance à un morphisme $\psi\colon K\to L$ 
qui soit un isomorphisme de
complexes dualisants potentiels. En effet, cela montrera que si $\phi\colon
K\to L$ est un autre isomorphisme, alors $\phi=\lambda\cdot \psi$ où
$\lambda\colon X\to \Lambda^\times$ est une fonction localement constante.
Demander que $\phi$ soit compatible aux épinglages impliquant que
$\lambda=1$, on aura bien un unique isomorphisme $K\iso L$ de
complexes dualisants potentiels.

D'après la propriété d'unicité des complexes dualisants potentiels sur
$X'$, $Y$ et $Y'$, on a des isomorphismes de complexes dualisants
potentiels $p^!K\iso p^!L$ et $i^! K\iso i^!L$ induisant le
même isomorphisme $q^!K\iso q^!L_Y$. On a construit plus haut un
triangle distingué canonique :
\[q_\star q^\star K\vers{(+,-)} i_\star i^! K\oplus p_\star p^!
K\vers{(+,+)} K\to q_\star q^\star K[1]\;\text{.}\]
En lui appliquant $\R \SheafHom(-,L)$, on obtient un autre triangle
distingué :
\[\R\SheafHom(K,L)\vers{(+,+)} i_\star
\R\SheafHom(i^!K,i^!L)\oplus p_\star \R\SheafHom(p^!K,p^!L)
\vers {(+,-)} q_\star \R\SheafHom(q^!K,q^!L)\vers +
\]
L'énoncé du théorème~\ref{17-theoreme-existence-unicite-cdp} pour $X'$, $Y$ et
$Y'$ implique aussitôt que les objets de cohomologie 
${\mathscr H}^q\R\SheafHom(K,L)$ sont nuls pour $q<0$, et, compte
tenu de la suite exacte canonique de faisceaux :
\[0\to \Lambda\vers {(+,+)} i_\star\Lambda\oplus p_\star\Lambda\vers{(+,-)}
q_\star\Lambda\to 0\;\text{,}\]
que l'on a un isomorphisme privilégié $\Lambda\iso {\mathscr
H}^0\R\SheafHom(K,L)$. Le morphisme $K\to L$ correspondant induit
bien entendu les uniques isomorphismes de complexes dualisants potentiels
$i^!K\iso i^!L$ et $p^!K\iso p^!L$. Ce morphisme $K\to L$ est donc
compatible aux épinglages sur $Y$ et sur $U$ : c'est un isomorphisme de
complexes dualisants potentiels. Ceci achève la démonstration de la
proposition~\ref{17-proposition-existence-unicite-cdp-recollement}.

\subsection{Cas général}

Montrons le théorème~\ref{17-theoreme-existence-unicite-cdp} dans le cas
général. Pour le montrer pour tous les $\ZZn$-schémas noethériens
excellents munis d'une fonction de dimension, par récurrence noethérienne,
on peut supposer le résultat connu pour les schémas finis sur un fermé
d'intérieur vide de $X$. On peut supposer que $X$ est réduit. Notons
$p\colon X'\to X$ la normalisation de $X$. Le morphisme $p$ est fini
surjectif et induit un isomorphisme au-dessus de l'ouvert dense de
normalité $U$ du schéma excellent $X$. Posons $Y=(X-U)_\red$ et formons
le carré cartésien suivant :
\[
\xymatrix{Y'\ar[r]^{i'} \ar[d]^{p'}& X'\ar[d]^p \\ Y\ar[r]^i & X}
\]
Comme $X'$ est normal, la 
proposition~\ref{17-proposition-theoreme-existence-unicite-cdp-cas-normal}
montre que l'énoncé du théorème~\ref{17-theoreme-existence-unicite-cdp} est
connu pour $X'$. L'hypothèse de récurrence noethérienne montre que c'est
aussi le cas pour $Y$ et $Y'$. La 
proposition~\ref{17-proposition-existence-unicite-cdp-recollement} donne la
conclusion souhaitée pour $X$.

\section{Le théorème de dualité locale}
\label{17-section-dualite}

\subsection{Énoncé du théorème}
\begin{theoreme2}\label{17-theoreme-dualite-locale}
Soit $X$ un $\ZZn$-schéma noethérien excellent 
muni d'une fonction de dimension $\delta$.
Soit $K$ le complexe dualisant potentiel de $(X,\delta)$. Alors
\begin{itemize}
\item $K\in\detalebctf X$ ;
\item $K$ est de dimension quasi-injective finie si et seulement si $X$ est
de dimension de Krull finie ;
\item le foncteur $\Dual_K=\R\SheafHom(-,K)$ préserve $\detalebc X$ ;
\item pour tout $M\in\detalebc X$, le morphisme de bidualité $M\to
\Dual_K\Dual_K M$ est un isomorphisme.
\end{itemize}
En particulier, si $X$ est de dimension de Krull finie, $K$ est un complexe
dualisant au sens de \sga{5}{i}{1.7}.
\end{theoreme2}

\subsection{Constructibilité, tor-dimension, dimension quasi-injective}
\subsubsection{Changement de coefficients}

\begin{proposition3}\label{17-proposition-cdp-changement-coefficients}
Soit $\Lambda=\mathbf{Z}/n\mathbf{Z}$. Soit $m$ un diviseur de $n$. Soit
$\Lambda'=\mathbf{Z}/m\mathbf{Z}$. L'anneau $\Lambda'$ est une
$\Lambda$-algèbre. Soit $K\in\detaleplus X$ un complexe dualisant potentiel
sur $X$
relativement à l'anneau de coefficients $\Lambda$. Alors
$K'=\R\SheafHom_\Lambda(\Lambda',K)\in\detalepluscoefs X {\Lambda'}$ est
naturellement muni d'une structure de complexe dualisant potentiel
relativement à l'anneau de coefficients $\Lambda'$. De plus, si
$M\in\detalepluscoefs X {\Lambda'}$, on a un isomorphisme canonique dans
$\detaleplus X$ :
\[\R\SheafHom_\Lambda(M,K)\simeq
\R\SheafHom_{\Lambda'}(M,K')\;\text{.}\]
\end{proposition3}

Ceci résulte facilement de la commutation des foncteurs de cohomologie à
supports avec le foncteur $\R\SheafHom_\Lambda(\Lambda',-)$.

\subsubsection{Constructibilité, tor-finitude}

\begin{proposition3}\label{17-proposition-constructibilite-cdp}
Soit $K$ un complexe dualisant potentiel sur un $\ZZn$-schéma noethérien
excellent muni d'une fonction de dimension $\delta$. Alors $K\in\detalebc X$.
\end{proposition3}

On va utiliser le lemme suivant :

\begin{lemme3}\label{17-lemme-devissage-constructibilite}
Soit $X$ un $\ZZn$-schéma noethérien excellent. Soit $i\colon Z\to X$ une
immersion fermée. Soit $j\colon U\to X$ l'immersion ouverte complémentaire.
Soit $K\in\detaleplus X$. Les conditions suivantes sont équivalentes :
\begin{itemize}
\item $K\in\detalebc X$ ;
\item $j^\star K\in\detalebc U$ et $i^!K\in\detalebc Z$.
\end{itemize}
\end{lemme3}

Ce lemme résulte aussitôt du fait non trivial que le foncteur
$\R j_\star $
envoie $\detalebc U$ dans $\detalebc X$, cf.~\refext{XIII}{13-constructibilité}.

\bigskip

Démontrons la proposition. On peut supposer $X$ réduit. Comme $X$ est
excellent, $X$ admet un ouvert dense régulier. Notons $Z$ le sous-schéma
fermé $(X-U)_{\red}$ et $i\colon Z\to X$ son immersion fermée dans $X$.
Par récurrence noethérienne,
on peut supposer que le complexe dualisant potentiel $i^!K$ de $Z$ est dans
$\detalebc Z$. En vertu du lemme, on est ramené à montrer que $j^\star
K\in\detalebc U$. Autrement dit, on peut supposer que $X$ est régulier. On peut
supposer de plus que $X$ est connexe. Notons $\eta$ le point générique de
$X$. D'après la proposition~\ref{17-proposition-cdpotentiel-regulier} et le
théorème~\ref{17-theoreme-existence-unicite-cdp}, on a un
isomorphisme canonique $K\simeq\Lambda(\delta(\eta))[2\delta(\eta)]$.
Ainsi, $K$ appartient bien à $\detalebc X$, ce qui achève la démonstration de
la proposition.

\begin{proposition3}
Soit $X$ un schéma noethérien excellent muni d'une fonction de dimension
$\delta$. Le complexe dualisant potentiel de $(X,\delta)$ appartient à
$\detalebctf X$.
\end{proposition3}

On sait déjà que le complexe dualisant potentiel $K$ de $(X,\delta)$
appartient à $\detalebc
X$. Il s'agit donc d'obtenir un résultat de tor-finitude pour $K$. Pour
cela, on peut supposer que $\Lambda=\mathbf{Z}/\ell^\nu\mathbf{Z}$ où
$\ell$ est un nombre premier et $\nu\geq 1$. D'après la
proposition~\ref{17-proposition-cdp-changement-coefficients},
$\R\SheafHom_\Lambda(\mathbf{Z}/\ell\mathbf{Z},K)$ est un complexe
dualisant potentiel relativement à l'anneau de coefficients
$\mathbf{Z}/\ell\mathbf{Z}$ ; d'après la 
proposition~\ref{17-proposition-constructibilite-cdp},
cet objet appartient à $\detalebccoefs X {\mathbf{Z}/\ell\mathbf{Z}}$, le
critère du lemme suivant permet de conclure que $K$ appartient à
$\detalebtf X$.

\begin{lemme3}\label{17-lemme-tor-finitude-z-l-n}
Soit $X$ un schéma noethérien. Soit $\ell$ un nombre premier. Soit $\nu\geq
1$. On pose $\Lambda={\mathbf{Z}/\ell^\nu\mathbf{Z}}$.
Soit $K\in\detaleb X$. Les conditions
suivantes sont équivalentes :
\begin{itemize}
\item[(i)] $K\in\detalebtf X$ ;
\item[(ii)] $K\otimes^{\L}_{\Lambda} \mathbf{Z}/\ell\mathbf{Z}\in
\detalebcoefs X {\mathbf{Z}/\ell\mathbf{Z}}$ ;
\item[(iii)]
$\R\SheafHom_\Lambda(\mathbf{Z}/\ell\mathbf{Z},K)\in\detalebcoefs X
{\mathbf{Z}/\ell\mathbf{Z}}$.
\end{itemize}
Dans ce cas, on a un isomorphisme canonique
$\R\SheafHom_\Lambda(\mathbf{Z}/\ell\mathbf{Z},K)\simeq
K\otimes_{\Lambda}^{\L}\mathbf{Z}/\ell\mathbf{Z}$ dans $\detalebcoefs X
{\mathbf{Z}/\ell\mathbf{Z}}$.
\end{lemme3}

L'équivalence entre (i) et (ii) est simplement indiquée pour mémoire. Il
s'agit ici de montrer que les conditions (ii) et (iii) sont équivalentes.
On représente $K$ par un complexe borné. On peut considérer le complexe
complexe double dont les colonnes sont représentées ci-dessous (les degrés
horizontaux étant indiqués en exposant) :
\[\dots\vers {\ell^{\nu-1}} \overset {-1} K\vers {\ell}
\overset 0 K \vers {\ell^{\nu-1}} \overset 1 K\vers
{\ell}\overset 2 K \vers {\ell^{\nu-1}} \dots \]
On note $C$ le complexe simple de faisceaux de
$\Lambda$-modules sur $X$ associé à ce complexe
double\;\footnote{Ce complexe étant concentré sur un nombre fini de lignes,
il n'y a pas lieu de préciser si l'on définit le complexe simple en termes
de sommes ou de produits.}. De façon évidente, on a un triangle distingué
dans $\detale X$ :
\[
C\vers p K\otimes_{\Lambda}^{\L}\mathbf{Z}/\ell\mathbf{Z}\to
\R\SheafHom_\Lambda(\mathbf{Z}/\ell\mathbf{Z},K)\vers i
C[1]\;\text{.}
\]

On remarque que :
\begin{itemize}
\item ${\mathscr H}^qC\simeq {\mathscr H}^{q+2}C$ pour tout
$q\in\mathbf{Z}$ ;
\item $p$ induit un isomorphisme sur les objets de cohomologie
${\mathscr H}^q$ pour $q$ suffisamment petit ;
\item $i$ induit un isomorphisme sur les objets de cohomologie
${\mathscr H}^q$ pour $q$ suffisamment grand.
\end{itemize}

Comme
$K\otimes_{\Lambda}^{\L}\mathbf{Z}/\ell\mathbf{Z}\in\detalemoinscoefs X 
{\mathbf{Z}/\ell\mathbf{Z}}$ et 
$\R\SheafHom_\Lambda(\mathbf{Z}/\ell\mathbf{Z},K)\in \detalepluscoefs X 
{\mathbf{Z}/\ell\mathbf{Z}}$, il en résulte aussitôt que les conditions
suivantes sont équivalentes :
\begin{itemize}
\item $C\simeq 0$ ;
\item $\R\SheafHom_\Lambda(\mathbf{Z}/\ell\mathbf{Z},K)\in\detaleb X$ ;
\item $K\otimes_{\Lambda}^{\L}\mathbf{Z}/\ell\mathbf{Z}\in\detaleb X$.
\end{itemize}

Si on suppose que $K\in\detalebctf X$, 
pour obtenir l'isomorphisme 
$\R\SheafHom_\Lambda(\mathbf{Z}/\ell\mathbf{Z},K)\simeq
K\otimes_{\Lambda}^{\L}\mathbf{Z}/\ell\mathbf{Z}$ dans $\detalebcoefs X
{\mathbf{Z}/\ell\mathbf{Z}}$, on peut représenter $K$ par un complexe borné
de $\Lambda$-modules plats. On peut conclure en utilisant l'isomorphisme
évident $K/\ell K \iso {}_\ell K$ induit 
par la multiplication par $\ell^{\nu-1}$, où ${}_\ell K$ désigne le
noyau de la multiplication par $\ell$ sur $K$. En effet, on montre
facilement que pour tout complexe
borné de $\Lambda$-modules plats le morphisme évident ${}_\ell K\to
\R\SheafHom_\Lambda(\mathbf{Z}/\ell\mathbf{Z},K)$ est un isomorphisme
dans $\detalepluscoefs X {\mathbf{Z}/\ell\mathbf{Z}}$.

\subsubsection{Préservation de \texorpdfstring{$\detalebc X$}{Dbc(X,Lambda)}}

\begin{proposition3}\label{17-proposition-la-dualite-preserve-dbc}
Soit $X$ un $\ZZn$-schéma noethérien excellent muni d'une fonction de dimension
$\delta$. Notons $K_X$ le complexe dualisant potentiel de $(X,\delta)$.
Alors, le foncteur $\Dual_X=\R\SheafHom(-,K_X)$ préserve $\detalebc X$.
\end{proposition3}

Grâce à la proposition~\ref{17-proposition-cdp-changement-coefficients}, on
peut supposer que l'anneau de coefficients est $\mathbf{Z}/\ell\mathbf{Z}$,
avec $\ell$ un nombre premier. Il s'agit de montrer que pour tout faisceau
constructible de $\Lambda$-modules $\mathscr M$ sur $X$,
$\Dual_X \mathscr M\in\detalebc X$. D'après la 
proposition~\ref{17-proposition-constructibilite-cdp}, si $\mathscr M$ est
constant, on a bien $\Dual_X\mathscr M\in\detalebc X$. Si $p\colon Y\to X$ est un
morphisme étale, on a un isomorphisme évident $p^\star
\Dual_X\mathscr M\simeq \Dual_Y p^\star\mathscr M$ pour tout
faisceau constructible $\mathscr M$ sur $X$. Si $\mathscr M$ est localement
constant, en choisissant pour $p$ un morphisme étale surjectif tel que
$p^\star \mathscr M$ soit constant,
on obtient que $\Dual_X\mathscr M\in\detalebc X$ si $\mathscr M$
est localement constant.

On raisonne alors par récurrence noethérienne. Pour tout faisceau
constructible $\mathscr M$, il existe un ouvert dense $U$ de $X$ sur lequel
$\mathscr M$ est localement constant. On note $j\colon U\to X$ l'immersion
ouverte correspondante et $i\colon Z\to X$ une immersion fermée
complémentaire. D'après ce qui précède, on a $j^\star \Dual_X\mathscr
M=\Dual_U{\mathscr M}_{|U}\in\detalebc U$. Par ailleurs, $i^!\Dual_X\mathscr
M\simeq \Dual_Zi^\star \mathscr M$. Par hypothèse de récurrence
noethérienne, on obtient que $i^!\Dual_X\mathscr M$ appartient à
$\detalebc Z$. Le
lemme~\ref{17-lemme-devissage-constructibilite} permet de conclure que
$\Dual_X\mathscr M$ appartient à $\detalebc X$.

\medskip

La stabilité de $\detalebc X$ par $\Dual_X$ permet d'énoncer le résultat
important suivant :

\begin{proposition3}\label{17-proposition-comportement-dualite-morphisme-regulier}
Soit $p\colon X'\to X$ un morphisme régulier
entre $\ZZn$-schémas noethériens excellents. On suppose $X$ muni
d'une fonction de dimension $\delta_X$ et on munit $X'$ de la fonction de
dimension $\delta_{X'}$ définie dans 
la proposition~\ref{17-proposition-stabilite-cdp-morphisme-regulier}. Alors,
pour tout $L\in \detalebc X$, on a un isomorphisme canonique $p^\star\Dual_X
L\iso \Dual_{X'}p^\star L$ et le morphisme de bidualité
$p^\star L\to\Dual^2_{X'} p^\star L$
s'identifie à l'image par $p^\star$ du morphisme de
bidualité $L\to \Dual^2_X L$.
\end{proposition3}

Cela résulte aussitôt des
propositions~\ref{17-proposition-stabilite-cdp-morphisme-regulier},
\ref{17-proposition-rhom-morphismes-reguliers} et
\ref{17-proposition-la-dualite-preserve-dbc}.

\subsubsection{Dimension quasi-injective}

\begin{proposition3}\label{17-proposition-dimension-quasi-injective}
Soit $X$ un $\ZZn$-schéma noethérien excellent muni d'une fonction de
dimension $\delta$. Notons $K_X$ le complexe dualisant potentiel de
$(X,\delta)$. La dimension quasi-injective de $K_X$ est $-2\inf_{x\in
X}\delta(x)$. En particulier, elle est finie si et seulement si $X$ est de
dimension de Krull finie.
\end{proposition3}

Commençons par minorer la dimension quasi-injective de $K_X$. Soit $x\in
X$. On note $i\colon Z\to X$ l'inclusion du sous-schéma intègre de 
$X$ de point générique $x$. On a un isomorphisme canonique 
$i^!K_X\simeq i^\star \R\SheafHom(\Lambda_Z,K_X)$. Le complexe
$K_Z=i^!K_X$ est un complexe dualisant potentiel pour $(Z,\delta_{|Z})$. Par
conséquent $(K_Z)_{\overline x}\simeq \Lambda(\delta(x))[2\delta(x)]$. Il
en résulte que la dimension quasi-injective de $K_X$ est au moins
$-2\delta(x)$. On obtient ainsi la minoration
\[-2\inf_{x\in X}\delta(x)\leq \dimquasiinj K_X\;\text{.}\]

Montrons que cette inégalité est en fait une égalité si $X$ est de
dimension de Krull finie. On peut procéder par récurrence sur la dimension
de $X$. On peut en outre supposer que $X$ est local strictement
hensélien (réduit) de point fermé $x$.
Soit $\mathscr M$ un faisceau constructible de
$\Lambda$-modules sur $X$. Il s'agit de montrer que $\Dual_X\mathscr
M\in\detaletronque X {\leq -2\delta(x)}$.
Il existe un ouvert affine dense $U$ sur lequel
$\mathscr M$ est localement constant. Le schéma $X$ étant réduit et
excellent, quitte à rétrécir $X$, 
on peut supposer que $U$ est régulier. Notons $j\colon U\to X$ l'immersion
de $U$. Notons
$U_1,\dots,U_n$ les composantes connexes de $U$, et $\eta_1,\dots,\eta_n$
les points génériques de ces composantes. Le schéma $U$ étant régulier, on
connaît la structure du complexe dualisant potentiel
$K_U$ : pour tout $1\leq i\leq n$, on a un
isomorphisme canonique $K_{U_i}\simeq
\Lambda(\delta(\eta_i))[2\delta(\eta_i)]$. En particulier,
$K_{U^i}\in\detaletronque {U_i}{\leq -2\delta(\eta_i)}$. Le faisceau ${\mathscr
M}_{|U_i}$ étant localement constant, on obtient que $\Dual_{U_i}{\mathscr
M}_{|U_i}\in\detaletronque {U_i}{\leq -2\delta(\eta_i)}$. D'après le théorème de
Lefschetz affine appliqué aux immersions ouverts affines
$j_i\colon U_i\to X$, il
vient alors que $\R {j_i}_\star \Dual_{U_i}{\mathscr
M}_{|U_i}$ appartient à $\detaletronque X {\leq \dim
\overline{\left\{\eta_i\right\}}-2\delta(\eta_i)}$. Comme on a  
$\dim \overline{\left\{\eta_i\right\}}-2\delta(\eta_i)\leq -2\delta(x)$,
il vient que $\R j_\star \Dual_U{\mathscr M}_{|U}$ appartient à 
$\detaletronque X {\leq -2\delta(x)}$.

Notons $i\colon Z\to X$ une immersion fermée complémentaire à $j$. Grâce à
la récurrence sur la dimension, on sait que 
$\Dual_Zi^\star \mathscr  M\in\detaletronque Z {\leq -2\delta(x)}$. En utilisant le
triangle distingué canonique
\[i_\star \Dual_Zi^\star {\mathscr M}\to \Dual_X\mathscr M\to
\R j_\star \Dual_U{\mathscr M}_{|U}\;\text{,}\]
on obtient bien que $\Dual_X\mathscr M\in\detaletronque X {\leq -2\delta(x)}$, ce
qui achève la démonstration de la proposition.

\subsection{Le théorème en degré négatif ou nul}

\begin{proposition2}\label{17-proposition-dualite-locale-degre-negatif}
Soit $X$ un $\ZZn$-schéma noethérien excellent muni d'une fonction de
dimension $\delta$. Soit $\mathscr M$ un faisceau constructible de
$\Lambda$-modules sur $X$. Alors, le morphisme canonique est un
isomorphisme dans $\detalebc X$.
\[\mathscr M\iso \tau_{\leq 0}\Dual_X\Dual_X \mathscr M\;\text{.}\]
\end{proposition2}

Au cours de cette démonstration, on dira qu'un faisceau constructible de
$\Lambda$-modules $\mathscr M$ sur $X$ est faiblement réflexif si le
morphisme canonique $\mathscr M\to \tau_{\leq 0}\Dual_X\Dual_X \mathscr
M$ de la proposition est un isomorphisme.

D'après le théorème~\ref{17-theoreme-existence-unicite-cdp}, on sait que
$\Lambda$ est faiblement réflexif. Si $g\colon Z\to X$ est une immersion
fermée et $\mathscr N$ un faisceau constructible de $\Lambda$-modules sur
$Z$, il est clair que $\mathscr N$ est faiblement réflexif sur $Z$ si et
seulement si $g_\star \mathscr N$ est faiblement réflexif sur $X$. Plus
généralement, si $f\colon Y\to X$ est un morphisme fini et $\mathscr N$ un
faisceau constructible de $\Lambda$-modules sur $Y$, alors $\mathscr N$ est
faiblement réflexif si et seulement si $f_\star \mathscr N$
l'est (voir \sga{5}{i}{1.13}). Notons aussi qu'une utilisation
appropriée du lemme des cinq montre que si on a une suite
exacte courte $0\to {\mathscr M}' \to \mathscr M\to {\mathscr M}''\to 0$ de
faisceaux constructibles de $\Lambda$-modules, et que ${\mathscr M}''$ est
faiblement réflexif, alors $\mathscr M$ est faiblement réflexif si et
seulement si ${\mathscr M}'$ est faiblement réflexif.

Grâce à la stabilité par extension énoncée plus haut et à la
proposition~\ref{17-proposition-cdp-changement-coefficients}, on peut supposer
que $\Lambda=\mathbf{Z}/\ell\mathbf{Z}$ où $\ell$ est un nombre premier.
Des remarques précédentes, il résulte que si $f\colon Y\to X$ est fini et
que $U$ est un ouvert de $Y$, alors $f_\star \Lambda_U$ est faiblement
réflexif. La classe des faisceaux constructibles de $\Lambda$-modules
faiblement réflexifs sur $X$ étant stable par facteurs directs et
extensions, on peut conclure en utilisant le dévissage des faisceaux
constructibles de \sga{4}{ix}{5.8}.

\subsection{L'argument de \texorpdfstring{\sgaquatreetdemi{Th.~finitude}{4.3}}{[SGA 4\textonehalf{} [Th. finitude] 4.3]}}

\begin{definition2}
Soit $X$ un $\ZZn$-schéma 
noethérien excellent muni d'une fonction de dimension
$\delta$. Soit $\mathscr M$ un faisceau constructible de $\Lambda$-modules
sur $X$. On dit que $\mathscr M$ est réflexif si le morphisme de
bidualité $\mathscr M\to \Dual_X\Dual_X\mathscr M$ est un isomorphisme. On
dira que le morphisme de bidualité est un isomorphisme pour $X$ si tout
faisceau constructible de $\Lambda$-modules sur $X$ est réflexif.
\end{definition2}

\begin{proposition2}\label{17-proposition-argument-de-sga-4-et-demi}
Soit $d\geq 0$.
Si le morphisme de bidualité est un isomorphisme pour les $\ZZn$-schémas
excellents noethériens de dimension $\leq d$, alors il l'est aussi pour les
schémas de type fini sur de tels schémas.
\end{proposition2}

\begin{remarque2}
Dans \sgaquatreetdemi{Th.~finitude}{4.3}, un tel
isomorphisme de bidualité est construit pour les schémas de type
fini sur un $\ZZn$-schéma régulier de dimension au plus un (mais
non nécessairement excellent). Aux hypothèses d'excellence et de régularité
du schéma de base près, il s'agit essentiellement du cas $d=1$ de la
proposition. La démonstration qui suit reprend et généralise celle de 
\sgaquatreetdemi{Th.~finitude}{4.3}.
\end{remarque2}

Par récurrence sur $d$, on peut supposer que l'isomorphisme de bidualité
est un isomorphisme pour les schémas de type fini (et leurs hensélisés
stricts) sur des $\ZZn$-schémas noethériens excellents de dimension $<d$.
Pour montrer que l'isomorphisme de bidualité est un isomorphisme pour tout
schéma $Y$ de type fini sur $X$, il suffit évidemment de le faire pour
$Y=(\mathbf{P}^1)^n\times X$ pour tout $n\in\mathbf{N}$.

Montrons par récurrence sur $n$ que pour tout $\ZZn$-schéma noethérien
excellent $X$ de dimension au plus $d$ (muni d'une fonction de dimension
$\delta$), le morphisme de bidualité est un isomorphisme pour
$(\mathbf{P}^1)^n\times X$. L'hypothèse de la proposition règle le cas
$n=0$. Supposons la propriété établie jusqu'au cran $n-1$, avec $n\geq 1$.
Soit $X$ un
schéma noethérien excellent de dimension $d$ muni d'une fonction de
dimension $\delta$. Montrons que le morphisme de bidualité est un
isomorphisme pour $(\mathbf{P}^1)^n\times X$. On peut supposer que $X$ est
local strictement hensélien, de point fermé $x$. Soit $\mathscr M$ un
faisceau constructible sur $(\mathbf{P}^1)^n\times X$. Notons $C$ un cône
du morphisme de bidualité $\mathscr M\to \Dual\Dual \mathscr M$. On sait
déjà que $C\in\detalebc {(\mathbf{P}^1)^n\times X}$
(cf.~proposition~\ref{17-proposition-la-dualite-preserve-dbc}). Nous
allons dans un premier temps montrer que les faisceaux de cohomologie de
$C$ sont en gratte-ciel, c'est-à-dire supportés par des points fermés.
L'hypothèse de récurrence sur $d$ montre que le support de $C$ est contenu
dans le fermé
$(\mathbf{P}^1)^n\times x$. Posons $Y=(\mathbf{P}^1\times X)_{(y)}$ où
$y$ est le point générique de $\mathbf{P}^1\times x\subset
\mathbf{P}^1\times X$. On considère les $n$ projections canoniques
$(\mathbf{P}^1)^n\times X\to \mathbf{P}^1\times X$ et leur changement de
base $(\mathbf{P}^1)^{n-1}\times Y\to Y$ au-dessus de $Y$. Le schéma $Y$
étant de dimension $d$ et la 
proposition~\ref{17-proposition-comportement-dualite-morphisme-regulier} montrant en
particulier que \guil{la dualité commute aux localisations},
l'hypothèse de récurrence pour $n-1$ implique que si
$z\in (\mathbf{P}^1)^n\times X$ est tel que $C_{\overline z}\neq 0$, alors
les images de $z$ par les $n$ projections canoniques $(\mathbf{P}^1)^n\times
X\to \mathbf{P}^1\times X$ sont des points fermés (puisqu'elles sont
au-dessus de $x$ et différentes du point générique $y$ de
$\mathbf{P}^1\times x$).
Par conséquent, un tel point $z$ est un point fermé de
$(\mathbf{P}^1)^n\times x$. Bref, les faisceaux de cohomologie de $C$ sont
supportés par des points fermés. Il en résulte que si on note  $\pi\colon
(\mathbf{P}^1)^n\times X\to X$ le morphisme canonique, alors pour montrer
que $C\simeq 0$, il suffit de montrer que $\R\pi_\star C\simeq 0$.
D'après \sgaquatreetdemi{Th.~finitude}{4.4},
$\R\pi_\star C$ s'identifie
au cône du morphisme de bidualité $\R\pi_\star \mathscr
M\to\Dual_X\Dual_X \R\pi_\star \mathscr M$, qui est un isomorphisme
par hypothèse. Par conséquent, $C\simeq 0$, ce qui achève la démonstration
de la proposition.

\subsection{Fin de la démonstration}

Démontrons le théorème~\ref{17-theoreme-dualite-locale}. Compte tenu des
résultats antérieurs, il ne reste plus qu'à montrer que le morphisme de
bidualité est un isomorphisme pour tout $\ZZn$-schéma noethérien excellent
$X$ muni d'une fonction de dimension $\delta$. Comme il suffit d'obtenir la
conclusion pour les hensélisés stricts de $X$, on peut supposer que $X$ est
de dimension de Krull finie $d$. On va raisonner par récurrence sur $d$.

\begin{definition2}
Soit $\mathscr M$ un faisceau constructible de $\Lambda$-modules sur un
$\ZZn$-schéma noethérien excellent $X$ muni d'une fonction de dimension
$\delta$. Pour tout $q\geq 1$, on dit que $\mathscr M$ vérifie la propriété
$\propD_q$ si le faisceau de cohomologie ${\mathscr H}^q(\Dual_X^2\mathscr
M)$ est nul.
\end{definition2}

D'après la proposition~\ref{17-proposition-dualite-locale-degre-negatif}, 
$\mathscr M$ est réflexif si et seulement s'il vérifie la propriété
$\propD_q$ pour tout $q\geq 1$.

La fonction de dimension $\delta$ sert à formuler la
propriété $\propD_q$. Pourtant, elle n'en dépend pas. En effet, si $\delta$
et $\delta'$ sont deux fonctions de dimension sur $X$ (connexe), il existe
un entier relatif $k$ tel que $\delta'=\delta+k$. Le complexe dualisant
potentiel $K_{X,\delta'}$ s'identifie canoniquement à
$K_{X,\delta}(k)[2k]$ : les foncteurs de bidualité $\Dual_{X,\delta}^2$ et
$\Dual_{X,\delta'}^2$ sont canoniquement isomorphes. Il n'y a donc pas lieu
de mentionner la fonction de dimension dans la notation $\Dual_X^2$, et les
propriétés $\propD_q$ définies relativement à $\delta$ et $\delta'$
sont équivalentes.

En outre, les propriétés $\propD_q$ sont clairement locales pour la
topologie étale. Comme les schémas noethériens excellents admettent
localement pour la topologie étale des fonctions de dimension, on peut leur
donner un sens même en l'absence d'une fonction de dimension globale. Par
recollement, on peut même donner un sens aux faisceaux de cohomologie
${\mathscr H}^q(\Dual_X^2\mathscr M)$.

\medskip

Soit $d\geq 0$. On suppose que l'isomorphisme de bidualité est un
isomorphisme pour tout $\ZZn$-schéma noethérien excellent de dimension au
plus $d-1$ muni d'une fonction de dimension. 

Nous allons montrer par récurrence sur $q\geq 1$ que tout faisceau
constructible de $\Lambda$-modules $\mathscr M$ sur un schéma noethérien
excellent $X$ de dimension $\leq d$ vérifie la propriété $\propD_q$.

Soit $q\geq 1$. On suppose que pour tout $1\leq q'<q$, tout faisceau de
$\Lambda$-modules sur un schéma noethérien excellent $X$ de dimension $\leq
d$ vérifie la propriété $\propD_{q'}$.

\begin{lemme2}\label{17-lemme-propriete-d-extensions-et-sous-objets}
Les entiers $d$ et $q$ ayant été fixés comme ci-dessus, la propriété
$\propD_q$ pour les faisceaux constructibles de $\Lambda$-modules sur les
schémas noethériens excellents de dimension au plus $d$ est stable par
extensions et sous-objets.
\end{lemme2}

En effet, si on une suite
exacte courte $0\to{\mathscr M}'\to\mathscr M\to {\mathscr M}''\to 0$
de faisceaux constructibles de $\Lambda$-modules sur un
tel schéma $X$, l'hypothèse de récurrence si $q\geq 2$ ou la
proposition~\ref{17-proposition-dualite-locale-degre-negatif} si $q=1$
implique que l'on a une suite exacte de faisceaux :
\[0\to \mathscr \HH^q(\Dual_X^2{\mathscr M}')\to 
\mathscr \HH^q(\Dual_X^2{\mathscr M})\to
\mathscr \HH^q(\Dual_X^2{\mathscr M}'')\]
La propriété $\propD_q$ est donc stable par extensions et sous-objets.

\begin{lemme2}\label{17-lemme-propriete-d-et-morphismes-finis}
Soit $p\colon Y\to X$ un morphisme fini entre schémas noethériens
excellents. Soit $\mathscr M$ un faisceau constructible de
$\Lambda$-modules sur $Y$. Alors, $\mathscr M$ vérifie la propriété
$\propD_q$ si et seulement si $p_\star\mathscr M$ la vérifie.
\end{lemme2}

Ceci résulte aussitôt de l'isomorphisme canonique
$p_\star {\mathscr H}^q(\Dual_Y^2\mathscr M)\simeq {\mathscr
H}^q(\Dual_X^2p_\star \mathscr M)$ (voir \sga{5}{i}{1.12~(a)})
et de la conservativité du foncteur $p_\star$.

\begin{lemme2}\label{17-lemme-devissage-standard-normal}
Soit $X$ un schéma noethérien excellent. 
Soit $\mathscr C$ une sous-catégorie strictement pleine de la catégorie
$\Cons(X,\Lambda)$ des faisceaux constructibles de $\Lambda$-modules sur
$X$, stable par facteurs directs et extensions. On suppose que pour tout
morphisme fini $p\colon Y\to X$, toute immersion ouverte $j\colon U\to Y$
avec $Y$ normal intègre
et tout nombre premier $\ell$ divisant $n$, on a $p_\star j_!
\mathbf{Z}/\ell\mathbf{Z}\in\mathscr C$. Alors, $\mathscr
C=\Cons(X,\Lambda)$.
\end{lemme2}

Il s'agit d'une variante facile de \sga{4}{ix}{5.8}.

\begin{lemme2}
Les entiers $d$ et $q$ ayant été fixés comme ci-dessus, 
si la propriété $\propD_q$ est satisfaite par le faisceau constant $\Lambda$
sur les schémas noethériens excellents normaux (strictement henséliens)
de dimension au plus $d$, alors la propriété $\propD_q$ est satisfaite par
tout faisceau constructible de $\Lambda$-modules sur un schéma noethérien
excellent de dimension au plus $d$.
\end{lemme2}

D'après le lemme~\ref{17-lemme-devissage-standard-normal}, il suffit d'établir
la propriété $\propD_q$ pour un faisceau de la forme $p_\star
j_!\mathbf{Z}/\ell\mathbf{Z}$ avec $\ell$ un diviseur premier de $n$, $p$
un morphisme fini et $j$ une immersion ouverte entre schémas normaux
intègres. D'après le 
lemme~\ref{17-lemme-propriete-d-et-morphismes-finis}, il suffit d'établir la
propriété $\propD_q$ pour $j_!\mathbf{Z}/\ell\mathbf{Z}$ avec $j\colon U\to
Y$ une immersion ouverte, avec $Y$ normal intègre. D'après la 
stabilité par sous-objet et extensions de la propriété $\propD_q$
(cf. lemme~\ref{17-lemme-propriete-d-extensions-et-sous-objets}), il suffit de
traiter le cas du faisceau $j_!\Lambda$, qui est lui-même un sous-faisceau
du faisceau constant $\Lambda$ sur le schéma normal $Y$, ce qui achève la
démonstration du lemme.

\bigskip

On est ainsi ramené à montrer la propriété $\propD_q$ pour le faisceau
constant $\Lambda$ sur les schémas noethériens excellents normaux $X$ de
dimension $d$. On peut supposer
$X$ local strictement hensélien de point fermé $x$ et de point générique
$\eta$. Si $d\leq 1$, $X$ est régulier, et alors, si on choisit la
fonction de dimension $\delta$ sur $X$ telle que $\delta(\eta)=0$,
le complexe dualisant potentiel associé $K_X$ sur $X$ est le
faisceau constant $\Lambda$, et alors il est évident que $\Lambda$ vérifie
la propriété $\propD_q$ puisque l'on a tautologiquement
$\Dual_X\Dual_X\Lambda\simeq \Lambda$. On
peut donc supposer que $d\geq 2$. En appliquant le lemme suivant à la
complétion $\widehat{X}\to X$, on voit qu'on peut supposer que $X$ est
complet :

\begin{lemme2}\label{17-lemme-propriete-d-et-morphismes-reguliers}
Soit $q\geq 1$.
Soit $p\colon X'\to X$ un morphisme régulier
entre $\ZZn$-schémas noethériens excellents. Soit $\mathscr M$
un faisceau constructible de $\Lambda$-modules sur $X$.
On suppose que $p$ est surjectif.
Alors, le faisceau $\mathscr M$ vérifie la propriété $\propD_q$ si et
seulement si $p^\star\mathscr M$ la vérifie.
\end{lemme2}

Ceci résulte aussitôt de la 
proposition~\ref{17-proposition-comportement-dualite-morphisme-regulier}.

\bigskip

Il nous reste à montrer que si $X$ est un schéma local strictement
hensélien noethérien normal complet de dimension $d\geq 2$, alors
le faisceau constant $\Lambda$ sur $X$ vérifie la propriété $\propD_q$.
D'après le théorème d'algébrisation
partielle \refext{V}{5-algébrisation partielle mixte raffinée géométrique},
il existe un morphisme $p\colon X'\to X$ fini surjectif tel que
\begin{itemize}
\item le schéma $X'$ soit normal ;
\item il existe un schéma local noethérien complet $Y$ de
dimension $<d$, un morphisme de type fini $Z\to Y$, un point géométrique
$\overline z\to Z$ et un isomorphisme $X'\simeq \widehat{Z_{(\overline z)}}$.
\end{itemize}

Le schéma $Y$ est noethérien excellent et de dimension $<d$. L'hypothèse
de récurrence sur $d$ et la
proposition~\ref{17-proposition-argument-de-sga-4-et-demi} impliquent que le
morphisme de bidualité est un isomorphisme pour $Z_{(\overline z)}$.  En
particulier, le faisceau constant $\Lambda$ sur $Z_{(\overline z)}$ est
réflexif. Appliqué au morphisme de complétion $X'\to Z_{(\overline z)}$,
la proposition~\ref{17-proposition-comportement-dualite-morphisme-regulier}
montre que le faisceau constant $\Lambda$ sur $X'$ est réflexif.
En particulier, le faisceau constant $\Lambda$ sur $X'$ satisfait
la propriété $\propD_q$.
D'après le lemme~\ref{17-lemme-propriete-d-et-morphismes-finis}, on peut en
déduire que le faisceau
$p_\star\Lambda$ sur $X$ vérifie la propriété $\propD_q$. Le morphisme fini
$p$ étant surjectif, le morphisme canonique $\Lambda\to p_\star\Lambda$ est
un monomorphisme de faisceaux. La propriété $\propD_q$ étant stable par
sous-objets (cf. lemme~\ref{17-lemme-propriete-d-extensions-et-sous-objets}),
le faisceau constant $\Lambda$ sur $X$ vérifie bien la propriété
$\propD_q$, ce qui achève la démonstration du théorème de dualité locale.

\section{Anneaux de coefficients généraux}
\label{17-section-coefficients-generaux}
\subsection{Énoncés}

\begin{definition2}
Soit $X$ un schéma noethérien. Soit $A$ un anneau noethérien. On appelle
complexe dualisant sur $\detalebccoefs X A$ (resp. $\detalebctfcoefs X A$)
un objet $K\in\detalebccoefs X A$ tel
que le foncteur $\Dual_K=\R\SheafHom(K,-)$ préserve $\detalebccoefs X A$
(resp. $\detalebctfcoefs X A$)
et que pour tout $M\in\detalebccoefs X A$ (resp. $M\in\detalebctfcoefs X A$),
le morphisme de bidualité $M\to \Dual^2_K M$ soit un isomorphisme.
\end{definition2}

Cette section vise à établir les deux théorèmes suivants :

\begin{theoreme2}\label{17-theoreme-dualite-d-b-c-coefficients}
Soit $A$ une $\Lambda$-algèbre noethérienne. Soit $X$ un schéma
noethérien. S'il en existe, les complexes dualisants sur $\detalebccoefs X A$
sont uniques au produit tensoriel près avec des objets inversibles.
Soit $R\in\dtotale(A)$ un complexe ponctuellement
dualisant fort au sens de \cite[page~120]{GDBC@Conrad}\;\footnote{On rappelle
que cela signifie ici que $R$ appartient à
$\dbc(A)$ et que pour tout $x\in \SP(A)$,
$R_{(x)}\in\dtotale(A_{(x)})$ est un complexe dualisant pour $A_{(x)}$ au sens
de \cite[page~258]{RD@Hartshorne}, ce qui signifie que $R_{(x)}$ est de
dimension injective finie et que le foncteur
$\R\Hom(-,R_{(x)})$ induit une involution de
$\dbc(A_{(x)})$. D'après,
\cite[lemma~3.1.5]{GDBC@Conrad}, il revient au même de demander que le foncteur
$\R\Hom_A(-,R)$ induise une involution de
$\dbc(A)$.}.
Soit $K$ un complexe dualisant sur $\detalebccoefs X \Lambda$. Alors,
$R\overset{\L}\otimes_{\Lambda} K$ est un complexe dualisant sur
$\detalebccoefs
X A$.
\end{theoreme2}

\begin{theoreme2}\label{17-theoreme-dualite-d-b-ctf-coefficients}
Soit $A$ une $\Lambda$-algèbre noethérienne. Soit $X$ un schéma
noethérien. S'il en existe, les
complexes dualisants sur $\detalebctfcoefs X A$ sont uniques au produit
tensoriel près avec des objets inversibles.
Soit $K$ un complexe dualisant sur $\detalebccoefs X \Lambda$.
Alors, $A\overset{\L}\otimes_{\Lambda} K$ est un complexe dualisant sur
$\detalebctfcoefs X A$.
\end{theoreme2}

\subsection{Systèmes locaux}
\label{17-subsection-systemes-locaux}

\begin{definition2}
Soit $X$ un schéma noethérien. On appelle système local (d'ensembles)
sur $X$
un faisceau d'ensembles sur $X_\et$ isomorphe à une limite inductive
filtrante de faisceaux représentés par des revêtements étales finis de $X$.
Un système local fini est un système local représenté par un revêtement
étale fini.
\end{definition2}

\begin{proposition2}\label{17-proposition-systemes-locaux-ind-objets-dans-rev}
Soit $X$ un schéma noethérien. La catégorie des systèmes locaux sur $X$ est
équivalente à la catégorie $\Ind(\Rev(X))$ des ind-objets dans la catégorie
$\Rev(X)$ des revêtements étales finis de $X$.
\end{proposition2}

Le foncteur qui à un revêtement étale fini $Y\to X$ associe le faisceau
d'ensembles sur $X_\et$ représenté par $Y$ est pleinement fidèle. En
utilisant \sga{1}{v}{8.7.5 a)}, on en déduit, par passage à la limite
inductive, un foncteur pleinement fidèle de la catégorie $\Ind(\Rev(X))$
vers celle des faisceaux d'ensembles sur $X_\et$. Par définition, l'image
essentielle de ce foncteur est la catégorie des systèmes locaux.

\begin{proposition2}
Soit $X$ un schéma noethérien connexe. Soit $\overline x$ un point
géométrique de $X$. Le foncteur qui à un système local $\mathscr F$ associe
la fibre ${\mathscr F}_{\overline x}$ est naturellement muni d'une action
de $\pi_1(X,\overline x)$ et définit une équivalence entre la
catégorie des systèmes locaux sur $X$ et
la catégorie $\pi_1(X,\overline x)-\Ens$ des ensembles sur lesquels le
groupe profini $\pi_1(X,\overline x)$ agit continûment.
Autrement dit, la catégorie des
systèmes locaux d'ensembles sur $X$
s'identifie à la catégorie des faisceaux
d'ensembles sur le topos classifiant du
groupe profini $\pi_1(X,\overline x)$.
\end{proposition2}

Notons $\Ensf$ la catégorie des ensembles finis.
D'après \sga{1}{v}{7}, le foncteur $\Rev(X)\to \Ensf$ qui à $Y$ associe
l'ensemble sous-jacent au schéma $Y_{\overline x}$ s'enrichit d'une action
du groupe $\pi_1(X,\overline x)$ pour définir une équivalence de catégories 
$\Rev(X)\iso \pi_1(X,\overline x)-\Ensf$ où $\pi_1(X,\overline x)-\Ensf$
est la catégorie des ensembles finis (discrets) munis d'une action continue
du groupe profini $\pi_1(X,\overline x)$. En passant cette équivalence aux
ind-objets, on obtient une équivalence $\Ind(\Rev(X))\iso
\pi_1(X,\overline x)-\Ens$, ce qui permet de conclure d'après la
proposition~\ref{17-proposition-systemes-locaux-ind-objets-dans-rev}.

\bigskip

À partir de la définition des systèmes locaux d'ensembles, on peut définir
les systèmes locaux de groupes abéliens, de groupes abéliens de torsion, de
modules, etc.

\begin{proposition2}\label{17-proposition-stabilites-systemes-locaux}
Pour tout schéma noethérien $X$, la catégorie des systèmes locaux
d'ensembles (resp. de groupes abéliens) sur $X$
admet des limites inductives et des limites projectives finies et le
foncteur d'inclusion de la catégorie des systèmes locaux d'ensembles (resp.
de groupes abéliens) sur $X$ dans la
catégorie des faisceaux d'ensembles (resp. de groupes abéliens)
sur $X_\et$ y commute.

Soit $p\colon Y\to X$ un morphisme entre schémas noethériens.
Si $\mathscr F$ est un système local sur $X$, alors $p^\star \mathscr F$
est un système local sur $Y$ ; la réciproque est vraie si $p$ est un
revêtement étale fini surjectif.

Si $\mathscr G$ est un système local sur $Y$ et $p$ un revêtement étale
fini, alors $p_\star \mathscr G$ est un système local sur $X$.
\end{proposition2}

Le seul fait non trivial à montrer ici est que si $p$ est un revêtement
étale fini surjectif et $\mathscr F$ un faisceau d'ensembles sur $X_\et$
tel que $p^\star \mathscr F$ soit un système local, alors $\mathscr F$ est
un système local. On a un isomorphisme évident entre $\mathscr F$ et
l'égalisateur des deux morphismes évidents $p_\star p^\star \mathscr F\to
p_\star p^\star p_\star p^\star \mathscr F$ déduits du couple de foncteurs
adjoints $(p^\star,p_\star)$. D'après les autres propriétés triviales de
stabilité des systèmes locaux, $p_\star p^\star \mathscr F$ et $p_\star
p^\star p_\star p^\star \mathscr F$ sont des systèmes locaux ;
l'égalisateur de deux morphismes entre systèmes locaux est encore un
système local.

\begin{proposition2}\label{17-proposition-devissage-systemes-locaux}
Soit $X$ un schéma noethérien. La catégorie abélienne des systèmes locaux
de groupes abéliens de torsion sur $X$ est stable par extension dans la
catégorie des faisceaux de groupes abéliens sur $X_\et$. Plus précisément,
si $0\to \mathscr F'\to\mathscr F\to\mathscr F''\to 0$ est une suite exacte
de faisceaux de groupes abéliens sur $X$ telle que $\mathscr F'$ et
$\mathscr F''$ soient des systèmes locaux et que $\mathscr F'$
soit de torsion, alors $\mathscr F$ est un système local de groupes
abéliens.
\end{proposition2}

Cette proposition résulte aussitôt du lemme suivant :

\begin{lemme2}
Soit $X$ un schéma noethérien. Soit $\mathscr G$ un faisceau de groupes
abéliens de torsion agissant librement sur un faisceau d'ensembles
$\mathscr T$. On note $p\colon \mathscr T\to\mathscr T/\mathscr G$ le
morphisme quotient. Si $\mathscr G$ et $\mathscr T/\mathscr G$ sont des
systèmes locaux, alors $\mathscr T$ aussi.
\end{lemme2}

On peut supposer que $X$ est connexe.
Notons $\mathscr Y=\mathscr T/\mathscr G$. On peut écrire $\mathscr Y$
comme réunion de sous-systèmes locaux finis $\mathscr Y'$. Pour chacun de
ces $\mathscr Y'$, on peut considérer $p^{-1}(\mathscr Y')$ : il s'agit
d'un faisceau d'ensembles sur $X_\et$ sur lequel $\mathscr G$ agit
librement avec $\mathscr Y'$ pour quotient. Le faisceau $\mathscr T$ étant
réunion des sous-faisceaux $p^{-1}(\mathscr Y')$ associés, pour montrer que
$\mathscr T$ est un système local, il suffit de montrer que pour tout
sous-système local fini $\mathscr Y'$ de $\mathscr Y$, $p^{-1}(\mathscr
Y')$ est un système local. Bref, on peut supposer que $\mathscr Y$ est un
système local fini.

On suppose ainsi que $\mathscr Y$ est représenté par un revêtement étale
fini $q\colon Y\to X$. On peut identifier le faisceau $\mathscr T$ sur $X$
muni du morphisme $p\colon \mathscr T\to \mathscr Y$ à un faisceau étale
$\tilde{\mathscr T}$ sur $Y$ \emph{via}
un isomorphisme canonique $\mathscr T\simeq q_!\tilde{\mathscr T}$.
L'action libre de $\mathscr G$ sur $\mathscr T$ correspond à une action
libre du système local $q^\star \mathscr G$ sur $\tilde{\mathscr T}$.
D'après la proposition~\ref{17-proposition-stabilites-systemes-locaux}, pour
montrer que $\mathscr T$ est un système local, il suffit de montrer que
$\tilde{\mathscr T}$ en est un. Bref, on peut supposer que $Y=X$.

On s'est ramené à la situation où $\mathscr T/\mathscr G$ est l'objet final
de la catégorie des faisceaux sur $X$, c'est-à-dire que $\mathscr T$ est un
torseur sous $\mathscr G$. Si $\mathscr G$ est un système local fini, alors
$\mathscr T$ est représentable par un revêtement étale fini et est donc un
système local ; nous allons nous ramener à ce cas-là.

La classe d'isomorphisme du $\mathscr G$-torseur $\mathscr T$ est définie
par un élément dans l'ensemble $\HH^1(X_\et,\mathscr G)$. Comme
$\HH^1(X_\et, -)$ commute aux limites inductives filtrantes, il existe un
sous-système local de groupes abéliens finis
$\mathscr G'$ de $\mathscr G$ (supposé de torsion),
un $\mathscr G'$-torseur $\mathscr T'$ et un $\mathscr G$-isomorphisme
$\mathscr T\simeq \mathscr G\otimes_{\mathscr G'}\mathscr T'$ où l'on a
noté $\otimes$ le foncteur d'extension du groupe structural (cf.
\cite[proposition~1.3.6, Chapitre~III]{CNA@Giraud}). L'extension du groupe
structural commutant aux limites inductives filtrantes, $\mathscr T$
s'identifie à la limite inductive des $\mathscr G''\otimes_{\mathscr
G'}\mathscr T'$ pour $\mathscr G''$ parcourant l'ensemble ordonné des
sous-systèmes locaux de groupes abéliens finis de $\mathscr G$ contenant
$\mathscr G'$. D'après ce qui précède, $\mathscr G''\otimes_{\mathscr
G'}\mathscr T'$ est un système local d'ensembles sur $X$.  Par passage à la
limite inductive, $\mathscr T$ est bien un système local.

\bigskip

Le résultat de l'exercice suivant montre que l'hypothèse \guil{de torsion}
est bien nécessaire dans la
proposition~\ref{17-proposition-devissage-systemes-locaux}, et que par
ailleurs, un faisceau qui est un système local localement pour la topologie
étale n'est pas forcément un système local.

\begin{exercice2}
Soit $A$ le sous-anneau de $\mathbf{C}[t]$ formé des
polynômes $f$ tels que $f(0)=f(1)$ : le schéma $C=\SP(A)$ correspondant
est obtenu en identifiant $0$ et $1$ dans la droite affine complexe
$\mathbf{A}^1_{\mathbf{C}}$.

\begin{itemize}
\item Montrer que $C$ est isomorphe à la cubique plane d'équation $x^3-y^2+xy=0$
dans le plan affine complexe $\SP(\mathbf{C}[x,y])$
(envoyer $x$ et $y$ respectivement sur $t(t-1)$ et $t^2(t-1))$.

\item Montrer que $C$ admet un unique point singulier $O$.

\item Montrer que le morphisme évident $p\colon
\mathbf{A}^1_{\mathbf{C}}\to C$ est le normalisé de $C$ et que le
sous-schéma fermé (réduit) $p^{-1}(O)$ de $\mathbf{A}^1_{\mathbf{C}}$
est $\{0,1\}$.

\item Construire un isomorphisme $\HH^1_\et(C,\mathbf{Z})\simeq \mathbf{Z}$.

\item Montrer qu'il existe un faisceau de groupes abéliens sur
$C_\et$ tel que :
\begin{enumerate}
\item[(i)] $\mathscr F$ soit extension de deux systèmes locaux, et s'insère
plus précisément dans une suite exacte courte 
$0\to\mathbf{Z}\to\mathscr F\to\mathbf{Z}\to 0$ ;
\item[(i')] Localement pour la topologie étale, $\mathscr F$ soit un système
local ;
\item[(ii)] $\mathscr F$ ne soit pas un système local.
\end{enumerate}
\end{itemize}
\end{exercice2}

\subsection{Partitions galoisiennes}

\begin{definition2}
Une partition galoisienne d'un schéma noethérien $X$ consiste en la donnée
d'une partition finie de $X$ par des sous-schémas (localement fermés)
réduits connexes $(S_i)_{i\in I}$ et d'un revêtement étale galoisien
$S'_i\to S_i$ pour tout $i\in I$.
\end{definition2}

\begin{definition2}
Soit $p\colon Y\to X$ un revêtement fini étale galoisien entre schémas
noethériens. Soit $A$ une
$\Lambda$-algèbre.
On dit qu'un système local de $A$-modules $\mathscr F$ sur $X$ est 
rendu ind-unipotent par $Y$, si pour un point géométrique $\overline y$ de
$Y$ (et donc pour tous) au-dessus d'un point géométrique $\overline x$ de
$X$, le $\pi_1(X,\overline x)$-module discret ${\mathscr F}_{\overline x}$
est ind-unipotent pour le sous-groupe distingué $\pi_1(Y,\overline y)$
(cf.~sous-section~\ref{17-subsection-modules-ind-unipotents-sous-groupe-distingue}),
autrement dit que $\mathscr F$ est limite inductive filtrante de faisceaux
localement constants extensions successives de faisceaux dont l'image
inverse par $p$ soit un faisceau constant.
\end{definition2}

\begin{definition2}
Soit $X$ un schéma noethérien muni d'une partition galoisienne $\mathscr
P=(S'_i\to S_i)_{i\in I}$. Soit $A$ une $\Lambda$-algèbre.
On dit d'un faisceau de
$A$-modules sur $X$ qu'il est faiblement constructible par rapport à
$\mathscr P$ si pour tout $i\in I$, sa restriction à $S_i$ est un système
local rendu ind-unipotent par $S'_i$. On note $\FCons^{\mathscr P}(X,A)$ la
sous-catégorie pleine de la catégorie des faisceaux de $A$-modules sur $X$
formée des faisceaux faiblement constructibles pour $\mathscr P$.
Si $\mathscr P'$ est une deuxième partition galoisienne, on dit que
$\mathscr P'$ raffine $\mathscr P$ si on a l'inclusion $\FCons^{\mathscr
P}(X,\Lambda)\subset \FCons^{\mathscr P'}(X,\Lambda)$ (et donc aussi
$\FCons^{\mathscr P}(X,A)\subset \FCons^{\mathscr
P'}(X,A)$ pour toute $\Lambda$-algèbre $A$).
\end{definition2}

\begin{proposition2}
Soit $X$ un schéma noethérien muni d'une partition galoisienne. Soit $A$
une $\Lambda$-algèbre.
La catégorie $\FCons^{\mathscr P}(X,A)$ est abélienne ; son foncteur
d'inclusion dans la catégorie des faisceaux de $A$-modules sur $X$ est
exact et commute aux limites inductives. $\FCons^{\mathscr P}(X,A)$ est 
stable par extensions dans la catégorie des faisceaux de $A$-modules
sur $X$.
\end{proposition2}

Ceci résulte aussitôt des propriétés des modules ind-unipotents pour un
sous-groupe
(cf.~proposition~\ref{17-proposition-proprietes-unipotents-sous-groupe-distingue})
et des propriétés générales des systèmes locaux
(cf.~sous-section~\ref{17-subsection-systemes-locaux}).

\begin{definition2}
On dit d'une partition galoisienne $\mathscr P=(S'_i\to S_i)_{i\in I}$ sur
un schéma noethérien $X$
qu'elle est dirigée si on a muni $I$ d'un ordre total tel que, soit $I$ est
vide, soit, si on note $i_0$ le plus petit élément de $I$, $S_{i_0}$ est
ouvert et, récursivement, $(S'_i\to S_i)_{i\in I-\{i_0\}}$ est une
partition galoisienne dirigée du fermé réduit $X-S_{i_0}$. On dit qu'une
partition galoisienne est dirigeable s'il existe un ordre total sur
l'ensemble d'indices qui en fasse une partition galoisienne dirigée.
\end{definition2}

\begin{proposition2}\label{17-proposition-raffinements-dirigeables}
Toute partition galoisienne d'un schéma noethérien est raffinée par une
partition galoisienne dirigeable.
\end{proposition2}

\begin{lemme2}
Soit $X'\to X$ un revêtement étale galoisien. Soit $(S_i\to X)_{i\in I}$
une partition de $X$ par un nombre fini de sous-schémas réduits connexes.
On note $S'_i$ une composante connexe du produit fibré $S_i\times_X X'$.
Alors, $(S'_i\to S_i)_{i\in I}$ est une partition galoisienne de $X$ qui
raffine la partition galoisienne $(X'\to X)$.
\end{lemme2}

Cela résulte aussitôt de la théorie de Galois.

\bigskip

Montrons la proposition par récurrence noethérienne sur $X$. Soit $\mathscr
P=(S'_i\to
S_i)_{i\in I}$ une partition galoisienne d'un schéma noethérien non vide
$X$. Tout d'abord, montrons que, quitte à 
raffiner $\mathscr P$, on peut supposer qu'il
existe un indice $i_0\in I$ tel que $S_{i_0}$ soit un ouvert. En effet, si
on choisit un $i_0\in I$ tel que $S_{i_0}$ contienne un point
maximal de $X$, $S_{i_0}$ contient un ouvert non vide $U$ de
$S_{i_0}$. Notons $V_1,\dots,V_n$ les composantes connexes du
fermé réduit $S_{i_0}-U$ de $S_{i_0}$.
D'après le lemme, il existe une partition
galoisienne $\mathscr Q=
(U'\to U,V'_1\to V_1,\dots,V'_n\to V_n)$ de $S_{i_0}$ qui
raffine la partition galoisienne $(S'_{i_0}\to S_{i_0})$ de $S_{i_0}$.
Quitte à remplacer la partition galoisienne initiale de $X$ par son
raffinement $\mathscr Q\cup \mathscr P'$ avec
$\mathscr P'=(S'_i\to S_i)_{i\in I-\{i_0\}}$, on peut effectivement
supposer que $S_{i_0}$ est ouvert.

On peut appliquer l'hypothèse de
récurrence noethérienne à la partition galoisienne $\mathscr P'$ 
de $X-S_{i_0}$ pour en obtenir un raffinement $\mathscr P''$ indexé par un
certain ensemble totalement ordonné $J$ qui fasse de $\mathscr P''$ une
partition galoisienne dirigée de $X-S_{i_0}$. La partition galoisienne
$(S'_{i_0}\to S_{i_0})\cup \mathscr P''$ raffine $\mathscr P$, et si on
étend l'ordre sur $J$ en un ordre sur la réunion disjointe $\{i_0\}\sqcup
J$ de façon à faire de $i_0$ le plus petit élément, on a obtenu une
partition dirigée.

\begin{exercice2}\label{17-exercice-partitions-ensemble-filtrant}
Montrer que si $\mathscr P$ et $\mathscr P'$ sont deux partitions
galoisiennes d'un schéma noethérien $X$, il existe une partition
galoisienne (dirigeable) raffinant à la fois $\mathscr P$ et $\mathscr P'$.
\end{exercice2}

\subsection{Dévissages}

Le but de cette sous-section est d'établir le résultat suivant :

\begin{proposition2}\label{17-proposition-devissage-coefficients}
Soit $A$ une $\Lambda$-algèbre noethérienne.
Soit $X$ un schéma noethérien.
Soit $\mathscr T$ une sous-catégorie triangulée strictement pleine de
$\detalebccoefs X A$ stable par facteurs directs. On suppose que pour tout
nombre premier $\ell$ divisant $n$ et tout $A/\ell A$-module $N$, le
foncteur $N\overset{\L}{\otimes}_{\mathbf{Z}/\ell\mathbf{Z}}-\colon
\detalebccoefs X
{\mathbf{Z}/\ell\mathbf{Z}}\to \detalebccoefs X A$ prend ses valeurs dans
$\mathscr T$. Alors, $\mathscr T=\detalebccoefs X A$.
\end{proposition2}

La démonstration de cette proposition est repoussée à la fin de cette
sous-section.

\begin{proposition2}\label{17-proposition-devissage-partitions-galoisiennes-trivial}
Soit $A$ une $\Lambda$-algèbre.
Soit $\mathscr P=(S'_i\to S_i)_{i\in I}$ avec
$I=\{1,\dots,N\}$ une
partition galoisienne dirigée d'un schéma noethérien $X$. 
Notons $k_i\colon S_i\to X$ l'immersion canonique pour tout $i\in I$. Pour
tout $i\in I$, le foncteur ${k_i}_!\colon \FCons^{(S'_i\to S_i)}(S_i,A)\to
\FCons^{\mathscr P}(X,A)$ est pleinement fidèle. Tout objet $\mathscr F$ de
$\FCons^{\mathscr P}(X,A)$ admet une filtration croissante (fonctorielle)
$(\Fil_n \mathscr F)_{n\in\mathbf{Z}}$ telle que le $\Fil_{-1}\mathscr F=0$,
$\Fil_{N}\mathscr F=\mathscr F$ et que pour tout $1\leq i\leq n$, le
quotient $\Fil_i\mathscr F/\Fil_{i-1}\mathscr F$ soit dans
${k_i}_!\FCons^{S'_i\to S_i}(S_i,A)$.
\end{proposition2}

C'est trivial.

\begin{proposition2}\label{17-proposition-objets-simples-faiblement-constructibles}
Soit $A$ une $\Lambda$-algèbre.
Soit $\mathscr P=(S'_i\to S_i)_{i\in I}$ une partition
galoisienne dirigeable d'un schéma noethérien $X$. Notons $k_i\colon S_i\to
X$ l'immersion canonique pour tout $i\in I$.
Si $A$ est noethérien, alors $\FCons^{\mathscr P}(X,A)$ est une catégorie
abélienne localement noethérienne (cf.~\cite[pages~325--326]{CA@Gabriel}).
Si $A$ est artinien, $\FCons^{\mathscr
P}(X,A)$ est localement finie et admet un nombre fini d'objets simples ;
plus précisément, si on note $k_i\colon S_i\to X$ les inclusions
canoniques, $W_i$ un ensemble (fini) représentatif des objets simples de la
catégorie des $\Gal(S'_i/S_i)$-modules (\emph{via} le choix d'un point
géométrique de $S'_i$, on identifie ces objets à des systèmes locaux sur
$S_i$ trivialisés par $S'_i$), alors les objets ${k_i}_!\mathscr F$ pour
$i\in I$ et $\mathscr F\in W_i$ forment un ensemble représentatif des
objets simples de $\FCons^{\mathscr P}(X,A)$.
\end{proposition2}

Supposons $A$ noethérien.
Montrons que $\FCons^{\mathscr P}(X,A)$ est localement noethérienne. On sait
déjà que cette catégorie abélienne admet des limites inductives filtrantes
et que celles-ci sont exactes. Il s'agit de montrer que tout objet de
$\FCons^{\mathscr P}(X,A)$ est limite inductive d'objets noethériens (ou
plus précisément, mais cela revient au même,
\guil{réunion} de ses sous-objets noethériens). Si
l'ensemble d'indice $I$ de $\mathscr P$ est vide, c'est trivial. Sinon, on
peut choisir un ordre total sur $I$ qui fasse de $\mathscr P$ une partition
galoisienne dirigée, et noter $i_0$ le plus petit élément de $I$. Notons
$j\colon S_{i_0}\to X$ l'immersion (ouverte) correspondante et $k\colon
X-S_{i_0}\to X$ l'immersion du fermé réduit complémentaire. Pour tout objet
$\mathscr F\in \FCons^{\mathscr P}(X,A)$, on a une suite exacte courte :
\[0\to j_!j^\star\mathscr F\to \mathscr F\vers\pi k_\star k^\star \mathscr F\to
0\;\text{.}\]
En raisonnant par récurrence sur
le cardinal de $I$, on peut supposer que $k^\star\mathscr F$ est
\guil{réunion} de ses sous-objets noethériens dans $\FCons^{\mathscr
P'}(X-S_{i_0},A)$ avec $\mathscr P'=\mathscr P-(S'_{i_0}\to S_{i_0})$. Bien
entendu, un objet $\mathscr R\in \FCons^{\mathscr P'}(X-S_{i_0},A)$ est
noethérien si et seulement si $k_\star\mathscr R$ l'est dans
$\FCons^{\mathscr P}(X,A)$. On dispose
donc d'un système inductif $({\mathscr H}_a)_{a\in A}$ indexé par un ensemble
ordonné filtrant $A$ de sous-objets noethériens de $k_\star k^\star
\mathscr F$ tel que $\mathscr F$ soit la réunion des sous-objets
$\pi^{-1}(\mathscr H_a)$ de $\mathscr F$ pour $a\in A$.
Si chacun des $\pi^{-1}(\mathscr H_a)$ est réunion de ses sous-objets
noethériens, alors $\mathscr F$ aussi. Ceci permet de supposer que
$k^\star \mathscr F$ est noethérien. Concernant $j^\star\mathscr F$, en
utilisant que 
la catégorie des $A[G]$-modules discrets (avec $G$ groupe profini) est
localement noethérienne, on obtient que l'objet $j^\star\mathscr F$ de
$\FCons^{S'_{i_0}\to S_{i_0}}(S_{i_0},A)$ est réunion de ses sous-objets
noethériens ; on en déduit aussitôt que $j_!j^\star\mathscr F$ est aussi
réunion de ses sous-objets noethériens dans $\FCons^{\mathscr P}(X,A)$. Pour
conclure que $\mathscr F$ est réunion de ses sous-objets noethériens, on
utilise le lemme suivant :

\begin{lemme2}
Soit $A$ une $\Lambda$-algèbre noethérienne.
Soit $\mathscr P$ une partition galoisienne dirigeable
d'un schéma noethérien $X$.
Soit $0\to \mathscr H\to\mathscr F\to \mathscr G\to 0$ une suite exacte
courte dans $\FCons^{\mathscr P}(X,A)$. On suppose que $\mathscr G$ est un
objet noethérien de $\FCons^{\mathscr P}(X,A)$ et que
$\mathscr H$ est réunion de ses sous-objets noethériens.
Alors, $\mathscr F$ est aussi réunion de ses sous-objets noethériens.
\end{lemme2}

Il est évident que $\mathscr G$ est un faisceau de $A$-modules
constructible sur $X$. D'après \sga{4}{ix}{2.7.3}, le foncteur
$\Ext^1(\mathscr G,-)$ de la catégorie des faisceaux de $A$-modules sur $X$
vers celle des $A$-modules commute aux limites inductives filtrantes. La
suite exacte donnée définissant un élément dans $\Ext^1(\mathscr G,\mathscr
H)$ et $\mathscr H$ s'écrivant comme une limite inductive filtrante de ses
sous-objets noethériens, il existe un sous-objet noethérien ${\mathscr
H}'$ de $\mathscr H$, une suite exacte courte $0\to \mathscr H'\to \mathscr
F'\to\mathscr G\to 0$ et un diagramme commutatif de la forme suivante, où
le carré de gauche est cocartésien :
\[
\xymatrix{0\ar[r] & \ar[d]\mathscr H'\ar[r] & \mathscr F'\ar[r]\ar[d]
 & \mathscr G\ar@{=}[d] \ar[r]& 0 \\
0\ar[r] & \mathscr H\ar[r] & \mathscr F\ar[r] &
\mathscr G \ar[r]& 0\;\text{,}}
\]
À vrai dire, on peut remplacer $\mathscr H'$ par tout sous-objet (noethérien)
$\mathscr H''$ de $\mathscr H$ contenant $\mathscr H'$, et $\mathscr F$
s'identifie à la réunion des sous-objets $\mathscr F''$ ainsi définis. Pour
conclure, il suffit de montrer qu'un tel $\mathscr F''$ est un objet 
noethérien de $\FCons^{\mathscr P}(X,A)$, ce qui est évident puisqu'il est
extension de deux objets noethériens $\mathscr H''$ et $\mathscr G$.

\bigskip

Supposons $A$ artinien. Il s'agit de trouver un ensemble fini d'objets
simples à partir desquels tous les objets noethériens s'obtiennent par
extensions successives. Compte tenu du dévissage de la
proposition~\ref{17-proposition-devissage-partitions-galoisiennes-trivial}, on
peut supposer que $\mathscr P$ est constitué d'un unique revêtement
galoisien $X'\to X$. Choisissons un point géométrique $\overline x'$ de
$X'$ au-dessus d'un point géométrique $\overline x$ de $X$. Notons
$G=\pi_1(X,x)$ et $H=\pi_1(X',x)$. Le groupe profini $H$ s'identifie à un
sous-groupe distingué fermé de $G$. Notons $K=G/H$ le groupe fini quotient.
L'anneau $A[K]$ est évidemment artinien à gauche, on peut en noter $W$ un
ensemble représentatif fini d'objets simples ; considérés comme des
$A[G]$-modules discrets, les éléments de $W$ sont encore simples.

La catégorie $\FCons^{\mathscr P}(X,A)$ s'identifie à la catégorie des
$G$-modules discrets ind-unipotents pour $H$. Tout objet de cette catégorie
s'écrit comme une réunion de sous-objets de type fini unipotents pour $H$,
et de tels sous-objets se dévissent eux-mêmes en extensions successives
d'objets sur lesquels $H$ agit trivialement, ces derniers s'identifiant à
des $A[K]$-modules de type fini, ils se dévissent en extensions successives
d'éléments de $W$.

\begin{definition2}
Soit $A$ une $\Lambda$-algèbre noethérienne.
Soit $X$ un schéma noethérien. Soit
$\mathscr P$ une partition galoisienne. On note $\Cons^{\mathscr P}(X,A)$
la sous-catégorie pleine de $\FCons^{\mathscr P}(X,A)$ formée des faisceaux
de $A$-modules constructibles (ce qui revient ici à dire que
les fibres sont des $A$-modules de type fini). Il s'agit bien entendu de la
sous-catégorie abélienne des objets noethériens dans $\FCons^{\mathscr
P}(X,A)$.
\end{definition2}

\begin{proposition2}\label{17-proposition-devissage-gabber}
Soit $A$ une $\Lambda$-algèbre noethérienne. Soit $X$ un schéma noethérien.
Soit $\mathscr P$ une partition galoisienne dirigeable de $X$. Soit $\mathscr
F$ un objet de $\Cons^{\mathscr P}(X,A)$. 
Il existe une filtration finie de $\mathscr F$ dans $\Cons^{\mathscr P}(X,A)$
dont les quotients successifs soient des facteurs directs d'objets de la
forme $M\otimes_{\mathbf{Z}/\ell\mathbf{Z}} {\mathscr F}_0$ où $\ell$ est
un nombre premier divisant $n$, $M$ un $A/\ell A$-module de type fini et
${\mathscr F}_0$ un objet de $\Cons^{\mathscr
P}(X,\mathbf{Z}/\ell\mathbf{Z})$.
\end{proposition2}

D'après la
proposition~\ref{17-proposition-devissage-partitions-galoisiennes-trivial}, on
peut supposer que $\mathscr P$ est constituée d'un unique revêtement
galoisien $X'\to X$. En reprenant les notations utilisées dans la
démonstration du cas artinien de la
proposition~\ref{17-proposition-objets-simples-faiblement-constructibles}, on
peut identifier $\mathscr F$ à un $A[G]$-module discret unipotent pour le
sous-groupe distingué fermé $H$. Cette propriété d'unipotence permet de
supposer que $H$ agit trivialement, de sorte qu'on se retrouve avec une
action du groupe fini $K=G/H=\Gal(X'/X)$ (bref, on peut supposer que 
$\mathscr F$ est un système local de $A$-modules trivialisé par $X'$). Le
lemme suivant appliqué à l'algèbre de groupe $B=A[K]$ permet de conclure.

\begin{lemme2}
Soit $A$ une $\Lambda$-algèbre noethérienne. Soit $B$ une $\Lambda$-algèbre
finie non nécessairement commutative. Tout $A\otimes_\Lambda B$-module (à
gauche) de type fini admet une filtration finie dont les quotients successifs
soient des facteurs directs de $A\otimes_\Lambda B$-modules de la forme
$N\otimes_{\mathbf{F}_\ell} L$ où $\ell$ est un nombre premier divisant
$n$, $N$ un $A/\ell A$-module de type fini et $L$ un $B/\ell B$-module simple.
\end{lemme2}

En premier lieu, comme $A\otimes_\Lambda B$ est évidemment noethérien à
gauche, on peut procéder à une récurrence noethérienne ; il suffit donc de
montrer que tout $A\otimes_\Lambda B$-module non nul admet un sous-module
non nul facteur direct d'un module de la forme $N\otimes_{\mathbf{F}_\ell}
L$ avec $N$ un $A/\ell A$-module de type fini, $L$ un $B/\ell B$-module
simple et $\ell$ un nombre premier divisant $n$. Ceci permet de supposer
que $B$ est un anneau semi-simple. En effet, pour tout $B$-module non nul
$M$, l'annulateur du radical de Jacobson $\mathscr N$ de $B$ dans $M$
est un sous-$B$-module non nul de $M$ (cf.~\cite[§4]{GTM131@Lam}) ;
si $M$ est un $A\otimes_{\Lambda} B$-module non nul, l'annulateur de
$\mathscr N$ dans $M$ s'identifie donc à un $A\otimes_{\Lambda}(B/\mathscr
N)$-module non nul et l'anneau $B/\mathscr N$ est bien semi-simple.

En deuxième lieu, l'énoncé du lemme est vrai pour un produit
$B=B_1\times\dots\times B_k$ d'anneaux si et seulement s'il est vrai pour
chacun des $B_i$ et l'énoncé est aussi invariant par équivalence de Morita
(cf.~\cite[§18]{GTM189@Lam}) puisqu'on peut le formuler intrinsèquement
en termes de la catégories des $B$-modules. Compte tenu du théorème
d'Artin-Wedderburn de structure des anneaux semi-simples
(cf.~\cite[3.5]{GTM131@Lam}), on peut donc supposer que $B$ est un
corps fini, \emph{a priori} non commutatif, mais effectivement commutatif
en vertu du théorème de Wedderburn (cf.~\cite[13.1]{GTM131@Lam}).

En troisième lieu, l'énoncé est vrai dans le cas particulier auquel on
s'est ramené ci-dessus. Soit $B=L$ une extension finie de
$\mathbf{F}_\ell$, pour un certain nombre premier $\ell$ divisant $n$. On
peut supposer que $\ell$ annule $A$.
Soit $M$ un
$A\otimes_{\mathbf{F}_\ell} L$-module. L'extension $L/\mathbf{F}_\ell$ est
galoisienne, notons $G$ son groupe de Galois.
On considère le $A\otimes_{\mathbf{F}_\ell} L$-module
$M'=\oplus_{\sigma\in G} \sigma^\star M$ où $\sigma^\star M$ est le
$A\otimes_{\mathbf{F}_\ell} L$-module obtenu par image inverse par
l'automorphisme de $\SP(A\otimes_{\mathbf{F}_\ell} L)$ induit par
$\sigma^{-1}\colon L\to L$ (on peut identifier $\sigma^\star M$ au
$A$-module sur $M$ sur lequel on ferait agir $\lambda\in L$ par la
multiplication d'origine par $\sigma(\lambda)$).
Le $A\otimes_{\mathbf{F}_\ell} L$-module $M$ est évidemment un facteur
direct de $M'$. On va montrer que $M'$ est de la forme
$N\otimes_{\mathbf{F}_\ell} L$ où $N$ est un $A$-module. On observe que
pour tout $\sigma\in G$, on dispose d'un isomorphisme évident de
$A\otimes_{\mathbf{F}_\ell}L$-modules $M'\simeq \sigma^\star M'$. Ces
isomorphismes vérifient tautologiquement les relations de cocycles qui
permettent de munir le $A\otimes_{\mathbf{F}_\ell} L$-module $M'$ d'une
donnée de descente relative au morphisme de schémas
$\SP(A\otimes_{\mathbf{F}_\ell} L)\to\SP(A)$ (qui est un revêtement étale
surjectif). Par descente fidèlement plate (cf.~\sga{1}{viii}{1.1}),
il existe un $A$-module $N$ tel
que $N\otimes_{\mathbf{F}_\ell} L$ s'identifie à $M'$ de façon compatible
avec les isomorphismes $M'\simeq \sigma^\star M'$. En particulier, on a un
isomorphisme de la forme voulue entre les $A\otimes_{\mathbf{F}_\ell}
L$-modules $M'$ et $N\otimes_{\mathbf{F}_\ell} L$, ce qui permet de
conclure.

\bigskip

Nous sommes maintenant en mesure de démontrer la
proposition~\ref{17-proposition-devissage-coefficients}. Soit $\mathscr T$ une
telle sous-catégorie triangulée de $\detalebccoefs X A$. Soit $\mathscr F$ un
faisceau constructible de $A$-modules sur $X$. Il s'agit de montrer que
$\mathscr F$ appartient à $\mathscr T$. Il existe évidemment une partition
galoisienne $\mathscr P$ telle que $\mathscr F$ appartienne à
$\Cons^{\mathscr P}(X,A)$. D'après la
proposition~\ref{17-proposition-raffinements-dirigeables}, on peut supposer
que $\mathscr P$ est dirigeable. On peut alors appliquer la
proposition~\ref{17-proposition-devissage-gabber} pour conclure que $\mathscr
F$ appartient à $\mathscr T$.

\subsection{Complexes dualisants sur \texorpdfstring{$\detalebccoefs X A$}{Dbc(X,A)}}
\subsubsection{Unicité}

\begin{proposition3}\label{17-proposition-unicite-complexe-dualisant}
Soit $X$ un schéma noethérien. Soit $A$ un anneau noethérien. Si
$K$ et $K'$ sont deux complexes dualisants sur $\detalebccoefs X A$ (resp.
$\detalebctfcoefs X A$), alors il
existe un objet inversible $L\in\detalebctfcoefs X A$ 
(cf.~proposition~\ref{17-proposition-objets-inversibles-d-t-a} pour plus de
précisions) tel que $K'$ soit isomorphe à $L\otimes^{\L}_{A}K$.
\end{proposition3}

\begin{lemme3}\label{17-lemme-dual-dans-d-b-implique-dans-d-b}
Soit $X$ un schéma noethérien. Soit $A$ un anneau noethérien. On suppose
que $K\in\detalebccoefs X A$ est dualisant. Pour tout
$F\in\detalemoinsccoefs X A$, si $\Dual_K F\in\detalebccoefs X A$, alors
$F\in\detalebccoefs X A$.
\end{lemme3}

Commençons par montrer que l'on peut supposer que $\Dual_K F=0$. Pour tout
$F\in\detalecoefs X A$, notons $\varepsilon_F\colon F\to \Dual_K\Dual_K F$ le
morphisme de bidualité. Le morphisme composé
\[\xymatrix{\Dual_K F\ar[r]^-{\varepsilon_{\Dual_K F}} &\Dual_K\Dual_K\Dual_K
F\ar[r]^-{\Dual_K(\varepsilon_F)}& \Dual_K F}\]
est l'identité de $\Dual_K F$. Comme $\Dual_K F$ est dans $\detalebccoefs X A$,
le fait que $K$ soit dualisant implique que $\varepsilon_{\Dual_K F}$ est
un isomorphisme. Par conséquent, $\Dual_K(\varepsilon_F)$ est un
isomorphisme. Quitte à remplacer $F$ par un cône de $\varepsilon_F$, on
peut supposer que $\Dual_K F=0$.

Par l'absurde, supposons que $F$ soit non nul. Il existe un morphisme
non nul $p\colon F\to F'$ avec $F'\in\detalebccoefs X A$ (par exemple, le
morphisme canonique $F\to \tau_{\geq n}F$ pour un entier $n$ bien choisi).
Considérons le carré commutatif :
\[
\xymatrix{F \ar[r]^p \ar[d]^{\varepsilon_F}& F'\ar[d]^\sim_{\varepsilon_{F'}}\\
\Dual^2_K F\ar[r]^{\Dual^2_K(p)} & \Dual^2_KF'}
\]
D'un côté, $\Dual^2_K F$ est nul, donc $\varepsilon_{F'}\circ p=0$, mais de
l'autre, $\varepsilon_{F'}$ est un isomorphisme, d'où $p=0$, ce qui conduit
à une contradiction.

\bigskip

Établissons la proposition~\ref{17-proposition-unicite-complexe-dualisant}
dans le cas des complexes dualisants sur $\detalebccoefs X A$, la démonstration
qui suit vaudra aussi pour $\detalebctfcoefs X A$ (à ceci près qu'il
ne sera plus alors nécessaire de recourir au lemme ci-dessus).
On pose $Y=\Dual_KK'\in\detalebccoefs X A$. Comme
$K$ est dualisant, on a aussi un isomorphisme privilégié $K'=\Dual_K Y$.
Pour tout $Z\in\detalebccoefs X A$, on a un isomorphisme fonctoriel
$\Dual_{K'}Z\simeq \Dual_K(Z\overset{\L}\otimes_A{Y})$ dans $\detalecoefs
X A$. Grâce au lemme, on en déduit que pour tout $Z\in\detalebccoefs X A$, on a
$Z\overset{\L}\otimes_A Y\in\detalebccoefs X A$. On a ainsi un triangle
commutatif de catégories et de foncteurs (à isomorphismes près de
foncteurs) :
\[
\xymatrix{\ar[rd]_{\Dual_{K'}}\detalebccoefs X A
\ar[r]^-{-\overset{\L}\otimes_A Y} & \ar[d]^{\Dual_K} \detalebccoefs X A\\
& (\detalebccoefs X A)^\opp }
\]
Comme $\Dual_K$ et $\Dual_{K'}$ sont des équivalences, le foncteur
$-\overset{\L}\otimes_A Y$ aussi. En particulier, $Y$ est un complexe
inversible. Notons $Y'$ l'inverse de $Y$. On a un isomorphisme de
foncteurs $\R\SheafHom(Y,-)\simeq Y'\overset{\L}\otimes_{A}-$
(si un foncteur est une équivalence, son adjoint à droite est un
quasi-inverse). Comme $K'=\Dual_K Y$, on peut en déduire que $K'\simeq
Y'\overset{\L}\otimes_A K$.

\subsubsection{Réduction au cas \texorpdfstring{$\Lambda=\mathbf{Z}/\ell\mathbf{Z}$}{Lambda=Z/lZ}}

\begin{proposition3}\label{17-proposition-reduction-complexe-dualisant-coh-ferme}
Soit $A$ un anneau noethérien. Soit $R\in \dtotale(A)$ un complexe
ponctuellement dualisant fort. Soit $J$ un idéal de $A$. On pose $A'=A/J$
et $R'=\R\Hom_A(A',R)\in\dtotale(A')$. Alors, $R'$ est un complexe
ponctuellement dualisant fort. Si on note $\Dual$ (resp. $\Dual'$) le
foncteur $\R\Hom_A(-,R)$ (resp. $\R\Hom_{A'}(-,R')$) sur
$\dtotale(A)$ (resp. $\dtotale(A')$) et $\oub\colon
\dtotale(A')\to\dtotale(A)$ le foncteur de \guil{restriction des scalaires},
on a un isomorphisme canonique :
\[\oub\circ \Dual'\simeq \Dual\circ \oub\;\text{.}\]
\end{proposition3}

C'est essentiellement trivial.

\begin{proposition3}\label{17-proposition-complexes-dualisants-quotient-par-ideal}
Soit $X$ un schéma noethérien. Soit $K\in\detalebctfcoefs X \Lambda$.
Soit $A$ une
$\Lambda$-algèbre noethérienne. Soit $J$ un idéal de $A$. Soit
$R\in\dtotale(A)$ un complexe ponctuellement dualisant fort. On pose
$A'=A/J$ et $R'=\R\Hom_A(A',R)\in\dtotale(A')$. On note
$K_R=K\overset{\L}\otimes_\Lambda R\in\detalebccoefs X A$ et $K_{R'}=
K\overset{\L}\otimes_\Lambda {R'}\in\detalebccoefs X {A'}$. On note
$\oub\colon \detalebccoefs X {A'}\to\detalebccoefs X A$ le foncteur conservatif
évident. Alors, pour tout $M\in\detalecoefs X {A'}$, on a un isomorphisme
canonique dans $\detalecoefs X A$ :
\[\oub(\R\SheafHom_{A'}(M,K_{R'}))\simeq \R\SheafHom_A(\oub
M,K_R)\;\text{.}\]
De plus, si $K_R$ est un complexe dualisant sur $\detalebccoefs X A$, alors
$K_{R'}$ en est un sur $\detalebccoefs X {A'}$ et la réciproque est vraie si
$J$ est nilpotent.
\end{proposition3}

Les autres assertions en étant des conséquences faciles, il s'agit de
montrer que l'on a un isomorphisme canonique
$\R\SheafHom_A(A',K_R)\simeq \oub(K_{R'})$ dans $\detalebccoefs X A$, ce
qui résulte de la proposition~\ref{17-proposition-coefficients-universels-2}.

\begin{corollaire3}\label{17-corollaire-reduction-theoreme-dualite-d-b-c-coefs-z-l}
Pour démontrer le théorème~\ref{17-theoreme-dualite-d-b-c-coefficients}, on
peut supposer que $A=\mathbf{Z}/\ell\mathbf{Z}$ où $\ell$ est un nombre
premier.
\end{corollaire3}

Étant entendu que la propriété d'unicité des complexes dualisants a déjà
été obtenu, il est évident que pour démontrer le
théorème~\ref{17-theoreme-dualite-d-b-c-coefficients}, on peut supposer que
$\Lambda=\mathbf{Z}/\ell^\nu\mathbf{Z}$ où $\ell$ est un nombre premier et
$\nu\geq 1$. Posons $A'=A/\ell A$. Notons $R\in\dtotale(A)$ un complexe
ponctuellement dualisant fort. D'après la
proposition~\ref{17-proposition-reduction-complexe-dualisant-coh-ferme}, le
complexe $R'=\R\SheafHom_A(A',R)\in\dtotale(A')$ en est un pour $A'$.
Appliquant dans un premier temps la
proposition~\ref{17-proposition-complexes-dualisants-quotient-par-ideal} au
cas où $\Lambda\to A$ est $\Lambda\to \mathbf{Z}/\ell\mathbf{Z}$, nous
obtenons que $K''=\R\SheafHom_\Lambda(\mathbf{Z}/\ell\mathbf{Z},K)$
est un complexe dualisant sur $\detalebccoefs X {\mathbf{Z}/\ell\mathbf{Z}}$.
D'après le lemme~\ref{17-lemme-tor-finitude-z-l-n}, on a aussi un isomorphisme
$K''\simeq K\overset{\L}\otimes_\Lambda\mathbf{Z}/\ell\mathbf{Z}$
dans $\detalebccoefs X {\mathbf{Z}/\ell\mathbf{Z}}$. Appliquons le
théorème~\ref{17-theoreme-dualite-d-b-c-coefficients} dans le cas de la
$\mathbf{Z}/\ell\mathbf{Z}$-algèbre $A'$ : on obtient que
$K''\overset{\L}\otimes_{\mathbf{Z}/\ell\mathbf{Z}} R'$ est un complexe
dualisant sur $\detalebccoefs X {A'}$. Cet objet
$K''\overset{\L}\otimes_{\mathbf{Z}/\ell\mathbf{Z}} R'$ s'identifie
aussi à $K\overset{\L}\otimes_\Lambda R'=K_{R'}$. D'après la
proposition~\ref{17-proposition-complexes-dualisants-quotient-par-ideal}, il
vient que $K\overset{\L}\otimes_\Lambda R=K_R$ est un complexe
dualisant sur $\detalebccoefs X A$, ce qui achève la démonstration de ce
corollaire.

\subsubsection{Démonstration du théorème~\ref{17-theoreme-dualite-d-b-c-coefficients}}

L'énoncé d'unicité des complexes dualisants sur $\detalebccoefs X A$
a déjà été obtenu,
cf.~proposition~\ref{17-proposition-unicite-complexe-dualisant}.
Pour l'énoncé d'existence, d'après le
corollaire~\ref{17-corollaire-reduction-theoreme-dualite-d-b-c-coefs-z-l}, on
peut supposer que $\Lambda=\mathbf{Z}/\ell\mathbf{Z}$ avec $\ell$ un nombre
premier. On se donne $K$ un complexe dualisant sur $\detalebccoefs X
{\mathbf{Z}/\ell\mathbf{Z}}$. Soit $A$ une $\mathbf{Z}/\ell\mathbf{Z}$-algèbre
noethérienne et $R\in\dtotale(A)$ un complexe ponctuellement dualisant
fort. Notons $\Dual_X=\R\SheafHom_{\Lambda}(-,K)$ le foncteur de
dualité sur $\detalebccoefs X {\Lambda}$ induit par $K$ et $\Dual_A$ celui
induit par $R$ sur $\dbc(A)$. Notons $\Dual_{X,A}$ le
foncteur $\R\SheafHom_A(-,K_R)$ sur $\detalecoefs X A$ où
$K_R=K\overset{\L}\otimes_A R$. La
proposition~\ref{17-proposition-coefficients-universels-3} ($\Lambda$ est un
corps) montre que l'on a
un isomorphisme canonique, pour tout $N\in\dbc(A)$ et
$\mathscr F\in\detalebccoefs X {\mathbf{Z}/\ell\mathbf{Z}}$ :
\[\Dual_{X,A}(N
\overset{\L}\otimes_{\mathbf{Z}/\ell\mathbf{Z}} \mathscr F)\simeq
(\Dual_A
N)\overset{\L}\otimes_{\mathbf{Z}/\ell\mathbf{Z}}(\Dual_X\mathscr F)\;
\text{.}\]
Par hypothèse, $\Dual_A N\in\dbc(A)$ et
$\Dual_X\mathscr F\in\detalebccoefs X {\mathbf{Z}/\ell\mathbf{Z}}$, ce qui
permet de déduire que $\Dual_{X,A}(N
\overset{\L}\otimes_{\mathbf{Z}/\ell\mathbf{Z}} \mathscr F)$
appartient à $\detalebccoefs X A$, puis que $\Dual_{X,A}$ préserve
$\detalebccoefs
X A$ grâce au dévissage de la
proposition~\ref{17-proposition-devissage-coefficients}.
Avec les mêmes notations, les morphismes de bidualité $N\to\Dual_A^2N$
et $\mathscr F\to\Dual_X^2\mathscr F$ sont des isomorphismes, le morphisme
de bidualité $N
\overset{\L}\otimes_{\mathbf{Z}/\ell\mathbf{Z}} \mathscr
F\to\Dual_{X,A}^2(
N
\overset{\L}\otimes_{\mathbf{Z}/\ell\mathbf{Z}} \mathscr F)$ en est
donc un aussi ; le même dévissage permet de conclure que $\Dual_{X,A}$
définit
une involution de $\detalebccoefs X A$, c'est-à-dire que $K_R$ est un complexe
dualisant sur $\detalebccoefs X A$, ce qui achève la démonstration du
théorème~\ref{17-theoreme-dualite-d-b-c-coefficients}.

\subsection{Complexes dualisants sur \texorpdfstring{$\detalebctfcoefs X A$}{Dbctf(X,A)}}

L'assertion d'unicité des complexes dualisants sur $\detalebctfcoefs X A$ a
déjà été énoncée dans la
proposition~\ref{17-proposition-unicite-complexe-dualisant}. L'essentiel de
cette sous-section vise à établir le théorème suivant, dont on va déduire
dans quelques lignes le
théorème~\ref{17-theoreme-dualite-d-b-ctf-coefficients} :

\begin{theoreme2}\label{17-theoreme-stabilite-d-b-ctf-par-rhom}
Soit $X$ un schéma noethérien. On suppose qu'il existe un complexe
dualisant sur $\detalebc X$. Soit $A$ une $\Lambda$-algèbre noethérienne.
Pour tous $K$ et $L$ objets de $\detalebctfcoefs X A$, l'objet
$\R\SheafHom_A(K,L)$ appartient à $\detalebctfcoefs X A$ et pour tout
$M\in\dplus(A)$, le morphisme canonique est un isomorphisme :
\[M\overset{\L}\otimes_A \R\SheafHom_A(K,L)\iso
\R\SheafHom_A(K,M\overset{\L}\otimes_A L)\;\text{.}\]
Si $A'$ est une $A$-algèbre,
que $K$ et $L$ sont des objets de
$\detalebctfcoefs X A$, alors on a un isomorphisme canonique dans
$\detalebcoefs X
{A'}$ :
\[A'\overset{\L}\otimes_A\R\SheafHom_A(K,L)\iso
\R\SheafHom_{A'}(A'\overset{\L}\otimes_A K,
A'\overset{\L}\otimes_A L)\;\text{.}\]
\end{theoreme2}

Montrons que l'on peut déduire le
théorème~\ref{17-theoreme-dualite-d-b-ctf-coefficients} de 
ce théorème~\ref{17-theoreme-stabilite-d-b-ctf-par-rhom} et du
théorème~\ref{17-theoreme-dualite-d-b-c-coefficients} qui a déjà été établi.
Commençons par un lemme :

\begin{lemme2}
Soit $A$ un anneau noethérien. Soit $K\in\dmc(A)$.
Alors $K$ est nul si et seulement si pour 
tout idéal maximal $\mathfrak m$ de $A$, l'objet
$A/\mathfrak m\overset{\L}\otimes_A K$ est nul. 
\end{lemme2}

On suppose que $K$ n'est pas nul. On veut montrer qu'il existe un idéal
maximal $\mathfrak m$ de $A$ tel que $A/\mathfrak
m\overset{\L}\otimes_A K$ ne soit pas nul.
Soit $q$ le plus grand entier tel que
$\HH^q(K)$ soit non nul. On peut supposer que $K$ est un complexe formé de
$A$-modules projectifs de type fini et nuls en degrés strictement plus
grands que $q$. Par construction du produit tensoriel dérivé, on a une
surjection $\HH^q(A/\mathfrak m\overset{\L}\otimes_A K)\to
\HH^q(K)/\mathfrak m \HH^q(K)$ pour tout idéal maximal $\mathfrak m$ de $A$. 
Pour conclure, il suffit donc de montrer qu'il existe un idéal maximal
$\mathfrak m$ tel que $\HH^q(K)\neq \mathfrak m \HH^q(K)$ ou encore, d'après le
lemme de Nakayama, que $\HH^q(K)\otimes_A A_{\mathfrak m}\neq 0$. Le support
de $\HH^q(K)$ est un fermé non vide de $\SP(A)$ (cf.
\egazero{i}{1.7}), il contient un point fermé que l'on
identifie à un idéal maximal $\mathfrak m$ de $A$,
et cet idéal maximal vérifie la condition voulue.

\bigskip

Soit $K$ un
complexe dualisant sur $\detalebc X$ et $A$ une $\Lambda$-algèbre noethérienne.
Grâce au lemme~\ref{17-lemme-tor-finitude-z-l-n}, il vient que $K$ appartient
à $\detalebctf X$ et donc que
$K_A=A\overset{\L}\otimes_\Lambda K$ appartient
à $\detalebctfcoefs X A$. D'après le
théorème~\ref{17-theoreme-stabilite-d-b-ctf-par-rhom}, le foncteur
$\Dual_A=\R\SheafHom_A(-,K_A)$ préserve $\detalebctfcoefs X A$. Il reste à
montrer que le morphisme de bidualité $L\to \Dual_A^2 L$ est un
isomorphisme pour tout $L\in\detalebctfcoefs X A$.
D'après le lemme, il suffit de
montrer qu'après produit tensoriel dérivé avec $A/\mathfrak m$, le
morphisme $L\to\Dual_A^2 L$ induit un isomorphisme. D'après le
théorème~\ref{17-theoreme-stabilite-d-b-ctf-par-rhom}, le foncteur de dualité
considéré commute au changement d'anneau, ainsi, après produit tensoriel
avec $A/\mathfrak m$, on obtient le morphisme de bidualité pour
$A/\mathfrak m\overset{\L}\otimes_A L$ dans $\detalebctfcoefs X
{A/\mathfrak m}$. Bref, on peut supposer que l'anneau $A$ est un corps.
Dans ce cas, on peut conclure en utilisant le
théorème~\ref{17-theoreme-dualite-d-b-c-coefficients}.

\begin{proposition2}\label{17-proposition-r-j-star-preserve-d-b-c-si-cplx-dualisant}
Soit $X$ un schéma noethérien. On suppose qu'il existe un complexe
dualisant sur $\detalebccoefs X \Lambda$. Alors, pour toute immersion $j$ d'un
ouvert $U$ de $X$, le foncteur $\R j_\star$ envoie
$\detalebc U $ dans $\detalebc X$.
\end{proposition2}

Soit $K$ un complexe dualisant sur $\detalebc X$. Pour des raisons
évidentes, $j^\star K$ est un complexe dualisant sur $\detalebc U$.
On note $\Dual_X$ (resp. $\Dual_U$) les dualités induites par $K$
et $j^\star K$ sur les catégories triangulées $\detalebc X $ (resp.
$\detalebc U $). On a un isomorphisme canonique $\Dual_X\circ
j_!\simeq \R j_\star \Dual_U$. On en déduit que pour tout
$M\in\detalebc U $, $\R j_\star M\simeq \Dual_Xj_!\Dual_U M$,
ce qui permet de conclure que $\R j_\star M$ appartient à $\detalebc X$.

\begin{definition2}
Si $X$ est un schéma noethérien, $A$ une $\Lambda$-algèbre noethérienne
et $\mathscr P$
une partition galoisienne dirigeable de $X$, on note
${\detalecoefs X A}^{\mathscr P}$ la sous-catégorie triangulée de
$\detalecoefs X A$
dont les objets de cohomologie sont dans $\FCons^{\mathscr P}(X,A)$. On
définit de même les variantes ${\detalebcoefs X A}^{\mathscr P}$,
${\detalebccoefs X A}^{\mathscr P}$, etc.
\end{definition2}

\begin{proposition2}\label{17-proposition-r-j-star-preserve-d-b-partitionne}
Soit $j\colon U\to X$ une immersion ouverte entre schémas noethériens. On
suppose que $\R j_\star$ applique $\detalebc U$ dans $\detalebc X$. Pour
toute partition galoisienne dirigeable $\mathscr P$ de $U$, il existe une
partition galoisienne $\mathscr P'$ de $X$ et un entier $c$ tel que
$\R j_\star$ envoie ${\detaleb U}^{\mathscr P}$ dans
${\detaleb X}^{\mathscr P'}$ et que pour tout $q>c$ et
$\mathscr F\in\FCons^{\mathscr P}(U,\Lambda)$, on ait $\R^q j_\star
\mathscr F=0$.
\end{proposition2}

On sait que la catégorie $\FCons^{\mathscr P}(U,\Lambda)$ est localement
finie et admet même un nombre fini d'objets simples (et ceux-ci sont des
faisceaux constructibles). Comme $\R j_\star$ applique $\detalebc U$ dans
$\detalebc X$, on peut choisir un entier $c$ et une partition galoisienne
$\mathscr P'$ de $X$ tels que pour tout objet simple $\mathscr F$ de
$\FCons^{\mathscr P}(U,\Lambda)$, $\R j_\star \mathscr F$
appartienne à ${\detaleb X}^{\mathscr P'}$ et ait des objets de cohomologie
nuls en les degrés strictement supérieurs à $c$. Ce résultat s'étend par
dévissage aux objets noethériens de $\FCons^{\mathscr P}(U,\Lambda)$ puis à
cette catégorie toute entière du fait de la commutation des foncteurs
$\R^q j_\star$ aux limites inductives filtrantes.

\begin{corollaire2}\label{17-corollaire-i-bang-preserve-d-b-partitionne}
Soit $i\colon Z\to X$ une immersion ouverte entre schémas noethériens. On
suppose que $i^!$ applique $\detalebc X$ dans $\detalebc Z$. Pour toute partition
galoisienne dirigeable $\mathscr P$ sur $X$, il existe une partition
galoisienne $\mathscr P'$ sur $Z$ et un entier $c$ tel que
$i^!$ envoie ${\detaleb X}^{\mathscr P}$ dans ${\detaleb Z}^{\mathscr P'}$ et que
pour tout $q>c$ et $\mathscr F\in\FCons^{\mathscr P}(X,\Lambda)$, on ait
$\mathscr H^q(i^!\mathscr F)=0$.
\end{corollaire2}

Si on note $j\colon U \to X$ l'immersion ouverte complémentaire,
l'hypothèse sur $i^!$ énoncée ici équivaut à celle exigée sur
$\R j_\star$ dans la
proposition~\ref{17-proposition-r-j-star-preserve-d-b-partitionne}. Quitte à
raffiner $\mathscr P$, on peut supposer que les constituants de $\mathscr
P$ sont soit au-dessus de $U$, soit au-dessus de $Z$. On peut ainsi écrire
$\mathscr P=\mathscr P_U\cup \mathscr P_Z$ où $\mathscr P_U$ et $\mathscr
P_Z$ sont des partitions galoisiennes de $U$ et $Z$ respectivement. On
applique la proposition~\ref{17-proposition-r-j-star-preserve-d-b-partitionne}
à $\mathscr P_U$. On obtient une partition galoisienne $\mathscr P''$ de
$X$ telle que $\R j_\star $ applique ${\detaleb U}^{\mathscr
P_U}$ dans ${\detaleb X}^{\mathscr P''}$ et un entier naturel  $c'$
tel que pour tout
$\mathscr F\in\FCons^{\mathscr P_U}(U,\Lambda)$, on ait $\R^q j_\star
\mathscr F=0$  pour $q>c'$. Quitte à raffiner $\mathscr P''$, on peut
supposer que $\mathscr P''=\mathscr P''_Z\cup \mathscr P''_U$ comme
ci-dessus. Quitte à raffiner $\mathscr P''_Z$, on peut supposer que cette
partition galoisienne de $Z$ raffine $\mathscr P_Z$. En utilisant le
triangle distingué
\[i^!K\to i^\star K\to i^\star\R j_\star j^\star K\vers +\]
pour tout $K\in\detaleplus X$, on obtient aussitôt que $\mathscr P'=\mathscr
P''_Z$ et $c=c'+1$ conviennent.

\begin{proposition2}\label{17-proposition-r-j-star-preserve-d-b-ctf-coefficients}
Soit $j\colon U\to X$ une immersion ouverte entre schémas noethériens. On
suppose que $\R j_\star$ envoie $\detalebc U$ dans $\detalebc X$. Alors, pour
toute $\Lambda$-algèbre noethérienne $A$, $\R j_\star$ envoie
$\detalebccoefs U A$ (resp. $\detalebctfcoefs U A$)
dans $\detalebccoefs X A$ (resp. $\detalebctfcoefs X A$).
En outre, pour tout $Y\in\detalebctfcoefs U A$ et $M\in\dplus(A)$, le morphisme
canonique $M\overset{\L}\otimes_A \R j_\star Y\to \R
j_\star (M\overset{\L}\otimes_A Y)$ est un isomorphisme.
Par ailleurs, si $i\colon Z \to X$ est une immersion fermée complémentaire
à $j$, alors $i^!$ envoie $\detalebccoefs X A$ (resp. $\detalebctfcoefs X A$) dans
$\detalebccoefs Z A$ (resp. $\detalebctfcoefs Z A$), et pour tout
$M\in\dplus(A)$
et $Y\in\detalebccoefs X A$, le morphisme canonique
$M\overset{\L}\otimes_A i^!Y\to i^!(M\overset{\L}\otimes_A Y)$
est un isomorphisme.
\end{proposition2}

L'énoncé sur $i^!$ se déduit aussitôt de celui sur $\R j_\star$, on
se concentre donc sur celui-là. Pour montrer que $\R j_\star$ envoie
$\detalebccoefs U A$ dans $\detalebccoefs X A$, on peut supposer par un dévissage
évident que $\Lambda=\mathbf{Z}/\ell\mathbf{Z}$, où $\ell$ est un nombre
premier. Par conséquent, pour tout $\Lambda$-module $N$ et tout objet
$Y\in\detalebccoefs U \Lambda$, on a un isomorphisme $\R j_\star
(N\overset{\L}\otimes_\Lambda Y)\simeq
N\overset{\L}\otimes_\Lambda \R j_\star Y$ (cf.
proposition~\ref{17-proposition-coefficients-universels-r-j-star}). Si $N$ a
une structure de $A$-module de type fini, comme on sait que $\R
j_\star Y$ appartient à $\detalebc X$, on peut en déduire que $\R j_\star
(N\overset{\L}\otimes_\Lambda Y)$ appartient à $\detalebccoefs X A$.
D'après le dévissage de la
proposition~\ref{17-proposition-devissage-coefficients}, il vient que
$\R j_\star$ envoie $\detalebccoefs U A$ dans $\detalebccoefs X A$.

Montrons maintenant que $\R j_\star$ envoie $\detalebctfcoefs U A$ dans
$\detalebctfcoefs X A$. Compte tenu du résultat précédent, il suffit de montrer
que si $\mathscr F$ est un faisceau de $A$-modules plat et constructible
sur $U$, alors $\R j_\star \mathscr F\in\detalebtfcoefs X A$,
c'est-à-dire que $M\overset{\L}\otimes_A\R j_\star \mathscr F$
est borné indépendamment du $A$-module $M$. D'après la
proposition~\ref{17-proposition-coefficients-universels-r-j-star}, il suffit
de montrer que $\R j_\star(\mathscr F\otimes_A M)$ est borné
indépendamment du $A$-module $M$, et si tel est le cas, la formule des
coefficients universels énoncée ici sera satisfaite. Il existe une
partition galoisienne dirigeable $\mathscr P$ de $U$ tel que $\mathscr F$
appartienne à $\FCons^{\mathscr P}(U,A)$. Pour tout $A$-module $M$,
$\mathscr F\otimes_A M$ est un objet de $\FCons^{\mathscr P}(U,A)$, ainsi,
il suffit de montrer qu'il existe un entier $c$ tel que pour tout objet
$\mathscr G$ de $\FCons^{\mathscr P}(U,A)$, on ait $\R^q j_\star
\mathscr G=0$ pour $q>c$, ce qui résulte de la
proposition~\ref{17-proposition-r-j-star-preserve-d-b-partitionne}.

\begin{definition2}
Soit $X$ un schéma noethérien. Soit $A$ une $\Lambda$-algèbre noethérienne.
Soit $K\in\detalebctfcoefs X A$. On dira que $K$ vérifie la condition (B) si
pour toute partition galoisienne $\mathscr P$ de $X$, il existe un entier
$c$ et une partition galoisienne $\mathscr P'$ de $X$ tels que pour
tout $L\in{\detalebcoefs X A}^{\mathscr P}$, $\R\SheafHom_A(K,L)$
appartienne à ${\detalebcoefs X A}^{\mathscr P'}$, que si on suppose que
$\mathscr H^qL=0$ pour $q>0$, alors 
$\mathscr H^q\R\SheafHom_A(K,L)=0$ pour $q>c$ et enfin, que si $L$
appartient à ${\detalebccoefs X A}^{\mathscr P}$, alors
$\R\SheafHom_A(K,L)$ appartient à ${\detalebccoefs X A}^{\mathscr P'}$.
\end{definition2}

\begin{proposition2}\label{17-proposition-tout-objet-de-d-b-ctf-verifie-b}
Soit $X$ un schéma noethérien tel qu'il existe un complexe dualisant sur
$\detalebc X$. Alors, tout objet $K\in\detalebctfcoefs X A$ vérifie la condition
(B).
\end{proposition2}

\begin{lemme2}
Soit $X$ un schéma noethérien tel qu'il existe un complexe dualisant sur
$\detalebc X$. Soit $i\colon Z\to X$ une immersion fermée. Soit $j\colon U\to
X$ l'immersion ouverte complémentaire. Soit $K\in\detalebctfcoefs X A$. On
suppose que $i^\star K$ et $j^\star K$ satisfont la condition (B). Alors, $K$
satisfait la condition (B).
\end{lemme2}

Pour tout $L\in\detalecoefs X A$, on a un
triangle distingué dans $\detalecoefs X A$ :
\[i_\star \R\SheafHom_A(i^\star K,i^!
L)\to\R\SheafHom_A(K,L)\to\R j_\star
\R\SheafHom_A(j^\star K,j^\star L)\vers +\]
Grâce au résultat de l'exercice~\ref{17-exercice-partitions-ensemble-filtrant},
on peut combiner d'une part le résultat sur $i^!$ du
corollaire~\ref{17-corollaire-i-bang-preserve-d-b-partitionne} et la condition
(B) pour $i^\star K$ et d'autre part la 
proposition~\ref{17-proposition-r-j-star-preserve-d-b-partitionne} concernant
$\R j_\star$ et la condition (B) pour $j^\star K$ pour obtenir que
$K$ vérifie la condition (B).

\bigskip

Démontrons la 
proposition~\ref{17-proposition-tout-objet-de-d-b-ctf-verifie-b}. La condition (B) définit une sous-catégorie
triangulée de $\detalebctfcoefs X A$. Pour montrer la proposition, il suffit de
montrer que si $K$ est un faisceau de $A$-modules plat et constructibles,
alors $K$ satisfait la condition (B).
L'existence d'un complexe dualisant étant une
condition préservée par passage à un sous-schéma, le lemme précédent
fournit un moyen de dévisser la situation pour se ramener au cas où $K$ est
localement constant. On est ramené au lemme suivant :

\begin{lemme2}
Soit $X$ un schéma noethérien. Soit $\mathscr F$ un faisceau de $A$-modules
constructible, plat et localement constant. Alors, $\mathscr F$ satisfait
la propriété (B).
\end{lemme2}

Tout d'abord, pour tout $L\in\detalebcoefs X A$, si
$\mathscr H^qL=0$ pour $q>0$, alors pour tout $q>0$, $\mathscr
H^q(\R\SheafHom_A(\mathscr F,L))=0$, et si $L\in\detalebccoefs X A$, alors
$\R\SheafHom_A(K,L)\in\detalebccoefs X A$.
Il reste donc à montrer que si $\mathscr P$ est une partition galoisienne
de $X$, il existe une partition galoisienne $\mathscr P'$ de $X$ telle que
pour tout $L\in{\detalebcoefs X A}^{\mathscr P}$, alors
$\R\SheafHom_A(\mathscr F,L)$ appartient à
$L\in{\detalebcoefs X A}^{\mathscr P'}$. On peut supposer que $\mathscr P$ est
constitué d'un unique revêtement étale galoisien $X'\to X$. On choisit
un revêtement étale galoisien $X''\to X$ tel que l'image inverse de
$\mathscr F$ sur $X''$ soit un faisceau constant, puis un
revêtement galoisien $X'''\to X$ coiffant $X'$ et $X''$. On voit aussitôt
que la partition galoisienne $\mathscr P'=(X'''\to X)$ de $X$ convient.

\bigskip

Démontrons le théorème~\ref{17-theoreme-stabilite-d-b-ctf-par-rhom}. Soit
$K\in\detalebctfcoefs X A$. D'après la
proposition~\ref{17-proposition-tout-objet-de-d-b-ctf-verifie-b}, $K$ vérifie
la condition (B). Soit $L\in\detalebctfcoefs X A$. Il existe une partition
galoisienne $\mathscr P$ de $X$ telle que $L$ appartienne à 
${\detalebcoefs X A}^{\mathscr P}$. Pour tout $A$-module $M$, l'objet
$M\overset{\L}\otimes_A L$ appartient encore à cette catégorie (et
est borné indépendamment de $M$). Il résulte de la condition (B) de $K$ que
$\R\SheafHom_A(K,L)$ appartient à $\detalebccoefs X A$ et qu'il
existe une partition galoisienne $\mathscr P'$ telle que
$\R\SheafHom_A(K,M\overset{\L}\otimes_A L)$ soit un objet de
${\detalebcoefs X A}^{\mathscr P'}$ borné indépendamment du $A$-module $M$.
La proposition~\ref{17-proposition-coefficients-universels-1} permet de
déduire que $\R\SheafHom_A(K,L)$ appartient à $\detalebtfcoefs X A$ et
que pour tout $M\in\dplus(A)$, le morphisme canonique
$M\overset{\L}{\otimes}_A \R\SheafHom_A(K,L)\to
\R\SheafHom_A(K,M\overset{\L}{\otimes}_A L)$ est un
isomorphisme. On déduit aussitôt de cette formule la compatibilité au
changement d'anneau $A\to A'$ pour toute $A$-algèbre $A'$.

\subsection{Élimination d'hypothèses noethériennes sur \texorpdfstring{$A$}{A}}

\begin{definition2}\label{17-definition-c-perfection}
Soit $A$ un anneau commutatif.
Soit $X$ un schéma noethérien. On dit d'un complexe
$K\in \detalecoefs X A$ qu'il est $\derivec$-parfait s'il existe une partition
finie $(U_i)_{i\in I}$ de $X$ par des sous-schémas (réduits) telle que pour
tout $i\in I$, $K_{|U_i}\in \detalecoefs {U_i} A$ soit un complexe parfait (cf.
\sga{6}{i}{4.8}). On note $\detalebcparfcoefs X A$ la sous-catégorie triangulée de
$\detalecoefs X A$ formée des complexes $\derivec$-parfaits.
\end{definition2}

Bien entendu, pour tout morphisme d'anneaux $A\to A'$, le foncteur
$A'\overset{\L}\otimes_A -\colon \detalecoefs X A\to \detalecoefs X {A'}$
induit un foncteur $\detalebcparfcoefs X A\to \detalebcparfcoefs X {A'}$.
En outre, si $A$ est un anneau noethérien,
$\detalebcparfcoefs X A=\detalebctfcoefs X A$.

\begin{theoreme2}\label{17-theoreme-dualite-d-b-c-parf-coefficients}
Soit $A$ une $\Lambda$-algèbre commutative. Soit $X$ un schéma
noethérien. S'il en existe, les
complexes dualisants sur $\detalebcparfcoefs X A$ sont uniques au produit
tensoriel près avec des objets inversibles.
Soit $K$ un complexe dualisant sur $\detalebccoefs X \Lambda$.
Alors, $A\overset{\L}\otimes_{\Lambda} K$ est un complexe dualisant
sur $\detalebcparfcoefs X A$.
En outre, le bifoncteur $\R\SheafHom_A$ préserve $\detalebcparfcoefs X A$
et commute à tout changement d'anneau $A\to A'$.
\end{theoreme2}

Le théorème~\ref{17-theoreme-stabilite-d-b-ctf-par-rhom} énonce
que si $B$ est une $\Lambda$-algèbre noethérienne, que $K$ et $L$
sont deux objets de $\detalebcparfcoefs X B$, alors pour toute $B$-algèbre
$A$, on a un isomorphisme canonique
\[\R\SheafHom_A(A\overset{\L}\otimes_B K,
A\overset{\L}\otimes_B
L)\simeq A\overset{\L}\otimes_B 
\R\SheafHom_B(K,L)\;\text{.}\]
Comme l'objet de droite appartient à $\detalebcparfcoefs X A$, il s'agit d'un
isomorphisme dans $\detalebcparfcoefs X A$.
Bref, compte tenu du
théorème~\ref{17-theoreme-stabilite-d-b-ctf-par-rhom} et
du théorème~\ref{17-theoreme-dualite-d-b-ctf-coefficients} (dont l'énoncé
d'unicité des complexes dualisants vaut aussi pour $\detalebcparfcoefs X A$ avec
la même démonstration), le théorème ci-dessus est ramené au lemme suivant :

\begin{lemme2}\label{17-lemme-approximation-noetherienne}
Soit $A$ un anneau commutatif.
Soit $X$ un schéma noethérien. Pour tout objet $K$ de
$\detalebcparfcoefs
X A$, il existe un sous-anneau noethérien (et même de type fini sur
$\mathbf{Z}$) $B$ de $A$ et $K'\in \detalebcparfcoefs X B$ tel que les objets
$K$ et $A\overset{\L}\otimes_B K'$ de $\detalebcparfcoefs X A$ soient
isomorphes.
\end{lemme2}

Comme il nécessite un examen plus attentif de la notion de
$\derivec$-perfection, on repousse la démonstration de ce lemme à la fin de
cette sous-section.

\begin{lemme2}
Soit $A$ un anneau commutatif. Soit $\mathscr F$ un faisceau de $A$-modules
sur $X_\et$. Les conditions suivantes sont équivalentes :
\begin{itemize}
\item[(i)] Il existe une partition $(U_i)_{i\in I}$ de $X$ par des
sous-schémas réduits tels que pour tout $i\in I$,
$\mathscr F_{|U_i}$ soit localement constant et que pour tout point
géométrique $\overline x$ de $X$, le $A$-module $\mathscr F_{\overline x}$
soit projectif de type fini ;
\item[(ii)] Le faisceau de $A$-modules $\mathscr F$ est
constructible\;\footnote{On utilise la définition donnée dans
\sga{4}{ix}{2.3} même si $A$ n'est pas noethérien, et non pas la définition
suggérée en note à cet endroit ; un faisceau constructible pour cette autre
définition est ce que pourrions appeler un faisceau
$\derivec$-pseudo-cohérent.} et pour tout point géométrique $\overline x$ de
$X$, le $A$-module $\mathscr F_{\overline x}$ est projectif  de type fini ;
\item[(iii)] Le faisceau de $A$-modules $\mathscr F$ est plat et
constructible.
\end{itemize}
\end{lemme2}

Par définition des faisceaux constructibles, on a évidemment l'équivalence
$\text{(i)}\Longleftrightarrow\text{(ii)}$. Si $\mathscr F$ est
constructible, les fibres $F_{\overline x}$ sont des $A$-modules de
présentation finie, il est alors équivalent d'exiger que ces modules soient
plats ou projectifs de type fini, ce qui montre l'équivalence
$\text{(ii)}\Longleftrightarrow\text{(iii)}$. 

On note $\mathscr C$ la catégorie fibrée au-dessus du site $X_\et$ qui à
$U\in X_\et$ fait correspondre la catégorie des faisceaux de $A$-modules
sur $U_\et$ et $\mathscr C_\derivec$ la sous-$X_\et$-catégorie de $\mathscr C$
formée des faisceaux de $A$-modules plats et constructibles. Nous allons
utiliser la terminologie de \sga{6}{i}{1.2}.
Il est évident qu'un objet de $\mathscr C$ qui est localement dans
$\mathscr C_\derivec$ est dans $\mathscr C_\derivec$ et que $\mathscr C_\derivec$ est stable par
noyau d'épimorphisme.
D'après \sga{4}{ix}{2.7}, un faisceau de $A$-modules sur $X_\et$ est
constructible si et seulement s'il est isomorphe au conoyau d'un morphisme
$A_V\to A_U$ pour $U$ et $V$ deux $X$-schémas étales et de présentation
finie. Les faisceaux $A_U$ pour $U$ étale et de présentation finis
sur $X$ sont donc évidemment plats et constructibles.
Il résulte de ces résultats qu'un objet de $\mathscr C_X$ est de
$\mathscr C_\derivec$-type fini si et seulement s'il est engendré par un nombre
fini de sections (ce qui revient à demander qu'il soit de $\mathscr
C_{0X}$-type fini) et qu'un objet de $\mathscr C_X$ est de $\mathscr
C_\derivec$-présentation finie si et seulement s'il est constructible (ce qui
revient encore à demander qu'il soit de $\mathscr C_{0X}$-présentation
finie). Il est par ailleurs évident que $\mathscr C_\derivec$ est quasi-relevable
dans $\mathscr C$ et même que $\mathscr C_{0X}$ est quasi-relevable dans
$\mathscr C_X$. Les catégories $\mathscr C_\derivec$ et $\mathscr C$ vérifient
donc les hypothèses de \sga{6}{i}{2.0} et de \sga{6}{ii}{1.1} (mais en général
pas de \sga{6}{i}{4.0}). On dispose donc d'une notion de complexe
pseudo-cohérent (relativement à $\mathscr C_\derivec)$ et celle-ci peut-être
définie de façon globale :

\begin{definition2}
Soit $A$ un anneau. Soit $X$ un schéma noethérien.
Soit $K\in C(X_\et,A)$. On dit que $K$ est strictement
$\derivec$-pseudo-cohérent
(resp. strictement $\derivec$-parfait) si
$K$ est un complexe borné supérieurement (resp. borné)
formé de faisceaux de $A$-modules plats et constructibles.
\end{definition2}

\begin{definition2}
Soit $A$ un anneau. Soit $X$ un schéma noethérien.
Soit $K\in \detalecoefs X A$. On dit que $K$ est
$\derivec$-pseudo-cohérent s'il est
isomorphe à l'image dans $\detalecoefs X A$ d'un complexe strictement
$\derivec$-pseudo-cohérent\;\footnote{En raison de l'existence de
\guil{résolutions globales}, la définition globale donnée ici
équivaut à la définition locale de \sga{6}{i}{2.3}.}.
Les objets $\derivec$-pseudo-cohérents forment une
sous-catégorie triangulée anonyme de $\detalecoefs X A$.
\end{definition2}

\begin{remarque2}
La $X_\et$-catégorie $\mathscr C$ contient aussi la $X_\et$-catégorie
$\mathscr C_0$ des faisceaux de $A$-modules facteurs directs de faisceaux
libres de $A$-modules. Bien entendu, $\mathscr C_0$ est contenue dans
$\mathscr C_{\derivec}$. Ainsi, les notions de stricte pseudo-cohérence (resp.
stricte perfection) définies relativement à $\mathscr C_0$ impliquent
les notions correspondantes relativement à $\mathscr C_\derivec$.
\end{remarque2}

\begin{proposition2}\label{17-proposition-critete-c-perfection}
Soit $A$ un anneau. Soit $X$ un schéma noethérien.
Soit $K\in \detalecoefs X A$. Alors, les conditions suivantes sont équivalentes
:
\begin{itemize}
\item[(i)] $K$ est $\derivec$-parfait ;
\item[(ii)] $K$ est $\derivec$-pseudo-cohérent et de tor-dimension finie ;
\item[(iii)] $K$ est isomorphe à l'image dans $\detalecoefs X A$ d'un complexe
strictement $\derivec$-parfait.
\end{itemize}
\end{proposition2}

\begin{lemme2}
Soit $A$ un anneau. Soit $X$ un schéma noethérien. Soit $(U_i)_{i\in I}$
une partition finie de $X$ par des sous-schémas réduits. Soit $K\in
\detalecoefs X
A$. On suppose que pour tout $i\in I$, la restriction de $K$ à $U_i$
est $\derivec$-pseudo-cohérente. Alors, $K$ est $\derivec$-pseudo-cohérent.
\end{lemme2}

Par les arguments habituels, on se ramène au cas où la partition de $X$ est
constituée d'un ouvert $U$ et d'un fermé $Z$. On note $i\colon Z\to X$ et
$j\colon U\to X$ les immersions correspondantes. On dispose d'un triangle
distingué dans $\detalecoefs X A$ :
\[j_!j^\star K\to K\to i_\star i^\star K\vers +\]
On sait que $j^\star K$ et $i^\star K$ sont $\derivec$-pseudo-cohérents.
Il est
évident que $j_!$ et $i_\star$ préservent la notion de faisceau plat et
constructible ; au niveau des catégories triangulées, ces foncteurs
préservent donc évidemment la notion de $\derivec$-pseudo-cohérence. Les objets
$j_!j^\star K$ et $i_\star i^\star K$ sont $\derivec$-pseudo-cohérent ;
il en résulte que $K$ aussi est $\derivec$-pseudo-cohérent.

\begin{lemme2}
Soit $A$ un anneau. Soit $X$ un schéma noethérien. Soit $K\in\detalecoefs X A$.
Si $K$ est $\derivec$-parfait, alors $K$ est $\derivec$-pseudo-cohérent.
\end{lemme2}

D'après le lemme précédent, quitte à passer à un recouvrement fini par des
localement fermés convenables, on peut supposer que $K$ est parfait.
Par conséquent, $K$ est pseudo-cohérent (sous-entendu relativement à la
$X_\et$-catégorie $\mathscr C_0$) ; \emph{a fortiori}, $K$ est
$\derivec$-pseudo-cohérent.

\bigskip

Démontrons la proposition~\ref{17-proposition-critete-c-perfection}.
L'implication $\text{(iii)}\Longrightarrow\text{(i)}$ est évidente.
L'implication $\text{(i)}\Longrightarrow\text{(ii)}$ résulte
essentiellement du lemme précédent ; il reste cependant à vérifier que si
$K$ est $\derivec$-parfait, alors il est de tor-dimension finie. Supposons
donc que $K$ est $\derivec$-parfait. Il existe un recouvrement fini
$(U_i)_{i\in I}$ de $X$ par des localement fermés tels que $K_{|U_i}$ soit
parfait pour tout $i\in I$. Pour obtenir (ii) pour $K$, il suffit de
montrer que le complexe parfait $K_{|U_i}$ est de tor-dimension finie, ce
qui résulte aussitôt de \sga{6}{i}{5.8.1}. Il reste à établir l'implication
$\text{(ii)}\Longrightarrow\text{(iii)}$, la plus intéressante pour nous.
Soit $K\in \detalecoefs X A$ un complexe $\derivec$-pseudo-cohérent
et de tor-dimension finie.
Par définition de la $\derivec$-pseudo-cohérence,
on peut remplacer si besoin est $K$ par un
complexe borné supérieurement formé de faisceaux de $A$-modules plats et
constructibles.
Pour tout $n\in\mathbf{Z}$, on peut considérer la troncature canonique
$\tau_{\leq n}K$ de $K$, si on note $Z^n$ le noyau de $K^n\to K^{n+1}$, il
s'agit du sous-complexe suivant de $K$ :
\[\dots \to K^{n-2} \to K^{n-1} \to Z^n \to 0\to \dots\;\text{.}\]
Comme $K$ est de tor-dimension finie, il existe un entier $a\in\mathbf{Z}$
tel que pour tout $A$-module $M$, $\tau_{\leq a}(K\otimes_A M)$ soit
acyclique. En appliquant ceci avec $M=A$, on obtient une résolution plate
de $Z^a$ :
\[\dots\to K^{a-3}\to K^{a-2}\to K^{a-1}\to Z^a \to 0\;\text{.}\]
Ensuite, on obtient aussitôt que pour un $A$-module $M$ quelconque, cette
suite reste exacte après passage au produit tensoriel avec $M$. Par suite,
pour tout $i>0$, $\SheafTor_i^A(Z^a,M)=0$ pour tout $A$-module $M$,
ce qui implique que $Z^a$ est un faisceau de $A$-modules plat. Par
ailleurs, $Z^a$ est le conoyau du morphisme $K^{a-2}\to K^{a-1}$, donc
$Z^a$ est un faisceau constructible. Ainsi, $Z^a$ est plat et constructible.
Le complexe $K$ est quasi-isomorphe au complexe
strictement $\derivec$-parfait
\[\dots \to 0\to Z^a \to K^a \to K^{a+1} \to \dots\;\text{,}\]
ainsi $K$ vérifie la condition (iii).

\bigskip

Nous sommes maintenant en mesure de démontrer le
lemme~\ref{17-lemme-approximation-noetherienne}. Compte tenu de la
proposition~\ref{17-proposition-critete-c-perfection}, il s'agit de montrer
que si $K\in C(X_\et,A)$ est un complexe strictement $\derivec$-parfait,
alors $A$ contient un sous-anneau $A'$ de type fini sur $\mathbf{Z}$ tel
qu'il existe un complexe strictement $\derivec$-parfait $K'\in C(X_et,A')$
et un isomorphisme $K\simeq A\otimes_{A'} K'$. Ceci résulte aussitôt du
lemme suivant :

\begin{lemme2}
Soit $(A_\alpha)_{\alpha\in I}$ un système inductif d'anneaux commutatifs
indexé par un ensemble ordonné filtrant $I$. On en note $A$ la limite
inductive. Soit $X$ un schéma noethérien. Alors, la donnée d'un faisceau de
$A$-modules constructible (resp. plat et constructible) sur $X$ équivaut à
la donnée d'un faisceau de $A_\alpha$-modules constructible (resp. plat et
constructible) sur $X$ pour $\alpha$ assez grand.
\end{lemme2}

Conformément aux grands principes de \ega{iv}{8},
ceci signifie d'une part que si $F$ est un faisceau de $A$-modules
constructible sur $X$, il existe $\alpha\in I$ et $F_\alpha$ un faisceau de
$A_\alpha$-modules constructible sur $X$ tel que $F$ soit isomorphe à
$A\otimes_{A_\alpha} F_\alpha$ et d'autre part que si $\alpha\in I$ et que
$F_\alpha$ et $G_\alpha$ sont deux faisceaux de $A_\alpha$-modules
constructibles sur $X$, si on note
$F_\beta=A_\beta\otimes_{A_\alpha}F_\alpha$ et
$G_\beta=A_\beta\otimes_{A_\alpha}G_\alpha$ pour tout $\beta\geq \alpha$ et
$F=A\otimes_{A_\alpha}F_\alpha$ et $G=A\otimes_{A_\alpha}G_\alpha$, alors
l'application canonique
\[\underset{\beta\geq \alpha}\colim \Hom_{A_\beta}(F_\beta,G_\beta)\to
\Hom_A(F,G)\]
est un isomorphisme. En outre, si $\alpha\in I$ et que $F_\alpha$ est un
faisceau de $A_\alpha$-modules constructible sur $X$, alors
$F=A\otimes_{A_\alpha}F_\alpha$ est plat si et seulement si pour $\beta\geq
\alpha$ assez grand, $F_\beta=A_{\beta}\otimes_{A_\alpha} F_\alpha$ est
plat.

L'énoncé dans le cas non respé résulte de la description des faisceaux de
$B$-modules constructibles (pour tout anneau commutatif $B$) comme conoyau
d'une flèche $B_V\to B_U$ où $U$ et $V$ sont étales et de présentation
finie sur $X$, et du fait que le foncteur $\HH^0(V,-)$ de la catégorie des
faisceaux de groupes abéliens sur $X_\et$ vers celle des groupes abéliens
commute aux limites inductives filtrantes \sga{4}{vii}{3.3}.

Il reste à montrer que si $F_\alpha$ est un faisceau de $A_\alpha$-modules
constructible tel que, avec les notations ci-dessus, $F$ soit $A$-plat,
alors pour $\beta\geq \alpha$ assez grand, $F_\beta$ est $A_\beta$-plat. En
utilisant une décomposition de $X$ en réunion de localement fermés connexes
au-dessus desquels $F_\alpha$ soit localement constant trivialisé par un
revêtement étale, on peut supposer que $X$ est connexe et que
$F_\alpha$ est localement constant et trivialisé par un revêtement étale
$Y\to X$. Pour vérifier la platitude de $F_\beta$, il suffit de l'obtenir
pour \emph{une} fibre $(F_\beta)_{\overline x}$ ; on est ainsi ramené au
lemme suivant :

\begin{lemme2}
Soit $(A_\alpha)_{\alpha\in I}$ un système inductif d'anneaux commutatifs
indexé par un ensemble ordonné filtrant $I$. On en note $A$ la limite
inductive. Soit $X$ un schéma noethérien. Soit $\alpha\in I$, soit
$M_\alpha$ un $A_\alpha$-module de présentation finie. On suppose que
$M=A\otimes_{A_\alpha} M_\alpha$ est $A$-plat. Alors, il existe $\beta\geq
\alpha$ tel que $M_\beta=A_\beta\otimes_{A_\alpha} M_\alpha$ soit
$A_\alpha$-plat.
\end{lemme2}

Les modules considérés étant de présentation finie, le module $M$ (resp.
$M_\beta$) est plat si et seulement s'il est facteur direct d'un module
libre de type fini. On peut conclure en utilisant convenablement
\ega{iv}{8.5.2}.


\ifx\danslelivre\undefined
\appendix
\fi

\section{Produits tensoriels de complexes non bornés}
\subsection{\texorpdfstring{$K$-platitude}{K-platitude}}

\begin{definition2}
Soit $(\mathscr T,\mathscr A)$ un topos annelé en anneaux commutatifs. On
note $C(\mathscr T,\mathscr A)$ (resp. $K(\mathscr T,\mathscr A)$,
$\dtotale(\mathscr T,\mathscr A)$) la catégorie des complexes de $\mathscr
A$-Modules (resp. la catégorie homotopique correspondance, la catégorie
dérivée associée). On dispose d'un bifoncteur $\otimes_{\mathscr A}$ sur
$C(\mathscr T,\mathscr A)$ induisant un bifoncteur sur $K(\mathscr
T,\mathscr A)$\;\footnote{Plus précisément, si $K$ et $L$ sont deux
complexes de $\mathscr A$-Modules, on peut définir un bicomplexe dont la
composante de bidegré $(p,q)$ est $K^p\otimes_{\mathscr A}L^q$ et alors
$K\otimes_{\mathscr A} L$ est le complexe simple (défini en termes
de sommes) associé à ce bicomplexe. Nous n'imposerons pas ici de
convention de signes au lecteur.}.
Soit $K\in C(\mathscr T,\mathscr A)$. On dit que $K$ est
$K$-plat si le foncteur triangulé $K\otimes_{\mathscr A}-\colon K(\mathscr
T,\mathscr A)\to K(\mathscr T,\mathscr A)$ préserve les
quasi-isomorphismes, autrement dit que pour tout complexe acyclique $L$ de
$C(\mathscr T,\mathscr A)$, le complexe $K\otimes_{\mathscr A} L$ est
acyclique (cf. \cite[definition~5.1]{RUC@Spaltenstein}).
\end{definition2}

La sous-catégorie pleine de $C(\mathscr T,\mathscr A)$ formée des complexes
$K$-plats est stable par limites inductives filtrantes, sommes directes et
facteurs directs. En outre, la sous-catégorie pleine de $K(\mathscr
T,\mathscr A)$ correspondante est une sous-catégorie triangulée de
$K(\mathscr T,\mathscr A)$.

\begin{proposition2}
Soit $(\mathscr T,\mathscr A)$ un topos annelé en anneaux commutatifs. Soit
$K$ un complexe borné supérieurement formé de $\mathscr A$-modules plats
(cf. \sga{4}{v}{1}). Alors, $K$ est $K$-plat.
\end{proposition2}

En utilisant les foncteurs de troncature bête et la stabilité par limites
inductives filtrantes de la $K$-platitude, on se ramène au cas où $K$ est
borné. Comme la $K$-platitude définit une
sous-catégorie triangulée de $K(\mathscr T,\mathscr A)$, on peut procéder à
un dévissage utilisant encore les troncatures bêtes pour se ramener au cas
où $K^q=0$ pour $q\neq 0$. On est alors ramené à montrer que si $F$ est un
$\mathscr A$-Module plat et $L\in C(\mathscr T,\mathscr A)$ un complexe
acyclique, alors $F\otimes_{\mathscr A}L$ est acyclique, ce qui résulte
aussitôt de la définition de la platitude.

\subsection{Résolutions \texorpdfstring{$K$-plates}{K-plates}}
\subsubsection{Définition du produit tensoriel dérivé}

\begin{theoreme3}\label{17-theoreme-existence-resolutions-k-plates}
Soit $(\mathscr T,\mathscr A)$ un topos annelé en anneaux commutatifs. Il
existe un foncteur $\rho\colon C(\mathscr T,\mathscr A)\to C(\mathscr
T,\mathscr A)$ et une transformation naturelle $\rho K\to K$ pour $K\in
C(\mathscr T,\mathscr A)$ telle que :
\begin{itemize}
\item pour tout $K\in C(\mathscr T,\mathscr A)$, $\rho K$ soit $K$-plat ;
\item pour tout $K\in C(\mathscr T,\mathscr A)$, le morphisme $\rho K\to K$
soit un quasi-isomorphisme ;
\item le foncteur $\rho$ commute aux limites inductives filtrantes.
\end{itemize}
\end{theoreme3}

Ce théorème sera démontré plus bas. Déduisons-en aussitôt la proposition
triviale suivante, qui constitue notre définition du produit tensoriel sur
$\dtotale(\mathscr T,\mathscr A)$ :

\begin{proposition3}\label{17-proposition-existence-produit-tensoriel-derive}
Soit $(\mathscr T,\mathscr A)$ un topos annelé en anneaux commutatifs.
Le foncteur dérivé total à gauche de $\otimes_{\mathscr A}\colon C(\mathscr
T,\mathscr A)\times C(\mathscr T,\mathscr A)\to C(\mathscr T,\mathscr A)$
existe. Plus précisément, pour tous $K$ et $L$ dans $C(\mathscr T,\mathscr
A)$, on note $K\overset{\L}\otimes_{\mathscr A} L=(\rho
K)\otimes_{\mathscr A}(\rho L)\in C(\mathscr T,\mathscr A)$ ; ce bifoncteur
$\overset{\L}\otimes_{\mathscr A}$ commute aux limites inductives
filtrantes en chaque argument et, préservant les quasi-isomorphismes, il
induit un bifoncteur du même nom $\dtotale(\mathscr T,\mathscr A)\times
\dtotale(\mathscr
T,\mathscr A)\to \dtotale(\mathscr T,\mathscr A)$ ; la transformation naturelle
évidente $\overset{\L}\otimes_{\mathscr A}\to\otimes_{\mathscr A}$ 
fait de $\overset{\L}\otimes_{\mathscr A}$ le foncteur dérivé total
à gauche de $\otimes_{\mathscr A}$
(cf.~\cite[remark~7.4, Chapter~~II]{SHT@Goerss-Jardine} pour une définition
des foncteurs dérivés totaux en termes d'extensions de Kan). En outre,
si $K$ et $L$ sont deux objets de $C(\mathscr T,\mathscr A)$ dont l'un
au moins est $K$-plat, alors le morphisme canonique
$K\overset{\L}\otimes_{\mathscr A}L\to K\otimes_{\mathscr A}L$ est un
isomorphisme dans $\dtotale(\mathscr T,\mathscr A)$.
\end{proposition3}

\subsubsection{Modules sur un anneau}

On se place ici dans le cas particulier où le topos $\mathscr T$ est
ponctuel. On peut identifier les faisceaux de $\mathscr A$-modules à des
$A$-modules pour un anneau $A$. Le lemme suivant
démontre le théorème~\ref{17-theoreme-existence-resolutions-k-plates} dans ce
cas particulier.

\begin{lemme3}\label{17-lemme-existence-resolutions-k-plates-topos-ponctuel}
Pour tout anneau commutatif $A$, on peut définir un foncteur $\rho_A$ et
une transformation naturelle $\rho_A\to\Id$ de foncteurs de la
catégorie $C(A)$ des
complexes de $A$-modules dans elle-même telle que $\rho_A$ commute aux
limites inductives filtrantes, préserve les monomorphismes,
que pour tout $K\in C(A)$, le morphisme de
complexes $\rho_A(K)\to K$ soit un quasi-isomorphisme, que pour tout
entier relatif $n$, $\rho_A(K)^n$ soit un $A$-module libre, et que
$\rho_A(K)$ soit la limite inductive filtrante de ses sous-complexes
bornés formés de $A$-modules libres (en particulier, $\rho_A(K)$ est
$K$-plat).

On peut définir de telles résolutions $K$-plates $\rho_A$ pour tout anneau
commutatif $A$ de sorte que si $A\to A'$ est un morphisme d'anneaux, on ait
un morphisme fonctoriel $\rho_A(K)\to\rho_{A'}(K)$ dans $C(A)$ pour
$K\in C(A')$, ce morphisme vérifiant une compatibilité évidente
à la composition des morphismes d'anneaux.
\end{lemme3}

On note $G$ le foncteur adjoint à gauche du foncteur d'oubli $\mathrm{oub}$
de la catégorie des $A$-modules vers celle des ensembles pointés. On pose
$F=G\circ \mathrm{oub}$. Si $M$ est un $A$-module, $FM$ est le quotient du
$A$-module libre de base l'ensemble $M$ par le sous-module libre de rang
$1$ engendré par le zéro de $M$. Le morphisme d'adjonction $FM\to M$ est un
épimorphisme.
Pour tout
$A$-module $M$, on note $Z_0$ le noyau de cet épimorphisme $FM\to M$
et on pose $(F'M)_0=FM$.
Ensuite, de façon évidente, pour tout entier naturel $n\geq 1$,
on peut définir par récurrence un objet $(F'M)_n=FZ_{n-1}$,
un morphisme $d_n\colon (F'M)_n\to (F'M)_{n-1}$ et le noyau $Z_n=\ker d_n$.
Il est évident que l'on définit ainsi un complexe $F'M$ concentré en degrés
négatifs ou nuls, muni d'une augmentation $F'M\to M$ qui soit un
quasi-isomorphisme. Comme $F$ préserve les monomorphismes, on voit que $F'$
préserve aussi les monomorphismes.

Soit $K\in C(A)$. Le foncteur $F'$ défini ci-dessus n'est pas additif (à
moins que $A=0$), mais il est tel que $F'(0)=0$. Ainsi, si on applique
terme à terme le foncteur $F'$ aux objets $K_n$ pour tout $n\in\mathbf{Z}$,
on obtient un complexe double dans la catégorie des $A$-modules. On note
$\rho_A(K)$ le complexe simple associé (défini en termes de sommes). On
dispose bien entendu d'un morphisme d'augmentation $\rho_A(K)\to K$. Il
est évident que $\rho_A$ commute aux limites inductives filtrantes,
préserve les monomorphismes et que pour tout entier $n$, $\rho_A(K)^n$
soit une $A$-module libre. Si $K$ est borné supérieurement, le fait que
pour tout $n\in\mathbf{Z}$, le morphisme $\rho_A(K^n)\to K^n$ soit un
quasi-isomorphisme implique, par passage au complexe simple, que
$\rho_A(K)\to K$ est un quasi-isomorphisme. Comme tout complexe de
$A$-module peut s'écrire comme une limite inductive filtrante de
sous-complexes bornés supérieurement, il vient que pour tout $K\in C(A)$,
le morphisme $\rho_A(K)\to K$ est un quasi-isomorphisme. Comme
$\rho_A$ préserve les monomorphismes et commute aux limites inductives
filtrantes, pour montrer que $\rho_A(K)$ est une limite inductive
filtrante de ses sous-complexes bornés formés de $A$-modules libres, on
peut supposer que $K$ est borné supérieurement : le résultat est alors
trivial puisque, dans ce cas, $\rho_A(K)$ est un complexe borné
supérieurement formé de $A$-modules libres.

La dernière assertion concernant le changement d'anneau étant évidente, on
peut considérer que le lemme a été démontré.

\subsubsection{Préfaisceaux de Modules}
\label{17-subsubsection-produit-tensoriel-derive-prefaisceaux}

On suppose maintenant que $\mathscr T$ est le topos des préfaisceaux sur
une petite catégorie $\mathscr C$. Le faisceau d'anneaux $\mathscr A$ est
un préfaisceau d'anneaux commutatifs sur $\mathscr C$.

Soit $K\in C(\mathscr T,\mathscr A)$. Pour tout objet $U$ de $\mathscr C$,
$K(U)$ s'identifie à un objet de $C(\mathscr A(U))$. On applique la
construction du
lemme~\ref{17-lemme-existence-resolutions-k-plates-topos-ponctuel} à l'anneau
$\mathscr A(U)$. On pose $(\rho K)(U)=\rho_{\mathscr A(U)}(K(U))\in
C(\mathscr A(U))$. Si $V\to U$ est un morphisme dans $\mathscr C$, on
définit un $\mathscr A(U)$-morphisme $(\rho K)(U)\to (\rho K)(V)$ de la
façon suivante :
\[\rho_{\mathscr A(U)}(K(U))\to \rho_{\mathscr A(U)}(K(V))\to
\rho_{\mathscr A(V)}(K(V))\]
où le morphisme de gauche est induit par la structure de complexes de 
préfaisceaux de $\mathscr A$-modules sur $K$ et la flèche de droite par la
compatibilité de la construction du
lemme~\ref{17-lemme-existence-resolutions-k-plates-topos-ponctuel} au
changement d'anneau. D'après ce lemme, ces morphismes de transition
définissent une structure de préfaisceau sur $\rho K$. Ainsi, on a défini
un objet $\rho K\in C(\mathscr T,\mathscr A)$ et il est muni d'un morphisme
fonctoriel $\rho K \to K$.

Par construction, on peut vérifier les vertus présumées de $\rho$ terme à
terme ; ainsi, ce foncteur $\rho$ permet d'établir le
théorème~\ref{17-theoreme-existence-resolutions-k-plates} dans le cas où le
topos est un topos de préfaisceaux.

\subsubsection{Faisceaux de Modules}
\label{17-subsubsection-produit-tensoriel-derive-faisceaux}

\begin{proposition3}\label{17-proposition-stabilite-k-platitude-par-faisceau-associe}
Soit $\mathscr T$ le topos des faisceaux sur un site dont la catégorie
sous-jacente est notée $\mathscr C$. On note $\mathscr T'$ le topos des
préfaisceaux d'ensembles sur $\mathscr C$. Soit $\mathscr A'$ un préfaisceau
d'anneaux commutatifs sur $\mathscr C$. On note $\mathscr A$ le faisceau
d'anneaux $a\mathscr A'$ sur $\mathscr T$ associé à $\mathscr A'$.
Si $K\in C(\mathscr
T',\mathscr A')$ est $K$-plat, alors $aK\in C(\mathscr T,\mathscr A)$
est $K$-plat.
\end{proposition3}

\begin{lemme3}
Avec les notations de la
proposition~\ref{17-proposition-stabilite-k-platitude-par-faisceau-associe},
si $K$ et $L$ sont des objets de $C(\mathscr T',\mathscr A')$ tel que $aK$
soit nul dans $\dtotale(\mathscr T,\mathscr A)$, alors
$a(K\overset{\L}\otimes_{\mathscr A'} L)$
(où le produit tensoriel dérivé au-dessus de
$\mathscr A'$ est celui défini plus haut dans le cas des préfaisceaux) est
nul dans $\dtotale(\mathscr T,\mathscr A)$.
\end{lemme3}

La compatibilité du produit tensoriel dérivé sur $\dtotale(\mathscr T',\mathscr
A')$ avec les limites inductives filtrantes calculées dans $C(\mathscr
T',\mathscr A')$ permet de supposer que $L$ est un complexe borné
supérieurement (utiliser les troncatures canoniques). Quitte à remplacer
$L$ par la résolution $K$-plate $\rho_{\mathscr A'}L$ sus-définie, on peut
ensuite supposer que $L$ est borné supérieurement et constitué de $\mathscr
A'$-Modules plats. En utilisant les troncatures bêtes sur ce complexe $L$,
on peut supposer de plus que $L$ est borné. Finalement, en utilisant les
triangles distingués dans $K(\mathscr T',\mathscr A')$ donnés par les
troncatures bêtes, on peut finalement supposer que $L$ est formé d'un
unique $\mathscr A'$-Module plat placé en degré $0$.

Bref, il faut montrer que si $L$ est un $\mathscr A'$-Module plat et que
$K\in C(\mathscr T',\mathscr A')$ est tel que $aK$ soit un complexe
acyclique dans $C(\mathscr T,\mathscr A)$, alors le complexe de faisceaux 
$a(K\otimes_{\mathscr A'} L)$ est acyclique. On a un isomorphisme
$a(K\otimes_{\mathscr A'} L)\simeq aK\otimes_{\mathscr A}aL$, donc pour
conclure que $a(K\otimes_{\mathscr A'} L)$ est acyclique, il suffit de
montrer que le $\mathscr A$-Module $aL$ est plat, ce qui est vrai d'après
\sga{4}{v}{1.7.1}.

\bigskip

Montrons la
proposition~\ref{17-proposition-stabilite-k-platitude-par-faisceau-associe}.
Notons $\mathscr Z$ la sous-catégorie triangulée de $\dtotale(\mathscr T',\mathscr
A')$ formée des complexes qui sont annulés par le foncteur faisceau associé
$\dtotale(\mathscr T',\mathscr A')\to \dtotale(\mathscr T,\mathscr A)$. Le foncteur
induit $\dtotale(\mathscr T',\mathscr A')/\mathscr Z\to \dtotale(\mathscr T,\mathscr A)$
est évidemment une équivalence de catégories triangulées. D'après le lemme,
le bifoncteur $\overset{\L}\otimes_{\mathscr A'}$ sur $\dtotale(\mathscr
T',\mathscr A')$ passe au quotient par $\mathscr Z$ pour définit un
bifoncteur sur $\dtotale(\mathscr T,\mathscr A)$. La proposition en résulte
aussitôt. En effet, soit $K\in C(\mathscr T',\mathscr A')$ $K$-plat, soit
$L\in C(\mathscr T,\mathscr A)$ tel que $L$ soit nul dans $\dtotale(\mathscr
T,\mathscr A)$. On peut identifier $L$ à un objet $L'$ de $C(\mathscr
T',\mathscr A')$ et cet objet $L'$ appartient à la sous-catégorie
triangulée $\mathscr Z$ de $\dtotale(\mathscr T',\mathscr A')$. Le lemme montre
que $K\overset{\L}\otimes_{\mathscr A'} L'$ appartient à $\mathscr
Z$, autrement dit, $K$ étant $K$-plat, que $K\otimes_{\mathscr A'} L'$
appartient à $\mathscr Z$, c'est-à-dire que le complexe 
de faisceaux $aK\otimes_{\mathscr A} L$ est acyclique, ce qui montre que
$aK\in C(\mathscr T,\mathscr A)$ est $K$-plat.

\begin{corollaire3}
Avec les notations de la
proposition~\ref{17-proposition-stabilite-k-platitude-par-faisceau-associe},
le foncteur $\rho_{\mathscr A}$ qui à un complexe $K\in C(\mathscr
T,\mathscr A)$ associe $a\rho_{\mathscr A'} K'$ où $K'$ est le $\mathscr
A'$-Module défini par $K$ est un foncteur de résolution $K$-plate sur
$C(\mathscr T,\mathscr A$) vérifiant les conditions du
théorème~\ref{17-theoreme-existence-resolutions-k-plates}.
Ainsi, on dispose d'un bifoncteur $\overset{\L}\otimes_{\mathscr  A}$
sur $\dtotale(\mathscr T,\mathscr A)$ et d'un isomorphisme bifonctoriel
\[aK\overset{\L}\otimes_{\mathscr A}aL\simeq
a(K\overset{\L}\otimes_{\mathscr A'}L)\]
dans $\dtotale(\mathscr T,\mathscr A)$ pour tous $K$ et $L$ dans $\dtotale(\mathscr
T',\mathscr A')$.
\end{corollaire3}

Avec l'énoncé de ce corollaire s'achève la démonstration du
théorème~\ref{17-theoreme-existence-resolutions-k-plates}.

\subsection{Compléments}
\subsubsection{Homomorphismes internes}

\begin{definition3}
Soit $(\mathscr T,\mathscr A)$ un topos annelé en anneaux commutatifs. On
note $\SheafHom_{\mathscr A}$ le foncteur adjoint à droite du foncteur
produit tensoriel $\otimes_{\mathscr A}$ sur $C(\mathscr T,\mathscr A)$,
c'est-à-dire que pour $X$, $Y$ et $Z$ des objets de $C(\mathscr T,\mathscr
A)$, on a un isomorphisme canonique de groupes abéliens :
\[\Hom_{C(\mathscr T,\mathscr A)}(X\otimes_{\mathscr A} Y,Z)\simeq 
\Hom_{C(\mathscr T,\mathscr A)}(X,\SheafHom_{\mathscr A}(Y,Z))\;\text{.}\]
\end{definition3}

On rappelle qu'un objet $L\in C(\mathscr T,\mathscr A)$ est $K$-injectif si
pour tout complexe acyclique $K\in C(\mathscr T,\mathscr A)$, le complexe
de groupes abéliens $\Hom^\bullet_{\mathscr A}(K,L)$ (complexe simple
associé au complexe double $(\Hom_{\mathscr A}(K^{-p},L^q))_{p,q}$) est
acyclique et qu'il existe des foncteurs de résolution
$K$-injectives\;\footnote{Ce résultat est énoncé dans
\cite[§4]{RUC@Spaltenstein} dans le cadre des espaces topologiques annelés,
mais la démonstration peut être étendue au cas des topos annelés. Le
principe de la démonstration est similaire à celui utilisé par Grothendieck
pour montrer l'existence de suffisamment d'injectifs dans les catégories de
faisceaux dans \cite[§I.10]{Tohoku@Grothendieck} ; dans le contexte de l'algèbre
homotopique, cet argument est connu sous le nom de \guil{raisonnement du
petit objet}. On peut trouver une démonstration pour le cas qui nous
intéresse dans \cite{Complexes@Hovey}.}.

La proposition suivante, indiquée pour mémoire, est essentiellement
triviale :

\begin{proposition3}
Soit $(\mathscr T,\mathscr A)$ un topos annelé en anneaux commutatifs.
Le foncteur $\SheafHom_{\mathscr A}$ sur $C(\mathscr T,\mathscr A)$ admet
un foncteur dérivé total à droite $\R\SheafHom_{\mathscr A}\colon
\dtotale(\mathscr T,\mathscr A)^\opp
\times\dtotale(\mathscr T,\mathscr A)\to \dtotale(\mathscr T,\mathscr
A)$, adjoint à droite de $\overset{\L}\otimes_{\mathscr A}$. En
outre, si $K$ et $L$ sont des objets de $C(\mathscr T,\mathscr A)$, avec
$L$ $K$-injectif, alors le morphisme canonique 
\[\SheafHom_{\mathscr A}(K,L)\to\R\SheafHom_{\mathscr A}(K,L)\]
est un isomorphisme dans $\dtotale(\mathscr T,\mathscr A)$ ; si on suppose
de plus que $K$ est $K$-plat, alors $\SheafHom_{\mathscr A}(K,L)$ est
$K$-injectif. Enfin, si $K$, $L$ et $M$ sont trois objets de
$\dtotale(\mathscr T,\mathscr A)$, on a un isomorphisme fonctoriel dans
$\dtotale(\mathscr T,\mathscr A)$ \guil{cher à Cartan} :
\[\R\SheafHom_{\mathscr A}(K\overset{\L}\otimes_{\mathscr A} L,M)
\simeq \R\SheafHom_{\mathscr A}(K,
\R\SheafHom_{\mathscr A}(L,M))\;\text{.}\]
\end{proposition3}

\subsubsection{Compatibilité aux images inverses}

\begin{definition3}
Soit $u\colon (\mathscr T,\mathscr A)\to (\mathscr T',\mathscr A')$ un
morphisme de topos annelés en anneaux commutatifs. Le foncteur
$u^\star\colon C(\mathscr T,\mathscr A)\to C(\mathscr T',\mathscr A')$
admet un foncteur dérivé total à gauche que l'on note $\L u^\star$
(construit en appliquant $u^\star$ à une résolution $K$-plate).
\end{definition3}

\begin{proposition3}\label{17-proposition-u-star-preserve-k-platitude}
Soit $u\colon (\mathscr T,\mathscr A)\to (\mathscr T',\mathscr A')$ un
morphisme de topos annelés en anneaux commutatifs. Alors, pour tout $K\in
C(\mathscr T',\mathscr A')$ $K$-plat, $u^\star K\in C(\mathscr T,\mathscr
A)$ est $K$-plat.
\end{proposition3}

Pour démontrer cette proposition, il est commode d'introduire une nouvelle
définition :

\begin{definition3}
Soit $(\mathscr T,\mathscr A)$ un topos annelé en anneaux commutatifs.
Soit $K\in C(\mathscr T,\mathscr A)$. On dit que $K$ est
universellement acyclique si pour tout $\mathscr A$-module $L$, le complexe
$K\otimes_{\mathscr A} L$ est acyclique.
\end{definition3}

\begin{lemme3}
Soit $(\mathscr T,\mathscr A)$ un topos annelé en anneaux commutatifs.
Soit $K\in C(\mathscr T,\mathscr A)$. Les conditions suivantes sont
équivalentes :
\begin{itemize}
\item[(i)] $K$ est universellement acyclique ;
\item[(i')] pour tout $L\in C(\mathscr T,\mathscr A)$, le complexe
$K\otimes_{\mathscr A} L$ est acyclique ;
\item[(ii)] $K$ est $K$-plat et acyclique.
\end{itemize}
\end{lemme3}

L'équivalence entre (i) et (i') est facile : supposant (i), on montre
d'abord que la condition énoncée dans (i') est vérifiée par les complexes
$L$ bornés (dévissage), puis par tous les complexes (passage à la limite
inductive utilisant les troncatures bêtes et canoniques). L'implication
$\text{(i)}\Longrightarrow\text{(ii)}$ est triviale. L'implication
$\text{(ii)}\Longrightarrow\text{(i)}$ résulte aussitôt de la
proposition~\ref{17-proposition-existence-produit-tensoriel-derive}.

\medskip

De ce lemme, on déduit le suivant :

\begin{lemme3}
Soit $(\mathscr T,\mathscr A)$ un topos annelé en anneaux commutatifs.
Soit $K\in C(\mathscr T,\mathscr A)$. Si $K'\to K$ est un
quasi-isomorphisme, avec $K'$ $K$-plat, alors $K$ est $K$-plat si et
seulement si le cône de $K'\to K$ est universellement acyclique.
\end{lemme3}

\begin{lemme3}
Soit $u\colon (\mathscr T,\mathscr A)\to (\mathscr T',\mathscr A')$ un
morphisme de topos annelés en anneaux commutatifs.
Si on note $\rho_{\mathscr A'}$ \guil{le}
foncteur de résolution $K$-plate construit dans le 
paragraphe~\ref{17-subsubsection-produit-tensoriel-derive-faisceaux}, alors
pour tout $K\in C(\mathscr T',\mathscr A')$, l'objet $u^\star\rho_{\mathscr
A'}K\in C(\mathscr T,\mathscr A)$ est $K$-plat.
\end{lemme3}

Compte tenu du procédé de construction de $\rho_{\mathscr A'}$ par
extension aux complexes de faisceaux d'une construction pour les faisceaux
par passage au complexe simple, on peut supposer que $K$ est constitué
d'une unique faisceau placé en degré $0$. Il s'agit de montrer que pour
tout $\mathscr A'$-Modules $F$, le complexe (borné supérieurement)
$u^\star \rho_{\mathscr A'}F$ est constitué de $\mathscr A$-Modules plats.
En revenant à la construction du
lemme~\ref{17-lemme-existence-resolutions-k-plates-topos-ponctuel}, on peut
observer pour tout $n\in\mathbf{Z}$, $(\rho_{\mathscr A'}F)^n$ est un
faisceau de $\mathscr A'$-Modules libres sur un faisceau d'ensembles
pointés ; cette propriété étant évidemment préservée après application de
$u^\star$, on obtient que $(u^\star \rho_{\mathscr A'}F)^n$ est un
$\mathscr A$-Module plat.

\bigskip

Les lemmes précédents réduisent la
proposition~\ref{17-proposition-u-star-preserve-k-platitude} au lemme suivant :

\begin{lemme3}
Soit $u\colon (\mathscr T,\mathscr A)\to (\mathscr T',\mathscr A')$ un
morphisme de topos annelés en anneaux commutatifs. Si $K\in
C(\mathscr T',\mathscr A')$ est universellement acyclique, alors
$u^\star K\in C(\mathscr T,\mathscr A)$ est universellement acyclique.
\end{lemme3}

On peut démontrer ce lemme en suivant la méthode des limites inductives
locales utilisée dans \sga{4}{v}{8.2.9} ; en présence d'une famille
conservative de foncteurs fibres, il est possible de faire plus simple.

\begin{corollaire3}
Si $u$ et $v$ sont des morphismes composables de topos annelés en anneaux
commutatifs, alors on a un isomorphisme de foncteurs
$\L v^\star\circ \L u^\star\iso \L (u\circ v)^\star$.
\end{corollaire3}

\begin{corollaire3}
Soit $u\colon (\mathscr T,\mathscr A)\to (\mathscr T',\mathscr A')$ un
morphisme de topos annelés en anneaux commutatifs.
Pour tous $K$ et $L$ objets de
$\dtotale(\mathscr T',\mathscr A')$, on a un isomorphisme canonique 
\[\L u^\star K\overset{\L}\otimes_{\mathscr A}\L u^\star
L\simeq \L u^\star (K\overset{\L}\otimes_{\mathscr A'}L)\]
dans $\dtotale(\mathscr T,\mathscr A)$.
\end{corollaire3}

\section{Complexes inversibles}

\begin{proposition}\label{17-proposition-objets-inversibles-d-a}
Soit $A$ un anneau commutatif. Soit $X\in \dtotale(A)$. Soit $Y\in
\dtotale(A)$. Soit $X\overset{\L}{\otimes}_AY\simeq A$ un isomorphisme
dans $\dtotale(A)$. Alors, il existe une fonction localement constante
$k\colon \SP(A)\to \mathbf{Z}$, un $A$-module inversible $L$ et des
isomorphismes $X\simeq L[k]$ et $Y\simeq L^\vee [-k]$ (le foncteur de
décalage $[k]$ étant défini de façon évidente).
\end{proposition}

Ce résultat apparaît dans \cite[lemma~3.3, Chapter~V]{RD@Hartshorne}
sous des hypothèses supplémentaires disant que $A$ est noethérien et que
$X$ et $Y$ sont dans $\dmc(A)$. Cette version implique
évidemment une autre où on demande à $A$ d'être noethérien et à $X$ et $Y$
d'être des complexes parfaits. Il est alors facile de supprimer l'hypothèse
noethérienne (cf.~\ega{iv}{8}).
Bref, pour achever la démonstration, il suffit de montrer que
dans les conditions de la proposition ci-dessus, le complexe $X$ (et donc
$Y$ par symétrie des rôles) est un complexe parfait. On dit d'un objet $X$
de $\dtotale(A)$ qu'il est de présentation finie si le foncteur
$\Hom_{\dtotale(A)}(X,-)$ de $\dtotale(A)$ vers la catégorie des groupes
abéliens commute aux sommes directes (infinies). On peut montrer que
$X\in\dtotale(A)$ est un complexe parfait si et seulement s'il est de
présentation finie\;\footnote{Il s'agit d'un bon exercice. Toutefois, on
peut aussi obtenir ce critère en utilisant des principes généraux. L'objet
$A$ de $\dtotale(A)$ est un générateur de présentation finie (\emph{i.e.} le
foncteur cohomologique $\Hom_{\dtotale(A)}(A,-)$ commute aux sommes
directes et est conservatif) ; la sous-catégorie triangulée de $\dtotale(A)$
formée des objets de présentation finie est donc l'enveloppe
pseudo-abélienne de la sous-catégorie triangulée engendrée par $A$,
c'est-à-dire la sous-catégorie des complexes parfaits : combiner
\cite[proposition~8.4.1]{TC@Neeman}, \cite[lemma~4.4.5]{TC@Neeman} et
\cite[remark~4.2.6]{TC@Neeman}.}. De l'isomorphisme
$X\overset{\L}{\otimes}_A Y\simeq A$, on tire un isomorphisme de
foncteurs $\R\SheafHom(X,-)\simeq
Y\overset{\L}{\otimes}_A-\colon \dtotale(A)\to\dtotale(A)$. Comme le
foncteur $Y\overset{\L}{\otimes}_A-$ commute évidemment
aux sommes directes,
c'est aussi le cas du foncteur $\Hom_{\dtotale(A)}(X,-)$, ce qui
montre que $X$ est de présentation finie : il s'agit d'un complexe parfait.

\begin{proposition}\label{17-proposition-objets-inversibles-d-t-a}
Soit $(\mathscr T,\mathscr A)$ un topos annelé en anneaux commutatifs. Soit
$X\in\dtotale(\mathscr T,\mathscr A)$.
Soit $Y\in\dtotale(\mathscr T,\mathscr A)$. Soit
$X\overset{\L}{\otimes}_{\mathscr A}Y\simeq \mathscr A$ un
isomorphisme dans $\dtotale(\mathscr T,\mathscr A)$.
Alors, l'objet final de $\mathscr T$ est recouvert par des objets $U$
tels que $X_{|U}$ et $Y_{|U}$ puissent être 
induits par des complexes parfaits $X'$ et $Y'$ de $\dtotale(\mathscr A(U))$
et que l'on ait un isomorphisme $X'\overset{\L}{\otimes}_{\mathscr
A(U)}Y'\simeq \mathscr A(U)$ compatible à l'isomorphisme donné ; la
forme des complexes $X'$ et $Y'$ est donc connue grâce à la
proposition~\ref{17-proposition-objets-inversibles-d-a}.
\end{proposition}

Démontrons la proposition~\ref{17-proposition-objets-inversibles-d-t-a}. On
peut supposer que $\mathscr T$ est le topos des faisceaux sur un site
$\mathscr C$. Notons $\mathscr T'$ le topos des préfaisceaux sur la
catégorie sous-jacente au site $\mathscr
C$, elle aussi notée $\mathscr C$,
et $\mathscr A'$ le préfaisceau d'anneaux sur $\mathscr C$ défini par
$\mathscr A$. On peut identifier $X\overset{\L}\otimes_{\mathscr A}
Y$ au faisceau associé à $\rho_{\mathscr A'}X\otimes_{\mathscr
A'}\rho_{\mathscr A'}Y$ où $\rho_{\mathscr A'}$ est le foncteur de
résolution $K$-plate défini dans le
paragraphe~\ref{17-subsubsection-produit-tensoriel-derive-prefaisceaux}.
Ainsi, l'isomorphisme donné $\mathscr A\iso
X\overset{\L}\otimes_{\mathscr A}Y$ peut être représenté localement
par un $0$-cocycle $s$ de $(\rho_{\mathscr A'}X\otimes_{\mathscr
A'}\rho_{\mathscr A'}Y)(U)=\rho_{\mathscr A(U)}(X(U))\otimes_{\mathscr
A(U)}\rho_{\mathscr A(U)}(Y(U))$.
Soit $U$ un objet de $\mathscr C$ sur lequel une telle description
est possible. Par construction de $\rho_{\mathscr A(U)}$, il vient qu'il
existe des sous-complexes bornés formés de $\mathscr A(U)$-modules libres de
type fini $X'$ de $\rho_{\mathscr A(U)}(X(U))$ et $Y'$ de $\rho_{\mathscr
A(U)}(X(U))$ tel que $X'\otimes_{\mathscr A(U)} Y'$ contienne $s$. Notons
$\mathscr T_{|U}$ le topos des faisceaux sur $\mathscr T$ au-dessus de $U$
(un site sous-jacent est donné par la catégorie $\mathscr C/U$). On dispose
d'un morphisme de topos annelés évident $\pi$ de $(\mathscr T_{|U},\mathscr
A_{|U})$ vers le topos ponctuel muni de l'anneau $\mathscr A(U)$. Posons
$X''=\pi^\star X'$ et $Y''=\pi^\star Y'$ : ce sont des complexes strictement
parfaits sur $\mathscr T_{|U}$. On dispose de morphismes évidents $X''\to
X_{|U}$ et $Y''\to Y_{|U}$ dans $\dtotale(\mathscr T_{|U},\mathscr A_{|U})$
et d'après ce qui précède, on a un morphisme $\mathscr A_{|U}\to
X''\otimes_{\mathscr A_{|U}} Y''$ factorisant l'isomorphisme $\mathscr
A_{|U}\simeq X_{|U}\overset{\L}\otimes_{\mathscr A_{|U}}Y_{|U}$.
Bref, $\mathscr A_{|U}\simeq X_{|U}
\overset{\L}\otimes_{\mathscr A_{|U}}Y_{|U}$ est un
facteur direct de $X''\otimes_{\mathscr A_{|U}} Y''$. Par ailleurs,
l'isomorphisme $\mathscr A_{|U}\iso X_{|U}
\overset{\L}\otimes_{\mathscr A_{|U}}Y_{|U}$ se factorise aussi par
$X''\otimes_{\mathscr A_{|U}} Y_{|U}$, donc $\mathscr A_{|U}$
est un facteur direct de
$X''\otimes_{\mathscr A_{|U}} Y_{|U}$ ; en tensorisant ce fait avec
$X_{|U}$, on obtient que $X_{|U}$ est un facteur direct de $X''$. Soit
$p\colon X''\to X''$ un projecteur dont l'image soit isomorphe à $X_{|U}$.
Comme $X'$ et $Y'$ sont (strictement) parfaits, quitte à remplacer $U$ par
une famille d'objets le recouvrant, on peut supposer que $p$ provient d'un
projecteur $\tilde p$ sur $X'$ dans $\dtotale(\mathscr A(U))$. L'image de
$\tilde p$ est un complexe parfait $\tilde X$ dans $\dtotale(\mathscr
A(U))$ induisant $X_{|U}$. De même, on peut supposer que $Y_{|U}$ est induit
par un complexe parfait de $\mathscr A(U)$-modules $\tilde Y$. Quitte à
raffiner le recouvrement de $U$, on peut donc supposer que l'isomorphisme
$\mathscr A_{|U}\iso X_{|U}\otimes_{\mathscr A_{|U}}Y_{|U}$ est induit
par un morphisme $\tilde s\colon \mathscr A(U)\to
\tilde X\otimes_{\mathscr A(U)}\tilde Y$. Notons
$C$ un cône de $\tilde s$. Le complexe $C$ est parfait et vérifie
$\pi^\star C\simeq 0$. Il en résulte que quitte à raffiner le recouvrement
de $U$, on peut supposer que $C$ est acyclique, c'est-à-dire que $\tilde s$
est un isomorphisme.

\section{Coefficients universels}
\subsection{Énoncés pour \texorpdfstring{$\R\SheafHom$}{RHom}}

\begin{proposition2}\label{17-proposition-coefficients-universels-1}
Soit $Z$ un schéma noethérien. Soit $A$ un anneau commutatif noethérien.
Soit $X\in \detalemoinsccoefs Z A$. Soit $Y\in\detalebtfcoefs Z A$. On suppose ou
bien que $A$ est un corps ou bien que
$\R\SheafHom(X,M\overset{\L}{\otimes}_{A} Y)$ est borné
indépendamment du $A$-module $M$. Alors, pour tout $M\in\dplus(A)$, le
morphisme canonique
\[M\overset{\L}{\otimes}_A \R\SheafHom(X,Y)\to
\R\SheafHom(X,M\overset{\L}{\otimes}_A Y)\]
est un isomorphisme. En outre, si on est dans le cas où
$\R\SheafHom(X,M\overset{\L}{\otimes}_{A} Y)$ est borné
indépendamment du $A$-module $M$, alors
$\R\SheafHom(X,Y)\in\detalebtfcoefs Z A$.
\end{proposition2}

On note $u_M\colon M\overset{\L}{\otimes}_A \R\SheafHom(X,Y)\to
\R\SheafHom(X,M\overset{\L}{\otimes}_A Y)$ le morphisme
canonique, pour tout $M\in\detalepluscoefs Z A$. L'hypothèse selon laquelle
$X\in\detalemoinsccoefs Z A$ implique que les foncteurs $\SheafExt^q(X,-)$
commutent aux limites inductives filtrantes de faisceaux de $A$-modules sur
$Z$. En particulier, le foncteur $\R\SheafHom(X,-)$ de la catégorie
$\detalepluscoefs Z A$ dans elle-même commute aux sommes directes
représentables. On en déduit que $u_M$ est un isomorphisme si $M$ est un
$A$-module libre. Par suite, si $M$ est un complexe borné de $A$-modules
libres, alors $u_M$ est un isomorphisme ; si $A$ est un corps, on peut
conclure que $u_M$ est un isomorphisme pour tout $M\in \detalepluscoefs Z A$.

On se place dorénavant dans le cas où où
$\R\SheafHom(X,M\overset{\L}{\otimes}_{A} Y)$ est borné
indépendamment du $A$-module $M$.
Soit $M\in\db(A)$. Montrons que $u_M$ est un isomorphisme.
Pour tout entier relatif $q$, il existe un
morphisme $P\to M$ avec $P$ un complexe borné de $A$-modules libres tel que
si on note $C$ un cône de ce morphisme, alors $C\leq q$ (de telles
inégalités sont à comprendre relativement à la $t$-structure canonique).
On considère le carré commutatif :
\[
\xymatrix{
P\overset{\L}{\otimes}_A \R\SheafHom(X,Y)
\ar[d]\ar[r]^{u_P}_{\sim}&
\R\SheafHom(X,P\overset{\L}{\otimes}_A Y) \ar[d] \\
M\overset{\L}{\otimes}_A \R\SheafHom(X,Y)\ar[r]^{u_M} &
\R\SheafHom(X,M\overset{\L}{\otimes}_A Y)
}\]
Notons $N$ le plus petit entier naturel tel que 
$\R\SheafHom(X,V\overset{\L}{\otimes}_A Y)$ soit $\leq N$ pour
tout $A$-module $V$. Par dévissage, on obtient que les cônes des flèches
verticales du diagramme ci-dessus sont $\leq q+N$. Les deux flèches
verticales induisent donc des isomorphismes sur les objets de cohomologie
${\mathscr H}^i$ pour $i\geq q+n+2$. La flèche du haut étant un
isomorphisme, la flèche du bas induit des isomorphismes sur les 
${\mathscr H}^i$ pour $i\geq q+n+2$. Ce fait étant vérifié pour tout
$q\in\mathbf{Z}$, le morphisme du bas est bien un isomorphisme.

On déduit aussitôt de ce qui précède que
$\R\SheafHom(X,Y)\in\detalebtfcoefs Z A$. On sait par ailleurs que
$Y\in\detalebtfcoefs Z A$ et que $X\in \detalemoinscoefs Z A$. Par conséquent, il
existe un entier relatif $N$ tel que si $M\in\dtotale^{\geq c}(A)$ pour un
certain entier $c$, alors la source et le but de $u_M$ sont $\geq c+N$.
En raisonnant comme ci-dessus, on déduit du fait que $u_M$ soit un
isomorphisme pour tout $M\in\db(A)$ que ce résultat vaut en
fait pour tout $M\in\dplus(A)$.

\begin{proposition2}\label{17-proposition-coefficients-universels-2}
Soit $Z$ un schéma. Soit $A$ un anneau commutatif. Soit
$Y\in\detalebtfcoefs Z
A$. Soit $M$ un complexe pseudo-cohérent dans $\dtotale(A)$. Soit
$N\in\dplus(A)$. Alors, le morphisme canonique 
\[ \R\Hom(M,N)\overset{\L}{\otimes}_A Y\to
\R\SheafHom(M,N\overset{\L}{\otimes}_A Y) \]
est un isomorphisme dans $\detalepluscoefs Z A$.
\end{proposition2}

Pour tout $M\in\dtotale(A)$, notons $v_M$ le morphisme canonique \[
\R\Hom(M,N)\overset{\L}{\otimes}_A Y\to
\R\SheafHom(M,N\overset{\L}{\otimes}_A Y)\;\text{.} \] Bien
entendu, si $M\in\dtotale(A)$ est un complexe parfait, $v_M$ est un
isomorphisme. Pour tout $M\in\dtotale(A)$ pseudo-cohérent
et tout entier relatif $q$, il
existe un morphisme $P\to M$ avec $P$ parfait dont le cône soit $\leq q$.
Pour pouvoir déduire que $v_M$ est un isomorphisme pour tout
$M\in\dtotale(A)$ du cas particulier où $M$ est supposé parfait, il suffit
donc de montrer qu'il existe un entier relatif $c$ tel que pour tout entier
$q$ et tout $M\in\dtotale^{\leq q}(A)$, la source et le but de $v_M$
soient $\geq -q+c$. Notons $a$ un entier tel que $N\geq a$ et $b$ un entier
tel que pour tout $A$-module $V$, $V\overset{\L}{\otimes}_A Y\geq b$.
On vérifie aussitôt que $c=a+b$ convient, ce qui achève la démonstration de
la proposition.

\begin{proposition2}\label{17-proposition-coefficients-universels-3}
Soit $Z$ un schéma noethérien. Soit $A$ un anneau commutatif noethérien.
Soit $X\in \detalemoinsccoefs Z A$. Soit $Y\in\detalebtfcoefs Z A$. On suppose ou
bien que $A$ est un corps ou bien que
$\R\SheafHom(X,N\overset{\L}{\otimes}_{A} Y)$ est borné
indépendamment du $A$-module $N$. Soit $B$ une $A$-algèbre.
Alors, pour tout complexe pseudo-cohérent
$M\in\dtotale(B)$ et pour tout $N\in\dplus(B)$, le morphisme canonique
\[\R\Hom_B(M,N)\overset{\L}{\otimes}_A \R\SheafHom_A(X,Y)
\to \R\SheafHom_B(M\overset{\L}{\otimes}_A X,
N\overset{\L}{\otimes}_A Y)\]
est un isomorphisme dans $\detalecoefs Z B$.
\end{proposition2}

Les hypothèses font que $\R\Hom_B(M,N)$ appartient à $\dplus(B)$.
Si $A$ n'est pas un corps, 
on peut appliquer le résultat de la
proposition~\ref{17-proposition-coefficients-universels-2} en remplaçant
respectivement $A$, $X$, $M$ et $N$ par $B$, $M$, $N$ et
$B\overset{\L}{\otimes}_A \R\SheafHom_A(X,Y)$ (on notera que
$\R\SheafHom_A(X,Y)\in\detalebtfcoefs Z A$ d'après la
proposition~\ref{17-proposition-coefficients-universels-1}). On obtient ainsi
un isomorphisme canonique
\[\R\Hom_B(M,N)\overset{\L}{\otimes}_A
\R\SheafHom_A(X,Y)\iso
\R\SheafHom_B(M,N\overset{\L}{\otimes}_A\R\SheafHom_A(X,Y))\;\text{.}\]
Si $A$ est un corps, on établit cet isomorphisme par un argument similaire
en se ramenant au cas où $M$ est un complexe parfait de $B$-modules.
En appliquant maintenant la
proposition~\ref{17-proposition-coefficients-universels-1}, on obtient un
nouvel isomorphisme :
\[\R\SheafHom_B(M,N\overset{\L}{\otimes}_A\R\SheafHom_A(X,Y))\iso
\R\SheafHom_B(M,\R\SheafHom_A(X,N\overset{\L}{\otimes}_A
Y))\;\text{.}\]
Enfin, on utilise l'isomorphisme cher à Cartan :
\[\R\SheafHom_B(M,\R\SheafHom_A(X,N\overset{\L}{\otimes}_A
Y))\simeq
\R\SheafHom_B(M\overset{\L}{\otimes}_A X,
N\overset{\L}{\otimes}_A Y)\;\text{.}\]
L'isomorphisme voulu est l'isomorphisme obtenu en composant les différents
isomorphismes canoniques ci-dessus.

\subsection{Conséquences pour \texorpdfstring{$\R j_\star$}{Rj*} et \texorpdfstring{$i^!$}{i!}}

\begin{proposition2}\label{17-proposition-coefficients-universels-r-j-star}
Soit $j\colon U\to X$ une immersion ouverte entre schémas noethériens. Soit
$A$ un anneau commutatif noethérien. Soit $Y\in\detalebtfcoefs U A$. On suppose
ou bien que $A$ est un corps ou bien que $\R j_\star
(M\overset{\L}\otimes_A Y)$ est borné indépendamment du $A$-module
$M$. Alors, pour tout $M\in\dplus(A)$, le morphisme canonique
\[M\overset{\L}{\otimes}_A\R j_\star Y\to
\R j_\star (M\overset{\L}{\otimes}_A Y)\]
est un isomorphisme dans $\detalecoefs X A$. Si $\R j_\star
(M\overset{\L}\otimes_A Y)$ est borné indépendamment du $A$-module
$M$, alors $\R j_\star Y$ appartient à $\detalebtfcoefs X A$.
\end{proposition2}

Ceci résulte de la
proposition~\ref{17-proposition-coefficients-universels-1}, compte tenu de la
formule $\R j_\star Y\simeq \R\SheafHom(j_!A,j_!Y)$ pour tout
$Y\in\detalecoefs U A$.

\begin{proposition2}
Soit $i\colon Z\to X$ une immersion fermée entre schémas noethériens. Soit
$A$ un anneau commutatif noethérien. Soit $Y\in\detalebtfcoefs Z A$. On suppose
ou bien que $A$ est un corps ou bien que $i^!
(M\overset{\L}\otimes_A Y)$ est borné indépendamment du $A$-module
$M$. Alors, pour tout $M\in\dplus(A)$, le morphisme canonique
\[M\overset{\L}{\otimes}_A i^! Y\to i^!
(M\overset{\L}{\otimes}_A Y)\]
est un isomorphisme dans $\detalecoefs X A$. Si $i^!
(M\overset{\L}\otimes_A Y)$ est borné indépendamment du $A$-module
$M$, alors $i^! Y$ appartient à $\detalebtfcoefs Z A$.
\end{proposition2}

Cette fois-ci, on utilise la formule $i^!Y\simeq \R\SheafHom(i_\star
A,i_\star Y)$.

\section{Modules ind-unipotents}
\subsection{Définitions}
Soit $G$ un groupe topologique. Soit $A$ un anneau. On appellera ici 
$A[G]$-module un $A[G]$-module (à gauche)
au sens où on l'entend usuellement en
considérant $G$ comme groupe discret. Un $A[G]$-module discret est un
$A[G]$-module dans lequel le stabilisateur de tout élément est ouvert.
Un $A[G]$-module est dit trivial si le groupe $G$ agit
trivialement sur lui ; un tel $A[G]$-module est discret.
On note $I_G$ l'idéal d'augmentation de $A[G]$ ; un
$A[G]$-module est trivial si et seulement s'il est annulé par $I_G$.

La catégorie des $A[G]$-modules discrets est une catégorie abélienne de
Grothendieck ; le foncteur d'inclusion de cette
catégorie dans celle des $A[G]$-modules est exact, commute aux limites
inductives, mais en général pas aux limites projectives. La catégorie des
$A[G]$-modules discrets est stable par sous-quotient (mais en général pas
par extensions) dans celle des $A[G]$-modules.

\begin{definition2}
Un $A[G]$-module discret est unipotent s'il admet une filtration finie dont
les quotients successifs sont des $A[G]$-modules triviaux. L'ordre
d'unipotence d'un $A[G]$-module unipotent est la plus petite longueur d'une
telle filtration ; ainsi, un $A[G]$-module trivial non nul est unipotent
d'ordre $1$. Un $A[G]$-module discret est ind-unipotent si tous ses
sous-$A[G]$-modules de type fini sont unipotents.
\end{definition2}

\begin{definition2}\label{17-definition-filtration-unipotence}
Soit $M$ un $A[G]$-module discret. On définit par récurrence
une filtration croissante
$(\Fil_n M)_{n\in\mathbf{Z}}$ de $M$ par des sous $A[G]$-modules, de façon
à ce que $\Fil_0M=0$ et que pour tout entier $n\in\mathbf{N}$, on ait un
isomorphisme canonique $\Fil_{n+1}M/\Fil_nM\simeq \HH^0(G,M/\Fil_nM)\subset
M/\Fil_n M$.
\end{definition2}

La filtration $\Fil_\bullet M$ est évidemment fonctorielle en $M$ au sens
où si $f\colon M\to M'$ est un morphisme de $A[G]$-modules,
$f(\Fil_nM)\subset \Fil_n(M')$ pour tout $n\in\mathbf{Z}$.

\subsection{Propriétés}

\begin{proposition2}\label{17-proposition-proprietes-unipotents}
Les notions d'unipotence et d'ind-unipotence des $A[G]$-modules discrets
jouissent des propriétés suivantes :
\begin{itemize}
\item[(i)] Le caractère unipotent des
$A[G]$-modules discrets est stable par sous-quotients
et extensions ; l'indice d'unipotence décroît par passage à un
sous-quotient et est sous-additif vis-à-vis des extensions ;
\item[(ii)] Le caractère ind-unipotent des $A[G]$-modules discrets est
stable par sous-quotients, et, si $G$ est profini, par extensions ;
\item[(iii)] La catégorie des $A[G]$-modules discrets unipotents
(resp. ind-unipotents) est abélienne, le foncteur d'inclusion
dans celle des $A[G]$-modules discrets est exact ;
\item[(iv)] Un $A[G]$-module discret $M$ est
unipotent si et seulement s'il existe un entier naturel $n$ tel que
$\Fil_nM=M$ ; dans ce cas, l'indice d'unipotence de $M$ est le plus petit
de ces entiers $n$ ;
\item[(v)] Le caractère unipotent (resp. ind-unipotent)
d'un $A[G]$-module discret ne dépend pas de l'anneau
des coefficients $A$ et cette notion n'est pas altérée non plus
si on considère $G$ comme groupe discret ;
\item[(vi)] Pour tout entier naturel $n$ et tout $A[G]$-module
discret, $\Fil_n M$ est l'annulateur de $I_G^n$ dans $M$ ;
\item[(vii)] Pour tout entier naturel $n$,
un $A[G]$-module discret est unipotent d'ordre au plus $n$
si et seulement s'il est annulé par l'idéal $I_G^n$. Un $A[G]$-module
discret $M$ est ind-unipotent si et seulement si tout élément de $M$ est
annulé par une puissance de $I_G$, c'est-à-dire que
$M=\cup_{n\in\mathbf{N}}\Fil_n M$ ;
\item[(viii)] Un $A[G]$-module unipotent est ind-unipotent ;
\item[(ix)] Un $A[G]$-module de type fini est unipotent si et seulement
s'il est ind-unipotent ;
\item[(x)] Une somme directe (resp. une limite inductive)
de $A[G]$-modules ind-unipotents est un $A[G]$-module ind-unipotent.
\item[(xi)] Les $A[G]$-modules (discrets) ind-unipotents sont exactement
les $A[G]$-modules limites inductives filtrantes de $A[G]$-modules
(discrets) unipotents.
\end{itemize}
\end{proposition2}

Montrons (i).
Si $(M,F_\bullet M)$ est un $A[G]$-module discret filtré, tout sous-objet
$M'$ (resp. sous-quotient $M''$) de $M$ est naturellement muni d'une
filtration $F_\bullet M'$ (resp. $F_\bullet M''$) telle que pour tout
$n\in\mathbf{Z}$, le morphisme induit au niveau des gradués $\Gr_n M'\to
\Gr_n M$ (resp. $\Gr_n M\to \Gr_n M''$) soit injectif (resp. surjectif) ;
comme tout sous-quotient d'un $A[G]$-module trivial est trivial, il en
résulte aussitôt que le caractère unipotent d'un $A[G]$-module est stable
par sous-quotient et que l'indice d'unipotence décroît dans cette
opération. Il est évident sur la définition que le caractère unipotent des
$A[G]$-modules discrets est stable par extensions, et que l'indice est
sous-additif vis-à-vis d'elles. (viii) est une conséquence triviale de (i).

Montrons la propriété (iv). Soit $M$ un $A[G]$-module discret.
Si $\Fil_n M=M$, la filtration $(\Fil_iM)_{0\leq i\leq n}$ est une
filtration finie de $M$ de longueur $n$ dont les quotients successifs sont
triviaux ; ainsi, $M$ est unipotent d'ordre au plus $n$. Inversement,
soit $M$ un $A[G]$-module muni d'une filtration $(F_i M)_{i\geq 0}$ telle
que $F_0M=0$ et dont les quotients successifs soient triviaux. Une
récurrence évidente sur $i\in\mathbf{N}$ montre que l'on a une inclusion
$F_i M\subset \Fil_i M$, de sorte que si un entier naturel $n$ est tel que
$F_n M=M$, alors $\Fil_n M=M$.

Dans le cas unipotent, la propriété (v) est une conséquence de la propriété
(iv) : la filtration $\Fil_\bullet M$ est la même que l'on considère $M$
comme $A[G]$-module ou comme $\mathbf{Z}[G]$-module, et que l'on considère
$G$ comme un authentique groupe topologique ou comme un groupe discret.
Quand ce ne sera pas pertinent, on ne précisera donc pas systématiquement
dans la suite que les $A[G]$-modules sont discrets.

Concernant la propriété (vi), il est facile de montrer par récurrence sur
$n\in\mathbf{N}$ que $\Fil_n M$ est l'annulateur de $I_G^n$ dans $M$.

Considérons (vii). On déduit trivialement de (iv) et (vi) qu'un
$A[G]$-module discret est unipotent d'ordre au plus $n$ si et seulement
s'il est annulé par $I_G^n$. Montrons l'autre partie de (vii). Soit $M$ un
$A[G]$-module discret ind-unipotent. Pour tout $m\in M$, le
sous-$A[G]$-module de $M$ engendré par $m$ est unipotent, d'après ce qu'on
vient de montrer, ce module est annulé par une puissance de $I_G$, en
particulier, $m$ est annulé par $I_G^n$ pour un certain entier naturel $n$.
Inversement, supposons que tout élément de $M$ soit annulé par un puissance
de $I_G$. Si on considère un sous-$A[G]$-module de type fini $N$ de $M$, en
appliquant l'hypothèse aux éléments d'un ensemble fini de générateurs de
$N$, on obtient que $N$ est annulé par $I_G^n$ pour un certain entier
naturel $n$, et donc que $N$ est unipotent. Ceci achève la démonstration de
(vii). On sait donc qu'un $A[G]$-module discret $M$ est ind-unipotent si et
seulement si $M=\cup_{n\in\mathbf{N}}\Fil_n M$. Comme la filtration $\Fil_n
M$ ne dépend pas de la topologie de $G$ ni de l'anneau des coefficients
$A$, on peut obtenir (v). La propriété (ix) est une conséquence immédiate
de (vii). La stabilité par sommes directes et sous-quotients du caractère
ind-unipotent en résulte aussi, ce qui établit (x).

Supposons $G$ profini
et montrons la stabilité par extensions des $A[G]$-modules discrets
ind-unipotents, ce qui achèvera la démonstration de (ii).
D'après (v), on peut supposer que $A=\mathbf{Z}$. On se
donne une suite exacte courte $0\to M'\to M\to M''\to 0$ de
$\mathbf{Z}[G]$-modules discrets avec $M'$ et $M''$ ind-unipotents.
Soit $N$ un sous-$\mathbf{Z}[G]$-module
de type fini de $M$. Il s'agit de montrer que $N$ est unipotent. On peut
considérer l'image $N''$ de $N$ dans $M''$, et $N'$ le noyau de la
projection $N\to N''$. Le $\mathbf{Z}[G]$-module $N'$ (resp. $N''$) est un
sous-$\mathbf{Z}[G]$-module de $M'$ (resp. $M''$). Par passage à des
sous-objets, $N'$ et $N''$ sont ind-unipotents. Le groupe $G$ étant profini
et $N$ un $\mathbf{Z}[G]$-module discret, $N$ est un groupe abélien de type
fini ; les groupes abéliens $N'$ et $N''$ qui en sont des sous-quotients
sont eux aussi de type fini. D'après (ix), $N'$ et $N''$ sont
unipotents ; d'après (i), $N$ est unipotent, ce qui achève la
démonstration de (ii).

(iii) résulte aussitôt de (i) et (ii). Enfin, (xi) découle de (vii), (viii)
et (x).

\begin{remarque2}
Si $H\to G$ est un morphisme de groupes topologiques et $M$ un
$A[G]$-module unipotent (resp. ind-unipotent), alors, en tant
que $A[H]$-module, $M$ est unipotent (resp. ind-unipotent).
\end{remarque2}

\subsection{Modules ind-unipotents pour un sous-groupe distingué}
\label{17-subsection-modules-ind-unipotents-sous-groupe-distingue}

Soit $A$ un anneau. Soit $G$ un groupe topologique. Soit $H$ un sous-groupe
distingué (fermé) de $G$. Si $M$ est un $A[G]$-module discret, on peut se
demander si en tant que $A[H]$-module, $M$ est unipotent (resp.
ind-unipotent).

Dans cette situation, on note $\Fil_\bullet M$ la filtration de la
définition~\ref{17-definition-filtration-unipotence} pour le groupe
$H$. Le fait que $H$ soit distingué dans $G$ permet d'obtenir aussitôt que
cette filtration $\Fil_\bullet M$ est constituée de sous-$A[G]$-modules de
$M$. De cette remarque et de la 
proposition~\ref{17-proposition-proprietes-unipotents}~(iv),
on tire :

\begin{proposition2}
Soit $M$ un $A[G]$-module. Les conditions suivantes sont équivalentes :
\begin{itemize}
\item en tant que $A[H]$-module, $M$ est unipotent ;
\item il existe une filtration finie de $M$ par des sous-$A[G]$-modules
telle que $H$ agisse trivialement sur les quotients successifs.
\end{itemize}
\end{proposition2}

\begin{proposition2}
Soit $M$ un $A[G]$-module. Les conditions suivantes sont équivalentes :
\begin{itemize}
\item en tant que $A[H]$-module, $M$ est ind-unipotent ;
\item tout sous-$A[G]$-module de type fini de $M$ est unipotent pour
$H$.
\end{itemize}
\end{proposition2}

Il s'agit de montrer que si un sous-$A[H]$-module $N$ de $M$ de type fini
est unipotent, alors le sous-$A[G]$-module (de type fini) de $M$ engendré
par $N$ est unipotent. Il suffit en fait d'établir le lemme suivant :

\begin{lemme2}
Soit $M$ un $A[H]$-module unipotent.
En tant que $A[H]$-module, $A[G]\otimes_{A[H]} M$ est unipotent.
\end{lemme2}

On peut conclure facilement en observant simplement que le $A[H]$-module
$A[G]\otimes_{A[H]} M$ s'identifie à une somme directe de copies de $M$.

\bigskip

Maintenant que les définitions d'unipotence et de ind-unipotences pour le
sous-groupe distingué $H$ sont clarifiées, on peut énoncer la proposition
suivante :

\begin{proposition2}\label{17-proposition-proprietes-unipotents-sous-groupe-distingue}
Les énoncés de la proposition~\ref{17-proposition-proprietes-unipotents}
restent vrais si on
remplace \guil{unipotent} par \guil{unipotent pour $H$},
\guil{ind-unipotent} par \guil{ind-unipotent pour $H$} et $I_G$ par $I_H$,
et que, comme ci-dessus, on définit la filtration $\Fil_\bullet$
relativement au groupe $H$.
\end{proposition2}

Compte tenu des clarifications faites ci-dessus, c'est trivial.


\ifx\danslelivre\undefined
\bibliography{../fichiers-configuration-latex/bibliographie}
\bibliographystyle{../fichiers-configuration-latex/smfalpha}


\begin{center}
XVIII. Cohomological dimension: First results\\
Luc Illusie
\end{center}

\tableofcontents

\else
\setboolean{english}{true}
\setboolean{francais}{false}
\selectlanguage{english}
\chapter{Cohomological dimension: First results}
\addtocontents{toc}{par Luc Illusie}
\begin{center}
Luc Illusie
\end{center}
\fi

In this chapter we establish Gabber's bound on cohomological dimension stated in the introduction
(in the comments on the proof of the finiteness theorem).

\section{Statement of the main result and applications}

\begin{theorem}[Gabber]\label{18A-GABBER'S BOUND}
Let $X$ be a strictly local, noetherian scheme of dimension $d >0$, and let $\ell$ be a prime
number invertible on $X$. Then, for any open subset $U$ of $X$, we have
\begin{equation}\label{18A-Gabber's bound1} \mathrm{cd}_{\ell}(U) \le 2d-1.  \end{equation}
\end{theorem}
Recall that, for a scheme $S$, $\mathrm{cd}_{\ell}(S)$ ($\ell$-\textit{cohomological
dimension} of $S$) denotes the infinimum of the integers $n$ such that for all $\ell$-torsion abelian
sheaves $F$ on $S$, and all $i > n$, $H^i(S,F) = 0$.  \begin{corollary}\label{18A-fraction field}
Let $X = \SP \,A$ be as in \ref{18A-GABBER'S BOUND}, and assume $A$ is a domain, with fraction field
$K$. Then \begin{equation}\label{18A-fraction field1} \mathrm{cd}_{\ell}(K) \le 2d-1.
\end{equation} \end{corollary} Indeed, it suffices to show that if $F$ is a finitely generated
$\mathbf{F}_{\ell}$-module over $\eta = \SP\, K$, then $H^i(\eta,F) = 0$ for $i > 2d-1$. But $\eta$
is a filtering projective limit of affine open subsets $U_{\alpha}$ of $X$, $F$ is induced from a
locally constant constructible $\mathbf{F}_{\ell}$-sheaf $F_{\alpha_0}$ on $U_{\alpha_0}$, and
$H^i(\eta,F) = \varinjlim H^i(U_{\alpha},F_{\alpha})$, where $F_{\alpha} = F_{\alpha_0} |
U_{\alpha}$ for $\alpha \ge \alpha_0$ (\sga{4}{vii}{5.7}).

\begin{remark}\label{18A-fraction field
sharp} (a) The proof shows that if (\ref{18A-Gabber's bound1}) holds for one $X$, and if $X$ is
integral, then (\ref{18A-fraction field1}) holds for $X$.
 
 (b) If $X$ is an integral noetherian scheme of dimension $d$, with generic point $\SP \, K$, and
$\ell$ is a prime number invertible on $X$, then $\mathrm{cd}_{\ell}(K) \ge d$ (\sga{4}{x}{2.5}).
Gabber can prove that under the assumptions of \ref{18A-fraction field} one has
$\mathrm{cd}_{\ell}(K) = d$.
\end{remark}

\begin{corollary}\label{18A-cd(f)} Let $Y$ be a
noetherian scheme of finite dimension, $f : X \ra Y$ a morphism of finite type, and $\ell$ a prime
number invertible on $Y$. Then $$ \mathrm{cd}_{\ell}(Rf_*) < \infty, $$ i. e. there exists an
integer $N$ such that for any $\ell$-torsion abelian sheaf $F$ on $X$, $R^qf_*F = 0$ for $q > N$.
\end{corollary}

\textit{Proof of \ref{18A-cd(f)}}. We may assume $Y$ affine. Covering $X$ by finitely many open
affine subsets $U_i$ ($0 \le i \le  n)$, and using the alternate \v Cech spectral sequence $$
E^{pq}_1 = \oplus R^q(f|U_{i_0 \cdots i_p})_*(F|U_{i_0 \cdots i_p}) \Rightarrow R^{p+q}f_*F, $$
where $U_{i_0 \cdots i_p} = U_{i_0} \cap \cdots \cap U_{i_p}$, we may assume $f$ separated.
Repeating the procedure, we may assume $X$ affine. Choose an immersion $X \ra \mathbf{P}^n_Y$, and
replace $\mathbf{P}^n_Y$ by the closure of the image of $X$. We get a commutative diagram $$
\xymatrix{X \ar[d]^f \ar[r]^j & \overline{X} \ar[dl]^g \\ Y} $$ with $j$ open and $g$ projective of
relative dimension $\le n$. By the proper base change theorem we have $\mathrm{cd}_{\ell}(Rg_*) \le
2n$. By the Leray spectral sequence $R^pg_*R^qj_*F \Rightarrow R^{p+q}f_*F$ it thus suffices to
prove \ref{18A-cd(f)} for $j$, in other words, we may assume that $f$ is an open immersion. Let $d$
be the dimension of $Y$. Let $x$ be a geometric point of $X$, with image $y = f(x)$ in $Y$, and let
$U = Y_{(y)} \times_Y X$ be the corresponding open subset of the strictly local scheme $Y_{(y)}$ (of
dimension $\le d$, that we may assume to be $>0$). Then $$ R^if_*(F)_y = H^i(U,F) $$ (where we still
denote by $F$ its inverse image on $U$). The conclusion follows from \ref{18A-GABBER'S BOUND}.
\begin{remark} \label{18A-remark-d-2d-1} (a) Under the assumptions of \ref{18A-GABBER'S BOUND}, if
$X$ is quasi-excellent and $U$ is affine, then by Gabber's affine Lefschetz theorem
(\refext{XV}{15-prop-devissage-classique}) we have $\mathrm{cd}_{\ell}(U) \le d$. 

(b) Gabber can show that, under the assumptions of \ref{18A-GABBER'S BOUND}, one has
$\mathrm{cd}_{\ell}(U) \ge d$ if $U$ is not empty and does not contain the closed point and that for
each $n$ such that $d \le n \le 2d-1$, there exists a pair $(X,U)$ as in \ref{18A-GABBER'S BOUND},
with $U$ affine, such that $\mathrm{cd}_{\ell}(U) = n$ (by (a), for $n > d$, $X$ is not
quasi-excellent). These results, as well as the one mentioned in \ref{18A-fraction field sharp}, are
proved in the second part of this chapter, \emph{raffinements et compléments}, by Fabrice Orgogozo.
\end{remark}

\section{Proof of \ref{18A-GABBER'S BOUND}} \begin{lemma}\label{18A-punctured spectrum} Let $X$ be
as in \ref{18A-GABBER'S BOUND}, and let $\{x\}$ be the closed point of $X$. Then (\ref{18A-Gabber's
bound1}) holds for $U = X - \{x\}$.  \end{lemma} \textit{Proof.} It suffices to show that for any
constructible $\mathbf{F}_{\ell}$-sheaf $F$ on $U$, $H^i(U,F) = 0$ for $i \ge 2d$. Let $\widehat X$
be the completion of $X$ at $\{x\}$ and set $\widehat U := \widehat X \times_X U = \widehat X -
\{x\}$. Let $\widehat F$ be the inverse image of $F$ on $\widehat U$. By Gabber's formal base change
theorem (\cite[6.6.4]{TUBULAR@Fujiwara}), the natural map $$ H^i(U,F) \ra H^i(\widehat U, \widehat
F) $$ is an isomorphism for all $i$. Therefore we may assume $X$ complete, and in particular,
excellent. Let $(f_1,\cdots,f_d)$ be a system of parameters of $X$, and let $U_i = X_{f_i}$, so that
$U = \cup_{1 \le i \le d}U_i$. Consider the (alternate) \v Cech spectral sequence $$ E^{pq}_1 =
\oplus H^q(U_{i_0 \cdots i_p},F) \Rightarrow H^{p+q}(U,F), $$ with $U_{i_0 \cdots i_p} = U_{i_0}
\cap \cdots \cap U_{i_p}$ as above.  By definition, $E^{pq}_1 = 0$ for $p \ge d$. On the other hand,
as $X$ is excellent, by Gabber's affine Lefschetz theorem
(\refext{XV}{15-prop-devissage-classique})), we have $E^{pq}_1 = 0$ for $q \ge d+1$. Therefore
$E^{pq}_1 = 0$ for $p+q \ge 2d$, hence $H^i(U,F) = 0$ for $i \ge 2d$.
\begin{lemma}\label{18A-punctual bound} Let $X$ be a noetherian scheme of finite dimension, $Y$ a
closed subset, $\ell$ a prime number invertible on $X$. Then, for any $\ell$-torsion sheaf $F$ on
$X$, $$ H^i_Y(X,F) = 0 $$ for $$ i > \mathrm{sup}_{x \in Y} (\mathrm{cd}_{\ell}(k(x)) + 2
\mathrm{dim} \OO_{X,x}).  $$ In particular, $$ \mathrm{cd}_{\ell}(X) \le \mathrm{sup}_{x \in
X} (\mathrm{cd}_{\ell}(k(x)) + 2 \mathrm{dim} \OO_{X,x}).  $$ \end{lemma} \textit{Proof.}
For $p \ge 0$, let $\Phi^p$ be the set of closed subsets of $Y$ of codimension $\ge p$ in $X$. We
have $\Phi^p = \emptyset$ for $p > \mathrm{codim}(Y)$.  Consider the (biregular) coniveau spectral
sequence of the filtration $(\Phi^p)$ cf. \cite[10.1]{BrauerIII@Grothendieck}),
\begin{equation}\label{18A-coniveau} E^{pq}_1 = H^{p+q}_{\Phi^p/\Phi^{p+1}}(X,F) \Rightarrow
H^{p+q}_Y(X,F).  \end{equation} We have $$ E^{pq}_1 = \oplus_{x \in Y^{(p)}}
H^{p+q}_{\{x\}}(X_x,F|X_x), $$ where $Y^{(p)}$ denotes the set of points of $Y$ of codimension $p$
in $X$, and $X_x = \SP \, \OO_{X,x}$.  For $x \in Y^{(p)}$ (i. e.
$\mathrm{dim}\OO_{X,x} = p$), let $\overline{x}$ be a geometric point above $x$. Consider
the diagram $$ \xymatrix{\{\overline{x}\} \ar[d] \ar[r]^{i_{\overline{x}}} & X_{(\overline{x})}
\ar[d] & \overline{U} \ar[l]_{\overline{j}} \ar[d] \\ \{x\} \ar[r]^{i_x} & X_x & U \ar[l]_{j}}, $$
where $U = X - \{x\}$, $\overline{U} = X_{(\overline{x})}- \{\overline{x}\}$. We have $$
R\Gamma_{\{x\}}(X_x,F|X_x) = R\Gamma(\{x\},Ri_x^!(F|X_x)).  $$ The stalk of $Ri_x^!(F|X_x)$ at
$\overline{x}$ is $$ Ri_x^!(F|X_x)_{\overline{x}} = Ri^!_{\overline{x}}(F| X_{(\overline{x})}), $$
as $(Rj_*(F|U))_{\overline{x}} = R{\overline{j}}_*(F|\overline{U})_{\overline{x}}$.  We thus have a
spectral sequence \begin{equation}\label{18A-Galois} E^{rs}_2 = H^r(k(x),R^si_{\overline{x}}^!(F|
X_{(\overline{x})})) \Rightarrow H^{r+s}_{\{x\}}(X_x,F|X_x) , \end{equation} It suffices to show
that, in the initial term of (\ref{18A-coniveau}), \begin{equation}\label{18A-E1}
H^{p+q}_{\{x\}}(X_x,F|X_x) = 0 \end{equation} for $p+q > \mathrm{cd}_{\ell}(k(x)) + 2p$. If $p = 0$,
then $Ri^!_{\overline{x}}(F| X_{(\overline{x})}) = F_{\overline{x}}$, and, in (\ref{18A-Galois}),
$E^{rs}_2 = 0$ for $s >0$, $E^{r0}_2 = 0$ for $r > \mathrm{cd}_{\ell}(k(x))$, so (\ref{18A-E1}) is
true in this case. Assume $p >0$. We have \begin{equation}\label{18A-local cohomology}
R^si_{\overline{x}}^!F = H^{s-1}(\overline{U},F|\overline{U}) \end{equation} for $s \ge 2$, where,
as above, $\overline{U} = X_{(\overline{x})} - \{\overline{x}\}$. By \ref{18A-punctured spectrum},
$H^{s-1}(\overline{U},F| \overline{U}) = 0$ for $s-1 \ge 2p$, hence, by (\ref{18A-local
cohomology}), $R^si_{\overline{x}}^!(F|X_{(\overline{x})}) = 0$ and $E^{rs}_2 = 0$ for $s \ge 2p
+1$. If $r+s \ge \mathrm{cd}_{\ell}(k(x)) + 2p + 1$ and $s \le 2p$, then $r >
\mathrm{cd}_{\ell}(k(x))$, hence $E^{rs}_2 = 0$ as well. Therefore, by (\ref{18A-Galois}),
(\ref{18A-E1}) holds, which finishes the proof. 

\medskip \textit{Proof of \ref{18A-GABBER'S BOUND}}. We prove \ref{18A-GABBER'S BOUND} by induction
on $d$. For $n \ge 0$ consider the assertion

\medskip $G_n$ : \textit{For every strictly local, noetherian scheme $X$ of dimension $n$, all open
subsets $U$ of $X$ and any prime number $\ell$ invertible on $X$, we have $\mathrm{cd}_{\ell}(U) \le
\mathrm{sup}(0,2n-1)$}. 

\medskip Let $d >0$. Assume $G_n$ holds for $n < d$, and let us prove $G_d$. Let $X$ be as in
\ref{18A-GABBER'S BOUND}. If $(X_i)_{1\le i \le r}$ are the reduced irreducible components of $X$
and $U_i = U \times_X X_i$, we have $\mathrm{cd}_{\ell}(U) \le \mathrm{sup}(
\mathrm{cd}_{\ell}(U_i))$, hence we may assume $X$ integral. Let $x$ be the closed point of $X$, and
$U = X - \{x\}$ the punctured spectrum. Let $j : V \ra U$ be a nonempty open subset of $U$, and $F$
be a constructible $\mathbf{F}_{\ell}$-sheaf on $V$. As $F = j^*j_!F$, by \ref{18A-punctured
spectrum} it suffices to show that, for any constructible $\mathbf{F}_{\ell}$-sheaf $L$ on $U$, the
restriction map $$ H^i(U,L) \ra H^i(V,j^*L) \leqno (*) $$ is an isomorphism for $i \ge 2d$. Let $Y =
U - V$. Consider the exact sequence $$ H^i_Y(U,L) \ra H^i(U,L) \ra H^i(V,j^*L) \ra H^{i+1}_Y(U,L).
$$ By \ref{18A-punctual bound}, we have $H^i_Y(U,L) = 0$ for $i > \mathrm{sup}_{y \in Y}
(\mathrm{cd}_{\ell}(k(y)) + 2 \mathrm{dim} \OO_{X,y})$. For $y \in Y$, denote by $Z$ the
closed, integral subscheme of $X$ defined by the closure of $\{y\}$ \textit{in} $X$. As $X$ is
integral and $V$ nonempty, $Z$ is a strictly local scheme of dimension $n < d$, with generic point
$y$.    By \ref{18A-fraction field sharp} (a) and $G_n$ (inductive assumption), we have
$\mathrm{cd}_{\ell}(k(y)) \le 2n -1$. We have $2n - 1 + 2\mathrm{dim}\OO_{X,y} \le 2d-1$.
Hence, for $i \ge 2d$, $H^i_Y(U,L) = H^{i+1}_Y(U,L) = 0$, and (*) is an isomorphism, which finishes
the proof.  \begin{remark} Gabber has an alternate proof of \ref{18A-GABBER'S BOUND}, based on the
theory of Zariski-Riemann spaces. By \ref{18A-punctual bound}, it suffices to show \ref{18A-fraction
field}. Here is a sketch, pasted from an e-mail of Gabber to Illusie of 2007, Aug. 15 :

\medskip "For $Y \ra X$ proper birational with special fiber $Y_0$, consider $i: Y_0 \ra Y$ and $j
: \eta \ra Y$, $\eta$ the generic point. We have by proper base change a spectral sequence $$
H^p(Y_0,i^*R^q j_*F) \rightarrow H^{p+q}(\eta,F) $$ for $F$ an $\ell$-torsion Galois module.  We
take the direct limit and get a spectral sequence involving cohomologies on the étale topos of
$ZRS_0$ defined as the limit of étale topoi of $Y_0$ or viewing $ZRS_0$ as a locally ringed topos
and applying a universal construction in the book of M.Hakim. The limit of the $R^q j_* F$ is $R^q
(\eta \ra ZRS)_* F$.  By a classical result of Abhyankar, also proved in Appendix 2 of the book of
Zariski-Samuel Vol. II, if $R$ is a noetherian local domain of dimension $d$ and $V$ a valuation
ring of $\mathrm{Frac}(R)$ dominating $R$, the sum of the rational rank and the residue
transcendence degree is at most $d$.  For a strictly henselian valuation ring  $V$ with residue
characteristic exponent $p$ and value group $\Gamma$, the absolute Galois group of
$\mathrm{Frac}(V)$ is an extension of the tame part (product for $\ell$ prime not equal to $p$ of
$\mathrm{Hom}(\Gamma, \mathbf{Z}_{\ell}(1)))$ by a $p$-group, so the $\ell$-cohomological dimension
is the dimension of $\Gamma$ tensored with the prime field $\mathbf{F}_{\ell}$, which is at most the
dimension of $\Gamma$ tensored with the rationals.  If $A$ is an $\ell$-torsion sheaf on the étale
topos of $ZRS_0$, let $\delta(A)$ be the sup of transcendence degrees of points where the stalk is
non-zero. I claim that $H^n(ZRS_0,A)$ vanishes for $n> 2 \delta(A)$. One reduces it to the finite
type case (passage to the limit \sga{4}{vi}{8.7.4})  using that the $\delta$ of the direct image of
$A$ to $Y_0$ is at most $\delta(A)$.  In $Y_0$ the transcendence degrees over the closed point of
$X$ are at most $d-1$ by the dimension inequality. Summing up, for the limit spectral sequence the
$q$-th  direct image sheaf  restricted to the special fiber has $\delta$ at most
$\mathrm{min}(d-1,d-q)$, giving vanishing for certain $E_2^{p,q}$ and the result." \end{remark}

\ifx\danslelivre\undefined
\bibliography{../fichiers-configuration-latex/bibliographie}
\bibliographystyle{../fichiers-configuration-latex/smfalpha}


\begin{center}
XVII-B. Dimension cohomologique : raffinements et compléments\\
Fabrice Orgogozo
\end{center}

\version
\tableofcontents

\else
\setboolean{english}{false}
\setboolean{francais}{true}
\selectlanguage{francais}
\chapter{Dimension cohomologique : raffinements et compléments}
\addtocontents{toc}{par Fabrice Orgogozo}
\begin{center}
Fabrice Orgogozo
\end{center}
\fi

\newcommand{\initial}{\mathop{\mathrm{in}}}

Le premier objectif de cette seconde partie de l'exposé est de
démontrer que pour tout nombre premier $ℓ$ il existe un ouvert
\emph{affine} d'un schéma nœthérien strictement hensélien régulier de
dimension $2$ dont la $ℓ$-dimension cohomologique est égale à $3$.
Outre les ingrédients cohomologiques — pureté, morphisme de Gysin et
comparaison à la complétion —, on utilise une construction dont le
principe est dû à Nagata : utilisant des « dilatations formelles »
construit un schéma nœthérien strictement hensélien $X$ de
dimension $2$, de complété $\chap{X}$ \emph{régulier de
dimension $2$}, et une courbe \emph{irréductible} $C$ dans $X$
devenant le $ℓ$-ième multiple d'un diviseur régulier dans $\chap{X}$.
Cette construction est ensuite étendue au cas, plus délicat, de la
dimension supérieure. À partir de là, on construit
aisément des schémas dont l'existence a été annoncée
dans la première partie (\refext{XVIII-A}{18A-remark-d-2d-1}).
Pour vérifier que leur dimension cohomologique est bien celle
attendue, on fait appel à une majoration
assez générale établie sans hypothèse d'excellence.
Enfin, on termine par une minoration de la dimension
cohomologique d'un ouvert (non nécessairement
affine) du spectre épointé d'un schéma nœthérien intègre
strictement local.

\section{Préliminaires}

\subsection{Dilatations formelles}
\label{18B-dilatations-formelles}

Soient $R$ un anneau, $π$ un élément non diviseur de $0$
et $f$ un élément de la complétion $π$-adique $\chap{R}$
de $R$. Pour tout $n ≥ 0$, choisissons un $f_n ∈ R$ tel que
$f ≡ f_n$ modulo $π^n$.

\begin{définition2}
On note $𝒟𝒾𝓁^f_π R$ la sous-$R$-algèbre de $R[π^{-1},F]$, où $F$ est une
indéterminée, colimite des $R$-algèbres $R[\frac{F-f_n}{π^n}]$.
\end{définition2}

On notera également $F$ l'image de cette variable dans
$𝒟𝒾𝓁^f_π R$.

\begin{remarque2}
Notons que les $R$-algèbres considérées sont toutes isomorphes
à une algèbre de polynômes en une variable sur $R$.
\end{remarque2}

\subsubsection{}On vérifie immédiatement les faits suivants :

\begin{enumerate}
\item la construction ne dépend pas des choix des $f_n$,
et ne dépend de l'élément $π$ qu'à travers l'idéal qu'il
engendre ;
\item les morphismes $R →  𝒟𝒾𝓁_π^f R$
et $𝒟𝒾𝓁_π^f R → \chap{R}$, $\frac{F-f_n}{π^n}↦ \frac{f-f_n}{π^n}$,
induisent des \emph{isomorphismes} sur
la complétion $π$-adique.
\end{enumerate}

\subsection{Platitude et nœthérianité}

\subsubsection{}Rappelons que si un morphisme $A → B$ est fidèlement plat,
$A$ est nœthérien si $B$ l'est (\ega{IV₂}{2.2.14}).
Pour vérifier la platitude, il est parfois commode d'utiliser
le critère suivant (\cite[II.1.4.2.1]{Platification@Raynaud-Gruson}).

\begin{proposition2}
\label{18B-critère-platitude}
Soient $M$ un $R$-module, et $π ∈ R$. On suppose que $π$ n'est diviseur de zéro
ni dans $R$ ni dans $M$.  Alors, $M$ est \emph{plat} sur $R$ si et seulement si
$M/π$ l'est sur $R/π$ et $M[π^{-1}]$ l'est sur $R[π^{-1}]$.
\end{proposition2}

\begin{remarque2}
Pour démontrer la nœthérianité des anneaux considérés
ci-après, on pourrait également utiliser le critère
de Cohen rappelé en \refext{XIX}{19-copremier},
en vérifiant notamment que les idéaux de hauteur $1$
sont principaux.
\end{remarque2}

\subsection{Gonflements}
\label{18B-gonflement}
\subsubsection{}Soit $A$ un anneau local nœthérien, d'idéal maximal $𝔪$.
Suivant \BourbakiAC{IX, appendice, §2}, on note $A]t[$
et on appelle « gonflement (élémentaire) »
de $A$ le localisé de l'anneau de polynômes $A[t]$ en l'idéal premier $𝔪A[t]$.
C'est un anneau local nœthérien.
(Il est noté $A(t)$, par analogie avec les
fractions rationnelles, dans \cite[p. 17—18]{Nag62} ;
voir aussi  \cite[p. 138]{KKK@Matsumura}.)
Plus généralement, on peut considérer un ensemble arbitraire
de variables $t_{e ∈ E}$ et définir l'anneau $A]t_e, e ∈ E[$,
localisé de $A[t_e, e ∈ E]$ en l'idéal premier engendré
par $𝔪$. Rappelons le fait suivant (\BourbakiAC{IX, appendice, prop. 2
et corollaire}).

\begin{proposition2}
L'anneau $A]t_e, e ∈ E[$ est local nœthérien de même dimension
que $A$.
\end{proposition2}

\subsubsection{}Notons que le cas d'un nombre fini de variables est
très élémentaire et que le cas général
résulte du lemme \ega{0$_{III}$}{10.3.1.3},
reproduit en \refext{XIX}{19-ega-noeth},
par passage à la (co)limite.
Pour une autre démonstration, voir également
\BourbakiAC{III.§5, exercice 7}.

\subsubsection{}
Notons que si le corps résiduel de $A$ est $κ$, celui
de $A]t_e, e ∈ E[$ est canoniquement isomorphe à son
extension transcendante pure $κ(t_e, e ∈ E)$.
De plus, on montre que si $F$ est un sous-ensemble de $E$,
le morphisme $A]t_e, e ∈ F[ → A]t_e, e ∈ E[$ est \emph{fidèlement plat}.
(Voir \BourbakiAC{IX, appendice, prop. 2} pour une démonstration
dans le cas  $F=∅$, auquel on se ramène immédiatement.)

\section{Construction de Nagata en dimension $2$, application
cohomologique}
\label{18B-Nagata-Gabber-dimension-2}

\subsection{Dilatation relativement à une série transcendante}

\subsubsection{}Soit $W$ un anneau de valuation discrète,
d'idéal maximal $𝔪_W=(π)$, corps résiduel $k$, corps des
fractions $K$ et de complété $\chap{W}$. On note
$\chap{K}$ le corps des fractions de $\chap{W}$.
Supposons qu'il existe un élément $φ∈ \chap{W}$ \emph{transcendant}
sur $K$. C'est le cas si $W$ est dénombrable ou, par exemple,
lorsque $W=k[[t]]$ auquel cas $φ=∑_n t^{n!}$ convient.
Par translation, on peut supposer que $φ$ appartient à l'idéal maximal de $\chap{W}$.

\subsubsection{}
Fixons un entier $ℓ ≥ 1$ et considérons l'élément $f=(y - φ)^ℓ$ du complété $π$-adique
$\chap{W[y]}$ de $W[y]$. Notons que ce complété s'injecte dans la complétion totale
$\chap{W}[[y]]$ et que $f$ appartient au sous-anneau $\chap{W}[y]$
de $\chap{W[y]}$ et $\chap{W}[[y]]$. En conséquence, le morphisme
canonique $𝒟𝒾𝓁_{π}^{f} W[y] → \chap{W[y]}$ se factorise
en un morphisme $𝒟𝒾𝓁_{π}^{f} W[y] → \chap{W}[y]$.

\begin{proposition2}
Le morphisme $𝒟𝒾𝓁_{π}^{f} W[y] → \chap{W}[y]$
est \emph{plat}.
\end{proposition2}

\begin{démo}
Notons pour simplifier $𝒟$ la source de ce morphisme.
D'après le critère de platitude rappelé ci-dessus,
il suffit de montrer que le morphisme
$𝒟[π^{-1}] → \chap{W}[y][π^{-1}]$ est plat
car $𝒟 → \chap{W}[y]$ induit visiblement un isomorphisme
modulo $π$. Or, lorsque l'on inverse $π$, les $W[y]$-algèbres
dont $𝒟$ est par définition la colimite sont
toutes isomorphes à un anneau de polynômes $K[y,F]$ où,
rappelons-le, $K$ est le corps des fractions $W[π^{-1}]$
de $W$. Par conséquent, il nous faut montrer
que le morphisme $K[y,F] → \chap{W}[y][π^{-1}]$,
composé des morphismes

\begin{minipage}[c]{.46\linewidth}
\begin{center}
\begin{equation*}
\begin{cases}
K[y,F] → K[y, F′]\\
F ↦  {F′}^ℓ
\end{cases}
\end{equation*}
\end{center}
\end{minipage} \hfill
\begin{minipage}[c]{.46\linewidth}
\begin{center}
\begin{equation*}
\begin{cases}
K[y, F′] → \chap{K}[y] \\
F′↦ y-φ
\end{cases}
\end{equation*}
\end{center}
\end{minipage}
La platitude du premier est évidente.
Pour le second, on se ramène par translation et changement
de base à montrer que le morphisme $K[F′] → \chap{K}$, $F′↦ φ$,
est plat. Il se factorise en le composé du passage aux fractions
$K[F′] → K(F′)$ avec l'\emph{injection} $K(F′) ↪ \chap{K}$
déduite de $φ$. Chacun de ces morphismes est plat.
\end{démo}

\begin{remarque2}
La construction de l'anneau de dilatation $𝒟𝒾𝓁_{π}^{f} W[y]$
est inspirée de celle de \cite{example@Nagata}, qui considère
le cas $ℓ=2$. (Voir aussi \cite[appendice, E4.1]{Nag62}
et \cite{Noetherian@HRW}.)
\end{remarque2}

\subsection{Le diviseur $C=V(F)$}

\subsubsection{} On conserve les notations précédentes.
Soit $A$ le localisé de l'anneau de dilatation $𝒟$ en l'idéal
premier image inverse de l'idéal $(π,y) ⊆ \chap{W}[y]$.
C'est un anneau nœthérien par fidèle platitude, de corps
résiduel $k$.

\begin{lemme2}
L'anneau $A$ satisfait les propriétés suivantes :
\begin{enumerate}
\item il est régulier et la suite $(π,F)$ est régulière ;
\item son quotient $A/F$ est intègre ;
\item l'intersection schématique du fermé $V(F)$ avec le spectre épointé
de $A$ est un schéma régulier.
\end{enumerate}
\end{lemme2}

\begin{démo} 
(i) Il suffit d'établir les deux énoncés pour $𝒟$.
Pour se faire, on peut compléter $π$-adiquement
(cf. par exemple \BourbakiAC{X, §4, nº2, cor. 3}).
Le complété de $𝒟$ est isomorphe à $\chap{W[y]}$
de sorte que la régularité de l'anneau et de
la suite $(π,F)$ — c'est-à-dire l'injectivité de
la multiplication par $π$, et par $F$ modulo $π$— sont évidents.
(ii) Si on inverse $π$, l'anneau $A$ devient une localisation de l'anneau de polynômes $K[y,F]$.
La restriction du diviseur à cet ouvert est intègre.
La multiplication par $π$ dans le quotient $A/F$
étant injective d'après (i), l'intégrité de $A/F$
résulte de celle de $A[π^{-1}]/F$.
Le quotient est non nul car $F ∈ 𝔪_A$.
(iii) Sur l'ouvert complémentaire de $V(π)$, l'élément $F$
est une indéterminée de sorte que le résultat est clair.
D'autre part, l'intersection du complémentaire de $V(y)$
avec le diviseur est contenu dans le complémentaire de $V(π)$
car l'équation est — dans le complété $π$-adique — de la forme $(y-φ)^ℓ$, avec $φ
∈ (π)$. Ceci suffit pour conclure.
\end{démo}

\subsubsection{}Notons le fait suivant, trivial mais crucial :
par construction, le diviseur $V(F)$ devient
$V((y - φ)^ℓ)$ dans le complété $\chap{A}$ de $A$
relativement à son idéal maximal.

\subsection{Hensélisation}
Pour simplifier les notations, on suppose
dorénavant le corps $k$ séparablement clos.

\subsubsection{}
Soient $A^\hs$ le hensélisé de $A$ en son point fermé,
$\chap{A}$ le complété de $A$ (ainsi que de $A^\hs$),
et notons $X=\Spec(A^\hs)$ et $\chap{X}=\Spec(\chap{A})$
leurs spectres respectifs, ainsi que $⋆$ et $\chap{⋆}$ les points fermés.
Comme $\Spec(A)$, le schéma $X$ est intègre.
De plus, le diviseur $C$ d'équation $F=0$ dans $X$
est réduit, cette propriété étant
également conservée par hensélisation.

\subsubsection{}
Vérifions que le diviseur $C$ est \emph{irréductible}.
D'après le théorème de comparaison de Elkik,
le morphisme $π₀(\chap{C}-\chap{⋆}) → π₀(C - ⋆)$ est une
bijection. Or, $\chap{C}-\chap{⋆}$ est connexe :
dans un anneau local régulier $B$, le spectre
$\Spec(B/g^ℓ)$ est irréductible pour tout $g ∈ 𝔪_B-𝔪_B²$.
Ainsi, l'ouvert $C - ⋆$ de $C$ est connexe et, finalement,
$C$ est irréductible.

\subsection{Application cohomologique}
\emph{On suppose dorénavant l'entier $ℓ$ inversible sur $X$.}

\subsubsection{}
Notons pour simplifier $G$ l'élément $(y-φ)^ℓ$ de $\chap{A}$,
de sorte que — par construction — on a l'égalité $F=G^ℓ$ dans $\chap{A}$.
Notons $U$ l'ouvert affine $X-C$ du schéma strictement local $X$.
Soit $j$ l'immersion ouverte de $U$ dans le spectre
épointé $X-⋆$ et $i$ l'immersion fermée $C-⋆ ↪ X-⋆$.

\subsubsection{}Le triangle
\[
i_* i^! → \mathrm{Id} → 𝖱 j_* j^*
\]
sur $X-⋆$ induit la suite exacte
\[
H³_{C-⋆}(X-⋆,𝐙/ℓ) → H³(X-⋆,𝐙/ℓ) → H³(U,𝐙/ℓ) →
H⁴_{C-⋆}(X-⋆,𝐙/ℓ)=H²(C-⋆,i^! 𝐙/ℓ[2]).
\]
Par pureté (\refext{XVI}{16-theoreme-purete-absolue}), le faisceau de cohomologie locale $i^! 𝐙/ℓ$ est
constant concentré en degré $2$. Or le groupe de cohomologie $H²(C-⋆,𝐙/ℓ)$ est
nul : la cohomologie du corps des fractions d'un anneau $B$
strictement local intègre de dimension $1$ est nul en
degré $≥2$. (On se ramène au cas bien connu d'un anneau de valuation
discrète en observant que le normalisé
de $B$ dans son corps des fractions est un
anneau nœthérien, de Dedekind, et strictement local
car colimite locale d'anneaux strictement locaux.)
Il en résulte que la flèche de restriction
$H³(X-⋆,𝐙/ℓ) → H³(U,𝐙/ℓ)$ est surjective ; nous allons
voir que c'est un isomorphisme.

\subsubsection{}Comme rappelé ci-dessus, le
morphisme $H³_{C-⋆}(X-⋆,𝐙/ℓ) → H³(X-⋆,𝐙/ℓ)$
s'identifie, par pureté, au morphisme
de Gysin $\mathrm{Gys}(f):H¹(C-⋆,𝐙/ℓ) → H³(X-⋆,𝐙/ℓ)$.
Il résulte de la commutativité du diagramme
\begin{center}
 \begin{tikzpicture}[auto]
  \matrix(diag)[matrix of math nodes,column sep=8em,row sep=5ex]
  {|(D)|  H¹(C-⋆,𝐙/ℓ)   & |(X)| H³(X-⋆,𝐙/ℓ)  \\
   |(Dc)|  H¹(\chap{C}-⋆,𝐙/ℓ) & |(Xc)| H³(\chap{X}-⋆,𝐙/ℓ),  \\};
  \draw[->] (D) -- (Dc);
  \draw[->] (D) -- node{$\mathrm{Gys}(F)$} (X);
  \draw[->] (Dc) -- node{$\mathrm{Gys}(F)=\mathrm{Gys}(G^ℓ)$} (Xc);
  \draw[->] (X) -- (Xc);
 \end{tikzpicture}
\end{center}
de l'égalité $\mathrm{Gys}(F)=ℓ ⋅ \mathrm{Gys}(G)$
et enfin du fait que les flèches verticales sont
des isomorphismes (comparaison à la complétion,
\cite[6.6.4]{TUBULAR@Fujiwara}) que le morphisme $\mathrm{Gys}(F)$ est nul.
(La commutativité du diagramme résulte par exemple
de la définition \refext{XVI}{16-definition-cl-i}
et de \refext{XVI}{16-proposition-classes-generalisees-cb}.)
Ainsi, le morphisme de restriction induit un
isomorphisme
\[
H³(X-⋆,𝐙/ℓ) ⥲ H³(U,𝐙/ℓ).
\]
Or le terme de gauche est non nul, à nouveau par pureté.
Le schéma affine $U$ est donc de dimension cohomologique $>2$. CQFD.

\section{Séries formelles de Gabber, application cohomologique}
\label{18B-Gabber-Nagata-dimension-supérieure}

On étend la construction précédente en dimension
arbitraire $≥ 2$.

\subsection{Une série formelle et sa décomposition}

\subsubsection{}
\label{18B-décomposition-série-formelle}
Soit $A$ un anneau commutatif.
Rappelons que l'application $A$-linéaire
\[
∑_{i=1}^n A[x_{j ≠ i}][[x_i]] → A[[x₁,…,x_n]]
\]
somme des injections canoniques est \emph{surjective} :
si $G ∈  A[[x₁,…,x_n]]$, on peut par exemple regrouper pour chaque $i ∈ [1,n]$
les termes $ax₁^{β₁}\cdots x_n^{β_n} ∈ A[[x₁,…,x_n]]$ de $G$
pour lesquels $β_i=\max_{j ∈ [1,n]} β_j$ et $β₁,…,β_{i-1}<β_i$.
(Cette dernière condition n'est là que pour définir $i$
de façon non ambiguë ; tout autre choix conviendrait.)
La somme $g_i$ de ces termes appartient à $A[x_{j ≠ i}][[x_i]]$, et
$G=g₁+\cdots+g_n$.

\subsubsection{}On fixe maintenant deux entiers non nuls $n$ et $ℓ$
et on considère
\[
S=\big(y + ∑_{i=1}^n ∑_{α=1}^{∞} t_{i α} x_i^α\big)^ℓ ∈ 𝐙[y,t_{i
∈ [1,n], α ≥ 1}][[x₁,…,x_n]].
\]
Il résulte de l'observation précédente que l'on peut écrire cette
série sous la forme
\[
y^ℓ + f₁+\cdots+f_n
\]
où chaque $f_i$ est une série formelle en $x_i$, à
coefficients polynomiaux en les autres variables.

\subsubsection{}
\label{18B-définition-fi}
Afin que la proposition de platitude ci-dessous soit vraie,
on procède de façon légèrement différente pour définir les
séries formelles $f_i ∈ 𝐙[y,t_{jα},x_{k ≠ i}][[x_i]]$ telles
que $S-y^ℓ=∑_{i=1}^n f_i$. Écrivons
\[
S=y^ℓ + ∑_{α=1}^∞ (∑_{i=1}^n t_{i α}x_i^α)^ℓ +
\big(\text{élément de degré }<ℓ\text{ en les }t_{jβ}\big).
\]
Soient $i ∈ [1,n]$ et $α ≥ 1$ des indices.
On considère les termes $ax₁^{β₁}\cdots x_n^{β_n}$ de $(∑_{i=1}^n t_{i α}x_i^α)^ℓ$
pour lesquels $i$ est le plus grand indice tel que $β_i ≠ 0$,
c'est-à-dire les termes de $(∑_{i=1}^n t_{i α}x_i^α)^ℓ$
qui sont dans $𝐙[t_{1 α},…,t_{n α},x₁,…,x_i]$ mais pas
dans $𝐙[t_{1 α},…,t_{n α},x₁,…,x_{i-1}]$.
À $i$ fixé, la somme sur $α$ de ces termes
est un élément $f_{i,=ℓ}$ de $𝐙[t_{jβ},x_{k ≠ i}][[x_i]]$.
Par construction, on a l'égalité
$ ∑_{α=1}^∞ (∑_{i=1}^n t_{i α}x_i^α)^ℓ=∑_{i=1}^n f_{i,=ℓ}$.
Enfin, on décompose le terme restant,
$S-y^ℓ-∑_{i=1}^n f_{i,=ℓ}$, en une somme $∑_{i=1}^n f_{i,<ℓ}$
où chaque $f_{i,<ℓ}$ appartient à $𝐙[y,t_{jβ},x_{k ≠ i}][[x_i]]$,
en procédant par exemple comme en \ref{18B-décomposition-série-formelle}.
On pose alors $f_i=f_{i,=ℓ}+f_{i,<ℓ}$ ; chacun de ces
éléments est de degré total en les $t_{j β}$ inférieur ou égal à $ℓ$,
et appartient donc également à $𝐙[y,t_{j,β ≠ α}][[x₁,…,x_n]][t_{1
α},…,t_{n α}]$ pour chaque $α ≥ 1$.

\begin{proposition2}
\label{18B-platitude-série-formelle-universelle}
Fixons $α ≥ 1$. Notons $T_i=t_{i α}$ pour chaque $i ∈ [1,n]$
et $R_α$ l'anneau $𝐙[y,t_{j,β ≠ α}][[x₁,…,x_n]]$.
Le morphisme
\[
R_α[F₁,…,F_n] → R_α[T₁,…,T_n]
\]
\[
F_i↦ f_i
\]
est \emph{libre}, donc \emph{plat}, au-dessus de l'ouvert
$x₁\cdots x_n ≠ 0$ de $\Spec(R)$.
\end{proposition2}

\begin{démo}
Par construction, chaque $f_i$ est une somme $f_{i,=ℓ} + f_{i,<ℓ}$, où
$f_{i,=ℓ}$ (resp. $f_{i,<ℓ}$) est un polynôme dans $R_α[T₁,…,T_n]$ de
degré total égal (resp. strictement inférieur) à $ℓ$.  De plus,
$f_{i,=ℓ}$ est, comme polynôme en $T_i$, de la forme $x_i^{αℓ} T_i^ℓ +
∑_{m<ℓ} c_m T_i^m$ où les $c_m$ appartiennent à $R_α[T_{j<i}]$.
Munissons les monômes de $R_α[T₁,…,T_n]$ de l'ordre lexicographique
gradué suivant : 
$T₁^{d₁} \cdots T_n^{d_n} \mathrel{\preceq} T₁^{d₁′} \cdots T_d^{d_n′}$
si et seulement si
$\sum d_i < \sum d′_i$ ou $\sum d_i = \sum
d′_i$ et $d_i < d'_i$ pour le plus grand $i$ tel que
$d_i \neq d'_i$.
Il est clair que le \emph{terme de tête} $\initial_{\preceq}(g)$
pour cet ordre d'un polynôme $g=f₁^{q₁}\cdots f_n^{q_n}$
en les $f_i$ est $T₁^{q₁ ℓ}\cdots T_n^{q_n ℓ}$, à multiplication près
par un monôme en les $x_i$.
Il en résulte immédiatement qu'en inversant les $x_i$,
l'anneau $R_α[T₁,…,T_n]$ est libre sur $R_α[F₁,…,F_n]$
de base les monômes $T₁^{r₁}\cdots T_d^{r_d}$,
avec $0 ≤ r_i < ℓ$.
\end{démo}

\subsubsection{}Soit $A$ un anneau local nœthérien régulier
de dimension $n$ dont on note $x₁,…,x_n$ un système
régulier de paramètres. On pourra penser par
exemple au localisé $k[x₁,…,x_n]_{(x₁,…,x_n)}$ d'un anneau
de polynômes sur un corps.
On considère le gonflement $A]\sous{t}[$ défini
en \ref{18B-gonflement}, où $\sous{t}$
est l'ensemble des variables $\{t_{i α}: i ∈ [1,…,n], α
∈ 𝐍_{≥1}\}$.

\subsection{Construction d'un anneau local régulier pathologique}
Soit $i ∈ [1,…,n]$. Notons encore $f_i$
l'image de la série formelle à coefficients entiers
considérée en \ref{18B-définition-fi}
dans le complété $x_i$-adique de $A]\sous{t}[ [y]$,
et $𝒟_i$ l'anneau de dilatation
$𝒟𝒾𝓁_{x_i}^{f_i} A]\sous{t}[ [y]$.
Le produit tensoriel $𝒫$ de ces $A]\sous{t}[ [y]$-algèbres s'envoie naturellement
dans le complété $\chap{A]\sous{t}[} [[y]]$,
où la première complétion est faite relativement à l'idéal
maximal $(x₁,…,x_n)$ de $A]\sous{t}[$.
On note $𝒟$ le localisé de $𝒫$ en l'image de l'idéal maximal $(x₁,…,x_n,y)$ de
$\Spec(\chap{A]\sous{t}[} [[y]])$.

\begin{proposition2}
Le morphisme $𝒟 → \chap{A]\sous{t}[} [[y]]$ est fidèlement \emph{plat}.
\end{proposition2}

Il en résulte que l'anneau $𝒟$ est local \emph{nœthérien}, régulier.

\begin{remarque2}
Notons qu'il est clair que $𝒟$ est « quasi-régulier ».
En effet, le gradué de $𝒫$ relativement à l'idéal $(x₁,…,x_n,y)$ est
une algèbre symétrique : pour chaque entier $r$,
le morphisme $A]\sous{t}[ [y] → 𝒫$ induit un isomorphisme
modulo $(x₁,…,x_n)^r$.
\end{remarque2}

\begin{démo}
Il suffit de montrer que le morphisme $𝒫 → \chap{A]\sous{t}[} [[y]]$ est
plat. D'après le critère de platitude rappelé
précédemment (\ref{18B-critère-platitude}), il suffit de montrer la platitude
sur l'ouvert $x₁\cdots x_n ≠ 0$.
En effet, le cas où seuls certains $x_i$
sont nuls se ramène à ce cas particulier :
en tensorisant par $A]\sous{t}[ [y]/x_i$ le morphisme
$𝒫 →  \chap{A]\sous{t}[} [[y]]$, on obtient une flèche
du même type définie par l'anneau $A/x_i$
de dimension $n-1$ et des séries qui coïncident
avec l'évaluation en $x_i=0$
des $f₁,…,\chap{f_i},…,f_n$.
Pour chaque $i$, $𝒟_i[x_i^{-1}]$ est une algèbre
de polynômes $A]\sous{t}[ [y,F_i][x_i^{-1}]$
de sorte que le morphisme dont on souhaite montrer la platitude
est
\[
A]\sous{t}[ [y,F₁,…,F_n] [\frac{1}{x₁ \cdots x_n}] →
\chap{A]\sous{t}[} [[y]] [\frac{1}{x₁ \cdots x_n}],
\]
\[
F_i↦ f_i.
\]
Il suffit de montrer que pour chaque sous-ensemble
\emph{fini} $𝒯$ des variables $\sous{t}$, le morphisme
\[
A]t ∈ 𝒯[ [y,F₁,…,F_n] [\frac{1}{x₁ \cdots x_n}]
→ \chap{A]t ∈ 𝒯[} [[y]] [\frac{1}{x₁ \cdots x_n}]
\]
est plat. Quitte à agrandir un tel ensemble $𝒯$,
on se ramène au cas où $𝒯$ est \emph{cofini}, de complémentaire
des variables $t_{1 α},…,t_{n α}$ pour un indice $α ≥ 1$
quelconque. Posant alors $R_α=A]t_{i,β≠α}[$
et $R′=A]\sous{t}[$, il suffit de montrer
que le morphisme $R_α[y,F₁,…,F_n]→ R′[[y]]$
est plat au-dessus de l'ouvert $x₁\cdots x_n ≠ 0$.
Un dernier dévissage nous ramène à montrer
la platitude du morphisme $R_α[y,F₁,…,F_n] →
R_α[y,t_{1 α},…,t_{n α}]$, au-dessus du même ouvert.
Ce dernier point résulte de la proposition
\ref{18B-platitude-série-formelle-universelle}.
\end{démo}

\begin{proposition2}
Le diviseur $C=V(y^ℓ+ F₁+\cdots+F_n)$ de $\Spec(𝒟)$
est \emph{régulier} hors du point fermé.
\end{proposition2}

Dans cet énoncé, on note abusivement $F_i$
l'image dans $𝒟$ de l'élément de $𝒟_i$ correspondant
à $f_i$ (cf. \ref{18B-dilatations-formelles}).

\begin{démo}
Il suffit de montrer que pour chaque sous-ensemble
strict $E$ de $[1,n]$, l'intersection schématique
de $C$ avec le sous-schéma
$X_E=\{x_i =0, i ∈ E ; x_i ≠ 0, i ∉ E\}$
de $\Spec(𝒟)$ est un diviseur régulier de $X_E$.
Si $E=∅$, cela résulte du fait que le schéma $X_∅$
est le localisé d'une algèbre de polynômes
en les $y,F₁,…,F_n$. Le cas général se ramène
aisément à ce cas particulier. (Remarquons que si
$x_i=0$, il en est de même de $F_i$.)
\end{démo}

\begin{proposition2}
L'image inverse de $C$ dans un localisé
strict de $\Spec(𝒟)$ est \emph{irréductible}.
\end{proposition2}

\begin{démo}
Même argument qu'en dimension $2$.
\end{démo}

\begin{corollaire2}
\label{18B-ouvert-dim-d-dimcoh-2d-1}
Pour tout entier $d ≥ 1$, il existe schéma nœthérien
strictement local régulier $X$ de dimension $d$ et un ouvert affine $U=D(f)$
de $ℓ$-dimension cohomologique $2d-1$.
\end{corollaire2}

\begin{démo}
La même démonstration qu'en dimension $2$
nous permet de minorer la dimension cohomologique
par $2d-1$. D'après \refext{XVIII-A}{18A-Gabber's bound1},
c'est une égalité.
\end{démo}

\section{Dimension cohomologique : majoration d'une « fibre de
Milnor générique »}
\label{18B-majoration-Milnor}

\subsection{Énoncé}

\begin{théorème2}
\label{18B-majoration Milnor générique}
Soit $R → R′$ un morphisme local essentiellement de type fini
d'anneaux nœthériens strictement locaux intègres.
Notons $K$ le corps des fractions de $R$.
Alors, pour tout nombre premier $ℓ$ inversible sur $R$,
on a la majoration
\[
\cd_ℓ(R′ ⊗_R K) ≤ \dim(R′),
\]
où le terme de gauche désigne la $ℓ$-dimension cohomologique
étale du spectre de l'anneau $R′ ⊗_R K$ et le terme
de droite désigne la dimension de Krull de $R′$.
\end{théorème2}

\subsubsection{}
\label{18B-essentiellement-tf}
Dans cet énoncé, l'hypothèse de finitude sur $f$
signifie que ce morphisme est localement de la forme
$\Spec(B) → \Spec(A)$, avec $B$ une colimite de $A$-algèbres
de type fini à morphismes de transition étales.

\begin{corollaire2}
Soit $R$ un anneau strictement local nœthérien intègre
de corps des fractions $K$ et soit $ℓ$ un nombre premier
inversible sur $R$. Alors, on a la majoration
\[
\cd_ℓ(K) ≤ \dim(R).
\]
\end{corollaire2}

\begin{remarque2}
Réciproquement, on peut montrer que si, sous les hypothèses du
corollaire, $U$ est un ouvert non vide strict de $\Spec(R)$, alors
$\mathrm{cd}_ℓ(U) ≥ \dim(R)$ et que, lorsque $R′ ⊗_R K ≠ 0$,
la majoration du théorème est une égalité.
Une façon de procéder, due à O. Gabber (non publié),
est d'utiliser une variante de la méthode (également due à O. Gabber) exposée
dans \cite[§6.1]{pdim@GO} et reposant sur une
« astuce quadratique ». Pour une autre méthode,
voir §\ref{18B-minoration-ouvert-spectre-épointé}, \emph{infra}.
Rappelons que la minoration « limite »
$\mathrm{cd}_ℓ(K) ≥ \dim(R)$ est élémentaire : on procède par
spécialisations successives en codimension $1$ (voir \sga{4}{X}{2.4}).
\end{remarque2}

\subsection{Démonstration}
\subsubsection{}On procède par récurrence sur $d′=\dim(R′)$ et
l'on se ramène au cas excellent.

\subsubsection{}
Notons $X=\Spec(R)$, $Y=\Spec(R′)$ et respectivement
 $X^⋆$ et $Y^⋆$ les spectres épointés.
Considérons l'ouvert $Y_⋆=Y×_X X^⋆$ de $Y^⋆$,
$V =\Spec(R′ ⊗_R K)$ la fibre générique de $Y → X$
et enfin $j$ le morphisme $V ↪ Y^⋆$.
Il résulte de l'hypothèse de récurrence que pour
chaque faisceau $ℱ$ de $𝐙/ℓ$-modules sur $V$,
le complexe $Φ_ℱ=𝖱 j_* ℱ$ appartient à $𝖣^{≤ \mathrm{cod}}(Y^⋆)$,
où $\mathrm{cod}$ est la fonction de perversité $y↦\dim 𝒪_{Y,y}$.
(Ceci est encore vrai avec $y↦d′-\dim \sur{\{y\}}$.)
On veut montrer que $H^r(V,ℱ)=H^r(Y^⋆,Φ_ℱ)$ est nul pour $r>d′$.
Fixons un tel $r$ et une classe $c ∈ H^r(Y^⋆,Φ_ℱ)$.

\subsubsection{}
On suppose $d ≥ 2$, et on choisit un système de paramètres $x₁,…,x_d$
pour l'anneau strictement local $R$.
Soit $Z=Y\}t₁,…,t_{d-1}\{$ le « gonflement étale »,
hensélisé strict de $𝐀^{d-1}_Y=Y[t₁,…,t_{d-1}]$ en un point générique
géométrique de la fibre spéciale sur $Y$.
L'« hyperplan » $H=V(t₁x₁+\cdots+t_{d-1}x_{d-1}+x_d)$ de $Z$
est de codimension $1$, essentiellement lisse au-dessus de $Y_⋆$.
Considérons le triangle distingué
\[
𝖱Γ_{H_⋆}(Z_⋆,Φ_ℱ) → 𝖱Γ(Z_⋆,Φ_ℱ) → 𝖱Γ(Z_⋆-H_⋆,Φ_ℱ) →,
\]
où l'on note $Z_⋆$ le produit fibré $Z ×_Y Y_⋆$ et, abusivement, $Φ_ℱ$
ses diverses images inverses. Soit $i$ l'immersion fermée $H_⋆ ↪ Z_⋆$,
où $H_⋆=H×_Z Z_⋆$.  Le morphisme $𝐙/ℓ → i^! 𝐙/ℓ(1)[2]$ de complexes
sur $H_⋆$ est un isomorphisme par pureté relative.  On montre par dévissage
qu'il en est de même de la flèche ${Φ_ℱ}_{|H_⋆} →
i^!{Φ_ℱ}_{|Z_⋆}(1)[2]$ obtenue par tensorisation à partir de la
précédente.  On utilise ici le fait que la restriction
de $Φ_ℱ$ à $Z_⋆$ provient de la base $Y_⋆$.  On en
tire le morceau de suite exacte :
\[
H^{r-2}(H_⋆,Φ_ℱ)(-1) → H^r(Z_⋆,Φ_ℱ) → H^r(Z_⋆-H_⋆,Φ_ℱ).
\]
Notons que $Z_⋆-H_⋆=Z-H$ car $H$ contient la fibre spéciale
de $Z → X$.
Soient $\chap{Y}$ le complété ($𝔪_{R′}$-adique) de $Y$
et $\gtilde{Z}$ un hensélisé strict du produit
fibré $Z×_Y \chap{Y}$.
Notons que le morphisme $\gtilde{Z} → Z$ est un morphisme local
entre schémas strictement locaux induisant
un isomorphisme sur la complétion le long
de la fibre spéciale sur $Y$.
Il résulte donc du théorème de comparaison
de Fujiwara-Gabber (\cite[6.6.4]{TUBULAR@Fujiwara}) que le morphisme
$H^r(Z_⋆,Φ_ℱ) → H^r(\gtilde{Z}_⋆,Φ_ℱ)$
est un \emph{isomorphisme} pour chaque $r$.
Il en est de même de $H^{r-2}(H_⋆,Φ_ℱ) →
H^{r-2}(\gtilde{H}_⋆, Φ_ℱ)$,
où $\gtilde{H}_⋆=H_⋆ ×_{Z_⋆} \gtilde{Z}_⋆$.
Le schéma $\gtilde{Z}_⋆-\gtilde{H}_⋆$ est un ouvert
affine, qui coïncide avec $\gtilde{Z}-\gtilde{H}$,
d'un schéma strictement local essentiellement de type fini
sur le schéma local nœthérien complet $\chap{Y}$.
Comme l'appartenance de $Φ_ℱ$ à $𝖣^{≤\mathrm{cod}}$
est préservée par complétion,
il résulte du théorème de Lefschetz affine
(\refext{XV}{15-prop-principale-locale}),
dans le cas excellent, que le groupe
de cohomologie $H^r(\gtilde{Z}_⋆-\gtilde{H}_⋆,Φ_ℱ)$
est nul pour chaque $r> \dim(Z)=\dim(Y)=d′$. En conséquence,
le morphisme $H^r(\gtilde{Z}_⋆,Φ) → H^r(\gtilde{Z}_⋆-\gtilde{H}_⋆,Φ)$
est nul pour les mêmes $r$. De ce fait, des théorèmes de comparaisons
susmentionnés et de la compatibilité du morphisme de Gysin
à la complétion, il résulte formellement que toute classe $c ∈ H^r(Z_⋆,Φ_ℱ)$
provient d'une classe dans $H^{r-2}(H_⋆,Φ_ℱ)(-1)$
et est donc tuée par restriction à $Z_⋆-H_⋆$.

\subsubsection{}
Il existe donc un voisinage étale $e:W → 𝐀^{d-1}_Y$ dont l'image
rencontre la fibre spéciale sur $Y$ tel que la classe $c ∈ H^r(Y^⋆,Φ_ℱ)$
soit tuée par restriction à $W-H_W$. Notons $k′$
le corps résiduel de $Y$, et $k$ celui de $X$.
L'ensemble $k^{d-1}$ est dense dans $𝐀^{d-1}_{k′}$, car $k$
est infini. Il en résulte qu'il existe une section $σ: Y →
𝐀^{d-1}_Y$, correspondant à des spécialisations des $t_i$
à valeurs dans $R$, telle que $W_σ$ ait une fibre spéciale
sur $Y$ non vide. Le schéma $Y$ étant strictement local, on relève
cette section en $Y → W_σ → W$. La classe de cohomologie $c$
est donc nulle sur $Y-H_Y$, où $H_Y$ est maintenant une hypersurface
d'équation $x_d+t₁x₁+\cdots+t_{d-1}x_{d-1}$ à coefficients $t_i$ dans $R$.
Cet ouvert affine $Y-H_Y$ contient la fibre générique $V=Y ⊗_R K$
car l'élément $x_d+t₁x₁+\cdots+t_{d-1}x_{d-1} ∈ R$ est non nul,
les $x_i$ constituant un système de paramètres de $R$.
Finalement la restriction de $c ∈ H^r(Y^⋆,Φ_ℱ)$ à $H^r(V,Φ_ℱ)=H^r(V,ℱ)$,
qui est la classe dont on est parti, est nulle. CQFD.

\section{Majoration : amélioration}
\label{18B-majoration-générale}

\subsection{Énoncé}

\subsubsection{}Soit $f:Y → X$ un morphisme entre espaces topologiques
sobres non vides.
On note
\[
\mathrm{dim.cat}(f)=\sup \{n ∈ 𝐍 : ∃ y₀ ⤳  y₁ ⤳  \cdots ⤳ y_n,
f(y₀) ≠ f(y₁) ≠ \cdots ≠ f(y_n) \} ∈ 𝐍 ∪ \{∞\}
\]
la \emph{dimension caténaire} de $f$, où chaque $y_i ⤳ y_{i+1}$ est une spécialisation.

\subsubsection{}
Par construction, $\mathrm{dim.cat}(f:Y → X)$ est majorée
par les dimensions de $X$ et de $Y$ avec égalité par exemple
lorsque $f$ est l'identité.
Plus généralement, lorsque $f$ est un morphisme \emph{générisant} (\ega{I′}{3.9.1})
— comme c'est le cas d'un morphisme plat de schémas —
la dimension caténaire coïncide avec la dimension de
l'image.

\begin{remarque2}
Si $f$ est un morphisme dominant essentiellement de type
fini (au sens du \ref{18B-essentiellement-tf}) entre schémas nœthériens intègres, on peut montrer
que la dimension caténaire de $f$
est la dimension de l'image d'une platification
de $f$.
\end{remarque2}

\begin{théorème2}
Soit $f:Y → X$ un morphisme essentiellement de type fini
entre schémas nœthériens strictement locaux et soit $V$ un
ouvert affine de $Y$. Alors, pour tout nombre premier $ℓ$
inversible sur $X$, on a la majoration
\[ \mathrm{cd}_ℓ(V) ≤
\dim(Y)+\max\big(0,\mathrm{dim.cat}(f)-1\big).  \] En
particulier, si $\dim(X) ≥ 1$, $\mathrm{cd}_ℓ(V) ≤
\dim(Y)+\dim(X)-1$.
\end{théorème2}

\begin{corollaire2}
Soit $f:Y → X$ un morphisme essentiellement de type fini
entre schémas nœthériens strictement locaux, où $\dim(X) ≥ 1$,
et soit $V$ un ouvert affine de $Y$. Alors, pour tout nombre premier $ℓ$
inversible sur $X$, on a la majoration
\[ \mathrm{cd}_ℓ(V) ≤
\dim(Y)+\dim(X)-1.\]
\end{corollaire2}

\begin{corollaire2}
\label{18B-n-entre-d-et-2d-1}
Soit $d ≥ 1$ un entier et soit $n$ un entier dans
l'intervalle fermé $[d,2d-1]$. Il existe un schéma
nœthérien strictement local $X$, régulier de dimension $d$,
et un ouvert affine $U$ de ce schéma tel pour tout nombre
premier $ℓ$ inversible sur $X$ on ait l'égalité
\[
\mathrm{cd}_ℓ(U)=n.
\]
\end{corollaire2}

\begin{démo}
Il suffit de montrer que pour tout entier $d ≥ 1$, et tout entier $r ≥ 0$,
il existe un schéma nœthérien strictement local régulier $Y$ de
dimension $d+r$ et un ouvert affine $V$ de $Y$
de $ℓ$-dimension cohomologique égale
à $2d+r -1$. Soient $X$ et $f$ comme en
\ref{18B-ouvert-dim-d-dimcoh-2d-1} :
l'ouvert affine $U=X[f^{-1}]$ est
de dimension $d$, $ℓ$-dimension cohomologique $δ=2d-1$.
Plus précisément, il résulte de la démonstration
qu'il existe une classe non nulle dans $H^δ(U,𝐙/ℓ)$.
Considérons maintenant $Y=X[T]_{(0)}$ un hensélisé stricte
de la droite affine sur $X$ en l'origine de la fibre
spéciale, $g=fT ∈ Γ(Y,𝒪_Y)$ et $V$ l'ouvert affine
$Y[g^{-1}]$. On a $\dim(Y)=d+1$. Par pureté cohomologique,
on vérifie immédiatement que le groupe de cohomologie $H^{d+1}(V,𝐙/ℓ)$
est également non nul. Par récurrence, on obtient une
paire $(Y,V)$ comme ci-dessus telle que $\mathrm{cd}_ℓ(V) ≥ 2d+r-1$.
D'après le corollaire précédent, on a également la majoration
$\mathrm{cd}_ℓ(V) ≤ 2d+r-1$, d'où l'égalité.
\end{démo}

\subsection{Démonstration}

\subsubsection{}On procède par récurrence sur la dimension caténaire de $f$.

\subsubsection{$\mathrm{dim.cat}(f)=0$} Cette égalité
se produit si et seulement si $Y$ est contenu dans la fibre spéciale. Le théorème
est donc connu dans ce cas : on est sur un corps
donc dans une situation excellente.

\subsubsection{$\mathrm{dim.cat}(f)=1$}
On peut supposer les schémas $X$ et $Y$ réduits.
Quitte à procéder par récurrence sur la dimension de $Y$, on peut également
supposer $Y$ irréductible : si $Y$ est la réunion
de deux fermés stricts $Y₁$ et $Y₂$, considérer
par exemple le morphisme $π:Y₁ ⨿ Y₂ → Y$ et
la suite exacte $0 → ℱ → π^* ℱ → i_*ℋ → 0$, où $i$
est l'immersion fermée $Y₁ ∩ Y₂ ↪ Y$ et $ℋ$ un faisceau sur
cette intersection. Quitte à remplacer $X$
par l'adhérence de l'image de $f$, on peut également supposer
la base intègre et $f$ dominant.
Soit $η$ (resp. $s$) le point générique (resp. fermé) de $X$
et $η′$ (resp. $s′$) le point générique (resp. fermé) de $Y$.
Comme tout point $y$ de $Y$ s'insère dans
une suite de spécialisations $η′ ⤳ y ⤳ s′$, d'image
$η ⤳ f(y) ⤳ s$, il résulte de l'hypothèse $\mathrm{dim.cat}(f)=1$
que $f(y)=η$ ou $f(y)=s$. Soient $\chap{X}$ le complété de $X$ et $\gtilde{Y}$ un hensélisé
strict du produit fibré $Y×_X \chap{X}$. C'est un schéma
strictement local de dimension $\dim(Y)$ et
\emph{excellent} car essentiellement de type fini sur le schéma
local nœthérien complet — donc excellent — $\chap{X}$.
On note $\gtilde{V}$ l'ouvert $V×_Y \gtilde{Y}$.
Il résulte de \sga{4½}{[Finitude]}{1.9} qu'au-dessus de $η$,
et donc au-dessus de $X-\{s\}$, la formation des images directes par $j:V ↪ Y$ commute au changement
de base $\gtilde{Y} → Y$.
En d'autres termes, si $V′=V ∪ (Y-Y_s)$ et $j′$ désigne
l'immersion intermédiaire $V ↪ V′$,
la formation de $𝖱 j′_*$ commute à $\gtilde{Y} → Y$.
Il en est de même pour $j″: V′ → Y$ d'après le théorème
de comparaison à la complétion de Fujiwara-Gabber
(\cite[6.6.4]{TUBULAR@Fujiwara}).
On utilise ici le fait que si $F$ est un fermé de $Y$
inclus dans la fibre spéciale alors les complétés
de $Y$ le long de $F$ et celui de l'hensélisé
de $\gtilde{Y}$ (en le point correspondant au point
fermé de $Y$) sont naturellement isomorphes.
Finalement, le foncteur $𝖱Γ(V)=𝖱Γ(Y) ∘ 𝖱j_*$ s'identifie
au foncteur $𝖱Γ(\gtilde{Y}) ∘ 𝖱 \gtilde{j}_*=𝖱Γ(\gtilde{V})$,
appliqué à l'image inverse.
Ainsi on a l'inégalité $\mathrm{cd}_ℓ(V) ≤ \mathrm{cd}_ℓ(\gtilde{V})$.
Le schéma $\gtilde{Y}$ est quasi-excellent, le terme de
droite est donc justiciable du théorème de Lefschetz affine.
L'inégalité $\mathrm{cd}_ℓ(V) ≤ \dim(Y)$ résulte alors
de l'égalité $\dim(Y)=\dim(Y′)$.

\begin{remarque2}
Lorsque $\dim(R)=1$, on a vu en \refext{XIII}{13-compléments annulation}
que le théorème peut également se démontrer
par normalisation.
\end{remarque2}

\subsubsection{$\mathrm{dim.cat}(f)>1$}
Notons à nouveau $j$ l'immersion ouverte $V ↪ Y$.
Par restriction à la fibre générique, on a un isomorphisme
$𝖱Γ(V_η,ℱ)=𝖱Γ(Y_η,(R j_* ℱ)_{|Y_η})$.
Si $y$ est un point géométrique de $Y$ localisé en $Y_η$,
la fibre $(R^q j_* ℱ)_y$ est nulle dès lors que $q>\dim(𝒪_{Y_η,y})$.
Cela résulte par passage à la limite du théorème d'Artin
pour les schémas affines de type fini sur un corps et du fait trivial
que $Y_{(y)} → Y$ se factorise à travers $Y_η$.
Soient $q ≥ 0$ un entier, $𝒢$ un sous-faisceau constructible de $(R^q j_* ℱ)_{|Y_η}$
et $S$ l'adhérence dans $Y$ de son support.
D'après ce qui précède, on a la majoration
$\codim(S_η,Y_η) ≥  q$. Il en résulte que $\codim(S,Y) ≥ q$,
et ce sans hypothèse de caténarité sur les schémas.
De la suite spectrale de composition des foncteurs
et du théorème \ref{18B-majoration Milnor générique},
on déduit que le groupe de cohomologie $H^n(Y_η,𝒢)$ est nul lorsque
$n>\dim(Y)$.

Considérons la flèche d'adjonction $ℱ → k_*k^* ℱ$,
où $k$ est l'immersion $V_η ↪ V$, et $𝒦$ son noyau.
Par construction, la restriction de $𝒦$ à $V_η$ est nulle.
La dimension de l'adhérence du support de $𝒦$ est donc
au plus $\dim(Y)-1$. Il résulte donc de l'hypothèse de récurrence
que le résultat d'annulation désiré est connu pour $𝒦$.
Procédant de même pour le conoyau de l'adjonction précédente,
on se ramène à démontrer l'annulation du groupe
$H^p(V,R⁰k_*k^* ℱ)$ pour $p ≥ \dim(Y)+\mathrm{dim.cat}(f) (>\dim(Y))$.
Compte tenu du résultat d'annulation précédemment établi
pour $𝖱Γ(V_η,ℱ)=𝖱Γ(V,R k_* k^* ℱ)$ et de la suite spectrale de Leray
$E₂^{p,q}=H^p(V,R^q k_* k^* ℱ) ⇒ H^{p+q}(V_η,ℱ)$,
il suffit de montrer que pour chaque tel $q$,
les groupes $H^{p-q-1}(V,R^q k_*k^* ℱ)$ sont nuls pour $q>0$.
Fixons $q$. Soit $y$ un point géométrique de $Y$ tel
que la fibre de $R^q k_*k^* ℱ$ en $y$ soit non nulle et $x$ le point
image de $Y$ dans $X$. Le schéma $η×_X X_{(x)}$ se décompose
en un coproduit de spectres de corps $η_α$ ; de même, le produit fibré
$Y_{(y)}×_X η$, dont on considère la cohomologie, est isomorphe au coproduit
des $Y_{(y)}×_{X_{(x)}} η_α$.
D'après \emph{op. cit.} (\ref{18B-majoration Milnor générique}),
ces derniers n'ont de cohomologie qu'en degré $q ≤ \dim Y_{(y)} ≤
\dim(Y)-\dim \sur{\{y\}}$. Il en résulte que la dimension
du support de chacun des sous-faisceaux constructibles de $R^q k_*k^*
ℱ$ est au plus $\dim(Y)-q$. De plus, la dimension caténaire du
morphisme $f$ restreint à un tel support est au plus
$\mathrm{dim.cat}(f)-1$. Il résulte donc de l'hypothèse
de récurrence que les groupes
$H^{p-q-1}(V,R^q k_*k^* ℱ)$ sont nuls lorsque
$p-q-1 ≥ (\dim(Y)-q)+(\mathrm{dim.cat.}(f)-1)$. CQFD.

\section{Dimension cohomologique d'un ouvert du spectre épointé : minoration}
\label{18B-minoration-ouvert-spectre-épointé}

\begin{théorème}
Soient $X$ un schéma intègre strictement local nœthérien de
dimension $d$ et $Ω$ un ouvert non vide du spectre épointé. Alors, pour tout nombre
entier $n>1$ inversible sur $X$, on a
\[
\mathrm{cd}_n(Ω) ≥ d.
\]
\end{théorème}

La démonstration occupe les paragraphes suivants.

\subsection{Construction combinatoire locale}
\label{18B-construction-combinatoire-locale}
\subsubsection{Notations}
Soit $X$ un schéma strictement local nœthérien régulier de
dimension $d ≥ 2$ et soit $t₁,…,t_{d-1},t_d$ un système régulier de
paramètres. Pour des raisons qui apparaîtront ultérieurement, on
note également $π$ l'élément $t_d$. Pour chaque $1 ≤ i ≤ d-1$, on note $H_i$ le
diviseur régulier $V(t_i)$ ; pour $i=d$, on pose
$H_d=V(t₁+\cdots+t_{d-1}-π)$. Enfin on note $U$ l'ouvert
affine $X[π^{-1}]$, $k$ l'immersion ouverte $U ↪ X$
et $j$ l'immersion ouverte $U-⋃_{i=1}^d H_i ↪ U$.
On fixe un entier $n$ inversible sur $X$ et on pose $Λ=𝐙/n$.

\subsubsection{}
\label{18B-cohomologie-via-pureté}
Soit $P$ une partie de $\{1,…,d\}$.
Notons $H_P$ l'intersection $⋂_{p ∈ P} H_p$, et désignons
par $H′_P$ l'intersection $H_P ∩ U$, ouverte dans $X$ et $H_P$
mais fermée dans $U$, et $k_P$ l'immersion ouverte $H′_P ↪ H_P$.
Pour chaque entier $q$, le groupe de cohomologie $H^q(U,Λ_{H′_P})$ est isomorphe
au groupe $H^q(H_P,Λ_{H′_P})$. Comme $H_P-H′_P$ est le diviseur
régulier défini par $π$ dans $H_P$, il résulte de la pureté
cohomologique que $H^q(U,Λ_{H′_P})$ est nul pour $q>1$,
isomorphe à $Λ$ pour $q=0$ et de rang $1$, engendré par
la classe de Kummer de $π$ pour $q=1$. Ceci vaut
également pour $P=∅$, avec la convention évidente
que $H_∅=X$ et $H′_∅=U$.

\subsubsection{}
\label{18B-qis}
Considérons maintenant le quasi-isomorphisme
\[
j_! Λ ⥲ \big( Λ → ⨁_{1 ≤ i ≤ d} Λ_{H′_i} → \cdots → ⨁_{|P|=d-1} Λ_{H′_P}
→ 0\big)
\]
de faisceaux de $Λ$-modules sur $U$,
où le premier terme du complexe de droite est placé en degré $0$.
(On utilise ici le fait que $H′_{P}=∅$ si $|P|=d$
car $H_P$ est alors le point fermé de $X$.)
Les différentielles sont des sommes, avec des signes, de flèches de restriction.
À cette résolution est associée — via la filtration « stupide » —
la suite spectrale
\[
E₁^{p,q}=H^q(U, ⨁_{|P|=p} Λ_{H′_P}) ⇒ H^{p+q}(U,j_! Λ).
\]
D'après les observations du paragraphe précédent,
toute classe de $E₁^{d-1,1}=⨁_{|P|=d-1} H¹(U,Λ_{H′_P})$
n'appartenant pas à l'image de $E₁^{d-2,1}$
survit dans l'aboutissement $H^d(U,j_! Λ)$.

\subsubsection{}
\label{18B-classe-locale-non-nulle}
De même que le faisceau $j_! Λ$ est isomorphe au produit tensoriel
\[
(U-H′₁ ↪ U)_! Λ ⊗ \cdots ⊗ (U - H′_d ↪ U)_! Λ),\]
le complexe quasi-isomorphe à $j_! Λ$ ci-dessus
est isomorphe au produit tensoriel des complexes $(Λ → Λ_{H′_i})$, $1 ≤ i ≤ d$, où la flèche
est l'unité de l'adjonction, isomorphes respectivement aux $(U-H′_i ↪ U)_! Λ$.
De cette observation, jointe à (\ref{18B-cohomologie-via-pureté}),
on en déduit que le complexe $E₁^{∙,1}$
est isomorphe à la troncation naïve $σ_{≤d-1}\big((Λ
\dessusdessous{\Id}{→} Λ)^{⊗d}\big)$
obtenue en remplaçant le $d$-ième terme, isomorphe à $Λ$, par zéro.
Il est bien connu que ce complexe (« de Koszul ») produit tensoriel
est acyclique (avant troncation), cf. par exemple \ega{III₁}{§1.1}.
(L'exactitude résulte également du quasi-isomorphisme ci-dessus, appliqué
à d'autres fermés.)
En particulier, l'image de la différentielle
$E₁^{d-2,1} → E₁^{d-1,1}$ est naturellement le noyau d'une forme linéaire
non nulle (explicite) sur $E₁^{d-1,1}$ et n'est donc pas $E₁^{d-1,1}$ tout entier.
Il existe donc des sommes directes de classes de Kummer de $π$
qui survivent dans $H^d(U,j_! Λ)$.

\subsection{Éclatement et normalisation partielle}

\subsubsection{}
\label{18B-stratégie-éclatement-normalisation-partielle}
Soit maintenant $X=\Spec(R)$ un schéma strictement local nœthérien intègre
de dimension $d ≥ 2$, de point fermé $x$, et soit $Ω$ un ouvert non vide
strict de $X$. Nous allons montrer qu'après éclatement
et « normalisation partielle » l'ouvert $Ω$ est
— localement et « modulo des nilpotents » —
un schéma régulier du type du schéma $U$ considéré ci-dessus.
Ceci permet de produire une classe de cohomologie non nulle de degré $d$
sur $Ω$.

\subsubsection{}Soit $Y=\Ecl_x(X)$. Notons $j$ l'immersion ouverte
de $Y-Y_x$ dans $Y$.
On désigne par $\chap{X}$ le complété du schéma local $X$,
par $Y′$ le produit fibré $Y×_X \chap{X}$ et par $Y′_{\red}$ la réduction de $Y′$.
Notons que le schéma $Y′$ est \emph{excellent} car $\chap{X}$
l'est.

\subsubsection{}Soit $𝒪_Y^{\japonaismath{の}}$ la normalisation de $𝒪_Y$
dans $j_* 𝒪_{Y-Y_x}$. On définit de même
$𝒪_{Y′}^{\japonaismath{の}}$ et $𝒪_{Y′_\red}^{\japonaismath{の}}$.
La $𝒪_Y$-algèbre $𝒪_Y^{\japonaismath{の}}$ est colimite (filtrante) de ses
sous-$𝒪_Y$-algèbres de type fini $𝒜_λ$.

\begin{proposition2}
Le foncteur envoyant une sous-$𝒪_Y$-algèbre $ℬ$
de $j_* 𝒪_{Y-Y_x}$ sur l'image de $(Y′ → Y)^*ℬ$ par l'adjonction
$(Y′ → Y)^*j_*𝒪_{Y-Y_x} → j′_* 𝒪_{Y′ -Y_x}$
induit une \emph{bijection} entre
les sous-$𝒪_Y$-algèbres finies de $j_* 𝒪_{Y-Y_x}$
et les sous-$𝒪_{Y′}$-algèbres finies de $j′_* 𝒪_{Y′-Y_x}$.
De plus, les algèbres $𝒪_Y^{\japonaismath{の}}$ et
$𝒪_{Y′}^{\japonaismath{の}}$ se correspondent par ce foncteur.
\end{proposition2}

\begin{démo}
On se ramène au cas où $Y=\Spec(A)$ et $Y_x=V(t)$.
Il suffit de montrer que si $A$ est une $R$-algèbre intègre
et $A′=A ⊗_R \chap{R}$ alors le morphisme
$A[t^{-1}]/A →  A′[t^{-1}]/A′$ est un isomorphisme
et que les normalisations se correspondent.
Notons que les anneaux $A$ et $A′$ ont même complétion $t$-adique.
Le premier point résulte alors du fait que si $M$ est un $A$-module
dont chaque élément est tué par une puissance de $t$, on a $M ⥲ M ⊗_A A′$.
Enfin, soit $(f′/t^r)^d + a′₁ (f′/t^r)^{d-1}+ \cdots +a′_d =0$
une relation de dépendance intégrale où $f′$ et les $a′_i$
appartiennent à $A′$.
Soit $N$ un entier assez grand. Écrivons $f′ = f+  t^N g′$,
$a′_i=a_i+ t^N b′_i$ où $f$ et les $a_i$ appartiennent à $A$.
La relation précédente devient
$(f/t^r)^d + a₁ (f/t^r)^{d-1}+ \cdots +a_d ∈ A′ ∩ A[t^{-1}]=A$.
Il en résulte que l'élément $f′/t^r=f/t^r + t^{N-r}g′$
est, modulo un élément de $A′$, dans l'image
de $A^{\!\japonaismath{の}}$.
\end{démo}

\subsubsection{}Par excellence, l'algèbre $𝒪_{Y′_\red}^{\japonaismath{の}}$
est finie sur $𝒪_{Y′_{\red }}$.
L'anneau $𝒪_{Y′} ^{\japonaismath{の}}$ en est l'« image inverse »
par la surjection naturelle. 
De ces observations et de la proposition précédente, on déduit
qu'il existe un indice $λ$ tel que, si $Z=\Spec(𝒜_λ)$
et $Z′= Z_{\chap{X}}$, alors $Z′_{\red}$ est intégralement clos
dans $Z′_{\red}-E′_{\red}$. Notons que $Z$ et $Y$ sont isomorphes hors
de $E$.

\subsubsection{}Soit $E$ une composante de dimension $d-1$ de $Y_x$
et soit $e$ un point maximal de $E′_{\red}$ dans $Z′_\red$.
Le localisé en $e$ est un anneau de valuation discrète : c'est un
anneau local réduit de dimension $1$, intégralement clos dans le
complémentaire du point fermé. Par excellence de $Z′$
il existe un ouvert dense de $E′_\red$ le long duquel $Z′_\red$
est régulier. On peut également supposer que $E′_\red$
est régulier sur cet ouvert. (Pour ce dernier
point il suffit de constater que $E′_\red$ est de type fini sur un corps.)
Soit $U′_\red$ un ouvert de $Z′_\red$
induisant l'ouvert de $E′_\red$ ci-dessus et $U$ un
ouvert de $Z$ induisant l'ouvert correspondant de $E$.
(Le morphisme $Z′ → Z$ est un isomorphisme sur $E$.)
On a $U′ ⊆ U×_Z Z′$. 
Ci-dessous, on s'autorise à rétrécir
les ouverts $U$ et $U′$, sous réserve
qu'ils contiennent tous les points maximaux de $E$.
On suppose de plus que $U ∩ Y_x = U ∩ E$.

\subsubsection{}On note $t$ une équation de $E$ dans $U$ et $π$
une équation de $E′_{\red}$ dans $U′_{\red}$
de sorte qu'il existe une unité $u$ et un entier $e$
tels que l'on ait l'égalité $t=u × π^e$ sur $U′_{\red}$.
L'existence d'un relèvement montre que l'on peut supposer
l'équation $π$ définie sur $U′$. Vérifions que l'\emph{on peut également
supposer $π$ définie sur $U$}. Les schémas $Z′$ et $Z$ ayant même
complétion $t$-adique, il suffit d'observer que si $a$ est une
fonction sur $U′$, on a l'égalité d'idéaux
$(π)=(π+at²)$, du moins lorsque
$1+ uaπ^{2e-1} ∈ 𝐆_m(U′)$, ce que l'on peut supposer quitte à
restreindre $U′$.

\subsubsection{}Soit $Ω$ un ouvert non vide de $X-\{x\}$
(cf. \ref{18B-stratégie-éclatement-normalisation-partielle}).
On note également $Ω$ ses images inverses dans $Y$ et $Z$ ;
elles lui sont isomorphes. Sur un voisinage ouvert des points maximaux
de $E$, l'ouvert $Ω$ coïncide avec le complémentaire de $E$.
Génériquement sur $E$, on a donc $Ω=Z[π^{-1}]$.
Soit $z$ un point fermé de $Z$ appartenant à un tel ouvert
ainsi qu'à l'ouvert $U$. Soient $t₁,…,t_{d-1}$ des fonctions de $𝒪_{Z ′_{\red},z}$
constituant, avec $π$, un système régulier de paramètres. On peut
supposer qu'elles s'étendent à $U′_\red$.
Utilisant à nouveau le fait que le morphisme $Z′ → Z$ est
un isomorphisme au-dessus de $E$, on peut également
supposer qu'elles proviennent de $U$, quitte
à les changer modulo $π$.
Pour chaque $i<d$, considérons l'adhérence schématique $H_i ⊆ Z$
de l'hypersurface $V(t_i)$ dans $U$  et $H_d$
l'adhérence schématique de $V(t₁+\cdots+t_{d-1}-π)$.
On note $H_i′=H_i ∩ Ω$.

\subsubsection{Stratégie}
On va construire une classe non nulle dans
le groupe de cohomologie $H^d(Ω,j_! Λ_{Ω-⋃₁^d H′_i})$,
où $j$ est l'immersion ouverte $Ω-⋃₁^d H′_i ↪ Ω$
et $Λ=𝐙/n$ avec $n$ inversible sur $X$.
\emph{Localement}, ces groupes de cohomologie sont invariants
par passage à la complétion de la base $X$
(et bien sûr au schéma réduit) de sorte que l'on va pouvoir utiliser les
calculs de \ref{18B-construction-combinatoire-locale}.
Il faut cependant prendre garde ici au fait que
l'intersection $⋂₁^d H′_i$ n'est pas nécessairement vide,
contrairement au cas local précédemment étudié : l'analogue
du complexe \ref{18B-qis} a donc un terme de plus (en degré $d$).
Malgré tout, on va relever à $Ω$ une classe de degré $d$
« locale » — c'est-à-dire du schéma $Ω×_Z Z_{(\sur{z})}$
(ou plutôt l'analogue sur $Z′_\red$) — à coefficients dans $j_!Λ$.

\subsubsection{}Considérons à nouveau la suite spectrale du
\ref{18B-qis} :
\[
E₁^{p,q}=H^q(Ω, Λ_p) ⇒ H^{p+q}(Ω,j_! Λ),
\]
où l'on note $Λ_p$ la somme directe des $Λ_{H′_P}$ avec $|P|=p$.
Pour chaque $P ⊆ [1,d]$ de cardinal $d-1$, l'intersection $H_P$
des hypersurfaces correspondantes de $Z$ est propre sur $X$.
Par construction, elle est aussi \emph{quasi-finie} au voisinage
du point fermé $z$. Il en résulte par le théorème de changement de
base propre pour le $π₀$ (ou le théorème principal de Zariski)
— d'après lequel sa décomposition en
composantes connexes se lit sur la fibre spéciale —
que chaque $H_P$ se décompose en le coproduit d'un schéma local
fini sur $X$ et d'un schéma ne rencontrant pas $z$.
Ainsi, $H^1(Ω,Λ_{d-1})$ a un facteur direct
isomorphe à $H¹(Z_{(z)}[π^{-1}],Λ_{d-1})$
et, par comparaison à la complétion, à $H¹({Z′_\red}_{(z)}[π^{-1}],Λ_{d-1})$.
De plus, la différentielle $E₁^{d-1,1} → E₁^{d,1}$
envoie ce facteur sur $0$ dans $H¹(Ω,Λ_d)$ car
l'intersection de $H_{[1,d]}$ avec $Ω$ est vide au voisinage de $z$.
En effet, si $t₁=t₂=\cdots=t_{d-1}=t₁+\cdots+t_{d-1}-π=0$
alors $π=0$ ; d'autre part, au voisinage de $z$, $Ω=Z[π^{-1}]$.
Ainsi, toute classe de cohomologie du facteur
direct « local » induit une classe dans $H^d(Ω,j_!Λ)$
qui relève la classe correspondante dans
$H^d({Z′_\red}_{(z)}[π^{-1}],j_!Λ)$.
On a vu (\ref{18B-classe-locale-non-nulle}) qu'il existe de telles
classes non nulles.
En conséquence, $H^d(Ω,j_! Λ) ≠ 0$ et
finalement $\cd_n(Ω) ≥ d$. CQFD.

\ifx\danslelivre\undefined
\bibliography{../fichiers-configuration-latex/bibliographie}
\bibliographystyle{../fichiers-configuration-latex/style-GTG}


\author{Yves Laszlo}

\title{XIX. Un contre-exemple}
\date{exposés oraux des  11 et 18 juin 2008}

\setcounter{tocdepth}{3}

\begin{abstract}L'exposé est destiné à construire, suivant Gabber, un exemple d'immersion ouverte
$j:\ U\rightarrow X$ de schémas noethériens telle que $R^1j_*\ZZ/2\ZZ$ ne soit pas constructible.
Ceci donne un contre-exemple à une version du théorème de constructibilité de Gabber
(\refext{XIII}{13-constructibilité}) sans hypothèse de quasi-excellence. Cette construction donne
aussi un exemple de diviseur $D$, dont $U$ est le complémentaire en l'occurrence, dans une surface
régulière localement à croisements normaux (au sens de de Jong) mais pas globalement (\ref{19-dcn}).
\end{abstract}

\maketitle

\tableofcontents

\else
\setboolean{english}{false}
\setboolean{francais}{true}
\selectlanguage{francais}
\chapter{Exemples et contre-exemples}
\addtocontents{toc}{par Yves Laszlo}
\begin{center}
Yves Laszlo
\end{center}
\fi

\section{Introduction}

L'exposé est destiné à construire, suivant Gabber, un exemple d'immersion ouverte $j:\ U\rightarrow
X$ de schémas noethériens telle que $R^1j_*\ZZ/2\ZZ$ ne soit pas constructible. Ceci montre que
l'hypothèse de quasi-excellence du théorème de constructibilité de Gabber
(\refext{XIII}{13-constructibilité}) est indispensable. D'un point de vue géométrique, la
construction est intéressante : $U$ est le complémentaire d'un diviseur $D$ dans une surface
régulière $X$  mais possède une infinité de points doubles ordinaires ; en particulier, son lieu
régulier n'est pas ouvert ce qui lui interdit d'être quasi-excellent. Ce diviseur est un exemple de
diviseur dans une surface régulière localement à croisements normaux (au sens de de Jong) mais pas
globalement (\ref{19-dcn}).

\section{La construction}\label{19-cex}

On va construire une immersion ouverte $j:U\hookrightarrow X$ de schémas noethériens tel que
$R^1j_*(\ZZ/2\ZZ)$ ne soit pas constructible. En particulier, $X$ n'est pas quasi-excellent d'après
le théorème de constructibilité de Gabber (\refext{XIII}{13-constructibilité}). Cette dernière
propriété est d'ailleurs évidente par construction puisque $X$ contient un fermé intègre de
dimension $1$ ayant une infinités de points doubles ! Suivant \textbf{SGA 4}, on note ici
$k\{\underline{x}\}$ l'hensélisé à l'origine d'un anneau de polynômes $k[\underline{x}]$. On choisit
un corps parfait $k$ de caractéristique différente de $2$, au plus dénombrable, qui a la propriété
suivante : toute extension finie de $k$ admet un élément qui n'est pas un carré. On note $\bar k$
une clôture algébrique de $k$.  Par exemple, on peut prendre pour $k$ un corps fini ou un corps de
nombres. Dans la suite, on note $\Lambda=\ZZ/2\ZZ$.

 On commence par regarder le plan $\mathbf{A}^2=\Spec(k[x,y])$ privé des points génériques des
courbes  irréductibles ne coupant pas la droite $\Delta=\Spec(k[x])$ d'équation $y=0$. Ces courbes
sont exactement les courbes irréductibles d'équation $u(1+yg(x,y)),u\in k^*$. On pose donc
$$A_0=(1+yk[x,y])^{-1}k[x,y].$$ Le morphisme de localisation $k[x,y]\ra A_0$ identifie $\Spec(A_0)$
au sous-ensemble du plan $\mathbf{A}^2=\Spec(k[x,y])$ cherché. Les points de $\Spec(A_0)$ sont de
trois sortes \begin{itemize} \item Le point générique de $\mathbf{A}^2$ ; \item Les points
génériques des courbes irréductibles du plan qui rencontrent $\Delta$ ; \item Les points de
$\Spec(A_0)$ fermés dans $\mathbf{A}^2$ (qui sont les points fermés de $\Delta$ comme on va le
voir).  \end{itemize}

Notons qu'un point générique d'une courbe $C$ qui coupe $\Delta$ se spécialise dans $\Spec(A_0)$ sur
n'importe quel point fermé de $C\cap \Delta$ et donc n'est pas fermé dans $\Spec(A_0)$.

Par ailleurs, un point de $\Specmax(A_0)$ est donc défini par $(\bar x,\bar y)\in\bar k$. Si $\bar
y$ est non nul, étant algébrique sur $k$, son inverse est dans $k[\bar y]$ ce qui entraîne $(\bar
x,\bar y)\not\in\Spec(A_0)(\bar k)$.  L'immersion fermée
$$\Spec(k[x])=\Spec(A_0/yA_0)\hookrightarrow\Spec(A_0)$$ induit donc un homéomorphisme
$\Specmax(k[x])\iso\Specmax(A_0)$.

Si $\xi\in\Specmax(A_0)$, on note $\pi_\xi\in k[x]$ le générateur unitaire des polynômes nuls en
$\xi$ et on choisit une racine de $\pi_\xi$ dans $\bar k$ définissant un point géométrique $\bar\xi$
au dessus de $\xi$. On le voit comme un élément de $A_0$ via le plongement tautologique
$k[x]\hookrightarrow A_0$. Le couple $(\pi_\xi,y)$ est un système de coordonnées locales de $A_0$ en
$\xi$, \emph{i.e.} on a un isomorphisme\footnote{\label{19-abus-coord}On devrait plutôt dire que le
morphisme $k[X,Y]_{(0,0)}\ra A_{0,\xi}$ qui envoie $X$ sur $\pi_\xi$ et $Y$ sur $y$ induit un unique
isomorphisme $k(\xi)\{\pi_\xi,y\}\iso A_{0,\xi}^{\h}$.} \begin{equation}\label{19-isom-coord}
k(\xi)\{\pi_\xi,y\}\iso A_{0,\xi}^{\h}.  \end{equation}

Comme $k$ est dénombrable, $\Specmax(A_0)$ est infini dénombrable comme $\Specmax(k[x])$. Soit
$[i],i\geq 0$ la suite des idéaux maximaux de $A_0$.

\begin{lemme}\label{19-lemme-la-constr} Il existe une tour d'extensions quadratiques $$A_n\subset
A_{n+1}=A_n[\sqrt{f_n}],\ n\geq 0$$ et une suite de points $\xi_n\in\Specmax(A_n)$ d'images
$\zeta_n=\xi_n\cap A_0$ dans $\Spec(A_0)$ telle que \begin{itemize} \item Pour tout $m\leq n$ et
tout $\m\in\Spec(A_n)$ au dessus de $[m]$, l'image de $f_n$ dans $k(\m)$ n'est pas un carré.  \item
Pour tout $m< n$ et tout $\m\in\Spec(A_n)$ au dessus de $\xi_m$, l'image de $f_n$ dans $k(\m)$ n'est
pas un carré.  \item L'extension $A_n/A_0$ est étale en $\xi_n$ et $$f_n\equiv\pi_n^2+y\mod\xi_n^3$$
où $\pi_n=\pi_{\zeta_n}$.  \end{itemize} \end{lemme}

\begin{proof} Soit $n\geq 0$. Supposons construits $A_m,f_m,\xi_m,\zeta_m,0\leq m<n$ vérifiant les
propriétés précédentes et construisons $f_n\in A_n$ et $\xi_n\in\Spec(A_n)$. Si $n>0$, on pose
$A_n=A_{n-1}[\sqrt{f_{n-1}}]$ (si $n=0$ l'anneau $A_n$ est déjà construit).

Observons tout de suite que, la caractéristique de $k$ étant différente de $2$, l'extension
$A_n/A_0$ est génériquement étale.  Comme elle est finie et que $\Specmax(A_0)=\Specmax(k[x])$, le
lieu $L_1$ des points de $\Specmax(A_n)$ où $A_n/A_0$ n'est pas étale est de cardinal fini et de
complémentaire $L_1^c$ infini. De même, l'ensemble des points $L_2$ de $\Spec(A_n)$ au dessus des
$[m],m\leq n$ et des $\zeta_m,m<n$ est fini.

Choisissons alors $\xi_n\in L_1^c-L_2$ et notons $\zeta_n=\xi_n\cap A_0$ son image dans
$\Spec(A_0)$. Par construction, on a $$\zeta_n\not\in\{[m],m\leq n,\zeta_m,m<n\}\text{ et
}A_n/A_0\text{ est étale en }\xi_n.$$

Soit $\m\in L_2$. Comme $k(\m)$ est une extension finie de $k$ (théorème des zéros), on peut choisir
$a_\m\in k(\m)$ qui n'est pas un carré. Comme $\xi_n\not\in L_2$, les idéaux de $L_2\cup\{\xi_n^3\}$
sont deux à deux étrangers. Le lemme chinois assure qu'il existe $f_n\in A$ tel que $$f_n\equiv
a_\m\mod\m\text{ pour tout }\m\in L_2\text{ et }f_n\equiv \pi_n^2+y\mod\xi_n^3.$$

De tels $f_n,\xi_n,\zeta_n$ conviennent.\end{proof}

On pose $A=\limind\ A_n$. Soit $n\in\NN$. Comme l'extension d'anneaux intègres $A_n\hookrightarrow
A$ est entière, on peut choisir un point géométrique $\overline{\xi_{n,\infty}}$ de $\Spec(A)$ au
dessus de $\xi_n$. Il définit donc des points géométriques $\overline{\xi_n}$ (resp.
$\overline{\zeta_n}$) au dessus de $\xi_n$ (resp. $\zeta_n$).

Par construction, l'inclusion $A_{n+1}\hookrightarrow A$ définit un isomorphisme (cf. la
note~\ref{19-abus-coord}) $$\bar k\{\pi_n,y,z\}/(z^2-f_n)\iso A^{\hs}_{\overline{\xi_{n,\infty}}}$$
compatible à (\ref{19-isom-coord}), \emph{i.e.} tel que le diagramme $$\xymatrix{\bar
k\{\pi_n,y,z\}/(z^2-f_n)\ar[r]^-{\thicksim}& A^{\hs}_{\overline{\xi_{n,\infty}}}\\ \bar
k\{\pi_n,y\}\ar[r]^-{\thicksim}\ar[u]& A_{0,\overline{\zeta_n}}^{\hs}\ar[u] }$$ commute (rappelons
que l'hensélisation stricte commute aux limites inductives filtrantes, cf. \ega{IV}{4}{18.8.18}).

Avec ces préparatifs, on peut énoncer le résultat principal.

\begin{proposition}\label{19-prop-cexgabber} Soit $j$  l'immersion ouverte
$\Spec(A[1/y])\hookrightarrow \Spec(A)$ et $\eta$ le point générique de $D=V(y)$.  \begin{itemize}
\item $A$ est noethérien.  \item La dimension (sur $\FF_2$) de
$(R^1j_*\Lambda)_{\overline{\xi_{n,\infty}}}$ est $2$, alors que la dimension de
$(R^1j_*\Lambda)_{\bar\eta}$  est $1$.  \item En particulier, $R^1j_*\Lambda$ n'est pas
constructible.  \end{itemize} \end{proposition}

\begin{remarque}\label{19-rem-cexgabber}Notons que le diviseur (intègre)  $D=V(y)$ de la surface
régulière $\Spec(A)$ admet chaque $\xi_{n,\infty}$ comme point double (ordinaire). Il n'est donc pas
quasi-excellent puisque son lieu régulier (ou normal, c'est la même chose ici) n'est pas ouvert.
Quitte à éventuellement localiser $A$, on obtient alors un contre-exemple à la constructibilité avec
un schéma ambiant régulier (mais certes pas excellent) ! \end{remarque}

Le point \textbf{iii} découle immédiatement des points 1) et 2). Passons à la preuve des deux
premiers points.  Le reste de l'exposé est destiné à prouver les points \textbf{i} et \textbf{ii},
seuls points restant à montrer.

\section{Noethérianité de $A$} On va adapter (cf. proposition~\ref{19-crit-noeth}) à la situation
(en l'utilisant)  le critère usuel de noethérianité des limites inductives  que l'on rappelle :

\begin{theoreme}[\egazero{III}{10.3.1.3}]\label{19-ega-noeth}Soit $(A_i,\m_i)$ un système inductif
filtrant d'anneaux locaux noethériens. On suppose que tous les $A_i$ noethériens et que les
morphismes de transitions sont locaux et plats. Alors, si pour tout $i\leq j$, on a $\m_iA_j=\m_i$,
alors $\limind A_i$ est noethérien.\end{theoreme}

On utilisera sans le rappeler ensuite le critère de noethérianité de Cohen (\cite[3.4]{Nag62}) :

\begin{proposition}[Cohen]\label{19-copremier} Un anneau est noethérien si et seulement si tout
idéal premier est de type fini.\end{proposition}

Soit $A_i,i\geq 0$ un système inductif d'anneaux et $A_\infty=\limind A_i$. On suppose que les
morphismes $A_i\ra A_{i+1}$ sont finis et injectifs et que chaque $A_i$ est noethérien (ou, ce qui
revient au même, que $A_0$ est noethérien). En particulier, $\Spec(A_{i+1})\ra\Spec(A_i)$ est fini
et surjectif et  $\Spec(A_\infty)\ra\Spec(A_i)$ est entier et surjectif pour tout $i$ ce qu'on
utilisera sans plus de précaution. Leurs  fibres sont de dimension nulle.  Pour $\p\in\Spec(A_0)$,
on note $\tilde\star_\p$ la propriété

\textbf{Propriété $\tilde\star_\p$:}\emph{ Il existe $i$ tel que pour tout $j\geq i$ et tout
$\q\in\Spec(A_j)$ au dessus de $\p$, l'idéal $\q A_\infty$ est premier.}

\begin{proposition}\label{19-crit-noeth0} $A_\infty$ est noethérien si et seulement si  $A_0$ est
noethérien et tout idéal premier $\p$ de $A_0$ vérifie la
propriété $\tilde\star_\p$.
\end{proposition}

\begin{proof}On note $f:\Spec(A_\infty)\ra\Spec(A_0)$. Prouvons le premier point.

\underline{Suffisance}. Soit $\q_\infty\in\Spec(A_\infty)$ et $\p$ son image dans $\Spec(A_0)$.
Montrons que $\q_\infty$ est de type fini. Choisissons $i$ comme dans $\tilde\star_\p$ et soit
$\q=\q_\infty\cap A_i$. On a d'une part $\q A_\infty\subset\q_\infty$ et, d'autre part
$$\q_\infty\cap A_i\subset\q A_\infty\subset\q_\infty$$ ce qui assure l'égalité $\p=\q_\infty\cap
A_0=\q A_\infty\cap A_0$ de sorte que $\q A_\infty$ se spécialise sur  $\q_\infty$ dans $f^{-1}(\p)$
qui est de dimension $0$. On a donc $\q_\infty=\q A_\infty$ ce qui prouve que $\q_\infty$ est de
type fini comme $\q$ et on invoque~\ref{19-copremier}.

\underline{Nécessité}. Supposons $A_\infty$ noethérien et soit $\p\in\Spec(A_0)$. La fibre
$$f^{-1}(\p)=\Spec(A_\infty\otimes_{A_0}\kappa(\p))$$ est noethérienne de dimension nulle, donc de
cardinal fini. Comme $A_\infty$ est noethérien, on peut donc supposer que tous les idéaux premiers
de $f^{-1}(\p)$ sont engendrés par des éléments de $A_i$ pour $i$ convenable. Soit alors
$\q\in\Spec(A_j), j\geq i$ au dessus de $\p$.  Soit $\q'\in\Spec(A_\infty)$ au dessus de $\q$. Comme
$\q'$ est engendré par $\q'\cap A_i$, il l'est par $\q'\cap A_j=\q$, de sorte $\q A_\infty=\q'$ qui
est donc premier.\end{proof}

\begin{proposition}\label{19-crit-noeth} On garde les hypothèses et les notations de
\ref{19-crit-noeth0}. Si de plus les extensions $A_{i+1}/A_{i}$ sont \emph{plates}, $A_\infty$ est
noethérien si et seulement si tout idéal \emph{maximal} $\m$ de $A_0$ vérifie la propriété
$\tilde\star_\m$.\end{proposition}

\begin{proof}La nécessité découle du premier point.  D'après~\ref{19-crit-noeth0}, il suffit de
prouver la suffisance. Soit donc $\p\in\Spec(A_0)$ et montrons que $\p$ vérifie $\star_\p$.

\begin{lemme}Sous les conditions de la proposition, la propriété $\star_\p$ est équivalente à la
\textbf{Propriété $\star_\p$ : } \emph{Il existe $i$ tel que pour tout $l\geq j\geq i$ et tout
$\q\in\Spec(A_j)$ au dessus de $\p$, l'idéal $\q A_{l}$ est premier.} \end{lemme} \begin{proof}
Supposons $\star_\p$ vérifiée. Soit alors $\q\in\Spec(A_j)$ au dessus de $\p$.  On déduit déjà
$1\not\in\q A_\infty$. De plus, si $xy\in\q A_\infty$, il existe $l\geq j$, tel que $x,y\in A_l$.
Donc, $xy\in\q A_\infty\cap A_l=\q A_l$ (fidèle platitude de $A_\infty/A_l$) et donc par exemple
$x\in\q A_l\subset\q A_\infty$. On a donc $\star_\p\Rightarrow\tilde\star_\p$. L'autre implication
est tautologique.\end{proof}

La clef est de constater que la condition $\star_\p$ ne dépend que des fibres schématiques de $f_i :
\Spec(A_i)\ra\Spec(A_0)$ et donc est \emph{invariante par changement de base}, ce qui va permettre
de se ramener au cas local pour appliquer (\ref{19-ega-noeth}). Précisons.

\begin{lemme}\label{19-star-fib} Soit $\p\in\Spec(A_0)$. Les deux propriétés suivantes sont
équivalentes.  \begin{itemize} \item  La propriété $\star_\p$ est satisfaite.  \item Il existe $i$
tel que pour tout $l\geq j\geq i$ le morphisme induit $$\phi : f_l^{-1}(\p)\ra f_{j}^{-1}(\p)$$
entre les fibres schématiques soit bijectif à fibres réduites.  \end{itemize} \end{lemme}
\begin{proof} Supposons  $\star_\p$  vérifiée et choisissons $i\leq j\leq l$ comme dans $\star_\p$.
Soit $\q\in f_j^{-1}(\p)$ et posons $A=A_j/\q$ et $B=A_l/\q A_l$.  La fibre schématique
$\phi^{-1}(\q)$ est le $k(\q)=\Frac(A)$-schéma $\Spec(B\otimes_A \Frac(A))$. Comme $\q A_l$ est
premier, $B$ est intègre et donc $B\otimes_A \Frac(A)$ également (la tensorisation par $\Frac(A)$
est une localisation). Comme $A_l/A_j$ est finie, l'extension $B/A$ est finie de sorte que
$B\otimes_A \Frac(A)$ est à la fois de dimension finie sur $\Frac(A)$ et intègre, donc c'est un
corps, ce qui entraîne la seconde condition.

Inversement, supposons que $\phi^{-1}(\q)$ soit le spectre d'un corps. Autrement dit,  on suppose
que $B\otimes_A\text{Frac}(A)$ est intègre et on veut montrer que $B$ est intègre. Mais on a
\begin{sslemme} Soient $\phi : A\ra B$ un morphisme plat d'anneaux noethériens. Supposons $A$
intègre. Alors, $B$ est intègre si et seulement si la fibre générique $B\otimes_A\Frac(A)$ est
intègre.\end{sslemme} \begin{proof} Seul la suffisance demande explication.  Comme $\phi$ est plat,
il envoie idéaux associés (resp. minimaux) sur l'unique point associé de $\Spec(A)$, son point
générique ($A$ est intègre). Comme la fibre générique $B\otimes_A\text{Frac}(A)$ est un localisé,
les points associés de la fibre générique s'identifient par localisation aux points associés de
$\Spec(B)$. Comme cette fibre générique  est intègre, $\Spec(B)$ n'a qu'un point générique $\eta$
(ce qui assure l'irréductibilité de $\Spec(B)$), dont l'anneau local $B_\eta$ est un localisé de
$B\otimes_A\text{Frac}(A)$. Ainsi $B$ est sans point immergé et génériquement réduit ($S_1$ et
$R_0$) donc réduit d'où finalement intègre.  \end{proof} \end{proof} \medskip

Choisissons $\m$ maximal dans $A_0$ contenant $\p$ et choisissons $i$ comme dans $\star_\m$. Soit
$\m_i\in\Specmax(A_i)$ au dessus de $\m$ (($\Spec(A_i)\ra\Spec(A_0)$ fini et surjectif) et relevons
$\p\subset\m$ en $\p_i\subset\m_i$ dans $A_i$ (platitude).

Par construction, l'idéal $\m_iA_j$ est premier pour tout $j\geq i$. Mais on a $\m_iA_j\cap
A_i=\m_i$ ($A_j/A_i$ est fidèlement plate) de sorte que $\m_iA_j$ est maximal ($A_j/A_i$ est finie).
Soit alors $A'=(A_i)_{\m_i}$ et changeons de base par $A_0\ra A'$. La fibre spéciale du morphisme
fini $$\Spec(A'_j)\ra\Spec(A'_i),j\geq i$$ est le spectre du \emph{corps} $A_j/\m_iA_j=A_j/\m_j$ de
sorte que $A'_j,i\geq j$ est \emph{local} et que $A'_i\ra A'_j$ est local. Comme les fibres de
$\Spec(A'_j)\ra \Spec(A'_0)$ sont des fibres de $\Spec(A_j)\ra\Spec(A_0)$, le système inductif
$(A'_j)$ vérifie $\star_{\m'}$ pour tout $\m'\in\Specmax(A'_0)$ (\ref{19-star-fib}) de sorte que
$\limind A'_j$ est noethérien (\ref{19-ega-noeth}). D'après la proposition~\ref{19-crit-noeth0}, la
condition $\star_{\p'}$ est vérifiée pour tout $\p'\in\Spec(A'_0)$. Comme la condition $\star$ ne
dépend que des fibres (\ref{19-star-fib}), on déduit que $\star_{\q}$ est vérifiée pour tout $\q$
dans l'image de $\Spec(A'_i)\ra\Spec(A_0)$. Or, par construction $\p$ est dans l'image, ce qui
termine la preuve du point \textbf{i} de la proposition \ref{19-prop-cexgabber}.\end{proof} \medskip

\section{Étude des points doubles} Reste à prouver le point \textbf{ii} de la
proposition~\ref{19-prop-cexgabber}.

Comme au voisinage de $\eta$ l'immersion fermée $D=V(y)\hookrightarrow\Spec(A)$ est une immersion
d'un diviseur régulier dans une schéma régulier, le théorème de pureté
(\refext{XVI}{16-corollaire-purete-absolue-dcn}) assure que la dimension de
$(R^1j_*\Lambda)_{\bar\eta}$ est $1$.

Pour alléger les notations, on pose $x=\pi_n, f=f_n$ et $R=\bar k\{x,y,z\}/(z^2-f)$ et on rappelle
la congruence $$f=x^2+y+O(3)$$ où $O(n)$ désigne un polynôme de degré $\geq n$. On doit démontrer
l'énoncé suivant

\begin{proposition} La dimension de $H^1(\Spec(R[1/y],\Lambda)$ est $2$.\end{proposition}

D'après le théorème de pureté,\refext{XVI}{16-corollaire-purete-absolue-dcn}, la proposition est
conséquence du lemme géométrique bien connu suivant

\begin{lemme}[Lemme de Morse]\label{19-Morse} Il existe un système de coordonnée $X,y,Z$  de $\bar
k\{x,y,z\}$ tel que  $z^2-f=Z^2-X^2-y$.\end{lemme}

\begin{proof} On a $$z^2-f=z^2-x^2-y+O(3).$$ Posons $x_1=\sqrt{-1}x,x_2$. Dans ces coordonnées, on a
$$F:=z^2-f+y=\sum_{i,j}g_{i,j}x_ix_j$$ où $(g_{i,j})\in M_2(R)$ est une matrice symétrique de taille
2 à coefficients dans $R$. On a $$g=(g_{i,j})=\Id+O(1).$$ Considérons alors
$M=\mathbf{A}^{\frac{2(2+1)}{2}}$ paramétrant les matrices carrées symétriques. La différentielle de
$$\mu:\ \left\{ \begin{array}{ccc} M & \rightarrow& M \\ P &\mapsto & {}^tPP \\ \end{array}
\right.$$  en $P=\Id$ est le double de l'application identique, donc $\mu$ est étale en ce point. Le
point $g\in M(R)$ se relève donc en $P\in M(R)$. Autrement dit, on peut écrire $g={}^tPP$.  Le
morphisme ${}^t(x_i)\mapsto P{}^t(x_i)={}^t(y,X_1,Z_1)$ est étale et définit donc des coordonnées de
$\bar k\{x,y,z\}$. Dans ces coordonnées, on a $$Z^2-f+y=X_1^2+Z_1^2$$ et il suffit de poser
$X=\sqrt{-1}X_1, Z=Z_1$ pour conclure.\end{proof}

\medskip Rappelons que la caractéristique de $k$ est différente de $2$.  Ceci termine la preuve de
la proposition~\ref{19-prop-cexgabber}.

\section{$D$ est localement mais pas globalement un diviseur à croisements
normaux}\label{19-sec-dcn}

Commençons par une définition. Dans cette section $D$ désigne un diviseur effectif d'un schéma
régulier $X$ et $j:\ U=X-D\hookrightarrow X$ l'immersion ouverte du complémentaire.

\begin{definition} On conserve les notations précédentes.  \begin{itemize} \item On dit que $D$ est
localement un diviseur à croisements normaux (en abrégé, \emph{localement dcn}) si pour tout $x\in
D$, le localisé de Zariski $\Spec(\OO_{D,x})$ est un diviseur à croisements normaux de
$\Spec(\OO_{X,x})$.  \item Supposons $D$ localement dcn. On note $\varepsilon(x)$ le nombre de
branches analytiques de $\Spec(\OO_{D,x})$ et $\zeta(x)$ son nombre de composantes irréductibles. La
fonction $\varepsilon:\ x\mapsto\varepsilon(x)$ (resp. $\zeta:\ x\mapsto\zeta(x)$) est appelée
fonction de comptage analytique (resp. fonction de comptage Zariski).  \end{itemize}
\end{definition}

Avec les notations précédentes, si $\bar x$ est un point géométrique au dessus de $x\in D$ avec $D$
localement à croisements normaux, l'hensélisé strict $D_{(\bar x)}$ est un diviseur à croisements
normaux stricts de  $D_{(\bar x)}$. On a alors la caractérisation suivante :

\begin{lemme}\label{19-carac-fonction-dcn}  Avec les notations précédentes, supposons de plus que
$D$ est localement dcn et $\Lambda=\ZZ/\ell\ZZ$ avec $\ell$ un nombre premier inversible sur $X$.
Alors, les propositions suivantes sont équivalentes.  \begin{itemize} \item $R^1j_*\Lambda$ est
constructibles ; \item $Rj_*\Lambda$ est constructible ; \item la fonction de comptage analytique
$\varepsilon$ est constructible.  \end{itemize} \end{lemme}

\begin{proof}D'après le théorème de pureté (\refext{XVI}{16-corollaire-purete-absolue-dcn}), la
fibre $(Rj_*\Lambda)_{\bar x}$ est l'algèbre extérieure sur $$(R^1j_*\Lambda)_{\bar
x}=\Lambda^{\varepsilon(x)}.$$ Le lemme en découle immédiatement grâce à la caractérisation des
faisceaux constructibles à fibres finies (\sga{4}{IX}{prop. 2.13 (iii)}).\end{proof}

L'intérêt de ce lemme réside dans la proposition  suivante.

\begin{proposition} Avec les notations précédentes, supposons de plus que $D$ est localement dcn.
Alors, $\varepsilon$ est constructible si et seulement si $D$ est un diviseur à croisement
normaux.\end{proposition}

\begin{proof} Grâce au théorème de pureté (\refext{XVI}{16-corollaire-purete-absolue-dcn}) et au
lemme précédent, il suffit de prouver la partie directe. Supposons donc $\varepsilon$ constructible
et montrons que $D$ est à croisements normaux. Soit $\bar x$ un point géométrique au dessus de $x\in
D$. Puisque $D_{(\bar x)}$ est un diviseur à croisements normaux stricts, il existe un voisinage
étale $\pi:\ X'\rightarrow X$ de $x$ dans $X$, tel que le diviseur $D'=\pi^{-1}(D)$ est la somme de
diviseurs $D'_i$ qui sont réguliers en $x$ et qui se coupent transversalement en $x'$, image de
$\bar x$ dans $X'$. La fonction de comptage analytique $\varepsilon'$ de $D'$ est la somme des
fonctions de comptage analytiques $\varepsilon'_i$. Comme $\varepsilon'$   ne dépend que de
l'hensélisé strict, on a donc $$\varepsilon'=\varepsilon\circ\pi=\sum\varepsilon'_i.$$ En
particulier, $\varepsilon'$ est constructible comme $\varepsilon$. La fonction de comptage Zariski
$\zeta'$ de $D'$ certainement constructible de sorte que la différence $\varepsilon'-\zeta'$ l'est
aussi. Par hypothèse, $\varepsilon'-\zeta'$ s'annule sur $\spec{\OO_{X',x'}}$, donc sur l'ensemble
des générisations de $x'$. Comme elle est constructible, elle est nulle sur un voisinage ouvert $U'$
(Zariski) de $x'$. Comme $\epsilon'_i\geq\zeta_i'$, on a $\epsilon'_i=\zeta_i'$ sur $U'$ de sorte
que, quitte à restreindre $U'$, chaque diviseur $D_i$ est régulier sur $U'$. En se restreignant au
localisé strict de chaque point de $U'$, sur lequel on sait que $D'$ est un diviseur à croisements
normaux, on obtient que les $D_i$ se coupent transversalement de sorte que  la restriction de $D'$ à
$U'$ est un diviseur à croisements normaux stricts.\end{proof}
 
\begin{remarque} L'argument précédent appliqué à $\zeta$ assure que si le localisé Zariski de $D$ en
tout point est un diviseur à croisements normaux strict alors $D$ est un diviseur à croisement
normaux strict.\end{remarque}

Avec les notations de la proposition~\ref{19-prop-cexgabber} , on a donc obtenu le résultat suivant.

\begin{corollaire}\label{19-dcn} Le diviseur $D$ de la surface régulière $\spec{A}$ est localement à
croisements normaux mais pas globalement.\end{corollaire}

\ifx\danslelivre\undefined
\bibliography{../fichiers-configuration-latex/bibliographie}
\bibliographystyle{../fichiers-configuration-latex/smfalpha}


\begin{center} XX. Rigidité\\ Yves Laszlo et Alban Moreau \end{center}

\date{exposés oraux des 22 juin 2007 et 11 juin 2008, notes du 29 août 2008}

\setcounter{tocdepth}{2}
\tableofcontents

\else
\setboolean{english}{false}
\setboolean{francais}{true}
\selectlanguage{francais}
\chapter{Rigidité}
\addtocontents{toc}{par Yves Laszlo et Alban Moreau}
\begin{center} Yves Laszlo et Alban Moreau
\end{center}
\fi

\section{Introduction} Le but de cet exposé est de démontrer les deux résultats techniques
\ref{20-prop-rig} (comparaison des torseurs sur l'ouvert complémentaire $\Spec{A}-V(I)$ défini par
un couple hensélien non nécessairement noethérien et l'ouvert correspondant du complété formel le
long de $V(I)$) et \ref{20rigfort} (rigidité de la ramification). Ils permettront dans l'exposé
suivant de montrer l'énoncé de finitude suivant (\refext{XXI}{21-local}):

\begin{theoreme0}[Gabber]   Soit $A$ un anneau strictement local de dimension $2$.  On suppose que
$A$ est normal, excellent, et on note $X'=\Spec(A)-\{\mathfrak{m}_A\}$  son spectre épointé.  Alors,
pour tout groupe fini $G$,   l'ensemble $H^1(X';G)$ est fini.\end{theoreme0}

Ce résultat est la clef pour démontrer le résultat de finitude général suivant
(\refext{XXI}{21-groupes}) :

\begin{theoreme0}[Gabber]    Soit $f:Y \rightarrow X$ un morphisme de type fini   entre schémas
quasi-excellents.  Soit $\mathbb L$ un ensemble de nombres premiers inversibles sur $X$.  Pour tout
faisceau constructible de groupes $F$ sur $Y_\et$    de $\mathbb L$-torsion,  le faisceau $\R^1
f_*(F)$ sur $X_\et$ est constructible. \end{theoreme0}

Par des techniques d'ultrafiltres, chères aux théoriciens des modèles, on est ramené à étudier des
revêtements étales de spectres épointés d'anneaux \textit{non noethériens}, ce qui explique qu'on
soit contraint de démontrer les énoncés techniques hors de tout cadre noethérien.

\begin{remarque0}\label{precau-champ} Soit $X$ un schéma. On considérera des champs en groupoïdes
$\C$ sur $X_\et$ (on dira simplement champs).  En général, la catégorie fibrée $\C$ n'est pas
scindée de sorte que si $x,y$ sont deux objets de $\C(S)$ où $S\ra X$ est étale, il faut quelques
précautions pour parler du faisceau  $\underline{\Hom}(x,y)$ sur $S_\et$. Précisément, suivant
\cite[I.2.6.3.1]{CNA@Giraud}, on considère l'équivalence de catégories fibrées $\C\ra{\lf}\C$ entre
$\C$ et la catégorie libre ${\lf}\C$ engendrée par $\C$, catégorie libre qui elle est scindée. On
définit alors $$\underline{\Hom}(x,y)(S')=\Hom_{{\lf}C(S')}({\lf}x',{\lf}y')$$ où
${\lf}x',{\lf}y'$ sont les images inverses par le morphisme étale $S'\ra S$ de $\lf x,\lf y$ dans
$\lf\C(S')$. Bien entendu (\cite[I.2.6.3.2 (1)]{CNA@Giraud}, $\lf$ induit une bijection
$$\Hom_{\C(S)}(x,y)\iso H^0(S,\underline{\Hom}(x,y)).$$ Ces remarques justifient qu'on puisse
si besoin supposer  sans dommage que les champs que l'on considérera sont scindés.  \end{remarque0}

\section{Lemme de rigidité} Soit $(A,I)$ un couple hensélien (\egalong{IV}{4}{18.5.5}) non
nécessairement noethérien, \textbf{avec $I$ de type fini}\footnote{Cette hypothèse sera utilisée
pour comparer les gradués $I$-adiques de $A$ et de son complété $\chap{A}$ (\ac{III}{2}{12}{})}.
Soit $U$ un ouvert de $X=\Spec(A)$ contenant $\Spec(A)-V(I)$. On note $\chap{A}$ le complété\footnote{On dira
simplement complété pour séparé complété.} $I$-adique de $A$ et $\chap{U}$ l'image inverse de $U$
par le morphisme de complétion $\pi:\chap{X}=\Spec(\chap{A})\ra X$.
On suppose pour simplifier $U$ quasi-compact (cf. \ref{20-U-non-qc}).

\subsection{Énoncés}\label{20-sec-enonc} Rappelons \sga{4}{IX}{1.5} qu'un faisceau en groupes $\F$
sur $X$  est ind-fini si pour tout ouvert étale $u:\ U\rightarrow X$ avec $U$ quasi-compact, le
groupe $\F(u)$ est limite inductives filtrante de ses sous-groupes d'indice fini. On dit alors qu'un
champ en groupoïdes $\C$ sur $X$ est ind-fini si pour tout si pour tout ouvert étale $u:\
U\rightarrow X$ avec $U$ quasi-compact et tout $x_u\in\C(u)$, le faisceau en groupes
$\pi_1(\C,x_u)=\Aut_\C(x_u)$ est ind-fini.

Le but de cette de cette section est de démontrer le théorème de rigidité suivant.

\begin{theoreme2}[Théorème de rigidité de Gabber] \label{20-prop-rig} Soit $\F$ un faisceau
d'ensembles sur $U_\et$. Alors on a \begin{itemize} \item[i)] la flèche naturelle $H^0(U,\F)\ra
H^0(\chap{U},\pi^*\F)$ est bijective ; \item[ii)] si $\F$ est de plus un faisceau en groupes
ind-fini, la flèche naturelle $H^1(U,\F)\ra H^1(\chap{U},\pi^*\F)$ est bijective.  \end{itemize}
\end{theoreme2}

Les deux énoncés du théorème précédent sont conséquence du théorème suivant, apparemment plus fort,
forme champêtre du théorème de rigidité\footnote{Les champs (ind-finis) en groupoïdes discrets
s'identifient aux faisceaux d'ensembles : on dira parfois un \textit{champ discret}.}.

\begin{theoreme2}[Théorème de rigidité de Gabber, forme champêtre] \label{20-prop-rig-champ} Soit
$\C$ un champ en groupoïdes ind-fini sur $U_\et$. Alors, la flèche naturelle
$\gamma(\C):\Gamma(U,\C)\ra \Gamma(\chap{U},\pi^*\C)$ est une équivalence.\end{theoreme2}

\begin{remarque2} En fait, le théorème~de rigidité \ref{20-prop-rig} est \textit{a priori}
équivalent à la version champêtre \ref{20-prop-rig-champ}. C'est ce qui ressort par exemple de
l'énoncé~\ref{20-prop-red-champ}. Mais, formellement, on n'a pas besoin de démontrer cela à ce
stade.\end{remarque2}

\begin{remarque2}
\label{20-U-non-qc}
 Les résultats précédents sont également valables lorsque $U$ n'est
pas nécessairement quasi-compact. Cela résulte du fait que
la catégorie des sections d'un champ sur $U$ est
équivalente à la $2$-limite projective des sections
sur les ouverts quasi-compacts de $U$ contenant $\Spec(A)-V(I)$.
L'hypothèse de quasi-compacité est utilisée dans un argument
d'éclatement ci-dessous (cf. \ref{20-fin-preuvebis}).
\end{remarque2}

\subsection{Réduction au cas constant}

Le résultat est le suivant

\begin{proposition2}\label{20-reduc-const}Supposons que pour tout $U$ comme plus haut,
\begin{itemize} \item[i)] pour tout ensemble fini $F$, la flèche $H^0(U,F)\ra H^0( \chap{U},F)$ est
bijective. Alors, \ref{20-prop-rig} i) est vrai, c'est-à-dire le théorème de rigidité
\ref{20-prop-rig-champ} est vrai pour les champs discrets.  \item[ii)]  pour tout groupe fini $G$,
la flèche $\Tors(U,G)\ra\Tors(\chap{U},G)$ est une équivalence et \ref{20-prop-rig} i) est vrai.
Alors le théorème de rigidité \ref{20-prop-rig-champ} est vrai.  \end{itemize} \end{proposition2}
\begin{proof} D'après \sga{4}{XII}{prop.6.5} (resp. \ref{20-prop-red-champ}), il suffit pour prouver
\ref{20-prop-rig} i) (resp. \ref{20-prop-rig-champ}) de prouver que pour tout $U'\ra U$ fini et tout
ensemble fini $F$ (resp. groupe fini $G$), la flèche \begin{equation}\label{20-lemme-sga}
H^0(U',F)\ra H^0(\chap{U'},F)\textup{ (resp. } \Tors(U',G)\ra\Tors(\chap{U'},G)) \end{equation}

 est bijective (resp. une équivalence) où $\chap{U'}=\chap{U}\times_UU'$.

\begin{lemme2} Il existe un schéma schéma affine $\Spec(B)$ et un diagramme cartésien
$$\xymatrix{U'\ar[r]\ar[d]\ar@{}[rd]|{\square}&\Spec(B)\ar[d]\\U\ar[r]&\Spec(A)}$$ où  $B$ est fini
sur $A$. Le morphisme $U'\ra \Spec(B)$ s'identifie à l'immersion ouverte $U_B\hookrightarrow
\Spec(B)$. De plus, $U_B$ contient $\Spec(B)-V(IB)$.  \end{lemme2}

\begin{proof} Comme $U'\ra U$ est fini, il est projectif (\ega{II}{6.1.11}). Comme $U$ est
quasi-compact, l'immersion ouverte $U\hookrightarrow X$ est quasi-affine (\ega{II}{5.1.1}), donc
quasi-projective de sorte que le composé $f:U'\rightarrow U\ra X$ est quasi-projectif
(\ega{II}{5.3.4}). Comme $X=\Spec(A)$ est affine, $\OO_X$ est certainement ample (cf. la définition
ou \ega{II}{5.1.2}). Les hypothèses du  théorème principal de Zariski (\egalong{IV}{3}{8.12.8}) sont
donc vérifiées. Il existe donc $X'\ra X$ fini de sorte que $f$ se factorise en $U'\hookrightarrow
X'\ra X$ où $U'\hookrightarrow X'$ immersion ouverte et $X'\ra X$ fini. L'adhérence schématique de
$U'$ dans $X'$ est fermée dans $X'$ : elle s'écrit donc $\Spec(B)$ où $B$ est fini sur $A$. On a
donc un diagramme commutatif $$\xymatrix{U'\ar[r]\ar[rd]&U_B\ar[d]\ar@{^(->}[r]&\Spec(B)\ar[d]\\
&U\ar@{^(->}[r]&\Spec(A)}$$ où les flèches non horizontales sont finies. La flèche $U'\ra U_B$ est
donc propre. Comme c'est aussi une immersion ouverte d'image dense, c'est un isomorphisme.  L'ouvert
$U$ contenant $\Spec(A)-V(I)$, on déduit que $U'=U_B$ contient
$\Spec(B)-V(IB)=(\Spec(A)-V(I))_B$.\end{proof}

D'après le lemme, la flèche (\ref{20-lemme-sga}) s'identifie à
\begin{equation}\label{20-lemme-sga-1} H^0(U_B,F)\ra H^0(\chap{U}_B,F)(\textup{ resp. }
\Tors(U_B,G)\ra\Tors(\chap{U}_B,G)) \end{equation}

 (où $?_B$ est l'extension des scalaires du $A$-schéma $?$ à $\Spec(B)$).  Il s'agit donc de montrer
que (\ref{20-lemme-sga-1}) est bijectif (resp. une équivalence).

Par définition, on a $$\chap{U}_B=\pi_C^{-1}(U)$$ où $\pi_C$ est la projection naturelle $$\pi_C:\
\Spec(C)\ra\Spec(A),\text{ avec } C=\chap{A}\otimes_AB.$$

Dans le cas noethérien, $C$ est le complété $IB$-adique $\chap{B}=C $ de $B$ ce qui prouve la
proposition dans ce cas -- appliquer l'hypothèse \ref{20-reduc-const} i) à $\F$ constant de valeur
$F$ sur $U_B$ --. Dans le cas général, la flèche $C\ra\chap{B}$ n'est pas en général un
isomorphisme.

\begin{lemme2}\label{20-lemme-hens}Avec les notations précédentes, on a \begin{itemize} \item [i)]
Soit $(A_n,I_n)$ un système projectif de couples henséliens. Le couple
$(A_\infty,I_\infty)=(\varprojlim{A_n},\varprojlim{I_n})$ est hensélien.  \item [ii)] Le complété
$I$-adique de $A$ est hensélien.  \item [iii)] Les couples $(B,IB)$ et $(C,IC)$ sont henséliens et
ont même complété $I$-adique.  \end{itemize} \end{lemme2} \begin{proof} Soit $P$ un polynôme
unitaire de $A_\infty[x]$ et $\bar a\in A_\infty/I_\infty$ une racine simple (c'est-à-dire  telle
que $P'(\bar a)$ inversible dans $A_\infty/I_\infty$). L'image $\bar a_n$ de $\bar a$ dans $A_n/I_n$
est une racine simple de $P$. Elle se relève donc de façon unique en une racine $a_n\in A_n$ de $P$
d'après le lemme de Hensel. Comme $I_{n+1}$ s'envoie dans $I_n$, par unicité des relèvements,
l'image de $a_{n+1}$ dans $A_n$ est égale à $a_n$ de sorte que la suite $a=(a_n)\in A_\infty$ est le
relèvement cherché de $\bar a$ ce qui prouve i).

Puisque $A\ra A/I^n$ est notoirement entier, les couples $(A/I^n, IA/I^n)$ sont henséliens de sorte
que ii) découle de i).

Par associativité du produit tensoriel, le morphisme naturel $B/I^nB\rightarrow C/I^nC$ s'identifie
au morphisme naturel $A/I^nA\rightarrow\chap{A}/I^n\chap{A}$. Comme ce dernier est un isomorphisme
(\ac{III}{2}{2}{, prop. 15 et cor. 2}), $B$ et $C$ ont même complété $I$-adique. iii) suit alors de
ii) car un couple fini sur un hensélien est hensélien.  \end{proof}

On a donc $\widehat{U_B}=\widehat{U_C}$. D'après le lemme précédent, sous les hypothèses de
\ref{20-reduc-const} i) (resp. ii)), la flèche naturelle $$H^0(U_B,F)\ra H^0(\widehat{U_B},F)=
H^0(\widehat{U_C},F)\leftarrow H^0({U_C},F)=H^0(\chap{U}_B,F)$$

(resp. $$\Tors(U_B,G)\ra \Tors(\widehat{U_B},G)= \Tors(\widehat{U_C},G)\leftarrow
\Tors({U_C},G)=\Tors(\chap{U}_B,G))$$ est alors une bijection (resp. équivalence), ce qu'on
voulait.\end{proof}

\subsection{Réduction au cas strictement hensélien} Résumons les notations dans le diagramme
cartésien suivant

$$\xymatrix{
\chap{U}\ar[r]^{\pi}\ar@{^(->}[d]_{{\chap{j}}}\ar@{}[rd]|\square\ar[r]^{\pi}&U\ar@{^(->}[d]^{j}\\
\chap{X}\ar[r]^{\pi}& X }$$ avec $U$ quasi-compact contenant $\Spec(A)-V(I)$.  Montrons le résultat
suivant.

\begin{proposition2}\label{20-reduchensel} Supposons que pour tout $U$ comme plus haut,
\begin{itemize} \item[i)] pour tout ensemble fini $F$, la flèche $H^0(U,F)\ra H^0(\chap{U} ,F)$ est
bijective si $A$ est de plus strictement local. Alors, \ref{20-prop-rig} i) est vrai (que $A$ soit
strictement local ou non).  \item[ii)]  pour tout groupe fini $G$, la flèche
$\Tors(U,G)\ra\Tors(\chap{U},G)$ est une équivalence si $A$ est de plus strictement hensélien et
\ref{20-prop-rig} i) est vrai. Alors, \ref{20-prop-rig-champ} est vrai (que $A$ soit strictement
local ou non).  \end{itemize} \end{proposition2}

\begin{proof} Commençons par un lemme.

\begin{lemme2}\label{20-cbase1} Supposons que pour tout $U$ comme plus haut, \begin{itemize}
\item[i)] pour tout ensemble fini $F$, la flèche $H^0(U,F)\ra H^0(\chap{U},F)$ est bijective si $A$
est de plus strictement local. Alors, la flèche de changement de base $$\gamma:\pi^*j_*F\ra
{\chap{j}}_*\pi^*F={\chap{j}}_*F$$ est un isomorphisme (que $A$ soit strictement local ou non).
\item[ii)]  pour tout groupe fini $G$, la flèche $\Tors(U,G)\ra\Tors(\chap{U},G)$ est une
équivalence si $A$ est de plus strictement hensélien et \ref{20-prop-rig} i) est vrai. Alors, la
flèche de changement de base $$\gamma:\pi^*j_*\underline{\Tors}(U,G)\ra
{\chap{j}}_*\pi^*\underline{\Tors}(U,G)
={\chap{j}}_*\underline{\Tors}(\chap{U},G),$$
où l'égalité résulte de \cite[III.2.1.5.7]{CNA@Giraud},
est une équivalence (que $A$ soit
strictement local ou non).
\end{itemize} \end{lemme2}
\begin{proof}

Les formules $j^*j_*=\Id\text{ et }\chap{j}^*\chap{j}_*=\Id$ assurent qu'on a
$$\chap{j}^*\pi^*j_*=\pi^*j^*j_*=\pi^*=\chap{j}^*\chap{j}_*\pi^*$$ de sorte que l'image inverse sur
$\chap{U}$ de la flèche de changement de base

\begin{equation}\label{20-CTCU} \pi^*j_*\C\ra {\chap{j}}_*\pi^*\C \end{equation} est une équivalence
pour tout champ en groupoïdes $\C$.

Soit $\hat x$ un point géométrique de $\chap{X}$ d'image le point géométrique $x=\pi\circ\hat x$ de
$X$ et montrons que la fibre de la flèche de changement de base  (\ref{20-CTCU}) en $\hat x$ est une
équivalence. D'après ce qui précède, on peut supposer $\hat x\not\in\chap{U}$. En particulier, $x\in
V(I)$.

Soit $A^{\hs}$ (resp. $X_{(x)}$)  l'hensélisé strict de $A$ (resp. $X$) en $x$ et $\chap{A}^{\hs}$
(resp. $\chap{X}_{(\hat x)}$) celui de $\chap{A}$ (resp. $\chap{X}$) en $\hat x$. On a un diagramme
commutatif où les flèches sont les flèches de fonctorialité, complétion ou stricte hensélisation

$$\xymatrix{\widehat{\chap{X}_{(\hat x)}}\ar[r]\ar[d]&\chap{X}_{(\hat
x)}\ar[r]\ar[d]&\chap{X}\ar[d]\\ {\chap{X}_{(\hat x)}}\ar[r]&X_{(x)}\ar[r]&X}$$

On note alors

$$\xymatrix{\widehat{\chap{U}_{(\hat x)}}\ar[r]\ar[d]&\chap{U}_{(\hat
x)}\ar[r]\ar[d]&\chap{U}\ar[d]\\ {\chap{U}_{(\hat x)}}\ar[r]&U_{(x)}\ar[r]&U}$$ l'image inverse du
diagramme par l'immersion ouverte $U\rightarrow X$. En particulier, $U_{(x)}$ (resp.
$\chap{U}_{(x)}$) désigne l'image inverse de l'hensélisé strict $X_{(x)}$ (resp. $\chap{X}_{(x)}$)
de $X$ (resp. $\chap{X}$) en $x$ (resp. $\hat x$) par $j$ (resp. ${\chap{j}}$). Comme $U$ est
quasi-compact, il en est de même des ouverts $U_{(x)},\chap{U}_{(\hat x)}$ de
$X_{(x)},\chap{X}_{(\hat x)}$ .

Les morphismes $j,{\chap{j}}$ étant cohérents, dans le cas i), la fibre $\gamma_{\hat x}$
s'identifie à la flèche naturelle $$H^0(U_{(x)},F)\ra H^0(\chap{U}_{(\hat x)},F)$$ tandis que dans
le cas ii) elle s'identifie à $$\Tors(U_{(x)},G)\ra \Tors(\chap{U}_{(\hat x)},G).$$

  On déduit que les flèches naturelles $$H^0(U_{(x)},F)\ra H^0(\widehat{ U_{(x)}},F)\text{ et }
H^0(\chap{U}_{(\hat x)},F)\ra H^0(\widehat{\chap{U}_{(\hat{x})}},F)$$ sont bijectives dans le cas i)
et que les flèches $$\Tors(U_{(x)},G)\ra \Tors(\widehat{U_{(x)}},G)\text{ et }\Tors(\chap{U}_{(\hat
x)},G) \ra \Tors(\widehat{\chap{U}_{(\hat{x})}},G)$$ sont des équivalences dans le cas ii). Il
suffit donc de voir que la flèche naturelle \begin{equation}\label{20iso-complet}
\widehat{\chap{U}_{(\hat{x})}} \ra\widehat{ U_{(x)}} \end{equation} est bijective, ou encore que
$$A^{\hs}\text{ et }\chap{A}^{\hs}\text{ ont même $I$-complété.}$$ Puisque l'anneau local
$\chap{A}^{\hs}$ est hensélien, il est a fortiori $I$-hensélien (exercice). Utilisant
(\ref{20-lemme-hens}), on constate que le $I$-complété $\widehat{\chap{A}^{\hs}}$ est hensélien.
Comme son corps résiduel est celui de $A^{\hs}$, il est strictement hensélien. La flèche $\chap{A}
\ra\widehat{A^{\hs}}$ induit donc une flèche $\chap{A}^{\hs}\ra\widehat{A^{\hs}}$ et donc, par
$I$-complétion, une flèche $$\widehat{\chap{A}^{\hs}}\ra\widehat{A^{\hs}}.\leqno (*)$$ Par ailleurs,
la flèche de complétion $A\ra\chap{A}$ induit par hensélisation stricte puis complétion une flèche
$$\widehat{A^{\hs}}\ra\widehat{\chap{A}^{\hs}}.  \leqno (**)$$ Les flèches (*) et (**) sont inverses
l'une de l'autre, d'où le lemme.\end{proof}

On a le diagramme commutatif à carré cartésien

$$\xymatrix{
\chap{U}\ar[r]^{\pi}\ar@{^(->}[d]_{{\chap{j}}}\ar@{}[rd]|\square\ar[r]^{\pi}&U\ar@{^(->}[d]^{j}\\
\chap{X}\ar[r]^{\pi}& X\\ X_I=\Spec(A/I)\ar@{^(->}[u]\ar@{^(->}[ru]& }$$

Comme on l'a observé, les paires $(A,I)$ et $(\chap{A},I\chap{A})$ sont henséliennes.  La flèche
$H^0(\chap{X},\C)\ra H^0(X_I,\C_{|X_I})$ est donc une équivalence pour tout champ ind-fini $\C$ sur
$X_\et$ d'après \cite[théorème 1']{Gab94}.

On déduit d'une part $$H^0(U,F)=H^0(X,j_*F)=H^0(X_I,(j_*F)_{|X_I})$$ et, d'autre part
$$H^0(\chap{U},F)=H^0(\chap{X},{\chap{j}}_*F)\stackrel{\ref{20-cbase1}}{=}H^0(\chap{X},\pi^*
j_*F)=H^0(X_I,(\pi^* j_*F)_{|X_I})$$ ce dernier n'étant autre que $H^0(X_I,(j_*F)_{|X_I})$ (bien
entendu l'isomorphisme induit $$H^0(U,F)\iso H^0(\chap{U},F)$$ est la restriction).

De même, on a
$$H^0(U,{\underline{\Tors}}(U,G))=H^0(X,j_*{\underline{\Tors}}(U,G))=H^0(X_I,j_*{\underline{\Tors}}(U,G)_{|X_I}))$$
et, d'autre part

\begin{eqnarray*} H^0(\chap{U},{\underline{\Tors}}(\chap{U},G))&=&
H^0(\chap{X},{\chap{j}}_*{\underline{\Tors}}(\chap{U},G)) \\ &\stackrel{\ref{20-cbase1}}{=}&
H^0(\chap{X},\pi^* j_*{\underline{\Tors}}(U,G))\\ &=&H^0(X_I,\pi^*
j_*{\underline{\Tors}}(U,G)_{|X_I})) \end{eqnarray*}

ce dernier n'étant autre que $H^0(X_I,j_*{\underline{\Tors}}(U,G)_{|X_I}))$ (bien entendu
l'isomorphisme induisant $$H^0(U,{\underline{\Tors}}(U,G))\iso
H^0(\chap{U},{\underline{\Tors}}(\chap{U},G)).$$ Reste à invoquer \ref{20-reduc-const}.\end{proof}

\subsection{Fin de la preuve de \ref{20-prop-rig-champ}} D'après~\ref{20-reduchensel}, pour prouver
\ref{20-prop-rig-champ}, il suffit de prouver l'énoncé suivant

\begin{proposition2}\label{20-fin-preuve}Supposons $A$ strictement hensélien (et
$I\subset\text{rad}(A)$) et soit $U$ comme plus haut.  \begin{itemize} \item[i)] pour tout ensemble
fini $F$, la flèche $H^0(U,F)\ra H^0(\chap{U},F)$ est bijective.  \item[ii)]  pour tout groupe fini
$G$, la flèche $\Tors(U,G)\ra\Tors(\chap{U},G)$ est une équivalence.\end{itemize} \end{proposition2}

La formule $\pi^*\underline{\Tors}(U,G)=\underline{\Tors}(\chap{U},G)$ (\cite{CNA@Giraud},
III.2.1.5.7) permet de réécrire \ref{20-fin-preuve} sous la forme suivante

\begin{proposition2}\label{20-fin-preuvebis}Supposons $A$ strictement hensélien (et
$I\subset\text{rad}(A)$) et soit $U$ comme plus haut. Désignons par $\C$ le champ discret $F_U$ ou
bien $\underline{\Tors}(U,G)$. Alors, la flèche $H^0(U,\C)\ra H^0(\chap{U},\pi^*\C)$ est une
équivalence.  \end{proposition2}

\begin{proof} On va se ramener par éclatement au cas ou l'idéal $J$ définissant le complémentaire de
$U$ est principal .

Pour tout idéal ${\tilde I}$ d'un anneau ${\tilde A}$, on note $$\Ecl_{{\tilde I}}({\tilde
A})=\Proj(\oplus_{n\geq 0}\tilde I^n)$$ l'éclatement de ${\tilde I}$ dans $\Spec(\tilde A)$.
Si ${\tilde I}$ est de type fini, le morphisme structural $e:\Ecl_{{\tilde I}}({\tilde
A})\ra\Spec(\tilde A)$ est projectif, en particulier propre. Si $S$ est un $A$-schéma, on note
$\chap{S}^J$ le complété de $S$ le long de $J$.

On suppose donc $A$ strictement hensélien de corps résiduel $k$ et $\F=F_U$ comme plus haut. On a
déjà observé que $\chap{A}$ était aussi strictement hensélien. Il suit en particulier que l'ensemble
des sections globales de tout faisceau étale sur $X$ ou $\chap{X}$ s'identifie à sa fibre spéciale,
ce qu'on utilisera sans plus de précaution.

Comme $U$ est quasi-compact, il existe un idéal $J$ de type fini tel que $U=\Spec(A)-V(J)$. Comme
$U$ contient $\Spec{A}-V(I)$ et que $I$ est de type fini, on peut supposer $I\subset J$.  Soit
$$Y=\Ecl_J(A)\text{ et }Y'=\Ecl_J(\chap{A}).$$ (On aurait dû écrire $\Ecl_{J\chap{A}}(\chap{A})$ pour
$\Ecl_J(\chap{A})$).  Pour des raisons de cohérences, on notera simplement $X'$ le complété
$\chap{X}=\Spec(\chap{A})$ (resp. $U'$ sa restriction $\chap{U}=\pi^{-1}(U)$ à $U$).

\begin{sous-lemme2}\label{20-compa-comple}Soient $n,m$ des entiers $\geq 0$. Le morphisme de
complétion définit des isomorphismes $$A/I^mJ^n\simeq \chap{A}/I^mJ^n\chap{A}\text{ et }A/J^n\simeq
\chap{A}/J^n\chap{A}$$ induisant un isomorphisme \begin{equation*} J^n/I^mJ^n\simeq
J^n\chap{A}/I^mJ^n\chap{A}.  \end{equation*} \end{sous-lemme2}

\begin{proof}Comme $I$ est de type fini, le morphisme de complétion induit des isomorphismes
$$A/I^{m+n}\simeq \chap{A}/I^{m+n}\chap{A}\text{ et } A/I^{n}\simeq \chap{A}/I^{n}\chap{A}$$ d'après
\ac{III}{2}{12}{, cor. 2 de la prop. 16}. Mais comme $J$ contient $I$, on a $$I^{m+n}\subset
I^mJ^n\text{ et }I^n\subset J^n,$$ de sorte que les changements de base $$A/I^{m+n}\rightarrow
A/I^mJ^n\text{ et }A/I^{n}\rightarrow A/J^n$$ donnent alors des isomorphismes  $$A/I^mJ^n\simeq
\chap{A}/I^mJ^n\chap{A}\text{ et }A/J^n\simeq \chap{A}/J^n\chap{A}$$ qui donnent
\ref{20-compa-comple}.  \end{proof} La flèche naturelle $Y'\ra Y$ est donc un isomorphisme au dessus
de $\Spec(A/I)\subset X$ car elle est induite par le morphisme gradué $$\oplus J^nA/IJ^{n}A\ra
J^n\chap{A}/IJ^{n}\chap{A}$$ qui est un isomorphisme. On identifiera ces restrictions par la suite.
En particulier, le morphisme $p_s:Y'_s\ra Y_s$ entre fibres spéciales (c'est-à-dire  au dessus du
point fermé de $s\in \Spec(A/I)\subset X$)  est un isomorphisme grâce auquel nous les identifierons.
Regardons le solide commutatif

$$\xymatrix{ \chap{U}=U'\ar@{^(->}[r]^{j'}\ar[dd]^{p}&Y'\ar[dd]^p\ar[rrd]^{e}&Y_s\ar@{=}[dd]
\ar[rrd]\ar@{_(->}[l]_{i'}\ar@{}[ldd]| \square\\
&&&\chap{X}=X'\ar[dd]^\pi&\Spec(k)\ar@{_(->}[l]\ar@{=}[dd]\\
U\ar@{^(->}[r]^{j}&Y\ar[rrd]^-e&Y_s\ar[rrd]\ar@{_(->}[l]_i\\ &&&X&\Spec(k)\ar@{_(->}[l]}$$

Admettons pour un temps le résultat suivant.

\begin{lemme2}\label{20-commut-princ} Soit $\C=F_U$
(resp. $\C=\underline{\Tors}(U,G)$).  Alors, la flèche de changement de base $$\gamma :p^*j_{*}\C\ra
j'_{*}p^*\C$$ est bijective (resp. une équivalence).\end{lemme2}

Déduisons alors l'isomorphisme cherché $$H^0(U,{\C})\iso H^0(U',{\C})=H^0(U',p^*{\C})$$ grâce au
théorème de changement de base propre d'Artin-Grothendieck (\cite{CNA@Giraud} dans le cas noethérien
et théorème \ref{20-th-cgt-base-propre} dans le cas général) appliqué aux faces inférieure et
supérieure du diagramme précédent.  On a en effet un diagramme commutatif où toutes les flèches sont
les flèches naturelles (obtenues par adjonction)

$$\xymatrix{H^0(U,{\C})\ar@{=}[r]\ar[d]^\alpha&H^0(Y,j_{*}{\C})\ar[r]^-{b}\ar[d]^a
&H^0(Y_s,i^*j_{*}{\C})\ar[d]\ar@{=}[r]&H^0(Y_s,i'^*p^*j_*{\C})\ar[ld]^c\\
H^0(U',p^*{\C})\ar@{=}[r]&H^0(Y',j'_{*}p^*{\C})\ar[r]^-{d} &H^0(Y_s,i'^{*}j'_{*}p^*{\C})}$$

Les flèches $b,d$ dont bijectives grâce au théorème de changement de base propre
(\ref{20-th-cgt-base-propre}) tandis que $c$ une bijection grâce à (\ref{20-commut-princ}). Il suit
que $a$ et $\alpha$ sont des bijections.

\begin{proof}\textit{Preuve du lemme \ref{20-commut-princ}:} Soit $x'$ un point géométrique de $Y'$
d'image $x$ dans $Y$. On peut supposer $x'\in V(J)$. Soit $B$ l'hensélisé (strict) de $Y$ en $x$ et
$B'$ celui de $Y'$ en $x$. On doit étudier la flèche $$H^0(\Spec(B)-V(JB),{\C})\ra
H^0(\Spec(B')-V(JB'),{\C})\leqno (*).$$ Observons que par définition de l'éclatement, $JB$ (resp.
$JB'$) est un idéal principal engendré par un élément non diviseur de zéro et non inversible $t\in
B$ (resp. $t'\in B'$) (équation locale du diviseur exceptionnel). Par ailleurs, les couples $(B,JB)$
et $(B',JB')$ sont henséliens car $B,B'$ sont locaux henséliens (exercice). Les isomorphismes
$$J^nA/J^{n+m}A\iso J^n\chap{A}/J^{n+m}\chap{A},n,m\geq 0$$ assurent  que $B$ et $B'$ ont même
complété $J$-adique $\chap{B}=\chap{B}^J$.

On utilise alors les généralisations des résultats d'Elkik \cite{Elk73} -- et donc de
Ferrand-Raynaud pour le $\pi_0$ -- au cas principal non noethérien de \cite{ART@Gabber-Ramero}.
Précisément, le théorème 5.4.37 appliqué au $B[t^{-1}]$-groupoïde discret
$F_B=\Spec(B[t^{-1}])\times F$ assure qu'on a
$$H^0(\Spec(B[t^{-1}]),F)=\pi_0(F_B)=\pi_0(F_{\chap{B}})=H^0(\Spec(\chap{B}[t^{-1}]),F)$$ et de même
en remplaçant $B,t$ par $B',t'$. Comme $B$ et $B'$ ont même complété $J$-adique, on a donc
$$H^0(\Spec(B[t^{-1}]),F)=H^0(\Spec(B'[t'^{-1}]),F),$$ ce qu'on voulait. Dans le cas
$\C=\underline{\Tors}(U,G)$, on déduit du cas discret que (*) est pleinement fidèle. Soit alors
$\chap{P}$ un revêtement galoisien de groupe $G$ sur $\chap{U}=\Spec(\chap{B})-V(J\chap{B})$.
D'après le théorème 5.4.53 de  \cite{ART@Gabber-Ramero}, il provient d'un (unique) revêtement $P$ de
$U$. La pleine fidélité de (*) assure que le groupe d'automorphismes de $P$ est $G$. Dire que $P$
est galoisien de groupe $G$, c'est dire que la flèche canonique $$\phi: P\times G\ra P\times_UP$$
est un isomorphisme. On peut voir cette flèche comme un morphisme de revêtements étales de $U$.
Après image inverse sur $\chap{U}$, elle s'identifie à la flèche analogue $$\chap{P}\times G\ra
\chap{P}\times_{\chap{U}}\chap{P}$$ qui est un isomorphisme (de revêtements étales de $\chap{P}$
donc de revêtements étales de $\chap{U}$) par hypothèse. La pleine fidélité de (*) assure que $\phi$
est un isomorphisme de sorte que $P$ est bien galoisien de groupe $G$. On a donc obtenu que le
foncteur naturel entre les catégories de $G$-revêtements galoisiens sur $U$ et $\chap{U}$ sont
équivalentes. Il en est donc de même pour le foncteur les catégories de $G$-revêtements galoisiens
sur $U'$ et $\chap{U'}$. On conclut en se souvenant de l'égalité $\chap{U}=\chap{U'}$. \end{proof}

\end{proof}

\begin{remarques2}\label{20CTC} Le théorème~\ref{20-prop-rig-champ} entraîne immédiatement que la
flèche de changement de base $$\pi^*j_*\C\ra\chap{j}_*\pi^*\C$$ est une équivalence. En effet, on
l'a déjà vu sur $\chap{U}$ (\ref{20-CTCU}). Si $\hat{x}\not\in\chap{U}$, on a déjà observé dans la
preuve de \ref{20-cbase1} que $U_{(\bar x)}$ et $\chap{U}_{(\overline{\hat x})}$ avaient même
complété $I$-adique de sorte que deux applications de \ref{20-prop-rig-champ} assurent que la fibre
de $$\pi^*j_*\C\ra\chap{j}_*\pi^*\C$$ en $\overline{\hat x}$ est une
équivalence.\end{remarques2}

\section{Rigidité de la ramification}

\subsection{La condition $c_2$}

Rappelons (\egalong{IV}{4}{18.6.7}) que l'hensélisé $A^\h$ d'un anneau semi-local $A$ est le produit
des hensélisés des localisés de $A$ en ses idéaux maximaux. Pour tout anneau noethérien, on note
$A^\nu$ son normalisé, à savoir la clôture intégrale de $A$ dans l'anneau total $K(A)$ des fractions
de $A_\red$. Puisque $K(A)$ est le produit des $K(A/\p)$ où $\p$ décrit les points maximaux de
$\Spec(A)$, le normalisé de $A$ est le produit des normalisés des $A/\p$. Si $A^\nu$ n'est en
général pas noethérien (\cite[exemple 5 de l'appendice]{Nag62}), il est en revanche quasi-fini sur
$A$ (\cite[V.33.10]{Nag62}). En particulier, si $A$ est local, $A^\nu$ est semi-local. Si $A$ est
local noethérien, $A^\nu$ est semi-local de sorte que son hensélisé est bien défini. On a alors
(comparer avec \cite[43.20 et exercice 43.21]{Nag62})

\begin{lemme}\label{20hnu=nuh}Soit $A$ un anneau local noethérien.  \begin{itemize} \item La flèche
canonique $A^\h\ra (A^\nu)^\h$ induit un isomorphisme $(A^\h)^\nu\iso (A^\nu)^\h$.  \item Cette
bijection induit une bijection canonique $\p\mapsto\p^*$ entre les points maximaux $\p$ de
$\Spec(A^h)$ et les points fermés $\p^*$ de $\Spec(A^\nu)$  de telle sorte que les anneaux intègres
$(A^h/\p)^\nu$ et $(A^\nu_{\p^*})^h$  sont (canoniquement) isomorphes.  \end{itemize} \end{lemme}

\begin{proof}D'après \egalong{IV}{4}{18.6.8}, le morphisme canonique $A^\nu\otimes_AA^h\ra
(A^\nu)^\h$ est un isomorphisme. Le morphisme canonique  $A\ra A^\h$ étant ind-étale, il est normal.
D'après \egalong{IV}{2}{6.14.4}, le morphisme canonique $A^\h\ra A^\nu\otimes_AA^h$ identifie
$A^\nu\otimes_AA^h$ à la fermeture intégrale de $A^\h$ dans $A^\h\otimes_AK(A)$. Si maintenant,
$A\ra B$ est étale, la fibre au point maximal $\p\in\Spec(A)$ s'identifie à $\Spec(K(B))$. En
passant à la limite, on déduit l'égalité $A^\h\otimes_AK(A)=K(A^\h)$ de sorte que
$A^\nu\otimes_AA^h$ s'identifie à la fermeture intégrale de $A^\h$ dans $A^\h\otimes_AK(A)=K(A^\h)$
et donc $(A^\h)^\nu\iso A^\nu\otimes_AA^h$. La composition $$(A^\h)^\nu\iso A^\nu\otimes_AA^h\iso
(A^\nu)^\h$$ est l'isomorphisme annoncé. Pour le second point, on observe d'une part que le spectre
du normalisé de $A^\h$ est la somme disjointe des normalisés de ses composantes irréductibles
\begin{equation}\label{20dec-norm} \Spec((A^\h)^\nu)=\coprod_{\p\textup{\ point
maximal}}\Spec((A^\h/\p)^\nu), \end{equation} chaque fermé $\Spec((A^\h/\p)^\nu)$ étant intègre
(puisque local et normal) de sorte que \ref{20dec-norm} est une la décomposition en composantes
irréductibles de $\Spec((A^\h)^\nu)$.  D'autre part, par définition de l'hensélisé d'un anneau
semi-local, on a \begin{equation}\label{20dec-hens} \Spec((A^\nu)^\h)=\coprod_{\p^*\textup{\ point
fermé}}\Spec((A^\nu_{\p^*})^\h).  \end{equation} Or, $(A^\nu_{\p^*})^\h$ est local et normal (comme
$A^\nu_{\p^*}$), donc intègre, prouvant que \ref{20dec-hens} est la décomposition en composantes
irréductibles de $\Spec((A^\nu)^\h)$. Le lemme suit.  \end{proof}

\begin{proposition}\label{20propc2} Soit $Z$ un sous-schéma fermé d'un schéma noethérien $X$.  Les
conditions suivantes sont équivalentes : \begin{itemize} \item [(i)] Soit $p:X^\nu\ra X$ le
morphisme de normalisation. Alors, $p^{-1}(Z)$ est de codimension $\geq 2$ dans $X^\nu$.  \item
[(ii)] Pour tout $z\in Z$, toutes les composantes irréductibles de $\Spec(\OO_{X,z}^\h)$ sont de
dimension $\geq 2$.  \item [(iibis)] Pour tout $z\in Z$, toutes les composantes irréductibles de
$\Spec(\OO_{X,z}^\hs)$ sont de dimension $\geq 2$.  \item [(iii)] Pour tout $z\in Z$, toutes les
composantes irréductibles de $\Spec(\widehat{\OO_{X,z}})$ sont de dimension $\geq 2$.

\end{itemize} \end{proposition}

\begin{proof}Notons $A=\OO_{X,z}$ pour $z\in Z$. Notons d'abord que le morphisme $A^\h\ra A^\hs$ est
injectif, entier et fidèlement plat. Ceci prouve que le morphisme $h:\Spec(A^\hs)\ra\Spec(A)$
vérifie $\dim(\overline{h(x)})=\dim(\overline{x})$ et induit une surjection au niveau des points
maximaux, ce qui prouve l'équivalence de (\textit{ii}) et (\textit{iibis}).

Un anneau intègre et son normalisé ainsi qu'un anneau local et son hensélisé, ont même dimension.
Conservant les notations de \ref{20hnu=nuh}, on a donc $$\dim A^h/\p=\dim A^\nu_{\p^*}.$$  Or, dire
$\codim p^{-1}(Z)\geq 2$, c'est dire $\dim A^\nu_{\p^*}\geq 2$ lorsque $\p^*$ décrit les points
fermés de  $$\Spec(A^\nu)=p^{-1}(\Spec(\OO_{X,z}))$$ lorsque $z$ décrit $Z$. Ceci revient donc à dire
que toutes les composantes irréductibles $\Spec(A^h/\p)$ de $\Spec(A^h)$ sont de dimension $\geq 2$
prouvant l'équivalence de (\textit{i}) et (\textit{ii}).

Pour montrer l'équivalence de (\textit{i}) et (\textit{iii}), on peut supposer que $X=\Spec(A)$ est
local hensélien et que $Z$ est réduit à son point fermé.

Prouvons d'abord que (\textit{iii}) implique (\textit{ii}). Soit $Y$ une composante irréductible de
$X$. Le morphisme de complétion $c:\chap{X}\ra X$ étant  fidèlement plat, $\chap{Y}=c^{-1}(Y)$ est
une réunion de composantes irréductibles de $\chap{X}$ de sorte qu'on a $\dim(\chap{Y})\geq 2$.
Comme $Y$ est local noethérien, on a $\dim(Y)=\dim(\chap{Y})\geq 2$.

Prouvons la réciproque. Quitte à se restreindre à une composante irréductible (réduite), on peut
supposer $X$ intègre de dimension $\geq 2$. Soit $\hat x$ (resp. $x$) le point fermé de $\chap{X}$
(resp. $X$) (ce n'est pas une composante irréductible de $\chap{X}$ qui est de dimension $\geq 2$).
Si une des composantes de $\chap{X}$ était de dimension $\leq 1$, elle serait de dimension $1$ (car
$\{\hat x\}$ n'est pas une composante) et donc son point générique serait un point isolé de
$\chap{X}-\{\hat x\}$ de sorte que $\chap{X}-\{\hat x\}$ serait disconnexe (étant de dimension $\geq
2$). Or, d'après \cite[corollaire 4.4]{Fer-Ray70}, la flèche $$\pi_0(\chap{X}-\{\hat
x\})=\pi_0(c^{-1}(X-\{x\}))\ra\pi_0(X-\{x\})$$ est bijective. Or, comme $X$ est intègre de dimension
$\geq 2$, l'ouvert  $X-\{x\}$ est intègre donc connexe. \end{proof}

\begin{definition}\label{20defc2} Avec les notations de \ref{20propc2}, si $Z$ vérifie les
conditions équivalentes de  \ref{20propc2}, on dit que $Z$ est $c_2$ dans $X$.  \end{definition}

\begin{remarque}\label{20remc2} Si $X$ est intègre et excellent, $Z$ est $c_2$ si et seulement si
$X-Z$ contient tous les points de codimension $\leq 1$. En effet, comme le morphisme de
normalisation est fini et $X, X^\nu$  caténaires, on a $\dim\OO_{X^\nu,z^\nu}=\dim\OO_{X,p(z^\nu)}$
pour tout $z^\nu\in p^{-1}(Z)$.\end{remarque}

\begin{proposition}\label{20platc2} Soit $f:X'\ra X$ un morphisme plat de schémas noethériens et $Z$
un fermé de $X$. Alors, si $Z$ est $c_2$ dans $X$, son image inverse $Z'=f^{-1}(Z)$ est $c_2$ dans
$X'$. En particulier, la condition $c_2$ est invariante par localisation Zariski ou
étale.\end{proposition}

\begin{proof}Soit $z'\in Z'$ d'image $z=f(z')\in Z$. On suppose donc (\ref{20propc2}) que toutes les
composantes de $A=\widehat{\OO_{X,z}}$ sont de dimension $\geq 2$ et on veut prouver que toutes les
composantes de $B=\widehat{\OO_{X',z'}}$ sont de dimension $\geq 2$. On peut donc supposer que $f$
est morphisme local de schémas noethériens, locaux et complets. Comme $f$ est plat, toute composante
de $X'$  domine une composante $X_0$ de $X$ et est une composante de $f^{-1}(X_0)$. On peut donc
supposer $X$ intègre de dimension $>1$, de point fermé $z$. D'après \sga{2}{VIII}{2.3}, le
$A$-module $\OO(X-z)$ est de type fini. Comme $B$ est plat sur $A$, on déduit que
$B\otimes_A\OO(X-z)=\OO(X'-f^{-1}(z)$ est de type fini sur $B$. Comme $B$ est noethérien, le sous
$B$-module $\OO(X'-z')$ de $\OO(X'-f^{-1}(z))$ est de type fini. Mais si une des composantes $X'_0$ de
$X'$ était de dimension $1$, le complémentaire $X'_0-z'$ serait réduit au point générique $\eta$ de
$X'_0$ qui serait isolé dans $X'-z'$ isolé de sorte que $\OO(X'_0-z')$ serait un sous $B$-module de
$\OO(X'-z')$, donc de type fini ($B$ est noethérien). A fortiori, $\OO(X'_0-z')$ serait de type fini
comme $\OO(X'_0)$-module, ce qui contredit \sga{2}{VIII}{2.3} puisque $X'_0$ est de dimension
$1$.\end{proof}

\section{Théorème de rigidité de la ramification I : forme faible } Nous allons commencer par
démontrer une variante du changement de base lisse qui est cruciale dans la preuve du théorème de
rigidité~\ref{20rigfaible}.  \subsection{Variante du théorème de changement de base lisse}

Soit $G$ un groupe fini. On va démontrer une variante du théorème de changement de base lisse
\sga{4}{XVI}{1.2} pour les faisceaux de $G$-torseurs sans hypothèse sur le cardinal de $G$, mais en
se restreignant au cas d'immersions ouvertes.

\begin{theoreme2}[Gabber]\label{20-cbl}Considérons un diagramme cartésien
$$\xymatrix{U'\ar@{^(->}[r]^{j'}\ar[d]\ar@{}[rd]|\square&X'\ar[d]^p\\ U\ar@{^(->}[r]^{j}&X}$$
Supposons $X$ excellent normal, $p:X'\ra X$ lisse et $j:U\ra X$  immersion ouverte telle que $U$
contient tous les points de codimension $\leq 1$. Alors, le morphisme de changement de base
$\Phi:p^*j_*\tors(U,G)\ra j'_*\tors(U',G)$ est une équivalence.  \end{theoreme2}

\begin{proof} D'après le théorème de changement de base lisse  pour les faisceaux d'ensembles
\sga{4}{XVI}{1.2}, $\Phi$ est pleinement fidèle. Il suffit de prouver l'essentielle surjectivité.
Soit $x'$ un point géométrique de $X'$ d'image $x=p(x')$. Passant aux fibres, on est ramené à
prouver que la flèche d'image inverse des torseurs $H^1(U_{(x)},G)\ra H^1(U'_{(x')},G)$ est
bijective, avec de plus $x'$ \emph{fermé} dans sa fibre \sga{4}{VIII}{3.13 b)}. La stricte
hensélisation préserve la normalité et la codimension (platitude).  Les propriétés de permanence des
anneaux excellents (cf. \refext{I}{1-henselisation-anneaux-excellents})
assurent donc qu'on peut supposer $X=\Spec(A), X'=\Spec(A')$
avec $A=\OO_{X,x}^\hs, A'=\OO_{X',x'}^\hs$ strictement locaux, normaux et excellents. Comme $p$ est
lisse, le choix de coordonnées locales $t_1,\cdots,t_n$ de $X'$ en $x'$ définit un $A$-isomorphisme
$A\{t_1,\cdots,t_n\}\iso A'$ où comme d'habitude $A\{t_1,\cdots,t_n\}$ désigne l'hensélisé strict de
$A[t_1,\cdots,t_n]$ à l'origine. Une récurrence évidente permet de supposer $n=1$. On s'est ramené à
la situation $$\xymatrix{{U'}\ar@{}[rd]|\square\ar[d]^p\ar@{^(->}[r]&X'\ar[d]^p\\
U\ar@{^(->}[r]^-{j}\ar@/^1pc/[u]_\sigma&X\ar@/_1pc/[u]_\sigma}$$ avec $A$ strictement local, normal
et excellent et $\sigma$ la section de $p$ définie par l'immersion fermée d'équation $t=0$. Comme
$X,X'$ sont locaux et normaux, ils sont intègres. Les ouverts non vides de $X,X'$ sont donc intègres
et donc connexes. Le composé $$\pi_1(U)\xrightarrow{\sigma_*}\pi_1({U'})\xrightarrow{p_*}\pi_1(U)$$
étant l'identité, il suffit de prouver que $\sigma_*$ est surjectif. Soit alors ${V'}$ un revêtement
étale connexe de ${U'}$. On doit prouver que sa restriction $V\ra U$ au fermé
$U\stackrel{\sigma}{\hookrightarrow} U'$ d'équation $t=0$ est connexe.

Comme ${X'}$ est excellent, la clôture intégrale ${Y'}$ de ${X'}$ dans ${V'}$  est finie sur ${X'}$,
normale et intègre (comme $X'$). Comme $X'$ est hensélien, il en de même de $Y'$ qui est donc une
union disjointe de ses composants locaux. Comme $Y'$ est intègre, $Y'$ est local.  Soit $D\subset
{Y'}$ le diviseur de Cartier d'équation $t=0$ :  $D$ est connexe, puisque fermé dans un schéma
local.

On a donc un diagramme commutatif à carrés cartésiens et où les flèches verticales sont finies (et
dominantes).

$$\xymatrix{V'\ar@{^(->}[r]\ar[d]&Y'\ar[d]&D\ar@{_(->}[l]\ar[d]&V\ar@{_(->}[l]\ar[d]\\
U'\ar@{^(->}[r]&X'&X\ar@{_(->}[l]_-\sigma&U\ar@{_(->}[l]}$$

Soit $x'$ un point de $D-V$, d'image $x$ dans $X-U$. Comme $D\ra X$ est fini, on a
$\dim\overline{\{x'\}}=\dim\overline{\{x\}}$ et $\dim(D)=\dim(X)$. Comme $X$ et $X'$ sont caténaires
(ils sont même excellents), on en déduit l'égalité $\dim\OO_{D,x'}=\dim\OO_{X,x}$ ce qui assure que
l'ouvert $V$ dans $D$ contient tous les points de codimension $1$ dans $D$ (de même que  le
complémentaire de $U$ dans $X$ contient tous les points de codimension $1$ dans $X$). D'après le
lemme \refext{XXI}{21-ref-deg-connexe} appliqué au diviseur de Cartier connexe du schéma normal,
excellent  ${Y'}$, le schéma $V$ est connexe.\end{proof}

\subsection{Énoncé et réductions}

\begin{theoreme2}[Rigidité de la ramification] \label{20rigfaible} Soient $X,X'$ des schémas
noethériens, $Z\subset X$ un sous-schéma fermé, $U\stackrel{{\tilde j}}{\hookrightarrow}X$ l'ouvert
complémentaire et $X'\xrightarrow{\pi} X$ un morphisme plat. Notons $U'\stackrel{{\tilde
j}'}{\hookrightarrow}X'$  l'immersion ouverte $U'=\pi^{-1}(U)\hookrightarrow X'$.  On suppose que
$\pi$ est régulier au dessus de $Z$. Soit $\C$ un champ en groupoïdes sur $U_\et$. Alors, la flèche
de changement de base $$\phi(\C):\pi^*{\tilde j}_*\C\ra {\tilde j}'_*\pi^*\C$$ est une équivalence
dans les deux cas suivants :
\begin{itemize} \item[(i)] $\C$ est discret (c'est-à-dire $\C$
équivalent à un faisceau d'ensembles).  \item[(ii)] $Z$ est $c_2$ et $\C=\tors(U,G)$ avec $G$ un
groupe (ordinaire) fini.
\end{itemize}
\end{theoreme2}

En considérant les fibres, on peut supposer que $\pi$ est un morphisme local de schémas strictement
locaux (la condition $c_2$ ne dépendant que des hensélisés stricts aux points de $Z$).

Soient $x,x'$ les points fermés respectifs de $X,X'$. Par récurrence sur la dimension de $X'$, on
peut supposer que $\phi(\C)_{\bar y'}$ est une équivalence en tout point géométrique $\bar y'$ de
$X'-\{x'\}$ et il suffit de prouver que $\phi(\C)_{x'}$ est une équivalence. On peut de plus
supposer $x\in Z$ (sinon $U=X$ et c'est terminé). Par hypothèse, la fibre spéciale $F=\pi^{-1}(x)$
de $\pi$ est géométriquement régulière.

On a un diagramme commutatif à « carrés » cartésiens (avec des notations un peu abusives)

$$\xymatrix{U'\ar@{^(->}[d]_{j'}\ar@/_3pc/[dd]_{\tilde
j'}\ar[r]&U\ar@{^(->}[d]^{j}\ar@/^3pc/[dd]^{\tilde j}\\
X'-F\ar@{^(->}[d]\ar[r]&X-\{x\}\ar@{^(->}[d]\\ X'\ar[r]^\pi&X}$$

Par hypothèse de récurrence, la flèche de changement de base associée au carré supérieur est une
équivalence de sorte qu'on a une équivalence $\pi^*j_*\C\iso j'_*\pi^*\C$ sur $X'-F$.  Comme $X,X'$
sont strictement henséliens, la flèche de changement de base $\phi(j_*\C)_{x'}$ \begin{eqnarray*}
H^0(X,\tilde j_*\C)&=& H^0(X-\{x\},j_*\C) \\ &\xrightarrow{\pi^*}& H^0(X'-F,\pi^*j_*\C) \\ &=&
H^0(X'-F,j'_*\pi^*\C) \\ &=&H^0(X',\tilde j'_*\pi^*\C) \end{eqnarray*}

s'identifie à la flèche d'image inverse

\begin{equation}\label{20bcphi} \pi^*:\ H^0(X-\{x\},j_*\C) \ra H^0(X'-F,\pi^*j_*\C).  \end{equation}

\begin{lemme2} On peut supposer que $\pi$ est un morphisme essentiellement lisse de schémas
strictement locaux et excellents.  \end{lemme2}

\begin{proof}Notons $\chap{X}$ le complété de $X$ le long de son point fermé et $\widehat{X'}$ le
complété de $X'$ le long de $F$.  Pour tout $S$-espace $\E$ sur $S_\et$ avec $S=X,X'$, on note
$\chap{\E}$ son image inverse sur $\chap{S}$.  On a un diagramme \emph{commutatif}
$$\xymatrix{\widehat{X'}\ar[r]^{\gamma'}\ar[d]^{\hat\pi}&X'\ar[d]^\pi\\\chap{X}\ar[r]^\gamma&X}$$ où
$\gamma,\gamma'$ sont les morphismes de complétion, donc sont plats, et $\hat\pi$ est plat comme
$\pi$ et est un morphisme local de schémas noethériens. Sa fibre spéciale est encore $F$ de sorte
qu'elle est géométriquement régulière. Ainsi, $\hat\pi$ est formellement lisse
(\egalong{IV}{4}{19.7.1}) et donc régulier (\cite{Localisation@Andre}) puisque $\chap{X}$ est local
noethérien complet donc excellent. D'après~\ref{20platc2}, $\chap{Z}=\chap{X}-\chap{U}$ est encore
$c_2$.  D'après le théorème de rigidité de Gabber (\ref{20-prop-rig-champ}) appliqué aux paires
henséliennes $(X,x)$ et $(X',F)$, il suffit, pour prouver que le foncteur~\ref{20bcphi} est une
équivalence, de prouver que le foncteur

\begin{equation}\label{20bcphic} \hat\pi^*:\ H^0(\chap{X}-\{x\},\widehat{j_*\C}) \ra
H^0(\chap{X'}-F,\widehat\pi^*\widehat{j_*\C}) \end{equation} est une équivalence. Mais on a d'une
part $$\widehat{j_*\C}={\chap{j}}_*\chap{\C}$$ d'après la version faisceautique (\ref{20CTC}) du
théorème de rigidité de Gabber,  et, d'autre part $$\chap{\C}=\tors(\chap{U},G)$$ d'après
\cite[III.2.1.5.7]{CNA@Giraud}. Pour prouver que (\ref{20bcphi}) est une équivalence, il suffit donc
de prouver que (\ref{20bcphic}) est une équivalence dans le cas où $X$ est complet, donc excellent
et $\pi$ un morphisme local régulier.

 D'après le théorème de Popescu (\cite{NPD@Swan}), le morphisme régulier $\pi$ est limite projective
filtrante de morphismes locaux essentiellement lisses $\pi_i: X'_i\ra X$. Notons que les $X'_i$ sont
strictement locaux et excellents comme $X$. Comme les $X'_i$ sont cohérents,  le foncteur section
globale commute à la limite projective au sens de \sga{4}{VII}{5.7} de sorte qu'il suffit de prouver
le théorème pour les $\pi_i$.\end{proof}

\subsection{Preuve de~\ref{20rigfaible}} On suppose donc que $\pi$ est un morphisme local
essentiellement lisse de schémas excellents.

Si $\C$ est discret, on invoque le théorème de changement de base par un morphisme (essentiellement)
lisse  pour conclure \sga{4}{XVI}{1.2} : le (\text{i}) du théorème~\ref{20rigfaible} est prouvé.

Supposons donc $\C=\tors(U,G)$. On doit donc prouver pour conclure la preuve  du
théorème~\ref{20rigfaible} la variante suivante du théorème de changement de base lisse de Gabber
(\ref{20-cbl}).

\begin{proposition2}\label{20-cblvar} Considérons un diagramme cartésien
$$\xymatrix{U'\ar@{^(->}[r]^{j'}\ar[d]\ar@{}[rd]|\square&X'\ar[d]^\pi\\ U\ar@{^(->}[r]^{j}&X}$$ où
$\pi$ est un morphisme essentiellement lisse de schémas excellents strictement locaux. On suppose
que le fermé complémentaire $Z=X-U$ est $c_2$ (c'est-à-dire sous ces hypothèses,  que $U$ contient
les points de codimension $1$ (\ref{20remc2})). Alors, le morphisme
$$\pi^*:H^0(X-\{x\},j_*\tors(U,G))\ra H^0(X'-\pi^{-1}\{x\},\pi^*j_*\tors(U,G))$$ est une
équivalence.\end{proposition2}

\begin{proof}Comme $X$ est excellent, le morphisme de normalisation $p:X^\nu\ra X$ est fini. Son
image est donc fermée. Comme $p$ est (ensemblistement) dominant, $p$ est surjectif. Comme $p$ est
surjectif, le foncteur $$j_*\tors(U,G)\ra
j_*p_*p^*\tors(U,G)\stackrel{\hbox{\cite[III.2.1.5.7]{CNA@Giraud}}}{=}p_*j^\nu_*\tors(U^\nu,G)$$ est
fidèle\footnote{On note $\E\mapsto\E^\nu$ le foncteur d'image inverse par $p$.}.  D'après
\ref{20-pertinence-eff} et  le  théorème~\ref{20rigfaible}, (\text{i}), il suffit de prouver que la
flèche

\begin{equation}\label{20imainv}\pi^*:H^0(X-\{x\},p_*j^\nu_*\tors(U^\nu,G))\ra
H^0(X'-\pi^{-1}\{x\},\pi^*p_*j^\nu_*\tors(U^\nu,G))\end{equation} est une équivalence.

Considérons le diagramme cartésien

$$\xymatrix{{X'^\nu}\ar[r]^-{\pi^\nu}\ar[d]^{p'}\ar@{}[rd]|\square&{X}^\nu\ar[d]^p\\
X'\ar[r]^{\pi}&X}$$

Comme $p$ est fini (donc propre), on a $\pi^*p_*=p'_*\pi^{\nu*}$ de sorte que \ref{20imainv}
s'identifie à la flèche d'image inverse

\begin{equation}\label{20imainvbis}\pi^*:H^0(X^\nu-\{x\}^\nu,j^\nu_*\tors(U^\nu,G))\ra
H^0(X'^\nu-(\pi^{\nu})^{-1}\{x\}^\nu,\pi^{\nu*}j^\nu_*\tors(U^\nu,G)).\end{equation}

Notons que, la condition $c_2$ ne dépendant que du normalisé, le complémentaire $Z^\nu$ de $U^\nu$
est encore $c_2$ dans $X^\nu$, et $U^\nu$ contient tous les points de codimension $1$.  On invoque
alors le théorème de changement de base lisse de Gabber \ref{20-cbl}.\end{proof}

\subsection{Comparaison à la complétion : cas des coefficients
abéliens dans le cas non nécessairement nœthérien}

\begin{quote}
Le paragraphe suivant est une \emph{esquisse} de démonstration
de l'analogue du théorème \ref{20rigfaible} pour les
coefficients abéliens. Le cas des schémas nœthériens est traité
dans \cite{TUBULAR@Fujiwara}. Nous reproduisons ici fidèlement
une lettre d'Ofer Gabber aux éditeurs (20 juin 2012).
\end{quote}

\vskip.5cm

Let $(A,I) \to (A',I')$ be a map of henselian pairs with $I$
finitely generated, $I' = I A'$, $\widehat{A} \iso
\widehat{A'}$ ($I$-adic completions). $X = \SP(A)$, $X' =
\SP(A')$, $\pi : X' \ra X$, $U = X - V(I)$, $U' = X' -
V(I')$, $j : U \to X$, $j' : U' \ra X'$.

\begin{quote}
\textbf{CTC} : For every torsion abelian  sheaf $F$ on $U$,
the base change arrow $\pi^*R^qj_*F \to R^qj'_*\pi^*F$ is an
isomorphism for all $q$.
\end{quote}

\textit{Analogue of \ref{20rigfaible}} (notations as there) : If $F$
is a sheaf of $\ZZ/n\ZZ$-modules on $U$ where
$n >0$ is invertible on $X$, then $\pi^*R^qj_*F \ra R^qj'_*
\pi^*F$ are isomorphisms.

This is reduced to CTC by the same argument.

\textit{Sketch of proof of CTC using Zariski-Riemann spaces}
: For comparing stalks we may assume $A$, $A'$ strictly
henselian and $I$ a proper ideal, and we want $$ H^q(U,F)
\iso H^q(U',F). \leqno (*) $$ We call a finitely generated
ideal $J \subset A$ containing a power of $I$
\textit{admissible}. We consider the admissible blow-ups
$\Bl_{J}(X)$ which form a cofiltered category using
$X$-scheme morphisms. In general there can be more than one
$X$-morphism between two admissible blow-ups but if we
restrict ourselves to $J$'s with $V(J) = V(I)$ (set
theoretically) (so that $U$ is schematically dense in the
blow-up), there is at most one. Define $J \le J'$ iff there
is an $X$-morphism $\Bl_{J'}(X) \ra
\Bl_J(X)$. This is a filtered preorder. When $V(J)
\subset V(J')$, $J \le J'$ is equivalent to the condition
that for some $n >0$ and ideal $K$, $J'^n = JK$. Thus we
have an isomorphism of the preordered set of admissible
$J$'s of full support in $A$ and the corresponding set for
$A'$. Let $ZRS_I(X) = \varprojlim \Bl_J(X)$ (a
locally ringed space). For the closed point $s$ of $X$ we
can consider the special fiber $ZRS_I(X)_s$ and its étale
topos, which for our purposes may be defined as the
projective limit of the étale topoi
$(\Bl_J(X)_{s})_{\et}$ as in
[\textbf{SGA 4}]. It has enough points by Deligne's theorem.
The points are given by "geometric points" of $ZRS_I(X)_s$
(i. e. a point and a choice of a separable closure of the
residue field). For every admissible $J$ we have $$ j_J : U
\hookrightarrow \Bl_J(X) $$ giving a spectral
sequence (using proper base change) $$
H^p(\Bl_J(X)_s, R^qj_{J*}F) \Rightarrow
H^{p+q}(U,F). \leqno (**) $$ We pass to the limit using the
general theory of \sgacourt{4}{vi}. We get a spectral sequence
$(**)_{\mathrm{lim}}$ involving cohomology on
$(ZRS_J(X)_s)_{\et}$. Since the latter topos is the
same for $X'$, to show $(*)$ we use the morphism of the limit
spectral sequence to reduce to stalks of the limits of the
$R^qj_{J*}$ sheaves.

Using the study in \cite{TUBULAR@Fujiwara} of the local rings of ZRS's
and their henselizations, one reduces $(*)$ to the case of
local rings at geometric points of the special fibers of
ZRS's.  Thus we are reduced to the case $A$, $A'$ are
henselian and $I$-valuative. Say $I = (\varphi)$. Then
$A[\varphi^{-1}]$ is a henselian local ring with maximal
ideal corresponding to $P = ⋂ I^n$, and $A/P$ is a
henselian valuation ring whose valuation topology is the
$\varphi$-adic one. In this case to prove $(*)$ one reduces to
the corresponding statement for $\mathrm{Frac}(A/P) \to
\mathrm{Frac}(A'/P')$. In fact for $K \to K'$ a dense
embedding of henselian valued fields, if we choose separable
closures $K_{\mathrm{sep}}$, $K'_{\mathrm{sep}}$ and a map
between them we have $\mathrm{Gal}(K'_{\mathrm{sep}}/K')
\iso \mathrm{Gal}(K_{\mathrm{sep}}/K)$, using forms of
Krasner's lemma (cf. \accourt{VI}{8, exercices 12, 14 a}).

Note : For admissible $J$, $\Bl_J(X') \ra
\Bl_J(X)$ gives an isomorphism on $I$-adic
completions (as in the discussion in the proof of
\ref{20-commut-princ}) as for every $m$ the map $$ \oplus_n
J^n \to \oplus_n J^nA' $$ is an isomorphism mod $I^m$.

\section{Rigidité de la ramification II : forme forte}

\subsection{Générisations étales immédiates} On note $X_{(x)},X^\h_{(x)}$ et $\chap{X}_{(x)}$ les
localisés, hensélisés et complétés respectivement de $X$ en $x$. On note $\overline{\{y\}}$
l'adhérence de $y$ dans $X$ munie de sa structure réduite. L'hensélisation et la complétion
commutent aux immersions fermées de sorte que $\overline{\{y\}}^\h,\chap{\overline{\{y\}}}$
respectivement coïncident avec l'image inverse de $\overline{\{y\}}$ par les morphismes
d'hensélisation, complétion respectivement.  Rappelons
(\refext{XIV}{14-definition-specialisation-immediate}) qu'une générisation $y\in X$ d'un point $x$
d'un schéma $X$ est une \textit{générisation étale immédiate} de $x$ si l'hensélisé strict en $\bar
x$ de l'adhérence de $y$ a une composante irréductible de dimension $1$.

\begin{lemme2}\label{20et-imm} Soit $y$ une générisation de $x$. Notons $c:\chap{X}_{(x)}\ra
X_{(x)}$ le morphisme de complétion. Alors, $y$ est une générisation étale immédiate de $x$ si et
seulement si l'un des points maximaux de $c^{-1}(y)$ est une générisation étale immédiate du point
fermé de $\chap{X}_{(x)}$.\end{lemme2}

 \begin{proof} Notons pour simplifier $Y=\overline{\{y\}}$. Observons d'abord qu'un des trois
schémas $Y_{(x)},Y^\h_{(x)}$ et $\chap{Y}_{(x)}$ possède un point maximal de dimension nulle si et
seulement si chacun est réduit (ensemblistement) à son point fermé. On peut donc exclure ce cas. Le
morphisme $Y_{(\bar x)}\ra Y^\h_{(x)}$ est fidèlement plat et entier. Donc,  l'hensélisé strict
possède un point maximal de dimension $1$\footnote{On devrait plutôt dire point maximal dont
l'adhérence est de dimension $1$.} si et seulement si l'hensélisé $Y^\h_{(x)}$ possède un point
maximal de dimension $1$. D'après (\ref{20propc2}), $Y^\h_{(x)}$ possède un point maximal de
dimension $1$ si et seulement si $\chap{Y}_{(x)}(\ref{20propc2})$ possède un point maximal de
dimension $1$. Par platitude de $c$, il s'envoie nécessairement sur $y$, le point générique de
$Y$.\end{proof}

On peut caractériser agréablement les générisations étales immédiates.

\begin{lemme2} Soit $f:X_{(\bar x)}\ra X_{(x)}$ le morphisme d'hensélisation strict. Les
générisations étales immédiates de $x$ sont les images $y=f(y')$ des $y'\in X_{(\bar x)}$ tels que
$\dim\overline{\{y'\}}=1$.\end{lemme2}

\begin{proof} Soit $y'\in X_{(\bar x)}$ tel que $\dim\overline{\{y'\}}=1$. L'image $y=f(y')$ est une
générisation stricte de $x$ (car par exemple les fibres de $f$ sont discrètes). Pour cette même
raison, $\overline{\{y'\}}$ est une composante de $f^{-1}(\overline{\{y\}})=\overline{\{y\}}_{(\bar
x)}$. Inversement, si $y$ est une générisation étale immédiate de $x$, le point générique $y'$ d'une
composante de dimension $1$ de $\overline{\{y\}}_{(\bar x)}$ s'envoie sur $y$ (platitude de $f$) et
son adhérence est de dimension $1$.\end{proof}

\begin{exemple2}Prenons l'exemple du pincement de \egalong{IV}{2}{5.6.11}. En conservant les
notations de \textit{loc. cit.}, l'anneau pincé $C$ est local noethérien de dimension $2$ et son
normalisé a deux idéaux maximaux de hauteur $1,2$ respectivement. D'après~\ref{20hnu=nuh},
l'hensélisé de $C$ a deux composantes irréductibles de dimension $1$ et $2$ de points génériques
$c,c'$. Comme dans la preuve de \refext{XIV}{14-prop-Laszlo-2}, , ceci assure l'existence de $\bar
c$ (au dessus de $c$) dans l'hensélisé strict de $C$ dont l'adhérence est de dimension $1$ et donc
que le point générique de $\Spec(C)$ est une générisation étale immédiate de son point
fermé.\end{exemple2}

\subsection{Couples associés et condition (*)}\label{20condstar} Commençons par une définition.

\begin{definition2}\label{20defass}Soit $x$ un point d'un schéma $X$. Choisissons une clôture
séparable de $k(x)$ définissant un point géométrique $\bar x$ de $X$.  \begin{itemize} \item[(i)]
Soit $G$ un schéma en groupes sur $X$. On définit les sections locales de $G$ à support dans $\bar
x$ par la formule $$H^0_{\bar x}(G)=\Ker(H^0(X_{(\bar x)},G)\ra H^0(X_{(\bar x)}-\{\bar x\},G)).$$
\item[(ii)] Soit $\C$ un champ (en groupoïdes)  sur $X$ et $p$ un nombre premier. On dit que $(x,p)$
est associé de $\C$ et on écrit $(x,p)\in \Ass(\C)$ si il existe $\sigma\in\C_{\bar x}$ tel que
$H^0_{\bar x}(\aut(\sigma))$ ait de la $p$-torsion.  \item[(iii)] Soit $\C$ un champ ind-fini (en
groupoïdes) sur un ouvert $U$ de $X$. On dit que $\C$ vérifie la condition (*) si pour tout $x\in
X-U$ de caractéristique $p>0$, il n'existe pas de générisation étale immédiate $y$ de $x$ telle que
$(p,y)\in\Ass(\C)$.  \end{itemize} \end{definition2}

Remarquons que la condition $(x,p)$ associé ne dépend pas du choix de la clôture séparable de
$k(x)$.

\begin{exemple2}\label{20exetoile} Supposons $X$ normal et $G$ groupe fini. Soit $U$ un ouvert de
$X$. Alors, $(x,p)$ est associé de $\C=\tors(U,G)$ si et seulement si $p|\card(G)$ et $x$ est un
point maximal de $U$. En effet, l'unique objet de $\C_{\bar x}$ est le torseur trivial $\sigma$ et
$\aut(\sigma)=G$. Or, $X_{(\bar x)}-\{\bar x\}$ est connexe (resp. vide) si $x$ non maximal (resp.
maximal). Ainsi, on a $H^0_{\bar x}(\aut(\sigma))=\{1\}$ (resp. $H^0_{\bar x}(\aut(\sigma))=G$). On
déduit que $\C$ vérifie (*) si et seulement si $U$ contient tous les points de codimension $1$ dont
la caractéristique divise l'ordre de $G$.  \end{exemple2}

\begin{lemme2}\label{20ass-comp} Soit $f:X\ra Y$ un morphisme plat de schémas noethériens, $x\in X$
d'image $y=f(x)$ dans $Y$  et $\C$ un champ en groupoïdes sur $Y$. Alors, $(x,p)\in\Ass(f^*\C)$ si
et seulement si $(y,p)\in\Ass(\C)$ et $x\in\Max(f^{-1}(y))$.  \end{lemme2}

\begin{proof} Choisissons un point géométrique $\bar x$ au dessus de $x$, qui définit $\bar y$ au
dessus de $y$.

Supposons $(x,p)\in\Ass(f^*\C)$.  Comme la flèche $(f^*\C)_{\bar x}\ra\C_{\bar y}$ est une
équivalence, il existe $\sigma\in \C_{\bar y}$ et $g\in\Aut(\sigma)$ tel que $f^*g$ est
d'ordre de $p$ et de  support $\{\bar x\}$. Notons $F$ l'hensélisé strict de $f^{-1}(y)$ en $\bar
x$. C'est aussi la fibre de l'hensélisé $\varphi: X_{(\bar x)}\ra Y_{(\bar y)}$ de $f$ au dessus de
$\bar y$. Si $F$ n'était pas réduit à $\bar x$, un des points de $F$ ne serait pas dans le support
de $f^*g$ de sorte que $f^*g$ serait l'identité en ce point. Mais $f^*g$ est constant sur
$F=\varphi^{-1}(\bar y)$ de sorte que $f^*g$ serait l'identité également en $\bar x\in F$, ce qui
n'est pas. Donc, $F$ est réduit à $\bar x$ de sorte que $\dim \OO_{f^{-1}(y),x}=0$ (puisqu'un anneau
local a même dimension que son hensélisé strict) et $x\in\Max(f^{-1}(y))$. De plus, $g$ est trivial
sur $\varphi(X_{(\bar x)}-\{\bar x\})=Y_{(\bar y)}-\{\bar y\}$ (fidèle platitude de $\varphi$) ce
qui assure $(y,p)\in\Ass(\C)$.

Inversement, supposons  $(y,p)\in\Ass(\C)$ et $x\in\Max(f^{-1}(y))$. On a donc un automorphisme $g$
d'ordre $p$ de $\sigma\in\C_{\bar y}$ de support $\{\bar y\}$. Le support de $\varphi^*g$ est la
fibre $\phi^{-1}(\bar y)=f^{-1}(y)_{(\bar x)}$. Comme $x$ est maximal dans $f^{-1}(y)$, on déduit
(dimension) que le schéma local $\phi^{-1}(\bar y)$ est de dimension nulle donc réduit à $\bar x$,
ce qu'on voulait.  \end{proof}

\begin{corollaire2}\label{20corass-comp} Soit $(X,x)$ un schéma local noethérien hensélien,
$U=X-\{x\}$ l'ouvert complémentaire du point fermé et $c:\chap{X}\ra X$ le morphisme de complétion.
Alors, le champ en groupoïdes $\C$ sur $U$ vérifie (*) si et seulement $\chap{\C}=c^*\C$ vérifie
(*).\end{corollaire2}

\begin{proof} On note encore $x$ le point fermé de $\chap{X}$ et on choisit un point géométrique
$\bar x$ au dessus de $x$.

Supposons que $\chap{\C}$ vérifie (*). Soit $(y,p=\car(x))\in\Ass(\C)$ et notons $Y$ l'adhérence de
$y$ dans le localisé (Zariski) de $X$ en $x$. Il s'agit de montrer que toutes les composantes de
$Y_{(\bar x)}$ sont de dimension $\geq 2$, ou encore (\ref{20propc2}) que toutes les composantes de
$\chap{Y}=\chap{Y}_{(x)}^h$ sont dimension $\geq 2$, c'est-à-dire (\ref{20propc2} à nouveau), que
toutes les composantes de $\chap{Y}_{(x)}$ sont de dimension $\geq 2$. Mais c'est bien le cas car,
d'après le lemme~\ref{20ass-comp}, on a $(\hat y,p)\in\Ass(\chap{\C})$.

Inversement, supposons que $\C$ vérifie (*). Soit donc $(\hat y,p=\car(x))\in\Ass(\chap{\C})$ et
soit $y=c(\hat y)$ et supposons que l'adhérence de $\hat y$ est de dimension $1$.  Si une des
composantes  $Y=\overline{\{y\}}$ était de dimension $>1$, toutes les composantes de
$\chap{Y}=c^{-1}(Y)$ seraient de dimension $>1$ (\ref{20propc2}). D'après le lemme~\ref{20ass-comp},
on sait que $\hat y$ est maximal dans $c^{-1}(y)$ donc dans $\chap{Y}$ de sorte que
$\dim\overline{\{\hat y\}}>1$, une contradiction. \end{proof}

\subsection{Le théorème de rigidité de la ramification}

\begin{theoreme2}[Rigidité de la ramification II]\label{20rigfort} Soit $\pi:X'\ra X$ un morphisme
plat de schémas noethériens, régulier au dessus d'un sous-schéma fermé $Z\subset X$.  Soit
$j:U=X-Z\hookrightarrow X$ l'immersion ouverte du complémentaire de $Z$ et $\C$ un champ ind-fini
sur $U$ vérifiant la condition (*). Alors, la flèche de changement de base $\pi^*j_*\C\ra
j'_*\pi'^*\C$ est une équivalence.\end{theoreme2}

\begin{proof} D'après le théorème de rigidité de la ramification I (\ref{20rigfaible}), le théorème
est vrai dans le cas discret de sorte que $\pi^*j_*\C\ra j'_*\pi'^*\C$ est toujours pleinement
fidèle.  Comme dans la preuve de \ref{20rigfaible}, on peut supposer $X,X'$ strictement locaux de
point fermés $x,x'$ et $\pi$ morphisme local. Par récurrence sur la dimension de $X$, on peut
supposer que le changement de base par $\pi$ est une équivalence pour l'immersion $U\hookrightarrow
X-\{x\}$ de sorte qu'on peut supposer $U=X-\{x\}$. Comme dans la preuve de \ref{20rigfaible} et en
utilisant l'invariance par complétion de la condition (*) (\ref{20corass-comp}), on peut supposer de
plus $X$ complet et $\pi$ morphisme essentiellement lisse et local et il s'agit de démontrer que la
flèche $$\pi^*:H^0(X-\{x\},\C)\ra H^0(X'-\pi^{-1}\{x\},\C')$$ est essentiellement surjective
(puisque en tout point $y'\not=x'$ la fibre du changement de base est une équivalence par
hypothèse).

Soit donc $\sigma'$ un objet de $H^0(X'-\pi^{-1}\{x\},\C')$. La condition (*) étant stable par
passage aux sous-gerbes (maximales), on peut comme dans la preuve de \ref{20-red-gerb} en
considérant la sous-gerbe maximale de $\C'$ engendrée par $\sigma'$,  supposer de plus que $\C$ est
une gerbe. Comme $\C$ est ind-finie, on peut supposer que $\C$ est \emph{constructible}.

Pour tout point maximal $y\in U$, notons $i_y$ le morphisme canonique $$i_y:\Spec(k(y))\ra
U=X-\{x\}.$$

 Soit $$\Psi : \C\ra\D:=\prod_{y\in \Max(U)}i_{y*}i_{y}^*\C$$ le morphisme déduit des morphismes
d'adjonction. La catégorie fibre de $i_{y*}i_{y}^*\C$ sur un ouvert étale $V\ra X$ s'identifie aux
sections rationnelles de $\C$ définies au voisinage (Zariski) des points maximaux de $V$ au dessus
de $y$. On déduit que $\Psi$ est conservatif et couvrant (étale localement surjectif sur les flèches
et les objets). Pour toute section $\tau\in H^0(U,\D)$ (vu comme un morphisme de $U$-espaces $\tau
:U\ra\D$), le champ des relèvements $K(\tau)=U\times_{\D}\C$ associé  est  donc une gerbe
(\cite[IV.2.5.4]{CNA@Giraud}), évidemment constructible.

Il suffit alors (exercice) de vérifier que la flèche de changement de base est essentiellement
surjective pour \begin{itemize}

  \item [1)] les gerbes $\G=i_{y*}i_{y}^*\C$; \item [2)] la gerbe des relèvements
$K(\tau)=U\times_{\D}\C$ associée à $\tau\in H^0(U,\D)$.  \end{itemize}

\subsubsection{Premier cas : changement de base pour
$\G=i_{y*}i_{y}^*\C$}

Supposons donc $\G=i_{y*}i_{y}^*\C$. Quitte à changer $X,U$ en $\overline{\{y\}},
U\cap\overline{\{y\}}$, on peut supposer que $X$ est irréductible de point générique $y$.

Si la dimension de $X$ est  $1$, on a $U=\{y\}$ et $\G_{\bar y}=\tors(\Spec(\overline{k(y)}),G)$
avec  $p=\car(y)$ ne divisant pas l'ordre de $G$ (cf. l'argument dans l'exemple~\ref{20exetoile}).
On invoque alors le changement de base par un morphisme lisse usuel (\cite[VII.2.1.2]{CNA@Giraud}).

On suppose donc que la dimension de $X$ est $>1$.

Choisissons une clôture séparable $k(y)\hookrightarrow k_y$ et notons ${j_y}:\Spec(k_y)\ra U$ est le
morphisme canonique. On a $$j_y^*\C=\tors(\Spec(k_y),G)$$ où $G$ est un groupe fini constant. Comme
$$i_{y*}i_y^*\C\ra {j_y}_*j_y^*\C$$ est fidèle, on peut (\ref{20-pertinence-eff}) remplacer $\G$ par
${j_y}_*j_y^*\C={j_y}_*\tors(\Spec(k_y),G)$.

\begin{lemme2} On a $R^1j_{y*}G=\{*\}$ et
${j_y}_*\tors(\Spec(k_y),G)=\tors(U,{j_y}_*G)$.\end{lemme2} \begin{proof} La seconde égalité découle
de la première et de la formule (\cite[V.3.1.5]{CNA@Giraud})
$$\pi_0({j_y}_*\tors(\Spec(k_y),G))=R^1j_{y*}G.$$ Soit $\tilde A$ l'hensélisé strict de $X=\Spec(A)$
en un point géométrique $\bar\xi$ de $X$. C'est une limite inductive filtrante d'algèbres $A_i$ de
type finies qui sont génériquement étales. On déduit que $j_y^{-1}(\Spec(\tilde A))$ est le spectre
de la limite inductive filtrante des algèbres étales $B_i=k_y\otimes_{k(y)}A_i$ qui sont donc
scindées puisque $k_y$ est séparablement clos. Ainsi, les schémas considérés étant cohérents, on a
$$(R^1j_{y*}G)_{\bar\xi}=H^1(j_y^{-1}(\Spec(\tilde A)),G)=\varinjlim
H^1(\Spec(B_i),G)=\{*\}.$$\end{proof}

En terme de module galoisien, $j_{y*}G$ est l'induite (continue) $\Hom_c(\Gamma,G)$ où
$\Gamma$ est le groupe profini $\Gal(k_y/k(y))$. En écrivant $\Gamma=\limind
\Gal(K_\alpha/k(y))$ où $(K_\alpha/k(y))_\alpha$ est le système inductif des sous-extensions
galoisiennes finies de $k_y/k(y)$, on trouve $$j_{y*}G=\limind j_{\alpha*}G$$ où
$j_\alpha:\Spec(K_\alpha)\ra\Spec(k(y))\ra U$ est le morphisme canonique. Comme $U, U'$ sont
noethériens donc cohérents, on a $$H^0(U,\tors(U,\limind j_{\alpha*}G))=\limind H^0(U,\tors(U,
j_{\alpha*}G))$$ et $$H^0(U',\pi^*\tors(U,\limind j_{\alpha*}G))=H^0(U',\tors(U',\pi^*\limind
j_{\alpha*}G))=\limind H^0(U,\tors(U', \pi^*j_{\alpha*}G))$$ de sorte qu'on est réduit pour le cas
1) à étudier le changement de base pour la gerbe $\G_\alpha=\tors(U,j_{\alpha*}G)$.

Soit $p:W\ra X$ la normalisation de $X$ dans $\Spec(K_\alpha)\ra X$ : c'est un morphisme fini (car
$X$ est excellent) et surjectif de sorte que  $W$ est semi-local et hensélien (comme $X$). On déduit
que $W$ est la réunion disjointe de ses hensélisés aux points fermés. Comme $W$ est intègre,  $W$
est  strictement local : on note $w$ son point fermé $w$. De plus, $W$ est normal, donc
géométriquement unibranche de sorte que $j_{\alpha *}G=p_*G_{|W-w}$. Comme $R^1p_*G$ est trivial
($p$ est fini), on déduit l'égalité $$\G_\alpha=p_*\tors(W-\{w\},G)$$ comme dans la preuve du
lemme~\ref{20pdirtors} \textit{infra}.

 En utilisant le changement de base propre pour $p$, on est ramené à prouver que la flèche
$$\Tors(W-\{w\},G)\ra\Tors(W'-\pi^{-1}\{w\},G)$$ est une équivalence. Deux cas se
présentent.

Rappelons qu'on a supposé que la dimension de $X$ (ou $W$, c'est la même chose) est  $>1$. Dans ce
cas, $\{w\}$ est $c_2$ dans $W$ et on invoque la variante du théorème de changement de base lisse de
Gabber (\ref{20-cblvar}).

\subsubsection{Deuxième cas : changement de base pour la gerbe des
relèvements $K=K(\tau)$}

Comme le morphisme $K(\tau)\ra \C$ est visiblement fidèle (\cite[IV.2.5.2]{CNA@Giraud}), la gerbe
des relèvements vérifie (*) comme $\C$. Comme $\Psi_{\bar y}$ est une équivalence pour tout point
maximal $y\in U$, on déduit que $K_{\bar y}$ est la gerbe triviale (autrement dit équivalente au
faisceau d'ensembles ponctuels) en tous ces points. Par hypothèse de récurrence, il suffit pour
achever la preuve de prouver le lemme suivant.

\begin{lemme2} Il existe une immersion fermée
$i:F\subset X$ nulle part dense telle que $K=i_*i^*K$.
\end{lemme2}

\begin{proof}Il suffit de prouver que pour tout $y$ maximal, il existe un ouvert de Zariski
contenant $y$ sur lequel $K$ est triviale. Par construction, il existe un voisinage étale $V\ra X$
de $y$ et $\sigma\in K(V)$. Comme $\underline{\Aut}(\sigma)$ est un faisceau constructible de
$V_\et$, l'isomorphisme $$\{\Id\}_{\bar y}\iso\underline{\Aut}(\sigma)_{\bar y}$$ provient
d'un isomorphisme $$\{\Id\}_{W}\iso\underline{\Aut}(\sigma)_W$$ sur un voisinage étale $W\ra
V\ra X$ de $y$. Quitte à localiser, on peut supposer que $W\ra X$ est un revêtement galoisien de son
image $U\subset X$. La section $\sigma$ descend sur $U$ et n'a pas d'automorphisme par construction,
ce qu'on voulait.  \end{proof} \end{proof}

\section{Appendice 1 : sorites champêtres} D'après \ref{20-prop-rig} i), on sait déjà que le
foncteur $H^0(U,\C)\ra H^0(\chap{U},p^*\C)$ est pleinement fidèle que $\C$ soit ind-fini ou non.

Soit $\C$ un champ ind-fini sur $Y$. On cherche des conditions assurant que l'hypothèse

\begin{hypothese}\label{20-hypo-champ}Soit $f:X\ra Y$ un morphisme de schémas. On suppose que pour
tout faisceau d'ensembles $\F$ sur $Y$, la flèche $H^0(Y,\F)\ra H^0(X,f^*\F)$ est
bijective.\end{hypothese}

entraîne que la conclusion

\begin{conclusion}\label{20-concl-champ}  La flèche $\phi: H^0(Y,\C)\ra H^0(X,f^*\C)$ est une
équivalence de catégories.\end{conclusion}

est vraie, autrement dit assurant que l'assertion

\begin{assertion}\label{20-cond-champ}On a l'implication
$\ref{20-hypo-champ}\Rightarrow\ref{20-concl-champ}$.  \end{assertion}

est vraie. On sait déjà que \ref{20-hypo-champ} entraîne que $\phi$ est pleinement fidèle (cf.
\ref{20-sec-enonc}).

\subsection{Premières réductions}

Commençons par un lemme formel :

\begin{lemme2}\label{20-pertinence-eff} Soit $f:X\ra Y$ un morphisme de schémas et $\C_1\ra\C_2$ un
morphisme de champs sur $Y$ qu'on suppose fidèle. Si \ref{20-cond-champ} est vraie pour $\C_2$,
alors \ref{20-cond-champ} est vraie pour $\C_1$.\end{lemme2}

\begin{proof} On a déjà observé (\ref{20-sec-enonc}) que $\phi$ est pleinement fidèle. Soit donc
$$c_1^X\in H^0(X,f^*\C_1)=\Hom_X(X,f^*\C_1)$$ dont on cherche un antécédent dans $H^0(Y,\C_1)$. Son
image $$c_2^X\in H^0(X,f^*\C_2)$$ a un antécédent (à isomorphisme près) $$c_2^Y\in H^0(Y,\C_2).$$

Le couple $(c_1^X,c_2^X=f^*c_2^Y)$ définit une section de la gerbe des relèvements
$$K(f^*c_2^X)=X\times_{f^*\C_1}f^*\C_2.$$

D'après \cite[V.1.4.3]{CNA@Giraud}, la gerbe $K(f^*c_2^X)=X\times_{f^*\C_1}f^*\C_2$ s'identifie à
$$f^*K(c_2^X)=f^*(Y\times_{\C_1}\C_2).$$ Or, $K(c_2^X)$ est un faisceau d'ensembles car $\C_1\ra
\C_2$ est fidèle. Donc, $(c_1^X,c_2^X=f^*c_2^Y)\in H^0(X,f^*K(c_2^X))$ a un unique antécédent de la
forme $(c_1^Y,c_2^Y)$ et $c_1^Y$ est bien l'antécédent cherché. \end{proof}

\begin{lemme2}\label{20-red-gerb} Si \ref{20-cond-champ} est vrai pour toute gerbe (resp. toute
gerbe ind-finie), alors \ref{20-cond-champ} est vrai pour tout champ (resp. tout champ
ind-fini).\end{lemme2}

\begin{proof} Soit $t\in H^0(X,f^*\C)$ et $\gamma_t\subset f^*\C$ la sous-gerbe maximale engendrée
par $t$ dans $f^*\C$ (\cite{CNA@Giraud}, III.2.1.3.2)\footnote{Dans \textit{loc. cit.}, $\pi_0(\C)$
est noté $\mathrm{Ger}(\C)$, qui n'est pas actuellement la notation standard.}.  Ceci définit une
section $\tau\in\pi_0(f^*\C)$ du faisceau d'ensembles $\pi_0(f^*\C)$ des sous-gerbes maximales de
$f^*\C$ (\textit{loc. cit.}, 2.1.4). D'après \textit{loc. cit.}, 2.1.5, la flèche naturelle
$$\pi_0(f^*\C)\ra f^*\pi_0(\C)$$ est bijective. Mais, par hypothèse, la flèche $$H^0(Y,\pi_0(\C))\ra
H^0(X,f^*\pi_0(\C))=H^0(X,\pi_0(f^*\C))$$ est bijective de sorte qu'il existe une (unique)
sous-gerbe (maximale) $\gamma\subset\C$ telle que $f^*\gamma=\gamma_t$, qui sera ind-finie si $\C$
l'est.  L'image dans $H^0(Y,\C)$ de l'antécédent de $t\in H^0(X,f^*\gamma)$ dans $H^0(Y,\gamma)$ est
l'antécédent cherché.\end{proof}

\subsection{Réduction au cas d'un champ de torseurs sous un groupe fini constant} Admettons pour un
instant le résultat suivant, généralisation au cas des champs de la résolution flasque de Godement.

\begin{lemme2}[Lemme d'effacement]\label{20-lemme-efface} Soit $\gamma$ une gerbe ind-finie sur un
schéma cohérent $X$. Il existe un groupe ind-fini $\G$ sur $X$ et un foncteur fidèle
$\gamma\hookrightarrow \underline{\Tors}(X,\G)$.\end{lemme2}

On peut alors prouver le critère suivant

\begin{proposition2}\label{20-prop-red-champ} Soit $f:X\ra Y$ un morphisme de schémas cohérents.  On
suppose que pour tout faisceau d'ensembles $\F$ sur $Y$, la flèche $H^0(Y,\F)\ra H^0(X,f^*\F)$ est
bijective (\ref{20-hypo-champ}).  On suppose en outre  que pour tout morphisme fini $p:Y'\ra Y$
induisant $f':X'=X\times_Y Y'\ra Y'$ et tout groupe fini constant $G$, la flèche
$\Tors(Y',G)\ra\Tors(X',G)$ est une équivalence.  Alors, pour tout champ ind-fini $\C$ sur $Y$, la
flèche $H^0(Y,\C)\ra H^0(X,f^*\C)$ est une équivalence.\end{proposition2}

\begin{proof} Seule l'essentielle surjectivité pose problème.  Les lemmes d'effacement,
\ref{20-pertinence-eff} et \ref{20-red-gerb} permettent de supposer que $\C=\underline{\Tors}(Y,\G)$
où $\G$  est un groupe ind-fini sur $Y$. Comme $X,Y$ sont cohérents, la cohomologie non abélienne
commute aux limites inductives filtrantes \sga{4}{VII}{remarque VII.5.14}. Comme $\G$ est ind-fini,
il est limite inductive filtrante de faisceaux en groupes constructibles \sga{4}{IX}{2.7.2} : on
peut donc supposer $\G$ constructible. Puisque $Y$ est cohérent,  il existe (\textit{loc. cit.},
2.14) une famille finie de morphismes finis $p_i:Y_i\ra Y$ et des groupes finis constants $G_i$ tels
que $\G$ se plonge dans le produit $\prod p_{i*}G_i$. On a donc un morphisme fidèle
$$\underline{\Tors}(Y,\G)\hookrightarrow \underline{\Tors}(Y,\prod{p_{i*}}G_i)=
\prod\underline{\Tors}(Y,{p_{i*}}G_i)$$ grâce à \cite{CNA@Giraud}, III.2.4.4.

Utilisant à nouveau \ref{20-pertinence-eff}, on peut supposer $$\C=\underline{\Tors}(Y,{p_*}G)$$
avec $G$ groupe fini constant et $p:Y'\ra Y$ fini.

\begin{lemme2}\label{20pdirtors}On a
$\underline{\Tors}(Y,{p_*}G)=p_*\underline{\Tors}(Y',G)$.\end{lemme2}

\begin{proof}Comme $p$ est fini, $R^1p_*G$ est trivial. Mais
$\pi_0(p_*\underline{\Tors}(Y',G))=R^1p_*G$ (\cite{CNA@Giraud}, V.3.1.9.1) de sorte que
$p_*\underline{\Tors}(Y',G)$ est une gerbe, visiblement neutre et vaut donc nécessairement
$\underline{\Tors}(Y,{p_*}G)$.\end{proof}

Le théorème de changement de base propre pour les faisceaux (trivial dans ce cas) assure  qu'on a
$f^*p_*G=p'_*f'^*G=p'_*G$. La flèche $$H^0(Y,p_*\underline{\Tors}(Y',G))\ra
H^0(X',f^*p_*\underline{\Tors}(Y',G))$$ s'identifie alors à la flèche naturelle

$$\begin{array}{llll} \Tors(Y',G)&=& H^0(Y,p_*\underline{\Tors}(Y',G))\\
&\ra&H^0(X,f^*p_*\underline{\Tors}(Y',G)) \\ &{=}& H^0(X,\underline{\Tors}(X,f^*p_*G))&\
(\text{d'après \ref{20pdirtors} et \cite[III.2.1.5.7]{CNA@Giraud}}) \\ &=&
H^0(X,\underline{\Tors}(X,p'_*G)) \\ &=& H^0(X',p'_*\underline{\Tors}(X',G)) \\ &{=}&
H^0(X',\underline{\Tors}(X',G)) \\ &=& \Tors(X',G) \end{array}$$ qui est bijective par
hypothèse.\end{proof}

\subsection{Preuve du lemme d'effacement} Soit $X$ un schéma cohérent.

\begin{lemme2}\label{cres-alg-clos}Il existe un schéma affine $X'$, un morphisme quasi-compact et
surjectif $f:X'\ra X$ tel que pour tout $x'\in X'$, le corps résiduel $k(x')$  est la clôture
algébrique du corps résiduel $k(f(x))$. \end{lemme2}

\begin{proof} Comme $X$ est quasi-compact, on peut recouvrir $X$ par un nombre fini d'ouverts
affines $X_i$. Le morphisme $f:\sqcup X_i\ra X$ est surjectif et quasi-compact ($X$ est
quasi-séparé). Comme toutes les extensions résiduelles sont des isomorphismes, on peut donc supposer
$X=\Spec(A)$ affine quitte à changer $X$ en  $\sqcup X_i$.

\begin{sslemme2}\label{20-c-co} Soient $I$ l'ensemble des polynômes unitaires (non constants) de
$A[X]$ et $$T'(A)=A[X_P, P\in I]/(P(X_P))$$  et $$f:X'=\Spec {T'(A)}\ra X=\Spec A.$$ Le morphisme
$f$ est surjectif et, pour tout $\xi\in X'$, le corps résiduel $k(\xi)$ est la clôture algébrique de
$k(f(\xi))$.\end{sslemme2}

\begin{proof} Soit $x\in X$ et $\overline{k(x)}$ une clôture algébrique de $k(x)$. La fibre
schématique $f^{-1}(x)$ est le spectre de $$B=k(x)[X_P, P\in I]/(\tilde P(X_P))$$ où $\tilde P$
désigne l'image de $P\in I$ par le morphisme de localisation des coefficients $$A[X]\rightarrow
k(x)[X].$$ Le choix de racines $x_P\in \overline{k(x)}$ pour tous les polynômes $P$ de $I$ définit
un point de $f^{-1}(x)$ assurant la surjectivité de $f$.

Soit alors $\xi\in f^{-1}(x)$ : c'est un point fermé car $f$ est entier et $k(\xi)$ est algébrique
sur $k(x)$. Soit $Q$ un polynôme unitaire de $k(x)$ de degré $d>0$. Il existe $a\in A$ d'image
$a(x)$ non nulle dans $k(x)$, un polynôme unitaire $P\in I$ et un entier $n>0$ tel que
$$s^{nd}Q(X)=P(s^nX).$$ On déduit que l'image $X_P/s^n$ dans $k(\xi)$ est une racine de $Q$. Comme
$k(\xi)$ est algébrique sur $k(x)$, ceci assure que $k(\xi)$ est une clôture algébrique de $k(x)$.
(C'est un exercice (facile) de théorie de Galois ou \acourt{V}{10, exercice 20}.) \end{proof} Comme
$f$ est quasi-compact puisqu'affine, le lemme est prouvé.  \end{proof}

Dans un second temps, rappelons la construction de la topologie constructible sur $X$ (cf.
\egalong{IV}{3}{}).

\begin{remarque2} Pour notre propos, on ne l'utilisera en fait que pour le schéma affine
$X'$.\end{remarque2}

On construit l'espace topologique $X^{\co}$ dont l'ensemble sous-jacent coïncide avec l'ensemble
sous-jacent $|X|$  de $X$ mais dont les ouverts (resp. fermés) sont les parties ind (resp.
pro)-constructibles, à savoir les réunions (resp. intersections) de parties constructibles.  Comme
$X$ est cohérent, $X$ est un espace topologique compact, totalement discontinu
(\egalong{IV}{1}{1.9.15}). De plus, la cohérence de $X$ entraîne que les parties constructibles sont
alors les réunions finies d'intersection $U\cap(X-V)$ avec $U,V$ ouverts quasi-compacts
(\egalong{III}{1}{9.1.3} et \egalong{IV}{1}{1.2.7}). Le complémentaire de $U\cap(X-V)$ étant
$(X-U)\cup V$, il est donc également ouvert dans $X^{\co}$ de sorte que $X$ admet une base d'ouverts
compacts.

L'identité de $|X|$ induit une application continue $X^{\co}\ra X$ puisqu'un fermé est
pro-constructible. Si $X=\Spec(A)$ est affine, $X^{\co}$ est naturellement homéomorphe au spectre
d'une certaine $A$-algèbre $T(A)$ pour un certain endo-foncteur $T$ de la catégorie des anneaux
(\cite{Oli67}, proposition 5). Cet homéomorphisme est compatible à la localisation de sorte que ces
structures schématiques se recollent munissant $X^{\co}$ d'une structure naturelle de $X$-schéma
relativement affine compatible avec l'application continue (identique !) $X^{\co}\ra X$. Si $x\in
|X|$, on a $\OO_{x}=k(x)$.

On définit alors $$f:X^c\ra X$$ comme le composé $$f:X^c=(X')^{\co}\ra X'\ra X.$$

Par construction, $X^c$ est compact -en particulier cohérent- (réduit), totalement discontinu et
admet une base de voisinages ouverts-compacts (qui sont donc ouverts-fermés puisque $X^c$ est
compact donc topologiquement séparé).  Ses corps résiduels sont algébriquement clos et $f$ est
quasi-compact\footnote{On appliquera ici cette construction à un ouvert (quasi-compact) d'un schéma
affine, donc à un schéma séparé de sorte que $f$ sera même affine dans ce cas.} et surjective (comme
composé de morphismes quasi-compacts) .

Rappelons que le faisceau vide sur un espace topologique est le faisceau associé au préfaisceau de
valeur constante $\varnothing$.  L'ensemble de ses sections sur tout ouvert non vide est
$\varnothing$ et est réduit à un point sur l'ouvert vide.

\begin{lemme2}\label{20-co-sec} Tout morphisme étale $f:Y\ra X^c$ est Zariski localement trivial.
En particulier, le morphisme canonique de topos $(\epsilon^{-1},\epsilon_*): X^c_\et\ra
X^c_{\textup{Zar}}$ est une équivalence d'inverse $(\epsilon^{-1},\epsilon_*)$. Tout faisceau ayant
des sections localement a des sections globales. De plus, tout torseur sur $X^c$ est trivial et
toute gerbe est neutre.\end{lemme2}

\begin{proof} Soit $y\in Y$ d'image $x\in X^c$. Comme $f$ est quasi-fini et $k(x)$ algébriquement
clos, l'inclusion $k(x)\hookrightarrow k(y)$ est une égalité. Le morphisme composé
$$\Spec(\OO_x)=\Spec(k(x))=\Spec(k(y))\ra Y$$ se prolonge au voisinage de $f$ en une section locale
de $f$ ce qui prouve le premier point.

Soit $\F$ un faisceau (Zariski ou étale, c'est la même chose) sur $X^c$. Il est donc Zariski
localement trivial. Comme $X^c$ est compact, on peut trouver un recouvrement fini par des ouverts
compacts $U_i$ sur lesquels $\F$ a une section. On montre par récurrence sur le nombre d'ouverts
qu'on peut raffiner ce recouvrement en un recouvrement fini $V_j$ par des ouverts compacts
disjoints. Comme pour tout $j$ il existe $i$ tel que $V_j\subset U_i$, le faisceau $\F$ a des
sections locales sur chaque $V_j$. Ces ouverts étant disjoints, ces sections se recollent en une
section globale. Le reste suit car tout torseur (resp. toute gerbe) sur $X^c$ a des sections
localement.\end{proof}

La preuve du lemme d'effacement est alors facile. Soit $\gamma$ une gerbe ind-finie sur $X$. Le
foncteur d'adjonction $$\gamma\ra f_*f^* \gamma$$ est fidèle car $f$ est surjectif. La gerbe
$f^*\gamma$ (\cite{CNA@Giraud}, III.2.1.5.6) est neutre et ind-finie (\ref{20-co-sec}) de sorte
qu'elle est équivalente à $\underline{\Tors}(X^c,G^c)$ pour un ind-groupe $G^c$ convenable. Par
ailleurs, comme on a déjà vu (\cite{CNA@Giraud}, V.3.1.9.1), le faisceau des sous-gerbes maximales
$\pi_0(f_*\underline{\Tors}(X^c,G^c))$ s'identifie à $R^1f_*G^c$.

\begin{lemme2} Le faisceau $R^1f_*G^c$ est trivial.\end{lemme2} \begin{proof} Soit $U\ra X$ un
morphisme étale. Comme dans la preuve de \ref{20-co-sec}, $X^c_U\ra X^c$ est étale et donc un
isomorphisme local pour la topologie de Zariski car les corps résiduels de $X^c$ sont algébriquement
clos. On en déduit que les topos étale et Zariski de $X^c_U$ sont équivalents. Tout $G^c$-torseur
étale sur $X_U^c$ provient donc d'un torseur Zariski.  Comme $X^c_U\ra X^c$ est un isomorphisme
local, $X_U^c$ a une base d'ouverts-fermés et donc est séparé. Comme $X^c_U$ est quasi-compact
(puisque $f$ est quasi-compact et $U$ quasi-compact comme $X$), l'espace topologique sous-jacent de
$X^c_U$ est de plus compact. On déduit comme dans \ref{20-co-sec} que tout torseur sur $X_U^c$ est
trivial $H^1(X_U^c,G^c)=\{*\}$. En passant à la limite (on n'utilise pas ici la cohérence de $f$),
on trouve que les fibres de $R^1f_*G_c$ sont triviales.\end{proof}

D'après le lemme, $\pi_0(f_*\underline{\Tors}(X^c,G^c))$ est le faisceau ponctuel ce qui assure que
$f_*\underline{\Tors}(X^c,G^c)$ est une gerbe. Comme elle a une section, elle est neutre de groupe
$G=f_*G^c$ et s'identifie donc à $\underline{\Tors}(X,G)$. Mais $f$ est quasi-compact de sorte que
$G$ est ind-fini comme $G^c$ \sga{4}{IX}{1.6}. La preuve du lemme d'effacement est complète.

\section{Appendice 2: théorème de changement de base propre  d'Artin-Grothendieck pour les champs
ind-finis sur des schémas non noethériens}

On va prouver l'énoncé suivant

\begin{theoreme}\label{20-th-cgt-base-propre}Considérons un
diagramme cartésien $$\xymatrix{X'\ar[r]^{g'}\ar[d]_{f'}\ar@{}[rd]|\square&X\ar[d]^f\\
Y'\ar[r]^g&Y}$$ avec $f$ propre. Alors, pour tout champ ind-fini $\C$, la flèche de changement de
base $g^*f_*\C\ra f'_*g'^*\C$ est une équivalence.
\end{theoreme}

\begin{proof} Notons que ce théorème est connu dans le cas discret \sga{4}{XII}{5.1} et si $Y$ est
localement noethérien (\cite{CNA@Giraud}, VII.2.2.2). La preuve qui suit est une adaptation de la
preuve de ce dernier énoncé.

\textit{Preuve du théorème : } D'après \cite{CNA@Giraud}, VII.2.2.5, il suffit de prouver l'énoncé
suivant : soit $X$ propre sur $S$ local hensélien et $i:X_0\hookrightarrow X$ l'immersion de la
fibre fermée. Alors, pour tout champ ind-fini $\C$ sur $X$, la flèche \begin{equation*} \gamma:
H^0(X,\C)\ra H^0(X_0,i^*\C) \end{equation*} est une équivalence. Notons que $X/S$ étant propre, il
est cohérent donc $X$ est cohérent comme $S$ \sga{4}{VI}{2.5}.  Ainsi, $i$ est un morphisme cohérent
de schémas cohérents. Si $\C$ est discret, le théorème est une conséquence immédiate du théorème de
changement de base propre pour les faisceaux d'ensembles \sga{4}{XII}{5.1 (i)}.  On en déduit que
$\gamma$ est pleinement fidèle. D'après \ref{20-prop-red-champ}, il suffit de montrer que pour tout
morphisme fini $X'\ra X$ (induisant une immersion fermée $X'_0\hookrightarrow X'$) et tout groupe
fini $G$, la flèche \begin{equation}\label{20-rest-ferm} \gamma: \Tors(X',G)\ra \Tors(X'_0,G)
\end{equation} est une équivalence (et en fait est essentiellement surjective puisqu'on sait déjà
qu'elle est pleinement fidèle). On applique alors \sga{4}{XII}{5.5 (ii)} au morphisme propre $X'\ra
S$ pour conclure.\end{proof}

\section{Appendice 3: sorites sur les gerbes}

On montre que toute gerbe ind-finie sur $X$ noethérien est limite inductive filtrante de ses
sous-gerbes constructibles (comparer avec \egalong{IV}{3}{IX.2.9} et \egalong{IV}{3}{IX.2.2}).

\subsection{Image d'un morphisme de champs} Soit $\varphi:\C\ra\C'$ un morphisme (cartésien) de
champs sur $X_\et$.  On définit l'image $\Im(\varphi)=\varphi(\C)$ comme la catégorie ayant pour
objets ceux de $\C$ et telles que
$\Hom_{\varphi(\C)}(g_1,g_2)=\Hom_{\C'}(\varphi(g_1),\varphi(g_2))$, la structure de catégorie
fibrée étant déduites des structures (compatibles) de $\C$ et $\C'$. Notons que $\Im(\varphi)$ est
naturellement équivalente à la sous-catégorie pleine  $\C'_\varphi$ de $\C'$ dont les objets sont
les images des objets de $\C$. On identifiera $\varphi(\C)$ et $\C'_\varphi$.

\begin{lemme2}\label{20-lemme-im} Soit $\bar x$ un point géométrique de $X$ et $f:W\ra X$ un
morphisme de schémas. Alors, \begin{itemize} \item la flèche naturelle $\Im(\varphi)_{\bar
x}\ra\Im(\varphi_{\bar x})$ est une équivalence ; \item on a une équivalence de champs canoniques
$f^*\Im(\varphi)\iso\Im(f^*\varphi)$.  \end{itemize} \end{lemme2}

\begin{proof} Les objets de $\Im(\varphi)_{\bar x}\ra\Im(\varphi_{\bar x})$ coïncident avec ceux de
$\C_{\bar x}$ et la flèche naturelle est simplement l'identité. Construisons l'inverse de la flèche.
Soit alors $a,b\in\bar C_{\bar x}$ et $\psi\in\Hom(\phi(a),\phi(b))$ qui provient de
$\Psi_S\in\Hom_S(\phi(\alpha),\phi(\beta))$. Mais $\Psi_S$ peut-être vu comme une flèche de
$\Im(\varphi)(S)$ : on prend son germe en $\bar x$ pour définir l'inverse (qui ne dépend pas des
choix). On vérifie que ceci définit l'inverse cherché.

Passons au second point et définissons la flèche. Par adjonction, on doit définir une flèche
(cartésienne) $$\Im(\varphi)\ra f_*\Im(f^*\phi).$$

Soit $S\ra X$ étale. Les objets du membre de gauche sont les objets de $\C(S)$ tandis que ceux de
droites sont ceux de $f^*\C(f^{-1}(S))=f_*f^*\C(S)$. La flèche d'adjonction $\C\ra f_*f^*\C$ permet
alors de définir  la flèche $x\mapsto f^*(x)$ cherchée au niveau des objets. Soient alors $x,y$ des
objets de $\C(S)$ et $$g\in\Hom_S(\varphi(x),\varphi(y))=\Hom_{\Im(\varphi)(S)}(x,y)$$ C'est donc
une section sur $S$ de $\underline{\Hom}(\varphi(x),\varphi(y))$ qui fournit (par image
inverse) une section sur $f^{-1}(S)$ de $\underline{\Hom}(f^*\varphi(x),f^*\varphi(y))$
(\cite[II.3.2.8.1 (4)]{CNA@Giraud}, donc une flèche de
$$\Hom_{f^{-1}(S)}(f^*\varphi(x),f^*\varphi(y))=\Hom_{f^{-1}(S)}((\varphi(f^*x),\varphi(f^*y))=
\Hom_{f_*\Im(f^*\varphi)(S)}(x,y).$$ Le foncteur ainsi défini est visiblement cartésien (comme
$\varphi$).  Le premier point assure que les fibres de ce foncteur sont des équivalences, ce qui
achève de prouver le lemme.\end{proof}

\subsection{Groupoïdes libres}Soit $\Gamma=\xymatrix{E\ar@<1ex>[r]^s\ar@<-1ex>[r]_b&V}$ un graphe
orienté ($E$ est l'ensemble des arêtes, $V$ l'ensemble des sommets, $b,s$ les applications « but,
source »). On associe (voir \cite{Ber95}) le groupoïde libre $L(\Gamma)$ qu'on peut décrire comme
suit. Soit $E^{\pm}$ l'ensemble $E^\pm=\{e^{\pm},e\in E\}$ union disjointe de deux copies de $E$ :
ses objets sont les sommets et les morphismes entre $v,v'\in V$ sont les mots (réduits)
$$e_1^{\pm}\cdots e_n^{\pm} \text{ avec } b(e_i)=s(e_{i+1}) (i=1,\cdots, n-1),  s(e_1)=v,
b(e_n)=v'.$$

\begin{remarque2} Il est bien connu que  $L(\Gamma)$ est le groupoïde fondamental
$\Pi_1(\Gamma_\RR)$ de la réalisation géométrique $\Gamma_\RR$ de $\Gamma$.\end{remarque2}

Par construction, les foncteurs de $L(\Gamma)$ dans un groupoïde $G$ s'identifient naturellement aux
familles $$(g_v)\in \Ob(G)^V, (\gamma_e)\in\mathrm{Fl}(G)^E\text{ telles que }\gamma_e\in
\Hom_G(g_{s(e)},g_{b(e)}).$$ Si on préfère, $L$ est l'adjoint à gauche du foncteur d'oubli
$\text{Groupoïdes}\ra\text{Graphes}.$

La construction se globalise de la manière suivante. Considérons un diagramme de $X$-schémas étales

$$\Gamma_X:\xymatrix{E\ar@<1ex>[rr]^s\ar@<-1ex>[rr]_b\ar[rd]&&V\ar[ld]\\&X}$$

Par fonctorialité de la construction $L$, on définit un préchamp sur $X_\et$ par la formule
$S\mapsto L(\Gamma_X(S))$ dont le champ associé est noté $L(\Gamma_X)$ : c'est le groupoïde libre
engendré par le graphe orienté. Par construction $L(\Gamma_X)$ a des sections locales si et
seulement si $V\ra X$ est surjectif. Les sections locales sont localement isomorphes si et seulement
si pour tout point géométrique $\bar x\ra X$, le graphe  $\xymatrix{E_{\bar
x}\ar@<1ex>[r]\ar@<-1ex>[r]&V_{\bar x}}$ est connexe. Par construction, on dispose de deux sections
tautologiques $$g\in \Ob(\Gamma_X)(V),\gamma\in\Hom_{(\Gamma_X)(E)}(s^*g,b^*g)$$  définies par
l'identité de $V$ et de $E$ respectivement. Soit $\G$ un groupoïde sur $X_\et$. On a alors la
propriété d'adjonction suivante : la flèche qui à un foncteur \emph{cartésien}
$$\varphi:L(\Gamma_X)\ra\G$$ associe
$$\varphi(g)\in\G(V),\varphi(\gamma)\in\Hom_{\G(E)}(s^*\varphi(g),b^*\varphi(g))$$ est bijective.

\subsection{Constructibilité de sous-gerbes} Considérons un hyper-recouvrement\footnote{La
terminologie est abusive : manque la section diagonale $V\ra E$ pour avoir un hyper-recouvrement
(tronqué).} de $X$, à savoir un diagramme de $X$-schémas étales (de type fini)

$$H_X:\xymatrix{E\ar@<1ex>[rr]^s\ar@<-1ex>[rr]_b\ar[rd]&&V\ar[ld]\\&X}$$ où les flèches $V\ra X$ et
$(s,b):E\ra V\times_XV$ sont surjectives. Pour tout point géométrique $\bar x\ra X$, le graphe
$H_{\bar x}$ est connexe de sorte que $L(H_X)$ est une gerbe. Soit $$g\in
\G(V),\gamma\in\Hom_{\G(E)}(s^*\varphi(g),b^*\varphi(g))$$ définissant un morphisme cartésien
$\varphi :L(H_X)\ra\G$.

\begin{definition2} Une gerbe $\G$ sur $X_\et$ est dite constructible si pour toute section locale
$\sigma\in\G(S), S\ra\X$ ouvert étale, le faisceau en groupes $\underline{\Aut}(\sigma)$ sur $S_\et$
est constructible.\end{definition2}

\begin{lemme2} Avec les notations précédentes, supposons $\G$ ind-finie. Alors, l'image
$I=\Im(\varphi)$ est constructible.\end{lemme2}

\begin{proof} Comme la formation de l'image et de $L$ commutent à l'image inverse, on peut procéder
par récurrence noethérienne. Il suffit donc de prouver que $I$ est constructible sur un ouvert non
vide de $X$. La constructibilité se testant après n'importe quel changement de base surjectif
localement de présentation finie,\sga{4}{IX}{2.8}, on peut supposer que $V,E$ sont des
\emph{revêtements} étales de $X$ complètement décomposés, autrement dit que $H_X$ est un graphe
constant $\Gamma$ de sorte que $\sigma$ s'identifie à $n=\card(V)$ sections $\sigma_i\in\G(X)$ deux
à deux isomorphes. Ainsi, $L(H_X)$ est le champ en groupoïdes constant $L(\Gamma)_X$. Soit $N$ le
noyau du morphisme de groupoïdes (abstraits) $$N=\Ker(L(\Gamma)\ra \G(X)).$$ Par adjonction,
$L(H_X)=L(\Gamma)_X\ra \G$ se factorise à travers la projection $L(\Gamma)_X\ra (L(\Gamma)/N)_X$.
Mais, comme $L(\Gamma)/N$ est un groupoïde \emph{fini}, il est constructible de sorte que $I$ l'est
aussi, comme quotient de $(L(\Gamma)/N)_X$.\end{proof}

\begin{proposition2}\label{20prop-constr} Soit $\pi:X\ra Y$ un morphisme de schémas noethériens,
$\G$ une gerbe ind-finie sur $Y_\et$ et $\sigma\in H^0(X,\pi^*\G)$. Il existe  une sous-gerbe
constructible $\G_1$ de $\G$ telle que $\sigma\in H^0(X,\pi^*\G_1)\subset H^0(X,\pi^*\G)$.
\end{proposition2}

\begin{proof} La formule $\pi^*\G_{x}=\G_{\pi(x)}$ assure que localement $\sigma$ provient par image
inverse d'une section locale de $\G$. Comme $Y$ est quasi-compact, on peut trouver $V\ra Y$ étale
(surjectif de type fini) et $\tau\in \G(V)$ telles que $\pi^*\tau$ et $s$ coïncident localement sur
$X$. En considérant $E\ra V\times_YV$ étale convenable définissant un isomorphisme $\gamma:
s^*\tau\ra b^*\tau$, il suffit d'après le lemme précédent de poser
$\G_1=\Im(L(V,E)\xrightarrow{(\tau,\gamma)}\G)$.\end{proof}

\begin{remarque2} On contourne ici l'absence de sorites sur les limites inductives. L'énoncé devrait
être en deux parties : d'abord qu'une gerbe ind-finie sur un schéma noethérien est limite inductive
filtrante de ses sous gerbes constructibles, ce qui est pour l'essentiel le contenu du lemme
précédent, ensuite que sur un schéma cohérent, le foncteur sections globales commute aux limites
inductives filtrantes.\end{remarque2}

\ifx\danslelivre\undefined \bibliography{../fichiers-configuration-latex/bibliographie}
\bibliographystyle{../fichiers-configuration-latex/style-GTG}


\begin{center}
XXI. Le théorème de finitude pour les coefficients non abéliens\\
Frédéric Déglise
\end{center}

\else
\setboolean{english}{false}
\setboolean{francais}{true}
\selectlanguage{francais}
\chapter{Le théorème de finitude pour les coefficients
non abéliens}
\addtocontents{toc}{par Frédéric Déglise}
\begin{center}
Frédéric Déglise
\end{center}
\fi

\begin{flushright}
À la mémoire de mon oncle Olivier.
\end{flushright}

\section{Introduction}

Le but de l'exposé est de démontrer les théorèmes suivants,
qui généralisent le théorème d'Artin (cf \sga 4 {XIV}{th. 1.1})
dans le cas ensembliste et non abélien~:
\begin{theoreme}[Gabber] \label{21-ensembles}
Soit $f:Y \rightarrow X$ un morphisme de type fini
 entre schémas noethériens.

Pour tout faisceau constructible $F$ sur $Y_\et$, 
 le faisceau $f_*F$ est constructible.
\end{theoreme}

\begin{theoreme}[Gabber] \label{21-groupes}
Soit $f:Y \rightarrow X$ un morphisme de type fini
 entre schémas quasi-excellents.
 
Soit $\mathbb L$ un ensemble de nombres premiers inversibles sur $X$.

Pour tout faisceau constructible de groupes $F$ sur $Y_\et$ 
 de $\mathbb L$-torsion,
le faisceau $\R^1 f_*(F)$ sur $X_\et$ est constructible.
\end{theoreme}

\begin{theoreme}[Gabber] \label{21-groupes2}
Soit $X$ un schéma normal excellent,
 $Z \subset X$ une partie fermée de codimension supérieure à $2$.
Notons $j:U \rightarrow X$ l'immersion ouverte du complémentaire.

Pour tout groupe fini $G$,
 le faisceau $\R^1 j_*(G_U)$ est constructible.
\end{theoreme}
 
\begin{theoreme}[Gabber] \label{21-local}
Soit $A$ un anneau strictement local de dimension $2$.
On suppose que $A$ est normal, excellent, et on note $X'=\Spec(A)-\{\mathfrak{m}_A\}$
 son spectre épointé.

Alors, pour tout groupe fini $G$,
 l'ensemble $H^1(X';G)$ est fini.
\end{theoreme}

Le théorème \ref{21-ensembles} est prouvé dans la section \ref{21-sec-ensembles}.
Ce théorème est utilisé par les suivants dans le cas où $X$ est quasi-excellent.
Ce cas est beaucoup plus simple, comme nous le dégageons dans la démonstration.

Le théorème \ref{21-groupes} est réduit -- en trois étapes --
 au théorème \ref{21-groupes2} dans la section \ref{21-sec-groupes}.
Toutefois, le lecteur attentif notera que ce dernier théorème 
 n'est pas un simple cas particulier car il n'est pas nécessaire
 de faire d'hypothèse sur le cardinal du groupe $G$.

Le théorème \ref{21-groupes2} est réduit au théorème \ref{21-local}
 dans la section \ref{21-sec-groupes2}.
Cette réduction apparaît en \ref{21-sec-red-21-groupes2}
 et utilise deux lemmes qui ont été établis auparavant 
 (lemmes \ref{21-principeG} et \ref{21-bertini}).

Le dernier théorème est bien un cas particulier de \ref{21-groupes2}.
Toutefois, nous avons choisi de le dégager dans cette introduction
 à la fois comme une résultat important
  et comme un point clé.
Il est démontré dans la section \ref{21-sec-local}
 suivant un raisonnement par l'absurde qui utilise
 la méthode des ultrafiltres (voir \ref{21-sec-ultra-filtres} pour des rappels).

\paragraph*{Notations et conventions}
\begin{itemize}
\item Quand une topologie
 sur un schéma est sous-entendue, il s'agit de la topologie étale.
\item
Étant donné un ensemble $D$ (resp. un groupe $G$), on notera parfois $D$ (resp. $G$)
 pour le faisceau étale constant induit sur un schéma $X$ lorsque $X$ est clair
 d'après le contexte. Si l'on veut préciser $X$,
 on note ce faisceau $D_X$ (resp. $G_X$), suivant l'usage.
\item
Quand on parle de la normalisation d'un schéma $X$,
 il s'agit du morphisme canonique
$$
X'=\sqcup_{i \in I} X'_i \rightarrow X
$$
où $I$ désigne l'ensemble des composante irréductibles de $X$
 et $X'_i$ désigne le schéma normalisé de la composante irréductible de $X$ correspondant
 à $i$, munie de sa structure de sous-schéma réduit.
On dit aussi que $X'$ est le schéma normalisé associé à $X$.
\end{itemize}

\section{Image directe de faisceaux d'ensembles constructibles}
\label{21-sec-ensembles}

Dans le cas où $X$ est quasi-excellent,
la preuve est une application de résultats déjà connus
 (cf \sgacourt 4 {IX}). 
Nous commençons par exposer la démonstration dans ce
cas, puis dans le cas général. 
Toutefois, l'étape de réduction exposée dans la section qui suit
est valable dans les deux cas.

\subsection{Réduction du théorème}

On commence par réduire le théorème \ref{21-ensembles} à l'assertion suivante:
\begin{enumerate}
\item[\ppE]
Soit $D$ un ensemble fini et
$j:U \rightarrow X$ une immersion ouverte entre schémas noethériens.
Alors, le faisceau d'ensembles $j_*(D_U)$ est constructible.
\end{enumerate}
Considérons les hypothèses du théorème \ref{21-ensembles}.
D'après \sga 4 {IX} {2.14}, on peut trouver un monomorphisme
$$
F \rightarrow \prod_{i=1}^n \pi_{i*}(C_i)=Q
$$
pour des morphismes finis $\pi_i:Y_i \rightarrow Y$ et 
 des faisceaux constants finis $C_i$ sur $Y_i$.
Comme un sous-faisceau d'un faisceau constructible est constructible
(\sga 4 {IX} {2.9(ii)}),
il suffit de montrer que $f_*(Q)$ est constructible.
On est donc ramené au cas de $(f\pi_i)_*(C_i)$ pour tout $i$,
ce qui montre qu'on peut supposer $F=D_Y$ pour
un ensemble fini $D$.

Notons que dans ce cas, le théorème est local en $Y$.
En effet, si l'on se donne un recouvrement étale
$ \pi:W \rightarrow Y$ ($Y$ est noethérien),
 le morphisme d'adjonction 
$$
D_Y \rightarrow \pi_*\pi^*(D_Y)=\pi_*(D_W)
$$ 
est un monomorphisme. En lui appliquant $f_*$,
 on en déduit un monomorphisme
$$
f_*(D_Y) \rightarrow (f \pi)_*(D_W).
$$
Il suffit donc de montrer que le membre de droite est 
constructible (\sga 4 {IX} {prop. 2.9(ii)} à nouveau).
Notamment, on peut donc supposer que $Y$ est affine.

Alors, $f$ est séparé de type fini. 
On peut donc considérer une factorisation
$$
Y \xrightarrow j \bar X \xrightarrow{\bar f} X
$$
de $f$ telle que $j$ est une immersion ouverte et $\bar f$
un morphisme propre. Le résultat est connu pour $\bar f$ 
(cf \sga 4 {XIV} {th. 1.1}) donc on est réduit au cas 
de l'immersion ouverte $j$, c'est-à-dire à l'assertion
\ppE.

\remarque \label{21-rem-notes}
\begin{enumerate}
\item Si l'on suppose que $X$ est quasi-excellent,
 le schéma $\bar X$ qui apparaît dans la réduction ci-dessus
 est encore excellent puisque $\bar f$ est de type fini.
\item Dans cette réduction,
 on a vu que \ppE est locale en $U$.
\end{enumerate}

\subsection{Cas où $X$ est quasi-excellent}

Notons le lemme facile suivant:
\begin{lemme}\label{21-lm-normalisation}
Considérons un carré cartésien de schémas noethériens
$$
\xymatrix@=10pt{
U'\ar^{j'}[r]\ar_{q}[d] & X'\ar^{p}[d] \\
U\ar^j[r] & X
}
$$
tel que $j$ est une immersion ouverte
 et $p$ un morphisme fini surjectif.
Alors, pour tout ensemble fini $D$,
 si $l_*(D_V)$ est constructible, $j'_*(D_{U'})$ est constructible.
\end{lemme}
\begin{proof}
Par hypothèse, $q$ est surjectif. On en déduit que le morphisme
 d'adjonction
$$
D_U \rightarrow q_*q^*(D_U)=q_*(D_{V'})
$$
est un monomorphisme.
Appliquant $j_*$, on en déduit un monomorphisme 
$$j_*(D_U) \rightarrow p_*(j'_*(D_{U'})).
$$
Puisque $p$ est fini,
 $p_*$ préserve la constructibilité d'après \sga 4 {IX}{prop. 2.14(i)}.
Le lemme en résulte
puisqu'un sous-faisceau d'un faisceau d'ensembles constructible
 est constructible
 (\sga 4 {IX} {2.9(ii)}).
\end{proof}

Avant de passer à la preuve dans le cas général,
 notons que la démonstration du théorème \ref{21-ensembles}
 dans le cas où $X$ est quasi-excellent est plus simple.
Grâce à la remarque précédente, on se réduit à l'assertion
\ppE dans le cas où $X$ est quasi-excellent.
Puisque cette assertion est locale en $X$,
on peut supposer que $X$ est excellent.
Dès lors, la normalisation $p:X' \rightarrow X$ de $X$ est finie. 
Ainsi, le lemme précédent appliqué au carré cartésien évident
nous ramène au cas où $X$ est normal.

Seul le cas où $X$ est connexe nous intéresse. Alors, d'après
 \sga 4 {IX} {lem. 2.14.1}, $j_*(D_U)=D_X$, ce qui conclut.

\subsection{Cas général}

Considérons les hypothèse de l'assertion \ppE.

Notons que cette assertion est trivial pour $X=\emptyset$.
On peut donc raisonner par induction noethérienne sur $X$.
On suppose plus précisément que l'hypothèse d'induction suivante 
est vérifiée:
\begin{enumerate}
\item[($\mathscr H$)] Pour tout fermé strict $Z$ de $X$,
 pour tout morphisme fini $Z' \rightarrow Z$
 et pour toute immersion ouverte $l:V' \rightarrow Z'$,
 $l_*(D_{V'})$ est constructible.
\end{enumerate}
L'assertion étant locale en $X$,
 on peut supposer que $X$ est le spectre d'un anneau
 noethérien réduit $A$.
Utilisant le lemme \ref{21-lm-normalisation} -- en prenant 
pour $X'$ la somme disjointe des composantes irréductibles de $X$ --
on peut supposer aussi que $A$ est intègre.
On a déjà vu que \ppE est aussi locale en $U$
 (point ii de la remarque \ref{21-rem-notes}). 
On peut donc se ramener au cas où $U=\Spec(A_f)$ pour un élément
$f \in A$.

Soit $A'$ la clôture intégrale de $A$ dans son corps des fractions.
On pose $X'=\Spec(A')$, $Z=\Spec(A/(f))$
 et on considère le diagramme formé de carrés cartésiens:
$$
\xymatrix@R=15pt@C=25pt{
&& Z''\ar_h[d]\ar@/^14pt/^\pi[dd] \\
U'\ar^{j'}[r]\ar_q[d] & X'\ar_p[d] & Z'\ar_{i'}[l]\ar_r[d] \\
U\ar^j[r] & X & Z\ar_-i[l]
}
$$
tel que $i$ et $j$ sont les immersions évidentes,
 $p$ (resp. $h$) est la normalisation de $X$
 (resp. $Z'$).

\noindent {\it Première étape (fibres génériques de $p$)}: \\
Considérons les points génériques $z'_1,...,z'_n$ de $Z'$,
 et posons $z_r=p(z'_r)$. Bien que $p$ ne soit pas nécessairement
 fini, 
 $p^{-1}(z_r)$ est fini et l'extension résiduelle $\kappa(z'_r)/\kappa(z_r)$
 est finie (voir \cite[chap. V, th. 33.10]{Nag62}).
D'après le lemme \ref{21-lm-normalisation},
 on peut toujours remplacer $A$ par une extension finie
  $A \subset B \subset A'$.
L'hypothèse ($\mathscr H$) est en effet encore vérifiée pour $Y=\Spec(B)$.
Dès lors, on peut supposer que les conditions suivantes sont vérifiées:
\begin{enumerate}
\item[(h1)] Pour tout indice $r$, $p^{-1}(\{z_r\})=\{z'_r\}$.
\item[(h2)] Pour tout indice $r$, $\kappa(z'_r)/\kappa(z_r)$ est triviale.
\end{enumerate}
Notons $A_i$ (resp. $A'_i$) l'anneau localisé de $X$ en $z_i$ (resp. $X'$ en $z'_i$).
Alors, $A'_i$ est un anneau de valuation discrète.
Du fait que l'extension induite $A'_i/A_i$ est entière, on déduit que $A_i$
est de dimension $1$, ce qui implique que $z_i$ est un point générique du diviseur
$Z$ de $X$. Comme $q$ est surjectif, on déduit de (h1) que $z_1,...,z_n$ est l'ensemble
des points génériques de $Z$.

\noindent {\it Deuxième étape (restriction à un ouvert de $Z$)}: \\
Puisque trivialement $j^*j_*(D_U)=D_U$,
 il suffit de montrer que $i^*j_*(D_U)$ est constructible.
Notons les faits suivants:
\begin{enumerate}
\item $q$ surjectif: $D_U \rightarrow q_*q^*(D_U)=q_*(D_{U'})$
 est un monomorphisme.
\item $X'$ normal, $j'$ dominante: $j'_*(D_{U'})=D_{X'}$.\footnote{Comme
 on l'a déjà vu, c'est \sga 4 {IX} {lem. 2.14.1}}
\item $h$ surjectif: $D_{Z'} \rightarrow h_*h^*(D_{Z'})=h_*(D_{Z''})$
 est un monomorphisme.
\end{enumerate}
On déduit de {\bf i} et {\bf ii} un monomorphisme
$$
j_*(D_U) \rightarrow j_*q_*(D_{U'})=p_*(D_{X'}).
$$
Notons que $p$ est pro-fini. Le théorème de changement de base propre
 s'étend à ce cas, ce qui donne la relation: $i^*p_*=r_*{i'}^*$.
Si on applique $i^*$ au monomorphisme précédent,
 on déduit de cette relation et de {\bf iii} un monomorphisme composé:
$$
\sigma:i^*j_*(D_U) \rightarrow i^*p_*(D_{X'})=r_*(D_{Z'}) \rightarrow \pi_*(D_{Z''}).
$$
Considérons maintenant une immersion ouverte dense $l:V \rightarrow Z$
 ainsi que le carré cartésien suivant:
$$
\xymatrix@=15pt{
Z''\ar_\pi[d] & V''\ar_/-2pt/{l'}[l]\ar^{\pi_V}[d] \\
Z & V.\ar_/-2pt/{l}[l]
}
$$
On en déduit un diagramme commutatif de faisceaux d'ensembles:
$$
\xymatrix@R=15pt@C=20pt{
i^*j_*(D_U)\ar_\sigma[d]\ar^-\alpha[r] & l_*l^*i^*j_*(D_U)\ar^{l_*l^*\sigma}[d] \\
\pi_*(D_{Z''})\ar^-\beta[r] & l_*l^*\pi_*(D_{Z''}).
}
$$
Or le morphisme $\beta$ -- induit par le morphisme de coünité l'adjonction $(l^*,l_*)$ --
 est un isomorphisme.
En effet, $l$ étant une immersion ouverte $l^*\pi_*=\pi_{V*}{l'}^*$ par changement de base.
 On obtient donc l'identification: $l_*l^*\pi_*(D_{Z''})=\pi_*l'_*(D_{V''})$.
De plus, à travers cette identification, le morphisme $\beta$ est
l'image\footnote{On
 le vérifie facilement en revenant à la définition du morphisme de changement de base
 à l'aide des adjonctions $(l\*,l_*)$ et $({l'}^*,l'_*)$ et en utilisant que la composée suivante de morphismes unités/coünités est l'identité:
$$
l_* \rightarrow l_*l^*l_* \rightarrow l_*.
$$
}
par $\pi_*$ du morphisme de coünité pour l'adjonction $({l'}^*,l'_*)$:
$$
D_{Z''} \rightarrow l'_*l'(D_{Z''})=l'_*(D_{V''}).
$$
Ce dernier est un isomorphisme puisque $Z''$ est normal et $l'$ dense
 (voir \sga 4 {IX}{lem. 2.14.1}).
Puisque $\sigma$ est un monomorphisme,
 on en déduit que $\alpha$ est un monomorphisme. \\
Or, d'après ($\mathscr H$), le faisceau d'ensembles $l_*(D_V)$ est constructible.
Pour conclure, il suffit donc (d'après \sga 4 {IX} {2.9(ii)})
 de trouver un ouvert $V$ de $Z$ tel que
\begin{equation} \label{21-eq-immergées}
j_*(D_U)|_V=D_V.
\end{equation}

\noindent {\it Troisième étape (composantes immergées de $Z$)}: \\
Considérons la réunion $T$ des composantes immergées de $Z$ dans $X$.
Alors, l'ouvert $V=Z-T$ satisfait la relation \eqref{21-eq-immergées}.

Il s'agit de démontrer que pour tout point géométrique $\bar x$ de $V$,
 la fibre du morphisme canonique $D_V \rightarrow j_*(D_U)|_V$
 au point $\bar x$ est un isomorphisme.
Autrement dit, le morphisme induit
\begin{equation} \label{21-eq-immergéesP}
\pi_0(X_{(\bar x)}) \rightarrow \pi_0(X_{(\bar x)}-\{\bar x\})
\end{equation}
est un isomorphisme.

Soit $x$ le point de $X=\Spec A$ correspondant à $\bar x$.
Si $x$ n'est pas un point générique, 
puisque $V$ n'a pas de composante immergée,
$$
\mathrm{Prof}_x(A/(f)) \geq 1 \Rightarrow \mathrm{Prof}_x(A) \geq 2.
$$
D'après le théorème de Hartshorne (\sga 2 {III} {3.6}),
 le morphisme \eqref{21-eq-immergéesP} est donc un isomorphisme. \\
Supposons que $x$ est un point générique. D'après (h1), il existe
un indice $i$ tel que $x=z_i$. Or, d'après (h2), le morphisme
entier birationnel
$$
X'_{(z'_i)}=\Spec(A'_i) \rightarrow \Spec(A_i)=X_{(z_i)}
$$
est radiciel. C'est donc un homéomorphisme universel.
On en déduit que le morphisme $X'_{(\bar z'_i)} \rightarrow X_{(\bar z_i)}$
est encore un homéomorphisme, où $\bar z'_i$ est le point géométrique
correspondant à la clôture séparable de $\kappa(z_i)$ définie par $\bar z_i$.
Dès lors, \eqref{21-eq-immergéesP} est un isomorphisme 
puisque la propriété correspondante est vraie
 pour le schéma normal $X'_{(\bar z'_i)}$. Ceci conclut.

\section{Image directe dérivée de faisceaux de groupes constructibles}
\label{21-sec-groupes}

\subsection{Réduction au cas d'un faisceau constant.}

\label{21-constant}

\begin{lemme} \label{21-lm-sous-fx}
Soit $f:Y \rightarrow X$ une morphisme de type fini
 et $u:F \rightarrow F'$ un monomorphisme de faisceaux en groupes 
  constructibles sur $Y_\et$.
  
Alors, $\R^1 f_*(F')$ constructible implique
 $\R^1 f_*(F)$ constructible.
\end{lemme}
\begin{proof}
Posons $C=F'/F$ vu comme faisceau pointé, 
 constructible sur $Y_{\et}$ par hypothèse.
On considère la suite exacte de faisceaux pointés
 (cf \sga 4 {XII}{3.1})
$$
f_*(F') \rightarrow f_*(C) \rightarrow \R^1 f_*(F)
 \xrightarrow v \R^1 f_*(F').
$$
Supposant que $\R^1 f_*(F')$ est constructible,
on peut trouver une famille génératrice de sections locales 
$(e_1,...,e_n)$ de $\R^1 f_*(F')$.
Soit $\Phi_i$ le faisceau fibre de $v$ en $e_i$, 
défini par le diagramme cartésien de faisceaux sur $X_\et$~:
$$
\xymatrix@=10pt{
\Phi_i\ar[r]\ar[d] & V_i\ar^{e_i}[d] \\
\R^1 f_*(F)\ar^v[r] & \R^1 f_*(F').
}
$$
Utilisant le critère de constructibilité par les fibres et les
spécialisations (cf \sga 4 {IX}{prop. 2.13(ii)}),
on voit aisément qu'il suffit de montrer que le faisceau
$\Phi_i$ est constructible pour tout $i$.

Pour montrer cela, il suffit de se rappeler de l'interprétation de 
$\Phi_i$ en termes de $F$-objets tordus.
Soit $x$ un point de $X_\et$. 
Il existe un voisinage étale $V$ de $x$
dans $X$ et une section $e$ de $\Phi_i$ sur $V/X$.
Quitte à restreindre le voisinage $V$, $e$ provient d'un $F$-torseur
sur $V \times_Y X$. On peut alors tordre par $e$ la suite exacte 
ci-dessus restreinte à $V$, 
chaque faisceau étant muni d'une action de $f_*(F)$~:
$$
f_*(F')^e \rightarrow f_*(C)^e \rightarrow \R^1 f_*(F)^e
 \xrightarrow{v^e} \R^1 f_*(F')^e.
$$
La suite obtenue est encore exacte et 
$\Phi_i|_V$ s'identifie avec le noyau de $v^e$.
On en déduit $\Phi_i|_V \simeq f_*(C)^e/f_*(F')^e$.
Or ce faisceau est constructible d'après le théorème 
\ref{21-ensembles}.
\end{proof}

Considérons les hypothèses du théorème \ref{21-groupes}.
D'après \sga 4 {IX}{prop. 2.14}, on peut trouver un monomorphisme
de groupes
$$
F \rightarrow F'=\prod_{i=1}^r \pi_{i*}(G_i)
$$
pour des morphismes finis $\pi_i:U_i \rightarrow U$ 
 et des groupes finis $G_i$ pour $i=1,...,r$.
D'après le lemme précédent, on est réduit au cas de $F'$.
Puisque $\pi_{i*}$ est exact,
on est donc ramené au cas du morphisme $f \circ \pi_i$
et du faisceau constant sur $U_i$ de groupe $G_i$
 pour chaque indice $i$.

\subsection{Réduction au cas d'une immersion ouverte}

Le lemme clé dans cette étape de réduction est 
le suivant~:
\begin{lemme}
\label{21-devissage1}
Soit $G$ un groupe fini.
Considérons un diagramme commutatif
$$
\xymatrix@R=10pt@C=20pt{
 & {Y'}\ar^{g}[rd] & \\
Y \ar^{h}[ru]\ar|f[rr] && X
}
$$
de morphismes de type fini.
Les conditions suivantes sont vérifiées~:
\begin{enumerate}
\item[(i)] Si $h$ est surjectif, \\
$\R^1 f_*(G_Y)$ constructible implique $\R^1 g_*(G_{Y'})$ constructible.
\item[(ii)] Si $g$ est propre, \\
 $\R^1 h_*(G_Y)$ constructible implique $\R^1 f_*(G_Y)$ constructible.
\end{enumerate}
\end{lemme}
\begin{proof}
Rappelons que l'on dispose de la suite exacte de faisceaux
d'ensembles sur $X_\et$ (cf \sga 4 {XII}{prop. 3.2})~:
$$
* \rightarrow \R^1 g_*(h_* G_Y)
 \xrightarrow u \R^1 f_*(G_Y) \xrightarrow v g_*\R^1 h_*(G_Y).
$$

Considérons l'assertion (i). Puisque $h$ est surjectif,
 le morphisme d'adjonction 
$$
G_{Y'} \rightarrow h_* h^* G_{Y'}=h_* G_Y
$$ 
est un monomorphisme.
Appliquant le lemme \ref{21-lm-sous-fx}, il suffit de 
montrer que $\R^1 g_*(h_* G_Y)$ est constructible.
On peut alors conclure puisque le morphisme
$u:\R^1 g_*(h_* G_Y) \rightarrow \R^1 f_*(G_Y)$
est un monomorphisme.

Considérons maintenant l'assertion (ii).
Puisque $g$ est propre, 
 le théorème de changement de base propre \sga 4{XIV}{th. 1.1}
  conjugué avec le théorème \ref{21-ensembles} montre que
   la source de $u$ est constructible.
Par hypothèse et une nouvelle application du théorème \ref{21-ensembles},
le but de $v$ est constructible.

Il suffit alors de raisonner comme dans la démonstration de 
\ref{21-lm-sous-fx} sur les fibres du morphisme $v$ associées à une famille
finie de sections locales de $g_*\R^1 h_*(G_Y)$ qui est génératrice.
Chacune de ses fibres est localement vide ou isomorphe au faisceau tordu 
$\R^1 g_*\big((h_* G_Y)^e\big)$ pour une de ces sections locales $e$.
Comme ce faisceau est toujours constructible, on peut conclure.
\end{proof}

Considérons maintenant les hypothèses du théorème \ref{21-groupes},
dans le cas $F=G_Y$.
Puisque $Y$ est noethérien, 
 il existe un recouvrement Zariski $\pi:W \rightarrow Y$
  tel que $W$ est affine. 
D'après l'assertion (i) du lemme ci-dessus,
il suffit de montrer le théorème pour $f \circ \pi$.
On peut donc supposer que $Y$ est affine.

Le morphisme $f:Y \rightarrow X$ est alors quasi-projectif.
On peut donc considérer une factorisation 
$Y \xrightarrow j Y' \xrightarrow g X$ de $f$
tel que $g$ est projectif et $j$ est une immersion ouverte.
D'après l'assertion (ii) du lemme ci-dessus, nous sommes 
réduit au cas de l'immersion ouverte $j$.

\subsection{Réduction au théorème \ref{21-groupes2} (\emph{i.e.} la codimension $2$)}

Grâce aux deux étapes de réduction précédentes,
 nous sommes ramenés au cas d'une immersion ouverte
  $j:U \rightarrow X$ et d'un faisceau constant sur $U$
   de groupe $G$.
Dans cette étape de réduction, 
on considère la codimension du complémentaire $Z$ de $U$ dans $X$.

Pour montrer que $\R^1 j_*(G_U)$ est constructible
 on peut raisonner localement sur $X$.
On peut donc supposer que $X$ est excellent.
Considérons la normalisation $p:X' \rightarrow X$
de $X$ ainsi que le carré cartésien~:
$$
\xymatrix@=10pt{
U'\ar^{j'}[r]\ar_q[d] & X'\ar^p[d] \\
U\ar^j[r] & X.
}
$$
Puisque $q$ est surjectif, le morphisme d'adjonction
$G_U \rightarrow q_*q^*(G_U)$ est un monomorphisme.
D'après le lemme \ref{21-lm-sous-fx},
 il suffit donc de montrer que $\R^1 j_*(q_* G_{U'})$
 est constructible.
Or $p$ et $q$ étant finis,
 $\R^1 j_*(q_* G_{U'})=p_* \R^1 j'_*(G_{U'})$.
On peut donc supposer que $X$ est normal.

Dès lors, $X$ est somme disjointe de ses composantes
irréductibles, et on peut donc le supposer intègre.
Notons $K$ son corps des fonctions.
Pour une extension finie $L/K$,
 on peut considérer le schéma normalisé $X'$ de $X$ dans $L$,
 ainsi que le diagramme cartésien suivant~:
$$
\xymatrix@=10pt{
U'\ar^{j'}[r]\ar_q[d] & X'\ar^p[d] \\
U\ar^j[r] & X.
}
$$
Par adjonction, on obtient un morphisme canonique
$$
\phi^{(L)}_{U/X}:\R^1 j_*(G_U) \rightarrow p_* \R^1 j'_*(G_{U'})
$$
induit par le morphisme qui à un $G$-revêtement d'un
 schéma étale $V/U$ associe son pullback sur
  $V'=V \times_U U'$.
\begin{lemme}
\label{21-principeG}
Considérons les hypothèses et notations précédentes.
Alors, les conditions suivantes sont équivalentes~:
\begin{enumerate}
\item[(i)] $\R^1 j_*(G_U)$ est constructible.
\item[(ii)] Il existe une extension finie séparable $L/K$ 
telle que $\phi^{(L)}_{U/X}$ est trivial.
\end{enumerate}
\end{lemme}
\begin{proof}
$(ii) \Rightarrow (i)$~: Soit $L/K$ une extension finie
telle que $\phi^{(L)}_{U/X}$ est trivial.
Avec les notations qui précèdent,
on pose $C=q_*(G_{U'})/G_U$, faisceau sur $U_\et$ pointé 
 de manière évidente.
On peut alors former la suite exacte de faisceaux pointés
 (cf \sga 4 {XII}{3.2})~:
$$
j_*q_*(G_{U'}) \rightarrow j_*(C)
 \rightarrow \R^1 j_*(G_U)
  \xrightarrow{(1)} \R^1 j_*\left(q_*(G_{U'})\right).
$$
Notons que, puisque $p$ est fini,
$\R^1 j_*\left( q_*(G_{U'}) \right)=p_* \R^1 j'_*(G_{U'})$ et le morphisme
(1) s'identifie au morphisme $\phi^{(L)}_{U/X}$.
Par hypothèse, la suite exacte ci-dessus implique donc que
$\R^1 j_*(G)=j_*(C)/j_*q_*(G_{U'})$ ce qui conclut
 d'après le théorème \ref{21-ensembles}. \\
$(i) \Rightarrow (ii)$~: Considérons une famille $(e_1,...,e_n)$ 
de sections de $\R^1 j_*(G_U)$ qui soit génératrice.
Pour tout indice $i$, $e_i$ correspond à un revêtement
$P_i \rightarrow V_i \times_X U$ pour un schéma étale $V_i/X$.
On peut supposer que $P_i$ est connexe de corps des fonctions $L_i$.
Le morphisme $P_i \rightarrow X$ correspond donc à une extension
finie séparable $L_i/K$.
On considère la clôture normale $L$ d'une extension de $K$
 composée des $L_i/K$.
Par définition, $\phi^{(L)}_{U/X}(e_i)=*$ pour tout entier $i$
 et le résultat suit puisque $(e_1,...,e_n)$ est génératrice.
\end{proof}

\begin{remarque} \label{21-rem-principeG}
Considérons les notations qui précèdent le lemme.
En termes de $G$-revêtements, la trivialité du morphisme
$\phi^{(L)}_{U/X}$ s'interprète comme suit~:
\begin{enumerate}
\item[(i)] Pour tout point géométrique $\bar s$ de $Z$
  et pour tout $G$-revêtement 
   $\pi:P \rightarrow X_{(\bar s)}^h-Z_{(\bar s)}^h$,
le revêtement 
$\phi^{(L)}_{U/X,s}(P)$ est trivial.
\item[(ii)] Pour tout schéma étale $V/X$ et tout 
$G$-revêtement $\pi:P \rightarrow V-Z_V$,
 il existe un recouvrement étale $W/V'$ tel que $\pi|_{W-Z_W}$
 s'étend à $W$ (l'extension est alors unique puisque $X$ est normal).
\end{enumerate}
\end{remarque}

Dans le cas des immersions ouvertes, on peut renforcer
le lemme \ref{21-devissage1} comme suit~:
\begin{lemme}
\label{21-devissage2}
Considérons un diagramme commutatif
$$
\xymatrix@R=10pt@C=20pt{
 & V\ar^{k}[rd] & \\
U \ar^{h}[ru]\ar|j[rr] && X
}
$$
d'immersions ouvertes tel que $X$ est intègre normal et $h$ est dominante.
Alors, les conditions suivantes sont vérifiées~:
\begin{enumerate}
\item[(i)] Si $\R^1 j_*(G)$ est constructible, alors $\R^1 k_*(G)$ est constructible.
\item[(ii)] Si $\R^1 h_*(G)$ est constructible
 et pour tout morphisme fini surjectif $X' \rightarrow X$, $j'=j \times_X X'$,
  le faisceau $\R^1 j'_*(G)$ est constructible,
 alors $\R^1 j_*(G)$ est constructible.
\end{enumerate}
\end{lemme}
\begin{proof}
Remarquons que l'hypothèse sur $h$ et $V$ entraîne que 
 $h_*(G_U)=G_V$ (cf \sga 4 {IX}{lem. 2.14.1}).
 
L'assertion (i) résulte donc simplement du fait que
le morphisme canonique $\R^1 k_*(G_V) \rightarrow \R^1 j_*(G_U)$
est toujours un monomorphisme.

Considérons les hypothèses de l'assertion (ii).
D'après le lemme précédent appliqué à $h$,
 on peut trouver une extension finie $L/K$ telle que
  $\phi_{U/V}^{(L)}=*$.
Soit $X'$ le schéma normalisé de $X$ dans $L/K$,
 $k':V' \rightarrow X'$ le pullback de $k$ sur $X'$.
Appliquant le lemme précédent à $k'$, on peut trouver
une extension finie $E/L$ telle que $\phi_{V'/X'}^{(E)}=*$.

On note $X''$ le normalisé de $X$ dans $E/K$,
 $X'' \xrightarrow{p'} X' \xrightarrow p X$ les morphismes
  canoniques. On note $h'$, $j'$, $k'$ 
   (resp. $h''$, $j''$, $k''$) 
     les pullback respectifs
    de $h$, $j$, $k$ sur $X'$ (resp. $X''$).
On peut alors conclure grâce au diagramme commutatif
 suivant~:
$$
\xymatrix{
& \R^1j_*(G)\ar[r]\ar^{\phi^{(L)}_{U/X}}[d]\ar@{-->}[ld]
 & k_*\R^1 h_*(G)\ar^{k_*\phi^{(L)}_{U/V}}[d] \\
p_*\R^1k'_*(G)\ar[r]\ar_{p_*\phi^{(E)}_{V'/X'}}[d]
 & p_* \R^1j'_*(G)\ar[r]\ar^{p_*\phi^{(E)}_{U'/X'}}[d]
 & (pk')_*\R^1 h'_*(G) \\
(pp')_*\R^1k''_*(G)\ar[r] & (pp')_* \R^1j''_*(G).
 &
}
$$
La flèche pointillée existe du fait de l'exactitude de
la suite horizontale du milieu et de $\phi_{U/V}^{(L)}=*$.
On conclut puisque $\phi^{(E)}_{U/X}=
 (p_* \phi^{(E)}_{U'/X'}) \circ \phi_{U/X}^{(L)}$.
\end{proof}

\begin{lemme}
Soit $X$ un schéma régulier
 et $Z$ un sous-schéma fermé régulier de $X$.
Notons $j:U \rightarrow X$ l'immersion ouverte complémentaire.

Soit $n$ l'ordre de $G$.
Alors, si $n$ est inversible sur $X$,
 $\R^1 j_*(G_U)$ est constructible.
\end{lemme}
\begin{proof}
L'assertion est évidemment locale en $X_\et$.
On peut donc supposer que $X$ est local strictement hensélien
 et régulier, spectre d'un anneau noté $A$.
 
Si $Z$ est de codimension supérieure à $2$ dans $X$,
 il résulte du théorème de pureté de Zariski-Nagata
  (cf \sga 1 {X} {3.1, 3.3}) que $\R^1 j_*(G)=*$.

En codimension $1$, $Z$ admet un paramètre régulier $f \in A$.
On considère le schéma strictement local 
$$
X'=\Spec(A[t]/t^n-f).
$$
Le lemme d'Abhyankar absolu (cf \sga 1 {XIII}{prop. 5.2})
montre alors précisément que
 pour tout revêtement $E \xrightarrow \pi U$
  principal galoisien de groupe $G$,
le revêtement $\pi \times_X X':E' \rightarrow U'$ se prolonge
 à $X'$. 
Il est donc trivial et le lemme \ref{21-principeG} 
 accompagné de la remarque \ref{21-rem-principeG}
  permet de conclure.
\end{proof}

Revenons au cas général d'une immersion ouverte
 $j:U \rightarrow X$ de fermé complémentaire $Z$,
  $X$ étant supposé normal. Supposons que l'ordre de $G$
   est inversible sur $X$.

Soit $T$ la réunion des lieux singuliers de $X$ et $Z$.
Posons $V=X-T$ et $W=X-(Z \cup T)$ et
 considérons les immersions ouvertes correspondantes~:
$$
\xymatrix@R=-2pt@C=20pt{
& V\ar^k[rd] & \\
W\ar^h[ru]\ar_l[rd]\ar|\nu[rr] && X \\
& U\ar_j[ru] & 
}
$$
D'après le lemme précédent, $\R^1 h_*(G)$ est constructible.
D'après le lemme \ref{21-devissage2},
 on est donc ramené à prouver que pour tout morphisme fini 
  surjectif $X' \rightarrow X$, $k'=k \times_X X'$,
   le faisceau $\R^1 k'_*(G)$ est constructible.
Cette dernière assertion est bien impliqué par le théorème 
\ref{21-groupes2}.

\section{Cas de codimension $2$ sans hypothèse sur la torsion}
\label{21-sec-groupes2}

\subsection{Résolution des singularités}

\begin{lemme} \label{21-trivial-kernel}
Soit $X$ un schéma normal connexe,
$Z \subset X$ une partie fermée de codimension $2$
 et $j:U \rightarrow X$ l'immersion ouverte complémentaire.

Supposons $\mathrm{dim}(Z)>0$.
Alors, il existe une partie fermée $T \subset Z$ 
de codimension $3$ dans $X$ telle que pour tous points géométriques 
$\bar s$ et $\bar t$ de $Z-T$
 et toute spécialisation
  $\eta:X_{(\bar t)}^h \rightarrow X_{(\bar s)}^h$,
le noyau du morphisme de spécialisation
$$
\eta^*:\R^1 j_*(G)_{\bar s} \rightarrow \R^1 j_*(G)_{\bar t}
$$
est trivial.
\end{lemme}
\begin{proof}
Cette assertion ne dépend que d'un voisinage ouvert de $Z$ dans $X$.
Quitte à enlever une partie nulle part dense de $Z$,
on peut donc supposer que le lieu singulier de $X$ est inclus dans $Z$.
De même, on peut supposer que les composantes irréductibles de $Z$ sont 
disjointes, et de là, que $Z$ est irréductible.

D'après \cite{D2DS@Lipman},
 on peut résoudre la singularité de $X$ au point générique de $Z$
  par une suite d'éclatements et de normalisations.
Donc, quitte à retirer de nouveau une partie fermée nulle part dense de $Z$, 
on peut supposer qu'il existe un diagramme formé
de carrés cartésiens
$$
\xymatrix@=22pt{
Z'\ar[r]\ar[d]
 & X'\ar_p[d] & U'\ar^q[d]\ar|f[ld]\ar_/-2pt/{j'}[l] \\
Z\ar[r] & X & U\ar_/-5pt/j[l]
}
$$
tel que $X'$ est régulier, $Z$' est un diviseur dans $X'$,
toute composante irréductible de $Z'$ domine 
 $Z$ et $q$ est un isomorphisme.

On en déduit donc $\R^1 j_*(G)=\R^1 f_*(G)$.
Comme $X'$ est régulier et $U'$ dominant dans $X'$,
on obtient $j'_*(G_{U'})=G_{X'}$ (cf \sga 4 {IX}{lem. 2.14.1}).
On obtient donc un monomorphisme canonique 
$\rho:\R^1 p_*(G) \rightarrow \R^1 f_*(G)$.

Considérant les notations du lemme (où l'on a supposé $T=\emptyset$), 
on obtient donc un diagramme commutatif d'ensembles pointés
$$
\xymatrix{
\R^1 p_*(G)_{\bar s}\ar^{\rho_s}[r]\ar[d]
 & \R^1 f_*(G)_{\bar s}\ar^{\eta^*}[d] \\
\R^1 p_*(G)_{\bar t}\ar[r]
 & \R^1 f_*(G)_{\bar t}.
}
$$
D'après le théorème de changement de base propre appliqué à $p$, 
 la composée $\eta^* \circ \rho_s$ est un monomorphisme.
Il nous suffit donc de vérifier que le noyau de $\eta^*$
 est inclus dans l'image de $\rho_s$.

Quitte à tirer la situation sur $X_{(\bar s)}^h$,
on peut supposer que $X=X_{(\bar s)}^h$ pour simplifier les
notations.
On se donne donc un $G$-revêtement principal connexe
 $P \rightarrow U'$ qui est trivial
  sur $U' \times_X X_{(\bar t)}^h$.
On peut supposer que $P$ est connexe.
Si $F$ désigne le corps des fonctions de $X'$, 
 ce revêtement correspond à une extension séparable $E/F$.
Soit $\bar P$ la clôture normale de $X'$ dans $E/F$.
Bien sûr, $\bar P \times_{X'} U'=P$ puisque $P$ est 
  la clôture normale de $U'$ dans $E/F$.
Donc $\bar P/X'$ est non ramifié au-dessus de $U'$.
Du fait que $P \times_X X_{(\bar t)}^h$ est trivial,
 il suit que $\bar P$ est non ramifié en tous points 
  génériques  de $Z'$
   (puisque ceux-ci dominent $Z_{(\bar t)}^h$).
Autrement dit, $\bar P$ est non ramifié en codimension $1$.
Mais d'après le théorème de Zariski-Nagata (cf \sga 1 {X}{th. 3.1}),
 $\bar P$ est non ramifié en codimension supérieure à $2$
  car $X'$ est régulier.
Ainsi, $\bar P$ est le revêtement étale cherché.
\end{proof}

\subsection{Un argument « à la Lefschetz »}

Pour cet argument, nous utiliserons le lemme suivant,
 qui est une application des résultats liés à « la méthode de
 Lefschetz » de \sga 2 {X}{\S 2}.

\begin{lemme}\label{21-ref-deg-connexe}
Soit $X$ un schéma normal excellent connexe,
 $D$ un diviseur de Cartier effectif connexe dans $X$,
  et $Z \subset D$ une partie fermée de codimension supérieure à $2$.
Alors, $D-Z$ est connexe.
\end{lemme}

\begin{proof}
Il suffit de montrer que pour tout point $s$ de $X$, le schéma 
 $(D-Z) \times_X X_{(s)}$ est connexe.
On peut même supposer que $s$ est un point générique de $Z$.
On peut donc supposer que $X$ est local, de dimension supérieure à $3$
 et que $Z$ est son point fermé.
Considérons le complété $\hat X$ du schéma local $X$.
Puisque $X$ est excellent, $\hat X$ est encore normal.
Puisque le morphisme $\hat X \rightarrow X$ est un épimorphisme universel,
 il suffit de montrer que $(D-Z) \times_X \hat X$ est connexe.
On peut donc supposer en outre que $X$ est complet.
 
Alors, le spectre épointé $X'=X-Z$ est normal, connexe et 
de dimension supérieure à $2$. 
Il résulte du {\it critère de normalité de Serre}
(cf \cite[23.8]{CRT@Matsumura})
 que pour tout point fermé $x$ de $X'$,
$$\prof(\mathscr O_{X',x}) \geq 2.$$
Dès lors, d'après \sga 2 X {ex. 2.1} (voir aussi plus directement
 \sga 2 {IX}{prop. 1.4}), on obtient un isomorphisme canonique
$$
\Gamma(X',\mathscr O) \simeq
 \Gamma\left(\widehat{X'}^D,\mathscr O\right)
$$
où $\widehat{X'}^D$ désigne le complété formel de $X'$ le long de $D$,
 ce qui conclut.
\end{proof}

\begin{lemme} \label{21-bertini}
Soit $X$ un schéma normal excellent,
 $D$ un diviseur principal dans $X$
  et $Z \subset D$ une partie fermée non vide de codimension 
   supérieure à $2$.
On pose $U=X-Z$, $V=D-Z$.

Considérons le carré cartésien formé des immersions
évidentes
$$
\xymatrix@=10pt{
V\ar^{j'}[r]\ar_{i_U}[d] & D\ar^i[d] \\
U\ar^j[r] & X.
}
$$
Alors, le morphisme de changement de base associé
$$
i^* \R^1 j_*(G) \rightarrow \R^1 j'_*(G)
$$
est un monomorphisme.
\end{lemme}
\begin{proof}
On peut supposer que $X$ et $D$ sont locaux strictement
henséliens. On doit montrer que le morphisme 
de restriction
$$
H^1(U;G) \xrightarrow{i_U^*} H^1(V;G)
$$
est injectif.

Remarquons que le lemme précédent nous montre déjà que $V$
 est connexe.
Soit $P$ et $P'$ deux $G$-torseurs sur $U$ qui coïncident
sur $V$. On considère le faisceau 
$L=\underline{\mathrm{Isom}}_G(P,P')$
des $G$-isomorphismes de $P$ dans $P'$ sur $U_\et$.
on doit montrer qu'il admet une section sur $U$.

Comme $L$ est localement constant constructible,
 il est représentable par un $U$-schéma étale fini noté $U'$.
Posons $V'=U' \times_U V$. 
Par hypothèse, $V'/V$ admet une section.

On va montrer que pour toute composante connexe $U'_0$ de $U'$
telle que $U'_0 \times_U V/V$ admet une section, il
existe une section de $U'_0/U$ ce qui suffira pour conclure.

Quitte à remplacer $U'$ par $U'_0$, on peut supposer 
pour montrer cela que $U'$ est connexe.
Les corps de fonctions de $U'$ et $U$ 
définissent une extension finie séparable $L/K$.
Soit $X'$ le schéma normalisé de $X$ dans $L/K$. 
On pose encore $D'=X' \times_X D$ et $Z'=X' \times_X Z$.

Notons que $X'$ est normal excellent et connexe. 
Dès lors, le lemme précédent implique que
  $V'=D'-Z'$ est connexe.
Par hypothèse, le $V$-schéma étale $V'$ admet une section, 
 donc $V'=V$. 
Il en résulte que le revêtement étale $U'/U$ est de degré $1$ 
au-dessus de $V$, donc $U'=U$.
\end{proof}

\subsection{Réduction au théorème \ref{21-local}}
\label{21-sec-red-21-groupes2}

On peut supposer que $X$ est affine.
Il alors limite projective de schémas affines de type fini sur
$\Spec(\ZZ)$. Comme il suffit de démontrer le théorème \ref{21-groupes2}
pour chacun des termes de cette limite projective, on est réduit
au cas où $X$ est de dimension finie.

On raisonne par induction sur la dimension de $X$.
Le cas où $X$ est de dimension $2$ résulte du théorème \ref{21-local}.

On se place dans le cas où $\dim(X)>2$.

Si $\codim(Z)>2$,
 on peut trouver un diviseur principal connexe $D \xrightarrow i X$
  qui contient $Z$.
Soit $j':D-Z \rightarrow D$ le morphisme induit. \\
D'après le lemme \ref{21-bertini}, 
on obtient un monomorphisme $i^* \R^1 j_*(G) \rightarrow \R^1 j'_*(G)$.
Par hypothèse de récurrence, $\R^1 j'_*(G)$ est constructible.
On peut donc conclure car $\R^1 j_*(G)=i_*i^*\R^1 j_*(G)$.

Plaçons nous dans le cas critique où $\codim(Z)=2$,
 ce qui entraîne $\dim(Z)>0$.

D'après hypothèse d'induction,
le théorème est connu pour le schéma semi-localisé de $X$
aux points génériques de $Z$ qui sont de codimension $2$
dans $X$. Il existe donc, d'après le lemme \ref{21-principeG},
une extension finie $L$ du corps des fonctions $K$ de $X$ 
telle que 
pour tout point générique $\eta$ de $Z$ de codimension $2$ dans $X$,
$\phi_{U/X,\eta}^{(L)}$ est trivial.

Considérons $X'$ la normalisation de $X$ dans $L/K$ et
$p:X' \rightarrow X$ sa projection.
On pose $j'=j \times_X X'$ et $Z'=Z \times_X X'$.
D'après le lemme \ref{21-trivial-kernel},
il existe une partie fermée $T' \subset Z'$ de codimension supérieure à $3$
dans $X'$ tel que les noyaux des flèches de spécialisations de
$\R^1 j'_*(G)$ aux points de $Z'-T'$ soient triviaux.

Soit $T$ la réunion des composantes irréductibles de $Z$ de codimension 
supérieure à $3$ et du fermé $p(T')$ dans $Z$.
Considérons un point géométrique $\bar s$ de $Z-T$.
Il existe un point générique géométrique $\bar t$ de $Z-T$ et
une spécialisation $\eta:X_{(\bar t)}^h \rightarrow X_{(\bar s)}^h$.
Considérant les fibres de $\phi^{(L)}_{U/X}$, on obtient le diagramme
suivant~:
$$
\xymatrix@R=20pt@C=34pt{
\R^1 j_*(G)_{\bar s}\ar^-{\phi^{(L)}_{U/X,\bar s}}[r]\ar[d]
 & [p_* \R^1 j'_*(G)]_{\bar s}\ar^{\eta^*}[d] \\
\R^1 j_*(G)_{\bar t}\ar^-{\phi^{(L)}_{U/X,\bar t}}[r]
 & [p_* \R^1 j'_*(G)]_{\bar t}.
}
$$
D'après le choix de $L/K$,
la composée $\eta^* \circ \phi^{(L)}_{U/X,\bar s}$ est triviale.
Par ailleurs, puisque $p$ est fini, $\eta^*$ a un noyau trivial.
On en déduit que $\phi_{U/X,\bar s}^{(L)}$ est trivial.
D'après le lemme \ref{21-principeG}, $\R^1 h_*(G)$ est constructible
 pour l'immersion ouverte $h:X-Z \rightarrow X-T$.
D'après le lemme \ref{21-devissage2},
 on est donc réduit à montrer que pour tout morphisme fini surjectif 
  $X' \rightarrow X$, $\R^1 k'_*(G)$ est constructible
   pour l'immersion ouverte $k':X'-T' \rightarrow X'$.

Or on peut trouver un diviseur principal $D \xrightarrow i X'$
 qui contient $T'$.
On pose $k''=k' \times_{X'} D$.

D'après le lemme \ref{21-bertini}, 
on obtient un monomorphisme $$i^* \R^1 k'_*(G) \rightarrow \R^1 k''_*(G).$$
On peut donc à nouveau conclure d'après l'hypothèse d'induction appliquée
à $k''$ et du fait que $\R^1 k'_*(G)=i_*i^* \R^1 k'_*(G)$.

\section{Revêtements principaux d'une surface strictement locale épointée}
\label{21-sec-local}

\subsection{Mise en place}

D'après le théorème de rigidité de Gabber (cf \refext{XX}{20-prop-rig}),
 on peut supposer que $A$ est complet.
Soit $X=\Spec(A)$ et $X'=\Spec(A)-\{\mathfrak{m}_A\}$.

D'après le théorème de Cohen-Gabber (cf \refext{IV}{4-Cohen-Gabber égale}),
 il existe un sous-anneau régulier $R \subset A$
 tel que $A/R$ est finie génériquement étale.
On note $m$ le degré générique de $A/R$
Notons tout de suite le fait suivant qui résulte de l'algèbre
commutative standard~:
\begin{lemme} \label{21-facile}
Soit $R$ un anneau local régulier de dimension $2$,
 $A/R$ une algèbre finie tel que $A$ est normal.
Soit $m$ le degré générique de $A/R$.

Alors $A$ est un $R$-module libre de rang $m$.
\end{lemme}
En effet, puisque $R$ est local régulier de dimension $2$,
 par définition $\prof(R)=2$.
Par ailleurs, puisque $A$ est local normal de dimension $2$,
 il résulte du critère de Serre que $\prof(A)=2$ (cf \cite[ex.
17.3]{CRT@Matsumura}).
D'après le théorème d'Auslander-Buchsbaum (cf \cite[th.
19.1]{CRT@Matsumura}),
$$\prof(A)+\mathrm{dim\ proj}(A)=\prof(R).$$ Il en résulte que 
le $R$-module $A$
 est nécessairement projectif, donc libre.

Faisant abstraction du groupe $G$,
 on fixe un entier $n>0$ et on montre que l'ensemble
  des revêtements étales de $X'$ de degré $n$ est fini.

On raisonne par l'absurde.
Considérons une suite $(P'_i \rightarrow X')_{i \in \NN}$
 de revêtements étales de degré $n$ telle que pour
  tout $i \neq j$, $P'_i$ est non $X'$-isomorphe à $P'_j$.

Soit $K$ le corps des fractions de $A$.
Pour tout entier $i$, $P'_i/X'$ correspond
 à une extension finie séparable $L_i/K$.
On note $B_i$ la clôture intégrale de $A$ dans $L_i$,
 $P_i=\Spec(B_i)$.
Remarquons par ailleurs que d'après le lemme précédent,
 $B_i/R$ est libre de rang $nm$.

\subsubsection{Questions de discriminant}

Rappelons la définition suivante~:
\begin{definition}
Soit $B/A$ une algèbre finie libre de rang $n$.
Soit $\mathscr B=(e_i)_{1 \leq i \leq n}$ une base de $B/A$.
Le déterminant de la matrice 
 $\left(\mathrm{Tr}_{B/A}(e_ie_j)\right)_{1 \leq i,j \leq n}$
est appelé le discriminant de $B/A$ relativement à $\mathscr B$.
Sa classe dans le monoïde multiplicatif $A/(A^\times)^2$
 est indépendante de $\mathscr B$.
On la note $\disc_{B/A}$.
\end{definition}
Par abus, on considérera la classe $\disc_{B/A}$ 
 comme un élément de $A$.
Rappelons que $B/A$ est étale si et seulement si
 $\disc_{B/A}$ est inversible dans $A$.
Par la suite, nous aurons besoin de la formule suivante
(cf \cite[2.1.4]{lm@Ramero})~:
Soit $B/A$ et $C/B$ deux algèbres finies libres.
Soit $n$ le rang de $C/B$.
Alors,
$$
\disc_{C/A}=\disc_{B/A}^n.N_{B/A}(\disc_{C/B}).
$$

Revenant à la situation du numéro précédent,
 on considère un idéal $\p$ de hauteur $1$ de $R$.
Soit $A_\p$ (resp. $B_{i,\p}$) l'anneau semi-localisé de $A$
(resp. $B_i$) correspondant à la fibre au-dessus de $\p$.

Notons que $A_\p$ est normal de dimension $1$.
Il résulte du lemme \ref{21-facile} que $B_{i,\p}/A_\p$
 est libre de rang $n$.
D'après la formule rappelée précédemment,
$$
\disc_{B_{i,\p}/R_\p}
=\disc_{A_\p/R_\p}^n.N_{A_\p/R_\p}(\disc_{B_{i,\p}/A_\p}).
$$
Or $A/R$ (resp. $B_i/R$) est génériquement étale
 et $B_{i,\p}/A_\p$ est étale.
On déduit de la relation précédente que l'élément 
$(\disc_{B_{i}/R})(\disc_{A/R}^n)^{-1}$
de $K^\times$ appartient à $R^\times_\p$. 
Comme ceci est valable pour tout $\p$
et que $R$ est normal, on en déduit~:
\begin{equation} \label{21-disc}
\frac{\disc_{B_{i}/R}}{\disc_{A/R}^n} \in R^\times
\end{equation}
 
\subsection{Lemme clé}

\subsubsection{Ultraproduits} \label{21-sec-ultra-filtres}

\begin{definition}
Soit $I$ un ensemble.
Un ultrafiltre $\F$ sur $I$ est la donnée d'un ensemble
 de parties de $I$ ordonné par inclusion vérifiant les propriétés
 suivantes~:
\begin{enumerate}
\item $\forall F \in \F, \ \forall G \in P(I), \ F \subset G \Rightarrow \ G \in \F$.
\item $\forall F,G \in \F, \ F \cap G \in \F$.
\item $\forall F \in P(I)$, $F \in \F$ ou bien $I\backslash F \in \F$.
\item $X \in \F$ et $\emptyset \notin \F$.
\end{enumerate}
\end{definition}

\begin{exemple}
Soit $a$ un élément de $I$. Alors, l'ensemble des parties de $I$
 contenant $a$ est un ultrafiltre $\F$ de $I$.
Dans ce cas, on dit que $\F$ est \textit{principal}.
\end{exemple}

D'après le lemme de Zorn, il existe des ultrafiltres
 non principaux sur un ensemble infini $I$.

\begin{definition}
Soit $\F$ un ultrafiltre sur un ensemble $I$
 et $\C$ une catégorie admettant des limites inductives filtrantes
  et des produits.

Soit $(X_i)_{i \in I}$ une famille d'objets de $\C$.
Le système inductif $(\prod_{i \in F} X_i)_{F \in \F}$ est filtrant.
On définit l'ultraproduit de $(X_i)_{i \in I}$
 suivant $\F$ comme la limite inductive de ce système~:
$$
\prod_{i\in I/\F} X_i
 =\colim_{F \in \F} \left( \prod_{i \in I} X_i \right).
$$

Si $(X_i)_i$ est la famille constante de valeur un objet $X$,
 on note $X^{I/\F}$ son ultraproduit, appelé l'ultrapuissance
  de $X$ suivant $\F$.
On dispose toujours de l'application diagonale
 $X \rightarrow X^{I/\F}$.
\end{definition}

Nous utiliserons cette notion dans le cas des anneaux ou des modules.
Nous aurons besoin du lemme élémentaire suivant~:
\begin{lemme} \label{21-lm-ultra-prod}
Soit $I$ un ensemble et $\F$ un ultrafiltre sur $I$.

Considérons une famille $(A_i)_{i \in I}$ d'anneaux.
On pose $A_\infty=\prod_{i \in I/\F} A_i$.
\begin{enumerate}
\item[(i)] Si pour tout $i \in I$,
 $A_i$ est intègre (resp. un corps),
  il en est de même de $A_\infty$.
\item[(ii)]  Si pour tout $i \in I$,
 $A_i$ est local d'idéal maximal $\mathfrak{m}_i$,
  $A_\infty$ est local d'idéal maximal $\prod_{i \in I/\F} \mathfrak{m}_i$.
\end{enumerate}

Considérons une famille d'algèbres $(B_i/A_i)_{i \in I}$,
 $B_\infty/A_\infty$ son ultraproduit suivant $\F$~:
\begin{enumerate}
\item[(iii)] Si pour tout $i \in I$,
 $B_i/A_i$ est locale (resp. libre de rang $m$),
  il en est de même de $B_\infty/A_\infty$.
\end{enumerate} 

Considérons un anneau $A$, et $A_\infty=A^{I/F}$ son ultrapuissance.
\begin{enumerate}
\item[(iv)] Si $M$ est un $A$-module de présentation finie, \\
 $M \otimes_A A^{I/\F}=M^{I/\F}$.
\item[(v)] Si $A$ est cohérent, l'application diagonale
 $A \rightarrow A^{I/\F}$ est plate.
\end{enumerate}
\end{lemme}

\subsubsection{Énoncé et conséquence du lemme clé}

Revenons à la situation qui nous occupe.
Soit $\F$ un ultrafiltre non principal sur $\NN$.
On note $B_\infty$ (resp. $A_\infty$, $R_\infty$)
l'anneau local obtenu par ultraproduit suivant $\F$
de $(B_i)_{i \in \NN}$ (resp. $A$, $R$).
Soit $\widehat{B_\infty}$, $\widehat{A_\infty}$,
 $\widehat{R_\infty}$ leurs complétions respectives.
On en déduit les tours d'anneaux locaux suivantes~:
$$
\xymatrix@R=12pt@C=30pt{
B_i\ar[r] & B_\infty\ar[r]
  & \widehat{B_\infty}
 \\
A\ar[u]\ar[r] & A_\infty\ar^\varphi[u]\ar[r]
  & \widehat{A_\infty}\ar^{\hat \varphi}[u] \\
R\ar[u]\ar[r] & R_\infty\ar[u]\ar[r]
   & \widehat{R_\infty}.\ar[u]
}
$$
Considérant le spectre des anneaux locaux sur les deux premières lignes, 
on obtient le diagramme commutatif suivant~:
$$
\xymatrix@R=12pt@C=30pt{
P_i\ar[d] & P_\infty\ar[l]\ar^p[d]
 & \hat P_\infty\ar[l]\ar^{\hat p}[d] \\
X & X_\infty\ar[l] & \hat X_\infty\ar[l]
}
$$
Considérant encore le spectre épointé des anneaux locaux $A$, $A_\infty$
 et $\widehat{A_\infty}$, on obtient par changement de base le diagramme
 suivant~:
$$
\xymatrix@R=12pt@C=30pt{
P_i'\ar[d] & P_\infty'\ar[l]\ar^{p'}[d]
 & \hat P_\infty'\ar[l]\ar^{\hat p'}[d] \\
X' & X_\infty'\ar[l] & \hat X_\infty'\ar[l]
}
$$\begin{lemme}
\begin{enumerate}
\item[(1)] Le schéma $\hat X_\infty$ est noethérien.
\item[(2)] Le morphisme $\hat X_\infty \rightarrow X$ est plat
 à fibre géométriques régulières.
\item[(3)] Le morphisme $\hat X_\infty \rightarrow X_\infty$
 est une immersion fermée surjective sur les points fermés
  de $X'_\infty$.
\item[(4)] Le morphisme fini
  $\hat p:\hat P_\infty \rightarrow \hat X_\infty$
   est étale au-dessus de $\hat X'_\infty$.
\item[(5)] Le morphisme fini
  $p:P_\infty \rightarrow X_\infty$
   est étale au-dessus de $X'_\infty$.
\end{enumerate}
\end{lemme}
Montrons pourquoi ce lemme permet de terminer la démonstration.

D'après la variante C2 du théorème de changement de base lisse
 (cf \refext{XX}{20rigfaible}, cas (ii))
 appliquée au morphisme $X \rightarrow \hat X_\infty$
  et à l'ouvert $\hat X_\infty'$,
  il existe un revêtement étale $Q \rightarrow X'$ tel que
   $\hat P'_\infty \simeq Q \times_{X'} \hat X'_\infty$.

D'après le théorème de rigidité de Gabber (cf \cite{})
appliqué à $A_\infty$, le morphisme
$$
H^1(X'_\infty,G) \rightarrow H^1(\hat X'_\infty,G)
$$
est un isomorphisme. 
Donc $P'_\infty  \simeq Q \times_{X'} X'_\infty$.

Rappelons que le morphisme $p'$
est la limite projective cofiltrante
suivant les éléments $F$ de l'ultrafiltre $\F$
des morphismes canoniques
$$
p'_F:\bigsqcup_{i \in F} P'_i \rightarrow \bigsqcup_{i \in F} X'.
$$

Mais alors, les données étant de présentation finie sur $X'$,
 il existe un élément $F$ de $\F$ tel que
  $\bigsqcup_{i \in F} P'_i \simeq \bigsqcup_{i \in F} Q$.
Ceci constitue une contradiction puisque, $\F$ étant non principal,
 $F$ contient au moins deux éléments.

\subsubsection{Preuve du lemme clé}~

\medskip

\noindent {\it Assertion (1)}~:
Notons que l'idéal de l'anneau local $A_\infty$ est de type
fini. Dès lors, le complété $\widehat{A_\infty}$ est noethérien
d'après \isoegazero I {7.2.7, 7.2.8}.
Il en est de même pour $\widehat{R_\infty}$ et $\widehat{B_\infty}$.

\medskip

\noindent {\it Assertion (2)}~:
D'après le critère de platitude par fibres (cf \isoega 4 {11.3.10}),
 il suffit de montrer que pour tout entier $l>0$,
  $A/\mathfrak{m}_A^l \rightarrow \widehat{A_\infty}/\mathfrak{m}_{\widehat{A_\infty}}^l$.
Or le membre de gauche n'est autre que $(A/\mathfrak{m}_A^l)^{I/\F}$
 qui est plat sur $A/\mathfrak{m}_A^l$ car ce dernier est cohérent
d'après la propriété (v) du lemme \ref{21-lm-ultra-prod}.
Notons que l'extension résiduelle $\kappa_A^{I/\F}/\kappa_A$ 
de $A_\infty/A$ est séparable. En effet,
 pour toute extension finie $L/\kappa_A$,
 $L \otimes_{\kappa_A} \kappa_A^{I/\F}=L^{I/\F}$ est un corps
  d'après le point (i) du lemme \ref{21-lm-ultra-prod}.
On peut alors conclure par application du 
{\it théorème de localisation de la lissité formelle}
d'André (cf \cite{Localisation@Andre}).

\medskip

\noindent {\it Assertion (3)}~:
Pour cette assertion,
il suffit d'appliquer le lemme suivant
à l'idéal maximal de $A_\infty$.
\begin{lemme}
Soit $A$ un anneau et $I$ un idéal de type fini tel que
 $I \subset \mathrm{rad}(A)$. 
Soit $\hat A$ le complété $I$-adique de $A$. \\
Alors le morphisme induit
 $\Spec(\hat A)-V(I.\hat A) \rightarrow \Spec(A)-V(I)$
est surjectif sur les points fermés.
\end{lemme}
Pour ce lemme, on doit montrer que pour tout fermé 
$Z \subset \Spec(A)$,
$$
Z \subset V(I) \Leftrightarrow f^{-1}(Z) \subset V(I.\hat A).
$$
Comme $I$ est de type fini, ceci équivaut à montrer que pour 
tout idéal $J \subset A$,
$$
\exists n>0 \mid I^n \subset J
  \Leftrightarrow
\exists n>0 \mid I^n.\hat A \subset J.\hat A.
$$
Or cela résulte facilement du lemme de Nakayama
puisque pour tout $m>0$, $I^m \subset \mathrm{rad}(A)$.

\medskip

\noindent {\it Assertion (4)}~: 
Remarquons d'abord que 
$\widehat{A_\infty}$ est normal
car $A$ est normal et les fibres géométriques de
 $A \rightarrow  \widehat{A_\infty}$ sont régulières
(cf \isoega 4 {6.14.1}).
Or, $\widehat{B_\infty}/\widehat{R_\infty}$ est libre (de rang $nm$),
donc sans torsion. Comme $\widehat{A_\infty}/\widehat{R_\infty}$
est entière, on en déduit que $\widehat{B_\infty}/\widehat{A_\infty}$
est sans torsion. 
Ainsi, $\widehat{B_\infty}/\widehat{A_\infty}$ est plat.

Or, des relations \eqref{21-disc} pour tout $i \in \NN$,
on déduit
$$
\frac{\disc_{\widehat{B_\infty}/\widehat{R_\infty}}}
 {\disc_{\widehat{A_\infty}/\widehat{R_\infty}}^n} \in \widehat R_\infty^\times
$$
Si $\p$ est un idéal de hauteur $1$ de $\widehat{R_\infty}$,
l'extension $\widehat{B_\infty})_\p/\widehat{A_\infty})_\p$ est sans torsion,
donc libre. 
La relation précédente montre que
$\disc_{\widehat{B_\infty})_\p/\widehat{A_\infty})_\p}$ est inversible
ce qui prouve que $\widehat{B_\infty})_\p/\widehat{A_\infty})_\p$ est 
étale, ce qui démontre (4).

\medskip

\noindent {\it Assertion (5)}~:
D'après le critère de platitude par fibre (cf \isoega 4 {11.3.10}),
pour tout point $x \in \hat X'_\infty$, $P_\infty)_x/X_\infty)_x$ est plat.
Comme dans le point (4), la relation \eqref{21-disc} permet de montrer que
$P_\infty)_x/X_\infty)_x$ est étale. La propriété (5) résulte donc du point 
(3).

\ifx\danslelivre\undefined
\bibliography{../fichiers-configuration-latex/bibliographie}
\bibliographystyle{../fichiers-configuration-latex/smfalpha}

%% file: lemme-proposition-theoreme
\newtheorem{theoreme}[subsection]{Th\'eor\`eme}

\newtheorem{corollaire}[subsection]{Corollaire}

\theoremstyle{definition}

\theoremstyle{definition}

\newtheorem{sous-lemme}[subsection]{Sous-lemme}
\newtheorem{sous-lemme2}[subsubsection]{Sous-lemme}
\newtheorem{sous-lemme3}[paragraph]{Sous-lemme}


%% file: numerotation
\renewcommand{\theenumi}{\textbf{\roman{enumi}}} 

\renewcommand\theequation{\thesection.\textbf{\alph{equation}}} 

\makeatletter
\@addtoreset{equation}{subsection}
\makeatother

\renewcommand\thesubsubsection{\textbf{\thesubsection.\arabic{subsubsection}}}
\renewcommand\thesubsection{\textbf{\thesection.\arabic{subsection}}}
\renewcommand\theparagraph{\textbf{\thesubsubsection.\arabic{paragraph}}}

\ifx\danslelivre\undefined
\else
\renewcommand\thechapter {\Roman{chapter}}
\fi

%% file: macros
\newcommand{\francaisouanglais}[2]{\ifthenelse{\boolean{francais}}{#1}{}\ifthenelse{\boolean{english}}{#2}{}}

\ifx\cacherXXX\undefined  \else  \fi 

\newcommand\refext[2]{exp.~\textbf{#1}, \ref{#2}}	

\newcommand{\version}{\begin{center}\small{version du \VCDate\ \`a \VCHM\,\textsc{tu}}\, {\tiny (\texttt{\VCRevision})\\ \ifx\VCModified=0 \VCModifiedText \fi}\end{center}}

\newcommand{\SP}{\mathrm{Spec}}		

\renewcommand{\L}{\mathrm{L}}		

\newcommand{\ZZ}{\mathbf{Z}}

\newcommand{\FF}{\mathbf{F}}

\renewcommand{\AA}{\mathbf{A}}

\def\commutatif{\ar@{}[rd]|{\circlearrowleft}}
\def\cartesien{\ar@{}[rd]|{\square}}

\newcommand{\ra}{\rightarrow}

\renewcommand{\lim}{{\mathrm{lim}}} 

\def\sga#1#2#3{[{\bf $\mathbf{SGA\,{#1}}$}~{\sc #2}~#3]} 	
\def\sgacourt#1#2{[{\bf $\mathbf{SGA\,{#1}}$}~{\sc #2}]}	
\def\sgaquatreetdemi#1#2{[{\bf $\mathbf{SGA\,{4\frac12}}$}~{[#1]}~#2]} 
\def\sgaquatreetdemicourt#1{[{\bf $\mathbf{SGA\,{4\frac12}}$}~{[#1]}]} 

\def\.{,\dots,}